%% file: orbifolds1.tex
\documentclass[12pt,psamsfonts,leqno,oneside,letterpaper]{amsbook}
\input newbookmacros.tex

\begin{document}

\author{Peter B. Shalen}
\address{Department of Mathematics, Statistics, and Computer Science (M/C 249)\\  University of Illinois at Chicago\\
 851 S. Morgan St.\\
  Chicago, IL 60607-7045} \email{shalen@math.uic.edu}
\thanks{Partially supported by NSF grant DMS-1207720}

%n

\abstractcomment{
The things below the stars are most likely non-problems, but I'm not quite ready to remove them.

$$\qquad*\qquad*\qquad*\qquad*\qquad*\qquad*\qquad*\qquad*$$

It has just occurred to me (7/14/17) that with my current def. of strong \simple ity, I could probably prove that a covering of a strongly \simple\ orbifold is strongly \simple\ by using Meeks-Simon-Yau for the manifold case. Something tp think about...

I want to check that I have not assumed anywhere that when $\oldXi$ is a $2$-orbifold and $x$ is a point of a $1$-dim. component of $\fraks_\oldXi$, the group associated to $x$ is generated by a reflection alone. This may once have been an issue in  the proof of Lemma \ref{affect}, which in the earlier draft kplus.tex contained a note about an issue raised by Shawn, but the second assertion, the part affected by this, is gone now. So...

I'm having doubts again about whether I've argued correctly about
$I$-fibrations. The vertical boundary of an $i$-fibered orbifold can
have components that are annular orbifolds with underlying surface a
disk and two singular points of order $2$. Is everything I've said
consistent with that example? The issue has come up again in
connection with Proposition \ref{when vertical}.The issue is understanding the case where
the annular orbifolds in the hypothesis of that prop. have singularities.

Here's a more specific point: if the base of an $I$-fibration is a
$2$-manifold, have I assumed that the total space is a $3$-manifold? I
don't think it's true. But I may have argued that one could reduce
certain facts to the manifold case by removing singularities from the
base. When one does this, the base becomes a $2$-manifold, but it's a
mistake to assume that the total space becomes a $3$-manifold.

What I believe at the moment is that if the base is a manifold and the total space is orientable, then the total space is a manifold. This illustrates how confused I get. In general, I somehow ought to
make sure I have never overlooked the case where the total space is an
orientable orbifold but the base space is a non-orientable orbifold,
which is a major source of complication. I've certainly been careful
about it in some places, but I wonder whether I have ignored it in others.

Every time I look at a missingref, I find mistakes in the text that
are independent of the missingref. Stupid things like using a letter
that denotes one object when I mean to refer to another object. I'm
not sure how I can catch all these things.

 I have been concerned about confusion
  between the calligraphic and fraktur fonts, specifically in the proof of Lemma \ref{pre-modification}. I tried
  mathscr in place of mathcal, and it's not at all clear that it's an improvement.

}

\title{Volume and Homology for Hyperbolic $3$-Orbifolds}
%\\ {To my parents, Marcia Heatter Shalen and Robert Ellis Shalen}

\begin{abstract}
 Let $\oldOmega$ be a
closed, orientable, hyperbolic 3-orbifold 
such that $\pi_1(\oldOmega)$
contains no hyperbolic triangle group. We show that strict upper bounds of 
$0.07625$, $0.1525$ and $0.22875$ for $\vol\oldOmega$ imply respective upper bounds of
$23$, $43$ and $79$ for $\dim H_1(\oldOmega;\FF_2 )$. Stronger results hold if 
we assume that the singular set $\fraks_{\oldOmega}$ is a link; specifically, under this assumption,  strict upper bounds of $0.305$, $0.4575$, $0.61$, $0.7625$ and $0.915$ for $\vol\oldOmega$ imply respective upper bounds of $7$, $13$, $14$, $28$  and $29$ for $\dim H_1(\oldOmega;\FF_2 )$. 
Irreducibility assumptions on the underlying manifold $|\oldOmega|$ of $\oldOmega$, and of the underlying manifolds of certain coverings of $\oldOmega$, also give stronger results. Specifically, if one assumes that $|\oldOmega|$ is irreducible and that $|\toldOmega|$ is irreducible for every two-sheeted cover or $((\ZZ/2\ZZ)\times(\ZZ/2\ZZ))$-cover $\toldOmega$ of $\oldOmega$, but does not assume that $\fraks_{\oldOmega}$ is a link, then a strict upper bound of 
$0.22875$ for $\vol \oldOmega$ implies an upper bound of $18$ for $\dim H_1(\oldOmega;\FF_2 )$. If one assumes that $|\oldOmega|$ is irreducible, that $|\toldOmega|$ is irreducible for every two-sheeted cover $\toldOmega$ of $\oldOmega$, and that $\fraks_{\oldOmega}$ is a link, then an upper bound of $0.61$, $1.22$ or $1.72$ for $\vol\oldOmega$ implies a respective bound of $7$, $11$ or $15$ for 
$\dim H_1(\oldOmega;\FF_2 )$. 

These upper bounds on $\dim H_1(\oldOmega;\FF_2 )$ for an orbifold $\oldOmega$ whose volume is subject to a suitable upper bound are deduced from upper bounds on $\dim H_1(|\oldOmega|;\FF_2 )$ for an orbifold $\oldOmega$ whose volume is subject to a suitable upper bound. Specifically, strict upper bounds of $0.305$, $0.61$ and $0.915$ for $\vol \oldOmega$ imply respective bounds of $4$, $ 9$ and $18$ for $\dim H_1(|\oldOmega|;\FF_2 )$. If $\fraks_{\oldOmega}$ is a link, a
strict upper bound of $0.305$ for $\vol \oldOmega$ implies an upper bound of $3$ for $\dim H_1(|\oldOmega|;\FF_2 )$, while for an integer $m$ with $2\le m\le6$, a strict upper bound of $0.305m$ for $\vol \oldOmega$ implies an upper bounds of $8(m-1)$ for $\dim H_1(|\oldOmega|;\FF_2 )$. 
If $|\oldOmega|$ is irreducible, and $\fraks_{\oldOmega}$ is not assumed to be a link,  a
strict upper bound of $0.915$ for $\vol \oldOmega$ implies an upper bound of $3$ for $\dim H_1(|\oldOmega|;\FF_2 )$.
If $|\oldOmega|$ is irreducible, and $\fraks_{\oldOmega}$ is  assumed to be a link,  then 
strict upper bounds of $1.22$, $1.83$ and $3.44$ for $\vol \oldOmega$ imply respective upper bounds of $3$, $5$ and $7$ for $\dim H_1(|\oldOmega|;\FF_2 )$.
\end{abstract}

\maketitle

\section{Introduction}

The theme of this monograph, which supersedes \cite{prelim}, is that certain explicit upper bounds on
the volume (denoted $\vol \oldOmega$) of a closed hyperbolic $3$-orbifold $\oldOmega$ impose
potentially useful upper bounds on the dimension of
$H_1(\oldOmega;\FF_2 )$. Since preliminary announcements of results of
this kind were made in \cite{arithmetic} and elsewhere, we would like
to emphasize that some of the results in this paper, unlike the preliminary
versions announced in \cite{arithmetic} and the results proved in \cite{prelim}, apply to orbifolds whose
singular sets are not assumed to be links. This makes them potentially
applicable to the problem of enumerating arithmetic lattices in
$\piggletwo(\CC)$ that was discussed in \cite{arithmetic} and will be
reviewed later in this Introduction. (On the other hand, for the case
where the singular set is a link, the results are numerically weaker
than those announced tentatively in \cite{arithmetic}.)

We emphasize that $H_1(\oldOmega;\FF_2 )$ denotes the first {\it
  orbifold} homology group of $\oldOmega$ with coefficients in
$\FF_2 $ (the field with two elements);
it is obtained from the orbifold fundamental group
$\pi_1(\oldOmega)$, denoted $\pi_1^{\rm orb}(\oldOmega)$ by some
authors, by abelianizing and tensoring with $\FF_2 $. We systematically
distinguish between an orbifold $\oldOmega$ and its underlying
topological space $|\oldOmega|$.
The singular set of  $\oldOmega$, which will be denoted
$\fraks_\oldOmega\subset|\oldOmega|$, contributes certain elements of
finite order, which are typically non-trivial, to $\pi_1(\oldOmega)$
(cf. \ref{orbifolds introduced}), and $\pi_1(|\oldOmega|)$ is
canonically identified with the quotient of $\pi_1(\oldOmega)$ by the
normal closure of the set of such elements of finite order. Thus
$H_1(|\oldOmega|;\FF_2 )$ may be regarded as a quotient of $H_1(\oldOmega;\FF_2 )$.

In the case where
$\oldOmega$ is an orientable $3$-orbifold, $|\oldOmega|$ is a
$3$-manifold, and each component of $\fraks_\oldOmega$ is either a simple
closed curve or a trivalent graph in $|\oldOmega|$. To say
that $\fraks_\oldOmega$ is a link means that all its components are
simple closed curves. 

It is a standard consequence of the Margulis Lemma \cite[Chapter
D]{bp} that an upper bound on the volume of 
%a closed hyperbolic $3$-orbifold
 $\oldOmega$ imposes some upper bound on the rank of
 $\pi_1(\oldOmega)$, and hence on the dimension of
 $H_1(\oldOmega;\FF_p)$ for any prime $p$. In \cite{rankfour}, \cite
 {last}, \cite{lastplusone}, and  \cite{fourfree}, for the case of a
 closed, orientable hyperbolic $3$-manifold $M$, relatively small upper
 bounds for $\vol M$ are shown to imply explicit bounds for $\dim
 H_1(M;\FF_2 )$ that are close to being sharp. Like the bounds given in
 \cite{rankfour}, \cite {last}, \cite{lastplusone}, and
 \cite{fourfree} for the manifold case, the bounds given in this
 monograph, for the case of hyperbolic $3$-orbifolds that are not assumed to be manifolds and are subject to relatively small upper bounds on volume, will improve the naive bounds by  orders of magnitude when they apply.

The approach used here to bounding the dimension of the $\FF_2 $-vector
space $H_1(\oldOmega;\FF_2 )$, where $\oldOmega$ is a closed,
orientable hyperbolic $3$-orbifold whose volume is subject to a given
bound, is to begin by finding bounds for  $\dim 
H_1(|\oldOmega|;\FF_2 )$ and $\dim H_1(|\toldOmega|;\FF_2 )$,
where $\toldOmega$ is an arbitrary two-sheeted covering orbifold or
$((\ZZ/2\ZZ) \times(\ZZ/2\ZZ ))$-covering orbifold of $\oldOmega$. 
%In the case where the singular set of $\oldOmega$ is a link, 
These bounds can be parlayed into bounds for $\dim H_1(|\oldOmega|;\FF_2 )$
by means of the following result, Proposition A, which is  proved in the body of this
monograph as Proposition  \ref{my little sony}. The proof is an
application of the Smith inequalities, but is rather involved and uses
the interesting combinatorial result stated as Proposition \ref{vector
  thing}.

\begin{propA}
Let $\oldOmega$ be a 
closed, orientable 3-orbifold.
Then either
\begin{enumerate}[(i)] 
\item $\oldOmega$ is covered with degree at most $2$ by some
orbifold $\toldOmega$ such that
$$\dim H_1(\oldOmega;\FF_2 )\le1+ \dim
H_1(|{\toldOmega}|;\FF_2 )+\dim H_1(|{\oldOmega}|;\FF_2 ),$$
or 
\item there exists a $((\ZZ/2\ZZ) \times(\ZZ/2\ZZ) )$-covering  $\toldOmega$ of $\oldOmega$ such that
%by some
%orbifold $\toldOmega$ such that
$$\dim H_1(\oldOmega;\FF_2 )\le 3+\dim H_1(|\oldOmega|;\FF_2 )+4\dim H_1(|\toldOmega|;\FF_2 ).$$
Furthermore, if $\fraks_\oldOmega$ is a link, then (i) holds.
\end{enumerate}
\end{propA}

(It will be noted that the  statement above
differs from that of  Proposition  \ref{my little sony} in  the
notation used.  Most of the statements in this introduction differ from
their counterparts in the body  of the monograph in that the latter
use notational conventions that are not established in the
introduction; for the same reason, notation in the sketches of proofs
given here differs from the notation used in the body of the
monograph. We also note that it has proved more convenient to state
some of the results in the body of the monograph in terms of the PL
category, but in this introduction we translate them into the smooth
category; cf. \ref{scarfuss}.

We will establish the following result, which provides bounds for
$\dim H_1(|\oldOmega|;\FF_2 )$ under certain volume bounds, and is proved in the
body of the monograph as Theorem \ref{manifold homology bound}; its
content is expressed in a table after the statement of Theorem
\ref{manifold homology bound}. In the statement of the theorem (and in
the table), the quantity $\lambda_\oldOmega$ is defined to be $2$ if
the  $\fraks_\oldOmega$ is a (possibly empty)
link in $|\oldOmega|$, and to be $1$ otherwise (see \ref{lambda thing}). The statement of Theorem B  also involves the notion of a {\it negative turnover}, which is
defined (see \ref{wuzza turnover}) to be a compact, orientable
$2$-orbifold which has negative Euler characteristic, has a $2$-sphere
as its underlying manifold, and has exactly three singular
points.

\begin{TheoremB} Let $\oldOmega$ be a closed,
orientable hyperbolic $3$-orbifold. Suppose that $\vol\oldOmega<0.915\lambda_\oldOmega$. (Thus we are assuming that either $\oldOmega$ has volume strictly less than
$
%\voct/2=
1.83$ and has a link as its singular set, or $\oldOmega$ has volume strictly less than $0.915$.) \abstractcomment{\tiny It was $3.44$. It seems at the moment that
  $\voct/2=1.83\ldots$ is the best I can do. Fix the proof and the
  application, where $1.72$ will be replaced by $\voct/4=0.9159\ldots$.} 
Suppose
 that $\oldOmega$ contains no  embedded negative turnovers. Then
$$\dim H_1(|\oldOmega|;\FF_2 )
\le\max\bigg(3,\bigg\lfloor\frac {\vol(\oldOmega)}{0.305}\bigg\rfloor+(3-\lambda_\oldOmega)
\max\bigg( (3\lambda_\oldOmega+1)
\bigg\lfloor
\frac{\vol(\oldOmega)}{0.305\lambda_\oldOmega}\bigg\rfloor
,\bigg\lfloor\frac{\vol(\oldOmega)}{0.305\lambda_\oldOmega}\bigg\rfloor+2\bigg)\bigg).
$$
%(In particular, if $\lambda_\oldOmega=2$, i.e. if $\fraks_{\oldOmega}$ is a link and $\vol\oldOmega<1.83$, we have $\dim
%H_1(|\oldOmega|;\FF_2 )=h(|\oldOmega|)\le40$; and if $\lambda_\oldOmega=1$, i.e. if $\fraks_{\oldOmega}$ is not a link and $\vol\oldOmega<0.915$, we have $\dim
%H_1(|\oldOmega|;\FF_2 ) =h(|\oldOmega|)\le18$.)
 \end{TheoremB}

Proposition A and Theorem B yield the following results, which provide bounds for $\dim H_1(\oldOmega;\FF_2 )$ 
under certain volume bounds, and are proved in the body of the
monograph as Theorem \ref{orbimain} and \ref{other orbimain}; their
content is expressed in a table after the statement of Theorem
\ref{other orbimain}. A {\it hyperbolic triangle group} is a group of
the form $\langle x,y\,|\,x^p=y^q=(xy)^r=1\rangle$ with $1/p\,+\,1/q\,+\,1/r<1$.

\begin{TheoremC}Let $\oldOmega$ be a 
closed, orientable, hyperbolic 3-orbifold such that $V:=\vol(\oldOmega)$ is strictly less than
$0.915$. Suppose that $\fraks_{\oldOmega}$ is a link, and that $\pi_1(\oldOmega)$
contains no hyperbolic triangle group. Then 
$$\begin{aligned}
\dim H_1(\oldOmega;\FF_2 )
\le1 &+\max\bigg(3,\bigg\lfloor\frac {V}{0.1525}\bigg\rfloor+
\max\bigg( 7
\bigg\lfloor
\frac{V}{0.305}\bigg\rfloor
,\bigg\lfloor\frac{V}{0.305}\bigg\rfloor+2\bigg)\bigg)\\&+
\max\bigg(3,\bigg\lfloor\frac {V}{0.305}\bigg\rfloor+
\max\bigg( 7
\bigg\lfloor
\frac{V}{0.61}\bigg\rfloor
,\bigg\lfloor\frac{V}{0.61}\bigg\rfloor+2\bigg)\bigg).
\end{aligned}
$$
In particular, $\dim H_1(\oldOmega;\FF_2 )\le29$. 
\end{TheoremC}

\begin{TheoremD}
Let $\oldOmega$ be a 
closed, orientable, hyperbolic 3-orbifold such that $V:=\vol(\oldOmega)$ is strictly less than
$0.22875$. Suppose that  $\pi_1(\oldOmega)$
contains no hyperbolic triangle group. Then
$$\dim H_1(\oldOmega;\FF_2 )\le7+4
  \max\bigg( 9
\bigg\lfloor
\frac{V}{0.07625}\bigg\rfloor
,3\bigg\lfloor\frac{V}{0.07625}\bigg\rfloor+4\bigg) .
  $$
%.305 4V
\end{TheoremD}

While triangle subgroups are excluded in the statements of Theorems C
and D, it is likely that similar results for the complementary case in
which $\pi_1(\oldOmega)$ does contain triangle subgroups can be
obtained by separate arguments based on the results of
\cite{rafalski-turnover}. As we mention below, the exclusion of
triangle subgroups is harmless in the potential applications to
arithmetic groups, which were the author's primary motivation. 

As by-products of the proofs of the results stated above, we will also
obtain the following three results, Propositions E, F, and G, which give much stronger conclusions
under certain topological restrictions; these are stated in the body
of the monograph as  Propositions \ref{lost corollary},
\ref{orbifirst}, and \ref{orbinext}. Propositions F and G will be
deduced from Proposition E using Proposition A.

\begin{propE}
Let $\oldOmega$ be a closed,
orientable, hyperbolic $3$-orbifold which contains no embedded negative
turnovers. 
Suppose that  $M:=|\oldOmega|$ of
$\oldOmega$  is irreducible. 
Then the following assertions are true.
\begin{itemize}
\item If the singular set of $\oldOmega$ is a link and $\vol(\oldOmega)\le3.44$, then $\dim H_1(M;\FF_2 )\le7$. 
\item If 
the singular set of $\oldOmega$ is a link and $\vol(\oldOmega)\le1.22$, then $\dim H_1(M;\FF_2 )\le3$.
\item If 
the singular set of $\oldOmega$ is a link and  $\vol(\oldOmega)<1.83$, then $\dim H_1(M;\FF_2 )\le6$.
\item If $\vol(\oldOmega)<0.915$, then $\dim H_1(M;\FF_2 )\le3$. 
\end{itemize}
\end{propE}

\begin{propF} Let $\oldOmega$ be a 
closed, orientable, hyperbolic 3-orbifold whose singular set 
is a link, and such that $\pi_1(\oldOmega)$
contains no hyperbolic triangle group. Suppose that $|\oldOmega|$ is irreducible, and that $|\toldOmega|$ is irreducible for every two-sheeted (orbifold) covering $\toldOmega|$ of $\oldOmega$. If $\vol \oldOmega\le1.72$ then
$\dim H_1(\oldOmega;\FF_2 )\le
13$. Furthermore, if
$\vol \oldOmega\le1.22$ then
$\dim H_1(\oldOmega;\FF_2 )\le
11$, and if $\vol \oldOmega\le0.61$ then
$\dim H_1(\oldOmega;\FF_2 )\le
7$. 
\end{propF}

\begin{propG}Let $\Mh$ be a 
closed, orientable, hyperbolic 3-orbifold such that 
$\pi_1(\Mh)$
contains no hyperbolic triangle group. Suppose that $|\Mh|$ is irreducible, and that $|\tMh|$ is irreducible for every two-sheeted (orbifold) covering $\tMh$ of $\Mh$ and for every $(\ZZ/2\ZZ \times\ZZ/2\ZZ )$-covering $\tMh$ of $\Mh$. If
$\vol \Mh<0.22875$, then $\dim H_1(\Mh;\FF_2 )\le
18$.
\end{propG}

It should be pointed out that the quantitative results stated above
are different from the versions tentatively stated in the expository
article \cite{arithmetic}. This confirms the statement made in
\cite{arithmetic} that ``the exact bound may be slightly different
when the paper is finished.'' 

Indeed, for the case where the singular set of a given orientable
hyperbolic $3$-orbifold is
assumed to be a link, which is the only case for which quantitative
results are discussed in \cite{arithmetic}, the results established here are weaker than
those stated in \cite{arithmetic}. However, on the positive side, we
wish to repeat that, unlike the statements given in \cite{arithmetic}
and \cite{prelim}, 
Theorem D and the final assertion of Proposition E do not require the hypothesis that the singular set be a
link. The reason for the improvement is the discovery of the proof of
Proposition A in its present form; when \cite{arithmetic} was written,
we knew the proof of Proposition A only in the case where the singular
set is a link, which is much easier than the general case.

As was mentioned above, the author's personal interest in the problems addressed in this
monograph  arises in part from their connection with the theory of
arithmetic groups. It was established in \cite{borel}, and is
explained in \cite{arithmetic}, that the most difficult step in
listing the arithmetic lattices of at most a given covolume in ${\rm
  PGL}(2,\CC)$ is to bound the dimension of the first homology with
coefficients in $\FF_2 $ of
certain such lattices (namely those which, in the notation explained
in \cite{arithmetic}, have the form $\Gamma_\frakO$ for some maximal
order $\frakO$ in a quaternion algebra). This is a special case of the
problem of bounding $\dim H_1(\oldOmega;\FF_2 )$, where $\oldOmega$ is
an orientable hyperbolic $3$-orbifold. For this application, it is
pointed out in \cite{arithmetic} that one can restrict attention to the case where $\oldOmega$
is closed and $\pi_1(\Mh)$ contains no hyperbolic triangle
groups.  At least in principle,  Theorem D should be directly applicable to this problem when the bound on the covolume  is sufficiently small. This illustrates the advantage of Theorem D over the results stated in \cite{arithmetic}, because the orbifolds that arise in this particular application almost never have a link as their singular set.

We would also like to mention that while Propositions E and G do not seem particularly natural from the purely geometric point
of view, these results, or the stronger Proposition \ref{new get lost}
which underlies them, have the potential of being useful for applications to the
question about enumerating arithmetic groups that we have mentioned,
if they are combined with the $\log(2k-1)$ Theorem (\cite{accs},
\cite{rankfour}), the three-dimensional version of Selberg's
eigenvalue bound, and Proposition A. As
this application is a bit speculative at the moment, we shall not give
details here.

{\bf The philosophy of the proofs.}
The starting point for the proofs of Theorem B and Proposition E  is
the well-known fact (essentially contained in Proposition \ref{hepcat}
of this monograph) that if the  irreducible manifold $M=|\oldOmega|$
itself admits a hyperbolic structure, then  $\vol
M\le\vol\oldOmega$. In this case, the conclusion of Proposition E
(which is stronger than that of Theorem B) can be deduced from the results of 
\cite{rankfour}, \cite {last}, and  \cite{fourfree}.
If  $M$ does not admit a hyperbolic structure, then it follows from
Perelman's geometrization theorem \cite{bbmbp}, \cite{Cao-Zhu},
\cite{kleiner-lott}, \cite{Morgan-Tian} that $M$ is either a small
Seifert fibered space, in which case $\dim H_1(M;\FF_2 )$ can be shown to be
at most $3$, or $M$ contains a surface $T$ which is either an
essential sphere or an  incompressible torus. One can
choose $T$ within its isotopy class so that 
$T=|\oldTheta|$, where
$\oldTheta$ is some incompressible (two-sided) $2$-suborbifold of
$\oldOmega$. (A number of elementary notions of orbifold theory
that are mentioned in this introduction,
including the notion of orbifold isotopy, 
the notion of
an incompressible suborbifold, and
the operation of splitting an
 orbifold along a two-sided
  suborbifold, are explained in
  some detail in the body of the monograph. We must also mention that 
 the arguments being sketched here are done in the PL category in the
 body of the monograph,
but the distinction between the smooth and PL categories will be ignored in this introduction.)
The challenge then becomes to prove that if $\oldOmega$ contains an incompressible $2$-orbifold whose underlying surface is a torus, then certain upper bounds for $\vol\oldOmega$ imply certain upper bounds for $\dim H_1(M;\FF_2 )$; or, contrapositively, that certain lower bounds for $\dim H_1(M;\FF_2 )$, together with the existence of such a $2$-suborbifold, imply certain lower bounds for  $\vol\oldOmega$.

For the case of a hyperbolic $3$-manifold, it was shown in \cite{ast}
that certain extrinsic topological invariants of an incompressible
surface give lower bounds on the hyperbolic volume of the
$3$-manifold. This was exploited in \cite{hodad}, \cite{after-hodad},
\cite{last}, and \cite{lastplusone}, to relate volume to the dimension
of the mod-$2$ homology. To meet the challenge described above, it is
necessary to adapt the results of \cite{ast} to the context of
orbifolds. This was first done in \cite{atkinson}. In this monograph
we need a more systematic version of the orbifold analogue of the
results of \cite{ast}. The latter results are stated in terms of the
(relative) characteristic submanifold of the manifold obtained by
splitting a hyperbolic $3$-manifold along an incompressible
surface. We therefore need a systematic version of the theory of the
(relative) characteristic suborbifold of a $3$-orbifold. After
preliminary material in Chapter \ref{prelim chapter} and Section \ref{fibration section},
% \redrealmissingref{List them. The list will depend on whether I split off a section on fibrations from Section \ref{characteristic section}}, 
the characteristic orbifold theory is developed in Section
\ref{characteristic section}, taking the main result of
\cite{bonahon-siebenmann} as a starting point. It will be seen that a
compact, orientable $3$-orbifold $\oldPsi$ which is \simple, in an orbifold sense to be
defined in the body of the monograph, has a characteristic
suborbifold which is well defined up to isotopy. 

 The orbifold analogue of the results of \cite{ast} is then carried
 out in Section \ref{darts section}. If $\oldTheta$ is an
 incompressible $2$-orbifold in a closed,
 orientable, hyperbolic $3$-manifold  $\oldOmega$ is, the formalism developed in  Section \ref{darts section} allows us to bound $\vol\oldOmega$ below by a certain invariant of
the (possibly disconnected) orbifold  obtained by splitting $\oldOmega$
along $\oldTheta$. We denote the latter orbifold by
$\oldOmega\cut\oldTheta$, and the invariant, which is additive over components, by $\volorb$.
%this invariant, which is additive over components, will be denoted by
%$\volorb(\oldOmega')$. 
(The acronym  stands for Agol, Storm, Thurston and
Atkinson, the authors of \cite{ast} and \cite{atkinson}.

For
any component $\oldPsi$ of $\oldOmega\cut\oldTheta$, the
quantity $\volorb(\oldPsi)$ 
 is bounded below by a positive constant multiple of
$-\chi(\overline{\oldPsi-\oldSigma_\oldPsi})$, where
  $\oldSigma_\oldPsi$ denotes the characteristic suborbifold of
  $\oldPsi$, and
$\chi$ denotes  orbifold Euler
  characteristic. (All $3$-orbifolds that will be encountered in this
  introduction have non-positive orbifold Euler characteristic.)
% the (orbifold) Euler characteristic of the complement  in $\oldPsi$
% of the  characteristic suborbifold of $\oldPsi$. 
For the case of a component $\oldPsi$ whose interior admits a
hyperbolic structure, an alternative approach to estimating the
invariant is to bound it below by $\vol(\inter\oldPsi)$; these two
ways of producing lower bounds turn out both to be useful, and to
complement each other. These lower bounds for the invariant are in
turn related to mod-$2$ homology, using topological arguments and the results of 
\cite{rankfour}, \cite {last}, and  \cite{fourfree}.

A major source of difficulty is that if $\oldTheta$ is an {\it arbitrary} $2$-suborbifold  of $\oldOmega$ whose components are incompressible, and such that the components of $|\oldTheta|$ are incompressible tori in $|\oldOmega|$, then $\oldTheta$ does not necessarily yield a strictly positive lower bound for volume by the arguments that we have described, and complicated combinatorial arguments are needed to replace a given such suborbifold by one that is better adapted to the purpose.

{\bf Foundational material.} 

It will be apparent from the discussion above that the monograph
contains an immense amount of foundational
 background material. We have mentioned that Sections \ref{characteristic section} and
 \ref{darts section} are devoted to the development of the
 characteristic suborbifold theory and the properties of the invariant
 $\volorb$. Both Section \ref{fibration section} and the bulk of Chapter
 \ref{prelim chapter} consist of background for Sections
 \ref{characteristic section}. Chapter
\ref{higher chapter}, whose role in the proofs will be indicated in
the sketch below, is also foundational material. The amount of such material partly accounts for the unusual
 length of the monograph, although the complexity of the
 central arguments in Chapters \ref{irreducible chapter} and
 \ref{general chapter} is a factor as well. 

{\bf Sketch of the proof of Proposition E.} 
In a sense that we shall explain below, the proof of Proposition E is
a major component of the proof of Theorem B, and we shall therefore
begin our more detailed discussion of the methods of this monograph with a sketch of
the proof of Proposition E, which occupies Sections \ref{tori
  section}---\ref{irr-M section}. (The important final section of Chapter
\ref{irreducible chapter}, Section \ref{A and C}, contains the proof of Propositions A from the Smith inequalities, and the
deduction of Propositions F and G from Propositions A and E.)

%I had the sentence ``If $\oldOmega$ is a closed, orientable,
%hyperbolic $3$-orbifold, we use a certain $2$-suborbifold $\oldTheta$
%of $\oldOmega$, whose components are incompressible, and such that the
%components of $|\oldTheta|$ are incompressible tori in $|\oldOmega|$. (If
%$|\oldOmega|$ is a hyperbolic manifold or a small Seifert fibered space we
%take $\oldTheta=\emptyset$.)'' I think that all goes out, because I've
%basically explained it above.}
%isotopy characteristic suborbifold
In this sketch we will focus
%To illustrate how the phllosophy described above plays out in this
%proof, \redcomment{That phrase doesn't fit well with the sentence before
%  it} we briefly discuss the
on the case in which $M$ is a graph manifold (but is not a small Seifert
fibered space); the case of a graph manifold may be thought of as the most difficult case,
because cutting $M$ along a system of incompressible tori cannot yield
any components whose interiors admit hyperbolic structures, and the
results of \cite{rankfour}, \cite {last}, and  \cite{fourfree}
therefore yield no information. In this
case, it can be shown (see Lemma \ref{just right}) that there is a compact, connected
$3$-dimensional submanifold $K$ of $M$ such that the components of
$\partial K$ are
incompressible tori, such that $\overline{M-K}$ is connected, and $\dim H_1(K;\FF_2 )$ and $\dim
H_1(\overline{M-K};\FF_2 )$ are both at least $m$, where $m\ge2$ is a
certain integer that increases monotonically with a prescribed lower
bound on 
$\dim H_1(M;\FF_2 )$. Among all
submanifolds satisfying these topological conditions, we choose $K$ so as to minimize the weight of $\partial K$, where the {\it weight} of a
subset of the underlying space of an orbifold is  defined to be the cardinality of its intersection with the
singular set. By symmetry we may assume that 
$\dim H_1(K;\FF_2 )\ge\dim
H_1(\overline{M-K};\FF_2 )$. The minimality can be used to show that
$\partial K=|\oldTheta|$ for some incompressible suborbifold
$\oldTheta$ of $\oldOmega$, and we may identify $K$ with $|\oldPsi|$
for some component $\oldPsi$ of $\oldOmega\cut\oldTheta$. We wish to bound $\volorb\oldPsi$ from
below by bounding 
$-\chi(\overline{\oldPsi-\oldSigma_\oldPsi})$ from below. 
%The strategy
%is to show that 
If $-\chi(\overline{\oldPsi-\oldSigma_\oldPsi})$ is
smaller than desired, which we think of as meaning that
$\oldSigma_\oldPsi $ is ``large,'' we can find a solid torus $J\subset
K$, constructed as the union of a submanifold of
$|\oldSigma_\oldPsi|$ with certain $3$-balls in
$\overline{\oldPsi-\oldSigma_\oldPsi}$, which intersects $\partial K$
is one or two annuli, homotopically non-trivial in $K$ and having strictly positive
weight. It can then be shown that some component $K_1$ of  $\overline{K-J}$ satisfies the same topological
conditions as $K$, but that $\partial K_1$  has smaller weight than
$\partial K$, a contradiction to
minimality. The details of the argument sketched here occupy much of
the technical sections \ref{seek 'n' find} and \ref{underlying tori section}.

{\bf Sketch of the proof of Theorem B.}
The proof of Theorem B, and the easy arguments needed to deduce
Theorems C and D from Proposition A and Theorem B, will be given in
Chapter \ref{general chapter} (after some preliminary material in
Chapter \ref{higher chapter}, the role of which will be briefly
explained below). As Theorem B differs from Proposition E in that the
underlying manifold $M:=|\oldOmega|$ of $\oldOmega$ is not assumed to be irreducible,
it is natural to begin the proof of Theorem B by considering a
$2$-manifold $\cals$ in $M$ whose components are $2$-spheres, such
that the components of the manifold obtained by splitting $M$ along
$\cals$ and capping off the boundary components with $3$-balls are
irreducible. Such a system of spheres $\cals$ is called an {\it
complete system} in the body of the monograph. A complete
system $\cals$ may be chosen within its isotopy class in such a way
that it is the
underlying surface of an incompressible $2$-suborbifold; this is
expressed by saying that $\cals$ is {\it admissible}. The $2$-orbifold
whose underlying surface is the complete
admissible system $\tcals$ will
%indeed, it will be seen from the sketch of
%the proof to be given here that it
 play the same role in
the proof  of Theorem B that the $2$-orbifold denoted by $\oldTheta$ above played in
the sketch of the proof of Proposition E, and will likewise be denoted
by $\oldTheta$ in the following sketch. 
%Let us denote by $\oldOmega'$
%the $3$-orbifold obtained from $\oldOmega$ by splitting $\oldOmega$
%along the suborbifold $\oldTheta$, by $M$ and $M'$ the manifolds
%$|\oldOmega|$ and $|\oldOmega'|$, and 
We will denote by $\rho$ the canonical map
from $|\oldOmega|$ to $|\oldOmega\cut\oldTheta|$.  

The orbifold $\oldOmega\cut\oldTheta $, although in general it is neither closed nor connected, satisfies a
condition which is  related to the hypothesis of Proposition E: its components
are topologically \simple\ (in an orbifold sense that is explained in
the body of the monograph), and values of the invariant $\volorb$ for the
components of $\oldOmega\cut\oldTheta )$ (see
above) are
%subject to
%a certain upper bound. \redcomment{Is there some reason I'm not just
  %saying
 bounded above in terms of $\vol\oldOmega$. which is in turn subject to
   certain bounds according to the hypotheses of Theorem B. Using
   \simple ity, and using  much the same arguments used to prove Proposition
   E, we obtain good upper bounds for $\dim H_1(|\oldOmega\cut\oldTheta |;\FF_2 )$
in terms of the values of $\volorb$ for the components of
$\oldOmega\cut\oldTheta$, and hence in terms of $\vol\oldOmega$; indeed,
these bounds are
% fact that the upper bounds for $\volorb(\oldOmega\cut\oldTheta )$ imply these good upper
%bounds for $\dim H_1(|\oldOmega\cut\oldTheta |;\FF_2 )$ is 
encoded in Proposition \ref{new get
  lost}, which is the main ingredient in the proof of Proposition
E. It should be noted that the need to prove Proposition \ref{new get
  lost}, rather than just the essentially weaker Proposition E,
contributes to the complexity of the technical arguments in Sections \ref{seek 'n' find} and \ref{underlying tori section}.

An elementary result, Proposition \ref{dual homology}, may be used to
show that upper bounds for \linebreak $\dim H_1(|\oldOmega\cut\oldTheta |;\FF_2 )$ and for the number of
components of $\cals$ give an upper bound for \linebreak $\dim H_1(|\oldOmega|;\FF_2 )$. Hence
in order to obtain an upper bound for $\dim H_1(|\oldOmega|;\FF_2 )$, it suffices to
give an upper bound for the number of components of $\cals$ in terms
of upper bounds for 
$\vol(\oldOmega )$. 
Our
method for doing this is based on the same general principles which
were invoked in the proof of Proposition
E: $\vol(\oldOmega )$ is bounded below by
$\volorb(\oldOmega\cut\oldTheta )$, or more generally by
$\volorb(\oldOmega\cut{\frakZ} )$ where $\frakZ$ is any union of
components of $\oldTheta$; and 
% for any component $\oldPsi$ of 
%$\oldOmega\cut\oldTheta$ or
 %$\oldOmega\cut{\frakZ}$, the quantity
$\volorb(\oldOmega\cut{\frakZ})$
%\oldPsi )$ 
is in turn 
%$\volorb(\oldOmega\cut\frakZ )$) is in turn
 bounded below by
a positive constant multiple of
$-\chi(\overline{{\oldOmega\cut\frakZ}
%(\oldOmega\cut\oldTheta)
-\oldSigma_{{\oldOmega\cut\frakZ}}})$. 
(The  characteristic
suborbifold $\oldSigma:=\oldSigma_{\oldOmega\cut\oldTheta}$ of $\oldOmega\cut\oldTheta $, is defined to be the union of the
characteristic suborbifolds of the components
of $\oldOmega\cut\oldTheta $ in the sense mentioned above.) 
For the purposes of the following discussion we will write $X_\frakZ$
as shorthand for $-\chi(\overline{{\oldOmega\cut\frakZ}
%(\oldOmega\cut\oldTheta)
-\oldSigma_{{\oldOmega\cut\frakZ}}})$. 
The challenge---which occupies the bulk of Chapters
\ref{higher chapter} and
\ref{general chapter}---then becomes to bound the number of components
of $\cals$ in terms of $X:=\max X_\frakZ$, where $\frakZ$ ranges over
all unions of components of $\oldTheta$.
%the sum of the quantities
%$-\chi(\overline{{\oldOmega\cut\frakZ}
%(\oldOmega\cut\oldTheta)
%-\oldSigma_{{\oldOmega\cut\frakZ}}})$.
%, where
%$\oldPsi$ ranges over the components of 
%$\oldOmega\cut\oldTheta$ or
 %$\oldOmega\cut{\frakZ}$.
%For the purpose of the discussion below let us denote this sum by $P_\frakZ$.
%the number of non-fully booked components of
%$\tcals$.

%use of the  characteristic
%suborbifold $\oldSigma:=\oldSigma_{\oldOmega\cut\oldTheta}$ of $\oldOmega\cut\oldTheta $, which we define to be the union of the
%characteristic suborbifolds of the components
%of $\oldOmega\cut\oldTheta $ in the sense mentioned above. 
A component $\tS$ of
$\tcals:=\partial
|\oldOmega\cut\oldTheta|$, which may be thought of as a ``side'' of one of the
spheres making up the system $\cals$, will be said to be {\it fully
  booked} if every component of $\overline{\cals\setminus\oldSigma}$
is the underlying set of a suborbifold of $\partial\oldOmega\cut\oldTheta $ whose
orbifold Euler characteristic is $0$. 
%Using the relationship between
%the invariant $\volorb$ and the characteristic suborbifold that was
%mentioned above, 
The quantity $X_\oldTheta\le X$
%$-\chi(\overline{{\oldOmega\cut\oldTheta}
%(\oldOmega\cut\oldTheta)
%-\oldSigma_{{\oldOmega\cut\oldTheta}}})$
is bounded below by an
expression in which each non-fully booked component of $\tcals$
makes a strictly positive contribution. Thus one way to bound the
number of components of $\cals$ from above in terms of $X$
%$-\chi(\overline{{\oldOmega\cut\oldTheta}
%(\oldOmega\cut\oldTheta)
%-\oldSigma_{{\oldOmega\cut\oldTheta}}})$
is to bound it from above in
terms of the number of non-fully booked components of $\tcals$.
%\oldTheta' \frakZ

{\bf Sketch of the proof of Theorem B, continued: a special case.}
In order to give a hint about how this is done, we first consider a
very special case. Let us define a {\it strong belt}
for a component $\tS$ of $\tcals$ to be a component $G$ of
$\tS\cap|\oldSigma|$ which is an annulus whose  weight (in the sense
defined above) is $0$, and has the property that the weights of the two components of $\tS-G$
are equal. (We use the word ``strong'' here because ``belts'' of a more
general kind are considered in the body of the monograph.) If there is
a strong
belt for $\tS$, it is unique once $\oldSigma$ has been fixed in its
  orbifold isotopy class, and we will then say that $\tS$
is {\it strongly belted.} The special case that we will consider is the one in
which (a) every fully booked component of $\tcals$ is strongly belted;  (b) for
every strongly belted component $\tS$ of $\tcals$, the component $L$ of
$|\oldSigma|$ containing the strong belt of $\tS$ is a solid torus which
meets  $\fraks'$ either in a core curve for $L$
or in the empty set, and meets every component of $\tcals$ either in a
strong belt for that component, which is a homotopically non-trivial
annulus in $L$, or in the empty set; and (c) for each
component of $\cals$ whose sides $\tS,\tS'\subset\tcals'$ are both
strongly belted, the images under $\rho$ of the strong belts for $\tS$ and $\tS'$ are isotopic in
$S\setminus\fraks$.

Suppose that (a), (b) and (c) hold. By (c) we may assume after modifying $\oldSigma$
within its isotopy class  that for every
component $S$ of $\cals$ whose sides $\tS$ and $\tS'$ are both
strongly belted,
the belts for $\tS$ and $\tS'$ have  the same image under $\rho$. If we denote by $\calw$
the image under $\rho$
of the union of all components of $|\oldSigma|$ containing strong belts for
components of $\tcals$, it then follows from (b) that $\calw$ is a
submanifold of $M$ whose boundary components are tori. Using the
\simple ity of $\oldOmega$ one can then deduce a topological
description of $\calw$: each of its components
is a
solid torus, whose intersection with $\fraks$ is either the empty set
or a core curve of the solid torus. We may write $\calw\cap\cals$ as a
disjoint union $\calb_1\cup\calb_2$, where $\calb_m$ is the union of all annuli in $\cals$
which are images under $\rho$ of strong belts of exactly $m$ sides of the components of
$\cals$ containing them. Then 
%$\calb_2$ is properly embedded in
%$\calw$, and
 $\calb_1\subset\partial\calw$. The components of
$\calb_2$
%is  annuli \redcomment{I want to say annuli in components of
  %$\cals$ that are images of strong belts for both sides, which is a
  %little awkward} 
are $\pi_1$-injective properly embedded annuli in
the solid torus components of $\calw$ and therefore separate the
components of $\calw$ that contain them. The system of annuli $\calb_2$
therefore defines a dual graph $\calf$ which is a forest in the sense that
each of its components is a tree; the edges of $\calf$ are  in
bijective correspondence with components of $\calb_2$, and the
vertices with components of $\calw-\calb_2$. In a tree, more than half the
vertices are of valence at most $2$. The forest $\calf$ has no
vertices of valence $0$. Using the topological description of $\calw$
given above, and properties of the characteristic suborbifold, one can
show that all but at most one of the valence-$1$ vertices in a given
component of $\calf$ correspond to components of $\calw-\calb_2$ that
contain components of $\calb_1$; each component of $\calb_1$
corresponds to an annulus contained in a component of $\cals$ having
at least one side which is not strongly belted, and is therefore not
fully booked by (a). One can
also use properties of the characteristic suborbifold to
show that any edge incident to a valence-$2$ vertex of $\calf$
corresponds to a component of $\calb_2$ which is 
contained in a component having
at least one side which is not fully booked (although it is
necessarily strongly belted). These observations are easily combined
to give the required upper bound for the number of components
of $\cals$ in terms of the number of non-fully booked components of
$\tcals$.

{\bf Sketch of the proof of Theorem B, concluded: some hints about the
  general case.}
The argument that we have sketched above is given in detail in the proofs of
Proposition \ref{new enchanted forest}, \ref{three or better}, and
\ref{alt twice times lemma}, but the context is  more
technical. In place of Condition (a) the argument assumes that the
system $\cals$  is ``\dandy'' in a sense defined in Section
\ref{semi-dandy section}. \Dandy\ systems are a class of admissible
systems of spheres that satisfy a weaker version of
Condition (a): a fully booked component of such a system is either
strongly belted is of one of several other special types. The
advantage of working with \dandy\ systems is that a \dandy\ complete
system of spheres exists in every closed orientable hyperbolic
$3$-orbifold: this is established by combinatorial arguments in
Section \ref{dandy existence section}. These arguments are based on
modifications of systems of spheres that similar in
spirit to the modifications of systems of tori that were mentioned in
the sketch of the proof of Proposition E given above.

In Propositions \ref{new enchanted forest}---\ref{alt twice times
  lemma}, Condition (c) above is replaced by the condition that
$\tcals$ has no ``\centralclashsphere s'' in the sense defined in
Section \ref{clash section}. A component $S$ of $\cals$ which violates
Condition (c), in the sense that both sides of $S$ are belted, but
the images of their belts under $\rho$ are non-isotopic in $S$, is an
example of a \centralclashsphere, but in the general setting it is
only one of several types that are allowed by the definition. In
general one does not expect to find a complete \dandy\ system of
spheres without \centralclashsphere s in $M=|\oldOmega|$ for a general closed orientable
hyperbolic $3$-orbifold $\oldOmega$. However, it can be shown (see
Propositions \ref{dandy subsystem}, \ref{fine and dandies exist} and \ref{no wonder 1}) that if $\cals$ is a
\dandy\ system in $M$ having at least two components, and if $S$ is a
\centralclashsphere\ of $\cals$ whose sides are
both fully booked, then
$\cals':=\cals-S$ is also a \dandy\ system; and that if $\oldTheta$ and
$\oldTheta'$ denote the suborbifolds of $\oldOmega$ with
$|\oldTheta|=\cals$ and $|\oldTheta'|=\cals'$, and if we set
$\oldPsi=\oldOmega\cut\oldTheta$ and
$\oldPsi'=\oldOmega\cut{\oldTheta'}$, then  in the notation introduced 
% explained
above we have
$X_{\oldTheta'}
% $-\chi(\overline{\oldPsi'-\oldSigma_{\oldPsi'}})
\ge X_\oldTheta
%-\chi(\overline{\oldPsi-\oldSigma_\oldPsi})
+c$ for some universal
  constant $c$.
Hence if a system $\cals^{(m)}$ is obtained from $\cals$ by repeating
$m$ times the operation of removing a  \centralclashsphere, then
$\cals^{(m)}$ is  a \dandy\ system; and  if 
$\oldTheta^{(m)}$ denote the suborbifold with
$|\oldTheta^{(m)}|=\cals^{(m)}$, 
%
%and if we set
 %$\oldPsi^{(m)}=\oldOmega\cut{\oldTheta^{(m)}}$, 
then 
$mc\le X_\oldTheta
%-\chi(\overline{\oldPsi-\oldSigma_{\oldPsi})}
+mc\le
X_{\oldTheta^{(m)}}\le X$
%$-\chi(\overline{\oldPsi^{(m)}-\oldSigma_{\oldPsi^{(m)}}})
%\ge X_\oldTheta
%-\chi(\overline{\oldPsi-\oldSigma_{\oldPsi})}
%+mc$. 
(cf. Corollary \ref{toss 'em
    out}). 
%The results cited here are not actually stated in terms of
%  \centralclashsphere s, but in terms of the closely related notion of
  %``\clashspheres''; however, in view of other technical results proved
  %in Section \ref{clash section}, this does not affect the truth of the
  %ststements made here.) 
%is dealt with} 
% But the properties
  %of $\volorb$ discussed above imply that
  %$-\chi(\overline{\oldPsi^{(m)}-\oldSigma_{\oldPsi^{(m)}}})$ is
    %bounded above by a constant multiple of
    %$\volorb(\oldPsi^{(m)}\le\vol\oldOmega$. 
Thus an upper bound
    on the quantity $X$ 
%\redcomment{
%I just realized I should have been
%      talking about the max of $X_\frakZ$ where $\frakZ$ ranges over
  %    all unions of components of $\oldTheta$. Fix this, and
 %fix other
   %   direct refs to
%$\vol\oldOmega$ below}
%  \redcomment{I have been talking as if I always
%      use an upper bound on the number of fully booked components, but
  %    this is one place where I used more than that. Another is the
    %  very last step, where I use that my sphere is weight $5$. The
  %    discussion above where I claim to be doing everything in terms
  %    of the number of fully booked components needs to be corrected}
 gives an upper bound $b$ on the number of times the
    operation can be repeated. Hence there is a \dandy\ system
    $\cals^*\subset\cals$ such that $\cals-\cals^*$ has at most $b$
    components and $\cals^*$ has no \centralclashsphere. This shows
    that in order  to find an upper bound for the number of components
    of a \dandy\ system in $\oldOmega$, it suffices
to find an upper bound for the number of components
    of a \dandy\ system $\cals^*$ without \centralclashsphere s, as the required
    bound will then be found by adding $b$.

We have argued that  a \dandy\ system $\cals^*$ without
\centralclashsphere s satisfies the counterpart of Conditions (a) and
(c) above. In order to bound the number of components
    of $\cals^*$, by carrying out an argument like the one sketched
above with $\cals^*$ playing the role of $\cals$,  we need the analogue of Condition (b). It turns out
    that for a \dandy\ system, the appropriate analogue of Condition
    (b) is always true if the system is {\it balanced} in the sense
    that its components all have the same weight. 

The hypothesis that $\oldOmega$ contains no embedded negative
turnovers implies that every component of $\cals^*$ has weight at
least $4$. If all the components of $\cals^*$ have weight exactly $4$,
then $\cals^*$ is balanced, and the argument can be carried out
following the sketch above. If $\cals^*$ has a component $S_0$ of weight at
least $5$, we can consider $S_0$ as forming a \dandy\ system in its
own right, and it turns out that a variant of the arguments sketched
above, using that $S_0$ has weight at least $5$, gives a lower bound
for $X_{S_0}$, and therefore for $X$, 
%$\vol\oldOmega$ 
which is strong enough to contradict the
hypotheses of Theorem B. 

{\bf The role of Chapter 4.}
In this  sketch of the proof of Theorem B, we have
mentioned the role of
Propositions \ref{fine and dandies exist} and \ref{no
  wonder 1}. It will be apparent from the discussion above that the
content of these results involves the relationship between the
characteristic suborbifold $\oldSigma_\oldPsi$ of
$\oldPsi=\oldOmega\cut\oldTheta$, where $\oldTheta$ is an
incompressible $2$-suborbifold of $\oldOmega$, and the characteristic
suborbifold 
$\oldSigma_{\oldPsi'}$ of $\oldPsi'=\oldOmega\cut{\oldTheta'}$, where
$\oldTheta'$ is a union of components of
$\oldTheta$. (In
the situation considered above, $\cals$ and $\cals'=\cals-S$ are the
underlying surfaces of $\oldTheta$ and $\oldTheta'$.) The problem
of describing $\oldSigma_{\oldPsi'}$ in terms of $\oldSigma_\oldPsi$
is analogous to the problem, addressed in \cite{bcsz} and elsewhere,
of describing ``reduced
    homotopies'' in a $3$-manifold equipped with an
    incompressible surface.  Chapter 4 is devoted to developing the
    machinery needed to do this; this constitutes a partial
    generalization to $3$-orbifolds  of  the machinery developed in \cite{bcsz}.
The 
canonical sequence $(V_n)_{n\ge1}$ of isotopy classes
of suborbifolds of $\toldTheta:=
\partial\oldPsi$
 defined in Section \ref{higher section} is analogous to the sequence
denoted $(\Phi_k^{\pm})_{k\ge1}$ in \cite{bcsz}; the terms of the latter
sequence are
subsurfaces, defined up to isotopy, of a bounded essential surface in a
$3$-manifold with torus boundary. The statement of Propositions \ref{fine and dandies exist} and \ref{no
  wonder 1} involve the notion of a ``\clashsphere,'' which is closely
related to the notion of a \centralclashsphere\ but is defined
in terms of the suborbifolds $V_n$ defined in Chapter 4.
The
argument given to cover Case I of the proof of Proposition
\ref{no wonder 1}, which requires the full machinery of Chapter 4,
seems indispensable for  establishing an inequality of the form $X_{\oldTheta'}
\ge X_\oldTheta
%-\chi(\overline{\oldPsi-\oldSigma_\oldPsi})
+c$ as in the discussion above.

{\bf Acknowledgments.}
I am grateful to Ian Agol, Francis Bonahon, Ted Chinburg, Marc Culler,
Benson Farb, Dave Futer, Tom Goodwillie, Tracy Hall, Christian Lange, Chris Leininger, Ben Linowitz, Shawn Rafalski, Matt Stover, and  John Voight for valuable discussions and encouragement. More specifically, Agol told me how to prove
  Proposition \ref{hepcat}; Rafalski read some of my first efforts to write about orbifolds and corrected a number of errors;  Hall and Goodwillie helped lead me to the proof of Proposition \ref {what i need?}; and Hall helped me find the correct statement and proof of Proposition \ref{vector thing}.  Some of the work presented here was done during visits to the Technion and the IAS, and I am grateful to these institutions
for their hospitality. Some of the work was also supported by NSF grant DMS-1207720.
%certain tion B
%\oldOmega' M' \Mh' \Mh h(M) \tMh \pl \chi $T $P\fraks\lambda

{\bf Dedication.} This work is dedicated to the memory of my parents,  Marcia Heatter Shalen and Robert Ellis Shalen.

\chapter{Preliminaries}\label{prelim chapter}

After introducing some very general conventions in Section \ref{prelim
  section}, we establish in Section \ref{manifolds section} a number
of technical results, very similar to results established in
\cite{Waldhausen} and \cite{FHS}, that will be needed later in the
monograph, especially in Section \ref{darts section} where we study
the basic properties of the invariant $\volorb$ which was discussed in
the Introduction. In Section \ref{orbifolds section} we establish some
conventions and very basic facts regarding orbifolds. Section
\ref{Orb3 section} is devoted more specifically to $3$-orbifolds, and
includes many orbifold analogues, not readily available in the
literature, of basic standard facts about
$3$-manifolds.

\section{General conventions}\label{prelim section}

\abstractcomment{\tiny A lot is tied up with the properties of the characteristic
  suborbifold. }

\Number
The set of all non-negative integers will be denoted $\NN$.
The cardinality of a finite set $S$ will be denoted $\card S$.
If $A$ and $B$ are subsets of a set, $A\setminus B$ will denote the
set of all elements of $A$ which do not belong to $B$. In the special
case where $B\subset A$, we will often use the alternative notation $A-B$.
A disjoint union of sets $A$ and $B$ will be denoted $A\discup B$: that is, writing $X=A\discup B$ means that $A\cap B=\emptyset$ and that $X=A\cup B$.
\EndNumber

\Number
If $X$ is a topological space, we will
denote by $\calc(X)$ the set of all connected components of $X$.  I
will set $\compnum(X)=\card \calc(X)$.
\EndNumber

\Number\label{injectify}
A map $f:X\to Y$ between path-connected spaces will
be termed $\pi_1$-injective if $f_\sharp:\pi_1(X)\to
\pi_1(Y)$ is injective. (Here and elsewhere, base points are
suppressed from the notation in cases where the choice of a base point
does not affect the truth of a statement.)

In general, a (continuous) map $f:X\to Y$ between arbitrary spaces $X$ will be
termed $\pi_1$-injective if each path component $C$ of $X$, the map
$f|C$ is a
%from  is
$\pi_1$-injective map from $C$ to the path component of $Y$ containing
$f(C)$.

A subset $A$ of a space $X$ is termed $\pi_1$-injective if the
inclusion map from $A$ to $X$ is $\pi_1$-injective.
\EndNumber

\Number\label{in and out soda} If $X$ is a topological space having
the homotopy type of a finite CW complex, the Euler characteristic
of $X$ will be denoted by $\chi(X)$, and we will set
\Equation\label{in and out equation}
\chibar(X)=-\chi(X).
\EndEquation
 We will denote by $h(X)$ the dimension of the
singular homology $H_1(X;\FF_2 )$ as an $\FF_2 $-vector space. 
%If $Y$ is
%a subset of $X$, we will denote by $h(Y\to X)$ the dimension of the
%image of the inclusion map $H_1(Y;\FF_2 )\to H_1(X;\FF_2 )$.
\EndNumber

\Number\label{search engine}
As the overview of our methods given in the introduction indicates, our argument depend heavily on the interaction between properties of an orientable $3$-orbifold and properties of its underlying $3$-manifold. Subsections \ref{nbhd stuff}---\ref{fhs-prop} are devoted to conventions and results concerning $3$-manifolds that will be used in this monograph. While some of the concepts involved will be generalized to orbifolds later, emphasizing the manifold case here will help pave the way for using the very rich literature on $3$-manifolds that is available.
\EndNumber

\section{Manifolds}\label{manifolds section}

\Number\label{nbhd stuff}
Let $M$ be a (topological, PL or smooth) compact manifold. A (respectively topological, PL or smooth) submanifold $\cals$ of $M$ will be termed {\it two-sided} if it has a neighborhood $N$ in $M$ such that the pair $(N,\cals)$ is (topologically, piecewise-linearly or smoothly) homeomorphic to $(\cals\times[-1,1],\cals\times\{0\})$. Such a neighborhood $N$ is called a {\it collar neighborhood} of $\cals$. Note that two-sidedness implies that $\cals$ is properly embedded in $M$ and has codimension $1$.
If $\cals$ is two-sided,  we will denote by $M\cut\cals$ the compact space obtained by splitting
$M$ along $\cals$. If $M$ is a topological or PL manifold, or if $M$ is smooth and $\cals$ is closed, then $M\cut\cals$ inherits the structure of a topological, PL or smooth manifold respectively. 
We will denote
by $\rho_\cals$ the natural surjection from $M\cut\cals$ to $M$. For
each component $S$ of $\cals$,  the set $\rho_\cals^{-1}(S)$
is the
union of two components of $\tcals:=\rho_\cals^{-1}(\cals)$. We will call these
components the {\it sides} of $S$. 
We will denote by $\grock_\cals$ the unique
involution $\grock$ of $\tcals$ which interchanges the sides of every component of $\cals$ and satisfies
$(\rho|\tcals)\circ\grock=\rho|\tcals$. 

We will regard $M\cut\cals$ as a completion of
$M-\cals$; in particular, $M-\cals$ is identified with $M\cut\cals-\rho_\cals^{-1}(\cals)$ (which is  the interior of
$M\cut\cals$ in the case where $M$ is closed). If $X$ is a union of components of $M-\cals$, its closure in
$M\cut\cals$ will be denoted $\hatX$. Note that $X\mapsto\hatX$ is a
bijection between the components of $X-\cals$ and those of
$X\cut\cals$. 
\EndNumber

\Number\label{dual graph}
By a {\it graph} we will mean a CW complex of dimension at most $1$. The $0$-cells and $1$-cells of a graph will be called its {\it vertices} and {\it edges} respectively.

Recall that if $\cals$ is a two-sided submanifold of a manifold $M$, 
%We conclude this section by reviewing the definition of
 the {\it dual graph} $\calg$ of $\cals$ in $M$ is defined as follows. The vertices of $\calg$ are in
bijective correspondence with the components of $M-\cals$; the
edges of $\calg$ are in bijective correspondence with the components
of $\cals$,
 %We also fix, for each 
% a biTo each endpoint $c$ of the edge corresponding to a 
%component $S$ of $\cals$, a bijective corresondence between the endpoints of $ we assign a component   
and a vertex $v$ of $\calg$ is incident to an edge $e$ if
and only if the component of $\cals$ corresponding to $e$ is contained
in the closure of the component of $M-\cals$ corresponding to
$v$. In particular, if a component $S$ of $\cals$ is contained in the closure of only component $C$ of $M-\cals$ (because both frontier components of a collar neighborhood of $S$ are contained in $C$), then the edge corresponding to $S$ is a loop. The dual graph of $\cals$ in $M$ is homeomorphic to a retract of $M$.
\EndNumber

\Proposition\label{dual homology}
Let $\cals$ be a two-sided submanifold of a compact (possibly disconnected) manifold $M$, and let $\beta$ denote the first betti number of 
 the dual graph of $\cals$ in $M$. Then we have $h(M)\le h(M-\cals)+\beta$, where $h(M)$ and $h(M-\cals)$ are  defined by \ref{in and out soda}. 
\EndProposition

\Proof
Let $\calg$ denote the dual graph of $\cals$ in $M$.
We argue by induction on $\compnum(\cals)$. If $\compnum(\cals)=0$, so that $\cals=\emptyset$, then $\calg$ has no edges and therefore $\beta=0$. In this case we have  $h(M)=h(M-\cals)=h(M-\cals)+\beta$. Now suppose that $\compnum(\cals)=n>0$, and that the assertion is true whenever $M$ and $\cals$ are replaced by a compact manifold $M'$ and a two-sided submanifold $\cals'$ of $M'$ with $n-1$ components. Fix a component $S$ of $\cals$, and set $M'=M\cut S$. Then $\cals':=\cals-S$ is two-sided in $M'$, so that $h(M')\le h(M'-\cals')+\beta'$, where $\beta'$ is the first betti number of the dual graph $\calg'$ of $\cals'$ in $M'$. We may write $\calg'=\calg-e$, where $e$ is the (open) edge of $\calg$ corresponding to $S$. Note also that $h(M'-\cals')=h(M-\cals)$.

If $S$ separates $M$ then $e$ separates $\calg$, so that $\beta'=\beta$. But in this case the Mayer-Vietoris theorem gives $h(M)\le h(M')$, so that $h(M)\le h(M'-\cals')+\beta'=h(M-\cals)+\beta$.

If $S$ does not separate $M$ then $e$ does not separate $\calg$, so that $\beta'=\beta-1$. But in this case the Mayer-Vietoris theorem gives $h(M)\le h(M')+1$, so that $h(M)\le h(M'-\cals')+\beta'+1=h(M-\cals)+\beta$.
\EndProof

\Number\label{esso}
The standard $1$-sphere $\SSS^1\subset\CC$ inherits a PL structure from $\RR$ via the covering map $t\mapsto e^{2\pi it}$, since the deck transformations $t\mapsto t+n$ are PL. This PL structure on $\SSS^1$ gives rise to a PL structure on $\DD^2$ by coning. When we work in the PL category, $\SSS^1$ and $\DD^2$ will always be understood to have these PL structures. Note that the self-homeomorphism of $\SSS^1$ or $\DD^2$ defined by an arbitrary element of $\OO(2)$ is then piecewise linear.
\EndNumber

\Number\label{just manifolds}
From this point on, all statements and arguments about manifolds are to be interpreted in the PL category except where another category is specified.
\EndNumber

\Number\label{great day}
In large part we will follow the conventions of \cite{hempel}
regarding $3$-manifolds. We will depart slightly from these
conventions in our use of the term ``irreducible'': we define a
$3$-manifold $M$ to be {\it irreducible} if $M$ is connected, every
(tame) $2$-sphere in $M$ bounds a $3$-ball in $M$, and $M$ contains no
homeomorphic copy of $\RR\PP^2\times[-1,1]$. Thus $M$ is ``connected
and $P^2$-irreducible'' in the sense of \cite{hempel}. The reason for
using the term ``irreducible'' in this stronger sense is that, unlike
the more classical definition, it coincides in the manifold case with
the definition of an irreducible orbifold to be given in Subsection \ref{oops}.

We will use the word ``incompressible'' only in the context of closed surfaces. A closed (possibly disconnected) surface $F$ in a $3$-manifold $M$ will be termed {\it incompressible} if $F$ is two-sided and $\pi_1$-injective in $M$, and has no component which bounds a ball in $M$. This is also consistent with the definition to be given below for orbifolds. 

As in \cite{hempel}, an irreducible $3$-manifold $M$ will be termed
{\it boundary-irreducible} if $\partial M$ is
 $\pi_1$-injective in $M$; by Dehn's lemma and the loop theorem, this is equivalent to the condition that
for every properly embedded disk
$D\subset M$ there is a disk $E\subset\partial M$ such that $\partial
E=\partial D$. We will say that an orientable $3$-manifold $M$ is {\it acylindrical} if it is
irreducible and boundary-irreducible, and every properly embedded
$\pi_1$-injective annulus in $M$ is
boundary-parallel. 
%{This will presumably be subsumed in the
  %later def. for orbifolds, as should be mentioned there. But the
%  def. has to be here, because it is used in Lemma \ref{torus goes to
  %  cylinder} and its apps in Sections \ref{underlying tori section} and \ref{irr-M section}. I think %I will add a result to Section \ref{tori section} about the relationship between acyclindricity and ``atorality'' for manifolds.

%Speaking of things that are subsumed, I think one of the very few mentions of characteristic submanifold in the paper is in the proof of Lemma \ref{graphology}, and that it could be defined there by specializing the def. of the characteristic suborbifold. Decide how to handle such things.}

A {\it graph manifold} is a closed, irreducible, orientable $3$-manifold
$M$ which contains a $\pi_1$-injective $2$-dimensional submanifold
$\calt$ such that every component of $\calt$ is a torus and every
component of $M\cut\calt$ is a Seifert fibered space. 

When $A$ is an annulus contained in the boundary of a solid torus $J$,
we will define the {\it winding number} of $A$ in $J$ to be the order
of the cyclic group $H_1(J,A;\ZZ)$ if this cyclic group is finite, and
to be $0$ if the cyclic group is infinite.
\EndNumber

\Definition\label{P-stuff}
If $X$ is a $3$-manifold, we will denote by $\plusX$ the $3$-manifold
obtained from $X$ by attaching a ball to each component of $\partial
X$ which is a $2$-sphere. We will say that $X$ is {\it $+$-irreducible}
if $\plusX$ is irreducible. We will say that $X$ is a {\it $3$-sphere-with-holes}
%, 
%or a {\it solid torus-with-holes}, \redmissingref{see comment about the terminology in the statement of Lemma \ref{uneeda}}
 if $\plusX$ is a $3$-sphere.
% or a solid torus, respectively.
\EndDefinition

\Number\label{plus-contained}
Note that if $K$ is a compact submanifold of a compact $3$-manifold
$N$, and if every component of $\partial K$ is either a component of
$\partial N$ or a surface of positive genus contained in $\inter N$,
then $K^+$ is naturally identified with a submanifold of $N^+$.
\EndNumber

The following slightly stronger version of \cite[Corollary
5.5]{Waldhausen} will be needed in Section \ref{darts section}:

\Proposition\label{stronger waldhausen}
Let $M$ be a closed, irreducible, orientable $3$-manifold. Let $T$ be
a closed, orientable $2$-manifold, no component of which is a sphere,
and let $i:T\to M$ and $j:T\to M$
be $\pi_1$-injective embeddings. Suppose that $i$ and $j$ are homotopic, and that the
components of $i(T)$ are pairwise non-parallel in $M$. Then $i$ and
$j$ are isotopic.
\EndProposition

\Proof
We argue by induction on $\compnum(T)$. If $\compnum(T)=0$ the
assertion is
  trivial, and if $\compnum(T)=1$ it is \cite[Corollary
5.5]{Waldhausen}. Now suppose that an integer $n\ge1$ is given and
that the result is true whenever
$\compnum(T)=n$. Let $M$, $T$, $i$ and $j$ be given,
satisfying the hypothesis, with
%For the induction step, suppose
$\compnum(T)=n+1$. We may assume without loss of generality that $T\subset M$ and that $i:T\to M$ is
the inclusion map. Thus $T$ is $\pi_1$-injective. Choose a  component $V$ of $T$, and set
$T'=T-V$. Since $j:T\to M$ is homotopic to the inclusion, so is
$j':=j|T':T'\to M$. Furthermore, $T'$ is in particular
$\pi_1$-injective, and $\compnum(T')=n$. Hence the induction
hypothesis implies that 
$j|T'$ is isotopic to the inclusion map. Thus after modifying 
$j$ within its isotopy class we may assume that $j|T'$ is the
inclusion. Now $j|V$ is homotopic to the inclusion $i_V:V\to M$. Since
$j$ and the inclusion map $i$ are embeddings, both $V$ and $j_V(V)$
are disjoint from $T'$.

We will show that $i_V$ and $j_V$ are isotopic by an ambient isotopy
of $M$ which is constant on $T'$. This will immediately imply that $i$
and $j$ are isotopic, and complete the induction.

Set $W=j_V(V)$.
We may assume after an isotopy that $W$ and $V$ intersect
transverally. We may further assumed that $j$ has been chosen within
its isotopy class rel $T'$ so as
to minimize $\compnum(W\cap V)$, subject to the
condition that $W$ and $V$ intersect transversally. We now claim:
\Equation\label{nothing there}
W\cap V=\emptyset.
\EndEquation

To prove (\ref{nothing there}), suppose that $W\cap V\ne\emptyset$. Since $j_V$ and the inclusion $i_V$ are homotopic in
$M$ and $\pi_1$-injective, it follows from
\cite[Proposition 5.4]{Waldhausen} that there
exist connected subsurfaces $A\subset V$ and $B\subset W$, and a compact submanifold
$X$ of $M$, such that
$\partial A\ne\emptyset$, $\partial X=A\cup B$,
%$\xi (A\times\{0\})=A$, $\xi ((A\times\{1\})\cup((\partial
%A)\times[0,1]))=B$, 
and the pair $(X,A)$ is homeomorphic to $(A\times[0,1],A\times\{0\})$.
%$A\cap W=\partial A$. 
Since $V$ and $W$
are disjoint from $T'$, we have $T'\cap\partial X=\emptyset$; hence
every component of $T'$ is either contained in $\inter X$ or disjoint
from $X$.
Note that $X$ is a handlebody since $\partial A\ne\emptyset$. If some
component $U$ of $T'$ is contained in $\inter X$, then the inclusion
homomorphism $\pi_1(U)\to\pi_1(X)$ is an injection from a
positive-genus surface group to a free group, which is
impossible. Hence $T'\cap X=\emptyset$. On the other hand, the
properties of $A$, $B$ and $X$  listed above imply that there is
an isotopy $(h_t)_{0\le t\le 1}$ of $M$, constant outside an
arbitrarily small neighborhood of $X$, such that $h_0$ is the identity
and $\compnum(h_1(W)\cap V)<\compnum(W\cap V)$. Since $T'\cap
X=\emptyset$, we may take $(h_t)$ to be constant on $T'$. This implies
that $j_V^1:=h_1\circ j$ is isotopic rel $T'$ to $j$, and that
$W_1:=j_V^1(V)$ satisfies $\compnum(W_1\cap V)<\compnum(W\cap V)$. This
contradicts the minimality of $\compnum(W\cap V)$, and (\ref{nothing
  there}) is proved.

Since $i_V$ and $j_V$ are homotopic and are $\pi_1$-injective, and $i_V(V)$ and $j_V(V)$ are disjoint by (\ref{nothing
  there}),  it
follows from \cite[Proposition 5.4]{Waldhausen}
that 
there
exists a compact submanifold
$Y$ of $M$, and a homeomorphism $\eta :V\times[0,1]\to Y$, such that
$\eta(v,0)=i_V(v)=v$ and $\eta(v,1)=j_V(v)$ for every $v\in V$. 
Since $V$ and $W$
are disjoint from $T'$, we have $T'\cap\partial Y=\emptyset$; hence
every component of $T'$ is either contained in $\inter Y$ or disjoint
from $Y$.
If some
component $U$ of $T'$ is contained in $\inter Y$, then since $U$ is
incompressible and $Y$ is homeomorphic to $V\times[0,1]$, it follows
from \cite[Proposition 3.1]{Waldhausen} that $U$
is parallel to each of the boundary components of $Y$, and in
particular to $V$. This contradicts the hypothesis that no two
components of $T=i(T)$ are parallel. Hence $T'\cap Y=\emptyset$. On the other hand, the
properties of $Y$ and $\eta $ stated above imply that there is
an isotopy $(h_t)_{0\le t\le 1}$ of $M$, constant outside an
arbitrarily small neighborhood of $Y$, such that $h_0$ is the identity
and $h_1\circ i_V=j_V$. Since $T'\cap
Y=\emptyset$, we may take $(h_t)$ to be constant on $T'$. This implies
that $j_V$ is isotopic rel $T'$ to $i_V$, as required.
\EndProof

The following result will also be needed in Section \ref{darts section}:

\Proposition\label{snuff}
Let $M$ be a closed, irreducible, orientable $3$-manifold. Let $V$ and $W$ be closed, orientable surfaces of strictly positive genus, and let $i_V:V\to M$ and $i_W:W\to M$ be $\pi_1$-injective embeddings such that $i_V(V)\cap i_W(W)=\emptyset$. Then  for any embeddings $f:V\to M$ and $g:W\to M$ homotopic to $i_V$ and $i_W$ respectively, such that $f(V)$ and $g(W)$ meet transversally, either (a) $f(V)\cap g(W)=\emptyset$, or (b) there
exist connected subsurfaces $A\subset f(V)$ and $B\subset g( W)$, and a compact submanifold
$X$ of $M$, 
%and a homeomorphism $\xi :A\times[0,1]\to X$, 
such that
$\partial A=\partial B\ne\emptyset$,  $\partial X=A\cup B$, and the pair $(X,A)$ is homeomorphic to
$(A\times[0,1],A\times\{0\})$.
%, and $\xi ((A\times\{1\})\cup((\partial
%A)\times[0,1]))=B$.
\EndProposition

\Proof
According to \cite[Corollary 5.5]{Waldhausen}, $g$ is isotopic to $i_W$. Hence after modifying $f$ and $g$ by a single ambient isotopy, we may assume that $g=i_W$. We may also assume that $V$ and $W$ are subsurfaces of $M$ and that $i_V$ and $g=i_W$ are the inclusion maps. The hypothesis then gives $V\cap W=\emptyset$. Then $f(V)$ meets $W$ transversally, and after a small isotopy, constant on $W$, we may assume that it also meets $V$ transversally. We may also suppose that among all embeddings in its isotopy class rel $W$, having the property that $f(V)$ meets $V$ transversally, $f$ has been chosen so as to minimize $\compnum(f(V)\cap V)$. Set $V'=f(V)$. Since $V'$ and $V$ are homotopic by hypothesis, it follows from \cite[Prop 5.4]{Waldhausen} that there
exist (not necessarily proper) connected subsurfaces $A_0\subset V'$ and $C\subset V$, a compact submanifold
$X_0$ of $M$, such that $\partial X_0=A_0\cup C$,
%$\xi_0 (A_0\times\{0\})=A_0$, $\xi_0 ((A_0\times\{1\})\cup((\partial
%A_0)\times[0,1]))=C$, 
the pair $(X_0,A_0)$ is homeomorphic to $(A_0\times[0,1],A_0\times\{0\})$, and
$A_0\cap V=\partial A_0$. We claim:

\Claim\label{fdjt}
Either $W\cap A_0\ne\emptyset$, or Alternative (a) of the conclusion of the lemma holds.
\EndClaim

To prove \ref{fdjt}, assume that $W\cap A_0=\emptyset$. In the case where $A_0=V'$, it follows that $W\cap V'=\emptyset$, which is Alternative (a). Now suppose that $A_0$ is a proper subsurface of $V'$, so that $\partial A_0\ne\emptyset$. The properties of $A$, $C$ and $X_0$ stated above then imply that $X_0$ is a handlebody, and that there is an embedding $f_1:V\to M$, isotopic to $f$ by an ambient isotopy which is constant on an arbitrarily small neighborhood of $X_0$, such that $\compnum(f(V_1)\cap V)<\compnum(V'\cap V)=\compnum(f(V)\cap V)$. The assumption $W\cap A_0=\emptyset$, together with the fact $V\cap W=\emptyset$, implies that $W$ is disjoint from $\partial X_0$. Since $W$ is a closed $\pi_1$-injective orientable surface of positive genus, it cannot be contained in the handlebody $X_0$. Hence $W\cap X_0=\emptyset$. We may therefore take the ambient isotopy between $f$ and $f_1$ to be constant on $W$. But then the inequality $\compnum(f(V_1)\cap V)<\compnum(f(V)\cap V)$ contradicts the minimality property of $f$. Thus \ref{fdjt} is proved.

Next, we claim:
\Claim\label{folderol}
If some component of $V'\cap W$ is homotopically trivial in $M$, then Alternative (b) of the conclusion holds.
\EndClaim

To prove \ref{folderol}, note that any homotopically trivial component of $V'\cap W$ must bound a disk in $V'$, since $V'$ is incompressible. Among all disks in $V'$ bounded by components of $V'\cap W$ choose one, $A$, which is minimal with respect to inclusion. Since $W$ is also incompressible, $\partial A$ bounds a disk $B\subset W$. The minimality of $D$ implies that $A\cap B=\partial A$, so that $A\cup B$ is a $2$-sphere; by irreducibility, $A\cup B$ bounds a $3$-ball $X\subset M$. Now the pair $(X,A)$ is homeomorphic to $(A\times[0,1],A\times\{0\})$, and hence Alternative (b) holds. Thus \ref{folderol} is established.

In view of \ref{fdjt} and \ref{folderol}, we may assume that  $W\cap A_0\ne\emptyset$ and that every component of $V'\cap W$ is homotopically non-trivial in $M$. Since $W\cap A_0\ne\emptyset$, there is a component $B$ of $W\cap X_0$ with $\partial B\ne\emptyset$. We have $B\cap C\subset W\cap V=\emptyset$, and hence $\partial B\subset A_0$. In particular, each component of $\partial B$ is a component of $V'\cap W$ and is therefore homotopically non-trivial in $M$, and in particular in $W$; this implies that $B$ is $\pi_1$-injective in the incompressible surface $W$, and is therefore $\pi_1$-injective in $M$. Thus $B$ is a properly embedded connected, $\pi_1$-injective surface in $X_0$, with $\partial B\subset A_0$. Since the pair $(X_0,A_0)$ is homeomorphic to $(A_0\times[0,1],A_0\times\{0\}$, it now follows from
\cite[Proposition 3.1]{Waldhausen} that $B$ is parallel in $X_0$ to a subsurface $A$ of $A_0$. This means that there is a submanifold $X$ of $X_0$ such that $\partial X=A\cup B$ and $(X,A)$ is homeomorphic to $(A\times[0,1],A\times\{0\})$. We have $\partial A=\partial B\ne\emptyset$. This gives Alternative (b) of the conclusion.

 \EndProof

\Proposition\label{final assertion}
Let $L$ be a connected, compact, $3$-dimensional
submanifold of an irreducible, orientable $3$-manifold $M$. Suppose
that $L$ is $\pi_1$-injective in $M$, that every component of $\partial
L$ is a torus, and that $L$ is not a solid torus. Then every component of
$\partial L$ is $\pi_1$-injective in $M$.
\EndProposition

\Proof Consider any component $T$ of 
$\partial L$. Since $L$ is
$\pi_1$-injective in $M$, it suffices to show that $T$ is
$\pi_1$-injective in $L$. If it is not, there is a properly embedded
disk $D\subset L$ whose boundary does not bound a disk in $T$. Let $Y$
denote a regular neighborhood of $D$ in $L$. Then $Q:=\overline{L-Y}$ is a $3$-manifold, and the component $S$ of $\partial Q$ containing $T\setminus Y$ is a $2$-sphere. Since $M$ is irreducible, $S$ must bound a ball $B\subset
M$. We have either $B\supset Q$ or $B\cap
Q=S$. If $B\supset Q$, then $Z=B\cup Y$ is a solid torus containing
$L$ and having boundary $T$. Since $L$ is $\pi_1$-injective in $M$ and is contained in the solid torus $Z$, it has a
cyclic fundamental group. The torus $T$ is one component of $\partial 
L$. Since $\pi_1(L)$ is cyclic, $\partial L$ cannot have a second
component of positive genus. Hence $L=Z$, which contradicts the hypothesis that $L$ is not a solid torus. On the other hand, if $B\cap
Q=S$, then $T\subset B$, so that the inclusion
homomorphism $\pi_1(T)\to\pi_1(M)$ is trivial. Since $L$ is
$\pi_1$-injective in $M$, the inclusion
homomorphism $\pi_1(T)\to\pi_1(L)$ is trivial. This is impossible,
because Poincar\'e-Lefschetz duality implies that the inclusion
homomorphism $H_1(T;\QQ)\to\pi_1(L;\QQ)$ is non-trivial. 
\EndProof

We will need the following result, the proof of which we will extract from \cite{FHS}:

\Proposition\label{fhs-prop}
Let $M$ be an orientable Riemannian $3$-manifold, let $V$ and $W$ be closed orientable $2$-manifolds, and let $f_0:V\to M$ and $g:W\to M$ be smooth embeddings, each of which has least area in its homotopy class. Suppose that for every embedding $f:V\to M$, homotopic to $f_0$,
%, and every embedding $g_0:W\to M$ homotopic to $g$ 
such that $f(V)$ and $g(W)$ meet transversally, either (a) $f(V)\cap g(W)=\emptyset$, or (b) there exist connected subsurfaces $A\subset f( V)$ and $B\subset g(W)$, and a compact submanifold
$X$ of $M$, 
%and a homeomorphism $\xi :A\times[0,1]\to X$, 
such that
$\partial A\ne\emptyset$, $\partial X=A\cup B$, and the pair $(X,A)$ is homeomorphic to $(A\times[0,1],A\times\{0\})$.
%$\xi (A\times\{0\})=A$, and $\xi ((A\times\{1\})\cup((\partial
%A)\times[0,1]))=B$.
% and $X\cap(f(V)\cup g(W))=\partial X$. 
Then $f_0(V)\cap g(W)=\emptyset$.
\EndProposition

\Proof
This is implicit in the proof of  \cite[Lemma 1.3]{FHS}. In the language of \cite{FHS}, Alternative (b) of the hypothesis of Proposition \ref{fhs-prop} is expressed by saying that
there is a product region between $f(V)$ and $g(W)$. In the special case where $f_0(V)$ and $g(W)$ meet transversally, we may apply the hypothesis of Proposition \ref{fhs-prop}, taking $f=f_0$, to deduce that either there is a product region between $f_0(V)$ and $g(W)$, or $f_0(V)\cap g(W)=\emptyset$. But according to \cite[Lemma 1.2]{FHS}, there cannot exist a product region between two subsurfaces of a Riemannian $3$-manifold which are the images of smooth embeddings of compact surfaces, each of which has least area in its homotopy class. Hence in this case we must have $f_0(V)\cap g(W)=\emptyset$, as required.

The proof of \cite[Lemma 1.2]{FHS} depends on the observation that if $f:V\to M$ and $g:W\to M$ are smooth embeddings of closed orientable surfaces in an orientable Riemannian $3$-manifold, and if $A$, $B$ and $X$  have the properties stated in Alternative (b), then we may define piecewise smooth embeddings $f':V\to M$ and $g':W\to M$ which agree with $f$ and $g$ on $\overline{V-f^{-1}(A)}$ and $\overline{W-f^{-1}(B)}$ respectively, and such that $f'|f^{-1}(A)$ and $g'|g^{-1}(B)$ are homeomorphisms of their respective domains onto $B$ and $A$, and are homotopic in $X$, rel $f^{-1}(\partial A)$ and $g^{-1}(\partial B)$, to $f|f^{-1}(A)$ and $g|g^{-1}(B)$ respectively. For the purpose of this proof, this construction of a pair of piecewise smooth embeddings $(f',g')$ from a pair of smooth embeddings $(f,g)$, involving a product region between $f(V)$ and $g(W)$, will be called a {\it swap}. If, keeping the same assumptions and notation, $f'':V\to M$ and $g'':W\to M$ are smooth embeddings which are homotopic $f'$ and $g'$ relative to annular neighborhoods of $f^{-1}(A)$ and $g^{-1}(B)$ respectively, and satisfy $f''(V)\cap g''(W)=(f(V)\cap g(W))-f(\partial A)$, we will say that $(f'',g'')$ is obtained from $(f,g)$ by a {\it smoothed swap}. If there is a product region between smooth embeddings $f$ and $g$, then there is a pair $(f'',g'')$ obtained from $(f,g)$ by a smoothed swap such that $\area f''(V)+\area g''(W)<\area f(V)+\area g(W)$; hence either $\area f''(V)<\area f(V)$ or $\area g''(W)<\area g(W)$, so that $f$ and $g$ cannot both have least area in their respective homotopy classes.

If $f_0(V)$ and $g(W)$ do not intersect transversally, it is shown  in the proof of  \cite[Lemma 1.3]{FHS} that there exists a number $\epsilon>0$ with the following property: $f_0$ may be $C^1$-approximated arbitrarily well by an embedding $f$ such that (1) $f$ is homotopic to $f_0$; (2) $f(V)$ and $g(W)$ meet transversally; (3)  if $f_0(V)\cap g(W)\ne\emptyset$ then $f(V)\cap g(W)\ne\emptyset$; and (4) if there is a product region between $f(V)$ and $g(V)$, there is a pair $(f'',g'')$, obtained from $(f,g)$ by a smoothed swap, such that 
%are smoothed andor any \redcomment{I need to explain how $f''$ and $g''$ are related to $f$ and $g$, and make it consistent with the application of (4) below; this will certainly require rewriting the next sentence. This is not in satisfactory shape at all} we have 
$\area f''(V)+\area g''(W)\le\area f(V)+\area g(W)-\epsilon$. By taking $f$ to be a good enough $C^1$-approximation to $f_0$ we can guarantee that
% $h$ is homotopic to $f$, that 
$\area f(V)<\area f_0(V)+\epsilon$. From (1), (2) and the hypothesis, it follows that either there is a product region between $f(V)$ and $g(W)$, or $f(V)\cap g(W)=\emptyset$. 
If there is a product region between $f(V)$ and $g(W)$, then (4) 
gives a pair $(f'',g'')$, obtained from $(f,g)$ by a smoothed swap, such that 
$\area f''(V)+\area g''(W)\le\area f(V)+\area g(W)-\epsilon<\area f_0(V)+\area g(W)$. Hence either $\area f''(V)<\area f_0(V)$, or $\area g''(W)<\area g(W)$; in either case we have a contradiction to the hypothesis that $f_0$ and $g$ have least area in their homotopy classes. We must therefore have $f(V)\cap g(W)=\emptyset$, which by (3) implies 
$f_0(V)\cap g(W)=\emptyset$, as required.
\EndProof
%(b)

\section{Orbifolds}\label{orbifolds section}

\Number\label{orbifolds introduced}
General references for orbifolds include \cite{bmp}, \cite{chk} and \cite{kapovich}. Although smooth orbifolds are emphasized in these books, the definition of orbifold goes through without change in the topological or PL category. (In reading the definition in the PL category, one should bear in mind that an orthogonal action of a finite group on a Euclidean space is in particular a PL action.) 

The material from here to the end of Lemma \ref{oh yeah he tweets} is
meant to be interpretable in each of the three categories, except
where a restriction on category is specified. Of course in the smooth
category the term``homeomorphism'' (of orbifolds) should be understood
to mean ``diffeomorphism.''

\abstractcomment{
One point to remember is that \cite{illman} applies only to very good orbifolds; this should be OK for the apps if I bear it in mind. Another point is that \cite{lange} does  not smooth pairs, only (very good) orbifolds; again I think it will OK for the apps if I bear it in mind. Do I need to think about what stratification, and the other issues in the def of general position (for example), mean in the PL context? They should be simpler, and I'm inclined to leave the issue implicit.

%I had the following passage at one point: ``
%In passages dealing with manifolds, the default category will be the
%{\it topological category} in which the objects are locally Euclidean
%metrizable spaces, submanifolds and embeddings are locally flat,
%ambient isotopies are continuous families of homeomorphisms,  and so
%forth. Many of the results that are quoted depend for their proofs on
%passing to either the piecewise linear or smooth category; as we do
%not consider manifolds of dimension greater than $3$, the details of
%this transition are standard. If $K$ is a closed subset of a manifold
%$M$, a {\it regular neighborhood} of $K$ in $M$ will mean a subset of
%$M$ which is a regular neighborhood of $K$ with respect to some PL
%structure on $M$; it exists if and only if $K$ is a PL subset of $M$
%with respect to some PL structure. A {\it hyperbolic structure} on a
%manifold $M$ is defined by a homeomorphism between $M$ and some
%hyperbolic manifold of finite volume.'' This was followed by a
%missingref: ``Decide whether I can live with this explanation. The answer is no, actually.''

%Another issue that's on my mind is that I've learned my defs. of solid torid orbifold and discal $3$-orbifold are different from the ones in \cite{bmp}. Decide what to do about this. A related point is that the classification of \torifold s can probably be deduced from the proof of ``preoccupani'' (I think it is), but not from the present version of the statement. This is closely related to the issue that will be dealt with in connection with the preoccupaninotes.

}

It will always be understood that an orbifold may have a boundary,
except where we specify otherwise.

Orbifolds will be denoted by capital fraktur letters ($\frakA, \frakB, \frakC,\ldots$).

The underlying space of an orbifold $\oldPsi$ will be denoted by $|\oldPsi|$. If $\oldPsi$ is PL, and if either $\dim\oldPsi\le2$, or $\dim\oldPsi=3$ and $\oldPsi$ is orientable, then $|\oldPsi|$ inherits the structure of a PL manifold of the same dimension as $\oldPsi$. 

By a {\it point} of an orbifold $\oldPsi$ we will mean simply a point
of $|\oldPsi|$. Every open subset of $|\oldPsi|$ is the underlying set
of a unique suborbifold of $\oldPsi$; such a suborbifold will be
referred to as an open subset of $\oldPsi$. We will say that a suborbifold $\frakV$
is a neighborhood of a point $v\in\oldPsi$ if $|\frakV|$ is a neighborhood of
$v$. 

Recall that an orbifold $\oldPsi$ is defined by specifying an underlying space
and a maximal family of mutually compatible chart maps (cf. \cite[Subsection 2.1.1]{bmp}). From our point
of view, a chart map is a map
%For every $v\in\oldPsi$ there exist an open neighborhood $W$ of
%$v$ in $\oldPsi$ and a chart map 
$\phi: U\to W$, where $W$ is an open subset of $|\oldPsi|
$ 
and $U\subset\RR^n$ is either (a) an open ball
about $0$ 
%(if $v\in\inter\oldPsi$) 
or (b) the intersection of an open ball
with a closed half-space $H$ whose boundary contains $0$, 
%(if
%$v\in\partial\oldPsi$), and
%$\phi(0)=v$, 
having the property that for some
% By definition there is a 
finite group
$G\subset \OO(n)$, which leaves $H$ invariant in  case (b)
and is arbitrary in case (a), 
%such that
$\phi$ induces a homeomorphism from $ |U/G|$ (the topological orbit
space defined by the natural action of $G$ on $U$) onto
$|W|$. The group $G$ is uniquely determined by the chart map
$\phi$ (see \cite[p. 443]{bonahon-siebenmann}). For every point $v$ of $\oldPsi$ there is a chart map $\phi$ such
that $\phi(0)=v$. For every $v\in\oldPsi$, either (a) holds for every 
 chart map $\phi$ with $\phi(0)=v$, in which case by definition we
 have $v\in\inter\Psi$, or (b) holds for every 
 chart map $\phi$ with $\phi(0)=v$, in which case by definition we
 have $v\in\partial\Psi$. The group $G$ is then determined  up to conjugacy in
$\OO(n)$ by the point $v$, and will be denoted
$G_v$. We will refer to the order of
$G_v$ as the {\it order} of
$v$. 
%\frakU Q v

For any integer $n\ge1$, we will denote by 
%$\DD^{n}$   the standard
%half-disk
%$\DD^{n}_+=\{(x_1,\ldots,x_n)\in\DD^{n}:x_n\ge0\}$. We will also
%denote by
$U^n=\inter\DD^n$ the open unit disk, and by $U^n_+$ the space
$\{(x_1,\ldots,x_n)\in U^{n}:x_n\ge0\}$. If 
$v$ is a point of an $n$-orbifold $\oldPsi$, we will set
%$\Delta_v=\DD^{n}$ and
 $U_v=U^{n}$ if $v\in\inter\oldPsi$, and
%$\Delta_v=\DD^{n}_+$ and
 $U_v=U^{n}_+$ if $v\in\partial\oldPsi$. Then any chart map
$\phi$ for $\oldPsi$ with $\phi(0)=v$ has as domain the image of $
U_v$ under some similarity transformation $S$ fixing $0$, and may be regarded as the composition of $S^{-1}$
with a chart map whose domain is
$U_v$.

Note that a chart map $\phi$ for an orbifold $\oldPsi$ is by
definition a map from a suitable set
$U\subset\RR^n$ to an open subset $W\subset|\oldPsi|$, and that this
point of view is needed to define an orbifold structure. However, once
the orbifold $\oldPsi$ has been specified, any chart $\phi:U\to W$
defines an orbifold covering map $\alpha:U\to \frakU$, where $\frakU$ is
the unique suborbifold of $\oldPsi$ with $|\frakU|=W$, and this
covering map factors through an orbifold homeomorphism of $U/G$ onto
$\frakU$. An orbifold covering map $\alpha$ which is defined in this
way from a chart map will be called a {\it post-chart map}.

The singular set of an orbifold $\oldPsi$ will be denoted
$\fraks_\oldPsi$. We regard it as a subset of $|\oldPsi|$. It consists
of all points $v\in\oldPsi$ such that $G_v\ne\{1\}$.

If $\oldPsi$ is an $n$-orbifold then $|\oldPsi|$ has a canonical
stratification \cite{kawasaki}, in which $\fraks_\oldPsi$ is the union
of all strata of dimension strictly less than $n$. If $v$ is a point
of $\oldPsi$, then up to conjugacy in $\OO(n)$, the group $G_v$ depends only on the
stratum $\sigma$ containing $v$, and may be denoted $G_\sigma$. The
order of $G_\sigma$ will be called the order of the stratum $\sigma$. If $n=3$ and $\oldPsi$ is orientable, then each component of $\fraks_\oldPsi$ is either a simple closed curve or (the
underlying space of) a graph whose vertices in $|\inter\oldPsi|$ have
valence $3$, and whose intersection with $|\partial\oldPsi|$ consists
of valence-$1$ vertices. In the former case, the given
component of $\fraks_\oldPsi$ is a single stratum. In the latter case,
each vertex lying in $|\inter\oldPsi|$  is a stratum, and the union of
each edge with those of its endpoints that lie in $|\partial\oldPsi|$
is a stratum.

If 
$n=2$ and $\oldPsi$ is orientable, then $G_v$ is cyclic for every $v\in\fraks_\oldPsi$. 
If 
$n=3$ and $\oldPsi$ is orientable, then $G_v$ is cyclic for every point $v$ lying in a one-dimensional stratum of $\oldPsi$.

The distinction between an orbifold $\oldPsi$ and its underlying space $|\oldPsi|$ will be rigidly observed. For example, if $\oldPsi$ is path-connected and $\star\in|\oldPsi|$ is a base point, $\pi_1(\oldPsi,\star)$ will denote the orbifold fundamental group of $\oldPsi$ based at $\star$ (denoted $\pi_1^{\rm orb}(\oldPsi,\star)$ by some authors). In contrast, $\pi_1(|\oldPsi|,\star)$ of course denotes the fundamental group of the underlying space $|\oldPsi|$ based at $\star$. 
(As in the case of spaces, we, will often suppress base points from the notation for the orbifold fundamental group in statements whose truth is independent of the choice of base point.) 
Similarly, $\partial\oldPsi$ will denote the orbifold boundary of $\oldPsi$, which is itself an orbifold.
If $\oldPsi$ is an orbifold of dimension $n\le3$ such that $|\oldPsi|$ is a topological manifold, then $\partial|\oldPsi|=|\partial\oldPsi|\cup\fraks^{n-1}_\oldPsi$, where $\fraks^{n-1}_\oldPsi$ denotes the union of all $(n-1)$-dimensional components of $\fraks_\oldPsi$. In the case of an orientable orbifold $\oldPsi$ of dimension $n\le3$, we have $\fraks^{n-1}_\oldPsi=\emptyset$ and hence $\partial|\oldPsi|=|\partial\oldPsi|$.

The Seifert van Kampen theorem
for orbifolds, which is proved in \cite[Section 2.2]{bmp}, will often
be used without being mentioned explicitly.

An orbifold, not necessarily connected, will be said to be {\it good}
(or, respectively, {\it very good}) if it admits an  (orbifold)
covering space (or, respectively, a  finite-sheeted (orbifold) covering space) which is a manifold. It is a standard consequence of the ``Selberg lemma'' that a compact hyperbolic orbifold is very good, and this fact will often be used implicitly. Note that a suborbifold of a very good orbifold is
  very good. Note also that an orbifold with finitely many components is very
good if and only if its components are very good.

If $v$ is a point of an orbifold $\oldPsi$, and $\phi: U\to\frakU$ is
a post-chart map  with $\phi(0)=v\in\frakU\subset\oldPsi$, then the image of
the inclusion homomorphism $\pi_1(\frakU,v)\to\pi_1(\oldPsi,v)$ is a
subgroup of $\pi_1(\oldPsi,v)$ which depends only on $\oldPsi$ and
$v$, not on the choice of post-chart map. This subgroup is the image of
$G_v$ under a homomorphism which is canonical modulo inner
automorphisms of $G_v$, and is injective if $\oldPsi$ is good. Using
this isomorphism we will often identify $G_v$ with a subgroup of
$\pi_1(\oldPsi,v)$ in the case where $\oldPsi$ is good. It follows that if
$\sigma$ is a stratum of $\oldPsi$, there is a homomorphism from
$G_\sigma$ to $\pi_1(\oldPsi)$ which is canonical modulo inner
automorphisms of $\pi_1(\oldPsi)$, and that if $\oldPsi$ is good then
$G_\sigma$ is identified with a subgroup of $\pi_1(\oldPsi)$ up to conjugacy.
%Q

An orbifold will be said to have {\it finite type} if it is homeomorphic to $\oldUpsilon-\frakE$, where $\oldUpsilon$ is a compact orbifold with $\fraks_\oldUpsilon\subset|\inter\oldUpsilon|$, and $\frakE$ is a union of components of $\partial\oldUpsilon$. Note that in particular,
according to this definition, a finite-type orbifold has compact boundary. 

An orbifold $\oldPsi$ of finite type has a well-defined (orbifold) Euler characteristic, which we
will denote by
%Likewise,
% in the case where $\oldPsi$ is compact \redproofreadingnote{That's not general enough, because I talk about Euler characteristics for finite-type $2$-orbifolds later. One option is to extend the notion of ``finite type,'' which is given elsewhere for $2$-orbifolds, to arbitrary dimensions: 
$\chi(\oldPsi)$. It
% the of $\oldPsi$, 
%when it is defined
is not in general equal to $\chi(|\oldPsi|)$. When $\oldPsi$ has finite type we will also set
$$\chibar(\oldPsi)=-\chi(\oldPsi)$$
in analogy with (\ref{in and out equation}).

Let $\oldPsi_1$ and $\oldPsi_2$ be orbifolds, and set
$n_i=\dim\oldPsi_i$ for $i=1,2$. We define an {\it immersion} (or a
{\it submersion}) from $\oldPsi_1$ to $\oldPsi_2$ to be a map $\frakf$
(of sets) from $|\oldPsi_1|$ to $|\oldPsi_2|$ such that, for every
$x\in\oldPsi_1$,  there exist
 chart maps $\phi_1:U_1\to W_1$ and $\phi_2:U_2\to W_2$ for
 $\oldPsi_1$ and $\oldPsi_2$ respectively, with $x\in W_1$,
%an open neighborhood $W_i$ of $x_i$ in $\oldXi$ which is a suborbifold of $\oldXi$, and a neighborhood $\frakV$ of $\frakf(x)$ in $\oldUpsilon$ which is a suborbifold of $\oldUpsilon$, 
%$\frakV$ of $\oldXi$ and $\oldUpsilon$, which are of $x$ and $y$ respectively, 
%chart maps $\phi_\frakU: U\to\frakU$ and $\phi_\frakV: V\to\frakV$,
%where $ U\subset\RR^m$ and $ V\subset\RR^n$ are open, 
and an injective
(or, respectively, surjective)  affine map $L:\RR^{n_1}\to\RR^{n_2}$,
such that $L(U_1)\subset U_2$
and
$\frakf\circ\phi_1=\phi_2\circ (L|U_1)$. (It follows that
$\frakf(W_1)\subset W_2$, and in particular that
 $\frakf(x)\in W_2$.)
The only ``maps'' between orbifolds that will be considered in this
monograph are immersions and submersions. An orbifold homeomorphism, or
more generally an orbifold covering map, is at once an immersion and a
submersion. An {\it embedding} of $\oldPsi_1$ in 
$\oldPsi_2$ is defined to be a composition
of a homeomorphism of $\oldPsi_1$ onto a suborbifold $\frakZ$ of
$\oldPsi_2$ with the inclusion $\frakZ\to\oldPsi_2$; any embedding is an injective immersion, but the converse is false.
%\psi\phi

Note that an immersion or submersion $\frakf:\oldXi\to\oldUpsilon$ is
in particular a continuous map from $|\oldXi|$ to $|\oldUpsilon|$;
thus if we ignore the orbifold structure, the immersion or submersion
$\frakf$ defines a continuous map of topological spaces, which we
denote by 
 $|\frakf|:|\oldXi|\to|\oldUpsilon|$. 

In this notation, if $\phi$ is a chart map for an orbifold
$\oldPsi$, 
%where $U\subset\RR^n$ and $W$ is an open subset of
%$\oldPsi$, 
and if
% $\frakU$ is the suborbifold of $\oldPsi$ with
%$|\frakU|=W$ and
 $\alpha$ denotes the post-chart map defined by
$\phi$, we have $|\alpha|=\phi$.

%\redcomment{That does not read well. Try to improve it.}
\nonessentialproofreadingnote{If
  I had defined general maps, not just immersions and
  submersions---and it's hard to see why I shouldn't do so---then we
  would have a genuine functor.}
\EndNumber

\begin{notationremarks}\label{obd}
Suppose that $\oldPsi$ is an orbifold, that $X$ is a topological space, and
that $f:X\to|\oldPsi|$ is a map with the property that for some orbifold $\oldXi$ we have $|\oldXi|=X$, and $f=|\frakf|$ for some orbifold immersion $\frakf:\oldXi\to\oldPsi$. In this situation, the orbifold $\oldXi$ is uniquely determined by $X$, $\oldPsi$ and $f$. In situations where it is clear from the context which $\oldPsi$ and $f$ are involved, $\oldXi$ will be denoted by $\obd(X)$. 

The reason why such a suborbifold $\oldXi$ and such a map $\frakf$ exist will usually be immediate from the context, and will be left implicit. 
 The most common situation in which the convention will be used is the
 one in which $X$ is given as a submanifold of $|\oldPsi|$ for some
 orbifold $\oldPsi$, in which case $f$ is understood to be the
 inclusion map $X\to|\oldPsi|$; in this case, $\obd(X)$ is defined if
 and only if $X=|\oldXi|$ for some suborbifold $\oldXi$ of $\oldPsi$,
 and if it is, we have $\obd(X)=\oldXi$. In the case where $\oldPsi$
 is orientable and $m:=\dim\oldPsi\le3$, a sufficient condition for
 $\oldXi$ and $\frakf$ to exist, and hence for $\obd(X)$ to be
 defined, is that $X\subset|\oldPsi|$ be a manifold which has
 dimension strictly less than $m$ and is in general position (see
 \ref{gen pos}) with respect to  $\fraks_\oldPsi$. Note also that
 $\obd(X)$ is always defined if $X$ is an open subset of $|\oldPsi|$;
 and that
 if $X$ is a closed subset of $|\oldPsi|$ such that $\obd(\Fr_{|\oldPsi|}X)$ is defined, then $\obd(X)$ is defined.

%Another situation in which the convention will be used is the one in which $X$ is given as a (not necessarily proper) subset of a covering space $\toldPsi$ of $\oldPsi$. In this case $f$ will be understood to be the restriction of the covering projection to $X$. In this situation, $\obd(X)$ is defined if and only if $X=|\oldXi|$ for some suborbifold $\oldXi$ of $\toldPsi$. \redproofreadingnote{Is this really ever used?  is the following sentence supposed to go?} 

Another situation in which $\oldXi$ and $\frakf$ exist  (and a particularly important one) is the one in
which $\oldTheta$ is a closed, two-sided (see \ref{two sides}) $2$-suborbifold of the
interior of an orientable PL $3$-orbifold $\oldOmega$, so that
$|\oldTheta|$ is a two-sided PL $2$-submanifold of $\inter|\oldOmega|$,
and we take $\oldPsi=\oldOmega$, $X=|\oldOmega|\cut{|\oldTheta|}$, and
$f=\rho_{|\oldTheta|}$. In this case 
$\obd(X)$ and $\frakf$ are always defined.
The PL orbifold
$\obd(X)=\obd(|\oldOmega|\cut{|\oldTheta|})$ will be denoted by
$\oldOmega\cut\oldTheta$, and the PL immersion $\frakf$ such that
$|\frakf|=f=\rho_{|\oldTheta|}$ will be denoted by
$\rho_\oldTheta$. \end{notationremarks}

\Number\label{unions and such}
%\Number\label{unions and such}
 If $\oldXi$ and $\oldUpsilon$ are
suborbifolds of an orbifold $\oldPsi$ we will set
$\oldXi\cup\oldUpsilon=\obd(|\oldXi|\cup|\oldUpsilon|)$, provided that
$|\oldXi|\cup|\oldUpsilon|$ is the underlying set of some suborbifold of
$\oldPsi$; more generally, if $(\oldXi_i)_{i\in I}$ is a
family of suborbifolds of $\oldPsi$ indexed by some set $I$, we will
set $\bigcup_{i\in I}\oldXi_i=\obd(\bigcup_{i\in
  I}|\oldXi_i|)$ provided that the right hand side is
defined. Likewise, for orbifolds $\oldXi$ and $\oldUpsilon$ of $\oldPsi$ we
will set $\oldXi\cap\oldUpsilon=\obd(|\oldXi|\cap|\oldUpsilon|)$ and
$\oldXi\setminus\oldUpsilon=\obd(|\oldXi|\setminus|\oldUpsilon|)$, provided that
the respective right hand sides are defined, and when
$\oldUpsilon\subset\oldXi$ and $\oldXi\setminus\oldUpsilon$ is defined we will
write $\oldXi-\oldUpsilon=\oldXi\setminus\oldUpsilon$. If
$\frakf:\oldXi\to\oldUpsilon$ is an immersion or submersion of
orbifolds, and $\oldXi$ is a suborbifold of $\oldXi$ (or
$\oldUpsilon$) such that 
$|\frakf|(|\oldXi|)$ (or, respectively,
$|\frakf|^{-1} (|\oldXi|)$) is the underlying subspace of a suborbifold of
$\oldUpsilon$ (or, respectively. $\oldXi$), we will set
$\frakf(\oldXi)=\obd(|\frakf|(|\oldXi|))$ (or, respectively,
$\frakf^{-1}(\oldXi)=\obd(|\frakf|^{-1}(|\oldXi|))$.
\nonessentialproofreadingnote{I've done a good deal more back-and-forth between suborbifolds and submanifolds than is really necessary, and I won't be able to eliminate this completely but should try to minimize it.}
 \EndNumber

\Number\label{isotopies}
By an {\it isotopy} of an orbifold $\oldPsi$ we mean a family 
$(\frakh_t)_{0\le t\le1}$ of self-homeomorphisms of $\oldPsi$ such that
(i) the map $(x,t)\mapsto\frakh_t(x)$ from $|\oldPsi\times[0,1]|$
to $|\oldPsi|$ is a submersion 
from $\oldPsi\times[0,1]$ to $\oldPsi$, and (ii)
$\frakh_0$
is the identity. We will say that $(\frakh_t)_{0\le t\le1}$ is {\it
  constant} on a suborbifold $\oldUpsilon$ of $\oldPsi $ if
$\frakh_t|\oldUpsilon$ is the identity for each $t$. We will use the orbifold analogues of standard
language for isotopy of manifolds; for example, two suborbifolds
$\oldXi,\oldXi'$ of $\oldPsi$ will be said to be (ambiently) isotopic, rel a suborbifold $\oldUpsilon$ of $\oldPsi$, if there is an isotopy 
$(\frakh)_{0\le t\le1}$ of $\oldPsi$, constant on $\oldUpsilon$, such that
$\frakh_1(\oldXi)=\oldXi'$. In the case where $\oldUpsilon=\emptyset$ we will simply say that $\oldXi$ and $\oldXi'$ are isotopic. We will say that $\oldXi$ and $\oldXi'$ are
{\it non-ambiently isotopic} if they are isotopic when regarded as
suborbifolds of the orbifold $\oldPsi^*\supset\oldPsi$ obtained from
the disjoint union of $\oldPsi$ with $(\partial\oldPsi)\times[0,1]$ by
gluing the suborbifolds $\partial\oldPsi\subset\oldPsi$ and
$(\partial\oldPsi)\times\{0\}\subset (\partial\oldPsi)\times[0,1]$ via
the homeomorphism $x\to(x,0)$.
\EndNumber

\Number\label{pairs}
A (topological, PL or smooth) {\it orbifold pair} is an ordered pair $(\oldUpsilon,\oldXi)$, where $\oldUpsilon$ is a
(topological, PL or smooth)  orbifold and $\oldXi$ is a (topological,
PL or smooth) suborbifold of $\oldUpsilon$. By a {\it smooth (or PL)
  structure} on a topological orbifold pair $(\oldUpsilon,\oldXi)$ we
mean a smooth (or PL)
  structure on $\oldUpsilon$ for which $\oldXi$ is a smooth (or PL)
suborbifold. A smooth structure and a PL structure on  $(\oldUpsilon,\oldXi)$ will be said to be {\it
  compatible} if they are compatible as a smooth structure and a PL
structure on  $\oldUpsilon$.

An {\it immersion} (or
{\it embedding} or {\it submersion})
% (or PL homeomorphism or diffeomorphism) from a topological (or PL or smooth) 
of an orbifold pair  $(\oldUpsilon,\oldXi)$ in another such pair
 $(\oldUpsilon',\oldXi')$ is an immersion (or, respectively, an embedding or submersion) of $j:\oldUpsilon\to
 \oldUpsilon'$ such that $j(\oldXi)\subset \oldXi'$, and such that
 $j|\oldXi\to\oldXi'$ is an immersion (or, respectively, an embedding or submersion). 
A {\it homeomorphism} of orbifold pairs (called a diffeomorphism in
the smooth category) is an immersion of pairs having an inverse which
is an immersion of pairs; equivalent, a homeomorphism from a  pair  $(\oldUpsilon,\oldXi)$ to a pair
 $(\oldUpsilon',\oldXi')$ is a homeomorphism from $\oldUpsilon$ to
 $\oldUpsilon'$ that carries $\oldXi$ onto $\oldXi'$. An {\it isotopy}
 of a pair $(\oldUpsilon,\oldXi)$ is an isotopy $(h_t)_{0\le t\le1}$
 of $\oldUpsilon$ such that $h_t(\oldXi)=\oldXi$ for every $t\in[0,1]$.

Two smooth (or
 PL) structures on a topological orbifold pair $(\oldUpsilon,\oldXi)$
 are {\it equivalent} if there is a self-homeomorphism of
 $(\oldUpsilon,\oldXi)$ carrying one structure to the other.

At one point we will also need to consider  PL  {\it orbifold triads}.
Such a triad is by definition an ordered triple $(\oldUpsilon,\oldXi_1,\oldXi_2)$, where $\oldUpsilon$ is a
 PL  orbifold and $\oldXi_i$ is a  PL suborbifold of $\oldUpsilon$ for
 $i-1,2$. A PL {\it  homeomorphism} between PL 
 orbifold triads  $(\oldUpsilon,\oldXi_1,\oldXi_2)$ and
 $(\oldUpsilon',\oldXi'_1,\oldXi'_2)$  is defined to be a
 homeomorphism from $\oldUpsilon$ to $\oldUpsilon'$ that carries
 $\oldXi_i$ onto $\oldXi'_i$ for $i=1,2$.
\EndNumber
%U\tfrakU V \tfrakV \frakZ (a) W L \alpha\oldLambda S

\Number\label{two sides}
Let $\oldOmega$ be a (topological, PL or smooth) compact orbifold. Generalizing the definition given for manifolds in \ref{nbhd stuff}, we define a (respectively topological, PL or smooth) suborbifold $\oldTheta$ of $M$ to be  {\it two-sided} if it has a neighborhood $\frakH$ in $\oldOmega$ (again called a {\it collar neighborhood}) such that the pair $(\frakH,\oldTheta)$ is (topologically, piecewise-linearly or smoothly) homeomorphic to $(\oldTheta\times[-1,1],\oldTheta\times\{0\})$. Note that this implies that $\oldTheta$ is properly embedded in $\oldOmega$ (in the sense that $\oldTheta\cap\partial\oldOmega=\partial\oldTheta$) and has codimension $1$.
\EndNumber

\Notation\label{wuzza weight}
If $\oldPsi$ is an orbifold, and $S$ is a subset of $|\oldPsi|$ such that $S\cap\fraks_\oldPsi$ is finite, the number $\card (S\cap\fraks_\oldPsi)$ will be called the {\it weight} of $S$ in $\oldPsi$, and will be denoted $\wt_\oldPsi S$, or simply by $\wt S$ when $\oldPsi$ is understood.
\EndNotation

\Number\label{gen pos}
Let $\oldPsi$ be an orbifold of dimension  $m\le3$, and let $\fraks^{(0)}_\oldPsi$ denote the union of all zero-dimensional strata of $\fraks_\oldPsi$ (see \ref{orbifolds introduced}). Then $|\oldPsi|-\fraks^{(0)}_\oldPsi$ is an $m$-manifold, and every positive-dimensional stratum of $\oldPsi$ is a submanifold of $|\oldPsi|-\fraks^{(0)}_\oldPsi$. 
A manifold $H\subset|\oldPsi|$, having   dimension strictly less than $m$, will be said to be {\it in general position with respect to  $\fraks_\oldPsi$} if $H$ is disjoint from $\fraks^{(0)}_\oldPsi$, and intersects every positive-dimensional stratum of $\fraks_\oldPsi$ transversally in the manifold $|\oldPsi|-\fraks^{(0)}_\oldPsi$. 
\EndNumber

\Number\label{new weak}
Let $X$ be a compact  PL subset of $|\oldPsi|$, where $\oldPsi$ is an
orientable (PL) orbifold of dimension at most $3$. We define a {\it
  strong regular neighborhood} of $X$ in $\oldPsi$ to be a suborbifold
$\frakC$ of $\oldPsi$ with the property that $|\frakC|$ is the second-derived
neighborhood of $X$ with respect to some triangulation $\calt$ of
$|\oldPsi|$, compatible with its PL structure, such that both $X$ and
$\fraks_\oldPsi$ are underlying sets of subcomplexes of $\calt$. This
definition implies that $|\frakC|$ is unique up to a PL isotopy
$(h_t)_{0\le t\le 1}$ of the manifold $|\oldPsi|$, constant on $X$,
such that $h_t(\fraks_\oldPsi)=\fraks_\oldPsi$ for each $t$. Hence the
regular neighborhood $\frakC$ of $X$ is unique up to PL orbifold
isotopy.

We define a {\it weak regular neighborhood} of $X$ in $\oldPsi$ to be a suborbifold $\frakC$ of $\oldPsi$ such that (i) $\frakC$ is a neighborhood of $X$ in $\oldPsi$, and (ii) if $\oldPsi^*$ is the $3$-orbifold obtained from the disjoint union $\oldPsi\discup((\partial\oldPsi)\times[0,1])$ by gluing $\partial\oldPsi\subset\oldPsi$ to $(\partial\oldPsi)\times\{0\}\subset(\partial\oldPsi)\times[0,1]$ via the homeomorphism $t\mapsto(t,0)$, then a strong regular neighborhood of $\frakC\subset\oldPsi\subset\oldPsi^*$ in $\oldPsi^*$ is also a strong regular neighborhood of $X$ in $\oldPsi^*$.
When $\partial\oldPsi=\emptyset$, the notions of strong and weak regular neighborhoods coincide, and we will simply use the term {\it regular neighborhood} in that case.

If $\oldXi$ is a PL suborbifold of $\oldPsi$, we define a strong (or
weak) regular neighborhood of $\oldXi$ in $\oldPsi$ to be a strong
(or, respectively weak) regular neighborhood of $|\oldXi|$.

\EndNumber
%\frakC\frakC\oldUpsilon

\Number
For our purposes, if $\oldPsi$ is an orbifold, $H_1(\oldPsi)$ will denote the abelianization of the orbifold fundamental group $\pi_1(\oldPsi)$, and for any abelian group $A$ we will define $H_1(\oldPsi,A)$ to be $H_1(\oldPsi)\otimes A$. We will not define (or use) higher-dimensional homology for orbifolds.
\EndNumber

\DefinitionsRemarks\label{annular etc}
As is standard,  an $n$-orbifold will be called {\it discal} if it is (orbifold-)homeomorphic to the quotient of $\DD^n$ by an orthogonal action of a finite group. It follows from the Orbifold Theorem \cite{blp}, \cite{chk} that a $3$-orbifold $\oldPsi$, such that $\fraks_\oldPsi$ has no isolated points, is discal if and only if $\oldPsi$ it admits $\DD^3$ as a(n orbifold) 
covering by an $n$-ball; this fact will often be used without an explicit reference. A $2$-orbifold $\oldTheta$ will be termed {\it
  spherical} if $\oldTheta$ is closed and connected and $\chi(\oldTheta)>0$; if $\oldTheta$ is very good, this is equivalent to the condition that $\oldTheta$ admits $\SSS^2$ as a covering.
A $2$-orbifold $\oldTheta$ will be termed {\it annular} or {\it toric} if $\oldTheta$ admits $\SSS^1\times[0,1]$ or $\TTT^2$, respectively, as a covering; this is equivalent to the condition that $\oldTheta$ is connected, that $\chi(\oldTheta)=0$,  and that $\partial\oldTheta$ is, respectively, non-empty or empty. 
A {\it \torifold} is defined to be a $3$-orbifold which is
covered by a solid torus. (We will give two other characterizations of  \torifold s in \ref{three-way equivalence} below.)
\EndDefinitionsRemarks

\Number\label{other injectify}We will use the analogues for orbifolds
of the conventions of \ref{injectify}. A connected suborbifold
$\frakU$ of a connected orbifold $\frakZ$ will be
termed {\it $\pi_1$-injective} if the inclusion homomorphism $\pi_1(\frakU)\to
\pi_1(\frakZ)$ is injective. 
In general, a suborbifold $\frakU$ of an orbifold $\frakZ$ will be
termed $\pi_1$-injective if each component of $\frakZ$ is
$\pi_1$-injective in the component of $\frakU$ containing
it. A closed (possibly disconnected) $2$-dimensional suborbifold $\frakU$ of a
$3$-orbifold $\frakZ$ will be termed {\it incompressible} if $\frakU$
is  contained in $\inter\frakZ$, is two-sided, is $\pi_1$-injective in $\frakZ$, and has no component which bounds a discal $3$-suborbifold of $\frakZ$. 
\EndNumber

\Definition\label{praxis}
Let  $\otheroldLambda$ be an orientable $2$-orbifold without boundary.
A $2$-dimensional suborbifold $\oldUpsilon$ of $\otheroldLambda$ will be
termed {\it taut} if (i) $\oldUpsilon$ is a closed subset of
$\otheroldLambda$ having compact boundary, (ii) $\partial\oldUpsilon$ is $\pi_1$-injective in
$\otheroldLambda$, and (iii) no component of
$\partial\oldUpsilon$ is contained in a suborbifold of $\otheroldLambda$ which is a $2$-manifold homeomorphic to
$\SSS^1\times[0,\infty)$.  

Note that in the case where $\otheroldLambda$ is a closed $2$-orbifold, Condition (i) is equivalent to saying that $\oldUpsilon$ is compact, and Condition (iii) holds automatically. Thus in the case where $\otheroldLambda$ is closed, a taut suborbifold of $\oldLambda$ is simply a compact suborbifold $\oldUpsilon$ such that $\partial\oldUpsilon$ is $\pi_1$-injective in
$\otheroldLambda$. Many of the assertions in the monograph in which tautness plays a role
% statements in this section, and all the applications in subsequent sections of the results of this section.
will involve a closed $2$-orbifold. However, the more general definition will come up in Section \ref{vegematic section}.

Note that the orientability of $\otheroldLambda$ implies that for any suborbifold $\oldUpsilon$, the subspace $|\oldUpsilon|$ is a submanifold of $|\otheroldLambda|$ whose boundary is disjoint from $\fraks_\otheroldLambda$. In particular, $\partial\oldUpsilon$ is a $1$-manifold, which is compact and $\pi_1$-injective in $\otheroldLambda$ if $\oldUpsilon$ is taut.
\EndDefinition

\Definition\label{just irreducible}
A $3$-orbifold $\oldPsi$ is said to be {\it irreducible} if $\oldPsi$
is connected and every two-sided spherical $2$-suborbifold of
$\oldPsi$ is the boundary of a discal $3$-suborbifold of
$\oldPsi$. As mentioned in Subsection \ref{great day}, this
generalizes our definition for the case of a $3$-manifold.
\EndDefinition

\Number\label{doubling}
Let $\oldPsi$ be a topological orbifold, and let $\oldXi$ be a
codimension-$0$ topological
suborbifold of $\partial\oldPsi$ (so that $\oldXi$ has codimension $1$ in $\oldPsi$).
The {\it double} of $\oldPsi$ along $\oldXi$, denoted ${\rm D}_\oldXi\oldPsi$, is the
closed orbifold obtained from the disjoint union of two copies of
$\oldPsi$ by gluing  together the copies of $\oldXi$ contained in
their boundaries via the identity homeomorphism. The definition of
${\rm D}_\oldXi\oldPsi$ gives two canonical embeddings of $\oldPsi$ in ${\rm D}_\oldXi\oldPsi$.
If $\oldPsi$ has a PL structure, and the suborbifold $\oldXi$ is PL,
then ${\rm D}_\oldXi\oldPsi$ inherits a PL structure, and the
canonical embeddings are PL.

If $\oldPsi$ is a topological or PL orbifold, we will write
${\rm D}\oldPsi$ for ${\rm D}_{\partial\oldPsi}\oldPsi$, which is a closed
topological or PL
orbifold.

Note that this definition of the double of a topological or PL
orbifold includes the case where
$\oldPsi$ is a topological or PL manifold;  in this case ${\rm D}\oldPsi$ is
a closed topological or PL
manifold.

\EndNumber

\Number\label{smooth doubling}
Defining doubling in the smooth category is a little trickier. The
only case that will be needed in this monograph  is the one in which 
$\oldPsi$ is a smooth, orientable orbifold of dimension at most $3$, and  $\oldXi$ is a
codimension-$0$ smooth
suborbifold of $\partial\oldPsi$. In this case we automatically have
$|\partial\oldXi|\cap\fraks_\oldPsi=\emptyset$. Hence for some
manifold neighborhood $U$ of $|\partial\oldXi|$ in $|\oldPsi|$, we have 
$U\cap\fraks_\oldPsi=\emptyset$, and so the manifold $U$ is a
suborbifold of  $\oldPsi$. Let $U'$ denote the smooth manifold-with-corner
obtained from $U$ by
%We may give $U$
%the structure of a smooth manifold by
introducing a corner along
 $|\partial\oldXi|$ (see \cite[Theorem 7.5.6]{mukherjee}); thus $U$ and $U'$ have the same underlying set.  We set $Y=U'\cap|\oldXi|$. The smooth manifold-with-corner structure of
 $U'$ and the smooth orbifold structure of $\oldPsi-\partial\oldXi$ restrict to the
 same smooth manifold structure on $U'-|\partial\oldXi|$. The smooth
 manifold-with-corner structure on $U'$ gives rise to the structure
 of a smooth manifold on ${\rm D}_Y U'$, and the smooth orbifold
 structure of $\oldPsi-\partial\oldXi$ gives rise to the structure of
a smooth orbifold on ${\rm D}_{\inter\oldXi}(\oldPsi-\partial\oldXi)$. We
may regard the topological orbifold
${\rm D}_\oldXi\oldPsi$ as the union of
${\rm D}_{\inter\oldXi}(\oldPsi-\partial\oldXi)$ with 
${\rm D}_{Y} U'$. The smooth orbifold
 structure of ${\rm D}_{\inter\oldXi}(\oldPsi-\partial\oldXi)$ and
 the smooth manifold
 structure of ${\rm D}_{Y} U'$ restrict to the same smooth manifold structure on their
intersection, ${\rm D}_{Y-\partial|\oldXi|}(U'-|\partial\oldXi|)$. Hence
these smooth structures define a smooth structure on
${\rm D}_\oldXi\oldPsi$. The canonical embeddings of
$\oldPsi$ in ${\rm D}_\oldXi\oldPsi$ are only piecewise smooth, but either of them restricts to a
smooth embedding of $\oldXi$ in ${\rm D}_\oldXi\oldPsi$.
\EndNumber

\Number\label{silvering}
Note that if $\oldPsi$ is a topological or PL orbifold, or a smooth
orientable orbifold of dimension at most $3$,  and $\oldXi$ is
respectively a topological, PL or smooth codimension-$0$
suborbifold of $\partial\oldPsi$, 
the orbifold ${\rm D}_\oldXi\oldPsi$ has a canonical (topological, PL or smooth) involution which maps each copy of 
$\oldPsi$ onto the other, via the identity homeomorphism. The quotient of
${\rm D}_\oldXi\oldPsi$ by this involution is a (topological, PL or smooth) orbifold, which
%$\oldPsi'$ such that
%$|\oldPsi'|=|\oldPsi|$ and
%$\fraks_{\oldPsi'}=\fraks_\oldPsi\cup|\oldXi|\subset|\oldPsi|$. The orbifold $\oldPsi'$
is said to be obtained from $\oldPsi$ by {\it silvering} the suborbifold $\oldXi$ of $\partial\oldPsi$. It will be denoted
$\silv_\oldXi\oldPsi$, or by $\silv\oldPsi$ in the special case where
$\oldXi=\partial\oldPsi$. (Cf. \cite[p. 65, Remark]{other-bs}.) Note that
$\silv_\oldXi\oldPsi$ is always non-orientable if
$\oldXi\ne\emptyset$. There is a canonical immersion of
$\oldPsi$ in $\silv_\oldXi\oldPsi$, 
% denote the natural
which may be defined as the composition of either of the canonical embeddings $\oldPsi\to
{\rm D}_\oldXi\oldPsi $ with the natural degree-$2$ covering map
$p:{\rm D}_\oldXi\oldPsi \to\silv_\oldXi\oldPsi$. This canonical immersion
will be denoted $\mu_{\oldPsi,\oldXi}$, or simply by $\mu_\oldPsi$ in
the case where $\oldXi=\partial\oldPsi$. If the pair
$(\oldPsi,\oldXi)$ is PL, then $\mu_{\oldPsi,\oldXi}$ is a PL
immersion. If 
$(\oldPsi,\oldXi)$ is smooth, then $\mu_{\oldPsi,\oldXi}$ is only a
piecewise smooth immersion, since the canonical embeddings of
$\oldPsi$ in
${\rm D}_\oldXi\oldPsi$ are only piecewise smooth; however,
$\mu_{\oldPsi,\oldXi}$ restricts to a smooth embedding of $\oldPsi$ in $\silv_{\oldXi}\oldPsi$.

Note that
$|\mu_{\oldPsi,\oldXi}|:|\oldPsi|\to|\silv_\oldXi\oldPsi|$ is a
homeomorphism of topological spaces, and that
$|\mu_{\oldPsi,\oldXi}|^{-1}(\fraks_{\silv_\oldXi\oldPsi})=\fraks_\oldPsi\cup|\oldXi|\subset|\oldPsi|$. Note
also that, in any of the three categories, $\mu_{\oldPsi,\oldXi}|(\oldPsi-\oldXi):
\oldPsi-\oldXi\to (\silv_\oldXi\oldPsi)-\mu_{\oldPsi,\oldXi} (\oldXi)$ is an
orbifold homeomorphism; this homeomorphism will be denoted $\mu^*_{\oldPsi,\oldXi}$.

As the definitions of doubling and silvering in the smooth category
apply only to orientable orbifolds of dimension at most $3$, it will be
understood that the smooth versions of Subsections \ref{first point}
and \ref{third point}, of Lemma \ref{shabbos}, and of Proposition
\ref{oh yeah he tweets} apply only when the orbifolds in question are
orientable and have dimension at most $3$. (In the case of 
Proposition
\ref{oh yeah he tweets} this is  only slightly stronger than the general hypotheses.)

 We define $[[0,1]$ to be the $1$-orbifold $\silv_{\{0\}}[0,1]$, and we define $[[0,1]]$ to be the $1$-orbifold $\silv_{\{0,1\}}[0,1]$.
\EndNumber
%\eta\alpha

\Number\label{gyrot}
We have observed that  if $(\oldPsi,\oldXi)$ is a topological, PL, or
smooth orbifold pair such that $\oldPsi$ is orientable, $\oldXi\subset\partial\oldPsi$, and
$\dim\oldXi=\dim\partial\oldPsi\le2$, then 
$\mu_{\oldPsi,\oldXi}|\oldXi:\oldXi\to \silv_\oldXi\oldPsi$ is
respectively a topological, PL, or smooth embedding. In particular,
$\mu_{\oldPsi,\oldXi}(\oldXi)$ is, respectively, a topological, PL, or
smooth suborbifold of
${\rm D}_\oldXi\oldPsi$. 
%Hence 
%$\mu_{\oldPsi,\oldXi}(
%\oldXi)$ 
%under the canonical
%immersion $\oldPsi\to\silv_\oldXi\oldPsi$ 
%is a suborbifold of
%$\silv_\oldXi\oldPsi$. 
Note also that the orientability of $\oldPsi$
implies that every singular stratum of $\oldPsi$ has codimension at
least $2$ in $|\oldPsi|$. Hence $|\mu_{\oldPsi,\oldXi}(
\oldXi)|$ is the union of the closures of the codimension-$1$ strata of
$\silv_\oldXi\oldPsi$.
\EndNumber

\Number\label{first point}
Let $\oldPsi$ be an orbifold,  let $\oldXi$ be a
codimension-$0$ 
suborbifold of $\partial\oldPsi$, and let $\oldUpsilon$ be a
suborbifold of $\oldPsi$ such that $\oldXi\cap\oldUpsilon$ is a
codimension-$0$ 
suborbifold of $\oldXi$. Then ${\rm D}_{\oldXi\cap\oldUpsilon}\oldUpsilon$ is
canonically identified with a suborbifold of ${\rm
  D}_\oldXi\oldPsi$. This is immediate in the topological and PL
categories. In the
  smooth category, identifying ${\rm D}_{\oldXi\cap\oldUpsilon}\oldUpsilon$  with a suborbifold of ${\rm
  D}_\oldXi\oldPsi$ requires carrying out the construction described
in \ref{smooth doubling} for both  $\oldPsi$ and its suborbifold
$\oldUpsilon$. Thus if $U$ is a manifold neighborhood  of
$|\partial\oldXi|$ in $|\oldPsi|$, we must introduce
%   \redcomment{Explanation that we must introduce 
corners both
  for $U$, along $|\partial\oldXi|$, and for $U\cap|\oldUpsilon|$,
  along $|\partial(\oldUpsilon\cap\oldXi)|$, in such a way that
  $U\cap|\oldUpsilon|$ becomes a submanifold-with-corners of $U$. 
% and for a manifold neighborhood of
  %$|\partial(\oldUpsilon\cap\oldXi)|$ in $|\oldUpsilon|$. 
This can be done even if $\oldUpsilon$ has non-empty intersection with $\partial\oldXi$,
  because $\partial\oldUpsilon$ and $\partial\oldPsi$ must be tangent
  along their intersection since $\oldUpsilon$ is a smooth suborbifold
  of $\oldPsi$.

In any of the three categories, when we
pass to the quotient by the canonical involution of ${\rm D}_\oldXi\oldPsi$, the identification of ${\rm D}_{\oldXi\cap\oldUpsilon}\oldUpsilon$ with a suborbifold of ${\rm
  D}_\oldXi\oldPsi$ gives rise to a
canonical identification of $\silv_{\oldXi\cap\oldUpsilon}\oldUpsilon$ with a
suborbifold of $\silv_\oldXi\oldPsi$. Under this identification we
have
$\silv_{\oldXi\cap\oldUpsilon}\oldUpsilon=\mu_{\oldPsi,\oldXi}(\oldUpsilon)$.
\EndNumber
%\oldXi_0

\Number\label{second point}
Let $\oldPsi$ be an orientable orbifold of dimension at most $3$, and let $\oldXi$ be a
codimension-$0$ 
suborbifold of $\partial\oldPsi$. Then for every suborbifold $\frakZ$ of
$\silv_\oldXi\oldPsi$,  no component of which is contained in
$\mu_{\oldPsi,\oldXi}(\oldXi)$, there is a
suborbifold
$\oldUpsilon$  of $\oldPsi$ such that $\oldXi\cap\oldUpsilon$ is a
codimension-$0$ 
suborbifold of $\oldXi$, and such that under the
canonical identification mentioned in \ref{first point}, we have 
$\frakZ=\silv_{\oldXi\cap\oldUpsilon}\oldUpsilon$. To prove this, note
that since the quotient map $p:{\rm D}_\oldXi\oldPsi\to\silv_\oldXi\oldPsi$
is a covering, $\tfrakZ:=p^{-1}(\frakZ)$ is a suborbifold of ${\rm D}_\oldXi\oldPsi$,
invariant under the canonical involution $\delta$ of
${\rm D}_\oldXi\oldPsi$. Thus $\delta|\tfrakZ$ is an  involution of the
orbifold $\tfrakZ$, no component of $\tfrakZ$ is pointwise fixed by
$\delta$, and no component of the complement of
$\Fix(\delta|\tfrakZ)=p^{-1}(\oldXi\cap\frakZ)$ relative to $\tfrakZ$
is $\delta$-invariant;  since $\dim\tfrakZ\le\dim {\rm D}_\oldXi\oldPsi\le3$, it
follows that $p^{-1}(\oldXi\cap\frakZ)$ is a
two-sided suborbifold of $\tfrakZ$. This implies that if $\oldPsi_1$
is one of the two canonical copies of $\oldPsi$ contained in
${\rm D}_\oldXi\oldPsi$,  then $\oldUpsilon_1:=\oldPsi_1\cap\tfrakZ$ is a
suborbifold of $\oldPsi_1$, and $\toldXi:=p^{-1}(\oldXi\cap \frakZ)$ is a
codimension-$0$ suborbifold of $\partial\oldUpsilon_1$. The orbifold
$\oldPsi_1$ comes equipped with a canonical homeomorphism
$\iota:\oldPsi\to\oldPsi_1$ such that
$p\circ\iota=\mu_{\oldPsi,\oldXi}$. It now follows that
%$p|\oldPsi_1:\oldPsi_1\to\oldPsi$ is a homeomorphism, it follows that 
$\oldUpsilon:=\iota(\oldUpsilon_1)$ is a
suborbifold of $\oldPsi$, and that $\iota(\toldXi)=\oldXi\cap\oldUpsilon$ is a
codimension-$0$ suborbifold of $\partial\oldUpsilon$. The construction
implies that $\mu_{\oldPsi,\oldXi}(\oldUpsilon)=\frakZ $, so that by
\ref{first point} the suborbifold of $\silv_\oldXi\oldPi$
canonically identified with $\silv_{\oldXi\cap\oldUpsilon}\oldUpsilon$
is $\frakZ$.
\EndNumber

\Number\label{third point}
Let $\oldPsi_1$ and $\oldPsi_2$ be orbifolds,  let $\oldXi_i$
be a
codimension-$0$ 
suborbifold of $\partial\oldPsi_i$ for $i=1,2$, and let
$f:(\oldPsi_1, \oldXi_1)\to(\oldPsi_2, \oldXi_2)$ be an immersion (or
a submersion) of pairs. For $i=1,2$, the orbifold ${\rm D}_{\oldXi_i}\oldPsi_i$
is by definition the union of two copies $\oldPsi_i^{(1)}$ and
$\oldPsi_i^{(2)}$ of $\oldPsi_i$, glued along $\oldXi_i$ via the
identity. We denote by ${\rm D}f:{\rm D}_{\oldXi_1}\oldPsi_i\to
{\rm D}_{\oldXi_1}\oldPsi_2$ the immersion (or, respectively, submersion)
defined to map $\oldPsi_1^{(j)}$ to $\oldPsi_2^{(j)}$, via the map
$f$, for $j=1,2$. If $\delta_i$ denotes the canonical involution of
${\rm D}_{\oldXi_i}\oldPsi_i$ for $i=1,2$, we have $\delta_2\circ
{\rm D}f={\rm D}f\circ\delta_1$, and hence ${\rm D}f$ induces an immersion (or
respectively a submersion) $\silv
f:\silv_{\oldXi_1}\oldPsi_1\to\silv_{\oldXi_2}\oldPsi_2$. Note that,
regarding $f$ as a submersion (or immersion) of $\oldPsi_1$ in $\oldPsi_2$, we have
\Equation\label{functor freak}
(\silv
f)\circ\mu_{\oldPsi_1,\oldXi_1}=\mu_{\oldPsi_2,\oldXi_2}\circ
f.
\EndEquation
 Furthermore, since $|\mu_{\oldPsi_i,\oldXi_i}|:
|\oldPsi_i|\to |\silv_{\oldXi_i}\oldPsi_i|$ is a homeomorphism for
$i=1,2$ (see \ref{silvering}), $\silv f$ is the unique immersion (or
respectively, submersion) for which (\ref{functor freak}) holds.

If the pairs $(\oldPsi_i,\oldXi_i)$ (for $i=1,2$) are smooth, and $f$
is a smooth immersion, then ${\rm D}f$ is a smooth immersion. (Here it
is understood that the doubling is carried out in the smooth category,
as in \ref{smooth
  doubling}.) As the canonical covering maps 
${\rm D}_{\oldXi_i}\oldPsi_i\to\silv_{\oldXi_1}\oldPsi_i$ are smooth,
it follows that $\silv f$ is a smooth immersion. This fact is {\it not}
readily deduced from (\ref{functor freak}), as the immersions
$\mu_{\oldPsi_i,\oldXi_i}$ are not smooth.

If $f:(\oldPsi_1, \oldXi_1)\to(\oldPsi_2, \oldXi_2)$ is a
homeomorphism of pairs, then $\silv f$ is a homeomorphism, as it
admits the inverse $\silv(f^{-1})$. More
generally, we claim that if $f$ is an embedding of pairs, and if $f^{-1}(\oldXi_2)=\oldXi_1$, then $\silv f$ is an
embedding. 
To prove this, note
that since $f$ is an embedding we may write $f=i\circ h$, where $h$ is
a homeomorphism of $\oldPsi_1$ onto a suborbifold $\oldPsi_1'$ of
$\oldPsi_2$, and $i:\oldPsi_1'\to\oldPsi_2$ is the inclusion. Set
$\oldXi_1'=h(\oldXi_1)$. Since $f^{-1}(\oldXi_2)=\oldXi_1$, we have
$\oldPsi_1'\cap\oldXi_2=\oldXi_1'$. Hence by \ref{first
    point}, $\silv_{\oldXi_1'}\oldPsi_1'$ is identified with a
  suborbifold of $\silv_{\oldXi_2}\oldPsi_2$; the map $\silv i:
  \silv_{\oldXi_1'}\oldPsi_1'\to\silv_{\oldXi_2}\oldPsi_2$ is then the
  inclusion. Since $h$ is a homeomorphism and $h(\oldXi_1)=\oldXi_1'$,
  the map $\silv
  h:\silv_{\oldXi_1}\oldPsi_1\to\silv_{\oldXi_1'}\oldPsi_1'$ is a homeomorphism.
But since $f=i\circ h$, we have
 $\silv f=(\silv i)\circ(\silv h)$; this exhibits $\silv f$ as the
composition of a homeomorphism with an inclusion map, and shows that
$f$ is indeed an embedding.
\EndNumber

Proposition \ref{oh yeah he tweets} below provides a converse to the
observations made in \ref{third point}. The following lemma is needed
for the proof of Proposition \ref{oh yeah he tweets}.

\Lemma\label{shabbos}
Let $\oldPsi_1$ and $\oldPsi_2$ be orientable $n$-orbifolds, where $n$
is an integer with $0<n\le3$. Let 
$f:\inter\oldPsi_1\to\inter\oldPsi_2$ be a homeomorphism, and suppose
that there is
 a homeomorphism
$h:|\oldPsi_1|\to|\oldPsi_2|$ such that
$|f|=h\big|\inter|\oldPsi_1|$. Then $f$ extends to a homeomorphism
$\overline{f}:\oldPsi_1\to\oldPsi_2$, and $|\overline{f}|=h$. 
\EndLemma

\Proof
We first show that there is a homeomorphism
$\overline{f}:\oldPsi_1\to\oldPsi_2$ such that $|\overline{f}|$ is
equal to the given homeomorphism $h$. To prove this, it suffices to
show:
\Claim\label{pesant}
For every point $v_1\in\oldPsi_1$, there is an open
neighborhood $ W_1$
% and $ W_2$ 
of
$v_1$
% and $v_2:=h(v_2)$ 
in $\oldPsi_1$ with the property that, if we set $ W_2=h( W_1)$,
% and $\oldXi_2$ respectively,   
there exist chart maps (see \ref{orbifolds introduced}) 
$\phi_i: U_i\to W_i$
 for $i=1,2$, and an affine automorphism $L$ of $\RR^n$ such that $L (U_1)= U_2$ and $|\phi_2|\circ(L|U_1)=h\circ|\phi_1|$.
\EndClaim

To prove \ref{pesant}, let
$v_1\in\oldPsi_1$ be given, and set
$v_2=h(v_1)$.  Note that \ref{pesant} is immediate in the case where
$v_i\notin\fraks_{\oldPsi_i}$ for
$i=1,2$. Note also that, in view of the orientability of the
$\oldPsi_i$, we have $
h(|\inter\oldPsi_1|)=h(\inter|\oldPsi_1|)=\inter|\oldPsi_2|=|\inter\oldPsi_2|$;
if the mutually equivalent conditions $v_1\in\inter\oldPsi_1$ and
$v_2\in\inter\oldPsi_2$ do hold, then in view of the existence of a homeomorphism
$f:\inter\oldPsi_1\to\inter\oldPsi_2$ with
$|f|=h\big|\inter|\oldPsi_1|$, the conclusion of \ref{pesant} holds. If the mutually homeomorphic orientable
orbifolds $\oldPsi_1$ have dimension at most $2$, we have
$\fraks_{\oldPsi_i}\subset\inter\oldPsi_i$ for $i=1,2$, and the
conclusion follows in this case. There remains the case in which
$\dim\oldPsi_i=3$ for $i=1,2$ while
$v_i\in\partial\oldPsi_i$ for $i=1,2$ and
$v_i\in\fraks_{\oldPsi_i}$ for some $i\in\{1,2\}$. After possibly
interchanging the roles of $\oldPsi_1$ and $\oldPsi_2$ and replacing
$h$ and $f$ by their inverses, we may assume
$v_1\in\fraks_{\oldPsi_1}$.

In this case, since $\oldPsi_1$ is an orientable $3$-orbifold and
$v_1\in\partial\oldPsi_1$, the stratum of $\fraks_{\oldPsi_1}$
containing $v$ must be a closed or half-open arc $A$ having $v$ as an
endpoint, and $G_{A}$ is cyclic of some order $n>1$. 
%Furthermore,
%$\inter A$ is  a stratum of $\inter\oldPsi_1$, and $G_{\inter A}$ is
%cyclic of order $n$. 
Let $B_1\subset|\oldPsi_1|$ be a $3$-ball such that
(i) $\gamma_1:=B_1\cap\fraks_{\oldPsi_1}$ is an arc contained in ${A}$,
properly embedded and unknotted in $B_1$, (ii) $\gamma_1\cap\partial
|\oldPsi_1|=\{v_1\}$, and (iii) $D_1:=B_1\cap\partial|\oldPsi_1|$ is a disk
having $v_1$ as an interior point. Since $\gamma_1\subset A$, the arc
$\gamma_1 =B_1\cap\fraks_{\oldPsi_1}=\fraks_{\obd(B_1)}$ is a single
  stratum of $\obd(B_1)$, and $G_{\gamma_1}$ is cyclic of order $n$. In particular, $\gamma_1^*:=\gamma_1-\{v_1\}=\gamma_1\cap|\inter\oldPsi_1|$ is a
  single stratum of $\obd(B_1)\cap\inter\oldPsi_1$ and is equal to 
$\fraks_{\obd(B_1)\cap\inter\oldPsi_1}$; and $G_{\gamma_1^*}$ is
cyclic of order $n$.

Set $B_2=h(B_1)$, $D_2=h(D_1)$, $\gamma_2=h(\gamma_1)$, and $\gamma_2^*=\gamma_2-\{v_2\}=h(\gamma_1^*)$. Since the homeomorphism $f:\inter\oldPsi_1\to\inter\oldPsi_2$ satisfies
$|f|=h\big|\inter|\oldPsi_1|$, we have
$
%\fraks_{\obd(B_2)}=
B_2\cap\fraks_{\oldPsi_2}\cap(\inter|\oldPsi_2|)=h(B_1\cap(\inter|\oldPsi_1|)\cap\fraks_{\oldPsi_1})=h(\gamma_1^*)=\gamma_2^*$; and
$\gamma_2^*$ is  a stratum of $\obd(B_2)\cap\inter\oldPsi_2$, with
 $G_{\gamma_2^*}$
cyclic of order $n$. Since $\gamma_2^*\subset\fraks_{\oldPsi_2}$,
we have $\gamma_2=\overline
{\gamma_2^*}\subset\fraks_{\oldPsi_2}$, and in particular
$v_2\in\fraks_{\oldPsi_2}$. But the
orientability of $\oldPsi_2$ implies that
$\fraks_{\oldPsi_2}\cap\partial|\oldPsi_2|$ is finite, and we may
therefore choose $B_1$ small enough so that
$B_2\cap\fraks_{\oldPsi_2}\cap\partial|\oldPsi_2|=\{v_2\}$. It now follows that $\fraks_{\obd(B_2)}=B_2\cap
\fraks_{\oldPsi_2}=\gamma_2$. Since $v_2\in\partial\oldPsi_2$, the
orientability of $\oldPsi_2$ implies that $\{v_2\}$ cannot be a
stratum of $\obd(B_2)$; hence $\gamma_2$ is a single stratum of
 $\obd(B_2)$, and $G_{\gamma_2}$ is a cyclic group of order $n$.

For $i=1,2$, set $
W_i=D_i\cup\inter B_i\subset B_i$, and let $\delta_i$ denote the
half-open arc $\gamma_i\cap W_i$. Let $p_i:\tW_i\to W_i$ denote the $n$-fold cyclic
branched cover of $W_i$ branched over $\gamma_i$. Since
$\gamma_2=h(\gamma_1)$, there is a
homeomorphism $\thh:\tW_1\to\tW_2$ such that
$p_2\circ\thh=(h|W_1)\circ p_1$. But since $W_1$ is a
$3$-ball, $D_1$ is a disk in $\partial W_1$, and $\gamma_1=W_1\cap
\fraks_{\oldPsi_1}$ is an arc meeting $D_1$ in a single endpoint,
there exists (in the notation of \ref{orbifolds introduced}) a homeomorphism $j:U^3_+\to W_1$, mapping
$U^2\times\{0\}\subset U^3_+$ onto $\inter D_1$, such that the deck
transformations of the branched covering $p_1\circ j:U^3_+\to W_1$
are orthogonal transformations.
 The branched covering $p_2\circ
\thh\circ j:U^3_+\to W_2$ has the same group of deck
transformations as $p_1\circ j$.  The conclusion of \ref{pesant} now
holds if we set 
%W_i=D_i\cup\inter W_i\subset W_i$ for $i=1,2$,  set
$U_1=U_2=U^2_+$, define $L$ to be the identity map of $U^2_+$, and set
$\phi_1=p_1\circ j|U^2_+$ and $\phi_2=p_2\circ
\thh\circ j|U^2_+$.
%B
%are orthogonal transformations. Then $p_2\circ j$ Since \redcomment{Fix from here to the end of the paragraph to
  %make it match the def. of a chart map. Use this if needed:
%  %$\tD_i:=p_i^{-1}(D_i)$ } Since for $i=1,2$ we have $\fraks_{\obd(B_i)}=B_i\cap
%\fraks_{\oldPsi_i}=\gamma_i$, and  $\gamma_i$ is a single stratum of
%$\obd(B_i)$ with $G_{\gamma_i}$  cyclic of order
%$n$, we may write $p_i=|\phi^0_i|$, where $\phi^0_i$ is a regular orbifold
%covering map from the manifold $\tB_i$ to the orbifold
%$\obd(B_i)$. The conclusion of \ref{pesant} now follows upon setting
%$ W_i=\obd(D_i\cup\inter B_i)\subset \obd(B_i)$,
%$ U_i=p_i^{-1}(| W_i|)\subset\tB_i$, and
%$\phi_i=\phi^0_i| U_i: U_i\to W_i$, for $i=1,2$; and
%setting $\tf=\thh| U_1: U_1\to U_2$. 

We have observed that \ref{pesant} implies the existence of a  homeomorphism
$\overline{f}:\oldPsi_1\to\oldPsi_2$ such that
$|\overline{f}|=h$. Since $f$ and
$\overline{f}|\inter\oldPsi_1$ are orbifold homeomorphisms from
$\inter\oldPsi_1$ to $\inter\oldPsi_2$ with
$|\overline{f}\big|\inter\oldPsi_1|=h\big||\inter\oldPsi_1|=|f|$,
we have $\overline{f}|\inter\oldPsi_1=f$, so that $\overline{f}$
extends $f$ as asserted.
\EndProof
%(i cyclic\phi\alpha \gamma j \th \tf^0\oldXi\frakU\obd\alpha 
%post-vQ \delta \DD \ref

\Proposition\label{oh yeah he tweets}
Let $n$ be an integer with $0<n\le3$, and let $\oldPsi_1$ and
$\oldPsi_2$ be   $n$-orbifolds. Suppose that $\oldPsi_2$ is
orientable, and that $\oldPsi_1$ has no $(n-1)$-dimensional stratum.
For $i=1,2$, let $\oldXi_i$ be an ($n-1$)-suborbifold of $\partial\oldPsi_i$.
Suppose that
$j:\silv_{\oldXi_1}\oldPsi_1\to\silv_{\oldXi_2}\oldPsi_2$ is an
embedding. Then there is a unique embedding of pairs $j^0:(\oldPsi_1,\oldXi_1)\to(\oldPsi_2,\oldXi_2)$ such
that
% $j^0(\oldXi_1)\subset\oldXi_2$ and 
$\silv j^0=j$. We have $(j^0)^{-1}(\oldXi_2)=\oldXi_1$. If $j$ is a homeomorphism then $j^0$ is a
homeomorphism of pairs. If
$(\oldPsi_1,\oldXi_1)=(\oldPsi_2,\oldXi_2)$, and $j$ is a homeomorphism isotopic to the
identity, then there is an isotopy $(j^0_t)_{0\le t\le1}$ of the pair 
$(\oldPsi_1,\oldXi_1)$ (see \ref{pairs}) such that $j^0=j^0_1$.
\EndProposition
%j'_0

\Proof
For $i=1,2$, set $\mu_i=
\mu_{\oldPsi_i,\oldXi_i}:
 \oldPsi_i\to\silv_{\oldXi_i}\oldPsi_i$,  
$\oldXi_i'=\mu(\oldXi_i)\subset\silv_{\oldXi_i}\oldPsi_i$, 
 and $\mu_i^*=
\mu^*_{\oldPsi_i,\oldXi_i}:
 \oldPsi_i-\oldXi_i\to(\silv_{\oldXi_i}\oldPsi_i)-\oldXi_i'$. 

We prove the uniqueness assertion first. Suppose that 
$j^0,j^1:(\oldPsi_1, \oldXi_1)\to(\oldPsi_2, \oldXi_2)$ are 
 embeddings such that 
%$j^i(\oldXi_1)\subset\oldXi_2$ and 
$\silv j^i=j$ for $i=0,1$. 
Since 
%$|\mu_1|$ and
$|\mu_2|$ is a homeomorphism, and since \ref{third point} gives
$\mu_2\circ j^0=j\circ\mu_1=\mu_2\circ j^1$, we have
$|j^0|=|j^1|$. Hence the orbifold embeddings $j^0$ and $j^1$ must
coincide, and uniqueness is established.

We next prove the existence assertion in the case where $j$
is a homeomorphism. For $i=1,2$, the hypothesis implies that
$\oldPsi_i$ has no $(n-1)$-dimensional strata. Hence $|\oldXi_i'|$ is the closure of the union of the
  $(n-1)$-dimensional strata of $\silv_{\oldXi_i}\oldPsi_i$.  The
  homeomorphism $j$ therefore maps $\oldXi_1'$ onto $\oldXi_2'$, and thus
  restricts to a homeomorphism of $(\silv_{\oldXi_1}\oldPsi_1)-\oldXi_1'$
  onto $(\silv_{\oldXi_2}\oldPsi_2)-\oldXi_2'$. Since $\mu_1^*$ and $\mu_2^*$
   are homeomorphisms   by \ref{silvering}, we have a homeomorphism
  $f:=(\mu_2^*)^{-1}\circ j\circ\mu_1^*$ from $ \oldPsi_1-\oldXi_1$
  onto $ \oldPsi_2-\oldXi_2$. In particular $ \oldPsi_1-\oldXi_1$ is
  orientable; this implies in particular that $\inter\oldPsi_1$ is
  orientable, and hence that $\oldPsi_1$ is orientable. (It may be
  worth pointing out that the orientability of $\inter\oldPsi_1$ is
  stronger than the assertion that $\obd(\inter|\oldPsi_1|)$ is
  orientable, which would {\it not} imply orientability for $\oldPsi_1$.)

Since $|\mu_1|$ and $|\mu_2|$
   are also homeomorphisms   by \ref{silvering}, we have a homeomorphism
  $h:=|\mu_2|^{-1}\circ| j|\circ|\mu_1|$ from $ |\oldPsi_1|$
  onto $ |\oldPsi_2|$. Since 
$j(\oldXi_1')=\oldXi_2'$, we have $h(|\oldXi_1|)=|\oldXi_2|$.
We have $|f|=h\big||\oldPsi_1-\oldXi_1|$; in
  particular, $|f|$ and $h$ agree on $\inter|\oldPsi_1|$. As
  $\oldPsi_1$ and $\oldPsi_2$ are orientable, it now follows from Lemma \ref{shabbos} that 
 $f$ extends to a homeomorphism
$j^0:\oldPsi_1\to\oldPsi_2$, and that $|\overline{j^0}|=h$. Hence
$j^0(\oldXi_1)=\oldXi_2$, so that $j^0$ is a homeomorphism from
$(\oldPsi_1,\oldXi_1)$ to $(\oldPsi_2,\oldXi_2)$.
It follows from this construction that $\mu_2\circ j^0$ and
$j\circ\mu_1$ agree on the dense suborbifold $\oldPsi_1-\oldXi_1$ of
$\oldPsi_1$, and hence that
$\mu_2\circ j^0=j\circ\mu_1$. By \ref{silvering} this implies that
$\silv j^0=j$, and the existence assertion is established in this case.

%, \redcomment{Why? This hasn't been made clear at
  %all. Come to think of it, I don't understand the last sentence of
  %the next paragraph either} and it
%follows formally from the constructions that
%that $\mu_2\circ j^0=j\circ\mu_1$. This proves the existence
%assertion in this case.

To prove the existence assertion in general, suppose that $j:\silv_{\oldXi_1}\oldPsi_1\to\silv_{\oldXi_2}\oldPsi_2$ is an
embedding, and write $j=\iota\circ j'$, where $j'$ is a homeomorphism
of $\silv_{\oldXi_1}\oldPsi_1$ onto a suborbifold $\frakZ$ of
$\silv_{\oldXi_2}\oldPsi_2$, and $\iota:\frakZ\to
\silv_{\oldXi_2}\oldPsi_2$ is the inclusion. According to \ref{second
  point} we may write
$\frakZ=\silv_{\oldXi_2'}\oldPsi_2'$, where 
$\oldPsi_2'$ is
a
suborbifold  of $\oldPsi_2$ such that $\oldXi_2':=\oldXi_2\cap\oldPsi_2'$ is a
codimension-$0$ 
suborbifold of $\partial\oldPsi_2'$. By the case of the lemma already
proved, the homeomorphism $j': \silv_{\oldXi_1}\oldPsi_1\to
\silv_{\oldXi_2'}\oldPsi_2'$ is equal to $\silv((j')^0)$ for some
homeomorphism  $(j')^0:(\oldPsi_1, \oldXi_1)\to(\oldPsi_2',
\oldXi_2')$.  Now if $\iota^0:(\oldPsi_2',\oldXi_2')\to(\oldPsi_2,\oldXi_2)$ denotes the
inclusion, we have $\iota=\silv\iota^0$. Since $(j')^0$ is a
homeomorphism, $j^0:=\iota^0\circ(j')^0$ is an embedding; and $\silv
j^0=\silv(\iota^0)\circ\silv((j')^0)=\iota\circ j'=j$. Furthermore,
since $\oldXi_2':=\oldXi_2\cap\oldPsi_2'$, we have 
$(\iota^0)^{-1}(\oldXi_2)=\oldXi_2'$, and since $(j')^0$ is a
homeomorphism, it then follows that
$(j^0)^{-1}(\oldXi_2)=\oldXi_1$.
This completes
the existence proof.

%Since $j$ is an orbifold embedding, and since
%$|\mu_1(\oldXi_1)|\subset\fraks_{\frakZ_1}$ is ($n-1$)-dimensional,
%$|j(\mu_1(\oldXi_1))|=|\mu_2(j^0(\oldXi_1))|$ is contained in the closure of the union of the
%($n-1$)-dimensional strata of $\frakZ_2$, which is in turn equal to
%$|\mu_2(\oldXi_2)|$ since $\oldPsi_2$ is orientable. This shows that
%$j^0(\oldXi_1)\subset\oldXi_2$.

Finally, suppose that $(\oldPsi_1,\oldXi_1)=(\oldPsi_2,\oldXi_2)$, and
that $j$ is a homeomorphism isotopic to the
identity. Fix an isotopy $(j_t)_{0\le t\le1}$, where
$j_t:\frakZ_1\to\frakZ_2$ is a homeomorphism for each $t$, while $j_0$
is the identity and $j_1=j$. By the assertions already proved, for
each $t\in[0,1]$ there is  a unique self-homeomorphism $j^0_t$ of $(\oldPsi_1,\oldXi_1)$ such
that $\silv j^0_t=j_t$. Then $(j^0_t)_{0\le t\le1}$ is an
isotopy of the pair
$(\oldPsi_1,\oldXi_1)$, and by the uniqueness assertion already
proved, $j^0_0$ is the identity and $j^0_1=j^0$. This proves the final assertion.
\EndProof
%\oldUpsilon\eta j' two \ref

%\Remark
%It is not difficult to show that the embedding $j^0$ given by the
%first assertion of Lemma \ref{oh yeah he tweets} has the property that
%$(j^0)^{-1}(\oldXi_2)=\oldXi_1$. Adding this conclusion (which will
%not be needed in the applications) would make
%first assertion of Lemma \ref{oh yeah he tweets} a converse to the
%fact pointed out  in \ref{third point}. 
%\EndRemark

\Number\label{scarfuss}
According to the main theorem of \cite{illman}, every very good smooth orbifold $\oldOmega$ admits a PL structure compatible with its smooth structure, and this PL structure is unique up to PL orbifold homeomorphism. We shall denote by $\oldOmega\pl$ the orbifold $\oldOmega$ equipped with this PL structure. The main result of \cite{lange} implies that every very good PL orbifold of dimension $n\le4$ is PL homeomorphic to $\oldOmega\pl$ for some smooth $n$-orbifold $\oldOmega$. However, as there is no result in the literature guaranteeing, even for $n=3$, that every PL suborbifold of $\oldOmega\pl$ is ambiently homeomorphic to a smooth suborbifold of $\oldOmega$, a little care will be required in passing between the two categories. (A more subtle question about the relationship between the categories is addressed by Proposition \ref{fibration-category}.)
\EndNumber

The following result extends the results quoted above from \cite{illman} and \cite{lange}.

\Proposition\label{yoga}
Let $(\oldUpsilon,\oldXi)$  be an orbifold pair, where $\oldUpsilon$ is compact, orientable, and very good, $n:=\dim\oldUpsilon\le3$, and $\oldXi\subset\partial\oldUpsilon$ is compact and has dimension $n-1$. Then every smooth structure on $(\oldUpsilon,\oldXi)$ is compatible with a PL structure, which is unique up to equivalence. Furthermore,
every PL structure on $(\oldUpsilon,\oldXi)$ is compatible with a smooth structure, which is unique up to equivalence. \EndProposition

\Proof
To prove the first assertion, first suppose that
$(\oldUpsilon,\oldXi)$ is equipped with a smooth structure. Since
$\oldUpsilon$ is orientable and $n\le3$, the orbifold
$\silv_\oldXi\oldUpsilon$ acquires a smooth structure according to the
discussion in \ref{smooth doubling} and \ref{silvering}. The main
result of \cite{illman} implies that $\silv_\oldXi\oldUpsilon$ has a
PL structure compatible with its smooth structure. Set
$\mu=
\mu_{\oldUpsilon,\oldXi}:\oldUpsilon\to\silv_\oldXi\oldUpsilon$  (see \ref{silvering}).  The PL structure that we have assigned to
$\silv_\oldXi\oldUpsilon$ may be pulled back via
 $\mu$ to define a PL structure on $\oldUpsilon$. 
Furthermore, if ${\rm D}_\oldXi\oldUpsilon $ is equipped with the PL structure obtained by pulling back the PL structure on 
$\silv_\oldXi\oldUpsilon$ via $p$, then the non-trivial deck transformation of the covering $p:{\rm D}_\oldXi\oldUpsilon \to\silv_\oldXi\oldUpsilon$ is a PL 
involution, and hence its fixed point set is a PL suborbifold of
${\rm D}_\oldXi\oldUpsilon $; this implies that 
$\oldXi$ is a PL suborbifold of $\oldUpsilon$, so that the orbifold pair $(\oldUpsilon,\oldXi)$ acquires a PL structure.

On the other hand, since the PL structure on $\silv_\oldXi\oldUpsilon$ is compatible with its smooth structure, the pulled back PL structure on $\oldUpsilon$ is compatible with the smooth structure obtained by pulling back the smooth structure on
$\silv_\oldXi\oldUpsilon$, which is simply the given smooth structure on $\oldUpsilon$. Hence the PL structure on
$(\oldUpsilon,\oldXi)$ is compatible with its given smooth structure. This proves the existence part of the first assertion.
%since the fixed point set of the canonical involution

 To prove the uniqueness part of the first assertion, suppose that
 $(\oldUpsilon,\oldXi)$ has two PL structures compatible with its
 given smooth structure. Equipping the pair $(\oldUpsilon,\oldXi)$
 with these two PL structures gives two PL orbifold pairs
 $(\oldUpsilon_1,\oldXi_1)$ and
 $(\oldUpsilon_2,\oldXi_2)$. The assumption that the given   PL structures on
 $(\oldUpsilon,\oldXi)$ are compatible with the same smooth structure
 may be paraphrased as saying that
% there are has two PL structures compatible with its
% given smooth structure, 
there are smooth structures on $(\oldUpsilon_1,\oldXi_1)$ and
 $(\oldUpsilon_2,\oldXi_2)$, compatible with their PL structures,
 under which they are diffeomorphic. 
Hence there are smooth structures on $\silv_{\oldXi_1}\oldUpsilon_1$
and $\silv_{\oldXi_2}\oldUpsilon_2$, compatible with the PL structures that
 they inherit from $(\oldUpsilon_1,\oldXi_1)$ and
 $(\oldUpsilon_2,\oldXi_2)$, under  which they are diffeomorphic.

It therefore follows from
%\redcomment{Fix from here.} These give rise to two PL structures on
% $\silv_\oldXi\oldUpsilon$, which are compatible with the smooth
% structure on $\silv_\oldXi\oldUpsilon$ determined by the smooth
 %structure on $(\oldUpsilon,\oldXi)$, and are therefore equivalent by
 the main result of \cite{illman} that there is a PL homeomorphism
 $h:\silv_{\oldXi_1}\oldUpsilon_1\to
 \silv_{\oldXi_1}\oldUpsilon_2$. Applying the PL version of Proposition \ref{oh yeah
   he tweets}, letting the $\oldUpsilon_i$ play the roles of the
 $\oldPsi_i$ and defining the $\oldXi_i$ as above, we obtain a PL homeomorphism 
%Since $\oldUpsilon$ is
 %orientable, $|\oldXi|$ is the union of the closures of the
 %$2$-dimensional components of $\fraks_\oldUpsilon$; hence $|h|$ maps
 %$|\oldXi|$ onto itself. On the other hand, 
 %$\mu:\oldUpsilon\to\silv_\oldXi\oldUpsilon$  is an
 %immersion of topological orbifolds, and $|\mu|$ is the identity
 %map of $|\oldUpsilon|=|\silv_\oldXi\oldUpsilon|$. Hence there is a
 %self-immersion  
$h^0:(\oldUpsilon_1,\oldXi_1)\to (\oldUpsilon_2,\oldXi_2)$ such
 that  $\silv h^0=h$. 
%According to the
 %definition of $\silv h^0$, this implies that $\mu\circ h^0=h\circ\mu$.
%Since $h$ and $|\mu|$ are
 %homeomrophisms, $h^0$ is a self-homeomorphism of $\oldUpsilon$.
%\redcomment{That sentence is not
  % clear to me at all. Furthermore, it makes no sense to say
   %$|\mu|$ is the identity; it's a homeo, as I pointed out in
   %\ref{silvering}.} 
%But if $\oldUpsilon$ is equipped with either of its two given PL
%structures, and if $\silv_\oldXi\oldUpsilon$ is assigned  the
%corresponding PL structure, then $\mu$ becomes a PL immersion. It
%therefore follows from the equality $\mu\circ h^0=h\circ\mu$ that $h^0$ carries one of the
%given PL structures on
% $\oldUpsilon$ to the other. Since
%$h(\oldXi)=\oldXi$, we may regard $h^0$ as a self-homeomorphism of
 %the
%pair $(\oldUpsilon,\oldXi)$ 
%carrying one of its given PL structures 
%to
%the other. \redcomment{This should be rewritten to avoid saying that
%  $\mu$ is PL, because when the argument is turned around to prove the
  %second assertion we will not have the info that $\mu$ is smooth. But
%now it seems to be important to make sure that Prop. \ref{oh yeah he
  %tweets} works as stated in the smooth category, which is something I
%have not thought out clearly.} 
This gives the required equivalence between the given PL 
structures, and completes the proof of the first
assertion. 
%h'

The proof of the second assertion is the same as the proof of the
first, except that the roles of PL and smooth structures are
interchanged,  the main theorem of \cite{lange} is used in place of
the main theorem of \cite{illman}, and the smooth version of
Proposition \ref{oh yeah he tweets} is used in place of the PL version.
\EndProof
%\alpha\eta\oldUpsilon \eta tweets\mu point} tweets} shabbos} \ref

\Number\label{categorille}For the rest of this monograph we will adopt
the convention, generalizing the convention stated in \ref{just
  manifolds},  that all statements and arguments about orbifolds are
to be interpreted in the PL category except where another category is
specified. The reason for emphasizing the PL category is that orbifold
fibrations, which will be defined and considered in the next chapter, and form the
basis for the characteristic suborbifold theory, would not
behave well in the smooth category if both the base and fiber
have non-empty boundary; in this case the total space of the fibration
would be an orbifold with corners. We have not attempted to define
orbifolds with corners in this monograph, although manifolds with
corners make a few transitory appearances, both above and below the
present passage.

Smooth orbifolds will indeed be considered at a number of points in
the monograph; in particular, as hyperbolic orbifolds are smooth by
definition, statements involving hyperbolic orbifolds (including some
of the main results of the monograph) should be interpreted in the
smooth category. This makes it necessary to make the transition
between the smooth and PL categories at a number of
points. Proposition \ref{fibration-category} is one of the more subtle
instances of this.

The results of \cite{bonahon-siebenmann} will play an important role
in this monograph. While these results (which are largely concerned
with fibrations, but almost exclusively with those that have closed fibers) are proved in the smooth category in \cite{bonahon-siebenmann}, the statements and proofs go through without change in the PL category, and it is the PL versions that will be quoted in the proofs of  Propositions \ref{when vertical} and \ref{new characteristic}. Indeed, the proofs of the PL versions of the results of \cite{bonahon-siebenmann} are simpler than the proofs of the smooth versions, in that some of the arguments involve extending an isotopy of the 2-sphere to an isotopy of  the 3-ball, and the orbifold analogue of this. In the PL category this follows from Alexander's coning trick, whereas in the smooth category it requires a deep  theorem due to Cerf. I am indebted to Francis Bonahon for these observations.

\EndNumber

\DefinitionsRemarks\label{wuzza turnover}
A $2$-orbifold of finite type (see \ref{orbifolds introduced}) will be
termed {\it negative} if each of its 
components has negative Euler characteristic. Note that the empty
$2$-orbifold is negative.

We define a {\it negative turnover} to be a negative, compact, orientable $2$-orbifold $\oldTheta$ such that $|\oldTheta|$ is a $2$-sphere and  $\card \fraks_\oldTheta=3$. 
\EndDefinitionsRemarks

\Number\label{2-dim case}
Let $\otheroldLambda$ be a finite-type  $2$-orbifold. Suppose that $\fraks_\otheroldLambda$ is finite (which is in particular true if $\otheroldLambda$ is orientable). Let $x_1,\ldots,x_n$ denote the distinct points of $\fraks_\otheroldLambda$, and let $p_i$ denote the order of the singular point $x_i$ for $i=1,\ldots,n$. Then we have $\chi(\otheroldLambda)=\chi(|\otheroldLambda|)-\sum_{i=1}^n(1-1/p_i)$. In particular we have $\chi(\otheroldLambda)\le\chi(|\otheroldLambda|)$. 

Another consequence of the formula for $\chi(\otheroldLambda)$ given above, which will be used many times in this monograph, often without an explicit reference, is that if $\otheroldLambda$ is an orientable annular orbifold, then either $|\otheroldLambda|$ is an annulus and $\fraks_\otheroldLambda=\emptyset$, or $|\otheroldLambda|$ is a disk and $\fraks_\otheroldLambda$ consists of two points, both of order $2$.
\EndNumber

\Number\label{cobound}An immediate consequence of the description of orientable annular orbifolds given in \ref{2-dim case} is that If $\frakB$ is an orientable annular $2$-orbifold, and if
$C\subset|\frakB|-\fraks_\frakB$ is a simple closed curve such that
$\omega(C)$ is $\pi_1$-injective in $\frakB$, then there is a
weight-$0$ annulus $A\subset|\frakB|$ having $C$ as a boundary
component, and having its other boundary component contained in
$\partial|\frakB|$. 
\EndNumber

The following lemma is preparation for the proof of Proposition \ref{at least a sixth}, which will be needed in Section \ref{dandy section}.

\Lemma\label{at least a lemma}
Let $\oldUpsilon$ be a negative (see \ref{wuzza turnover}), compact, connected $2$-orbifold which is  not a negative turnover. Then we have
%\compnum(|\oldUpsilon|)\le6
$\chibar(\oldUpsilon )\ge1/6$.
Furthermore, if $\chibar(\oldUpsilon )<1/3$, then either (a) $|\oldUpsilon|$ is a disk, and
$\fraks_\oldUpsilon$ consists of exactly two points, exactly one of which has
order $2$, or (b) $|\oldUpsilon|$ is a sphere, and
$\fraks_\oldUpsilon$ consists of exactly four points, exactly one of which has
order greater than $2$.
\EndLemma

\Proof
We must prove that either (i) $\chibar(\oldUpsilon)\ge1/3$, or (ii) 
$|\oldUpsilon|$ is a disk, 
$\fraks_\oldUpsilon$ consists of exactly two singular points, exactly one of which has
order $2$, and $\chibar(\oldUpsilon )\ge1/6$, or (iii) $|\oldUpsilon|$ is a sphere, and
$\fraks_\oldUpsilon$ consists of exactly four points, exactly one of which has
order greater than $2$, and $\chibar(\oldUpsilon )\ge1/6$.

Set $m=\card\fraks_\oldUpsilon\ge0$, and set $\fraks_\oldUpsilon=\{x_1,\ldots,x_m\}$. Let $p_i\ge2$ denote the order of
the singular point $x_i$. Let $g$ denote the genus of
$|\oldUpsilon|$, and set $b=\compnum(\partial|\oldUpsilon|)$.
%, so that the right hand sides of (\ref{freep your feet}) and
%(\ref{mud your mississip}) are respectively
%$\max(1,2b)/6\le\max(1,b)/3$ and $\max(1,b)/3$. 
We have (cf. \ref{2-dim case})
\Equation\label{fancy delancey}
\chibar(\oldUpsilon )=\chibar(|\oldUpsilon |)+\sum_{i=1}^m(1-1/p_i)
=(2g+b-2)+\sum_{i=1}^m(1-1/p_i).
\EndEquation
If $g\ge2$, or if $g=1$ and $b\ge1$, then (\ref{fancy delancey})
implies that $\chibar(\oldUpsilon )\ge 2g+b-2\ge1$, so that (i) holds. 
If $g=1$ and $b=0$, then (\ref{fancy delancey})
implies that $\chibar(\oldUpsilon )=\sum_{i=1}^m(1-1/p_i)$. Since by
hypothesis we have $\chibar(\oldUpsilon )>0$, we must have $m\ge1$; and
since each term in the sum $
\sum_{i=1}^m(1-1/p_i)$ is at least $1/2$, we have $\chibar(\oldUpsilon
)\ge1/2$, so that (i) holds.

For the rest of the proof we will assume that $g=0$, so that
(\ref{fancy delancey}) becomes
\Equation\label{hott mott}
\chibar(\oldUpsilon )=
b-2+\sum_{i=1}^m(1-1/p_i).
\EndEquation
If $b>2$, (\ref{hott mott}) implies that $\chibar(\oldUpsilon )\ge
b-2\ge1$, so that (i) holds.
If $b=2$, (\ref{hott mott}) gives 
$\chibar(\oldUpsilon)=\sum_{i=1}^m(1-1/p_i)$. Since
$\chibar(\oldUpsilon)>0$ we must then have $m>0$, and since each of the
terms $1-1/p_i$ is at least $1/2$, we have $\chibar(\oldUpsilon
)\ge1/2$, so that (i) holds.
If $b=1$ then (\ref{hott mott}) becomes 
$\chibar(\oldUpsilon)=-1+\sum_{i=1}^m(1-1/p_i)$.
Since $\chibar(\oldUpsilon )>0$, and since each term
$1-1/p_i$ is strictly less than $1$, we must have 
$m\ge2$. In the case where $b=1$ 
and $m\ge3$, since each of the terms $1-1/p_i$ is at least $1/2$, we
have $\chibar(\oldUpsilon)\ge1/2$, so that (i) holds.

If $b=1$ 
and $m=2$, then since $\chibar(\oldUpsilon )>0$,
it follows from (\ref{hott mott}) that either $p_1$ or $p_2$ is
at least $3$, and hence that $\chibar(\oldUpsilon )\ge-1+1/2+2/3=1
/6$. We distinguish two subcases. If neither $p_1$ nor $p_2$
is equal to $2$, then $p_1$ and $p_2$ are each at least $3$, and hence
$\chibar(\oldUpsilon )\ge-1+2/3+2/3=1
/3$; thus (i) holds in this subcase. Now consider the subcase in which either $p_1$ or $p_2$ is equal to
$2$. Then $\fraks_\oldUpsilon$ consists of two points, exactly one of which has
order $2$. The surface
$|\oldUpsilon|$ is a disk since $g=0$ and $b=1$. Thus (ii) holds in this subcase.

Now suppose that $b=0$. Then (\ref{hott mott}) becomes 
$\chibar(\oldUpsilon)=-2+\sum_{i=1}^m(1-1/p_i)$. Since $\chibar(\oldUpsilon )>0$, and since each term
$1-1/p_i$ is strictly less than $1$, we must have 
$m\ge3$. But if $m=3$, then since $b=g=0$, the orbifold $\oldUpsilon$ is a
negative turnover, a contradiction to the hypothesis. Hence
$m\ge4$. 
If $m>4$, then since every term $1-1/p_i$ is at least $1/2$, we have
$\chibar(\oldUpsilon )\ge1/2$, and (i) holds.
If $m=4$, then since $\chibar(\oldUpsilon )>0$,
one of the $p_i$ is
at least $3$, and hence $\chibar(\oldUpsilon )\ge-2+1/2 +1/2 +1/2 +2/3
=1/6$.
We distinguish two subcases. If at least
two of the $p_i$ are greater than $2$, then we have
$\chibar(\oldUpsilon )\ge-2+1/2+1/2+2/3+2/3=1
/3$; thus (i) holds in this subcase. Now consider the subcase in which only one of the $p_1$  is greater than
$2$; that is, $\fraks_\oldUpsilon$ consists of four points, exactly one of which has
order greater than $2$. The surface 
$|\oldUpsilon|$ is a sphere since $b=g=0$. Thus (iii) holds in this subcase. 
\EndProof

\Proposition\label{at least a sixth}
Let $\oldGamma$ be a negative, compact $2$-orbifold, no component of which is a negative turnover. Then we have
\Equation\label{freep your feet}
\compnum(|\oldGamma|)\le6
\chibar(\oldGamma ).
\EndEquation
If in addition we assume that 
%for every
%integer $d>1$, 
the number of points of order $2$ in $\fraks_\oldGamma$
% the $2$-orbifold
%$\overline{(\partial\oldPsi)-\oldPhi(\oldPsi)}$ 
is even, then we have
\Equation\label{mud your mississip}
\compnum(|\oldGamma|)\le2\lfloor
3\chibar(\oldGamma )\rfloor.
\EndEquation
\EndProposition

\Proof
First note that for any component $\oldUpsilon$
of $\oldGamma$, it follows from the first assertion of Lemma \ref{at least a lemma} that
$\chibar(\oldUpsilon )\ge1/6$.
Summing over the components of $\oldGamma$, we deduce that $\chibar(\oldGamma)\ge
\compnum(|\oldGamma|)/6$, which is equivalent to the first
assertion of the proposition.

To prove the second assertion, first consider the case in which
$\oldGamma$ has a component $\oldUpsilon_0$ with the property that
$\chibar(\oldUpsilon)\ge1/3$. Let $\oldUpsilon_1,\ldots,\oldUpsilon_n$ denote
the remaining components of $\oldGamma$, where $n\ge0$. According to the
first assertion of Lemma \ref{at least a lemma}, we have
$\chibar(\oldUpsilon _i)\ge1/6$ for $i=1,\ldots,n$. Hence
$$\chibar(\oldGamma)=\chibar(\oldUpsilon _0)+\sum_{i=1}^n\chibar(\oldUpsilon _i)\ge\frac13+n\cdot\frac16=\frac{n+2}6,$$
which implies that
$$\compnum(|\oldGamma|)=n+1\le6\chibar(\oldGamma)-1.$$
Since $\compnum(|\oldGamma|)$ is an integer, it follows that
$$\compnum(|\oldGamma|)\le\lfloor6\chibar(\oldGamma)\rfloor-1=\lfloor2\cdot3\chibar(\oldGamma)\rfloor-1\le2\lfloor
3\chibar(\oldGamma )\rfloor,$$
which gives the second assertion in this case.

There remains the case in which $\chibar(\oldUpsilon)<1/3$ for every
component $\oldUpsilon$ of $\oldGamma$. In this case, it follows from the
second assertion of Lemma \ref{at least a lemma} that for every
component $\oldUpsilon$ of $\oldGamma$, either (a) the surface $|\oldUpsilon|$ is a disk, and
$\fraks_\oldUpsilon$ consists of two points, exactly one of which has
order $2$, or (b) $|\oldUpsilon|$ is a sphere, and
$\fraks_\oldUpsilon$ consists of four points, exactly one of which has
order greater than $2$.  In particular,  $\compnum(|\oldGamma|)$ is
congruent modulo $2$ to the number of points of order $2$ in
$\fraks_\oldGamma$. According to the hypothesis of the second assertion of the
present proposition, the number of point of order $2$ in
$\fraks_\oldGamma$ is even. Hence  $\compnum(|\oldGamma|)$ is
even. Let us write $\compnum(|\oldGamma|)=2k$, where $k$ is a positive integer, and let $\oldUpsilon_1,\ldots,\oldUpsilon_{2k}$ denote
the components of $\oldGamma$. Then according to the
first assertion of Lemma \ref{at least a lemma}, we have
$\chibar(\oldUpsilon _i)\ge1/6$ for $i=1,\ldots,2k$. Hence
$$\chibar(\oldGamma)=\sum_{i=1}^{2k}\chibar(\oldUpsilon _i)\ge2k\cdot\frac16=\frac{k}3,$$
which implies that
$k\le3 \chibar(\oldGamma)$. Since $k$ is an integer, we have $k\le\lfloor3
\chibar(\oldGamma)\rfloor$. Hence
$$\compnum(|\oldGamma|)=2k\le2\lfloor
3\chibar(\oldGamma )\rfloor,$$
which gives the second assertion in this remaining case.
\EndProof
%\eta $T

\section{Some results on $3$-orbifolds}\label{Orb3 section}

\Lemma\label{cover-irreducible}
Suppose that $\oldLambda$ is a connected, orientable $3$-orbifold with non-empty boundary, and that some regular covering of $\oldLambda$ is irreducible. Then $\oldLambda$ is irreducible.
\EndLemma

\Proof
Let $p:\toldLambda\to\oldLambda$ be an irreducible regular covering of $\oldLambda$.
Let $\oldPi$ be any two-sided spherical $2$-suborbifold of $\inter\oldPi$. Then every component $\toldPi$ of $p^{-1}(\oldPi)$ is two-sided and spherical, and hence bounds a discal $3$-suborbifold of $\toldLambda$; this discal $3$-suborbifold is unique since $\partial\toldLambda\ne\emptyset$, and will be denoted by $\tfrakB_{\toldPi}$. If $\toldPi$ and $\toldPi'$ are distinct components of $p^{-1}(\oldPi)$, there is a deck transformation $t$ carrying $\toldPi$ to $\toldPi'$; since $\toldPi$ and $t(\toldPi)=\toldPi'$ are disjoint, and $\partial\toldLambda\ne\emptyset$, we must have either (a) $t(\tfrakB_{\toldPi})\subset\inter\tfrakB_{\toldPi}$, (b) $\tfrakB_{\toldPi}\subset\inter t(\tfrakB_{\toldPi})$, or (c) $t(\tfrakB_{\toldPi})\cap\tfrakB_{\toldPi}=\emptyset$. If (a) or (b) holds, there are infinitely many deck transformations mapping the compact suborbifold $\tfrakB_{\toldPi}$ of $\toldLambda$ into itself, which is impossible. Hence the discal orbifolds
$\tfrakB_{\toldPi}$ and $\tfrakB_{\toldPi'}=t(\tfrakB_{\toldPi})$ are disjoint. Since this holds for any two distinct components $\toldPi$ and $\toldPi'$ of $p^{-1}(\oldPi)$, there is a suborbifold $\frakB$ of $\oldLambda$ such that $p^{-1}(\frakB)=\bigcup_{\toldPi\in\calc(p^{-1}(\oldPi))}\tfrakB_{\toldPi}$; we have $\partial\frakB=\oldPi$, and for any component $\toldPi$ of $p^{-1}(\oldPi)$, the orbifold $\frakB$ is homeomorphic to the quotient of $\tfrakB_{\toldPi}$ by its stabilizer in the group of deck transformations. Since $\tfrakB_{\toldPi}$ is discal, its quotient $\frakB$ by a finite group action is discal according to \ref{annular etc}.
\EndProof
%\toldPi'$t\tau\frakt

\Proposition\label{kinda dumb}
Let $\oldPsi$ be a very good $3$-orbifold, 
and let $\oldTheta$ be a $2$-suborbifold of $\oldPsi$ which is either contained in  $\partial\oldPsi$  or two-sided in $\oldPsi$. Then the following conditions are equivalent:
\begin{enumerate}[(a)]\item$\oldTheta$ is $\pi_1$-injective in $\oldPsi$; 
\item For every discal $2$-suborbifold $\frakD$ of $\oldPsi$ such that $\frakD\cap\oldTheta=\partial\frak D$, there is a discal $2$-suborbifold $\frakE$ of $\oldTheta$ with $\partial\frakE=\partial\frakD$.
\end{enumerate}
If in addition we assume that $\oldPsi$ is orientable and that $\oldTheta$ has no spherical component, then (a) and (b) are equivalent to:
\begin{enumerate}[(c)]
\item For every orientable discal $2$-suborbifold $\frakD$ of $\oldPsi$ such that $\frakD\cap\oldTheta=\partial\frakD$, there is a discal $2$-suborbifold $\frakE$ of $\oldTheta$ with $\partial\frakE=\partial\frakD$.
\end{enumerate}
\EndProposition

\Proof
We first give the proof in the case where $\oldTheta\subset\partial
\oldPsi$. Since $\oldPsi$ is very good, we may fix a regular finite-sheeted cover $p:\toldPsi\to\oldPsi$ such that $\toldPsi$ is a $3$-manifold. Let $G$ denote the group of deck transformations of this covering.

Let us show that (a) implies (b).
Suppose that (a) holds, and that $\frakD$ is a discal $2$-suborbifold of $\oldPsi$  such that $\frakD\cap\oldTheta=\partial\frakD$. Then every component of $p^{-1}(\frakD)$ is a disk whose boundary lies in $p^{-1}(\oldTheta)$. Since $\oldTheta$ is $\pi_1$-injective in $\oldPsi$, in particular $p^{-1}(\oldTheta)$ is $\pi_1$-injective in $\toldPsi$, and so each component of $p^{-1}(\frakD)$ bounds a disk in $p^{-1}(\oldTheta)$. Among all disks in $p^{-1}(\oldTheta)$ bounded by components of $p^{-1}(\frakD)$, choose one, $D_0$, which is minimal with respect to inclusion. If  $g\in G$ is given, $g(D_0)\cup D_0$ cannot be a component of $p^{-1}(\oldTheta)$, since $\oldTheta$ has no spherical component; in view of the minimality of $D_0$ it follows that either $g(D_0)=D_0$ or $g(D_0)\cap D_0=\emptyset$. Since one of these alternatives holds for each $g\in G$, the orbifold $\frakE:=p(D_0)$ is homeomorphic to the quotient of $D_0$ by its stabilizer in $G$, and is therefore a discal suborbifold of $\oldTheta$ with $\partial\frakE=\frakD$. This establishes (b).

To show that (b) implies (a), we will suppose that (a) does not hold, and produce a discal orbifold $\frakD$ violating (b).
Since (a) does not hold, 
%$\oldTheta$ is not $\pi_1$-injective in $\oldPsi$. Hence 
$p^{-1}(\oldTheta)$ is not 
$\pi_1$-injective in $\toldPsi$. It therefore follows from the PL version of the equivariant loop theorem 
\cite[Theorem 4.3]{jaco-rubinstein} that there is a non-empty, properly embedded, $G$-invariant submanifold $\cald$ of $\toldPsi$, each of whose components is a disk, such that $\partial\cald\subset p^{-1}(\inter\oldTheta)$, and no component of $\partial\cald$ bounds a disk in $p^{-1}(\oldTheta)$. If $D_1$ is a component of $\cald$ then $\frakD:=p(D_1)$ is orbifold-homeomorphic to the quotient of $D_1$ by its stabilizer in $G$. Hence $\frakD$ is a discal suborbifold of $\oldPsi$. 
To show that $\frakD$ violates (b), it is enough to show that there is no  discal $2$-suborbifold $\frakE$ of $\oldTheta$ with $\partial\frakE=\partial\frakD$. If $\frakE$ is such a discal $2$-suborbifold, then $\tfrakE:=p^{-1}(\frakE)$ is a disjoint union of disks, and $\partial\tfrakE=p^{-1}(\partial\frakD)$. In particular the component $\partial D_1$ of $\partial\cald$ bounds a disk in $p^{-1}(\oldTheta)$, a contradiction. Thus we have shown that (a) and (b) are equivalent (in the case $\oldTheta\subset\partial \oldPsi$).

To prove the second assertion (in the case $\oldTheta\subset\partial \oldPsi$), it suffices to show, under the additional assumption that $\oldPsi$ is orientable and that $\oldTheta$ has no spherical component, that if
$\frakD$  is a discal $2$-suborbifold of $\oldPsi$ with
$\frakD\cap\oldTheta=\partial\frakD$, such that no discal
$2$-suborbifold of $\oldTheta$ has the same boundary as  $\frakD$,
then there is an orientable discal $2$-suborbifold $\frakD'$ of
$\oldPsi$ with $\frakD'\cap\oldTheta=\partial\frakD'$, such that no
discal $2$-suborbifold of $\oldTheta$ has the same boundary as
$\frakD'$. If $\frakD$ is orientable, we need only set
$\frakD'=\frakD$. Now suppose that $\frakD$ is non-orientable. Then
$\frakD\cap\oldTheta=\partial\frakD$ is homeomorphic to
$[[0,1]]$. After modifying $\frakD$ by a small isotopy we may assume
that $\frakD\cap\partial\oldPsi=\partial\frakD$. Since $\oldTheta$ is
a suborbifold of $\partial\oldPsi$, we have
$\fraks_{\partial\frakD}\subset\inter\oldTheta$. Hence after another
small isotopy we may assume that
$\partial\frakD\subset\inter\oldTheta$. 
% an arcThen  take and if there is no discal $2$-suborbifold of $\oldTheta$ having the same boundary as  $\frakD$,  $\oldTheta$ with $\partial\frakE=\partial\frakD$. For every orientable discal $2$-suborbifold $\frakD$ of $\oldPsi$ such that $\frakD\cap\oldTheta=\partial\frakD$, there is a discal $2$-suborbifold $\frakE$ of $\oldTheta$ with $\partial\frakE=\partial\frakD$. that  It is trivial that (b) implies (c). To prove the converse (in the case
%$\oldTheta\subset\partial \oldPsi$), suppose that (c) holds, and that
%$\frakD$ is a discal $2$-suborbifold of $\oldPsi$ such that
%$\frakD\cap\oldTheta=\partial\frakD$. We must show that there %is a
%discal $2$-suborbifold $\frakE$ of $\oldTheta$ with
%$\partial\frakE=\partial\frakD$. Since $\frakD$ is discal, %$|\frakD|$
%is a disk, and $\fraks_\frakD$ is either (1) the empty set, (2) a
%single point of $\inter |\frakD|$, or (3) an arc
%$A\subset\partial|\frakD|$ whose interior is an order-$2$ stratum of
%$\fraks_\frakD$. If (1) or (2) holds, or if (3) holds and the stratum of $\fraks_\oldPsi$ containing $\inter A$ is two-dimensional, then $|\frakD|$ is in general position with respect to  $\fraks_\oldPsi$, and so the required conclusion follows from (c). 
%There remains the subcase in which (3) holds and  $\inter A$ is itself a stratum of $\fraks_\oldPsi$. In this subcase, since $\oldTheta$ is a suborbifold of $\partial\oldPsi$, we have $\partial A\subset\inter|\oldTheta|$.  We may therefore assume, after modifying $\frakD$ by a small (orbifold-)isotopy, that $\partial\frakD=\frakD\cap\partial\oldPsi\subset\inter\oldTheta$.
Let $\frakR$ denote a strong regular neighborhood  of $\frakD$ in $\oldPsi$, which we
%, and we note that $\Fr_\oldPsi\frakR$  is in general position with respect to  $\fraks_\oldPsi$.
%After possibly modifying $\frakD$ by a small (orbifold) isotopy, we may assume that
%$\frakD\cap\partial\oldPsi=\partial\frak D\subset\inter\oldTheta$. If $\frakR$ is a strong regular neighborhood of $\frakD$ in $\oldPsi$, then $\Fr_\oldPsi\frakR$  is in general position with respect to  $\fraks_\oldPsi$. 
 choose to be small enough so that
% the strong regular neighborhood  $\frakR$  of $\frakD$ to be sufficiently small, we can guarantee that 
$\frakB:=\frakR\cap\partial\oldPsi\subset\oldTheta$. Then $\frakB$ is annular and $|\frakB|$ is a disk. Since $\frakD$ is discal, the $3$-orbifold $\frakR$ is discal, and the orientability of $\oldPsi$ implies that $|\partial\frakR|$ is a $2$-sphere. Set $\frakD'=\Fr\frakR$. Since $|\partial\frakR|$ is a sphere and $|\frakB|$ is a disk, $|\frakD'|=|\partial\frakR|-\inter|\frakB|$ is connected. The orientability of $\oldPsi$ implies that $\frakD'$ is orientable. We have $\chi(\frakD')=\chi(\partial\frakR)-\chi(\frakB)=\chi(\partial\frakR)=\chi(\frakR)/2>0$, so that $\frakD'$ is discal. 

It remains only to show that no discal $2$-suborbifold of $\oldTheta$ has the same boundary as  $\frakD'$. Suppose
 that $\partial\frakD'$ does bounds a discal suborbifold $\frakE'$ of $\oldTheta$. Since $\frakB$ is annular, we have $\frakE'\ne\frakB$. Hence
% suborbifold of $\oldTheta$, whose boundary is $\partial\frakD'=\partial\frakE'$.  We must therefore have 
$\frakB\cap\frakE'=\partial\frakB=\partial\frakE'=\partial\frakD'$, and $\frakB\cup\frakE'$ is a component of $\oldTheta$. We have $\chi(\frakB\cup\frakE')=\chi(\frakE')>0$, a contradiction to the hypothesis that $\oldTheta$ has no spherical components. This completes the proof of the proposition in the case $\oldTheta\subset\partial \oldPsi$.

Now consider the case in which $\oldTheta$ is   two-sided. Set
$\oldPsi'=\oldPsi\cut\oldTheta$, and
$\oldTheta'=\rho_\oldTheta^{-1}(\oldTheta)\subset\partial\oldPsi'$. It
follows from the Seifert-van Kampen theorem for orbifolds (see
\cite[Section 2.2]{bmp}) that $\oldTheta$ is $\pi_1$-injective in
$\oldPsi$ if and only if $\oldTheta'$ is $\pi_1$-injective in
$\oldPsi'$. It is  clear that each of the conditions (a), (b) and (c) for $\oldPsi$ and $\oldTheta$ is equivalent to the same condition with $\oldPsi$ and $\oldTheta$ replaced by $\oldPsi'$ and $\oldTheta'$. Furthermore, if 
$\oldPsi$ is orientable and that $\oldTheta$ has no spherical component, then $\oldPsi'$ is orientable and $\oldTheta'$ has no spherical component.  Hence the assertion of the proposition in this case follows by applying the case already proved, with $\oldPsi'$ and $\oldTheta'$ playing the respective roles of $\oldPsi$ and $\oldTheta$.
\EndProof
%%\tau $%g\obd \frakR (1)$D\frakB

\Corollary\label{injective hamentash}
Let $\oldPsi$ be a very good, orientable $3$-orbifold, and let $S\subset|\oldPsi|$ be a $2$-sphere of weight $3$, in general position with respect to $\fraks_\oldPsi$, such that $\chi(\obd(S)\le0$. Then $\obd(S)$ is $\pi_1$-injective in $\oldPsi$.
\EndCorollary

\Proof
We may assume after a small non-ambient isotopy that $S\subset\inter|\oldPsi|$. Set $\oldTheta=\obd(S)$. Since $\oldPsi$ and the closed surface $S\subset\inter|\oldPsi|$ are orientable, $\oldTheta$ is two-sided in $\oldPsi$.
According to Proposition \ref{kinda dumb} (more specifically the implication (c)$\Rightarrow$(a)), it suffices to show that if $\frakD$ is an orientable discal $2$-suborbifold of $\oldPsi$ such that $\frakD\cap\oldTheta=\partial\frak D$, then $\partial\frakD$ is the boundary of a discal $2$-suborbifold of $\oldTheta$. Since $\frakD$ is orientable,  $\partial\frakD$ is a closed $1$-submanifold of $\oldTheta$. But since $|\oldTheta|$ is a $2$-sphere and $\fraks_\oldTheta$ has cardinality $3$, every closed $1$-submanifold of $\oldTheta$ bounds a discal suborbifold of $\oldTheta$. 
\EndProof

\NotationRemark\label{lambda thing}
If $\oldPsi$ is a compact, orientable $3$-orbifold, we will define an
integer $\lambda_\oldPsi$ by setting $\lambda_\oldPsi=2$ if every
component of $\fraks_\oldPsi$ is an arc or a simple closed curve, and
$\lambda_\oldPsi=1$ otherwise.

Note that if $\oldOmega$ is a closed, orientable $3$-orbifold, and $\cals\subset|\oldOmega|$ is a $2$-manifold which is in general position (see \ref{gen pos}) with respect to $\fraks_\oldOmega$ (so that $\obd(\cals)\subset\oldOmega$ is defined by \ref{obd}), then we have $\lambda_{\oldOmega\cut{\obd(\cals)}}=\lambda_\oldOmega$. It follows that for every component $\oldPsi$ of $\oldOmega\cut{\obd(\cals)}$ we have $\lambda_\oldPsi\ge\lambda_\oldOmega$.
\EndNotationRemark

\Definition\label{parallel def}
Let $\oldPsi$ be an orientable $3$-orbifold, and let $\oldXi$ be a
$2$-suborbifold of $\oldPsi$. Let $\frakB,\frakB'\subset\oldPsi$ be
$2$-suborbifolds whose boundaries are contained in $\oldXi$. We will say
that $\frakB$ and $\frakB'$ are {\it parallel in the pair $(\oldPsi,\oldXi)$}
(or simply {\it parallel in $\oldPsi$} in the case where
$\oldXi=\partial\oldPsi$, or  {\it parallel} when
$\oldXi=\partial\oldPsi$ and it is understood which orbifold $\oldPsi$
is involved) if there is
an embedding $j:\frakB\times[0,1]\to\oldPsi$ 
such that
$j(\frakB\times\{0\})=\frakB$, 
$j(\frakB\times\{1\})=\frakB'$,
and
$(\partial\frakB)\times[0,1]\subset j^{-1}(\oldXi)\subset
\partial(\frakB\times[0,1])$.

In all cases where these terms are used, $\oldXi$ will be either a suborbifold of $\partial\oldPsi$ or a two-sided suborbifold of $\oldPsi$. Note that in the case where $\oldXi\subset\partial\oldPsi$, the condition $(\partial\frakB)\times[0,1]\subset j^{-1}(\oldXi)\subset
\partial(\frakB\times[0,1])$ which appears in the definition is equivalent to the simper condition $j((\partial\frakB)\times[0,1])\subset\oldXi$. 
\EndDefinition

\DefinitionsRemarks\label{oops}
A $3$-orbifold $\oldPsi$ will be termed {\it weakly \simple} if it is very good,  compact,
irreducible, and non-discal, and every $\pi_1$-injective two-sided toric suborbifold
of $\inter\oldPsi$ is parallel in $\oldPsi$  to a component of $\partial\oldPsi$. An $3$-orbifold $\oldPsi$ will be termed {\it strongly
  \simple} if (I) $\oldPsi$ is  compact and non-discal, (II) some finite-sheeted regular covering of $\oldPsi$ is an irreducible $3$-manifold, and (III) $\pi_1(\oldPsi)$ has no free abelian subgroup of rank $2$.

Note that if $\oldPsi$ is orientable and strongly \simple, Condition (II) in the definition implies that $\oldPsi$ is very good. According to Lemma \ref{cover-irreducible}, Condition (II) also implies that $\oldPsi$ is irreducible; since $\oldPsi$ is non-discal, it follows that no component of $\partial\oldPsi$ is spherical. Note also that Condition (III) implies that $\oldPsi$ contains no $\pi_1$-injective toric suborbifold whatever; hence a strongly \simple\ orientable $3$-orbifold is weakly \simple.

If $\Mh$ is a closed, orientable hyperbolic $3$-orbifold, then $(\Mh)\pl$ is strongly \simple. This fact will often be used implicitly.
\EndDefinitionsRemarks

\Proposition\label{almost obvious}
Let $\oldPsi$ be a strongly \simple, orientable  $3$-orbifold 
containing
no embedded negative turnovers (see \ref{wuzza turnover}). Suppose that $\lambda_\oldPsi=2$. Then $|\oldPsi|$ contains no weight-$3$
sphere which is in general position with respect to $\fraks_\oldPsi$.
\EndProposition

\Proof
Suppose that $S\subset|\oldPsi|$ is a weight-$3$ sphere in general position with respect to $\fraks_\oldPsi$. First suppose that $\chi(\obd(S))\le0$. By
Corollary \ref{injective hamentash}, $\omega(S)$ is $\pi_1$-injective in $\oldPsi$. Since $\oldPsi$ and the closed surface $S\subset\inter|\oldPsi|$ are orientable, $\oldTheta$ is two-sided in $\oldPsi$. (The hypothesis in Corollary \ref{injective hamentash} that $\oldPsi$ is very good follows from the strong \simple ity of $\oldPsi$.)  If $\chi(\obd(S))=0$ then $\obd(S)$ is a $\pi_1$-injective two-sided toric suborbifold of $\oldPsi$; as observed in \ref{oops}, this contradicts the strong \simple ity of $\oldPsi$. If $\chi(\obd(S))<0$ then $\obd(S)$ is an embedded negative turnover, a contradiction to the hypothesis. 

There remains the possibility that $\chi(\obd(S)>0$, so that $\obd(S)$ is spherical. Since a strongly \simple\ orientable $3$-orbifold is irreducible (see \ref{oops}), $\obd(S)$ is the boundary of a discal $3$-suborbifold $\frakB$ of $\oldPsi$. Since $\lambda_\oldPsi=2$, every component of $|\frakB|\cap\fraks_\oldPsi$ is a properly embedded arc or simple closed curve in $|\frakB|$. If $m$ denotes the number of arc components of $|\frakB|\cap\fraks_\oldPsi$, we have $\wt S=2m$, a contradiction since $\wt S=3$.
\EndProof

\begin{definitionremark}\label{what, no soap?}
Generalizing the definition given for manifolds in \ref{great day}, we define a $3$-orbifold $\oldPsi$ to be {\it
boundary-irreducible} if $\partial\oldPsi$ is $\pi_1$-injective in $\oldPsi$. If $\oldPsi$ is very good, it follows from Proposition \ref{kinda dumb} that $\oldPsi$ is boundary-irreducible if and only if for every  discal $2$-orbifold  $\oldDelta\subset\oldPsi$ with $\oldDelta\cap\partial\oldPsi=\partial\oldDelta$, there is a discal $2$-suborbifold of $\partial\oldPsi$ having the same boundary as $\oldDelta$. If $\oldPsi$ is orientable, has no spherical boundary components, and is very good (as is the case if $\oldPsi$ is strongly \simple, cf. \ref{oops}), it follows from Proposition \ref{kinda dumb} that $\oldPsi$ is boundary-irreducible if and only if for every  orientable discal $2$-orbifold  $\oldDelta\subset\oldPsi$ with $\oldDelta\cap\partial\oldPsi=\partial\oldDelta$, there is a discal $2$-suborbifold of $\partial\oldPsi$ having the same boundary as $\oldDelta$.

\nonessentialproofreadingnote{In the presence of
irreducibility, this implies that every  properly embedded discal suborbifold is boundary-parallel.  In the one place where that fact has come up, I've just given the argument (it's
% been used, and whether it should be stated here are quoted, or just pointed out as  its comes up. (Added 7/19/17: It just occurred to me that in the one place I remember using this, perhaps 
in the proof of Prop \ref{when vertical}.)
%, I used the fact that I was in an $I$-fibred gadget to deduce that a discal $3$-suborbifold is a ball with at most an unknotted arc as singular set; this seems unnecessary because the sphere bounding the ball contains at most two singular points by construction.) In either case it would be good to have cross-refs. to this subsection. 
A related point is that I had been talking as if every properly embedded discal suborbifold of an orientable $3$-orbifold were two-sided. This is of course wrong, because of the example in which $\partial|\oldDelta|$ is made up of an arc in $\partial|\oldPsi|$ and an arc in $\fraks_\oldPsi$. The corresponding statement for annular suborbifolds is also false, because a M\"obius band (which is a manifold and in particular an orbifold) can be properly embedded in an orientable $3$-manifold. That misconception was the reason why I didn't have the term `two-sided'' in the statement just made. At the moment I think this note can probably be ignored...}
\end{definitionremark}

\Number\label{boundary is negative}
Let $\oldPsi$ be a compact, orientable $3$-orbifold. If $\oldPsi$ is
weakly \simple, then in particular it is irreducible and non-discal;
this implies that every boundary component of $\oldPsi$ has
non-positive Euler characteristic.

If $\oldPsi$ is strongly \simple\ and boundary-irreducible, then in fact every boundary
component of $\oldPsi$ has strictly negative Euler
characteristic. This is because any boundary component of $\oldPsi$
having Euler characteristic $0$ would be toric by definition and would
be $\pi_1$-injective by boundary irreducibility; according to an observation in \ref{oops}, this would contradict
strong \simple ity.
\EndNumber

\Definition
An orbifold $\oldPsi$ will be termed {\it componentwise irreducible},
or {\it componentwise boundary-irreducible}, if each of its components
is, respectively, irreducible or boundary-irreducible. We will say
that $\oldPsi$ is {\it componentwise weakly (or strongly) \simple} if
$\oldPsi$ is compact and each of its components is weakly (or,
respectively, strongly) \simple.
\EndDefinition

\Lemma\label{oops lemma} 
Let $\oldPsi$ be a compact orientable $3$-orbifold, and let $\oldTheta$ be a two-sided $2$-suborbifold of $\oldPsi$ whose components  have non-positive Euler characteristic and are $\pi_1$-injective. 
\begin{itemize}
\item If $\oldPsi$ is componentwise irreducible, or componentwise strongly \simple, then $\oldPsi\cut\oldTheta$ is, respectively, componentwise irreducible or componentwise strongly \simple.
\item If $\oldTheta$ is closed and $\oldPsi$ is componentwise boundary-irreducible, then
$\oldPsi\cut\oldTheta$ is componentwise
boundary-irreducible. 
\end{itemize}
\EndLemma

\Proof
We may assume that $\oldTheta$ is connected, since the general case will follow from the connected case by induction on the number of components of $\oldTheta$. 

If $\oldPsi_0$ denotes the component of $\oldPsi$ containing $\oldTheta$, every component of $\oldPsi\cut\oldTheta$ is either a component of $(\oldPsi_0)\cut\oldTheta$, or a component of $\oldPsi$ distinct from $\oldPsi_0$; hence we may assume that $\oldPsi$ is connected.

First suppose that $\oldPsi$ is  irreducible. To prove that $\oldPsi\cut\oldTheta$ is componentwise irreducible, it suffices to prove that every two-sided spherical $2$-suborbifold $\frakV$ of $\oldPsi-\oldTheta$ bounds a discal $3$-suborbifold of $\oldPsi-\oldTheta$. Since $\oldPsi$ is  irreducible, $\frakV$ bounds a discal $3$-suborbifold $\frakB$ of $\oldPsi$. Since $\frakV\cap\oldTheta=\emptyset$, we must have either $\frakB\subset\oldPsi-\oldTheta$ or $\frakB\supset\oldTheta$. If the latter alternative holds, the $\pi_1$-injectivity of $\oldTheta$ implies that $\pi_1(\oldTheta)$ is isomorphic to a subgroup of $\pi_1(\frakB)$; this is impossible, because the hypothesis $\chi(\oldTheta)\le0$ implies that $\pi_1(\oldTheta)$ is infinite,
% and that $\oldTheta$ is $\pi_1$-injective in $\oldPsi$ and hence in $\frakB$, 
whereas the discality of $\frakB$ implies that $\pi_1(\frakB)$ is finite. Thus $\oldPsi\cut\oldTheta$ is indeed componentwise irreducible.

Now suppose that $\oldPsi$ is  strongly \simple. 
Since $\oldTheta$ is $\pi_1$-injective in $\oldPsi$, the immersion $\rho_\oldTheta$ is $\pi_1$-injective; since the  strong \simple ity of $\oldPsi$ implies that $\pi_1(\oldPsi)$ has no rank-$2$ free abelian subgroup (see Definition \ref{oops}), it follows that no component of $\oldPsi\cut\oldTheta$ has a fundamental group with a rank-$2$ free abelian subgroup, as required by Condition (III) of Definition \ref{oops}. 

Next note that the  strong \simple ity of $\oldPsi$ gives a finite-sheeted regular covering $p:\toldTheta\to\oldTheta$ such that $\oldPsi$ is an irreducible $3$-manifold. Note that $\toldTheta:=p^{-1}(\oldTheta)$ is two-sided and $\pi_1$-injective in $\toldPsi$. Applying the assertion proved above, with $\toldPsi$ and $\toldTheta$ playing the respective roles of $\oldPsi$ and $\oldTheta$, we deduce that $\toldPsi\cut{\toldTheta}$ is componentwise irreducible. Every component of $\oldPsi\cut\oldTheta$ admits some component of $\toldPsi\cut{\toldTheta}$ as a regular covering, and thus satisfies Condition (II) of Definition \ref{oops}. 

To show that $\oldPsi\cut\oldTheta$ is componentwise strongly \simple, it remains only to prove that it has no discal component. If $\frakB$ is any component of $\oldPsi\cut\oldTheta$, then $\partial\frakB$ has a $\pi_1$-injective $2$-suborbifold $\frakV$ homeomorphic to some component of $\oldTheta$. The hypothesis gives $\chi(\frakV)\le0$, so that $\pi_1(\frakV)$ is infinite. Hence $\pi_1(\frakB)$ is infinite, and $\frakB$ cannot be discal.

Finally, suppose that $\oldTheta$ is closed and that $\oldPsi$ is  boundary-irreducible. To prove that $\oldPsi\cut\oldTheta$ is componentwise boundary-irreducible, we must show that $\partial(\oldPsi\cut\oldTheta)$ is $\pi_1$-injective in $\oldPsi\cut\oldTheta$. But since $\oldTheta$ is closed, each
component of $\partial(\oldPsi\cut\oldTheta)$ is either a component of $\partial\oldPsi$, which by the  boundary-irreducibility of $\oldPsi$ is $\pi_1$-injective in $\oldPsi$ and hence in $\partial(\oldPsi\cut\oldTheta)$, or a component of $\rho_\oldTheta^{-1}(\oldTheta)$. A component of the latter type is also $\pi_1$-injective since $\oldTheta$ is $\pi_1$-injective in $\oldPsi$.
\EndProof
%\obd

%\Lemma\label{lambda lemma}
%Let $\oldPsi$ be a compact, orientable $3$-orbifold, and let $K$ be a
%submanifold of $N:=|\oldPsi|$ whose frontier is a closed $2$-submanifold
%of $\inter N$, transverse to
%$\fraks_\oldPsi$. Then \redcomment{I don't see how I can conclude that
  %$\wt\Fr_N K$ is divisible by $\lambda_\oldPsi$ in this situation,
  %although it would be trivial if $K$ were contained in $\inter
 % N$. The problem is that only the weaker condition holds in the very
  %first app (proof of Lemma \ref{modification}.}
%\EndLemma

%\Proof
%\redproofsummary{This is really trivial under the assumption
%  $K\subset\inter N$, but I don't see how it can be true without
  %that. The solution will be to change the statement to a variant of the trivial one, in which we have the boundary in place of the frontier, and adapt the application by introducing $\wt^*$ later. Sort this out.}
%\EndProof

\Definitions\label{acylindrical def}
Let $\oldPsi$ be an orientable, componentwise irreducible $3$-orbifold, and let $\oldXi$ be a
$\pi_1$-injective $2$-suborbifold of $\partial\oldPsi$. 
An annular $2$-orbifold $\oldPi\subset\oldPsi$ 
will be called {\it
  essential in the pair $(\oldPsi,\oldXi)$} (or simply   {\it
  essential (in $\oldPsi$)} in the case where $\oldPsi$ is componentwise
  boundary-irreducible and $\oldXi=\partial\oldPsi$) if (1)
$\oldPi$ is two-sided in $\oldPsi$ and has its boundary contained in $\oldXi$, (2) $\oldPi$ is
$\pi_1$-injective in $\oldPsi$, and (3) $\oldPi$ is {\it not} parallel
 in the pair $(\oldPsi,\oldXi)$ either to a suborbifold of $\oldXi$ or to a component of
$\overline{\partial\oldPsi-\oldXi}$. 
We define an {\it acylindrical pair} to be an ordered pair 
$(\oldPsi,\oldXi)$, where $\oldPsi$ is an orientable $3$-orbifold, $\oldXi$ is a
$\pi_1$-injective $2$-suborbifold of $\partial\oldPsi$, and
$\oldPsi$ contains no annular $2$-orbifold which is    essential in
the pair
$(\oldPsi,\oldXi)$. This definition generalizes the definition of an acylindrical $3$-manifold given in \ref{great day} in the sense that, given a an orientable, irreducible, boundary-irreducible $3$-manifold $M$ which is not a ball, $M$ is acylindrical if and only if the manifold pair $(M,\partial M)$, regarded as an orbifold pair, is acylindrical.
%\redproofreadingnote{I have replaced ``strongly \simple'' in that def. by
%  ``irreducible.'' This is consistent with the decision that
  %acylindrical pairs should not be def'd to be strongly \simple. Check
  %that everything makes sense.
%\redmissingref{
%Now I'm thinking the main reason for including \simple ity in the def. of an acylindrical pair may have to do with the proof of Prop. \ref{new characteristic}, specifically Claim \ref{new-claim}. And I'm not sure that was done correctly! Certainly I never argued that the acylindrical pieces were not covered by Seifert fibered manifolds, for example; it wasn't even in the def. when I wrote that. I'm guessing that the proof will go through without putting \simple ity in the def. of an acylindrical pair. The specific fact that these pieces contain no $\pi_1$-injective tori should follow from the fact that the {\it ambient} orbifold is \simple. Check it. If T I'm right, a nice by-product will be that the corrected def. of acylindrical pair really will generalize the def. of an acylindrical mfd.}
%I have done a lot of revisions to get rid of the presence of strong \simple ity in the def. of an acylindrical pair. This comes up, for example, in Lemma \ref{butthurt} and its app. I should be especially careful to make sure everything is consistent.
%}
\EndDefinitions
%\oldPsi\oldLambda

\Proposition\label{what essential means}
Let $\oldPsi$ be a very good, orientable, componentwise irreducible $3$-orbifold, let $\oldXi$ be a
$\pi_1$-injective, non-spherical $2$-suborbifold of $\partial\oldPsi$, and suppose that
$\oldPi$ is an   essential  annular $2$-orbifold 
in the pair $(\oldPsi,\oldXi)$. Let $D\subset|\oldPsi|$ be a
weight-$0$ disk such that $\beta:=D\cap|\oldPi|$ is an arc in
$\partial D$, and $D\cap\partial|\oldPsi|= \overline{\partial
  D-\beta}\subset\oldXi$. Then $\beta$ is the frontier in $|\oldPi|$
of a weight-$0$ disk.
\EndProposition

\Proof
Since $\oldPi$ is annular and orientable, $|\oldPi|$ is either a
weight-$0$ annulus, or a weight-$2$ disk such that the two points of
$|\oldPi|\cap\fraks_\oldPsi$ are both of order $2$. Hence if we assume
that $\beta$ is not the frontier in $|\oldPi|$
of a weight-$0$ disk, and if $\frakB$ denotes a strong regular neighborhood
of $\omega(\beta)$ in $\oldPi$, then $|\overline{\oldPi-\frakB}|$ is
either a weight-$0$ disk or a disjoint union of two weight-$1$
disks. Thus $\overline{\oldPi-\frakB}$ has either one or two
components, and each of its components is an orientable discal
$2$-orbifold. Set $R=|\frakB|$ and $
J=\calc(\overline{\oldPi-\frakB})$ (so that $\card J\le2$).

Let $Z$ be a weight-$0$ ball in $|\oldPsi|$ such that (i) $D\subset
Z$, (ii)
$Z\cap|\oldPi|=R$,
% is a regular neighborhood of $\beta$ in
%$|\oldPi|$, 
(iii) $R\subset\partial Z$, (iv)
$Q:=Z\cap\partial|\oldPsi|$ is a disk contained in $|\oldXi|$, and (v)
$Q\cup R$ is a strong regular neighborhood of $\partial\beta$ in
$|\partial\oldPi|$ (and thus consists of two
arcs). Then $Q\cup R$ is an
annulus in $\partial Z$, and hence $\cald:=\overline{\partial
  Z-(Q\cup R)}$ is a disjoint union of two weight-$0$ disks.

%Set $\frakB=\omega(R)$.
For each $\frakV\in  J$, let 
$\oldTheta_\frakV$ denote the union of $\frakV$ with the component or
components of $\obd(\cald)$ that meet $\frakV$. Then $\oldTheta_\frakV$ is a
two-sided $2$-orbifold in $\oldPsi$ whose boundary is
contained in $\oldXi$; since $\frakV$ is discal, and since each
component of $\cald$ is a weight-$0$ disk meeting $|\frakV|$ in an arc
or the empty set, $\oldTheta_\frakV$ is also discal. Since $\oldXi$
is $\pi_1$-injective, it follows from 
Proposition \ref{kinda dumb} (more specifically the implication (a)$\Rightarrow$(b)), applied with $\oldTheta_\frakV$ playing
the role of $\frakD$, that there is a discal $2$-suborbifold
$\frakH_\frakV$ of $\oldXi$ with
$\partial\frakH_\frakV=\partial\oldTheta_\frakV$. 
Then $\oldTheta_\frakV\cup\frakH_\frakV$
is a spherical $2$-orbifold, and is non-ambiently isotopic to a two-sided suborbifold of $\oldPsi$ since $\oldTheta$ is two-sided and $\frakH_\frakV\subset\partial\oldPsi$. Since $\oldPsi$ is irreducible, $\oldTheta_\frakV\cup\frakH_\frakV$
 bounds a discal $3$-orbifold
$\oldGamma_\frakV\subset\oldPsi$. 

lf $\frakV$ and $\frakV'$ are distinct elements of $J$, we have $\oldTheta_\frakV\cap\oldTheta_{\frakV'}=\emptyset$. Hence
$\oldTheta:=\bigcup_{\frakV\in J}\oldTheta_\frakV$ is a $2$-orbifold
whose components are the orbifolds $\oldTheta_\frakV$ for $\frakV\in
J$. By construction we have $\oldTheta\cap\inter Z=\emptyset$, and
hence $Z$ is contained in the closure of some component $U$ of
$|\oldPsi-\oldTheta|$. For each $\frakV\in J$, we have $Z\cap|\oldTheta_\frakV|\supset Z\cap|\frakV|\ne\emptyset$, and hence $|\oldTheta_\frakV|\subset \overline{U}$. It follows that $\oldTheta\subset \Fr_\oldPsi\frakU$, where $\frakU=\obd(\overline{U})$. 
This in turn implies that $|\oldPi|=R\cup\overline{|\oldPi-\frakB|}\subset Z\cup|\oldTheta|\subset|\frakU|$, so that $\oldPi\subset\frakU$.

For each $\frakV\in J$, we have 
$\Fr\oldGamma_\frakV=\oldTheta_\frakV\subset\Fr\frakU$; hence either $\frakU\subset\oldGamma_\frakV$ or $\frakU\cap\oldGamma_\frakV=\oldTheta_\frakV$. If
$\frakU\subset\oldGamma_\frakV$, then in particular
$\oldPi\subset\oldGamma_\frakV$, and since $\oldPi$ is
$\pi_1$-injective in $\oldPsi$ by essentiality, it is in particular $\pi_1$-injective in $\oldGamma_\frakV$; this is impossible, since the annularity of $\oldPi$ implies that $\pi_1(\oldPi)$ is infinite, while the discality of $\oldGamma_\frakV$ implies that $\pi_1(\oldGamma_\frakV)$ is finite. Hence $\frakU\cap\oldGamma_\frakV=\oldTheta_\frakV$ for each $\frakV\in J$. It follows that $\oldGamma_\frakV$ is a component of $\overline{\oldPsi-\frakU}$ for each $\frakV\in J$. Furthermore, if $\frakV$ and $\frakV'$ are distinct elements of $J$, then since the frontiers 
$\oldTheta_\frakV$ and $\oldTheta_{\frakV'}$ of
$\oldGamma_\frakV$ and $\oldGamma_{\frakV'}$ are disjoint, $\oldGamma_\frakV$ and $\oldGamma_{\frakV'}$ are distinct components of
$\overline{\oldPsi-\frakU}$ and are therefore disjoint.
Hence
$\oldGamma:=\bigcup_{\frakV\in J}\oldGamma_\frakV$ is a $2$-orbifold whose components are the orbifolds $\oldGamma_\frakV$ for $\frakV\in J$, and $\Fr\oldGamma=\oldTheta$.

Since $\frakU\cap\oldGamma_\frakV=\oldTheta_\frakV$ for each $\frakV\in J$, we have
$\frakU\cap\oldGamma=\oldTheta$. Since $\frakZ:=\obd(Z)\subset\frakU$, it follows that 
$\frakZ\cap\oldGamma=\frakZ\cap\oldTheta=\obd(\cald)$. Thus if we set $\oldLambda=\frakZ\cup\oldGamma$, we have $\Fr_\oldPsi\oldLambda=\overline {\frakZ-\obd(\cald)}\cup\overline{\oldGamma-\obd(\cald)}=\frakB\cup\overline{\oldPi-\frakB}=\oldPi$.

%nnd final paragraph to fit with the revisions I've made. The regular neighborhood argument is confusing and should be replaced by something more direct. The contradiction from ``$\oldPi$ is
  %  boundary-parallel'' should really come from ``$\oldPi$ is
    %parallel in the pair $(\oldPsi,\oldPi)$ to a suborbifold of $\oldPsi$.''
%A crucial observation will be that $Z\cap\oldGamma_\frakV=\frakV$ for each $\frakV$. 

%Set $W= |\oldPi|\cup Z$.
%For each $\frakV\in J$, since $W$ is
%connected and
%$\Fr_{|\oldPsi|}|\oldGamma_\frakV|=|\oldTheta_\frakV|\subset\Fr_\oldPsi
%W$, we have either $W\subset|\oldGamma_\frakV|$ or
%$W\cap|\oldGamma_\frakV|=|\oldTheta_\frakV|$. If for some $\frakV$ we
%have $W\subset|\oldGamma_\frakV|$, then since $\oldPi\subset\omega(W)$
%is annular and $\pi_1$-injective, $\pi_1(\oldGamma_\frakV)$ must be infinite,
%a contradiction since $\oldGamma_\frakV$ is discal. Hence
%$W\cap|\oldGamma_\frakV|=|\oldTheta_\frakV|$ for each 
%$\frakV\in J$. 

Now for each 
$\frakV\in J$, since the 
%orientable $3$-orbifold $\oldGamma_\frakV$ is
%discal, and the 
$2$-orbifolds $\oldTheta_\frakV$ and
$\overline{\partial\oldGamma_\frakV-\oldTheta_\frakV}=\frakH_\frakV$
are discal, 
$\fraks_{\partial\oldGamma_\frakV}$
consists of at most
two points. Hence the discal $3$-orbifold $\oldGamma_\frakH$ is homeomorphic to the quotient of to the
quotient of $\DD^3$ by an orthogonal action of a (possibly trivial)
finite cyclic  group. It follows that 
the pair
$(\oldGamma_\frakV,\frakH_\frakV)$ is homeomorphic to
$(\frakH_\frakV\times[0,1],\frakH_\frakV\times\{0\})$. If we set $\frakH=\bigcup_{\frakV\in J}\frakH_\frakV=\oldGamma\cap\partial\oldPsi$, it
 follows that
$(\oldGamma,\frakH)$ is homeomorphic  to $(\frakH\times[0,1],\frakH\times\{0\})$.
Since each component of $\cald$ is a disk in $|\oldTheta|$ meeting
$|\partial\oldTheta|=|\partial\frakH|$ in an arc, we deduce that the
triad $(\oldGamma,\frakH,\obd(\cald))$ is homeomorphic
(see \ref{pairs})
 to $(\frakH\times[0,1],\frakH\times\{0\},\obd(\alpha)\times[0,1])$, where $\alpha=\cald\cap|\partial\frakH|\subset Q$. On the other hand, the triad 
$(\frakZ,\obd(Q),\obd(\cald))$ is homeomorphic to $(\obd(Q)\times[0,1],\obd(Q)\times\{0\},\obd(\alpha)\times[0,1])$. Since $\oldLambda=
\frakZ\cup\oldGamma$ and $\frakZ\cap\oldGamma=\obd(\cald)$, it now follows that
the pair $(\oldLambda,\frakC)$, where
$\frakC=\oldLambda\cap\partial\oldPsi=\oldLambda\cap\oldXi$, is homeomorphic to $(\frakC\times[0,1],\frakC\times\{0\})$. Hence $\oldPi=\Fr\oldLambda$ is
    parallel in the pair $(\oldPsi,\oldXi)$ to a suborbifold of $\oldXi$,
a contradiction to the hypothesis that $\oldPi$ is essential. 
\EndProof
%\frakE\frakH

\Number\label{standard torifold}
Let $q\ge1$ be an integer. There is a action of the cyclic group
$\langle x\,|\,x^q=1\rangle$ on the solid torus $\DD^2\times \SSS^1 $
defined by $x\cdot (e^{2\pi i/q}z,w)=(z,w)$, where $\DD^2$ and
$\SSS^1$ are
identified with subsets of the complex plane. This action is at once smooth in the standard smooth structure on $\DD^2\times \SSS^1$, and PL with respect to have the standard product PL structure on $\DD^2\times \SSS^1$ (see \ref{esso}). The quotient of $\DD^2\times \SSS^1 $ by this action, which inherits both a smooth and a PL structure from $\DD^2\times \SSS^1$, will be denoted $\frakJ_q$. Up to a homeomorphism which is at once smooth and PL, we may identify $|\frakJ_q|$ with a solid torus in such a way that $\fraks_{\frakJ_q}$ is empty if $q=1$, and is a core curve of the solid torus $|\frakJ_q|$, having order $q$, if $q>1$. 

There is also an action, which is again at once smooth and PL, of the dihedral group $\langle x,t\,|\,x^q=1,t^2=1,txt=x^{-1}\rangle$ on $\DD^2\times \SSS^1 $, defined by $x\cdot (z,w)=(e^{2\pi i/q}z,w)$ and $t\cdot(z,w)=(\overline{z},\overline{ w})$, where the bars denote complex conjugation.  
We will let $\frakJ'_q$ denote the quotient of $\DD^2\times \SSS^1 $
by this action, which again inherits both a smooth and a PL structure
from $\DD^2\times \SSS^1$. We may regard $\frakJ'_q$ as the quotient
of $\frakJ_q$ by the involution $\tau_q$ induced by
$(z,w)\mapsto(\overline{z},\overline{w})$; furthermore, if 
we
identify $|\frakJ_q|$ with the standard solid torus
$\DD^1\times\SSS^1$ as above, $|\tau_q|$ is itself given by $(z,w)\mapsto(\overline{z},\overline{w})$, where bars denote complex conjugation in $\DD^2\subset\CC$ or $\SSS^1\subset\CC$. 
Furthermore,
there is a PL homeomorphism $h:|\frakJ'_q|\to \DD^3$ such that (i) $h(\fraks_{\frakJ'_q})$ contains two parallel line segments $\ell_q^1,\ell_q^2$ whose endpoints are in $\SSS^2=\partial \DD^3$, (ii) $s_q:=\overline {h(\fraks_{\frakJ'_q})-(\ell_q^1\cup\ell_q^2)}$ is either the empty set or a line segment having one endpoint in $\inter\ell_q^i$ for each $i\in\{1,2\}$, and (iii) for $i=1,2$, each point of $h^{-1}(\ell_q^i\setminus s_q)$ has order $2$.

We will say that an action of a finite group $G$ on $\DD^2\times \SSS^1 $ is {\it standard} if either $G$ is a cyclic group of some order $q\ge1$ and the action is the first one described above, or  $G$ is a dihedral group of order $2q$ for some $q\ge1$ and the action is the second one described above. 
\EndNumber
%\frakJ \tau t_q

\Lemma\label{prepre}
Let $\oldPsi$ be a very good, irreducible, orientable $3$-orbifold, 
%set 
%$N=|\oldPsi|$, 
and let $\oldTheta\subset\inter\oldPsi$  be a two-sided toric
$2$-suborbifold of $\oldPsi$ which is not $\pi_1$-injective.
%, transverse to $\fraks_\oldPsi$, such that
%$\chi(\omega(T))=0$. 
Then either (a) $\oldTheta$ bounds a $3$-dimensional
suborbifold  of $\inter\oldPsi$ which is homeomorphic to the quotient of $\DD^2\times \SSS^1 $ by the standard action of a finite group, 
or (b)
$\oldTheta$ is contained in the interior of a $3$-dimensional discal suborbifold
 of
$\inter \oldPsi$.
\EndLemma

\Proof
Set $N=|\oldPsi|$ and $T=|\oldTheta|$.
Since $\oldTheta$ is connected and non-spherical, and
%$\oldPsi$ is strongly \simple, the  toric $2$-suborbifold
 $\oldTheta$ is not
  $\pi_1$-injective,
% (see \ref{oops}). Hence
it follows from Proposition \ref{kinda dumb} (more specifically the implication (c)$\Rightarrow$(a)) that there is an orientable discal $2$-suborbifold $\frakD$ of $\oldPsi$, with $\frakD\cap\oldTheta=\partial\frakD$, such that
%$D:=|\frakD|$ is in general position with respect to $\fraks_\oldPsi$, and 
$\partial\frakD$ does not bound a discal suborbifold of $\oldTheta$. This means that
  $D:=|\frakD|\subset N$ is a
disk
such that (1) $D\cap T=\partial D$, 
(2) $D$ 
is
in general position with respect
to
 $\fraks_\oldPsi$
and
meets it
 in at most one point, and
(3) 
any disk in $T$  bounded by
$\gamma:=\partial D$ must meet $\fraks_\oldPsi$
 in at least two points. 
If $D\cap\fraks_\oldPsi$ consists of a single point, let $q$ denote
the order of this point in $\fraks_\oldPsi$; if $D\cap\fraks_\oldPsi=\emptyset$, set $q=1$. 

We claim:
\Claim\label{well, in that case}
Either (i) $T$ is a torus, $\wt T=0$  (see \ref{wuzza weight}), and  $\gamma$ is a non-separating curve in $T$, or (ii) $T$ is a sphere, $\wt T=4$, each point of $\fraks_\oldPsi$ has order $2$, and  $\gamma$ separates $T$ into weight-$2$ disks.
\EndClaim

To prove \ref{well, in that case}, set $n=\wt\oldTheta$ and $\fraks_\oldTheta=\{x_1,\ldots,x_n\}$, and let $p_i$ denote the order of the singular point $x_i$ for $i=1,\ldots,n$.
Since $\oldTheta$ is toric and orientable, it follows from \ref{2-dim case} that 
%Let $\otheroldLambda$ be a finite-type, orientable $2$-orbifold, let $x_1,\ldots,x_n$ denote %the distinct points of $\fraks_\otheroldLambda$, and let $p_i$ denote the order of the singular %point $x_i$ for $i=1,\ldots,n$. Then we have
\Equation\label{hey, i don't make the rules}
0= \chi(\oldTheta)=\chi(T)-\sum_{i=1}^n(1-1/p_i). 
\EndEquation
In particular we have $\chi(T)\ge0$, so that $T$ is a torus or a sphere. If $T$ is a torus, it follows from (\ref{hey, i don't make the rules}) that $n=0$; Condition (3) above then implies that $\gamma$ is a non-separating curve in $T$, and Alternative (i) of \ref{well, in that case} holds. If $T$ is a sphere, then $\gamma$ separates $T$ into two disks $\Delta_1$ and $\Delta_2$; Condition (3) implies that each $\Delta_i$ has weight at least $2$. Hence $n=\wt(\Delta_1)+\wt(\Delta_2)\ge4$. Since $\chi(T)=2$, and each term $1-1/p_i$ in (\ref{hey, i don't make the rules}) is at least $1/2$, it now follows that $n=4$ and that $1-1/p_i=1/2$, i.e. $p_i=2$, for each $i$. Hence $\wt \Delta_1=\wt \Delta_2=2$, and Alternative (ii) holds. This completes the proof of \ref{well, in that case}.

%. \redcomment{Revise the argument below to conclude, a posteriori, that   $\card(T\cap\fraks_\oldPsi)=4$ and that $\fraks_\oldPsi$ consists of points of order $2$.} In Case (i), it follows from Condition (3) above that $\gamma$ is a non-separating curve in $T$. In Case (ii), it follows from Condition (3) that $\gamma$ separates $T$ into two disks, each of which contains exactly two points of $\fraks_\oldPsi$, each of order $2$. 

Since $D\cap T=\partial D=\gamma$, there is a ball $E$, with $D\subset E\subset N$, such that
$E\cap T$ is an annulus $A$, which is contained in $\partial E$ and is a 
regular neighborhood of $\gamma $ in $T$. 
We may choose $E$ in such a way that there exists a homeomorphism $ h:D\times[-1,1]\to E$ such that $ h(D\times\{0\})=D$, $ h((\partial D)\times[-1,1])=A$, and $ h((D\cap\fraks_\oldPsi)\times[-1,1])=E\cap\fraks_\oldPsi$. Set $D_i= h(D\times\{i\})\subset\partial E$ for $i=\pm1$. 
Set 
%$K=\rho(\tK)$ and 
$F=\overline{T-A}\cup\overline{\partial E-A}$. Since $\oldTheta=\obd(T)$ is two-sided, so is $\obd(F)$.
%\rho(\tF)=\partial K$.

If Alternative (i) of \ref{well, in that case} holds, $F$ is a
  $2$-sphere.
%, and hence so is $F$. 
Furthermore, in this case
  $\fraks_{\obd(F)}=F\cap\fraks_\oldPsi$ is empty if $q=1$, and consists of two points of order $q$ if $q>1$; furthermore, in the latter case, one point of $\fraks_{\obd(F)}$ is contained in $D_{1}$, and one in $D_{-1}$.

If Alternative (ii) of \ref{well, in that case} holds,
%In Case (ii), 
$F$ is a disjoint union of two
$2$-spheres $S_1\supset D_1$ and $S_{-1}\supset D_{-1}$.
%, and hence so is $F$. 
In this case, for $i=\pm1$, the intersection $S_i\cap\fraks_\oldPsi$ consists of two points of order $2$ in $S_i-D_i$, if $q=1$; and if $q>1$ it consists of three points, two of which lie in $S_i-D_i$ and have order $2$, and one of which lies in $\inter D_i$ and has order $q$.

Hence in any event, for each component $S$ of $F$, the two-sided $2$-orbifold $\obd(S)$ is spherical. Since $\oldPsi$ is irreducible, $\omega(S)$ bounds a $3$-dimensional discal suborbifold of $\oldPsi$. For each component $S$ of $F$ we choose a discal  $3$-suborbifold $\oldUpsilon_S$ of $\oldPsi$ with $\partial\oldUpsilon_S=\omega(S)$. In particular, for each component $S$ of $F$, the manifold $|\oldUpsilon_S|$ is a $3$-ball, and the group
%  $\fraks_{\obd(F)}=F\cap\fraks_\oldPsi$ either is empty or consists of two points. Since
  %the hyperbolic $3$-orbifold $\oldPsi$ is in particular irreducible,
  %$F$ must bound a ball $B\subset N$ such
  %that $B\cap\fraks_\oldPsi$ is either the empty set or an unknotted arc
 $\pi_1(\oldUpsilon_S)$ is finite. 
The closures of the components of $N-S$ are $|\oldUpsilon_S|$ and
$\overline{N-|\oldUpsilon_S|}$. 

For each component $S$ of $F$, set $Y_S={(T\cup E)-S}$.
The construction gives that $Y_S=A\cup(\inter E)\cup(F-S)$; and that $F-S$ is empty if Alternative (i) of \ref{well, in that case} holds, and is a component of $F$ if Alternative (ii) of \ref{well, in that case} holds. In either case it follows that $Y_S$ is connected. Since $Y_S$ is disjoint from $S$, it must be
%Since $T\cap S$ is a non-separating annulus
%in $T$ in Case (i) and is  a disk in Case (ii),
%$F$ is both the boundary of $|\oldUpsilon_S|$ and a boundary component of $K$, 
%$(\partial|\oldUpsilon_S|)\cap T=S\cap T\subset F\cap T=
%$T\setminus S$ is connected, and is therefore 
contained either in $|\oldUpsilon_S|$ or in
$\overline{N-|\oldUpsilon_S|}$.  Hence $T\cup E=Y_S\cup S$ is contained in either $|\oldUpsilon_S|$
or in
$\overline{N-|\oldUpsilon_S|}$. Since one of these alternatives holds for every component of $S$, we have either
$T\cup E\subset\overline{N-|\oldUpsilon_S|}$ for every component $S$ of $F$, or
$T\cup E\subset|\oldUpsilon_S|$ for some component $S$ of $F$; that is:
\Claim\label{la nouvelle monique}
Either (A) $(T\cup E)\cap|\oldUpsilon_S|=S$ for every component $S$ of $F$, or (B) 
$T\cup E\subset|\oldUpsilon_{S_0}|$ for some component ${S_0}$ of $F$.
\EndClaim

Consider the case in which Alternative (B) of \ref{la nouvelle monique}
holds. Define $\frakB$ to be a strong regular neighborhood of
$\oldUpsilon_{S_0}$ in $\oldPsi$.
%, and set $B=|\frakB|$.
%chosen so that $B\cap\fraks_\oldPsi$ is a regular neighborhood of
%$\fraks_\oldUpsilon=|\oldUpsilon|\cap\fraks_\oldPsi$ in
%$\fraks_\oldPsi$. \redmissingref{It should say that $\obd(B)$ (if I
%call it that) is a (strong or weak?) regular neighborhood of
%$\oldUpsilon$ in $\oldPsi$, which is probably stronger than what I've
%said and will be more economical when this is cleaned up.} 
Then $\frakB$ is (orbifold-)homeomorphic to $\oldUpsilon_{S_0}$ and
is therefore discal. We have $\oldTheta\subset\oldUpsilon_{S_0}\subset\inter
\frakB$. This gives Alternative (b) of the conclusion of the
proposition.

The rest of the proof will be devoted to the case in which Alternative (A) of \ref{la nouvelle monique} holds; in this case, we will show that Alternative (a) of the conclusion of the proposition holds. 

Set $Z=\bigcup_{S\in\calc(F)}|\oldUpsilon_S|$ and
%\Equation\label{nouveau doucement}
$J=E\cup Z$.
It follows from Alternative (A) of \ref{la nouvelle monique} that $J$ is obtained, up to homeomorphism, from $T\cup E$ by gluing a ball to each component of $F$ along the boundary of the ball. Hence $J$ is a manifold whose boundary is $T$. 
In view of the existence of the homeomorphism $ h:D\times[-1,1]\to E$  with the properties stated above, we can now deduce that $\frakJ:=\obd(J)$ is homeomorphic to an orbifold obtained from $\frakZ:=\obd(Z)$ by gluing together the suborbifolds $\oldDelta_1:=\omega(D_{1})$ and $\oldDelta_{-1}:=\omega(D_{-1})$ of $\partial\oldUpsilon$ by a homeomorphism. It now suffices to prove that $\frakJ:=\omega(J)$ is (orbifold-)homeomorphic to  the quotient of $\DD^2\times \SSS^1 $  by the standard action of a finite group.

Consider first the subcase in which Alternative (i) of \ref{well, in that case} holds. In this subcase, we have seen that $F$ is a single $2$-sphere; and that $\fraks_{\obd(F)}$ is empty if $q=1$, and consists of two points of order $q$ if $q>1$. Furthermore, we have seen that in the latter case, one point of $\fraks_{\obd(F)}$ is contained in $D_1$, and one in $D_{-1}$. Hence $\frakZ$ is a single discal $3$-orbifold; $Z$ is a ball; and $\fraks_\frakZ$ is empty if $q=1$, and is a single unknotted arc of order $q$, with one endpoint in $D_1$ and one in $D_{-1}$, if $q>1$. It follows that $\frakJ$ is (orbifold-)homeomorphic to the orbifold $\frakJ_q$ described in \ref{standard torifold}, and is therefore homeomorphic to the quotient of $\DD^2\times \SSS^1 $ by a standard action of a cyclic group of order $q$. Thus Alternative (a) of the conclusion holds.

Now consider the subcase in which Alternative (ii) of \ref{well, in that case} holds.  
In this case we have seen that $F$ is a disjoint union of two
$2$-spheres $S_1\supset D_1$ and $S_{-1}\supset D_{-1}$.
%, and hence so is $F$. 
Furthermore, we have seen that if $q=1$ then for $i=\pm1$ we have $S_i\cap\fraks_\oldPsi=\{x_i,x_i'\}$ for some points $x_i,x_i'\in S_i-D_i$ having order $2$; and that if $q>1$ we have $S_i\cap\fraks_\oldPsi=\{x_i,x_i',x_i''\}$ for some points $x_i,x_i'\in S_i-D_i$ having order $2$, and some point $x_i''\in\inter D_i$ having order $q$.
Hence $\frakZ$ is a disjoint union of two discal $3$-orbifolds $\oldUpsilon_1$ and $\oldUpsilon_{-1}$. Here $Y_i:=|\oldUpsilon_i|$ is a ball for $i=\pm1$. If $q=1$ then $\fraks_{\oldUpsilon_i}$ is an unknotted arc in $Z_i$, having order $q$, and having endpoints $x_i$ and $x_i'$. If $q>1$ then $t_i:=\fraks_{\oldUpsilon_i}$ is a cone on $\{x_i,x_i',x_i''\}$; moreover, $t_i$ is unknotted in the sense that it is contained in a properly embedded disk in $Z_i$. The interiors of the arcs joining $x_i$, $x_i'$ and $x_i''$ to the cone point have orders $2$, $2$ and $q$ respectively. 

It follows that $J=|\frakJ|$ is (orbifold-)homeomorphic to the orbifold $\frakJ_q'$ described in \ref{standard torifold}, and is therefore homeomorphic to the quotient of $\DD^2\times \SSS^1 $  by a standard action of a dihedral group of order $2q$. Thus Alternative (a) of the conclusion holds.
\EndProof
%$h \eta\tau t_

\Proposition\label{three-way equivalence}
Let $\oldLambda$ be an orientable $3$-orbifold. The following conditions are mutually equivalent:
\begin{enumerate}
\item $\oldLambda$ is a \torifold;
\item $\oldLambda$ has a toric boundary component and is  strongly \simple;
\item $\oldLambda$ is homeomorphic to the quotient of $\DD^2\times \SSS^1 $ by the standard action of a finite group.
%either (A) $|\oldLambda|$ is a solid torus and $\fraks_\oldLambda$ is either the empty set or a core curve of $|\oldLambda|$, or (B) 
%each point of $\fraks_\oldLambda$ has order $2$, and 
%there is a homeomorphism $h:|\oldLambda|\to \DD^3$ such that (i) $h(\fraks_\oldLambda)$ contains two parallel line segments $\ell_1,\ell_2$ whose endpoints are in $\SSS^2=\partial \DD^3$, (ii) $s:=\overline {h(\fraks_\oldLambda)-(\ell_1\cup\ell_2)}$ is either the empty set or a line segment having one endpoint in $\inter\ell_i$ for each $i\in\{1,2\}$, and (iii) for $i=1,2$, each point of $h^{-1}(\ell_i\setminus s)$ has order $2$.   
\end{enumerate}
\EndProposition

\Proof
It is trivial that (3) implies (1).

If (1) holds then $\oldLambda$ is finitely covered by a solid torus $J$. Some finite-sheeted covering of $J$ is a regular covering of $\oldLambda$, and is a solid torus since $J$ is one; hence we may assume without loss of generality that $J$ is a regular covering of $\oldLambda$. Since $\partial J$ is a single torus, and the inclusion homomorphism $\pi_1(\partial J)\to\pi_1(J)$ has infinite kernel, $\partial\oldLambda$ is a single toric orbifold. To show that (2) holds, it now suffices to prove that $\oldLambda$ is strongly \simple. Since
% $\oldLambda$ is covered by the compact manifold $J$, it is 
the solid torus $J$ is an irreducible $3$-manifold, Condition (II) of Definition \ref{oops} holds with $\oldLambda$ playing the role of $\oldPsi$. Since $J$ is compact, so is $\oldLambda$; and since $\pi_1(J)$ is infinite cyclic, $\pi_1(\oldLambda)$ is infinite, so that $\oldLambda$ is non-discal, and $\pi_1(\oldLambda)$ has no rank-$2$ free abelian subgroup. Hence Conditions (I) and (III) of Definition \ref{oops} hold as well. This
completes the proof that (1) implies (2).

Now suppose that (2) holds, and let $\oldTheta$ be a toric boundary component of $\oldLambda$. Let $\oldLambda'$ denote the $3$-orbifold obtained from the disjoint union  $\oldLambda\discup((\partial\oldLambda)\times[0,1])$ by gluing $\partial\oldLambda\subset\oldLambda$ to $(\partial\oldLambda)\times\{0\}\subset(\partial\oldLambda)\times[0,1]$ via the homeomorphism $x\mapsto(x,1)$. Then $\oldLambda'$ is homeomorphic to $\oldLambda$, and is therefore strongly \simple. In particular $\oldLambda'$ is irreducible, and the toric $2$-orbifold $\oldTheta$ is not $\pi_1$-injective. Thus  the hypotheses of Lemma \ref{prepre} hold with $\oldLambda'$ 
%and $\partial\oldLambda$ 
playing the role of $\oldPsi$, and with $\oldTheta$ chosen as above. Hence either (a) $\partial\oldLambda$ bounds a $3$-dimensional
suborbifold   $\oldLambda''$  of $\inter\oldLambda'$ which is homeomorphic to the quotient of  $\DD^2\times \SSS^1 $  by the standard action of a finite group, 
or (b)
$\oldTheta$ is contained in the interior of a $3$-dimensional discal suborbifold $\frakB$
 of
$\inter N$. If (a) holds, we must have either $\oldLambda''=\oldLambda$ or $\oldLambda''=\oldTheta\times[0,1]$; the latter alternative is impossible, because the solid toric $3$-orbifold $\oldLambda''$ has connected boundary. Hence $\oldLambda''=\oldLambda$, so that (3) holds. If (b) holds, then $\pi_1(\frakB)$ is finite. But the construction of $\oldLambda'$ implies that $\oldLambda$ is  $\pi_1$-injective in $\oldLambda'$, and hence $\pi_1(\oldLambda)$ is finite. But if $\toldLambda$ is a finite-sheeted manifold covering of $\oldLambda$ (which exists by Condition (II) of Definition \ref{oops}), then $\partial\toldLambda$ has a torus component since $\oldTheta$ is toric; it follows that $H_1(\toldLambda;\QQ)\ne0$, which contradicts the finiteness of $\pi_1(\oldLambda)$. Hence (a) cannot occur, and we have shown that (2) implies (3).
\EndProof
%\obd

\Proposition\label{preoccupani}
Let $\oldPsi$ be a componentwise strongly \simple, orientable $3$-orbifold, and let $\oldTheta\subset\inter\oldPsi$  be a two-sided toric
$2$-suborbifold. Then either (a) $\oldTheta$ bounds a \torifold\ contained in $\inter\oldPsi$, or (b)
$\oldTheta$ is contained in the interior of a discal $3$-suborbifold of
$\inter \oldPsi$.\abstractcomment{The rest of this statement has been removed, and can be found in irreducible1.tex.}
\EndProposition

\Proof
We may assume without loss of generality that $\oldPsi$ is connected, and is therefore strongly \simple. In particular $\oldLambda'$ is irreducible, and the two-sided toric $2$-orbifold $\oldTheta$ is not $\pi_1$-injective. The conclusion is now an immediate consequence of Lemma \ref{prepre} and (the trivial part of) Proposition \ref{three-way equivalence}. \EndProof

\Corollary\label{preoccucorollary}
Let $\oldPsi$ be a componentwise strongly \simple, orientable $3$-orbifold, and let $\oldTheta\subset\oldPsi$  be an  orientable toric
$2$-suborbifold. Suppose that the image of the inclusion homomorphism $\pi_1(\oldTheta)\to\pi_1(\inter \oldPsi)$ is infinite. Then $\oldTheta$ bounds a \torifold\ contained in $\oldPsi$.
\EndCorollary

\Proof
After modifying $\oldTheta$ by a small non-ambient orbifold isotopy, we may assume $\oldTheta\subset\inter\oldPsi$. Since $\oldPsi$ and $\oldTheta$ Then $\oldPsi$ are orientable, it now follows that $\oldTheta$ is two-sided; thus $\oldPsi$ and $\oldTheta$ satisfy the hypothesis of Proposition \ref{preoccupani}. But since   the inclusion homomorphism $\pi_1(\oldTheta)\to\pi_1(\oldPsi)$ has infinite image, Alternative (b) of the conclusion of Proposition \ref{preoccupani} cannot hold. Hence Alternative (a) must hold.
\EndProof

\tinymissingref{\tiny In this version I am reworking Proposition \ref{at least a
    sixth}, which I think was screwed up. For the old version, see
  eclipse5.tex. I fixed the statement and app. in eclipse6.tex. Here I
will try to add a lemma to make the proof go smoothly.}

The following technical lemma will be needed in the proof of  Prop. \ref{new enchanted forest}.

\Lemma\label{occupani}
Let $\Mh$ be a closed, orientable hyperbolic $3$-orbifold, and set $\oldOmega=(\Mh)\pl$ and
$M=|\oldOmega|$. Let $W$ be a connected $3$-submanifold of $M$,
transverse to $\fraks_\oldOmega$, and let $T$ be a component of
$\partial W$ with $\chi(\obd(T))=0$. Let $\frakZ$ be a negative, compact, connected $2$-orbifold with connected boundary, and let $f:\frakZ\to\obd(\overline{M-W})$ be an (orbifold) immersion such that $f^{-1}(\partial W)=f^{-1}(T)=\partial\frakZ$. Suppose that $f$, regarded as an (orbifold) immersion of $\frakZ$ in $\oldOmega$, is $\pi_1$-injective.
% Some simple closed
%curve $c\subset T$, disjoint from $\fraks_\oldOmega$, is the boundary of
%a properly embedded connected $2$-manifold $X$ in $ \overline{M-W}$, transverse to  $\fraks_\oldOmega$,
%such that $\obd(X)$ is negative and is $\pi_1$-injective in
%$\oldOmega$. 
Then
$W$  is contained in a
 submanifold $J$
of $M$ such that $\obd(J)$ is a \torifold, and
$\partial J=T$.
\EndLemma

\Proof
Note that $\oldTheta:=\obd(T)$ satisfies the hypotheses of Proposition
\ref{preoccupani}, with $\oldOmega$ and $M$ playing the roles of
$\oldPsi$ and $N$. Hence one of the alternatives (a) or (b) of the conclusion of Proposition \ref{preoccupani} must hold in this context.

First suppose that Alternative (a) holds, i.e. that $T$ bounds a $3$-dimensional submanifold $J$ of $M$ such that $\obd(J)$ is a \torifold.
Since $T$ is a component of $\partial W$, and $W$ is connected, we have either $W\subset J$ or $W\cap J=T$. If $W\subset J$ then the conclusion of the lemma holds. If $W\cap J=T$, then---since a \torifold\ has connected boundary---$J$ is the component of $\overline{M-W}$ containing $T$, and hence $|f|(|\frakZ|)\subset J$. Since $f$, regarded as an immersion of $\frakZ$ in $\oldOmega$, is $\pi_1$-injective, in particular $f:\frakZ\to\obd(J)$ is $\pi_1$-injective; but the \torifold\ $\obd(J)$ has virtually cyclic fundamental group, and hence $\pi_1(\frakZ)$ is virtually cyclic. This is impossible because the $2$-orbifold $\frakZ$ is negative. Thus the subcase $W\cap J=T$ cannot occur.

Now suppose that
Alternative (b) holds; that is, $T$ is contained in the interior a $3$-dimensional submanifold $B$ of $M$, transverse to $\fraks_\oldPsi$, such that $\obd(B)$ is a discal orbifold. Since $\pi_1(\obd(B))$ is finite, the inclusion homomorphism $\pi_1(\obd(T))\to\pi_1(\oldOmega)$ then has finite image. 

On the other hand,  since $\frakZ$ is negative, the $1$-manifold $\partial\frakZ$ is $\pi_1$-injective in the $2$-orbifold $\frakZ$. Since
$f$, regarded as an immersion of $\frakZ$ in $\oldOmega$, is in turn $\pi_1$-injective,
it follows that $f|\partial\frakZ$, regarded as an immersion of $\partial\frakZ$ in $\oldOmega$, is  $\pi_1$-injective.
Since $\pi_1(\partial\frakZ)\cong\ZZ$ and $\partial\frakZ\subset T$, the $\pi_1$-injectivity of
$f|\partial\frakZ:\partial\frakZ\to\oldOmega$ 
% $\obd(c)$  in $\oldOmega$ 
implies that the inclusion homomorphism $\pi_1(\obd(T))\to\pi_1(\oldOmega)$ has infinite image.
This is a contradiction. Thus Alternative (b) cannot occur, and the proof of the lemma is complete.
\EndProof

%U\frakB

\Proposition\label{silver irreducible}
Let $\oldPsi$ be a very good, compact, orientable $3$-orbifold, and let $\oldXi$ be a
%$\pi_1$-injective, 
 $2$-suborbifold of $\partial\oldPsi$ having no spherical component. Then 
$\silv_{\oldXi}\oldPsi$ (see \ref{silvering}) is irreducible if and only if (1)
$\oldPsi$ is irreducible and (2) $\oldXi$ is $\pi_1$-injective in
$\oldPsi$.
\EndProposition

\Proof
Set $\mu=
\mu_{\oldPsi,\oldXi}
:\oldPsi\to
 \silv_\oldXi\oldPsi$ and 
$\mu^*=\mu^*_{\oldPsi,\oldXi}:
\oldPsi-\oldXi\to (\silv_\oldXi\oldPsi)-\mu (\oldXi)$
(see
 \ref{silvering}).

First suppose that (1) and (2) hold. If $\oldGamma$
is a two-sided spherical $2$-suborbifold of
$\silv_{\oldXi}\oldPsi$, then since $\oldGamma$ is two-sided and
closed,
% and is not contained in $\oldXi$ since $\oldXi$ has no
%spherical component. It therefore
it follows from \ref{second point} that we may
write $\oldGamma=\silv\frakD$ for some two-sided  $2$-suborbifold
$\frakD$ of $\oldPsi$ with $\partial\frakD
%=\frakD\cap\partial\oldPsi
\subset\oldXi$.   We have
$\chi(\frakD)=\chi(\oldGamma)>0$, so that $\frakD$ is either spherical
or discal. If $\frakD$ is spherical then $\frakD\subset\inter\oldPsi$;
 since $\oldPsi$ is irreducible according to (1), $\frakD$ bounds a discal
 $3$-suborbifold of $\inter\oldPsi$. 
%we may identify
The discal orbifold bounded by $\frakD$ is mapped 
%via 
by the 
homeomorphism
$\mu^*$
%|(\oldPsi-\oldXi):
%\oldPsi-\oldXi\to (\silv_\oldXi\oldPsi)-\mu(\oldXi)$
% is an
%orbifold 
 onto a
 discal suborbifold of $\silv_{\oldXi}\oldPsi$ bounded by
 $\oldGamma$. 

Now suppose that $\frakD$ is discal. Since $\oldXi$ is $\pi_1$-injective in $\oldPsi$ according to (2), it follows from Proposition \ref{kinda dumb} (more specifically the implication (a)$\Rightarrow$(b))
that there is a discal suborbifold $\oldPhi$ of
$\oldXi$ with $\partial\oldPhi=\partial\frakD$. Thus $\frakV:=\oldPhi\cup\frakD$ is a spherical
suborbifold of $\oldPsi$, non-ambiently isotopic  to a two-sided suborbifold of $\inter\oldPsi$. Since $\oldPsi$ is irreducible, $\frakV$
bounds a discal $3$-suborbifold $\frakH$ of $\oldPsi$. But since
$\oldPhi$ and $\frakD$ are discal and orientable,
$\fraks_\frakH\cap\partial\frakH=\fraks_\frakV=\fraks_\oldPhi\cup\fraks_\frakD$ consists of at most
two points. Hence the discal $3$-orbifold $\frakH$ is homeomorphic to the quotient of  $\DD^3$ by an orthogonal action of a (possibly trivial)
finite cyclic  group. It follows that the pair $(\frakH,\oldPhi)$ is homeomorphic to
$(\oldPhi\times[0,1],\oldPhi\times\{0\})$ (where $|\oldPhi|$ is a disk of
weight at most $1$).
Now
$\frakK:=\silv_\oldPhi\frakH$ is a suborbifold of
$\silv_{\oldXi}\oldPsi$. Our description of the pair $(\frakH,\oldPhi)$
 implies that
%Since $\frakH$, $\oldPhi$ and
%$\frakD=\overline{\partial\frakH-\oldPhi}$ are discal and orientable,
%),  and hence
 $\frakK$ is discal. We have $\partial\frakK=\oldGamma$.
Thus in either case $\oldGamma$ is the
boundary of a discal $3$-suborbifold of $\silv_{\oldXi}\oldPsi$; this proves that $\silv_{\oldXi}\oldPsi$ is irreducible.

Conversely, suppose that $\silv_{\oldXi}\oldPsi$ is irreducible. If
$\oldPi$ is a two-sided spherical $2$-suborbifold of $\inter\oldPsi$, then
$\oldPi$ is mapped by
the homeomorphism
$\mu^*$
%|(\oldPsi-\oldXi):
%\oldPsi-\oldXi\to (\silv_\oldXi\oldPsi)-\mu(\oldXi)$
% is an
%orbifold 
 onto a
% discal suborbifold of $\silv_{\oldXi}\oldPsi$ bounded by
 %$\oldGamma$. 
%may be identified with a
 two-sided spherical $2$-suborbifold of
$\silv_{\oldXi}\oldPsi$, which bounds a discal $3$-suborbifold
$\frakK$ of $\silv_{\oldXi}\oldPsi$ by irreducibility; we have
$\frakK\subset\inter \silv_\oldXi\oldPsi
\subset (\silv_\oldXi\oldPsi)-\mu(\oldXi)$, and $(\mu^*)^{-1}(\frakK)$
% the inverse
%image of 
%we may identify
%under
%the homeomorphism
%$\mu|(\oldPsi-\oldXi)$
is  a discal $3$-suborbifold  of $\oldPsi$ bounded
by $\oldPi$. This establishes (1). To prove (2), according to
Proposition \ref{kinda dumb} (more specifically the implication (c)$\Rightarrow$(a), it suffices to show that if $\frakD$ is
an orientable discal $2$-suborbifold of $\oldPsi$, such that   $\frakC :=\partial
\frakD=\frakD\cap\oldXi$, then $\frakC $ bounds a discal
suborbifold of $\oldXi$. If we are given such a suborbifold $\frakD$,
%the general position property of $|\frakD|$ implies
%$\dim\fraks_\frakD\le\dim\fraks_\oldPsi-1$; but
%$\dim\fraks_\oldPsi\le1$ by the orientability of $\oldPsi$, and hence
%$\fraks_\frakD$ is a finite set. This implies that $\frakD$ is
%orientable. 
we may assume after a small isotopy that $\frakD$ is
two-sided in $\oldPsi$. Then $\silv\frakD$ is a two-sided spherical $2$-suborbifold of 
$\silv_{\oldXi}\oldPsi$ (cf. \ref{first point}), which bounds a discal $3$-suborbifold
$\frakK$ of $\silv_{\oldXi}\oldPsi$ by irreducibility. 
According to \ref{second point},  we may write
$\frakK=
\silv_\oldPhi\frakH$, where  $\frakH$ is some suborbifold of $\oldPsi$
and $\oldPhi:=\oldXi\cap\frakH$ is a suborbifold of $\oldXi$. We have
$\partial\oldPhi=\frakC$, and it suffices to prove that $\oldPhi$ is discal.
%Since
%$\frakK$ is discal, and $|\mu_{\frakH,\oldXi}|:|\frakH|\to|\frakK|$ is
%a homeomorphism by \ref{silvering}, $|\frakH|$ is a $3$-ball; and
%since $\frakD$ is discal, $|\frakD|\subset\partial|\frakH|$ is a
%disk. Hence $|\oldPhi|=\overline{|\frakH|-|\frakD|}$ is connected. 

Set
$\oldPhi'=\mu(\oldPhi)$, so that $\oldPhi$ and $\oldPhi'$ are
homeomorphic. According to \ref{gyrot}, applied with $\frakH$ and
$\oldPhi$ playing
the respective roles of $\oldPsi$ and $\oldXi$, we have that
(I) $\oldPhi'$ is a suborbifold of 
 $\frakK=\silv_\oldXi\frakH$, and (II)  
$|\oldPhi'|$ is the union of the closures of the two-dimensional strata
of  $\frakK$. But since $\frakK$ is discal, it may be identified
homeomorphically with $\DD^3/G$ for some  $G\le\OO(3)$. If $p:\DD^3\to
\frakK$ denotes the orbit map, It follows
from (II) that $p^{-1}(\oldPhi')$ is the union of all two-dimensional fixed
point sets of elements of $G$; each such fixed-point set is the
intersection of $\DD^3$ with a plane in $\RR^3$ containing $0$. But it
follows from (I) that $p^{-1}(\oldPhi')$ is a manifold, and as
$\oldPhi'=\mu(\oldPhi)\supset\mu(\frakC)\ne\emptyset$, the union defining $p^{-1}(\oldPhi')$ must
consist of exactly one term. Hence $p^{-1}(\oldPhi')$ is a disk, so
that $\oldPhi'$ is discal, and hence so is $\oldPhi$. Thus
(2) is established.
\EndProof
%(i)(ii)(a)(b)(1)(2) C \frakV
%\oldPhi\cup\frakD\oldUpsilon\eta
%canonical \eta \frakE (I) (II){second point} tweets

\Proposition\label{silver acylindrical}
Let $\oldLambda$ 
be an orientable, componentwise strongly \simple\ $3$-orbifold, and let $\oldXi$ be a 
$\pi_1$-injective $2$-suborbifold of $\partial\oldLambda$. Then
a two-sided  annular $2$-orbifold $\oldPi\subset\oldLambda$, with $\oldPi\cap\partial\oldLambda=\partial\oldPi\subset \oldXi$,
is
  essential in the pair $(\oldLambda,\oldXi)$ if and only if the  two-sided toric orbifold $\silv\oldPi$, which is canonically identified with a suborbifold of
$\silv_\oldXi\oldLambda$ (see \ref{first point}),
  is incompressible in
$\silv_\oldXi\oldLambda$ and is not parallel in $\silv_\oldXi\oldLambda$ (see \ref{parallel def}) to
a boundary component of
$\silv_\oldXi\oldLambda$. \nonessentialproofreadingnote{ I think
  strong \simple ity is needed only for the {\it second} assertion. Of
  course on the face of it it costs nothing to make the assumption for
  the whole statement, but the question arises whether leaving it out
  could simplify the logic in the apps.} Furthermore, 
$\silv_\oldXi\oldLambda$ is componentwise weakly \simple\ if and only if the pair $(\oldLambda,\oldXi)$
is acylindrical (see \ref{acylindrical def}).
\EndProposition

\Proof 
We may assume without loss of generality that $\oldLambda$ is
connected (and therefore strongly \simple).
We set $\mu=
\mu_{\oldLambda,\oldXi}:\oldLambda\to\silv_\oldXi\oldLambda$
and 
$\mu^*=\mu^*_{\oldLambda,\oldXi}:\oldLambda-\oldXi\to
(\silv_\oldXi\oldLambda)-\mu (\oldXi)$. 
According to  \ref{first point}) under the
 canonical identification of $\silv\oldPi$ with a suborbifold of
$\silv_\oldXi\oldLambda$, we have $\silv\oldPi=\mu(\oldPi)$.

To prove the first assertion, first suppose that $\oldPi$
is
essential in the pair $(\oldLambda,\oldXi)$. To show that $\silv\oldPi=\mu(\oldPi)$ is incompressible (i.e. $\pi_1$-injective) in
$\silv_\oldXi\oldLambda$, we apply Proposition \ref{kinda dumb}, more specifically the implication (b)$\Rightarrow$(a),
letting $\silv_\oldXi\oldLambda$ and $\mu(\oldPi)$ play the roles of
$\oldPsi$ and $\oldTheta$. (The strong \simple ity of
$\oldLambda$ implies that $\oldLambda$ is very good (see \ref{oops}) and hence that
$\silv_\oldXi\oldLambda$ is very good, as required for the application
of  Proposition \ref{kinda
  dumb}.) 
%The non-sphericity of $\mu(\oldPi)$, also required for the application of the final assertion
%of Proposition \ref{kinda
 % dumb}, follows from the annularity of $\oldPi$.)
% that $\oldPsi$ is very good and that $\oldTheta$ has no
%spherical component follow from the strong \simple ity of
%$\oldLambda$, cf. \ref{oops}. \redcomment{Does it? I'm confused because
  %we are talking about icompressibility in $\silv_\oldXi\oldLambda$,
  %not in $\oldLambda$.} 
To prove $\pi_1$-injectivity, it suffices to show that Condition (b) of Proposition \ref{kinda dumb} holds. Thus we consider an arbitrary  discal orbifold $\frakD\subset
\silv_\oldXi\oldLambda$ such that  
$\frakD\cap\mu(\oldPi)=\partial\frakD$.
%, and (ii) $|\frakD|$ meets $\fraks_{\silv_{\oldXi}\oldLambda}$ in general position in $|{\silv_{\oldXi}\oldLambda}|$. 
We are required to show that $\partial\frakD$ is
the boundary of some discal suborbifold of $\mu(\oldPi)$. Consider first the case in which
$|\partial\frakD|\subset\inter|\mu(\oldPi)|=|\mu|(\inter|\oldPi|)$. Then
$\frakD\subset (\silv_\oldXi\oldLambda)-\mu (\oldXi)$, and
$(\mu^*)^{-1}(\frakD)$ is
% $\eta|(\oldLambda-\oldXi):
%\oldLambda-\oldXi\to (\silv_\oldXi\oldLambda)-\eta(\oldXi)$ (see
%\ref{silvering}), with 
a discal suborbifold of $\oldLambda-\oldXi$ whose
boundary is contained in $\inter\oldPi$. Since $\oldPi$ is essential,
it is in particular $\pi_1$-injective in $\oldLambda$. It therefore follows from Proposition \ref{kinda dumb}, this time applied with $\oldLambda$ and $\oldPi$ playing the roles of $\oldPsi$ and $\oldTheta$, that $\partial((\mu^*)^{-1}(\frakD))=\partial\frakE$ for some discal
suborbifold $\frakE$ of $\oldPi$. Since 
$\partial((\mu^*)^{-1}(\frakD))\subset\inter\oldPi$, we have $\frakE\subset\inter\oldPi$, and $\mu^*(\frakE)$ is
a discal suborbifold of $\mu(\oldPi)$ whose boundary is $\partial\frakD$; this gives the required conclusion in this case. Now consider the case in which $|\partial\frakD|$ meets
$\partial|\mu(\oldPi)|\subset|\mu(\oldXi)|\subset\fraks_{\silv_\oldXi\oldLambda}$. 
In particular we then have $|\partial\frakD|\cap\fraks_\frakD\ne\emptyset$; since
$\frakD$ is discal, it follows that $|\frakD|$ is a disk, and that
$\fraks_\frakD$ consists of either (i) a single one-dimensional
stratum, which is a closed arc, or (ii) two one-dimensional strata,
each of which is a half-open arc, and a single zero-dimensional
stratum, which is contained in the closure of each of the
one-dimensional strata. In either case, $\alpha_0:=\fraks_\frakD$  is
an arc contained in $\partial|\frakD|$, and is the closure of the
union of the one-dimensional strata of $\frakD$; furthermore, each
one-dimensional stratum of $\frakD$ contains an endpoint of the arc
$\alpha_0$. It now follows that $\beta_0:=|\partial
\frakD|=\overline{\partial|\frakD|-\alpha_0}$ is also an arc, and that
the arcs $ \alpha_0$ and $\beta_0$ have the same endpoints. Since
$\frakD\cap\mu(\oldPi)=\partial\frakD$, we have
$|\frakD|\cap|\mu(\oldPi)|=\beta_0$, and hence
$\partial\beta_0\subset|\mu(\partial\oldPi)|$.

If $\sigma$ is any one-dimensional stratum of $\frakD$, then since
$\frakD$ is a suborbifold of $\silv_\oldXi\oldPsi$, some stratum
$\tau$ of $\silv_\oldXi\oldPsi$ contains $\sigma$. We have seen that
any one-dimensional stratum of $\frakD$ contains an endpoint of
$\alpha_0$, so that
$\emptyset\ne\sigma\cap\partial\alpha_0=\sigma\cap\partial\beta_0\subset\tau\cap|\mu(\partial\oldPi)|\subset\tau\cap|\mu(\oldXi)|$. Since $|\mu(\oldXi)|$ is a union of strata of $\silv_\oldXi\oldPsi$ by \ref{silvering}, it follows that $|\mu(\oldXi)|\supset\tau\supset\sigma$. Thus every  one-dimensional stratum of $\frakD$ is contained in $|\mu(\oldXi)|$. Since
 $\alpha_0:=\fraks_\frakD$  is the closure of the union of the one-dimensional strata of $\frakD$, it follows that $\alpha_0\subset|\mu(\oldXi)|$.

Recalling from \ref{silvering} that  $|\mu|:|\oldLambda|\to|\silv_\oldXi\oldLambda|$ is a
homeomorphism, we set $D=|\mu|^{-1}|\frakD|$, $\alpha=|\mu|^{-1}|\alpha_0|$, and $\beta=|\mu|^{-1}|\beta_0|$. Then $\beta=D\cap|\oldPi|$ is an arc in
$\partial D$, and $\alpha=\overline{\partial
  D-\beta}$. Since $\alpha_0\subset|\mu(\oldXi)|$ we have $\alpha\subset|\oldXi|$. 

We wish to apply Proposition \ref{what essential means}, which
requires showing that $D$ has weight $0$. For this purpose, first note
that
$|\mu|(D\cap\fraks_\oldPsi)\subset|\frakD\cap\fraks_{\silv\oldPsi}|=|\fraks_\frakD|=\alpha_0\subset|\mu(\oldXi)|$,
so that $D\cap\fraks_\oldPsi\subset|\oldXi|\subset|\partial\oldPsi|$. Since $\oldPsi$ is
orientable, it follows that any point of intersection of the properly
embedded disk $D\subset|\oldPsi|$ with $\fraks_\oldPsi$ would be an
endpoint of some closed or half-open arc which is a one-dimensional
stratum of $\oldPsi$; in this case $|\mu|(D)$ would not be the
underlying set of a suborbifold
of $\silv_\oldXi\oldPsi$. Hence 
$D\cap\fraks_\oldPsi=\emptyset$, i.e. $\wt D=0$. It now follows from Proposition \ref{what essential means} that $\beta$ is the
frontier of a weight-$0$ disk  $E\subset|\oldPi|$.  
(The requirement in Proposition \ref{what essential means} that the
ambient orbifold be very good, and that the suborbifold $\oldXi$ be
non-spherical, holds by virtue of the strong \simple ity of
$\oldLambda$, cf. \ref{oops}.) 
%But this implies that 
Hence 
$\partial\frakD=\obd(\beta_0)=\mu(\obd(\beta))$ is
the boundary of the discal suborbifold $\mu(\obd(E))$ of $\mu(\oldPi)$, and Condition (b) of Proposition \ref{kinda dumb} is verified in this case as well. This  shows that  $\mu(\oldPi)$ is incompressible in
$\silv_\oldXi\oldLambda$. 
%again have a contradiction. 
%\redproofreadingnote{To check that this app. of Proposition \ref{what essential means} matches the statement, I think I need to know that $\alpha\subset\oldXi$. This seems a little tricky, because a priori a non-orientable $2$-suborbifold of $\oldLambda$ could have its $1$-dim. singular set contained in the singular locus of $\inter\oldLambda$. That won't happen here because of the way we showed $\frakD$ has singularities on $\partial|\frakD|$, but it shows we have to be careful.} 

Now suppose that  $\silv\oldPi$ is parallel in $\silv_\oldXi\oldLambda$ to a boundary component
% $\frakB$
 of  $\silv_\oldXi\oldLambda$. Then by definition there exists an embedding $j:(\silv\oldPi)\times[0,1]\to\silv_\oldXi\oldLambda$ such that
%$j((\partial\oldPi)\times[0,1])\subset\oldXi$, 
$j((\silv\oldPi)\times\{0\})=\silv\oldPi$, and
$j((\silv\oldPi)\times\{1\})\subset\partial(\silv_\oldXi\oldLambda)$. Applying
Proposition \ref{oh yeah he tweets}, with $\oldPi\times[0,1]$,
$(\partial\oldPi)\times[0,1]$, $\oldLambda$, and $\oldXi$ playing the
respective roles of $\oldPsi_1$, $\oldXi_1$, $\oldPsi_2$, and
$\oldXi_2$, we obtain an embedding 
$j^0:(\oldPi\times[0,1], (\partial\oldPi)\times[0,1])
\to(\oldLambda
,\oldXi)
$ such
that 
%$j^0((\partial\oldPi)\times[0,1])\subset\oldXi$ and 
$\silv j^0=j$. Hence
$\mu\circ j^0=j\circ(\mu_\oldPi\times\id_{[0,1]})$. Since
$j((\silv\oldPi)\times\{0\})=\silv\oldPi$,
% and $\mu\circ
%j^0=j\circ(\mu_\oldPi\times\id)$, 
we have $j^0(\oldPi\times\{0\})=\oldPi$. 
%, then $|\mu_\oldPi|(|\partial\oldPi|)\times[0,1]$ is a topological $2$-manifold contained in
%$\fraks_{(\silv\oldPi)\times[0,1]}$; hence $|j|(|\mu_\oldPi|(|\partial\oldPi|)\times[0,1])$ is contained in the union of the closures of the two-dimensional strata of $\silv_\oldXi\oldLambda$, which is $|\mu|(|\oldXi|)$. On the other hand, $|\mu_\oldPi|(\inter|\oldPi|)\times[0,1]$ is an open subset of $|(\silv\oldPi)\times[0,1]|$ whose intersection with  $\fraks_{\oldPi\times[0,1]}$ is at most one-dimensional, and hence $|j|(\inter|\oldPi|)\times[0,1])\cap|\oldXi|=\emptyset$.  It follows that $|j|^{-1}(|\oldXi|)=
%|\mu_\oldPi|(|\partial\oldPi|)\times[0,1]$, 
%and therefore that there is an embedding $j^0:\oldPi\times[0,1]\to\oldLambda$ such that $|j^0|=|j|$. We then have
%$j^0((\partial\oldPi)\times[0,1])\subset\oldXi$ and
%$j^0(\oldPi\times\{0\})=\oldPi$. 
Furthermore, since
$j((\silv\oldPi)\times\{1\}) \subset\partial(\silv_\oldXi\oldLambda)=\mu(\overline{(\partial\oldLambda)-\oldXi})$, we have 
$j^0(\oldPi\times\{1\})
\subset\overline{(\partial\oldLambda)-\oldXi}$. This shows that
$\oldPi$ is parallel in the pair $(\oldLambda,\oldXi)$ to a suborbifold of $\overline{(\partial\oldLambda)-\oldXi}$, a contradiction to the essentiality of $\oldPi$. Hence 
$\silv\oldPi$ is not parallel in $\silv_\oldXi\oldLambda$ to a boundary component
% $\frakB$
 of  $\silv_\oldXi\oldLambda$.

Conversely, suppose that
$\silv\oldPi$
  is incompressible in
$\silv_\oldXi\oldLambda$ and is not parallel to
a boundary component of $\silv_\oldXi\oldLambda$. Note that  
$\mu_\oldPi:\oldPi\to\silv\oldPi$ is $\pi_1$-injective, since it is the
composition of the $\pi_1$-injective inclusion $\oldPi\to {\rm
  D}\oldPi$ 
with the natural covering map ${\rm D}\oldPi\to\silv\oldPi$. By
incompressibility the inclusion map
$i:\silv\oldPi\to\silv_\oldXi\oldLambda$ is also $\pi_1$-injective,
and hence so is $i\circ\mu_\oldPi:\oldPi\to\silv_\oldXi\oldLambda$. But we
have $i\circ\mu_\oldPi=\mu\circ h$, where $h:\oldPi\to\oldLambda$ is the inclusion. Hence $h$ is $\pi_1$-injective,
i.e. $\oldPi$ is $\pi_1$-injective in $\oldLambda$.

To show that $\oldPi $ is essential in the pair $(\oldLambda,\oldXi)$, it remains to show that it is not %boundary-parallel. Suppose to the contrary that 
%$\oldPi$ is
 parallel in $(\oldLambda,\oldXi)$ either to a suborbifold of $\oldXi$ or to a component of
$\overline{\partial\oldLambda-\oldXi}$. Assume to the contrary that there is
an embedding 
$k:\oldPi\times[0,1]\to\oldLambda$ such that 
$k((\partial\oldPi)\times[0,1])\subset\oldXi$,
$k(\oldPi\times\{0\})=\oldPi$, and
$k(\oldPi\times\{1\})$ is contained either in $\oldXi$ or in $\overline{(\partial\oldLambda)-\oldXi}$.
If $k(\oldPi\times\{1\})\subset\oldXi$, and if we set
$\frakJ=k(\oldPi\times[0,1])$ and $\frakB=k(\oldPi\times\{1\})$, then
${\rm D}_\frakB\frakJ\subset {\rm D}_\oldXi\oldLambda$ is a \torifold\ whose
boundary is ${\rm D}\oldPi$; this is impossible, because the
incompressibility of $\silv\oldPi$ in $\silv_\oldXi\oldLambda$ implies
that the 
%kernel of the
 inclusion homomorphism ${\rm D}\oldPi\to
{\rm D}_\oldXi\oldLambda$ is injective.
%has order
%at most $2$. \redcomment{That looks backwards!}
Now consider the case in which
$k(\oldPi\times\{1\})\subset\overline{\partial\oldLambda-\oldXi}$. 
%I'd written that ``another ref. to lemma \ref{oh yeah he tweets} should be
%given below, and will replace some glop.'' Not sure now whether that
%is quite right.
%  Decide whether the proof of Prop. \ref{what essential means} needs
  %another proofreading. } 
Then, since $k((\inter\oldPi)\times[0,1))\subset\inter\oldLambda$ and 
$k((\partial\oldPi)\times[0,1])\subset\oldXi$, we have
$k^{-1}(\oldXi)=(\partial\oldPi)\times[0,1]$; hence  \ref{third
  point}, applied to the embedding of pairs 
$k:(\oldPi\times[0,1], (\partial\oldPi)\times[0,1])\to(\oldLambda,\oldXi)$,
implies that
$k':=\silv k :(\silv\oldPi)\times[0,1]\to\silv_\oldXi\oldLambda$
is an 
embedding. Since \ref{third point} also gives $\mu\circ
k=k'\circ(\mu_\oldPi\times\id_{[0,1]})$, we have
%$k'((\partial\oldPi)\times[0,1])\subset\oldXi$,
$k'((\silv\oldPi)\times\{0\})=\silv\oldPi$ and
% since $|\overline{\partial\oldLambda-\oldXi}|=|\partial\silv_\oldXi\oldLambda|$, we have
$k'((\silv\oldPi)\times\{1\})\subset\mu(\overline{\partial\oldLambda-\oldXi})=
\partial\silv_\oldXi\oldLambda$. This shows that $\silv\oldPi$ is parallel in $\silv_\oldXi\oldLambda$ to a suborbifold of $\partial(\silv_\oldXi\oldLambda)$, which must be a component of $\partial(\silv_\oldXi\oldLambda)$ since $\silv\oldPi$ is closed; this is a contradiction. This completes the proof that $\oldPi$ is essential in $(\oldLambda,\oldXi)$, and the first assertion of the proposition is established.

% $\overline{(\partial\oldLambda)-\oldXi}$.

%there is an embedding
%$k:\oldPi\times[0,1]\to\oldLambda$ such that 
%$|j^0|=|j|$. We then have
%$k((\partial\oldPi)\times[0,1])\subset\oldXi$ and
%$j^0(\oldPi)\times\{0\})=\silv\oldPi$. Furthermore, since
%$|j|(|\oldPi|\times\{1\}) \subset|\partial(\silv_\oldXi\oldLambda)|=\overline{(|\partial\oldLambda)|-|\oldXi|}$, we have 
%$j^0(\oldPi\times\{1\}) \subset\overline{(\partial\oldLambda)-\oldXi}$. 

To prove the second assertion, first suppose that
$\silv_\oldXi\oldLambda$ is weakly \simple. If $\oldPi$ is an orientable annular
orbifold which is essential in the pair $(\oldLambda,\oldXi)$, then the first
assertion implies that the 
toric  orbifold
 $\silv\oldPi$
  is incompressible in
$\silv_\oldXi\oldLambda$ and is not parallel to
a boundary component of $\silv_\oldXi\oldLambda$; this contradicts the
weak
\simple ity of $\silv_\oldXi\oldLambda$. Hence $(\oldLambda,\oldXi)$
is acylindrical. 

Conversely, suppose that  $(\oldLambda,\oldXi)$
is acylindrical. Since $\oldLambda$ is strongly \simple, it is irreducible, and the suborbifold $\oldXi$ of $\partial\oldLambda$ has no spherical component (see \ref{oops}); as $\oldXi$ is $\pi_1$-injective, it then follows from Proposition \ref{silver irreducible} that $\silv_\oldXi\oldLambda$ is
irreducible. The strong \simple ity of $\oldLambda$ implies that
$\oldLambda$ is
very good and non-discal, and hence 
$\silv_\oldXi\oldLambda$ has the same properties. 
Now suppose that $\oldTheta$ is an incompressible toric
suborbifold of $\silv_\oldXi\oldLambda$.
% which  is not parallel to
%a boundary component of $\silv_\oldXi\oldLambda$. 
Since $\oldTheta$ is in
particular a two-sided  closed suborbifold of
$\silv_\oldXi\oldLambda$, it follows from \ref{second point} that we may
write $\oldTheta=\silv\oldPi$ for some two-sided suborbifold
$\oldPi$ of $\oldLambda$ with $\partial\oldPi\subset\oldXi$. We have
$\chi(\oldPi)=\chi(\oldTheta)=0$, so that $\oldPi$ is either toric or
annular. If $\oldPi$ is toric, then 
the
homeomorphism $\mu^*$
% $\eta|(\oldLambda-\oldXi):
%\oldLambda-\oldXi\to (\silv_\oldXi\oldLambda)-\eta(\oldXi)$
maps $\oldPi$ onto $\oldTheta$, and the
incompressibility of $\oldTheta$ in $\silv_\oldXi\oldLambda$ implies that
$\oldPi$ is incompressible in $\oldLambda$, a contradiction to the
strong \simple ity of $\oldLambda$ see \ref{oops}. Hence $\oldPi$ is annular. The acylindricity of the pair
$(\oldLambda,\oldXi)$  implies that $\oldPi$ is not
essential in the pair $(\oldLambda,\oldXi)$. It therefore follows from the first
assertion of the present proposition that $\oldTheta$ is  parallel to
a boundary component of $\silv_\oldXi\oldLambda$. This shows that $\silv_\oldXi\oldLambda$ is
weakly \simple. 
\EndProof
%\oldPsi\obd***\eta canonical\xi\oldUpsilon \frakZ \oldPsi identif \xi
%\eta \xi canonical j' {second point} (a) (i \fraks_\frakD \partial\frakD\mu^*\silv\oldPsi

\abstractcomment{\tiny This lemma is apparently not used. For introduction and proof summary, see cravat.tex.
%The following graph-theoretical result will be needed in Section \ref{structure section}.

\Lemma\label{graph lemma}
Let $\calt$ be a finite graph which has no isolated vertices, su y ch that each component of $\calt$ is a tree. Let $\calv$ and $\cale$ denote respectively the vertex set and edge  set of $\calt$. Let $\calv_0$ be a subset of $\calv$ such that (a) each component of $\calt$ contains at most one vertex in $\calv_0$, and (b) each vertex in $\calv_0$ has valence $1$ in $\calt$. Let $\cale_0$ be a subset of $\cale$ such that each endpoint of each  edge in $\cale$ either has valence at least $3$ or belongs to $\calv_0$. Then $\card \cale_0<(\card \cale)/2$.
\EndLemma
}

%\Number\label{move it}
%\EndNumber strongly

\abstractcomment{\tiny What else?}

\chapter{Fibrations, characteristic suborbifolds, and
  invariants}\label{fibration chapter}

The importance of the notion of a characteristic suborbifold, and of
the invariant $\volorb$, was stressed in the Introduction. The central
concept underlying the definition of the characteristic suborbifold is
that of a fibration (with $1$-dimensional fibers) of an orbifold, which we study
in Section \ref{fibration section}. The theory of the characteristic
suborbifold itself is developed in Section \ref{characteristic
  section}, taking as a starting point the results proved in
\cite{bonahon-siebenmann} for the case of a closed orbifold, and
following the suggestions made there for extending to the bounded
case, which is the case that we need in the present work. The
invariant $\volorb$ is defined and studied in Section \ref{darts
  section}, as are several closely related invariants which we have
found convenient for our arguments. A crucial point is to relate these
invariants to the theory of the characteristic suborbifold; this
connection was mentioned in the Introduction.

\section{Fibrations of orbifolds}\label{fibration section}

%The material in Subsections \ref{can sub} and \ref{fibered stuff}
%---\ref{fibration-category}, all statements
%will be understood to apply to the smooth, PL or topological category.

\Number\label{can sub}
Let $n$ be an integer with $1\le n\le3$, and let $ U$ denote either $ U^{n-1}$ or 
$ U^{n-1}_+$. Let $J$ be a compact
$1$-manifold, which we equip with a metric in such a way that each
component of $J$ is isometric to either $\SSS^1$ or $[0.1]$. Suppose that
$G$ is a finite subgroup of 
 $\Isom( U)\times\Isom(J)$
%^2
(where $\Isom( U)$ is either $\OO(n-1)$ or a subgroup of
$\OO(n-1)$ having order at most $2$, while $\Isom(J)$ is a product of
copies of $\OO(2)$ and groups of order $2$). 
 Let $X\le\Isom( U)$ denote the 
 image of $G$ under the projection to the first factor. Then the
 projection $ U\times J\to  U$ induces a canonical
 map of sets $( U\times J)/G\to  U/X$.

According to \ref{esso}, $U$ and $ U\times J$ have canonical PL
manifold structures, and the actions  of $X$ and $G$  on $U$ and
$ U\times J$ are PL, so that $( U\times J)/G$ and $  U/X$ acquire
structure of PL orbifolds. The canonical map $( U\times J)/G\to  U/X$ is
then a PL immersion. If we have either $U=U^{n-1}$ or $\partial
J=\emptyset$ then $U$ and $ U\times J$ have canonical smooth
manifold structures,  the actions  of $X$ and $G$  on $U$ and
$ U\times J$ are smooth, and the canonical map $( U\times J)/G\to
U/X$ is a smooth immersion of smooth orbifolds. (In the smooth category we avoid the case in which  $U=U^{n-1}_+$ and $\partial
J\ne\emptyset$, as $U\times J$ would then have only the structure of a
manifold with corners.)

\EndNumber
%X

\DefinitionsRemarks\label{fibered stuff} (Cf.
\cite{bonahon-siebenmann}, \cite{other-bs}.) 
Suppose that $n$ is an integer with $1\le n\le3$,  that
$\oldLambda$ is a compact PL or smooth $n$-orbifold, 
%, a connected $1$-orbifold $J$, 
 that $\frakB$ is respectively a compact PL or smooth
$(n-1)$-orbifold, and that $J$ is a compact $1$-manifold. In the case
where $\oldLambda$ and $\frakB$ are smooth, assume that either
$\partial\frakB=\emptyset$ or that
$\partial J=\emptyset$. Let $q:|\oldLambda|\to|\frakB|$ be a continuous map. Let a
point $v\in\frakB$ be given. We define a (PL or smooth) {\it
  local $J$-standardization of $q$ at $v$} to be a 
quadruple $(\frakV,G,\psi,\zeta)$, where (i) $\frakV$ is an open
neighborhood of $v$ in $\frakB$; 
(ii) $G$ is a finite subgroup   of 
 $\Isom( U_v)\times\Isom(J)$, where $U_v$ is defined as in
 \ref{orbifolds introduced}, and $J$ is equipped with a metric under
 which each of its components is isometric to either $\SSS^1$ or $[0,1]$; (iii) $\psi:( U_v\times J)/G \to \obd(q^{-1}(|\frakV|))$
 and $\zeta:\frakV\to  U_v/X$  are, in the respective cases, PL
 orbifold homeomorphisms or diffeomorphisms (cf. \ref{can sub}), where $X\le\Isom( U_v)$ denotes the
 image of $G$ under the projection to the first factor; 
and (iv)
 $|\zeta|\circ q\circ|\psi|=|\beta|$, where $\beta:( U_v\times
 J)/G\to  U_v/X$  denotes the canonical submersion described in
 \ref{can sub}. Note that if $(\frakV,G,\psi,\zeta)$ is a   local
 $J$-standardization of $q$ at $v$, there is a
  local $J$-standardization
$(\frakV',G',\psi',\zeta')$  of $q$ at $v$ such that $\frakV'\subset\frakV$ and
 $\zeta'(0)=v$. 

If $s=(\frakV,G,\psi,\zeta)$ is a local $J$-standardization of
$q$ at $v$, we will denote by $p_1^s: \Isom( U_v)\times\Isom(J)\to \Isom( U_v)$
the projection to the first factor, and by $X_s$ the group
$p_1^s(G)\le\Isom( U_v)$. 
Composing the quotient map
$ U_v\to| U_v/X_s|$ with $|\zeta|^{-1}: | U_v/X_s|\to
|\frakV|$, we obtain a chart map $U_v\to |\frakV|$
for $\frakB$, which will be denoted $\phi_s$. 
We will 
denote by
$\pi_s :  U_v\times J\to( U_v\times J)/G$ the quotient
 submersion, 
and by
$P_s$ the covering map $\psi\circ \pi_s:  U_v\times J\to \obd(q^{-1}(|\frakV|))$.
Note that the fibers of $q|q^{-1}(\obd(|\frakV|))$ are the images under $|P_s|$
 of the $1$-manifolds of the form $\{x\}\times J$ for $x\in
 U_v$. 
%\OO(n-1) \DD^{n-1}

By a PL or smooth {\it (orbifold) fibration} 
$q:\oldLambda\to\frakB$ (or a PL or smooth {\it orbifold) fibration} of
$\oldLambda$ with base $\frakB$) we mean a continuous map 
$q:|\oldLambda|\to|\frakB|$ such that for every
$v\in\frakB$, there exists respectively a PL or smooth local $J$-standardization of $q$
at $v$, where $J$ is some compact $1$-manifold, and where, in the
smooth case, either $\partial J=\emptyset$ or $\partial\frakB=\emptyset$. We
may always choose this local standardization in such a way that $\psi(0)=v$.
%(\obd(q^{-1}(|\frakV|)))$
%is the diameter
 %$[-1,1]\times\{0\}\subset\DD^2_+= U_v$. \frakJ

Note that any PL or smooth orbifold fibration is in particular a PL or
smooth orbifold submersion, respectively.

If $q:\oldLambda\to\frakB$
is a PL or smooth orbifold fibration, then for each $w\in\frakB$ the {\it fiber}
$q^{-1}(w)$ is respectively a PL or smooth compact $1$-suborbifold of $\oldLambda$. If $J$ is a
compact $1$-manifold and
$(\frakV,G,\psi,\zeta)$ is a
 local $J$-standardization of $q$ at a point $v\in\frakB$, then
for every point $w\in\frakV$ we have
$\compnum(q^{-1}(w))=\compnum(q^{-1}(\frakV))$. In particular, the
function $w\to \compnum(q^{-1}(w))$ is locally constant on
$\frakB$. Hence if $\frakB$ is connected and the fiber of $q$ over some point of $\frakB$
s connected, then  the fiber of $q$ over every point of $\frakB$
s connected. 

Note that in the smooth category we have not defined orbifold
fibrations with non-closed fibers over non-closed bases, nor will we
refer to such fibrations in the smooth category anywhere in the
monograph. (According to an observation made in \ref{can sub},
defining such fibrations would require the notion of a smooth orbifold with corners.)

A PL or smooth  orbifold fibration $q:\oldLambda\to\frakB$ will be called an {\it
  $\SSS^1$-fibration} if every fiber is homeomorphic to a
quotient of $\SSS^1$ (and therefore to $\SSS^1$ or to $[[0,1]]$). It will be called an {\it
  $I$-fibration} if every fiber is homeomorphic to  
quotient of $[0,1]$ (and therefore to $[0,1]$ or to $[[0,1]$); note
that according to our conventions, the notion of a smooth
$I$-fibration makes sense only when the base is a closed orbifold. Note
that for any arbitrary fibration $q$, if $J$ is a
compact $1$-manifold and
$(\frakV,G,\psi,\zeta)$ is a
 local $J$-standardization of $q$ at a point $v\in\frakB$, then
the fiber over every point $w\in\frakV$ is homeomorphic to the
quotient of $J$ by its stabilizer in $G$; it follows that if $\frakB$
is connected, and if some (and hence every) fiber is connected, then $q$
is either an $\SSS^1$-fibration or an $I$-fibration.

A PL or smooth $I$-fibration $q:\oldLambda\to\frakB$ will be termed {\it trivial} if there exists a homeomorphism $t:\frakB\times[0,1]\to\oldLambda$ such that $q\circ t:\frakB\times[0,1]\to\frakB$ is the projection to the first factor. 

If a compact, connected, PL or smooth $2$-orbifold $\oldLambda$
admits a (respectively PL or smooth) orbifold
fibration whose base is a closed $1$-orbifold, then
$\chi(\oldLambda)=0$. Hence if $\oldLambda$ admits an $I$-fibration
with a closed base,
then $\oldLambda$ is an annular orbifold (see \ref{annular etc}).

If $\oldLambda$ is a compact, orientable PL $3$-orbifold equipped with
a PL $I$-fibration $q:\oldLambda\to\frakB$, where $\frakB$ is a $2$-orbifold, 
we define the {\it vertical boundary} $\partialv\oldLambda$ to be the
suborbifold $q^{-1}(\partial\frakB)$, of $\oldLambda$, and we define the {\it horizontal boundary}
$\partialh\oldLambda$ to be the suborbifold of $\oldLambda$ which is
the union (in the sense of \ref{unions and such}) of the (orbifold-)boundaries of
all the fibers. 
%We regard $\partialv\oldLambda$ and
%$\partialv\oldLambda$ as suborbifolds of $\oldLambda$. 
We have
$\partial\oldLambda=\partialh\oldLambda\cup\partialv\oldLambda$, and
$\partial(\partialh\oldLambda)=\partial(\partialv\oldLambda)=\partialh\oldLambda\cap\partialv\oldLambda$. If
$\frakB$ (or equivalently $\oldLambda$) is connected, then
$\partialh\oldLambda$ has at most two components, because each
component contains an endpoint of each fiber. 

Note that $q$ restricts to an $I$-fibration of $\partialv\oldLambda$ whose base is $\partial\frakB$; hence $\partialv\oldLambda$ is annular. Note also that $q|\partialh\oldLambda:\partialh\oldLambda\to\frakB$ is a degree-$2$ orbifold covering. In particular, if the $I$-fibration $q$ is non-trivial, then $\partialh\oldLambda$ is connected, and the image of the inclusion homomorphism $\pi_1(\partialh\oldLambda)\to\pi_1(\oldLambda)$ has index $2$.
\EndDefinitionsRemarks
%\zeta \alpha \phi X

Note that, according to the convention posited in \ref{categorille},
the PL category will be the default category of orbifolds for the rest
of this section. (In Proposition \ref{fibration-category} and in some
of the succeeding proofs, the smooth
category is explicitly considered.)

\Lemma\label{one way to look at it}
Let $q:\oldPsi\to\frakB$ be an $\SSS^1$-fibration of 
orbifolds, where $\dim\oldPsi=3$ and $\dim\frakB=2$. Let
$\frakJ$ be a $1$-suborbifold of
$\partial\frakB$, and set $\oldXi=q^{-1}(\frakJ)$, so that
$q:(\oldPsi,\oldXi)\to(\frakB,\frakJ)$ is a submersion of pairs. Then the
submersion $\silv q:\silv_\oldXi\oldPsi\to\silv_\frakJ\frakB$ defined as in
\ref{third point} is an $\SSS^1$-fibration. 
\EndLemma

\Proof
Set $\mu=\mu_{\oldPsi,\oldXi}$ and
$\nu=\mu_{\frakB,\frakJ}$. We must show
that for each point $v\in\silv\frakB$ there exists a local $\SSS^1$-standardization  $(\frakV,G,\psi,\zeta)$ of $\silv q$ at $v$.
If
$v\in\frakB-\nu(\frakJ)$, the existence of such a
local $\SSS^1$-standardization is immediate from the hypothesis that $q$ is an
$\SSS^1$-fibration (and the fact that $\mu^*_{\frakB,\frakJ}$ and $\mu^*_{\oldPsi,\oldXi}$ are
homeomorphisms, see \ref{silvering}). Now suppose that
$v\in\nu(\frakJ)$. Write $v=\nu(v^0)$ for some $v^0\in\frakJ$. In the notation of
\ref{fibered stuff} we then have $ U_v= U^2_+$. 
Since $q$ is an
$\SSS^1$-fibration, there exists a local $\SSS^1$-standardization
$s^0=(\frakV^0,G^0,\psi^0,\zeta^0)$ of $q$ 
%such that $\frakV^0$ is
%a neighborhood of 
at $v^0$; we may choose $s^0$ in such a way that
$\psi^0(0)=v^0$. Then 
$\zeta^0(\frakV^0\cap\partial\frakB)$
is the quotient by $X^0:=X^{s^0}$ of 
 $(-1,1)\times\{0\}\subset U^2_+= U_v$.
If $v^0\in\inter\frakJ$ (in which case $v^0$, regarded as a point of
$|\frakB|$, 
may or may not lie in $\partial|\frakJ|$), we may choose $s^0$ so that
$\frakV^0\cap\partial\frakB\subset\frakJ$, which implies that  
$\zeta^0(\frakV^0\cap\frakJ)=((-1,1)\times\{0\})/X^0$.
If $v^0\in\partial\frakJ$, we may choose $s^0$ so that
$\frakV^0\cap\partial\frakB$ is an open arc, and 
$\frakV^0\cap\frakJ$ is a half-open subarc of
$\frakV^0\cap\partial\frakB$ having $v^0$ as its  endpoint. In the latter
case we have $X^0=\{1\}$, so that $\zeta^0$ is a homeomorphism from
$\frakV^0$ to $ U^2_+$; after possibly modifying $\zeta^0$ by
postcomposition with the symmetry $(x,y)\mapsto(-x,y)$ of $ U^2_+$,
we may assume that
$\zeta^0(\frakV^0\cap\frakJ)=[0,1)\times\{0\}$. 
%$\zeta^0(\frakV^0\cap\frakJ)=((-1,1)\times\{0\})/X^0$.
%$\frakV^0\cap\partial\frakB\subset\frakJ$
%which implies that  
Let $A$ denote the preimage of $\zeta^0(\frakV^0\cap\frakJ)$ under the
quotient map $ U^2_+\to U^2_+/X^0$, so that $A$ is equal to
$(-1,1)\times\{0\}$ if $v^0\in\inter\frakJ$, and to $[0,1)\times\{0\}$
if $v^0\in\partial\frakJ$. In either case, $A$ is
$X^0$-invariant.

The $X^0$-invariance of $A$ implies that
${\rm D}g$ is a well-defined
self-homeomorphism of 
${\rm D}_{A} U^2_+$ for every
$g\in X^0$, and ${\rm D}G^0:=\{{\rm D}g:g\in G^0\}$ is a group of
self-homeomorphisms of   ${\rm D}_{A} U^2_+$, isomorphic to $G^0$, and normalized
(in fact centralized)  by the canonical involution $\delta_1$ of
${\rm D}_{A} U^2_+$. This implies that ${\rm D}X^0$  and $\delta_1$ generate a finite group $X$ of
self-homeomorphisms of ${\rm D}_{A} U^2_+$. If
$v^0\in\inter\frakJ$, so that
$A=(-1,1)\times\{0\}$, then under the canonical identification of
${\rm D}_{A} U^2_+$ with $ U^2$ we have $X\le\Isom( U^2)=\OO(2)$.
If
$v^0\in\partial\frakJ$, so that
$A=[0,1)\times\{0\}$, then 
${\rm D}_{A} U^2_+$  may be identified homeomorphically with $ U^2_+$
in such a way that  $X\le\Isom( U^2_+)$. But we have $v=\nu(v^0)\in\inter\silv\frakB$
  if $v^0\in\inter\frakJ$, and $v=\nu(v^0)\in\partial\silv\frakB$
  if $v^0\in\partial\frakJ$. Thus in either
case,
${\rm D}_{A} U^2_+$  is identified with $ U_v$, and
 $X$  with a subgroup of $\Isom( U_v)$.

 Thus $({\rm D}_{A} U^2_+)\times \SSS^1=
{\rm D}_{A\times \SSS^1}( U^2_+\times \SSS^1)
$
is homeomorphically identified with $ U_v\times \SSS^1$, and the
group $G$ generated by ${\rm D}G^0$ and $\delta_1\times\{\id\}$ is
identified with a subgroup of $\Isom( U_v)\times \OO(2)$ in such a
way that $p_1(G)=X$, where 
$p_1:\Isom( U_v)\times\OO(2)\to
\Isom( U_v)$ denotes the projection to the first factor.

% We may identify ${\rm D}_{A}\DD^2_+$  homeomorphically with $\DD^2_+$ in such a way that $X$
%becomes a subgroup of $\OO(2)$; the manifold

% such that $\frakV^0\cap\frakJ$
%is a connected $1$-orbifold,  and
%$\frakV^0\subset\frakU$. Property (b) of $\eta$, together with the
%connectedness of $\frakV^0\cap\frakJ$, implies that $\zeta^{-1}(v^0)$
%consists of a single point which lies in  each component of
%$\zeta^{-1}(\frakJ)$, and hence that
%$A:=\zeta^{-1}(\frakJ)\subset\partial \DD^2_+$ is connected. But Property
%(c) of $\frakU$ implies that  $\frakV^0\cap\frakJ$ is not a closed orbifold; since
%$\zeta|A:A\to\frakV^0\cap\frakJ$ is a covering map, it now follows that
%$A$ is properly contained in $\partial \DD^2_+$. Hence $A$ is an arc.

%We may suppose $\frakV^0$ to be chosen so that
%$\frakJ':=\frakJ\cap\frakV^0$ is a single arc. J

In view of the definition of a local $\SSS^1$-standardization, we have
$P_{s^0}^{-1}(q^{-1}(\frakV^0)\cap\oldXi)=P_{s^0}^{-1}(q^{-1}(\frakV^0\cap\frakJ))=A\times
\SSS^1$ (where $P_{s^0}$
is defined as in \ref{fibered stuff}).
We therefore have 
$q^{-1}(\frakV^0)\cap\oldXi=
P_{s}(A\times \SSS^1)=\psi^0((A\times \SSS^1)/G^0)$.
Thus $\psi^0$ may
be regarded as a homeomorphism between the pairs 
$(( U^2_+\times \SSS^1)/G^0,(A\times \SSS^1)/G^0)$ and
$(q^{-1}(\frakV^0) ,q^{-1}(\frakV^0)\cap\oldXi)$, and in particular
${q^{-1}(\frakV^0)\cap\oldXi}$ is a codimension-$0$ suborbifold of the boundary
of $q^{-1}(\frakV^0)$. We therefore have a homeomorphism
$\silv\psi^0:\silv_{(A\times \SSS^1)/G^0}(( U^2_+\times \SSS^1)/G^0)\to\silv_{q^{-1}(\frakV^0)\cap\oldXi}(q^{-1}(\frakV^0))$. 
 
%Since 
%$P_{s^0}^{-1}(q^{-1}(\frakV^0)\cap\oldXi)=A\times
%\SSS^1$, the group $G^0\le\OO(2)\times\OO(2)$, regarded as a
%group of homeomorphisms of $ \DD^2_+\times \SSS^1$, leaves $A\times \SSS^1$ is invariant. It follows that
%$X^0\le\OO(2)$, regarded as a
%group of homeomorphisms of $ \DD^2_+$, leaves $A$  
%invariant 
%(and therefore has order at most $2$). Hence 

%\redcomment{Fix from here, using some
  %of the \%-ed out stuff, to get the group $G$, its action on
  %$\DD^2\times \SSS^1$, etc. The rest
  %of the present sentence is wrong---I have an extra generator, $\delta=(\delta_1,\id)$} if $p_1: {\rm
  %O}(2)\times\OO(2)\to\OO(2)$ denotes the projection to the
%first factor, we have $p_1(G)=p_1(G^0)$. 

Our identification of
$ U_v\times
\SSS^1$ with
${\rm D}_{A\times \SSS^1}( U^2_+\times
\SSS^1)$ gives rise to an identification of $( U_v\times
\SSS^1)/G$ with
$({\rm D}_{A\times \SSS^1}( U^2_+\times
\SSS^1))/G$, which is in turn  canonically identified with 
$\silv_{(A\times \SSS^1)/G^0}(( U^2_+\times \SSS^1)/G^0)$,
the domain of the homeomorphism $\silv\psi^0$.
%the latter actiontherefore gives rise to an action of $G^0$ on
%$\silv_{A\times \SSS^1}(\DD^2\times \SSS^1)$ via doubling: any element
%$g\in X^0$ acts via the
%self-homeomorphism ${\rm D}g$ of   $\silv_{A\times \SSS^1}(\DD^2\times \SSS^1)$. This (faithful) action gives an identification of $G^0$ with a
%group of self-homeomorphisms of $\silv_{A\times \SSS^1}(\DD^2\times
%\SSS^1)$, which is normalized by the canonical involution $\delta$ of
%$\silv_{A\times \SSS^1}(\DD^2\times
%\SSS^1$. Hence $G^0$ and $\delta$ generate a finite group $G$ of
%self-homeomorphisms of $\silv_{A\times \SSS^1}(\DD^2\times \SSS^1)$. 
 %We may identify $\silv_{A\times \SSS^1}(\DD^2\times \SSS^1)$  homeomorphically with $\DD^2\times \SSS^1$ in such a way that $G$
%becomes a subgroup of $\OO(2)\times\OO(2)$; \redcomment{The rest
  %of this sentence is wrong---I have an extra generator. Decide what
%  matters here} if $p_1: {\rm
  %O}(2)\times\OO(2)\to\OO(2)$ denotes the projection to the
%first factor, we have $p_1(G)=p_1(G^0)$. 
%The orbifold $(\DD^2\times
%\SSS^1)/G
%=(\silv_{A\times \SSS^1}(\DD^2\times
 %\SSS^1))/G$ is also canonically identified with 
%$\silv_{(A\times \SSS^1)/G^0}((\DD^2\times \SSS^1)/G^0)$,
%the domain of the homeomorphism $\silv\psi^0$.
On the other hand, according to \ref{first point}, the
image $\silv_{q^{-1}(\frakV^0)\cap\oldXi}(q^{-1}(\frakV^0))$ of $\silv\psi^0$
is canonically identified with a suborbifold of $\silv_\oldXi\oldLambda$, and we
have $\silv_{q^{-1}(\frakV^0)\cap\oldXi}(q^{-1}(\frakV^0))=\mu
(q^{-1}(\frakV^0))$. If we now set $\frakV=\nu(\frakV^0)$, it follows that
$\silv_{q^{-1}(\frakV^0)\cap\oldXi}(q^{-1}(\frakV^0))=(\silv
q)^{-1}(\frakV)$. We may therefore regard
$\silv\psi^0$ as a homeomorphism $\psi:( U_v\times
\SSS^1)/G\to (\silv q)^{-1}(\frakV)$. 

Note also that the homeomorphic identification of 
${\rm D}_{A} U^2_+$   with $ U_v$
 induces a homeomorphic identification of
$\silv_{A/X^0}( U^2_+/X^0)$ with $ U_v/X$, and that $\frakV=\nu(\frakV^0)$
is identified with $\silv_{\frakV^0\cap\frakJ}\frakV^0$ by \ref{first
  point}. We may therefore regard $\zeta:=\silv\zeta^0$ as a homeomorphism from
$\frakV$ to $ U_v/X$. The constructions of $\frakV$, $G$,
$\psi$ and $\zeta$ that we have now given imply that
$(\frakV,G,\psi,\zeta)$ is the required local $\SSS^1$-standardization
of $\silv q$ at $v$. 
\EndProof
%^0_0 (a) \frakJ\eta v\frakV' \psi' G'' U^2\DD U p_1 _0
% U\zeta G_1 p_1 X 

\Lemma\label{another way to look at it}
Let $q:\oldLambda\to\frakB$ be an $I$-fibration of
orbifolds, where $1\le \dim\oldLambda\le3$.  Set $\oldXi=\partialh\oldLambda$ and  $\mu=\mu_{\oldLambda,\oldXi}$. Then there is a unique $\SSS^1$-fibration
 $q':\silv_\oldXi\oldLambda\to\frakB$ such that $q'\circ
 \mu=q$. 
\EndLemma

\Proof
Set $n=\dim\oldLambda$, so that $\dim\frakB=n-1$. Since $|\mu|:
|\oldLambda|\to|\silv_\oldXi\oldLambda|$ is a homeomorphism by
\ref{silvering}, we may define a map of sets $Q:
|\silv_\oldXi\oldLambda|\to|\frakB|$ by $Q=|q|\circ|\mu|^{-1}$. We are
required to show that $Q$ defines an $\SSS^1$-fibration  
$\silv_\oldXi\oldLambda\to\frakB$.
%, i.e. we must show 
%that for each point $v\in\frakB$ there exist a neighborhood $\frakV$ of $v$ in
%$\frakB$, a suborbifold
%$\oldUpsilon$ of $\silv_\oldXi\oldLambda$ such that
%$|\oldUpsilon|=(Q)^{-1}(|\frakV|)$, 
%, a connected $1$-orbifold $J$, 
%and a homeomorphism $\psi:
%( U^{{n-1}}\times
%\SSS^1)/G
%\to
%\oldUpsilon$,
%where 
%$G$ is a finite subgroup of 
 %$\OO({n-1})\times\Isom(J)$,
%such that
%if $\pi:
%\DD^{n-1}\times
%J
%\to
%(\DD^{n-1}\times
%\SSS^1)/G$
%denotes the orbit map, then
%$|\psi(\pi(\{x\}\times \SSS^1))|$ 
%in $\psi:r^{-1}(\frakV)\to(\DD^{n-1}\times
%J)/G$
%is a fiber of $Q$
%for each $x\in \DD^{n-1}$. \redcomment{This is one of the passages that
%  need to be rewritten in light of the new def. of fibration (and
  %should involve a ref. to the notion of a local strandardization).}

Consider any point  $v\in\frakB$. Set $ U= U_v$ (see \ref{orbifolds introduced}). 
%For this purpose, note that 
Since $q$ is an $I$-fibration,
% and both
%$\oldLambda$ and $\frakB$ are orientable, 
there exists a local $[0,1]$-standardization $(\frakV,G,\psi,\zeta)$
of $q$ at $v$.   It follows from
\ref{fibered stuff}  
that
 the fibers of $q|q^{-1}(\frakV)$ are the images under
 $P_s$ of the arcs of the form $\{x\}\times [0,1]$ for $x\in
  U $. But
$q^{-1}(\frakV)\cap\oldXi=q^{-1}(\frakV)\cap\partialh\oldLambda$ is
the set of endpoints of fibers of $q|q^{-1}(\frakV)$, and hence
$q^{-1}(\frakV)\cap\oldXi=
P_{s}( U \times\{0,1\})=\psi(( U \times\{0,1\})/G)$. Thus we may regard $\psi$ as a homeomorphism between the pairs 
$(( U \times[0,1])/G,( U \times
\{0,1\})/G)$ and
$(q^{-1}(\frakV) ,q^{-1}(\frakV)\cap\oldXi)$, and in particular
${q^{-1}(\frakV)\cap\oldXi}$ is a codimension-$0$ suborbifold of the boundary
of $(q^{-1}(\frakV))$. We therefore have a homeomorphism
$\silv\psi:\silv_{( U \times\{0,1\})/G}(( U \times
[0,1])/G)\to\silv_{q^{-1}(\frakV)\cap\oldXi}(q^{-1}(\frakV))$. 
 
The canonical action of $G\le\Isom( U)\times\Isom([0,1])$
on $  U \times[0,1]$ gives rise to an action of $G$ on
${\rm D}_{ U \times\{0,1\}}( U \times [0,1])$ via doubling: any element
$g\in G$ acts via the
self-homeomorphism ${\rm D}g$ of   ${\rm D}_{ U \times\{0,1\}}( U \times
[0,1])$. This (faithful) action gives an identification of $G$ with a
group of self-homeomorphisms of ${\rm D}_{ U \times\{0,1\}}( U \times
[0,1])$, which is normalized by the canonical involution $\delta$ of
${\rm D}_{ U \times\{0,1\}}( U \times
[0,1])$. Hence $G$ and $\delta$ generate a finite group $G'$ of
self-homeomorphisms of ${\rm D}_{ U \times\{0,1\}}( U \times
[0,1])$. We may identify ${\rm D}_{ U \times\{0,1\}}( U \times
[0,1])$  homeomorphically with $ U \times \SSS^1$ in such a way that $G'$
becomes a subgroup of $\Isom( U)\times\OO(2)$; if $p_1: \Isom( U)\times\OO(2)\to\Isom( U)$ denotes the projection to the
first factor, we have $p_1(G')=p_1(G)$. The orbifold $( U \times
\SSS^1)/G'
=(\silv_{ U \times\{0,1\}}( U \times
[0,1]))/G'$ is also canonically identified with $\silv_{( U \times\{0,1\})/G}(( U \times
[0,1])/G)$, the domain of the homeomorphism $\silv\psi$. On the other hand, according to \ref{first point}, the
image $\silv_{q^{-1}(\frakV)\cap\oldXi}(q^{-1}(\frakV))$ of $\silv\psi$
is canonically identified with a suborbifold of $\silv_\oldXi\oldLambda$, and we
have $\silv_{q^{-1}(\frakV)\cap\oldXi}(q^{-1}(\frakV))=\mu
(q^{-1}(\frakV))$, so that
$|\silv_{q^{-1}(\frakV)\cap\oldXi}(q^{-1}(\frakV))|=
Q^{-1}(|\frakV|)$. Thus $\obd(Q^{-1}(|\frakV|))$ is well defined and
is canonically identified with
$\silv_{q^{-1}(\frakV)\cap\oldXi}(q^{-1}(\frakV))$. We may therefore regard
$\silv\psi$ as a homeomorphism $\psi':( U \times
\SSS^1)/G'\to \obd(Q^{-1}(|\frakV|))$. Note also that since
$p_1(G')=p_1(G)$, we may regard $\zeta$ as a homeomorphism from
$\frakV$ to $ U /p_1(G')$. It now follows from the construction that
$(\frakV,G',\psi',\zeta)$ is a local $\SSS^1$-standardization
of $Q$ at $v$. This proves that $Q$ defines an $\SSS^1$-fibration.
\EndProof
%\oldUpsilon\oldLambda\oldPsi Q q'_0 q_0' ^0 \eta 2\DD^{n-1} \zeta

\Proposition\label{fibration-category}Let $q:\oldLambda\to\frakB$ be a
smooth $\SSS^1$-fibration \abstractcomment{Do I need this for
  $I$-fibrations too? I don't think so} of a smooth, orientable,
compact, very good $3$-orbifold $\oldLambda$ over a smooth, compact,
very good $2$-orbifold $\frakB$. 
%Let $|\frakB|$ be equipped with a distance function that defines its
%topology, and let $\epsilon>0$ be given. 
Then there exist PL structures on $\oldLambda$ and $\frakB$, compatible with their smooth structures, and a sequence $(q_n')_{n\ge1}$, where $q_n':\oldLambda\to\frakB$ is a PL $\SSS^1$-fibration for each $n\ge1$, such that $(|q_n
'|)_{n\ge1}$ converges uniformly to $|q|$, and  for each compact PL subset $K$ of $|q|^{-1}(|\frakB|-\fraks_\frakB)$, the sequence $(|q_n'|\big|K)$ converges to $q|K$ in the $C^1$ sense.
\EndProposition

Note that in the statement of Proposition \ref{fibration-category},
the uniform convergence of $(|q'_n|)$ to $|q|$ implies that for each
compact PL subset $K$ of $|q|^{-1}(|\frakB|-\fraks_\frakB)$ we have
$|q_n'|(K)\subset |\frakB|-\fraks_\frakB$ for sufficiently large $n$,
so that the final assertion makes sense.

\Proof[Proof of Proposition \ref{fibration-category}]
We may assume without loss of generality that $\frakB$ is connected.

Since $q$ is an $\SSS^1$-fibration and $\oldLambda$ and $\frakB$ are very good, there exist finite-sheeted regular coverings $p_\oldLambda:\tL\to\oldLambda$ and $p_\frakB:\tB\to\frakB$ such that $\tL$ and $\tB$ are connected manifolds, and a locally trivial fibration $\tq:\tL\to\tB$ whose fiber is a (possibly disconnected) closed $1$-manifold, such that $p_\frakB\circ\tq=q\circ p_\oldLambda$. Since the fiber of $\tq$ is closed, $\tq$ is a boundary-preserving map. If $\calg$ and $\calx $ denote the covering groups of the coverings $p_\oldLambda$ and $p_\frakB$ respectively, $\tq$ is equivariant in the sense that there is a homomorphism $\rho:\calg\to \calx $ such that $\tq\circ g=\rho(g)\circ\tq$ for every $g\in \calg$. It follows from the main result of \cite{illman} that $\tL$ and $\tB$ have triangulations, compatible with their smooth structures, and invariant under the actions of $\calg$ and $\calx $ respectively. For the rest of the proof, the manifolds $\tL$ and $\tB$ will be understood to be equipped with the PL structures defined by these triangulations; note that the orbifolds $\oldLambda$ and $\frakB$ inherit PL structures, compatible with their smooth structures, and that $p_\oldLambda$ and $p_\frakB$ are then PL maps. After replacing the triangulation of $\tL$ by its first barycentric subdivision, we may assume:
\Claim\label{stable is fixed}
Each simplex of $\tL$ is pointwise fixed by its stabilizer under $\calg$.
\EndClaim

Let us fix a distance function $h$ on (the total space of) the tangent bundle $T\tB$ which determines its topology.  For any subset $K$ of $\tL$ which is a union of closed simplices, any piecewise smooth maps $f,g:K\to\tB$, and any $\delta>0$, we will say that $f$ and $g$ are {\it $C^1$ $\delta$-close} if for every closed simplex $\sigma\subset K$, every point $x\in\sigma$ and every tangent vector $w$ to $\sigma$ at $x$, we have $h(d(f|\sigma)(w),d(g|\sigma)(w))<\delta$.

For $k=-1,0,1,2,3$, let $\tL^{(k)}$ denote the union of all simplices of dimension at most $k$ in $\tL$. Note that since $\tL^{(k)}$ is $\calg$-invariant, it makes sense to speak of equivariant maps from $\tL$ to $\tB$. By induction, for $-1\le k\le3$, we will show:

\Claim\label{by induction}
For any $\delta>0$ there is an equivariant PL map $r^{(k)}:\tL^{(k)}\to\tB$ such that (1) $r^{(k)}$ is $\delta$-close to $\tq|\tL^{(k)}$, and (2) $r^{(k)}(\tL^{(k)}\cap\partial\tL)\subset\partial\tB$.
\EndClaim

Since \ref{by induction} is trivial for $k=-1$, we need only show that if $k$ is given with $0\le k\le3$, and if \ref{by induction} is true with $k-1$ in place of $k$, then it is true for the given $k$. Let $\delta'$ be a positive number less than $\delta$, which for the moment will be otherwise arbitrary; we will impose a finite number of smallness conditions on $\delta'$ in the course of the proof of \ref{by induction}. Let $r^{(k-1)}:\tL^{(k-1)}\to\tB$ be an equivariant PL map such that $r^{(k-1)}$ is  $\delta'$-close to $\tq|\tL^{(k-1)}$, and $r^{(k-1)}(\tL^{(k-1)}\cap\partial\tL)\subset\partial\tB$. 

Fix a complete set of orbit representatives $\cald$ for the action of $\calg$ on the set of open $k$-simplices of $\tL$. For each $\Delta\in\cald$, the stabilizer $\calg_\Delta$ fixes $\overline\Delta$ pointwise by \ref{stable is fixed}. Since $\tq$ and $r^{(k-1)}$ are equivariant, the sets $\tq(\overline{\Delta})$ and $r^{(k-1)}(\partial\overline{\Delta})$ are contained in $
%E_\Delta:=
\Fix(\rho(\calg_\Delta))$. If $\Delta$ is an element of $\cald$ such that $\Delta\subset\partial\tL$, then since $\tq$ is boundary-preserving, we have $\tq(\overline{\Delta})\subset\partial\tB$; furthermore, in this case we have $
r^{(k-1)}(\partial\overline{\Delta})\subset
r^{(k-1)}(\tL^{(k-1)}\cap\partial\tL)\subset\partial\tB$. Thus if we set
$E_\Delta=
\Fix(\rho(\calg_\Delta))$ for each $\Delta\in\cald$ such that $\Delta\not\subset\partial\tL$, 
and $E_\Delta=
\partial\tB\cap\Fix(\rho(\calg_\Delta))$ for each $\Delta\in\cald$ such that $\Delta\subset\partial\tL$, then
$\tq(\overline{\Delta})$ and $r^{(k-1)}(\partial\overline{\Delta})$ are contained in $E_\Delta$ for each $\Delta\in\cald$. 
 Note that $\Fix(\rho(\calg_\Delta))$ is a PL subset of $\tB$ for each
 $\Delta$, since the action of $\calx $ on $\tB$ is piecewise linear;
 since $\partial\tB$ is a PL set, it follows that  $E_\Delta$ is a PL subset of $\tB$ for each $\Delta\in\cald$. Note also that %since the action of $\calx $ on $\tB$ is piecewise linear; and that 
$r^{(k-1)}|\partial\overline{\Delta}$ is  $\delta'$-close to 
$\tq|\partial\overline{\Delta}$. If we choose $\delta'$ sufficiently small, it follows that $r^{(k-1)}|\partial\overline{\Delta}$ may be extended to a PL map $r_\Delta:\overline\Delta\to E_\Delta\subset\tB$, and that $\tq_\Delta$ may be taken to be arbitrarily $C^1$-close to $\tq|\overline{\Delta}$.

If $\Delta$ is any open  $k$-simplex of $\tL$, we may choose $g\in \calg$ so that $\Delta_0:=g^{-1}(\Delta)\in\cald$, and define  a PL map $ r_\Delta:\overline{\Delta}\to\tB$ by
$ r_\Delta=\rho(g)\circ r_{\Delta_0}\circ g^{-1}$. Since 
$r_{\Delta_0}(\overline{\Delta_0})\subset E_{\Delta_0}\subset \Fix(\rho(\calg_{\Delta_0}))$, the map $r_\Delta$ does not depend on the choice of $g$, and this definition of $r_\Delta$ specializes to the earlier one when $\Delta\in\cald$. Note  that if $\Delta$ is any $k$-simplex contained in $\partial\tL$, then choosing $g$ and defining $\Delta_0$ as above, we have $\Delta_0\subset\partial\tL$, and hence $E_\Delta\subset
 \partial\tB$; it follows that $ r_\Delta(\overline{\Delta})\subset\rho(g)(E_\Delta)\subset\partial\tB$.
Note also that since, for $\Delta\in\cald$, we may take $r_\Delta$ to
be arbitrarily  $C^1$-close to $r|\overline{\Delta}$ if $\delta'$ is
small enough, we may choose $\delta'$ so that $r_\Delta$ is $C^1$
$\delta$-close to $r|\overline{\Delta}$ for each $k$-simplex $\Delta$
of $\tL$.

We may extend $r^{(k-1)}$ to a PL map $r^{(k)}:\tL^{(k)}\to\tB$ by setting $r^{(k)}|\overline{\Delta}=r_{\Delta}$ for every
open $k$-simplex $\Delta$ of $\tL$. This map is equivariant by construction. Since
$r^{(k-1)}(\tL^{(k-1)}\cap\partial\tL)\subset\partial\tB$,  and
$ r_\Delta(\overline{\Delta})\subset\partial\tB$ for each $k$-simplex $\Delta\subset\partial\tL$, we have
$r^{(k)}(\tL^{(k)}\cap\partial\tL)\subset\partial\tB$.
Since $r_\Delta$ is $C^1$ $\delta$-close to $r|\overline{\Delta}$ for each $k$-simplex $\Delta$ of $\tL$, and 
 $r^{(k-1)}$ is $\delta$-close to $r|\tL^{(k-1)}$, the map  $r^{(k)}$ is $\delta$-close to $r|\tL^{(k)}$. This proves \ref{by induction}.

It follows from the case $k=3$ of \ref{by induction} that there is a
sequence $(\tq'_n)_{n\ge1}$ of boundary-preserving equivariant PL maps
from $\tL$ to $\tB$ which converges in the $C^1$ sense to
$\tq$. It follows from equivariance that for each $n\ge1$ there is a
unique PL map $Q_n:|\oldLambda|\to|\frakB|$ such that
$|p_\frakB|\circ\tq'_n=Q_n\circ |p_\oldLambda|$. 
%On the other hand, 

The $C^1$-convergence of $\tq'_n$ to $\tq$ implies in particular:
\Claim\label{when I first put this uniform on}
The sequence $(Q_n)$ converges uniformly to $|q|$, and  for each
compact PL subset $K$ of $|q|^{-1}(|\frakB|-\fraks_\frakB)$, the
sequence $(Q_n|K)$ converges to $q|K$ in the $C^1$
sense. 
\EndClaim

Since the smooth locally trivial fibration $\tq$ is in particular a
submersion, the $C^1$-convergence of $\tq'_n$ to $\tq$ also implies that
 $\tq'_n$ is a PL submersion  for sufficiently large $n$. Hence after
 truncating the sequence $(q'_n)$, we may assume:
%Since $\tq'_n$ is also boundary-preserving, we deduce:
\Claim\label{veriest dunce}
For every $n$, the boundary-preserving map
$\tq'_n:\tL\to\tB$ is a PL submersion.
\EndClaim

%\redcomment{The original statement here was about a locally trivial
  %fibration, which made it a PL version of Ehresmann's theorem. I'm
  %going to try to incorporate that in the next step.

%Set $B=|\frakB|-\fraks_\frakB$. Since $B$ is connected, it follows
%from \ref{veriest dunce} that the map
%$|q'_n|\big||q'_n|^{-1}(B):|q'_n|^{-1}(B)\to B$ is a locally trivial
%fibration for sufficiently large $n$, \redcomment{That seems wrong, since
  %we can guarantee immersedness only over compact subsets of $B$}  and that its fiber has the form $|\frakC|$ for some orbifold quotient of the fiber of $\tq'_n$; in particular, the fiber of $|q'_n|\big||q'_n|^{-1}(B)$ is a $1$-manifold (so that each of its components is homeomorphic to $\SSS^1$ or to $[0,1]$). We claim:

Now we claim:

\Claim\label{the rest of the story}
For every $n$, the PL map $Q_{n}:|\oldLambda|\to|\frakB|$
defines a PL orbifold fibration of $\oldLambda$ over $\frakB$ whose fibers are closed.
\EndClaim

%To prove \ref{the rest of the story}, note that, according to \ref{veriest dunce} and \ref{peripatetics}, for sufficiently large $n$, we have that (a) $\tq'_{n}:\tL\to\tB$ is a locally trivial fibration whose fiber is a closed $1$-manifold; and (b) the locally trivial fibration $|q'_{n}|\big||q'_{n}|^{-1}(B):|q'_{n}|^{-1}(B)\to B$ has connected  fiber. We will establish \ref{the rest of the story} by showing that $q':=q'_{n}:\oldLambda\to\frakB$ is an $\SSS^1$-fibration when (a) and (b) hold. 

To prove \ref{the rest of the story}, let $n$ be an arbitrary positive
integer and set $\tq'=\tq_n'$ and $Q=Q_n$. Then $|p_\frakB|\circ\tq'=Q\circ |p_\oldLambda|$. Let an arbitrary point $v\in\frakB$ be
given, and choose a point $\tv\in p_\frakB^{-1}(v)$. Fix a chart $\phi$ for $\frakB$ with domain $U_v$, such that
$\phi(0)=v$, let $\alpha$ denote the post-chart map defined by $\phi$,
and set $\frakV=\alpha(U_v)\subset\frakB$. Let 
$\tV$ denote the  component of $p_\frakB^{-1}(\frakV)$ containing $\tv$. The chart $\phi$ may be chosen so
that $\tV$ is an arbitrarily small neighborhood of $\tv$. There is a homeomorphism $\iota$
 of $U_v$ onto $\tV$  such that
$\alpha=(p_\frakB|\tV)\circ\iota$. 
%is identified with $\phi$; 
In
particular, if $X$ denotes the group associated with the chart
$\phi$, so that $\alpha$ induces
a homeomorphism of $U_v/X$ onto $\frakV$, then $\iota$ conjugates
$X$ onto the
stabilizer of $\tV$
in $\calx$, which we denote by $X_0$. 

Now set $N=(\tq')^{-1}(\tV)$ and $J=(\tq')^{-1}(\tv)$. Since 
$\tq':\tL\to\tB$ is a  boundary-preserving PL submersion, we may
assume, by taking $\tV$ to be a sufficiently small
neighborhood of $\tv$, that there is a PL homeomorphic 
identification of $N$ with $\tV\times J$ such that
$\tq'|N:N\to\frakV$ is the projection to the first factor. The
stabilizer of $N$ in $\calg$ is $G_0:=\rho^{-1}(X_0)$. Since $\rho$
is a surjective homomorphism, we have $\rho(G_0)=X_0$.

Let $\Gamma $ denote the stabilizer of $J=(\tq')^{-1}(\frakV)$ in
$G_0$. Since $\Gamma $ is a finite group of homeomorphisms of the closed
$1$-manifold $J$, we may fix a PL metric on $J$ for which each of its components
is isometric to $\SSS^1$, and $\Gamma \le\Isom(J)$. 
On the other hand, if $\tV$ is a sufficiently small
neighborhood of $\tv$, then the product structure on $N$ may be chosen
so that for each $u\in J$ and each $\gamma\in \Gamma $ we have 
$\gamma(\tV\times\{u\})=
\tV\times\gamma(\{u\})$. It now follows that $G_0\le
X_0\times\Gamma\le X_0\times\Isom(J)$, where $X_0\times\Isom(J)$, is
regarded in the natural way as a group of self-homeomorphisms of
$\tV\times J=N$. Furthermore, since $\rho(G_0)=X_0$, the
projection of $X_0\times\Isom(J)$ to its first factor maps $G_0\le
X_0\times\Isom(J)$ onto $X_0$.

The homeomorphism $\iota^{-1}\times\id$ of $N=\tV\times J$ onto $ U_v\times J$
conjugates $G_0$ onto a group $G\le X\times\Isom(J)$ whose projection
to the first factor is $X$. On the other hand, since 
  $|p_\frakB|\circ\tq'=Q\circ |p_\oldLambda|$, the restriction of
  $p_\oldLambda$ to $N\subset\oldLambda$ induces a homeomorphism of
  $N/G_0$ onto $\obd(Q^{-1}(|\frakV|)$. Precomposing the latter homeomorphism
  with the homeomorphism from $(U_v\times J)/G$ to $N/G_0$ induced by
  $\iota\times\id$, we obtain a homeomorphism $\psi: (U_v\times
  J)/G\to \obd(Q^{-1}(|\frakV|)$. If we now define $\zeta:\frakV\to
  U_v/X$ to be the inverse of the homeomorphism induced by $\alpha$, it follows from the constructions that
  $(\frakV, G,\psi,\zeta)$ is a
local $J$-standardization for $Q$ at
$v$. This proves \ref{the rest of the story}.
%\tfrakV
 
%, the action of $X$ pulls back via $\iota$ to an action on $U_v$
%is the group associated with the chart $\phi$, 
%\redcomment{I want to say that horizontal sections map to horizontal
  %sections, which allows us to identify $G_0$ with a subgroup of
  %$X_0\times{\rm Homeo}(J)$. Then there is a homeo of $U_v\times J$ onto
  %$N$ defined as $\iota\times\id$, which induces a homeo of
  %$(U_v\times J)/G$ onto $N/G_0$, where $G$ is the conjugate of $G_0$
  %via this homeo. Notation needs work, obviously. But $N/G_0$ is
  %identified with $Q^{-1}(\frakV)$ via the identity
  %$|p_\frakB|\circ\tq'=Q\circ |p_\oldLambda|$. Then the homeo of
  %$(U_v\times J)/G$ onto $N/G_0$ becomes $\psi$ in the local $J$-standardization for $Q$ at
%$v$. The $\frakV$ and $G$ are the ones already defined, and the $\zeta$ comes from $\iota$ or
%$\phi$. For junk, see not-connected.tex.}

%q'_n q_n' \oldLambda \frakB  $p_\frakB$ p_\oldLambda

For each $n\ge1$, we denote by $q'_n:\oldLambda\to\frakB$ the PL
orbifold fibration with closed fibers given by \ref{the rest of the story}.
We claim:

\Claim\label{peripatetics}
For sufficiently large $n$, the orbifold fibration $q'_n:\oldLambda\to \frakB$
is an $\SSS^1$-fibration.
\EndClaim

To prove \ref{peripatetics}, we choose a point $v_0\in|\frakB|-\fraks_\frakB$,
%Actually I
  %think any point in $\frakB$ may work. Fix up this paragraph, which
  %is a ``work in progress''}, and choose a
  %local $J$-standardization $s=(\frakV,G,\psi,\zeta)$ of $q$ at $v_0$, for some closed
  %$1$-manifold  $J$.  is a to be a 
%quadruple 
and fix closed disk neighborhoods $D_0$ and $D_1$ of $v_0$ in $B$, with
$D_0\subset\inter D_1$. By \ref{when I first put this 
  uniform on}, $(Q_n)$ converges uniformly to $|q|$; hence for sufficiently large $n$, say
for $n\ge n_0$, we have $|q'_n|^{-1}(v_0)\subset |q|^{-1}(D_0)\subset
%$, and
%$
|q'_n|^{-1}(D_1)$.
%\supset |q|^{-1}(D_0)$. 
But $|q|^{-1}(D_0)$ is connected since $q$ is an
$\SSS^1$-fibration. Hence for any $n\ge n_0$, the set $|q|^{-1}(D_0)$
is contained in a component   $P_n$ of
$|q'_n|^{-1}(D_1)$, and in particular we have $|q'_n|^{-1}(v_0)\subset
P_n$.  Thus $|q'_n|^{-1}(v_0)$ is the fiber of the map
$|q'_n|\big|P_n:P_n\to D_1$. Since $D_1\cap\fraks_\frakB=\emptyset$, and
$q'_n$ is an orbifold fibration, the
map 
$|q'_n|\big|P_n:P_n\to D_1$ is a 
 locally trivial
fibration. Since this fibration has a
connected total space and a simply connected base, it follows from the
exact homotopy sequence of this fibration that the fiber
$|q'_n|^{-1}(v_0)$ is connected. Thus for $n\ge n_0$, the
orbifold fibration
 $q'_{n}:\oldLambda\to\frakB$, whose fibers are closed, has a  connected
 fiber
$(q'_n)^{-1}(v_0)$.  Since $\frakB$ is connected, it now follows from
\ref{fibered stuff} that
$q'_{n}:\oldLambda\to\frakB$ is an $\SSS^1$-fibration for $n\ge n_0$. This proves \ref{peripatetics}.
%v D P

According to \ref{the rest of the story}, we may assume after
truncating the sequence $(q'_n)$ that $q'_n$ is an $\SSS^1$-fibration
for every $n$. This fact, together with \ref{when I first put this
  uniform on}, implies the conclusion of the proposition.
%The $C^1$-convergence of
%$\tq'_n$ to $\tq$ implies that
%for each
%compact PL subset $K$ of $|q'_n|^{-1}(|\frakB|-\fraks_\frakB)$, the
%sequence $(q_n'|K)_{n\ge1}$, converges to $q|K$ in the $C^1$ sense. As $q'_n$ is an $\SSS^1$-fibration
%for every $n$,
%this gives the conclusion of the proposition.
\EndProof
%\obd canoonical H \epsilon distance $v \zeta $t s N suborbifold
%G_1\DD^2 \toldLambda \tfrakB _{n} truncat G^\oldLambda G^\frakB \iota
%X X_0 G G_0 X ---> Xfuzz X_0 ---> X Xfuzz ---> X_0
% G ---> Gfuzz G_0 ---> G Gfuzz ---> G_0 K $P $u\gamma\Gamma $P K
% {first point} {second point} {third point} {oh yeah

\Lemma\label{ztimesz-lemma}
Let $\oldPsi$ be a weakly \simple\ $3$-orbifold such that every
component of $\partial\oldPsi$ is toric, and $\pi_1(\oldPsi)$ is infinite. Suppose that $\oldPsi$ admits no (piecewise linear) $\SSS^1$-fibration, and that no component of $\fraks_\oldPsi$ is $0$-dimensional. Then for every rank-$2$ free abelian subgroup $H$ of $\pi_1(\oldPsi)$, there is a component $\frakK$ of $\partial\oldPsi$ such that $H$ is contained in a conjugate of the image of the inclusion homomorphism $\pi_1(\frakK)\to\pi_1(\oldPsi)$. Furthermore, no finite-sheeted cover of $\oldPsi$ admits an $\SSS^1$-fibration.
%over a very good $2$-orbifold.
\EndLemma

\Proof
The weakly \simple\ orbifold $\oldPsi$ is by definition very
good. Hence, by the main result of \cite{lange}, we may write
$\oldPsi=\oldOmega \pl$ for some smooth orbifold $\oldOmega $. Since
$\oldPsi$ has only toric boundary components, and has no
$0$-dimensional components in its singular set, the same is true of
$\oldOmega $. Since $\oldPsi$ admits no PL $\SSS ^1 $-fibration, it
follows from Proposition \ref{fibration-category} that $\oldOmega $
admits no smooth $\SSS ^1 $-fibration  over a very good
$2$-orbifold. If  $\oldOmega $
were to admits a smooth $\SSS ^1 $-fibration  over a 
$2$-orbifold which is not very good, it would be a $3$-manifold
homeomorphic to a lens space, a contradiction to the hypothesis that
$\pi_1(\oldOmega)$ is infinite. Hence $\oldOmega $
admits no smooth $\SSS ^1 $-fibration.

If $\oldTheta$ is any smooth, closed $2$-suborbifold of  $\inter\oldOmega $, then by the first assertion of the main theorem of \cite{illman}, there is a PL structure on the very good orbifold $\oldOmega $, compatible with its smooth structure, for which $\oldTheta$ is a PL suborbifold. The second assertion of the main theorem of
\cite{illman} implies that $\oldOmega $, equipped with this PL structure, is PL homeomorphic to $\oldPsi$. This shows that  any smooth, closed $2$-suborbifold of  $\inter\oldOmega $ is (topologically) ambiently homeomorphic to a PL surface in $\oldOmega $. Since $\oldPsi$ is weakly \simple\ (and in particular irreducible), it now follows that (A) every $\pi_1$-injective, smooth, two-sided toric suborbifold $\oldTheta$ of $\inter\oldOmega $ is topologically boundary-parallel in $\oldOmega $, in the sense that $\oldTheta$ is the frontier of a suborbifold topologically homeomorphic to $\oldTheta\times[0,1]$; and (B) every two-sided smooth spherical  $2$-suborbifold $\oldTheta$ of $\inter\oldOmega $ bounds a topological discal $3$-suborbifold of $\oldOmega $. 

Since $\oldOmega $  admits no smooth $\SSS^1$-fibration, has only toric boundary components, has no $0$-dimensional components in its singular set, and has the properties (A) and (B) just stated, it follows from the Orbifold Theorem \cite{blp}, \cite{chk} (in the case where $\fraks_\oldOmega \ne\emptyset$) or from Perelman's geometrization theorem \cite{bbmbp}, \cite{Cao-Zhu}, \cite{kleiner-lott}, \cite{Morgan-Tian} (in the case where $\fraks_\oldOmega =\emptyset$) that
$\inter\oldOmega $ admits  a hyperbolic metric of finite volume. Hence  for every rank-$2$ free abelian subgroup $H$ of $\pi_1(\oldOmega )$, there is a component $\frakK$ of $\partial\oldOmega $ such that $H$ is contained in a conjugate of the image of the inclusion homomorphism $\pi_1(\frakK)\to\pi_1(\oldOmega )$. The first conclusion follows. 

To prove the second assertion, assume that
some finite-sheeted cover $\toldOmega$ of $\oldPsi$ admits a (PL) $\SSS^1$-fibration over a compact $2$-orbifold $\frakB$. The
finiteness of $\vol\oldPsi$ implies that $\pi_1(\frakB)$ has no
abelian subgroup of finite index; hence $\chi(\frakB)<0$. This implies
that $\frakB$ has a finite-sheeted cover $B$ which is an orientable
$2$-manifold with $\chi(B)\le-2$. There is a finite-sheeted cover
$M$ of $\toldOmega$ admitting a locally trivial fibration
$p:M\to B$ with fiber $\SSS^1$; in particular $M$ is a
manifold. Since $\chi(B)\le-2$, there is a simple closed curve
$C\subset\inter B$ which is not boundary-parallel in $B$. But then
the image of the inclusion homomorphism $\pi_1(p^{-1}(C))\to\pi_1(M)$
is a rank-$2$ free abelian subgroup of $\pi_1(M)$ which is not carried
by any component of $\partial M$. This contradicts the first
assertion. 
\EndProof
%Psi N \oldPsi ---> \frog   \oldOmega ---> \oldPsi   \frog ---> \oldOmega

\Lemma\label{affect}
Let $\oldLambda$ be a compact, orientable $3$-orbifold, let $\frakB$ be a
compact ${2}$-orbifold, and let $q:\oldLambda\to\frakB$ be an
$I$-fibration. Then $|q|:|\oldLambda|\to|\frakB|$ is a homotopy
equivalence. 
\EndLemma

\Proof
For each $x\in|\frakB|$ we have $|q|^{-1}(x)=|\frakI|$ for some fiber $\frakI$ of $q$. As observed in \ref{fibered stuff}, $\frakI$ is homeomorphic to either $[0,1]$ or $[[0,1]$, and hence $|\frakI|$ is a topological arc. In particular $|\frakI|$ is contractible and locally contractible; hence by the theorem of \cite{smale}, $|q|$ induces isomorphisms between homotopy groups in all dimensions. As $|\oldLambda|$ and $|\frakB|$ are triangulable, it follows from Whitehead's Theorem \cite[Theorem 4.5]{hatcherbook} that $|q|$ is a homotopy equivalence.
\EndProof

\Proposition\label{one for the books}
Let $\oldUpsilon$ be a very good, compact, connected, orientable $3$-orbifold such that  no component of $\partial\oldUpsilon$ is spherical. Let $\frakZ$ be a compact, connected  $2$-suborbifold of
$\partial\oldUpsilon$ with $\chi(\frakZ)\le0$. Suppose that the inclusion
homomorphism $\pi_1(\frakZ)\to\pi_1(\oldUpsilon)$ is an
isomorphism. Then   $\oldUpsilon$ admits a trivial $I$-fibration  under which $\frakZ$  is a component of $\partialh\toldLambda$.
%$(\oldLambda,\frakZ)$ 
%is
%(orbifold-)homeomorphic to $(\frakZ\times[0,1],\frakZ\times\{0\})$.
\EndProposition

\Proof
First consider the case in which $\oldUpsilon$ is a $3$-manifold. To emphasize that $\oldUpsilon$ and $\frakZ$ are manifolds, I will denote them by $W$ and $Z$ in this case. I claim:

\Claim\label{madigan}
$W$ is irreducible.
\EndClaim

To prove \ref{madigan}, suppose that $S\subset\inter W$ is a $2$-sphere. The hypothesis implies that the inclusion homomorphism $H_1(Z,\ZZ)\to H_1(W,\ZZ)$ is surjective; hence if $\alpha$ is any element of $H_1(Z,\ZZ)$, there is a (possibly singular) closed curve contained in $Z$ and representing $\alpha$. Since the curve is in particular disjoint from $S$, the intersection number of $\alpha$ with $[S]$ is $0$. Since this is true of every $\alpha\in H_1(W,\ZZ)$, the sphere $S$ separates $W$. 

We may therefore write $W=X\cup Y$, where $X$ and $Y$ are compact $3$-submanifolds of $W$ and $Z\cap Y=S$. We may label them in such a way that $Z\subset X$. The simple connectivity of $S$ implies that the inclusion homomorphisms from $\pi_1(X)$ and $\pi_1(Y)$ to $\pi_1(W)$ are injective; and that if we choose a base point in $S$ and identify $\pi_1(X)$ and $\pi_1(Y)$ isomorphically with their images under the inclusion homomorphisms, we have $\pi_1(W)=\pi_1(X)\star\pi_1(Y)$. In particular, $\pi_1(X)\cap\pi_1(Y)=\{1\}$. But by hypothesis the inclusion homomorphism $\pi_1(Z)\to\pi_1(W)$ is surjective; in particular
the inclusion homomorphism $\pi_1(X)\to\pi_1(W)$ is surjective, so that under the identifications we have $\pi_1(W)=\pi_1(X)$. Hence $\pi_1(Y)=\pi_1(X)\cap\pi_1(Y)=\{1\}$, i.e. $Y$ is simply connected. It then follows from a standard application of Poincar\'e-Lefschetz duality that every component of $\partial Y$ is a $2$-sphere. Since by hypothesis no component of $\partial W$ is a $2$-sphere, we have $\partial Y=S$. But it follows from the Poincar\'e Conjecture \cite{morgan-poincare} that a compact, simply connected $3$-manifold whose boundary is a single sphere is homeomorphic to a $3$-ball. Thus \ref{madigan} is proved.

According to \cite[Theorem 10.2]{hempel}, if $W$ is a compact, irreducible, orientable $3$-manifold, and $Z$ is a compact, connected $2$-manifold of
$\partial\oldUpsilon$ such that the inclusion
homomorphism $\pi_1(Z)\to\pi_1(W)$ is an
isomorphism, then the pair  $(W,Z)$ is homeomorphic to $(Z\times[0,1],Z\times\{0\})$. With \ref{madigan}, this completes the proof of the proposition in the manifold case.

For the proof in the general case, note that since $\oldUpsilon$ is
very good, there is a regular finite-sheeted cover $p:\toldUpsilon\to\oldUpsilon$ such that $\toldUpsilon$ is a $3$-manifold. Since $\frakZ$ is connected and the inclusion homomorphism $\pi_1(\frakZ)\to\pi_1(\oldUpsilon)$ is an isomorphism, $\tfrakZ:=p^{-1}(\frakZ)$ is also connected, and the inclusion homomorphism $\pi_1(\tfrakZ)\to\pi_1(\toldUpsilon)$ is also an isomorphism. By the manifold case of the proposition, which has already been proved, we may identify $\toldUpsilon$ with $\tfrakZ\times[0,1]$ in such a way that $\tfrakZ$ is identified with $\tfrakZ\times\{0\}$. Then the action of the covering group $G$ of the regular covering $p:\toldUpsilon\to\oldUpsilon$ is identified with an action on $\tfrakZ\times[0,1]$ that leaves $\tfrakZ\times\{0\}$ invariant.

Let $\frakR$ be a $G$-invariant regular neighborhood of $\tfrakZ$ in
$\partial\toldUpsilon$. By uniqueness of regular neighborhoods (see
\ref{new weak}), we may suppose the identification of $\toldUpsilon$ with $\tfrakZ\times[0,1]$ to be chosen so that $\frakR=(\tfrakZ\times\{0\})\cup((\partial\tfrakZ)\times[0,1])$. Then $\tfrakZ\times\{0\}$ and $\tfrakZ\times\{1\}$ are both $G$-invariant. Furthermore, $\frakZ_1:=p(\tfrakZ\times\{1\})$ is a suborbifold of $\partial\oldUpsilon$ disjoint from $\frakZ$.

We will now apply
 \cite[Theorem 8.1]{meeks-scott}: ``If $F$ is a compact surface, not $\SSS^2$ or $P^2$, and if $G$ is a finite group acting smoothly on $F \times [0,1]$ so as to preserve $F\times\{0,1\}$, then the action of $G$ is conjugate to an action which preserves the product structure.'' (Here it is understood that $F\times[0,1]$ is equipped with the smooth structure obtained from the product of the smooth structures on $F$ and $[0,1]$ by straightening the angle (see \cite[Theorem 7.5.3]{mukherjee}); it is readily shown that an action on $F\times\{0,1\}$,  which preserves the product structure is smooth in terms of this smooth structure.)

In order to apply  \cite[Theorem 8.1]{meeks-scott} in the present situation, we must make
the transition between the smooth and PL categories. First note that by Proposition \ref{yoga}, the pair 
$(\oldUpsilon,\frakZ\cup\frakZ_1)$  admits a smooth structure compatible with its PL structure. The pair 
$(\toldUpsilon,p^{-1}(\frakZ\cup\frakZ_1))=
(\tfrakZ\times[0,1],\tfrakZ\times\{0,1\})$ inherits a smooth structure
from  $(\oldUpsilon,\frakZ\cup\frakZ_1)$, which is compatible with its
own PL structure. By the equivalence of the PL and smooth categories
for $3$-manifolds, the pair
$(\tfrakZ\times[0,1],\tfrakZ\times\{0,1\})$ has only one smooth
structure compatible with its PL structure, and hence the smooth
structure which it inherits from $(\oldUpsilon,\frakZ\cup\frakZ_1)$ is
a product smooth structure (in the sense explained above). We may
therefore apply   \cite[Theorem 8.1]{meeks-scott},
% with $\tfrakZ$ playing the role of $F$, to deduce that the pair 
%$(\oldUpsilon,\frakZ\cup\frakZ_1)$, with the smooth structure that has %been assigned to it, is diffeomorphic to 
%$(\frakZ\times[0,1],\frakZ\times\{0,1\})$.
%
%We apply this result,
 taking $F=\tfrakZ$, and using the action of $G$ on $\tfrakZ\times[0,1]$ constructed above. Since $\tfrakZ\times\{0\}$ and $\tfrakZ\times\{1\}$ are both $G$-invariant, in particular their union is, and hence the action is conjugate to  an action which preserves the product structure. But since  each of the suborbifolds $\tfrakZ\times\{0\}$ and $\tfrakZ\times\{1\}$ is $G$-invariant, the conjugated action has the form $g\cdot(x,t)=(g\cdot x,t)$ for some action of $G$ on $\tfrakZ$. If we define  $\frakZ'$ to be the quotient of $\tfrakZ$ by this action of $G$ on $\frakZ$, it now follows that the pair $(\oldUpsilon,\frakZ\discup\frakZ_1)$, with the smooth structure that we have assigned to it, is diffeomorphic to $(\frakZ'\times[0,1],\frakZ'\times\{0,1\})$. But by
Proposition \ref{yoga}, $(\oldUpsilon,\frakZ\discup\frakZ_1)$ has, up to equivalence, only one PL structure compatible with its smooth structure; hence $(\oldUpsilon,\frakZ\discup\frakZ_1)$, with its original PL structure, is PL homeomorphic to $(\frakZ'\times[0,1],\frakZ'\times\{0,1\})$.
In particular, $(\oldUpsilon,\frakZ)$ is PL homeomorphic to
$(\frakZ'\times[0,1],\frakZ'\times\{0\})$, and therefore  to $(\frakZ\times[0,1],\frakZ\times\{0\})$.
\EndProof

\Proposition\label{robert strange}
Let $\oldLambda$ be a compact, orientable $3$-orbifold equipped with
an $I$-fibration over a
$2$-orbifold. Let $\frakZ$ be a $\pi_1$-injective, connected,
 two-sided
 $2$-suborbifold
of $\oldLambda$ such that $\chi(\frakZ)\le0$ and
$\partial\frakZ\subset\inter\partialv\oldLambda$. Then $\frakZ$ is
parallel in the pair $(\oldLambda,\partialv\oldLambda)$ (see \ref{parallel def}) either to 
some component of $\partialh\oldLambda$ (and hence to each component of the horizontal boundary of the component of $\oldLambda$ containing $\frakZ$), or to an annular suborbifold of $\partialv\oldLambda$. 
%\redmissingref{I need to add the def. of horizontal if it's not
%there. It means transverse to the fibers.} 
\EndProposition

\Proof 
We first observe that if $\oldLambda_0$ denotes the component of $\oldLambda$ containing $\frakZ$, and if
$\partialh\oldLambda_0$ has more than one component, then the $I$-fibration of $\oldLambda_0$ is trivial; hence in any event, if $\frakZ$ is parallel in   $(\oldLambda,\partialv\oldLambda)$  to one component of
    $\partialh\oldLambda$, it is parallel
in   $(\oldLambda_0,\partialv\oldLambda_0)$  to each component of
    $\partialh\oldLambda_0$.
 This justifies the parenthetical phrase in the conclusion of the proposition.

To prove the proposition, note that after replacing $\oldLambda$ by the component containing $\frakZ$, we
may assume that $\oldLambda$ is connected. The base of the
fibration of $\oldLambda$, which will be denoted by $\frakC$, is then
connected. Since $\frakZ$ has non-positive Euler characteristic and is $\pi_1$-injective,
$\pi_1(\oldLambda)$ is infinite. It follows that $\pi_1(\frakC)$ is
infinite, and hence that there is a finite-sheeted covering $\tfrakC$ of $\frakC$ which is an orientable $2$-manifold. Pulling back the $I$-fibration of $\oldLambda$ to $\tfrakC$, we obtain a finite-sheeted covering $\toldLambda$ of $\oldLambda$ which admits an $I$-fibration with base $\tfrakC$. Since $\toldLambda$ and $\tfrakC$ are both orientable, the fibration of $\toldLambda$ over $\tfrakC$ is a locally trivial fibration with fiber $[0,1]$; hence $\toldLambda$ is a $3$-manifold. This shows that  $\oldLambda$ is very good. Moreover, since $\pi_1(\tfrakC)$ is infinite, $\toldLambda$ is irreducible, and hence $\oldLambda$ is irreducible by Lemma \ref{cover-irreducible}.

Fix a component $\frakV$ of $\partialh\oldLambda$. Let $E$ denote the subgroup of $\pi_1(\oldLambda)$ defined, up to
conjugacy, as the image of the inclusion homomorphism
$\pi_1(\frakV)\to\pi_1(\oldLambda)$. Then $E$ has index at most
$2$ in $\pi_1(\oldLambda)$. 
If $\frakZ$ does not separate $\oldLambda$, then the image $A$ of the inclusion homomorphism
$\pi_1(\oldLambda-\frakZ)\to\pi_1(\oldLambda)$, which is also defined up to
conjugacy, is contained in the kernel of a homomorphism from
$\pi_1(\oldLambda)$ onto $\ZZ$, and hence has infinite index in
$\pi_1(\oldLambda)$; since $E$ is contained in a conjugate of $A $,
this is a contradiction. Hence $\frakZ$ separates $\oldLambda$. Let
$\oldUpsilon_0$ and $\oldUpsilon_1$ denote the closures of the components of
$\oldLambda-\frakZ$, indexed in such a way that $\frakV\subset\oldUpsilon_0$.
Fix a basepoint $\star\in\frakZ$, let $C$ denote the image
of the inclusion homomorphism
$\pi_1(\frakZ,\star)\to\pi_1(\oldLambda,\star)$, and for $i=0,1$, let $A_i$ denote the image
of the inclusion homomorphism
$\pi_1(\oldUpsilon_i,\star)\to\pi_1(\oldLambda,\star)$. Then
$\pi_1(\oldLambda,\star)$ may be
written as a free product with amalgamation
$A_0\star_CA_1$. But since $E$ is conjugate to a subgroup of $A_0$, the
index of $A_0$ in $\pi_1(\oldLambda)$ is at most $2$. In a non-trivial
free product with amalgamation, each factor has infinite index; hence
for some $i_0\in\{1,2\}$ we must have
$A_{i_0}=C$, i.e. the inclusion homomorphism
$\pi_1(\frakZ)\to\pi_1(\oldUpsilon_{i_0})$ is an isomorphism. Set
$\oldUpsilon=\oldUpsilon_{i_0}$, and set
$\oldUpsilon'=\overline{\oldLambda-\oldUpsilon}=\oldUpsilon_{3-i_0}$.

Thus the inclusion homomorphism
$\pi_1(\frakZ)\to\pi_1(\oldUpsilon)$ is an isomorphism. Furthermore,
$\oldUpsilon$ is very good since $\oldLambda$ has been seen to be very
good. Since $\oldLambda$ has been seen to be irreducible,
$\oldUpsilon$ is irreducible by Lemma \ref{oops lemma}. If some
component of $\partial\oldUpsilon$ were spherical, then by
irreducibility $\oldUpsilon$ would be discal, and
$\pi_1(\frakZ)\cong\pi_1(\oldUpsilon)$ would be finite; this is a
contradiction, since by hypothesis we have $\chi(\frakZ)\le0$. Hence $\partial\oldUpsilon$ has no spherical component.
We may now apply Proposition \ref{one for the books}
% since $\frakZ$ is non-discal, non-spherical and $\pi_1$-injective.
to deduce:
\Claim\label{watty piper}
The pair $(\oldUpsilon,\frakZ)$ is (orbifold-)
homeomorphic to $(\frakZ\times[0,1],\frakZ\times\{0\})$. 
\EndClaim

Set $\frakE=(\partial\oldUpsilon)-\inter\frakZ=\oldUpsilon\cap\partial\oldLambda$. Since $\frakZ$ is connected and $\pi_1$-injective in $\oldLambda$, it follows from
\ref{watty piper} that $\frakE$ is also connected and $\pi_1$-injective in $\oldLambda$. We have
$\partial\frakE=\partial\frakZ\subset\inter(\partialv\oldLambda)$. If
$\oldXi$ is any component of $\frakE\cap\partialv\oldLambda$, and
$\frakB$ is the (annular) component of $\partialv\oldLambda$
containing $\oldXi$, then each component of $\partial\oldXi$ either is
a component of $\partial\frakB$ or is contained in $\inter\frakB$. It
therefore follows from \ref{cobound} that for any component  $\oldXi$ of $\frakE\cap\partialv\oldLambda$, we have either (i) $\oldXi\subset\inter(\partialv\oldLambda)$, (ii) $\oldXi$ is a component of $\partialv\oldLambda$, or (iii) $|\oldXi|$ is a weight-$0$ annulus having one boundary component in $\inter(\partialv\oldLambda)$ and one in $\partial(\partialv\oldLambda)$.  Furthermore, in each case, $\oldXi$ is annular.

If Alternative (i) holds for some component $\oldXi$ of $\frakE\cap\partialv\oldLambda$, then since $\frakE$ is connected we have $\frakE\subset\partialv\oldLambda$. It then follows from \ref{watty piper} that
$\frakZ$ is
parallel in the pair $(\oldLambda,\partialv\oldLambda)$ to an annular suborbifold of $\partialv\oldLambda$; thus the conclusion of the proposition holds in this case. If Alternative (iii) holds for every component $\oldXi$ 
of $\frakE\cap\partialv\oldLambda$, then $\frakE\cap\partialv\oldLambda$ is a strong regular neighborhood of $\partial\frakE$ in the connected $2$-orbifold $\frakE$, and the complement of this regular neighborhood in $\frakE$ is the interior of a component $\frakZ'$ of $\partialh\oldLambda$. In view of \ref{watty piper}, it then follows that $\frakZ$ is
parallel in the pair $(\oldLambda,\partialv\oldLambda)$ to $\frakZ'$, and the conclusion holds in this case as well.

The rest of the proof will be devoted to the remaining case, in which $\frakE\cap\partialv\oldLambda$ has a component $\oldXi_0$ such that 
%Alternative (i) does not hold for any component $\oldXi$ of, and there is at least one of $\frakE\cap\partialv\oldLambda$ such that Alternative (iii) fails to hold with $\oldXi=\oldXi_0$.
%In particular, 
Alternative (ii) holds with $\oldXi=\oldXi_0$. Since
$\partial\frakE\subset\inter(\partialv\oldLambda)$, the component
$\oldXi_0$ of $\partial\oldLambda$ is contained in
$\inter\frakE$. Since the component $\oldXi_0$ of
$\partialv\oldLambda$ must meet each component of
$\partialh\oldLambda$, we have
$\partialh\oldLambda\subset\frakE$. Hence
$\frakE_0:=\oldXi_0\cup\partialh\oldLambda$ is connected and is
contained in $\frakE$. Since $\partial\frakE_0$ is a union of components of
$\partial(\partialh\oldLambda)$, it is disjoint from
$\partial\frakE=\partial\frakZ\subset\inter(\partialv\oldLambda)$. Hence $\frakE_0\subset\inter\frakE$.

Since $\frakZ$ is $\pi_1$-injective in $\oldLambda$, it follows from
\ref{watty piper} that $\frakE$ is $\pi_1$-injective in
$\oldLambda$. But $\partial\frakE_0$ is a union of components of
$\partial(\partialh\oldLambda)$, and is therefore $\pi_1$-injective in
$\oldLambda$. This implies that $\frakE_0$ is taut, and therefore
$\pi_1$-injective, in $\frakE$. Hence in this case we have:
\Claim\label{sam hillbilly}
$\frakE_0$ is $\pi_1$-injective in $\oldLambda$.
\EndClaim

We claim :
\Claim\label{climb}
The inclusion
homomorphism $\pi_1(\oldXi_0)\to\pi_1(\oldLambda)$ is an
isomorphism. Furthermore, if  $\frakD$ is a component of $\partialh\oldLambda$, the manifold
$|\frakD|$ is a weight-$0$ annulus, and $\frakD\cap\oldXi_0$ is
connected, i.e. is a single component of $\partial\frakD$.
\EndClaim

To prove \ref{climb}, let $d$
denote the index in $\pi_1(\oldLambda)$ of the image of the inclusion homomorphism
$\pi_1(\frakD)\to\pi_1(\oldLambda)$, so that $d=1$ if the
$I$-fibration of $\oldLambda$ is trivial, and $d=2$ otherwise. Let
$d_1$ and $d_2$ denote, respectively, the indices in
$\pi_1(\frakE_0)$  and $\pi_1(\oldLambda)$ of the images of the inclusion homomorphisms
$\pi_1(\frakD)\to\pi_1(\frakE_0)$ and
$\pi_1(\frakE_0)\to\pi_1(\oldLambda)$. It follows from \ref{sam
  hillbilly} that $d_1d_2=d$; in particular, $d_1$ and $d_2$ (which a
priori could be infinite) are at most $2$.

Consider the subcase in which the $I$-fibration
of $\oldLambda$ is trivial. In this  subcase it is immediate that $\frakD\cap\oldXi_0$ 
 is a single component of $\partial\frakD$.
If $\frakD'$ denotes the
component of $\partialh\oldLambda$ distinct from $\frakD$, we have
$\frakE_0=\frakD\cup\oldXi_0\cup\frakD'$. Using a basepoint in $\oldXi_0$,
we may write $\pi_1(\frakE_0)$ as a free product with amalgamation
$W\star_U W'$, where $U$, $W$ and  $W'$ respectively denote the images
of the injective inclusion homomorphisms from
$\pi_1(\oldXi_0)$, $\pi_1(\oldXi_0\cup\frakD)$ and $\pi_1(\oldXi_0\cup\frakD')$
to $\pi_1(\frakE_0)$. In this subcase $|\oldXi_0|$ is a weight-$0$ annulus, so that $U$
is cyclic; and $\frakD$ is homeomorphic to $\frakD'$, so that $W$ and
$W'$ are isomorphic. If $|\frakD|$ is not a weight-$0$ annulus, the factors $W$
and $W'$ are non-cyclic, and hence $\pi_1(\frakE_0)=W\star_UW'$ is a
non-trivial free product with amalgamation; this implies that the
index $d_1$ of $W$ in $\pi_1(\frakE_0)$ is infinite, a
contradiction. Hence $|\frakD|$ is a weight-$0$ annulus. The
weight-$0$ annuli $|\oldXi_0|$ and $|\frakD|$ share a boundary curve, so
that $\oldXi_0$ and $\frakD$ are isotopic in
$\partial\oldLambda$. Since $d=1$ in this subcase, i.e. the inclusion
homomorphism $\pi_1(\frakD)\to\pi_1(\oldLambda)$ is an isomorphism,  the inclusion
homomorphism $\pi_1(\oldXi_0)\to\pi_1(\oldLambda)$ is also an isomorphism. This  proves
\ref{climb} in this subcase.

Now consider the subcase in which the $I$-fibration
of $\oldLambda$ is non-trivial. In this  subcase we have  $d=2$; furthermore, we have
$\partialh\oldLambda=\frakD$, so that
$\frakE_0=\frakD\cup\oldXi_0$ and
$\frakD\cap\oldXi_0=\partial\oldXi_0$. Since the component $\oldXi_0$
of $\partialv\oldLambda$ is annular, $|\oldXi_0|$ is either a
weight-$0$ annulus or a weight-$2$ disk in which both singular points
have order $2$. If $|\oldXi_0|$ is  a
weight-$0$ annulus, then the intersection of the connected orbifolds
$\frakD$ and $\oldXi_0 $ has two components, and hence the image of
the inclusion homomorphism $\pi_1(\frakD)\to\pi_1(\frakE_0)$ has
infinite index in $\pi_1(\frakE_0)=\pi_1(\frakD\cup\oldXi_0)$; this
says that $d_1=\infty$, a contradiction. Hence $|\oldXi_0|$ must be a weight-$2$ disk in which both singular points
have order $2$. In particular, $\frakD\cap\oldXi_0=\partial\oldXi_0$
is a single
component of $\partial\frakD$.

Using a basepoint in $\partial\oldXi_0$,
we may therefore write $\pi_1(\frakE_0)$ as a free product with
amalgamation $X\star_V T$, where $X$, $V$ and  $T$ respectively denote the images
of the injective inclusion homomorphisms from
$\pi_1(\oldXi_0)$, $\pi_1(\partial\oldXi_0)$ and $\pi_1(\frakD)$
to $\pi_1(\frakE_0)$. The infinite cyclic group $V$ has index $2$ in
the infinite dihedral group $X$. If $V$ is a proper subgroup of $T$,
then
%case $|\oldXi_0|$ is a weight-$0$ annulus, so that $U$
%is cyclic; and $\frakD$ is homeomorphic to $\frakD'$, so that $W$ and
%$W'$ are isomorphci. If $|\frakD|$ is not a weight-$0$ annulus, the factors $W$
%and $W'$ are non-cyclic, and hence
$\pi_1(\frakE_0)=X\star_VT$ is a
non-trivial free product with amalgamation; this implies that the
index $d_1$ of $T$ in $\pi_1(\frakE_0)$ is infinite, a
contradiction. Hence $V=T$, which implies that $|\frakD|$ is a
weight-$0$ annulus.  The equality $V=T$, with the definitions of $d$
and $T$, also implies that
$d_1=|\pi_1(\frakE_0):T|=|X\star_VT:T|=|X:V|=2$, so that
$d_2=d/d_1=1$. In view of the definition of $d_2$, this
% $X=\pi_1(\frakE_0)$, which
means that the inclusion homomorphism $\pi_1(\frakE_0)\to\pi_1(\oldLambda)$
is an isomorphism. But since $V=T$ we have
$\pi_1(\frakE_0)=X\star_VT=X$, so that the inclusion homomorphism
$\pi_1(\oldXi_0)\to\pi_1(\frakE_0)$
is an isomorphism, and hence so is the inclusion homomorphism $\pi_1(\oldXi_0)\to\pi_1(\oldLambda)$.  Thus 
\ref{climb} is proved in all subcases. 

It follows from the second assertion of \ref{climb} that
$\frakE_0$ is a regular neighborhood of $\oldXi_0$. Hence, by the first
assertion of \ref{climb},  the inclusion
homomorphism $\pi_1(\frakE_0)\to\pi_1(\oldLambda)$ is an
isomorphism.  Using Proposition \ref{one for the books}, we
now deduce:

\Claim\label{chili if you want it}
The pair $(\oldLambda, \frakE_0)$ is homeomorphic to
$(\oldXi_0\times[0,1], 
\oldXi_0\times\{0\} )$.
\EndClaim

Since $\oldXi_0$ is annular, it follows from \ref{chili if you want it} that $\oldLambda$ is a
\torifold, and in particular it is
boundary-reducible. Since $\frakE$ is $\pi_1$-injective in
$\oldLambda$ we must have $\frakE\ne\partial\oldLambda$, and hence $\partial\frakE\ne\emptyset$. It also follows from
\ref{chili if you want it} that
$\oldXi_1:=\overline{(\partial\oldLambda)-\frakE_0}$ is homeomorphic
to $\oldXi_0$, and in particular is connected. Hence $\oldXi_1$ is a
single component of $\partialv\oldLambda$. Since
$\frakE_0\subset\inter\frakE$, we have
$\emptyset\ne\partial\frakE\subset\inter\oldXi_1$. Since $\oldXi_1$ is
annular, it follows from \ref{cobound} that every component of
$\partial\frakE$ cobounds an annulus with a component of
$\partial\oldXi_1$. Since $\frakE$ is connected, it now follows that
%\Claim\label{funny color}
%The orbifold
$\frakE$ is a regular neighborhood of
$\frakE_0$ in $\partial\oldLambda$.
%\EndClaim
With
\ref{chili if you
  want it}, this implies that
the pair $(\oldLambda, \frakE)$ is homeomorphic to
$(\frakE\times[0,1], 
\frakE\times\{0\} )$.
%\redcomment{Fix from this point. I should replace \ref{funny color} by the
%statement that $\frakE$ is a regular neighborhood of $\frakE_0$. With
%the new version of  \ref{chili if you
 % want it}, this will give that
%the pair $(\oldLambda, \frakE)$ is homeomorphic to
%$(\frakE\times[0,1], 
%\frakE\times\{0\} )$.
But it follows from \ref{watty piper} that
$(\oldLambda, \frakE)$ is homeomorphic to $(\oldUpsilon',\frakZ)$
(where as above we set
$\oldUpsilon'=\overline{\oldLambda-\oldUpsilon}$). Hence
$(\oldUpsilon',\frakZ)$ is homeomorphic to $(\frakZ\times[0,1], 
\frakZ\times\{0\} )$, i.e. 
% The product
% structure of the latter pair implies (clarify this) that
$\frakZ$ is
parallel to a suborbifold of $\oldUpsilon'\cap\partial\oldLambda=
\overline{(\partial\oldLambda)-\frakE}\subset\oldXi_1\subset\partialv\oldLambda$. 
This gives the conclusion of the proposition in this case.
\EndProof
%\toldLambda

\Corollary\label{covering annular} 
Let $p:\toldPsi\to \oldPsi$ be a covering map
of orientable $3$-orbifolds that are 
componentwise strongly \simple. Let $\oldXi$ be a $\pi_1$-injective suborbifold of $\partial\oldPsi$. Then for any essential  annular $2$-orbifold $\oldPi$ in the pair $(\oldPsi,\oldXi)$, every component of $p^{-1}(\oldPi)$ is
essential in the pair $(\toldPsi,p^{-1}(\oldXi))$. Furthermore, if
$\frakR,\frakR'\subset\oldPsi$ are
 annular $2$-orbifolds which are both essential in
$(\oldPsi,\oldXi)$ and are not parallel, then no component of $p^{-1}(\frakR)$
is parallel in $(\toldPsi,p^{-1}(\oldXi))$ to any component of $p^{-1}(\frakR')$.
\EndCorollary

\Proof
Set $\toldXi=p^{-1}(\oldXi)$.

To prove the first assertion, suppose that
 $\oldPi$ is an
essential orientable annular $2$-orbifold in the pair $(\oldPsi,\oldXi)$. Then in particular $\oldPi$ is $\pi_1$-injective in $\oldPsi$;
and since $p$ is an orbifold covering map,
$(p|\toldPi)_\sharp:\pi_1(\toldPi)\to\pi_1(\oldPi)$ is injective for each component $\toldPi$ of  $p^{-1}(\oldPi)$. It follows
that $p^{-1}(\oldPi)$ is $\pi_1$-injective in $\toldPsi$. To establish the first assertion, it remains to show:

\Claim\label{go kurds}
No
component of $p^{-1}(\oldPi)$ is parallel
 in the pair $(\toldPsi,\toldXi)$ either to a suborbifold of $\toldXi$ or to a component of
$\overline{\partial\toldPsi-\toldXi}$. 
\EndClaim

Assume that \ref{go kurds} is false. 
Then there is a connected $3$-suborbifold $\toldLambda$ of $\toldPsi$ such that (i) $\Fr_{\toldPsi} \toldLambda$ is a single component of $p^{-1}(\oldPi)$, and (ii) $\toldLambda$ admits a trivial $I$-fibration under which $\Fr_{\toldPsi}\toldLambda$ is a component of $\partialh\toldLambda$,  the other component of $\partialh\toldLambda$ is contained either in  $\toldXi$ or in 
$\overline{\partial\toldPsi-\toldXi}$, and
$\partialv\toldLambda\subset\toldXi$.
We may assume that, among all $3$-suborbifolds of $\toldLambda$ for which (i) and (ii) hold, $\toldLambda$ is minimal with respect to inclusion. 
Set $\toldPi_0=\Fr_{\toldPsi}\toldLambda$, so that $\toldPi_0$ is a  component of $p^{-1}(\oldPi)$.
Let $\tfrakZ_0$ denote the  component of $\partialh\toldLambda$ that is contained either in  $\toldXi$ or in 
$\overline{\partial\toldPsi-\toldXi}$.

%As a preliminary to the proof of \ref{go kurds}, 
Let us set $\frakC=\overline{(\partial\toldLambda)-\toldPi_0}=(\partialv\toldLambda)\cup\tfrakZ_0=\toldLambda\cap\partial\toldPsi$.
Since the component $\toldPi_0$ of $\partialh\toldLambda$ is annular, $\toldLambda$ is a \torifold, and in particular $\partial\toldLambda$ is toric. Since $p^{-1}(\oldPi)$ is $\pi_1$-injective in $\toldPsi$, the annular orbifold $\toldPi_0$ is $\pi_1$-injective in $\toldPsi$ and hence in the orientable toric orbifold $\partial\toldLambda$; it then follows that $\frakC=\overline{
(\partial\toldLambda)-\toldPi_0}$ is annular.

Since
$\toldPi_0=\Fr_{\toldPsi}\toldLambda$ is a  component of $p^{-1}(\oldPi)$, either (a) there is a component $\toldPi_1$ of $p^{-1}(\oldPi)$ contained in $\toldLambda$ and disjoint from $\toldPi_0$, or (b)
$\toldLambda\cap p^{-1}(\oldPi)=\toldPi_0$. In each case we will obtain a contradiction.

First suppose that (a) holds. 
We have $\partial\toldPi_1\subset\inter\frakC$.  The
$\pi_1$-injectivity of $p^{-1}(\oldPi)$  in $\toldPsi$  implies that
each boundary curve of $\toldPi_1$ is $\pi_1$-injective in the
orientable annular orbifold $\frakC$. It then follows from
\ref{cobound} that each component of $\partial\toldPi_1$ cobounds a
weight-$0$ annulus in $\frakC$ with some component of
$\partial\frakC$. Hence $\toldPi_1$ is annular, and after modifying the $I$-fibration of $\toldLambda$ by an isotopy, we may assume that $\partial\toldPi_1\subset\inter\partialv\toldLambda$. It then follows from Proposition \ref{robert strange} that $\toldPi_1$ is parallel in the pair $(\toldLambda,\partialv\toldLambda)$ either to $\tfrakZ_0$ or to a suborbifold of $\partialv\toldLambda$. Hence there is a suborbifold $\toldLambda_1$ of $\toldLambda$, admitting a trivial $I$-fibration, such that
$\partialv\toldLambda_1\subset \partialv\toldLambda\subset \toldXi$,
while $\toldPi_1=\Fr_{\toldPsi}\toldLambda_1$ is one component of $\partialh\toldLambda_1$,  and the remaining component $\tfrakZ_1$ of $\partialh\toldLambda_1$ is either  $\tfrakZ_0$ or  a suborbifold of $\partialv\toldLambda\subset\toldXi$. Thus $\tfrakZ_1$ 
is contained either in  $\toldXi$ or in 
$\overline{\partial\toldPsi-\toldXi}$.
% or in $\partialv\toldLambda\subset\overline{\toldLambda-\toldXi}$. \oldGamma
Since  $\toldPi_1= \Fr_\oldPsi\toldLambda_1$ is disjoint from $\toldPi_0= \Fr_\oldPsi\toldLambda$,
%$\partialh\toldLambda_0=\Fr_{\toldPsi}\toldLambda_1$, 
the suborbifold $\toldLambda_1$ of $\toldLambda$ is proper. This contradicts the minimality of $\toldLambda$.
%, and \ref{sassafrass} is proved.
%frakK

 Now suppose that (b) holds. Then $\toldLambda$ is the closure of a component of $\toldPsi-p^{-1}(\oldPi)$. Hence we may write $\toldLambda=p^{-1}(\oldLambda)$, where $\oldLambda$ is the closure of some component of $\oldPsi-\oldPi$. Thus $q:=p|\toldLambda:\toldLambda\to\oldLambda$ is a(n orbifold) covering. Since (b) holds, we have $q^{-1}(\oldPi)=\toldPi$. Since $\oldLambda$ is covered by the \torifold\ $\toldLambda$, it follows from the definition that $\oldLambda$ is itself a \torifold. Since 
$\toldPi$ is connected, and since the inclusion homomorphism $\pi_1(\toldPi)\to\pi_1(\toldLambda)$ is an isomorphism by (ii), the inclusion homomorphism $\pi_1(\oldPi)\to\pi_1(\oldLambda)$ is also an isomorphism. 
It now follows from Proposition \ref{one for the books}, applied with $\oldLambda$ and $\oldPi$ playing the roles of $\oldUpsilon$ and $\frakZ$, that $\oldLambda$ admits a trivial $I$-fibration under which $\oldPi$ is a component of $\partialh\oldLambda$. (The hypotheses in Proposition \ref{one for the books} that $\oldUpsilon$ is very good and has no spherical boundary components hold here because $\oldLambda$ is a \torifold.)
Let $\frakZ$ denote the component of $\partialh\oldLambda$ distinct from $\oldPi$, so that $\oldLambda\cap\partial\oldPsi=\frakZ\cup\partialv\oldPsi$.

%Let
%$\tfrakZ_0$ denote the component of $\partialh\toldLambda$ distinct from $\toldPi$,
%under the $I$-fibration of $\toldLambda$ given by  Condition (ii). 
According to Condition (ii) we have $\partialv\toldLambda\subset\toldXi$, and $\tfrakZ_0$ is contained either in $\toldXi$ or in $\partial\toldPsi-\inter\toldXi$. If 
$\tfrakZ_0\subset\toldXi$ then $\toldLambda\cap\partial\toldPsi\subset\toldXi$; hence
$\oldLambda\cap\partial\oldPsi\subset\oldXi$, i.e. $\partialv\oldLambda\cup\frakZ\subset\oldXi$. Thus in this subcase
$\oldPi$ is parallel in the pair $(\oldPsi,\oldXi)$ to a suborbifold of $\oldXi$,  a contradiction to the essentiality of $\oldPi$.

If
$\tfrakZ_0\subset(\partial\toldPsi)-\inter\toldXi$, then $\tfrakZ_0=\toldLambda\cap((\partial\toldPsi)-\inter\toldXi)$. But $\tfrakZ_0$ is connected, and the inclusion homomorphism $\pi_1(\tfrakZ_0)\to\pi_1(\toldLambda)$ is an isomorphism. Hence 
$\frakV:=\oldLambda\cap\overline{(\partial\oldPsi)-\oldXi}$ is connected, and the inclusion homomorphism $\pi_1(\frakV)\to\pi_1(\oldLambda)$ is an isomorphism. But we have $\frakV\subset\oldLambda\cap\partial\oldPsi=(\partialv\oldLambda)\cup\frakZ$. The orbifold $(\partialv\oldLambda)\cup\frakZ$ is connected, and the inclusion homomorphism $\pi_1((\partialv\oldLambda)\cup\frakZ)\to\pi_1(\oldLambda)$ is an isomorphism. 
The inclusion homomorphism $\pi_1(\frakV)\to\pi_1((\partialv\oldLambda)\cup\frakZ)$ is therefore also an isomorphism, so that $(\partialv\oldLambda)\cup\frakZ$ is a 
regular neighborhood of $\frakV$ in $\partial\oldPsi$.
As $(\partialv\oldLambda)\cup\frakZ$ is also a
regular neighborhood of $\frakZ$ in $\partial\oldPsi$, it follows that 
$\frakV$ and $\frakZ $ are isotopic in $(\partialv\oldLambda)\cup\frakZ$. Hence after modifying the $I$-fibration of $\oldLambda$ we may assume that 
$\frakV=\frakZ $. This means that
$\oldPi$ is parallel in the pair $(\oldPsi,\oldXi)$ to the component $\frakV$ of
$\overline{\partial\oldPsi-\oldXi}$, and we again have a contradiction to the essentiality of $\oldPi$. This completes the proof of \ref{go kurds} and thus establishes the first assertion of the corollary.
%\oldLambda

We now turn to the proof of the second assertion. Suppose that $\frakR,\frakR'\subset\oldPsi$ are two-sided
 annular $2$-orbifolds which are both essential in
$\oldPsi$ and are not parallel. Suppose that some component of $p^{-1}(\frakR)$
is parallel in $\toldPsi$ to some component  of $p^{-1}(\frakR')$. Then there is a connected $3$-suborbifold $\tfrakK$ of $\toldPsi$ admitting a trivial $I$-fibration such that (i) $\partialv\tfrakK\subset\toldXi$, and  (ii) $\partialh\tfrakK$ has one component contained in $p^{-1}(\frakR)$
and one contained in $p^{-1}(\frakR')$. Among all 
connected $3$-suborbifolds $\tfrakK$ of $\toldPsi$ that admit a trivial $I$-fibration for which (i) and (ii) hold, choose one, say $\tfrakK_0$, which is minimal with respect to inclusion. Let $\tfrakR_0\subset p^{-1}(\frakR)$ and $\tfrakR_0'\subset p^{-1}(\frakR')$ denote the components of
 $\partialh\tfrakK_0$. We claim:

\Claim\label{sassafrass}
We have $\tfrakK_0\cap p^{-1}(\frakR\cup\frakR')=\tfrakR_0\cup\tfrakR_0'$.
\EndClaim

To prove \ref{sassafrass}, suppose that $\tfrakK_0\cap
p^{-1}(\frakR\cup\frakR')$ has a component distinct from $\tfrakR_0$
and $\tfrakR_0'$. By symmetry we may assume this component is
contained in $p^{-1}(\frakR)$, and we shall denote it by
$\tfrakR_1$. Since $\frakR$ is essential in $\oldPsi$, it follows from
the first assertion of the present corollary, which has already been
proved, that $\tfrakR_1$ is essential in $\toldPsi$. In particular
$\tfrakR_1$ is $\pi_1$-injective in $\toldPsi$ and hence in
$\tfrakK_0$. Since $\frakR$ and $\frakR'$ are annular, so is
$\tfrakR_1$. On the other hand, since $\tfrakR_1$ is two-sided in $\toldPsi$, contained in $\tfrakK_0$ and disjoint from $\tfrakR_0\cup\tfrakR_0'=\partialh\tfrakK_0$, we have $\partial\tfrakR_1\subset\inter\partialv\tfrakK_0$. It therefore follows from Proposition \ref{robert strange} that $\tfrakR_1$ is parallel in the pair $(\tfrakK_0,\partialv\tfrakK_0)$ either to $\tfrakR_0'$ or to a suborbifold of $\partialv\tfrakK_0$. But the latter alternative would imply that $\tfrakR_1$ is parallel in $\toldPsi$ to a suborbifold of $\toldXi$, a contradiction to the essentiality of $\tfrakR_1$. Hence 
$\tfrakR_1$ is parallel to $\tfrakR_0'$ in the pair $(\tfrakK_0,\partialv\tfrakK_0)$. This means that there is a connected $3$-suborbifold $\tfrakK_1$ 
of $\toldPsi$ admitting a trivial $I$-fibration such that (i) $\partialv\tfrakK_1\subset\tfrakR$, and  the components of 
$\partialh\tfrakK_1$ are $\tfrakR_1\subset p^{-1}(\frakR)$
and  $\tfrakR_0'\subset p^{-1}(\frakR')$. Since the component $\tfrakR_1$ of $\Fr_{\toldPsi}\tfrakK_1$ is disjoint from $\partialh\tfrakK_0=\Fr_{\toldPsi}\tfrakK_1$, the suborbifold $\tfrakK_1$ of $\tfrakK_0$ is proper. This contradicts the minimality of $\tfrakK_0$, and \ref{sassafrass} is proved.

It follows from \ref{sassafrass} that $\tfrakK_0$ is the closure of a component of $\toldPsi-p^{-1}(\frakR\cup\frakR')$. Since $p:\toldPsi\to\oldPsi$ is a covering map, it follows that there is a suborbifold $\frakK_0$ of $\oldPsi$, which is the closure of a component of $\oldPsi-(\frakR\cup\frakR')$, such that $\tfrakK_0$ is a component of $p^{-1}(\frakK_0)$, and that $p_0=p|\tfrakK_0:\tfrakK_0\to\frakK_0$ is a covering map. Since $\tfrakR_0=p_0^{-1}(\frakR)$ and $\tfrakR_0'= p_0^{-1}(\frakR')$ are the components of $\Fr_{\toldPsi}\tfrakK_0$, the components of $\Fr_\oldPsi\frakK_0$ are $\frakR$ and $\frakR'$. 

We have $p_0^{-1}(\frakK_0\cap\partial\oldPsi)
=\tfrakK_0\cap\partial\toldPsi=\partialv\tfrakK_0\subset\toldXi$.
%whose components are annular orbifolds. 
Thus each component of $\frakK_0\cap\partial\oldPsi$ is covered by an annular orbifold and is therefore annular. On the other hand, if $\frakB$ is a component of $\frakK_0\cap\partial\oldPsi$, and if we choose a component $\tfrakB$ of $p_0^{-1}(\frakB)$, then the triviality of the $I$-fibration of $\tfrakK_0$ implies that the component $\tfrakB$ of $\partialv\tfrakK_0$ has one boundary component in $\partial\tfrakR_0$ and one in $\partial\tfrakR_0'$. Hence $\frakB$ has boundary components in both $\partial\frakR$ and $\partial\frakR'$. In particular the orientable annular orbifold $\frakB$ has two distinct boundary components, so that $|\frakB|$ is a weight-$0$ annulus with one boundary component in $|\partial\frakR|$ and one in $|\partial\frakR'|$.

Since $p_0^{-1}(\frakR)=\tfrakR_0$ is connected, and since the inclusion homomorphism $\pi_1(\tfrakR_0)\to\pi_1(\tfrakK_0)$ is an isomorphism in view of the triviality of the $I$-fibration of $\tfrakK_0$, the inclusion homomorphism $\pi_1(\frakR)\to\pi_1(\frakK_0)$ is also an isomorphism. 
Furthermore, $\frakK_0$ is very good since the componentwise strongly \simple\ $3$-orbifold $\oldPsi$ is automatically very good (see \ref{oops}). Since strong \simple ity also implies that $\oldPsi$ is irreducible, $\frakK_0$ is irreducible by Lemma \ref{oops lemma}. If some component of $\partial\frakK_0$ were spherical, then by irreducibility $\frakK_0$ would be discal, and so $\pi_1(\frakK_0)$ would be finite; but $\frakK_0$ contains the annular orbifold $\frakR$, which is essential and therefore $\pi_1$-injective in $\oldPsi$, and hence $\pi_1(\frakK_0)$ is infinite. This shows that $\partial\oldLambda$ has no spherical component. 
It now follows from Proposition \ref{one for the books} that $\frakK_0$ admits a trivial $I$-fibration  under which $\frakR$  is a component of $\partialh\frakK_0$.

Since $\partial\frakR\subset\partial(\frakK_0\cap\partial\oldPsi)$, and since each component of $|\frakK_0\cap\partial\oldPsi|$ is a weight-$0$ annulus having exactly one boundary component in common with $|\frakR|$, we may choose the trivial $I$-fibration of $\frakK_0$ in such a way that $\partialv\frakK_0=\frakK_0\cap\partial\oldPsi$. Hence the component of $\partialh\frakK_0$ distinct from $\frakR$ is $\overline{(\partial\frakK_0)-(\frakR\cup(\frakK_0\cap\partial\oldPsi))}=\frakR'$. It now follows that $\frakR$ and $\frakR'$ are parallel in $\oldPsi$, a contradiction. Thus the second assertion of the corollary is proved.
 \EndProof
%\obd\frakK\frakB\partial\toldPsi\frakR\tfrakR oldPi\toldPi\oldUpsilon\oldGamma\tfrakZ

\Corollary\label{haken annuli}
For every strongly \simple, orientable, boundary-irreducible $3$-orbifold $\oldPsi$,
there is a natural number $N$ with the following property: if $\frakC\subset\oldPsi$ is any two-sided $2$-orbifold
such that each component of $\frakC$ is an essential annular suborbifold of $\oldPsi$, and no two components of $\frakC$ are parallel in $\oldPsi$, then $\compnum(\frakC)\le N$.
\EndCorollary

\Proof
If we replace the hypothesis that $\oldPsi$ is a strongly \simple, orientable, boundary-irreducible $3$-orbifold  by the hypothesis that $\oldPsi$ is an irreducible $3$-manifold, this is a special case of Haken's finiteness theorem (see \cite[Lemma 13.2]{hempel}). %\redcommentD{But  he assumes irreducibility. Just checked that. The underlying issues here include whether-if strong \simple ity is retaioned a hypothesis---whether it is preserved by passing to the covering. This kind of thing is ultimately really irrelevant, because strong \simple ity is probably used only to guarantee very good-ness. But I need to know whether, once the statement and apps have been corrected (which will probably involve removing the hyp. of strong \simple ity), the reference works }. 
To prove the corollary for a strongly \simple, orientable, boundary-irreducible $3$-orbifold $\oldPsi$, we first note that in view of the definition of  strong \simple ity (\ref{oops}), we may fix a finite-sheeted covering $p:\toldPsi\to\oldPsi$ such that $\toldPsi$ is an irreducible $3$-manifold. We may then fix a natural number $N$ with the  property that if $\tfrakC\subset\toldPsi$ is any two-sided $2$-manifold
such that each component of $\tfrakC$ is an essential annulus in $\toldPsi$, and no two components of $\tfrakC$ are parallel in $\toldPsi$, then $\compnum(\tfrakC)\le N$. Now let
$\frakC\subset\oldPsi$ be a two-sided $2$-orbifold
such that each component of $\frakC$ is an essential annular suborbifold of $\oldPsi$, and no two components of $\frakC$ are parallel in $\oldPsi$. Let $\frakC_1,\ldots,\frakC_n$ denote the components of $\frakC$. For $i=1,\ldots,n$, choose a component $\tfrakC_i$ of $p^{-1}(C_i)$. By Corollary \ref{covering annular}, each $\tfrakC_i$ is an essential annulus in $\toldPsi$, and for any two distinct indices $i,j\in\{1,\ldots,n\}$, the annuli $\tfrakC_i$ and $\tfrakC_j$ are non-parallel in $\toldPsi$. Hence $n\le N$.
\EndProof
%\obd

%\toldXi\tfrakA\obd

\Corollary\label{still essential}
Let $\oldPsi$ be an orientable $3$-orbifold which is componentwise strongly \simple\ and componentwise boundary-irreducible. Let $\oldPsi'$ be a
suborbifold of  $\oldPsi$ whose frontier components are essential annular suborbifolds of $\oldPsi$ (so that
$\oldPsi'$ is componentwise strongly \simple\ by Lemma \ref{oops lemma}, and  $\oldPsi'\cap\partial\oldPsi$ is
$\pi_1$-injective in $\partial\oldPsi$). Suppose that $\oldPi$ is an essential annular orbifold in the pair $(\oldPsi', \oldPsi'\cap\partial\oldPsi)$. Then $\oldPi$ is essential in $\oldPsi$.
\EndCorollary

\Proof
Since $\oldPi$ is two-sided in $\oldPsi'$, and since $\partial\oldPi\subset\partial\oldPsi$, we deduce that $\oldPi$ is two-sided in $\oldPsi$. Furthermore, $\oldPi$ is $\pi_1$-injective  in $\oldPsi$, since $\oldPi$ is $\pi_1$-injective in  $\oldPsi$ and since $\Fr_\oldPsi\oldPsi'$ is 
$\pi_1$-injective in $\oldPsi$. Now suppose that $\oldPi$ is parallel in $\oldPsi$ to a suborbifold $\frakC$ of $\partial\oldPsi$. Then there is a connected suborbifold $\oldLambda$ of $\oldPsi$ admitting a trivial $I$-fibration under which $\oldPi$ and $\frakC$ are the components of $\partialh\oldLambda$, and $\partialv\oldLambda\subset\partial\oldPsi$. If $\oldLambda\subset\oldPsi'$, then  $\frakC\subset\oldPsi'\cap\partial\oldPsi$, and $\oldPi$ is parallel to $\frakC$ in the pair $(\oldPsi', \oldPsi'\cap\partial\oldPsi)$; this contradicts the hypothesis that $\oldPi$ is essential in this pair. There remains the case in which $\oldLambda\not\subset\oldPsi'$. Then $\frakB:=(\Fr_\oldPsi\oldPsi')\cap\oldLambda$ is a non-empty union of components of $\Fr_\oldPsi\oldPsi'$. Furthermore, $\frakB$ is two-sided in $\oldLambda$, and $\partial\frakB\subset\inter\partialv\oldLambda$. The hypothesis that 
the components of
$\Fr_\oldPsi\oldPsi'$ are essential annular suborbifolds of $\oldPsi$
implies in particular that $\frakB$ is $\pi_1$-injective in $\oldPsi$
and therefore in $\oldLambda$, and that the components of $\frakB$ are
annular. Hence if we fix a component $\frakZ$ of $\frakB$, it follows from Proposition \ref{robert strange} that $\frakZ$ is
parallel in the pair $(\oldLambda,\partialv\oldLambda)$ either to $\frakC$ or to an annular suborbifold of $\partialv\oldLambda$. In either subcase, $\frakZ$ is in particular parallel in $\oldPsi$ to a suborbifold of $\partial\oldPsi$. This contradicts the hypothesis that the components of $\Fr_\oldPsi\oldPsi'$ are essential annular suborbifolds of $\oldPsi$.
\EndProof

\Definition\label{S-pair def}
Let $\oldLambda$ be a compact, orientable $3$-orbifold, and let $\oldXi$ be a suborbifold of $\partial\oldLambda$. We shall say that an orbifold fibration of $\oldLambda$ is {\it compatible with $\oldXi$} if either (i) the fibration is an $I$-fibration in which $\oldXi=\partialh\oldLambda $, or (ii) the fibration is an $\SSS^1$-fibration in which $\oldXi$ is saturated.
We define an {\it \spair} to be an ordered pair $(\oldLambda,\oldXi)$ %\redmissingref{Decide whether $\oldXi$ is a good letter here}, 
where $\oldLambda$ is a compact, orientable (possibly disconnected) $3$-orbifold, $\oldXi\subset\partial\oldLambda$ is a $2$-orbifold, and $\oldLambda$ admits an orbifold fibration which is compatible with $\oldXi$.
%either
%Let $\oldPsi$ be a $3$-orbifold, and let $\oldLambda$  be a connected suborbifold of $\oldPsi$ which is equipped with an $I$-fibration or $\SSS^1$-fibration over a $2$-orbifold. We will say that $\oldLambda $ is {\it standardly embedded} in $\oldPsi$ if either 
%$\oldLambda$ admits an orbifold fibration such that either in which 
We will say that an \spair\ is {\it\pagelike} if Alternative (i) of the definition of compatibility holds, and is {\it\bindinglike} if Alternative (ii) holds; these are not mutually exclusive conditions. We will say that a \pagelike\ \spair\ $(\oldLambda,\oldXi)$, with $\oldLambda$ connected, is {\it\untwisted} if the $I$-fibration in Condition (i) can be taken to be a trivial $I$-fibration; otherwise we will say that it is twisted. Thus, in the case where $\oldLambda$ is connected, the \pagelike\ \spair\ $(\oldLambda,\oldXi)$ is twisted if and only if $\oldXi$ has two components, and is untwisted if and only if $\oldXi$ is connected.
\EndDefinition

\Proposition\label{what i need?}
Let $\oldLambda$ be an orientable \torifold, and let $\frakE$ be a $\pi_1$-injective $2$-suborbifold of $\partial\oldLambda$, each component of which is annular. Then $(\oldLambda,\frakE)$ is a \bindinglike\ \spair. 
\EndProposition

\Proof
In this proof, $\DD^2$ and $\SSS^1$ will be identified with subsets of
the complex plane.

For any integer $q\ge1$, let $\frakJ_q$ be defined as in \ref{standard torifold}. Recall that $\frakJ_q$ has both a natural smooth structure and a natural PL structure; these are mutually compatible. Recall also that, up to a homeomorphism which is at once smooth and PL, $|\frakJ_q|$ may be identified with
%e $3$-orbifold such that $|\frakJ_q|=
$\DD^2\times \SSS^1$, in such a way that $\fraks_{\frakJ_q}$ is the curve $\{0\}\times \SSS^1$ and has order $q$ (in the sense of \ref{orbifolds introduced}) if $q>1$, and $\fraks_{\frakJ_1}=\emptyset$ if $q=1$. 
%Let $\frakJ_1$ denote the $3$-orbifold such that $|\frakJ_1|=\DD^2\times \SSS^1$ and $\fraks_{\frakJ_1}=\emptyset$. 
For any integer $r>1$, let $\frakD_r$ denote the $2$-orbifold such that $|\frakD_r|=\DD^2$, and such that $\fraks_{\frakD_r}=\{0\}$, and $0$ has order $r$. Let $\frakD_1$ denote the $2$-orbifold such that $|\frakD_1|=\DD^2$ and $\fraks_{\frakD_1}=\emptyset$.
For any $q\ge1$, and for any relatively prime integers $m$ and $n$
with $m\ne0$, we may define a smooth $\SSS^1$-fibration
$p_{m,n,q}:\frakJ_q\to\frakD_{qm}$ by  
$p_{m,n,q}(z,w)=z^mw^n$. 

For the purpose of this proof, for any integer $m>0$, we define an {\it $m$-admissible system of arcs} to be a set $K\subset \SSS^1$ which is a finite union of closed arcs, and has the property that the map $z\mapsto z^m$ from $K$ to $\SSS^1$ is injective. 
For any $q\ge1$, for any relatively prime integers $m$ and $n$,
% with $|n|<m$, 
and for any $m$-admissible system of arcs $K$, set
$\frakE_{m,n,q,K}=\{(ut^n,t^{-m}):u\in K, t\in \SSS^1\}\subset \SSS^1\times \SSS^1=\partial \frakJ_q$. %Added 6/27/17: I've just supplied the minus sign in the exponent. I have also added the def. of an admissible system to make the following stuff work. Think about all this one more time. Added 7/11/17: Actually to guarantee that the map $(t,u)\mapsto(t^n,ut^{-m})$ is one-to-one, we need the map $z\mapsto z^n$ from $K$ to $\SSS^1$ is injective, not $z\mapsto z^m$. So maybe the right condition is ``$n$-admissibility.'' That would suggest we want $n>0$, not $m>0$, and $|m|<n$, not $|n|<m$.  But I haven't yet got the justification for (either version of) the last inequality straight yet. Later, same day: OK, now I see I have to modify the identfication of the domain with $\DD^2\times \SSS^1$ by doing a Dehn twist. But that should probably be reflected in the argument below. } 
Then 
$\frakE_{m,n,q,K}$ is  saturated in the fibration $p_{m,n,q}$, each component of $|\frakE_{m,n,q,K}|$ is a smooth annulus disjoint from $\fraks_{\frakJ_q}$, and $|\frakE_{m,n,q,K}|$ has the same number of components as $K$ (because the $m$-admissibility of $K$ implies that the map $(u,t)\mapsto(ut^n,t^{-m})$ is one-to-one on $K\times \SSS^1$).

Suppose that $\oldLambda$ and $\frakE$ satisfy the hypotheses of the proposition. (According to our conventions, $\oldLambda$ and its suborbifold $\frakE$ are PL.)
According to Proposition \ref{three-way equivalence}, $\oldLambda$ may be identified with the quotient of $\DD^2\times \SSS^1$ by a standard action of a finite group. 
In view of the description given in \ref{standard torifold} of the quotient of $\DD^2\times \SSS^1$ by a standard action, 
%$\frakJ_q$ here and $\frakJ_q$ in \ref{standard torifold}. Are they the same thing??} 
it follows that up to PL homeomorphism we have either (i) $\oldLambda=\frakJ_q$ for some $q\ge1$, or (ii) there is a $q\ge1$ such that $\oldLambda$ is the quotient of $\frakJ_q$ by the involution $(z,w)\mapsto(\overline{z},\overline{w})$, where bars denote complex conjugation. If  (i) holds, the hypothesis that the components of $\frakE\subset\partial\oldLambda=\SSS^1\times \SSS^1$ are annular and $\pi_1$-injective in $\oldLambda$ implies that $\frakE$ is piecewise smoothly isotopic to $\frakE_{m,n,q,K}$ for some relatively prime integers $m$ and $n$ 
%with $|n|<m$ 
and some $m$-admissible system of arcs $K$. 
Now since $\partial \frakJ_q$ is disjoint from $\fraks_{\frakJ_q}$, it follows from Proposition \ref{fibration-category} that there is a PL $\SSS^1$-fibration $p_{m,n,q}\PL$ of $\frakJ_q$ such that $p_{m,n,q}\PL|\partial\frakJ_q$ is arbitrarily close in the $C^1$ sense to $p_{m,n,q}$. In particular we can choose $p_{m,n,q}\PL$ so that $\frakE_{m,n,q,K}$, which is  saturated in the fibration $p_{m,n,q}$, is piecewise smoothly isotopic to a PL submanifold $\frakE_{m,n,q,K}\PL$ which is saturated in $p_{m,n,q}\PL$. 
Thus $\frakE$ is piecewise smoothly, and hence piecewise linearly, isotopic to $\frakE_{m,n,q,K}\PL$, and is therefore
%, we may define an $\SSS^1$-fibration $p_{m,n,q}:\frakJ_q\to\frakD_{qm}$ by $p_{m,n,q}(z,w)=z^mw^n$. We have $L\subset \SSS^1=\partial \DD^2=\partial|\frakD|$, and $|p|^{-1}(L)=|\frakE|$; hence $\frakE=p^{-1}(\omega(L))$, so that $\frakE$ 
saturated in some PL $\SSS^1$-fibration of $\frakJ_q=\oldLambda$. This means that 
%Consider first the case in which (i) holds. The hypothesis that the components of $\frakE\subset\partial\oldLambda=\SSS^1\times \SSS^1$ are annular and $\pi_1$-injective in $\oldLambda$ implies that $\frakE$ is isotopic to $\{(t^n,ut^m):t\in \SSS^1,u\in K\}$, where $m$ and $n$ are relatively prime integers with $m\ne0$, and $K$ is a finite union of disjoint closed arcs in $\SSS^1$. Now if 
 %$\fraks_\oldLambda=\{0\}\times \SSS^1$, let $q>1$ denote the order of the singular curve $\{0\}\times \SSS^1$, and if
 %$\fraks_\oldLambda=\emptyset$, set $q=1$. Let $\frakD$ denote the $2$-orbifold defined by setting $|\frakD|=\DD^2$ and $\fraks_\frakD=\{0\}$, and taking the singular point $0$ to have order $qm$. We may define an $\SSS^1$-fibration $p_{m,n}:\oldLambda\to\frakD$ by $p(z,w)=z^mw^n$. We have $L\subset \SSS^1=\partial \DD^2=\partial|\frakD|$, and $|p|^{-1}(L)=|\frakE|$; hence $\frakE=p^{-1}(\omega(L))$, so that $\frakE$ is saturated in the $\SSS^1$-fibration. 
%This proves that 
$(\oldLambda,\frakE)$ is a \bindinglike\ \spair. This proves the  proposition in the case where (i) holds.

We now turn to the case in which (ii) holds. We have $\wt|\partial\oldLambda|=4$ in this case, and every point of $\fraks_{\partial\oldLambda}$ has order $2$. 
%Each component of $|\frakE|$ is either a weight-$0$ annulus or a weight-$2$ disk. 
According to \ref{standard torifold}, $\oldLambda$ is the quotient of $\frakJ_q$ by the involution $\tau_q$
given by $|\tau_q|(z,w)=(\overline{z},\overline{w})$. Hence the quotient map $\sigma_q:\frakJ_q\to\oldLambda$ is a degree-two (orbifold) covering map, $\tau_q$ is its non-trivial deck transformation, and $\oldLambda$ inherits a smooth structure from $\frakJ_q$, compatible with its PL structure. We claim:
\Claim\label{behind the annulus}
If $A$ is any smooth  annulus in $|\partial\oldLambda|\setminus\fraks_\oldLambda$ such that $\omega(A)$ is $\pi_1$-injective in $\oldLambda$, then some self-diffeomorphism of $\oldLambda$ carries $A$ onto an annulus $A'$ such that
$\sigma_q^{-1}(\obd(A'))=\frakE_{m,n,q,K}$ for some relatively prime integers
$m,n$
% with $|n|<m$
 and some admissible system of arcs $K$ having exactly
two components, which are  interchanged by complex
conjugation. 
%Furthermore, the image of the
%inclusion homomorphism $\pi_1(\omega(A))\to\pi_1(\oldLambda)$ has index
%$2mq$ in $\pi_1(\oldLambda)$. 
\EndClaim

To prove \ref{behind the annulus}, let $E:\RR^2\to \SSS^1\times \SSS^1=\partial\frakJ_q$ denote the map defined by $E(x,y)=(e^{2\pi ix},e^{2\pi iy})$. Then $\sigma_q\circ E:\RR^2\to\partial\oldLambda$ is a smooth orbifold covering map whose group of deck transformations is the group $\Gamma$ generated by the integer translations of $\RR^2$ and the involution $(x,y)\mapsto(-x,-y)$. Thus $\partial\oldLambda-\fraks_{\partial\oldLambda}$ may be identified with the four-punctured sphere $(\RR^2-\Lambda)/\Gamma$, where $\Lambda=(1/2)\ZZ^2$ is the set of all fixed points of non-trivial elements of $\Gamma$. It is a standard fact in two-dimensional topology (cf. \cite[p. 243, Figure 7]{jones-reid}) that every simple closed curve $C$ in
$(\RR^2-\Lambda)/\Gamma$ which separates $(\RR^2-\Lambda)/\Gamma$ into two twice-punctured disks is isotopic to the image in $(\RR^2-\Lambda)/\Gamma$ of a line of rational or infinite slope in $\RR^2$ disjoint from $\Lambda$; the slope of this line, an element of $\QQ\cup\{\infty\}$ which is uniquely determined by $C$, will be referred to as the {\it slope} of $C$. Now if $A$ satisfies the hypotheses of
 \ref{behind the annulus}, the $\pi_1$-injectivity of $A$ implies that a core curve of $A$ separates two points of $\fraks_{\partial\oldLambda}$ from the other two, and therefore has a well-defined slope. The $\pi_1$-injectivity of $A$ also implies that the slope of a core curve is
non-zero; we will write it as $m/n$, where $m$ and $n$ are relatively prime integers and $m>0$. Then $A$ is isotopic in $\partial\oldLambda-\fraks_{\partial\oldLambda}=(\RR^2-\Lambda)/\Gamma$ to an annulus $A_0'$ such that
$\sigma_q^{-1}(\obd(A_0'))=
\frakE_{m,n,q,K}$ for some admissible system of arcs $K$ having exactly
two components, which are  interchanged by complex
conjugation. 
%n_0
%Now write $n_0=dm+n$, where $d,n\in\ZZ$ and $0\le n<m$. The self-diffeomorphism $\xi:(z,w)\mapsto(zw^{-d},w)$ of $\SSS^1\times \SSS^1=\partial\frakJ_q$ maps $\frakE_{m,n_0,q,K}$ onto $\frakE_{m,n,q,K}$ and commutes with $\tau_q|\partial\frakJ_q$. Hence $\xi$ induces a self-homeomorphism of $\oldLambda$ which maps $A_0'$ onto an annulus $A'$ such that
%$\sigma_q^{-1}(\obd(A'))=
%\frakE_{m,n,q,K}$. 
This establishes  \ref{behind the annulus}.
%A_1

Let us now prove the proposition in the subcase in which $|\frakE|$ is a single  weight-$0$ annulus. In this subcase, choose a smooth annulus $\frakE\smooth$ which is ambiently, piecewise-smoothly isotopic to $\frakE$. In view of \ref{behind the annulus}, some self-diffeomorphism of $\oldLambda$ maps $\frakE\smooth$ onto an annulus $\frakE'\smooth$ such that
% applied with $A=|\frakE|$, we may assume without loss of generality that 
%$\frakE$ is isotopic in $\partial\oldLambda$ to a suborbifold $\frakE'$ such that 
$\sigma_q^{-1}(\frakE'\smooth)=\frakE_{m,n,q,K}$ for some relatively
prime integers $m,n$
% with $|n|<m$ 
and
some $m$-admissible system of arcs $K$ having exactly two components.
The definitions of $\tau_q$ and $p_{m,n,q}$ imply that
$|p_{m,n,q}|\,\circ\,|\tau_q|=|c\circ p_{m,n,q}|$, where $c$ denotes the
involution of $\frakD_{mq}$ such that $|c|$ is
complex conjugation on $\DD^2=|\frakD_{mq}|$. Hence $|p_{m,n,q}|$
induces a map of  spaces from
$|\oldLambda|$ to $|\frakD_{mq}/c|$; this map defines
a smooth
$\SSS^1$-fibration 
$\overline{p}:\oldLambda\to\frakD_{mq}/c$.  Since $\sigma_q^{-1}(\frakE'\smooth)=\frakE_{m,n,q,K}$ is saturated in the fibration $p:=p_{m,n,q}$, the annulus $\frakE'\smooth$ is saturated in $\overline{p}$.

Now let $V$ be a compact PL $2$-submanifold of $\partial\oldLambda$
which is a neighborhood of $\frakE'\smooth$ in $\partial\oldLambda$, and is disjoint from
$\fraks_{\frakJ_q}$.  It follows from Proposition
\ref{fibration-category} that there is a PL $\SSS^1$-fibration
$\overp\pl$ of $\frakJ_q$ such that $|\overp\pl|$ is arbitrarily close
in the uniform sense to $|\overp|$, and $\overp\pl|V$ is arbitrarily close in the $C^1$ sense to $\overp|V$. In particular we can choose $\overp\pl$ so that $\frakE'\smooth$, which is contained in $\inter V$ and is saturated in the fibration $\overp$, is piecewise smoothly isotopic to a PL submanifold $\frakE'$ which is  contained in $\inter V$ and is saturated in $\overp$. Thus $\frakE$ is piecewise smoothly, and hence piecewise linearly, isotopic to $\frakE'$, and is therefore
%, we may define an $\SSS^1$-fibration $p_{m,n,q}:\frakJ_q\to\frakD_{qm}$ by $p_{m,n,q}(z,w)=z^mw^n$. We have $L\subset \SSS^1=\partial \DD^2=\partial|\frakD|$, and $|p|^{-1}(L)=|\frakE|$; hence $\frakE=p^{-1}(\omega(L))$, so that $\frakE$ 
saturated in some PL $\SSS^1$-fibration of $\frakJ_q=\oldLambda$. 
Hence the pair $(\oldLambda,\frakE)$ is a \bindinglike\ \spair, and the  conclusion is established in this subcase.

We next consider the more general subcase in which every component of
$|\frakE|$ is a  weight-$0$ annulus. The $\pi_1$-injectivity of
$\frakE$ implies that each component of $|\frakE|$ separates two
points of $\fraks_{\partial\oldLambda}$ from the other two. The core
curves of the components of $|\frakE|$ are therefore all parallel, and hence there is a single $\pi_1$-injective weight-$0$ annulus $B$ containing $|\frakE|$. By the subcase proved above, letting $\frakB:=\omega(B)$ play the role of $\frakE$, the suborbifold $\frakB$ is saturated in some $\SSS^1$-fibration $p$ of $\oldLambda$. The restriction of $p$ to $\frakB$ is an $\SSS^1$-fibration of the annulus $\frakB$, and since the components of $\frakE$ are $\pi_1$-injective in $\oldLambda$ and hence in $\frakB$, the suborbifold $\frakE$ is isotopic in $\frakB$ to a suborbifold which is saturated in the restricted fibration of $\frakB$, and is therefore saturated in the fibration $p$ of $\oldLambda$. This shows that $(\oldLambda,\frakE)$ is a \bindinglike\ \spair, and the  conclusion of the proposition is established in this more general subcase.

We can now prove the  conclusion whenever (ii) holds. Let $\frakC$ denote a strong regular neighborhood of $\partial\frakE$ relative to $\partial\oldLambda-\inter\frakE$. Then $\frakC$ is $\pi_1$-injective in $\oldLambda$, and each component of $|\frakC|$ is a weight-$0$ annulus. Applying the subcase just proved, with $\frakC$ playing the role of $\frakE$, we obtain an $\SSS^1$-fibration of $\partial\oldLambda$ in which $\frakC$ is saturated. Since $\frakE$ is a union of components of $\partial\oldLambda-\inter\frakC$, it follows that $\frakE$ is saturated in this fibration, and the proof  is complete.
\EndProof
%\tau$t\obd \xi J\times {1}\times J \times {1} \times JK \tau
%\overline $c<

%\redproofreadingnote{I moved  the next two results from the section after this one, and worked a lot on the order of results in this section to try to make it make sense. Give it another think.}

\Proposition\label{unique fibration}
Let $(\oldLambda,\oldXi)$ be a \pagelike\ \spair\ such that each
component of $\oldLambda$ has strictly negative Euler
characteristic. Let $\frakB_1$ and $\frakB_2$ be $2$-orbifolds, and
suppose that for $i=1,2$ we are given an $I$-fibration
$q_i:\oldLambda\to\frakB_i$ which is compatible with $\oldXi$. Then
there exist a homeomorphism $\alpha:\frakB_1\to\frakB_2$, and an
isotopy $(\frakg_t)_{0\le t\le1}$ of the pair  $(\oldLambda,\oldXi)$ (see  \ref{pairs}), such that $q_2\circ\frakg_1=\alpha\circ q_1$.
\EndProposition

\Proof
Set $\oldLambda^{0}=\silv_\oldXi\oldLambda$, and set
$\mu=
\mu_{\oldLambda,\oldXi}
:\oldLambda\to\oldLambda^{0}$.  For
$i=1,2$, according to Lemma \ref{another way to look at it}, there is a unique $\SSS^1$-fibration
$q_i^{0}:\oldLambda^{0}\to\frakB_i$ such that $ q_i^{0}
\circ\mu=q_i$.  Since $q_1^{0}$ and $q_2^{0}$ are both $\SSS^1$-fibrations of $\oldLambda^{0}$ over bases of negative Euler characteristic, it follows from \cite[ p. 446, Fibration Uniqueness Theorem 2]{bonahon-siebenmann} that the fibrations $q_1^{0}$ and $q_2^{0}$ are isotopic: more precisely, 
there exist a homeomorphism $\alpha:\frakB_1\to\frakB_2$,  and a homeomorphism
 $\frakg^{0}:\oldLambda^{0}\to\oldLambda^{0}$ that is isotopic to the identity,  such that
$q_2^{0}\circ\frakg^{0}=\alpha\circ q_1^{0}$. It follows from Lemma \ref{oh
  yeah he tweets}, applied with $\oldLambda$ playing the roles of both
$\oldPsi_1$ and $\oldPsi_2$, and with $\oldXi$ playing the roles of both
$\oldXi_1$ and $\oldXi_2$,
%Now  $|\oldXi|$ is the union of the $2$-dimensional strata of
%$\oldLambda^{0}$. For $i=1,2$, and for $0\le t\le 1$, since
%$\frakg_i^{0}:\oldLambda^{0}\to\oldLambda^{0}$ is a homeomorphism,
%$|\frakg_i^{0}|$ leaves $|\oldXi|$ invariant, and
that  
%$t\in[0,1]$ 
there is an isotopy $(\frakg_t)_{0\le t\le1}$ of the pair $(\oldLambda,\oldXi)$ such
that $\silv\frakg_1=\frakg^0$. According to \ref{third point}, we have
$\frakg_1\circ\mu
=
\mu\circ
\frakg^{0}$.
%Since 
%$|\mu|$ and $\frakg_t^{0}$ are homeomorphisms, $\frakg_t$ is a homeomorphism.
It now follows that 
%$(\frakg_t)_{0\le t\le1}$ is an isotopy of $\oldLambda$,
% that $\frakg_t(\oldXi)=\oldXi$ for each $t\in[0,1]$, 
%and that 
$q_2\circ \frakg_1=\alpha\circ q_1$.
\EndProof
%\eta\alpha tweets \oldUpsilon ' \eta {oh ^{0}

\Definition\label{strong equivalence}
Let $\otheroldLambda$ be a $2$-orbifold. Suppose that $\oldXi_1$ and $\oldXi_2$ are suborbifolds of $\otheroldLambda$, and that $\iota_i$ is an involution of $\oldXi_i$ for $i=1,2$. We will say that $\iota_1$ and $\iota_2$ are {\it strongly equivalent (in $\otheroldLambda$)} if there is an isotopy $(\frakh_t)_{0\le t\le1}$ of $\oldLambda$ such that $h_1(\oldXi_1)=\oldXi_2$ and $\iota_2=\frakh_1\circ \iota_1\circ \frakh_1^{-1}$.

Note that strong equivalence is indeed an equivalence relation on the set of all involutions of suborbifolds of $\oldXi$. Note also that if the involutions $\iota_i$ of $\oldXi_i$ are strongly equivalent, then in particular $\oldXi_1$ and $\oldXi_2$ are isotopic. 
\EndDefinition

%\redcomment{Figure out where to put this, and rewrite so as to emphasize that it's about $[\iota]$ being well defined:} two strongly equivalent involutions are in particular isotopic as self-homeomorphisms of $\oldXi$.

\Number\label{more about strong equivalence}
Let $\otheroldLambda$ be a negative $2$-orbifold, let $\oldXi$ be a
negative, taut suborbifold of $\otheroldLambda$, and let $\iota_1$ and
$\iota_2$ be involutions of $\oldXi_i$. If $\iota_1$ and $\iota_2$ are
strongly equivalent in $\otheroldLambda$, then there is a
self-homeomorphism $h$ of $\oldXi$, isotopic in $\oldLambda$ to the
identity map of $\oldXi$, such that and $\iota_2=\frakh\circ
\iota_1\circ \frakh^{-1}$. It then follows from Corollary \ref{i guess} % either rewrite the following corollary so it's about strong equivalence {\it in $} 
that $\frakh$ is isotopic in $\oldXi$ to the identity map. This shows that if two involutions of $\oldXi$ are strongly equivalent in $\otheroldLambda$ then they are strongly equivalent in $\oldXi$. The converse is trivial.
\EndNumber

\Corollary\label{flubadub o'connor}
Let $(\oldLambda,\oldXi)$ be a \pagelike\ \spair\ such that each component of $\oldLambda$ has strictly negative Euler characteristic. Let $\frakB_1$ and $\frakB_2$ be $2$-orbifolds, and suppose that for $i=1,2$ we are given an $I$-fibration $q_i:\oldLambda\to\frakB_i$ which is compatible with $\oldXi$. For $i=1,2$, let $\iota_i$ denote the non-trivial deck transformation of the covering map $q_i|\oldXi:\oldXi\to\frakB_i$. Then $\iota_1$ and $\iota_2$ are strongly equivalent in $\oldXi$.
\EndCorollary

\Proof
According to Proposition \ref{unique fibration}, there exist a
homeomorphism $\alpha:\frakB_1\to\frakB_2$, and an isotopy
$(\frakg_t)_{0\le t\le1}$ of the pair $(\oldLambda,\oldXi)$, such that
$q_2\circ\frakg_1=\alpha\circ q_1$. If we set
$\frakh_t=\frakg_t|\oldXi$ for $0\le t\le1$, then  $(\frakh_t)_{0\le
  t\le 1}$ is an isotopy of $\oldXi$. If we set
 $p_1=\alpha\circ(q_1|\oldXi)$ and $p_2=q_2|\oldXi$, we have $p_2\circ \frakh_1=p_1$, i.e. $\frakh_1$ is an equivalence of coverings between $p_1:\oldXi\to\frakB_2$ and $p_2:\oldXi\to\frakB_2$. Since $\iota_i$ is the unique non-trivial deck transformsation of $p_i$ for $i=1,2$, it follows that $\frakh_1\circ\iota_1\circ \frakh_1^{-1}=\iota_2$.
\EndProof
%\alpha\eta\oldLambda(i

\begin{remarksdefinitions}\label{Do I need it?}
Let $\oldPsi$ be an  orientable $3$-orbifold which is componentwise strongly \simple\ and
componentwise boundary-irreducible. If  $\oldLambda$ is a $3$-suborbifold of $\oldPsi$ such that
$(\oldLambda,\oldLambda\cap\partial\oldPsi)$ is an \spair, then by
definition $\oldLambda$ admits an orbifold
fibration compatible with $\oldLambda\cap\partial\oldPsi$. The definitions imply that, with respect to such a fibration, the suborbifold 
$\Fr_\oldPsi\oldLambda$ of $\oldLambda$ is saturated; and if the \spair\ $(\oldLambda,\oldLambda\cap\partial\oldPsi)$ is \pagelike\ we may take the fibration to be an $I$-fibration with $\partialh\oldLambda=\oldLambda\cap\partial\oldPsi$, in which case $\Fr_\oldPsi\oldLambda=\partialv\oldLambda$. In any event, the saturation of $\Fr_\oldPsi\oldLambda$ implies that each of its components is an annular or toric orbifold, and each of its components must be annular if $(\oldLambda,\oldLambda\cap\partial\oldPsi)$ is \pagelike\ (see \ref{fibered stuff}). 

We define an {\it \Ssuborbifold} of $\oldPsi$ 
to be a $3$-suborbifold $\oldLambda$ of $\oldPsi$ such that
the
components of $\Fr_\oldPsi\oldLambda$ are annular suborbifolds of $\oldPsi$
which are essential in $\oldPsi$, and
$(\oldLambda,\oldLambda\cap\partial\oldPsi)$ is an \spair. Note that this definition allows $\oldLambda$ to be disconnected; note also that $\oldLambda$ is an \Ssuborbifold\ of $\oldPsi$ if and only if each of its components is an \Ssuborbifold\ of $\oldPsi$. An \Ssuborbifold\
$\oldLambda$ will be called {\it\pagelike} or {\it\bindinglike} if
$(\oldLambda,\oldLambda\cap\partial\oldPsi)$ is a \pagelike\ or \bindinglike\ 
\spair, respectively.  Similarly, a connected \pagelike\ \Ssuborbifold\ $\oldLambda$ will be termed
{\it\untwisted} or {\it\twisted}   if the \pagelike\ \spair\ 
$(\oldLambda,\oldLambda\cap\partial\oldPsi)$ is 
\untwisted\ or \twisted,  respectively. 

A $3$-suborbifold $\oldLambda$ of $\oldPsi$
will be called an {\it \Asuborbifold} of $\oldPsi$ if the
components of $\Fr_\oldPsi\oldLambda$ are essential annular suborbifolds of $\oldPsi$,
%which are essential in the pair $(\oldPsi,\partial\oldPsi)$, 
and
$(\oldLambda,\oldLambda\cap\partial\oldPsi)$ is an acylindrical pair (see 
\ref{acylindrical def}).
\end{remarksdefinitions}

\Lemma\label{when a tore a fold}
Let $\oldLambda$ be a connected suborbifold of a componentwise strongly \simple,
componentwise boundary-irreducible, orientable $3$-orbifold $\oldPsi$. Then $\oldLambda$
is a \bindinglike\ \Ssuborbifold\ of $\oldPsi$ if and only if (i)
$\oldLambda$ is a \torifold\, and (ii) every component of
$\Fr_\oldPsi\oldLambda$ is an  essential annular suborbifold of
$\oldPsi$. 
\EndLemma

\Proof
Suppose that $\oldLambda$
is a \bindinglike\ \Ssuborbifold\ of $\oldPsi$. Consider first the case in which $\partial\oldLambda\ne\emptyset$, and fix a component $\oldTheta$ of $\partial\oldLambda$. Since $\oldLambda$ is a \bindinglike\ \Ssuborbifold, $\oldTheta$ admits an $\SSS^1$-fibration over a closed $1$-orbifold, and is therefore toric. 
%But $\oldLambda$ is $\pi_1$-injective in $\oldPsi$
Since the components of $\Fr\oldLambda$ are essential, and in particular $\pi_1$-injective, in the strongly \simple\ orbifold $\oldPsi$,
%; hence $\oldTheta$ is not $\pi_1$-injective in $\oldLambda$. On the other hand, 
it follows from Lemma \ref{oops lemma} that $\oldLambda$ is strongly \simple. Since $\oldTheta$ is toric, $\oldLambda$ satisfies Condition (2) of Proposition \ref{three-way equivalence}, and according to the latter proposition it is a \torifold. This is Condition (i) of the statement of the present lemma, and
Condition (ii) is immediate from the definition of an \Ssuborbifold.

Now consider the case in which $\oldLambda$ is closed. Then $\oldLambda$ is a component of $\oldPsi$ and is therefore strongly \simple. According to Condition (II) of Definition \ref{oops}, some finite-sheeted covering space $\toldLambda$ of $\oldLambda$ is an irreducible $3$-manifold. But by definition the \bindinglike\ \Ssuborbifold\ $\oldLambda$ admits an $\SSS^1$-fibration, and $\toldLambda$ inherits an $\SSS^1$-fibration $q:\toldLambda\to\frakB$ for some closed $2$-orbifold $\frakB$. Since $\oldLambda$ is irreducible, $\frakB$ is infinite and therefore has an element of infinite order. This implies that $\pi_1(\toldLambda)$ has a rank-$2$ free abelian subgroup, a contradiction to Condition (II) of Definition \ref{oops}. Hence this case cannot occur, and we have proved that if $\oldLambda$ is a \bindinglike\ \Ssuborbifold\ then (i) and (ii) hold.

Conversely, suppose that (i) and (ii) hold. It follows from (ii) that the components of $\Fr\oldPsi$ are $\pi_1$-injective in $\oldPsi$ as well as being annular. Hence Proposition \ref{what i need?} implies that $(\oldPsi,\Fr\oldPsi)$ is a \bindinglike\ \spair. As (ii) includes the assertion that the components of $\Fr\oldPsi$ are essential in $\oldPsi$, it now follows that $\oldLambda$ is a \bindinglike\ \Ssuborbifold\ of $\oldPsi$.
\EndProof
%\obd

\Proposition\label{doesn't exist}
Let $\oldPsi$ be a compact, orientable $3$-orbifold which is strongly \simple\ and
boundary-irreducible,  let $\oldLambda$ be a connected \Ssuborbifold\ of $\oldPsi$, and set $\oldXi=\oldLambda\cap\partial\oldPsi$. Let $\oldXi_0$ be a component of $\oldXi$. Then the inclusion homomorphism $H_1(|\oldXi_0|;\QQ)\to H_1(|\oldLambda|;\QQ)$ is surjective. Furthermore, if $\oldLambda$ is pagelike, then the image of the inclusion homomorphism $\pi_1(|\oldXi_0|)\to\pi_1(|\oldLambda|)$ has index at most $2$ in $\pi_1(|\oldLambda|)$.
\EndProposition

\Proof
We prove the second assertion first. Consider the commutative diagram
$$
\begin{matrix}
\pi_1(\oldXi_0)&\longrightarrow&\pi_1(\oldLambda)\\
\downarrow&\empty&\downarrow\\
\pi_1(|\oldXi_0|)&\longrightarrow&\pi_1(|\oldLambda|)
\end{matrix}
$$
in which the horizontal maps are induced by inclusion, and each vertical map is the canonical surjection from the fundamental group of a connected orbifold to the fundamental group of its underlying space. Let us fix an $I$-fibration $q:\oldLambda\to\frakB$, where $\frakB$ is some $2$-orbifold, under which $\oldXi$ is the horizontal boundary of $\oldLambda$. Then $q|\oldXi_0:\oldXi_0\to\frakB$ is a covering map of some degree $d\le2$, and hence the  image of the top horizontal map in the diagram has index  $d\le2$ in $\pi_1(\oldLambda)$. Since the right vertical map is surjective, it follows that the image of the inclusion
$\pi_1(|\oldXi_0|)\to\pi_1(|\oldLambda|)$ has index at most $2$ in $\pi_1(|\oldLambda|)$, as required.

In the case where the \spair\ $(\oldLambda,\oldXi)$ is \pagelike, the first assertion follows immediately from the second. To prove the first assertion in the case where $(\oldLambda,\oldXi)$ is \bindinglike, note that by  Lemma \ref{when a tore a fold}, 
$\oldLambda$ is a \torifold\, and every component of
$\Fr_\oldPsi\oldLambda$ is an  essential annular suborbifold of
$\oldPsi$. It follows that the component $\oldXi_0$ of $\oldXi=\overline{\partial\oldLambda-\Fr_\oldPsi\oldLambda}$ has Euler characteristic $0$. Now since $\oldLambda$ is a torifold, $\partial\oldLambda$ is connected; furthermore, we have $\Fr_\oldPsi\oldLambda\ne\emptyset$, as otherwise the \torifold\ $\oldLambda$ would be a component of $\oldPsi$, a contradiction to the componentwise boundary-irreducibility of $\oldPsi$. It follows that $\partial\oldXi_0\ne\emptyset$, so that $\oldXi_0$ is annular. The essentiality of the components of $\Fr_\oldPsi\oldLambda$ implies that $\partial\oldXi_0$ is $\pi_1$-injective in $\oldLambda$, and hence $\oldXi_0$ is $\pi_1$-injective in $\oldLambda$. On the other hand, since $\oldLambda$ is a \torifold, it follows from Proposition \ref{three-way equivalence} that either $|\oldLambda|$ is a $3$-ball or $\oldLambda$ is a solid torus. If $|\oldLambda|$ is a $3$-ball, the surjectivity of the inclusion homomorphism $H_1(|\oldXi_0|;\QQ)\to H_1(|\oldLambda|;\QQ)$ is trivial. If $\oldLambda$ is a solid torus, and in particular a manifold, then the orientable annular orbifold $\oldXi_0$ must be an annulus, and the $\pi_1$-injectivity of $\oldXi_0$ in $\oldLambda$ implies that the image of the inclusion homomorphism $\pi_1(\oldXi_0)\to \pi_1(\oldLambda)$ has finite index in $\pi_1(\oldLambda)$; in particular,
the inclusion homomorphism $H_1(\oldXi_0;\QQ)\to H_1(\oldLambda;\QQ)$ is surjective.
\EndProof
%\obd

\Proposition\label{when vertical}
Let $\oldLambda$ be an orientable
$3$-orbifold equipped with an $I$-fibration over a negative $2$-orbifold,
and let $\frakC\subset\oldLambda$ be a two-sided $2$-orbifold with $\partial\frakC\subset\partialh\oldLambda$. Suppose that every component of $\frakC$ is annular and is $\pi_1$-injective in $\oldLambda$. Suppose also that no component of $\frakC$ is parallel in the 
 pair $(\oldLambda,\partialh\oldLambda)$ to a suborbifold of
 $\partialh\oldLambda$. 
Then $\frakC$ is isotopic, via an isotopy of the pair
$(\oldLambda,\partialh\oldLambda)$ (see \ref{pairs}), to a saturated annular suborbifold.
\EndProposition

\Proof
First consider the case in which $\frakC$ is connected (and therefore
annular). After replacing $\oldLambda$ by its component containing
$\frakC$, we may assume that $\oldLambda$ is itself connected. The hypothesis implies that  $\chi(\oldLambda)<0$.  Since, by hypothesis,
the annular orbifold $\frakC$ is $\pi_1$-injective in $\oldLambda$, and no component of $\frakC$ is parallel in the 
 pair $(\oldLambda,\partialh\oldLambda)$ to a suborbifold of
 $\partialh\oldLambda$, the definition of essentiality implies that
 either (a) $\frakC$ is essential in $(\oldLambda,\partialh\oldLambda)$, or (b) 
$\frakC$ is parallel in the 
 pair $(\oldLambda,\partialh\oldLambda)$ to a component of $\partialv\oldLambda$. If (b) holds then the conclusion of the proposition is immediate. Now suppose that (a) holds.

Set $\oldXi=\partialh\oldLambda$ and
$\oldPsi=\silv_\oldXi\oldLambda$. Set
$\mu=\mu
_{\oldLambda,\oldXi}
:\oldLambda\to\oldPsi$.
Let $\oldTheta$ denote the two-sided toric suborbifold
$\silv\frakC$ of $\inter\oldPsi$. Since $\frakA$ is essential, it
follows from Proposition \ref{silver acylindrical} that $\oldTheta$ is
$\pi_1$-injective in $\oldPsi$ and has no component which is parallel
to a component of $\partial\oldPsi$.

Let us write the given $I$-fibration of $\oldLambda$ as
$q:\oldLambda\to\frakB$, where $\frakB$ is some $2$-orbifold.
According to Lemma \ref{another way to look at it}, there is a unique $\SSS^1$-fibration
$q':\oldPsi\to\frakB$ such that $ q' \circ\mu=q$.
%The $I$-fibration of $\oldLambda$
%gives rise  to an $\SSS^1$-fibration of $\oldPsi$ in which the
%fibers are obtained by silvering the fibers of $\oldLambda$. The base
%of the $I$-fibration of $\oldLambda$, which we will denote by
%$\frakB$, is also the base of the $\SSS^1$-fibration of $\oldPsi$. 
It follows from the PL version of 
 \cite[Verticalization Theorem 4]{bonahon-siebenmann} (see
 \ref{categorille}) that $\oldTheta$ is isotopic in $\oldPsi$ to a
 suborbifold $\oldTheta'$ which is either ``vertical,'' in the sense
 of being saturated in the $\SSS^1$-fibration $q'$ of $\oldPsi$, or
 ``horizontal,'' in the sense of being transverse to the fibers of $q'$. 
If $\oldTheta'$ is horizontal then it is a(n orbifold)
covering of $\frakB$; if $d$ denotes the degree of the covering, the
toricity of $\oldTheta$ gives
$0=\chi(\oldTheta)=\chi(\oldTheta')=d\chi(\frakB)$, so that
$\chi(\frakB)=0$. But by hypothesis we have
$0>\chi(\oldLambda)=\chi(\frakB)$, a contradiction. Hence $\oldTheta'$
must be vertical, i.e. saturated. We may therefore write
$\oldTheta'=\silv\frakC'$ for some saturated suborbifold $\frakC'$ of
$\oldLambda$. 

Let
$\frakh: \oldPsi\to \oldPsi$ be a homeomorphism,
isotopic to the identity,  such that
$\frakh(\oldTheta)=\oldTheta'$. It follows from Lemma \ref{oh
  yeah he tweets}, applied with $\oldLambda$ playing the roles of both
$\oldPsi_1$ and $\oldPsi_2$, and with $\oldXi$ playing the roles of both
$\oldXi_1$ and $\oldXi_2$,
that  there is an
isotopy $(\frakg_t)_{0\le t\le 1}$ of the pair $(\oldLambda,\oldXi)\to(\oldLambda,\oldXi)$ such that
$\silv\frakg_1=\frakh$. 
According to \ref{third point} we have
$\mu\circ\frakg_1=\frakh\circ\mu$, and it now follows that 
%$(\frakg_t)_{0\le t\le1}$ is an isotopy of $\oldLambda$,
%and that
 $\frakg_1(\frakC)=\frakC'$. 
%\redcomment{Fix from here to make it consistent with my
  %silvering conventions. This will be another app. of the revised
  %version of Prop \ref{oh yeah he tweets}. Let $\oldUpsilon\subset\oldPsi$ denote the image of $\oldXi$
%under the canonical immersion $\oldLambda\to\oldPsi$ (see \ref{silvering}). Since
%$\oldLambda$ is orientable, $|\oldUpsilon|\subset|\oldPsi|$ is the union of the closures of
%the two-dimensional strata of $\oldPsi$. Hence an isotopy of  $\oldTheta=\silv\frakC$ is isotopic in
%$\oldPsi$ carrying $\oldTheta$ to $\oldTheta'$ must leave $\oldUpsilon%$
%invariant, and therefore defines an isotopy of $\oldLambda$ which
%leaves $\oldXi$ invariant and carries $\frakC$ onto $\frakC'$. % of $\oldLambda$ are
%isotopic. 
%\redcomment{That requires a little more argument. Also, it's incorrectly stated,
  %since the statement involves an isotopy that leaves the horizontal
  %boundary invariant. I have corrected the statement in this regard,
  %and should fix apps. too. And of course the rest of the proof.} 
 This completes the proof in the case in which $\frakC$ is connected.

To prove the proposition in the general case, we use induction on
$\compnum(\frakC)$. If $\compnum(\frakC) =0$ the assertion is trivial,
and for $\compnum(\frakC)=1$ it has already been proved. Now suppose
that $n>1$ is given, and that the assertion has already been proved
for the case  $\compnum(\frakC) =n-1$. Suppose that the hypotheses
hold and that $\compnum(\frakC)=n$. Choose a component $\frakC_0$ of
$\frakC$. According to the induction hypothesis we may assume, after
modifying $\frakC$ by an isotopy of
$\oldLambda$ that leaves $\partialh\oldLambda$ invariant, that $\frakC':=\frakC-\frakC_0$
is saturated in the fibration of $\oldLambda$. If
$q:\oldLambda\to\frakB$ denotes the $I$-fibration of $\oldLambda$,
where $\frakB$ is some $2$-orbifold, the saturation of $\frakC'$
implies that $\frakC'=p^{-1}(\frakK)$ for some two-sided $1$-orbifold
$\frakK\subset\frakB$. If we set
$\oldLambda_0=\oldLambda\cut{\frakC'}$ and
$\frakB_0=\frakB\cut\frakK$, and
$\rho=\rho_{\frakC'}:\oldLambda_0\to\oldLambda$ and
$\sigma=\rho_{\frakK}:\frakB_0\to\frakB$, there is a unique
$I$-fibration $q_0:\oldLambda_0\to\frakB_0$ such that $\sigma\circ
q_0=q\circ\rho$. The horizontal boundary of $\oldLambda_0$ with
respect to this $I$-fibration contains the boundary of the annular
suborbifold $\rho^{-1}(\frakC_0)$ of $\oldLambda_0$. By the case of
the proposition that has already been proved, $\frakC$ is isotopic,
via an isotopy of the pair $(\oldLambda_0,\partialh\oldLambda_0)$, to a saturated annular suborbifold $C_1$ of
$\oldLambda_0$. Now $\frakC'\cup\rho(C_1)$ is saturated in the
fibration of $\oldLambda_0$ and is isotopic to $\frakC$ via an isotopy
of $\oldLambda$ that leaves $\partialh\oldLambda$ invariant. This completes the induction.
\EndProof
%\obd \frakK \oldUpsilon\oldXi_1 \oldPsi\oldLambda'\oldUpsilon \eta
%{oh yeah he tweets}

\Corollary\label{vertical corollary}
Let $(\oldLambda,\oldXi)$ be an \spair\ such that $\oldLambda$ is strongly \simple,
and let $\frakC$ be an
% two-sided properly embeddedprop $2$-suborbifold of
%$\oldLambda$, each of whose components is
 essential  annular orbifold in 
the pair $(\oldLambda,\oldXi)$. Then $\oldLambda$ admits an orbifold fibration compatible with $\oldXi$ in which $\frakC$ is saturated.
\EndCorollary

\Proof
In the case where
$\chi(\oldLambda)<0$, the \spair\ $(\oldLambda,\oldXi)$ must be
\pagelike, since an orbifold admitting an $\SSS^1$-fibration has Euler
characteristic $0$; and if $\frakB$ denotes the base of an $I$-fibration of
$\oldLambda$ compatible with $\oldXi$, we have
$\chi(\frakB)=\chi(\oldLambda)<0$. Thus  in this case the assertion
follows from Proposition \ref{when vertical}. 

Now suppose that
$\chi(\oldLambda)\ge0$. Since $\frakC$ is annular, the (possibly
disconnected) $3$-orbifold $\toldLambda:=\oldLambda\cut\frakC$ also
has non-negative Euler characteristic. By Lemma \ref{oops lemma},
$\toldLambda$ is componentwise strongly \simple; in particular,
$\toldLambda$ is componentwise irreducible and very good and has no
discal component, so $\partial\toldLambda$ is very good and has no
spherical component. Since
$\chi(\partial\toldLambda)=2\chi(\toldLambda)\ge0$, it now follows
that every component of $\partial\toldLambda$ is toric. Since by
construction every component of $\oldLambda$ has non-empty boundary,
it follows from
Proposition \ref{three-way equivalence}, each component of
$\toldLambda$ is a \torifold. 

%Now set $\oldXi^*=(\partial\oldLambda)-\inter\oldXi$. 

Set
$\tfrakC=\rho_\frakC^{-1}(\frakC)$ and
$\toldXi=\rho_\frakC^{-1}(\oldXi)$ (see \ref{nbhd stuff}),  so that
the components of $\tfrakC$ and $\toldXi$ are annular and
$\pi_1$-injective in $\toldLambda$. Set
$\toldXi^*=(\partial\toldLambda)-\inter\toldXi$, so that $\toldXi^*$ and
$\tfrakC$ are disjoint suborbifolds of $\partial\toldLambda$. 
Since $\oldLambda$ is connected,
each component of $\toldLambda$ contains a component of
$\toldXi$. Hence each component of $\toldXi^*$ is a proper suborbifold
of the boundary of some \torifold\ component of $\toldLambda$. Since
the components of $(\partial\toldLambda)-\inter\toldXi^*=\toldXi$ are
$\pi_1$-injective annular suborbifolds of $\toldLambda$, it
now follows that each component of $\toldXi^*$ is also a
$\pi_1$-injective annular suborbifold of $\toldLambda$.   Proposition \ref{what i need?}, applied to the
components of $\toldLambda$, shows that $(\toldLambda,\toldXi^*\discup\tfrakC)$ is a \bindinglike\ \spair, i.e. $\toldLambda$ has an $\SSS^1$-fibration in which $\tfrakC$ and $\toldXi^*$ are saturated. Since the $\SSS^1$-fibration of an orientable annular $2$-orbifold is unique up to isotopy, we may choose the $\SSS^1$-fibration of $\toldLambda$ in such a way that $\grock_\frakC$ (see \ref{nbhd stuff})  interchanges the induced fibrations of the components of $\tfrakC$. This implies that
$\oldLambda$ has an $\SSS^1$-fibration in which $\frakC$ and $\oldXi^*:=p
(\partial\oldLambda)-\inter\oldXi
=\rho_\frakC(\toldXi^*)$ are saturated. Since $\oldXi^*$ is saturated, so is
$\oldXi$, and the conclusion follows in this case.
\EndProof
%\frakA\frakE\tfrakE\tfrakC\tfrakA\oldLambda'\discup\cup

\Proposition\label{new what they look like}
Let $\oldLambda$ be a compact, connected, orientable $3$-orbifold, and let $\oldXi$ be a compact $2$-suborbifold of $\partial\oldLambda$. Then $\silv_\oldXi\oldLambda$ admits an $\SSS^1$-fibration over a $2$-orbifold if and only if $(\oldLambda,\oldXi)$ is an \spair.
\EndProposition

\Proof 
Set $\oldOmega=\silv_\oldXi\oldLambda$, and set $\mu=\mu_{\oldLambda,\oldXi}:\oldLambda\to
\oldOmega$. (see \ref{silvering}). Set $\oldXi'=\mu(\oldXi)\subset\oldOmega$; by
\ref{gyrot}, $\oldXi'$ is a suborbifold of $\oldOmega$. Set 
$\mu^*=\mu^*_{\oldLambda,\oldXi}:\oldLambda-\oldXi\to\oldOmega-\oldXi'$.

First suppose that  $(\oldLambda,\oldXi)$ is an \spair. By definition this means that there is a fibration $q:\oldLambda\to\frakB$ of $\oldLambda$
over a $2$-orbifold $\frakB$ such that either (i) $q$ is an
$I$-fibration and $\oldXi=\partialh\oldLambda$, or (ii) $\oldLambda$ is an
$\SSS^1$-fibration and $\oldXi$ is saturated. If (i) holds, then
according to Lemma \ref{another way to look at it} there is a fibration $q':
\oldOmega\to\frakB$ such that $q'\circ\mu= q$. If (ii) holds, we
may write $\oldXi=q^{-1}(\frakJ)$ for some $1$-orbifold
$\frakJ\subset\partial\frakB$, and 
according to Lemma \ref{one way to look at it}, $\silv
q:\oldOmega\to\silv_\frakJ\frakB$ is an $\SSS^1$-fibration.

Conversely, suppose that $r$ is an $\SSS^1$-fibration of $\oldOmega$ over a
$2$-orbifold $\oldGamma$. 
According to the definition of an
$\SSS^1$-fibration (see \ref{fibered stuff}), for each  point
$v\in\oldGamma$ there is a local $\SSS^1$-standardization
of $r$ at
$v$. For any local $\SSS^1$-standardization
$s=(\frakV,G,\psi,\zeta)$ of $r$ at $v$, we will set 
 $Z_s=P_s^{-1}(\oldXi'\cap r^{-1}(\frakV))$ (where $P_s$ is defined as
 in \ref{fibered stuff}).

%, and make sure the following points are dealt with: `` Also make sure $G_2$ can't be dihedral''} 
%that
%According to \redproofreadingnote{This is
%  sort of in \cite{other-bs}. I think it differs  from the definition only in that is says the covering is regular. If I'm right, I can either add regularity to the def., or point out when I give the def. that it comes for free }, 

% the canonical two-sheeted covering map from $D_\oldXi\oldLambda$
%to $\oldOmega$ maps $\oldXi$ homeomorphically onto a $2$-orbifold
%$\oldXi'\subset\oldOmega$. \abstractcomment{I was thinking it's not a
%suborbifold, but now after studying defs. I think it is.} 
According to \ref{gyrot},
%the
%the
%definition of $\oldOmega=\silv_\oldXi\oldLambda$ implies 
%\Claim\label{looks good}
%Each point of $\oldOmega$ has a neighborhood $U\subset\oldOmega$ such that
%for every manifold covering $\tU$ of $U$, the pre-image in $\tU$ of 
%that 
$|\oldXi'|\subset|\oldOmega|$ is the
union of all closures of $2$-dimensional strata of $
%U\cap
\fraks_\oldOmega$. Hence if  $s=(\frakV,G,\psi,\zeta)$ is any local
$\SSS^1$-standardization of $r$ at an arbitrary point $v\in\oldGamma$, then $Z_s$ 
is
%of 
the
union of all $2$-dimensional fixed point sets of elements of $G$. Each
such $2$-dimensional fixed point set is equal
either to (A)
  $ U_v\times\{-c,c\}$ for some $c\in \SSS^1$, or (B) $\ell\times
  \SSS^1$, where either  $v\in\inter\oldGamma$ and $\ell$  is some diameter 
of $ U^2 = U_v$,  or
$v\in\partial\oldGamma$ and
$\ell=\{0\}\times[0,1)\subset U^2_+= U_v$. But $Z_s $ 
is a $2$-submanifold of $ U_v\times \SSS^1$ since $\oldXi'$ is a suborbifold of $\oldOmega$.
It
follows that $Z_s $ 
is either a (possibly empty) union of fixed point sets of the form
(A), or a single fixed point set of the form (B). This proves:

\Claim\label{better onion}
For every 
local
$\SSS^1$-standardization
$s=(\frakV,G,\psi,\zeta)$  of $r$ at an arbitrary point $v\in\oldGamma$, we have either
%$v\in\oldGamma$ there exist a neighborhood $\frakV$ of $v$ in
%$\oldGamma$ and  a regular covering map $P: U^2\times \SSS^1\to r^{-1}(\frakV)$ such that
%(a) $P$ maps $\{x\}\times \SSS^1$ onto a fiber of $r$ for each $x\in
% U^2$, and (b)
%$\frakZ $ either
 (I) $Z_s= \emptyset$, (II) $Z_s=
 U_v\times\Phi$ for some non-empty finite set $\Phi\subset \SSS^1$
with $-\Phi=\Phi$, or (III) $Z_s=\ell\times \SSS^1$, where either  $v\in\inter\oldGamma$ and $\ell$  is some diameter 
of $ U^2 = U_v$,  or
$v\in\partial\oldGamma$ and
$\ell=\{0\}\times[0,1)\subset U^2_+= U_v$. Furthermore, $ Z_s$ is the union of all two-dimensional
components of fixed point sets of elements of $G$; in particular,
$ Z_s$ is $G$-invariant.
%if (II) or (III) holds, there is a deck
%transformation of the covering $P: U^2\times \SSS^1\to r^{-1}(\frakV)$ which
%interchanges the components of $( U^2\times \SSS^1)-P^{-1}(\oldXi'\cap
%r^{-1}(\frakV))$. 
\EndClaim 
%\frakZ

Let $ T$ denote the set of all points $v\in|\oldGamma|$ such that
Alternative (I) or (III) of 
\ref{better onion} holds for some local
$\SSS^1$-standardization
$s=(\frakV,G,\psi,\zeta)$  of $r$ at $v$; and let $ T'$ denote the set of all points $v\in|\oldGamma|$ such that
Alternative (II)  of 
\ref{better onion} holds for some local
$\SSS^1$-standardization
$s=(\frakV,G,\psi,\zeta)$  of $r$ at $v$.
Then $ T\cup
 T'=|\oldGamma|$. It is clear that $ T$ and $ T'$ are open. 
Furthermore, the
fiber over every point of $ T$ is either contained in $\oldXi'$ or
disjoint from $\oldXi'$, while the fiber over every point of $ T'$
meets $\oldXi'$ in a non-empty finite set. Thus $ T\cap T'=\emptyset$. But $|\oldGamma|$ is connected since
$\oldLambda$ is connected, and hence either $|\oldGamma|= T$ or
$|\oldGamma|= T'$.

Consider first the case in which $|\oldGamma|=T'$. In this case we let
$q$ denote the submersion $r\circ\mu$. We will show that $q:\oldLambda\to\oldGamma$ is an $I$-fibration, and
that with respect to this fibration we have
$\partial_h\oldLambda=\oldXi$; this will imply that $(\oldLambda,\oldXi)$
is  a \pagelike\ \spair.

Let
$v\in|\oldGamma|$ be given. Since $|\oldGamma|=T'$, there is
 a local
$\SSS^1$-standardization $s=(\frakV,G,\psi,\zeta)$ of $r$ at $v$ for which
Alternative (II) of 
\ref{better onion} holds. If $\Phi\subset \SSS^1$ is the set given by (II), let us
denote by $\iota_c$, for each $c\in\Phi$, the unique
orientation-reversing isometry of $\SSS^1$ whose fixed points are $c$ and
$-c$; thus $\iota_{-c}=\iota_c$ for each $c\in\Phi$. It follows from
the last sentence of \ref{better onion} that $ U_v\times\Phi$ is
$G$-invariant, and that
%the $G$ acts on the set
%of components of $( U_v\times
%\SSS^1)-Z_s= U_v\times(\SSS^1-\Phi)$, and that 
the elements of $G$
having two-dimensional fixed point sets are precisely the
self-homeomorphisms of $ U_v\times \SSS^1$ having the
form $\id_{ U_v}\times\iota_c$ for some $c\in\Phi$. The set
$\{\iota_c:c\in\Phi\}$ generates a finite dihedral group
$L\le\OO(2)$ (where we regard a group of order $2$ as a
degenerate case of a dihedral group), and we have
$\{\id_{ U_v}\}\times L\le G$. The 
$G$-invariance of
$ U_v\times\Phi$ implies that
$\{\id_{ U_v}\}\times L$ is a normal subgroup of $G$.

Let $L_0$ denote the index-two cyclic
subgroup of $L$ consisting of orientation-preserving isometries of
$S^1$.  Choose a component of $\SSS^1-\Phi$,
and let $E$ denote its closure. Then 
$( U_v\times
\SSS^1)/(\{\id_{ U_v}\}\times L_0)$ is canonically identified with
${\rm D}_{ U_v\times\partial E}( U_v\times E)$, so that 
$( U_v\times
\SSS^1)/(\{\id_{ U_v}\}\times L)$ is canonically identified with
${\silv}_{ U_v\times\partial E}( U_v\times E)$. 

Let $G^0$ denote the stabilizer of $ U_v\times E$ in $G$. Note
that $\Isom(E)$ (which has order at most $2$) is canonically
identified with a subgroup of $\Isom(\SSS^1)=\OO(2)$, and that we have
$G^0\le(\Isom( U_v)\times\Isom(E))\cap G\le
\Isom( U_v)\times\OO(2)$. Now
$G^0$  is mapped isomorphically
onto $G/(\{\id_{ U_v}\}\times L)$ by the
quotient homomorphism $G\to G/(\{\id_{ U_v}\}\times L)$. Hence the
%maps the stthis latter
identification of
$( U_v\times
\SSS^1)/(\{\id_{ U_v}\}\times L)$  with
${\silv}_{ U_v\times\partial E}( U_v\times E)$ induces a canonical identification of $( U_v\times
\SSS^1)/G$ with 
${\silv}_{( U_v\times\partial E)/G^0}(( U_v\times E)/G^0)$. We
may therefore regard $\psi$ as an embedding of
${\silv}_{( U_v\times\partial E)/G^0}(( U_v\times E)/G^0)$ in
$\oldOmega=\silv_\oldXi\oldLambda$.

It follows from the last sentence of \ref{better onion} that for any $g\in G^0$, the set
$( U_v\times E)\cap\Fix g$ has no two-dimensional component. Hence
$( U_v\times E)/G^0$ has no two-dimensional stratum. As
$\oldLambda$ is orientable, we may now apply Proposition  \ref{oh yeah he
  tweets}, letting 
$\psi$, 
$( U_v\times E)/G^0$, $( U_v\times\partial E)/G^0$, 
$\oldLambda$ and $\oldXi$ play the respective roles of $j$,
$\oldPsi_1$, $\oldXi_1$, $\oldPsi_2$ and $\oldXi_2$, to obtain  an 
%embedding of
%${\silv}_{( U_v\times\partial E)/G^0}(( U_v\times E)/G^0)$ in
%$\oldOmega=\silv_\oldXi\oldLambda$. 
%there is a unique
 embedding of pairs $\psi^0:
(( U_v\times E)/G^0, ( U_v\times\partial E)/G^0)\to(\oldLambda,\oldXi)$
such
that
% $j^0(\oldXi_1)\subset\oldXi_2$ and 
$\silv \psi^0=\psi$. The latter equality, with the definition of $q$
and the fact that $s$ is a local $\SSS^1$-standardization, implies
that $|\psi^0(( U_v\times E)/G^0)|=|\mu|^{-1}(|\psi(
{\silv}_{( U_v\times\partial E)/G^0}(( U_v\times
E)/G^0))|)=|q|^{-1}(|\frakV|)$. Hence we may regard $\psi^0$ as a
homeomorphism of $( U_v\times E)/G^0$ onto $q^{-1}(\frakV)$.

Set $p_1=p_1^s:\Isom( U_v)\times\OO(2)\to\Isom( U_v)$, and
$X=X_s=p_1(G)$.
Since 
$\{\id_{ U_v}\}\times L\le\ker p_1$, and since
the
quotient homomorphism $G\to G/(\{\id_{ U_v}\}\times L)$
 maps $G^0$ isomorphically
onto $G/(\{\id_{ U_v}\}\times L)$, we have
$X=p_1(G^0)$.

Since $s$ is a local
$\SSS^1$-standardization of $r$, the submersion $\beta:=\zeta\circ r\circ\psi$ is the canonical submersion $( U_v\times
 \SSS^1)/G\to  U_v/X$, i.e. it is induced by the projection
 $ U_v\times \SSS^1\to  U_v$ (see \ref{can sub}). If we set
 $\beta^0=\beta\circ\mu_
{( U_v\times E)/G^0, ( U_v\times\partial E)/G^0}$,
the definitions of $\psi^0$
 and $q$ imply
  that 
%$\zeta\circ q\circ\psi^0:
%( U_v\times E)/G^0\to   U_v/X$ agrees on $\inter( U_v\times E)$ with
%the submersion $\gamma: ( U_v\times E)/G^0\to   U_v/X$ induced by
%the projection
% $ U_v\times J\to  U_v$. Since $\inter( U_v\times E)$ is dense in
 %$ U_v\times E$, 
%we have
 $\zeta\circ q\circ\psi^0=\beta^0$. But $E$, with the intrinsic metric inherited from $\SSS^1$ and scaled
to have length $1$, may be isometrically identified with $[0,1]$. The
 equality $X=p_1(G^0)$, together with the fact that  
$\beta: ( U_v\times
 \SSS^1)/G\to  U_v/X$ is the canonical submersion, 
then implies that $\beta^0: ( U_v\times E)/G^0\to
  U_v/X$ is the canonical submersion described in \ref{can sub},
 with $J=[0,1]$. 
This shows that $s^0:=(\frakV,G^0,\psi^0,\zeta)$ is
a local $[0,1]$-standardization of $q$ at $v$. 

In particular, it follows from
\ref{fibered stuff}  
that
% if
 %$\pi^0:  U_v\times E\to( U_v\times E)/G^0$ denotes the quotient
 %submersion, then
 the fibers of $q|q^{-1}(\frakV)$ are the images under
 $P_{s^0}$ of the arcs of the form $\{x\}\times E$ for $x\in
  U_v$. Hence if $H$ denotes the set of
endpoints of fibers of $q|q^{-1}(\frakV)$, we have  $|P_{s^0}|^{-1}(H)=
U_v\times\partial E= U_v\times (E\cap\Phi)$. But
Condition (II) gives $ U_v\times (E\cap\Phi)=( U_v\times E)\cap Z_s=( U_v\times E)\cap
P_s^{-1}(\oldXi'\cap r^{-1}(\frakV))=P_{s^0}^{-1}(\oldXi\cap
q^{-1}(\frakV))$. Hence $H=|\oldXi\cap
q^{-1}(\frakV)|$. 
Since these observations apply to every $v\in\oldGamma$,  it
now follows that $q:\oldLambda\to\oldGamma$ is an $I$-fibration, and
that with respect to this fibration we have
$\partial_h\oldLambda=\oldXi$, as required.

%\redcomment{Fix the following discussion of the horizontal boundary. The
  %mention of fibers that had preceded the next sentence (now just
  %below this comment) has ben \%-ed out, so I'll need to make the
  %transition from a local standardization to something about fibers.}
%If 
%$\pi^0:
% U_v\times
%E
%\to
%( U_v\times
%E)/G^0$
 %denotes the orbit map, then
%$\psi^0(\pi^0(\{x\}\times \SSS^1))=
%\mu^{-1}(\psi(\pi_s(\{x\}\times \SSS^1)))$ 
%in $\psi:r^{-1}(\frakV)\to( U_v\times
%\SSS^1)/G$
%s a fiber of $q$
%for each $x\in  U_v$. 
%Furthermore, if we
%set $P^0=\psi^0\circ\pi^0: U_v\times E\to q^{-1}(\frakV)$, we
  %  have  
%$(P^0)^{-1}(\oldXi\cap q^{-1}(\frakV))
%=( U_v\times E)\cap\mu^{-1}(P_s^{-1}(\oldXi'\cap
%r^{-1}(\frakV)))=( U_v\times E)\cap\mu^{-1}(Z_s)= U_v\times \partial E$.
%Since the fibers of $q|q^{-1}(\frakV)$ are the suborbifolds of the
%form $P^0(\{x\}\times \SSS^1)$ for $x\in  U_v$, the suborbifold
%$\oldXi\cap q^{-1}(\frakV)=P^0( U_v\times\partial E)$ is the set of
%endpoints of fibers of $q|q^{-1}(\frakV)$. 
%\redcomment{  The next step
%  should be to explain why
  %that's the set of endpoints of fibers (which now seems pretty
  %obvious in view of the description of fibers given above). One thing that
  %needs to be fixed above is that the def. of $r$ should precede the
  %sentence ``Let
%$v\in|\oldGamma|$ be given...'' The exact wording will depend on what
%the next couple of sentences look like.} 
%This shows that $(\oldLambda,\oldXi)$
%is  a \pagelike\ \spair\ in this case.  U_v  U^{2}

There remains the case in which $|\oldGamma|=T$. For the purpose of
the argument in this case, we set $ U^2_{++}=\{(x,y)\in
U^2:x,y\ge0\}\subset  U^2_{+}$. We fix a homeomorphism $k $
of $U^2_+$ onto $ U^2_{++}$ such that
$k ((-1,1)\times\{0\})=(\{0\}\times[0,1))\cup([0,1)\times\{0\})$, and
such that the order-$2$ groups $k^{-1}\,\Isom(U^2_{++})\,k$ and
$\Isom(U^2_+)$ coincide.

If $s=(\frakV,G,\psi,\zeta)$  is
a local
$\SSS^1$-standardization
of $r$ at a point $v\in\oldGamma$, satisfying Alternative (III) of
\ref{better onion}, we will denote by $\ell_s$ the set $\ell$ given by
Alternative (III), and by $\rho_s$  the reflection about
$\ell_s$. According to the last sentence of \ref{better onion},  
 $(\rho_s,\id)$ is the unique element of $G$ having a two-dimensional fixed
point set in $ U_v$, and $\Fix (\rho_s,\id)=\ell_s\times\SSS^1$ is
$G$-invariant. In particular
$\langle(\rho_s,\id)\rangle$ is a central subgroup of $G$.

%In this case we set
%$Y=|r(\oldXi')|$. 
We shall say that a local
$\SSS^1$-standardization
$s=(\frakV,G,\psi,\zeta)$  
of $r$ at a point $v\in\oldGamma$ is {\it normalized} if either
Alternative (I) of \ref{better onion} holds, or
$v\in\partial\oldGamma$ and Alternative (III) holds, or
$v\in\inter\oldGamma$ and Alternative (III) holds and
$\ell_s=(-1,1)\times\{0\}\subset U^2$. The assumption that 
$|\oldGamma|=T$ implies that $r$ admits a normalized local
$\SSS^1$-standardization at every point of $\oldGamma$. 

Consider any point $v\in\oldGamma$, and 
 any
normalized local
$\SSS^1$-standardization  $s=(\frakV,G,\psi,\zeta)$ for $r$ at $v$. If
$s$ satisfies Alternative (I) of \ref{better onion}, we have $
\oldXi'\cap r^{-1}(\frakV)
=P_s(Z_s)=\emptyset$, and in particular
$\oldXi'\cap r^{-1}(v)
=\emptyset$.
If
$s$ satisfies Alternative (III), we have 
$
\oldXi'\cap r^{-1}(\frakV)
=P_s(Z_s)=P_s(\ell_s \times \SSS^1)$; in particular we then have
$r^{-1}(v)=P_s(\{0\}\times\SSS^1)\subset\oldXi'$. This shows that
there is a (unique) set $Y\subset|\oldGamma|$ such that
$|r|^{-1}(Y)=|\oldXi'|$; and that  an arbitrary 
normalized local
$\SSS^1$-standardization  $s$ for $r$ at a point $v\in\oldGamma$
satisfies (I) if $v\notin Y$ and satisfies (III) if $v\in Y$.

For any point $v\in\oldGamma$, let us set
$\ell_v=(-1,1)\times\{0\}\subset U^2= U_v$ if
$v\in\inter\oldGamma$, and
$\ell_v=\{0\}\times[0,1)\subset U^2_+= U_v$ if
  $v\in\partial\oldGamma$. We denote by $\rho_v$ the involution of
  $ U_v$ given by reflection about
  $\ell_v$. We also define a set $K_v\subset U_v$ by
setting $K_v= U_v$ if $v\notin Y$, setting $K_v= U^2_+$ if $v\in
Y\setminus\partial\oldGamma$, and setting
$K_v= U^2_{++}$ if
$v\in Y\cap\partial\oldGamma$.
%if $ U_v= U^2$
%and $\ell_s =(-1,1)\times\{0\}$. 
Then $K_v$ is the closure of a 
component
of $ U_v-\ell_v $ whenever $v\in Y$.  Note that we have $\ell_s=\ell_v$ and $\rho_s=\rho_v$ whenever
  $s$ is  a normalized local
$\SSS^1$-standardization
of $r$ at $v$.

Let $v$ be any point of $\oldGamma$. We define a {\it normalized chart
  centered at $v$} to be a chart $\phi$ for $\oldGamma$ such that
(a) $\phi$ has domain $U_v$, so that the post-chart map defined by $\phi$
induces a homeomorphism of $U_v/X$ onto an open subset of $\oldGamma$
for some finite group $X=X_\phi\le\OO(2)$, (b)
$\phi(0)=v$, and (c) either (i) $\phi(U_v)\cap Y=\emptyset$, or (ii) $\ell_v$ is
 $X$-invariant,
$\phi^{-1}(Y)=\ell_v$, and $\rho_v\in X$. The discussion above shows that for every normalized local
$\SSS^1$-standardization  $s$ for $r$ at $v$, the chart $\phi_s$ (see
\ref{fibered stuff}) is a
normalized chart centered at $v$, and we have $X_{\phi_s}=X_s$
(where $X_s$ is defined by \ref{fibered stuff}). In particular, for each $v\in\oldGamma$ there exists a
normalized chart centered at $v$.

If $\phi$ is a
normalized chart centered at a point $v\in\oldGamma$, we will denote
by $X^1_\phi$ the stabilizer of $K_v$ in $X_\phi$. Note that
$\phi|K_v$ induces a homeomorphism of
 $K_v/X^1_\phi$ onto $|\phi(U_v)|$.

Now for each
$v\in\oldGamma$, we set $U'_v=U^2_+$ if $v\in
\partial\oldGamma$, or if $v$, regarded as a point of $|\oldGamma|$,
lies in $Y$. Otherwise, we set   $U'_v=U^2$. Thus we have $ U'_v=K_v$ except when
$v\in
Y\cap|\partial\oldGamma|$, in which case $K_v= U^2_{++}$ and
$ U'_v= U^2_+$. Let us define a homeomorphism $\kappa_v:U'_v\to K_v$
by setting $\kappa_v=k$ if $v\in
Y\cap|\partial\oldGamma|$, and taking $k$ to be the identity otherwise.
Whenever  $\phi$ is a
normalized chart centered at $v\in\oldGamma$, we  set
$\theta_\phi=(\phi|K_v)\circ\kappa_v:U'_v\to|\phi(U_v)|$. Our choice
of $k$ guarantees that $\kappa_v^{-1}X^1_\phi\kappa_v\le\Isom(U_v')$, and
$\theta_\phi$ induces a  homeomorphism of $|U'_v/(\kappa_v^{-1}X^1_\phi\kappa_v)|$ onto
$|\phi(U_v)|\subset|\oldGamma|$. 

If we let $\phi$ vary over all normalized charts centered at arbitrary
points of
$\oldGamma$, the homeomorphisms  $\theta_{\phi}$ satisfy the compatibility conditions
which is required for them to be an atlas for some orbifold structure
on the space $|\oldGamma|$; indeed, these follow directly from the
compatibility conditions that the $\phi$ satisfy as charts for $\oldGamma$.
This gives a new $2$-orbifold $\oldGamma^1$
with $|\oldGamma^1|=|\oldGamma|$.

The identity map of $|\oldGamma|$ may
be regarded as an orbifold immersion
$\xi:\oldGamma^1\to\oldGamma$. If $\oldPhi$ and $\oldPhi^1$ denote the
open subsets of $\oldGamma$ and $\oldGamma^1$ respectively such that
$|\oldPhi|=|\oldPhi^1|=|\oldGamma|-Y$, then  $\xi$ restricts to a
homeomorphism of $\oldPhi^1$ onto $\oldPhi$.
Note also that for any $v^1\in\oldGamma^1$ we have
%Also, this is the
%  place to point out that
 $ U_{v^1}= U'_{\xi(v_1)}$.

According to the construction we have given, 
for every normalized chart
$\phi$ 
centered at a point $v\in\oldGamma$,
the post-chart map defined by $\theta_{\phi}$ induces a
homeomorphism of $U'_v/(\kappa_v^{-1}X^1_{\phi}\kappa_v)$ onto
$\xi^{-1}({\phi}(U_v))$.

Now define a map of spaces $Q:
|\oldLambda|\to |\oldGamma^1|$ by $Q=|\xi|^{-1}|r|\circ|\mu|$. We
will show that $Q$  defines an $S^1$-fibration of $\oldLambda$ over
$\oldGamma^1$.  As we have shown that
$|r|^{-1}(Y)=|\oldXi'|$, the suborbifold  $\oldXi$ will automatically
be  saturated in this fibration; thus showing that $Q$  defines an $\SSS^1$-fibration will prove that $(\oldLambda,\oldXi)$
is  a \bindinglike\ \spair, thereby completing the proof of the
proposition in this remaining case.

We must show that for each
point $v^1\in\oldGamma^1$, there is a
 local
$\SSS^1$-standardization for $Q$ at $v^1$. Set
$v=\xi(v^1)\in\oldGamma$, so that $ U_{v^1}= U'_v$,
and fix a
normalized local
$\SSS^1$-standardization  $s=(\frakV,G,\psi,\zeta)$ for $r$ at
$v$. 

First consider the subcase
$v\notin Y$. 
In this subcase, we
have $r^{-1}(\frakV)\subset\oldOmega-\oldXi'$, since $s$ satisfies
(I); and we have
$ U'_v= U_v$. 
Since $\xi|\oldPhi^1:\oldPhi^1\to\oldPhi$ is a homeomorphism, $\xi$ maps
some neighborhood $\frakV^1$ of $v^1$ in $\oldGamma^1$
homeomorphically onto 
$\frakV$.
Since
$\mu^*:\oldLambda-\oldXi\to\oldOmega-\oldXi'$ is a homeomorphism by
\ref{silvering}, we may define a homeomorphism
$\psi^1=(\mu^*)^{-1}\circ\psi$ of 
$( U_v\times\SSS^1)/G=( U_{v_1}\times\SSS^1)/G$ onto
$(\mu^*)^{-1}(r^{-1}(\frakV))=
%a suborbifold $\oldPi$
%of
%$\oldLambda$. Since $s$ is a local $\SSS^1$-standardization for $r$, we %have
%$|\psi(( U_v\times \SSS^1)/G)|= r^{-1}(|\frakV|))$ and hence
%$|\oldPi|=
\obd(Q^{-1}(|\frakV^1|))$. The properties required to make 
%fact that
 %$s$ is a local $\SSS^1$-standardization for $r$ also implies that 
%$\zeta:\frakV\to  U_v/X_s$ is an 
 %orbifold homeomorphism, so that
 %$\zeta^1:=\zeta\circ(\xi|\frakV^1)$ is an orbifold homeomorphism
 %of $\frakV^1$ onto $ U_v/X_s= U_{v_1}/X_s$; and that
  %$|\zeta|\circ r\circ|\psi|=|\beta|$, where $\beta:( U_v\times
 %\SSS^1)/G\to  U_v/X$  denotes the canonical submersion, so that
 %$|\zeta^1|\circ Q\circ|\psi^1|=|\beta|$. This shows that, in this
 %subcase, 
$s^1:=
(\frakV^1,G,\psi^1,\zeta^1)$ a local
$\SSS^1$-standardization of $Q$ at $v^1$ follow formally from the
corresponding properties of the local
$\SSS^1$-standardization $s$ of $r$ at $v$. 

To construct the required local $\SSS^1$-standardization  in the
subcase $v\in Y$, recall that if $v\in\inter\oldGamma$ we have $ U_v= U^2$
and $\ell_v =(-1,1)\times\{0\}$, and that if $v\in\partial\oldGamma$
we have $ U_v= U^2_+$ and $\ell_v=\{0\}\times[0,1)$. In either case,
$\rho_s=\rho_v$ is the reflection of $ U_v$ about $\ell_s=\ell_v$.
Since $v\in Y$, we have
$K_v= U^2_+$ if $v\in\inter\oldGamma$, and
%Y\setminus\partial\oldGamma$, and setting
$K_v= U^2_{++}$ if
$v\in \partial\oldGamma$.
We may therefore canonically identify 
$ U_v\times
\SSS^1$ with ${\rm D}_{\ell_s \times S^1}(K_v\times\SSS^1)$ in such a way
that $(\rho_v,\id)$
is identified with the canonical involution of ${\rm
  D}_{\ell_s \times S^1}(K_v\times\SSS^1)$. Hence
$( U_v\times
\SSS^1)/\langle(\rho_v,\id)\rangle$ is canonically identified with
${\silv}_{\ell_s \times\SSS^1}(K_v\times \SSS^1)$.

Let $G^*$ denote the stabilizer of $K_v\times \SSS^1$ in $G$. 
Then $G^*$  is mapped isomorphically
onto $G/\langle(\rho_v,\id)\rangle$ by the
quotient homomorphism $G\to G/\langle(\rho_v,\id)\rangle$. Hence the
%maps the stthis latter
identification of
$( U_v\times
\SSS^1)/\langle(\rho_v,\id)\rangle$  with
${\silv}_{\ell_s \times\SSS^1}(K_v\times \SSS^1)$
induces a canonical identification of $( U_v\times
\SSS^1)/G$ with 
$(\silv_{\ell_s \times\SSS^1}(K_v\times \SSS^1))/G^*$.
We
may therefore regard $\psi$ as an embedding of
$(\silv_{\ell_s \times\SSS^1}(K_v\times \SSS^1))/G^*$
in
$\oldOmega=\silv_\oldXi\oldLambda$. 

Since $v\in Y$, Condition (III) of \ref{better onion} holds, and hence
$\ell_s \times \SSS^1$ is the
only set that occurs as the two-dimensional fixed point set of an
element of $G$. It follows that
 for any $g\in G^*$ the set
$(K_v\times \SSS^1)\cap\Fix g$ has no two-dimensional
component. Hence
$(K_v\times \SSS^1)/G^*$ has no two-dimensional stratum. As
$\oldLambda$ is orientable, we may now apply Proposition  \ref{oh yeah he
  tweets}, letting 
$\psi$, $(K_v\times \SSS^1)/G^*$, $(\ell_s \times \SSS^1)/G^*$,
$\oldLambda$ and $\oldXi$ play the respective roles of $j$,
$\oldPsi_1$, $\oldXi_1$, $\oldPsi_2$ and $\oldXi_2$, to obtain  an 
%embedding of
%${\silv}_{( U_v\times\partial E)/G^0}(( U_v\times E)/G^0)$ in
%$\oldOmega=\silv_\oldXi\oldLambda$. 
%there is a unique
 embedding of pairs $\psi^*:
((K_v\times \SSS^1)/G^*, (\ell_s \times \SSS^1)/G^*)\to(\oldLambda,\oldXi)$
such
that
% $j^0(\oldXi_1)\subset\oldXi_2$ and 
$\silv \psi^*=\psi$. The latter equality, with the definition of $q$
and the fact that $s$ is a local $\SSS^1$-standardization, implies
that $|\psi^*((K_v\times \SSS^1)/G^*)|=|\mu|^{-1}(|\psi(
{\silv}_{\ell_s \times\SSS^1}((K_v\times \SSS^1)/G^*))|)=|\mu|^{-1} (|q|^{-1}(|\frakV|))$. Hence we may regard $\psi^*$ as a
homeomorphism of $
(K_v\times \SSS^1)/G^*$
onto $\obd(|q\circ\mu|^{-1}(|\frakV|)$. If we set 
$\frakV^1=\xi^{-1}(\frakV)\subset\oldGamma^1$,
%$\frakV^1=\psi(\frakV)$, 
it
now follows that
$\psi^1:=\psi^*\circ(\kappa_v\times\id):
U'_v\times\SSS^1\to \obd(Q^{-1}(|\frakV^1|)$ is a homeomorphism.

Now set $G^1=(\kappa_v^{-1}\times\id)G^*(\kappa_v\times\id)$. Since
$\kappa_v$ is  either $k$ or the identity, and
$G^*\le\Isom(K_v)\times\OO(2)$, our choice of $k$
guarantees that in all cases we have $G^1
\le\Isom(U_v')\times\OO(2)=\Isom(U_{v^1})\times\OO(2)$.
According to an observation about general normalized charts made above, the post-chart map defined by $\theta_{\phi_s}$ induces a
homeomorphism of $U'_v/(\kappa_v^{-1}X^1_{\phi_s}\kappa_v)$ onto
$\xi^{-1}({\phi_s}(U_v))=\xi^{-1}(\frakV)=\frakV^1$, whose inverse we denote by $\zeta^1$.  If $p_1:\Isom(U'_v)\times\OO(2)$
denotes the projection to the first factor, we have $p_1(G^1)=
\kappa_v^{-1}X^1_{\phi_s}\kappa_v$, so that $\zeta^1$ is a homeomorphism
from $\frakV^1$ to $U'_v/p_1(G^1)$.
In this subcase, the properties required to make 
$s^1:=
(\frakV^1,G^1,\psi^1,\zeta^1)$ a local
$\SSS^1$-standardization of $Q$ at $v^1$ follow formally from the
corresponding properties of the local
$\SSS^1$-standardization $s$ of $r$ at $v$, together with the
constructions and
the fact that $\psi^1$ and $\zeta^1$ are homeomorphisms.
\EndProof
%\obd canonical \oldUpsilon \silv \eta\alpha $c \frakE \alpha \frakU \frakV \frakE\frakH
%$U \oldGamma \Phi \frakZ E \pi P\beta  \eta $h $f r^{-1}^0 _0' '_0
%\frakK\oldPsi T
%,\psi) A E J P^0\pi^0 {n-1} H G' G^0 r STOP SEARCHING FOR r q J Y Z W
%\xi\frakE \nu (a) \gamma \epsilon H L \Delta \zeta (1) \alpha \beta
%\epsilon \iota $j $h $i \xi \oldGamma' $k J $A \delta $T U \rho (a)
%WK \rho \ell \theta K \frakQ \frakZ \oldPi \frakP \kappa j k P $X Y
%\frakC $X \oldPhi \frakF \frakC \phi \psi'

\Proposition\label{butthurt}
Let 
$(\oldLambda,\oldXi)$ be an acylindrical pair which 
is not an \spair. Suppose that $\oldLambda$ is strongly \simple, 
%\redcomment{With the present def. of strong \simple ity this should be pretty immediate, so not needed as a hyp. Remove it, fix the proof, and check apps} 
and that every component of
$(\partial\oldLambda)-\inter\oldXi$ is an annular orbifold.
Let
$p:\toldLambda\to\oldLambda$ be
a
finite-sheeted covering of $\oldLambda$, and suppose that $\toldLambda$ is strongly \simple. Then
$(\toldLambda,p^{-1}(\oldXi))$ is acylindrical and is not an \spair. 
\EndProposition

\Proof
Since $(\oldLambda,\oldXi)$ is an acylindrical pair and $\oldLambda$ is strongly \simple, it follows
from Proposition \ref{silver acylindrical} that 
$\oldPsi:=\silv_\oldXi\oldLambda$ is weakly \simple. Since $(\oldLambda,\oldXi)$
is  not an \spair, it follows from Proposition \ref{new what they look like} that 
$\oldPsi$ does not admit an
$\SSS^1$-fibration. Furthermore, since every component of 
$\partial\oldLambda-\inter\oldXi$ is annular, every component of
$\partial\oldPsi$ has Euler characteristic $0$ and is therefore
toric. Since $\oldLambda$ is orientable by the definition of an
acylindrical pair, $\fraks_\oldPsi$ has no $0$-dimensional
components. It now follows from Lemma \ref{ztimesz-lemma} that for
every rank-$2$ free abelian subgroup $H$ of $\pi_1(\oldPsi)$, there
is a component $\frakK$ of $\partial\oldPsi$ such that $H$ is
contained in a conjugate of the image of the inclusion homomorphism
$\pi_1(\frakK)\to\pi_1(\oldPsi)$. (The hypothesis in Lemma
\ref{ztimesz-lemma} that $\pi_1(\oldPsi)$ is infinite follows from the
strong \simple ity of $\oldLambda$.)

If we set $\toldXi=p^{-1}(\oldXi)$, then $\toldPsi:=\silv_{\toldXi}\toldLambda$ is a finite-sheeted cover of $\oldPsi$.
% and hence $\inter\toldPsi$ admits a finite-volume hyperbolic metric. 
Lemma \ref{ztimesz-lemma} also
implies that the finite-sheeted cover  $\toldPsi$ of $\oldPsi$ admits no $\SSS^1$-fibration; hence by Proposition
\ref{new what they look like},
$(\toldLambda,\toldXi)$ is not an \spair.

To show that $(\toldLambda,\toldXi)$ is acylindrical, suppose to the contrary that there is an essential annular $2$-orbifold $\oldPi$ in the pair $(\toldLambda,\toldXi)$. Let $q:\toldLambda'\to\toldLambda$ be a finite-sheeted covering such that $\toldLambda'$ is an orientable manifold. Set $\toldXi'=q^{-1}(\toldXi)$. Then $\toldPsi':=\silv_{\toldXi'}\toldLambda'$ is a finite-sheeted covering of $\toldPsi$, and therefore of $\oldPsi$. %$\inter\toldPsi'$ admits a finite-volume hyperbolic metric. 
Since 
the manifold $M:={\rm D}_{\toldXi'}\toldLambda'$ is a two-sheeted cover of $\toldPsi'$, it follows that $M$ is a finite-sheeted cover of $\oldPsi$. Hence every rank-$2$ free abelian subgroup of $\pi_1(M)$ is carried up to conjugacy by some component of $\partial M$.

Now fix a component $A$ of $q^{-1}(\oldPi)$. Since $\toldLambda'$ is a manifold, $A$ is an annulus. Since $\oldPi$ is essential in $(\toldLambda,\toldXi)$, it follows from Corollary \ref{covering annular} that $A$ is essential in $(\toldLambda',\toldXi')$. In particular $A$ is $\pi_1$-injective in $\toldLambda$, and is not parallel in $(\toldLambda',\toldXi')$ to an annulus in $\toldXi'$. Hence the torus $T:={\rm D}A$ is $\pi_1$-injective in $M$. Thus the image of the inclusion $\pi_1(T)\to\pi_1(M)$ is a free abelian group of rank $2$, and is therefore carried up to conjugacy by some component of $\partial M$. This implies (via a standard application of \cite[Lemma 5.1]{Waldhausen}) that $T$ is boundary-parallel in $M$; that is, there is a submanifold $K$ of $M$, homeomorphic to $\TTT^2\times[0,1]$, with $\Fr_M K=T$. The canonical involution $\delta$ of $M={\rm D}_{\toldXi'}\toldLambda'$, which interchanges the two copies of $\toldLambda'$ in $M$, leaves $T={\rm D}A$ invariant; if $\delta$ were to interchange $K$ and $\overline{M-K}$ then $M$ would be homeomorphic to $\TTT^2\times[-1,1]$, a contradiction to finite volume. Hence $K$ is $\delta$-invariant. Thus for some submanifold $J$ of $\toldLambda'$ we have  $K={\rm D}_BJ$, where
$B=\toldXi'\cap J$. Since $K$ is connected, so is $J$.

The essentiality of $A$ implies that $\Fr_{\toldXi'}B=\partial A$ is $\pi_1$-injective in $\toldLambda'$, and in particular in $\toldXi$; hence $B$ is $\pi_1$-injective in $\toldXi$, therefore in $\toldLambda'$, and in particular in $J$. This implies that $J$ is $\pi_1$-injective in $K$, so that $\pi_1(J)$ is free abelian of rank at most $2$. But $\rank\pi_1(J)\ne2$ since $\oldLambda$ is strongly \simple, and $\rank\pi_1(J)\ne0$ since $J$ contains the $\pi_1$-injective torus $A$. Hence $\pi_1(J)$ is infinite cyclic, so that $J$ is a solid torus (see \cite[Theorem 5.2]{hempel}). Since $A$ is a $\pi_1$-injective annulus on the torus $T$, the image of the inclusion homomorphism $\pi_1(A)\to\pi_1(K)$ is a maximal cyclic subgroup of $\pi_1(K)$, and hence the inclusion homomorphism $\pi_1(A)\to\pi_1(J)$ is an isomorphism. If we set $A'=J\cap\overline{\partial\toldLambda'-\toldXi}$, then $T':={\rm D}A'$ is a boundary torus of $K$, and hence $A'$ is an annulus which is homotopically non-trivial in $K$ and therefore in $J$. It now follows that the triad $(J,A,A')$ homeomorphic to $(A,A\times\{1\},A\times\{0\})$, so that
$A$ is parallel in $(\toldLambda',\toldXi')$ to the component $A'$ of $\overline{\partial\toldLambda'-\toldXi}$. This contradicts the essentiality of $A$.
\EndProof
%\obd\oldUpsilon \oldOmega \oldPsi \toldOmega

The following definition and proposition will be used in Chapters \ref{higher chapter} and \ref{general chapter}. The material is included in this section because the proof of Proposition \ref{newer no disks lemma} quotes Proposition \ref{robert strange}.

\Definition\label{reduced intersection}
Let $\oldPsi$ be an orientable $3$-orbifold, and let $\frakP$ be an
orientable closed two-sided $2$-suborbifold of $\inter\oldPsi$. We will say that
a two-sided  $2$-suborbifold $\frakB$ of $\oldPsi$ has {\it reduced
  intersection} with $\frakP$ if $\frakB$ meets $\frakP$
transversally, and there is no compact $2$-suborbifold
$\frakC$ of $\frakP$
 such that (i) $\partial\frakC\ne\emptyset$, (ii) $\frakC\cap\frakB=\partial\frakC$, and (iii) $\frakC$ is
 parallel in the pair $(\oldPsi,\frakB)$
%relative to $\partial\frakC$ 
(see \ref{parallel def}) to
 a suborbifold of $\frakB$.
\EndDefinition

\Proposition\label{newer no disks lemma}
Let 
$\oldPsi$
be 
an  orientable $3$-orbifold which is componentwise strongly \simple, and let $\frakP$ be a closed incompressible $2$-suborbifold of $\inter\oldPsi$.
Then:
\begin{enumerate}
\item
Every two-sided   $\pi_1$-injective $2$-suborbifold $\frakB$ of
$\oldPsi$ is isotopic rel $\partial \frakB$  to a suborbifold which
has reduced intersection with $\frakP$.
\item 
If $\frakB$ is a two-sided  $\pi_1$-injective $2$-suborbifold 
$\oldPsi$ which
has reduced intersection with $\frakP$, then
% that $\omega(A)$ is $\pi_1$-injective \redmissingref{should mention $\pi_1$-injectivity in the app to Prop \ref{dandy subsystem}---that may be done already} in $\oldPsi_1$. Then there is a
%is a properly embedded weight-$0$ annulus 
no component of
$\frakB\cap\frakP$ bounds a discal suborbifold of $\frakB$.
%\redcomment{This condition will be replaced by one saying that $\frakB_0$
  %has reduced intersection with $\frakP$. An extra ``item'' will be
  %added saying that reduced intersection implies the property stated
  %here. The final item will have reduced intersection in its
  %hypothesis. The ``proof summaries'' in the middle of the following proof
  %indicated modifications needed for this changed statement.} 
\item If $\oldLambda$ is any connected \pagelike\ \Ssuborbifold\ of $\oldPsi$
  such that
  $\Fr\oldLambda$ has reduced intersection with $\frakP$, then
%bounds a disk in $\Fr\oldLambda$, then
%  $\oldLambda$ admits an $I$-fibration such that
  %$\oldLambda\cap\partial\oldPsi=\partialh\oldLambda $, and 
each component of 
%the
  %suborbifold
 $\oldLambda\cap\frakP$ 
%of $\oldLambda$ 
is
parallel in the pair $(\oldLambda,\Fr\oldLambda)$
to some (and hence to each) component of
    $\oldLambda\cap\partial\oldPsi$.
%, in which case there is no need to
  %  mention the fibration explicitly here. I'm leaning toward the
%    latter choice.}
\end{enumerate}
\EndProposition

\Proof
We may assume without loss of generality that $\oldPsi$ is connected, and therefore strongly \simple.

We prove Assertion (1) first. After a small isotopy, constant on
$\partial\frakB$, we may assume that $\frakB$ is transverse to
$\frakP$.  Among all two-sided suborbifolds of $\oldPsi$ that
are transverse to
$\frakP$ and isotopic rel boundary to $\frakB$,  let us choose one, $\frakB_0$,
that has the smallest possible number of components of intersection
with $\frakP$. We claim that $\frakB_0$ has reduced intersection with
$\frakP$; this will imply  Assertion (1).

To this end, suppose that $\frakC$ is a compact $2$-suborbifold
of $\frakP$
 such that (i) $\partial\frakC\ne\emptyset$, (ii) $\frakC\cap\frakB_0=\partial\frakC$, and (iii) $\frakC$ is
 parallel in the pair $(\oldPsi,\frakB_0)$ to
 a suborbifold of $\frakB_0$. Thus  there is an embedding $j:\frakC\times[0,1]\to\oldPsi$ such that
$\frakC':=j(\frakC\times\{1\})\subset\frakB_0$,  $j(\frakC\times\{0\})=\frakC$, and
$(\partial\frakC)\times[0,1]\subset j^{-1}(\frakB_0)\subset
\partial(\frakC\times[0,1])$. 
The existence of the embedding $j$ implies that there is 
an (orbifold) ambient isotopy $(\frakh_t)_{0\le t\le1}$ of $\oldPsi$, constant outside an arbitrarily small compact neighborhood $U$ of $\frakJ:=j(\frakC\times[0,1])$ in $\oldPsi$, such that $\frakh_0=\id_\oldPsi$ and $\frakh_1(\frakJ)\cap\frakP=\emptyset$. Set $\frakB_1=\frakh_1(\frakB)$. 
% into an arbitrarily small neighborhood $V$ of $ %:\oldPsi\times\oldPsi\to
%\frakP
% homeomorphism $\eta^+:\frakC\to\frakC'$ which is non-ambiently isotopic in $\frakJ:=j(\frakC\times[0,1])$ to the inclusion map $\frakC\to \frakJ$. We may define an orbifold homeomorphism $\eta^+:\frakB_0\to\frakB_1$ to agree with $\eta$ on $\frakC$ and to be the identity on $\frakB_0-\inter\frakC$. Then $\eta^+$ is isotopic to the inclusion map, so that $\frakB_1$ is isotopic to $\frakB$.
%Modifying $\frakB_1$ by a small isotopy, we
  %obtain an orbifold $\frakB_1'$, still properly embedded in $\oldPsi$ and isotopic to $\frakB$ rel $\partial\frakB$,
%  and having boundary $C\cup C'$, 
%such
The orbifold  $\frakB_1$ is two-sided in $\oldPsi$ and isotopic to $\frakB$ rel $\partial\frakB$. 
If we take the neighborhood $U$ of $\frakJ$ to be sufficiently small, we have
  $\frakB_1\cap\frakP=(\frakB_0\cap\frakP)\setminus\frakC$. Thus $\frakB_1\cap\frakP$ is a union of components of $\frakB_0\cap\frakP$, but does not contain any of the boundary components of $\frakC$. Since $\partial\frakC\ne\emptyset$, it follows that $\compnum(\frakB_1\cap\frakP)<\compnum(\frakB_0\cap\frakP)$, a contradiction to
    the defining property of $\frakB_0$. This proves Assertion
    (1).

To prove Assertion (2), suppose that $\frakB\subset\oldPsi$ is two-sided and  $\pi_1$-injective, and
has reduced intersection with $\frakP$.
Suppose that some component $\oldXi$ of $\frakB\cap\frakP$
bounds a discal suborbifold $\frakD$ of $\frakB$. 
In particular   $\frakD$ is a discal suborbifold of $\oldPsi$.
Since 
$\frakP$ is incompressible in $\oldPsi$, it follows from Proposition \ref{kinda dumb} (more specifically the implication (a)$\Rightarrow$(b), which applies since the strongly \simple\ $3$-orbifold is very good by \ref{oops}) that the curve $c:=|\oldXi|$ bounds a
discal suborbifold of $\frakP$. Since $\frakP$ is two-sided in $\oldPsi$ and therefore orientable, this means that $c$ bounds a
  disk $ F\subset|\frakP|$ such that
  $\wt F\le1$.  

Among all subdisks of $ F$
% weight-$0$
%disks in $\frakP$ 
that are bounded by components of $|\frakB\cap\frakP|$, choose one,
say $ F_0$, that is minimal with respect to inclusion. The
minimality of $ F_0$ implies that
$(\inter F_0)\cap|\frakB|=\emptyset$, i.e. $ F_0\cap|\frakB|=\partial F_0$. On the other hand,
$\omega( F_0)$ is a discal suborbifold of $\oldPsi$. Since
$c_0:=\partial F_0\subset|\frakB|$, and since $\frakB$ is
$\pi_1$-injective in $\oldPsi$, we have $c_0=|\partial\frakD_0|$  for
some discal suborbifold $\frakD_0$ of $\frakB$. Since $\frakD_0$ is
discal and orientable, $D_0:=|\frakD_0|$ is a disk of weight at most $1$. Since
$\inter F_0\cap D_0\subset \inter F_0\cap\frakB=\emptyset$, we
have $ F_0\cap D_0=c_0$. Hence $H:= F_0\cup D_0$ is a
$2$-sphere of weight at most $2$. Furthermore, $\frakH:=\omega(H)$ is
a good orbifold since it is a suborbifold of the very good orbifold
$\oldPsi$. It follows that $\frakH$ is spherical. Moreover, since the $2$-sphere $H$ is in particular an orientable $2$-submanifold of the interior of the orientable $3$-manifold $|\oldPsi|$, the suborbifold $\frakH$ of $\oldPsi$ is two-sided. Since $\oldPsi$
is strongly \simple\ and therefore irreducible, $\frakH$ bounds a discal $3$-suborbifold $\frakE$ of
$\oldPsi$. Since $\wt\frakH\le2$, the set $\fraks_\frakE$ is either a
single arc or the empty set; if $\fraks_\frakE$ is an arc, it has one
endpoint in $ F_0$ and one in $D_0$, and is unknotted in
$|\frakE|$ by the Smith Conjecture \cite{smithbook}.
Hence in any case $\omega( F_0)$ is parallel to $\frakD_0$ in the pair $(\oldPsi,\frakB)$. This
contradicts the assumption that $\frakB$ has reduced intersection with
$\frakP$, and Assertion (2) is  proved. 
%\redcomment{Have I made it clear
 % here that $\omega( F_0)$ meets $\frakB$ only in its boundary,
  %which is Condition (ii) of the def. of reduced intersection? It
  %seems important. I'm not proofreading this proof in detail because
  %rthere are too many issues.}

%\frakB_0

%$ F_0\subset F$\gamma
%we have $wt( F_0)\le\wt F\le1$; furthermore, $\inter F_0\cap D_0\subset$ hence if we set $\frakH=\frakD_0\cup\omega( F_0)$, then $|\frakH|$ is a $2$-sphere, and $\wt(\frakH)\le2$. 

 %$ F_0\cap\fraks_{\oldPsi}=\emptyset$, the curve $\gamma_0:=\partial D_0$ is homotopically trivial in $\oldPsi$. Now since we have observed that $C$ and $C'$ are both homotopically non-trivial in $\oldPsi$, and since $\frakB_0\cap\fraks_{\oldPsi}=\emptyset$, the orbifold $\omega(\frakB_0)$ is $\pi_1$-injective in $\oldPsi$. Hence 
%$\gamma_0$ is homotopically trivial in $\omega(\frakB_0)$, and therefore bounds a disk $D_0\subset \frakB_0$. The minimality of $ F_0$ implies that $\frakB\cap\inter
  % F_0=\emptyset$, and hence $\frakB_1:=(\frakB_0-\inter D_0)\cup F_0$ is an
  %annulus in $\oldPsi$, properly embedded in $\oldPsi$, 
  %and having boundary $C\cup C'$. Furthermore, since $\frakB_0$ and $ F_0$ have weight $0$, the annulus $\frakB_1$ also has weight $0$. 

To prove Assertion (3), suppose that $\oldLambda$ is a \pagelike\ \Ssuborbifold\ of $\oldPsi$
  such that
%  transverse to $\frakP$, such that no component of
  $\Fr\oldLambda$ has reduced intersection with $\frakP$. Let
  $\frakZ$ be a component of $\oldLambda\cap\frakP$. The
definition of an \Ssuborbifold\ implies that $\Fr_\oldPsi\oldLambda$ is
$\pi_1$-injective in $\oldPsi$, and Assertion (2) implies that no component of
$(\Fr_\oldPsi\oldLambda)\cap\frakP$ bounds a discal suborbifold of
$\Fr_\oldPsi\oldLambda$. Hence the $1$-manifold
$(\Fr_\oldPsi\oldLambda)\cap\frakP=\Fr_\frakP\frakZ$ is
$\pi_1$-injective in $\oldPsi$. This implies that $\frakZ$ is
non-discal and $\pi_1$-injective in $\frakP$; as $\frakP$ is in turn
incompressible in $\oldPsi$, it follows that $\frakZ$ is
$\pi_1$-injective in $\oldPsi$, and hence in $\oldLambda$. The strong
\simple ity of $\oldPsi$ implies that the incompressible suborbifold $\oldPi$
is negative, and hence $\chi(\frakZ)<0$. On the other hand, the hypothesis that 
$\oldLambda$ is a \pagelike\ \Ssuborbifold\ of $\oldPsi$ means that
$\oldLambda$ admits an $I$-fibration such that
$\oldLambda\cap\partial\oldPsi=\partialh\oldLambda $. Since
$\frakP\subset\inter\oldPsi$, we have
$\partial\frakZ\subset(\partial\oldLambda)\cap(\inter\oldPsi)=\inter(\partialv\oldLambda)$. Thus
$\oldLambda$ and $\frakZ$ satisfy the hypotheses of Prop.
\ref{robert strange}. It follows that
$\frakZ$ is parallel in the pair  $(\oldLambda,\partialv\oldLambda)$  either to some component of
    $\partialh\oldLambda=\oldLambda\cap\partial\oldPsi$, or to a
    suborbifold of $\partialv\oldLambda=\Fr\oldLambda$. The
    latter alternative cannot hold since   $\Fr\oldLambda$ has
    reduced intersection with $\frakP$, and hence $\frakZ$ is parallel in   $(\oldLambda,\partialv\oldLambda)$  to some, and hence to each, component of
    $\oldLambda\cap\partial\oldPsi$.
Thus
 Assertion (3) is
    proved. 
\EndProof
%\frakC\frakH\oldGamma\frakG\frakE\oldXi\frakX\gamma\eta \frakh \eta
%\sigma \zeta smooth $D E F \oldPhi \Delta \alpha \beta

\section{The characteristic suborbifold}\label{characteristic section}

The proof of the following result, which is the basis of our treatment of the characteristic suborbifold theory, was suggested by the discussion on p. 445 of \cite{bonahon-siebenmann}.

\Proposition\label{new characteristic} Let $\oldPsi$ be a $3$-orbifold
which is orientable, componentwise strongly \simple, and componentwise
boundary-irreducible. 
Then up to (orbifold) isotopy there
exists a unique two-sided $2$-suborbifold $\frakQ$ of $\oldPsi$
such that (1) each component of $\frakQ$ is  an essential annular suborbifold (see
\ref{acylindrical def}) of $\oldPsi$, (2) each component
of $\overline{\oldPsi-\frakH}$, where $\frakH$ is a
strong regular neighborhood of $\frakQ$ in $\oldPsi$, 
%the pair
%$(\oldLambda,\oldLambda\cap\partial\oldPsi)$ 
is either an \Ssuborbifold \ or an \Asuborbifold\ of $\oldPsi$ (see Definition \ref{Do I need it?}), and (3) Condition (2) becomes false if $\frakQ$
is replaced by the union of any proper subset of its components. 
\EndProposition

The two alternatives given in Condition (2) of the above statement are not mutually exclusive: a component $\oldLambda$ of 
$\overline{\oldPsi-\frakH}$ may be both an \Ssuborbifold\ and an \Asuborbifold.

\Proof[Proof of Proposition \ref{new characteristic}]
\nonessentialproofreadingnote{I had a note here saying ``This has been revised a lot, and I was confused
  about things like the strong \simple ity of $\oldLambda$. Be very
  careful about the proofreading.'' After the mega-proofreading, I looked at that note and noticed that the sentence that now begins the proof of \ref{new-claim} was at the beginning of the proof of the lemma, where it made no sense! What am I supposed to do about these things?}
We may assume without loss of generality that $\oldPsi$ is connected,
and therefore strongly \simple\ and boundary-irreducible. According to \ref{oops}, the strong \simple ity of $\oldPsi$ implies that $\partial\oldPsi$ has no spherical component.

Set $\oldOmega=\silv\oldPsi$. It follows from Proposition \ref{silver
  irreducible}, applied with $\partial\oldPsi$ playing the role of $\oldXi$, that $\oldOmega$
is irreducible; indeed, Conditions (1) and (2) of Proposition \ref{silver
  irreducible} follow respectively from the strong \simple ity (cf. \ref{oops}) and the boundary-irreducibility
of $\oldPsi$. It then
follows from the PL version of \cite[p. 444, Splitting Theorem
1]{bonahon-siebenmann} (see \ref{categorille}) that there 
exists, up to (orbifold) isotopy, a unique $2$-dimensional suborbifold
$\frakE$ of $\oldOmega$ such that (1\,$'$) each component of $\frakE$
is toric and is
incompressible in $\oldOmega$, (2\,$'$) each component 
%$\frakB$  
of
$\oldOmega\cut\frakE$ 
either is weakly \simple\ or admits
an $\SSS^1$-fibration (or both), and (3\,$'$)
Condition (2\,$'$) becomes false if $\frakE$ is replaced by the union of
any proper subset of its components. 

We claim:
\Claim\label{new-claim}
Let  $\frakQ$ be any two-sided  suborbifold of $\oldPsi$ (so that $\silv\frakQ$ is identified with a two-sided
suborbifold of $\oldOmega$ according to \ref{first point}).
%, and set
%$\silv\frakQ=\frakE$. 
Then Conditions (1), (2) and (3) of the statement of the present proposition hold for $\frakQ$ if and only if Conditions (1\,$'$), (2\,$'$) and (3\,$'$) hold when we set
$\frakE=\silv\frakQ$. 
\EndClaim

To prove \ref{new-claim}, 
fix a strong regular neighborhood $\frakH$ of $\frakQ$.
Since $\overline{\oldPsi-\frakH}$ is homeomorphic to
  $\oldPsi\cut\frakQ$, it follows from Lemma \ref{oops lemma} that  $\overline{\oldPsi-\frakH}$
is componentwise strongly \simple.

Set $\frakE=\silv\frakQ$ as in the
statement of \ref{new-claim}.  Note that since 
 $\oldPsi$ is strongly \simple, no component of $\frakQ$ is an incompressible
 toric suborbifold (see \ref{oops}). Hence a necessary condition for all components of $\frakE$ to
 be incompressible and toric is that all components of $\frakQ$ be
 annular. It then follows from 
%First note that 
%Let  $\frakQ$ denote the unique suborbifold of $\oldPsi$ such that
%$\silv\frakQ=\frakE$.  \redmissingref{Is the existence and uniqueness of
 % $\frakQ$ obvious, or do I need to say something or give a ref?} 
the first
  assertion of Proposition \ref{silver acylindrical},  applied with
  $\oldPsi$ and $\partial\oldPsi$
  playing the roles of $\oldLambda$ and $\oldXi$, that Condition (1\,$'$)
  for $\frakE$ is equivalent to Condition (1) for $\frakQ$. (The condition in
  Proposition \ref{silver acylindrical} that $\frakE$ not be parallel to
a boundary component of $\oldOmega$ holds vacuously since $\oldOmega$ is closed.)

The second step in the proof of \ref{new-claim} is to show that if Condition (1) holds for $\frakQ$, then
$\frakE$ satisfies Condition (2\,$'$) if and only if  $\frakQ$
satisfies Condition (2). For this purpose, first observe that
 if  $\oldLambda$ is any  component of
$\overline{\oldPsi-\frakH}$, then by \ref{first point},
$\silv_{\oldLambda\cap\partial\oldPsi}\oldLambda$ is (orbifold-)homeomorphic to 
%$\hat \frakB$ for
 some component
% $\frakB$ 
of
$\oldOmega\cut\frakE$; and conversely that, according to \ref{second point},
%for
 every component
% $\frakB$
 of
$\oldOmega\cut\frakE$
%, the orbifold $\hat\frakB$ 
is homeomorphic to
$\silv_{\oldLambda\cap\partial\oldPsi}\oldLambda$ for some
component
$\oldLambda$  of
$\overline{\oldPsi-\frakH}$. 
 Thus Condition (2\,$'$) is equivalent to the assertion that for each component  $\oldLambda$ of  
$\overline{\oldPsi-\frakH}$, the orbifold
$\silv_{\oldLambda\cap\partial\oldPsi}\oldLambda$ either is weakly
 \simple\ or admits
an $\SSS^1$-fibration. But since $\overline{\oldPsi-\frakH}$ is
  componentwise strongly \simple, $\oldLambda$ is strongly \simple, and
  we may therefore apply the second assertion of Proposition \ref{silver
  acylindrical}, with $\oldLambda\cap\partial\oldPsi$ playing the role
of $\oldXi$, to deduce that  $\silv_{\oldLambda\cap\partial\oldPsi}\oldLambda$ is
weakly \simple\ if and only if
the pair
$(\oldLambda,\oldLambda\cap\partial\oldPsi)$
is acylindrical.
%  the pair
%$(\oldLambda,\oldLambda\cap\partial\oldPsi)$ is
%acylindrical. 
%\redmissingref{I've spent a lot of time checking that the
  %term \simple\ is used in a way consistent with the def., which was
  %{\it not} involve boundary-irreducibility. I should keep checking
  %that everything is OK.} 
Likewise, Proposition
\ref{new what
    they look like} implies that
 $\silv_{\oldLambda\cap\partial\oldPsi}\oldLambda$ admits
an $\SSS^1$-fibration if and only if
$(\oldLambda,\oldLambda\cap\partial\oldPsi)$ is an \spair.
But Condition (1) for $\frakQ$ implies that  the
components of $\Fr_\oldPsi\oldLambda$ are essential annular suborbifolds of $\oldPsi$, and hence that
$(\oldLambda,\oldLambda\cap\partial\oldPsi)$
is an acylindrical pair or an \spair\ if and only if $\oldLambda$ is an \Asuborbifold\ or \Ssuborbifold\ of
$\oldPsi$, respectively.
This shows that Condition  (2\,$'$) for $\frakE$ is equivalent to Condition (2) for  $\frakQ$, provided that Condition (1) holds for $\frakQ$. 

The third and final step in the proof of \ref{new-claim} is to show
that if Condition (1) holds for $\frakQ$, then $\frakE$ satisfies Condition (3\,$'$) if and only if
$\frakQ$ satisfies Condition (3). For this purpose,
first observe that
 if $\frakQ'$ is the union
of a proper subset of the components of $\frakQ$, then
$\frakE':=\silv\frakQ'\subset\oldOmega$ is the union
of a proper subset of the components of $\frakE$; and conversely, that 
if $\frakE'$ denotes the union
of a proper subset of the components of $\frakE$, then we may write
$\frakE'=\silv\frakQ'$, where
$\frakQ'$ is the union
of some proper subset of the components of $\frakQ$. But if Condition (1) holds for $\frakQ$, and if 
$\frakQ'$ is the union
of a proper subset of the components of $\frakQ$, then the second step of the proof of \ref{new-claim}, applied with $\frakQ'$ in place of $\frakQ$, shows that 
$\frakE':=\silv\frakQ'$ satisfies Condition (2\,$'$) if and only if $\frakQ'$ satisfies Condition (2). This shows that $\frakE$ fails to satisfy Condition (3\,$'$) if and only if  $\frakQ$ fails to satisfy Condition (3), and completes the proof of \ref{new-claim}.

By the existence assertion of the PL version of \cite[p. 444, Splitting Theorem
1]{bonahon-siebenmann}, we may fix 
a $2$-suborbifold
$\frakE_0$ of $\oldOmega$ such that Conditions (1\,$'$)---(3\,$'$) hold with $\frakE_0$ in place of $\frakE$.
Since $\frakE_0$ is in particular a two-sided suborbifold of
$\oldOmega=\silv\oldPsi$, it follows from \ref{second point} that
$\silv\frakQ_0=\frakE_0$ for some two-sided $2$-suborbifold   $\frakQ_0$ of $\oldPsi$,
.   
By \ref{new-claim}, Conditions (1)---(3) hold with $\frakQ_0$ in place
of $\frakQ$. Now let $\frakQ$ be any two-sided  $2$-suborbifold of
$\oldPsi$ such that (1)---(3) hold, and set $\frakE=\silv\frakQ\subset\oldOmega$. By
\ref{first point} we may write $\frakE_0=\mu_\oldPsi(\frakQ_0)$ and $\frakE=\mu_\oldPsi(\frakQ)$. By \ref{new-claim}, Conditions (1\,$'$)---(3\,$'$) hold for $\frakE$. According to
the uniqueness assertion of the PL version of \cite[p. 444, Splitting Theorem
1]{bonahon-siebenmann}, $\frakE=\mu_\oldPsi(\frakQ)$ is
(orbifold-)isotopic to $\frakE_0=\mu_\oldPsi(\frakQ_0)$ in
$\oldOmega$. Let $h:\oldOmega\to\oldOmega$ be a homeomorphism such
that $h(\frakE)=\frakE_0$. According to Proposition \ref{oh yeah he tweets},
applied with $\oldPsi$ playing the role of both the $\oldPsi_i$, with
$\partial\oldPsi$ playing the role of both the $\oldXi_i$, and  with
$h$ playing the role of $j$, 
there is a homeomorphism $h^0:\oldPsi\to\oldPsi$, isotopic to the identity, such
that $\silv h^0=h$,  and therefore  $\mu_\oldPsi\circ h^0=h\circ\mu_\oldPsi$. It now follows that $h^0(\frakQ)=\frakQ_0$.
\EndProof

\DefinitionNotation\label{oldSigma def}
If $\oldPsi$ is any componentwise strongly \simple, componentwise boundary-irreducible, orientable $3$-orbifold, we will denote by $\frakQ(\oldPsi)$ the suborbifold  $\frakQ$ of $\oldPsi$ (defined up to isotopy) given by Proposition \ref{new characteristic}, and we will denote by $\frakH(\oldPsi) $ a strong regular neighborhood of $\frakQ(\oldPsi)$ in $\oldPsi$. 
%Thus every
%component of $\overline{\oldPsi- \frakQ} $ either is \simple\ or admits
%an $\SSS^1$-fibration. 
We will denote by
$\oldSigma_1(\oldPsi)$ the union of all
components  of $\overline{\oldPsi- \frakH(\oldPsi)} $ that are \Ssuborbifold s (see \ref{Do I need it?}) of $\oldPsi$,
%which admit fibrations under which they are standardly embedded in $\oldPsi$, 
and by $\oldSigma_2(\oldPsi)$ the union of all components of $\frakH(\oldPsi)$ that do not meet any component of $\oldSigma_1(\oldPsi)$.
We will define the (relative)
%{\it absolute
{\it characteristic suborbifold} of $\oldPsi$, denoted
$\oldSigma(\oldPsi)$, to be
$\oldSigma_1(\oldPsi)\cup\oldSigma_2(\oldPsi)$. 
%\redmissingref{Stray sentence: ``Some of these may also be \simple.''
%I think this is covered by the remark after Prop. \ref{new
%characteristic} that the alternatives are not mutually exclusive.}
Thus $\frakH(\oldPsi) $, $\oldSigma_1(\oldPsi)$, $\oldSigma_2(\oldPsi)$ and $\oldSigma(\oldPsi)$,  as well as $\frakQ(\oldPsi) $, are well-defined up to
(orbifold) isotopy. Note that each component of $\oldSigma(\oldPsi)$ is an \Ssuborbifold\ of $\oldPsi$.
\EndDefinitionNotation

\Number\label{tuesa day}
If 
$\oldPsi$ is a componentwise strongly \simple, componentwise
boundary-irreducible, orientable $3$-orbifold, we will set
$\oldPhi(\oldPsi)=\oldSigma(\oldPsi)\cap\partial\oldPsi$ and $\frakA(\oldPsi)=\Fr_{\oldPsi}\oldSigma(\oldPsi)$, so that $\oldPhi(\oldPsi)$ and $\frakA(\oldPsi)$ are
suborbifolds of $\partial|\oldSigma(\oldPsi)|$ whose union is $\partial\oldSigma(\oldPsi)$, and
$\oldPhi(\oldPsi)\cap\frakA(\oldPsi)=\partial \oldPhi(\oldPsi)=\partial\frakA(\oldPsi)$. (One may think of $\oldPhi(\oldPsi)$ as a ``characteristic $2$-orbifold'' by analogy with the notion of a ``characteristic surface'' developed in \cite{bcsz}; the analogy will be pursued further in section \ref{higher section}.) %\redproofreadingnote{I have eliminated essentially all non-canonical uses of $\frakA$, and have replaced $\cala(\oldPsi)$ and $\obd(\cala(\oldPsi)$ almost everywhere by $|\frakA(\oldPsi)|$ and $\frakA(\oldPsi)$. The exceptions are two proofs where $\cala$ appears with subscripts.}
It follows from Proposition \ref{new characteristic} and the definition of $\oldSigma(\oldPsi)$  
% the case where $\oldPsi$ is
%componentwise \simple, 
that every component of $\frakA(\oldPsi)$ is an 
annular $2$-orbifold, essential in $\oldPsi$. The essentiality of the
components of $\frakA(\oldPsi)$ implies that $\partial\oldPhi(\oldPsi)=
\partial\frakA(\oldPsi)$ is $\pi_1$-injective in $\oldPsi$, and
in particular in $\partial\oldPsi$. Thus  $\oldPhi(\oldPsi)$ is taut (see \ref{praxis}) in $\partial\oldPsi$ and is  $\pi_1$-injective in $\oldPsi$.

By \ref{boundary is negative}, the componentwise strong \simple ity and componentwise boundary-irreducibility of $\oldPsi$ imply that each component of
$\partial\oldPsi$ has strictly negative Euler characteristic. In view of the
$\pi_1$-injectivity of
$\partial\oldPhi(\oldPsi)$ in $\partial\oldPsi$, it then follows that
each component of $\oldPhi(\oldPsi)$ or $\overline{\partial\oldPsi-\oldPhi(\oldPsi)}$ has non-positive Euler
characteristic.

Suppose that $\oldLambda$ is a component of $\oldSigma(\oldPsi)$ or of $\overline{\oldPsi-\oldSigma(\oldPsi)}$. Since $\oldLambda$ is a compact $3$-orbifold, we have $\chi(\oldLambda)=\chi(\partial\oldLambda)/2$. On the other hand, since each component of $\Fr\oldLambda$ is a component of $\frakA(\oldPsi)$ and is therefore annular, we have $\chi(\Fr\oldLambda)=0$. It follows that
$\chi(\oldLambda\cap\partial\oldPsi)=\chi(\partial\oldLambda)$, and hence that
 $\chi(\oldLambda)=\chi(\oldLambda\cap\partial\oldPsi)/2$. But $\chi(\oldLambda\cap\partial\oldPsi)\le0$ since $\oldLambda\cap\partial\oldPsi$ is a union of components of $\oldPhi(\oldPsi)$ or of $\overline{\partial\oldPsi-\oldPhi(\oldPsi)}$. Hence $\chi(\oldLambda)\le0$, or equivalently $\chi(\partial\oldLambda)\le0$.

 We will denote by $\book(\oldPsi)$ the union of
$\oldSigma(\oldPsi)$ with all components $\frakC$ of $\frakH(\oldPsi)$ such that $\frakC\cap\oldSigma(\oldPsi)=\Fr_\oldPsi\frakC$.
(The notation is meant to suggest that $\book(\oldPsi)$ is an orbifold
version of a ``book of $I$-bundles.'' See \cite{hyperhaken},
\cite{canmac}, \cite{inject}.) We will
set $\kish(\oldPsi)=\overline{\oldPsi-\book(\oldPsi)}$. Thus $\book(\oldPsi)$ and
$\kish(\oldPsi)$ are suborbifolds of $\oldPsi$.

Note that $\overline{\oldPsi-\oldSigma}$ is the disjoint union of $\kish(\oldPsi)$ with the union of certain components of $\frakH$. Since the components of $\frakH$ are strong regular neighborhoods of two-sided annular suborbifolds of $\oldPsi$, it follows that $\chi(\overline{\oldPsi-\oldSigma})=\chi(\kish(\oldPsi))$.
\EndNumber
%injective essential

%\frakP\oldPi

\Lemma\label{i see a book}
Let $\oldPsi$ be a $3$-orbifold
which is orientable, componentwise strongly \simple, and componentwise
boundary-irreducible. Let $\frakB$ be a compact
suborbifold of $\oldPsi$. Then the following conditions are
equivalent:
\begin{itemize}
\item $\frakB$ is (orbifold-)isotopic to  $\book(\oldPsi)$ in
  $\oldPsi$.
\item There exist a two-sided $2$-suborbifold $\oldPi$ of
  $\oldPsi$, and a strong
regular neighborhood $\oldGamma$  of $\oldPi$ in $\oldPsi$,
such that (a) each component of $\oldPi$ is  an essential annular suborbifold of $\oldPsi$, (b) each component
$\oldLambda$ of $\overline{\oldPsi-\oldGamma}$
%either $\oldLambda$ admits a
%fibration under which it is standardly embedded in $\oldPsi$, or 
%the pair 
%$(\oldLambda,\oldLambda\cap\partial\oldPsi)$ 
is either an \Asuborbifold\ or an \Ssuborbifold\ of $\oldPsi$,
%acylindrical pair or an \spair, 
and (c)
$\frakB$ is the union of $\oldGamma$ with all components of 
$\overline{\oldPsi-\oldGamma}$ that are \Ssuborbifold s.
\end{itemize}
\EndLemma

\Proof
Let $\calp$ denote the set of all two-sided $2$-suborbifolds
$\oldPi$ of $\oldPsi$ such that every component of $\oldPi$ is   an
essential annular suborbifold of $\oldPsi$. Let $\calp_0$ denote the set of all
elements $\oldPi$ of $\calp$ such that 
%for
 each component
%$\oldLambda$ 
of $\overline{\oldPsi-\oldGamma}$, where $\oldGamma$ is a strong regular neighborhood of $\oldPi$ in $\oldPsi$,
%either $\oldLambda$ admits a
%fibration under which it is standardly embedded in $\oldPsi$, or 
%the pair 
%$(\oldLambda,\oldLambda\cap\partial\oldPsi)$ 
is either an \Asuborbifold\ or an \Ssuborbifold.
%acylindrical pair or an \spair. 
For each $\oldPi\in\calp$, if $\oldGamma$ is a strong regular neighborhood of $\oldPi$ in $\oldPsi$, let $\frakB_\oldPi$ denote the union of $\oldGamma$ with all components of 
$\overline{\oldPsi-\oldGamma}$ that are \Ssuborbifold s of $\oldPsi$; thus $\frakB_\oldPi$ is determined up to isotopy by $\oldPi$. We claim:

\Claim\label{untold us before}
If $\oldPi_0$ and $\oldPi$ are elements of $\calp_0$ and $\calp$ respectively, such that $\oldPi_0\subset\oldPi$, then $\oldPi\in\calp_0$, and $\frakB_{\oldPi}$ is isotopic to $\frakB_{\oldPi_0}$.
\EndClaim

To prove \ref{untold us before}, it is enough to consider the case in which $\oldPi-\oldPi_0$ has exactly one component, as the general case then follows by induction on the number of components of $\oldPi-\oldPi_0$. 
In this case, let $\frakC$ denote the essential annular suborbifold $\oldPi-\oldPi_0$ of $\oldPsi$. A strong regular neighborhood of $\oldPi$ in $\oldPsi$ may be written as $\oldGamma=\oldGamma_0\discup\oldGamma_\frakC$, where $\oldGamma_0$ and $\oldGamma_\frakC$ are strong regular neighborhoods of $\oldPi_0$ and $\frakC$ respectively. Let $\oldLambda$ denote the component of $\overline{\oldPsi-\oldGamma_0}$ containing  $\oldGamma_\frakC$. Since $\oldPi_0\in\calp_0$, the suborbifold
$\oldLambda$ of $\oldPsi$ is either an \Asuborbifold\ or an
\Ssuborbifold. 
%It follows from the definition of an \Asuborbifold\ that if
%$\oldLambda$ is an \Asuborbifold\ of $\oldPsi$, then $\frakC$ is non-essential in the pair $(\oldLambda,\oldLambda\cap\partial\oldPsi)$.
%, or $\oldLambda$ is an \Ssuborbifold\ of $\oldPsi$. \redproofreadingnote{Yes, that certainly follows: if $\oldLambda$ is an \Asuborbifold, the def. implies that $\frakC$ is not essential. Do I add a word of explanation or leave it as is?}

Consider first the subcase in which  $\frakC$ is essential in the pair
$(\oldLambda,\oldLambda\cap\partial\oldPsi)$. In
this subcase, 
it follows from the definition of an \Asuborbifold\ that
$\oldLambda$ is not an \Asuborbifold\ of $\oldPsi$, and is therefore an
\Ssuborbifold\ of $\oldPsi$. Furthermore,
 it follows from Lemma \ref{oops lemma} that $\oldLambda$ is strongly \simple. 
%\redproofreadingnote{Revise this, since Proposition \ref{when vertical} no longer mentions strong \simple ity. Actually the problem seems bigger than that. This \Ssuborbifold\ is not assumed to be \pagelike. Doesn't Prop. \ref{when vertical} apply only to \pagelike\Ssuborbifold s? Or does the proof of that prop. not match its statement? By the way, the condition $\chi<0$ which is now in the statement of Prop. \ref{when vertical} holds only for \Ssuborbifold s that are ``strictly \pagelike'' in the sense that they're not also \bindinglike.} 
Thus Corollary \ref{vertical corollary} implies that $\oldLambda$ admits an orbifold fibration, compatible with $\oldLambda\cap\partial\oldPsi$, in which $\frakC$ is saturated. Hence every component of $\overline{\oldLambda-\oldGamma_\frakC}$ is an \Ssuborbifold\ of $\oldPsi$. Since the components of $\overline{\oldPsi-\oldGamma}$ are precisely the components of $\overline{\oldPsi-\oldGamma_0}$ distinct from $\oldLambda$ and the components of $\overline{\oldLambda-\oldGamma_\frakC}$, it follows that in this subcase we have $\oldPi\in\calp_0$ and $\frakB_\oldPi=\frakB_{\oldPi_0}$.

Now consider the subcase in which $\frakC$ is non-essential in 
the pair $(\oldLambda,\oldLambda\cap\partial\oldPsi)$.  Since $\frakC$ is an essential annular orbifold in $\oldPsi$, in particular, $\frakC$ is $\pi_1$-injective in $\oldLambda$, and is not parallel in the pair
$(\oldLambda,\oldLambda\cap\partial\oldPsi)$ to a suborbifold of 
$\oldLambda\cap\partial\oldPsi$. Since $\frakC$ is not essential in $(\oldLambda,\oldLambda\cap\partial\oldPsi)$, it must be
parallel in the pair
$(\oldLambda,\oldLambda\cap\partial\oldPsi)$ to a component of $\Fr_\oldPsi\oldLambda$.
%The last couple of sentences seem to explain the role of the condition of essentiality in the statement. I had originally said $\pi_1$-injective instead of
  %essential, but it
%; this is one of the differences between the conditions
  %here and those in Prop. \ref{new characteristic}. Make sure this
  %works; if not, change the statement and the app. in the proof of
  % It now says essential, so I'm guessing it
  %didn't work. Check whether I've fixed the app. in the proof of Prop. \ref{Comment A}.} 
Hence $\overline{\oldLambda-\oldGamma_\frakC}$ has some component $\frakJ$ such that $\frakJ^+:=\frakJ\cup\oldGamma_\frakC$ is a strong regular neighborhood of $\frakC$ in $\oldPsi$. 
Thus $\frakJ^+$ is a \bindinglike\ \Ssuborbifold\ of $\oldPsi$. Furthermore,  $\oldLambda':=\overline{\oldLambda-\frakJ^+}$, is isotopic to
%$(\oldLambda',\oldLambda'\cap\partial\oldPsi)$ is homeomorphic to
$\oldLambda$,
%\oldLambda\cap\partial\oldPsi)$ 
and is therefore either
%an acylindrical pair or an \spair; that is, $\oldLambda$ is either
 an \Asuborbifold\ or an \Ssuborbifold\ of $\oldLambda$. Since the components of
$\overline{\oldPsi-\oldGamma}$ are precisely $\frakJ$, $\oldLambda'$,
and the components of $\overline{\oldPsi-\oldGamma_0}$ distinct from
$\oldLambda$, we have $\oldPi\in\calp_0$. If $\oldLambda$ is an
\Ssuborbifold\ of $\oldPsi$, we have
$\frakB_\oldPi=\frakB_{\oldPi_0}$. If $\oldLambda$ is not an
\Ssuborbifold\ of $\oldPsi$, we have
$\frakB_\oldPi=\frakB_{\oldPi_0}\cup\frakJ^+$; furthermore,
$\frakB_{\oldPi_0}\cap\frakJ^+$ is then a component of
$\Fr_\oldPsi\oldLambda$, and is the component of the horizontal
boundary of a trivial $I$-fibration of $\frakJ^+$ for which the
vertical boundary is $\frakJ^+\cap\partial\oldPsi$. Hence $\frakB_\oldPi$ is isotopic to $\frakB_{\oldPi_0}$ in this situation. This completes the proof of \ref{untold us before}.

Now set $\frakQ=\frakQ(\oldPsi)$, $\frakH=\frakH(\oldPsi)$, $\oldSigma_i=\oldSigma_i(\oldPsi)$ for $i=1,2$, and $\oldSigma=\oldSigma(\oldPsi)$ (see \ref{oldSigma def}). Note that according to \ref{oldSigma def} and Proposition \ref{new characteristic} we have $\frakQ\in\calp_0$. We claim:
\Claim\label{just as you thought}
$\frakB_{\frakQ}$ is isotopic to $\book(\oldPsi)$.
\EndClaim

To prove \ref{just as you thought}, first note that we may write
$\frakH=\oldSigma_2\discup\frakH_1\discup\frakH_2$, where for $m=1,2$
we denote by $\frakH_m$  the union of all components of $\frakH$ that
meet $\oldSigma_1$ in exactly $m$ frontier components. 
The definition of $\book(\oldPsi)$ gives that $\book(\oldPsi)= \oldSigma\cup\frakH_2=\oldSigma_1\cup\oldSigma_2\cup\frakH_2$.
On the other hand, the definitions imply that $\frakH_1$ is disjoint from $\oldSigma_2$ and $\frakH_2$, and that each component of $\frakH_1$ is the strong regular neighborhood of a two-sided suborbifold of $\oldPsi$ and meets $\oldSigma_1$ in exactly one frontier component. Hence $\frakB_\frakQ=\frakH\cup\oldSigma_1=\oldSigma_1\cup \oldSigma_2\cup\frakH_1\cup\frakH_2$ is isotopic to $\book(\oldPsi)$, as asserted by \ref{just as you thought}.
 %and therefore to $\frakB$. Since Conditions (a), (b) and (c) hold with $\frakQ_0$, $\frakH_0$ and $\frakB_0$ in place of $\frakQ$, $\frakH$ and $\frakB$, an isotopy from $\frakB_0$ to $\frakB$ gives  a two-sided, properly embedded $2$-suborbifold $\frakQ$ of
  %$\oldPsi$, and a
%regular neighborhood $\frakH$  of $\frakQ$, such that (a), (b) and (c) hold.
%(\oldPsi)

The assertion of the present lemma may be paraphrased as saying that a
compact suborbifold $\frakB$ of $\oldPsi$ is isotopic to
$\book(\oldPsi)$ if and only if it is isotopic to $\frakB_\frakP$ for
some $\frakP\in\calp_0$. If $\frakB$ is isotopic to $\book(\oldPsi)$,
then by \ref{just as you thought}, $\frakB$ is isotopic to
$\frakB_\frakQ$, and we have $\frakQ\in\calp_0$. Conversely, suppose
that is isotopic to $\frakB_\frakP$ for some $\frakP\in\calp_0$. Among
all elements of $\calp_0$ that are
% contained in $\oldPi$ (and are
%therefore
 unions of components of $\oldPi$, choose one, $\oldPi_0$,
which has the smallest number of components, and is therefore minimal with respect to inclusion. It then follows from \ref{oldSigma def} and Proposition \ref{new characteristic} that $\oldPi_0$ is isotopic to $\frakQ$. Hence by \ref{just as you thought}, $\frakB_{\oldPi_0}$ is isotopic to $\book(\oldPsi)$. But it follows from \ref{untold us before} that $\frakB_{\oldPi}$ is isotopic to $\frakB_{\oldPi_0}$, and hence to $\book(\oldPsi)$.
\EndProof

\Proposition\label{slide it}
Let $\oldPsi$ be a componentwise strongly \simple, componentwise boundary-irreducible orientable $3$-orbifold. 
%Let $\oldLambda$ be a \Ssuborbifold\ of $\oldPsi$ with $\chi(\oldLambda
Then every (possibly disconnected) \Ssuborbifold\ of $\oldPsi$ 
%whose frontier components are
%essential
 is isotopic to a suborbifold contained in
 $\oldSigma(\oldPsi)$. 
\EndProposition

\Proof
We may assume without loss of generality that $\oldPsi$ is connected, and is therefore strongly \simple\ and boundary-irreducible. Set $\oldSigma=\oldSigma(\oldPsi)$ and $\frakH=\frakH(\oldPsi)$. 

Let $\frakU$ be an \Ssuborbifold\ of $\oldPsi$. 
Let $\calx$ denote the set of all two-sided $2$-orbifolds
$\frakC\subset\oldPsi$ such that (i) $\frakC\cap\frakU=\emptyset$,
(ii) each component of $\frakC$ is an essential annular suborbifold of
$\oldPsi$, and (iii) no two components of $\frakC$ are parallel in
$\oldPsi$. We have $\calx\ne\emptyset$ since $\emptyset\in\calx$. According to Corollary \ref{haken annuli}, there is a natural number $N$ such that $\compnum(\frakC)\le N$ for every $\frakC\in\calx$. Hence $\calx$ has an element $\frakC_0$ which is maximal with respect to inclusion. Let $\frakK$ be a strong regular neighborhood of the two-sided $2$-orbifold $\frakC_0\discup\Fr_\oldPsi\frakU$ in $\oldPsi$.
% such that $\frakK\cap\frakU=\emptyset$.
We claim:
\Claim\label{carry me back}
Each component
%$\oldLambda$ 
of $\overline{\oldPsi-\frakK}$
%, the pair
%$(\oldLambda,\oldLambda\cap\partial\oldPsi)$ 
is either an \Ssuborbifold\  or an \Asuborbifold.
\EndClaim
To prove \ref{carry me back}, note that for any component
$\oldLambda$ of $\overline{\oldPsi-\frakK}$ we have either (a) $\oldLambda\subset\frakU$ or (b) $\oldLambda\cap\frakU=\emptyset$. If (a) holds, then since $\frakC_0\cap\frakU=\emptyset$, we may describe $\oldLambda$ as the closure of the suborbifold of $\oldPsi$
%the suborbifold
that is obtained  from some component $\oldLambda'$ of
$\frakU$ by removing a strong regular neighborhood of $\Fr\oldLambda'$;
thus $\oldLambda$ is isotopic to $\oldLambda'$, and hence 
$\oldLambda$ is an \Ssuborbifold.
%, i.e.
%$(\oldLambda,\oldLambda\cap\partial\oldPsi)$ is an \spair. 
We shall complete the proof of \ref{carry me back} by showing that if (b) holds then 
$\oldLambda$ is an \Asuborbifold.
%\frakX

Suppose that $\oldPi$ is an essential  annular orbifold in the pair
$(\oldLambda,\oldLambda\cap\partial\oldPsi)$. Since $\oldPi\subset\oldLambda\subset \overline{\oldPsi-\frakK}$ we have $\oldPi\cap\frakC_0=\emptyset$; since $\oldPi$ and $\frakC_0$ are two-sided $2$-suborbifolds
of $\oldPsi$, it follows that $\frakC_1:=\oldPi\cup\frakC_0$ is also a two-sided $2$-suborbifold of
$\oldPsi$. We have $\frakC_0\cap\frakU=\emptyset$ since $\frakC_0\in\calx$, and $\oldPi\cap\frakU=\emptyset$ since (b) holds and $\oldPi\subset\oldLambda$. Hence $\frakC_1\cap\frakU=\emptyset$. Since $\oldPi$ is essential in $(\oldLambda,\oldLambda\cap\partial\oldPsi)$, it follows from Corollary \ref{still essential}, with $\oldPsi$ and $\oldPi$ defined as above, and with $\oldLambda$ playing the role of $\oldPsi'$ in the corollary, that $\oldPi$ is essential in $\oldPsi$. Thus Conditions (i) and (ii) of the definition of $\calx$ hold with $\frakC_1$ in place of $\frakC$. But the maximality of $\frakC_0$ implies that $\frakC_1\notin\calx$. Hence Condition (iii) of the definition of $\calx$ cannot hold with $\frakC_1$ in place of $\frakC$; that is, $\frakC$ has two components that are parallel in
$\oldPsi$. Since Condition (iii) of the definition of $\calx$ does hold with $\frakC_0$ in place of $\frakC$, it follows that $\oldPi$ is parallel in
$\oldPsi$ to some component $\frakD$
%\redcomment{Name it, and make sure parallelism in $\oldPsi$ is the
%relevant thing} 
of $\frakC_0$.
There therefore exists 
a connected suborbifold $\oldUpsilon$ of $\oldPsi$ admitting a trivial $I$-fibration such that
% an embedding $j:\oldPi\times[0,1]\to\oldLambda$ such that
$\partialv\oldUpsilon\subset\partial\oldPsi$, and the components of
$\partialh\oldUpsilon$ are $\oldPi$ and $\frakD$. 
% an embedding $j:\oldPi\times[0,1]\to\oldLambda$ such that
%$j((\partial\oldPi)\times[0,1])\subset\partial\oldPsi$, $j(\oldPi\times\{0\})=\oldPi$, an%d
%$j(\oldPi\times\{1\})=\frakD$. Set
%$\oldUpsilon=j(\oldPi\times[0,1])$. 
Since $\frakD\subset\frakC_0$ and $\oldLambda\subset
\overline{\oldPsi-\frakK}$, we have $\frakD\cap\oldLambda=\emptyset$; since, in addition,
$\oldPi$ is properly embedded in
$\oldLambda$ and disjoint from $\Fr\oldLambda$, the intersection of $\oldUpsilon$ with
$\Fr\oldLambda$ is a non-empty union of components of
$\Fr\oldLambda$. The components of $\oldUpsilon\cap\Fr\oldLambda$ are disjoint from $\oldPi$ and $\frakD$, and hence  their boundaries are contained in $\inter\partialv\oldUpsilon$.  Furthermore, since the components of
$
\oldUpsilon\cap\Fr\oldLambda$
%\subset
%$\Fr\oldLambda\subset\Fr\frakK$ 
are frontier components of the \Ssuborbifold\ $\oldLambda$ of $\oldPsi$, they are essential annular suborbifolds of $\oldPsi$, and in particular they are $\pi_1$-injective in $\oldUpsilon$. 
%parallel to components of
%$\frakC_0\discup\Fr_\oldPsi\frakU$, they are essential in
%$\oldPsi$. 
%It therefore follows from Corollary \ref{still essential} that they are essential in $\oldUpsilon,
Hence if $\oldXi$ is any component of $\oldUpsilon\cap\Fr\oldLambda$, we may apply
%Using the $\pi_1$-injectivity of $\oldUpsilon\cap\Fr\oldLambda$,
%and applying
 Proposition \ref{robert strange}, with $\oldUpsilon$ and $\oldXi$
%an arbitrary component of $\oldUpsilon\cap\Fr\oldLambda$,  
playing the respective roles of $\oldLambda$ and $\frakZ$ in that proposition, 
%\redmissingref{To apply Prop. \ref{robert strange} we need to know 
%$\partial
%(\oldUpsilon\cap\Fr\oldLambda)
%\subset\inter\partialv\oldUpsilon$}
to deduce that each component of
$\oldUpsilon\cap\Fr\oldLambda$ is parallel in the pair  $(\oldUpsilon,\oldUpsilon\cap\partial\oldPsi)$  either to a suborbifold of $\oldUpsilon\cap\partial\oldPsi$ or to
$\oldPi$. 
But since the components of
$
\oldUpsilon\cap\Fr\oldLambda
$ are essential in
$\oldPsi$, none of them can be parallel to a suborbifold of
$\oldUpsilon\cap\partial\oldPsi$. 
Thus $\oldXi$
%for each component of
%$\oldUpsilon\cap\Fr\oldLambda$ 
is parallel in the pair $(\oldUpsilon,\oldUpsilon\cap\partial\oldPsi)$ 
to 
$\oldPi$; that is, there is a suborbifold $\oldUpsilon_\oldXi$ of $\oldUpsilon$ with an $I$-fibration such that $\oldPi$ and $\oldXi$ are the components of $\partialh\oldUpsilon_\oldXi$, and $\partialv\oldUpsilon_\oldXi\subset \oldUpsilon\cap\partial\oldPsi)$. Since $\Fr_\oldPsi\oldUpsilon_\oldXi=\oldPi\cup\oldXi$ for every $\oldXi\in\calc(\oldUpsilon\cap\Fr\oldLambda)$,  the set $\{\oldUpsilon_\oldXi:\oldXi\in\calc(\oldUpsilon\cap\Fr\oldLambda)\}$ is linearly ordered by inclusion. Let us write its least element as $\oldUpsilon_{\oldXi_0}$ for some $\oldXi_0\in\calc(\oldUpsilon\cap\Fr\oldLambda)$. Then $\oldUpsilon_{\oldXi_0}$ meets no components of $\Fr_\oldPsi\oldLambda$ other than $\oldXi_0$, and is therefore contained in $\oldLambda$.
% component $\oldUpsilon_0$ of
%$\oldUpsilon\cap\oldLambda$ containing $\oldPi$ admits an $I$-fibration for
%which $\partialv\oldUpsilon_0\subset\oldLambda\cap\partial\oldPsi$, and the
%components of $\partialh\oldUpsilon_0$ are $\oldPi$ and a component
%$\frakD_0$ of $\Fr\oldLambda$. 
This means that $\oldPi$ is parallel to
$\oldXi_0$ in $(\oldLambda, \oldLambda\cap\partial\oldPsi)$, a contradiction to the essentiality of
$\oldPi$ in
$(\oldLambda,\oldLambda\cap\partial\oldPsi)$. This shows that
$(\oldLambda,\oldLambda\cap\partial\oldPsi)$ is acylindrical when (b)
holds, and completes the proof of \ref{carry me back}.

It follows from \ref{carry me back} that there is a suborbifold
$\frakK_0$ of $\oldPsi$ such that (1) $\frakK_0$ is a union of
components of $\frakK$, (2) each component
 of $\overline{\oldPsi-\frakK_0}$
%, the pair
%$(\oldLambda,\oldLambda\cap\partial\oldPsi)$ 
is either an \Ssuborbifold\  or
an \Asuborbifold\ of $\oldPsi$, and (3) Condition (2) becomes false if
$\frakK_0$ is replaced by the union of any proper subset of its
components. According to Proposition \ref{new characteristic} and
\ref{oldSigma def}, $\frakK_0$ is isotopic to $\frakH$. But since $\frakK$ is a regular neighborhood of
  $\frakC_0\discup\Fr_\oldPsi\frakU$, the suborbifold $\frakK$ is isotopic to  a suborbifold of $\oldPsi$ disjoint from $\frakU$; in particular, $\frakK_0$ is isotopic to  a suborbifold of $\oldPsi$ disjoint from $\frakU$. Hence
we may suppose $\frakH$ to be chosen within its isotopy class in
such a way that $\frakH\cap\frakU=\emptyset$.

Now set $\oldSigma_1=\oldSigma_1(\oldPsi)$ and $\oldSigma^*=\overline{\oldPsi-(\frakH\cup\oldSigma_1)}$; thus $\oldSigma^*$ is the union of all components of
$\overline{\oldPsi-\frakH}$ which are not \Ssuborbifold s of
$\oldPsi$. Since 
$\frakH\cap\frakU=\emptyset$, we
may write $\frakU=\frakU_1\discup\frakU^*$, where
$\frakU_1\subset\oldSigma_1$ and 
$\frakU^*\subset\oldSigma^*$.
We claim:  
\Claim\label{new squeeze it} For every component $\frakB$ of $\Fr_\oldPsi\frakU^*$,
there is a unique connected suborbifold $\frakV_\frakB$ of $\oldSigma^*$
such that (A) $\Fr_\oldSigma\frakV_\frakB=\frakB$ and (B)
$\frakV_\frakB$ has a trivial $I$-fibration for which one component
of $\partialh\frakV_\frakB$ is
$\frakB$ and the other is a component of
$\Fr_\oldPsi\oldSigma^*$. 
%believe
  %the latter is correct and is better, but I should double-check
  %this. I refer to
  %$\Fr_\oldPsi\oldSigma^*\subset\Fr_\oldPsi\frakH$ in the
  %proof of uniqueness below.} 
\EndClaim

To prove \ref{new squeeze it}, let $\oldLambda$ denote the component of $\oldSigma^*$ containing $\frakB$. The definition of $\oldSigma^*$ implies that $\oldLambda$ is a component of $\overline{\oldPsi-(\frakH\cup\oldSigma_1)}$ and is not an \Ssuborbifold\ of $\oldPsi$; in view of the defining properties of $\oldSigma_1$ (see Proposition \ref{new characteristic} and
\ref{oldSigma def}),  $\oldLambda$ 
is an \Asuborbifold, which by definition means that $(\oldLambda,\oldLambda\cap \partial\oldPsi)$ is an acylindrical pair. Hence the two-sided annular orbifold  $\frakB$ is not essential in $(\oldLambda,\oldLambda\cap \partial\oldPsi)$. On the other hand, since the component $\frakB$ of $\frakB$ of $\Fr_\oldPsi\frakU^*$ is in particular a frontier component of the \Ssuborbifold\ $\frakU$ of $\oldPsi$, the annular suborbifold $\frakB$ of $\oldPsi$ is essential in $\oldPsi$; in particular, $\frakB$ is $\pi_1$-injective in $\oldPsi$, and in particular in $\oldLambda$. Hence $\frakB$ is parallel in $(\oldLambda,\oldLambda\cap\partial \oldPsi)$ 
to
a suborbifold of either $\oldLambda\cap\partial\oldPsi$ or
$\Fr_\oldPsi\oldLambda$. The former alternative is impossible, because $\frakB$, as a frontier component of the \Ssuborbifold\ $\frakU$, is essential in $\oldPsi$. Hence $\frakB$  is parallel 
in $(\oldLambda,\oldLambda\cap\partial \oldPsi)$ to
a suborbifold of 
$\Fr_\oldPsi\oldLambda$, and the existence of a suborbifold
$\frakV_\frakB$ with the asserted property follows. 
To prove uniqueness, 
note that if $\frakV$ and $\frakV'$ are distinct suborbifolds of
$\oldSigma^*$ such that (A) and (B) hold with $\frakV$ in place of $\frakV_\frakB$ and also with $\frakV'$ in place of $\frakV_\frakB$, then in particular $\frakV$ and $\frakV'$ are distinct suborbifolds of $\oldLambda$ with frontier $\frakB$, which implies that $\frakV\cap\frakV'=\frakB$ and that $\frakV\cup\frakV'=\oldLambda$. Condition (B) for $\frakV$ and $\frakV'$ then implies that 
$\oldLambda$ has a trivial $I$-fibration under which its horizontal boundary is contained in
$\Fr_\oldPsi\oldSigma^*$; since the components of $\Fr_\oldPsi\oldSigma^*$ are essential annuli, this implies that $\oldLambda$ is a \torifold\ whose frontier components are essential annuli. Hence by Lemma \ref{when a tore a fold}, $\oldLambda$ is a (\bindinglike) \Ssuborbifold\ of $\oldPsi$, a contradiction. Thus \ref{new squeeze it} is proved.

Now we claim: 
\Claim\label{barbara stanwick}
If $\frakB$ and $\frakB'$ are components of $\Fr_\oldPsi\frakU^*$, and if $\frakV_\frakB$ and $\frakV_{\frakB'}$ are defined by \ref{new squeeze it}, we have either $\frakV_\frakB\cap\frakV_{\frakB'}=\emptyset$, $\frakV_\frakB\subset\frakV_{\frakB'}$, or $\frakV_{\frakB'}\subset\frakV_{\frakB}$.
\EndClaim

In proving \ref{barbara stanwick}, we may assume that
$\frakB\ne\frakB'$. Note that if
$\frakV_\frakB$ and $\frakV_{\frakB'}$ are contained in distinct components of $\oldSigma^*$ then they are disjoint. Hence we may assume that they are contained in the same component $\oldLambda$ of $\oldSigma^*$. Next note that
 $\Fr_\oldLambda\frakV_{\frakB}=\frakB$ and
 $\Fr_\oldLambda\frakV_{\frakB'}=\frakB'$ are distinct components of
 $\Fr_\oldPsi\frakU^*$; in particular they are disjoint connected sets. Since the
 frontiers of $\frakV_\frakB$ and $\frakV_{\frakB'}$ relative to the
 connected space $\oldLambda$ are connected and disjoint, \ref{barbara
   stanwick} must hold unless
 $\frakV_\frakB\cup\frakV_{\frakB'}=\oldLambda$. 
Suppose that the latter equality holds.
Then $\frakB'\subset\frakV_{\frakB}$, and $\partial\frakB'\subset\inter{\frakV_{\frakB}\cap\partial\oldPsi}=\inter(\partialv\frakV_{\frakB})$, where 
$\frakV_\frakB$ is equipped with the 
trivial $I$-fibration
given by Condition (B) of \ref{new squeeze it}. Furthermore, ${\frakB'}$, as a frontier
component of the \Ssuborbifold\ $\frakU$, is an essential annular orbifold in
$\oldPsi$, and in particular is $\pi_1$-injective in $\frakV_\frakB$. We may therefore apply 
Proposition \ref{robert strange}, letting $\frakV_\frakB$ and $\frakB'$ play the respective roles of $\oldLambda$ and $\frakZ$ in that proposition, to deduce that ${\frakB'}$ is
parallel in the pair $(\frakV_\frakB,\frakV_\frakB\cap\partial\oldPsi)$ either to
the component $\oldGamma:=\frakV_\frakB\cap\Fr_\oldPsi\oldLambda$ of $\partialh\frakV_\frakB$, or to an annular suborbifold of $\partialv\frakV_\frakB$. 
The second of these alternatives cannot occur, since ${\frakB'}$ is 
%a
%component of $\Fr_\oldPsi\frakU^*\subset\Fr_\oldPsi\frakU$ and is therefore 
essential in
$\oldPsi$. Hence the first alternative must hold, so that some suborbifold $\frakJ$ of $\frakV_\frakB$ admits a trivial $I$-fibration in which its horizontal boundary components are ${\frakB'}$ and $\oldGamma$, and its vertical boundary is contained in $\frakV_\frakB\cap\partial\oldPsi$.
%the other is  a component of $\Fr_\oldPsi\oldSigma^*$. 
It follows that $\Fr_{\oldLambda}\frakJ=\frakB'=\Fr_\oldLambda\frakV_\frakB'$, and since $\oldGamma\subset\frakJ$ is disjoint from $\frakV_{\frakB'}$, we must have $\frakJ=\overline{\oldLambda-\frakV_{\frakB'}}$. But by 
Condition (B) of \ref{new squeeze it}, $\frakV_{\frakB'}$ also admits
a trivial $I$-fibration in which $\frakB'$ is one of its horizontal boundary
components, and $\frakV_{\frakB'}\cap\partial\oldPsi$ is its vertical boundary. It follows that $\oldLambda=\frakV_{\frakB'}\cup\frakJ$ admits
a trivial $I$-fibration in which its vertical boundary
is $\oldLambda\cap\partial\oldPsi$; since the frontier components of the component $\oldLambda$ of $\oldSigma$
%; since the components of $\Fr_\oldPsi\oldSigma^*$ 
are essential annuli, this implies that $\oldLambda$ is a \torifold\ whose frontier components are essential annuli. Hence by Lemma \ref{when a tore a fold}, $\oldLambda$ is a (\bindinglike) \Ssuborbifold\ of $\oldPsi$, a contradiction. Thus \ref{new squeeze it} is proved.

%, and hence that $\oldLambda$
%is an \Ssuborbifold\ of $\oldPsi$, a contradiction. This proves \ref{barbara stanwick}.
%(i) (ii) \frakD

Now set $\frakZ=\bigcup_{\frakB\in\calc(\Fr_{\frakU^*})}\frakV_\frakB\subset\oldSigma^*$, where the $\frakV_{\frakB}$ are defined by \ref{new squeeze it}. It follows from 
\ref{barbara stanwick} that $\frakZ$ is a {\it disjoint} union of sets of the form $\frakV_\frakB$, where $\frakB$ ranges over certain components of $\Fr_\oldPsi\frakU^*$. The defining properties of the $\frakV_\frakB$ given by \ref{new squeeze it} imply that each $\frakV_\frakB$ admits a trivial $I$-fibration for which $\frakV_\frakB\cap\Fr_\oldPsi\oldSigma^*$  is one component of $\partialh\frakV_\frakB$. Hence:
\Claim\label{byline baloney}
Each component $\frakV$ of $\frakZ$ admits a trivial $I$-fibration for
which $\frakV\cap\Fr_\oldPsi\oldSigma^*$  is one component of
$\partialh\frakV$. 
\EndClaim

On the other hand, we claim that
\Claim\label{bobaloo}
$\frakU^*\subset\frakZ$.
\EndClaim
To prove \ref{bobaloo}, consider an arbitrary component $\frakE$ of $\frakU^*$, and let $\oldLambda$ denote the component of $\oldSigma^*$ containing $\frakE$.
For each component $\frakB$ of
$\Fr_\oldPsi\frakE$ we have $\Fr\frakV_\frakB=\frakB$ by \ref{new squeeze it}. Hence either
$\frakV_\frakB\cap\frakE=\frakB$ or
$\frakV_\frakB\supset\frakE$. If
$\frakV_\frakB\cap\frakE=\frakB$ for every component $\frakB$ of
$\Fr\frakE$, then $\oldLambda$ is a strong regular neighborhood of $\frakE$
and is therefore an \Ssuborbifold\ of $\oldPsi$, 
%; that is,  $(\oldLambda,\oldLambda\cap\partial\oldPsi)$
%is an \spair, 
a contradiction to the definition of $\oldSigma^*$. Hence there is a component $\frakB_0$ of
$\Fr_\oldPsi\frakE$ such that $\frakV_{\frakB_0}\supset\frakE$; thus $\frakE\subset\frakV_{\frakB_0}\subset\frakZ$, and \ref{bobaloo} is proved.

Now let $\oldSigma^!$ denote the union of $\oldSigma_1$ with all those components of $\frakH$ whose frontiers are not entirely contained in $\oldSigma_1$. Then $\oldSigma\subset\oldSigma^!$, and each component of $\overline{\oldSigma^!-\oldSigma}$ is a component of $\frakH$ sharing exactly one frontier component with $\oldSigma$. Since $\frakH$ is a strong regular neighborhood of a two-sided suborbifold of $\oldPsi$, it follows that $\oldSigma^!$ is isotopic to $\oldSigma$. We claim:

\Claim\label{what am i doing}  $\oldSigma^!\cup\frakZ$ is isotopic to $\oldSigma^!$, and hence to $\oldSigma$.
\EndClaim

To prove \ref{what am i doing}, consider
an arbitrary component $\frakV$   of $\frakZ$. According to \ref{byline baloney}, $\frakV$  admits a trivial $I$-fibration for
which $\partialh^+\frakV:=\frakV\cap\Fr_\oldPsi\oldSigma^*$  is one component of
$\partialh\frakV$. Since $\oldSigma^*=\overline{\oldPsi-(\frakH\cup\oldSigma_1)}$, we have
$\partialh^+\frakV=\frakV\cap\Fr_\oldPsi(\frakH\cup\oldSigma_1)$. Since $\oldSigma_1$ is a union of components of $\overline{\oldPsi-\frakH}$, we have $\partialh^+\frakV\subset\Fr\frakH$; and since $\frakV\subset\frakZ\subset\oldSigma^*$, the component $\frakH_0$ of $\frakH$ containing $\partialh^+\frakV$ does not have its frontier contained entirely in $\oldSigma_1$; that is, $\frakH_0\subset\oldSigma^!$, and hence $\partialh^+\frakV\subset\Fr_\oldPsi\oldSigma^!$. Thus each component of  $\frakZ$ meets $\oldSigma^!$ precisely in one component of its horizontal boundary (with respect to its trivial $I$-fibration). This implies \ref{what am i doing}.

To prove the conclusion of the proposition, we need only note that
 we have $\frakU_1\subset\oldSigma_1\subset\oldSigma^!$ by the definitions, and $ \frakU^*\subset\frakZ$ by \ref{bobaloo}; hence $\frakU=\frakU_1\cup\frakU^*\subset\oldSigma^!\cup\frakZ$. Since $\oldSigma^!\cup\frakZ$ is isotopic to $\oldSigma$ by \ref{what am i doing}, the conclusion follows.

%It follows that $\frak
%Note also that $\Fr_\oldPsi\oldSigma^*$ is a union of components of $\Fr_\oldPsi
\EndProof
%\frakE\frakC\frakD\frakB\frakJ\frakZ\oldXi\frakX\frakH(\oldPsi)\oldSigma( canonical

\Lemma\label{no west for the bleary}
Let $\oldPsi$ be a  strongly \simple, componentwise boundary-irreducible, orientable $3$-orbifold.
If  $\oldDelta$ is a component of $\overline{\oldPsi-\oldSigma(\oldPsi)}$ which is an \Ssuborbifold\ of $\oldPsi$,
%dmits a fibration under which it is standardly embedded in $\oldPsi$, 
then 
$\oldDelta$ is a component of $\frakH(\oldPsi)$, and
every component of $\overline{\oldPsi- \frakH(\oldPsi)} $ which shares a frontier component with $\oldDelta$ is a component of $\oldSigma_1(\oldPsi)$.
\EndLemma

\Proof
Recall from Definition \ref{oldSigma def}
that $\frakH:=\frakH(\oldPsi)$ is 
%constructed from  
a strong regular neighborhood of 
the $2$-suborbifold $\frakQ(\oldPsi) $ in $\oldPsi$.
%, well-defined up to isotopy, having the properties stated in Proposition \ref{new characteristic}. 
Thus every component of $\frakH$ has two frontier components relative to $\oldPsi$. We have $\oldSigma:=\oldSigma(\oldPsi)=\oldSigma_1\cup\oldSigma_2$, where according to \ref{oldSigma def}, $\oldSigma_1:=\oldSigma_1(\oldPsi)$ is 
the union of all
components  of $\overline{\oldPsi- \frakH} $ that are \Ssuborbifold s of $\oldPsi$;
and $\oldSigma_2 :=\oldSigma_2(\oldPsi)$ is the union of all components of $\frakH $ that do not meet any component of $\oldSigma_1$. It follows that if $\oldDelta$ is any component of $\overline{\oldPsi-\oldSigma}$, then either (i) $\oldDelta$ is a component of $\frakH$, and every component of $\overline{\oldPsi- \frakH} $ which shares a frontier component with $\oldDelta$ is a component of $\oldSigma_1$, or (ii) there is a component $\oldLambda$ of $\overline{\oldPsi-\frakH}$, which is {\it not} a component of $\oldSigma_1$, such that $\oldDelta$ is the union of $\oldLambda$ with all components of $\frakH$ that share one frontier component with $\oldLambda$ and one with $\oldSigma_1$. If (ii) holds, then $\oldDelta$ is a strong regular neighborhood of $\oldLambda$, and $\oldLambda$ is {\it not} an \Ssuborbifold\ of $\oldPsi$;
%  admit a fibration under which it is standardly embedded in $\oldPsi$; 
therefore $\oldDelta$
is {\it not} an \Ssuborbifold\ of $\oldPsi$.
% fibration under which it is standardly embedded in $\oldPsi$. 
The conclusion follows.
\EndProof

\Proposition\label{lady edith}
Let $\oldPsi$ be a componentwise strongly \simple, componentwise boundary-irreducible, orientable $3$-orbifold, and let $V$ be a component of
$\partial|\oldPsi|-\inter|\oldPhi(\oldPsi)|$. Suppose that $\obd(V)$ is an annular orbifold. Then $V$ is a component of 
$\frakH(\oldPsi)\cap\partial|\oldPsi|$, and  $V$ is a weight-$0$ annulus.
%disjoint from $\fraks_\oldPsi$, \redmissingref{``weight-$0$ annulus'' is
%better; make corresponding changes in proof and apps.} 
Furthermore, exactly one component of
$\partial V$ is contained in a component $F$ of $|\oldPhi(\oldPsi)|$ such
that $\obd(F)$ is an annular orbifold. Finally, if $L$ denotes
the component of $|\oldSigma(\oldPsi)|$ containing $F$, then $\omega(L)$ 
is not a \pagelike\ \Ssuborbifold\ of $\oldPsi$ (see \ref{S-pair def}).
\EndProposition

\Proof
We may assume without loss of generality that $\oldPsi$ is connected, and is therefore strongly \simple\ and boundary-irreducible.

Set $N=|\oldPsi|$, $\oldSigma=\oldSigma(\oldPsi)$,
$\oldPhi=\oldPhi(\oldPsi)$ and
$\frakA=\frakA(\oldPsi)$. 
We have $\partial V\subset\partial|\oldPhi|=\partial|\frakA|$. Let $Z$
denote the component of $\overline{N-|\oldSigma|}$ containing $V$. It follows from Lemma \ref{oops lemma} that $Z$ is strongly \simple. By Proposition \ref{new characteristic} and Definition  \ref{oldSigma def}, 
%and
%\ref{tuesa day}, 
$\omega(Z)$ is an \Asuborbifold\ of $\oldPsi$, so that in particular the pair $(\omega(Z),\omega(Z\cap\partial N))$ is acylindrical. Let $V'\subset\partial Z$
denote the union of $V$ with all components of $|\frakA|$ that
share boundary components with $V$. Then $V'$ is connected since $V$
is. Furthermore, since $\omega(V)$ is annular, and each component of
$\frakA$ is annular by \ref{tuesa day}, we have
$\chi(\omega(V'))=0$. If we let $\oldUpsilon$ denote a strong
 regular neighborhood of $\omega(V')$ in $\omega(Z)$, and set
$V''=\Fr_Z|\oldUpsilon|$, then $\omega(V'')$  is
non-ambiently isotopic to $\omega(V')$, so that $V''$ is connected and $\chi(\omega(V''))=0$; thus $\omega(V'')$ is either annular or  toric. We claim:
\Claim\label{scott to be}
There is a \torifold\ $\frakK\subset\omega(Z)$ such that $\Fr_Z|\frakK|=V''$ and $\frakK\ne\oldUpsilon$.
\EndClaim

To prove \ref{scott to be}, first consider the case in which $\omega(V'')$ is annular. The orbifold $\omega(V')$ is homeomorphic to $\obd(V'')$ and is therefore annular. 
Since $V$ is annular we have $\partial V\ne\emptyset$, and hence $V'$ contains at least one component $A_0$ of $|\frakA|$. By \ref{tuesa day}, $\omega(A_0)$ is $\pi_1$-injective in $\oldPsi$, and in particular in $\omega(V')$. Since $\omega(A_0)$ and $\omega(V')$ are both annular, the image of the injective inclusion homomorphism $\pi_1(\omega(A_0))\to\pi_1(\omega(V' ))$ has index at most $2$ in $\pi_1(\omega(A_0))$; furthermore, if the image has index $2$ then 
$\pi_1(\omega(V'))$ is infinite cyclic and $\pi_1(\omega(A_0))$ is
infinite dihedral (see \ref{cobound}). If a homomorphism from the infinite dihedral group
$D_\infty$ to another group restricts to an injection on the index-$2$
cyclic subgroup of $D_\infty$, it must itself be injective. Since
$\omega(A_0)$ is $\pi_1$-injective in $\oldPsi$, it now follows that $\omega(V')$ is also $\pi_1$-injective in $\oldPsi$. 
Hence $\omega(V'')$ is $\pi_1$-injective in $\oldPsi$, and  in particular it is $\pi_1$-injective in $\omega(Z)$.

On the other hand, 
the definition of $V''$ gives $\partial V''\subset\partial N$, so that
$\partial(\omega(V''))\subset \omega(Z\cap\partial N)$; since
$(\omega(Z),\omega(Z\cap\partial N))$ is acylindrical, Definition \ref{acylindrical def} says that $\omega(V'')$ is not essential in the pair $(\omega(Z),\omega(Z\cap\partial N))$. Since $\omega(V'')$ is $\pi_1$-injective in $\omega(Z)$, and is not essential in 
$(\omega(Z),\omega(Z\cap\partial N))$, it follows from Definition \ref{acylindrical def} that $\omega(V'')$ is parallel 
in $(\omega(Z),\omega(Z\cap\partial N))$
to $\oldTheta$, where $\oldTheta$ is either (i) a suborbifold of $\omega(Z\cap\partial N)$ or (ii) a component of
$\omega(\overline{(\partial Z)\setminus(\partial N)}=\omega(\Fr_NZ)$. This in turn means that there is a suborbifold $\frakK$ of $\omega(Z)$ which can be equipped with a trivial $I$-fibration in such a way that $\omega(V'')$ and $\oldTheta$ are the components of $\partialh\frakK$, and $\partialv\frakK\subset \omega(Z\cap\partial N)$. In particular we have $\Fr_Z|\frakK|=V''$.
%, and $\frakK\cap\partial\oldPsi=\overline{(\partial\frakK)-\omega(Z)}$ is connected and meets $\omega(V'')$. \redmissingref{That last clause seems messed up, and I don't yet know what it's supposed to say. It's OK if I really know that $\frakK\cap\partial\oldPsi$ is connected and meets $\omega(V'')$, but the expression $\overline{(\partial\frakK)-\omega(Z)}$ is definitely wrong} 
Since the component $\omega(V'')$ of $\partialh\frakK$ is annular, 
$\frakK$ is a \torifold. 

To prove \ref{scott to be}, it remains to show that 
$\frakK\ne\oldUpsilon$.
If $\oldTheta$ satisfies (i) then 
$\overline{(\partial|\frakK|)-V''}$  is contained in $\inter(Z\cap\partial N)$, and is therefore disjoint from $\Fr_NZ$. By contrast, $\overline{(\partial|\oldUpsilon|)-V''}$ contains the component $A_0$ of $\Fr_NZ$, and therefore
$\frakK\ne\oldUpsilon$.

If $\oldTheta$ satisfies (ii), then the intersection of
$\overline{(\partial|\frakK|)-V''}$ with each component $X$ of $Z\cap\partial N$ is a strong regular neighborhood in $X$ of a common boundary component of $X$ and $|\oldTheta|$. In particular, $\overline{(\partial|\frakK|)-V''}$ contains no component of $Z\cap\partial N$. By contrast, $\overline{(\partial|\oldUpsilon|)-V''}$ contains the component $V$ of $Z\cap\partial N$,
 and therefore
$\frakK\ne\oldUpsilon$.
%we have ) or is a regular neighborhood in $\partial Z$ of a component of $\Fr_NZ=\overline{(\partial Z)\setminus(\partial N)}$ (and therefore does not contain .  Since 
%$\overline{(\partial|\oldUpsilon|)-V''}$ contains both $V\subset\partial N$ and $A_0\subset\Fr_NZ$, we must have $\frakK\ne\oldUpsilon$.
%has the component $V$, which is disjoint from $V''$. Since $|\frakK|\cap\partial N$ is connected and meets $\omega(V'')$, it follows that  $\frakK\ne\oldUpsilon$. 
This 
%shows that 
%$\frakK\ne\oldUpsilon$, and
 completes the proof of \ref{scott to be} in the case where $\omega(V'')$ is annular.

Now consider the case in which $\omega(V'')$ is toric. Note that $\omega(V'')$ is also two-sided, and hence orientable, since $V''=\Fr_Z|\oldUpsilon|$. Since 
the annular orbifold $\omega(V)\subset\omega(V')$
% which is impossible because 
is $\pi_1$-injective in $\oldPsi$, and since $\omega(V'')$ is  non-ambiently isotopic to $\omega(V')$,  the inclusion homomorphism $\pi_1(\omega(V''))\to\pi_1(\oldPsi)$ has infinite image. Since $\omega(Z)$ is strongly \simple, 
we may apply Corollary \ref{preoccucorollary}, with $\omega(Z)$ and $\obd(V'')$ playing the respective roles of $\oldPsi$ and $\oldTheta$, to deduce that 
%either (a) 
$V''$ bounds a $3$-dimensional submanifold $J$ of $Z$ such that $\omega(J)$ is a \torifold.
%, or (b) $V''$ is contained in the interior a $3$-dimensional submanifold $B$ of $Z$, transverse to $\fraks_\oldPsi$, such that $\omega(B)$ is a discal suborbifold. If (b) holds, the inclusion homomorphism $\pi_1(\omega(V''))\to\pi_1(\oldPsi)$ has finite image. ButHence (a) holds. 
In this case we define $\frakK$ to be the \torifold\ $\omega(J)$. We have $\partial|\frakK|=V''\subset\inter Z$, whereas $\partial|\oldUpsilon|$ contains the non-empty submanifold $V\subset\partial Z$; hence $\frakK\ne\oldUpsilon$ in this case as well, and the proof of \ref{scott to be} is complete.

If we fix a \torifold\ $\frakK$ with the properties given by \ref{scott to be}, then since  $\Fr_Z|\frakK|=V''=\Fr_Z|\oldUpsilon|$ and $\frakK\ne\oldUpsilon$, and since $Z$ is connected, we have $|\frakK|\cup|\oldUpsilon|=Z$ and $|\frakK|\cap|\oldUpsilon|=V''$. Since $\omega(V'')$ is the frontier of a strong regular neighborhood of $\omega(V')$ in $\omega(Z)$, it follows that $\omega(Z)$ is (orbifold-)homeomorphic to $\frakK$. Hence $\omega(Z)$ is a \torifold. But the components of $\Fr_\oldPsi\omega(Z)$ are components of $\frakA(\oldPsi)$, and are therefore essential annular suborbifolds of $\oldPsi$. In view of Lemma \ref{when a tore a fold}, it follows that:

\Claim\label{the real mcginty}
The orbifold $\omega(Z)$ is a \bindinglike\ \Ssuborbifold\ of $\oldPsi$.
\EndClaim

%\redmissingref{Decide whether to isolate this as a general principle earlier in the section.}

Since $Z$ is by definition a component of $\overline{N-|\oldSigma|}$, It follows from \ref{the real mcginty} and
% \redmissingref{this other fact that I've referred to several times---I think it should just say that if a \torifold\ sits in an orientable $3$-orbifold so that its frontier components are $\pi_1$-injective annular orbifolds, then it's an \Ssuborbifold. Actually the part about having a \bindinglike\ \spair has been proved, but being an \Ssuborbifold\ also involves essentiality, which seems to have been largely ignored in this mess of a proof} that $\omega(Z)$ is an \Ssuborbifold\ of $\oldPsi$. With 
Lemma \ref{no west for the bleary} that:

\Claim\label{blearier and blearier}
The \torifold\ $\omega(Z)$ is a component of $\frakH:=\frakH(\oldPsi)$, and
every component of $\overline{\oldPsi- \frakH} $ which shares a frontier component with $\omega(Z)$ is a component of $\oldSigma_1:=\oldSigma_1(\oldPsi)$.
\EndClaim

Since  $\omega(Z)$ is a component of $\frakH$ by \ref{blearier and blearier}, and $V$ is a component
of $Z\cap\partial N$, the surface $V$ is a component of
$|\frakH|\cap\partial N$. But since $\frakH$ is a strong
regular neighborhood of the two-sided $2$-suborbifold $\frakQ:=\frakQ(\oldPsi)$
of $\oldPsi$, each component of $|\frakH|\cap\partial N$ is a
weight-$0$ annulus. Thus
%Since $V$ is a component of $Z\cap\partial N$, it
%follows from that
 $V$ is a weight-$0$
annulus, and
%\redmissingref{At the moment I don't understand the logic of the
  %last couple of sentences, but I'm focused on adding the phrase
  %``weight-$0$ annulus.''} Thus
 the first sentence of the conclusion of the proposition is proved.

Now let $\oldSigma_Z$ denote the union of all components of $\oldSigma_1$ which share frontier components with $\omega(Z)$. By \ref{blearier and blearier}, we have $\Fr_\oldPsi\omega(Z)\subset\Fr_\oldPsi\oldSigma_Z$, and $\Fr_\oldPsi\omega(Z)$
% have the same frontier, and this frontier 
has two components. By definition each component of $\oldSigma_Z$ contains at least one component of $\Fr_\oldPsi\omega(Z)$, and hence $\oldSigma_Z$ has either one or two components. Since $\oldSigma_Z$ is a union of components of $\oldSigma_1$, each component of $\oldSigma_Z$ is an \Ssuborbifold\ of $\oldPsi$.
We claim:

\Claim\label{woozy}
The orbifold $\oldSigma_Z$ has exactly two components, of which exactly  one is a \bindinglike\ \Ssuborbifold\ of $\oldPsi$,
%$\SSS^1$-fibration under which it is standardly embedded in $\oldPsi$, 
and exactly one is a \pagelike\ \Ssuborbifold\ of $\oldPsi$ (see \ref{S-pair def}).
\EndClaim

To prove \ref{woozy}, first recall from \ref{oldSigma def} that every component of $\oldSigma$ is an \Ssuborbifold. Hence if \ref{woozy} is false, then either $\oldSigma_Z$ is connected, or it has two components which are either both \bindinglike\ \Ssuborbifold s\ of $\oldPsi$, or are  both \pagelike\ \Ssuborbifold s of $\oldPsi$. Since we have seen that $\oldSigma_Z$ has at most two components, it follows that 
%\redcomment{The next sentence is complete nonsense, and the one after it is not much better. In fact the statement of \ref{woozy} is also nonsense. This is one of the lemmas that need a huge proofreading according to an earlier comment.} In any event, 
either every component of $\oldSigma_Z$ is a \bindinglike\ \Ssuborbifold, or every component of $\oldSigma_Z$ is a \pagelike\ \Ssuborbifold.
Hence 
$\oldSigma_Z$ itself has a fibration $q:\oldSigma_Z\to\frakB$, over some $2$-orbifold $\frakB$, such that either $q$ is an $\SSS^1$-fibration and $\oldSigma_Z\cap\partial\oldPsi$ is saturated, or $q$ is an $I$-fibration and 
 $\oldSigma_Z\cap\partial\oldPsi=\partialh\oldLambda $. In either case, $\Fr_\oldPsi\oldSigma_Z$ is saturated in the fibration $q$, and in particular the union of components $\Fr\omega(Z)$ of $\Fr\oldSigma_Z$ is saturated.

Since $\frakH$ is a strong regular neighborhood of $\frakQ$ in $\oldPsi$, and $\omega(Z)$ is a component of $\frakH$ by \ref{blearier and blearier}, $\omega(Z)$ is a strong regular neighborhood of some component $\frakQ_Z$ of $\frakQ$. We may identify $\omega(Z)$ homeomorphically with $\frakQ_Z\times[0,1]$ in such a way that $\frakQ_i:=\frakQ_Z\times\{i\}$, for $i=0,1$, are the components of $\Fr_\oldPsi\omega(Z)$. For $i=0,1$, the fibration $q$ restricts to a fibration $q_i$ of $\frakQ_i$ whose base is the $1$-orbifold $\beta_i:=q(\frakQ_i)\subset\partial\frakB$. Up to isotopy, a given annular $2$-orbifold has at most one $\SSS^1$-fibration and at most one $I$-fibration. Hence the $1$-orbifolds $\beta_0$ and $\beta_1$ are homeomorphic, and moreover there is a fibration $r$ of $\omega(Z)$ over the $2$-orbifold $\frakE:=\beta_0\times[0,1]$ such that, under suitable homeomorphic identifications of $\beta_i$ with $\beta_0\times\{i\}\subset\frakE$ for $i=0,1$, 
we have $r|\frakQ_i=q_i$ for $i=0,1$. It follows that if $\frakB'$ denotes the  $2$-orbifold constructed from the disjoint union $\frakB\cup\frakE$ by identifying $\beta_i\subset\partial\frakB$ with $\beta_i=\beta_0\times\{i\}\subset\partial\frakE$ for $i=0,1$, then there is a fibration $q':\oldSigma_Z\cup\omega(Z)\to\frakB'$ which restricts to $q$ on $\oldSigma_Z$ and to $r$ on $\omega(Z)$. Hence $\oldSigma_Z\cup\omega(Z)$ is an \Ssuborbifold\ of $\oldPsi$, which is  \pagelike\ or \bindinglike\ according as $q$ is an $I$-fibration or an $\SSS^1$-fibration. This in turn implies that Condition (2) of Proposition \ref{new characteristic} continues to hold if $\frakQ$ is replaced by $\frakQ-\frakQ_Z$. But this is a contradiction, since $\frakQ$ satisfies Condition (3) of Proposition \ref{new characteristic}. This completes the proof of \ref{woozy}.

It follows from \ref{woozy} that $|\oldSigma_Z|$ has two distinct components, and that they may be labeled $L$ and $L'$ in such a way that $\omega(L)$ is a \bindinglike\ \Ssuborbifold\ of $\oldPsi$  but not a \pagelike\ \Ssuborbifold\ of $\oldPsi$, while
$\omega(L')$ is a \pagelike\ \Ssuborbifold\ of $\oldPsi$  but not a \bindinglike\ \Ssuborbifold\ of $\oldPsi$.
Since $\omega(L)$ and $\omega(L')$ are respectively a \bindinglike\ and a \pagelike\ suborbifold of $\oldPsi$, we may fix an $\SSS^1$-fibration and an $I$-fibration of $\omega(L)$ and $\omega(L')$ respectively such that $\omega(L)\cap\partial\oldPsi$ is saturated, and $\omega(L')\cap\partial\oldPsi=\partialh\omega(L')$.  Since $\omega(L')$ is not \bindinglike, it follows from Lemma \ref{when a tore a fold} that $\omega(L')$ is not a \torifold, and hence the base of the $I$-fibration of $\omega(L')$ is not annular. Since, by \ref{boundary is negative}, every component of $\partial\oldPsi$ has negative Euler characteristic, the base of the $I$-fibration of $\omega(L')$ is not toric; and by \ref{tuesa day}, we have $\chi(L)\le0$. Hence the base of the $I$-fibration of $\omega(L')$ has strictly negative Euler characteristic.
% $\chi(\omega(L))=0$ and $\chi(\omega(L'))<0$. According to \redmissingref{same cross-ref as above}, $\omega(L)$ has an $\SSS^1$-fibration, and $\omega(L')$ has an $I$-fibration over a $2$-orbifold of negative Euler characteristic, and under which they are standardly embedded in $\oldPsi$. Hence 
If $F$ and $F'$ respectively denote the components of $L\cap\partial
N$ and $L'\cap\partial N$ that meet $V$, the orbifolds $\omega(F)$ and
$\omega(F')$ are respectively a saturated suborbifold of $\omega(L)$
and a component of $\partialh\omega(L')$. It follows that
$\chi(\omega(F))=0$ and that $\chi(\omega(F'))<0$. Since $\obd(F)$ is not toric (again by virtue of \ref{boundary is negative}), it follows that $\omega(F)$
is annular, while $\omega(F')$ is not. This proves the second sentence of
the conclusion of the proposition. Since we chose $L$ in such a way
that $\omega(L)$ is not a \pagelike\ \Ssuborbifold\ $\oldPsi$, the
final sentence is proved as well.
\EndProof
%acylindrical

\Corollary\label{mangenfeffer}
Let $\oldPsi$ be a componentwise strongly \simple, componentwise boundary-irreducible, orientable $3$-orbifold. Let $\oldGamma$ be an annular
  component of $\oldPhi(\oldPsi)$ such that the component of $|\partial\oldPsi|$ containing $|\oldGamma|$ is a $2$-sphere. Then no component of
  $\partial\oldPsi-\inter\oldGamma$ can be annular.
%\redcomment{
%Now that I have  this statement, it occcurs to me that it could be applied in lots of other places as well. For example, doesn't it imply that in the complement of a belt for a sphere of weight $4$, each component has a singular point of order $>2$? Could this conceivably make some of the final results stronger?? At this point I really doubt it.}
\EndCorollary

\Proof
Let $S$ denote the component of $|\partial\oldPsi|$ containing $\oldGamma$. Suppose that some component $\frakE$ of $S-\inter\oldGamma$ is annular. Since $S$ is a $2$-sphere and $\oldGamma$ is connected, $E:=|\frakE|$ is a disk.
(Since $\frakE$ is annular,
$E\cap\fraks_\oldPsi$ must consist of two points, each of which has order $2$.) Set $C=\partial E$. Since $\oldGamma$ is a component of $\oldPhi(\oldPsi)$, there is a component $V$ of 
%_2
$\partial|\oldPsi|-\inter|\oldPhi|$
such that $C\subset\partial V$ and $V\subset E$. By
\ref{tuesa day},  $\obd(\partial V)$ is
$\pi_1$-injective in $\oldPsi$ and hence in $\obd(E)$. But since
$\obd(E)$ is annular and orientable, any suborbifold of $\obd(E)$ whose
frontier consists of $\pi_1$-injective simple closed curves in
$\inter E$ must itself be annular by \ref{cobound}.
Hence $\obd(V)$ is annular, and
it follows from Proposition \ref{lady edith} that $V$ is a weight-$0$
annulus. Let $C'$ denote the component of $\partial V$ distinct from
$C$. If $F$ denotes the component of $|\oldPhi|$ containing $C'$,
then $\obd(F)$ is again a suborbifold of $\obd(E)$ whose frontier
consists of $\pi_1$-injective simple closed curves in $\inter E$,
and must therefore be annular (see \ref{cobound}). But the component of $|\oldPhi|$
containing $C$ is $|\oldGamma|$, and $\oldGamma$ is also annular. Thus
both components of $\partial V$ are contained in underlying surfaces
of annular components of $\oldPhi $; this contradicts the uniqueness
assertion of Proposition \ref{lady edith}.
\EndProof

\Corollary\label{abu gnu}
Let $\oldPsi$ be a componentwise strongly \simple, componentwise boundary-irreducible, orientable $3$-orbifold. Then
$(\partial\oldPsi)\cap\book(\oldPsi)$ is the union of $\oldPhi(\oldPsi)$ with all annular components of
$\overline{(\partial\oldPsi)-\oldPhi(\oldPsi)}$.
\EndCorollary

\Proof
Let $\frakH_0$ denote the union of 
all components $\frakC$ of $\frakH(\oldPsi)$ such that $\frakC\cap\oldSigma(\oldPsi)=\Fr_\oldPsi\frakC$. 
By definition (see \ref{tuesa day}) we have $\book(\oldPsi)=\oldSigma(\oldPsi)\cup\frakH_0$. Hence it suffices to show that $\frakH_0\cap\partial\oldPsi=\frakE$, where $\frakE$ denotes the union of all annular components of
$\overline{(\partial\oldPsi)-\oldPhi(\oldPsi)}$.

A component $\frakC$ of $\frakH_0$ is in particular a component of $\frakH$, and therefore intersects $\partial\oldPsi$ in a suborbifold whose components are annular; since $\frakC\cap\oldSigma(\oldPsi)=\Fr_\oldPsi\frakC$, the annular orbifold $\frakC\cap\partial\oldPsi$ is a component of $\overline{(\partial\oldPsi)-\oldPhi(\oldPsi)}$. This shows that $\frakH_0\cap\partial\oldPsi\subset\frakE$. To prove the reverse inclusion, suppose that $\frakV$ is an annular component of $\overline{(\partial\oldPsi)-\oldPhi(\oldPsi)}$.
It follows from the first sentence of the conclusion of Proposition
\ref{lady edith} that $V:=|\frakV|$
%, that every component of $|\frakE|$ 
is a component
of $|\frakH(\oldPsi)|\cap\partial|\oldPsi|$. 
Since
$\Fr V\subset\oldPhi(\oldPsi)$, the component of $\frakH(\oldPsi)$ containing $V$ must have its frontier contained in $\oldSigma(\oldPsi)$, and is therefore a component of $\frakH_0$. This shows that $V\subset|\frakH_0|\cap\partial|\oldPsi|$, and the required inclusion 
$\frakE\subset
\frakH_0\cap\partial\oldPsi$ is thereby proved.
\EndProof
%\frakE\oldGamma

\Corollary\label{after abu gnu}
Let $\oldPsi$ be a componentwise strongly \simple, componentwise boundary-irreducible, orientable $3$-orbifold. Then every component of $\kish(\oldPsi)\cap\partial\oldPsi$ has strictly negative Euler characteristic.
\EndCorollary

\Proof
Let $\oldGamma$ be a component of  $\kish(\oldPsi)\cap\partial\oldPsi$. It follows from \ref{tuesa day} that $\chi(\oldGamma)\le0$. It follows from Corollary \ref{abu gnu} that $\oldGamma$ is not  annular. Since
$\oldPsi$ is  componentwise strongly \simple\ and componentwise boundary-irreducible, $\oldGamma$ cannot be  toric according to \ref{boundary is negative}.
\EndProof

\Corollary\label{less than nothing}
Let $\oldPsi$ be a $3$-orbifold
which is orientable, componentwise strongly \simple, and componentwise
boundary-irreducible. Then every component of $\kish\oldPsi$ either is closed or has strictly negative Euler characteristic.
\EndCorollary

\Proof
Let $\frakK$ be a component of $\kish\oldPsi$, and suppose that $\frakK$ is not closed, so that $\Fr\frakK\ne\empty$. Every component of $\Fr\frakK$ is a union of components of $\obd(\cala)$; as these components are annular, it follows that $\frakK\cap\partial\oldPsi\ne\emptyset$. Since each component of $\frakK\cap\partial\oldPsi$ has strictly negative Euler characteristic by Corollary \ref{less than nothing}, we have $\chi(\frakK\cap\partial\oldPsi)<0$. But it follows from \ref{tuesa day} that $\chi(\frakK)= \chi(\frakK\cap\partial\oldPsi)/2$.
\EndProof

\Lemma\label{why the role}
Let $\oldPsi$ be a componentwise  strongly \simple, componentwise boundary-irreducible, orientable $3$-orbifold, 
%set $N=|\oldPsi|$, 
and let $A$ be a component of $|\frakA(\oldPsi)|$.  Let $\oldLambda$ denote
the component of $\oldSigma(\oldPsi)$ such that $|\oldLambda|\supset A$.
Suppose that $\oldLambda$ is a \torifold, and that there is an annular suborbifold $\oldUpsilon$ of
$\partial\oldPsi$ such that $\partial |\oldUpsilon|=\partial A$. Then we have
$\Fr_{|\oldPsi|}|\oldLambda|=A$ and
$|\oldLambda|\cap\partial|\oldPsi|=|\oldUpsilon|$. Furthermore, the inclusion
homomorphism $\pi_1(\oldUpsilon)\to\pi_1(\oldLambda)$ is not surjective.
\EndLemma

\Proof
According to \ref{tuesa day}, $\obd(A)$ is an  annular
$2$-orbifold, essential in $\oldPsi$. Since $\oldUpsilon$ and $\obd(A)$ are  annular, we have $\chi(\oldUpsilon\cup\obd(A))=0$, so that
$\oldUpsilon\cup\obd(A)$ is toric.
Since 
 $\obd(A)$ is essential it is $\pi_1$-injective, and hence the inclusion homomorphism
$\pi_1(\oldUpsilon\cup\obd(A))\to\pi_1(\oldPsi)$ has infinite image. Furthermore, since $\obd(A)$ is two-sided in the orientable $3$-orbifold $\oldPsi$, and $\oldUpsilon\subset\partial\oldUpsilon$, the $2$-orbifold $\oldUpsilon\cup\obd(A)$ is orientable. It therefore follows from Corollary \ref{preoccucorollary} that
  $\oldUpsilon\cup\obd(A)$ bounds a \torifold\
  $\oldDelta\subset\oldPsi$. Thus $\Fr_N\oldDelta=\obd(A)$ and $\oldDelta\cap\partial\oldPsi=\oldUpsilon$.

Since  $\omega(A)$ is an essential annular suborbifold of $\oldPsi$, and
since $\oldDelta$ is a \torifold,   Lemma \ref{when a tore a fold} gives:

\Claim\label{now is the time}
 $\oldDelta$ is a \bindinglike\ \Ssuborbifold of $\oldPsi$.
\EndClaim

%The definition of $\oldLambda$ gives $A\subset\Fr\oldLambda$. Since $|\oldLambda|$ is connected \redmissingref{Is this important?} and $A=\Fr|\oldDelta|$, we have either (i) $|\oldLambda|\subset|\oldDelta|$ or (ii) $|\oldLambda|\cap|\oldDelta|=A$.

It follows from Definition \ref{oldSigma def} 
%and  Proposition \ref{new characteristic},  
that
$\oldSigma(\oldPsi)=\oldSigma_1\discup\oldSigma_2$, where $\oldSigma_i=\oldSigma_i(\oldPsi)$, and that every component of 
% every component of $
%We will denote by
$\oldSigma_1$ or $\oldSigma_2$ is respectively a 
component of $\overline{\oldPsi- \frakH} $ or
% which admit fibrations under which they are standardly embedded in $\oldPsi$, 
%and by $\oldSigma_2(\oldPsi)$ the union of all components of 
$\frakH$,
%(\oldPsi)$ that do not meet any component of $\oldSigma_1(\oldPsi)$.
where $\frakH= \frakH(\oldPsi)$.
Hence $\oldLambda$ is a component of $\oldSigma_1$ or $\oldSigma_2$, and the component $A$ of $\Fr_N|\oldLambda|$ is a component of $\Fr_N|\frakH|$. The component $\frakH_0$ of $\frakH$ containing $\obd(A)$ is by definition a strong regular neighborhood of a component $\frakQ_0$ of $\frakQ:=\frakQ(\oldPsi)$. Since $\obd(A)$ is a component of $\Fr_\oldPsi\frakH_0$, and $\Fr_N\oldDelta=\obd(A)$, we have:

\Claim\label{but not both}
Either
$\frakH_0\subset\oldDelta$ or $\frakH_0\cap\oldDelta=A$.
\EndClaim

Now we claim:

\Equation\label{nobody else}
(\frakH-\frakH_0)\cap\oldDelta=\emptyset.
\EndEquation

To prove (\ref{nobody else}), first note that since $\obd(A)=\Fr\oldDelta$ is a component of the frontier of the component $\frakH_0$ of $\frakH$, each component of $\frakH-\frakH_0$ is either disjoint from $\oldDelta$ or contained in $\inter\oldDelta$. Suppose that $\inter\oldDelta$ contains at least one component of $\frakH-\frakH_0$, and let $\frakH_1,\ldots\frakH_m$ denote the components of $\frakH-\frakH_0$ contained in $\inter\oldDelta$, where $m>0$. Then for $i=1,\ldots,m$ the orbifold $\frakH_i$ is a strong regular neighborhood of some component $\frakQ_i$ of $\frakQ$. Thus $\frakQ':=\frakQ-(\frakQ_1\cup\cdots\cup\frakQ_m)$ is the union of a proper subset of the components of $\frakQ$, and $\frakH':=\frakH-(\frakH_1\cup\cdots\cup\frakH_m)$ is a strong regular neighborhood of $\frakQ'$.

According to Definition \ref{oldSigma def}, $\frakQ$ satisfies Condition (2) of Proposition \ref{new characteristic}: each component 
of $\overline{\oldPsi-\frakH}$
%the pair
%$(\oldLambda,\oldLambda\cap\partial\oldPsi)$ 
is either an \Ssuborbifold \ or an \Asuborbifold\ of $\oldPsi$. Now consider an arbitrary component $\frakJ$ of $\overline{\oldPsi-\frakH'}$. It follows from \ref{but not both} and the definition of $\frakH'$ that either (a) $\frakJ=\oldDelta$, (b) $\frakH_0\subset\oldDelta$ and
 $\frakJ=\overline{\oldDelta-\frakH_0}$, or (c) $\frakJ$ is a component of $\overline{\oldPsi-\frakH}$ disjoint from $\inter\oldDelta$. If (a) or (b) holds then $\frakJ$ is isotopic to $\oldDelta$ and is therefore an \Ssuborbifold\ of $\oldPsi$ by \ref{now is the time}. If (c) holds then $\frakJ$ is either an \Ssuborbifold\ or an \Asuborbifold\ of $\oldPsi$ since $\frakQ$ satisfies Condition (2) of Proposition \ref{new characteristic}. This shows that Condition (2) of Proposition \ref{new characteristic} holds with $\frakQ'$ in place of $\frakQ$. But this is a contradiction, since $\frakQ$ satisfies $\frakQ$ satisfies Condition (3) of Proposition \ref{new characteristic} according to Definition \ref{oldSigma def}. This completes the proof of (\ref{nobody else}).

Next we claim:
\Claim\label{isadora}
In the notation of \ref{oldSigma def}, $\oldLambda$ is a component of $\oldSigma_1(\oldPsi)$.
\EndClaim

To prove \ref{isadora}, first note that $\oldLambda$ is by hypothesis a component of $\oldSigma(\oldPsi)$, and by \ref{oldSigma def} is therefore a component of either $\oldSigma_1(\oldPsi)$ or $\oldSigma_2(\oldPsi)$. If $\oldLambda$ is a component of  $\oldSigma_2(\oldPsi)$, then in particular it is a component of $\frakH$; since $\obd(A)\subset\oldLambda$, it follows that $\oldLambda=\frakH_0$. Now by \ref{but not both}, we have
either
$\frakH_0\subset\oldDelta$ or $\frakH_0\cap\oldDelta=A$. Let us define a suborbifold $\oldDelta^*$ of $\oldPsi$ by setting
$\oldDelta^*=\overline{\oldDelta-\frakH_0}$ if
$\frakH_0\subset\oldDelta$, and $\oldDelta^*=\oldDelta$ if $\frakH_0\cap\oldDelta=A$. It follows from \ref{nobody else} that $\inter\oldDelta^*$ is disjoint from $\frakH$; since the definition of $\oldDelta^*$ also implies that $\Fr_\oldPsi\oldDelta^*\subset\Fr_\oldPsi\frakH_0$, it follows that $\oldDelta^*$ is a component of $\overline{\oldPsi-\frakH}$. But $\oldDelta^*$ is isotopic in $\oldPsi$ to $\oldDelta$, and is therefore an \Ssuborbifold\ by \ref{now is the time}. Thus $\frakH_0$ shares a frontier component with a component of 
$\overline{\oldPsi-\frakH}$ which is an \Ssuborbifold\ of $\oldPsi$. By definition (see \ref{oldSigma def}) this means that $\frakH_0$ is not a component of $\oldSigma_2(\oldPsi)$. Thus the assumption that $\oldLambda$ is a component of $\oldSigma_2(\oldPsi)$ leads to a contradiction, and \ref{isadora} is proved.

%, or $\frakH=\overline{\frakJ^+-\frakH_0}$, where $\frakJ^+$ is a component of $\overline{\oldDelta\frakH

According to \ref{but not both}, we have either (i) $\frakH_0\cap\oldDelta=A$ or (ii)
$\frakH_0\subset\oldDelta$. We will complete the proof by showing that (i) implies the conclusion of the lemma, and that (ii)  leads to a contradiction.

 First suppose that (i) holds. In view of (\ref{nobody else}) it follows that $\frakH\cap\inter\oldDelta=\emptyset$. Since $\Fr_\oldPsi\oldDelta=\obd(A)\subset\Fr\frakH$, it follows that $\oldDelta$ is a component of $\overline{\oldPsi-\frakH}$. But it follows from \ref{isadora} and the definition of $\oldSigma_1(\oldPsi)$ (see \ref{oldSigma def}) that
$\oldLambda$ is also a component of $\overline{\oldPsi-\frakH}$. Since both $\oldLambda$ and $\oldDelta$ contain $\obd(A)$, we have $\oldLambda=\oldDelta$.
Hence $\Fr_N|\oldLambda|=\Fr_N|\oldDelta|=A$
and $|\oldLambda|\cap\partial|\oldPsi|=|\oldDelta|\cap\partial|\oldPsi|=|\oldUpsilon|$.
If the inclusion
homomorphism $\pi_1(\oldUpsilon)\to\pi_1(\oldLambda)$ 
is surjective, then it follows from  Proposition \ref{one for the books} 
applied with $\oldLambda$ and $\oldUpsilon$ playing the roles of $\oldUpsilon$ and $\frakZ$, 
that $(\oldLambda,\oldUpsilon)$ is homeomorphic to $(\oldUpsilon\times[0,1],\oldUpsilon\times\{0\})$.
(The hypotheses in Proposition \ref{one for the books} that the $3$-orbifold in question is very good and has no spherical boundary components hold here because $\oldLambda$ is a \torifold.) Hence  $\omega(A)$ is parallel in $\oldPsi$ to a suborbifold of $\partial\oldPsi$. Since $\omega(A)$ is essential by \ref{tuesa day}, this is a contradiction. This shows that the inclusion
homomorphism $\pi_1(\oldUpsilon)\to\pi_1(\oldLambda)$ is not surjective, and the conclusions of the lemma are established in the case where (i) holds.

Now suppose that (ii) holds. Set $\oldDelta^*=\overline{\oldDelta-\frakH_0}$. According to (\ref{nobody else}) we then have $\frakH\cap\inter\oldDelta^*=\emptyset$.  
Since $\Fr_\oldPsi\oldDelta^*$ is a component of $\Fr_\oldPsi\frakH_0$ (the component distinct from $\obd(A)$), it follows that $\oldDelta^*$ is a component of $\overline{\oldPsi-\frakH}$.  
Since we have seen that
$\oldLambda$ is also a component of $\overline{\oldPsi-\frakH}$ by \ref{isadora}, and since
$\oldLambda$ contains $A$ but $\oldDelta^*$ does not, we have $\oldLambda\cap\oldDelta^*=\emptyset$, and hence $\oldLambda\cap\oldDelta=\oldLambda\cap(\oldDelta^*\cup\frakH_0)=\oldLambda\cap\frakH_0=\obd(A)$.

By \ref{now is the time}, $\oldDelta$ is a \bindinglike\ \Ssuborbifold of $\oldPsi$. By hypothesis, $\oldLambda$ is a component of $\oldSigma(\oldPsi)$, so that its frontier components are essential annular suborbifolds of $\oldPsi$ (see \ref{tuesa day}). The hypothesis also gives that $\oldLambda$ is a\torifold, so that according to Lemma \ref{when a tore a fold}, $\oldLambda$ is also a \bindinglike\ \Ssuborbifold\ of $\oldPsi$. 
Thus there exist $\SSS^1$-fibrations $q:\oldDelta\to \frakB$ and $q':\oldLambda\to \frakB'$, where $\frakB$ and $\frakB'$ are $2$-orbifolds, with respect to which $\oldDelta\cap\partial\oldPsi$ and  $\oldLambda\cap\partial\oldPsi$ are saturated.

%These fibrations restrict to fibrations $q_A$ and $q_{A'}$ of $A$ and $A'$, whose bases are $1$-orbifolds $\beta$ and $\beta'$. We may homeomorphically identify $\frakH_0$ with $A\times[0,1]$ so that $A$ and $A'$ are identified with $A\times\{0\}$ and $A\times\{1\}$ respectively.
Up to isotopy, a given orientable annular $2$-orbifold has at most one $\SSS^1$-fibration. Hence 
%the $1$-orbifolds $\beta$ and $\beta$ are homeomorphic, and
% indeed can be homeomorphically identified with $\beta\times\{0\}$ and $\beta\times\{1\}$ in such a way that
  %moreover there is a fibration $r$ of $\frakH_0$ over the $2$-orbifold $\frakE:=\beta\times[0,1]$ such that, under suitable homeomorphic identifications of $\beta$ and $\beta'$ with $\beta\times\{0\}$ and $\beta\times\{1\}$, 
the fibrations $q$ and $q'$ may be chosen so that, under a suitable homeomorphic identification of the orbifolds $q(\obd(A))$ and $q'(\obd(A))$, we have $q|\obd(A)=q'|\obd(A)$. Hence if 
%we have $r|A=q$ and $r|A'=q'$. It follows that if 
$\frakB''$ denotes the  $2$-orbifold constructed from the disjoint union $\frakB\discup\frakB'$ by gluing $q(\obd(A))\subset\partial\frakB$ to $q(\obd(A)')\subset\partial\frakB'$ via the chosen identification, 
%and $\beta'\subset\partial\frakB'$ with $\beta\times\{1\}\subset\partial\frakE$, 
then there is a fibration of
$\oldDelta\cup\oldLambda$ over $\frakB''$
%$q':\oldSigma_X\cup\omega(X)\to\frakB'$ 
which restricts to $q$ and $q'$ on their respective domains.
% $\oldSigma_X$ and to $r$ on $\omega(X)$. 
With respect to this fibration we have
$(\oldDelta\cup\oldLambda)\cap\partial\oldPsi=\partialh(\oldDelta\cup\oldLambda)$. Thus
$\oldDelta\cup\oldLambda$ is a \pagelike\
\Ssuborbifold\ of $\oldPsi$. Since we have seen that $\oldDelta^*$ and $\oldLambda$ are components of $\overline{\oldPsi-\frakH^*}$, the \Ssuborbifold\ $\oldDelta\cup\oldLambda$ is a component of $\overline{\oldPsi-(\frakH-\frakH_0})$. The other components of $\overline{\oldPsi-(\frakH-\frakH_0})$ are components of $\overline{\oldPsi-\frakH}$, and are therefore \Ssuborbifold s and \Asuborbifold s. It now follows that Condition (2) of Proposition \ref{new characteristic} continues to hold if $\frakQ$ is replaced by $\frakQ-\frakQ_0$. But this is a contradiction, since $\frakQ$ satisfies Condition (3) of Proposition \ref{new characteristic}. Hence this case cannot occur, and the proof of the lemma is complete.
\EndProof
%\frakJ\overlineq_0comment

\Proposition\label{Comment A}If $p:\toldPsi\to \oldPsi$ is a covering map
of orientable $3$-orbifolds that are 
componentwise strongly \simple\ and componentwise boundary-irreducible, then $p^{-1}(\book(\oldPsi))$ is orbifold-isotopic to
$\book(\toldPsi)$. 
\EndProposition

\Proof
We may assume without loss of generality that $\oldPsi$ and $\toldPsi$
are connected (and are therefore strongly \simple\ and boundary-irreducible).
Set $\frakB=\book(\oldPsi)$. Then it follows from Lemma \ref{i see
  a book} that there exist a two-sided $2$-suborbifold $\oldPi$ of
  $\oldPsi$, and a strong
regular neighborhood $\oldGamma$  of $\oldPi$ in $\oldPsi$,
such that Conditions (a), (b) and (c) of that lemma hold. Since, by Condition (a),
 each component of $\oldPi$ is  an essential  annular suborbifold of
 $\oldPsi$, it follows from Corollary \ref{covering annular} that each component of the two-sided $2$-orbifold
 $\toldPi:=p^{-1}(\oldPi)$ is  annular and is essential in
 $\toldPsi$. 

Let
$\oldLambda$ be any component of $\overline{\oldPsi-\oldGamma}$, and
set $\toldLambda:=p^{-1}(\oldLambda)$. Then $\Fr_{\toldPsi}\toldLambda$ is a union
of components of $\toldPi$, which are essential  annular suborbifolds of $\toldPsi$. Since $\toldPsi$ is componentwise
strongly \simple, it follows from Lemma \ref{oops lemma} that
$\toldLambda$ is strongly \simple. On the other hand, by Condition (b), either $\oldLambda$ is an \Ssuborbifold\ of
%fibration under which it is standardly embedded in 
$\oldPsi$, which immediately implies that $\toldLambda$ is a \Ssuborbifold\ of
$\toldPsi$;
%so that $p^{-1}(\oldLambda)$ admits a
%fibration under which it is standardly embedded in $\toldPsi$;  
or $\oldLambda$ is an \Asuborbifold\ but not an \Ssuborbifold\ of
%admits no fibration under which it is
%standardly embedded in
 $\oldPsi$. In the latter case,
the pair
$(\oldLambda,\oldLambda\cap\partial\oldPsi)$ is acylindrical and is
not an \spair, and each component of $\overline{\partial\oldLambda\setminus\partial\oldPsi}=\Fr_\oldPsi{\oldLambda}
$ is annular. 
%admits no fibration under which it is
%standardly embedded in
% $\oldPsi$, in which 
%case
 In this case, in view of the strong \simple ity of $\toldLambda$, it follows from Proposition \ref{butthurt} that
$(\toldLambda,\toldLambda\cap\partial\toldPsi)$ is acylindrical and is not an \spair. Thus $\toldLambda$ is either an \Ssuborbifold\ or an \Asuborbifold\ of $\toldPsi$ 
for every component $\oldLambda$ of $\overline{\oldPsi-\oldGamma}$; and $\toldLambda$ is an \Ssuborbifold\ if and only if $\oldLambda$ is an \Ssuborbifold. 
Hence 
if we denote by $\toldGamma$ the strong regular neighborhood $p^{-1}(\oldGamma)$ of
$\toldPi$, then each component of $\overline{\toldPsi-\toldGamma}$ is
either an \Asuborbifold\ or an \Ssuborbifold\ of $\toldPsi$. Furthermore, since Condition (c) of Lemma \ref{i see
  a book} asserts that $\frakB$ is the union of $\oldGamma$ with all components of 
$\overline{\oldPsi-\oldGamma}$ that are \Ssuborbifold s of $\oldPsi$,
%admit
%fibrations under which they are standardly embedded in $\oldPsi$, 
it now follows that $\tfrakB:=p^{-1}(\frakB)$ is the union of $\toldGamma$ with all components of 
$\overline{\toldPsi-\toldGamma}$ that are \Ssuborbifold s of
%admit
%fibrations under which they are standardly embedded in
$\toldPsi$. Thus Conditions (a), (b) and (c) of Lemma \ref{i see a
  book} hold with $\toldPsi$, $\toldPi$, 
% of
  %$\oldPsi$, and a 
%regular neighborhood 
$\toldGamma$ and $\tfrakB$ in place of
$\oldPsi$, $\oldPi$, 
$\oldGamma$ and $\frakB$. It therefore follows from Lemma \ref{i see
  a book} that $\tfrakB$ is isotopic to $\book(\toldPsi)$.
\EndProof
%\obd

\Proposition\label{fortunately}
Let $\oldPsi$ be a componentwise strongly \simple, componentwise boundary-irreducible,
orientable $3$-orbifold. Then for every 
%component of
%$\fraks_\oldPsi$ is an arc or a simple closed curve. Then for every
integer $d>1$, the number of points of order $d$ in
$\fraks_{\overline{\partial\oldPsi-\oldPhi(\oldPsi)}}$ is divisible by $\lambda_\oldPsi$. In particular, $\card \fraks_{\overline{\partial\oldPsi-\oldPhi(\oldPsi)}}$ is divisible by $\lambda_\oldPsi$ (see \ref{lambda thing}).
\EndProposition

\Proof
The assertion is trivial if $\lambda_\oldPsi=1$. We shall therefore
assume that $\lambda_\oldPsi=2$; thus every component of $\fraks_\oldPsi$ is an arc or simple
closed curve, and therefore has a well-defined order greater than
$1$. 

Given an integer $d>1$, let $C$ denote the
union of all components of $\fraks_\oldPsi$ that have order
$d$. Set 
$X=\overline{|\oldPsi|-|\oldSigma(\oldPsi)|}$ and
$C'=C\cap X$, so that each component of $C'$ is an arc or a simple closed curve. Then $C\cap\partial X=\partial C'$, and hence $
%\nu:=
\card(C\cap \partial X)$ is even. We have $\partial X= |\frakA|\cup G$, where
$\frakA=\frakA(\oldPsi)$ and
$G=\overline{\partial\oldPsi-\oldPhi(\oldPsi)}$. Since $\frakA$
is an orientable orbifold, $|\frakA|\cap G=\partial |\frakA|$ is disjoint from
$\fraks_{\frakA}=\cals\cap\fraks_\oldPsi$,
and hence
%$\fraks_\oldPsi$, and hence 
from
 $C$; we therefore have 
$\card(C\cap \partial X)=\card(C\cap |\frakA|)+\card(C\cap G)$. Thus
$\card(C\cap |\frakA|)\equiv\card(C\cap G)\pmod2$. We are required to prove that
$\card(C\cap G)$ is divisible by $\lambda_\oldPsi=2$. Thus it suffices to show that
$\nu:=\card(C\cap |\frakA|)$, the number of points of order $d$ in
$\fraks_{\frakA}$ is even. According to \ref{tuesa day}, each component of $\frakA$ is an
orientable annular orbifold; hence each component of $|\frakA|$ is an
annulus containing no points of $\fraks_{\frakA}$, or a disk
containing two points of $\fraks_{\frakA}$, each of order $2$. It follows that $\nu=0$ if $d>2$, and that $\nu$ is twice the number of disk components of $|\frakA|$ if $d=2$. In either case, $\nu$ is even as required.
\EndProof

\Proposition\label{what does it say?}
Let $\oldPsi$ be a componentwise strongly \simple, componentwise boundary-irreducible, orientable $3$-orbifold, and let $L$ be a component of
$|\oldSigma(\oldPsi)|$. If for at least one component $F_0$ of
$L\cap\partial|\oldPsi|$ the $2$-orbifold $\omega(F_0)$ is annular, then
$\omega(L)$ is a \torifold, and for  every component $F$ of
$L\cap\partial|\oldPsi|$ the $2$-orbifold $\omega(F)$ is annular.
\EndProposition

\Proof
Set
$\oldLambda=\omega(L)$, and $\oldXi=\omega(L\cap\partial|\oldPsi|)$. 
According to \ref{oldSigma def}, the component $\oldLambda$ of
$\oldSigma(\oldPsi)$ is an \Ssuborbifold\ of $\oldPsi$.
%, i.e. $(\oldLambda,\oldXi)$ is an \spair. 
If $\oldLambda$ is a \bindinglike\ \Ssuborbifold,
%; thus either $\oldLambda$ is $I$-fibered over $\frakB$ and $\oldXi=\partialh\oldLambda$, or $\oldLambda$ is $\SSS^1$-fibered over $\oldLambda$ and $\oldXi$ is saturated in the fibration. It the latter case 
it follows from Lemma \ref{when a tore a fold} that $\omega(L)$ is a \torifold\ and that $\obd(F)$ is annular for  every component $F$ of
$L\cap\partial|\oldPsi|$, which is the conclusion in this case.

Now suppose that $\oldLambda$ is a \pagelike\ \Ssuborbifold; thus $\oldLambda$ is $I$-fibered over some $2$-orbifold $\frakB$, and $\oldXi=\partialh\oldLambda$.
Since the definition of $\oldSigma(\oldPsi)$ implies that each component of $\Fr\oldLambda$ is an essential annular suborbifold of $\oldPsi$, the suborbifold $\oldLambda$ is $\pi_1$-injective in $\oldPsi$; and since the strongly \simple\ orbifold $\oldPsi$ is definition very good by \ref{oops},  $\oldLambda$ is also very good. Let  $\toldLambda$ be a finite-sheeted manifold covering of $\oldLambda$.  Then $\toldLambda$ is equipped with an $I$-fibration over a surface $S$, so that 
every component of $\partialh\toldLambda$ is a covering of $S$ (of degree at
most $2$). Furthermore, $\partialh\toldLambda$ is a
manifold covering of $\partialh\oldLambda$. Since the component
$\obd(F_0)$ of $\partialh\oldLambda$ is annular, the component $X$ of
$\partialh\toldLambda$ covering $\omega(F_0)$ is an annulus. Since $X$
is also a covering of $S$, the surface $S$ is an annulus or a M\"obius
band. Hence $\toldLambda$ is a solid torus, and each component of its
horizontal boundary is an annulus. In view of the definitions of
\torifold\ and annular orbifold, the conclusion follows in this case as well.
\EndProof
%\obd

\Number\label{manifolds are different}We conclude this section by observing that the treatment of the characteristic suborbifold given here covers only the case in which the given orbifold is strongly \simple. In the proofs of Lemma \ref{graphology} and Proposition \ref{crust chastened}, where we consider the characteristic submanifolds of $3$-manifolds which, regarded as orbifolds, are not strongly \simple, we will refer directly to the ``classical'' sources
\cite{Jo} and\cite{js}.
\EndNumber

%_0

\section{ASTA volume for orbifolds}\label{darts section}

\Definition\label{hexe}
The {\it Gromov volume} of a closed, orientable topological $3$-manifold $M$ is
defined to be $\vtet\|M\|$, where $\vtet$ denotes the volume of a
regular ideal tetrahedron in $\HH^3$, and $\|M\|$ denotes the Gromov
norm \cite{gromov}, \cite[Chapter C]{bp}. We will denote the Gromov volume of $M$ by $\volG(M)$. 
\EndDefinition 

The relevance of $\volG(M)$ to the present monograph arises from Gromov's
result (see
\cite[Theorem C.4.2]{bp}) that for any closed
hyperbolic $3$-manifold $M$ we have $\volG M=\vol M$.

\begin{remarksdefinition}\label{sem ting}
According to \cite[Section 0.2]{gromov}, if $\tM$ is a $d$-sheeted cover
of a
compact, orientable topological $3$-manifold $M$, where $d$ is a positive integer,
we have 
\Equation\label{i'm silly}\volG(\tM)=d\volG(M).
\EndEquation

Now suppose that $\oldOmega$ is a closed, connected, very good topological orbifold. By definition, $\oldOmega$ admits a finite-sheeted covering which is a
manifold.

Suppose that $\toldOmega_1$ and $\toldOmega_2$ are finite-sheeted manifold
coverings of $\oldOmega$. Then $\toldOmega_1$ and $\toldOmega_2$ admit a common
covering $\toldOmega_3$. If $d_i$ denote the degree of the covering
$\toldOmega_i$ of $\oldOmega$, then for $i=1,2$ we have $d_i|d_3$, and
$\toldOmega_3$ is a $d_3/d_i$-fold covering of $\toldOmega_i$. By (\ref{i'm
  silly}) we have
$(d_3/d_1)\volG(\toldOmega_1)=\volG(\toldOmega_3)=(d_3/d_2)\volG(\toldOmega_2)$,
so that
\Equation\label{smoking fistule}
\volG(\toldOmega_1) /d_1=\volG(\toldOmega_2) /d_2.
\EndEquation
% there is a well-defined invariant $\volG(\oldOmega)$ given by
%setting $\volG(\oldOmega)=\volG(\toldOmega)/d$, where $\toldOmega$ is any
%finite-sheeted cover of $\oldOmega$ which is a manifold, and $d$ denotes
%the degree of the covering $\toldOmega$.

If $\oldOmega$ is a closed, connected, very good topological orbifold, then by
(\ref{smoking fistule})   there is a well-defined invariant $\volG(\oldOmega)$ given by
setting $\volG(\oldOmega)=\volG(\toldOmega)/d$, where $\toldOmega$ is any
finite-sheeted cover of $\oldOmega$ which is a manifold, and $d$ denotes
the degree of the covering $\toldOmega$.

More generally, if $\oldOmega$ is a possibly disconnected, very good, closed,
orientable $3$-orbifold, let us define
$$\volG(\oldOmega)=\sum_{C\in\calc(\oldOmega)}\volG(C).$$
\end{remarksdefinition}

\Notation\label{t-defs}
This subsection will contain definitions of a number of invariants which,
like the invariant $\volG$ discussed above, will turn out to be closely related to 
hyperbolic volume. It is convenient to collect the definitions of all  these
invariants here; their connection with hyperbolic volume will be
established in the course of the section. 
%Many of these invariants will be used in this paper, but there are a few whose importance will become apparent only in \cite{second}. 

%We

For any very good, compact, orientable (but possibly disconnected) topological $3$-orbifold $\oldPsi$, we will set
$$\volorb(\oldPsi)=\frac12\volG({\rm D}\oldPsi),$$
where ${\rm D}\oldPsi$ is defined as in \ref{doubling} (and is very good since $\oldPsi$ is very good and orientable). The acronym ASTA stands for Agol, Storm, Thurston and Atkinson; see \cite{ast} and \cite{atkinson}.
Note that if $\oldPsi$ is a disjoint union of suborbifolds $\oldPsi_1$ and $\oldPsi_2$, we have
$2\volorb(\oldPsi)=\volG({\rm D}\oldPsi)=\volG(\oldPsi_1)+\volG(\oldPsi_2)=2\volorb(\oldPsi_1)+2\volorb(\oldPsi_2)$, and hence $\volorb(\oldPsi)=\volorb(\oldPsi_1)+\volorb(\oldPsi_2)$. 

Note that, according to the convention posited in \ref{categorille}, the PL category will be the default category of orbifolds for the rest of this section.

For any very good, compact, orientable (but possibly disconnected) $3$-orbifold $\oldPsi$, we will set
$$\smock_0(\oldPsi)=\sup_\oldPi \volorb(\oldPsi\cut\oldPi),$$ 
where $\oldPi$ ranges
over all (possibly empty, possibly disconnected) incompressible $2$-suborbifolds of $\inter\oldPsi$. (Note that for any such $\oldPi$, the orbifold $\oldPsi\cut\oldPi$ is very good since $\oldPsi$ is.)
{\it A priori}, $\smock_0(\oldPsi)$ is an element of the extended real
number system. 
%\redcomment{Decide whether to retain this notation. The arguments in this section involving it seem circuitous, %but...}

Note that the definition of $\smock_0(\oldPsi)$ implies that
$\smock_0(\oldPsi)\ge \volorb(\oldPsi\cut\emptyset)$, i.e.
\Equation\label{more kitsch}
\smock_0(\oldPsi)\ge \volorb(\oldPsi).
\EndEquation

We will also set
%$$\smock(\oldPsi)= 2\bigg\lfloor\frac{\smock_0(\oldPsi)}{0.61} \bigg\rfloor
%\qquad
%\text{and}\qquad
$$\smock(\oldPsi)= \frac{\smock_0(\oldPsi)}{0.305}.
%\smock'(\oldPsi)= \bigg\lfloor\frac{12}{\voct}\smock_0(\oldPsi)\bigg\rfloor.
$$ 

If $\oldPsi_1,\ldots,\oldPsi_m$ are the components of $\oldPsi$, then every  incompressible
$2$-suborbifold of $\inter\oldPsi$ has the form $\oldPi=\oldPi_1\cup\cdots\cup\oldPi_m$, where $\oldPi_i$ is an  incompressible $2$-suborbifold of $\inter\oldPsi$; furthermore, we have $\volorb(\oldPsi\cut\oldPi)=\volorb((\oldPsi_1)\cut{\oldPi_1}) +\cdots+\volorb((\oldPsi_m)\cut{\oldPi_m})$. It follows that $\smock_0$ is additive over components, and hence so is $\smock$.

If $\oldPsi$ is a compact, orientable $3$-orbifold such that every boundary component of $N:=|\oldPsi|$ is a sphere, so that
$\plusN$ is closed, we will set
$$
%\theta(\oldPsi)=2\bigg\lfloor\frac{\volG(\plusN)}{0.61}\bigg\rfloor
%\qquad
%\text{and}\qquad
\theta(\oldPsi)=\frac{\volG(\plusN)}{0.305}.
$$

As a hint about why the curious-looking number 
  $0.305$ appears in these
definitions, we mention that the inequality $\voct/12>0.305$, where, as in \ref{voct def}, $\voct$ denotes the volume of a regular ideal octahedron in $\HH^3$,
% will be denoted
%. We have $\voct=3.6638\ldots$.
%ich is
is used in the proof of Corollary \ref{bloody hell} below (and is close to being an equality since $\voct=3.6638\ldots$), while
Theorem 1.1 of \cite{rankfour}, which implies that
  any closed, orientable hyperbolic $3$-manifold $M$ of volume at most
  $1.22=4\cdot0.305$ satisfies $h(M)\le3$, is used in the proof of Proposition \ref{new get
    lost}. (See \ref{in and out 
    soda} above for the definition of $h(M)$.) 

If a compact, orientable $3$-orbifold $\oldPsi$ is componentwise strongly \simple\ and boundary-irreducible,
we will set
$$\sigma(\oldPsi)=
12\chibar(\kish(\oldPsi))
,$$
where $\kish(\oldPsi)$ is defined by \ref{tuesa day}. 
(The appearance of the number $12$ in this definition is for convenience; it guarantees that certain lower bounds for $\sigma(\oldPsi)$ that appear in the course of the arguments in this monograph are integers rather than fractions.)

We will  set 
$$\sigma'(\oldPsi)=\lambda_\oldPsi\bigg\lfloor\frac{\sigma(\oldPsi)}{\lambda_\oldPsi}\bigg\rfloor.$$ 
%=
%2
%\bigg
%\lfloor6\chibar(\kish(\oldPsi))
%\bigg
%\rfloor\qquad
%\text{and}\qquad
%\sigma'(\oldPsi)=
%\bigg
%\lfloor
%12\chibar(\kish(\oldPsi))
%\bigg
%\rfloor,
%$$\sigma'(\oldPsi)=\sup_\oldPi\bigg\lfloor12\chibar(\kish(\oldPsi\cut\oldPi))\bigg\rfloor.$$ 
%\text
%{ if } \lambda_\oldPsi=2.
%$$
Note that we then have $\lambda_\oldPsi|\sigma'(\oldPsi)$, and $\sigma'(\oldPsi)\le\sigma(\oldPsi)$.

We will also set
$$\delta(\oldPsi)=\sup_\oldPi\sigma(\oldPsi\cut\oldPi),
%\qquad\text{and}\qquad
%\delta'(\oldPsi)=\sup_\oldPi\sigma'(\oldPsi\cut\oldPi),$$
%and
%$$\sigma'(\oldPsi)=\sup_\oldPi\bigg\lfloor12\chibar(\kish(\oldPsi\cut\oldPi))\bigg\rfloor,
$$ 
where $\oldPi$ ranges
over all (possibly empty and possibly disconnected) closed,  incompressible $2$-suborbifolds of $\inter\oldPsi$. In view of Lemma \ref{oops lemma}, the incompressibility of  $\oldPi$ guarantees that $\oldPsi\cut\oldPi$
is itself componentwise strongly \simple\ and boundary-irreducible, so that
$\sigma(\oldPsi\cut\oldPi)$ 
%and $\sigma'(\oldPsi\cut\oldPi)$ are
is defined.
Like $\smock_0(\oldPsi)$, the invariant $\delta(\oldPsi)$ is {\it a priori} an element of the extended real
number system.

Note that the definition of $\delta(\oldPsi)$ implies that 
$\delta(\oldPsi)\ge\sigma(\oldPsi\cut\emptyset)$, i.e.
%and 
%$\delta'(\oldPsi)\ge\sigma'(\oldPsi\cut\emptyset)$, i.e.
\Equation\label{kitsch}
\delta(\oldPsi)\ge\sigma(\oldPsi).
%\qquad\text{and}\qquad
%\delta'(\oldPsi)\ge\sigma'(\oldPsi).
\EndEquation

Now suppose that $\oldOmega$ is a closed, orientable, strongly \simple\
$3$-orbifold, let $\cals $ be a (possibly disconnected) $2$-dimensional submanifold of $M:=|\oldOmega|$, 
in general position with respect
to
$\fraks_\oldOmega$, and let $X$ be a union of components of
$M- \cals $. We will set
$t_\oldOmega(X)=\smock(\obd(\hatX))$.
%\qquad\text{and}\qquad t'_\oldOmega(X)=\smock'(\obd(\hatX))
%.$$
(See \ref{nbhd stuff} for the definition of $\hatX$. Note that since $\oldOmega$ is strongly \simple, and therefore very good by \ref{oops}, $\obd(\hatX)$ is also very good, and hence $\smock(\obd(\hatX))$ is defined.)
If  $\obd(\cals)$ is  incompressible, then since $\oldOmega$ is strongly \simple\ (see \ref{oops}), it follows from Lemma \ref{oops lemma} that the
components of $\obd(\hatX)$ are strongly \simple\ and boundary-irreducible. In this case we will set 
$s_\oldOmega(X)=\sigma(\obd(\hatX))$, $s'_\oldOmega(X)=\sigma'(\obd(\hatX))$,
and 
$y_\oldOmega(X)=\delta(\obd(\hatX))$. Since we have observed that $\sigma'(\oldPsi)\le\sigma(\oldPsi)$ for any compact, orientable $3$-orbifold $\oldPsi$ which is componentwise strongly \simple\ and boundary-irreducible,
we have $s'_\oldOmega(X)\le s_\oldOmega(X)$.
If every component of $\partial \hatX$ is a sphere, we will set
$q_\oldOmega(X)=\theta(\obd(\hatX))$.
\EndNotation

%The following lemma, which has a simple proof and fits naturally into
%the present discussion, is not quoted in the present paper but will be
%used in \cite{second}.

\Lemma\label{frobisher}
Let  $\oldOmega$ be any very good, compact, orientable (possibly disconnected) $3$-orbifold, and let $\oldTheta$ be a (possibly disconnected) incompressible, closed  $2$-suborbifold of $\inter\oldOmega$. Then 
$\smock_0(\oldOmega)\ge\smock_0(\oldOmega\cut\oldTheta)$,
and hence
$\smock(\oldOmega)\ge\smock(\oldOmega\cut\oldTheta)$.
\EndLemma

\Proof
Let $\oldPi$ be any (possibly empty, possibly disconnected) incompressible closed $2$-suborbifold of $\inter(\oldOmega\cut\oldTheta)$. Then $\oldTheta\cup\oldPi$ is identified with an incompressible $2$-suborbifold of $\inter\oldOmega$. Hence the definition of
$\smock_0(\oldOmega)$ implies that $\volorb(\oldOmega\cut{\oldTheta\cup\oldPi})\le \smock_0(\oldOmega)$. Since $\oldOmega\cut{\oldTheta\cup\oldPi}$ is homeomorphic to $(\oldOmega\cut\oldTheta)\cut\oldPi$, it follows that 
$\volorb((\oldOmega\cut\oldTheta)\cut\oldPi)\le \smock_0(\oldOmega)$. Since the latter inequality holds for every
incompressible closed $2$-suborbifold $\oldPi$ of $\inter(\oldOmega\cut\oldTheta)$, the definition of $\smock_0(\oldOmega\cut\oldTheta)$ gives $\smock_0(\oldOmega\cut\oldTheta) \le\smock_0(\oldOmega)$. The inequality $\smock(\oldOmega\cut\oldTheta) \le\smock(\oldOmega)$ then follows at once in view of the definition of the invariant $\zeta$.
\EndProof

\Theorem\label{darts theorem}
Let $\Mh$ be a closed, orientable, hyperbolic $3$-orbifold. Set $\oldOmega=(\Mh)\pl$, and let
$\oldTheta$ be a  (possibly disconnected) incompressible  $2$-suborbifold of
$\oldOmega$. Then
$$\vol(\Mh)\ge\volorb(\oldOmega\cut\oldTheta).$$
\EndTheorem

\Proof
We first give the proof under the following additional assumptions:
\begin{enumerate}
\item $\Mh$ is a closed, orientable, hyperbolic $3$-manifold (so that $\oldTheta$ is an incompressible $2$-submanifold of
$\oldOmega$);
\item no component of $\oldTheta$ is the boundary of a twisted $I$-bundle (over a closed, non-orientable surface) contained in the manifold $\oldOmega$; and
\item no two components of $\oldTheta$ are parallel in $\oldOmega$.
\end{enumerate}
For this part of the proof, to emphasize that $\Mh$ is assumed to be a manifold, I will denote it by $M$. Because $\volorb$ is a purely topological invariant, and because of the equivalence between the PL and smooth categories for $3$-manifolds,  
%we may fix a smooth structure on $M$ compatible with its PL structure, and 
it suffices to prove that if $T$ is a smooth submanifold of the hyperbolic manifold $M$ such that Conditions (1)---(3) hold with $M$ and $T$ in place of $\oldOmega$ and $\oldTheta$, then 
$\vol(M)\ge\volorb(M\cut T)$.  The proof of this will use the observations that Propositions \ref{stronger waldhausen} and \ref{snuff} above, and \cite[Corollary 5.5]{Waldhausen}, which are theorems in the PL category of $3$-manifolds, remain true without change in the smooth category; this also follows from the equivalence of categories mentioned above.
%We begin with the subcase in which no two components of $T$ are parallel. 
For each component $V$ of $T$, it follows from
  \cite[Theorem 3.1]{schoen-yau} that the inclusion map $i_V:V\to M$ is
  homotopic to a immersion  $j_V$ which has least area in its homotopy class. By  \cite[Theorem
  5.1]{FHS}, for each component $V$ of $T$, either $j_V$ is an
  embedding,  or there is  a one-sided, closed, connected surface $F_V\subset M$ such
  that $j_V$ is a two-sheeted covering map from $V$ onto $F_V$. If the
  latter alternative holds, then $j_V$, and hence $i_V$, is homotopic
  to a homeomorphism $j'_V$ of $V$ onto the boundary of a tubular
  neighborhood $N_V$ of $F_V$, which is a twisted $I$-bundle. Since
  $i_V$ and $j_V'$ are homotopic embeddings of $V$ into $M$, it
  follows from the smooth version of \cite[Corollary 5.5]{Waldhausen} that they are
  isotopic; this implies that $V$ is the boundary of a twisted
  $I$-bundle in $M$, a contradiction to Assumption (2). Hence for
  every component $V$ of $T$, the immersion $j_V$ is an embedding.

Now define an immersion $j:T\to M$ by setting $j|V=j_V$ for each component $V$ of $T$.  We claim:
\Claim\label{two by two}
The immersion $j$ is an embedding.
\EndClaim
We have already shown that $j_V$ is an embedding for each component
$V$ of $T$. Hence, to prove \ref{two by two}, it suffices to prove
that for any distinct components $V$ and $W$ of $T$ we have  $j_V(V)\cap j_W(W)=\emptyset$. 

%Set $V^*=j_V(V)$ and $W^*=j_W(W)$, and assume that $V^*\cap
%W^*\ne\emptyset$. Suppose for the moment that $V^*$ and $W^*$ intersect transversally. 
Since  the inclusions $i_V$ and $i_W$ are $\pi_1$-injective and have disjoint images, 
%and sinceare homotopic in
%$M$ and $\pi_1$-injective, 
it follows from the smooth version of Proposition \ref{snuff}
that for any embeddings $f:V\to M$ and $g:W\to M$ homotopic to $i_V$ and $i_W$ respectively, such that $f(V)$ and $g(W)$ meet transversally, either (a) $f(V)\cap g(W)=\emptyset$, or (b) there
exist connected subsurfaces $A\subset f(V)$ and $B\subset g( W)$, and a compact submanifold
$X$ of $M$, 
%and a homeomorphism $\xi :A\times[0,1]\to X$, 
such that
$\partial A\ne\emptyset$, $\partial X=A\cup B$, and the pair $(X,A)$ is homeomorphic to
$ (A\times[0,1],A\times\{0\})$.
%, and $\xi ((A\times\{1\})\cup((\partial
%A)\times[0,1]))=B$. 
Since  $j_V$ and $j_W$ are respectively homotopic to $i_V$ and $i_W$, and have least area within their homotopy classes, this shows that the hypothesis of Proposition \ref{fhs-prop} holds with $f_0=j_V$ and $g=j_W$. Thus Proposition \ref{fhs-prop} asserts that $j_V(V)\cap j_W(W)=\emptyset$, and \ref{two by two} is proved.
%''}$\pi_1$-injective and have disjoint images, 
%and sinceare homotopic in

Now since $j_V$ is homotopic to $i_V$ for each component $V$ of $T$,
the map $j$ is homotopic to $i$. Since $i$ is an embedding, and $j$ is
an embedding by \ref{two by two}, and since Assumption (3) holds, it
follows from the smooth version of Proposition \ref{stronger waldhausen} that $i$ and $j$
are isotopic. In particular, the (possibly disconnected) manifold
$M\cut T$ is homeomorphic to $X:=M\cut{j(T)}$. 

Since each of the embeddings $j_V$ is a least area immersion, the
surface $j(T)\subset M$ is a minimal surface. The manifold $X$ inherits a metric from the hyperbolic manifold
$M$, and with this metric, $X$ is bounded by a minimal
surface. According to Theorem 7.2 of \cite{ast}, we have
$\vol(X)\ge\volorb(X)$. 
%did not expect to see a volume on one side and a $\volorb$ on the
%other. The way this subcase is quoted below seems incorrect, because
%I have $\volorb$ on both sides below. OK, that part is OK, but the
%real problem is in the induction for the case where $T$ is
%disconnected. The problem is that the statement, $M$ is closed, but
%I'm applying it for the case where $M$ need not be closed in the
%course of the induction. Oops. Actually the induction doesn't even
%make sense, given that the inequality has $\vol$ on the left and not
%$\volorb$.} 
Since it is clear that $\vol(M)=\vol(X)$, and since we have seen that 
$M\cut T$ is homeomorphic to $X$, we have
$\vol(M) \ge\volorb(M\cut T)$, as required.
Thus the proof is complete under Assumptions (1)--(3).

%\cite

Next, we prove the theorem under the assumptions that (1) and (2)
hold, but (3) may not. (This argument, and the rest of the proof, will be done in the PL category.) Since (1) is still assumed to hold, we continue
to write $M$ in place of $\oldOmega$, and we write $T$ in place of $\oldTheta$ (so that $T$ now denotes a PL $2$-submanifold of $M$). Under Assumptions (1) and (2), the proof will
proceed by induction on $\compnum(T)$. The assertion is trivial if
$T=\emptyset$. Now suppose that $n\ge0$ is given and that the
assertion is true in the case where $\compnum(T)=n$. Suppose that $M$
and $T$ satisfy the hypotheses of the theorem and Assumptions (1) and
(2), and that $\compnum(T)=n+1$. If (3) also holds, the conclusion has
already been established. 

Now suppose that (3) does not hold,
i.e. that $T$ has two parallel components. Then there exists a
submanifold of $M$ which is homeomorphic to the product of a closed, connected
surface with $[0,1]$, and whose boundary components are components of
$T$. Among all such submanifolds of $T$ choose one, say $Y$, which is
minimal with respect to inclusion. If $\inter Y$ contains a component
$W$ of $T$, then since $W$ is incompressible in $M$ and hence in $Y$,
it follows from \cite[Proposition 3.1]{Waldhausen} that $W$ is parallel in $Y$
to each of the components of $\partial Y$; this contradicts the
minimality of $Y$. Hence $\inter Y$ contains no component of $T$, and
therefore $Y$ is the closure of a component of $M-T$.

Fix a component $V$ of $\partial Y$.
Set $T'=T-V$. Then the abstract disjoint union $Y\discup M\cut {T'}$ is homeomorphic to
$M\cut{T}$. 
By \ref{t-defs} it follows that 
%, and hence that
%${\rm D}(M\cut {T})$ is homeomorphic to the disjoint
%union ${\rm D}Y\discup {\rm D}(M\cut{T'})$, so that
$\volorb(M\cut {T})=\volorb(Y)+\volorb(M\cut{T'})$. Since
$Y$ is homeomorphic to $V\times[0,1]$, the manifold ${\rm D}Y$ is
homeomorphic to $V\times \SSS^1$. It then follows from \cite[Proposition C.3.4]{bp}) that $\volG({\rm D}Y)=0$, so that $\volorb(Y)=\volG({\rm D}Y)/2=0$. Thus we have 
%$2\volorb(M\cut {T})=\volG({\rm D}(M\cut{T'}))$, which implies
%that
 $\volorb(M\cut{{T'}})=\volorb(M\cut{{T}})$. Since $\compnum(T)=n$,
 the induction hypothesis gives $\vol
 M\ge\volorb(M\cut{T'})=\volorb(M\cut{T})$, and the induction is
 complete. This establishes the result under Assumptions (1) and (2).

We now turn to the proof of the theorem in the general case. Since $\Mh$ is hyperbolic, it is very good. Fix a finite-sheeted orbifold covering $\ph:\Nh\to\Mh$ where $\Nh$ is a hyperbolic $3$-manifold. Since the PL orbifold $\oldOmega$ has the same underlying topological orbifold as the hyperbolic orbifold $\Mh$, we may pull back the PL structure via $\ph$ to obtain a PL structure on $\Nh$. The manifold $\Nh$, equipped with this PL structure, will be denoted $N$, and the map $\ph$ may then be regarded as a PL covering map from $N$ to $\oldOmega$.

Set $\calq=p^{-1}(\oldTheta)$. Note that since $\oldTheta$ is incompressible in $\oldOmega$, the $2$-manifold $\calq$ is incompressible in $N$.  We claim: 
\Claim\label{even before that}
There is a two-sheeted covering $q:N'\to N$ such that $N'-q^{-1}(\calq)$ has no component whose closure  is a twisted $I$-bundle.
\EndClaim

To prove \ref{even before that}, let $\calx$ denote the set of all twisted $I$-bundles that are  closures of components of $N-\calq$. First consider the case in which the elements of $\calx$ are pairwise disjoint. In this case, let $Y$ denote the union of all elements of $\calx$, and set $Z=\overline{N-Y}$. Note that $Z$ is connected, since $N$ is connected and each element of $\calx$ has connected boundary. Since each component of $Y$ is a twisted $I$-bundle, there exist a compact (possibly disconnected) $2$-manifold $B$, and a $2$-sheeted covering map $q_Y:B\times[0,1]\to Y$, such that $q_Y$ maps $B\times\{i\}$ homeomorphically onto $\partial Y$ for $i=0,1$. Since $\partial Y=\partial Z$, we may regard $q_Y|(B\times\{i\})$ as a homeomorphism $r_i:B\times\{i\}\to\partial Z$. Let $Z_0$ and $Z_1$ be homeomorphic copies of $Z$, equipped with homeomorphisms $h_i:Z\to Z_i$. Then we may define a homeomorphism $\alpha$ from $\partial(B\times[0,1])=B\times\{0,1\}$ to $\partial(Z_1\discup Z_2)$ by setting $\alpha|(B\times\{i\})=h_i\circ r_i$ for $i=0,1$. In this case we define $N'$ to be the closed $3$-manifold obtained from the disjoint union $(B\times[0,1])\discup Z_1\discup Z_2$ by gluing $\partial(B\times[0,1])$ to $\partial(Z_1\discup Z_2)$ via the homeomorphism $\alpha$. Now we have a well-defined two-sheeted covering map $q:N'\to N$ given by setting $q|(B\times[0,1])=q_Y$ and $q|Z_i=h_i^{-1}$ for $i=0,1$. To prove the property of this covering stated in \ref{even before that}, note that the closure of each component of $N'-q^{-1}(\calq)$ either is equal to $Z$ or is a component of $B\times[0,1]$. A component of the latter type have disconnected boundary and is therefore not a twisted $I$-bundle. That $Z$ is not a twisted $I$-bundle follows from the definitions of $\calx$ and $Z$.

To complete the proof of \ref{even before that}, it remains to consider the case in which there are two elements of $\calx$ have non-empty intersection. If two such elements are denoted $\overline{X_1}$ and $\overline{X_2}$, where $X_1$ and $X_2$ are components of $N-\calq$, then since each of the twisted $I$-bundles $\overline{X_i}$ and $\overline{X_2}$ has connected boundary, and since $N$ is connected, we have $N=\overline{X_1}\cup\overline{X_2}$, and $\calx=\{\overline{X_1},\overline{X_2}\}$. We have $\overline{X_1}\cap\overline{X_2}=\partial \overline{X_1}=\partial\overline{ X_2}=\calq$.
For $j=1,2$, fix a
compact, connected $2$-manifold $B_j$, and a $2$-sheeted covering map $q_j:B_j\times[0,1]\to \overline{X_j}$, such that $q_j$ restricts to a homeomorphism 
$u_{ij}$ of $B_j\times\{i\}$ onto $\calq=\partial\overline{X_j} $ for $j=1,2$ and for $i=0,1$.
In this case we define $N'$ to be the closed $3$-manifold obtained from the disjoint union $(B_1\times[0,1])\discup( B_2\times[0,1])$ by gluing $B_1\times\{i\}$ to $B_2\times\{i\}$ via the homeomorphism 
$u_{i2}^{-1}\circ u_{i1}$, for $i=0$ and for $i=1$. 
%$\partial(Z_1\discup Z_2)$ via the homeomorphism $\alpha$. 
Now we have a well-defined two-sheeted covering map $q:N'\to N$ given by setting $q|(B_j\times[0,1])=q_j$ for $j=1,2$.
To prove the property of this covering stated in \ref{even before that}, note that the closures of the components of $N'-q^{-1}(\calq)$ are $B_1\times[0,1]$ and $B_2\times[0,1]$; these have disconnected boundary, and hence neither of them is a twisted $I$-bundle. 
Thus \ref{even before that} is proved.

Let $q:N'\to N$ be a two-sheeted cover with the properties stated in \ref{even before that}, and set $\calq'=q^{-1}(\calq)$. We claim

\Claim\label{new before general}
No component of $\calq'$ is the boundary of a twisted $I$-bundle in $N'$.
\EndClaim

To prove \ref{new before general}, suppose that some component  of $\calq'$ bounds a twisted $I$-bundle $R\subset N'$. We may assume that
$R$ is
minimal with respect to inclusion among all among all twisted $I$-bundles that are bounded by  components of $\calq'$. 
According to the property of $N'$ stated in \ref{even before that},
$R$ cannot be the closure of a component of $N'-\calq'$. Hence 
$\inter R$ contains some component $Q'$ of $\calq'$. Since $H_2(R,\ZZ)=0$, there is a compact submanifold $R_1$ of $\inter R$ with $\partial R_1=Q'$. Now fix a two-sheeted covering $t:\tR\to R$, where $\tR$ is a trivial $I$-bundle over a surface. 
Since 
$\calq$ is incompressible, 
$t^{-1}(\calq')=(pt)^{-1}(\calq)$ is also incompressible. It therefore follows from \cite[Prop. 3.1]{Waldhausen}, that each component of $\partial (t^{-1}(R_1))=t^{-1}(Q')\subset t^{-1}(\calq')$ is parallel to the boundary components of $R$. A second application of  \cite[Prop. 3.1]{Waldhausen} then shows that the two-sheeted covering $t^{-1}(R_1)$ of $R_1$ is a trivial $I$-bundle over a surface. Hence by  \cite[Prop. 3.1]{Waldhausen}, $R_1$ is a twisted $I$-bundle over a surface. But this contradicts the minimality of $R$,
and so
\ref{new before general} is proved.

Now since $N$ is a manifold, $N'$ is also a manifold. Since $N$ has the same underlying topological orbifold as $\Nh$, we may pull back the hyperbolic structure of $\Nh$ via the covering map $q$ to obtain a hyperbolic structure on $N'$. We will denote the manifold $N'$ equipped with this hyperbolic structure by $\Nh'$. By construction the PL structure of $N'$ is compatible with the hyperbolic structure of $\Nh'$, so that the PL manifold $Nh'\pl$ (whose PL structure is defined only up to PL homeomorphism) may be taken to be $N'$. Condition (1) then holds with $\Nh'$ and $\calq'$
%  a manifold cover, we claim
%that, under the general hypotheses of the theorem:
%\Claim\label{before general}
%There is a finite-sheeted cover $p:\toldOmega\to\oldOmega$ such that
%the conditions (1) and (2) hold with $\toldOmega$ and $\toldTheta:=p^{-1}(\oldTheta)$ 
in place
of $\Mh$ and $\oldTheta$. 
According to \ref{new before general},
Condition (2) also holds with $\Nh'$ and $\calq'$
in place
of $\oldOmega$ and $\oldTheta$. 
It therefore follows from the case
of the theorem already proved that
%7.1 of \cite{ast} 
\Equation\label{tea quation}
\vol(\Nh')\ge\volorb((N')\cut{\calq'}).
\EndEquation

%\EndClaim
%\redproofsummary{Use the hypothesis that $\oldOmega$ is very good to get
%  a manifold cover. To get a higher cover for which (2) holds, I may
  %need to distinguish the case where some component of $\oldTheta$
%  bounds {\it two} twisted $I$-bundles; in that case one can show it's
  %the only component, and the double cover that is fibered over $\SSS^1$
  %does the trick. If we're not in this case, there is a homomorphism
  %of $H_1$ onto $\FF_2 $ that maps each of the generators arising from
  %a twisted $I$-bundle to the generator; the kernel of this
  %homomorphism gives the required covering. I'm beginning to think
  %this is basically the same construction, and that I don't need two
  %cases.} This proves \ref{before general}.

%Now let $p:\toldOmega\to\oldOmega$  be the finite-sheeted covering
%given by \ref{before general}. 
%letsuch that
%$\toldOmega$ is a manifold \redcomment{we really need one in which (1)
  %and (2) both hold}, and set
%The $\toldTheta:=p^{-1}(\oldTheta)$ Since the conditions
%(1) and (2) hold with $\toldOmega$ and $\toldTheta$ in place of $M$
%and $T$, 
Now  if $d$ denotes the degree of the covering
$q\circ p:N'\to\oldOmega$, we have $\vol(\Nh')=d\cdot\vol(\Mh)$. On
the other hand, since $\calq'=(q\circ p)^{-1}(\oldTheta)$,
%$N'$ is a $d$-fold covering of $\oldOmega$
the manifold 
$(N')\cut{\calq'}$ is a $d$-fold covering of
% $\oldOmega$, 
%$\volG({\rm D}(
$\oldOmega\cut\oldTheta$, and
${\rm D}((N')\cut{\calq'})$ is therefore a $d$-fold covering of
% $\oldOmega$, 
%$\volG({\rm D}(
${\rm D}(\oldOmega\cut\oldTheta)$.
%${\rm D}N'$ is a $d$-fold covering of ${\rm D}\oldOmega$, 
Thus
the
definitions of $\volorb((N')\cut{\calq'})$,
$\volG({\rm D}(\oldOmega\cut\oldTheta))$, and $\volorb(\oldOmega\cut\oldTheta)$ give
$$\volorb((N')\cut{\calq'})=
\frac12\volG({\rm D} ((N')\cut{\calq'}))=
\frac d2\volG({\rm D} (\oldOmega\cut\oldTheta))=
d\cdot\volorb(\oldOmega\cut\oldTheta).$$ 
Hence the conclusion of the present theorem follows from (\ref{tea
  quation}) upon dividing both sides by $d$.
\EndProof
%\obd\plN'

\Corollary\label{smockollary}
If $\Mh$ is a
 closed, orientable, hyperbolic $3$-orbifold, and if we set $\oldOmega=(\Mh)\pl$, then $\vol\Mh=\smock_0(\oldOmega)$.
\EndCorollary

\Proof 
Since by definition we have $\smock_0(\oldOmega)=\sup_\oldTheta \volorb((\oldOmega)\cut\oldTheta)$, where $\oldTheta$ ranges
over all incompressible closed $2$-suborbifolds of $\oldOmega$,
the inequality $\vol\Mh\ge\smock_0(\oldOmega)$ follows
from Theorem \ref{darts theorem}.
% and the definition of
%$\smock_0(\oldOmega)$. 
To prove the inequality
$\vol\Mh\le\smock_0(\oldOmega)$, note that since $\Mh$ is hyperbolic, it admits a $d$-sheeted manifold
cover $\tMh$  for some integer $d>0$. According to
\cite[Theorem C.4.2]{bp} we have   $\volG \tMh=\vol \tMh$. In
view of Definition \ref{sem ting}, we then have
$\vol\Mh=(\vol\tMh)/d=(\volG\tMh)/d=\volG\Mh=\volG\oldOmega$. Since
$\oldOmega$ is closed, ${\rm D}\oldOmega$ is a disjoint union of two copies of
$\oldOmega$, so that the definition of the invariant $\volorb$ gives $\volorb(\oldOmega)=\volG(\oldOmega)=\vol\Mh$. But
(\ref{more kitsch}), applied with $\oldOmega$ playing the role
of $\oldPsi$, gives
$\volorb(\oldOmega)\le\smock_0(\oldOmega)$, and the conclusion follows.
\EndProof
%\obd oldOmega\plA

\Number\label{voct def}
The volume of a regular ideal octahedron in $\HH^3$ will be denoted
$\voct$. We have $\voct=3.6638\ldots$.
\EndNumber

\Theorem\label{i want my v-8}
If $\oldPsi$ is an orientable $3$-orbifold which is
componentwise strongly \simple\ and componentwise boundary-irreducible,
we have
$$\volorb(\oldPsi)\ge\voct\chibar(\kish(\oldPsi))$$
(where $\volorb(\oldPsi)$ is defined since $\oldPsi$ is very good by \ref{oops}).
\EndTheorem

\Proof
Since both sides are additive over components, we may assume that
$\oldPsi$ is connected, and therefore strongly \simple\ and boundary-irreducible.

First consider the case in which $\oldPsi=M$ is a strongly 
\simple, boundary-irreducible, orientable 
$3$-manifold. In this case, because of the equivalence between the PL and smooth categories for $3$-manifolds, $M$ has a smooth structure compatible with its PL structure, and up to topological isotopy, the smooth $2$-submanifolds of $M$ are the same as its PL $2$-submanifolds. This will make it unnecessary to distinguish between smooth and PL annuli in the following discussion.

Since the components of $\frakA(\oldPsi)$ are  annular orbifolds by \ref{tuesa day}, and since $\frakA(\oldPsi)=\frakA(M)$ is a $2$-manifold in this case,
the components of $\Fr_M(\calb)\subset \frakA(M)$ are annuli. Let
$E_1,\ldots,E_n$ denote the components of $\kish M=\overline{M-\book M}$. Then
$\cala_i:=\Fr_ME_i$ is a union of components of $\frakA(M)$ for
 $i=1,\ldots, n$. According to the manifold case of Lemma \ref{i see a book}, the pair $(E_i,E_i\cap\partial M)$ is acylindrical for $i=1,\ldots,n$. By the manifold case of Corollary \ref{less than nothing}, we have $\chi(E_i)<0$ for $i=1,\ldots,n$. It is a standard consequence of Thurston's geometrization theorem that if $(E,F)$ is an acylindrical manifold pair with $E$ compact and $\chi(E)<0$, and if every component of $\overline{\partial E-F}$ is an annulus, then $(\inter E)\cup(\inter F)\subset E$ has a finite-volume hyperbolic metric with totally geodesic boundary. Thus for $i=1,\ldots,n$ the manifold
 $E_i^-:=E_i-\cala_i$ has such a
metric. According to the three-dimensional case of Theorem 4.2
of \cite{miyamoto}, we have
$\vol E_i^-\ge\voct\chibar(E_i)$. Hence  ${\rm D}(E_i^-)=\inter {\rm D}_{E_i'\cap\partial M}E_i\subset {\rm D}_{E_i'\cap\partial M}\subset {\rm D}M$ is a hyperbolic manifold
of finite volume, and 
\Equation\label{john donne}
\vol {\rm D}(E_i^-)=2\vol (E_i^-)\ge2\voct\chibar(E_i)).
\EndEquation
But for each $i$, since the components of the $\cala_i$ are essential
annuli, and since ${\rm D}M$ is canonically identified with a two-sheeted
covering of $\silv M$, it follows from the manifold case of Proposition \ref{silver acylindrical} that the components of ${\rm D}\cala_i=\Fr_{{\rm D}M} {\rm D}(E_i^-)$ are incompressible tori in ${\rm D}M$. Hence
%It therefore follows from
\cite[Theorem 1]{soma}, together with the fact (see 
\cite[Section 6.5]{thurstonnotes}) that the volume of a finite-volume hyperbolic $3$-manifold is equal to the relative Gromov volume of its compact core, implies that
  $\volG({\rm D}M)\ge\vol E_1'+\cdots+\vol E_n'$ (an inequality which could be shown to be an equality). With
(\ref{john donne}), this gives
$$\volG({\rm D}M)\ge2\voct(\chibar(E_1)+\cdots+\chibar(E_n))=2\voct\chibar(\kish(M)),$$
so that
$$\volorb(M)=\frac12\volG({\rm D}M)\ge\voct\chibar(\kish(M)).$$
This completes the proof in the case where $\oldPsi$ is a manifold.

Now suppose that $\oldPsi$ is an arbitrary strongly \simple, boundary-irreducible,
orientable $3$-orbifold. According to Condition (II) of Definition \ref{oops}, $\oldPsi$ admits a finite-sheeted covering $p:\toldPsi\to\oldPsi$
such that $\toldPsi$ is an irreducible $3$-manifold. Since, by Condition (III) of Definition \ref{oops}, $\pi_1(\oldPsi)$ has no rank-$2$ free abelian subgroup,  $\pi_1(\toldPsi)$ also has no rank-$2$ free abelian subgroup. Furthermore, since Condition (III) of Definition \ref{oops} implies that $\oldPsi$ is not discal, $\toldPsi$ is also non-discal; and $\toldPsi$ is a degree-$1$ regular irreducible $3$-manifold covering of itself. Hence $\toldPsi$ is strongly \simple.
Since $\oldPsi$ is boundary-irreducible, its covering $\toldPsi$ is also boundary-irreducible. Hence by the case of the theorem for a  manifold, which was proved
above, we have 
\Equation\label{summer vacation}
\volorb(\toldPsi)\ge\voct\chibar(\kish(\toldPsi)).
\EndEquation

Let $d$ denote the degree of the
covering. Doubling $p$, we obtain a $d$-fold covering map ${\rm D}\toldPsi\to
{\rm D}\oldPsi$. Hence by Definition \ref{sem ting}, we have
$\volG({\rm D}\toldPsi)=d\volG({\rm D}\oldPsi)$. Dividing both sides of the latter
equality by $2$, and applying the definition of $\volorb$, we obtain 
\Equation\label{geese}
\volorb(\toldPsi)=d\cdot\volorb(\oldPsi).
\EndEquation
On the other hand, it follows from Proposition \ref{Comment A} that $\kish(\toldPsi)$ is (orbifold-)isotopic to
$p^{-1}(\kish(\oldPsi))$; thus $\kish(\toldPsi)$ is homeomorphic to a $d$-fold
covering of the (possibly disconnected) orbifold $\kish(\oldPsi)$, so
that 
\Equation\label{more geese}
\chibar(\kish(\toldPsi))=d\cdot \chibar(\kish(\oldPsi)).
\EndEquation
Combining (\ref{summer vacation}) with (\ref{geese}) and (\ref{more geese}), we obtain
$$d\cdot\volorb(\oldPsi)=\volorb(\toldPsi)
\ge\voct\chibar(\kish(\toldPsi))=d\cdot\voct\chibar(\kish(\oldPsi)),$$
which gives the conclusion. 
\EndProof
%\obd\oldUpsilon_0\silv

\Corollary\label{bloody hell}
If $\oldPsi$ is a compact, orientable $3$-orbifold which is
componentwise strongly \simple\ and componentwise boundary-irreducible (and hence very good by \ref{oops}),
we have
$$\smock(\oldPsi)\ge\delta(\oldPsi).$$
\EndCorollary

\Proof
Let $\oldPi$ be any (possibly empty and possibly disconnected) incompressible, closed $2$-suborbifold of $\inter\oldPsi$. According to Lemma \ref{oops lemma}, the components of
$\oldPsi\cut\oldPi$ are strongly \simple\ and boundary-irreducible, and so 
Theorem \ref{i want my v-8} gives
$\volorb(\oldPsi\cut\oldPi)\ge\voct\chibar(\kish(\oldPsi\cut\oldPi))$. Hence,
using that $\voct=3.6638\ldots>3.66$, and using the definitions of
$\smock_0(\oldPsi)$ and $\smock(\oldPsi)$, we obtain
$$\sigma(\oldPsi\cut\oldPi))=12\chibar(\kish(\oldPsi\cut\oldPi))\le\frac{12}{\voct} \volorb(\oldPsi\cut\oldPi)
\le\frac{12}{\voct} \smock_0(\oldPsi)<\frac {\smock_0(\oldPsi)}{0.305}
=\smock(\oldPsi).
$$
Since the inequality $\sigma(\oldPsi\cut\oldPi))<\smock(\oldPsi)$ holds for every choice of $\oldPi$, we have
$\delta(\oldPsi)=\sup_\oldPi \sigma(\oldPsi\cut\oldPi))\le\smock(\oldPsi)$.
\EndProof

%The following corollary will be used in \cite{second}.

\Corollary\label{lollapalooka}
Let $\Mh$ be a closed,
orientable hyperbolic $3$-orbifold. Set $\oldOmega=(\Mh)\pl$ and $M=|\oldOmega|$.
Let $\cals $ be a (possibly disconnected) $2$-submanifold  of $M$, 
in general position with respect
to
$\fraks_\oldOmega$, such that $\obd(\cals)$ is incompressible in $\oldOmega$,
and let $X$ be a union of components of
$M- \cals $. Then
$$s_{{\oldOmega}}(X)\le y_{{\oldOmega}}(X)\le t_\oldOmega(X).$$ 
\EndCorollary

\Proof
According to the definitions, this
%the inequality $y_{{\oldOmega}}(X)\le
%t_\oldOmega(X)$ 
means that $\sigma(\omega(\hatX))\le\delta(\omega(\hatX))\le
\zeta(\omega(\hatX))$.  
The inequality $\sigma(\omega(\hatX)\le\delta(\omega(\hatX)$ follows
from (\ref{kitsch}), and the inequality
$\delta(\omega(\hatX)\le
\zeta(\omega(\hatX))$
follows from Corollary \ref{bloody hell}. (The conditions that $\hatX$ is componentwise strongly \simple\ and componentwise boundary-irreducible, which are needed for the application of Corollary \ref{bloody hell}, follow from Lemma \ref{oops lemma}, in view of the strong \simple ity of $\oldOmega$ and the incompressibility of  $\obd(\cals)$.)
\EndProof

\Proposition\label{hepcat} Let $\oldPsi$ be a very good, compact, orientable
$3$-orbifold, and set $N=|\oldPsi|$. 
\begin{itemize} 
\item If $\oldPsi$ (or equivalently $N$) is closed, then $\volG(N)\le\volG(\oldPsi)$.
\item If every component of
$\partial N$ is a sphere (so that $\plusN$ is a  closed $3$-manifold), then
$\volG(\plusN)\le\volorb(\oldPsi)$.
\end{itemize}
\EndProposition

\Proof
It follows from the definitions that the two sides of each of the asserted
inequalities are additive over components. Hence we may assume that
$\oldPsi$, or equivalently $N$, is connected.

To prove the first assertion, fix a covering $p:M\to\oldPsi  $ such that $M$ is a manifold and $d:=\deg p<\infty$. Then the definition of $\volG\oldPsi$ gives $\volG(\oldPsi)=\volG(M)/d$.
On the other hand, since $\oldPsi$ is orientable, $|p|:M\to N$ is a branched covering of degree $d$, and hence by \cite[Section 0.2]{gromov} (or \cite[Remark C.3.3]{bp}) we have $\volG M\ge d\cdot\volG N$. Thus 
$\volG(N)\le\volG(M)/d =\volG(\oldPsi)$, as required.

To prove the second assertion, first consider the case in which $N$ is closed. We have $\plusN=N$ in this case. Furthermore, ${\rm D}\oldPsi$ is a disjoint union of two copies of
$\oldPsi$, and hence $\volorb(\oldPsi)=\volG(\oldPsi)$. In view of the first assertion, we now have
%If $p:M\to\oldPsi  $ is a covering such that $M$ is a manifold and $d:=\deg p<\infty$, then the definition of $\volG\oldPsi$ gives $\volG(\oldPsi)=\volG(M)/d$.
%On the other hand, since $\oldPsi$ is orientable, $|p|:M\to N$ is a branched covering of degree $d$, and hence by \cite[Section 0.2]{gromov} (or \cite[Remark C.3.3]{bp}) we have $\volG M\ge d\cdot\volG N$. Thus 
$\volG(\plusN)=\volG(N)\le\volG(\oldPsi)=\volorb(\oldPsi)$, as required.

% by \redrealmissingref{reference or
  %cross-reference}. 
Now suppose that $N$ has $n\ge1$
boundary components. Since the components of $\partial N$ are spheres,
${\rm D}N$ is homeomorphic to the connected sum of two copies of $\plusN$ and
$n-1$ copies of $\SSS^2\times \SSS^1$. Since Gromov volume is additive under
connected sum (see \cite[Theorem 1]{soma} or \cite[Section 3.5]{gromov}), and since $\volG(\SSS^2\times \SSS^1)=0$ by \cite[Proposition C.3.4]{bp}), it follows
that $\volG({\rm D}N)=2\volG(\plusN)$. Since
$|{\rm D}\oldPsi |={\rm D}N$ we may write this as 
\Equation\label{deep-sea fission}
\volG(|{\rm D}\oldPsi |)=2\volG(\plusN).
\EndEquation
Now the closed orbifold ${\rm D}\oldPsi$ is very good since $\oldPsi$ is good, and we may apply the first assertion with ${\rm D}\oldPsi$ in place of $\oldPsi$ to deduce that  $\volG(|{\rm D}\oldPsi|)\le\volG({\rm D}\oldPsi)$. Combining this with the definition of
$\volorb(\oldPsi)$ and 
%, \redrealmissingref{It said ``same reference or
%  cross-reference,'' meaning it's the same basic argt as in the closed case. Figure out how to avoid doing the argt twice}, and 
(\ref{deep-sea fission}) we find
$$2\volorb(\oldPsi)=\volG({\rm D}\oldPsi)
%\le \volG({\rm D}\oldPsi)/=2\volG(N),$$
\ge \volG(|{\rm D}\oldPsi |)=2\volG(\plusN),$$
from which the conclusion follows.
\EndProof
%\obd

\Corollary\label{bloody hep}
If $\oldPsi$ is a very good, compact, orientable $3$-orbifold such that every
component of $\partial|\oldPsi |$ is a sphere,
we have
$$\smock(\oldPsi)\ge\theta(\oldPsi).$$
\EndCorollary

\Proof
Set $N=|\oldPsi|$.
By (\ref{more kitsch}) 
 and Proposition \ref{hepcat} we have 
%$\volorb(\oldPsi)
$\smock_0(\oldPsi)\ge \volorb(\oldPsi) \ge \volG(\plusN)$.
By definition we have $\smock(\oldPsi)= \smock_0(\oldPsi)/0.305$, and hence 
$\smock(\oldPsi)\ge
\volG(\plusN)/0.305
%\smock_0(\oldPsi)
=\theta(\oldPsi)$. 
\EndProof

%The following corollary will be used in \cite{second}.

\Corollary\label{lollaheplooka}
Let $\Mh$ be a closed,
orientable hyperbolic $3$-orbifold. Set $\oldOmega=(\Mh)\pl$ and $M=|\oldOmega|$.
Let $\cals $ be a (possibly disconnected) $2$-submanifold of $M$, 
in general position with respect to
$\fraks_\oldOmega$, each of whose components is a sphere, and let $X$ be a union of components of
$M- \cals $.
Then
$$q_{{\oldOmega}}(X)\le t_\oldOmega(X).$$ 
\EndCorollary

\Proof
Using Corollary \ref{bloody hep} and the definitions, we find that
$q_{{\oldOmega}}(X)=
\theta(\omega(\hatX))\le
\smock(\omega(\hatX)
)=t_\oldOmega(X)$.
\EndProof

\chapter{Hyperbolic $3$-orbifolds whose underlying manifolds are
  irreducible}\label{irreducible chapter}

This chapter is devoted to the proofs of Propositions A, E, F and G of
the Introduction (Propositions \ref{my little sony}, \ref{lost
  corollary}, \ref{orbifirst} and \ref{orbinext}). We mentioned in the
Introduction that Proposition \ref{new get lost} is the main result
underlying Proposition E (Proposition \ref{lost corollary}), and will also be needed in the final section
for the proof of Theorem B. The technical preparation for the proof of Proposition \ref{new get lost}  occupy Sections \ref{seek 'n'
  find}, \ref{underlying tori section}. As was indicated in the
Introduction, this material involves studying incompressible
suborbifolds of a  $3$-orbifold whose underlying surfaces are
incompressible tori in the underlying $3$-manifold; the needed
background about incompressible tori in $3$-manifolds is done in
Section \ref{tori section}. The proofs
of Propositions \ref{new get lost} and  \ref{lost corollary} appear
in \ref{irr-M section}. In Section \ref{A and C} we give the
self-contained proof of Proposition A, and combine it with Proposition
E to deduce Propositions F and G. 

\section{Tori in $3$-manifolds}\label{tori section}

This section is exclusively concerned with $3$-manifolds. We have found it convenient to adopt classical terminology here for Seifert fibered spaces \cite{hempel}. Although these could be regarded as $3$-manifolds equipped with $\SSS^1$-fibrations over orbifolds, we regard them here as $3$-manifolds equipped with maps to $2$-manifolds, with the local structure described in \cite[Chapter 12]{hempel}.

\Lemma\label{torus goes to cylinder}
Let $M$ be a compact, orientable, irreducible $3$-manifold. Suppose
that every component of $\partial M$ is a torus, that $M$ admits no
Seifert fibration, and that every incompressible torus in $\inter M$
is boundary-parallel in $M$. Then $M$ is acylindrical. 
\EndLemma

\Proof
We may assume that $\partial M\ne\emptyset$, as otherwise the conclusion is vacuously true.
Suppose that $A$ is a $\pi_1$-injective, properly embedded annulus in $M$. We must show that $A$ is boundary-parallel. Let $N$ denote a regular neighborhood of the union of $A$ with the component or components of $\partial M$ that meet $A$. Since each component of $\partial M$ is a torus, $N$ admits a Seifert fibration. 

Let $T$ be any component of $\Fr_M N$. Since $N$ admits a Seifert fibration, $T$ is a torus. By the manifold case of Lemma \ref{prepre}, it now follows that either 
(a) $T$ bounds a solid torus in $\inter M$,
or (b)
$T$ is contained in the interior of a ball in
$\inter M$, or (c) $T$ is $\pi_1$-injective in $M$.
But since
$A\subset M$ is a properly embedded, $\pi_1$-injective annulus, $T$ must contain an annulus which is $\pi_1$-injective in $M$; in particular, the inclusion homomorphism $\pi_1(T)\to\pi_1(M)$ is non-trivial. Hence (a) cannot hold. If (c) holds then $T$ is boundary-parallel according to the hypothesis of the present lemma. This shows that every 
 component $T$ of $\Fr_M N$ is the frontier of a submanifold $K_T$ of $M$ which either is a solid torus or is homeomorphic to $\TTT^2\times[0,1]$. We fix such a $K_T$ for every  component $T$ of $\Fr_M N$.

For each  component $T$ of $\Fr_M N$, we must have either $K_T\supset N$ or $K_T\cap N=T$. Consider the case in which
$K_{T_0}\supset N$ for some component $t_0$ of $\Fr_MN$. The construction of $N$ implies that $N$ contains at least one boundary component of $M$, say $t_1$. Then $t_0$ and $t_1$ are distinct components of $\partial K_{T_0}$. Hence $K_{T_0}$ cannot be a solid torus, and must therefore be homeomorphic to $\TTT^2\times[0,1]$. In particular, $t_0$ and $t_1$ are the only components of $\partial K_{T_0}$. Since $A\subset N\subset K_{T_0}$ has its boundary contained in $\partial N$, we must have $\partial A$ contained in $t_1$. It now follows from \cite[Proposition 3.1]{Waldhausen} that $A$ is boundary-parallel in $K_{T_0}$, and hence in $M$; this is the required conclusion.

There remains the case in which
 $K_T\cap N=T$ for each  component $T$ of $\Fr_M N$. If we set $\calt=\calc(\Fr_MN)$, we have $M=N\cup\bigcup_{T\in\calt}K_T$. Fix a Seifert fibration of $N$, and let $\calt'\subset\calt$ denote the set of all components $T$ of $\Fr_MN$ such that either (i) $K_T$ is homeomorphic to $\TTT^2\times[0,1]$ or (ii) $K_T$ is a solid torus whose meridian curve is isotopic in $T$ to a fiber of $N$. Then the Seifert fibration of $N$ extends to a Seifert fibration of $M':=N\cup\bigcup_{T\in\calt'}K_T$. If $\calt'=\calt$ then $M=M'$, so that $M$ is Seifert fibered, a contradiction to the hypothesis. If $\calt'\ne\calt$ then $M$ is obtained from the Seifert fibered space $M'$ by attaching one or more solid tori along boundary components of $M'$ in such a way that each attached solid torus has a meridian curve which is attached along a fiber of $M'$. Since $\partial M\ne\emptyset$, it follows that $M$ is either  solid torus (which admits a Seifert fibration) or a reducible $3$-manifold. In either case we have a contradiction to the hypothesis.
\EndProof
%M'N

\Lemma\label{and four if by zazmobile}
Let $M$ be a closed, connected, orientable $3$-manifold equipped with a Seifert
fibration. If $h(M)\ge4$, then $M$ contains a saturated incompressible torus.
\EndLemma

\Proof
Fix a Seifert fibration $p:M\to G$, where $G$ is a compact, connected $2$-manifold.
let $n\ge0$ denote the number of singular fibers, and let
$x_1,\ldots,x_n\in G$ denote the images of the singular fibers under $p$. If $G$
is not a sphere or a projective plane, then $G$ contains a simple
closed curve $C$ which does not bound a disk in $G$, and $p^{-1}(C)$ is a saturated
incompressible torus. 

Now suppose that $G$ is a sphere or a projective plane. Fix disjoint disks
$D_1,\ldots,D_n\subset G$ such that $x_i\in D_i$ for each $i$. 
For $i=1,\ldots,n$, set $V_i=p^{-1}(D_i)$.
% and $t_i=\partial
%V_i=p^{-1}(\partial D_i)$, so that $t_1,\ldots,T_n$ are the boundary components of $M$.
Set
$G'=G-\bigcup_{i=1}^n\inter D_i$, and
$M'=p^{-1}(G')=\overline{M-(V_1\cup\cdots\cup V_n)}$. The map
$p':=p|M':M'\to G'$ is a fibration with fiber $\SSS^1$, and hence $h(M')\le
h(G)+1
%\le n+1
$. Since the $V_i$ are solid tori, the inclusion homomorphism $H_1(M';\FF_2 )\to
H_1(M;\FF_2 )$ is surjective, we have $h(M)\le h(M')\le h(G')+1$.

In the case where $G$ is a sphere we have $h(G')=\max(0, n-1)$, and hence $h(M)\le\max( n,1)$. The hypothesis then implies that $n\ge4$. 
%In particular
%$\partial G'\ne\emptyset$, and since $G$ and $M$ are both orientable,
%$p'$ must be a trivial fibration. It follows that $h(M')=n+1$. 
%Since in particular we have $n\ge2$, the inclusion homomorphism
%$j:H_1(T_1,\FF_2 )\to H_1(M';\FF_2 )$ is injective.
%Now set
%$V_i=p^{-1}(D_1)$. 
%Since $j$ is injective, and since $V$ is a solid
%torus, we have $h_1(M'\cup V)=h_1(M)-1=n$. But the inclusion
%homomorphism $H_1(M'\cup V;\FF_2 )\to H_1(M;\FF_2 )$ is surjective since
%$V_2,\ldots,V_n$ are solid tori; hence $h(M)\le h(M'\cup V)=n$. With
%the hypothesis we now deduce that $n\ge4$. 
Hence there is a simple
closed curve $C\subset G-\{x_1,\ldots,x_n\}$ such that each of the
disks bounded by $C$ contains at least two of the $x_i$. This implies
that $p^{-1}(C)$ is a saturated incompressible torus.

In the case where $G$ is a projective plane we have $h(G')=\max(1, n)$, and hence $h(M)\le\max(2, n+1)$. The hypothesis then implies that $n\ge3$. 
%In particular
%$\partial G'\ne\emptyset$, and since $G$ and $M$ are both orientable,
%$p'$ must be a trivial fibration. It follows that $h(M')=n+1$. 
%Since in particular we have $n\ge2$, the inclusion homomorphism
%$j:H_1(T_1,\FF_2 )\to H_1(M';\FF_2 )$ is injective.
%Now set
%$V_i=p^{-1}(D_1)$. 
%Since $j$ is injective, and since $V$ is a solid
%torus, we have $h_1(M'\cup V)=h_1(M)-1=n$. But the inclusion
%homomorphism $H_1(M'\cup V;\FF_2 )\to H_1(M;\FF_2 )$ is surjective since
%$V_2,\ldots,V_n$ are solid tori; hence $h(M)\le h(M'\cup V)=n$. With
%the hypothesis we now deduce that $n\ge4$. 
In this case, define $C\subset G$ to be a simple closed curve bounding a disk $E\subset G$ with $x_1,\ldots,x_n\in\inter E$. Since $n\ge3>2$, and since $G-\inter E$ is a M\"obius band,
%closed curve $C\subset G-\{x_1,\ldots,x_n\}$ such that each of the
%disks bounded by $C$ contains at least two of the $x_i$. This implies
$p^{-1}(C)$ is a saturated incompressible torus.
\EndProof

\Lemma\label{before graphology}
Let $\Sigma$ be a compact, orientable (but possibly disconnected)
$3$-manifold, and let $A\subset\partial\Sigma$ be a compact
$2$-manifold. Suppose that $\Sigma$ admits a Seifert fibration in
which $A$ is saturated. Let
$F$ be an incompressible, closed (possibly disconnected) $2$-dimensional submanifold of
$\inter\Sigma$. Then
either $\Sigma$ admits a Seifert fibration in which $F\cup A$ is saturated,
or some component of $F$ is a fiber in some fibration of a component
of $\Sigma$ over $\SSS^1$.
\EndLemma

\Proof
First consider the case in which $F$ and $\Sigma$ are both
connected. In this case, \cite[Lemma II.7.3]{js} asserts that either
(a) $F\cup A$ is saturated in some Seifert fibration of $\Sigma$, or (b)
$\Sigma$ has a finite-sheeted covering space homeomorphic to
$F\times[0,1]$. The proof shows that (b) may be replaced by the
stronger condition (b') $F$ is a fiber in some fibration of $\Sigma$
over $\SSS^1$ for which the monodromy homeomorphism is of finite
order. In particular, the conclusion of the lemma holds in this case.

If $F$ is connected but $\Sigma$ is not, the assertion of the lemma follows upon
applying the case discussed above, with the role of $\Sigma$ played by
the component of $\Sigma$ containing $F$.

To prove the result in general, we use induction on $\compnum(F)$. If
$\compnum(F)=1$, i.e. $F$ is connected, we are in the case discussed
above, and if $\compnum(F)=0$, i.e. $F=\emptyset$, the result is
trivial. Now suppose that $\compnum(F)=n>1$ and that the result is
true for surfaces with $n-1$ components. Assume that no component of
$F$ is a fiber in a fibration of a component of $\Sigma$ over $\SSS^1$. Fix a component $F_0$ of $F$. Since the lemma has been
proved in the case of a connected surface, there is a Seifert
fibration of $p:\Sigma\to E$, for some compact $2$-manifold $E$, in which $F_0$ is saturated. Hence $F$ is a
torus. Let $N$ denote a
regular neighborhood of $F_0$ which is saturated in the Seifert
fibration $p$. We may identify $N$
homeomorphically with $F_0\times[-1,1]$ in such a way that $F_0$ is
identified with $F_0\times\{0\}$. For $\epsilon=\pm1$, fix an annulus
$C_\epsilon\subset F_0\times\{\epsilon\}\subset\partial N$ which is
saturated in the Seifert fibration $p$. Note that the annuli $C_1$ and
$C_{-1}$ are homotopic in $N$. Set
$\Sigma'=\overline{\Sigma-N}$, and set $A'=A\cup C_1\cup
C_{-1}\subset\partial\Sigma'$, so that $A'$ is saturated with respect to
the Seifert fibration $p|\Sigma'$. Now set $F'=F-F_0$, so that
$\compnum(F')=n-1$. By the induction hypothesis, either (i) some
component $F_1'$
of $F'$ is a fiber in a fibration over $\SSS^1$ of a component $\Sigma_1'$ of $\Sigma'$
containing $F_1'$, or (ii) $F'\cup A'$ is saturated in some Seifert fibration
of $\Sigma'$. If (i) holds, then since 
no component of
$F$ is a fiber in a fibration of a component of $\Sigma$ over $\SSS^1$,
we must have $\Sigma_1'\subset\Sigma_0$, where $\Sigma_0$ denotes the component of $\Sigma$
containing $F_0$. This implies that
$\Sigma_1'\subset\overline\Sigma_0-N$, so that $\Sigma_1'$ is not
closed; this is impossible, because $F_1'$ is a closed surface, and
therefore cannot be the fiber in a fibration of a non-closed
$3$-manifold over $\SSS^1$. Hence (ii) must hold. Fix a Seifert fibration
$q:\Sigma'\to B$, for some compact $2$-manifold $B$, in which $F'$ and $A'$ are saturated. In particular the
annuli $C_1$ and $C_{-1}$ are saturated in the fibration $q$. Since
$C_1$ and $C_{-1}$ are homotopic in $N=F_0\times[-1,1]$, the Seifert
fibration $q|(F_0\times\{-1,1\})$ extends to a Seifert fibration of
$N$ in which $F_0=F_0\times\{0\}$ is saturated. It follows that $q$
extends to a Seifert fibration $\overline q$ of $\Sigma$ in which $F_0$ is
saturated. Since $F'$ and $A\subset A'$ are saturated in the Seifert fibration $q$, it now follows
that $F=F'\cup F_0$ and $A$ are saturated in the Seifert fibration
$\overline q$, and the induction is complete.
\EndProof
%\obd

\Lemma\label{graphology}
Let $M$ be  a closed graph manifold (see \ref{great day}). Let $X$ be a compact (possibly
empty) $2$-submanifold of $M$ whose boundary components are all
incompressible  tori in $M$. Suppose that $h(M)\ge \max(4,h(X))$, and let $m$ be an integer such that
$\max(2,h(X))\le m\le h(M)$. Then
there is a compact, connected submanifold $L$ of $ M$ such that
$X\subset \inter L$,
each component of
$\partial L$ is an incompressible  torus in $M$, and $h(L)=m$.
\EndLemma

\Proof
Let $\Sigma$ denote the characteristic submanifold of $M$. (By definition this means that $\Sigma=M$ if $M$ is Seifert-fibered, and that if $M$ is a non-Seifert-fibered graph manifold, and hence a Haken manifold, then $\Sigma$ is the characteristic submanifold as defined in \cite{js}; cf. \ref{manifolds are different}). Since $M$
is a graph manifold, each component  of $\overline{M-\Sigma}$ is
homeomorphic to $[0,1]\times \SSS^1\times  \SSS^1$. Since the
  components of $\partial X$ are incompressible tori,  we may assume after an isotopy that $\partial X\subset\inter\Sigma$. Since each
  component of
$\Sigma$ is a Seifert fibered space, $\Sigma$ itself admits a Seifert
fibration. Applying Lemma \ref{before graphology}, taking
$F=\partial X$ and $A=\emptyset$, we deduce that
either (i) some
component $F_0$ of $\partial X$ is a fiber in a fibration over $\SSS^1$
of some component $\Sigma_0$  of $\Sigma$, or (ii) $\partial X$ is saturated in some Seifert fibration
of $\Sigma$. If (i) holds, then since the components of $\partial X$
are tori, $\Sigma_0$ is closed and $h(\Sigma_0)\le3$. Since the graph
manifold $M$ is by definition irreducible and therefore connected, it
follows that $M=\Sigma_0$ and hence that $h(M)\le3$, a contradiction
to the hypothesis. Hence (ii) holds. Let us fix a
Seifert fibration
  $p:\Sigma\to B$, for some compact $2$-manifold $B$,
 in which $\partial X$ is saturated.

The images under $p$ of singular fibers of $\Sigma$ will be referred to as {\it singular points} of $B$.

Let $\calK$ denote the set of all
compact, connected, non-empty, $3$-dimensional submanifolds $K$ of $ M$ such that
(1) $X\subset \inter K$ and
(2) each component of
$\partial K$ is an incompressible saturated torus in $\inter\Sigma$.
% (3) $X\subset
%K$, and (4) $h(K)\le m$.

Set $\mu=\max(2,h(X))$. We claim that:
\Claim\label{oswald zdestiny}
There is an element $K$ of $\calk$ with $h(K)\le\mu$.
\EndClaim

Indeed, if $X\ne\emptyset$, then since $\partial X$ is saturated, $X$ admits a
regular neighborhood $K$ with  saturated boundary; it is clear that $K\in\calk$, and
$h(K)=h(X)\le\mu$ by the definition of $\mu$.
Hence \ref{oswald zdestiny} is true in this case. If
  $X=\emptyset$ and $\Sigma\ne M$, then
  $\partial\Sigma\ne\emptyset$. In this case, let us choose a
  saturated boundary-parallel torus
  $T\subset\inter\Sigma$, and let $K$ denote a saturated regular neighborhood
  of $T$ in $\inter\Sigma$. We clearly have $K\in\calk$, and we have  $h(K)=2\le \mu$. If
  $X=\emptyset$ and $\Sigma=M$, then since
  $h(M)\ge4$, Lemma \ref{and four if by zazmobile} gives a saturated incompressible torus in $M$. 
%either $B$ is a sphere containing at least four
  %singular points, or $B$ is a projective plane containing at least
  %three singular points, or $\chi(B)\le0$. In each of these cases,
  %there is a simple closed curve $C\subset B$ such that any disk in $B$
  %bounded by $C$ contains at least two singular points. In this case,
  If $K$ denotes a saturated neighborhood of such a torus, then  again we
  clearly have $K\in\calk$, and again we have  $h(K)=2\le \mu$. This proves \ref{oswald zdestiny}.

Now for each $K\in\calk$, since $\partial K$ is saturated in the
Seifert fibration of $\Sigma$, the $3$-dimensional submanifold
$\overline{\Sigma\setminus K}$ of $\Sigma$ is also saturated. Hence
%whose boundary contains $T$. The fact that
%$K\in\calk$ also implies that. Hence 
we may write $\overline{\Sigma\setminus K}=p^{-1}(R_K)$ for
a unique $2$-submanifold $R_K$ of $B$. 
Each component of
$\partial R_K$ is contained in either $\partial B$ or $p(\partial K)$.
Let $\alpha_1(K)$ denote the number of components
of $\partial B$ contained in $R_K$, let $\alpha_2(K)$ denote the number of singular points lying in $R_K$,
% =\card(\fraks_B\cap R_K)$ \redcomment{I don't even know what this means. What is , 
and let $\gamma(K)$ denote the sum of the squares of the first betti
numbers of the components of $R_K$. Set 
$\nu(K)=(\alpha_1(K)+\alpha_2(K),\gamma(K))\in\NN^2$. The
set $\NN^2$ of pairs of non-negative integers will be endowed with
the lexicographical order. We claim:

\Claim\label{sooey pig}
Suppose that $K$ is an element of $\calk$ with non-empty boundary, and
that $T$ is a component of
$\partial K$. Then $T=p^{-1}(\ell)$ for some simple closed curve
$\ell\subset\inter B$. 
%Furthermore, if 
% and let
%of $\ and hence $R_K\ne\emptyset$. Choose a
Furthermore, we have $\ell\subset\partial R_K$, and if $R_0$ denotes the component of $R_K$ whose boundary contains
$\ell$, then there exist
%In each case, I will construct 
a subsurface $\Delta$ of
$R_0$, and a regular neighborhood $K'$ of
%\begin{itemize}
$K\cup p^{-1}(\Delta)$ in $M$, such that
% whose
%boundary components are saturated tori in $\Sigma$,
%\end{itemize}
%$\Delta\cap p(K)$ is a $1$-dimensional
%submanifold of $p(\partial K)$, 
%and such that
 %for any sufficiently
%small regular neighborhood $K'$ of $K\cup p^{-1}(\Delta)$, we have
\begin{itemize}
%whose intersection with 
\item 
%$L:= and
$K'\in\calk$,
\item
$h(K')\le h(K)+1$, and
\item $\nu(K')<\nu(K)$.
\end{itemize}
\EndClaim

%. In each of the cases I will then set
%$J:=p^{-1}(\Delta)$ and $K'=K\cup J$, and show that $K'$ has the
%properties asserted in \ref{sooe
%. 

To prove \ref{sooey pig}, first note that the definition of $\calk$
implies that $T$ is a saturated torus in $\inter\Sigma$, and therefore has
the form $p^{-1}(\ell)$ for some simple closed curve
$\ell\subset\inter B$. Since $T$ is a component of $\partial K$, the
curve $\ell$ is a component of $\partial R_K$. Let $R_0$ denote the component of $R_K$ whose boundary contains
$\ell$. 

If $R_0$ were a disk containing at most one singular point, then
$p^{-1}(R_0)$ would be a solid torus with boundary $T$; this would
contradict the incompressibility of $T$ in $M$. Hence one of the
following cases must
occur: (i) $R_0$ is a disk containing exactly
two singular points; (ii) $R_0$ is a disk containing at
least three singular points, or $\chi(R_0)\le0$ and $R_0$ contains at
least one singular point; (iii) $R_0$ contains no singular points
and $\chi(R_0)<0$; (iv) $R_0$ is an annulus or M\"obius band containing no
singular points, and $\partial R_0\subset p(\partial K)$; or (v) $R_0$ is an annulus containing no
singular points, and the component of $\partial R_0$ distinct from
$\ell$ is contained in $\partial B$.

In Cases (i), (iv), and (v) we will take $\Delta=R_0$. In Case (ii) we
will take $\Delta$ to be a disk contained in $R_0$, meeting
$\partial R_0$ in an arc contained in $\ell$, and containing exactly
one singular point in its interior and none on its boundary. In Case
(iii) we will take $\Delta$ to be a regular neighborhood relative to $R_0$ of
a non-boundary-parallel  arc in $R_0$ which has both endpoints in
$\ell$. 

In each case, $\Delta$ is connected, and $\Delta\cap p(K)$ is a non-empty $1$-dimensional
submanifold of $p(\partial K)$; hence $K_0':=K\cup p^{-1}(\Delta)$ is a connected
$3$-manifold. Furthermore, since $K\in\calk$, each component of
$\partial K$ is an incompressible saturated torus in
$\inter\Sigma$. Hence each boundary component of $K_0'$ is a saturated torus in
$\Sigma$. 
%(In cases (i)---(iv) the boundary of $K_0'$ is contained in
%$\inter\Sigma$, but this is not true in Case (v).) 

It is also clear from the construction of $\Delta$ in each of the five cases that
no boundary component of $\overline{R_0-\Delta}$ bounds a disk in
$R_0$ containing at most one singular point. This implies that every
boundary component of $J:=p^{-1}(\overline{R_0-\Delta})$ is
$\pi_1$-injective in $Z:=p^{-1}({R_0})$. The components of
$\partial Z$ are incompressible in $M$, because each of them is a
component of either $\partial K$ or $\partial\Sigma$. Hence $Z$ is $\pi_1$-injective
in $M$; it follows that  $\partial J$ is incompressible in
$M$. But every boundary component of $K_0'$
is a boundary component of either
$K$ or $J$. Hence  $\partial K_0'$ is incompressible in $M$.

We will take $K'$ to
be a small regular neighborhood of $K_0'$ in $M$. Since the components of
$\partial K_0'$ are saturated tori in $\inter\Sigma$ and are incompressible in
$M$, we may choose  $K'$ so that its boundary components are also saturated tori in
$\inter\Sigma$ and are incompressible in $M$. Since in addition  $\inter K'\supset\inter K\supset X$, we have $K'\in\calk$. 

We must show that $h(K')\le h(K)+1$, or equivalently that $h(K_0')\le h(K)+1$. 
For this purpose it suffices to show that
$\dim H_1(K_0',K;\FF_2 )\le1$, or equivalently that $\dim H_1(L,Y;\FF_2 )
\le1$, where
$L=p^{-1}(\Delta)$ and $Y=p^{-1}(\Delta)\cap K=
p^{-1}(\Delta\cap\partial R_0)$. This is readily verified in each of
the cases (i)---(v). Indeed, in Case (i), $L $ is a
Seifert fibered space over the disk with two singular fibers and
$Y=\partial L$. 
In Case (ii), $L $ is a solid torus and
%Seifert fibered space over an annulus with one singular fiber and
$Y$ is a homotopically non-trivial annulus in
% component of 
$\partial L$. 
In Cases (iii) and (v), the pair $(L,Y) $ is
respectively homeomorphic to
$([0,1]\times [0,1]\times \SSS^1,[0,1]\times\{0,1\}\times \SSS^1)$ and
$(\SSS^1\times \SSS^1\times[0,1] ,\SSS^1\times \SSS^1\times\{0\})$,
%[0,1]\times\{0,1\}\times \SSS^1)$
%\SSS^1\times I\times\partial I)$ and $(\SSS^1\times
%\SSS^1 \times I,\SSS^1\times \SSS^1\times\{0\})$, 
while in Case (iv), $L$ is homeomorphic to an $\SSS^1$-bundle over an annulus or M\"obius band, and $Y=\partial L$.

%$(\SSS^1\times
%\SSS^1 \times I,\SSS^1\times \SSS^1\times\partial I)$ or to 

%KL

To prove \ref{sooey pig}, it remains to show that $\nu(K')<\nu(K)$. For this purpose, note that,
as a consequence of the definitions and the constructions, $\alpha_1(K)-\alpha_1(K')$ is equal
 to the
  number of components of $\partial B$ 
  contained in $\Delta$, while
$\alpha_2(K)-\alpha_2(K')$ is equal to 
 the number of points of singular points lying in $\Delta$.
This gives $\alpha_1(K')\le\alpha_1(K)$ and $\alpha_2(K')\le\alpha_2(K)$ for
  $i=1,2$, and at least one of these inequalities will be strict if
$\Delta$   contains either a
component of $\partial B$ or a singular point of $B$. It follows from the construction of $\Delta$
 that this is true in Cases (i), (ii), and (v). Thus in these cases we
 have $\alpha_1(K')+\alpha_2(K')<\alpha_1(K)+\alpha_2(K)$, and in
 particular $\nu(K')<\nu(K)$. In Cases (iii) and (iv) we have
$\alpha_1(K')+\alpha_2(K')\le\alpha_1(K)+\alpha_2(K)$, and we need to
prove that $\gamma(K')<\gamma(K)$. Since $R_{K'}$ is homeomorphic to $R_{K'_0}$, it suffices
to
prove that $\gamma(K'_0)<\gamma(K)$. 

In Case (iii), $R_{K'_0}$ is obtained from $R_K$ by removing a regular
neighborhood of a properly embedded arc in the component $R_0$ of
$R_K$; this arc is not boundary-parallel in $R_0$, and its endpoints
are in the same component of $\partial R_0$. If the arc does not
separate $R_0$, then $R_{K'_0}$ is obtained from $R_K$ by replacing the
component $R_0$ by a new component whose first betti number is one less
than that of $R_0$; hence $\gamma(K'_0)<\gamma(K)$ in this subcase. If the arc 
separates $R_0$, then $R_{K'_0}$ is obtained from $R_K$ by replacing the
component $R_0$ by two new components, whose first betti numbers are
strictly positive and add up to the first betti number of $R$. Hence
$\gamma(K'_0)<\gamma(K)$ in this subcase as well. In Case (iv), $R_{K'_0}$
is obtained from $R_K$ by discarding a component which is an annulus or a M\"obius band. Hence
$\gamma(K'_0)=\gamma(K)-1$ in this case. This completes the proof of
\ref{sooey pig}.
%K'

Now let $\calk^*$ denote the subset of $\calk$ consisting of all
elements $K$ such that $h(K)\le m$. Since $m\ge\mu$ by hypothesis, it
follows from \ref{oswald zdestiny} that $\calk^*\ne\emptyset$. Since $\NN^2$
is well ordered, it now follows  that there is an element $L\in \calk^*$ such that
$\nu(L)\le\nu(K)$ for every $K\in\calk^*$. By the definition of $\calk^*$
we have $h(L)\le m$. Suppose that $h(L)<m$. Since $h(M)\ge m$ by hypothesis, we then have $h(L)<h(M)$, so
that $L\ne M$ and therefore $\partial L\ne\emptyset$. Hence
 by \ref{sooey pig}, applied with $L$ playing the role of $K$, there is an
element $L'$ of $\calk$ such that $\nu(L')<\nu(L)$ and $h(L')\le
h(L)+1\le m$. This means that $h(L')\in\calk^*$, a
contradiction. This shows that $h(L)=m$. In view of the definition of $\calk$,
it follows that $L$ has the properties asserted in the lemma.
\EndProof

\Lemma\label{even easier} Let $Y$ be a compact, connected, orientable
$3$-manifold.  Let $\calt\subset\inter Y$ be a compact $2$-manifold
whose components are all of strictly positive genus, and suppose that
$Y-\calt$ has exactly two components, $B$ and $C$. Then $h(Y)\le
h(B)+h(C)-1$.
\EndLemma

\Proof
All homology groups considered in this proof will be understood to
have coefficients in $\FF_2 $. Let $m$ denote the number of components of
$\calt$. Consider the exact sequence
\Equation\label{yellow pig}
H_1(\calt)\to H_1(\overline B)\oplus
H_1(\overline C)\to H_1(Y)\to H_0(\calt)\to H_0(\overline B)\oplus
H_0(\overline C)\to H_0(Y)\to 0,
\EndEquation
which is a fragment of the Mayer-Vietoris sequence.
Let $V$ denote the image of the map $H_1(\calt)\to H_1(\overline B)\oplus
H_1(\overline C)$ in (\ref{yellow pig}). 
The sequence (\ref{yellow pig}) gives rise to an exact sequence
\Equation\label{purple pig}
0\to V\to H_1(\overline B)\oplus
H_1(\overline C)\to H_1(Y)\to H_0(\calt)\to H_0(\overline B)\oplus
H_0(\overline C)\to H_0(Y)\to 0.
\EndEquation
The exactness of (\ref{purple pig}) implies that
$$\begin{aligned}
0&=\dim V-\dim(H_1(\overline B)\oplus
H_1(\overline C))+\dim H_1(Y)-\dim H_0(\calt)+\dim (H_0(\overline B)\oplus
H_0(\overline C))-\dim H_0(Y)\\
&=\dim V-(h(B)+h(C))+h(Y)-m+2-1,
%\\
%&\ge h(Y)-(h(B)+h(C))+1,
\end{aligned}$$
%which gives the conclusion.
so that
\Equation\label{knute rockne}
h(Y)=(h(B)+h(C))+(m-\dim V)-1.
\EndEquation
Now let $i:H_1(\calt)\to H_1(\bar B)$ and $j:H_1(\calt)\to H_1(\bar
C)$ denote the inclusion homomorphisms. The intersection pairing on
$H_1(\calt)$ is non-singular because $\calt$ is a closed surface; but
since $\calt\subset\partial B$, this pairing is trivial on $\ker
i$. Hence $\dim\ker i\le(\dim H_1(\calt ))/2$, and therefore $\dim \image
i\ge (\dim H_1(\calt ))/2$.  Since each component of $\calt$ has strictly
positive genus, we have $\dim H_1(\calt )\ge 2m$, and so $\dim \image
i\ge m$. On the other hand, the map $H_1(\calt)\to H_1(\overline B)\oplus
H_1(\overline C)$ in (\ref{yellow pig}) is defined by
$x\mapsto(i(x),j(x))$, and so $\dim V\ge\dim\image i$. It follows that
\Equation\label{punkin}
\dim V\ge m.
\EndEquation
The conclusion of the lemma follows from (\ref{knute rockne}) and (\ref{punkin}).
\EndProof
%B

\Lemma\label{another goddam torus lemma}
Let $U$ and $V$ be compact, connected $3$-dimensional submanifolds of an
orientable $3$-manifold. Suppose that $\inter
U\cap\inter V=\emptyset$, and that every component of $U\cap V$ is a
torus. Then $h(U\cup V)\ge h(U)-1$.
\EndLemma

\Proof
All homology groups in this proof will be understood to have
coefficients in $\FF_2 $. Set $\calt= U\cap V$ and $Z=U\cup V$. We may
assume that $\calt\ne\emptyset$, so that $Z$ is connected. Then $\calt$ is a union of $m$ boundary tori of
$V$ for some integer $m\ge1$. The
intersection pairing on $H_1(\calt;\FF_2 )$ is nondegenerate, and the
kernel $K$ of the inclusion homomorphism $H_1(\calt)\to
H_1(V)$ is self-orthogonal with respect to this pairing. Hence
$\dim K\le (\dim H_1(\calt;\FF_2 ))/2=m$. The exactness of the
Mayer-Vietoris fragment
$$H_1(\calt)\longrightarrow H_1(U)\oplus H_1(V)\longrightarrow
H_1(Z)$$
implies that the kernel $L$ of the inclusion homomorphism $H_1(U)\to
H_1(Z)$ is the image of $K$ under the inclusion homomorphism $H_1(\calt)\to
H_1(U)$. Hence $\dim L\le m$. Now since $U$ and $Z$ are
connected, the homology exact
sequence of the pair $(Z,U)$ gives rise to an
exact sequence
$$0\longrightarrow L \longrightarrow H_1(U) \longrightarrow H_1(Z)
\longrightarrow H_1(Z,U) \to0,$$
which implies that 
\Equation\label{glass my ass}
h(U)-h(Z)=\dim L-\dim H_1(Z,U)\le m-\dim H_1(Z,U).
\EndEquation
%shows that 
On the other hand, by excision we have $H_1(Z,U)\cong H_1(V,\calt)$,
and the exact homology sequence
$$H_1(V,\calt)\longrightarrow H_0(\calt)\longrightarrow H_0(V)$$
shows that $\dim H_1(V,\calt)\ge\dim H_0(\calt)-\dim H_0(V)=m-1$. Thus
$\dim H_1(Z,U)\ge m-1$, which combined with (\ref{glass my ass}) gives $h(U)-h(Z)\le1$.
\EndProof
%UVXY

\Lemma\label{just right}
Let $m\ge2$ be an integer, and
let $M$ be  a closed graph manifold with $h(M)>\max(3,m(m-1))$. Then there is
a compact submanifold $P$ of $M$ such that
\begin{itemize}
\item each component of $\partial P$ is an incompressible  torus in
  $M$,
\item $P$ and $M-P$ are connected, and
\item $\min(h(P),h(M-P))\ge m$.
\end{itemize}
\EndLemma

\Proof
Note that the hypothesis implies that 
$h(M)\ge 4$ and that
$2\le m\le h(M)$. Hence by Lemma \ref{graphology}, applied with $X=\emptyset$, there is a compact, connected submanifold $L$ of $ M$ such that
(a) each component of
$\partial L$ is an incompressible  torus in $M$, and (b) $h(L)=m$. It
follows from the special case $T=\emptyset$ of \cite[Theorem V.2.1]{js}
that we may
  take $L$ to be ``maximal'' among all compact, connected submanifolds
  satisfying (a) and (b) in the sense that, if $L'$ is any
  compact, connected submanifold of $\inter M$, such that
  $L\subset\inter L'$,
each component of
$\partial L'$ is an incompressible  torus in $M$, and $h(L')=m$, then
$L'$ is a regular neighborhood of $L$. Since the hypothesis implies
that $h(M)>m$, we have $L\ne M$, so that $\partial L\ne\emptyset$.

It is a standard consequence of Poincar\'e-Lefschetz duality that
the total genus of the boundary of a compact, connected, orientable
$3$-manifold $L$ is at most $h(L)$. Since the boundary components of
the manifold $L$ we have chosen are tori, it follows that $\partial L$
has at most $m$ components. Thus if $k=\compnum(\partial L)$ and
$r=\compnum(M-L)$, we have $1\le r\le k\le m$. If 
$P_1,\ldots,P_r$ are the components of $\overline{M-L}$, we may
index the components of $\partial L$ as $S_1,\ldots,S_k$, in such a
way that $S_i\subset\partial P_i$ for $i=1,\ldots,r$. 

For
$p=0,\ldots,r$, let $L_p$ denote the connected $3$-manifold obtained from the
disjoint union $L\discup P_1\discup\cdots\discup P_p$ by gluing
$S_i\subset L$ to $S_i\subset P_i$, for $i=1,\ldots,p$, via the
identity map. Then $h(L_0)=h(L)=m$. For $1\le p\le r$, it follows from
Lemma \ref{even easier} that $h(L_p)\le h(L_{p-1})+h(P_p)-1$. Hence 
$$h(L_r)\le h(L_0)+\sum_{i=1}^r h(P_i)-r= m-r+\sum_{i=1}^r h(P_i).$$
Up to homeomorphism, the manifold $M$ may be obtained from $L_r$ by
gluing the boundary components $S_{r+1},\ldots,S_k$ of $L\subset L_r$
to boundary components of $P_1\discup\cdots\discup P_r\subset
L_r$. Hence
$$h(M)\le h(L_r)+k-r\le2 (m-r)+\sum_{i=1}^r h(P_i).$$
Since $h(M)>m(m-1)$ by the hypothesis of the lemma, we obtain
$m^2-3m<-2r+\sum_{i=1}^r h(P_i)=\sum_{i=1}^r (h(P_i)-2)$, so that
$$\sum_{i=1}^r (h(P_i)-2)> m(m-3)\ge r(m-3).$$
Hence $h(P_{i_0})-2> m-3$ for some
$i_0\in\{1,\ldots,r\}$. If we set $P=P_{i_0}$ it follows that $h(P)\ge
m$. Furthermore, $P=P_{i_0}$ is by definition connected, and
$\overline{M-P}=L\cup\bigcup_{i\ne i_0}P_i$ is connected because each of
the connected submanifolds $P_i$ meets the connected submanifold
$L$. It remains only to show that $h(M-P)\ge m$.

Assume that $h(\overline{M-P})<m$. We again apply Lemma \ref{graphology},
this time taking $X=\overline{M-P}$. This gives a compact, connected submanifold $L'$ of $ M$ such that
$\overline{M-P}\subset \inter L'$,
each component of
$\partial L'$ is an incompressible  torus in $M$, and $h(L')=m$. Since
$L\subset \overline{M-P}$, we in
particular have $L\subset\inter L'$. In view of our choice of $L$,
it follows that $L'$ is a regular neighborhood of $L$. Since
$L\subset \overline{M-P}\subset \inter L'$, and since  $\partial\overline{M-P}$ is incompressible, it follows that
$L'$ is also a regular neighborhood of $\overline{M-P}$. 
But this is
impossible because $h(L')=m>h(\overline{M-P})$.
\EndProof
%T\cals I P K characteristic

\section{Finding useful solid tori}\label{seek 'n' find}

The main result of this section, Lemma \ref{uneeda}, is the first step
in the proof of Proposition E (Proposition \ref{lost corollary}) that
was sketched in the introduction.
The following result, Lemma \ref{uneeda this first}, is
required for the proof of Lemma \ref{uneeda}.

\Lemma\label{uneeda this first}
Let $p:\oldUpsilon\to\frakB$ be a covering map of compact, connected $2$-orbifolds. Suppose that $\oldUpsilon$ is orientable, that $\chi(|\oldUpsilon|)=0$, and that $|p|: |\oldUpsilon| \to |\frakB|$ (see \ref{orbifolds introduced})   is $\pi_1$-injective. Then $\chi(|\frakB|)=0$.
 Furthermore, the index in $\pi_1(|\frakB|)$ of the image of $|p|_\sharp:\pi_1(|\oldUpsilon|)\to\pi_1(|\frakB|)$ is at most the degree of the orbifold covering map $p$. Finally, if $|\oldUpsilon|$ is an annulus, if every finite subgroup of $\pi_1(\frakB)$ is cyclic, and if the degree of $p$ is at least $2$, then $\card\fraks_\oldUpsilon\ne1$. 
\EndLemma

%\redproofreadingnote{As the assertion $\fraks_\frakB\subset|\inter\frakB|$ has  been removed, so has the followingn passage:
%``Note that asserting $\fraks_\frakB\subset|\inter\frakB|$, as we do in the conclusion of Lemma \ref{uneeda this first}, is not the same as asserting
%$\fraks_\frakB\subset \inter|\frakB|$. Indeed, the latter inclusion need not hold under the hypotheses of the lemma; for example, the orientable $2$-manifold $\oldUpsilon:=\SSS^1\times[-1,1]$ is a two-sheeted covering of the $2$-orbifold $\frakB:=\SSS^1\times[[0,1]$, but $\fraks_\frakB$ is one of the boundary components of the annulus $|\frakB|$.''}

\Proof[Proof of Lemma \ref{uneeda this first}]
Let $d$ denote the degree of the orbifold covering map $p$.

In the case where $|\oldUpsilon|$ is a torus, the $\pi_1$-injectivity of $|p|$ implies that $\pi_1(|\frakB|)$ has a subgroup isomorphic to $\ZZ\times\ZZ$, and therefore that the $2$-manifold $|\frakB|$ is a torus or Klein bottle. Hence $\chi(|\frakB|)=0$. Since $|\frakB|$ is closed, we have $\dim\fraks_\frakB\le0$, so that $|p|$ is a branched covering map; in particular $|p|$ is a degree-$d$ map of closed manifolds, and hence $|\pi_1(|\frakB|):|p|_\sharp(|\pi_1(|\oldUpsilon|)|\le d$.

The rest of the proof will be devoted to the case in which
$|\oldUpsilon|$ is an annulus. Since $\frakB$ is connected, the
covering map $p$ is surjective, and hence restricts to a surjective
covering map from $\partial\oldUpsilon$ to $\partial\frakB$. Since
$\partial\frakB$ is a closed $1$-orbifold, each of its components
is homeomorphic to $\SSS^1$ or $[[0,1]]$. If $\oldGamma$ is any  component of
$\partial\frakB$ is homeomorphic to $[[
0,1]]$, there is a
%n $r$ maps some
component $\toldGamma$ of $\partial\oldUpsilon$ such that $p(\toldGamma)=\oldGamma$. Since
$\frakB$ is orientable, $|\toldGamma|$ is a component of
$\partial|\oldUpsilon|$, and is therefore $\pi_1$-injective in the
annulus $|\oldUpsilon|$. The $\pi_1$-injectivity of $|p|$ then
implies that $|p|\big||\toldGamma|:|\toldGamma|\to|\oldGamma|$ is
$\pi_1$-injective;  since $|\toldGamma|$ is a
$1$-sphere it follows that $|\oldGamma|$ cannot be an arc, and hence
$\oldGamma$ cannot be homeomorphic to $[[0,1]]$. This shows that every component of
$\partial\frakB$ is (orbifold)-homeomorphic to $\SSS^1$; equivalently, $|\partial\frakB|$ is a closed $1$-manifold.

Let $\fraks^1_\frakB$ denote the union of all $1$-dimensional components of $\fraks_\frakB$. According to \ref{orbifolds introduced},  we have $\partial|\frakB|=|\partial\frakB|\cup \fraks^1_\frakB$. In the present situation, since
$|\partial\frakB|$ is a closed  $1$-manifold, and
$\partial|\frakB|$ is of course a closed  $1$-manifold,
$\fraks^1_\frakB$ must also be a closed $1$-manifold, and 
%we may write $\partial\frakB$ \redcomment{I think this should be
%  $\partial|\frakB|$, and I have to decide whether the sentence reads
%  OK---maybe I should be more explicit about the fact that
  $\partial|\frakB|$
must be
%  $\fraks^1_\frakB$ is a closed $1$-manifold} as
 a disjoint union $|\partial\frakB|\discup \fraks^1_\frakB$.
% of closed topological $1$-manifolds. 

Since 
%$p$ is a surjective
%orbifold covering and 
$\oldUpsilon$ is orientable, we have
$|p|(\partial|\oldUpsilon|)=|p(\partial\oldUpsilon)|=|\partial\frakB|$. Since 
$|\partial\frakB|\cap\fraks^1_\frakB=\emptyset$, it follows that
$C:=
|p|^{-1}(\fraks^1_\frakB)\subset\inter|\oldUpsilon|$. On the other
hand, since $p$ is an orbifold covering map, the map of spaces $|p|\big| C:
C\to
\fraks^1_\frakB$ is locally surjective. As we have seen that
$\fraks^1_\frakB$ is a topological $1$-manifold, it follows that
$C$
is (the underlying space of) a graph with no endpoints or isolated
vertices,
%It follows
%that for each point $x$ of $C$ there exist a finite (cyclic or
%dihedral) subgroup $G$ of $\OO(2)$ containing a reflection, an
%orientation-preserving (cyclic and possibly trivial) subgroup $G'$ of
%$G$, orbifold neighborhoods $\frakU'$ and $\frakU$ of $x$ and $p(x)$
%in $\oldUpsilon$ and $\frakB$ respectively, a neighborhood $V$ of $0$
%in $\RR^2$, and orbifold homeomorphisms $h:V/G\to \frakU$ and
%$h':V/G'\to \frakU'$, such that $p(\frakU')=\frakU$, and $h^{-1}\circ
%p\circ h'$ is the natural quotient map from $V/G'$ to $V/G$. In
%particular, $\frakU'$ may be chosen so that $C\cap|\frakU'|$ is
%homeomorphic to the intersection of $\DD^2$ with a finite union of
%lines through the origin in $\RR^2$.  Since this holds for every
%$x\in C$, the set $C$ is homeomorphic to a finite graph in which
%every vertex has even, positive valence. In particular,  $C$ is a
%finite graph, without endpoints or isolated vertices, 
contained in
the interior of the annulus $|\oldUpsilon|$. 
This implies:
\Claim\label{either or} Either (a) every component of $C$ is a homotopically non-trivial curve in $\inter|\oldUpsilon|$, or (b) some component of $(\inter|\oldUpsilon|)-C$ is an open disk.
\EndClaim

Let us set $W=|\frakB|-\fraks^1_\frakB$ and $\tW=|p|^{-1}(W)=|\oldUpsilon|-C$. The definition of $\fraks^1_\frakB$ implies that $\omega(W)$ has only isolated singular points, and hence that $|p|\big|\tW$ is a branched covering map of degree $d$ from the (possibly disconnected) $2$-manifold $\tW$ to $W$. Since $\fraks^1_\frakB\subset\partial|\frakB|$, we have $W\supset \inter|\frakB|$. Hence: 
%$|\frakB|-\fraks^1_\frakB$ meets $\fraks_\frakB$ only in isolated points, and hence $p|\Upsilon-C:\Upsilon-C\to |\oldUpsilon|-\fraks^1_\frakB$
%Now since $\fraks^1_\frakB\subset\partial|\frakB|$, we have $W\supset\inter|\frakB|$.
%On the other hand, since $|\oldUpsilon|$ is an annulus, each component of the closed $1$-manifold $C\subset\inter|\oldUpsilon|$ separates $|\oldUpsilon|$; hence if $X$ is any component of $\tW$, the $2$-submanifolds $X$ and $\overline{X}$ of $\oldUpsilon$ have the same interior. Hence:
\Claim\label{Is this it?}
For every component $X$ of $\tW$, the restriction of $p$ to $\inter X$ is a branched covering map of degree at most $d$ from $\inter X$ to $\inter|\frakB|$ (and is surjective since the $2$-manifold $|\frakB|$ is connected). 
\EndClaim
If Alternative (b) of \ref{either or} holds, so that $\tW$ has a component $X$ such that $\inter X$ is an open disk, then \ref{Is this it?} implies that
%If some component of $C$ is homotopically trivial in $|\oldUpsilon|$, then $|\oldUpsilon|-C$ has some component $X_0$ such that $\overline{X_0}$ is a disk. Then by \ref{Is this it?}, 
$|p|$ restricts to a (surjective) branched covering map from $\inter{X}$ to $\inter|\frakB|$, and it follows that $\inter|\frakB|$ is an open disk. This is impossible, since $|\oldUpsilon|$ is an annulus and $|p|\big||\oldUpsilon|:|\oldUpsilon|\to|\frakB|$ is $\pi_1$-injective. Hence Alternative (a) of \ref{either or} holds, i.e.
\Claim\label{It's either, not or}
Every component of $C$ is a homotopically non-trivial curve in the interior of the annulus $|\oldUpsilon|$.
\EndClaim

%Hence every component of $C$ is homotopically non-trivial in the annulus $|\oldUpsilon|$. 
It follows from \ref{It's either, not or} that every component of
$\tW$ is the interior of an annulus in $|\oldUpsilon|$. If we fix a
component $X_0$ of $\tW$, then by \ref{Is this it?}, $|p|$ restricts
to a branched covering map of degree at most $d$ from the open annulus
$\inter{X_0}$ to $\inter|\frakB|$, and it follows that $|\frakB|$ is a
annulus, a M\"obius band or a disk; but again, the $\pi_1$-injectivity
of $|p|$ implies that $|\frakB|$ is not a disk. Thus we have shown
that $|\frakB|$ is an annulus or a M\"obius band, and so
$\chi(|\frakB|)=0$. Furthermore, since the branched covering map
$|p|\big|\inter X_0: \inter X_0\to\inter|\frakB|$ has degree at most
$d$, the image of $(|p|\big|\inter {X_0})_\sharp$ 
%\pi_1(\inter{X_0})$ 
has index at most $d$ in $\pi_1(\inter|\frakB|)$. In particular, the image of 
$|p|_\sharp$ has index at most $d$ in $\pi_1(|\frakB|)$.
%_1

To prove the final assertion, suppose that every finite subgroup of $\pi_1(\frakB)$ is cyclic. Under this additional hypothesis, we claim:
\Equation\label{with additional}
 \fraks_\oldUpsilon\cap C=\emptyset.
\EndEquation
To prove (\ref{with additional}), suppose that $x$ is a point of $ \fraks_\oldUpsilon\cap C$. Set $y=p(x)$.
%Let $V$ be a neighborhood of $p(x)$ in $\frakB$ which is orbifold-homeomorphic to $\DD^2/G$ for some finite subgroup $G$ of $\OO(2)$. 
Since $x\in C$ we have $y\in \fraks^1_\frakB$, so that $G_y$ contains a reflection (where the notation $G_y$ is defined by \ref{orbifolds introduced}). But since $x\in\fraks_\oldUpsilon$ and $\oldUpsilon$ is orientable, $G_y$ must also contain a rotation. Hence $G_y$ is non-cyclic. But \ref{2-dim case} gives $\chi(\frakB)\le\chi(|\frakB|)=0$, which implies that $\frakB$ is a very good orbifold, and hence that $G_y$ is isomorphic to a subgroup of $\pi_1(\frakB)$ (see \ref{orbifolds introduced}).
This contradicts the condition that every finite subgroup of $\pi_1(\frakB)$ is cyclic, and (\ref{with additional}) is proved. 

Now assume further that $d\ge2$ and that $\card\fraks_\oldUpsilon=1$. In particular $\fraks_\oldUpsilon\ne\emptyset$. Since $\oldUpsilon$ is orientable we have $\fraks_\oldUpsilon\subset\inter|\oldUpsilon|=|\inter\oldUpsilon|$, which with (\ref{with additional}) gives $\emptyset\ne\fraks_\oldUpsilon\subset|\inter\oldUpsilon|-C$. %and so Hence $\emptyset\ne p(\fraks_\oldUpsilon)\subset \frakB-\fraks^1_\oldUpsilon$. 
Since $p$ is a covering map we have
$p(\inter\oldUpsilon)\subset\inter\frakB$ and
$p(\fraks_\oldUpsilon)\subset\fraks_\frakB$. Recalling that
$C=p^{-1}(\fraks^1_\frakB)$, we deduce that $\emptyset\ne
p(\fraks_\oldUpsilon)\subset\fraks_\frakB\cap(
|\inter\frakB|-\fraks^1_\frakB)$. But since we have seen that
$\fraks^1_\frakB$ is a closed $1$-manifold, we have
$|\inter\frakB|-\fraks^1_\frakB=\inter|\frakB|$. Hence
%$\fraks_\frakB\cap(
%|\inter\frakB|-\fraks^1_\frakB)=\fraks_\frakB\cap\inter|\frakB|$, and therefore
$\emptyset\ne
p(\fraks_\oldUpsilon)\subset\fraks_\frakB\cap\inter|\frakB|$. In particular,
\Equation\label{oykhie-doykh}
\fraks_\frakB\cap\inter|\frakB|\ne\emptyset.
\EndEquation
We have observed that $W\supset\inter|\frakB|$, so that $\inter|\frakB|=\inter W$, and that $\omega(W)$ has only isolated singular points.
%, so that $\obd(\inter(|\frakB|))$ is orientable. \redproofreadingnote{Orientability doesn't really seem to follow. But isolated singularities may be all that is needed to apply the formula.} 
 %Hence $p|\obd(\inter \tW):\obd(\inter\tW)\to\obd(\inter W)$ is a covering map of orientable $2$-orbifolds of finite type. 
Set $n=\card\fraks_{\omega(\inter W)}$; by (\ref{oykhie-doykh}) we
have $n\ge1$. If $\fraks_{\omega(\inter W)}=\{y_1,\ldots,y_n\}$, and
if $k_i$ denotes the order of $y_i$, then by \ref{2-dim case} we have
$\chibar(\omega(\inter W))=\chibar(\inter W)+\sum_{i=1}^n(1-1/k_i)$;
since we have seen that $\chi(|\frakB|)=0$, it follows that
$\chibar(\omega(\inter W))=\sum_{i=1}^n(1-1/k_i)\ge n/2\ge1/2$. Since
$p|\obd(\inter \tW):\obd(\inter\tW)\to\obd(\inter W)$ is a degree-$d$
(orbifold) covering, and $d\ge2$, we have
$\chibar(\obd(\inter\tW))=d\cdot\chibar(\obd(\inter W))\ge1$. But we
have assumed that $\card\fraks_\oldUpsilon=1$; and by (\ref{with
  additional}), and the orientability of $\oldUpsilon$, we have  $\fraks_\oldUpsilon\subset\inter\tW$.
Thus
$\card\fraks_{\obd(\inter\tW)}=1$. If $\ell$ denotes the order of the unique point of $\fraks_{\inter\obd
(\tW)}$, then \ref{2-dim case} gives $\chibar(\obd(\inter\tW))=\chibar(\inter \tW)+(1-1/\ell)$.  Since each component of $\tW$ is a half-open annulus by \ref{It's either, not or}, we obtain $\chibar(\obd(\inter \tW))=1-1/\ell<1$, a contradiction. This establishes the final assertion.
\EndProof
%\tW\cala\obd $r

\Lemma\label{uneeda}
Let $\frakK
$ be a strongly \simple, boundary-irreducible, orientable
$3$-orbifold. Set $K=|\frakK|$. Suppose that
$K$ is boundary-irreducible and $+$-irreducible (see
Definition \ref{P-stuff}), and that $K^+$ is not homeomorphic to
$\TTT^2\times[0,1]$ or to a twisted $I$-bundle over a Klein bottle (where $K^+$ is defined as in \ref{P-stuff}). Suppose that $X$ is 
a $\pi_1$-injective, connected subsurface of $\partial K^+$, with $X\cap\fraks_\frakK\ne\emptyset$ and  $\chi(X)=0$,  
%such that the inclusion homomorphism $\pi_1(X)\to\pi_1(K)$
%-injective in $\partial K$ 
%is non-trivial, \redcomment{Why doesn't it say $\pi_1$-injective? That's what I seem to have in the application}
%  original version, which can be found in book35.tex, made no
  %sense. This version matches the use of the hypothesis in the first
  %paragraph of the proof. Check the rest of the proof and the
  %apps.  } 
and that there exist
disks $G_1,\ldots,G_n\subset\inter X$, where $n\ge0$, such that
$\overline{X-(G_1\cup\ldots\cup G_n)}$ is a component of $|\oldPhi(\frakK)|$.
Then $X$ is an annulus, and there is a
$\pi_1$-injective solid torus $J\subset K^+$, 
%\redmissingref{It said ``fix all cross-refs!!'' What does that mean?} 
with $\partial J\subset K\subset K^+$, such that one of the following alternatives holds:
\begin{itemize}
\item We have $\partial J\cap\partial K=X\discup X'$ for some annulus $X'\subset\partial K^+ $; furthermore, each of the annuli $X$ and $X'$ has winding
  number $1$ in $J$ (see \ref{great day}) and has non-empty intersection with $\fraks_\frakK$, and $
\partial  J\cap\fraks_\frakK\subset \inter X\cup \inter X'$.
\item 
We have $\partial J\cap\partial K=X$; furthermore, $X$ has winding
  number $1$ or $2$ in $J$, and we have
  $\wt (\partial J) \ge\lambda_\frakK$ (see \ref{lambda thing})  and $\partial J\cap\fraks_\frakK\subset\inter X $. 
\end{itemize}
\EndLemma

\Proof
Since $X$ is $\pi_1$-injective in $\partial K^+$, which is a union of components of $\partial K$, and since $K$ is boundary-irreducible, $X$ is $\pi_1$-injective in $K$.

For $i=1,\ldots,n$, the simple closed curve $\partial\obd( G_i)$ is a
component of $\partial\oldPhi(\frakK)$ and is therefore
$\pi_1$-injective in $\frakK$ by \ref{tuesa day}. Hence the orientable
$2$-orbifold $\obd(G_i)$ is not discal. Since $G_i$ is a disk, it
follows that $\wt G_i\ge2$. We have $\wt X\ge\sum_{i=1}^n\wt G_i$, and hence
\Equation\label{read this one first}
\wt X\ge2n.
\EndEquation

Let  $X_0$ denote the component $\overline{X-(G_1\cup\ldots\cup G_n)}$ of $|\oldPhi(\frakK)|$.
Since $X$ is $\pi_1$-injective in $K$, since the surface $X$ has Euler characteristic $0$ and is therefore non-simply connected, 
%the
%inclusion homomorphism $\pi_1(X)\to\pi_1(K)$ is non-trivial, 
and since the inclusion homomorphism $\pi_1(X_0)\to\pi_1(X)$ is surjective, the inclusion homomorphism $\pi_1(X_0)\to\pi_1(K)$ is non-trivial.

%The non-triviality of the
%inclusion homomorphism $\pi_1(X)\to\pi_1(K)$ also implies that $\pi_1(X)\ne\{1\}$. Since in additional  $X$ is orientable, we have $\chi(X)\le0$. \redcomment{Both those facts now follow from the assumption $\chi(X)=0$, which may not have been there at an earlier stage. Rewrite accordingly.} 
Since $\chi(X)=0$ and $X\cap\fraks_\frakK\ne\emptyset$, it follows from \ref{2-dim case} that %so that $\omega(X)$ has at least one singular point, it now follows that
 $\chi(\omega(X))<0$. Hence if $n=0$, so that $X_0=X$, we have
 $\chi(\omega(X_0))=\chi(\omega(X))<0$. If $n>0$ then \ref{2-dim case} gives
 $\chi(\omega(X_0))\le\chi(X_0)=\chi(X)-n=-n<0$. Thus in any event we have
%In particular, we have
 $\chi(\obd(X_0))<0$. If $\oldLambda_0$ denotes the component of
$\oldSigma(\frakK)$ containing $\obd(X_0)$ (so that $\oldLambda_0$ is
by definition an \Ssuborbifold\ of $\frakK$), then by \ref{tuesa day} we have 
$\chi(\oldLambda_0)=\chi(\oldLambda_0\cap\partial\frakK)/2$; every component of  $\oldLambda_0\cap\partial\frakK$ is a component of $\oldPhi(\frakK)$ and therefore has non-positive Euler characteristic by \ref{tuesa day}. Hence $\chi(\oldLambda_0)\le\chi(\obd(X_0))/2<0$. In particular $\oldLambda_0$ is not a \torifold, so that by Lemma \ref{when a tore a fold}, $\oldLambda_0$ is not a \bindinglike\ \Ssuborbifold\ of $\frakK$. It follows that $\oldLambda_0$ is a \pagelike\ \Ssuborbifold\ of
%can be equipped with an $I$-fibration over a $2$-orbifold $\frakB$ in such a way that it is %standardly embedded in
 $\frakK$. 

We may therefore fix an $I$-fibration $q_0:\oldLambda_0\to\frakB_0$, where $\frakB_0$ is some compact, connected $2$-orbifold, such that
$\calf_0:=|\partialh\oldLambda_0|=|\oldLambda_0|\cap\partial K$. Note that $X_0$ is a component of $\calf_0$. 
 In the case where
the \pagelike\
 \Ssuborbifold\ $\oldLambda_0$ is twisted, 
$q_0$ is a non-trivial fibration, so that 
 $\calf_0$ is connected; hence we have $\calf_0=X_0$ in this case.
 In the case where
the \pagelike\
 \Ssuborbifold\ $\oldLambda_0$ is untwisted, 
$q_0$ is a trivial fibration, so that
% which has at most two components (see \ref{fibered stuff}). In the
 %case where
  $\calf_0$ has exactly two components; in this case the component $\calf_0$ distinct from
$X_0$ will be denoted $X_0'$. 
In this latter case the inclusion map from $X'_0$ to $K$ is homotopic in $K$ to a homeomorphism of $X_0'$ onto $X_0$; since the inclusion homomorphism
 $\pi_1(X_0)\to\pi_1(K)$ is non-trivial, it then follows that the
 inclusion homomorphism $\pi_1(X_0')\to\pi_1(K)$ is also non-trivial in this case.

% so that the images (defined up
 %to conjugacy) of the.$X_0'$ has boundary components that are joined
 %to the components of $\partial X$ by annulus components of $\Fr
 %\oldLambda_0$, and are therefore homotopically non-trivial in $K$. 
Hence in all cases:
\Claim\label{all cases cool}
For each component $F$ of $\calf_0$, the inclusion homomorphism $\pi_1(F)\to\pi_1(K)$ is non-trivial.
\EndClaim

Let us write $|\Fr \oldLambda_0|=\cala_0\discup\cala_1$, where
$\cala_0$ (respectively $\cala_1$) denotes the union of all components
$A$ of $|\Fr \oldLambda_0|$ such that the inclusion homomorphism
$\pi_1(A)\to\pi_1(K)$ is trivial (respectively,
non-trivial). Since
each component of $\Fr \oldLambda_0\subset\frakA(\frakK)$ is an orientable annular orbifold (see 
\ref{tuesa day}), each component of $\cala_i$ is an annulus or a
disk for $i=0,1$, and any annulus component of $\cala_i$ has weight $0$. It is obvious from the definition that $\cala_1$ has no disk components, and hence:
\Claim\label{lame claim}
Each component of $\cala_1$ is a weight-$0$ annulus. 
\EndClaim

If $A$ is any component of $\cala_1$, then $A$ is a component of
$|\partialv\oldLambda_0|$, and hence has at least one boundary
component contained in each component of
$|\partialh\oldLambda_0|=\calf_0$. In particular, at least one
component of $\partial A$ is  contained in $X_0$, and is therefore a
component of $\partial X_0$; we choose such a component and call it $C$. Since $A$ is an annulus by \ref{lame claim}, the definition of $\cala_1$ implies that  $C$ is homotopically non-trivial in $K$. This implies that $C$ cannot have the form $\partial G_i$ with $1\le i\le n$; hence $C$ is a component of $\partial X$. Since $C$ is  homotopically non-trivial in $K$, it is in particular homotopically non-trivial in the $2$-manifold $X$, and is therefore $\pi_1$-injective in $X$. Since  $X$ is $\pi_1$-injective in  $K$, it follows that $C$ is $\pi_1$-injective in $K$. As $C$ is a component of $\partial A$, the annulus $A$ is $\pi_1$-injective in $K$. This shows:
\Claim\label{omnibust}
 $\cala_1$ is $\pi_1$-injective in $K$.
\EndClaim

Since $\partial G_1,\ldots,\partial G_n$ are components of $\partial
X_0\subset\partial\calf_0$ that bound disks in $\partial K$, they are
components of $\partial\cala_0$. On the other hand, since $X$ is a
compact, orientable $2$-manifold with $\chi(X)=0$, any component of
$\partial X$ carries $\pi_1(X)$; since $X$ is non-simply connected, and is $\pi_1$-injective in $K$,
% is
%non-trivial by hypothesis, 
it follows that any component of $\partial
X$ is a component of $\partial\cala_1$. 
%\redcomment{In the argument I
%  just removed here, which can be found in book35.tex,
  %boundary-irreducibility is used. This raises the question of whether
  %I need boundary-irreducbility in the hypothesis.} 
Thus:

\Claim\label{mobsters}
We have $\partial X_0\cap\cala_0=\partial G_1\cup\cdots\cup\partial G_n$, and
$\partial X_0\cap\cala_1=\partial X$.
\EndClaim

(We will show below, in \ref{and an annulus too!}, that $\partial
X\ne\emptyset$, which by \ref{mobsters} implies that $\cala_1\ne\emptyset$;
however, at this stage of the proof we allow the possibility that
$\cala_1$ is empty.)

If $C$ is any component of $\partial\cala_0$, then $C$ is a simple
closed curve in $\partial\calf_0$ which is homotopically trivial in
$K$. The boundary-irreducibility of $K$ then implies that $C$ bounds a
disk in $\partial K$. A disk bounded by $C$ cannot contain a
component $F$ of $\calf_0$, since the inclusion homomorphism
$\pi_1(F)\to\pi_1(K)$ is non-trivial by \ref{all cases cool}. Hence the disk $D_C$ bounded by $C$ is unique, and satisfies $D_C\cap\calf_0=C$. Thus $D_C$ is a component of $\partial K-\inter\calf_0$.

According to \ref{tuesa day},
$\omega(\cala_0)\subset\frakA(\frakK)$ is an  essential
annular suborbifold of $\frakK$, and hence every component
$C$ of $\omega(\partial\cala_0)$ is $\pi_1$-injective in
$\frakK$. This shows:
\Claim\label{hit it}
For every component  $C$ of $\partial\cala_0$ we have $D_C\cap\fraks_\frakK\ne\emptyset$.
\EndClaim

For each component $A$ of $\cala_0$, we set $S_A=A\cup\bigcup _{C\in\calc(\partial A)}D_C$. If $A$ is an annulus, then $\partial A$ has two components
$C$ and $C'$; the disks $D_C$ and $D_{C'}$ are components of $\partial K-\inter\calf_0$ with distinct boundaries, and are therefore disjoint.
If $A$ is a disk, then of course $\partial A$ has only one component. Hence in any case,
$S_A$ is a $2$-sphere.
Since $K^+$ is irreducible, and since it is immediate from the hypothesis that $\partial K^+\ne\emptyset$, the sphere $S_A$ bounds a unique ball $E_A\subset K^+$ for each component $A$ of $\cala_0$, and we have $\Fr_{K^+}E_A=A$. 

%\redmissingref{A \%-ed out passage about $C$ and $C'$ has been removed here. It can be found in first-move.tex.}

Since $X_0\subset|\oldLambda_0|$, and since the inclusion homomorphism $\pi_1(X_0)\to\pi_1(K)$ is non-trivial while the inclusion homomorphism $\pi_1(K)\to\pi_1(K^+)$ is an isomorphism,  the inclusion homomorphism $\pi_1(|\oldLambda_0|)\to\pi_1(K^+)$ is non-trivial. Hence if $A$ is a component of $\cala_0$, the ball $E_A$  cannot contain $|\oldLambda_0|$. We must therefore have 
\Equation\label{lamb}
E_A\cap|\oldLambda_0|=A.
\EndEquation
 Thus $E_A$ is a component of $\overline{K^+-|\oldLambda_0|}$.

If $A$ and $A'$ are distinct components of $\cala_0$, then $E_A$ and
$E_{A'}$ are components of $\overline{K^+-|\oldLambda_0|}$, and are
distinct because their frontiers $A$ and $A'$ are distinct. Hence:
\Claim\label{piggy-wig}
 $(E_A)_{A\in\calc(\cala_0)}$ is a disjoint family.
\EndClaim
 We set $J=|\oldLambda_0|\cup\bigcup_{A\in\calc(\cala_0)}E_A\subset K^+$. It follows from (\ref{lamb}), \ref{piggy-wig} and the definitions that
\Equation\label{where's the beef?}
\Fr\nolimits_{K^+}J=\cala_1.
\EndEquation

It follows from (\ref{where's the beef?}) that
$\partial J=\cala_1\cup(J\cap\partial{K^+})$. Since $\cala_1$ is properly embedded in $K$, it now follows that
\Equation\label{sundried}
\partial J \cap\partial K\subset\partial K^+
\EndEquation
and that
\Equation\label{two for a dime}
\partial J\subset K.
\EndEquation

%It follows from \ref{omnibust} that  $\cala_1$ is $\pi_1$-injective in
Since $K^+$ is an irreducible
$3$-manifold  according to the hypothesis, every $2$-sphere in $\inter
J$ bounds a $3$-ball  $B\subset\inter K^+$; no component of $\cala_1$
can be contained in $B$,
as each component of $\cala_1$ meets $\partial
K^+$. In view of \ref{where's the beef?} it follows that $B\subset
K^+$. This shows:

\Claim\label{omni or bust}
$J$ is irreducible.
\EndClaim

We set $\calf=J\cap\partial K^+$. It follows from the definitions of
$J$ and of the $E_A$ that
\Equation\label{no end}
\calf=\calf_0\cup\bigcup_{C\in\calc(\partial\cala_0)}D_C.
\EndEquation

Now $X_0$ is a component of $\calf_0$, and according to \ref{mobsters},  $G_1,\ldots,G_n$ are the only disks that have the form $D_C$ for $C\in\calc(\partial\cala_0)$ and meet $X_0$. Since $X=X_0\cup\bigcup_{i=1}^nG_i$, it follows from (\ref{no end}) that
\Claim\label{change lobsters}
$X$ is a component of $\calf$.
\EndClaim

The next step involves two orientable $3$-orbifolds $\frakZ$ and $\frakZ'$, equipped
with $I$-fibrations 
%$r:\frakZ\to\oldDelta$
%and $r':\frakZ'\to\oldDelta'$ 
over $2$-orbifolds,  which are
constructed as follows. We
define $\frakZ$ to be the manifold $\DD^2\times[0,1]$, set
$\oldDelta=\DD^2$, and define an $I$-fibration $r:\frakZ\to\oldDelta$ to the projection to the first
factor. We define $\frakZ'$ to be the quotient of $\frakZ$ by the
involution $\tau:(z,t)\mapsto(c(z),1-t)$, where $c$ is the involution
$(x,y)\mapsto(x,-y)$ of $\DD^2$, and set $\Delta'=\DD^2/c$. We have $r\circ\tau=c\circ r$, and
hence induces a  map $r':\frakZ'\to\oldDelta'$, which is also an
$I$-fibration.
%defined by
%complex conjugation in $\DD^2$. Since the involution maps fibers to
%fibers, 
%$\frakZ'$ inherits an orbifold fibration from $\frakZ$; its
%base, which we denote by $\oldDelta'$, is the quotient of $\DD^2$ by the
%involution $z\mapsto\overline{z}$. 
Note that $|\oldDelta'|$ is a disk,
and that $\partial|\oldDelta'|$ is the union of the two arcs
$|\partial\oldDelta'|$ and $\fraks_\oldDelta$, which meet only in
their endpoints. Note also that $\partialv\frakZ$ is an annulus,
while $\partialv\frakZ'$ is an orbifold having two singular points of
order $2$, and a disk as underlying surface. 
%orbifold

For each component $A$ of $\cala_0$,
define $\frakZ_A$ to be a homeomorphic copy of $\frakZ$ if $A$ is an
annulus, and define  $\frakZ_A$ to be a homeomorphic copy of $\frakZ'$
if $A$ is a disk and $\card\fraks_{\omega(A)}=2$. Thus each $\frakZ_A$ has a fibration $r_A:\frakZ_A\to\oldDelta_A$, where
$\oldDelta_A$ is a $2$-orbifold homeomorphic to either $\oldDelta$ or
$\oldDelta'$, such that under suitable homeomorphic identifications of $\frakZ_A$ and $\oldDelta_A$ with $\frakZ$ and $\oldDelta$ or with $\frakZ'$ and $\oldDelta'$, the fibration $r_A$ is identified with $r$ or $r'$ respectively. 

The fibration of any given orientable annular $2$-orbifold is unique
up to fiber-preserving (orbifold) homeomorphism. For each component
$A$ of $\cala_0$, since $\Fr (\frakZ_A)$ and $\omega(A)$ are
homeomorphic by construction, it follows that there exist homeomorphisms
$\theta_A:\partialv\frakZ_A\to\omega(A)$ and
 $h_A:\partial\oldDelta_A\to q_0(\omega(A)) $
such that 
$q_0\circ\theta_A=h_A\circ (r_A|\partialv\frakZ_A)$. (It is worth bearing in mind that the homeomorphic $1$-manifolds  $|\partial\oldDelta_A|$ and $|\omega(q_0(A))| $ are arcs when $A$ is a disk, and are $1$-spheres when $A$ is an annulus.) Let $\oldLambda$ denote the connected $3$-orbifold obtained from the disjoint union of $\oldLambda_0$ with the $\frakZ_A$, where $A$ ranges over the components of $\cala_0$, by gluing $\partialv\frakZ_A$ to $\omega(A)\subset\oldLambda_0$ via $\theta_A$ for every $A$. Let $\frakB$ denote the $2$-orbifold obtained from the disjoint union of $\frakB_0$ with the $\oldDelta_A$, where $A$ ranges over the components of $\cala_0$, by gluing $\partial\oldDelta_A$ to $\omega(q_0(A))\subset\frakB_0$ via $h_A$ for every $A$. Then there is a well-defined $I$-fibration $q:\oldLambda\to\frakB$ which restricts to $q_0$ on $\oldLambda_0$, and to $r_A$ on $\frakZ_A$ for each $A$.

For each component $A$ of $\cala_0$, both $|\frakZ_A|$ and $E_A$ are $3$-balls, and $|\partialv\frakZ_A|\subset\partial |\frakZ_A|$ and $A\subset\partial E_A$ are either both annuli or both disks. Hence, the homeomorphism 
$|\theta_A|:|\partialv\frakZ_A|\to A$ 
%defined by
%$\theta_A:\partialv(\frakZ_A)\to\omega(A)$
may be extended to a homeomorphism from $|\frakZ_A|$ to $E_A$; we fix such an extension $t_A:|\frakZ_A|\to E_A$ for each $A$, and define a map $T:|\oldLambda|\to J$ to be the identity on $|\oldLambda_0| $ and to restrict to $t_A$ on $\frakZ_A$ for each $A$. Since
$(E_A)_{A\in\calc(\cala_0)}$ is a disjoint family by \ref{piggy-wig}, and since
$E_A\cap|\oldLambda_0|=A$ for each $A$ by (\ref{lamb}), the map $T:|\oldLambda|\to J$ is a homeomorphism. 

The orientability of $\frakK$ implies that $K=|\frakK|$ is orientable,
and hence that $K^+$ is orientable. In particular $J\subset K^+$ is
orientable. Since $T$ is a homeomorphism, the $3$-manifold
$|\oldLambda|$ is orientable. But the orientability of $\frakK$, $\frakZ$ and $\frakZ'$ 
implies that every component of $\fraks_\frakK$, and every component
of $\fraks_{\frakZ_A}$ for every $A\in\calc(\cala_0)$, is
one-dimensional; hence every component of $\fraks_\oldLambda$ is
one-dimensional. With the orientability of $|\oldLambda|$, this
implies that $\oldLambda$ is orientable.

It follows from (\ref{no end}) and the definition of $T$ that
\Equation\label{still no end}
T(|\partialh\oldLambda|)=\calf.
\EndEquation
From \ref{change lobsters} and (\ref{still no end}) it follows that
$T^{-1}(X))$ is a component of $|\partialh\oldLambda|$, so that
%we have
%$T(|\partialh\oldLambda|)=\calf$.
% $X=X_0\cup\bigcup_{i=1}^nG_i$ is a component of $T(\partialh\oldLambda)$. \redmissingref{This is superseded by a corresponding statement about $X$ being a component of $\calf$.}
%Hence:
\Claim\label{ypsilanti}
$\oldUpsilon:=\omega(T^{-1}(X))$ is a component of $\partialh\oldLambda$. 
\EndClaim

According to \ref{fibered stuff},
$q|\partialh\oldLambda:\partialh\oldLambda\to\frakB$ is a degree-$2$
orbifold covering. In particular, by \ref{ypsilanti},
$q|\oldUpsilon:\oldUpsilon\to\frakB$ is an orbifold covering of degree
at most $2$. On the other hand, the homeomorphism $T:|\oldLambda|\to
J$ maps $|\oldUpsilon|$ onto the surface $X$.  By hypothesis $X$ is
$\pi_1$-injective in $ K$ 
%\redcomment{That's the hypothesis now. I guess I should check the app. one last time, oy. See the \%-ed out stuff below for some crap I would have to do if I didn't have $\pi_1$-injectivity}
%  homomorphism is non-trivial, which implies $\pi_1$-injectivity if
  %$X$ is an annulus. (Well, I guess. Dealing with torsion could be
  %messy, and I'm not sure it's the best way.) 
%If $X$ is a torus, it
  %follows from boundary-irreducibility. This should also have been
  %point out at the beginning of the proof. Also, the organization of
  %this paragraph is confusing.}
% and hence in $K$ by
%boundary-irreducibility, $X$ is $\pi_1$-injective in $K$ 
and hence in
$K^+$, and in particular in $J$. 
%\redcomment{That last sentence really
%  makes no sense at all. The hypothesis says that the inclusion
%hom. $\pi_1(X)\to\pi_1(K)$ is non-trivial, not injective. I think the
%solution will be to weaken the hypothesis of Lemma \ref{uneeda this
 % first} to require only non-triviality, not injectivity, of the
%induced map on $\pi_1$. The proof of Lemma \ref{uneeda this first}
%will be unaffected in the case where the domain is an annulus. In the
%case where it's a torus, a different argument is needed. I'm thinking
%na\"ively that the point may simply be that the range of a branched
%covering of surfaces has Euler characteristic $\ge$ that of the
%domain; if the range had $\chi>0$ it would contradict the
%non-triviality of the homomorphism. A little care is required if the
%range is $RP^2$: in that case the orientability of the domain
%contradicts non-triviality of the homomorphism. Fix up this paragraph
%after Lemma \ref{uneeda this first} is revised.}
%I guess, 
%which I think doesn't happen here? Do we know the image is infinite,
%not just non-trivial?}
The hypothesis also gives
$\chi(X)=0$. Hence $|\oldUpsilon|$ is $\pi_1$-injective in
$|\oldLambda|$, and $\chi(|\oldUpsilon|)=0$. According to Lemma \ref{affect}, $|q|:|\oldLambda|\to|\frakB|$ is a
homotopy equivalence, and hence
$|q\big|\oldUpsilon|=|q|\big||\oldUpsilon|$ is
$\pi_1$-injective. As we have shown that $\oldLambda$ is orientable,
$\oldUpsilon\subset\partial\oldLambda$ is orientable. It
therefore follows from Lemma \ref{uneeda this first},  applied
with $q|\oldUpsilon$ playing the role of $p$, that $\chi(|\frakB|)=0$. Hence  $|\frakB|$  is a torus, a Klein bottle, an annulus
or a M\"obius band. 
Since $|q|:|\oldLambda|\to|\frakB|$ is a
homotopy equivalence, and $T:|\oldLambda|\to
J$ is a homeomorphism, $J$ is homotopy-equivalent to $|\frakB|$. 

First suppose that $|\frakB|$ is a torus or Klein
bottle.  Thus $J$ is homotopy-equivalent to a torus or Klein
bottle. Since the compact, orientable $3$-manifold $J$ is irreducible by \ref{omni or bust}, it follows from \cite[Theorem 10.6]{hempel} that $J$ is 
homeomorphic to $\TTT^2\times[0,1]$ or to a twisted $I$-bundle over a
Klein bottle.  On the other hand, since $|\frakB|$ is closed, $\frakB$ must also be closed and $\fraks_\frakB$ must be $0$-dimensional. 
Since
$\partialh\oldLambda$ is a covering space of $\frakB$, it follows
that $|\partialh\oldLambda|$ is a branched cover of $|\frakB|$, and
is therefore closed. Since the components of $\partialv\oldLambda$
are annular orbifolds, and are therefore not closed, it follows that
$\partialv \oldLambda=\emptyset$,
i.e. $\partial|\oldLambda|=|\partialh\oldLambda|$. Since (\ref{still no
  end}), with the definition of $\calf$, implies that 
$T(|\partialh\oldLambda|)\subset\partial {K^+}$, it now follows that $\partial J=T(\partial|\oldLambda|)\subset\partial {K^+}$, which with the connectedness of ${K^+}$ implies that $J={K^+}$. But then ${K^+}$ is homeomorphic to $\TTT^2\times[0,1]$ or to a twisted $I$-bundle over a Klein bottle, a contradiction to the hypothesis. 

Hence $|\frakB|$ is either an annulus or a M\"obius band. Since   $J$ is homotopy-equivalent to $|\frakB|$, it follows that $\pi_1(J)$ is cyclic. Since the compact, orientable $3$-manifold $J$ is irreducible by \ref{omni or bust}, it then follows from \cite[Theorem 5.2]{hempel} that
\Claim\label{i'll say it's solid}
$J$ is a solid torus.
\EndClaim

Since we have seen that $X$ is $\pi_1$-injective in $K $, it is $\pi_1$-injective in $J$. Furthermore, the hypothesis gives $\chi(X)=0$. In view of \ref{i'll say it's solid}, it follows that
\Claim\label{and an annulus too!}
$X$ is an annulus.
\EndClaim

Since $\pi_1(X)$ and $\pi_1(J)$ are infinite cyclic by \ref{i'll say it's solid} and \ref{and an annulus too!}, and $X$ is $\pi_1$-injective in $K$, it follows that

\Claim\label{stubbins}
$J$ is $\pi_1$-injective in $K$.
\EndClaim

Note that the properties of $X$ and $J$ stated in the lemma before the alternatives given in the two bullet points are covered by (\ref{two for a dime}), \ref{i'll say it's solid}, \ref{and an annulus too!} and \ref{stubbins}.

We have seen that
$q|\partialh\oldLambda:\partialh\oldLambda\to\frakB$ is a degree-$2$ orbifold covering. Furthermore, by \ref{ypsilanti}, $\oldUpsilon$ is one component of $\partialh\oldLambda$. Hence either
(a) $\partialh\oldLambda=\oldUpsilon$, and $q|\oldUpsilon$ is a $2$-sheeted orbifold covering map onto $\frakB$; or
(b) $\partialh\oldLambda$ has two components, $\oldUpsilon$ and a second component $\oldUpsilon'$, and $q|\oldUpsilon$ and $q|\oldUpsilon'$ are orbifold homeomorphisms onto $\frakB$.

Suppose that (a) holds. In this case we will show that the solid torus $J$ and the annulus $X$ satisfy the second alternative in the conclusion of the lemma. 
We have $\partial J\cap\partial K=J\cap\partial
K^+=\calf=T(|\partialh\oldLambda|)$, by (\ref{sundried}) and
(\ref{still no end}) (and the definition of $\calf$), while $\partialh\oldLambda=\oldUpsilon$ by (a). 
Hence $\partial J\cap\partial K =T(|\oldUpsilon|)$. But by definition (see \ref{ypsilanti}) we have $\oldUpsilon=\omega(T^{-1}(X))$, so that $X=T(|\oldUpsilon|)$. Thus $\partial J\cap\partial K=X$. 

Next note that $\partial J=(J\cap\partial K^+)\cup\Fr_{K^+}J=X\cup\cala_1$, in view of \ref{where's the beef?}. Since 
$\cala_1\cap\fraks_\frakK=\emptyset$ by \ref {lame claim}, it follows that 
$\partial J\cap\fraks_\frakK\subset \inter X$. 

%It follows that $\wt (\partial J)=\wt X=\card\fraks_{\omega(X)}$. Since (a) holds, the 
To show that $\wt (\partial J) \ge\lambda_\frakK$, or equivalently that $\wt X\ge\lambda_\frakK$, we distinguish several subcases. By hypothesis we have $X\cap\fraks_\frakK\ne\emptyset$, i.e. $\wt X\ge1$; hence the assertion is true if $\lambda_\frakK=1$. If $\lambda_\frakK=2$ and $n\ge1$, then by \ref{read this one first}
we have $\wt X\ge2n\ge2=\lambda_\frakK$. There remains the subcase in which $\lambda_\frakK=2$ and $n=0$. In this subcase, we must show $\wt X\ge2$, and since $\wt X\ge1$, it is enough to show that $\wt X\ne1$. Since $n=0$, it follows from \ref{mobsters} that 
$\partial X_0=\partial X$, and since $X_0\subset X$ it then follows
that $X=X_0$. It also follows from \ref{mobsters} that $\partial
X_0\cap\cala_0=\emptyset$ in this subcase; since 
each component of $\cala_0=\partialv\oldLambda_0$ is saturated in $\oldLambda_0$, and therefore meets
$X_0=\partialh\oldLambda_0$, it follows that $\cala_0=\emptyset$, so
that $\oldLambda=\oldLambda_0$ and $\frakB=\frakB_0$.    Since $T$ is
defined to be the identity on $|\oldLambda_0|$, it is the identity on
its entire domain $|\oldLambda|$. Hence $\obd(X)=
\obd(T^{-1}(X))=\oldUpsilon$ (cf. \ref{ypsilanti}). 
%Hence $\wt X=\card\fraks_\oldUpsilon$. 
To show that $\wt X\ne1$, i.e. that $\card\fraks_\oldUpsilon\ne1$, we
will apply the last sentence of  Lemma \ref{uneeda this first}, with
$q|\oldUpsilon$ playing the role of $p$.  We have already seen that
the general hypotheses of Lemma \ref{uneeda this first} hold with this
choice of $p$. By \ref{and an annulus too!}, $|\oldUpsilon|=X$ is an
annulus. Since (a) holds, $p=q|\oldUpsilon$ has degree $2$.  

Suppose that some finite subgroup of $\pi_1(\frakB)$ is
non-cyclic. Then there is a point $Q\in\frakB$ such that $G_Q$ is a
dihedral group of order $2\ell$ for some integer $\ell>1$. Since
$\oldLambda=\oldLambda_0\subset\frakK$ is orientable, it now follows
that the fiber of $Q$ under the $I$-fibration $q$ contains a point
$\tQ$ which is a valence-$3$ vertex of a graph component of
$\fraks_\oldLambda\subset\fraks_\frakK$; the three oriented edges
having $\tQ$ as a terminal point are of orders $2$, $2$ and
$\ell$. This is a contradiction since
$\lambda_\frakK=2$. Hence every finite subgroup of
% $\pi_1(\frakK)$ is
%cyclic. \redcomment{Does that require very good or something?} But
%$\pi_1(\frakB)$ is isomorphic to $\pi_1(\oldLambda)$ since
%$\oldLambda$ has an $I$-fibration over $\frakB$, and $\oldLambda$ is
%$\pi_1$-injective in $\frakK$ since it is an \Ssuborbifold\ of
%$\frakK$. Hence every finite subgroup of 
$\pi_1(\frakB)$ is cyclic. It therefore indeed follows from
 the last sentence of  Lemma \ref{uneeda this first}
that $\card\fraks_\oldUpsilon\ne1$, and the proof that 
$\wt (\partial J) \ge\lambda_\frakK$ is complete in the case where (a) holds.

% \redproofreadingnote{
%I had written ``This is the kind of crap that will disappear if I don't go back and forth between working with $J$ and $X$ on the
%  one hand and working with $|\oldLambda|$ and $|\oldUpsilon|$ on the
  %other.'' I also wrote ``Part of what makes this
  %read like a bit of a mess is all the back and forth
%between working with $J$ and $X$ on the
  %one hand and working with $|\oldLambda|$ and $|\oldUpsilon|$ on the
  %other. I may have
%resolved the formal conflicts in the notation, but I should be able to
%make it read better by being more consistent.'' I now think these
%comments are based on a misinterpretation: $T$ is a homeomorphism
%between $J$ and $|\oldLambda|$, but it does {\it not} give an orbifold
%homeo between $\omega(J)$ and $\oldLambda$. So the two things should
%be kept distinct. The argument can certainly use some polishing, but
%there is nothing this simple.
%}

To complete the proof of the second alternative of the conclusion in
the case where (a) holds, it remains to show that the winding number of $X$
in $J$ (which is non-zero by the $\pi_1$-injectivity of $X$) is at most $2$. Since the homeomorphism $T:|\oldLambda|\to
J$ maps $|\oldUpsilon|$ onto $X$, this is equivalent to showing that, if $P$ denotes
the image of the inclusion homomorphism
$\pi_1(|\oldUpsilon|)\to\pi_1(|\oldLambda|)$, then the index $|\pi_1(|\oldLambda|):P|$ is at most $2$. 
Since
$|q|:|\oldLambda|\to|\frakB|$ is a homotopy equivalence by Lemma \ref{affect}, we have
$|\pi_1(|\oldLambda|):P|=|\pi_1(|\frakB|):|p|_\sharp(\pi_1(|\oldUpsilon|))|$,
%$\pi_1(|\oldUpsilon|)\to\pi_1(|\oldLambda|)=
%$|p|_\sharp(\pi_1(|\oldUpsilon|))$ has index $1$ or $2$ in
%$\pi_1(|\frakB|)$, 
where $p$ again denotes the  orbifold covering
$q|\oldUpsilon:\oldUpsilon\to\frakB$. 
According to 
%the final assertion of
%Now
%$\partial|\frakB|$ is the union of $|\partial\frakB|$ with the
%$1$-dimensional components of $\fraks_\frakB$.
%But
Lemma \ref{uneeda this first}, the index $|\pi_1(|\frakB|):|p|_\sharp(\pi_1(|\oldUpsilon|))|$ is bounded above by the degree of the orbifold covering $p$, which is equal to $2$. This completes the proof of the conclusion in the case where Alternative (a) holds.
Now suppose that (b) holds. In this case we will show that the solid
torus $J$ and the annulus $X$ satisfy the first alternative in the
conclusion of the lemma. By definition (see \ref{ypsilanti}) we have
$\oldUpsilon=\obd(T^{-1}(X))$, so that $X=T(|\oldUpsilon|)$.  Set
$X'=T(|\oldUpsilon'|)$.  Since $q|\oldUpsilon$ and $q|\oldUpsilon'$
are orbifold homeomorphisms onto $\frakB$, the orbifolds $\oldUpsilon$
and $\oldUpsilon'$ are homeomorphic. In particular, $|\oldUpsilon|$
and $|\oldUpsilon'|$ are homeomorphic, and hence so are $X$ and
$X'$. Thus (in view of \ref{and an annulus too!}) $X'$ is an annulus.

Since $\oldUpsilon$ and $\oldUpsilon'$ are the components of $\partialh\oldLambda$, it follows from (\ref{still no end}) that $X$ and $X'$ are the components of $\calf$.
Since $X\cup X'=\calf=J\cap\partial {K^+}$, and since 
$\Fr_{K^+}J=\cala_1$ by \ref {where's the beef?}, each component of $\overline{\partial J- (X\cup X')}$ is a component of $\cala_1$, and hence by \ref{lame claim}
is a weight-$0$ annulus. It follows that
$\overline{\partial J- (X\cup X')}$ is a union of two weight-$0$
annuli. Hence $\partial  J\cap\fraks_\frakK\subset \inter X\cup \inter
X'$.
% and the annuli $X$ and $X'$ have the same winding number in
%$J$. \redcomment{The equality of winding numbers is true for a pretty direct reason, but this makes it
  %sound as if it follows from what comes before, which I don't see.} 
To determine the winding numbers of $X$ and $X'$ in $J$, note that
since $q|\oldUpsilon:\oldUpsilon\to\frakB$ and $q|\oldUpsilon':\oldUpsilon'\to\frakB$ are orbifold homeomorphisms, and since
$|q|:|\oldLambda|\to|\frakB|$ is a homotopy equivalence by the first
assertion of Lemma \ref{affect}, the inclusions
$|\oldUpsilon|\to|\oldLambda|$ and $|\oldUpsilon'|\to|\oldLambda|$
are homotopy equivalences; since  $X=T(|\oldUpsilon|)$ and
$X'=T(|\oldUpsilon'|)$, the inclusions $X\to J$ and $X'\to J$  are
homotopy equivalences, and hence $X$ and $X'$ both have   winding number  $1$ in $J$.

To show that  the first alternative in the conclusion of the lemma holds in this case, it remains only to show that each of the annuli $X$ and $X'$ has non-empty intersection with $\fraks_\frakK$. It follows directly from the hypothesis that $X\cap\fraks_\frakK\ne\emptyset$. To show that $X'\cap\fraks_\frakK\ne\emptyset$, we first note that by \ref{no end},
the component $X'$ of $\calf$ contains a component of $\calf_0$; since $X'\ne X$, this component of $\calf_0$ is distinct from $X_0$, and is therefore equal to $X_0'$ in the notation introduced above. By \ref{no end}, $X'$ is the union of $X_0'$ with all disks of the form $D_C$, where $C$ ranges over the components of $\partial\cala_0$ contained in $X_0'$. Consider first the subcase in which $\cala_0\ne\emptyset$, and choose a component $A$ of $\cala_0$. Since $A$ is in particular a component of $|\Fr \oldLambda_0|$, it is saturated in the fibration of $\oldLambda_0$, and therefore meets
%; since the base $\frakB_0$ of the $I$-fibration of $\oldLambda_0$ is connected, $|A|$ meets 
%every component of $ |\partialh\oldLambda_0|=\calf_0$ (see \ref{fibered stuff}). r I have %removed Lemma ``products left,'' which isn't quoted
  %anywhere for the time being. It can be found in no-beginning.tex.}
 every component of $ |\partialh\oldLambda_0|=\calf_0$ (see \ref{fibered stuff}). Hence some component
 $C_0$ of $\partial A\subset\partial\cala_0$ is contained in $X_0'$. We therefore have $D_{C_0}\subset X'$. But by \ref{hit it} we have
 $D_{C_0}\cap\fraks_\frakK\ne\emptyset$, and hence $X'\cap\fraks_\frakK\ne\emptyset$ in this subcase.
 
 Finally, consider the subcase in which $\cala_0=\emptyset$. We then have $X=X_0$, $X'=X_0'$, $\oldLambda=\oldLambda_0$ and $\frakB=\frakB_0$.    Since $T$ is defined to be the identity on $|\oldLambda_0|$, it is the identity on its entire domain $|\oldLambda|$. Hence $\obd(X)= \obd(T^{-1}(X))=\oldUpsilon$ and $\obd(X')= \obd(T^{-1}(X'))=\oldUpsilon'$.  It follows that $X\cap\fraks_\frakK=\fraks_\oldUpsilon$ and that $X'\cap\fraks_\frakK=\fraks_{\oldUpsilon'}$. We have observed that $\oldUpsilon$ and $\oldUpsilon'$ are homeomorphic orbifolds, and hence $\card(X\cap\fraks_\frakK)=\card(\fraks_\oldUpsilon)=\card(\fraks_{\oldUpsilon'})=\card(X'\cap\fraks_\frakK)$. Since we have $X\cap\fraks_\frakK\ne\emptyset$, it follows that $X'\cap\fraks_\frakK\ne\emptyset$ in this subcase as well.
 \EndProof
 %\cala\oldLambda\obd \lambda \eta \theta B Q $h\alpha\ell

\section{Tori in the underlying space of a
  $3$-orbifold}\label{underlying tori section}

In this section  we continue the technical preparation, sketched in
the Introduction, for the proof of
Proposition E (Proposition \ref{lost corollary}).

\begin{notationremarks}\label{wait star}
If $\oldPsi$ is an orbifold, and
$S$ is a subset of $|\oldPsi|$ such that
$S\cap\fraks_\oldPsi$ is finite, we will define a quantity
$\wt^*_\oldPsi S$ by setting $\wt^*_\oldPsi S=\wt_\oldPsi S$ if
$\wt_\oldPsi S$ is even or $\lambda_\oldPsi=1$, and $\wt^*_\oldPsi
S=\wt_\oldPsi S+1 $ if $\wt_\oldPsi S$ is odd and
$\lambda_\oldPsi=2$. (See \ref{wuzza weight} and \ref{lambda thing} for the definitions of $\wt_\oldPsi\in\NN$ and $\lambda_\oldPsi\in\{1,2\}$.) Note that $\wt^*_\oldPsi S$ is always divisible by $\lambda_\oldPsi$. Note also that if $S$ and $S'$ are subsets of $\oldPsi$ such that $S\cap\fraks_\oldPsi$ and $S'\cap\fraks_\oldPsi$ are both finite, and if $\wt_\oldPsi S\le\wt_\oldPsi S'$, then $\wt^*_\oldPsi S\le\wt^*_\oldPsi S'$. Hence
 if $\wt^*_\oldPsi S'<\wt^*_\oldPsi S$ then $\wt_\oldPsi S'<\wt_\oldPsi S$. Moreover, 
 if $\wt_\oldPsi S=\wt_\oldPsi S'$ then $\wt_\oldPsi^* S=\wt_\oldPsi^* S'$.

We will write $\wt^* S$ for $\wt^*_\oldPsi S$ when the orbifold $\oldPsi$ is understood.
\end{notationremarks}

\Lemma\label{pre-modification} Let $\oldPsi$ be a
compact, orientable, strongly \simple, boundary-irreducible 
$3$-orbifold containing no embedded negative
turnovers. 
%, such that each component of $\fraks_\oldPsi$ is either an arc
%or a simple closed curve. 
Set $N=|\oldPsi |$. Suppose that each component of
$\partial N$ is a sphere, and that $N$ is $+$-irreducible (see
Definition \ref{P-stuff}). Let $K$
 be a non-empty, proper, compact, connected, $3$-dimensional
 submanifold of $N$. 
Assume that $\calt:=\Fr_NK$ is contained in $\inter N$ and is 
in general position with respect to
$\fraks_\oldPsi$,  and  that its components are all incompressible tori in $N$ (so that $K^+$ is naturally identified with a submanifold of $N^+$ by \ref{plus-contained}).
% is irreducible and boundary-irreducible). 
Assume that $K^+$ is not homeomorphic to $\TTT^2\times[0,1]$ or to a twisted $I$-bundle over a Klein bottle.
Suppose that either (i)  $\obd(\calt)$ fails to be
incompressible in $\oldPsi$, or (ii)  $\obd(\calt)$ is
incompressible in $\oldPsi$ (so that $\obd(K)$ is boundary-irreducible and strongly \simple\ by Lemma \ref{oops lemma}, and hence $\kish(\obd(K))$ is defined in
view of \ref{tuesa day}) and
$$\chibar(\kish(\obd(K)))<\min\bigg(1,\frac14\wt^*\nolimits_\oldPsi(\calt)\bigg).$$
Then at least one of the following conditions holds:
\begin{enumerate}
\item There exist a disk $D\subset N$ with
$\partial D=D\cap \calt $, such that $D$ 
in general position with respect to
$\fraks_\oldPsi$,
and 
a disk $E\subset \calt $ such that $\partial E=\partial D$ and 
$\wt_\oldPsi( E)>\wt_\oldPsi(D)$; furthermore, if 
$\lambda_\oldPsi=2$,
then $\max(\wt_\oldPsi(
E)-\wt_\oldPsi(D),\wt_\oldPsi D)\ge2$. 
\item There is 
a solid torus $J\subset K^+$, $\pi_1$-injective in $N^+$,
 with $\partial J\subset
K\cap\inter N
\subset K\subset K^+$,
%\redmissingref{cross-ref for $+$ notation?} 
  such that 
  $\partial J\cap\partial K$ is a union of two disjoint annuli $X$ and $X'$ contained in $\calt$, each having winding
  number $1$ in $J$ (see \ref{great day}) and each having non-empty
  intersection with $\fraks_\oldPsi$, and 
  $\partial J\cap\fraks_\oldPsi\subset \inter X\cup \inter X'$.
\item 
There is 
a solid torus $J\subset K^+$, $\pi_1$-injective in $N+$, 
 with $\partial J\subset K\cap\inter N\subset K\subset K^+$,   such that 
  $\partial J\cap\partial K$ is an annulus $X\subset\calt$, having winding
  number $1$ or $2$ in $J$, and 
we have
  $\wt (\partial J) \ge\lambda_\oldPsi$ and $\partial J\cap\fraks_\oldPsi\subset \inter X $.
%such that
  %$\emptyset\ne (\partial J)\cap\fraks_\oldPsi\subset X $.
\end{enumerate}
\EndLemma

\Proof[Proof of Lemma \ref{pre-modification}]
Let us  consider the case of the lemma in which Alternative (i) of the
hypothesis holds. 
Since  $\obd(\calt )$ is not incompressible in $\oldPsi$, it
follows from Proposition \ref{kinda dumb}, more specifically the implication (c)$\Rightarrow$(a), that there is an orientable discal $2$-suborbifold $\frakD$ of $\oldPsi$ such that $\frakD\cap\obd(\calt)=\partial\frak D$, but such that there is no discal $2$-suborbifold $\frakE$ of $\obd(\calt)$ with $\partial\frakE=\partial\frakD$. (The hypothesis in the final assertion Proposition \ref{kinda dumb} that $\oldPsi$ is very good follows from the strong \simple ity of $\oldPsi$, while the hypothesis that $\oldTheta:=\obd(\calt)$ is non-spherical follows from the inequality $\chi(\obd(\calt))\le\chi(\calt)=0$.)
%general position  implies that $|\frakD|\cap\fraks_\oldPsi$ is %finite, so that $\frakD$ is orientable, and discality then implies that $\wt|\frakD|\le1_m$. 
Assume that such a $\frakD$ does exist. Then
there is a disk $D\subset N$, in general position with respect to $\fraks_\oldPsi$,  with
$\gamma:=\partial D=D\cap \calt $, such that $\wt_\oldPsi D\le1$ and such that $\gamma$ does not bound any disk of weight at most $1$ in $ \calt $. Since $\calt $ is incompressible by hypothesis,
there is a disk $E\subset \calt $ with $\partial E=\gamma$. Hence
$\wt E\ge2$. In particular we have
$\wt E>\wt D$. Furthermore, if $\wt (E)-\wt(D)=1$, then we must have
$\wt E=2$ and $\wt D=1$; this implies that $D\cup E$ is a $2$-sphere
of weight $3$ in $|\oldPsi|$, 
in general position with respect to $\fraks_\oldPsi$. If $\lambda_\oldPsi=2$, then since $\oldPsi$ 
%is
%\simple\ \redmissingref{and boundary-irreducible?} and 
contains
no embedded negative turnovers, it follows from Proposition \ref{almost obvious}
that $|\oldPsi|$ contains no weight-$3$
sphere in general position with respect to $\fraks_{\oldPsi}$, and hence  $\wt (E)-\wt(D)\ge2$. 
This shows that Alternative (1) of
the conclusion of the lemma holds in this case.
%T

For the rest of the proof, we will assume that Alternative (ii) of the hypothesis
holds (so that in particular $\obd(\calt)$ is
incompressible in $\oldPsi$). The strategy of the proof under this
assumption is to try to find a subsurface $X\subset K$ satisfying the
hypothesis of Lemma \ref{uneeda}, with $\frakK=\obd(K)$; when such a subsurface
 can be found, Alternative (2) or (3) of the
conclusion of the present lemma will be seen to hold. When the quest
for such a subsurface fails, Alternative (1) of the conclusion will turn out to hold.

Set 
$\frakU=\overline{\obd(K)-\oldSigma(\obd(K))}$. According to the definition of $\kish(\obd(K))$ (see \ref{tuesa day}), $\frakU$ is the disjoint union of $\kish(\obd(K))$ with certain components of $\frakH(\obd(K))$. Since the latter components have Euler characteristic $0$, we have $\chi(\frakU)=\chi(\kish(\obd(K)))$. Hence Alternative (ii) of the hypothesis gives
\Equation\label{bought}
\chibar(\frakU)<\min\bigg(1,\frac14\wt^*\nolimits(\calt)\bigg).
\EndEquation

We claim:
\Claim\label{a half or better}
If $\oldXi$ is any component of $\frakU\cap\obd(\partial K)$ such that
$\wt|\oldXi|$ is odd, we have $\chibar(\oldXi)\ge1/2$. 
%\redmissingref{Am I allowed to write $\wt\oldXi$ or do I have to write $\wt|\oldXi|$? Also, this should be isolated as a general lemma, because it is also used in improving Lemma \ref{modification}, replacing the inequality $\wt_\oldPsi(\calt _1)<
%\wt_\oldPsi\calt$ by $\wt_\oldPsi(\calt _1)\le
%\wt_\oldPsi\calt$ (which implies $\wt^*_\oldPsi(\calt _1)<
%\wt^*_\oldPsi\calt$). This improvement is needed later in the proof %of Proposition \ref{crust chastened}. See the comment beginning ``P%roblem here.''}
\EndClaim

To prove \ref{a half or better}, set $q=\wt|\oldXi|$, write
$\fraks_\oldXi=\{z_,\ldots,z_q\}$, and let $e_i$ denote the order of
$z_i$ for $i=1,\ldots,q$. Then by \ref{2-dim case} we have
$\chibar(\oldXi)=\chibar(|\oldXi|)+\sum_{i=1}^q(1-1/e_i)\ge\chibar(|\oldXi|)+q/2$. If
$|\oldXi|$ is not a disk or a sphere, we have $\chibar(|\oldXi|)\ge0$,
and $q\ge1$ since $q$ is odd; hence $\chibar(\oldXi)\ge1/2$. If
$|\oldXi|$ is a disk, we have $q\ne1$ since $\partial\oldXi\subset\partial\oldPhi(\omega(K))$
is $\pi_1$-injective in $\partial(\obd(K))$ by
\ref{tuesa day}; since $q$ is odd we have $q\ge3$, and hence
$\chibar(\oldXi)\ge-1+3/2=1/2$. Now suppose that $|\oldXi|$ is a sphere. Since
$\oldXi\subset\partial\obd(K)=\partial\oldPsi\cup\calt$, and the
components of $\calt$ are tori, $\oldXi$ is a component of
$\partial\oldPsi$. In this case, if $q=1$ then $\oldXi$ is a bad $2$-orbifold; this is impossible since the strongly \simple\
$3$-orbifold $\oldPsi$ is in particular very good (see \ref{oops}). Hence $q\ne1$. Next note that since
$\oldPsi$ is strongly \simple\ and boundary-irreducible, its boundary component $\oldXi$ has negative Euler characteristic by \ref{boundary is negative}; hence if $q=3$, then $\oldXi$ is a negative
turnover, a contradiction to the hypothesis. 
Since $q$ is odd we now have $q\ge5$, and hence $\chibar(\oldXi)\ge-2+5/2=1/2$. Thus \ref{a half or better} is proved.

Now we claim:
\Claim\label{cause i may have said so}
If $Q$ is a union of components of $\partial|\frakU|$, each of which has strictly positive
genus, we have
$$\wt^*\nolimits Q\le4\chibar (\frakU).$$
\EndClaim

To prove \ref{cause i may have said so}, set $n=\wt Q$,
and write $Q\cap\fraks_\oldPsi=\{x_1,\ldots,x_n\}$. Then $x_1,\ldots,x_n$
are the singular points of the $2$-orbifold $\obd(Q)$. Let $p_i$
denote the order of the singular point $x_i$. Observing that
$p_i\ge2$ for $i=1,\ldots,n$, and that $\chibar(Q)\ge0$ since each
component of $Q$ has
positive genus, and applying \ref{2-dim case}, we find that
$\chibar(\obd(Q))=\chibar(Q)+\sum_{i=1}^n(1-1/{p_i})\ge n/2$.
Hence $n\le2\chibar(\obd(Q))$, i.e.
\Equation\label{mcfoofus}
\wt Q\le2\chibar (\obd(Q)).
\EndEquation
 Since $\obd(Q)$ is a union of components of
$\partial\frakU$, and since every  component of
$\partial\frakU$ has non-positive Euler characteristic by
\ref{tuesa day}, it follows that 
$\wt Q\le2\chibar(\partial\frakU)$. But since $\frakU$ is a compact $3$-orbifold we have
$\chibar(\partial\frakU)=2\chibar(\frakU)$, and hence
%$n\le4\chibar(\frakU)$, i.e.
\Equation\label{almost cause i may have said so}
\wt Q\le4\chibar (\frakU).
\EndEquation
If $\lambda_\oldPsi=1$ or $\wt Q$ is even, the conclusion of
\ref{cause i may have said so} follows from (\ref{almost cause i may
  have said so}) in view of the definition of $\wt^*Q$. There remains
the case in which $\lambda_\oldPsi=2$ and $\wt Q$ is odd. In this
case, write $\partial|\frakU|=Q\discup Q'$, where $Q'$ is a union of
components of $\partial|\frakU|$. We have $\wt(|\partial\frakU|)=\wt
Q+\wt Q'$. Since $\lambda_\oldPsi=2$, each component of
$\fraks_\oldPsi\cap|\frakU|$ is an arc or a simple closed curve. Hence
$\wt\partial|\frakU|$ is equal to twice the number of arc components
of $\fraks_\oldPsi\cap|\frakU|$, and is therefore even. Since $\wt Q$
is odd, it now follows that $\wt Q'$ is odd. Each component of $Q'\cap\Fr\frakU$ is a component of $\Fr\frakU=\frakA(\obd(K)$, and is therefore annular (and orientable) by \ref{tuesa day}; in particular, each component of $Q'\cap\Fr\frakU$ has weight $0$ or $2$. Hence $\wt(Q'\cap\partial K)$ is odd. There therefore exists a component $\oldXi$ of $\obd(Q'\cap\partial K)$ such that $\wt|\oldXi|$ is odd. According to \ref{a half or better}, we have $\chibar(\oldXi)\ge1/2$. But 
by
\ref{tuesa day},
$\partial\oldXi\subset\partial\oldPhi(\omega(K))$
is $\pi_1$-injective in
$\obd(K)$ and hence in $\obd(Q')$; since every  component of
$\partial\frakU$ has non-positive Euler characteristic by
\ref{tuesa day}, it follows that $\chibar(\obd(Q'))\ge\chibar(\oldXi)\ge1/2$. Combining this with (\ref{mcfoofus}), we find $\chibar(\partial\frakU)=\chibar(\obd(Q'))+\chibar(\obd(Q)\ge(1+\wt Q)/2$. Since $\wt^* Q=1+\wt Q$ in this case, it follows that $\wt^*Q\le2\chibar(\partial\frakU)=4\chibar(\frakU)$, as required. This completes the proof of \ref{cause i may have said so}.
%\wt\partial\frakU=\wt\
%\wt\nolimits _\oldPsi (Q)\le2\chibar (\obd(Q)).Thus we have 
%the $\pi_1$-injectivity of
%\redproofsummary{ 
%%The proof can be extracted from the present proof of \ref{if not why not} (or more specifically %(\ref{hot poop}) 
%In the special case where $\lambda_\oldPsi=1$ or $\wt\caly_0$ is even, the proof is %essentially as follows: ``Since each component of $Q_0$ is
%a torus, and in particular has strictly positive genus, we may
%therefore apply \ref{where's the oakie?} to deduce that
%\Equation\label{hot poop}
%$\wt\caly_0\le\wt Q_0\le4\chibar (\frakU)$.''  .In the case where $\lambda_\oldPsi=2$ and %$\wt\caly_0$ is odd, it requires further
%In this case use Lemma \ref{fortunately} to show that of $\frakU\cap\partial$ has an odd-weight component not contained in $\caly_0$. It has $\chibar\ge1/2$ by \ref{a half or better}, and the result follows upon combining this with (\ref{almost cause i may have said so}). 
%\wt_

%as asserted in \ref{cause i may have said so}.

Let $\calx_0$ denote the union of all components of
$|\oldSigma(\obd(K))|\cap\calt=|\oldPhi(\obd(K))|\cap\calt$ (see \ref{tuesa day}) that
%are contained in $\calt$ but
 are not contained in disks in $\calt$. 
%Since each component of
%$\partial K$ is either a sphere or a component of $\calt$, we have
%$\calx_0\subset \calt$. 
Let $\calx\subset\calt$ denote the union of $\calx_0$ with all disks
in $\calt$ whose boundaries are contained in $\partial \calx_0$. Since the
components of $\calt$ are tori, each component of $\calx$ is either an
annulus or a component of $\calt$. Furthermore, no annulus component
of $\calx$ is contained in a disk in $\calt$; that is, the annulus component
of $\calx$ are homotopically non-trivial in the torus components of
$\calt$ containing them. Hence each component of
$\caly:=\overline{\calt-\calx}$ is also  either a homotopically non-trivial
annulus in $\calt$ or a component of $\calt$. Since  $\calt$ is incompressible in $N$, it now follows that the
submanifolds $\calx$ and $\caly$ are $\pi_1$-injective in $N$.

It follows from the definition of $\caly$ that every component of $(\inter
\caly)\cap|\oldSigma(\obd(K)|$ is contained in a disk in $\inter \caly$. Hence $(\inter
\caly)\cap|\oldSigma(\obd(K)|$ is contained in a (possibly empty) disjoint
union of disks
  $E_1,\ldots,E_m\subset\inter \caly$, where $\partial E_i\subset
 |\oldSigma(\obd(K)|$ whenever $1\le i\le m$. The definitions also imply
  that $\caly_0:=\overline{\caly-(E_1\cup\cdots\cup E_m)}$ is a union of components
    of $\calt\cap|\frakU|$. 
\abstractcomment{\tiny The expression $\calt\cap|\frakU|$ replaces
  $\overline{\calt\setminus|\oldPhi(\obd(K)|}$ in the version new-5point4.tex.}
%
%Since $\caly_0$ is a union of components
  %  of $\overline{\calt\setminus|\oldPhi(\obd(K)|}$, it is in particular  a union of component%s
  %  of $\overline{(\partial
    %  K)\setminus|\oldPhi(\obd(K)|}$. 
Hence $\obd(\caly_0)$ is a union of
    components of $\obd(\calt)\cap\frakU$, which are in particular components of $\frakU\cap\partial\obd(K)$.
%$\overline{\obd(\partial
  %    K)-\oldPhi(\obd(K)}$. 
Since each component of $\frakU\cap\partial\obd(K)$ has non-positive Euler
    characteristic by \ref{tuesa day}, it follows that
\Equation\label{gunman}
\chibar(
\obd(\calt)\cap\frakU))\ge\chibar(\obd(\caly_0)).
\EndEquation

On the other hand, since each component of $\caly$
is an annulus or a torus, we have $\chibar(\caly_0)=m$. By
\ref{2-dim case} it follows that
\Equation\label{by crossref}
\chibar(\obd(\caly_0))\ge m.
\EndEquation
Now by \ref{tuesa day} we have
    $\chibar(\frakU)
%=\chibar(\partial\frakU)/2=
=\chibar(\frakU\cap\omega(\partial
K))/2=(\chibar(\frakU\cap\omega(\calt))+\chibar(\frakU\cap\partial
\oldPsi))/2$. Since $\frakU\cap\partial \oldPsi$ is a union of
components of 
$\overline{\partial\obd(K)-\oldPhi(\obd(K))}$, it follows from
\ref{tuesa day} that
%\frakU\cap\partial \oldPsi By \redmissingref{standard reason} we have
$\chibar(\frakU\cap\partial \oldPsi)\ge0$, and hence
$\chibar(\frakU)
\ge\chibar(\obd(\calt)\cap\frakU)/2$.
%=\chibar(\overline{\obd(\partial
  %    K)-\oldPhi(\obd(K)})/2$. 
%\redcomment{That doesn't seem right. Should it be $\chibar(\partial\frakU)/2\ge\chibar(\obd(\calt)\cap\frakU)/2$? That may be OK, because the contributions of $\partial\oldPsi\cap\frakU$ to $\chibar$ should be non-negative, and it seems like enough for the app. just below.}
With (\ref{gunman}) and (\ref{by
      crossref}) this gives $\chibar(\frakU)\ge m/2$. Since
$    \chibar(\frakU)<1$ by (\ref{bought}), it follows that
$m<2$. Hence
\Equation\label{charmin}
m\le1.
\EndEquation

%\redcomment{Decide whether in \ref{if not why not} and \ref{cause i said
    %so}, the right inequalities are the ones in the visible version,
  %the ones in the \%ed out version, or neither.}

Next we claim:
\Claim \label{mabel syrble}
If a component $F$ of $|\partial\frakU|$ contains a component
of $\caly_0$, then $F$ is a torus.
\EndClaim

To prove \ref{mabel syrble}, let $Y_0$ be a component of $\caly_0$ contained
in $F$, and let $Y$ and $T$ respectively denote the components of $\caly$
and $\calt$ containing $Y_0$. Since $Y$ is a component of $\caly$, either $Y$ is an  annulus which is homotopically
non-trivial in $N$, or $Y=T$. If $Y=T$ then $Y_0$ is a genus-$1$ subsurface of
$F$. If $Y$ is a homotopically non-trivial annulus, then the
components of $\partial Y\subset Y_0\subset F$ are homotopically non-trivial
simple closed curves in $N$. Thus in any event, $F$ contains either a genus-$1$
subsurface or a simple closed curve which is homotopically non-trivial in $N$. Hence
$F$ cannot be a sphere.

If $F$ has genus at least two, then $\chibar(F)\ge2$. By \ref{2-dim case}, it follows that $\chibar(\obd(F))\ge2$. Since every
component of $ \partial\frakU$
%(\overline{K-|\oldSigma(\obd(K)|})$ 
has non-positive
Euler characteristic by \ref{tuesa day}, we have
$\chibar( \partial\frakU)
%(\overline{K-|\oldSigma(\obd(K)|}))
\ge\chibar(\obd(F))\ge2$, and hence $\chibar(\frakU)=\chibar( \partial\frakU)/2\ge1$.
%$$\begin{aligned}\chibar(\frakU)&=\chibar(\overline{\obd(K)-\oldSigma(\obd(K)})\\&\ge \%chibar(|\overline{\obd(K)-\oldSigma(\obd(K)}|)\\&=
%\chibar(\overline{K-|\oldSigma(\obd(K)|}))\\&=\frac12\chibar( \partial(\overline{K-|\oldSigma%(\obd(K)|}))\ge1.\end{aligned}.$$
This a
contradiction, since $\chibar(\frakU)<1$ by (\ref{bought}). Thus \ref{mabel syrble} is proved.

Now we claim:
\Claim\label{if not why not}
%If $\card (\caly\cap\fraks_\oldPsi)>\card ((\partial\frakU)\cap\fraks_\oldPsi)$ then
%the
If $\wt^*\caly>4\chibar (\frakU)$, then Alternative (1) of the
conclusion of the lemma holds. 
\EndClaim

To prove \ref{if not why not}, 
let $Q_0$ denote the union of all
components of $|\partial\frakU|$ that contain components of
$\caly_0$. According to \ref{mabel syrble},  each component of $Q_0$ is
a torus. Since in particular every component of $Q_0$ has strictly positive genus, it follows from 
\ref{cause i may have said so} that
%If $Q$ is a union of components of $\partial|\frakU|$, each of which has strictly positive
%genus, we have
%\Equation\label{yet another
%\wt^*\nolimits _\oldPsi (Q_0)\le4\chibar (\frakU).
%\EndEquation
%We claim:
\Equation\label{where's the oakie?}
\wt\nolimits^*Q_0\le4\chibar(\frakU).
\EndEquation
%first note that
On the other hand, since $\caly_0\subset Q_0$, we have $\wt\caly_0\le\wt Q_0$, which by \ref{wait star} implies $\wt^*\caly_0\le\wt^*Q_0$. With (\ref{where's the oakie?}), this implies $\wt^*\caly_0\le4\chibar(\frakU)$.
%; by
%let $Q_0$ denote the union of all
%components of $|\partial\frakU|$ that contain components of
%$\caly_0$. According to \ref{mabel syrble},  each component of $Q_0$ is
%a torus (and in particular has strictly positive genus). We may
%therefore apply 
%\ref{cause i may
  %have said so} it follows that
%to deduce that
%\Equation\label{hot poop}
%$\wt\caly_0\le\wt Q_0\le4\chibar (\frakU)$.  
%\EndEquation
Now if $m=0$, we have
$\caly=\caly_0$, and hence
%\subset |\partial\frakU|$, \redmissingref{Better to be consistent about writing $|\partial\frakU|$ or $\partial|\frakU|$} and (\ref{hot poop}) implies that
$\wt^*\caly\le4\chibar (\frakU)$.
Thus \ref{if not why not} is vacuously true when $m=0$. In view of
(\ref{charmin}), it remains only to consider the subcase $m=1$.

In this
subcase set $E=E_1$, so that $\caly_0=\overline{\caly-E}$. 
%In particular, $\caly_0$
%is a genus-one surface with connected boundary. Hence if $T$ denotes
%the component of 
Let us index the components of $\caly$ as $Y^0,\ldots,Y^r$, where $r\ge0$
and $E\subset Y^0$. 
Then 
$Y^0_0:=\overline{Y^0-E}$ is a component of $\caly_0$. In particular we have
$Y^0_0\subset\calt\cap|\frakU|\subset\partial|\frakU|$.
%setminus|\oldPhi(\obd(K)|}\subset\partial(\overline{K-|\oldSigma(\obd(K)|})$.
Let $T$ denote the (torus) component of $\calt$ containing $Y^0$.
Let $F$ denote the component of
$Q_0$ containing $Y^0_0$. Then $F$ is
a torus by \ref{mabel syrble}.

Let $\Delta$ denote the component of $\overline{F-Y^0_0}$ containing
$\alpha:=\partial E\subset\partial Y^0_0$. If $Y^0=T$, then
$Y^0_0$ is a genus-$1$ surface with one boundary component, and
hence $\Delta$ is a disk. If $Y^0$ is an annulus, then $Y^0_0$ is a
connected planar surface with three boundary components; two of
these, the components of $\partial Y^0$, are homotopically non-trivial
in $N$ and hence in the torus $F$. It follows that $\alpha$, the third
boundary component of $Y^0_0$, is homotopically trivial in $F$. Thus
in any event,
%so that
$\Delta$ is a disk.
% in this case as well.

Set $Z=Y^0_0\cup\Delta$.
%, and for $1<i\le r$ set $Z^i=Y^i$.
Then
the surfaces $Z,Y^1,\ldots, Y^r$ are contained in
$Q_0$. The surfaces $Y^1,\ldots,Y^r$ are pairwise disjoint
because they are distinct components of $\caly$. If $1\le i\le r$, we
have $Y^0\cap Y^i=\emptyset$ and hence $Z\cap
Y^i\subset\Delta\subset\inter Z$. Since
$Y^i$ is connected, and $Q_0$ is a closed $2$-manifold, we must have either $Y^i\subset\Delta$ or $Z\cap
Y^i=\emptyset$. But since $Y^i$ is a homotopically non-trivial annulus or an incompressible torus
in $N$, it cannot be contained in the disk
$\Delta$. This shows that the surfaces $Z,Y^1,\ldots, Y^r$ are
pairwise disjoint. Hence
\Equation\label{jooolie}
\wt\nolimits Z+\sum_{i=1}^r\wt\nolimits Y^i\le \wt Q_0.
\EndEquation
On the other hand, (\ref{where's the oakie?}) and the hypothesis of Claim
\ref{if not why not} give
$\wt^* Q_0\le4\chibar
(\frakU)<\wt^*\caly$, which by \ref{wait star} implies $\wt Q_0<\wt\caly$. Since $\wt\caly=\sum_{i=0}^r\wt Y^i$, it follows that
\Equation\label{not sloop}
\wt Q_0<\sum_{i=0}^r\wt Y^i.
\EndEquation
From (\ref{jooolie}) and (\ref{not sloop}) it follows that
$\wt Z< \wt Y^0$,
which may be rewritten as
$\wt Y^0_0+\wt \Delta<
\wt Y^0_0+\wt E.$
Hence 
\Equation\label{one last wisk}
\wt\Delta<
\wt E.
\EndEquation

Since the disk $\Delta$ is contained in $K$ and its boundary $\beta$ lies in
$\partial K$, we may modify $\Delta$ by a small isotopy, constant on
$\beta$, to obtain a disk $D$ which is properly embedded in $K$. In
particular we have $D\cap \calt=\beta$. We may choose the isotopy in such a
way that $D$ is in general position with respect to $\fraks_\oldPsi$ and
$\wt D=\wt\Delta$. With (\ref{one last
  wisk}), this gives $\wt D<
\wt E$. 
%If $\lambda_\oldPsi=1$ it follows that 
%$\wt
%(E)-\wt(D),\wt D)\ge \lambda_\oldPsi$, which implies Alternative (1) of the conclusion of the
%lemma. 
Now suppose that $\lambda_\oldPsi=2$ and that
$\max(\wt
(E)-\wt(D),\wt D)<2$. Since $\wt D<
\wt E$, we have either $\wt E=2$ and $\wt D=1$, or
$\wt E=1$ and $\wt D=0$; hence
$D\cup E$ is a $2$-sphere
of weight $1$ or $3$ in $|\oldPsi|$. If the sphere $D\cup E$ has weight $1$ then $\obd(D\cup E)$ is a bad $2$-suborbifold of $\oldPsi$; this is impossible since the strongly \simple\ $3$-orbifold
$\oldPsi$ is in particular very good (see \ref{oops}). If $\lambda_\oldPsi=2$, then since $\oldPsi$ 
%is
%\simple\ \redmissingref{and boundary-irreducible?} and 
contains
no embedded negative turnovers, it follows from Proposition \ref{almost obvious} that $|\oldPsi|$ contains no weight-$3$
sphere 
in general position with respect to $\fraks_\oldPsi$. This contradiction shows that if $\lambda_\oldPsi=2$ then $\max(\wt
(E)-\wt(D),\wt D)\ge2$. Thus Alternative (1) of the conclusion of the
lemma holds  in this situation, and Claim \ref{if not why not} is  established.

%ZY^

Next we claim:
\Claim\label{cornway}
Either $\calx\cap\fraks_\oldPsi\ne\emptyset$, or Alternative (1) of the
conclusion of the lemma holds. 
\EndClaim

To prove \ref{cornway}, assume that Alternative (1) of the
conclusion of the lemma does not hold.
Then by \ref{if not why not}, we have
\Equation\label{cause i said so}
%\card (Y\cap\fraks_\oldPsi)\le\card ((\partial\frakU)\cap\fraks_\oldPsi).
\wt\nolimits^* \caly\le4\chibar (\frakU).
\EndEquation

Note that
\Equation\label{me quation}
\wt\calx+\wt \caly=\wt\calt. 
\EndEquation
%It follows from (\ref{bought}) that
%\Equation\label{requation}
%\wt\nolimits^*\calt>4\chibar(\frakU).
%\EndEquation
If we assume that $\wt\calx=0$, then it follows from (\ref{me quation}) that 
$\wt \caly=\wt\calt$, which with \ref{wait star} implies that $\wt^* \caly=\wt^*\calt$. Combined with
(\ref{cause i said so}) this gives
$\wt^* \calt\le4\chibar (\frakU)$, which contradicts
(\ref{bought}). Hence we must have
% \redcomment{Fix the next paragraph and then remove the comment following it.}
%It follows from (\ref{cause i said so}), (\ref{me quation}), and
% (\ref{requation}) that
 $\wt\calx>0$, i.e. 
%that
$\calx\cap\fraks_\oldPsi\ne\emptyset$, and
\ref{cornway} is proved.

In view of \ref{cornway}, we will assume for the remainder of the proof
that $\calx\cap\fraks_\oldPsi\ne\emptyset$. Let us fix a component $X$ of
$\calx$ such that $X\cap\fraks_\oldPsi\ne\emptyset$. The definitions of 
$\calx_0$ and $\calx$ imply that $X$ contains a unique component $X_0$ of
$\calx_0$, and that there are disjoint disks $G_1,\ldots,G_n$ in $\inter X_0$ such that
$X_0=\overline{X-(G_1\cup \cdots\cup G_n)}$. The definition of
$\calx_0$ also implies that $X_0$ is a component of
$|\oldPhi(\obd(K))|$. As we have seen that $\calx$ is $\pi_1$-injective in $N$, in particular $X$ is $\pi_1$-injective in $K$.
%and that the
%inclusion homomorphism $\pi_1(X_0)\to\pi_1(K)$ is non-trivial
%(and hence so is the
%inclusion homomorphism $\pi_1(X)\to\pi_1(K)$). 
Since $X\subset K$,
we have $X\cap\fraks_{\obd(K)}=X\cap\fraks_\oldPsi\ne\emptyset$. Note
also that $K$ is boundary-irreducible and $+$-irreducible, since $N$
is $+$-irreducible, $\partial N$ consists of sphere components, and  $\calt=\partial K^+$ is
incompressible in $N$. We have seen that each component of $\calx$ is
either a (torus) component of $\calt$ or an
annulus; in particular, we have
% $X$ is
%$\pi_1$-injective in $\calt$, and therefore in $K$, and
 $\chi(X)=0$. Furthermore, by
hypothesis,  $K^+$ is not homeomorphic to $\TTT^2\times[0,1]$ or to a
twisted $I$-bundle over a Klein bottle. Thus
all the hypotheses of Lemma \ref{uneeda} hold with
$\frakK=\obd(K)$. Hence by Lemma 
\ref{uneeda}, $X$ is an annulus, and there is a
$\pi_1$-injective solid torus $J\subset K^+$, 
%\redmissingref{cross-ref for $+$ notation. Also, fix all
%cross-refs!!} 
with $\partial J\subset K\subset K^+$, such that one of the following alternatives holds:
\begin{itemize}
\item We have $\partial J\cap\partial K=X\discup X'$ for some annulus $X'\subset\partial K^+ $; furthermore, each of the annuli $X$ and $X'$ has winding
  number $1$ in $J$ and has non-empty
  intersection with $\fraks_\frakK$, and 
$\partial  J\cap\fraks_\frakK\subset \inter X\cup \inter X'$.
\item We have $\partial J\cap\partial K=X$; furthermore, $X$ has winding
  number $1$ or $2$ in $J$, and 
%is disjoint from $\fraks_\frakK$, and $
%(\partial  J)\cap\fraks_\frakK\subset \inter A_1\cup \inter A_2$.
  %is a single annulus  $A\subset\partial K^+ $, having winding
  %number at most $2$ in $J$, and
we have  $\wt (\partial J) \ge\lambda_\frakK$ and $\partial J\cap\fraks_\frakK\subset \inter X $.
%  $\emptyset\ne \partial J\cap\fraks_\frakK\subset X $.
\end{itemize}

Note also that since $\frakK$ is a suborbifold of $\oldPsi$, and $J\subset K$, we have $\partial J\cap\fraks_\oldPsi=\partial J\cap\fraks_\frakK$. Hence the two alternatives above imply respectively that $\partial  J\cap\fraks_\oldPsi\subset \inter X\cup \inter X'$ and that  
$\partial  J\cap\fraks_\oldPsi\subset \inter X$. Likewise, since $\frakK$ is a suborbifold of $\oldPsi$,
it follows from the definitions that $\lambda_\frakK\ge\lambda_\oldPsi$, so that when the second alternative stated above holds, we have $\wt(\partial J)\ge\lambda_\oldPsi$. Finally, since the components of $\calt$ are incompressible in $N$, the submanifold $K^+$ of $N^+$ is $\pi_1$-injective, and hence the $\pi_1$-injectivity of $J$ in $K^+$ implies that it is $\pi_1$-injective in $N^+$. Thus one of the alternatives (2), (3) of the
conclusion holds.
\EndProof
%\oldTheta

\Lemma\label{modification} Let $\oldPsi$ be a
compact, orientable, strongly \simple, boundary-irreducible
$3$-orbifold containing no embedded negative turnovers. 
Set $N=|\oldPsi |$. Suppose that each component of
$\partial N$ is a sphere, and that $N$ is $+$-irreducible. Let $K$
 be a non-empty, proper, compact, connected, $3$-dimensional
 submanifold of $N$, and set $\calt=\Fr_N K$.  Assume that $\calt$ is contained in $N$ and is in general position with respect to
$\fraks_\oldPsi$,  and  that its components are all incompressible tori in $N$ (so that $K^+$, which by \ref{plus-contained} is naturally identified with a submanifold of $N^+$,  is irreducible and boundary-irreducible), 
 and that   either (a) $K^+$ is
acylindrical, or (b) $\overline{N-K}$ is connected and
$h(K)\ge3$. 
%Set $p=\wt_\oldPsi^*\calt$.  
Suppose in addition that either (i)  $\obd(\calt )$ fails to be
incompressible in $\oldPsi$, or (ii)  $\obd(\calt )$ is
incompressible in $\oldPsi$ (so that $\obd(K)$ is boundary-irreducible and strongly \simple\ by Lemma \ref{oops lemma}, and hence $\kish(\obd(K))$ is defined in
view of \ref{tuesa day}) and
$\chibar(\kish(\obd(K)))<\lambda_\oldPsi/4$, or (ii$'$) $\wt_\oldPsi^*\calt \ge4$,  $\obd(\calt )$ is
incompressible in $\oldPsi$, and
$\chibar(\kish(\obd(K)))<1$. Then there is a compact,
connected $3$-manifold $K_1\subset N$ having the following properties:
\begin{itemize}
\item every component of $\calt _1:=\Fr_NK_1$ is an incompressible torus in
  $\inter N$, in general position with respect to
$\fraks_\oldPsi$;
\item $\wt_\oldPsi(\calt _1)<
\wt_\oldPsi\calt$; and if $\lambda_\oldPsi=2$, then $\max(\wt_\oldPsi(\calt )-
\wt_\oldPsi(\calt_1)
,\wt_\oldPsi(\calt_1))\ge2$; 
\item $h(K_1)\ge h(K)/2$; and
\item $h(N-K_1)\ge h(N-K)$.
\end{itemize}
Any such $K_1$ is a proper, non-empty submanifold of $N$.

Furthermore, if (a) holds, $K_1$ may be chosen so that
$K_1^+$ and $K^+$ (which by \ref{plus-contained} are identified with
  submanifolds of $\plusN$) are isotopic in $\plusN$; and if (b) holds, $K_1$ may be chosen so that $\overline{N- K_1} $ is connected.
\EndLemma

\Proof
We set $P=K^+$. As was observed in the statement of the lemma, since $\calt$ is closed and has no sphere components, it follows from \ref{plus-contained} that we may
identify $P$ with a submanifold of $N^+$. It was also observed in the statement of the lemma that the $+$-irreducibility of $P$, the incompressibility of  $\calt$ and the hypothesis that $N$ has only sphere components, imply that $P$ is irreducible and boundary-irreducible. (This implication, which is standard in $3$-manifold theory, is also included  in the manifold case of Lemma \ref{oops lemma}.) 

Note that the hypotheses imply that $K$ has at
least one frontier component, or equivalently that it has at
least one torus boundary component. In particular,
$h(K)$ and $h(N-K)$ are both strictly positive. Hence if $K_1$ is a
submanifold with the first four properties listed, we have $h(K_1)\ge
h(K)/2>0$ and $h(N-K_1)\ge
h(N-K)>0$, and therefore $K_1\ne\emptyset$ and $N-K_1\ne\emptyset$. Thus
the first four properties listed do imply that $K_1$ is a non-empty,
proper submanifold of $N$, as asserted.

We will now turn to the proof that there is a submanifold $K_1$ with
the stated properties. 

Since $K$ is a non-empty, proper submanifold of $N$, we have $\calt\ne\emptyset$.
If $\wt\calt=0$ then the components of $\obd(\calt)$ are incompressible toric suborbifolds of the
$3$-orbifold $\oldPsi$, a contradiction since $\oldPsi$ is strongly \simple\ (see \ref{oops}). Hence
$\wt\calt>0$, which by the definition of $\wt^*\calt$ implies $\wt^*\calt>0$. Since $\wt^*\calt$ is divisible by $\lambda_\oldPsi$ according to \ref{wait star}, we have $\wt^*\calt\ge\lambda_\oldPsi$. It follows that if
Alternative (ii) of the hypothesis holds, i.e. if 
 $\obd(\calt )$ is
incompressible in $\oldPsi$ and
$\chibar(\kish(\obd(K)))<\lambda_\oldPsi/4\le1/2$, then
$\chibar(\kish(\obd(K)))<\min(1,{\wt^*(\calt)/4})$; this is Alternative (ii) of
the hypothesis of Lemma \ref{pre-modification}. Likewise, it is
immediate that if
Alternative (ii$'$) of the hypothesis of the present lemma holds, i.e. if 
 $\obd(\calt )$ is
incompressible in $\oldPsi$ and if $\wt^*\calt\ge4$ and
$\chibar(\kish(\obd(K)))<1$,
then
$\chibar(\kish(\obd(K)))<\min(1,{\wt^*(\calt)/4})$. 
Thus either of the
Alternatives (ii) or (ii$'$) of the hypothesis of the present lemma implies 
Alternative (ii) of the hypothesis of Lemma \ref{pre-modification}.
Alternative (i) of the hypothesis of the present lemma is identical to
Alternative (i) of the hypothesis of Lemma
\ref{pre-modification}. Note also that if $K^+$ 
%satisfies the hypothesis
%of Lemma \ref{pre-modification} which asserts that it
%is not 
were homeomorphic to $\TTT^2\times[0,1]$ or to a twisted $I$-bundle over a
Klein bottle, then $K^+$ would fail to be acylindrical, and we would have
$h(K)=h(K^+)=2$, so that neither of the alternatives (a) or (b) of the
hypothesis of the present lemma would hold; thus $K^+$ is not
homeomorphic to either of these manifolds. Hence the hypotheses of the present
lemma imply those of Lemma
\ref{pre-modification}, so that under these hypotheses one of the alternatives (1)---(3) of the conclusion of Lemma
\ref{pre-modification} must hold.

Consider the case in which Alternative (1) of the conclusion of Lemma
\ref{pre-modification} holds: that is,
there is a disk $D\subset N$ with
$\partial D=D\cap \calt $, such that $D$ is in general position with respect to $\fraks_\oldPsi$ and 
there is a disk $E\subset \calt $ such that $\partial E=\partial D$ and 
$\wt E>\wt D$. Furthermore, if $\lambda_\oldPsi=2$,
we may suppose $D$ and $E$ to be chosen so that $\max(\wt(
E)-\wt(D),\wt D)\ge2$.

Since $N$ is $+$-irreducible, the sphere $D\cup
E$ bounds a ball in $\plusN$, and hence the surface $(\calt-E)\cup
D$ is isotopic to $\calt$ in $\plusN$, by an isotopy that is constant
on $\overline{\calt-E}$. Hence $\calt_1:=(\calt-E)\cup
D\subset\inter N$ bounds a submanifold $P_1$ of $\plusN$ which is isotopic to $P$ in $\plusN$, by an isotopy that is constant
on $\overline{\calt-E}$. Set $K_1=P_1\cap N$; since $\partial P_1=\calt_1\subset\inter N$, we have
$K_1^+=P_1$.

We have $\wt(\calt)-\wt(\calt_1)=\wt(E)-\wt(D)$.
Since $D$
is
in general position with respect 
to $\fraks_\oldPsi$ and
$\wt E>\wt D$,
%$D$ meets $\fraks_\oldPsi$ transversally in
%at most one point, and since
%$D\cap\fraks_\oldPsi\ge2$, 
the surface $\calt_1=\Fr_NK_1$ is 
in general position with respect to
$\fraks_\oldPsi$, and
$\wt\calt_1<\wt\calt$. Furthermore, if $\lambda_\oldPsi=2$, 
%so that 
%$\max(\wt(
%E)-\wt(D),\wt D)\ge2$,
then 
$\max(\wt(\calt )-
\wt(\calt_1)
,\wt(\calt_1))\ge \max(\wt(
E)-\wt(D),\wt D)\ge2$.
Since $P_1$ is isotopic to $P$,
the components of $
\Fr_{N}K_1=
\Fr_{N^+}P_1$ are incompressible tori. 
The existence of an isotopy between $P$ and $P_1$  implies that
$P=K^+$ is homeomorphic to $P_1=K_1^+$, and that
$\overline{N^+-P}=(\overline{N-K)}^+$ is homeomorphic to
$\overline{N^+-P_1}=(\overline{N-K_1)}^+$. Hence we have
  $h(K_1)=h(P_1)=h(P)=h(K)\ge h(K)/2$ and $h(N-K_1)=h(N^+-P_1)=h(N-P)=
  h(N-K)$. 
To verify the last sentence of the lemma in this case, note that 
%$K_1^+=P_1$ and $K^+=P$ are isotopic regardless of whether (a) or (b)
%holds; and that the existence of an isotopy between $P_1$ and $P$
%implies that $\overline{N-K_1}=\overline{N^+-P_1}$ is homeomorphic to
%$\overline{N-K}=\overline{N^+-P}$, 
$\compnum(N-K_1)=\compnum(N^+-P_1)=\compnum(N-P)=
  \compnum(N-K)$,
so that if (b) holds then $\overline{N-K_1}$ is
connected. Furthermore, our construction gives that $K_1^+=P_1$ is
isotopic to $K^+=P$ regardless of whether (a) or (b) holds.
Thus the proof of the lemma is complete in this case.

%L
%(Note
%that the hypothesis that one of the alternatives (a), (b) holds is not
%needed in this case. \redcomment{That sentence may create more confusion
%  than it resolves, in light of the last sentence of the Lemma.})

Now consider the case in which Alternative (2) of the conclusion of Lemma
\ref{pre-modification} holds: that is, 
there is 
a solid torus $J\subset P$, $\pi_1$-injective in $N^+$,
 with $\partial J\subset
K\cap\inter N
\subset K\subset P$, 
  such that 
  $\partial J\cap\partial K$ is a union of two disjoint annuli 
%
%there is 
%a $\pi_1$-injective solid torus $J\subset K^+$, 
 %with $\partial J\subset
%K\cap\inter N
%\subset K\subset K^+$, 
  %such that 
  %$(\partial J\cap(\partial K)$ is a union of two disjoint annuli
 $X$ and $X'$ contained in $\calt$, each having winding
  number $1$ in $J$, and we have
\Equation\label{harps}
X\cap\fraks_\oldPsi\ne\emptyset\quad\text{and}\quad X'\cap\fraks_\oldPsi\ne\emptyset,
\EndEquation
%and each having non-empty
%  intersection with $\fraks_\oldPsi$, 
and 
  $\partial J\cap\fraks_\oldPsi\subset \inter X\cup \inter X'$.
%there is a $\pi_1$-injective solid torus $J\subset K$ such that $J\cap\partial K$
  %is a union of two disjoint annuli 
%$A_1$ and $A_2$ contained in $\calt $, each having winding
  %number $1$ in $J$, and
  %such that 
%\Equation\label{harps}
%A_i\cap\fraks_\oldPsi\ne\emptyset\quad\text{for}\quad i=1,2
%\EndEquation
% and $
%  J\cap\fraks_\oldPsi\subset \inter A_1\cup \inter A_2$. 
Note that the set
$\Fr_{P}J=\overline{\partial J-(X\cup X')}$ is also a disjoint
union of two annuli, $A_1$ and $A_2$, which also have winding number
$1$ in $\calt$, and are properly embedded in $P$.  Since $\partial 
  J\cap\fraks_\oldPsi\subset \inter X\cup \inter X'$, we have 
\Equation\label{houris}
A_i\cap\fraks_\oldPsi
=\emptyset\quad\text{for}\quad i=1,2.
\EndEquation
Note also that since the annuli $X$, $X'$, and $A_1$ and $A_2$ have winding number
$1$ on the solid torus $J$, which is $\pi_1$-injective in $N^+$, they are themselves
$\pi_1$-injective in $N^+$.

Set $P_0=\overline{P-J}$. 
%Set $K_0=\overline{K-J}$.
 Then $P_0$ is a compact, possibly
disconnected, $3$-dimensional submanifold of $P\subset N$. Since $\partial J\subset\inter N$,  
every component of $K\cap\partial N$ is a sphere contained in either $\inter J$ or $\inter P_0$. %Hence every component of $N^+-\inter N$ is a ball contained in either $\inter J$ or $\inter P_0$. 
Thus if we set $K_0=P_0\cap N$, each component of $P_0$ is obtained from a component of $K_0$ by attaching balls along certain spheres in the boundary. Set
$\calt_0=\Fr_NK_0=\Fr_{N^+}P_0=\partial P_0$. We have $\calt_0=(\calt-(X\cup X'))\cup(A_1\cup A_2)$, so that
%It follows from
(\ref{harps}) and (\ref{houris}) imply
%that $\calt_0$ is transverse to
%$\fraks_\oldPsi$ and
% that 
\Equation\label{on travel}
\wt\calt_0\le (\wt\calt)-2.
\EndEquation

The manifold $K_1$ whose existence is asserted by the lemma will be
constructed as a suitable component of $K_0$. Before choosing a
suitable component, we will prove a number of
facts about $P_0$ and $\calt_0$ that will be useful in establishing the
properties of $K_1$.

Since $P=P_0\cup J$, and since
$P_0\cap J=\Fr_{P}J$ is the union of the disjoint annuli $A_1$ and $A_2$,
each of which has winding number $1$ in the solid torus $J$, we have
\Equation\label{erfer}
h(P)\le h(P_0)+1.
\EndEquation
%(where the second inequality is an equality if and only if $A_1$ and
%$A_2$ lie in the same component of $P_0$).

The following property of $\calt_0$ will also be needed:

\Claim\label{tony baloney}
Every component of $\calt_0$ is a (possibly compressible) torus in
  $\inter {N}$, 
in general position with respect to
$\fraks_\oldPsi$.
\EndClaim

To prove \ref{tony baloney}, we need only note that since the annuli
$X$ and $X'$ are $\pi_1$-injective in the (torus) components of
$\calt$ containing them, the components of $\overline{\calt-(X\cup
  X')}$ are tori and annuli. As the closed, orientable $2$-manifold $\calt_0$ is
obtained from the disjoint union of $\overline{\calt-(X\cup
  X')}$ with $A_1$ and $A_2$  by gluing certain boundary components
in pairs, its components have Euler characteristic $0$. Since $\calt$
is by hypothesis 
in general position with respect to $\fraks_\oldPsi$, it follows from
\ref{houris} that $\Fr_{N^+}P_0$ is in general position with respect  to
$\fraks_\oldPsi$. This establishes \ref{tony
  baloney}.

We will also need:
\Claim\label{da props}
Let $L$ be any component of $P_0$. Then:
\begin{itemize}
\item $L$ is $\pi_1$-injective in ${N^+}$;
\item every component of $\Fr_{N^+}L$ is a torus in
  $\inter {N}$, in general position with respect to
$\fraks_\oldPsi$;
\item ${\wt}(\Fr_{N^+}L)\le \wt(\calt)-2$;  and
\item if $L$ is not a solid torus then  $\Fr_{N^+}L$
  is incompressible in ${N^+}$.
\end{itemize}
\EndClaim

To prove the first assertion of \ref{da props}, note that the
frontier of $L$ in ${P}$ consists of one or both of the
$\pi_1$-injective annuli $A_1$, $A_2$. Hence
$L$ is itself $\pi_1$-injective in ${P}$. On the other hand, since the
components of $\Fr_{N^+}{P}=\calt$ are incompressible tori, ${P}$ is in turn
$\pi_1$-injective in ${N^+}$. Hence $L$ is $\pi_1$-injective in ${N^+}$.

Next note that every component of $\Fr_{N^+}L$
is a component of $\calt_0$; hence the assertion that every component of $\Fr_{N^+}L$ is a torus in
  $\inter {N^+}$, in general position with respect to
$\fraks_\oldPsi$, is an immediate consequence of \ref{tony baloney}. 

Since $\Fr_{N^+}L\subset\Fr_{N^+}P_0=\calt_0$, the inequality
$\wt(\Fr_{N^+}L)\le \wt(\calt)-2$ is an immediate
consequence of (\ref{on travel}).

Since ${N^+}$ is irreducible, and since we have seen that $L$ is
$\pi_1$-injective in $N^+$ and that every component of $\Fr_{N^+}
L$ is a torus, the fourth assertion of \ref{da props} follows from
Proposition \ref{final assertion}, applied with $M=N^+$, and with $L$ defined as above. Thus \ref{da props} is proved.

The rest of the proof of the lemma in this case is divided into subcases. Consider the
subcase in which Alternative (a) of the hypothesis holds. 
Since $P$ is
acylindrical in this subcase, the
mutually parallel annuli $A_1$ and $A_2$
cannot be essential in $P$.
Since these annuli are $\pi_1$-injective in
${N^+}$, 
%they do not meet sphere components of $\partial {P}$, and are therefore
%properly embedded in $P$. Since
they are in particular
$\pi_1$-injective in $P$;
% and are non-essential, 
they must therefore be
boundary-parallel in $P$. 
It follows that some component $P_1$
%$K_1^+$ 
of
$P_0$ is isotopic to $P$ in $N^+$. 
%We fix such a component of $P_0$
%and denote it by $P_1$. 
We set $K_1=P_1\cap N$.
It follows from \ref{da props},
taking $L=P_1$, that 
$\wt(\Fr_{N}K_1)=\wt(\Fr_{N^+}P_1)\le \wt(\calt)-2$.  
Since $P_1$ is isotopic to $P$,
the components of $
\Fr_{N}K_1=
\Fr_{N^+}P_1$ are incompressible tori; by \ref{tony baloney} they are
in general position with respect to $\fraks_\oldPsi$.
The existence of an isotopy between $P$ and $P_1$ also implies that
%\redmissingref{Maybe it's better to be more explicit here. I think we
  $h(K_1)=h(P_1)=h(P)=h(K)\ge h(K)/2$ and that $h(N-K_1)=h(N^+-P_1)=h(N-P)=
  h(N-K)$. Thus the conclusions of the lemma (including the last sentence, which in this subcase asserts that $P=K^+$ and $P_1=K_1^+$ are isotopic) hold with this choice of $K_1$.

Next consider the
subcase in which Alternative (b) of the hypothesis holds and $K_0$ is
connected. In this subcase we take $K_1=K_0$. According to (\ref{erfer})
we have $h(K_1)=h(P_0)= h(P)-1= h(K)-1$; since Alternative (b) of the hypothesis
gives $h(K)\ge3$, it follows that $h(K_1)\ge h(K)/2$, as asserted in
the conclusion of the lemma. It also follows that $h(K_1)\ge2$, so that
$K_1^+$ is not a solid torus. Hence by \ref{da props}, the components of
$\calt_1:=\Fr_{N} K_1=\Fr_{N^+} P_1$ are incompressible tori in
${N^+}$, in general position with respect to $\fraks_\oldPsi$. 

It follows from \ref{da props} that  
$\wt\calt_1\le \wt(\calt)-2$;
%Since $\Fr_{N^+}L\subset\Fr_{N^+}P_0=\calt_0$, it follows from  (\ref{on travel}) that
%the inequality
%$\wt(\Fr_{N^+}L)\le \wt\calt-2$. 
%This in turn trivially implies that $\Fr_NK_1$ has weight strictly less than $\wt\calt$, and that $\max(\wt(\calt)-
%\wt(\calt_1)
%,\wt(\calt_1))\ge
%\wt(\calt)-
%\wt(\calt_1)\ge
%2$, regardless of whether $\lambda_\oldPsi$ is $1$ or $2$; 
in particular, the second bullet point of the conclusion of the lemma holds in this subcase.

%It follows from \ref{not really final} that $h({N}-K_1)=h({N^+}-P_0)\ge
%h({N^+}-P) =h({N}-K)$. Furthermore, 
Since Alternative (b) holds, the manifold ${N}-K$
is connected, and hence so is ${N^+}-P$.
Since
$\overline{N^+-P_0}=\overline{N^+-P}\cup J$, and since
$\overline{N^+-P}\cap J$ is the union of the disjoint annuli $X$ and $X'$,
each of which has winding number $1$ in the solid torus $J$, the manifold 
%we have:
%\Claim\label{not really final}
${N^+}-P_0$ is also connected, and $h({N^+}-P_0)\ge
h({N^+}-P)$. This means that ${N}-K_1$ is connected, as required by the last sentence of the lemma when (b) holds. and that $h({N}-K_1)\ge
h({N}-K)$.
Thus all
the conclusions of the lemma are established in this subcase.

The remaining subcase of this case is the one in which Alternative (b) of the hypothesis holds and $K_0$ is
disconnected. In this subcase, $K_0$ has two components since $\Fr_PJ$
has two components $A_1$ and $A_2$. For $i=1,2$, let $K_i$ denote
the component of $K_0$ containing $A_i$. After re-indexing if
necessary we may assume that $h(K_1)\ge h(K_2)$. In this subcase,
$P_0$ has two components, and they may be indexed as $P_1$ and $P_2$,
where $P_i$ is the union of $K_i$ with certain ($3$-ball) components
of $\overline{N^+-N}$. Since the $A_i$
have winding number $1$ in $J$, the manifold $P$ is homeomorphic to
the union of $P_1$ and $P_2$ glued along an annulus. Hence $h(K)=h(P)\le
h(P_1)+h(P_2) =h(K_1)+h(K_2)\le 2h(K_1)$, i.e. $h(K_1)\ge h(K)/2$, as asserted in
the conclusion of the lemma. Since Alternative (b) of the hypothesis
gives $h(K)\ge3$, it follows that $h(K_1)\ge2$, so that
$K_1^+$ is not a solid torus. Hence by \ref{da props}, the components of
$\calt_1:=\Fr_NK_1=\partial K_1$ are incompressible tori in ${N^+}$, in general position with respect to $\fraks_\oldPsi$. It also follows from \ref{da props} that
$\wt\calt_1=\wt(\Fr_{N^+}P_1)\le \wt(\calt)-2$.
%that ${N^+}-L$ is connected and that $h({N^+}-L)\ge h({N^+}-K)$. 

To establish the conclusions of the lemma in this subcase, it remains to
prove that
$\overline{N-K_1}$ is connected and that $h(\overline{N-K_1})\ge h(\overline{N-K})$. 
Since Alternative (b) holds, 
the manifold $\overline{N-K}$
is connected, and hence so is $\overline{N^+-P}$. We have
$\overline{N^+-P_0}=\overline{N^+-P}\cup J$, where $J$ is connected, and
$\overline{N^+-P}\cap J=X\cup X'\ne\emptyset$; hence $\overline{N^+-P_0}$ is
connected, and therefore so is $\overline{N-K_0}$. Next note that
%has non-empty intersection with $\overline{N-P}$
%the manifold
%${N^+}-P={N}-K$ is connected. \redmissingref{I should check the equality 
%${N^+}-P={N}-K$.}
%$\overline{{N^+}-P_0}=\overline{{N^+}-K}\cup J$, and $\overline{{N^+}-K}\cap
%J=X\cup X'\ne\emptyset$. The solid torus $J$ is connected, and
%$\overline{{N^+}-K}$ is connected because (b) holds. Hence
%$\overline{{N^+}-P_0}$ is connected. 
%We have
$\overline{{N}-K_1}=\overline{{N}-K_0}\cup K_2$, where $K_2$ is
connected, and has non-empty intersection with $\overline{{N}-K_0}$
because it is a component of $K_0$. Hence $\overline{{N}-K_1}$ is
connected, as required by the last sentence of the lemma when (b) holds.

To estimate 
%$h({N^+}-K_1)=
$h(\overline{{N}-K_1})$, first note that each component of
$P_0=\overline{P-J}$ contains either $A_1$ or $A_2$.
% a component of $\overline{\calt-(X\cup X')}$. 
Since
$P_0=\overline{P-J}$ is disconnected in this subcase, 
$A_1$ and $A_2$ must lie in distinct components of $P_0$. If $X$ and
$X'$ are contained in distinct components of $\calt$, and if $T$
denotes the component of $\calt$ containing $X$, then the 
containing the annulus $\overline{T-X}\subset P_0$ meets both $A_1$ and
$A_2$, a contradiction.
%$\overline{\calt-(X\cup X')}$ is also disconnected. But each of the
%annuli $X$, $X'$ is $\pi_1$-injective in ${N^+}$, and therefore is
%non-separating in the (torus) component of $\calt$ containing it. 
Hence
%$P_0=\overline{P-J}$ is disconnected in this su, n
$X$ and $X'$
must lie on the same component of $\calt$. As disjoint, homotopically
non-trivial annuli on a torus, $X$ and $X'$ are homotopic in $\calt$ and hence
in $\overline{{N^+}-P}$. We have written $
%\overline{{N}-K_0}=
\overline{{N^+}-P_0}$ as the
union of $\overline{{N^+}-P}$ with the solid torus $J$, and the
intersection of $\overline{{N^+}-P}$ with $J$ consists of the two
annuli 
$X$ and $X'$, which have winding number $1$ in $J$. As these annuli
are homotopic in $\overline{{N^+}-P}$, we have
$h(\overline{{N^+}-P_0})=h(\overline{{N^+}-P})+1$. On the other hand, we have
$\overline{{N^+}-P_1}=\overline{{N^+}-P_0}\cup P_2$, and the components of
$ \overline{{N^+}-P_0}\cap P_2$ are components of $\Fr_{N^+}P_2$ and are
therefore tori by \ref{da props}. Applying Lemma \ref{another goddam
  torus lemma}, taking $U=\overline{{N^+}-P_0}$ (which has been seen to be
connected) and $V=P_2$, we
deduce that $h(\overline{{N^+}-P_1})\ge
h(\overline{{N^+}-P_0})-1$. Hence $h(\overline{{N^+}-P_1})\ge
h(\overline{{N^+}-P})$, or equivalently $h(\overline{{N}-K_1})\ge h(\overline{{N}-K})$. 
Thus
the conclusions of the lemma are established in this subcase.

Finally, consider the case in which Alternative (3) of the conclusion of Lemma
\ref{pre-modification} holds: that is, 
%there is a $\pi_1$-injective
%solid torus $J\subset K$ such that $J\cap\partial K$
  %is a single annulus  $A\subset\calt $, having winding
  %number $2$ in $J$, and
  %$\emptyset\ne(\partial J)\cap\fraks_\oldPsi\subset A$.
there is 
a solid torus $J\subset P$, $\pi_1$-injective in $N^+$,
 with $\partial J\subset K\cap\inter N\subset K\subset P$,   such that 
  $\partial J\cap\partial K$ is an annulus $X\subset\calt$, having winding
  number $1$ or $2$ in $J$, and we have
  $\wt (\partial J) \ge\lambda_\oldPsi$ and $\partial J\cap\fraks_\oldPsi\subset \inter X $.

In this case,
  $A:=\overline{\partial J-X}$ is a properly
 embedded,
  $\pi_1$-injective annulus in $P$, also having winding number $1$ or $2$ in $J$. If Alternative (a) of the
  hypothesis holds, then $A$ must be boundary-parallel in $P$; that
  is, some component $Z$ of $\overline{P-A}$ must be a solid torus in
  which $A$ has winding number $1$. Thus $\overline{P-Z}$ is isotopic
  to $P$.
If  $Z=\overline{P-J}$, then $P$ is homeomorphic to $J$
  and is therefore a solid torus, a contradiction to the
  incompressibility of $\calt$.
Hence we must have $Z=J$, so that $P$ is isotopic to $\overline{P-J}$. In this subcase, we set
  $K_1=(\overline{P-J})\cap N$. Setting $\calt_1=\Fr
  K_1=\partial(\overline{P-J})$, we have $\calt_1=(\calt-X)\cup A$, 
  so that $\wt\calt_1=\wt\calt-\wt X+\wt A$; furthermore, $\calt_1$ is
  in general position with respect to $\fraks_\oldPsi$, because $\calt_1$ is
  in general position with respect to $\fraks_\oldPsi$ and $\partial
  J\cap\fraks_\oldPsi\subset \inter X $. Since   $\wt (\partial J)
  \ge\lambda_\oldPsi$ and $\partial J\cap\fraks_\oldPsi\subset \inter X $, we have
%It follows from (\ref{i
  %  too}) that
 $\wt X\ge\lambda_\oldPsi$ and $\wt A=0$. Hence
  $\wt\calt_1\le\wt\calt-\lambda_\oldPsi$; in particular, $\wt\calt_1<\wt\calt$, and if $\lambda_\oldPsi=2$ then  
$\max(\wt(\calt )-
\wt(\calt_1)
,\wt(\calt_1))\ge
\wt(\calt )-
\wt(\calt_1)\ge
\lambda_\oldPsi=2$.
The existence of an isotopy between $P$ and $P_1$ also implies that
%\redmissingref{Maybe it's better to be more explicit here. I think we
  $h(K_1)=h(P_1)=h(P)=h(K)\ge h(K)/2$ and that $h(N-K_1)=h(N^+-P_1)=h(N-P)=
  h(N-K)$. Thus the conclusions of the lemma (including the last sentence, which in this subcase asserts that $P=K^+$ and $P_1=K_1^+$ are isotopic) hold with this choice of $K_1$.

There remains the subcase in which Alternative (b) holds, i.e. 
$\overline{{N}-K}$ is connected and
$h(K)\ge3$. 
In this subcase we will set $P_1=\overline{P-J}$ and $K_1={P_1}\cap N$. Then $P_1$ is connected since $P$ and $A=\Fr_PJ$ are connected; hence $K_1$ is connected. Since $P_1\cup J=P$, and
since $P_1\cap J=A$ is connected, we have $h(P)\le
h(P_1)+h(J)=h(P_1)+1$, so that $h(P_1)\ge h(P)-1=h(K)-1$. Since $h(K)\ge3$ it
follows that $h(P)\ge2$.

Note that $P$ is $\pi_1$-injective in ${N^+}$ since its frontier is incompressible, and that $P_1$ is $\pi_1$-injective in
$P$ since its frontier $A$ is $\pi_1$-injective. Hence $P_1$ is
$\pi_1$-injective in ${N^+}$, i.e. $K_1$ is
$\pi_1$-injective in ${N}$. Since $\calt_1:=\Fr_{N} K_1=\Fr_{N^+} P_1=\overline{\calt
-X}\cup A$, the surfaces $\calt$ and $\calt_1$ are homeomorphic; in
particular, the components of $\calt_1$ are tori. 
Furthermore, $\calt_1$ is
  in general position with respect to $\fraks_\oldPsi$, because $\calt_1$ is
  in general position with respect to $\fraks_\oldPsi$ and $\partial
  J\cap\fraks_\oldPsi\subset \inter X $.
Since
$h(P_1)\ge2$, the manifold $P_1$ is not a solid torus. It therefore
follows from Lemma \ref{final assertion}, applied with $N^+$ and $P_1$ playing the respective roles of  $M$ and $L$, that  $\calt_1$ is incompressible in $N^+$, and hence in $N$. This proves the first
  property of $K_1$ asserted in the conclusion of the present
  lemma. 
Next note that since   $\wt (\partial J) \ge\lambda_\oldPsi$ and $\partial J\cap\fraks_\oldPsi\subset X $,
%  $\emptyset\ne \partial J\cap\fraks_\oldPsi\subset X$, 
we have
  $\wt X\ge\lambda_\oldPsi$ and   $\wt A=0$. Hence 
$\wt\calt_1\le\wt(\calt)-\lambda_\oldPsi=\wt(\calt)-\lambda_\oldPsi$, which implies that
$\wt\calt_1<\wt\calt$, and that
$\max(\wt(\calt )-
\wt(\calt_1)
,\wt(\calt_1))\ge \wt(\calt )-
\wt(\calt_1)\ge2$ if $\lambda_\oldPsi=2$.
 The inequality 
$h(K_1)\ge h(K)/2$ holds because $h(K_1)\ge h(K)-1$ and
$h(K)\ge3$. Since ${N}-K$ is connected by Alternative (b), ${N^+}-P$ is connected. Since, in addition, $J$ is a solid torus and $\overline{{N^+}-P}\cap J=X$ is an annulus,
the manifold $\overline{{N^+}-P_1}=\overline{{N^+}-P}\cup J$ is connected, so that $N-K_1$ is connected, as required by the last sentence of the lemma when (b) holds; and
$h(\overline{{N^+}-P_1})$ is equal either to $h(\overline{{N^+}-P})$ or to
$h(\overline{{N^+}-P})+1$. In particular, $h({N}-K_1)=h({N^+}-P_1)\ge
h({N^+}-P)=h({N^+}-K)$. Thus all the conclusions are seen to hold in this final subcase.
\EndProof

%boundary \partial{pre-mod PPP K_0 N K^+K_1 XA_1A_2A_1'
%pUVWZ\wt P_1 M L^+ T
%p^* \wt transverse

\abstractcomment{\tiny
The inequality involving $\sigma$ in the hypothesis
could be weakened by replacing $\sigma$ by the corresponding thing
involving only the contribution of the tori, not the spheres. This may
mean all the results of this section remain true (and are a little
stronger) if $\sigma$ is replaced by that variant. I need to think
about whether this affects the final conclusions of the paper.

In the proof in the case where (3) holds, I could easily have gotten
$p-2$ as the upper bound for $\wt\calt_1$, since the
intersection with $J$ is an arc which must have two endpoints. The
only case where one doesn't get this is the case in which (1) holds,
but I had an idea that the non-existence of hyperbolic triangle groups might give
it there as well. This means that one may get much more info than I
thought without the assumption of no nodes.}

\Lemma\label{get a newer mop}
Let $m\ge2$ be an integer. Let $\oldPsi$ be a strongly \simple,
boundary-irreducible, orientable  $3$-orbifold containing no embedded
negative turnovers. 
%, such that each component of $\fraks_\oldPsi$ is either an arc
%or a simple closed curve. 
Set $N=|\oldPsi |$. Suppose that each
component of $\partial N$ is a sphere, that $\plusN$ is a graph manifold, and that $h(N)\ge\max(4m-4,m^2-m+1)$.
Then
there is a compact
submanifold $K$ of $N$ such that
\begin{itemize} 
\item each component of $\Fr_N K $ is an incompressible torus in
  $\inter N $, in general position with respect to $\fraks_\oldPsi$,
\item $K$ and $N -K$ are connected, 
\item $\min(h(K),h(N-K))\ge m$, 
\item  $\obd(\Fr_N K )$ is incompressible in $\oldPsi$ (so that $\obd(K)$ is boundary-irreducible and strongly \simple\ by Lemma \ref{oops lemma}, and hence $\kish(\obd(K))$ is defined in
view of \ref{tuesa day}), and
\item either \begin{enumerate}[(A)]
\item $\chibar(\kish(\obd(K)))\ge1$, or 
\item  $\chibar(\kish(\obd(K)))\ge\lambda_\oldPsi/4$ and
$\wt^*
(\Fr_NK)<4$. 
\end{enumerate}
\end{itemize}
\EndLemma

\Proof
Set $M=\plusN$. 
%We will first prove assertion (1) in the case where $M$
%is a graph manifold. 
Since $h(M)=h(N)\ge\max(4,m^2-m+1)$, the hypothesis of Lemma
\ref{just right} holds.
% with $m=2$. 
Hence there is a compact
submanifold $P$ of $M$ such that (I) each component of $\Fr_M P=\partial P$ is an incompressible torus in
  $M$, (II) $P$ and $M-P$ are connected, and
(III) $\min(h(P),h(M-P))\ge m$. 
In particular, (III) implies that $P$ is a proper, non-empty
submanifold of $M$, so that $\Fr_MP\ne\emptyset$.
After an isotopy we may assume that (IV)
$\Fr_MP $ is contained in $\inter  N$ and is in general position with respect to $\fraks_\oldPsi$. Among all compact
submanifolds of $M$ satisfying (I)---(IV), let us choose $P$ so
as to minimize the quantity   $\wt( \Fr_MP)$. 

Set $K=P\cap N$, so that $P=K^+\subset N^+=M$.
Set 
$\calt=\partial P=\Fr_MP=\Fr_NK$. Then it follows from (I)---(IV) that
each component of $\calt$ is an incompressible torus in
  $\inter N $,  in general position with respect to $\fraks_\oldPsi$; that $K$ and $N -K$ are connected; and that $\min(h(K),h(N-K))\ge m$.

In view of Condition (II), the hypotheses of Lemma \ref{even easier}
hold with $Y=M$ and with $\calt$ defined as above. Hence $h(M)\le
h(P)+h(\overline{M-P})-1$. Since $h(M)\ge 4m-4$, 
we have $h(P)+h(\overline{M-P})\ge4m-3$, so that at least one of the
integers $h(P)$ and $h(\overline{M-P})$ is at least $2m-1$. 
Since Conditions (I)---(IV), and the value of $\wt\calt$, are unaffected
when $P$ is replaced by $\overline{M-P}$, we may assume $P$ to have been chosen so that
$h(P)\ge2m-1$.
%, and since $m\ge2$ we then have $h(P)\ge m$. 
We will
show that with this choice of $P$, the suborbifold $\obd(\calt)$
is incompressible in $\oldPsi$, and  one of the Alternatives (A) or (B)
of the statement holds. This will imply the conclusion of the  lemma.
%K

This step is an application of Lemma \ref{modification}.  According to
the hypotheses of the present lemma, $\oldPsi$ is a compact, orientable,
strongly \simple\ $3$-orbifold containing no embedded negative turnovers, %each component of $\fraks_\oldPsi$ is either an arc
%or a simple closed curve, 
and each boundary component of $\partial N$
is a sphere. Since $M=\plusN$ is a graph manifold, it is by definition
irreducible, so that $N$ is $+$-irreducible. We have seen that $K\subset N$ is a compact,
connected $3$-manifold, that the components
of $\Fr_N K$ are all incompressible tori in $\inter N$, in general position with respect to
$\fraks_\oldPsi$, that $\overline{N-K}$ is connected,
% (because $\overline{M-K}$ is connected and $\calt\subset\inter M$),
and that $h(K)\ge2m-1$; since $m\ge2$, we have in particular that
$h(K)\ge3$. This gives Alternative (b) of the hypothesis of Lemma
\ref{modification}.

Now assume, with the aim of obtaining a contradiction, that
either some component of $\obd(\calt)$ is compressible in $\oldPsi$, or
that  $\obd(\calt)$ is incompressible in $\oldPsi$
but that both the Alternatives (A) or (B) of the conclusion of the present lemma are false. If
$\obd(\calt)$ has a compressible component in $\oldPsi$ then
Alternative (i) of Lemma \ref{modification} holds.
 If  $\obd(\calt)$ is incompressible in $\oldPsi$,
but (A) and (B) are
both false, then $\chibar(\kish(\obd(K)))<1$, and either
$\chibar(\kish(\obd(K)))<\lambda_\oldPsi/4$ or
$\wt^*\calt\ge4$. If 
$\chibar(\kish(\obd(K)))<\lambda_\oldPsi/4$ then Alternative (ii) of Lemma
\ref{modification} holds. If $\chibar(\kish(\obd(K)))<1$ and $\wt^*\calt\ge4$,
then 
%since $p$ is even by
%\redmissingref{replace by a $\lambda$-reference}, 
%and has been seen to be strictly positive,
%we have $p\ge4$, and 
then Alternative (ii$'$) of Lemma \ref{modification} holds. Thus in any
event, our assumption implies that one of the alternatives (i), (ii)
or (ii$'$)
of Lemma \ref{modification} holds, and hence that there is a compact,
connected $3$-manifold $K_1\subset N$ having the properties stated in
the conclusion of that lemma. According to the last sentence of Lemma
\ref{modification}, since Alternative (b) holds, we may choose
$K_1$ so that $\overline{N-K_1}$ is connected.

The conclusion of Lemma \ref{modification} gives 
$h(K_1)\ge h(K)/2$. Since $h(K)\ge2m-1$ it follows that $h(K_1)\ge
m$. The conclusion of Lemma \ref{modification} also gives  
%Furthermore, we have and hence $h(K_1)\ge2$ since $h(K_1)$ is
%an integer; and
$h(N-K_1)\ge h(N-K)$. Now since $K$ satisfies Condition (III)
above, we have $h(N-K)\ge m$; hence $h(N-K_1)\ge m$. 
Thus Condition (III) holds when $P$ is
replaced by $P_1:=K_1^+$. It is immediate from the conclusion of Lemma
\ref{modification} (including the connectedness of $\overline{N-K_1}$) that Conditions (I), (II) and
(IV) above also hold when $P$ is replaced by $P_1$. But since Lemma \ref{modification} also gives 
$\wt(\Fr_M P_1)=\wt(\Fr_N K_1)< \wt\calt$, this contradicts the minimality
of $\wt\calt$. 
\EndProof
%p>0 Kp^*

\Proposition\label{crust chastened}
Let $\oldPsi$ be a compact,  orientable, strongly \simple, boundary-irreducible
  $3$-orbifold containing no embedded negative turnovers.
%, such that each component of $\fraks_\oldPsi$ is either an arc
%or a simple closed curve. 
Set $N=|\oldPsi |$. Suppose that
$N$ is $+$-irreducible, and that each
component of $\partial N$ is a sphere, so that $\plusN$ is closed.
Then:
\begin{enumerate}
\item if $\plusN$ contains an incompressible torus and $h(N)\ge4$, we
  have $\delta(\oldPsi)\ge3\lambda_\oldPsi$ (see \ref{t-defs}); 
and
\item if $\plusN$ is a graph manifold, and if $h(N)\ge8$ and $\lambda_\oldPsi=2$, we have
  $\delta(\oldPsi)\ge12$.
\end{enumerate}
\EndProposition

\Proof
To prove Assertion (1), suppose, in addition to the general hypotheses of the proposition,  that
$h(N)\ge4$ and that $\plusN$ contains an incompressible torus. We claim:
\Claim\label{Henry F. Schricker}
There is a compact, connected submanifold $K$ of $N$, whose frontier
components are incompressible tori
%$\calt:=\Fr_NK$ is 
in general position with respect to $
\fraks_\oldPsi$, such that 
 $\oldPi:=\obd(\Fr_N K)$ is
incompressible in $\oldPsi$, and
$\chibar(\kish(\obd(K)))\ge\lambda_\oldPsi/4$.
\EndClaim

We will first prove \ref{Henry F. Schricker} in the case where
$\plusN$ is a graph
  manifold. In this case, the hypotheses of Lemma \ref{get a newer
    mop} hold with
$m=2$. Hence there is  a compact
submanifold $K$ of $N$ such that the conclusions of Lemma \ref{get a
  newer mop} hold (with $m=2$). 
In particular, $\Fr_NK$ is 
 in general position with respect to $
\fraks_\oldPsi$, the components of $\Fr_NK$ are incompressible tori in $\inter N$, and  $\obd(\Fr_N K)$ is incompressible in $\oldPsi$. Since one of the alternatives (A) or (B) of
Lemma \ref{get a newer mop} holds (and since $\lambda_\oldPsi\in\{1,2\}$), we have
%\Equation\label{John T. 69McCutcheon}
$\chibar(\kish(\obd(K)))\ge\lambda_\oldPsi/4$, and the proof of (\ref{Henry
  F. Schricker}) is complete in this case.

To prove  (\ref{Henry
  F. Schricker}) in the case where $\plusN$ is not a graph
manifold, we note that since $\plusN$ contains an incompressible torus and is irreducible, it is a Haken manifold. We let  $\Sigma$ denote the characteristic submanifold of $\plusN$ in the sense of \cite{js} (cf. \ref{manifolds are different}), and we note that in this case, according to the definition of a graph manifold, there is a component $L$ of $\overline{\plusN-\Sigma}$ which is not homeomorphic to $\TTT^2\times[0,1]$. Since $\partial L\subset\partial\Sigma$, each component of $\partial L$ is an incompressible torus in $\plusN$. 

If $L$ admits a Seifert fibration, it follows from the definition of the characteristic submanifold \cite[p. 138]{js} that the inclusion map $L\to\plusN$ is homotopic to a map of $L$ into $\Sigma$. Since $\partial L$ is incompressible, it follows that the identity map of $L$ is homotopic in $L$ to a map of $L$ into $\partial L$. This implies by \cite[Lemma 5.1]{Waldhausen} that $L$ is homeomorphic to $\TTT^2\times[0,1]$, a contradiction. Hence $L$ does not admit a Seifert fibration.

If $V$ is an incompressible torus in $\plusN$, the definition of the characteristic submanifold implies that the inclusion $V\to\plusN$ is homotopic to a map of $V$ into $\Sigma$. Since $\plusN$ contains at least one incompressible torus, we have $\Sigma\ne\emptyset$.

If $V$ is an incompressible torus in $\inter L$, then since the inclusion map $j:V\to\plusN$ is homotopic to a map of $V$ into $\Sigma$, and since $\partial L$ is incompressible, it follows that $j$ is homotopic in $L$ to a map of $V$ into $\partial L$, which by $\pi_1$-injectivity can be taken to be a covering map from $V$ to a component $V'$ of $\partial L$. It then follows by \cite[Lemma 5.1]{Waldhausen} that $V$ and $V'$ are parallel in $L$. This shows that every incompressible torus in $\inter L$ is boundary-parallel in $L$.

Since $L$ does not admit a Seifert fibration, and every incompressible torus in $\inter L$ is boundary-parallel in $L$, it follows from Lemma \ref{torus goes to cylinder} that $L$ is acylindrical. 
After an isotopy we may assume that $\partial L\subset\inter
N$. Then $(L\cap N)^+=L$. Note also that $L$ is a proper submanifold of $\plusN$, since it is a component of $\overline{\plusN-\Sigma}$ and since $\Sigma\ne\emptyset$. In particular, there exists a compact, proper submanifold
$K$ of $N$ such that the components of $\Fr_NK$ are incompressible
tori, and $K^+$ is acylindrical. We
may choose such a $K$ so that $\Fr_NK$ is contained in $\inter N$ and  in general position with respect to
$\fraks_\oldPsi$, and so that, for every
submanifold $K'$ of $N$ such that $(K')^+$ is isotopic to $K^+$ 
in $N^+$, and $\Fr_NK'$ is  contained in $\inter N$ and in general position with respect to
$\fraks_\oldPsi$, we have $\wt(\Fr_NK')\ge \wt(\Fr_N K)$. Set $\calt=\Fr K$ and $\oldPi=\obd(\calt)$.
 
It now suffices to prove that
 $\oldPi$ is
incompressible in $\oldPsi$, and that $\chibar(\kish(\obd(K)))\ge\lambda_\oldPsi/4$.
Suppose to the contrary that  $\oldPi$ fails to be
incompressible in $\oldPsi$, or that it is  incompressible but that $\chibar(\kish(\obd(K)))<\lambda_\oldPsi/4$.
This means that one of
the alternatives (i) or (ii) of the hypothesis of Lemma
\ref{modification} holds. Since $K^+$ is acylindrical, Alternative (a)
of \ref{modification} also holds. Hence there is a submanifold $K_1$
that satisfies the conclusions of Lemma \ref{modification}. In particular, $K_1^+$ is
isotopic to $K^+$ in $N^+$, and $\Fr K_1$ is 
contained in $\inter N$ and
 in general position with respect to $\fraks_\oldPsi$, and has weight strictly less than
$\wt\calt$. This contradiction to the minimality
of $\wt\calt$ completes the proof of \ref{Henry F. Schricker}.
%\Fr_NK

%\redmissingref{I took out a \%-ed out passage here. It can be found in right.tex.}

Now let $K$ and  $\oldPi$ be given by \ref{Henry F. Schricker}. Since $\oldPi$ is incompressible in
$\oldPsi$ by \ref{Henry F. Schricker}, the definitions of $\delta(\oldPsi)$ and $\sigma(\oldPsi\cut\oldPi)$ (see \ref{t-defs})
imply that 
%\Equation\label{Wilbur Shaw}
%\delta(\oldPsi)\ge\sigma(\oldPsi\cut\oldPi)=
%2\lfloor6\chibar(\kish(\oldPsi\cut\oldPi))\rfloor
%=
%2\lfloor6(
%\chibar(\kish(\obd(K)))+\chibar(\kish(\obd(\overline{N-K})))) \rfloor
%12\chibar(\kish(\oldPsi\cut\oldPi))=12(\chibar(\kish(\obd(K)))+\chibar(\kish(\obd(\overline{N-K})))).
%\EndEquation
%and
\Equation\label{mr greengrass}
\delta(\oldPsi)\ge\sigma(\oldPsi\cut\oldPi)=
12\chibar(\kish(\oldPsi\cut\oldPi))
=
12 (
\chibar(\kish(\obd(K)))+\chibar(\kish(\obd(\overline{N-K})))).
\EndEquation
Since 
$\chibar(\kish(\obd(K)))\ge\lambda_\oldPsi/4$ by \ref{Henry F. Schricker}, and
$\chibar(\kish(\obd(\overline{N-K})))\ge0$ by Corollary \ref{less than nothing}, it follows from (\ref{mr greengrass}) that 
$\delta(\oldPsi)\ge3\lambda_\oldPsi$, and Assertion (1) is proved. 
%If in addition we have $\lambda_\oldPsi=2$, then $\chibar(\kish(\obd(K)))\ge1/2$, which with (\ref{mr greengrass}) and the inequality
%$\chibar(\kish(\obd(\overline{N-K})))\ge0$ gives
%$\delta(\oldPsi)\ge6$, so that Assertion (2) is proved as well. 
%\bigg

Let us now turn to Assertion (2).
%and (4). 
Suppose, in addition to the general hypotheses of the proposition,  that
$\plusN$ is a graph manifold, that $\lambda_\oldPsi=2$, and that $h(N)\ge8$. Then the hypotheses of Lemma \ref{get a newer mop} hold with
$m=3$. Hence there is  a compact
submanifold $K$ of $N$ such that the conclusions of Lemma \ref{get a
  newer mop} hold with $m=3$. Thus $K$ and $\overline{N-K}$ are
connected, and
$\min(h(K),h(N-K))\ge 3$. Furthermore, $\calt:=\Fr_NK$ is 
 in general position with respect to $\fraks_\oldPsi$; the
 $2$-orbifold $\oldPi:=\obd(\calt)$ is
 incompressible in $\oldPsi$; and either (A) $\chibar(\kish(\obd(K)))\ge1$, or (B)  $\chibar(\kish(\obd(K)))\ge\lambda_\oldPsi/4=1/2$ and
$\wt^*\calt<4$. 
By
the
definitions of $\delta(\oldPsi)$ and $\sigma(\oldPsi\cut\oldPi)$, the inequality (\ref{mr greengrass}) holds
in this context.
If (A) holds, then since $\chibar(\kish(\obd(\overline{N-K})))\ge0$,
the right hand side of (\ref{mr greengrass}) is bounded below by $12$; hence
$\delta(\oldPsi)\ge12$, so that the conclusion of (2) is true when (A) holds.

Now suppose that (B) holds. We claim:
%and if it happens that
\Equation\label{football heh heh}
\chibar(\kish(\obd(\overline{N-K})))\ge1/2.
\EndEquation
% This
%will imply that the right hand side of (\ref{Wilbur
  %Shaw}) is again at least $12$, so that
%$\delta(\oldPsi)\ge12$ as required.

To prove (\ref{football heh heh}), suppose that 
%(B) holds but that
$\chibar(\kish(\obd(\overline{N-K})))<1/2$. Since (B) holds we have
$\wt^*\calt<4$; since in addition $\wt^*\calt$ is divisible by $\lambda_\oldPsi=2$ by \ref{wait star}, we have $\wt^*\calt\le2$. 
Set
$K^*=\overline{N-K}$. We will apply Lemma \ref{modification}, with
$K^*$ playing the role of $K$ in that lemma. We have observed that
$K^*$ is connected, and that the components of $\calt=\Fr_NK^*$ are
incompressible tori,  in general position with respect to $\fraks_\oldPsi$. Since
$K=\overline{N-K^*}$ is connected, and $h(K^*)\ge3$, Alternative (b) of
the hypothesis of Lemma \ref{modification} holds with $K^*$ playing the role of $K$ in that lemma. Since  $\oldPi:=\obd(\Fr_NK^*)$ is incompressible in $\oldPsi$, and since we
have assumed that $\chibar(\kish(\obd(\overline{N-K})))<1/2=\lambda_\oldPsi/4 $,
Alternative (ii) of the hypothesis of Lemma \ref{modification}
holds. Hence there is a compact, connected $3$-manifold $K_1^*\subset N$
such that the conclusions of Lemma \ref{modification} hold
%, with $p=2$
%and 
with $K_1^*$ and $K^*$ playing the roles of $K_1$ and $K$
respectively. In particular, $K_1^*$ is a proper, non-empty
submanifold of $N$, so that $\calt_1:=\Fr_NK_1^*\ne\emptyset$; every
component of $\Fr_N K_1^*$ is an incompressible toric suborbifold of $\inter N$;
and we have $\wt\calt_1<\wt\calt\le2$. Since $\lambda_\oldPsi=2$,
Lemma \ref{modification} also guarantees that $\max(\wt(\calt )-
\wt(\calt_1),\wt(\calt_1))\ge2$. 
Since $\wt\calt_1<2$, it follows that $\wt(\calt _1)\le
\wt(\calt)-2\le0$,
i.e. $\calt_1\cap\fraks_\oldPsi=\emptyset$. Thus any
component of $\obd(\calt_1)=\calt_1\ne\emptyset$ is an
incompressible torus in $\oldPsi$; this is a contradiction since $\oldPsi$
is strongly \simple (see \ref{oops}). This
  completes the proof of (\ref{football heh heh}). 

It follows from (\ref{football heh heh}), together with the inequality
  $\chibar(\kish(\obd(K)))\ge1/2$ which is contained in Condition (B),
  that the right hand side of (\ref {mr greengrass}) is again at least $12$, so that
$\delta(\oldPsi)\ge12$ as required.
\EndProof

\Lemma\label{get the basin}
Let $\oldPsi$ be a compact,  orientable, strongly \simple, boundary-irreducible
$3$-orbifold containing no embedded negative turnovers. Suppose that
$\lambda_\oldPsi=2$. 
% such that each component of $\fraks_\oldPsi$ is either an arc
%or a simple closed curve. 
Set $N=|\oldPsi |$. Suppose that each
component of $\partial N$ is a sphere,
%that $h(N)\ge8$, 
that $N$ is $+$-irreducible, and that
there
is an acylindrical manifold $Z\subset \plusN$ such that 
\begin{itemize}
\item $\Fr_N Z$ is a single torus which is
incompressible in $\plusN$, and 
\item$h(Z)\le h(N)-2$.
\end{itemize}
Then $\delta(\oldPsi)\ge12$. \abstractcomment{\tiny I hope I have fixed up the notation in the
  application to be correct and to fit with
  the statement.}
\EndLemma

\Proof
After an isotopy we may assume that $\partial Z\subset
\inter N$. Then $(Z\cap N)^+=Z$. In particular, there exists a submanifold
$W$ of $N$ such that $W^+$ is acylindrical (and in particular connected), $\Fr_NW$ is a single torus which is
contained in $\inter N$ and is incompressible in $\plusN$, and $h(W)\le h(N)-2$.

We
may choose such a $W$ so that $T:=\Fr_NW$ is in general position with respect to
$\fraks_\oldPsi$, and  so that, for every
submanifold $W'$ of $N$ such that $\Fr_NW$ is contained in $\inter N$ and in general position with respect to
$\fraks_\oldPsi$, and such that $(W')^+$ is isotopic to $W^+$ 
in $N^+$, we have $\wt(\Fr_NW')\ge \wt T$.

We will apply Lemma \ref{modification},
letting $W$ play the role of $K$. Since $W^+$ is acylindrical, Alternative (a) of Lemma
\ref{modification} holds. If one of the alternatives (i), (ii) or (ii$'$)  of Lemma
\ref{modification} also holds, then Lemma \ref{modification} gives a
submanifold $W_1$ of $N$ such that $W_1^+$ is isotopic to $W^+$ in  $\plusN$, and $\Fr_N W_1$ is contained in $\inter N$, is in general position with respect to $\fraks_\oldPsi$ and has weight strictly less than
$\wt T$, a contradiction to the
minimality of $\wt T$. Hence Alternatives  (i), (ii), and (ii$'$) of Lemma
\ref{modification} must all fail to hold. Since Alternative  (i)
fails to hold, $\obd(T)$ is
incompressible in $\oldPsi$. Since Alternative  (ii)
fails to hold, we have
\Equation\label{obadiah}
\chibar(\kish(\obd(W)))\ge1/2.
\EndEquation
Since Alternative  (ii$'$)
fails to hold, either
$\wt_\oldPsi^* T<4$ or 
$\chibar(\kish(\obd(W)))\ge1$. But by \ref{wait star},
$\wt^*T$ is divisible by $\lambda_\oldPsi=2$; the definition of $\wt^*T$ given in \ref{wait star} also implies that $\wt T\le\wt^*T$. Hence
%$\wt^*_\oldPsi S$ is always divisible by $\lambda_\oldPsi$ and $p$ is even by  \redmissingref{replace by a $\lambda$ reference}, we have
\Equation\label{jeremiah}
\chibar(\kish(\obd(W)))\ge1\quad \text{or}\quad \wt T\le2.
\EndEquation
Now since $\obd(T)$ is incompressible in
$\oldPsi$ the definitions of $\delta(\oldPsi)$ and $\sigma(\oldPsi\cut{\obd(T)})$ (see \ref{t-defs})
imply that 
\Equation\label{isaiah}
\delta(\oldPsi)\ge\sigma(\oldPsi\cut{\obd(T)})=12\chibar(\kish(\oldPsi\cut{\obd(T)}))
=12
(\chibar(\kish(\obd(W)))+\chibar(\kish(\obd(\overline{N-W})))).
\EndEquation
%\bigg
If the first alternative of (\ref{jeremiah}) holds, i.e. if
$\chibar(\kish(\obd(W)))\ge1$, then it follows from (\ref{isaiah})
that $\delta(\oldPsi)\ge12$, which is the conclusion of the lemma. 

It
remains to consider the case in which the second alternative of
(\ref{jeremiah}) holds, i.e. $\wt T\le2$.
In this case, we will make a second application of Lemma \ref{modification}, this time
taking $K=\overline {N-W}$. According to Lemma \ref{even easier}, we
have $h(N)\le h(W)+h(N-W)-1$, so that $h(N-W)\ge h(N)-h(W)+1$. Since $h(W)\le h(N)-2$, it follows that $h(N-W)\ge 3$. Since $W$ is connected, Alternative (b) of Lemma
\ref{modification} is now seen to hold with $K=\overline{N-W}$. If Alternative (ii)  of Lemma
\ref{modification} also holds with this choice of $K$, then since $\lambda_\oldPsi=2$, Lemma
\ref{modification} gives a proper, non-empty, compact 
submanifold $K_1$ of $N$, whose frontier components are incompressible
tori in $N$, such that $t_1:=\Fr_N
K_1$ has weight less than
$ \wt T$. 
Since 
 $\lambda_\oldPsi=2$, Lemma
\ref{modification} also guarantees that 
$\max(\wt T-\wt T_1,\wt T_1)
%\wt_\oldPsi(\calt_1)
%,\wt_\oldPsi(\calt_1))
\ge2$.
But $\wt T_1<\wt T\le2$, and hence $\wt T-\wt T_1\ge2$, i.e. $\wt T_1\le \wt T-2\le0$.
This means that
%and since
%$\wt(\Fr_N W_1)$ is even by
%\redmissingref{replace by a $\lambda$ reference}, we must then have 
$T_1\cap
K_1=\emptyset$, and so any component of $\obd(\calt_1)=\calt_1$ is an
incompressible toric suborbifold of $\oldPsi$. This is a contradiction since $\oldPsi$
is strongly \simple (see \ref{oops}).
Hence Alternative  (ii) of Lemma
\ref{modification} must fail to hold with $K=\overline{N-W}$. As we
have already seen that $\obd(T)=\obd(\Fr_N(\overline{N-W}))$ is incompressible, this
means that 
\Equation\label{obadobah}
\chibar(\kish(\obd(\overline{N-W})))\ge1/2.
\EndEquation
It now follows from (\ref{isaiah}),  (\ref{obadiah}) and (\ref{obadobah}) that $\delta(\oldPsi)\ge12$.
\EndProof
%\oldOmega\wt^*p

\section{Homology of underlying
  manifolds}\label{irr-M section}

\Proposition\label{agol-plus}
If an orientable hyperbolic manifold $M$
has at least two cusps then $\vol M\ge\voct$ (see \ref{voct def}).
\EndProposition

\Proof
We may of course assume that $\vol M<\infty$. If $M$ has exactly two
cusps, the result follows from \cite[Theorem 3.6]{twocusps}. If $M$
has more than two cusps, it is a standard consequence of Thurston's
hyperbolic Dehn filling theorem \cite[Chapter E]{bp} that there is a hyperbolic manifold $M'$ having
exactly two cusps which can be obtained from $M$ by Dehn filling, and
that $\vol(M')<\vol(M)$. Since $\vol(M')\ge\voct$, we have
$\vol(M')>\voct$ in this case.
\EndProof

\Proposition\label{new get lost}
Let $\oldPsi$ be a compact, strongly \simple, boundary-irreducible,
orientable $3$-orbifold containing no embedded negative
turnovers. Set $N=|\oldPsi |$. Suppose that every
component of $\partial N$ is a sphere, that $N$ is
$+$-irreducible, and that $\smock_0(\oldPsi)\le3.44$. (See \ref{t-defs}, and note that $\smock_0(\oldPsi)$ is well defined since the strongly \simple\ orbifold $\oldPsi$ is in particular very good.)
Then the following conclusions hold.
\begin{enumerate}
\item If $\lambda_\oldPsi=2$ then $h(N)\le7$.
\item If $\lambda_\oldPsi=2$ and $h(N)\ge6$ then either $\theta(\oldPsi)>10$ or
  $\delta(\oldPsi)\ge6$.
\item If $\lambda_\oldPsi=2$ and $h(N)\ge 4$ then either $\theta(\oldPsi)>4$ or $\delta(\oldPsi)\ge 6$.
\item If $h(N)\ge 4$ then either $\theta(\oldPsi)> 4$ or $\delta(\oldPsi)\ge3$.
\end{enumerate}
\EndProposition

\Proof
The hypothesis and the definition of $\smock(\oldPsi)$ (see \ref{t-defs}) give
\Equation\label{for future reference}
\smock(\oldPsi)= \frac{\smock_0(\oldPsi)}{0.305}
%\bigg\rfloor
\le 
\frac{3.44}{0.305}
%\cdot3.44
%\bigg\rfloor=10.
<11.
\EndEquation
%since $\voct=3.66\ldots$.
%(\oldPsi)\bigg\rfloor.$$ 

Set $M=\plusN$, so that $M$ is a closed,
orientable $3$-manifold. The $+$-irreducibility of $N$
means that $M$ is irreducible. 
%Since 
We have $h(M)=h(N)$.

By (\ref{more kitsch}) and the
hypothesis of the present proposition, we have
$\volorb(\oldPsi)\le
 \smock_0(\oldPsi)\le3.44$. Applying Proposition \ref{hepcat} (and recalling that $\oldPsi$ is very good by \ref{oops}), we find
\Equation\label{slob cat}
\volG M\le \volorb(\oldPsi)\le 3.44.
\EndEquation

Because of the equivalence between the PL and smooth categories for $3$-manifolds, $M$ has a smooth structure compatible with its PL structure, and up to topological isotopy, the smooth $2$-submanifolds of $M$ are the same as its PL $2$-submanifolds. This will make it unnecessary to distinguish between smooth and PL tori in the following discussion.

Consider the case in which $M$ (regarded as a smooth manifold) is
hyperbolic. In this case, we have $\vol M=\volG M$ by
\cite[Theorem C.4.2]{bp}, and (\ref{slob cat}) gives
$$
\vol M\le \volorb(\oldPsi)\le 3.44. 
$$
According to  \cite[Theorem 1.7]{fourfree},
  any closed, orientable hyperbolic $3$-manifold $M$ of volume at most
  $3.44$ satisfies $h(M)\le7$. This establishes Conclusion (1) in this
  case.

According to \cite[Theorem 1.1]{rankfour},
  any closed, orientable hyperbolic $3$-manifold $M$ of volume at most
  $1.22$ satisfies $h(M)\le3$. Hence if $h(N)\ge4$, we have 
$$\vol
  M>1.22.$$
 But since $\volG M=\vol M$ by
\cite[Theorem C.4.2]{bp}, the definition of
 $\theta(\oldPsi)$ (see \ref{t-defs}) gives
$\theta(\oldPsi)=
%\frac{
%6}{\voct}
\vol(M)/0.305$.
%\bigg\rfloor.$$
Hence
$\theta(\oldPsi)>
%2\bigg\lfloor
%\frac{
%6}{\voct}\cdot
1.22/{0.305}
%\bigg\rfloor=
=4$.
%\redcomment{This doesn't seem quite right. $\voct$ is just over $3.66$, so $1.22/\voct$ is just under $1/3$. That would seem to give a value just under $2$ for the argument of the floor, giving $2$, not $4$, for the floor value. Maybe I need to fiddle with the defs. For example, I wonder whether using $3.66$ in place of $\voct$ in the defs. would take care of the issue, but I'm not sure.}
This establishes conclusions (2) and (4) in this case.

According to \cite[Theorem 1.2]{lastplusone}, 
  any closed, orientable hyperbolic $3$-manifold $M$ of volume at most
  $3.08$ satisfies $h(M)\le5$. 
Hence if $h(N)\ge6$, we have $\vol
  M>3.08$, and therefore
$\theta(\oldPsi)=\vol(M)/{0.305}>
%6}{\voct}\cdot
3.08/0.305>10$.
%(since $\voct=3.66\ldots$). 
%\redmissingref{That value should have been given earlier, and see the comment above.} 
This establishes conclusion (3) in this case, and completes the proof
of the proposition in the case where $M$ is hyperbolic.

Now consider the case in which  $M$ is
not hyperbolic. 
Since $M$ is irreducible, it follows from Perelman's geometrization theorem
\cite{bbmbp}, \cite{Cao-Zhu}, \cite{kleiner-lott}, \cite{Morgan-Tian}, that there is a (possibly empty) $2$-manifold
$\calt\subset M$, each component of which is an incompressible torus,
such that for each component $C$ of $M-\calt$, either $\hatC$ is a Seifert
fibered space, or $C$ admits
a
hyperbolic metric of finite volume; in the latter case, we will say
more briefly that $C$ is hyperbolic. 

We claim:
\Claim\label{found a doughnut}
If $h(N)\ge4$ then $M$ contains at least one incompressible torus.
\EndClaim

Note that \ref{found a doughnut} is obvious if $\calt\ne\emptyset$. If
$\calt=\emptyset$, then since $M$ is not itself hyperbolic, $M$ must
be a Seifert fibered space. But since $h(N)\ge4$, it follows from Lemma \ref{and four if by zazmobile}
that $M$
contains an incompressible torus. Thus \ref{found a doughnut} is established.

If
%According to (\ref{four more years}) we have
$h(N)\ge4$, then by \ref{found a doughnut} and 
Assertion (1) of Proposition \ref{crust
  chastened}, we have
$\delta(\oldPsi)\ge3\lambda_\oldPsi$. 
 Thus in particular we have $\delta(\oldPsi)\ge3$; and if, in addition to assuming
%According to (\ref{four more years}) we have
$h(N)\ge4$, we assume $\lambda_\oldPsi=2$, then 
%\ref{found a doughnut} and 
%Proposition \ref{crust
  %chastened} imply
$\delta(\oldPsi)\ge6$. 
This establishes Conclusions (2), (3) and (4) in this case.
%, we need only note that i

To prove Conclusion (1) in this case, let $\ch$ denote the set of all hyperbolic components of
$M-\calt$. (The set $\ch$ may be empty.) Since $M$ is not itself
hyperbolic, we have $\partial \hatC\ne\emptyset$ for any
$C\in\ch$. 

It follows from \cite[Theorem 1]{soma}, together with the fact (see \cite[Section 6.5]{thurstonnotes}) that the volume of a finite-volume hyperbolic $3$-manifold is equal to the relative Gromov volume of its compact core, that $\sum_{C\in\ch}\vol( C)=\volG M$. With (\ref{slob cat}), this gives
\Equation\label{last cat standing}
\sum_{C\in\ch}\vol( C)
%=\volG M
\le 3.44.
\EndEquation
Since $\partial\hatC\ne\emptyset$ for each $C\in\ch$, each 
hyperbolic manifold in the collection $\ch$ has at least one cusp. Hence by
\cite[Theorem 1.1]{cao-m}, we have
\Equation\label{what do you want? four fewer years?}
\vol( C)\ge2\calv=2.029\ldots \quad\text{for each }C\in\ch,
\EndEquation
where $\calv$
is the volume of the ideal regular
tetrahedron in hyperbolic 3-space. It follows from
(\ref{last cat standing}) and (\ref{what do you want? four fewer
  years?}) that $\card \ch\le1$.

Let us now distinguish two subcases, according as $\ch=\emptyset$ or
$\card \ch=1$. In each of these subcases, we will show that Conclusion
(1) of the proposition holds.

In the subcase $\ch=\emptyset$, the manifold $M=\plusN$ is by definition a
graph manifold. 
To prove Conclusion (1) in this subcase, note that if
$\lambda_\oldPsi=2$ and $h(M)\ge8$, we may apply Assertion (2) of
Proposition \ref{crust chastened} to deduce that 
$\delta(\oldPsi)\ge12$. By Corollary \ref{bloody hell},  we therefore have
$\smock(\oldPsi)\ge12$,
%
%.$$ 
%
%\le2\bigg\lfloor\frac{6}{\voct}\volorb(\oldPsi)\bigg\rfloor\le
%2\bigg\lfloor
%\frac{6}{\voct}\cdot3.44\bigg\rfloor=10,$$ 
a contradiction to (\ref{for future reference}). Thus (1) holds in this subcase.

There remains the subcase in which $\card \ch=1$. Let $C_0$ denote the
unique element of $\ch$. Then $ C_0$ is a (connected)
finite-volume, orientable hyperbolic $3$-manifold with at least one
cusp. According to (\ref{last cat standing}) we have
$\vol( C_0)
%=\volG C_0
\le 3.44$. On the other hand, if $C_0$
had at least two cusps, we would have $\vol C_0\ge\voct=3.66\ldots$ by
Proposition \ref{agol-plus}. Hence $ C_0$ has exactly one cusp,
i.e. $\partial \hatC_0$ is a single torus.

According to \cite[Theorem 6.2]{hodad}, if $Y$ is a complete,
finite-volume, orientable hyperbolic 3-manifold having exactly one
cusp, and if
$h(Y)\ge6$, then $\vol(Y)>5.06$. Since $ C_0$ has exactly one cusp, and
since $\vol( C_0)\le3.44$ by (\ref{last cat standing}),
we have $h(C_0)\le5$. Thus if we assume that $\lambda_\oldPsi=2$ and $h(N)\ge7$, we have
$h(C_0)\le h(N)-2$. Since $C_0$ is acylindrical by Lemma \ref{torus
  goes to cylinder}, 
we may apply Lemma
\ref{get the basin}, with $Z=C_0$, to deduce that 
$\delta(\oldPsi)\ge12$. By Corollary \ref{bloody hell},  we therefore have
$\smock(\oldPsi)\ge12$,  a contradiction to (\ref{for future reference}).
%. By Corollary \ref{bloody hell},  we therefore have
%$\smock(\oldPsi)\ge12$. In view of the definition of $\smock(\oldPsi)$, this
%means that
%$$2\bigg\lfloor\frac{6}{\voct}\volorb(\oldPsi)\bigg\rfloor\ge12,$$
%i.e. that
 %$\volorb(\oldPsi)\ge3.6$, a
%contradiction to the hypothesis 
%$\volorb(\oldPsi)\le 3.44$. 
Hence in this subcase $\lambda_\oldPsi=2$ implies $h(M)\le6$, and in particular
Conclusion (1)
holds in this 
subcase as well.

\abstractcomment{\tiny Check carefully for editing errors. In particular check that
  it's clear that the hypotheses of Lemma \ref{get the basin} hold for
that last app. I'm confusing about the role of the boundary spheres
that constitute the difference between $N$ and $M$. How am I making
the transition? I guess the answer is contained in the proofs of
Proposition \ref{crust chastened} and Lemma \ref{get the basin}, which are about
$\plusX$ (or is it $\plusN$? I don't know what I was talking about
here exactly), but depend on earlier lemmas about $X$ (or is it $N$?).}

\EndProof
%\cite very

The following result was mentioned in the introduction as Proposition E:

\Proposition\label{lost corollary}
Let $\Mh$ be a closed,
orientable, hyperbolic $3$-orbifold such that $\oldOmega:=(\Mh)\pl$ contains no embedded negative
turnovers. 
% such that each component of $\fraks_\oldOmega$ is a simple closed curve.  
Suppose that $M:=|\oldOmega|$ is irreducible.
Then:
\begin{itemize}
\item If $\lambda_\oldOmega=2$ and $\vol(\Mh)\le3.44$, then $h(M)\le7$. 
\item If $\lambda_\oldOmega=2$ and  $\vol(\Mh)<1.83$, then $h(M)\le5$.
\item If  $\lambda_\oldOmega=2$ and $\vol(\Mh)\le1.22$, then $h(M)\le3$.
\item If $\vol(\Mh)<0.915$, then $h(M)\le3$. 
\end{itemize}
\EndProposition

\Proof Note that the hypothesis of each of the assertions implies that $\vol\Mh\le3.44$. According to Corollary \ref{smockollary} and the hypothesis, we then
have $\smock_0(\oldOmega)= \vol\Mh\le3.44$. Furthermore, the hyperbolicity of $\Mh$ implies that $\oldOmega$ is strongly \simple\ (see \ref{oops}). Thus the hypotheses of
Proposition \ref{new get lost} hold with $\oldOmega$ and $M$ playing the
respective roles of $\oldPsi$ and $N$. According to Assertion (1)
of Proposition \ref{new get lost}, if $\lambda_\oldOmega=2$, we have $h(M)\le7$, which
is the first conclusion of the present proposition. 

To prove the remaining
conclusions, note that 
%let us set 
%$m=\max(\theta(\oldOmega),\delta(\oldOmega)$.
Corollary \ref{bloody hep} gives 
$\smock(\oldOmega)\ge\theta(\oldOmega)$, and that
Corollary \ref{bloody hell} gives
$\smock(\oldOmega)\ge\delta(\oldOmega)$. Hence $\smock(\oldOmega)\ge\max(\theta(\oldOmega),\delta(\oldOmega))$.
In view of the definition of $\smock(\oldOmega)$ (see \ref{t-defs}), this
means that
$\smock_0(\oldOmega)/{0.305} \ge\max(\theta(\oldOmega),\delta(\oldOmega))$. Again using that
$\vol\Mh=\smock_0(\oldOmega\pl)$ by  Corollary \ref{smockollary}, we deduce:
\Equation\label{that's not an all}
%2\lfloor
\vol(\Mh) \ge0.305\max(\theta(\oldOmega),\delta(\oldOmega)).
\EndEquation

If  $\lambda_\oldPsi=2$ and  $\vol(\Mh)<1.83$, it follows from (\ref{that's not an all})
that $\max(\theta(\oldOmega),\delta(\oldOmega))<6$, which by
Assertion (2) of Proposition \ref{new get lost} implies that
$h(M)\le5$. This is the second assertion of the present proposition.

If  $\lambda_\oldPsi=2$ and  $\vol(\Mh)\le1.22$, it follows from (\ref{that's not an all})
that $\max(\theta(\oldOmega),\delta(\oldOmega))\le4$, which by
Assertion (3) of Proposition \ref{new get lost} implies that
$h(M)\le3$. This is the third assertion of the present proposition.

If $\vol(\Mh)<0.915$, it follows from (\ref{that's not an all})
that $\max(\theta(\oldOmega),\delta(\oldOmega))<3$, which by
Assertion (4) of Proposition \ref{new get lost} implies that
$h(M)\le3$. This is the fourth assertion of the present proposition.
\EndProof

\section{Hyperbolic volume and orbifold homology: a special case}\label{A and C}

%\Number\label{psidef}
%We define a real number $\psi>0$ by $\psi=\prod _{j=1}^\infty(1-2^{-j})$. Note that the infinite product converges absolutely because the series $\sum _{j=1}^\infty2^{-j}$ converges.
%\EndNumber

\Number\label{tsar}
Recall that if $v$ is a vertex in a finite graph (see \ref{dual graph}), the {\it valence} of $v$ is defined to be the number of  oriented edges whose terminal vertex is $v$. (The underlying edges of these oriented edges are not necessarily distinct, as the graph may contain loops.)

Let $\oldOmega$ be a closed, orientable $3$-orbifold. Then $M:=|\oldOmega|$ is a closed, orientable $3$-manifold, and $\Sigma:=\Sigma_\oldOmega$ is a union of strata of dimension at most $1$. Every component $C$ of $\Sigma$ is either (a) a single $1$-dimensional stratum which is a simple closed curve, or (b) the underlying set of a graph in which the vertices and edges are strata of $\Sigma$. We will denote by $\Phi=\Phi_\oldOmega$ the union of all components of $\Sigma$ satisfying (a), and
 by $\Psi=\Psi_\oldOmega$ the union of all components of $\Sigma$
 satisfying (b). Then $\Psi$ is itself the underlying set of a graph in which the vertices and edges are strata of $\Sigma$; these will be referred to as vertices and edges of $\Psi$.

Each vertex $v$ of $\Psi$ has valence $3$. If $\vec{e}_1$, $\vec{e}_2$
and $\vec{e}_3$ denote the  oriented edges of $\Psi $ whose terminal
vertex is $v$, and if $p_i$ denotes the order of the underlying edge
$e_i$ of $\vec{e}_i$, we have $1/p_1+1/p_2+1/p_3>1$. Furthermore, the
orientations of the $\vec{e}_i$ define 
generators $c_i$ of the cyclic groups $G_{e_i}$ 
%of $\pi_1(\oldOmega)$ 
for $i=1,2,3$ (see \ref{orbifolds introduced}),
such that
%may be chosen within their conjugacy classes, and a 
% may be chosen,  in such a way that
 $c_1c_2c_3=1$. The group $G_v$ is generated by $c_1$, $c_2$ and $c_3$ and is a $(p_1,p_2,p_3)$-triangle group. (Cf. \cite[Corollary 3.11]{blp}.)

For each $1$-dimensional stratum $\sigma$ of $\Sigma$, 
regarding $G_\sigma$ as a subgroup of $\pi_1(\oldOmega)$, well-defined
up to conjugacy, as in \ref{orbifolds introduced},
%r the natural homomorphism $\pi_1(\oldOmega)
we define an element $x_\sigma$ of $H_1(\oldOmega;\FF_2 )$ to be the image of a generator of  $G_\sigma$ under the natural homomorphism $\pi_1(\oldOmega)\to H_1(\oldOmega;\FF_2 )$; note that $x_\sigma$ is independent of 
%the choice of $G_e$ within its conjugacy class, and of 
the choice of a generator of $G_\sigma$. If $v$ is a vertex of $\Psi$,  if $\vec{e}_1$, $\vec{e}_2$ and $\vec{e}_3$ denote the oriented edges of $\Psi $ whose terminal vertex is $v$, and if $e_i$ denotes the  underlying edge of $\vec{e}_i$, we have $x_{e_1}+x_{e_2}+x_{e_3}=0$.

We denote by $\Sigma'=\Sigma'_\oldOmega$ the union of the closures of all strata $\sigma$ of $\Sigma$ such that $x_\sigma\ne0$. (A necessary condition for $x_\sigma$ to be non-zero is that $G_\sigma$ have even order.) Thus $\Sigma'\cap\Phi$ is a union of components of $\Phi$, while $\Sigma'\cap\Psi$ is (the underlying set of) a subgraph of (the graph defining) $\Psi$. 

If $v$ is a vertex of $\Psi$ lying in $\Sigma'\cap\Psi$,  
if $\vec{e}_1$, $\vec{e}_2$ and $\vec{e}_3$ denote the oriented edges
of $\Psi $ whose terminal vertex is $v$, and if $e_i$ denotes the
underlying edge of $\vec{e}_i$, then since
$x_{e_1}+x_{e_2}+x_{e_3}=0$, there cannot be exactly one index
$i\in\{1,2,3\}$ such that $x_{e_i}\ne0$; that is, $v$ cannot have
valence $1$ in $\Sigma'\cap\Psi$. Since $\Sigma'$ is by definition a
union of closures of $1$-dimensional strata of $\Sigma$, the valence
of $v$ in $\Sigma'\cap\Psi$ cannot be $0$. Hence every vertex of
$\Psi$ lying in $\Sigma'\cap\Psi$ has valence $2$ or $3$ in
$\Sigma'\cap\Psi$. It follows that  every component $C'$ of $\Sigma'$
is either (a$'$) a  simple closed curve which is either a component of
$\Phi$ or a subcomplex of $\Psi$, or (b$'$) the underlying set of a
graph in which the vertices are vertices of $\Sigma$, and each edge is
a union of vertices and edges of $\Sigma$; the vertices and edges of
this graph will be referred to as vertices and edges of $C'$. We will
denote by $\Phi'=\Phi'_\oldOmega$ the union of all components of
$\Sigma'$ satisfying (a$'$), and
 by $\Psi'=\Psi'_\oldOmega$ the union of all components of $\Sigma'$ satisfying (b$'$).

We define a {\it superstratum} of $\Sigma'$ to be a subset of $\Sigma'$ which is either a component of $\Phi'$ or an edge of $\Psi'$. It follows from this definition that every superstratum $s$ of $\Sigma'$ is either a component of $\Sigma'\cap\Phi$, or a topological open arc or simple closed curve which is a union of edges and vertices of $\Psi$ that are contained in $\Sigma'$. In the latter case, every vertex of $\Psi$ contained in $s$ has valence $2$ in $\Sigma'\cap\Psi$.
Note that in any event the superstrata of $\fraks'_\oldOmega$ are the components of the set of $1$-manifold points of $\fraks'_\oldOmega$.

We will associate with every superstratum $s$ of $\Sigma'$ a non-zero element $y_s$ of $H_1(\oldOmega;\FF_2 )$. In the case where $s$ is a component of $\Sigma'\cap\Phi$, the superstratum $s$ is in particular a stratum of $\Sigma$, and we set $y_s=x_s$; since $s\subset\Sigma'$, we have $x_s\ne0$. In the case where $s$ is a topological open arc or simple closed curve which is a union of edges and vertices of $\Psi$ that are contained in $\Sigma'$, we wish to define $y_s=x_e$, where $e$ is an arbitrary edge of $\Psi$ contained in $s$. To show that $y_s$ is well-defined we must show that for any two edges $e_1,e_2$ of $\Psi$ contained in $s$, we have $x_{e_1}=x_{e_2}$. Since $s$ is connected, it suffices to prove this under the additional assumption that $e_1$ and $e_2$ share a vertex $v$. Since $v$ lies in the superstratum $s$, it has valence $2$ in $\Sigma'\cap\Psi$. Hence we may label the oriented edges with terminal vertex  $v$ as
$\vec{e}_1$, $\vec{e}_2$ and $\vec{e}_3$, in such a way that  $e_i$ is the  underlying edge of $\vec{e}_i$ for $i=1,2$, and the underlying edge $e_3$ of $\vec e_3$ does not lie in $\Sigma'$. This means that $x_{e_3}=0$, and since $x_{e_1}+x_{e_2}+x_{e_3}=0$ it follows that $x_{e_1}=x_{e_2}$, as required. Thus $y_s$ is well-defined. Since $s\subset\Sigma'$ we have $x_e\ne0$ for any edge $e$ of $\Psi$ contained in $s$, and hence $y_s\ne0$.
\EndNumber

\Number\label{kernel}
Let $\oldOmega$ be a closed, orientable $3$-orbifold. We will denote by $K=K_\oldOmega$ the kernel of the natural homomorphism $H_1(\oldOmega;\FF_2 )\to H_1(|\oldOmega|;\FF_2 )$.

The kernel of the natural homomorphism $\pi_1(\oldOmega)\to\pi_1(|\oldOmega|)$ is
% The normal subgroup $K_\oldOmea$ of $\pi_1(\oldOmega)$ is
 the normal closure of the union of the subgroups $G_\sigma$, where $\sigma$ ranges over the $1$-dimensional strata of $\Sigma_\oldOmega$. Hence $K$ is spanned by the elements $x_\sigma$, where $\sigma$ again ranges over the $1$-dimensional strata of $\Sigma_\oldOmega$.
If a $1$-dimensional stratum $\sigma$ of $\Sigma$ is not contained in $\Sigma':=\Sigma'_\oldOmega$, then by definition we have $x_\sigma=0$. Hence  $K$ is spanned by the elements $x_\sigma$, where $\sigma$ now ranges over the $1$-dimensional strata of $\Sigma_\oldOmega$ that are contained in $\Sigma'$. Furthermore, for each
$1$-dimensional stratum $\sigma$ of $\Sigma_\oldOmega$ that is contained in $\Sigma'$, we have $x_\sigma=y_s$, where $s$ is the superstratum of $\Sigma'$ containing $s$. Hence $K$ is spanned by the elements $y_s$, where $s$ ranges over all superstrata of $\Sigma'$.
\EndNumber

\Notation\label{just on accounta}
Let $\oldOmega$ be a closed, orientable $3$-orbifold. For each non-zero linear function $\phi:H_1(\oldOmega;\FF_2 )\to\FF_2 $, we denote by $L_\phi$ the  PL subset of $|\oldOmega|$ defined to be the union of the closures of all superstrata $s$ of $\Sigma'_\oldOmega$ such that $\phi(y_s)=1$.
\EndNotation

\Lemma\label{sick of that shoot}
Let $\oldOmega$ be a closed, orientable $3$-orbifold, and let $\phi:H_1(\oldOmega;\FF_2 )\to\FF_2 $ be a   non-zero linear function. Then $L_\phi$ is a link in the $3$-manifold $|\oldOmega|$, i.e. a (PL) closed $1$-manifold. Furthermore, if $p:\toldOmega\to\oldOmega$ denotes the two-sheeted (orbifold) covering of $\oldOmega$ 
defined by the codimension-$1$ subspace
%index-$2$ subgroup of $\pi_1(\oldOmega)$ which is the
%preimage of
$\ker\phi$ of
%under the natural homomorphism from to
$H_1(\oldOmega;\FF_2 )$, then $|p|:|\toldOmega|\to|\oldOmega|$ is a branched covering of manifolds whose branch locus is $L_\phi$. Finally, 
we have
%the number of components of
 $\compnum(L_\phi)
%$ is bounded above by the
\le 1+h(|\toldOmega|)$.
%\dim H_1(|\oldOmega|;\FF_2 )$.
\EndLemma

\Proof
To prove the first assertion, it suffices to show that each of the sets $L_\phi\cap\Phi'_\oldOmega$ and $L_\phi\cap\Psi'_\oldOmega$ is a $1$-manifold. To show that $L_\phi\cap\Phi'_\oldOmega$ is a $1$-manifold, we need only observe that $L_\phi$ is by definition a union of superstrata of $\fraks'_\oldOmega$, and that the superstrata contained in $\Phi'_\oldOmega$ are the components of $\Phi'_\oldOmega$. To show that $L_\phi\cap\Psi'_\oldOmega$ is a $1$-manifold, it suffices to show that for every vertex $v$ of $\Psi'_\oldOmega$ (in the sense defined in \ref{tsar}) there are exactly two oriented edges of $L_\phi$ having $v$ as their terminal vertex. According to the definition of $\Psi'_\oldOmega$ there are exactly three oriented edges $\vec{e'_1}$, $\vec{e'_2}$ and $\vec{e'_3}$ of $\Psi'_\oldOmega$ having $v$ as their terminal vertex. For $i=1,2,3$ let $\vec{e_i}$ denote the unique oriented edge of $\fraks_\oldOmega$ which is contained in $\vec{e_i'}$, inherits its orientation from $\vec{e_i'}$, and has $v$ as its terminal vertex. Let $e_i$ and $e_i'$ denote the underlying edges of $\vec{e_i}$ and $\vec{e_i'}$ respectively. Then by \ref{tsar} we have $x_{e_1}+x_{e_2}+x_{e_3}=0$. But \ref{tsar} also gives that $y_{e_i'}=x_{e_i}$ for $i=1,2,3$. Thus we have $y_{e_1'}+y_{e_2'}+y_{e_3'}=0$ and therefore $\phi(y_{e_1'})+\phi(y_{e_2'})+\phi(y_{e_3'})=0$. Hence there are an even number of indices $i\in\{1,2,3\}$ for which $\phi(y_{e_i'})=1$, i.e. for which $e_i'\subset L_\phi$. But this number is strictly positive since $L_\phi$ is by definition a union of closures of  superstrata of $\Sigma'_\oldOmega$. Hence the number in question must be $2$, and the first assertion of the lemma is proved.

As a first step in proving the second assertion, first note that since
$p:\toldOmega\to\oldOmega$ is a  covering map of orbifolds,
$|p|\,\big|\,|p|^{-1}(|\oldOmega|-\fraks_\oldOmega): |p|^{-1}(|\oldOmega|-\fraks_\oldOmega)\to |\oldOmega|-\fraks_\oldOmega$ is a covering map of manifolds. Since $\dim\fraks_\oldOmega\le1$, it follows that $|p|:|\toldOmega|\to|\oldOmega|$ is a branched covering of manifolds whose branch locus is contained in $\fraks_\oldOmega$.

Now we claim:

\Claim\label{cohen pone}
Let $Q$ be a point of a $1$-dimensional stratum $\sigma$ of $\fraks_\oldOmega$. Then $Q$ belongs to the branch locus of $|p|$ if and only if $\sigma\subset L_\phi$.
\EndClaim

To prove \ref{cohen pone}, we choose a point $\tQ\in p^{-1}(Q)$, so
that $G_{\tQ}$ is canonically identified with a subgroup of
$G_Q$. Regarding $G_Q$ as a subgroup of $\OO(n)$, and using (orbifold)
chart maps that map $0$ to $Q$ and $\tQ$, we may identify some neighborhoods
$\frakV$ and $\tfrakV:=p^{-1}(\frakV)$ of $Q$ and $\tQ$ via orbifold
homeomorphisms with $U/G_Q$ and $U/G_{\tQ}$, where $U$ is a ball in
$\RR^n$, in such a way that $p:\tfrakV\to\frakV$ is the natural
projection from $U/G_{\tQ}$ to $U/G_Q$. It follows that for $p$ to be
a local homeomorphism at $\tQ$---or, equivalently, for $|p|$ to be a
local homeomorphism at $\tQ$---it is necessary and sufficient that the
$G_{\tQ}$ be the full group $G_Q$. On the other hand, $G_{\tQ}$ and
$G_Q$ are also canonically identified with subgroups of
$\pi_1(\toldOmega,\tQ)$ and $\pi_1(\oldOmega,Q)$; and if $\Delta$ denotes the image of $p_\sharp:\pi_1(\toldOmega,\tQ)\to\pi_1(\oldOmega,Q)$, we have $G_{\tQ}=G_Q\cap\Delta$. Hence $|p|$ is a local homeomorphism at $\tQ$ if and only if $G_Q\subset\Delta$.

Now let $\gamma$ denote a generator of the cyclic group $G_Q$, and
let $j$ denote the canonical surjection from $\pi_1(\oldOmega,Q)$ to $H_1(\oldOmega;\FF_2 )$. The definition of $x_\sigma$ (see \ref{tsar}) says that $j(\gamma)=x_\sigma$, and the definition of the covering $\toldOmega$ says that $\Delta=j^{-1}(\ker\phi)$. Thus the condition $G_Q\subset\Delta$ is equivalent to $\gamma\in\Delta$, and therefore to $x_\sigma\in\ker\phi$. This shows that $|p|$ is a local homeomorphism at $\tQ$ if and only if $\phi(x_\sigma)=0$. Hence $Q$ lies in the branch locus of $|p|$ if and only if $\phi(x_\sigma)=1$.

A necessary condition for the equality $\phi(x_\sigma)=1$ to hold is that $x_\sigma$ be non-zero, or equivalently (see \ref{tsar}) that $\sigma$ be contained in $\Sigma'_\oldOmega$. If $\sigma$ is contained in $\Sigma'_\oldOmega$, then by \ref{tsar} we have $x_\sigma=y_s$, where $s$ is the stratum of $\Sigma'_\oldOmega$ containing $\sigma$; in view of the definition of $L_\phi$ it then follows that $\phi(x_\sigma)=1$ if and only if $s\subset L_\phi$, i.e. if and only if $\sigma\subset L_\phi$. This completes the proof of  \ref{cohen pone}.

Since $|p|$ is a branched covering, its  branch locus is a closed subset of $|\oldOmega|$ whose components are all of codimension at most $2$; thus the branch locus has no isolated points. It therefore follows from \ref{cohen pone} that this branch locus is precisely $L_\phi$, and the second assertion of the lemma is proved.

%To finish the proof of the second assertion, just look at where branching occurs. }

To prove the third assertion, let us set $r=\compnum(L_\phi)$, and let
$t:\toldOmega\to\toldOmega$ denote the non-trivial deck transformation
of the two-sheeted orbifold covering
$p:\toldOmega\to\oldOmega$. Since $L_\phi$ is the  branch locus of the
branched covering $|p|$, the involution $|t|$ of the manifold $|\toldOmega|$ has fixed point set $|p|^{-1}(L_\phi)$, and $\compnum(\Fix|t|)=r$.

According to the Smith inequality (see \cite[p. 126, Theorem 4.1]{bredon}), we
have
\Equation\label{a shitte hier, a shitte her}
\sum_{d\ge0}\dim H_d(\Fix(|t|);\FF_2 )\le \sum_{d\ge0}\dim
H_d(|\toldOmega|;\FF_2 ).
\EndEquation
(In general this inequality holds if $|t|$ is replaced by an order-$p$
homeomorphism for a given prime $p$, and homology with coefficients in
$\FF_p$ is used; in this application we have $p=2$.)

Since $\Fix(|t|)$ is a disjoint union of $r$ simple closed curves,
the left hand side of (\ref{a shitte hier, a shitte her}) is $2r$. By
Poincar\'e duality, the right hand side of (\ref{a shitte hier, a
  shitte her}) is $2+2h(|\toldOmega|)$. Hence
$r\le1+h(|\toldOmega|)$. 
%t\Mh
%
%we have $\Fix(|t|)$is an application of the Smith inequalities that can be extracted from the proof of  Proposition \ref{boogie-woogie bugle boy}, and the proof of that proposition should be revised so as to quote the present result, among other things.
\EndProof
%\toldOmega Q j v x

\Lemma\label{boogie-woogie bugle boy}
Let $\oldOmega$ be a 
closed, orientable 3-orbifold such that $\Psi'_{\oldOmega}=\emptyset$.
Then $\oldOmega$ is covered with degree at most $2$ by some
orbifold $\toldOmega$ such that
$$\dim H_1(\oldOmega;\FF_2 )\le1+ h(|{\toldOmega}|)+h(|{\oldOmega}|).$$
%$$\ge-1.$$
\EndLemma

\Proof
Set $K=K_\oldOmega$
Since $K$ is the kernel of the surjection $H_1(\oldOmega;\FF_2 )\to
H_1(|\oldOmega|;\FF_2 )$, we have
%(\ref{freestyle}) gives
 $\dim H_1(\oldOmega;\FF_2 )=\dim H_1(|\oldOmega|;\FF_2 )+\dim K$, i.e.
\Equation\label{shells}
\dim H_1(\oldOmega;\FF_2 )=h(|\oldOmega|;\FF_2 )+\dim K.
\EndEquation

According to \ref{kernel}, $K$ is spanned by the elements $y_s$, where $s$ ranges over all superstrata of $\Sigma':=\Sigma'_\oldOmega$. Since $\Psi'_\oldOmega=\emptyset$,
we have $\Phi'_\oldOmega=\Sigma'$, so that the superstrata of $\Sigma'$ are simply its components. Hence $K$ has a basis $B$ consisting of vectors of the form $y_C$ where $C$ is a component of $\Sigma'$. Define a linear function $\phi^K:K\to \FF_2 $ by setting $\phi(y)=1$ for every $y\in B$, and extend $\phi^K$ arbitrarily to a linear function $\phi:H_1(\oldOmega;\FF_2 )\to\FF_2 $.

Since  the superstrata of $\Sigma'$ are simply its components, the definition of $L_\phi$ given in \ref{just on accounta} says that $L_\phi$ is the union of all components $C$ of $\Sigma'_\oldOmega$ such that $\phi(y_C)=1$.
 Since for every $y\in B$
we have
 $\phi(y)=1$, and $y=y_C$ for some component $C$ of $\Sigma'$, it follows that $\compnum(L_\phi )\ge\card B=\dim K$. With (\ref{shells}) this gives 
$\dim H_1(\oldOmega;\FF_2 )\le h(|\oldOmega|;\FF_2 )+\compnum(L_\phi)$. But according to Lemma \ref{sick of that shoot} we have
$\compnum(L_\phi)
\le 1+h(|\toldOmega|)$, where 
$\toldOmega$ denotes the two-sheeted (orbifold) covering of $\oldOmega$ 
defined by the codimension-$1$ subspace
%index-$2$ subgroup of $\pi_1(\oldOmega)$ which is the
%preimage of
$\ker\phi$ of
% under the natural homomorphism from to
$H_1(\oldOmega;\FF_2 )$. Hence
$\dim H_1(\oldOmega;\FF_2 )\le 1+h(|\oldOmega|;\FF_2 )+h(|\toldOmega|)$, as asserted by the proposition.
\EndProof
%C_TS\tau $t\oldOmega\Mh link\Mh hyperbolic

\Proposition\label{vector thing}
Let $V$ be a  finite dimensional vector space over $\FF_2 $. Let $I$ be a non-empty finite index set, and let $(z_i)_{i\in I}$ be an indexed family of non-zero elements of $V$. Then there is a codimension-$1$ subspace $W$ of $V$ such that 
$$\card\{i\in I:z_i\notin W\}>\frac12\card I.$$
\EndProposition

\Proof
Since $I\ne\emptyset$, and $z_i\ne0$ for each $i\in I$, we have $V\ne0$.
Set $n=\dim V>0$ and $m=\card I>0$. For $0\le k\le n-1$, we will show by induction on $k$ that there is a $k$-dimensional subspace $W_k$ of $V$ such that
\Equation\label{induction thing}
\card\{i\in I:z_i\notin W_k\}>(1-2^{k-n})m.
\EndEquation
%The definition of $\psi$ implies that the right hand side of (\ref{induction thing}) is greater than $\psi m$ for $k=0,\ldots,n-1$. Hence
For $k=n-1$, this provides an $(n-1)$-dimensional subspace $W_{n-1}$
of $V$ such that $$\card\{i\in I:z_i\notin
W_{n-1}\}>(1-2^{-1})m=\frac m2,$$ which gives
% (\ref{induction thing}) holds implies
 the conclusion of the proposition.

For $k=0$, if we set $W_0=\{0\}$, the hypothesis that $z_i\ne0$ for every $i\in I$ implies that 
$$\card\{i\in I:z_i\notin W_k\}=m>(1-2^{-n})m,$$
%\frac12\cdot\frac{(2^{n}-1)}{2^{n-1}-1}m,$$
which establishes the base case of the induction.
%j
%
%the left hand side of (\ref{induction thing}) is equal to $M$, so that the inequality (\ref{induction thing}) is an equality
Now suppose that $k$ is given with $0\le k<n-1$, and that $W_k$ is a $k$-dimensional subspace of $V$ for which (\ref{induction thing}) holds. Set 
$m_k=
\card\{i\in I:z_i\notin W_k\}$, so that $m_k >(1-2^{k-n})m$
%\ge 
%(2^{k-1}(2^{n-k}-1))m/(2^{n-1}-1)$
by (\ref{induction thing}).
We may regard $W_k$ as a subgroup of the additive group of $V$, and we denote by $\calx$ the set of all proper cosets of $W_k$ in $V$. For each $X\in\calx$, set $a_X=\card\{i\in I:z_i\in X\}$. We have $\{i\in I:z_i\notin W_k\}=\coprod_{X\in\calx}\{i\in I:z_i\in X\}$, so that $m_k=\sum_{X\in\calx}a_X$. Since $\card\calx=2^{n-k}-1$, it follows that for some $X_0\in\calx$ we have $a_{X_0}\le m_k/(2^{n-k}-1)$. 

Since $W_k$ is a $k$-dimensional subspace of the $\FF_2$-vector space $V$, and $X_0$ is a proper additive coset of $W_k$ in $V$, the set $W_{k+1}:=W_k\discup X_0$ is a $(k+1)$-dimensional subspace of $V$. Set $m_{k+1}=
\card\{i\in I:z_i\notin W_{k+1}\}$. We have $\{i\in I:z_i\notin W_{k+1}\}=\{i\in I:z_i\notin W_{k}\}-\{i\in I:z_i\in X_0\}$, and hence
$$
\begin{aligned}
m_{k+1}&=m_k-a_{X_0}\ge m_{k}-m_k/(2^{n-k}-1)=\bigg(1-\frac1{2^{n-k}-1}\bigg) m_k\\
&>\bigg(1-\frac1{2^{n-k}-1}\bigg)(1-2^{k-n})m
%\redcomment{\text{fix from
 %   here}}\qquad
%\frac{2^{k-1}(2^{n-k}-1)}{2^{n-1}-1}\cdot m=
%\frac{2^{k}(2^{n-k-1}-1)}{2^{n-1}-1}\cdot m,
=(1-2^{k-n+1})m,
\end{aligned}
$$
which shows that (\ref{induction thing}) holds with $k+1$ in place of $k$, and thus completes the induction.
%%( Hence $\sum_{X\in\calx}a_{X_0}\ge\redcomment{Finish this.}
\EndProof
%vz $W$

\Lemma\label{when lambda is one}
Let $\oldOmega$ be a 
closed, orientable 3-orbifold such that $\Psi'_\oldOmega\ne\emptyset$. Then there exists a $((\ZZ/2\ZZ)\times(\ZZ/2\ZZ))$-covering  $\toldOmega$ of $\oldOmega$ such that
%by some
%orbifold $\tMh$ such that
$\dim H_1(\oldOmega;\FF_2 )\le 3+h(|\oldOmega|)+4h(|\toldOmega|)$.
\EndLemma

%\redcomment{I've used $\oldOmega$ instead of $\Mh$ because I don't know what $\Mh$ or the hyperbolicity hypothesis was doing in Proposition \ref{boogie-woogie bugle boy}. I need to check the proof and apps of that one to see whether $\Mh$ should be removed both there and here. }

\Proof
Set $K=K_\oldOmega$. According to \ref{kernel}, $K$ is spanned by the elements $y_s$, where $s$ ranges over all superstrata of $\Sigma':=\Sigma'_\oldOmega$. 
Recall that by definition a superstratum of $\Sigma'$  is either a
component of $\Phi':=\Phi'_\oldOmega$ or an edge of
$\Psi':=\Psi'_\oldOmega$. 
Let $K_0$ denote the subspace of $K$ spanned by the elements $y_e$ where $e$ ranges over the edges of  $\Psi'$.

By \ref{tsar} we have $y_e\ne0$ for every edge $e$ of $\Psi'$, and since $\Psi'\ne\emptyset$ by hypothesis, we have $K_0\ne0$. 
Let us fix  a basis $B_0$ of $K_0$ 
 consisting of vectors of the form $y_e$ where $e$ is an edge of $\Psi'$. The quotient space $K/K_0$ is spanned by images of elements of $K$
having the form $y_s$ for some superstratum $s$ of $\Sigma'$. Hence
there is a set $B_1\subset K$, each of whose elements has the form
$y_s$ for some superstratum $s$ of $\Sigma'$, such that the quotient
map $K\to K/K_0$ maps $B_1$ bijectively onto a basis of $K/K_0$. In
particular no element of $B_1$ lies in $K_0$, and therefore no element
of $B_1$ has the form $y_e$ for an edge $e$ of $\Psi'$. Hence every
element of $B_1$ has the form $y_s$ for some component $s$ of
$\Phi'$. Since $B_0$ is a basis of $K_0$, and 
since the quotient map $K\to K/K_0$ maps $B_1$ bijectively onto a basis of $K/K_0$, we have a basis $B:=B_0\discup B_1$ of $K$.

We denote by $\phi_1^K:K\to\FF_2 $ the unique  $\FF_2$-linear function which maps each element of $B_0$ to $1$ and maps each element of $B_1$ to $0$. We fix, arbitrarily, a linear extension $\phi_1: H_1(\oldOmega;\FF_2 )\to\FF_2 $ of $\phi_1^K$. Then $L_{\phi_1}\cap\Psi'$ is by definition a subcomplex of $\Psi'$, and according to Lemma \ref{sick of that shoot}, 
$L_{\phi_1}\cap\Psi'$ is also a closed $1$-manifold. Hence if we
let $\calv$ denote the set of all vertices of $L_{\phi_1}\cap\Psi'$, then $n:=\card\calv$ is also equal to the number of edges of $L_{\phi_1}\cap\Psi'$. But by the definition of $L_{\phi_1}$, an edge $e$ of $\Psi'$ lies in $L_{\phi_1}\cap\Psi'$ if and only if $\phi_1(y_e)=1$, and the definition of $\phi_1$ guarantees that $\phi_1(y_e)=1$ whenever $y_e\in B_0$. Hence we have $n\ge\card B_0$, i.e.
\Equation\label{sunday bloody funday}
n\ge\dim K_0.
\EndEquation
In particular, since $K_0\ne0$, we have $n>0$, i.e. $\calv\ne\emptyset$.

The non-triviality of $K_0$ also implies that the linear function $\phi_1$ is not identically $0$, since it takes the value $1$ at every element of the basis $B_0$ of $K_0$.

Consider an arbitrary element $v$ of $\calv$. The definition of $\Psi'$ implies that $v$ has valence $3$ in 
$\Psi'$; but since $L_{\phi_1}\cap\Psi'$ is also a closed $1$-manifold, 
$v$ has valence $2$ in 
 $L_{\phi_1}\cap\Psi'$. Hence there is a unique edge $e_v$ of $\Psi'$ which is not contained in 
$L_{\phi_1}\cap\Psi'$
but has $v$ as an endpoint; and $e_v$ is not a loop.
We set $z_v=y_{e_v}$. According to \ref{tsar} we have $ z_v \ne 0$, and the definition of $K_0$ gives $z_v\in K_0$. 
Thus $(z_v)_{v\in \calv}$ is an indexed family of non-zero elements of the $\FF_2$-vector space $K_0$.

We may now apply Proposition
\ref{vector thing}, letting  $K_0$  play the role of $V$ in that proposition, and letting $\calv$ play the role of the non-empty index set $I$. This gives a codimension-$1$
subspace $W$ of $K_0$ such that the set $\calv^*:=\{v\in \calv:z_v\notin W\}$ satisfies $\card\calv^*>\card(\calv)/2$, i.e.
\Equation\label{dinner at eight}
\card\calv^*>n/2.
\EndEquation

Let $\phi_2^{K_0}:K_0\to\FF_2 $ denote the unique $\FF_2$-linear
function on $K_0$ whose kernel is $W$. Since the quotient map $K\to
K/K_0$ maps $B_1$ bijectively onto a basis of $K/K_0$, we may define a
linear extension $\phi_2^K:K\to\FF_2 $ by stipulating that $\phi_2^K$
maps every element of $B_1$ to $1$. We fix, arbitrarily, a linear
extension $\phi_2: H_1(\oldOmega;\FF_2 )\to\FF_2 $ of $\phi_2^K$.

Since $n>0$ it follows from \ref{dinner at eight} that $\calv^*\ne\emptyset$. For any $v\in\calv^*$ we have $z_v\notin W$ by definition, so that $\phi_2(z_v)=\phi_2^{K_0}(z_v)=1$. On the other hand, by construction we have $z_v=y_{e_v}\not\subset L_{\phi_1}$, so that $\phi_1(z_v)=0$. This shows both that $\phi_2\ne\phi_1$ and that $\phi_2$ is not identically $0$. Since we have also shown that  $\phi_1$ is not identically $0$, the linear map $\psi:H_1(\oldOmega;\FF_2 )\to(\FF_2 )^2$ defined by $\psi(y)=(\phi_1(y),\phi_2(y))$ is surjective. Hence $\ker\psi$ defines a $((\ZZ/2\ZZ)\times(\ZZ/2\ZZ))$-covering $p:\toldOmega\to\oldOmega$.

Since the linear function $\phi_1$ is not identically $0$, its kernel is a subspace of codimension  $1$ in $H_1(\oldOmega;\FF_2 )$. Hence
%the preimage of
$\ker\phi_1$
% under the natural homomorphism from $\pi_1(\oldOmega)$ to $H_1(\oldOmega;\FF_2 )$ has index  $2$ in $\pi_1(\oldOmega)$, and hence
defines a degree-$2$ orbifold covering 
$p_1:\toldOmega_1\to\oldOmega$. Note that since $\ker\psi$ is a codimension-$1$ subspace of $\ker\phi_1$, which is in turn a codimension-$1$ subspace of $H_1(\oldOmega;\FF_2 )$, we may write $p=p_2\circ p_1$ for some degree-$2$ covering map $p_2:\toldOmega\to\toldOmega_1$. Note also that, in view of the surjectivity of $\psi$, the linear function $\phi_2|(\ker\phi_1)$ is not identically $0$, so  that the linear function  
 $\tphi_2:=\phi_2\circ (p_1)_*:H_1(\toldOmega_1;\FF_2 )\to\FF_2 $ is not identically $0$; and that $p_2:\toldOmega\to\toldOmega_1$ is the degree-$2$ covering defined by the codimension-$1$ subspace $\ker\tphi_2$ of $H_1(\toldOmega_1;\FF_2 )$.

Let $t:\toldOmega_1\to\toldOmega_1$ denote the non-trivial deck transformation of the degree-$2$ orbifold covering $p_1:\toldOmega_1\to\oldOmega$. Since $t$ is an involution of the orbifold $\toldOmega_1$, the involution $|t|$ of $|\toldOmega_1|$ leaves $\fraks_{\toldOmega_1}$ invariant. Let $\sigma$ be any $1$-dimensional stratum of $\fraks_{\toldOmega_1}$. Then $t(\sigma)$ is also a stratum of $\fraks_{\toldOmega_1}$, and the outer automorphism $t_\sharp$ of $\pi_1(\toldOmega_1)$ carries the conjugacy class of the subgroup $G_\sigma$ of $\pi_1(\toldOmega_1)$ (see \ref{orbifolds introduced}) onto the conjugacy class of $G_{t(\sigma)}$. Hence the involution $t_*$ of $H_1(\toldOmega_1;\FF_2 )$ maps $x_\sigma$ to $x_{t(\sigma)}$. Since this is true for every $1$-dimensional stratum $\sigma$ of $\fraks_\oldOmega$, the set $\fraks'_{\toldOmega_1}$ is $t$-invariant. Since the superstrata of $\fraks'_{\toldOmega_1}$ are the components of the set of $1$-manifold points of $\fraks'_{\toldOmega_1}$, each superstratum of $\fraks'_{\toldOmega_1}$ is mapped by $t$ onto a superstratum of $\fraks'_{\toldOmega_1}$. Furthermore, if $s$ is any superstratum of $\fraks'_{\toldOmega_1}$, and if we choose any stratum $\sigma$ of $\fraks_{\toldOmega_1}$ with $\sigma\subset s$, then we have $t(s)\subset t(\sigma)$ and hence $t_*(y_s)=t_*(x_\sigma)=x_{t(\sigma)}=y_{t(s)}$. Thus we have shown:

\Claim\label{ice people}
For every superstratum $s$ of $\fraks'_{\toldOmega_1}$, the set $t(s)$ is a superstratum of $\fraks'_{\toldOmega_1}$, and we have $t_*(y_s)=y_{t(s)}$.
\EndClaim
%\toldOmega

Now suppose that $s$ is a superstratum of $\fraks'_{\toldOmega_1}$ contained in $L_{\tphi_2}$. By definition this means that $\tphi_2(y_s)=1$. But since $p_1\circ t=p_1$, we have $(p_1)_*\circ t_*=(p_1)_*$, and hence $\tphi_2\circ t_*=\phi_2\circ (p_1)_*\circ t_*=\phi_2\circ (p_1)_*=\tphi_2$. Using \ref{ice people}, we now find that $\tphi_2(y_{t(s)})=\tphi_2(t_*(y_s))=\tphi_2(y_s)=1$, so that $t(s)\subset L_{\tphi_2}$. This proves:

\Claim\label{you should be so lucky}
The deck transformation $t$ leaves $L_{\tphi_2}$ invariant.
\EndClaim

Now by Lemma \ref{sick of that shoot}, $L_{\phi_1}$ is the branch locus of the branched covering $|p_1|$. Hence:
\Equation\label{pete's a pie}
\Fix t=|p_1|^{-1}(L_{\phi_1}).
\EndEquation

We claim:
\Claim\label{or he was, anyway}
No component of $L_{\tphi_2}$ is contained in
$|p_1|^{-1}(L_{\phi_1})$.
\EndClaim

If \ref{or he was, anyway} is false, then in particular some
superstratum $s$ of  $\fraks'_{\toldOmega_1}$ 
is contained both in $L_{\tphi_2}$ and in
$|p_1|^{-1}(L_{\phi_1})$. Since $s\subset L_{\tphi_2}$ we have $\tphi_2(y_s)=1$.  On the other hand, if $\sigma$ is a stratum of $\fraks_\oldSigma$ contained in $s$, then $y_s=x_\sigma$ is the image of a generator $\gamma$ of $G_\sigma$ under the natural homomorphism $\pi_1(\toldOmega_1)\to H_1(\oldOmega_1;\FF_2 )$. Up to conjugacy in $\pi_1(\toldOmega_1)$, we may represent $\gamma$ by an oriented simple closed curve bounding a disk that meets $\sigma$ transversally in a single point, and is disjoint from the other strata of $\fraks_\oldSigma$. Since $\sigma\subset s\subset |p_1|^{-1}(L_{\phi_1})$, and
since
$L_{\phi_1}$ is the branch locus of the two-fold branched covering map $|p_1|$
by Lemma \ref{sick of that shoot}, the element $(p_1)_\sharp(\gamma)$, which is defined up to conjugacy in $\pi_1(\oldOmega)$, is a square in $\pi_1(\oldOmega)$. Hence $(p_1)_*(y_s)=0$. In particular $\tphi_2(y_s)=
\phi_2( (p_1)_*(y_s)=0$. This contradiction establishes \ref{or he was, anyway}.

Now suppose that $C$ is a component of $L_{\tphi_2}$ such that $C\cap |p_1|^{-1}(L_{\phi_1})\ne\emptyset$. Since $L_{\tphi_2}$ is a closed $1$-manifold by \ref{sick of that shoot}, $C$ is a simple closed curve. Since $|p_1|^{-1}(L_{\phi_1})=\Fix t$ by (\ref{pete's a pie}), and since $L_{\tphi_2}$ is $t$-invariant by \ref{you should be so lucky}, it follows that $C$ is $t$-invariant. But by \ref{or he was, anyway}, $C$ is not contained in $|p_1|^{-1}(L_{\phi_1})=\Fix t$. Hence $t|C$ is a non-trivial involution of the simple closed curve $C$, and therefore has at most two fixed points, i.e. $C$ meets $|p_1|^{-1}(L_{\phi_1})=\Fix t$ in at most two points. This proves:

\Claim\label{far too many}
Every component of $L_{\tphi_2}$ meets $|p_1|^{-1}(L_{\phi_1})$ in at most two points.
\EndClaim

%$1$-manifold

%Since $\tphi_2$ is the composition of $(p_1)_*$ with a linear function on $H_1(\oldOmega;\FF_2 )$, the set $L_{\tphi_2}$ is invariant under the deck transformation of the covering $\toldOmega_1$ of $\oldOmega$. On the other hand, no component of $L_{\tphi_2}$ is contained in the fixed point set of the deck transformation, which is $|p_1|^{-1}(L_{\phi_1})$. This is because $L_{\phi_1}$ is the branch locus of $p_1$ by Lemma \ref{sick of that shoot}, so that for any edge $e$ of $|p_1|^{-1}(L_{\phi_1})$, the map $(p_1)_*$ takes $x_e$ to $0$. (In fact it takes a generator of $G_e$ to a square.) So in particular $x_e$ is in the kernel of $\tphi_2$, which means $e$ is not in $L_{\tphi_2}$. This implies that any component of $L_{\tphi_2}$ contains at most two points of $p^{-1}(L_{\phi_1})$, because if it contains at least one it is invariant under the deck transformation, but the deck transformation is not the identity on it, so the restriction has at most two fixed points.

Set $\tcalv^*=|p_1|^{-1}(\calv^*)$. Since $\calv^*\subset\calv\subset 
L_{\phi_1}\cap\Psi'\subset L_{\phi_1}$, and by Lemma \ref{sick of that shoot} $L_{\phi_1}$ is the branch locus of the two-fold branched covering $|p_1|$, the map $p_1\big|\tcalv^*:\tcalv^*\to\calv^*$ is bijective, and hence $\card\tcalv^*=\card\calv^*$. We claim that
\Equation\label{imf}
\tcalv^*\subset L_{\tphi_2}.
\EndEquation

To prove (\ref{imf}), consider any point $\tv\in\tcalv^*$, and set $v=p_1(\tv)\in\calv^*$. Recall that the edge $e_v$ of $\Psi'$  is not contained in 
$L_{\phi_1}\cap\Psi'$
but has $v$ as an endpoint.
%, and that $e_v$ is not a loop. 
Hence
$v\in\overline{e_v}$ and $e_v\cap L_{\phi_1}=\emptyset$. 
%\redcomment{Why should that
  %follow? Can't the other endpoint be in $L_{\phi_1}$? Is this used?} 
It follows that there is a $1$-dimensional stratum $\sigma \subset
e_v$ of $\fraks_\oldOmega$  such that
$v\in\overline{\sigma}$ and $\sigma\cap L_{\phi_1}=\emptyset$. 
% $\overline{\sigma }\cap L_{\phi_1}=\{v\}$. 
Since $L_{\phi_1}$ is the branch locus of the  branched covering
$|p_1|$, it follows that there exist an open arc
$\alpha\subset|\toldOmega_1|$ with $\tv\in\overline\alpha$, 
%which is mapped homeomorphically onto
%$\overline{\sigma }$ by $|p_1|$, \redcomment{This doesn't quite follow
%  from what I said. Do I need the closed arc to map homeomorphically?
  %If so, I can probably deduce it from the fact, which I have \%-ed
  %out, that $e_v$ is not a loop}
%with $\tv\in\alpha$, 
and an open neighborhood $\tU$ of $\alpha$ in $|\toldOmega_1|$,
such that
%which is
$p_1$  maps $\tU$ homeomorphically  onto a neighborhood $U$ of $\sigma $ in $|\oldOmega|$. This implies that $\alpha$ is contained in a $1$-dimensional stratum $\tsigma$ of $\fraks_{\toldOmega_1}$.

Choose a point $\tQ\in\alpha$, and set $q=p(\tQ)\in \sigma $. Since $|p_1|\big|\tU:\tU\to U$  is a homeomorphism, $p_1$ is in particular a local homeomorphism at $\tQ$, and hence the homomorphism $(p_1)_\sharp
:\pi_1(\toldOmega_1)\to\pi_1(\oldOmega)
$, which is defined modulo inner automorphisms, carries the conjugacy class of $G_{\tQ}$ onto the conjugacy class of $G_Q$. But we have $G_{\tsigma}=G_{\tQ}$ and $G_{\sigma} =G_Q$, so that $(p_1)_\sharp$ carries the conjugacy class of $G_{\tsigma}$ onto the conjugacy class of $G_{\sigma} $. Hence $(p_1)_*(x_{\tsigma})=x_{\sigma} $. Since $\sigma \subset e_v\subset\Psi'$, we have $x_{\sigma} \ne0$, and therefore $x_{\tsigma}\ne0$; thus $\tsigma\subset\Psi'_{\toldOmega_1}$. 
 Let ${\ts}$ denote the superstratum of $\Psi'_{\toldOmega_1}$ containing $\tsigma$. Then $(p_1)_*(y_{{\ts}})=(p_1)_*(x_{\tsigma})=x_{\sigma}=y_{e_v}=z_v$. Thus we have $\tphi_2(y_{{\ts}})=\phi_2((p_1)_*(y_{{\ts}}))=\phi_2(z_v)$. But since $v\in \calv^*$, the definition of $\calv^*$ gives 
%$v\in\calv\subset\Psi^*$ $ and 
$z_v\notin W$. Furthermore, since $e_v\subset\Psi'$, we have $z_v=y_{e_v}\in K_0$. The definition of $\phi_2$ then gives $\phi_2(z_v)=\phi_2^{K_0}(z_v)$, and since $W=\ker\phi_2^{K_0}$ we have $\phi_2(z_v)=1$. This shows that  $\tphi_2(y_{{\ts}})=1$, so that $\ts\subset L_{\tphi_2}$. Since $\tv\in\overline{\alpha}\subset\overline{\tsigma}\subset\overline{\ts}$, it now follows that $\tv\in L_{\tphi_2}$, and (\ref{imf}) is proved.
%sigma {\ts_v}$ denote the superstratum of %$\Psi'_{\toldOmega_1}$ containing $\tsigma$
%$\inter\alpha$ is contained in a $1$-dimensional stratum $\tsigma$ of
%$\fraks_{\toldOmega_1}$. sigma_v \inter s_v

%
%{ and for each $\tv\in\tcalv^*$, if we set $v=p_1(\tv)\in\calv^*$, there's an edge $\ts_v$ of $p^{-1}Psi'$ lying above $s_v$. Since $v\in\calv^*:=\{v\in calv:z_v\notin W\}$, we have $(p_1)_*(x_{\ts_v})\notin W$. (I'm slightly confused about the $x$ vs. $y$ issue here, and the corresponding issue of strata vs. superstrata.) So $x_{\ts_v}\notin\ker\tphi_2$, i.e. $\ts_v\subset L_{\tphi_2}$. 

According to (\ref{imf}), each point of $\tcalv^*$ lies on some component of $ L_{\tphi_2}$. Since $\tcalv^*=|p_1|^{-1}(\calv^*)\subset 
|p_1|^{-1}(\calv)\subset
|p_1|^{-1}(L_{\phi_1})$, it follows from \ref{far too many} that every component of $L_{\tphi_2}$ contains at most two points of $\tcalv^*$. Hence the number of components of $L_{\tphi_2}$ that have non-empty intersection with $\tcalv^*$ is at least $\card(\tcalv^*)/2=\card(\calv^*)/2$. 
But by \ref{dinner at eight} and (\ref{sunday bloody funday}), we have $\card(\calv^*)> n/2\ge(\dim K_0)/2$. Hence:

\Claim\label{the boys are marching}
The number of components of $L_{\tphi_2}$ that have non-empty intersection with $\tcalv^*$ is strictly greater than $(\dim K_0)/4$.
\EndClaim

Now we claim:

\Claim\label{pharmacy}
The number of components of $L_{\tphi_2}$ that are contained in $|p_1|^{-1}(\Phi')$  is at least $\dim(K/K_0)$. 
\EndClaim

To prove \ref{pharmacy}, set $m=\dim(K/K_0)$. Recall that $B_1\subset K$ is mapped onto a basis of $K/K_0$ by the quotient map, so that $\card B_1=m$. Furthermore, we have seen that 
every element of $B_1$ has the form $y_s$ for some component $s$ of $\Phi'$. Hence if we denote by $\scrt$
the set of all components $s$ of $\Phi'$ such that $y_s\in B_1$, we have $\card\scrt\ge m$. 

Consider an arbitrary element $s$ of $\scrt$. Since $y_s\in B_1$, the construction of $\phi_1$ implies that $\phi_1(y_s)=0$. Hence $s$ is disjoint from the branch locus
$L_{\phi_1}$  of the branched covering $|p_1|$. 
% $|p_1|^{-1}(s)$ has two components, and that
Each component of $|p_1|^{-1}(s)$ is therefore a simple closed curve,
and has an open neighborhood $\tT$  which is mapped by $|p_1|$ onto an
open neighborhood $T$ of $s$; furthermore, $|p_1|$ restricts to a
covering map from $\tT$ to $T$. It follows that if  $\tsigma$ is any   $1$-dimensional stratum  of $\fraks_{\toldOmega_1}$ contained in $|p_1|^{-1}(s)$, then $|p_1|$ maps $\tsigma$  onto a $1$-dimensional stratum $\sigma\subset s$ of $\fraks_\oldOmega$. Moreover, if $\tQ$ is any point of $\tsigma$, and if we set $Q=p_1(\tQ)\in\sigma$, then since $|p_1|$ is a local homeomorphism at $\tQ$, the homomorphism $(p_1)_\sharp:\pi_1(\toldOmega_1)\to\pi_1(\oldOmega)$, which is well-defined modulo inner automorphisms, maps $G_{\tQ}$ onto $G_Q$, i.e. it maps $G_{\tsigma}$ onto $G_\sigma$. Hence $(p_1)_*:H_1(\toldOmega_1;\FF_2 )\to H_1(\oldOmega;\FF_2 )$ maps $x_{\tsigma}$ to $x_\sigma=y_s$. Since $\sigma\subset s\subset\Phi'$, we have $x_\sigma\ne0$, and hence $x_{\tsigma}\ne0$; that is, $\tsigma\subset\fraks'_{\toldOmega_1}$. We also have $\tphi_2(x_{\tsigma})=\phi_2((p_1)_*(x_{\tsigma}))=\phi_2(y_s)$. Since $y_s\in B_1$, the construction of $\phi_2$ gives $\phi_2(y_s)=\phi_2^K(y_s)=1$, so that $\tphi_2(x_{\tsigma})=1$. If we denote by $\ts$ the superstratum of $\fraks'_{\toldOmega_1}$ containing $\tsigma$, we have $y_{\ts}=x_{\tsigma}$, so that $\tphi_2(y_{\ts})=1$, i.e. $\ts\subset L_{\tphi_2}$. Since this holds for every $1$-dimensional stratum $\tsigma$ of $\fraks_{\toldOmega_1}$ contained in $|p_1|^{-1}(s)$, we deduce that $|p_1|^{-1}(s)\subset L_{\tphi_2}$.

Thus for every $s\in\scrt$, each component of the set $|p_1|^{-1}(s)$ 
%has two components,
%each of which
 is a simple closed curve contained in
$L_{\tphi_2}$. We have $|p_1|^{-1}(s)\ne\emptyset$ since the branched
covering map $|p_1|$ is surjective. 
%\redcomment{I'm not sure I've made it clear why they are
  %simple closed curves. Also, I'm not sure I see the relevance of the
  %phrase ``Since the elements of $\scrt$ are in particular components
  %$s$ of $\Phi'$'' below. Most importantly, the claim about 
  %two components seems to be wrong. This means the coefficient of $2$
  %in \ref{pharmacy} is wrong, but it doesn't affect the final result
  %of the lemma.} 
Since the components of $L_{\tphi_2}$ are simple
closed curves by Lemma \ref{sick of that shoot}, each component of
$|p_1|^{-1}(s)$ is a component
 of $L_{\tphi_2}$. Since
% the elements of $\scrt$ are in particular components $s$ of $\Phi'$, and since 
 $\card\scrt\ge m$, it follows that
there are at least $m$
components of $L_{\tphi_2}$  contained in $|p_1|^{-1}(\Phi')$. Thus  \ref{pharmacy} is proved.

Now since $\tcalv^*=|p_1|^{-1}(\calv^*)\subset|p_1|^{-1}(\calv)
\subset |p_1|^{-1}(\Psi')$, we have $\tcalv^*\cap|p_1|^{-1}(\Phi')=\emptyset$. Hence
no component of $L_{\tphi_2}$ can be contained in $|p_1|^{-1}(\Phi')$ and have non-empty intersection with
$\tcalv^*$. In view of \ref{the boys are marching} and \ref{pharmacy}, it follows that
\Equation\label{freestyle}
\compnum(L_{\tphi_2})>\frac{\dim K_0}4+\dim(K/K_0)\ge\frac{\dim K}4.
\EndEquation
%\dim B_0\frak

Since $K$ is the kernel of the surjection $H_1(\oldOmega;\FF_2 )\to H_1(|\oldOmega|;\FF_2 )$, (\ref{freestyle}) gives
\Equation\label{scales and fins}
\dim H_1(\oldOmega;\FF_2 )=\dim H_1(|\oldOmega|;\FF_2 )+\dim K<h(|\oldOmega|)+4\compnum(L_{\tphi_2}). 
\EndEquation

On the other hand, since
$p_2:\toldOmega\to\toldOmega_1$ is the two-sheeted (orbifold) covering of $\toldOmega_1$ defined by $\ker\tphi_2$,
%defined by the codimension-$1$
%index-$2$ subgroup of $\pi_1(\oldOmega)$ which is the
%preimage of
%subspace $\ker\phi$ of
% under the natural homomorphism from to
%$H_1(\oldOmega;\FF_2 )$. Then $\toldOmega$ is a degree-$4$ covering of $\oldOmega$.
% with the terminology elsewhere in the section, where I explicitly mention the index-$2$ subgroup of $\pi_1$}
we may apply Lemma \ref{sick of that shoot}, with $\toldOmega_1$, $\tphi_2$ and $p_2$ playing the respective roles of $\oldOmega$, $\phi$ and $p$,  to deduce that $$\compnum(L_{\tphi_2})\le1+h(|\toldOmega|).$$
Combining this with (\ref{scales and fins}), we obtain
$$\dim H_1(\oldOmega;\FF_2 )<4+h (|\oldOmega|)+4h(|\toldOmega|).$$
 This implies  the conclusion of the lemma. \EndProof
%\tphi_1\toldOmega_2$\alphaU v^*\ts\Mh hyperbolic V W T

The following result was stated in the introduction as Proposition A. 

\Proposition\label{my little sony}
Let $\oldOmega$ be a 
closed, orientable 3-orbifold.
Then either
\begin{enumerate}[(i)] 
\item $\oldOmega$ is covered with degree at most $2$ by some
orbifold $\toldOmega$ such that
$$\dim H_1(\oldOmega;\FF_2 )\le1+ h(|{\toldOmega}|)+h(|{\oldOmega}|),$$
or
\item there exists a $((\ZZ/2\ZZ)\times(\ZZ/2\ZZ))$-covering  $\toldOmega$ of $\oldOmega$ such that
%by some
%orbifold $\tMh$ such that
$$\dim H_1(\oldOmega;\FF_2 )\le 3+h(|\oldOmega|)+4h(|\toldOmega|).$$
Furthermore, if $\fraks_\oldOmega$ is a link, then (i) holds.
\end{enumerate}
\EndProposition

\Proof
If  $\Psi'_{\oldOmega}=\emptyset$, it follows from Lemma \ref{boogie-woogie bugle boy} that (i) holds. If  $\Psi'_{\oldOmega}\ne\emptyset$, it follows from Lemma  \ref{when lambda is one} that (ii) holds. This proves the first assertion.

Now suppose that $\fraks_\oldOmega$ is a link. Since
 $\Psi'_\oldOmega$ is contained in $\fraks_\oldOmega$ and is the
 underlying set of a trivalent graph, we have $\Psi'_\oldOmega=\emptyset$. It therefore follows from Lemma \ref{boogie-woogie bugle boy} that (i) holds. This proves the second assertion.
\EndProof

The following result was stated in the introduction as Proposition F.

\Proposition\label{orbifirst}Let $\Mh$ be a 
closed, orientable, hyperbolic 3-orbifold such that 
%of volume strictly less than
%$
%\voct/4=
%0.915$.
$\fraks_{\Mh}$ is a link, and such that $\pi_1(\Mh)$
contains no hyperbolic triangle group. Suppose that $|\Mh|$ is irreducible, and that $|\tMh|$ is irreducible for every two-sheeted (orbifold) covering $\tMh$ of $\Mh$. If $\vol \Mh\le1.72$ then
$\dim H_1(\Mh;\FF_2 )\le
13$. Furthermore, if
$\vol \Mh\le1.22$ then
$\dim H_1(\
Mh;\FF_2 )\le
11$, and if $\vol \Mh\le0.61$ then
$\dim H_1(\Mh;\FF_2 )\le
7$. 
\EndProposition

\Proof
Set $\oldOmega=(\Mh)\pl$. Then $\fraks_\oldOmega$ is a link in $|\oldOmega|$. 
%Since by \ref{tsar}, $\Psi'_\oldOmega$ is contained in $\fraks_\oldOmega$ and is the underlying set of a trivalent graph, we have $\Psi'_\oldOmega=\emptyset$. This is the ``useful'' thing, I think(??)
It therefore follows from the final sentence of Proposition \ref{my little sony}
%\ref{boogie-woogie bugle boy} 
%\redcomment{It's a lemma now, and this will have to be rewritten} 
that $\Mh$ is covered with degree at most $2$ by an
orbifold $\tMh$ such that 
\Equation\label{number me}
\dim H_1(\Mh;\FF_2 )\le1+ h(|{\tMh}|)+h(|{\Mh}|).
\EndEquation
%Since $\vol\Mh<0.915$, we have
%$\vol\tMh=2\vol\Mh<1.83$. 
Set $\oldOmega=(\Mh)\pl$ and  $\toldOmega=(\tMh)\pl$.
Since $\Mh$ has a link as its singular set,
so does $\tMh$; thus $\lambda_\oldOmega=\lambda_{\toldOmega}=2$. Since  $\pi_1(\Mh)$
contains no hyperbolic triangle group, it follows from Corollary \ref{injective hamentash}
that neither $\oldOmega$ nor $\toldOmega$ contains
any embedded negative turnover. The hypothesis implies that
$|\oldOmega|$ and $|\toldOmega|$ are both irreducible. If
$\vol\Mh\le1.72$, then 
%$\vol\Mh$ and 
$\vol\tMh=2\vol\Mh
%$ are both bounded above by $
\le3.44$; it then follows from Proposition  \ref{lost corollary} that
$h(|\Mh|)\le5$ and $h(|\tMh|)
%$ are both bounded above by $
\le7$. We therefore have $\dim H_1(\Mh;\FF_2 )\le13$ by (\ref{number me}). 
If $\vol\Mh\le1.22$, then  Proposition \ref{lost corollary} gives 
$h(|\Mh|)\le3$; and since 
$\vol\tMh\le2.44<3.44$,  Proposition \ref{lost corollary} also gives 
$h(|\tMh|)\le7$, so that $\dim H_1(\Mh;\FF_2 )\le11$ by (\ref{number me}). Likewise, if $\vol\Mh\le0.61$, then $\vol\Mh$ and $\vol\tMh=2\vol\Mh$ are both bounded above by $1.22$; it then follows from Proposition \ref{lost corollary} that $h(|\Mh|)$ and $h(|\tMh|)$ are both bounded above by $3$, so that $\dim H_1(\Mh;\FF_2 )\le7$ by (\ref{number me}). 
\EndProof
%C_TSCorollary oldOmega

The following result was stated in the introduction as Proposition G.

\Proposition\label{orbinext}Let $\Mh$ be a 
closed, orientable, hyperbolic 3-orbifold such that 
$\pi_1(\Mh)$
contains no hyperbolic triangle group. Suppose that $|\Mh|$ is irreducible, and that $|\tMh|$ is irreducible for every two-sheeted (orbifold) covering $\tMh$ of $\Mh$ and for every $((\ZZ/2\ZZ)\times (\ZZ/2\ZZ))$-covering $\tMh$ of $\Mh$. If
$\vol \Mh<0.22875$, then $\dim H_1(\Mh;\FF_2 )\le
18$.
\EndProposition

\Proof
Set $\oldOmega=(\Mh)\pl$. 
%In the case where  $\Psi_\oldOmega'=\emptyset$, Proposition \ref{boogie-woogie bugle boy} \redcomment{It's a lemma now, and this will have to be rewritten}  asserts that 
According to Proposition \ref{my little sony}, either $\oldOmega$ has a covering $\toldOmega$ of degree at most $2$ such that
$\dim H_1(\oldOmega;\FF_2 )\le1+ h(|{\toldOmega}|)+h(|{\oldOmega}|)$, or
%. In the case where $\Psi'_\oldOmega\ne\emptyset$, it follows from Proposition \ref{when lambda is one} \redcomment{It's a lemma now, and this will have to be rewritten} 
%that
 there is a $((\ZZ/2\ZZ)\times(\ZZ/2\ZZ))$ cover $\toldOmega$ of $\oldOmega$ such that
%orbifold $\tMh$ such that
$\dim H_1(\oldOmega;\FF_2 )\le 3+h(|\oldOmega|)+4h(|\toldOmega_2|)$.
Thus in either case, $\Mh$ has a regular covering $\tMh$, with covering group isomorphic to $(\FF_2 )^d$ for some $d\le2$, such that 
\Equation\label{forget me knot}
\dim H_1(\Mh;\FF_2 )\le 3+h(|\Mh|)+4h(|\tMh|).
\EndEquation

Since  $\pi_1(\Mh)$
contains no hyperbolic triangle group, it follows from Corollary \ref{injective hamentash}
that neither $\oldOmega$ nor $\toldOmega$ contains
any embedded negative turnover. The hypothesis implies that $|\oldOmega|$ and $|\toldOmega|$ are both irreducible. Since $\vol\Mh<0.22875$, and since the covering $\toldOmega$ of $\oldOmega$ have degree at most $4$, the volumes of $\Mh$ and $\tMh$ are both strictly bounded above by $4\times0.22875=0.915$. It now follows from Proposition  \ref{lost corollary} that $h(|\Mh|)$ and $h(|\tMh|)$ are both bounded above by $3$. We therefore have $\dim H_1(\Mh;\FF_2 )\le18$ by (\ref{number me}). \EndProof
%\Mh

% $\vol \Mh\le1.72$ then
%$\dim H_1(\Mh;\FF_2 )\le
%15$. 

%This should be followed by a new result that includes the $\lambda=1$ case. I think it should be stated as a separate prop. and also mentioned in the intro and cross-referenced, as this one is. I think it says that if $\vol \Mh<0.22875$, and $\Mh$ and all its four-fold covers are irreducible (but the singular set is not assumed to be a link) then $\dim H_1(\Mh;\FF_2 )\le
%18$. 

\chapter{Essential intersections and higher characteristic $2$-orbifolds}\label{higher chapter}

As was mentioned in the Introduction, the material of this chapter includes a 
partial analogue for $3$-orbifolds  of  the machinery developed for
    $3$-manifolds in \cite{bcsz} and elsewhere. This machinery is
    presented in Section \ref{higher section}, and depends on notions
    involving $2$-orbifolds that are developed in Section \ref{vegematic section}.

\section{Isotopies and essential intersections in $2$-orbifolds}\label{vegematic section}

In subsections \ref{vere's dot notation}---\ref{inchworm
  lemma}, we study isotopy classes of suborbifolds of a $2$-orbifold,
and develop the notion of the ``essential intersection'' of such
classes. In subsections \ref{before i guess}---\ref{more
  associativity}, we  study isotopy classes of
self-homeomorphisms of a $2$-orbifold, and the interaction between
such classes and essential intersections. We conclude the section with
a result, Proposition \ref{steinbeck}, which is of interest in its own
right and will be needed in Section \ref{clash section}.

\NotationRemarks
\label{vere's dot notation}
If $\otheroldLambda$ is an orientable $2$-orbifold, the manifold $|\otheroldLambda|-\fraks_\otheroldLambda$ will be denoted by $\dot\otheroldLambda$.

Note that if $\otheroldLambda$ has finite type then $\fraks_\otheroldLambda$ is finite, and hence $\dot\otheroldLambda$ has finite type.

If $\xi$ is an embedding of a $2$-orbifold $\oldUpsilon$ in an orientable $2$-orbifold $\otheroldLambda$, we have $\xi^{-1}(\fraks_\otheroldLambda)=\fraks_\oldUpsilon$, and $\xi(x)$ and $x$ have the same order for every $x\in\fraks_\oldUpsilon$. Hence $|\xi|\big|\dot\oldUpsilon:\dot\oldUpsilon\to\dot\otheroldLambda$ is an embedding of $2$-manifolds, which will be denoted $\dot\xi$.
\EndNotationRemarks

The following lemma, like several other results in this section, involves the general version of the definition of tautness that was given as Definition \ref{praxis}.

\Lemma\label{geventlach}
Let $\otheroldLambda$ be an orientable $2$-orbifold without boundary, and let
$\oldGamma$ be a taut  
finite-type $2$-suborbifold of $\otheroldLambda$. Then
$\dot\oldGamma$ is a taut,
finite-type submanifold of
$\dot\otheroldLambda$. Furthermore, is $\oldGamma$ is negative then $\dot\oldGamma$ is negative.
\EndLemma

\Proof
Since $\oldGamma$ has finite type, $\dot\oldGamma$ has finite type by \ref{vere's dot notation}. 
To show that $\dot\oldGamma$ is taut in $\otheroldLambda$, we reason from Definition \ref{praxis}. First note that $\dot\oldGamma$ is a closed subset of $\dot\otheroldLambda$ since $\oldGamma$ is a closed subset of $\otheroldLambda$. Next note that we have $\partial\dot\oldGamma=\partial\oldGamma$; since $\oldGamma$ is taut in $\otheroldLambda$, the $1$-manifold $\partial\dot\oldGamma=\partial\oldGamma$ is $\pi_1$-injective in $\otheroldLambda$ and hence in $\dot\otheroldLambda\subset\otheroldLambda$. Finally, suppose that some component  of $\partial\dot\oldGamma$ is contained in a $2$-manifold $\frakX\subset\otheroldLambda$ homeomorphic to $\SSS^1\times[0,\infty)$. We may write $\frakX=\dot\frakK$ for some component $\frakK$ of 
$\overline{\otheroldLambda-\oldGamma}$. Since $\dot\frakK$ is
homeomorphic to $\SSS^1\times[0,\infty)$, either (a) $|\frakK|$ is a disk
and $\card\fraks_\frakK=1$, or (b) $\frakK$ is itself a manifold homeomorphic to $\SSS^1\times[0,\infty)$. But (a) contradicts the $\pi_1$-injectivity of $\partial\oldGamma$ in $\otheroldLambda$ given by Condition (ii) of Definition \ref{praxis}, while (b) contradicts Condition (iii) of Definition \ref{praxis}. 

Now suppose that $\oldGamma$ is negative, and let $\oldGamma_0$ denote an arbitrary component of $\oldGamma_0$. Let $x_1,\ldots,x_p$ denote the points of $\oldGamma_0$ (where $p\ge0$), and let $d_i\ge2$ denote the order of $x_i$. Then $\chi(\dot\oldGamma_0)=\chi(\oldGamma_0)-\sum_{i=1}^p(1-1/p_i)\le\chi(\oldGamma_0)<0$; this shows that $\dot\oldGamma$ is negative. 
\EndProof
%KX

\Lemma\label{historic fact}
Let
$\oldUpsilon$ and $\otheroldLambda$ be orientable $2$-orbifolds. Then 
two proper embeddings
 $\xi$ and $\xi'$ of $\oldUpsilon$ in $\otheroldLambda$ are properly isotopic if and only if $\dot\xi$ and $\dot\xi'$ are properly isotopic embeddings of $\dot\oldUpsilon$ in $\dot\otheroldLambda$. Furthermore, if a proper embedding $\eta:\dot\oldUpsilon\to\dot\otheroldLambda$ is properly isotopic to $\dot\xi$ for some embedding $\xi:\oldUpsilon\to\otheroldLambda$, then $\eta=\dot\xi'$ for some embedding $\xi':\oldUpsilon\to\otheroldLambda$ (which is properly isotopic to $\xi$ by the first assertion).
\EndLemma

\Proof
First suppose that $\xi$ and $\xi'$ are properly isotopic proper embeddings
 of $\oldUpsilon$ in $\otheroldLambda$. If
$(\xi_t)_{0\le t\le 1}$ is a proper isotopy with $\xi_0=\xi$
and $\xi_1=\xi'$, then for each $t\in[0,1]$ we have
$\xi_t^{-1}(\fraks_\otheroldLambda)=\fraks_\oldUpsilon$; since
$\fraks_\otheroldLambda$ is discrete by orientability, $(\xi_t)_{0\le t\le
  1}$ is constant on $\fraks_\oldUpsilon$, and therefore restricts to a
proper isotopy $(\xi_t\big|(|\oldUpsilon|-\fraks_\oldUpsilon))_{0\le t\le 1}$,
where $\xi_t\big|(|\oldUpsilon|-\fraks_\oldUpsilon)$ is an embedding of
$|\oldUpsilon|-\fraks_\oldUpsilon$ in $|\otheroldLambda|-\fraks_\otheroldLambda$ for each
$t$. By definition we have
$\xi_0\big|(|\oldUpsilon|-\fraks_\oldUpsilon)=\dot\xi$ and
$\xi_1\big|(|\oldUpsilon|-\fraks_\oldUpsilon)=\dot\xi'$, so that $\dot\xi$
and $\dot\xi'$ are isotopic. 

To complete the proof of the lemma, it now suffices to show that if
$\xi$ is an embedding of $\oldUpsilon$ in $\otheroldLambda$, and $\eta$
is a proper embedding of $\dot\oldUpsilon$ in $\dot\otheroldLambda$ which
is properly isotopic to $\dot\xi$, then there is a unique embedding
$\xi'$ of $\oldUpsilon$ in $\otheroldLambda$ such that $\eta=\dot\xi'$,
and $\xi'$ is properly isotopic to $\xi$. To prove this, fix a
proper isotopy $(\eta_t)_{0\le t\le 1}$, where $\eta_t$ is a proper
embedding of $\dot\oldUpsilon$ in $\dot\otheroldLambda$ for each $t$, and
such that $\eta_0=\dot\xi$ and and $\eta_1=\eta$. For each
$t\in[0,1]$, since $\eta_t$ is proper, it admits a unique extension
to a proper embedding $\hateta_t:Y\to L$, where $Y$ and $L$ denote
the end compactifications of the manifolds $\dot\oldUpsilon$ and $\dot\otheroldLambda$. We
have $\hateta_t(Y-\dot\oldUpsilon)\subset L-\dot\otheroldLambda$ for each
$t$.  Since $\fraks_\oldUpsilon$ and $\fraks_\otheroldLambda$ are discrete
subsets of the manifolds $|\oldUpsilon|$ and $|\otheroldLambda|$
respectively, $Y$ and $L$ are canonically identified with the end
compactifications of $|\oldUpsilon|$ and $|\otheroldLambda|$; in particular,
$|\oldUpsilon|$ and $|\otheroldLambda|$ are canonically identified with
subsets of $Y$ and $L$, and under these identifications we have
$\hateta_0\big||\oldUpsilon|=|\xi|$. Since $L-\dot\otheroldLambda$ is discrete, the map $t\mapsto\hateta_t(y)$ from $[0,1]$ to
$L-\dot\otheroldLambda$ is constant for every $y\in Y-\dot\oldUpsilon$. Hence
$\hateta_t|(Y-\dot\oldUpsilon)=|\xi|\big|(Y-\dot\oldUpsilon)$ for each
$t$. In particular, for each $t$, we have 
$\hateta_t|(|\oldUpsilon|-\dot\oldUpsilon)=|\xi|\big|(|\oldUpsilon|-\dot\oldUpsilon)$. Hence
$\hateta_t(|\oldUpsilon|)\subset |\otheroldLambda|$ for each $t$; and $\hateta_t(y)=|\xi|(y)$ is a point of $\fraks_\otheroldLambda$
having the same order as $y$, for every $t$ and for every $y\in\fraks_\oldUpsilon$. It
follows that $\hateta_t\big||\oldUpsilon|=|\xi_t|$ for some orbifold embedding
$\xi_t:\oldUpsilon\to\otheroldLambda$. Now $(\xi_t)_{0\le t\le 1}$ is
an isotopy of embeddings of $\oldUpsilon$ in $\otheroldLambda$, and we have
$\xi_0=\xi$ and $\dot\xi_1=\eta$. Thus if we set
$\xi'=\xi_1$ then $\eta=\dot\xi'$,
and $\xi'$ is properly isotopic to $\xi$. Since
$\dot\oldUpsilon$ is dense in $\oldUpsilon$, the embedding $\xi'$
of $\oldUpsilon$ in $\otheroldLambda$ is the only one for which $\dot\xi'=\eta$.
\EndProof
%\fraks

%\tG\otheroldLambda'\otheroldLambda_0\oldGamma

\Number\label{chidef}Let $\otheroldLambda$ be a negative, finite-type, orientable $2$-orbifold without boundary. We will
denote by $\Theta(\otheroldLambda)$ the set of all finite-type (possibly empty) $2$-suborbifolds $\oldUpsilon$ of
$\otheroldLambda$ such that (i) $\oldUpsilon$ is taut and (ii) no component of $|\oldUpsilon|$ which is a weight-$0$ annulus shares a boundary component with a
 component of $|\otheroldLambda-\inter\oldUpsilon|$ which is a weight-$0$ annulus. 

We will denote by
$\barcaly(\otheroldLambda)$ the set of all (orbifold-)isotopy classes of elements of $\Theta(\otheroldLambda)$. If $\oldUpsilon\in\Theta(\otheroldLambda)$ is given, we will denote its isotopy class by $[\oldUpsilon]\in\barcaly(\otheroldLambda)$.
We will denote by $\Theta_-(\otheroldLambda)$ the set of all elements
$\oldUpsilon\in\Theta(\otheroldLambda)$ that are negative orbifolds, and by
$\barcaly_-(\otheroldLambda)$ the set of all elements of
$\barcaly(\otheroldLambda)$ that have the form $[\oldUpsilon]$ for some
$\oldUpsilon\in\Theta_-(\otheroldLambda)$. Note that any finite-type,
negative, taut suborbifold of $\otheroldLambda$ belongs to
$\Theta_-(\otheroldLambda)$, because such a suborbifold has no annular
components, and therefore automatically satisfies Condition (ii) in the definition of
$\Theta(\otheroldLambda)$.

Note also that any union of components of $\otheroldLambda$ is an element of $\Theta_-(\otheroldLambda)$; elements of this type will play a role in the arguments to be given in Section \ref{higher section} (see for example \ref{oldXi}).

We have a well-defined function $\chi$ on $\barcaly(\otheroldLambda)$ given
by $\chi([\oldUpsilon])=\chi(\oldUpsilon)$. We also set
$\chibar([\oldUpsilon])=-\chi([\oldUpsilon])=\chibar(\oldUpsilon)$. Furthermore,
we have functions $\newalpha$, $\beta$, and $\newgamma$ given on
$\barcaly(\otheroldLambda)$ by setting 
$\beta([\oldUpsilon])=\compnum(\partial|\oldUpsilon|)$ and $\newgamma(\oldUpsilon)=\card\fraks_\oldUpsilon$, and defining
$\newalpha([\oldUpsilon])$  to be the number of components of $|\oldUpsilon|$ that are weight-$0$ annuli.

Let us define a relation $\preceq$ on $\barcaly(\otheroldLambda)$ by writing $[\oldUpsilon_1]\preceq[\oldUpsilon_2]$ if $\oldUpsilon_1$ is isotopic in $\otheroldLambda$ to a suborbifold of $\oldUpsilon_2$; this is clearly independent of the choice of the $\oldUpsilon_i$ in their isotopy classes.
\EndNumber

\Lemma\label{more geventlach}
Let $\otheroldLambda$ be an orientable $2$-orbifold without boundary, and let
$\oldGamma$ be an element of $\Theta(\otheroldLambda)$. Then
$\dot\oldGamma$ is an element of
$\Theta(\dot\otheroldLambda)$. 
\EndLemma

\Proof
Since $\oldGamma$ has finite type it is automatic that $\dot\oldGamma$ has finite type. Since $\oldGamma$ is taut, it follows from Lemma \ref{geventlach} that $\dot\oldGamma$ is taut. Now suppose that some annular component $\frakZ$ of $\dot\oldGamma$ shares a boundary component with an
annular component $\frakZ'$ of $\dot\otheroldLambda-\inter\dot\oldGamma$. Since $\frakZ$ and $\frakZ'$ are in particular compact, they are components of $\oldGamma$ and $\otheroldLambda-\inter\oldGamma$ respectively. But this is impossible, since no annular component  of $\oldGamma$ shares a boundary component with an
annular component of $\otheroldLambda-\inter\oldGamma$.
\EndProof

\Lemma\label{old partial order}
Let $\otheroldLambda$ be a negative, finite-type, orientable $2$-orbifold  without boundary.
%, the relation
%$\preceq$ is a partial ordering on
%$\barcaly(\otheroldLambda)$. Furthermore, 
Let $Y_1,Y_2$ be elements of $\barcaly(\otheroldLambda)$ such that  $Y_1\preceq Y_2$. Then
$\chibar(Y_1)\le\chibar(Y_2)$ and $\newgamma(Y_1)\le\newgamma(Y_2)$.
If
$\chibar(Y_1)=\chibar(Y_2)$, then the
equivalence classes $Y_1$ and $Y_2$ are respectively represented 
%by
%elements
 %$\oldUpsilon_1,\oldUpsilon_2\in\Theta(\otheroldLambda)$
by elements $\oldUpsilon_1,\oldUpsilon_2$ of $\Theta(\otheroldLambda)$ such that $\oldUpsilon_1\subset\inter\oldUpsilon_2$ and every component of $\oldUpsilon_2-\inter\oldUpsilon_1$ is  annular.
Finally, if
$\chibar(Y_1)=\chibar(Y_2)$ and $\newgamma(Y_1)=\newgamma(Y_2)$, then the
equivalence classes $Y_1$ and $Y_2$ are respectively represented 
%by
%elements
 %$\oldUpsilon_1,\oldUpsilon_2\in\Theta(\otheroldLambda)$
by elements $\oldUpsilon_1,\oldUpsilon_2$ of $\Theta(\otheroldLambda)$ such that $\oldUpsilon_1\subset\inter\oldUpsilon_2$ and every component of $|\oldUpsilon_2-\inter\oldUpsilon_1|$ is a weight-$0$ annulus.
\EndLemma

\Proof
Since $Y_1\preceq Y_2$,
there are elements $\oldUpsilon_1$ and $\oldUpsilon_2$ of
$\Theta(\otheroldLambda)$ such that $\oldUpsilon_1\subset\inter\oldUpsilon_2$ and
$Y_i=[\oldUpsilon_i]$ for $i=1,2$. Since the $\oldUpsilon_i$ have
$\pi_1$-injective boundaries, so does 
 $\oldGamma:=\oldUpsilon_2-\inter\oldUpsilon_1$; in view of the negativity of $\otheroldLambda$, it follows that
 each component of $\oldGamma:=\oldUpsilon_2-\inter\oldUpsilon_1$ has non-positive Euler characteristic. In particular we have $\chibar(\oldGamma)\ge0$. Thus
\Equation\label{en passant}
\chibar(Y_2)-\chibar(Y_1)=\chibar(\oldUpsilon_2)-\chibar(\oldUpsilon_1)=\chibar(\oldGamma)\ge0,
\EndEquation
so that $\chibar(Y_1)\leq\chibar( Y_2)$. Note also that since $Y_1$ is a suborbifold of $Y_2$ we have $\fraks_{Y_1}\subset\fraks_{Y_2}$, and hence $\newgamma(Y_1)\le\newgamma(Y_2)$.

In the special case $\chibar(Y_2)=\chibar(Y_1)$, it follows from
(\ref{en passant}) that $\chi(\oldGamma)=0$; since each component
of $\oldGamma:=\oldUpsilon_2-\inter\oldUpsilon_1$ has non-positive
Euler characteristic, it follows that each component of $\oldGamma$
has Euler characteristic $0$. But since $\otheroldLambda$ is negative, no component of $\otheroldLambda$ has Euler characteristic $0$; hence %$\oldGamma$  has no component with empty boundary and Euler characteristic $0$; thus 
every component of $\oldGamma$ has
non-empty boundary. A connected, orientable $2$-orbifold that has non-empty
boundary and Euler characteristic $0$ is either an annular orbifold or
a $2$-manifold homeomorphic to $\SSS^1\times[0,\infty)$. If a component
of $\oldGamma$ were homeomorphic to $\SSS^1\times[0,\infty)$, then in particular some component of $\partial\oldUpsilon_1$  would bound a suborbifold of $\otheroldLambda$ which is a $2$-manifold  homeomorphic to $\SSS^1\times[0,\infty)$; since $\oldUpsilon_1$ is taut, this contradicts Condition (iii) of Definition
\ref{praxis}. Hence every component of
$\oldGamma$ is annular. If in addition we have
$\newgamma(Y_1)=\newgamma(Y_2)$, then $\card(\fraks_{\oldGamma})=\card(\fraks_{\oldUpsilon_2}-\fraks_{\oldUpsilon_1})=\newgamma(Y_2)-\newgamma(Y_2)=0$. Thus each component of $|\oldGamma|$ has weight $0$, and since the components of $\oldGamma$ are annular and orientable, the components of $|\oldGamma|$ must be weight-$0$ annuli.
\EndProof
%\alpha\gamma

\Proposition\label{new partial order}
For any negative, finite-type, orientable $2$-orbifold $\otheroldLambda$ without boundary, the relation
$\preceq$ is a partial ordering on
$\barcaly(\otheroldLambda)$. Furthermore, if $Y_1\preceq Y_2$ for given
elements $Y_1,Y_2\in\barcaly(\otheroldLambda)$, then in terms of the lexicographical order on $\ZZ^{4}$, we have 
$$(\chibar(Y_1),\newgamma(Y_1),\newalpha(Y_1),-\beta(Y_1))\le(\chibar(Y_2),\newgamma(Y_2),\newalpha(Y_2),-\beta(Y_2)),$$
%\chibar(Y_2),\beta$ and $\newgamma(Y_1)\le\newgamma(Y_2)$; and if $\chibar(Y_1)=\chibar(Y_2)$ and $\newgamma(Y_1)=\newgamma(Y_2)$, then the equivalence classes $Y_1$ and $Y_2$ are respectively represented by elements $\oldUpsilon_1,\oldUpsilon_2\in\Theta$ such that $\oldUpsilon_1\subset\inter\oldUpsilon_2$ and every component of $\oldUpsilon_1-\inter\oldUpsilon_2$ is an annular orbifold. Finally, if $Y_1\preceq Y_2$ and $\chibar(Y_1)=\chibar(Y_2)$, then $\beta(Y_1)\ge\beta(Y_2)$, 
with equality only if $Y_1=Y_2$. 
\EndProposition

\Proof
We first observe that the second sentence of the proposition implies
the first. Indeed, to prove the first sentence of the proposition we need only prove
antisymmetry, since reflexivity and transitivity are obvious. If
$Y_1,Y_2\in\barcaly(\otheroldLambda)$ satisfy $Y_1\preceq Y_2$, then the
second sentence of the present proposition implies that
$(\chibar(Y_1),\newgamma(Y_1),\newalpha(Y_1),-\beta(Y_1))\le(\chibar(Y_2),\newgamma(Y_2),\newalpha(Y_2),-\beta(Y_2))$,
in terms of the lexicographical order on $\ZZ^{4}$. If in addition we
have  
$Y_2\preceq Y_1$, then the second sentence implies that
$(\chibar(Y_2),\newgamma(Y_2),\newalpha(Y_2),-\beta(Y_2))\le(\chibar(Y_1),\newgamma(Y_1),\newalpha(Y_1),-\beta(Y_1))$. Since
the lexicographical order is antisymmetric, we have
$(\chibar(Y_1),\newgamma(Y_1),\newalpha(Y_1),-\beta(Y_1))=(\chibar(Y_2),\newgamma(Y_2),\newalpha(Y_2),-\beta(Y_2))$. The
second sentence then  implies that $Y_1=Y_2$, and the antisymmetry of
$\preceq$ follows.

To prove the second sentence, assume that
$Y_1,Y_2\in\Theta(\otheroldLambda)$ satisfy $Y_1\preceq Y_2$. Then according to Lemma \ref{old partial order},
we have $\chibar(Y_1)\le\chibar(Y_2)$ and $\newgamma(Y_1)\le\newgamma(Y_2)$.
For the rest of the proof, we will assume that $\chibar(Y_1)=\chibar(Y_2)$ and
$\newgamma(Y_1)=\newgamma(Y_2)$. We will prove that in terms of the
lexicographical ordering of $\ZZ^2$ we have
$(\newalpha(Y_1),-\beta(Y_1))\le(\newalpha(Y_2),-\beta(Y_2))$,  with
equality only if $Y_1=Y_2$. This will complete the proof.
%\alpha\gamma

Since $\chibar(Y_1)=\chibar(Y_2)$ and
$\newgamma(Y_1)=\newgamma(Y_2)$,  Lemma \ref{old partial order} gives
elements $\oldUpsilon_1,\oldUpsilon_2$ of $\Theta(\otheroldLambda)$ with
$[\oldUpsilon_i]=Y_i$ for $i=1,2$, such that
$\oldUpsilon_1\subset\inter\oldUpsilon_2$, and such that if we set
$\oldGamma:=\oldUpsilon_2-\inter\oldUpsilon_1$, then every component of $|\oldGamma|$ is a weight-$0$ annulus. We claim: 
\Claim\label{mitzi}
Every weight-$0$ annulus component of $|\oldUpsilon_1|$ is contained in
a weight-$0$ annulus component of $|\oldUpsilon_2|$, and every weight-$0$ annulus component
of $|\oldUpsilon_2|$ contains at most component of $|\oldUpsilon_1|$.
\EndClaim

To prove \ref{mitzi}, let $\frakB$ be any component of
$\oldUpsilon_1$ such that $|\frakB|$ is a weight-$0$ annulus, and let $\frakH$ denote the component of
$\oldUpsilon_2$ containing $\frakB$. Consider an arbitrary boundary
component $C$ of $|\frakB|$, and let
$\frakE$ denote the component of $\frakH\setminus\inter\oldUpsilon_1$
containing $\obd(C)$. Since $\frakE$ is in particular a component of
$\oldGamma$, the surface $|\frakE|$ is a weight-$0$ annulus. We have
$\partial\frakE\subset\partial\frakH\cup\partial\oldUpsilon_1$. If
$\partial\frakE\subset\partial\oldUpsilon_1$, then $\frakE$ is a
component of $\otheroldLambda-\inter\oldUpsilon_1$. 
Thus $|\frakB|$ and $|\frakE|$ are respectively weight-$0$ annulus components of $|\oldUpsilon_1|$ and $|\otheroldLambda-\inter\oldUpsilon_1|$, and they
share the
boundary component $\obd(C)$; 
since
$\oldUpsilon_1\in\Theta(\otheroldLambda)$,
%The annular component
%$\frakB$ of $\oldUpsilon_1\in\Theta(\otheroldLambda)$ then shares the
%boundary component $\obd(C)$ with the annular
%component  $\frakE$ of $\otheroldLambda-\inter\oldUpsilon_1$; 
this
contradicts Condition (ii) of the definition of $\Theta(\otheroldLambda)$ (see
\ref{chidef}). Hence $|\partial\frakE|$ must have a component $C'$ which
is contained in $|\partial\frakH|$. 
%Now the annular orbifold $\frakE$ has distinct
%boundary components  $\obd(C)$ and $\obd(C')$, and hence
%$|\frakE|$ must be a weight-$0$ annulus with
%$\partial|\frakE|=C\cup C'$. 
This shows that for every boundary component $C$ of
the weight-$0$ annulus $|\frakB|$, there is a weight-$0$ annulus having $C$ as one
boundary component, having its other boundary component contained in
$|\partial\frakH|$, and meeting $|\frakB|$ only in $C$. Hence $|\frakH|$ is
a weight-$0$
%n orbifold regular neighborhood of $\frakB$, and is therefore
annulus.

Now suppose that for some  component $\frakH$ of $\oldUpsilon_2$, the surface $|\frakH|$ is a weight-$0$ annulus containing two distinct components of $|\oldUpsilon_1|$. Since
$|\frakH|$ is a weight-$0$ annulus and $\oldUpsilon_1$ is taut, every component
of $|\frakH\cap\oldUpsilon_1|$ or of
$|\frakH\setminus\inter\oldUpsilon_1|$ is a weight-$0$ annulus. Since
$\frakH\cap\oldUpsilon_1$ has at least two components, and $\frakH$ is
connected, there is a component $\frakE$ of
$\frakH\setminus\inter\oldUpsilon_1$ that meets at least two
components of $\oldUpsilon_1$. Since the weight-$0$ annulus $|\frakE|$  has only
two boundary components, we must have
$\partial\frakE\subset\oldUpsilon_1$, and hence $\frakE$ is a
component of $\otheroldLambda-\inter\oldUpsilon_1$. If $\frakD$ denotes one of
the components of $\oldUpsilon_1$ that meet $\frakE$, then $|\frakD|$
and $|\frakH|$ are weight-$0$ annulus components of $|\oldUpsilon_1|$ and
$|\otheroldLambda-\inter\oldUpsilon_1|$ that share a boundary component, and
we again have a contradiction to Condition (ii) of the definition of
$\Theta(\otheroldLambda)$. This completes the proof of \ref{mitzi}.

It follows from \ref{mitzi} that $\newalpha(Y_1)\preceq\newalpha(Y_2)$. Now
assume that  $\newalpha(Y_1)=\newalpha(Y_2)$; we must show that
$\beta(Y_1)\ge\beta(Y_2)$, with equality only if $Y_1=Y_2$. From
\ref{mitzi} and the assumption that $\newalpha(Y_1)=\newalpha(Y_2)$ we
deduce:
\Claim\label{more mitzi}
Every annular component of $\oldUpsilon_2$ contains an annular
component of $\oldUpsilon_1$.
\EndClaim

Next note that 
%To prove the third sentence, again suppose that we have $Y_1\preceq
%Y_2$, $\chibar(Y_1)=\chibar(Y_2)$, and
%$\gamma(Y_1)=\gamma(Y_2)$. Then by the assertion just proved,
% each component of $\oldGamma$ is an annular orbifold. In particular, 
for each component $\oldDelta $ of $\oldGamma$, the weight-$0$ annulus $|\oldDelta|$ has two boundary components,
%an annular orbifold, the $2$-manifold $|\oldDelta |$ is an annulus or
%a disk. 
each  of which is contained in
either $\partial|\oldUpsilon_1|$ or $\partial|\oldUpsilon_2|$. For
$i=1,2$, let  $c_i(\oldDelta )$ denote the number of components of
$|\partial\oldDelta |$ contained in $\partial |\oldUpsilon_i|$. We have
$c_1(\oldDelta )+c_2(\oldDelta )=\compnum(\partial|\oldDelta
|)=2$.
If $c_1(\oldDelta)=0$ then $\partial |\oldDelta |\subset\partial|\oldUpsilon_2|$,
which implies that $\oldDelta $ is a component of $\oldUpsilon_2$;
since the annular orbifold $\oldDelta$ is a component of $\oldGamma$,
and is therefore by definition disjoint from $\inter\oldUpsilon_1$, this contradicts
\ref{more mitzi}. Hence $c_1(\oldDelta)\ge1$ for each component $\oldDelta $ of $\oldGamma$. It follows that
$c_1(\oldDelta)\ge c_2(\oldDelta)$ for each $\oldDelta$. On the other hand, since $\oldUpsilon_1\subset\inter\oldUpsilon_2$, we have $\partial\oldUpsilon_i\subset\partial\oldGamma$ for $i=1,2$, so that 
$\compnum(\partial|\oldUpsilon_i|)=\sum_{\oldDelta \in\calc(\oldGamma)}c_i(\oldDelta)$.
We now have
\Equation\label{the mitziest}
\beta(Y_1)=\compnum(\partial|\oldUpsilon_1|)=\sum_{\oldDelta \in\calc(\oldGamma)}c_1(\oldDelta)
\ge
\sum_{\oldDelta \in\calc(\oldGamma)}c_2(\oldDelta)=\compnum(\partial|\oldUpsilon_2|)
=
\beta(Y_2).
\EndEquation
This establishes the inequality $\beta(Y_1)
\ge
\beta(Y_2)$. 

Finally, assume that $\beta(Y_1)=\beta(Y_2)$. Then (\ref{the mitziest}), together with the inequality $c_1\ge c_2$,
implies that $c_1(\oldDelta )=c_2(\oldDelta )$ for each component $\oldDelta $ of $\oldGamma$. Since for each $\oldDelta $ we have %$c_1(\oldDelta )\ge1$ and
% $\compnum(\partial|\oldDelta |)=
$c_1(\oldDelta )+c_2(\oldDelta )=2$, we must have $c_1(\oldDelta )=c_2(\oldDelta )=1$.
% and $\compnum(\partial|\oldDelta |)=2$. Since $\oldDelta $ is an annular orbifold, $|\oldDelta |$ must be an
Thus the weight-$0$ annulus $|\oldDelta |$ has one boundary component in $\partial\oldUpsilon_1$ and one in $\partial\oldUpsilon_2$. As this holds for every component $\oldDelta $ of $\oldGamma$, it follows that $\oldUpsilon_2$ is an orbifold regular neighborhood of $\oldUpsilon_1$, and hence that $Y_1=Y_2$.
\EndProof

\Number\label{involute}
Let   $\otheroldLambda$ be a negative,  orientable $2$-orbifold without boundary.
It follows from the definition of
tautness (\ref{praxis}) that if $\oldUpsilon$ is a taut suborbifold of $\otheroldLambda$, then
$\otheroldLambda-\inter\oldUpsilon$ is also taut.
It then follows from the definition of
$\Theta(\otheroldLambda)$ (\ref{chidef}) that for any $\oldUpsilon\in\Theta(\otheroldLambda)$
we have
$\otheroldLambda-\inter\oldUpsilon\in\Theta(\otheroldLambda)$. 
%\redcomment{I'm not
  %sure this is right. If $\oldUpsilon$ is homeomorphic to
  %$\SSS^1\times[0,1)$, I think we may have $\oldUpsilon\in\Theta(\otheroldLambda)$
%but
%$\otheroldLambda-\inter\oldUpsilon\notin\Theta(\otheroldLambda)$, according to
%the way tautness is defined. Maybe the solution at this point is
%to assume that $\otheroldLambda$ is closed. 
The
assignment $\oldUpsilon\mapsto\otheroldLambda-\inter\oldUpsilon$ is an
involution of $\Theta(\otheroldLambda)$, and defines an involution
$[\oldUpsilon]\mapsto[\otheroldLambda-\inter\oldUpsilon]$ of
$\barcaly(\otheroldLambda)$. Note that this involution of
$\barcaly(\otheroldLambda)$ is order-reversing with respect to the partial
order $\preceq$.
\EndNumber

The statement of the following lemma involves the special case of the notation introduced in \ref{chidef} in which the given boundaryless $2$-orbifold is a $2$-manifold.

\Lemma\label{another stupidity}
Let $L$ be an orientable $2$-manifold without boundary, and let
$F$ be an element of $\Theta(L)$. Then no two annulus components of $F$ are isotopic.
\EndLemma

\Proof
Let $X$ and $X'$ be annulus components of $F$. If $X$ and $X'$ are isotopic, their core curves cobound an annulus by \cite[Lemma 2.4]{epstein}. Hence some component $Y$ of $\overline{L-(X\cup X')}$ is an annulus having one boundary component $C$ contained in $X$ and one boundary component $C'$ contained in $X'$. Then $C$ and $C'$ are components of $\partial F$; any other component of $\partial F$ which intersects $Y$ must be a simple closed curve contained in $\inter Y$, and is homotopically non-trivial since $F$ is taut (by Condition (i) in the definition of $\Theta(L)$). Hence the closure of the component of $Y\setminus\partial F$ containing $C$ is an annulus $Y_0$ having $C$ as one boundary component and having its other boundary component contained in $\partial F$; since $Y$ is a component of $\overline{L-(X\cup X')}$, the annulus $Y_0$ is a component of $\overline{L-F}$. But then the annulus component $X$ of $F$ shares a boundary component with the annulus component $Y_0$ of $\overline{L-F}$, a contradiction to Condition (ii) in the definition of $\Theta(L)$.
\EndProof

\begin{definitionremarks}\label{nondeg cross}
Let $\otheroldLambda$ be a finite-type orientable $2$-orbifold without boundary, and
let $C_1$ and $C_2$ be  $1$-submanifolds of
$|\otheroldLambda|-\fraks_\otheroldLambda$. We will say that $C_1$ and $C_2$
have a {\it degenerate crossing} if there is a  disk
$D\subset|\otheroldLambda|-\fraks_\otheroldLambda$ such that $\partial D$ has
the form $a_1\cup a_2$, where $a_i$ is an arc contained in $C_i$ for
$i=1,2$, and $a_1\cap a_2=\partial a_1=\partial a_2$. Such a disk will
be called a {\it disk of degeneracy} for $C_1$ and $C_2$.

If $C_1$ and $C_2$ have a degenerate crossing, there is a disk of
degeneracy $D_0$ which is minimal with respect to inclusion
among all disks of degeneracy for $C_1$ and $C_2$. The minimality of
$D_0$ implies that $D_0\cap(C_1\cup C_2)=\partial D_0$. It then
follows that $C_1$ is isotopic, by an ambient orbifold isotopy which is
supported on a small neighborhood of $D_0$, to a $1$-manifold $C_1'$
such that $\card(C_1'\cap C_2)<\card(C_1\cap C_2)$. Hence if $C_1$ and
$C_2$ are compact $1$-submanifolds of
$|\otheroldLambda|-\fraks_\otheroldLambda$ which intersect transversally, and have been chosen with their isotopy
classes in  $|\otheroldLambda|-\fraks_\otheroldLambda$ so as to minimize the number of their intersection points,
then they do not have a degenerate crossing.

Suppose that $\oldUpsilon_1$ and
$\oldUpsilon_2$ are taut finite-type $2$-suborbifolds of $\otheroldLambda$. Thus each $|\partial\oldUpsilon_i|$ is a submanifold of $|\otheroldLambda|$ whose boundary is disjoint from $\fraks_\otheroldLambda$, cf. Definition \ref{praxis}. We will say that $\oldUpsilon_1$ and
$\oldUpsilon_2$ are  in {\it standard position} if (1) the $1$-submanifolds $|\partial\oldUpsilon_1|$
and $|\partial\oldUpsilon_2|$ of  $|\otheroldLambda|-\fraks_\otheroldLambda$ meet transversally and have no degenerate
crossing, and (2) no component of
$|\overline{\otheroldLambda-(\oldUpsilon_1\cup\oldUpsilon_2)}|$ is a
weight-$0$ annulus in $|\otheroldLambda|$ having one boundary component
contained in $|\partial\oldUpsilon_1|$ and  one contained in
$|\partial\oldUpsilon_2|$. Note that if
$|\overline{\otheroldLambda-(\oldUpsilon_1\cup\oldUpsilon_2)}|$ does have a
component $A$ of the type ruled out by Condition (2), and if we set $\oldUpsilon_i'=\oldUpsilon_i\cup\obd(A)$ for $i=1,2$, then  $\oldUpsilon_i'$ is isotopic to $\oldUpsilon_i$, while $\card(\partial\oldUpsilon_1'\cap \partial\oldUpsilon_2')=\card(\partial\oldUpsilon_1\cap\partial\oldUpsilon_2)$, and $\compnum(\otheroldLambda-(\oldUpsilon_1'\cup\oldUpsilon_2'))<\compnum(\otheroldLambda-(\oldUpsilon_1\cup\oldUpsilon_2))$. Hence if
$\oldUpsilon_1$ and
$\oldUpsilon_2$ have been chosen with their isotopy
classes so as to minimize the element  
$(\card
(\partial\oldUpsilon_1\cap
\partial\oldUpsilon_2),\compnum(\otheroldLambda-(\oldUpsilon_1'\cup\oldUpsilon_2')))$
with respect to the lexicographical ordering of $\NN\times\NN$, then 
$\oldUpsilon_1$ and $\oldUpsilon_2$ of $\otheroldLambda$ are in
standard position.

In particular, any two taut finite-type suborbifolds of
$|\otheroldLambda|$ are orbifold-isotopic to suborbifolds that
are in standard position.
\end{definitionremarks}

\Number\label{dot standard}
Let $\otheroldLambda$ be a finite-type orientable $2$-orbifold without boundary. If $\oldXi$ is any $2$-suborbifold of $|\otheroldLambda$, then $|\oldXi|$ is a weight-$0$ disk or a weight-$0$ annulus if and only if $|\dot\oldXi|$ is, respectively, a weight-$0$ disk or a weight-$0$ annulus (in which case $\dot\oldXi=\oldXi$).
It follows that if $C_1$ and $C_2$ are compact $1$-submanifolds of $|\otheroldLambda|-\fraks_\otheroldLambda$, then $C_1$ and $C_2$ have a degenerate crossing in $\otheroldLambda$ if and only if $C_1$ and $C_2$, regarded as curves in $\dot\otheroldLambda$, have a degenerate crossing; and that $C_1$ and $C_2$ are in standard position in $\otheroldLambda$ if and only if $C_1$ and $C_2$, regarded as curves in $\dot\otheroldLambda$, are in standard position.
\EndNumber

The following result generalizes \cite[Lemma 2.5]{epstein}.

\Proposition\label{beats epstein}
Let $C_1$ and $C_2$ be disjoint, $\pi_1$-injective simple closed curves in an orientable $2$-manifold $L$. Suppose that $C_1'$ and $C_2'$ are simple closed curves in $L$ such that (i) $C_i'$ is homotopic to $C_i$ for $i=1,2$, and (ii) the intersection of $C_1'$ and $C_2'$ is non-empty and transverse. Then $C_1'$ and $C_2'$ have a degenerate crossing.
\EndProposition

\Proof
According to \cite[Theorem 2.1]{epstein}, $C_1$ and $C_1'$ are isotopic. Hence after modifying $C_1'$ and $C_2'$ by a (single) ambient isotopy, we may assume that $C_1'=C_1$. 
Let $A$ be an annulus neighborhood of $C_1$ disjoint from $C_2$. Since $C_2'$ meets $C_1$ transversally, we may choose $A$ in such a way that $A$ admits a homeomorphic identification with $\SSS^1\times[-1,1]$ under which $\SSS^1\times\{0\}=C$ and $A\cap C_2'=E\times[-1,1]$ for some finite set $E\subset \SSS^1$.
%Since $C_1\cap C_2'=C_1'\cap C_2'\ne\emptyset$, we may choose $A$ so that $(\partial A)\cap C_2'\ne\emptyset$. 
Let $G$ denote the component of $L-\inter A$ containing $C_2$.
Fix a basepoint  $x\in G\subset L$,
and let $p:(\tL,\tx)\to ( L,x)$ denote the based covering space of
$( L,x)$ defined by the subgroup of $\pi_1( L,x)$
which is the image of the inclusion homomorphism
$\pi_1( G,x)\to\pi_1( L,x)$. Let
$s:( G,x)\to(\tL,\tx)$ denote the based lift
of the inclusion map $( G,x)\to( L,x)$; then
$p$ maps  $\tGee:=s( G)$ homeomorphically onto
$ G$. The definition of the covering $\tL$ implies
that the inclusion homomorphism
$\pi_1(\tGee)\to\pi_1(\tL)$ is
surjective. Furthermore, $\partial\tGee$ is $\pi_1$-injective in $\tL$ since $C_1$ is $\pi_1$-injective in $ L$. Hence every component of  %$\tGee$ is a deformation retract of $\tL$, and every component of  
 $\tL-(\inter\tGee)$ is a half-open annulus whose boundary is a component of $\partial\tGee$. 

The curve $C_2$ is contained in $G$ and therefore admits a lift to $\tL$. Since $C_2$ and $C_2'$ are homotopic, $C_2'$ also admits a lift to $\tL$; that is, there is a simple closed curve $\tC_2'\subset \tL$ which is mapped homeomorphically onto $C_2'$ by $p$. Since $p$ maps $\tGee$ homeomorphically onto $G$, and $C_1\cap C_2'=C_1'\cap C_2'\ne\emptyset$, we cannot have $\tC_2'\subset\inter\tGee$.

Consider the case in which $\tC_2'\cap\tGee=\emptyset$. In this case, $\tC_2'$ is a $\pi_1$-injective simple closed curve contained in some component of  $\tL-(\inter\tGee)$, which is a half-open annulus whose boundary is a component of $\partial\tGee$. It follows that $\tC_2'$ is homotopic to a component of $\partial\tGee$. But the construction of $\tGee$ implies that each of its boundary components is mapped homeomorphically by $p$ onto a boundary component of an annulus neighborhood of $C_1$. Hence $C_2'$ and $C_1$ are homotopic. Since the intersection of $C_1$ and $C_2'$ is non-empty and transversal, it follows from \cite[Lemma 2.5]{epstein} that $C_1$ and $C_2'$ have a degenerate crossing.

There remains the case in which $\tC_2'\cap\partial\tGee\ne\emptyset$. In this case, in view of transversality, $\tC_2'\cap(\tL-\inter\tGee)$ is a non-empty disjoint union of properly embedded arcs in $\tL-\inter\tGee$. Choose one of these arcs, say $a$. Since the component of $\tL-\inter\tGee$ containing $a$ is a half-open annulus, $a$ is the frontier of a disk in $\tL-\inter\tGee$. This shows that $\tC_2'$ has a degenerate crossing with some component of $\partial\tGee$. In view of our choice of $A$, this implies that $\tC_2'$ has a degenerate crossing with some simple closed curve  which is mapped homeomorphically onto $C_1$ by $p$. In particular $p^{-1}(C_1)$ and $p^{-1}(C_2')$ have a degenerate crossing.

Among all disks of degeneracy for $p^{-1}(C_1)$ and $p^{-1}(C_2')$, choose one, say $D$, which is minimal with respect to inclusion. Write $\partial D=a_1\cup a_2$, where $a_i$ is an arc contained in $p^{-1}(C_i)$, and $a_1\cap a_2=\partial a_1=\partial a_2$. If $(\inter D)\cap p^{-1}(C_1)\ne\emptyset$, then some component of $D\cap p^{-1}(C_1)\ne\emptyset$ is a properly embedded arc $a_1'$ in $D$, which must have its endpoints in $a_2$; it follows that $a_1'$ is the frontier of a disk of degeneracy for
$p^{-1}(C_1)$ and $p^{-1}(C_2')$, properly contained in $D$. This contradicts minimality, and it follows that $D\cap p^{-1}(C_1)=a_1$. The same argument shows that $D\cap p^{-1}(C_2')=a_2$.

We claim:
\Claim\label{brownell}
 For each point $u\in\partial a_1=\partial a_2$, the map  $p|((\partial D)-\{u\})$ is one-to-one.
\EndClaim
To prove \ref{brownell}, suppose that $x$ and $y$ are distinct points of $(\partial D)-\{u\}$ such that $p(x)=p(y)$. First consider the subcase in which $x$ and $y$ both lie on $a_1$. Since  $x,y\ne u$, at most one of the points $x$ and $y$ can be an endpoint of $a_1$. By symmetry we may assume that $x$ is not an endpoint of $a_1$. If $c$ denotes the sub-arc  of $a_1$ having $x$ and $y$ as endpoints, then $p$ must map $c-\{y\}$ onto $C_1$. Since $C_1\cap C_2'\ne\emptyset$, it follows that $(c-\{y\})\cap p^{-1}(C_2')\ne\emptyset$, a contradiction since
$D\cap p^{-1}(C_2')=a_2$ is disjoint from $c-\{y\}$. We obtain the same contradiction if $x$ and $y$ both lie on $a_2$. This leaves the subcase in which one of $x$ and $y$ lies on $\inter a_1$ and the other on $\inter a_2$; by symmetry, we may assume that $x\in\inter a_1$ and $y\in\inter a_2$. Since $p(x)=p(y)\in C_2'$, we have $x\in(\inter a_1)\cap p^{-1}(C_2')$, a contradiction since 
$D\cap p^{-1}(C_2')=a_2$. This proves \ref{brownell}.

Now let  $u_1$ and $u_2$ denote the points of $\partial a_1=\partial a_2$. According to \ref{brownell}, $p|((\partial D)-\{u_i\})$ is one-to-one for $i=1,2$. Hence if $p|\partial D$ is not one-to-one, it follows from \ref{brownell}  that $p(a_i)=C_i'$ for $i=1,2$, and that $C_1=C_1'$ and $C_2'$ meet precisely in the point $p(u_1)=p(u_2)$. Since $C_1$ and $C_2'$ meet transversally, their homological intersection number is $1$, a contradiction since $C_2'$ is homotopic to $C_2$ and disjoint from $C_1$. Hence $p|\partial D$ is one-to-one. It now follows from \cite[Lemma 1.6]{epstein} that $p|D$ is one-to-one. This implies that $p(D)$ is a disk of degeneracy for $C_1=C_1'$ and $C_2$.
\EndProof

\Number\label{trimming}
Let $\otheroldLambda$ be a negative, orientable, finite-type $2$-orbifold without boundary. If $\oldUpsilon$ is a $\pi_1$-injective $2$-suborbifold of $\otheroldLambda$, let us define $\oldDelta=\oldDelta_\oldUpsilon^\otheroldLambda$ to be the union of all  components of $\oldUpsilon$ that either are discal or are manifolds homeomorphic to $\SSS^1\times\RR$, and let us define $\frakB=\frakB_\oldUpsilon^\otheroldLambda$ to be the union of all components $\frakC$ of $\oldUpsilon$ such that $|\frakC|$ is a weight-$0$ annulus. (We will denote  $\oldDelta_\oldUpsilon^\otheroldLambda$ and $\frakB_\oldUpsilon^\otheroldLambda$ by $\oldDelta_\oldUpsilon$ and $\frakB_\oldUpsilon$ when it is clear which ambient orbifold $\otheroldLambda$ is involved.) We will say that a component $\frakC$ of $\frakB$ is {\it redundant}  if there exist a component $\frakE$ of $\oldUpsilon-(\frakB\cup\oldDelta)$ and a weight-$0$ annulus $R\subset|\otheroldLambda-\inter(\frakC\cup\frakE)|$ having one boundary component in $\partial |\frakC|$ and one in $\partial|\frakE|$. We will say that components $\frakC,\frakC'$ of $\frakB$ are {\it similar} if there is a weight-$0$ annulus $R\subset|\otheroldLambda-\inter(\frakC\cup\frakC')|$ having one boundary component in $\partial |\frakC|$ and one in $\partial|\frakC'|$. Then similarity is an equivalence relation on the set of components of $\frakB$, and any component of $\frakB$ which is similar to a redundant component is redundant.  We define a {\it trimming} of $\oldUpsilon$ to be a suborbifold of $\oldUpsilon$ having the form $(\oldUpsilon-(\frakB\cup\oldDelta))\cup\frakB_0$, where $\frakB_0$ is a union of components of $\frakB$ which contains no redundant components, and contains exactly one representative of each similarity class of non-redundant components.

Although the trimming of $\oldUpsilon$ depends on the choice of a union of components $\frakB_0$ of $\frakB$ having the properties stated above, it follows from the definitions that the trimming is uniquely defined up to isotopy.

Note also that, according to the above definitions, if $\oldUpsilon'$ is a trimming of $\oldUpsilon$, then $\oldDelta_{\oldUpsilon'}$ is empty, $\frakB_{\oldUpsilon'}$ has no redundant components, and no two distinct components of $\frakB_{\oldUpsilon'}$ are similar. Furthermore, since $\oldUpsilon'$ is a union of components of the $\pi_1$s-injective suborbifold $\oldUpsilon$, the suborbifold $\oldUpsilon'$ of $\otheroldLambda$ is itself $\pi_1$-injective.
\EndNumber

\Lemma\label{dot trimming}
Let  $\oldUpsilon$ be a $\pi_1$-injective $2$-suborbifold of be a negative, orientable, finite-type $2$-orbifold without boundary $\otheroldLambda$. Then a $2$-suborbifold $\oldGamma$ is a trimming of $\oldUpsilon$ if and only if the $2$-suborbifold $\dot\oldGamma$ of $\dot\otheroldLambda$ is a trimming of $\dot\oldUpsilon$.
\EndLemma

\Proof
A $2$-suborbifold $\oldXi$ of $\otheroldLambda$ either is discal or is a manifold homeomorphic to $\SSS^1\times[0,\infty)$ if and only if the suborbifold $\dot\oldXi$ of $\dot\oldUpsilon$ either is discal or is a manifold homeomorphic to $\SSS^1\times[0,\infty)$. Hence we have 
$\oldDelta_{\dot\oldUpsilon}^{\dot\otheroldLambda}=
(\oldDelta_\oldUpsilon^\otheroldLambda)\dot\empty$.
Likewise, if $\oldXi$ is a suborbifold of $\otheroldLambda$, then $|\oldXi|$ is a weight-$0$ annulus if and only if $|\dot\oldXi|$ is a weight-$0$ annulus (in which case $\dot\oldXi=\oldXi$). Hence 
$\frakB_{\dot\oldUpsilon}^{\dot\otheroldLambda}=
(\frakB_\oldUpsilon^\otheroldLambda)\dot\empty$; a component $\frakC$ of $\frakB_{\dot\oldUpsilon}^{\dot\otheroldLambda}$ is redundant if and only if $\dot\frakB$ is a redundant component of $\frakB_{\dot\oldUpsilon}^{\dot\otheroldLambda}$; and two components $\frakC$ and $\frakC'$ of $\frakB_\oldUpsilon^\otheroldLambda$ are similar if and only if the components $\dot\frakC$ and $\dot\frakC'$ of $\frakB_{\dot\oldUpsilon}^{\dot\otheroldLambda}$ are similar. It follows that if $\oldGamma$ is a trimming of $\oldUpsilon$ then $\dot\oldGamma$ of $\dot\otheroldLambda$ is a trimming of $\dot\oldUpsilon$. The converse now follows in view of the uniqueness of a trimming up to isotopy.
\EndProof

\Proposition\label{lattice} 
Let $\otheroldLambda$ be a negative, orientable, finite-type $2$-orbifold without boundary.
\begin{enumerate} 
\item The set $\barcaly(\otheroldLambda)$,
ordered by $\preceq$, is a lattice: that is, any subset of
$\barcaly(\otheroldLambda)$ with cardinality $1$ or $2$ has a supremum and
an infimum. 
\item If $\oldUpsilon_1$ and $\oldUpsilon_2$ are
elements of $\Theta(\otheroldLambda)$ that  are in standard position, then $\oldUpsilon_1\cap \oldUpsilon_2$ is $\pi_1$-injective in $\otheroldLambda$, and thus admits a trimming by \ref{trimming}. Furthermore, if
$\oldGamma$ denotes the trimming of 
$\oldUpsilon_1\cap \oldUpsilon_2$,
then 
$\oldGamma\in\Theta(\otheroldLambda)$, and
$[\oldGamma]$ is an infimum for
$\{[\oldUpsilon_1],[\oldUpsilon_2]\}$ in $\barcaly(\otheroldLambda)$.
\item If $\oldUpsilon_1$ and $\oldUpsilon_2$ are
elements of $\Theta(\otheroldLambda)$ such that $\otheroldLambda-\inter\oldUpsilon_1$ and $\otheroldLambda-\inter\oldUpsilon_2$  are in standard position, then  $\otheroldLambda-\inter (\oldUpsilon_1\cup \oldUpsilon_2)$ is $\pi_1$-injective in $\otheroldLambda$, and thus admits a trimming by \ref{trimming}. Furthermore, if
 $\oldTheta$ denotes the trimming of 
 $\otheroldLambda-\inter (\oldUpsilon_1\cup \oldUpsilon_2)$,
% with all discal components
%of $\otheroldLambda-\inter(\oldUpsilon_1\cup \oldUpsilon_2)$, 
then 
$\otheroldLambda-\inter\oldTheta\in\Theta(\otheroldLambda)$, and
$[\otheroldLambda-\inter\oldTheta]$
%, then
 %$[\oldGamma]$ 
is a supremum for $\{[\oldUpsilon_1],[\oldUpsilon_2]\}$ in
$\barcaly(\otheroldLambda)$.
\end{enumerate}
\EndProposition

\Proof
We first prove Assertion (2). 
Note that each component of $\Fr_{\oldUpsilon_1}(\oldUpsilon_1\cap\oldUpsilon_2)$ is either a two-sided arc in  $\oldUpsilon_1$ or a component of $\partial\oldUpsilon_2$. Since $\oldUpsilon_2$ is taut, $\partial\oldUpsilon_2$ is $\pi_1$-injective in $\otheroldLambda$. Thus 
$\Fr_{\oldUpsilon_1}(\oldUpsilon_1\cap\oldUpsilon_2)$ is $\pi_1$-injective in $\oldUpsilon_1$, and hence $\oldUpsilon_1\cap\oldUpsilon_2$ is $\pi_1$-injective in $\oldUpsilon_1$; since $\oldUpsilon_1$ is taut, it follows that $\oldUpsilon_1\cap\oldUpsilon_2$ is $\pi_1$-injective in $\otheroldLambda$.

Now let $\oldGamma$ be a trimming of $\oldUpsilon_1\cap\oldUpsilon_2$.
To prove that
$\oldGamma\in\Theta(\otheroldLambda)$, 
first note that by an observation made in \ref{trimming},  $\oldDelta_{\oldGamma}$ is empty, $\frakB_{\oldGamma}$ has no redundant components, and no two distinct components of $\frakB_{\oldGamma}$ are similar.
In view of the definitions (and the negativity of $\otheroldLambda$, which guarantees that $\otheroldLambda$ is not toric), this means that $\oldGamma$ has no discal components and no manifold components homeomorphic to $\SSS^1\times\RR$, and that no weight-$0$ annulus component of $|\oldUpsilon|$ shares a boundary component with a weight-$0$
annulus component of $|\otheroldLambda-\inter\oldUpsilon|$. The latter fact is Condition (ii) in the definition (\ref{chidef}) of $\Theta(\otheroldLambda)$. Thus to show $\oldGamma\in\Theta(\otheroldLambda)$, it remains only to verify Condition (i) in the definition of $\Theta(\otheroldLambda)$, namely that $\oldGamma$ is taut.

To verify Condition (i) in the definition of tautness (see \ref{praxis}),
note that since the $\oldUpsilon_i$ are closed subsets of $\otheroldLambda$, so is $\oldUpsilon_1\cap\oldUpsilon_2$; and since the $\partial\oldUpsilon_i$ are compact, so is $\partial\oldGamma\subset\partial(\oldUpsilon_1\cap\oldUpsilon_2)\subset (\partial\oldUpsilon_1)\cup(\partial\oldUpsilon_2)$. 
To verify Condition (ii) of the definition of tautness, first note that  $\oldGamma$ is $\pi_1$-injective in $\otheroldLambda$ by \ref{trimming}. Since no component of $\oldGamma$ is discal, $\partial\oldGamma$ is $\pi_1$-injective in $\oldGamma$ and hence in $\otheroldLambda$. This is Condition (ii) of the definition of tautness.

To verify Condition (iii) of the definition of tautness, suppose that some component $\frakJ$ of $\partial\oldGamma$ is contained in a suborbifold 
$\oldXi$ 
of $\otheroldLambda$
% $\overline{\otheroldLambda-\oldUpsilon}$ 
which is a $2$-manifold homeomorphic to
$\SSS^1\times[0,\infty)$. Since we have seen that $\partial\oldGamma$ is $\pi_1$-injective in $\otheroldLambda$, in particular $\frakJ$ is $\pi_1$-injective in $\oldXi$, and hence $\frakJ$ is the boundary of a $2$-manifold $\oldXi'\subset\oldXi\subset\otheroldLambda$ which is homeomorphic to
$\SSS^1\times[0,\infty)$. Now since $\oldGamma$ is a trimming of $\oldUpsilon_1\cap\oldUpsilon_2$, we have $\oldGamma\subset\oldUpsilon_1\cap\oldUpsilon_2$, and in particular $\frakJ\subset\oldUpsilon_1\cap\oldUpsilon_2$. Hence if $i\in\{1,2\}$ is given, we have $\partial\oldXi'=\frakJ\subset\oldUpsilon_i$. This implies that either $\oldXi'\subset\oldUpsilon_i$, or $\oldXi'$ contains a component of $\partial\oldUpsilon_i$. But the latter alternative is ruled out by the tautness of $\oldUpsilon_i\in\Theta(\otheroldLambda)$; indeed, since $\oldXi'$  is homeomorphic to
$\SSS^1\times[0,\infty)$, this alternative would contradict Condition (iii) of the definition of tautness. The remaining possibility is that for $i=1,2$ we have $\oldXi'\subset\oldUpsilon_i$. This means that $\oldXi'\subset\oldUpsilon_1\cap\oldUpsilon_2$. On the other hand, since $\oldGamma$ is a trimming of $\oldUpsilon_1\cap\oldUpsilon_2$, it is in particular a union of components of $\oldUpsilon_1\cap\oldUpsilon_2$; thus the
%if $\oldGamma_0$ denotes the
 component  $\frakJ$ of $\partial\oldGamma$ is a component of $\partial(\oldUpsilon_1\cap\oldUpsilon_2)$. We now have $\oldXi'\subset\oldUpsilon_1\cap\oldUpsilon_2$ and
$\partial\oldXi'=\frakJ\subset\partial(\oldUpsilon_1\cap\oldUpsilon_2)$. Since $\oldXi'$ is connected, this implies that $\oldXi'$ is a component of $\oldUpsilon_1\cap\oldUpsilon_2$. Since $\oldGamma$ is a union of components of $\oldUpsilon_1\cap\oldUpsilon_2$, and $\partial\oldXi'=\frakJ\subset\oldGamma$, it now follows that $\oldXi'$ is a component of $\oldGamma$. But this is impossible, since we have seen that no component of $\oldGamma$ is a $2$-manifold homeomorphic to $\SSS^1\times[0,\infty)$. This contradiction establishes Condition (iii) of the definition of tautness, and completes the proof that $\oldGamma\in\Theta(\otheroldLambda)$.

We now turn to the proof that $[\oldGamma]$ is an infimum for
$\{[\oldUpsilon_1],[\oldUpsilon_2]\}$.
For $i=1,2$ we have $\oldGamma\subset
\oldUpsilon_1\cap\oldUpsilon_2\subset
\oldUpsilon_i$ and
hence $[\oldGamma]\preceq[\oldUpsilon_i]$. Given an
element $Y$  of $\barcaly(\otheroldLambda)$ such that $Y
\preceq[\oldUpsilon_i]$ for $i=1,2$, we need to prove
\Equation\label{late for dinner}
Y\preceq[\oldGamma].
\EndEquation

We begin by proving (\ref{late for dinner}) in the case where $\otheroldLambda$ is a
$2$-manifold. In this case we will set $L=\otheroldLambda$,
$F_i=\oldUpsilon_i$ for $i=1,2$, and $G=\oldGamma$: the Roman letters
are intended to emphasize that $\otheroldLambda$, $\oldGamma$ and the $\oldUpsilon_i$ are
  manifolds. Let
us write $Y=[ F]$, where
$ F\in\Theta( L)$. We may choose $ F$ within
its proper isotopy class so that $\partial F$ is transverse to
$\partial F_i$ for $i=1,2$, and is disjoint from
$\partial F_1\cap\partial F_2$. Furthermore,
among all elements of $\Theta( L)$ that represent the proper isotopy
class $Y$, and whose boundaries are transverse to $\partial F_i$ for
$i=1,2$ and disjoint from $\partial F_1\cap\partial F_2$, we may
choose $ F$ to minimize the quantity
$\card(\partial F\cap (\partial F_1\cup
\partial F_2))$. We claim:
\Claim\label{either one}
For $i=1,2$ we have $\partial F\cap \partial F_i=\emptyset$.
\EndClaim

To prove \ref{either one}, assume that  $\partial F\cap 
\partial F_{i_0}\ne\emptyset$ for some $i_0\in\{1,2\}$. Since $Y\preceq[ F_{i_0}]$,
the $1$-manifold $\partial F$ is isotopic to a $1$-manifold
disjoint from $\partial F_{i_0}$. It then follows from Proposition \ref{beats epstein} that 
$\partial F_{i_0}$ and $\partial F$ have a degenerate
crossing. Among all disks in $ L$ that
are disks of degeneracy for 
$\partial F_{i}$ and $\partial F$ for some $i\in\{1,2\}$, choose one, say
$D$, which is minimal with respect to inclusion. By symmetry we may
assume that $D$ is
a disk of degeneracy for 
$\partial F_{1}$ and $\partial F$. By the definition of
a disk of degeneracy, we may write $\partial D=a_1\cup a$, where $a_1$
and $a$ are arcs contained in $\partial F_1$ and $\partial F$ respectively, and $a\cap a_1=\partial a=\partial a_1$.
Since $ F_2$ is taut, each component of
$\partial F_2\cap D$ is a properly embedded arc in $D$. If
some component of $\partial F_2\cap D$ has both its endpoints
in $a_1$, then $\partial F_1$ and $\partial F_2$
have a degenerate crossing, which is a contradiction since
$ F_1$ and $ F_2$ are in standard position. If
some component $b$ of $\partial F_2\cap D$ has both its endpoints
in $a$, then $a$ is the frontier of a subdisk of $D$ which is a disk
of degeneracy for $\partial F_2$ and $\partial F$,
which contradicts the minimality of $D$. There remains the possibility
that every component of $\partial F_2\cap D$ has one endpoint
in $a_1$ and one in $a$. Then we have
$\card(a\cap\partial F_2)=\card(a_1\cap\partial F_2)$. Hence
the $1$-manifold $C':=(\partial F-a)\cup a_1$, which is
isotopic to $C:=\partial F$, satisfies
$\card(C'\cap\partial F_2)=\card(C\cap\partial F_2)$. The
$1$-manifold $C'$ shares the arc $a_1$ with $\partial F_1$,
but can be modified by a small isotopy to obtain a $1$-manifold $C''$, transverse to $\partial F_i$ for
$i=1,2$ and disjoint from $\partial F_1\cap\partial F_2$,
such that
$\card(C''\cap\partial F_1)=\card(C\cap\partial F_1)-2$
and
$\card(C''\cap\partial F_2)=\card(C'\cap\partial F_2)=\card(C\cap\partial F_2)$. If
$ F''$ is a submanifold of $ L$ that is properly isotopic to
$ F$ and has boundary $C''$, then we have $\card(\partial F'\cap (\partial F_1\cup
\partial F_2))=\card(\partial F\cap (\partial F_1\cup
\partial F_2))-2$, which contradicts the minimality property of
$ F$. Thus \ref{either one} is proved.

Since the trimming $G$ is by definition a union of components of $F_1\cap F_2$, we have $\partial G\subset\partial(F_1\cap F_2)\subset(\partial F_1)\cup(\partial F_2)$, so that by \ref{either one} we have 
\Equation\label{ryan in common}
(\partial F)\cap(\partial G)=\emptyset.
\EndEquation

Let $ F_-$ denote the union of all negative components of $ F$. We claim:
\Claim\label{new alabama}
Suppose that $i\in\{1,2\}$ is given, and that $ Z$ is a component 
of $ F_-$. Then there exists a compact submanifold $ X$ of $ Z$ such that (a) the closure of each component of $ Z-  X$ is an annulus meeting $\partial Z$ in a single component of $\partial Z$, and (b) $ X\subset F_i$.
\EndClaim

To prove \ref{new alabama}, we may assume by symmetry that $i=1$. Let $ L_0$ denote the component of
$ L$ containing $ Z$. Since
$[ F]=Y\preceq[ F_1]$, the inclusion map
$ Z\to L_0$ is isotopic in $ L_0$ to an embedding
$h: Z\to L_0$ with
$h( Z)\subset F_1$. Let $ Z_1\subset L_0$ denote the component
of $ F_1$ containing $h( Z)$. Let $z$ be a base point in
$ Z$, set $u=h(z)\in Z_1$, let
$j: Z_1\to L_0$ denote the inclusion map, and let $(\tL_0,\tu)$ denote the
based covering space of $( L_0,u)$ determined by the subgroup
$j_\sharp(\pi_1( Z_1,u))$ of $\pi_1( L_0,u)$. Then $j$ admits a
based lift 
$\tj:( Z_1,u)\to(\tL_0,\tu)$. Thus
$\tZ_1:=\tj( Z_1)$ is a subsurface of $\tL_0$ which is
mapped homeomorphically onto $ Z_1$ by the covering
projection $p:\tL_0\to L_0$. The construction implies
that the inclusion homomorphism $\tZ_1\to\tL_0$ is
surjective; since the tautness of $ Z_1$ implies that
$\partial\tZ_1$ is $\pi_1$-injective in $\tL$, it follows that each component of
$\tL_0-\inter\tZ_1$ is a half-open annulus.

The map $h:( Z,z)\to( L_0,u)$
admits a based lift to $(\tL_0,\tu)$ since
$h( Z)\subset Z_1$; and
since the inclusion map
$ Z\to L_0$ is (freely) homotopic to
$h$,
% $j$ is homotopic to $h$, 
it follows from the covering homotopy
property of covering spaces that this inclusion map admits a lift to
$ \tL_0$. Hence there is a 
%: Z\to\tL_0$. Since $h$ is an embedding, the
 subsurface $\tZ$ of $\tL_0$ which is mapped homeomorphically
 onto $ Z$ by $p$. Note that since $\partial Z\cap\partial Z_1=\emptyset$ by \ref{either one},
we have $\partial\tZ\cap \partial\tZ_1=\emptyset$.

Set $\tX=\tZ\cap\tZ_1$. Consider an arbitrary component
$A$ of the closure of $\tZ-\tX$. Since  $\partial\tZ\cap
\partial\tZ_1=\emptyset$, the set $A$ is a $2$-manifold, and
each component of $\partial A$ is either a component of
$\partial\tZ$, or a component
of $\partial\tZ_1$ contained in $\inter\tZ$. Since $ Z$ and
$ Z_1$ are taut in $ L$, the $1$-manifold
$\partial A$ is $\pi_1$-injective in $\tL_0$. On the other hand, we have
$A\subset\tL_0-\inter\tZ_1$, and we have observed that each
component of $\tL_0-\inter\tZ_1$ is a half-open
annulus. Thus $A$ is a
compact, connected subsurface of a half-open annulus component $B$  of
$\tL_0-\inter\tZ_1$, and $\partial A$ is
$\pi_1$-injective in $B$; hence $A$
must be a closed annulus. We cannot have $\partial A\subset\partial B$,
since $\partial B$ is a single simple closed curve. If both components of $\partial A$ are
components of $\partial\tZ$, then $A$ is a component of $\tZ$;
this is impossible because $A$ is an annulus, and $\tZ$, being homeomorphic to the
component $Z$ of $ F_-$, is negative. Hence $\partial A$ must have
one component contained in $\partial\tZ$ and one contained in
$\tZ_1$. 

This shows that the closure of each component of $\tZ -\tX$ is
an annulus meeting $\partial\tZ$ in a single component of
$\partial\tZ$. Now set $ X=p(\tX)\subset Z$. Since $p$
maps $\tZ$ homeomorphically onto $ Z$, each component of
$ Z - X$ is the homeomorphic image under $p$ of a component of
$\tZ -\tX$, and is therefore
an annulus meeting $\partial Z$ in a single component of
$\partial Z$. Since the definition of $\tX$ implies that
$\tX\subset\tZ_1$, we have $ X\subset
p(\tZ_1)= Z_1\subset F_1$. This completes the proof
of \ref{new alabama}.

It follows immediately from \ref{new alabama} that for $i=1,2$, there
exists a compact submanifold $ W_i$ of $ F_-$ such that
(a) the closure of each component of $ F_-- W_i$ is an
annulus meeting $\partial F_-$ in a single component of
$\partial F_-$, and (b)
$ W_i\subset F_i$. After removing an open boundary collar from each
$W_i$, we may assume that $ W_i\subset \inter F_-$ for
$i=1,2$. For each component $C$ of
$\partial F_-$, and each index $i\in\{1,2\}$, 
%let us define a
%set $A_i(C)\subset F_-$ as follows: if $C\subset W_i$, we
%set $A_i(C)=C$; otherwise, we define $A_i(C)$ to be 
the closure of the
component of $ F_-- W_i$ that contains $C$ is an annulus which we will
denote by $A_i(C)$. We have 
%is either equal to $C$, or is an annulus
% meeting 
$A_i(C)\cap \partial F_-=C$.
%precisely in $C$. 
For $i=1,2$ we have
\Equation\label{feldfusser}
 F_--\bigcup_{C\in\calc(\partial F_-)}A_i(C)=\inter W_i.
\EndEquation
 We
claim:
\Claim\label{one's in the other}
If $C$ and $C'$ are components of $\partial F_-$, not necessarily distinct, we have either (i) $A_1(C)\cap A_2(C')=\emptyset$, or (ii) $C=C'$, and either
$A_1(C)\subset A_2(C)$ or $A_2(C)\subset A_1(C)$.
\EndClaim

To prove \ref{one's in the other},
% first note that the assertion is
%trivial if $C\subset W_1$; for in
%that case we have $A_1(C)=C$, which implies $A_1(C)\subset A_2(C)=A_2(C')$ if $C=C'$, and $A_1(C)\cap A_2(C')=\emptyset$ otherwise. Likewise, it is trivial
%if $C'\subset W_2$. 
%, or if $C'$ is contained in either $ W_1$ or $ W_2$. 
%Now suppose that $C\not\subset W_1$ and $C'\not\subset W_2$,
%nor $C'$ is contained in either of the
%$ W_i$, 
%so that $A_1(C)$ and $A_2(C')$ are
% annuli; 
let $C_1$ denote the boundary component of $A_1(C)$ distinct from $C$, and  let $C_2'$ denote the boundary component of $A_2(C')$ distinct from $C'$. Then $C_1,C_2'\subset\inter F_-$
are components of $\partial F_1$ and $\partial F_2$ respectively. The curves $C_1$ and
$C_2'$ are respectively isotopic to $C$ and $C'$, which are either equal or disjoint. %particular, $C_1$ and $C_2$ are isotopic to and hence to each other; ]
Furthermore, since
$ F_1$ and $ F_2$ are in standard position, $C_1$
and $C_2$ have no degenerate crossings. Since $C_1$ and $C_2$ are isotopic to disjoint curves, it then follows from Proposition \ref{beats epstein} that $C_1$ and $C_2$ are disjoint.
Thus $A_1(C)$ and
$A_2(C')$ meet $\partial F_-$ in $C$ and $C'$ respectively, and their frontiers
are disjoint. If $C=C'$, it follows that one of the annuli  $A_1(C)$ and
$A_2(C')$ is contained in the other. If $C\ne C'$ it follows that either $A_1(C)\cap
A_2(C')=\emptyset$, or  $A_1(C)\cup A_2(C')$ is an annulus component of $ F_-$. Since $ F_-$ is negative by construction, it has no annulus components, and
\ref{one's in the other} is proved.

Set $N=\bigcup_{C\in\calc(\partial F_-)}(A_1(C)\cup A_2(C))$. It follows from \ref{one's in the other} that every component of $N$ has the form $A_i(C)$ for some $i\in\{1,2\}$ and some component $C$ of $\partial F_-$. In particular, each component  of $N$ is an annulus contained in $ F_-$, meeting $\partial F_-$ in a component of $\partial  F_-$. it follows that $ F_-':=\overline{ F_--N}$ is a subsurface of $ L$ properly isotopic to $ F_-$. On the other hand, it follows from (\ref{feldfusser}) and the definitions of $N$ and $ F_-'$ that
$$
 F_-'=\bigcap_{C\in\calc(\partial F_-)}\bigg(\overline{ F_--A_1(C)}\cap \overline{ F_--A_2(C)}\bigg)= W_1\cap W_2.$$
Since $ W_i\subset \inter F_i$ for $i=1,2$, this gives
%\Equation\label{snuff stuff}
$ F_-'\subset\inter( F_1\cap F_2)$.
%$\EndEquation
Since $ F_-'$ is properly isotopic to $ F_-$, it is taut and negative, so that the component of $ F_1\cap F_2$ containing any component of $F_-'$ is negative. Since $G$ is a trimming of $ F_1\cap F_2$ it contained every negative component of $ F_1\cap F_2$, Hence we have
%, and is therefore contained in the interior of the union of all non-discal components of $ F_1\cap F_2$; that is,
\Equation\label{murgadoodle}
 F_-'\subset\inter G.
\EndEquation

Note that we can choose the ambient proper isotopy from $ F_-$ to $ F_-'$ to be constant outside an arbitrarily small neighborhood of $ F_-$. In particular,  if $ \Fann$ denotes the union of all annulus components of $ F$, then the isotopy is constant on $\Fann$. Now since $\otheroldLambda$ is negative and $F$ is taut, each component of $F$ has non-positive Euler characteristic. Since $\otheroldLambda$ is negative by hypothesis, the submanifold $F$ has no boundaryless components. Furthermore, the tautness of $F$ implies that none of its components is a manifold homeomorphic to $\SSS^1\times\RR$. Hence $F=F_-\cup \Fann$. It now follows that:
\Claim\label{no move}
We have $ \Fann\cap  F_-'=\emptyset$, and $ F$ is properly isotopic to $F':= \Fann\cup F_-'.$
\EndClaim

We now claim:

\Claim\label{polychromatic fowl}
Let  $E$ be a  component of $\partial\Fann$, and let an index $i\in\{1,2\}$ be given. Suppose that $E\not\subset\inter F_i$. Then  there is an annulus $K$ having $E$ as one boundary component, and whose other boundary component $C:=(\partial K)-E$ satisfies $K\cap F_i=C$. In particular $C$ is a component of $\partial F_i$.

Furthermore, if $K$ is an annulus having the properties stated above, and if we define $C$ as above and set $i'=3-i$ (so that $i'$ is the index in $\{1,2\}$ distinct from $i$), then:

\begin{enumerate}[(a)]
%\EndClaim
%\redcomment{Prove this.}
%For each index $i\in\{1,2\}$, and each component  $E$ of $\Fann$ such that $E\not\subset\inter F_i$, we choose an annulus $K_E^{(i)}$ having the properties stated in Alternative (ii) of \ref{polychromatic fowl}. In this case we define $C_E^{(i)}$ and $B_E^{(i)}$ as in \ref{polychromatic fowl}. \redcomment{Explain why} $K_E^{(i)}$ may be chosen in such a way that 
%\Equation
%E\subset K_E^{(i)}.
%\EndEquation
%Next, we claim:
%\Claim\label{mustard}
%If an index $i\in\{1,2\}$ is given,  if $E$ is a component  of $\Fann$ such that $E\not\subset\inter F_i$, and if we set  $K=K_E^{(i)}$, then 
\item $K\cap \partial F_{i'}$ is a (possibly empty) union of pairwise disjoint simple closed curves contained in $\inter K$ and homotopically non-trivial in $K$; and
\item 
%Let an index $i\in\{1,2\}$ be given, and suppose that $E$ is a component  of $\Fann$ such that $E\not\subset\inter F_i$.   Then 
if  $K\cap F_{i'}\ne\emptyset$, then $C\subset F_{i'}$. 
\end{enumerate}
\EndClaim

To prove \ref{polychromatic fowl}, first note that by \ref{either one}, we have $E\cap \partial F_i\subset\partial F\cap\partial F_i=\emptyset$. Since $E$ is connected and $E\not\subset\inter F_i$, it follows that $E\cap F_i=\emptyset$. On the other hand, since $E\subset\Fann\subset F$, and since $[F]=Y\preceq[F_i]$, the curve $E$ is isotopic to a simple closed curve $E'\subset\inter F_i$, which we may take to meet $\partial F_i$ transversally. We have $E\cap E'\subset E\cap F_i=\emptyset$. It then follows from \cite[Lemma 2.4]{epstein} that $E$ and $E'$ are the boundary components of some annulus $A\subset L$. Since $E\cap F_i=\emptyset$ and $E'\subset\inter F_i$, we have $A\cap\partial F_i\ne\emptyset$, and $(\partial A)\cap(\partial F_i)=\emptyset$. Hence each component of $A\cap\partial F_i$ is a simple closed curve in $\inter A$, which is homotopically non-trivial in $L$---and hence in $A$---since $F_i$ is taut. Thus if $K$ denotes the closure of the component of $A\setminus\partial F_i$ containing $E$, then $K$ is an annulus having $E$ as one boundary curve; the other boundary curve $C$ of $K$ is a simple closed curve component of $A\cap \partial F_i$, and is therefore a component of $\partial F_i$. The definition of $K$ implies that $K\cap F_i=C$, and the first assertion is established.

To prove Assertion (a), first note that the boundary components of $K$ are $E$ and $C$. By \ref{either one}, we have $E\cap \partial F_{i'}\subset\partial F\cap\partial F_{i'}=\emptyset$.
The curve $C\subset\partial F_i$ meets $\partial F_{i'}$ transversally since $F_1$ and $F_2$ are in standard position. Hence each component of $K\cap \partial F_{i'}$ is either a simple closed curve in $\inter K$, or a properly embedded arc in $A$ whose boundary lies in $C$. Such an arc, say $a$, must be the frontier of a disk in $A$ whose boundary is the union of $a$ with an arc $a'\subset\partial A$. Hence the existence of such an arc implies that $\partial F_{i'}$ and $C\subset\partial F_i$ have a degenerate crossing, a contradiction since $F_1$ and $F_2$ are in standard position. Thus (a) is proved.

To prove Assertion (b), note that if $K\cap F_{i'}\ne\emptyset$, then either $K\subset F_{i'}$ or $K\cap\partial F_{i'}\ne\emptyset$. If $K\subset F_{i'}$ then in particular $C\subset F_{i'}$. Now suppose that $K\cap\partial F_{i'}\ne\emptyset$. By Assertion (a), $K\cap\partial F_{i'}$ is a disjoint union of homotopically non-trivial simple closed curves in $\inter K$. Hence the closure of the component of $K\setminus\partial F_{i'}$ containing $C$ is an annulus $K'$ whose boundary curves are $C$ and $E'$, where $E'$ is some component of $K\cap\partial F_{i'}$. We have $K'\cap\partial F_{i'}=E'\subset\partial K'$. Hence either $K'\subset F_{i'}$ or $K'\cap F_{i'}= E'$. If $K'\subset F_{i'}$, then in particular $C\subset F_{i'}$. If $K'\cap F_{i'}= E'$, then since $K\supset K'$ has its interior disjoint from $F_i$, the annulus $K'$ is contained in 
$\overline{L-(F_1\cup F_2)}$, and has one boundary component (namely $C$) contained in $F_i$, and one (namely $E'$) contained in $F_{i'}$. This contradicts Condition (2) of the definition of standard position (see \ref{nondeg cross}). Thus Assertion (b) is established, and the proof of \ref{polychromatic fowl} is complete.

We also claim:

\Claim\label{catch-up}
If $E$ is a component of $\partial\Fann$, then either (A) $E\subset\inter F_1\cap\inter F_2$, (B) $E$ is isotopic to a component of $\partial F_1$ contained in $\inter F_2$, or (C) $E$ is isotopic to a component of $\partial F_2$ contained in $\inter F_1$.
\EndClaim

To prove \ref{catch-up}, suppose that Alternative (A) does not hold. Then either $E\not\subset\inter F_1$ or $E\not\subset\inter F_2$. 
%Since the disjunction of the remaining alternatives is symmetric in $F_1$ and $F_2$ss, we may assume that $E\not\subset\inter F_1$.
Let us fix an index $i\in\{1,2\}$ such that $E\not\subset\inter F_i$. Then by \ref{polychromatic fowl},  there is an annulus $K$ having $E$ as one boundary component, and whose other boundary component $C:=(\partial K)-E$ satisfies $K\cap F_i=C$. In particular $C$ is a component of $\partial F_i$.

We distinguish two cases. Consider first the case in which $K\cap F_{i'}\ne\emptyset$. In this case, according to Assertion (b) of \ref{polychromatic fowl}, we have $C\subset F_{i'}$. Since $C$ is a component of $\partial F_i$, and since $\partial F_i$ and $\partial F_{i'}$ are in standard position and in particular intersect transversally, we must have $C\subset\inter F_{i'}$. Thus one of the alternatives (B) and (C) holds in this case.

Now consider the case in which $K\cap F_{i'}=\emptyset$. In particular we then have $E\not\subset\inter F_{i'}$. We may therefore apply \ref{polychromatic fowl},  with the roles of $i$ and $i'$ reversed, to obtain an annulus $K'$ having $E$ as one boundary component, and whose other boundary component $C':=(\partial K')-E$ satisfies $K'\cap F_{i'}=C'$. In particular $C'$ is a component of $\partial F_{i'}$.
We now distinguish two subcases. In the subcase where $K'\cap F_i\ne\emptyset$, we can argue exactly as in the first case, with the roles of $i$ and $i'$ reversed, and with $K'$ and $C'$ playing the roles of $K$ and $C$, to show that one of the alternatives (B) and (C) holds.

There remains the subcase in which $K'\cap F_{i}=\emptyset$. As we are in the second case, we also have $K\cap F_{i'}=\emptyset$. Since $K$ is disjoint from $F_{i'}$, and $C'$ is contained in $F_{i'}$, we have $K\cap C'=\emptyset$. The same argument shows that $K'\cap C=\emptyset$. Now since the annuli $K$ and $K'$ share the boundary curve $E$, and since $K$ is disjoint from $C'=(\partial K')-E$ and $K'$ is disjoint from $C=(\partial K)-E$, we must have $K\cap K'=E$. Hence $Q:=K\cup K'$ is an annulus with boundary curves $C$ and $C'$. Since $K\cap F_i=C$ and $K'\cap F_i=\emptyset$, we have $ Q\cap F_i=C$; similarly $ Q\cap F_{i'}=C'$. Thus the annulus $Q$ is contained in 
$\overline{L-(F_1\cup F_2)}$, and has one boundary component (namely $C$) contained in $F_i$, and one (namely $C'$) contained in $F_{i'}$. This contradicts Condition (2) of the definition of standard position. Thus this subcase cannot occur, and the proof of \ref{catch-up} is complete.

Next, we claim:

\Claim\label{hardly hurricanes}
Let $E$ be a component of $\partial F$. If $E$ is isotopic to a curve contained in a component of $F_1\cap F_2$ which is not a component of $G$, then $E$ is isotopic to a component of $\partial G$.
\EndClaim

To prove \ref{hardly hurricanes}, let $S$ be a curve  isotopic to $E$, and contained in a component $W$  of $F_1\cap F_2$ which is not a component of $G$. Since $G$ is a trimming of $F_1\cap F_2$, either $W$ is a disk, or $W$ is homeomorphic to $\SSS^1\times\RR$, or $W$ is a redundant annulus in $F_1\cap F_2$, or $W$ is an annulus similar to an annulus component of $G$ (see \ref{trimming}). If $W$ is a disk, then  $E$ is homotopically trivial in $L$, a contradiction to the tautness of $F$. If $W$ is homeomorphic to $\SSS^1\times\RR$, then $E$ is contained in a submanifold of $L$ homeomorphic to $\SSS^1\times\RR$, and we again have a contradiction to tautness. If $W$ is a redundant annulus in $F_1\cap F_2$, then the homotopically non-trivial curve $S\subset W$ is isotopic in $L$ to a boundary component of
$(F_1\cap F_2)-(\frakB_{F_1\cap F_2}\cup\oldDelta_{F_1\cap F_2})$, which is in particular a component of $\partial G$. If $W$ is an annulus similar to an annulus component $W_0$ of $G$, then $S$ is isotopic to a component of $\partial W_0$, which is in particular a component of $\partial G$. This completes the proof of \ref{hardly hurricanes}.

We now claim:

\Claim\label{hereford}
For every component $X$ of $\Fann$,  either (a) some component of $\partial X$ is contained in $\inter G$, or (b)  a core curve of the annulus $X$ is isotopic to some component of $\partial G$. Furthermore, if (b) holds and (a) does not, there is an annulus $U\subset L$  such that $U\cap G$ is one component of $\partial U$, and the other component of $\partial U$ is contained in $\partial X$.
\EndClaim

To prove \ref{hereford},  suppose that $X$ is a component of $\Fann$, and fix a component $E$ of $\partial X$. Then $E$ satisfies one of the alternatives (A)---(C) of \ref{catch-up}. First suppose that (A) holds. Let $W$ denote the component of $F_1\cap F_2$ whose interior contains $E$. If $W$ is a component of $G$, then (a) holds. If $W$ is not a component of $G$, then since $E$ is a component of $\partial F$, and is of course isotopic to itself, we may apply \ref{hardly hurricanes} to deduce that $E$ is isotopic to a component of $\partial G$; it follows that (b) holds. 

If one of the alternatives (B), (C) of \ref{catch-up} holds, then in particular $E$ is isotopic to a boundary component of some component $V$ of $F_1\cap F_2$. If $V$ is a component of $G$ then (b) holds. If $V$ is not a component of $G$, then since $E$ is a component of $\partial F$, it follows from \ref{hardly hurricanes} that $E$ is isotopic to a component of $\partial G$, so that again (b) holds. This establishes the first assertion of \ref{hereford}.

To prove the second assertion of \ref{hereford}, suppose that (b) holds and (a) does not. By \ref{ryan in common} we have $(\partial F)\cap(\partial G)=\emptyset$.  In particular, if $E$ is a boundary curve of $X$, we have $E\cap\partial G=\emptyset$. Since (a) does not hold, we have $E\not\subset \inter G$, and hence $E\cap G=\emptyset$. But (b) implies that $E$ is isotopic to a component $E'$ of $\partial G$, so that by \cite[Lemma 2.4]{epstein} there is an annulus $U_0\subset L$ with $\partial U_0=E\cup E'$. Since $\partial G$ is a $1$-manifold disjoint from $E$, the curve $E'$ is a component of $U_0\cap\partial G$, and any other components of $U_0\cap\partial G$ is a simple closed curve contained in $\inter U_0$, and is homotopically non-trivial since $G$ is taut. Hence the closure  of the component of $U_0\setminus\partial G$ containing $E$ is an annulus $U$ having $E$ as one boundary component, and its other boundary component $E''$ is contained in $\partial G$. The definition of $U$ guarantees that $U\cap\partial G=E'$, and since 
$E\cap G=\emptyset$ it follows that
$U\cap G=E'$. This completes the proof of \ref{hereford}.

Now write $\Fann=\Fannin\discup\Fannout$, where $\Fannin$ is the union of all components of $\Fann$ that have at least one boundary component contained in $\inter G$, and $\Fannout$ is the union of all components of $\Fann$ that have no such boundary component. According to \ref{hereford}, for each component $X$ of $\Fannout$ we may select  an annulus $U_X\subset L$  such that $U_X\cap G$ is one component of $\partial U_X$, which will be denoted $C_X$, and the other component of $\partial U_X$, which we will denote $C'_X$, is contained in $\partial X$. We claim:

\Claim\label{totalitously}
If $X$ and $X'$ are distinct components of $\Fannout$, we have $U_X\cap U_{X'}=\emptyset$, Furthermore, if $X$ is a component of $\Fannout$ and $X'$ is any component of $\Fann$, we have $U_X\cap X'=\emptyset$. Finally, for each component $X$ of $\Fannout$, 
%the annuli $X$ and $U_{X'}$ share a boundary com 
we have either $U_X\cap X=C'_X$, or $U_X\subset\inter X$ or $X\subset\inter U_X$.
\EndClaim

To prove \ref{totalitously}, define an index set $J$ by $J=(\calc(\Fannout)\times\{1\})\times (\calc(\Fann)\times\{2\})\subset \calc(\Fann)\times\{1,2\})$. Define an annulus $A_j\subset L$ for each index $j\in J$ by setting $A_{X,1}=U_X$ for each $X\in\calc(\Fannout)$, and $A_{X,2}=X$ for each $X\in\calc(\Fann)$.
Then \ref{totalitously} is equivalent to the assertion that for any two distinct indices $j,j'\in J$ we have $A_j\cap A_{j'}=\emptyset$, unless $j$ and $j'$ are equal, in some order, to $(X,1)$ and $(X,2)$ for some $X\in\calc(\Fannout)$; and that in the latter case we have either $A_j\cap A_{j'}=C'_X$, or $A_j\subset\inter A_{j'}$, or $A_{j'}\subset\inter A_j$.

As a preliminary to proving this, we observe that for any index $j=(X,i)\in J$, the annulus $A_j$ is isotopic to $X$. This is trivial if $i=2$ since we then have $A_j=X$. If $i=1$, we have $X\in\calc(\Fannout)$ and $A_j=U_X$; since the annuli $U_X$ and $X$ share the boundary curve $C'_X$, they are isotopic.

Now note that for any index $j$, each boundary component of $A_j$ is a component of either $\partial\Fann$ or $\partial G$. By \ref{ryan in common} we have $(\partial F)\cap(\partial G)=\emptyset$. Hence if $j,j'\in J$ are distinct indices,  $(\partial A_j)\cap(\partial A_{j'})$ is a union of common components of $A_j$ and $A_{j'}$. If $A_j\cap A_{j'}\ne\emptyset$, it follows that one of the annuli $A_j$, $A_{j'}$ contains a boundary component of the other. Since tautness implies that the boundary curves of $A_j$ and $A_j'$ are $\pi_1$-injective in $L$, this implies that $A_j$ and $A_{j'}$ are isotopic. But by the preliminary observation above, if we write $j=(X,i)$ and $j'=(X',i')$, then $A_j$ and $A_{j'}$ are respectively isotopic to $X$ and $X'$. Hence $X$ and $X'$ are isotopic. If $X\ne X'$ this contradicts  Lemma \ref{another stupidity} since $X$ and $X'$ are annulus components of $F\in\Theta(L)$. 
 Thus we can have $A_j\cap A_{j'}=\emptyset$ only if $X=X'$. In the latter case, $A_j$ and $A_{j'}$ are equal to $X$ and $U_X$ in some order, and therefore share the boundary curve $C'_X$. But we have observed that $(\partial A_j)\cap(\partial A_{j'})$ is a union of common boundary components of $A_j$ and $A_{j'}$. Hence the boundary curves $(\partial A_j)-C'_X$ and $(\partial A_{j'})-C'_X$ are either equal or disjoint. Since $L$ is negative and therefore not a torus, we must have either $A_j\cap A_{j'}=C'_X$, or $A_j\subset\inter A_{j'}$, or $A_{j'}\subset\inter A_j$. This proves \ref{totalitously}.

Now for each component $X$ of $\Fannout$, it follows from \ref{totalitously} that  $M_X:=X\cup U_X$ is an annulus having at least one boundary component in common with $X$, and that $M_{X}\cap M_{X'}=\emptyset$ for any distinct components $X$ and $X'$ of $\Fannout$. 
Hence there is a proper isotopy $(h_t)_{0\le t\le1}$ of $L$ (with $h_0$ equal to the identity) such that $h_t(X)=M_X$ for each component $X$ of $\Fannout$, and we may take $(h_t)$ to be constant outside an arbitrarily small neighborhood of $\bigcup_{X\in\calc(\Fannout)}M_X$. In particular, since \ref{totalitously} also gives that $M_X\cap{X'}=\emptyset$ for any component $X$ of $\Fannout$ and any component $X'$ of $\Fannin$, we may take $(h_t)$ to be constant on $\Fannin$. Note also that if $X$ is a component of $\Fann$, we have $X\cap F_-'
\subset
\Fann\cap  F_-'=\emptyset$ by \ref{no move}; and since
$F_-'\subset\inter G$ by
\ref{murgadoodle}, and $U_X\cap G=C_X\subset\partial G$, we have
$U_X\cap F_-'=\emptyset$. Thus $M_X\cap F_-'=\emptyset$. We may therefore take $(h_t)$ to be constant on $F_-'$. It now follows that $F':= \Fann\cup F_-'$, which by \ref{no move} is a disjoint union, is carried by $h_1$ onto $F''=(\bigcup_{X\in\calc(\Fannout)}M_X)\cup\Fannin\cup F_-'$, and that the latter union is also disjoint. Since $F$ is properly isotopic to $F'$ by \ref{murgadoodle}, $F$ is also properly isotopic to $F''$.

We may write $F''=\Fann''\discup F_-'$, where $\Fann''=\bigg(\coprod_{X\in\calc(\Fannout)}M_X\bigg)\discup\Fannin$.
Each annulus of the form $M_X$, where $X$ is a component of $\Fannout$, has $C_X\subset\partial G$ as either a boundary component or a core curve; and by definition, each component of $\Fannin$ is an annulus having one boundary curve in $\partial G$. Thus each component of $\Fann$ has a boundary component or a core curve contained in $G$. Hence there is a proper isotopy of $L$, constant outside an arbitrarily small neighborhood of $\Fann''$, carrying $\Fann''$ into $\inter G$. In particular we may take this isotopy to be constant on $F_-'$. It then carries $F''$ onto a surface $F'''\subset\inter G$. Since $F$ is properly isotopic to $F''$, it is also properly isotopic to $F'''$, and hence $Y=[F]\preceq G$. This
completes the proof of (\ref{late for dinner}) in the case where $\otheroldLambda$ is a $2$-manifold.

%\oldUpsilon Fann W \frak \old \oldXi \frakZ X
%\otheroldLambda\toldLambda V non-trivial

%\redcomment{I have checked the proof above locally, but I'm not really comfortable about why all the pieces are needed.}

To prove  (\ref{late for dinner}) in the general case, write $Y=[\oldUpsilon]$, where
$\oldUpsilon\in\Theta(\otheroldLambda)$. Since $Y
\preceq[\oldUpsilon_i]$ for $i=1,2$, there are elements $\oldXi_i$ of $\Theta(\otheroldLambda)$ for $i=1,2$ such that $\oldXi_i$ is properly isotopic to $\oldUpsilon$ in $\otheroldLambda$, and $\oldXi_i\subset\oldUpsilon_i$. By Lemma \ref{more geventlach}, $\dot\oldUpsilon$, $\dot\oldXi_1$ and $\dot\oldXi_2$ belong to 
$\Theta(\dot\otheroldLambda)$. For $i=1,2$, since $\oldXi_i\subset\oldUpsilon_i$, we have $\dot\oldXi_i\subset\dot\oldUpsilon_i$. But it follows from Lemma \ref{historic fact} that  $\dot\oldXi_i$ is properly isotopic to $\dot\oldUpsilon$ in $\dot\otheroldLambda$; hence if we set $\dot Y=[\dot\oldUpsilon]\in\barcaly(\dot\otheroldLambda)$, we have $\dot Y\preceq[\dot\oldUpsilon_i]$ for $i=1,2$. Furthermore, since $\oldUpsilon_1$ and $\oldUpsilon_2$ are in standard position, it follows from \ref{dot standard} that $\dot\oldUpsilon_1$ and $\dot\oldUpsilon_2$ are in standard position. Since $\oldGamma$ is a trimming of 
$\oldUpsilon_1\cap\oldUpsilon_2$, it follows from Lemma \ref{dot trimming} that $\dot\oldGamma$ is a trimming of $
(\oldUpsilon_1\cap\oldUpsilon_2)\dot\empty=
\dot\oldUpsilon_1\cap\dot\oldUpsilon_2$.

 Since  (\ref{late for dinner}) has already been proved in the manifold case, it follows that $\dot Y\preceq[\dot \oldGamma]$. By definition this means that $\dot\oldUpsilon$ is properly isotopic in $\dot\otheroldLambda$ to a taut submanifold $V$ of $\dot\oldGamma$. It now follows from Lemma \ref{historic fact} that $V=\dot\frakP$ for some taut suborbifold $\frakP$ of $\oldGamma$, and that $\frakP$ is properly isotopic to $\oldUpsilon$. Since $V\subset\dot\oldGamma$, we have $\frakP\subset\oldGamma$. This implies that $Y\preceq[\oldGamma]$, as asserted by  (\ref{late for dinner}). Thus (\ref{late for dinner}) is established in complete generality, and Assertion (2) is thereby proved.

To prove Assertion (3), set $\oldUpsilon_i'=\otheroldLambda-\inter\oldUpsilon_i$ for
$i=1,2$, so that $\oldUpsilon_1'$ and $\oldUpsilon_2'$ are in standard
position. By \ref{involute} we have
$\oldUpsilon_1',\oldUpsilon_2'\in\Theta(\otheroldLambda)$. According to Assertion (2), $\otheroldLambda-\inter (\oldUpsilon_1\cup \oldUpsilon_2)=
\oldUpsilon_1'\cap \oldUpsilon_2'$ is $\pi_1$-injective in $\otheroldLambda$, and if
 $\oldTheta$ is a trimming of 
$\otheroldLambda-\inter (\oldUpsilon_1\cup \oldUpsilon_2)=
\oldUpsilon_1'\cap \oldUpsilon_2'$,
then
$\oldTheta\in\Theta(\otheroldLambda)$ and 
$[\oldTheta]$ is an infimum for
$\{[\oldUpsilon_1'],[\oldUpsilon_2']\}$. Since $\oldTheta\in\Theta(\otheroldLambda)$, we have $[\otheroldLambda-\inter\oldTheta]\in\Theta(\otheroldLambda)$ by \ref{involute}. Furthermore, since, as observed in
\ref{involute}, the involution of
$\barcaly(\otheroldLambda)$ defined by 
$[\oldUpsilon]\mapsto[\otheroldLambda-\inter\oldUpsilon]$ is order-reversing with respect to the partial
order $\preceq$, it now follows that $[\otheroldLambda-\inter\oldTheta]$ is a supremum for
%$[\oldTheta]$ is an infimum for
$\{[\otheroldLambda-\inter\oldUpsilon_1'],[\otheroldLambda-\inter\oldUpsilon_2']\}=
\{[\oldUpsilon_1],[\oldUpsilon_2]\}$,
% But the definitions of
%$\oldTheta$ and $\oldTheta$ imply that
%$\oldTheta=\otheroldLambda-\inter\oldTheta$, 
and Assertion (3) is proved.

We now turn to the proof of Assertion (1). If elements $Y_1$ and $Y_2$ of $\barcaly(\otheroldLambda)$ are
given, then by \ref{nondeg cross} we may write $Y_i=[\oldUpsilon_i]$
for $i=1,2$, where 
%$\oldUpsilon_i\in\Theta_-(\otheroldLambda)$, and
$\oldUpsilon_1$ and $\oldUpsilon_2$ are elements of $\Theta(\otheroldLambda)$ that are in standard position.
%have no degenerate
%crossings. 
It follows from Assertion (2), which has already been proved, that if
$\oldGamma$ is
%and
%$\oldTheta$ are
 defined as in that assertion, then $\oldGamma\in\Theta(\otheroldLambda)$,
 and  $[\oldGamma]\in\barcaly(\otheroldLambda)$ is an infimum for
$\{Y_1,Y_2\}$. 
It also follows from
\ref{nondeg cross}, with \ref{involute}, that we may write $Y_i=[\oldUpsilon_i']$
for $i=1,2$, where $\oldUpsilon_1',\oldUpsilon_2'$ are elements of
%$\oldUpsilon_i\in\Theta_-(\otheroldLambda)$, and
%$\oldUpsilon_1$ and $\oldUpsilon_1$ are elements of
 $\Theta(\otheroldLambda)$ such that
%$\oldUpsilon_1':=
$\otheroldLambda-\inter\oldUpsilon_1'\in\Theta(\otheroldLambda)$ and
%$\oldUpsilon_2':=
$\otheroldLambda-\inter \oldUpsilon_2'\in\Theta(\otheroldLambda)$ 
%are elements of
%$\oldUpsilon_i\in\Theta_-(\otheroldLambda)$, and
%$\oldUpsilon_1$ and $\oldUpsilon_1$ are elements of
 %$\Theta(\otheroldLambda)$ that
 are in standard position.
%have no degenerate
%crossings. 
It follows from Assertion (3), which has already been proved, that if
%$\oldGamma$ is
%and
$\oldTheta$ is
 defined as in that assertion with $\oldUpsilon_i'$ playing the role of
 $\oldUpsilon_i$ for $i=1,2$, then $\otheroldLambda-\inter\oldTheta\in\Theta(\otheroldLambda)$, and
$[\otheroldLambda-\inter\oldTheta]
\in\barcaly(\otheroldLambda)$  is a supremum for
$\{Y_1,Y_2\}$.

\EndProof
%C'GPQKZXYW$X$UV$SHM isotop\oldTheta\frakK(i)(A)(1)(2)(3)C_X

\Corollary\label{negative lattice}
Let $\otheroldLambda$ be a negative, closed, orientable
$2$-orbifold. 
\begin{enumerate}
\item The set $\barcaly_-(\otheroldLambda)$,
ordered by $\preceq$, is a lattice.
\item If $\oldUpsilon_1$ and $\oldUpsilon_2$ are
elements of $\Theta_-(\otheroldLambda)$ that are in standard position,
and if 
$\oldGamma_-$ denotes the union of all negative components of 
$\oldUpsilon_1\cap \oldUpsilon_2$, then
$\oldGamma_-\in\Theta_-(\otheroldLambda)$, and
$[\oldGamma_-]$ is an infimum for
$\{[\oldUpsilon_1],[\oldUpsilon_2]\}$ in $\barcaly_-(\otheroldLambda)$.
\item If $\oldUpsilon_1$ and $\oldUpsilon_2$ are
elements of $\Theta_-(\otheroldLambda)$ such that  $\otheroldLambda-\inter\oldUpsilon_1$ and $\otheroldLambda-\inter\oldUpsilon_2$ are in standard position, and
$\frakK$ denotes the
union of $\oldUpsilon_1\cup \oldUpsilon_2$ with all components
of $\otheroldLambda-\inter(\oldUpsilon_1\cup \oldUpsilon_2)$ which are discal or are manifolds homeomorphic to $\SSS^1\times[0,\infty)$,
then $\frakK\in\Theta_-(\otheroldLambda)$, and
$[\frakK]$ is a supremum for
$\{[\oldUpsilon_1],[\oldUpsilon_2]\}$ in $\barcaly(\otheroldLambda)$. In particular,
$[\frakK]$ is a supremum for
$\{[\oldUpsilon_1],[\oldUpsilon_2]\}$ in $\barcaly_-(\otheroldLambda)$.
\end{enumerate}
\EndCorollary

\Proof
We first prove (3).
Set $\frakH=\otheroldLambda-\inter (\oldUpsilon_1\cup \oldUpsilon_2)$, and let
 $\oldTheta$ denote the trimming of 
 $\frakH$.
According to Assertion (3) of Proposition \ref{lattice}, we have 
$\otheroldLambda-\inter\oldTheta\in\Theta(\otheroldLambda)$, and
$[\otheroldLambda-\inter\oldTheta]$
is a supremum for $\{[\oldUpsilon_1],[\oldUpsilon_2]\}$ in
$\barcaly(\otheroldLambda)$.

We have $\otheroldLambda-\inter\frakH=\oldUpsilon_1\cup\oldUpsilon_2$. Hence every component of $\otheroldLambda-\inter\frakH$ contains a component of either $\oldUpsilon_1$ or $\oldUpsilon_2$.  For $i=1,2$, since $\oldUpsilon_i\in\Theta_-(\otheroldLambda)$, the components of $\oldUpsilon_i$ are negative and $\pi_1$-injective in $\otheroldLambda$. Hence if $\oldXi$ is any component of $\otheroldLambda-\inter\frakH$, the image of the inclusion homomorphism $\pi_1(\oldXi)\to\pi_1(\otheroldLambda)$ is non-solvable; it follows that $\otheroldLambda-\inter\frakH$ is a negative orbifold. In particular no component of $|\otheroldLambda-\inter\frakH|$ is a weight-$0$ annulus. This implies, in the terminology of \ref{trimming}, that $\frakB_\frakH$ has no redundant components, and that every similarity class of components of $\frakB_\frakH$ has at most one element. Hence the trimming $\oldTheta$ is obtained from $\frakH$ by removing all
components
which are discal or are manifolds homeomorphic to $\SSS^1\times[0,\infty)$. This shows that, in the notation of Assertion (3) of the present corollary, we have $\frakK=\otheroldLambda-\inter\oldTheta$. Thus
$\frakK\in\Theta(\otheroldLambda)$, and 
$[\frakK]$
is a supremum for $\{[\oldUpsilon_1],[\oldUpsilon_2]\}$ in
$\barcaly(\otheroldLambda)$. 

According to the definition of $\frakK$, each component of $\frakK$ contains a component of $\oldUpsilon_1\cup\oldUpsilon_2$. We have seen that, if $\oldXi$ is any component of $\oldUpsilon_1\cup\oldUpsilon_2=\otheroldLambda-\inter\frakH$, the inclusion homomorphism $\pi_1(\oldXi)\to\pi_1(\otheroldLambda)$ has non-solvable  image. Hence if $\frakC$ is any component of $\frakK$, the inclusion homomorphism $\pi_1(\frakC)\to\pi_1(\otheroldLambda)$ has non-solvable  image. This implies that $\frakK$ is negative. Since 
$\frakK\in\Theta(\otheroldLambda)$, it now follows that $\frakK\in\Theta_-(\otheroldLambda)$. Since $[\frakK]$ is a supremum for $\{[\oldUpsilon_1],[\oldUpsilon_2]\}$ in
$\barcaly(\otheroldLambda)$, and lies in $\barcaly_-(\otheroldLambda)$, it is a supremum for $\{[\oldUpsilon_1],[\oldUpsilon_2]\}$ in
$\barcaly_-(\otheroldLambda)$. This proves (3).

Let us now prove (2). According to Proposition \ref{lattice}, if
$\oldGamma$ denotes the
trimming
of $\oldUpsilon_1\cap \oldUpsilon_2$, then
$\oldGamma\in\Theta(\otheroldLambda)$, and $[\oldGamma]$ is an infimum for $\{[\oldUpsilon_1],[\oldUpsilon_2]\}$ in
$\barcaly(\otheroldLambda)$. On the other hand, it follows from the
definitions of $\oldGamma$ and $\oldGamma_-$ that $\oldGamma_-$ is the
union of all negative components of $\oldGamma$. Hence
$\oldGamma_-\in\Theta_-(\otheroldLambda)$, and the relation
$[\oldGamma_-]\preceq[\oldGamma]$ holds in $\barcaly(\otheroldLambda)$. But
in $\barcaly(\otheroldLambda)$, since $[\oldGamma]$  is an infimum for
$\{[\oldUpsilon_1],[\oldUpsilon_2]\}$, we also have
 $[\oldGamma]\preceq[\oldUpsilon_i]$ for $i=1,2$. Hence
 $[\oldGamma_-]\preceq[\oldUpsilon_i]$ for $i=1,2$. Now suppose that
 $Z\in\barcaly_-(\otheroldLambda)$ satisfies  $Z\preceq[\oldUpsilon_i]$ for
 $i=1,2$. Since 
$[\oldGamma]$  is an infimum for
$\{[\oldUpsilon_1],[\oldUpsilon_2]\}$ in $\barcaly(\otheroldLambda)$, we
have $Z\preceq[\oldGamma]$. Hence we may write $Z=[\frakZ]$, with
 $\frakZ\in\Theta_-(\otheroldLambda)$ and $\frakZ\subset\oldGamma$. But
 since $\frakZ$ is negative and taut, any component of $\oldGamma$
containing a component of $\frakZ$ must be negative. This shows that
$\frakZ\subset\oldGamma_-$, so that $Z\preceq[\oldGamma_-]$, and the
proof that $[\oldGamma_-]$  is an infimum for $\{[\oldUpsilon_1],[\oldUpsilon_2]\}$ in
$\barcaly_-(\otheroldLambda)$ is complete.

We now turn to the proof of Assertion (1). If elements $Y_1$ and $Y_2$ of $\barcaly_-(\otheroldLambda)$ are
given, then by \ref{nondeg cross} we may write $Y_i=[\oldUpsilon_i]$
for $i=1,2$, where 
%$\oldUpsilon_i\in\Theta_-(\otheroldLambda)$, and
$\oldUpsilon_1$ and $\oldUpsilon_1$ are elements of $\Theta_-(\otheroldLambda)$ that are in standard position.
%have no degenerate
%crossings. 
It follows from Assertion (2), which has already been proved, that if
$\oldGamma_-$ is
%and
%$\frakK$ are
 defined as in that assertion, then
 $\oldGamma_-\in\Theta_-(\otheroldLambda)$, and $[\oldGamma_-]$ is an infimum for
$\{Y_1,Y_2\}$. 
It also follows from
\ref{nondeg cross}, with \ref{involute}, that we may write $Y_i=[\oldUpsilon_i']$
for $i=1,2$, 
where $\oldUpsilon_1',\oldUpsilon_2'$ are elements of
 $\Theta_-(\otheroldLambda)$ such that
$\otheroldLambda-\inter \oldUpsilon_1'\in\Theta(\otheroldLambda)$ and
$\oldUpsilon_2':=\otheroldLambda-\inter \oldUpsilon_2'\in\Theta(\otheroldLambda)$ 
 are in standard position.
It follows from Assertion (3), which has already been proved, that if
%$\oldGamma$ is
%and
$\frakK$ is
 defined as in that assertion with $\oldUpsilon_i'$ playing the role of
 $\oldUpsilon_i$ for $i=1,2$, then 
$\frakK\in\Theta_-(\otheroldLambda)$, and
$[\frakK]\in\barcaly_-(\otheroldLambda)$ is a supremum for
$\{Y_1,Y_2\}$. 
%We now turn to the proof of (1). If elements $Y_1$ and $Y_2$ of $\barcaly_-(\otheroldLambda)$ are
%given, then by \ref{nondeg cross} we may write $Y_i=[\oldUpsilon_i]$
%for $i=1,2$, where 
%$\oldUpsilon_i\in\Theta_-(\otheroldLambda)$, and
%$\partial\oldUpsilon_1$ and $\partial\oldUpsilon_1$ are elements of $\Theta_-(\otheroldLambda)$ that are in standard position.
%have no degenerate
%crossings. 
%By the assertions already proved, if $\oldGamma_-$ and
%$\frakK$ are defined as in the statement, $[\oldGamma_-]$ and
%$[\frakK]$ are respectively an infimum and a supremum for
%$\{Y_1,Y_2\}$. \redcomment{This contains the same error of exposition
  %pointed out in the last paragraph of the proof of
  %Prop. \ref{lattice}. As in that proof, this works for the infimum,
  %but for the supremum we need to put the closure of the components,
  %not the given suborbifolds themselves, in standard position.} 
This shows that $\barcaly_-(\otheroldLambda)$ is a lattice.
\EndProof
%\oldDelta\frakD\oldXi\frakX\frakC

\Corollary\label{nafta} Let $\otheroldLambda$ be a negative, closed, orientable
$2$-orbifold, and let $\oldUpsilon_1$ and $\oldUpsilon_2$ be elements of $\Theta_-(\otheroldLambda)$. If $\partial|\oldUpsilon_1|\cap
\partial|\oldUpsilon_2|=\emptyset$, and if  $\oldGamma_-$
denotes the union of all negative components of
$\oldUpsilon_1\cap\oldUpsilon_2$, then
$\oldGamma_-\in\Theta_-(\otheroldLambda)$, and $[\oldGamma]$ is an infimum for
$\{[\oldUpsilon_1],[\oldUpsilon_2]\}$ in $\barcaly_-(\otheroldLambda)$.
If $|\oldUpsilon_1|\cap |\oldUpsilon_2|=\emptyset$, then
$\oldUpsilon_1\cup\oldUpsilon_2\in\Theta_-(\otheroldLambda)$, and
$[\oldUpsilon_1\cup\oldUpsilon_2]$ is a supremum for
$\{[\oldUpsilon_1],[\oldUpsilon_2]\}$ in $\barcaly(\otheroldLambda)$. 
In particular, $[\oldUpsilon_1\cup\oldUpsilon_2]$ is a supremum for
$\{[\oldUpsilon_1],[\oldUpsilon_2]\}$ in $\barcaly_-(\otheroldLambda)$. 
\EndCorollary

\Proof
To prove the first assertion, let $X$ denote the union of all components of $|\otheroldLambda-\inter(\oldUpsilon_1\cup\oldUpsilon_2)|$ that are weight-$0$ annuli have one boundary component in $\partial\oldUpsilon_1$ and one in $\oldUpsilon_2$. Set $\oldUpsilon_i'=\oldUpsilon_i\cup\obd(X)$ for $i=1,2$. We have $[\oldUpsilon_i']=[\oldUpsilon_i]$. Since $\partial|\oldUpsilon_1|\cap
\partial|\oldUpsilon_2|=\emptyset$, it follows from the definition (see \ref{nondeg cross}) that $\oldUpsilon_1'$ and $\oldUpsilon_2'$ are in standard position. We have $\oldUpsilon_1'\cap\oldUpsilon_2'=(\oldUpsilon_1\cap\oldUpsilon_2)\discup\obd(X)$, and the components of $\obd(X)$ are annular. Hence the union of all negative components of
$\oldUpsilon_1'\cap\oldUpsilon_2'$ is $\oldGamma_-$. It therefore follows from the second assertion of Corollary \ref{negative lattice} that $\oldGamma_-\in\Theta_-(\otheroldLambda)$, and that $[\oldGamma_-]$ is an infimum for $\{[\oldUpsilon_1'],[\oldUpsilon_2']\}=\{[\oldUpsilon_1],[\oldUpsilon_2]\}$ in $\barcaly_-(\otheroldLambda)$. 

To prove the second assertion, note that since $|\oldUpsilon_1|\cap |\oldUpsilon_2|=\emptyset$, it follows from the definition (see \ref{nondeg cross}) that $\otheroldLambda-\inter\oldUpsilon_1$ and $\otheroldLambda-\inter\oldUpsilon_2$ are in standard position. Note also that every boundary component of $\oldUpsilon_1\cup\oldUpsilon_2$ is a boundary component of one of the $\oldUpsilon_i$, and is therefore $\pi_1$-injective; hence no component of $\otheroldLambda-\inter(\oldUpsilon_1\cup\oldUpsilon_2)$ is discal. It therefore follows from the third assertion of Corollary \ref{negative lattice} that
$\oldUpsilon_1\cup\oldUpsilon_2\in\Theta_-(\otheroldLambda)$; and that $[
\oldUpsilon_1\cup\oldUpsilon_2]$ is a supremum for $\{[\oldUpsilon_1],[\oldUpsilon_2]\}$ in $\barcaly(\otheroldLambda)$, and in particular in $\barcaly_-(\otheroldLambda)$. 
\EndProof
%\oldGamma

\Number\label{components}
Let $\otheroldLambda$ be a negative, closed, orientable $2$-orbifold. If $Y$ is an element of $\barcaly_-(\otheroldLambda)$, if $\oldUpsilon$ is an
element of $\Theta_-(\otheroldLambda)$ such that $[\oldUpsilon]=Y$, and if
$\oldUpsilon_1,\ldots,\oldUpsilon_m$ denote the components of $\oldUpsilon$, it
is clear that the elements $[\oldUpsilon_1],\ldots,[\oldUpsilon_m]$
depend only on $Y$; they will be called the {\it components} of
$Y$. We will say that $Y$ is {\it connected} if it has exactly one
component, i.e. if it has the form $[\oldUpsilon]$ for some connected
$\oldUpsilon\in\Theta_-(
\otheroldLambda)$. 

Note that if $Y,Y'\in\barcaly_-(\otheroldLambda)$ satisfy $Y\preceq
Y'$, and $Y$ is connected, then $Y\preceq Z$ for some component $Z$ of $Y$.
\EndNumber

\Number \label{join and meet}
Let $\otheroldLambda$ be a negative, closed, orientable $2$-orbifold.
Since $\preceq$ is a partial order on $\barcaly(\otheroldLambda)$ (and
hence on $\barcaly_-(\otheroldLambda)$) by Proposition \ref{new partial order},
a subset of $\barcaly_-(\otheroldLambda)$ can have at most one infimum. Hence if
$Y_1,Y_2\in\barcaly_-(\otheroldLambda)$, Corollary \ref{negative lattice} guarantees that
the set $\{Y_1,Y_2\}$, which is of cardinality at most $2$, has a
unique infimum and a unique supremum in $\barcaly_-(\otheroldLambda)$. These will be denoted by
$Y_1\wedge Y_2$ and $Y_1\vee Y_2$, respectively. Standard general
principles about lattices guarantee that $\wedge$ and $\vee$ are
commutative and associative; in particular, $Y_1\wedge\cdots\wedge
Y_p$ and $Y_1\vee\cdots\vee Y_p$ are well defined for any
$Y_1,\ldots,Y_p\in\barcaly_-(\otheroldLambda)$. As is customary in lattice theory,
these will be referred to, respectively, as the {\it meet} and {\it
  join} of $Y_1,\ldots,Y_p$.

It follows by induction from the second assertion of Corollary \ref{nafta} that if 
$\oldUpsilon_1,\ldots,\oldUpsilon_p$ are pairwise disjoint elements of $\Theta_-(\otheroldLambda)$, then $\oldUpsilon_1\cup\cdots\cup\oldUpsilon_p\in\Theta_-(\otheroldLambda)$; and that 
$[\oldUpsilon_1\cup\cdots\cup\oldUpsilon_p]$
is a supremum of $\{[\oldUpsilon_1],\cdots,[\oldUpsilon_p]\} $ in $\barcaly_-(\otheroldLambda)$, so that in particular
$[\oldUpsilon_1\cup\cdots\cup\oldUpsilon_p]=[\oldUpsilon_1]\vee\cdots\vee[\oldUpsilon_p]$. 

It follows that if $Y$ is any element of $\barcaly_-(\otheroldLambda)$, and $Y_1,\ldots,Y_p\in\barcaly_-(\otheroldLambda)$ denote the components of $Y$, then $Y=Y_1\vee\cdots\vee\ Y_p$.
\EndNumber

\Proposition\label{mystery}
Let $\otheroldLambda$ be a negative, closed, orientable $2$-orbifold, and let
$Z$ and $Z'$ be elements of $\barcaly_-(\otheroldLambda)$. Suppose that for every component $Y$ of $Z$ we have $Y\preceq Z'$. Then $Z\preceq Z'$.
\EndProposition

\Proof
Let $Y_1,\ldots,Y_p$ denote the component of $Z$. According to
\ref{join and meet} we have $Z=Y_1\vee\cdots\vee Y_p$; that is, $Z$ is
the supremum of $\{Y_1,\ldots,Y_p\}$. Since $Y_i\preceq Z'$ for each
$i$, the definition of a supremum implies that $Z\preceq Z'$.
\EndProof

\Corollary\label{mystery corollary}
Let $\otheroldLambda$ be a negative, closed, orientable $2$-orbifold, let
$Z$ be an element of $\barcaly_-(\otheroldLambda)$, and let $Z'$ be an element of $\barcaly(\otheroldLambda$). Suppose that for every component $Y$ of $Z$ we have $Y\preceq Z'$. Then $Z\preceq Z'$.
\EndCorollary

\Proof
Write $Z'=[\frakZ']$ for some $\frakZ'\in\Theta(\otheroldLambda)$. Let
$\frakZ''$ denote the union of all negative components of
$\frakZ'$. Then $Z'':=[\frakZ'']\in\barcaly_-(\otheroldLambda)$. If $Y$ is
any component of $Z$, we have $Y\preceq Z'$ by hypothesis, and hence
$Y\preceq [\oldUpsilon']$ for some component $\oldUpsilon'$ of
$\frakZ'$ (see \ref{components}). Since $Z\in\barcaly_-(\otheroldLambda)$ we have $\chibar(Y)>0$, and
Lemma \ref{old partial order} gives $\chibar(\oldUpsilon')\ge\chibar(Y)>0$. Hence
$\oldUpsilon'$ is a component of $\frakZ''$, and so
$Y\preceq[\frakZ'']$. This shows that the hypothesis of Proposition
\ref{mystery} holds with $Z''$ in place of $Z'$; hence $Z\preceq
Z''\preceq Z'$.
\EndProof

\Proposition\label{ascending}
If $\otheroldLambda$ is a negative, closed, orientable $2$-orbifold, there exists no strictly monotone increasing sequence in the lattice $\barcaly_-(\otheroldLambda)$. Furthermore, every subset of $\barcaly_-(\otheroldLambda)$ has a unique supremum.
\EndProposition

\Proof
To prove the first assertion, suppose that $Y_1\prec Y_2\prec\cdots$
is a strictly increasing sequence in $\barcaly_-(\otheroldLambda)$. For each $i\ge0$, set 
$\epsilon_i=(\chibar(Y_i),-\beta(Y_i))\in\ZZ^2$, where  $\beta(Y_1))$ is  defined as in \ref{chidef}. Note that in the notation of \ref{chidef} we have $\alpha(Y_i)=\gamma(Y_i)=0$ for each $i$, since $\otheroldLambda$ is closed and negative and $Y_i\in \barcaly_-(\otheroldLambda)$. 
%\le(\chibar(Y_2),\gamma(Y_2),\alpha(Y_2),-\beta(Y_2))$, 
Since in
particular $Y_i\preceq Y_{i+1}$, it follows from Proposition
\ref{new partial order} that $(\epsilon_i)_{i\ge1}$ is a (weakly) monotone increasing sequence in the lexicographical order of $\ZZ^2$.

%d
%Since $\otheroldLambda$ is a closed orbifold, there is a natural number $K$ such that for any taut, annular suborbifold $\frakJ$ of $\otheroldLambda$ with $\compnum(\frakJ)>K$, some component of $\otheroldLambda-\inter\frakJ$ is a weight-$0$ annulus. 
Set $N=\chibar(\otheroldLambda)$. For every $i\ge1$, we have $Y_i\le [\otheroldLambda]$, so that
Lemma \ref{old partial order}
 implies that $\chibar(Y_i)\le\chibar(\otheroldLambda)= N$. We have %$\gamma(Y_i)\le\compnum(\otheroldLambda)\le N$, and  
$-\beta(Y_i)\le0\le N$. 
%Now suppose that for some $i$ we have $\alpha(Y_i)>K$. Write $Y_i=[\oldUpsilon_i]$ for some $\oldUpsilon_i\in\Theta(\otheroldLambda)$, and let $\oldUpsilon_i'$ denotes the union of all annular components of $\oldUpsilon_i$. The definition of $\alpha(Y_i)$ gives $\compnum(\oldUpsilon_i')=\alpha(Y_i)>K$, which by our choice of $K$ implies that some component $A$ of $|\otheroldLambda-\inter \oldUpsilon_i'|$ is a weight-$0$ annulus. 
%Since the components of $\oldUpsilon_i-\oldUpsilon_i'$ are negative, none of them can be contained in $\obd(A)$, and hence $A$ is a component of $|\otheroldLambda-\inter \oldUpsilon_i'|$. 
%But $\obd(A)$ is then an annular component of $\otheroldLambda-\inter \oldUpsilon_i'$ which shares a boundary component with an annular component of $\oldUpsilon_i$; since $\oldUpsilon_i\in\Theta(\otheroldLambda)$, this contradicts the definition of $\Theta(\otheroldLambda)$ (see \ref{chidef}). Hence we have $\alpha(Y_i)\le K\le N$ for every $i$. 
It now follows that for every $i$ we have $\epsilon_i\in P_N^2$, where $P_N$ denotes the set of all (possibly negative) integers that are less than or equal to $N$. But in the set $P_N^2$ with the lexicographical ordering, every non-empty subset has a greatest element, and hence every monotone increasing sequence is eventually constant. 
%, where $P_N$ denotes the set of all (possibly negative) integers that are less than or equal to $N$. 
%The set $\{0,\ldots,N\}^3$, with the lexicographical ordering, is a finite totally ordered set, and hence every monotone increasing sequence is eventually constant. 
In particular we have
$\epsilon_{i+1}=\epsilon_i$ for some $i\ge1$. Since $\alpha(Y_i)=\alpha(Y_{i+1})=\gamma(Y_i)=\gamma(Y_{i+1})=0$, it then follows from Proposition \ref{new partial order} that
$Y_{i+1}=Y_i$, a contradiction since the given sequence was assumed to
be strictly monotone increasing.

This proves the first assertion. Since $\barcaly_-(\otheroldLambda)$ is a
lattice by Corollary \ref{negative lattice}, the first assertion implies formally that every non-empty subset has a unique supremum. Since $[\emptyset]\in\barcaly_-(\otheroldLambda)$ is clearly a supremum for $\emptyset\subset\barcaly_-(\otheroldLambda)$, the second assertion follows.
\EndProof

\Proposition\label{all at once}
Let $\otheroldLambda$ is a negative, closed, orientable $2$-orbifold, and let $U_1,\ldots,U_N$ be elements of $\barcaly_-(\otheroldLambda)$. Suppose that for every pair of distinct indices $i,j\in\{1,\ldots,N\}$, there exist elements $\frakU$ and $\oldUpsilon$ of $\Theta_-(\otheroldLambda)$ such that $[\frakU]=U_i$, $[\oldUpsilon]=U_j$, and $\frakU\cap\oldUpsilon=\emptyset$. Then there are elements $\frakU_1,\ldots,\frakU_N$ of $\Theta_-(\otheroldLambda)$ such that $[\frakU_i]=U_i$ for $i=1,\ldots,N$ and $\frakU_i\cap\frakU_j=\emptyset$ for each pair of indices $i,j\in\{1,\ldots,N\}$.
\EndProposition

\Proof
We use induction on $N$. If $N\le2$ the assertion is trivial. Now suppose that $N>2$ and that the assertion is true with $N-1$ in place of $N$. Suppose that we are given elements $U_1,\ldots,U_N$ of $\barcaly_-(\otheroldLambda)$ that satisfy the hypothesis. The induction hypothesis gives elements $\frakU_1,\ldots,\frakU_{N-1}$ of $\Theta_-(\otheroldLambda)$ such that $[\frakU_i]=U_i$ for $i=1,\ldots,n-1$ and $\frakU_i\cap\frakU_j=\emptyset$ for each pair of indices $i,j\in\{1,\ldots,N-1\}$. By an observation made in \ref{join and meet}, it then follows that $\frakU_1\cup\cdots\cup\frakU_{N-1}\in\Theta_-(\otheroldLambda)$, and that $[\frakU_1\cup\cdots\cup\frakU_{N-1}]$ is a supremum of
$\{U_1,\ldots,U_{N-1}\}$ in $\barcaly_-(\otheroldLambda)$. 

Now let us fix an element $\frakU_N^0$ of $\Theta_-(\otheroldLambda)$ such that $[\frakU_N^0]=U_N$. By \ref{involute} we have $\otheroldLambda-\inter \frakU_N^0\in\Theta(\otheroldLambda)$. Let $\frakZ$ denote the union of all negative components of $\otheroldLambda-\inter \frakU_N^0$, so that $\frakZ\in\Theta_-(\otheroldLambda)$.

Suppose that an index $i\in\{1,\ldots,N-1\}$ is given. By the hypothesis,  there are elements $\frakU$ and $\oldUpsilon$ of  $\Theta_-(\otheroldLambda)$ such that $[\frakU]=U_i$, $[\oldUpsilon]=U_N$, and $\frakU\cap\oldUpsilon=\emptyset$. Hence $\frakU$ is isotopic to a suborbifold $\frakU'$ of $\otheroldLambda$ with $\frakU'\subset \otheroldLambda-\inter\frakU_N^0$. Since $\frakU'$ is taut and negative, the component of 
$\otheroldLambda-\inter\frakU_N^0$ containing $\frakU'$ must be negative. This means that $\frakU'\subset\frakZ$.
% Since $[\frakU]=U_i$ and $[\oldUpsilon]=U_N$, the suborbifolds $\frakU$ and $\frakU_i$ are isotopic, as are the suborbifolds $\oldUpsilon$ and $\frakU_N^0$. This implies that $\frakU_i$ is isotopic to a suborbifold disjoint from $\frakU_N^0$. This means that 
Hence the relation $U_i\preceq[\frakZ]$ holds in $\Theta(\otheroldLambda)$. Since this is true for every $i\in\{1,\ldots,N-1\}$, and since $[\frakU_1\cup\cdots\cup\frakU_{N-1}]$ is a supremum of
 $\{U_1,\ldots,U_{N-1}\}$ in $\barcaly_-(\otheroldLambda)$, we have $[\frakU_1\cup\cdots\cup\frakU_{N-1}]\preceq
[\frakZ]$. By definition this means that there is a suborbifold $\frakZ'$ of $\otheroldLambda$, isotopic to $\frakZ$, such that $\frakU_1\cup\cdots\cup\frakU_{N-1}\subset\inter \frakZ'$. Since $\frakZ\subset \otheroldLambda-\inter \frakU_N^0$, there is a suborbifold $\frakV$ of $\otheroldLambda$, isotopic to 
$ \otheroldLambda-\inter \frakU_N^0$, such that
$\frakU_1\cup\cdots\cup\frakU_{N-1}\subset\inter \frakV$.
 Then $\frakU_N:=\otheroldLambda-\inter\frakV$ is isotopic to $\frakU_N^0$, so that $\frakU_N\in\Theta_-(\otheroldLambda)$ and $[\frakU_N]=U_N$. Since $\frakU_1\cup\cdots\cup\frakU_{N-1}\subset\inter\frakV$, we have $\frakU_i\cap\frakU_N=\emptyset$ for $i=1,\ldots,N$.
\EndProof
%_n_{n-1} $N n-1\}\oldUpsilon Y\oldUpsilon'\oldLambda\frakZ\frakU'\frakW\oldUpsilon\frakV\frakU_N^*

\Number\label{distribute} 
If 
$\otheroldLambda$ is a negative, closed, orientable $2$-orbifold, and
$(Y_i)_{i\in S}$ is an indexed family of elements of
  $\barcaly_-(\otheroldLambda)$, we shall denote its supremum, which is
  well-defined by Proposition \ref{ascending}, by $\bigvee_{i\in
    S}Y_i$. Note that for any such family, and for any element $Z$ of
  $\barcaly_-(\otheroldLambda)$, we have the identity
$$(\bigvee_{i\in
    S}Y_i)\wedge Z=
\bigvee_{i\in
    S}(Y_i\wedge Z),$$
as this identity holds in any lattice in which every subset has a
supremum. Note also that in the case where the index set is finite and is given in the form $S=\{1,\ldots,p\}$, where $p$ is a positive integer, we have $\bigvee_{i\in S}Y_i=Y_1\vee\cdots\vee Y_p$ in the notation of \ref{join and meet}.
\EndNumber

\Lemma\label{babi also}
Let $\otheroldLambda$ be a negative, closed, orientable $2$-orbifold, let $Y$ and $T$ be elements of $\barcaly_-(\otheroldLambda)$, and let $T_1,\ldots,T_p$ denote the components of $T$. Then $\chibar(Y\wedge T)=\sum_{i=1}^p\chibar(Y\wedge T_i)$.
\EndLemma

\Proof
As we pointed out in \ref{join and meet}, we have $T=T_1\vee\cdots\vee T_p$. Hence by \ref{distribute}, we have $T\wedge Y=(T_1\wedge Y)\vee\cdots\vee (T_p\wedge Y)$.

Write $T=[\oldXi]$ for some $\oldXi\in\Theta_-(\otheroldLambda)$. Then the components of $\oldXi$ may be indexed as $\oldXi_1,\ldots,\oldXi_p$ in such a way that $[\oldXi_i]=T_i$ for $i=1,\ldots,p$. For each $i$, since $T_i\wedge
Y\preceq T_i$, there is an element $\frakZ_i$ of $\Theta_-({\otheroldLambda})$ such that $[\frakZ_i]=T_i\wedge
Y$ and $\frakZ_i\subset\oldXi_i$. In particular the $\frakZ_i$ are pairwise disjoint. According to \ref{join and meet} we therefore have
$[\frakZ_1\cup\cdots\cup\frakZ_p]
 =(T_1\wedge
Y)\vee\cdots\vee(T_p\wedge
Y)=T\wedge Y$. It follows that $\chibar(T\wedge Y)= %\chibar([\frakZ_1])+\cdots+\chibar(\frakZ_p)]=
\chibar([\frakZ_1])+\cdots+\chibar([\frakZ_p])=\chibar(T_1\wedge Y)+\cdots+\chibar(T_p\wedge Y)$.  
\EndProof
%X'')}\oldXi\eta $t

\Definition\label{meet me in st. louis}
Let $C_1$ and $C_2$ be $1$-dimensional submanifolds of the interior of a closed, orientable $2$-orbifold $\otheroldLambda$ (so that the $C_i$ are submanifolds of $|\otheroldLambda|$ disjoint from $\fraks_\otheroldLambda$). We will say that $C_1$ and $C_2$ {\it meet essentially} if for all $1$-manifolds $C_1',C_2'\subset\inter\otheroldLambda$ such that $C_i'$ is (orbifold-)isotopic to $C_i$ for $i=1,2$, we have $C_1'\cap C_2'\ne\emptyset$.
\EndDefinition

\Proposition\label{tpp} Let $\otheroldLambda$ be a negative, closed, orientable
 $2$-orbifold. If $\oldUpsilon_1$ and
$\oldUpsilon_2$ are elements of $\Theta_-(\otheroldLambda)$ such that 
%$Y_1$ and $Y_2$ are elements of $\barcaly_-(\otheroldLambda)$ such that   
$\partial\oldUpsilon_1$ and $\partial\oldUpsilon_2$ meet essentially,
%\ne\emptyset$ for every $\oldUpsilon_1\in Y_1$ and every $\oldUpsilon_2\in Y_2$, 
and if $Y=[\oldUpsilon_1]\wedge[\oldUpsilon_2]\in \barcaly_-(\otheroldLambda)$, then $\chibar(Y)<\min(\chibar(\oldUpsilon_1),\chibar(\oldUpsilon_2))$.
\EndProposition

\Proof
It follows from \ref{nondeg cross} that we may assume $\oldUpsilon_1$ and
$\oldUpsilon_2$ to have been chosen within their isotopy classes so as to be in standard position.
%that 
%$\partial\oldUpsilon_1$ and $\partial\oldUpsilon_2$ intersect
%transversally and have no
%nondegenerate crossing. 
Let $A_1,\ldots,A_m$ denote the
components of $\oldUpsilon_1\cap\partial \oldUpsilon_2$, where a
priori we have $m\ge0$, and each $A_i$ is an arc or a simple
closed curve. Since $\partial\oldUpsilon_1$ and
$\partial\oldUpsilon_2$ meet essentially, we have
$\partial\oldUpsilon_1\cap\partial\oldUpsilon_2\ne\emptyset$; hence
$m\ge1$, and at least one of the $A_i$ is an arc. We may
therefore assume the $A_i$ to be indexed so that $A_1$ is an
arc. For $i=0,\ldots,m$ let $\oldXi_i$ denote the $2$-orbifold
$(\oldUpsilon_1)\cut{A_1\cup\ldots\cup A_i}$. Then for $i=1,\ldots,m$,
the arc or curve $A_i$ is
canonically identified with an arc or curve in $\oldXi_{i-1}$, and
$\oldXi_i$ is canonically identified with $(\oldXi_{i-1})\cut{A_i}$.

If $\frakZ$ is any compact $2$-orbifold, we will define
$\chi_-(\frakZ)$ to be $\sum_{i=1}^n\max(0,\chibar(\frakZ_i))$, where
$\frakZ_1,\ldots,\frakZ_n$ denote the components of $\frakZ$. We
claim:
\Claim\label{just because}
For $i=1,\ldots,m$ we have
$\chi_-(\oldXi_i)\le\chi_-(\oldXi_{i-1})$. Furthermore, if $i$ is an
index with $1\le i\le m$ such that $A_i$ is an arc and is not the frontier of a
weight-$0$ disk in $\oldXi_{i-1}$, and such that the component of
$\oldXi_{i-1}$ containing $A_i$ has strictly negative Euler
characteristic, then $\chi_-(\oldXi_i)<\chi_-(\oldXi_{i-1})$.
\EndClaim

To prove \ref{just because}, let an index $i$ with $1\le i\le m$ be
given. Set $A=A_i$, let $\frakU$ denote the component of $\oldXi_{i-1}$ containing
$A$, and set $\frakU'=\frakU\cut A$. Note that $\frakU'$
has either one or two components. We have
$\chi_-(\oldXi_{i-1})-\chi_-(\oldXi_{i})=\chi_-(\frakU)-\chi_-(\frakU')$. Hence
it suffices to show that
$\chi_-(\frakU')\le\chi_-(\frakU)$; and that if $A$ is an arc
and is not the frontier of a
weight-$0$ disk in $\frakU$, and $\chi(\frakU)<0$, then $\chi_-(\frakU')<\chi_-(\frakU)$.

First consider the case in which $A$ is a simple closed
curve. Since $\oldUpsilon_2$ is taut,
% by the definition of
%$\Theta_-(\otheroldLambda)$, the curve 
$A$
%, regarded as a boundary
%component of $\oldUpsilon_2$, 
is $\pi_1$-injective in $\otheroldLambda$ and
hence in 
%$\oldUpsilon_1$; hence $A$ is also $\pi_1$-injective when regarded as a simple
%closed curve in
$\frakU$. It follows that each
component of $\frakU'$ has non-positive Euler characteristic, 
%and hence 
so that $\chi_-(\frakU')=\chibar(\frakU')$. Since $A$ is a
simple closed curve we have $\chibar(\frakU)=\chibar(\frakU')\ge0$,
and hence
$\chi_-(\frakU)=\chibar(\frakU)=\chibar(\frakU')=\chi_-(\frakU')$. Thus
the required inequality is an equality in this case.

Now consider the case in which $A$ is an arc, and $\frakU'$ has
no discal component. Then  each
component of $\frakU'$ has non-positive Euler characteristic, and
hence $\chi_-(\frakU')=\chibar(\frakU')$. Since $A$ is an
arc, we have $\chibar(\frakU)=1+\chibar(\frakU')>0$, and hence
$\chi_-(\frakU)=\chibar(\frakU)=1+\chibar(\frakU')>\chibar(\frakU')=\chi_-(\frakU')$.

Next consider the case in which $A$ is an arc and each
component of $\frakU'$ is discal. Then $\chi_-(\frakU')=0$, and since
$\chi_-(\frakU)$ is by definition non-negative, we have
$\chi_-(\frakU')\le\chi_-(\frakU)$. Furthermore, if $\chi(\frakU)<0$, we have
$\chi_-(\frakU)=\chibar(\frakU)>0$, so that $\chi_-(\frakU')<\chi_-(\frakU)$.

There remains only the case in which $A$ is an arc and
$\frakU'$ has two components, $\frakD$ and $\frakV$, where $\frakD$
is discal and $\frakV$ is not. Then $|\frakD|$ is a disk, and
$\fraks_\frakD$ consists of at most one point. If $\fraks_\frakD$
consists of one point, we denote the order of the point by $p>1$; if
$\fraks_\frakD=\emptyset$, we set $p=1$. Then
$\chi(\frakD)=1/p$. Since $A$ is an arc, we have
\Equation\label{ach du lieber}
\chibar(\frakU)=1+\chibar(\frakD)+\chibar(\frakV)=(1-1/p)+\chibar(\frakV)\ge\chibar(\frakV).
\EndEquation
Since $\frakV$ is not discal, we have $\chibar(\frakV)\ge0$ and hence
$\chibar(\frakU)\ge0$ by (\ref{ach du lieber}); thus
$\chi_-(\frakU)=\chibar(\frakU)$, while
$\chi_-(\frakU')=\chibar(\frakV)$. Now by (\ref{ach du lieber}), we
have
$\chi_-(\frakU)=\chibar(\frakU)\ge\chibar(\frakV)=\chi_-(\frakU')$. If
$A$ is not the frontier of a weight-$0$ disk in $\frakU$, then
$\fraks_\frakD\ne\emptyset$, and hence $p>1$. Thus (\ref{ach du
  lieber}) gives $\chibar(\frakU)>\chibar(\frakV)$, so that
$\chi_-(\frakU)=\chibar(\frakU)>\chibar(\frakV)=\chi_-(\frakU')$. This
completes the proof of \ref{just because}.

The orbifold $\oldXi_0=\oldUpsilon_1$ is negative since
$\oldUpsilon_1\in\Theta_-(\otheroldLambda)$, and $A_1$ is an arc by our
choice of indexing. Since $\oldUpsilon_1$ and $\oldUpsilon_2$ are in standard position, $\partial\oldUpsilon_1$ and
$\partial\oldUpsilon_2$ have no degenerate crossing, and therefore  $A_1$ is not the frontier of a
weight-$0$ disk in $\oldXi_{0}$.
%,  no component of
%$\oldXi_1=(\oldUpsilon_1)\cut{A_1}$ is a weight-$0$ disk. 
Hence
\ref{just because} gives
$\chi_-(\oldXi_1)<\chi_-(\oldXi_0)=\chibar(\oldUpsilon_1)$. Since
\ref{just because} also gives $\chi_-(\oldXi_i)\le\chi_-(\oldXi_{i-1})$
whenever $1<i\le m$, it follows that
$\chi_-(\oldXi_m)<\chibar(\oldUpsilon_1)$. But $\oldUpsilon_1\cap\oldUpsilon_2$
is canonically identified with a union of components of
$
%(\oldUpsilon_1)\cut{\oldUpsilon_1\cap\partial\oldUpsilon_1}=
\oldXi_m$. Hence 
\Equation\label{supoib}
\chi_-(\oldUpsilon_1\cap\oldUpsilon_2)\le
%\chi_-((\oldUpsilon_1)\cut{\oldUpsilon_1\cap\partial\oldUpsilon_1})=
\chi_-(\oldXi_m)<\chibar(\oldUpsilon_1).
\EndEquation

Now let $\oldGamma_-$ denote the union of all negative components of
$\oldUpsilon_1\cap\oldUpsilon_2$. Since $\oldUpsilon_1$ and $\oldUpsilon_2$ are in standard position,  it follows from Assertion (2) of Corollary
\ref{negative lattice} that $Y=[\oldGamma_-]$.  On
the other hand, it follows from the definitions that
$\chibar(\oldGamma_-)=\chi_-(\oldUpsilon_1\cap\oldUpsilon_2 )$. Using
(\ref{supoib}), we now obtain
$\chibar(Y)=\chibar(\oldGamma_-)=\chi_-(\oldUpsilon_1\cap\oldUpsilon_2
) <\chibar(\oldUpsilon_1)$. The same argument, with the roles of
$\oldUpsilon_1$ and $\oldUpsilon_2$ reversed, shows that
$\chibar(Y)<\chibar(\oldUpsilon_2)$, and the proposition is proved.
\EndProof
%$1$A\alpha

The following two special results will be needed in Section \ref{clash section}:

\Lemma\label{before inchworm}
Let $\otheroldLambda$ be a negative, closed, orientable $2$-orbifold, and let $Z$ and $Z'$ be elements of $\barcaly_-(\otheroldLambda)$ such that $Z'\preceq Z$ and $\chibar(Z)=\chibar(Z')$. Then there exist elements $\frakZ$ and $\frakZ'$ of $\Theta_-(\otheroldLambda)$ such that (i) $[\frakZ]=Z$, (ii)$[\frakZ']=Z'$, (iii)
$\frakZ'\subset\frakZ$, (iv) $\partial\frakZ\subset\partial\frakZ'$, and (v) each component of $\overline{\frakZ-\frakZ'}$ is an annular orbifold contained in $\inter\frakZ$.
\EndLemma

\Proof
Since $Z'\preceq Z$, there are elements $\frakZ_0$ and $\frakZ'$ of $\Theta_-(\otheroldLambda)$ such that $[\frakZ_0]=Z$, $[\frakZ']=Z'$, and $\frakZ'\subset\inter\frakZ_0$. By hypothesis we have $\chibar(Z)=\chibar (Z')$, and since $\otheroldLambda$ is closed and negative, and 
$Z, Z'\in\barcaly_-(\otheroldLambda)$, we have $\gamma(Z)=\gamma(Z')=0$ in the notation of \ref{chidef}. It therefore follows from Lemma \ref{old partial order} that every component of $\frakZ_0-\inter\frakZ'$ is an annular orbifold. In particular, if for each component $C$ of $|\partial\frakZ_0|$ we denote by $\frakB_C$ the component of $\frakZ_0-\inter\frakZ'$ containing $\obd(C)$, then $\frakB_C$ is an annular orbifold. Since $\frakZ_0$ is negative by the definition of $\Theta_-(\otheroldLambda)$, the annular orbifold $\frakB_C$ cannot be a component of $\frakZ_0$; it must therefore have a boundary component which is distinct from $C$ and is contained in $\inter\frakZ_0$. It follows that $|\frakB_C|$ is a weight-$0$ annulus, and that it has $C$ as one boundary component, while its other boundary component is contained in $|\partial\frakZ'|$. Hence the orbifold $\frakZ:=\overline{\frakZ_0-\bigcup_{C\in\calc(\partial\frakZ_0)}\frakB_C}$ is isotopic to $\frakZ_0$ in $\otheroldLambda$, so that $\frakZ\in\Theta_-(\otheroldLambda)$ and $ [\frakZ]=Z$; and we have $\frakZ'\subset\frakZ$, and $\partial\frakZ\subset\partial\frakZ'$. Furthermore, each component of $\overline{\frakZ-\frakZ'}$ is a component of $\frakZ_0-\inter\frakZ'$ contained in $\inter\frakZ$, and is therefore an annular orbifold contained in $\inter\frakZ$. 
\EndProof
%\inter

\Lemma\label{inchworm lemma}
Let $\otheroldLambda$ be a negative, closed, orientable $2$-orbifold, and let $Z_1$ and $Z_2$ be elements of $\barcaly_-(\otheroldLambda)$ such that $\chibar(Z_1\wedge Z_2)=\chibar(Z_1)=\chibar(Z_2)$. For $i=1,2$, let $\frakZ_i$ be an element of $\Theta_-(\otheroldLambda)$ such that $[\frakZ_i]=Z_i$, and let $\oldUpsilon_i$
denote the union of $\frakZ_i$ with all components of $\otheroldLambda-\inter\frakZ_i$ that are annular orbifolds. Then $[\oldUpsilon_1]=[\oldUpsilon_2]$.
\EndLemma

\Proof
First consider the special case in which $Z_2\preceq Z_1$. In this case the hypothesis asserts that $\chibar(Z_1)=\chibar (Z_2)$. Thus the hypotheses of Lemma \ref{before inchworm} are satisfied with $Z_1$ and $Z_2$ playing the respective roles of $Z$ and $Z'$. Hence there exist elements $\frakZ$ and $\frakZ'$ of $\Theta_-(\otheroldLambda)$ such that the conclusions of Lemma \ref{before inchworm} hold with $Z_1$ and $Z_2$ in place of $Z$ and $Z'$. On the other hand, it is clear that the truth of the conclusion of the present lemma depends only on the isotopy classes of $\frakZ_1$ and $\frakZ_2$; we may therefore assume that $\frakZ=\frakZ_1$ and that $\frakZ'=\frakZ_2$, i.e.
%properties asserted for the $\frakZ_i$ are independent of the way they are chosen within their isotopy classes, we may assume 
that 
(i) $[\frakZ_1]=Z_1$, (ii)$[\frakZ_2]=Z_2$, (iii)
$\frakZ_2\subset\frakZ_1$, (iv)$\partial\frakZ_1\subset\partial\frakZ_2$, and (v) each component of $\overline{\frakZ_1-\frakZ_2}$ is an annular orbifold contained in $\inter\frakZ_1$. 

It follows that every component of $\overline{\frakZ_1-\frakZ_2}$ is an annular component of $\otheroldLambda-\inter\frakZ_2$; and that the annular components of $\otheroldLambda-\inter\frakZ_2$ that are not components of $\overline{\frakZ_1-\frakZ_2}$ are precisely the annular components of $\otheroldLambda-\inter\frakZ_1$. 
%Hence  $\oldUpsilon_2$ may be described as the union of $\frakZ_1$ with all components of $\otheroldLambda-\inter\frakZ_1$ that are annular orbifolds. Since $\frakZ_1$ is isotopic to $\frakZ_1$, it follows from this description of 
In view of the definition of the $\oldUpsilon_i$, it now follows that $\oldUpsilon_2=\oldUpsilon_1$. This proves the lemma in the special case where $Z_2\preceq Z_1$.

To prove the lemma in general, set $Z_3=Z_1\wedge Z_2$, and note that by definition we have $Z_3\preceq Z_i$ for $i=1,2$. The hypothesis gives that $\chibar(Z_3)=\chibar(Z_i)$ for $i=1,2$. Hence if we denote by $\oldUpsilon_3$ the union of $\frakZ_3$ with all components of $\otheroldLambda-\inter\frakZ_3$ that are annular orbifolds, the special case of the lemma already proved, with $Z_i$ and $Z_3$ playing the roles of $Z_1$ and $Z_2$, shows that $[\oldUpsilon_i]=[\oldUpsilon_3]$ for $i=1,2$. It follows that  $[\oldUpsilon_1]=[\oldUpsilon_2]$.
\EndProof
%\oldUpsilon Y Z \frakZ X

Whereas the  section up to now has been primarily about suborbifolds and their equivalence classes, the rest of the section will be concerned with homeomorphisms between suborbifolds, and equivalence classes of such homeomorphisms.

\Proposition\label{before i guess}
Let $\otheroldLambda$ be an orientable $2$-orbifold without boundary,
 let $\otheroldLambda_0$ be a negative, taut,
finite-type $2$-suborbifold of $\otheroldLambda$, and
let $\oldUpsilon$ be a finite-type
orientable $2$-orbifold, no component of which is discal or is a manifold homeomorphic to $\SSS^1\times[0,\infty)$. 
Let $\xi,\xi':\oldUpsilon\to\otheroldLambda_0$ be proper
$\pi_1$-injective embeddings which are properly isotopic in
$\otheroldLambda$. Then either (A) $\xi$ and $\xi'$ are properly isotopic in
$\otheroldLambda_0$, or (B) some component $F$ of $|\oldUpsilon|$ is a weight-$0$ annulus, and a core curve of the weight-$0$ annulus $|\xi(F)|$ cobounds a weight-$0$ annulus in $|\otheroldLambda_0|$ with some boundary component of $|\partial\otheroldLambda|$.
\EndProposition

%\Remark
%The condition in the statement of Lemma \ref{i guess} that every component of %$\partial\otheroldLambda_0$ is contained in $\partial\otheroldLambda$ or in $\inter\otheroldLambda$
%is convenient for the proof. In the applications we will have the stronger condition %that $\otheroldLambda_0\subset\inter\otheroldLambda$.
%\EndRemark

\Proof
First consider the case in which $\otheroldLambda$ and $\oldUpsilon$ are
$2$-manifolds. In this case, to emphasize the fact that $\otheroldLambda$,
$\otheroldLambda_0$ and $\oldUpsilon$ are manifolds, we will write $L=\otheroldLambda$,
$G=\otheroldLambda_0$ and $Y=\oldUpsilon$.  We may assume without loss of generality that
$ L$ is connected
%, that $ G\subset\inter L$,
and that $\xi$ and $\xi'$ map $Y$ into
$\inter G$. Let $ G_1,\ldots, G_m$ denote the
components of $ G$, and set
$Y_i=\xi^{-1}( G_i)$ for $i=1,\ldots,m$. (Note
that the $Y_i$  may be disconnected or empty.) For
$i=1,\ldots,m$, fix a basepoint  $x_i\in G_i\subset L$,
and let $p_i:(\tL_i,\tx_i)\to ( L,x_i)$ denote the based covering space of
$( L,x_i)$ defined by the subgroup of $\pi_1( L,x_i)$
which is the image of the inclusion homomorphism
$\pi_1( G_i,x_i)\to\pi_1( L,x_i)$. Let
$s_i:( G_i,x_i)\to(\tL_i,\tx_i)$ denote the based lift
of the inclusion map $( G_i,x_i)\to( L,x_i)$; then
$p_i$ maps  $\tGee_i:=s_i( G_i)$ homeomorphically onto
$ G_i$. The definition of the covering $\tL_i$ implies
that the inclusion homomorphism
$\pi_1(\tGee_i)\to\pi_1(\tL_i)$ is
surjective. Furthermore, $\partial\tGee_i$ is $\pi_1$-injective in $\tL_i$ since $ G$ is taut in $ L$. Hence every component of  %$\tGee$ is a deformation retract of $\tL_i$, and every component of  
 $\tL_i-(\inter\tGee_i)$ is a half-open annulus whose boundary is a component of $\partial\tGee_i$.

Set $\txi_i=s_i\circ(\xi|Y_i)$, so that $\txi_i$ is
a lift of $\xi|Y_i$. Let
$H_i:Y_i\times[0,1]\to L$ be a non-ambient proper
isotopy from $\xi|Y_i$ to $\xi'|Y_i$; thus
$H_i$ is a proper map,
$h_t^i:t\mapsto H_i(x,t)$ is an embedding for each $t$, and we
have $h_0^i=\xi|Y_i$ and
$h_1^i=\xi'|Y_i$. 
By the covering homotopy property, $H_i$ admits a lift $\tH_i:Y_i\to\tL_i$ such that, if we denote by $\thh^i_t$ the map $x\mapsto H_i(x,t)$ for $0\le t\le1$, we have $\thh^i_0=\txi_i$. Since $h^i_t$ is an embedding for each $t$, so is $\thh^i_t$; that is, $\tH_i$ is an isotopy. Since the isotopy $H_i$ is proper, so is $\tH_i$. For any component $F$ of $Y_i$, since $\xi'(F)\subset\xi'(Y)\subset G$, the lift $\thh^i_1|F$ of $\xi'|F$ must map $F$ into some component of $p_i^{-1}( G)$. We distinguish two subcases.

{\bf Subcase I.} For some $i_0\in\{1,\ldots,m\}$ and some component $F_0$ of $Y_{i_0}$, the component of $p_{i_0}^{-1}(G)$ containing $\thh^{i_0}_1(F_0)$ is distinct from $\tGee_{i_0}$.

In this subcase, let $U\ne \tGee_{i_0}$ denote the component of $p_{i_0}^{-1}(G)$ containing $\thh^{i_0}_1(F_0)$. Then $U$ is contained in some component $V$ of $\tL_{i_0}-\tGee_{i_0}$. We have observed that every component of $\tL_{i_0}-\inter\tGee_{i_0}$ is a half-open annulus. Hence $V$ is an open annulus, and by $\pi_1$-injectivity, $\pi_1(F_0)$ is cyclic. Since $F_0$ is a component of $ Y$, the hypothesis implies that $F_0$ is not a disc and is not homeomorphic to $\SSS^1\times[0,\infty)$. Hence $F_0$ is a (closed) annulus.

Let $C$ denote a core curve of $F_0$. For $t=0,1$ set $\tC_t=\thh^{i_0}_t(C)$. Then $\tC_0=\txi_{i_0}(C)$ is a core curve of the annulus $\txi_{i_0}(F_0)$, so that $C_0:=p_{i_0}(\tC_0)=\xi(C)$ is a core curve of the annulus $\xi(F_0)=p_{i_0}(\txi_{i_0}(F_0))$; and the curves $\tC_0$ and $\tC_1$ are isotopic, so that by \cite[Lemma 2.4]{epstein} they cobound an annulus $\tA \subset\tL_{i_0}$. Since $\tC_0\subset\inter\tGee_{i_0}$ and $\tC_1\subset
V\subset\tL_{i_0}-\tGee_{i_0}$, we have $\tA \cap\partial\tGee_{i_0}\ne\emptyset$. By tautness, every component of 
$\tA \cap\partial\tGee_{i_0}$ is a homotopically non-trivial simple closed curve in $\tL_{i_0}$, and hence in $\tA $. If $\tB$ denotes the closure of the component of $\tA \setminus\partial\tGee_{i_0}$ containing $\tC_0$, it follows that $\tB$ is an annulus contained in $\tGee_{i_0}$ whose boundary components are $\tC_0$ and some component of $\partial\tGee_{i_0}$. Hence $B:=p_{i_0}(\tB)$ is an annulus contained in $G$ whose boundary components are $C_0$ and some component of $\partial G$. Thus Alternative (B) of the conclusion of the proposition holds in this subcase.
%AB

{\bf Subcase II.} For every $i\in\{1,\ldots,m\}$ and every component $F$ of $Y_i$, the component of $p_i^{-1}(G)$ containing $\thh^i_1(F)$ is equal to $\tGee_i$.

In this subcase, for $i=1,\ldots,m$, the lift $\thh^i_1$ of $\xi'|Y_i$  maps $Y_i$ into $\tGee_i$.  In particular it follows, upon composing with $p_i$, that $\xi'_i(Y_i)\subset G_i$. We will show that Alternative (A) of the conclusion of the proposition holds in this subcase.

Let $N_i$ be a regular neighborhood of $\partial\tGee_i$ in $\tGee_i$ which is  disjoint from $\txi(Y_i)$ and from $\txi'(Y_i)$. Set $\tGee_i^-=\overline{\tGee_i-N_i}$ for $i=1,\ldots,m$.
Since
every component of %$\tGee$ is a deformation retract of $\tL_i$, and every component of  
 $\tL_i-\inter\tGee_i$ is a half-open annulus whose boundary is a component of $\partial\tGee_i$, 
 there is a homeomorphism $\phi_i:\tL_i\to\inter\tGee_i$ which is the identity on $\tGee_i^-$. In particular $\phi_i$ is the identity on $\txi(Y_i)$ and on $\txi'(Y_i')$. Hence $\phi_i\circ\tH_i$ is an isotopy in $\tGee_i$ from $\txi_i$ to $\txi_i'$. Since $p_i$ maps $\tGee_i$ homeomorphically onto $ G_i$, we have an isotopy
$H_i':=p_i\circ\phi_i\circ\tH_i$ in $ G_i$ from $\xi|Y_i$ to
$\xi'|Y_i$. Since $G_1,\ldots,G_m$ are pairwise disjoint, we may now define an isotopy $H'$ from $\xi$ to $\xi'$ in $ G=G_1\cup\cdots\cup G_m$ by setting $H'|(Y_i\times[0,1])=H_i'$ for $i=1,\ldots,m$. In order to show that Alternative (A) holds in this subcase,
% Since such an isotopy exists for $i=1,\ldots,m$ (and
%), the maps $\xi$ and
%$\xi'$ are isotopic in $ G=G_1\cup\cdots\cup G_m$, and
%proof of the lemma in the case where  $\otheroldLambda$ and $\oldUpsilon$ are $2$-manifolds, 
it now suffices to show that the isotopy $H'$ is proper. For this
purpose, it is enough to show each $H_i'$ is proper; and since $p_i$
maps $\tGee_i$ homeomorphically to $G_i$, for each $i$, it suffices to
show that $\phi_i\circ\tH_i$ is proper. This will be established
using the notion of ends of a locally compact space. For any  locally compact space $X$, we will denote
the set of ends of $X$ by $\Ends(X)$. If $f:X\to X'$ is a proper map between locally compact spaces, it induces a map $f_\bends:\Ends(X)\to\Ends(X')$. 

Since $G$ is taut, the boundary of each $G_i$ is compact, and since $\tGee_i$ is homeomorphic to $G_i$ for $i=1,\ldots,m$, the $1$-manifold
%each $\tGee_i$ also have compact boundary. Hence 
$\Fr_{\tL_i}\tGee_i=\partial\tGee_i$ is compact for each $i$. We therefore have a natural identification of $\Ends(\tGee_i)$ with a subset of $\Ends(\tL_i)$. Since $N_i$ is compact,  $\Ends(\tGee_i)$ is also canonically identified with $\Ends(\tGee_i^-)$. Since $\txi_i(Y_i)\subset\tGee_i$, the image of  $(\txi_i)_\bends:\Ends(Y_i)\to\Ends(\tL_i)$ is contained in $\Ends(\tGee_i)\subset\Ends(\tL_i)$. Since $\tH_i$ is a proper homotopy and $\thh^i_0=\txi_i$, it follows that the image of $(\tH_i)_\bends:\Ends(Y_i\times[0,1])\to\Ends(\tL_i)$ is also contained in $\Ends(\tGee_i)=\Ends(\tGee_i^-)$. This implies that $Z:=\overline{(Y_i\times[0,1])-\tH_i^{-1}(\tGee_i^-)}$ is compact. Since $\phi_i$ is the identity on $\tGee^-$, we can now deduce that $\phi_i\circ\tH_i$ agrees with $\tH_i$ except on the compact set $Z$. Since the isotopy $\tH_i$ is proper, so is the isotopy $\phi_i\circ\tH_i$. 
This completes the 
proof of the lemma in the case where  $\otheroldLambda$ and $\oldUpsilon$ are $2$-manifolds.

To prove the lemma in the general case, note that by Lemma \ref{historic fact},  $\dot\xi$ and $\dot\xi'$ are  properly isotopic embeddings of 
$\dot\oldUpsilon$ in $\dot\otheroldLambda$. Since no component of $\oldUpsilon$ is a discal orbifold or a manifold homeomorphic to $\SSS^1\times[0,\infty)$,  no component of $\dot\oldUpsilon$ is a disk or a manifold homeomorphic to $\SSS^1\times[0,\infty)$. Furthermore, $\dot\otheroldLambda_0$ has finite type by \ref{vere's dot notation}, and it follows from Lemma \ref{geventlach} that $\dot\otheroldLambda_0$ is a negative, taut submanifold of $\dot\otheroldLambda$.  By the case of the present lemma
already proved, it follows that $\dot\xi$ and $\dot\xi'$ are
properly isotopic in $\dot\otheroldLambda_0$. Hence by
Lemma \ref{historic fact}, $\xi$ and $\xi'$ are
properly isotopic in $\otheroldLambda_0$, as
required.
\EndProof
%alpha \xi_i H_(H^i) Z X W U V C (i) (ii) Case case

In the special case of Proposition \ref{before i guess} in which $\oldUpsilon$ is negative, it is obvious that Alternative (B) of the conclusion cannot occur. Hence we obtain the

\Corollary\label{i guess}
Let $\otheroldLambda$ be an orientable $2$-orbifold without boundary, let $\oldUpsilon$ be a negative, finite-type,
orientable $2$-orbifold, and let $\otheroldLambda_0$ be a negative, taut,
finite-type $2$-suborbifold of $\otheroldLambda$. 
Let $\xi,\xi':\oldUpsilon\to\otheroldLambda_0$ be proper
$\pi_1$-injective embeddings which are properly isotopic in
$\otheroldLambda$. Then $\xi$ and $\xi'$ are properly isotopic in
$\otheroldLambda_0$. 
\NoProof
\EndCorollary

\Proposition\label{unique embedding}
Let $\otheroldLambda$ be a negative, closed, orientable $2$-orbifold. Let
$\oldUpsilon_1$ and $\oldUpsilon_2$ be elements of $\Theta_-(\otheroldLambda$), such that $[\oldUpsilon_1]\preceq[\oldUpsilon_2]$. Then there is an embedding $\xi:\oldUpsilon_1\to\oldUpsilon_2$, unique up to isotopy in $\oldUpsilon_2$, such that $\xi$, regarded as an embedding of $\oldUpsilon_1$ in $\otheroldLambda$, is isotopic to the inclusion  map $\oldUpsilon_1\to\otheroldLambda$.
\EndProposition

\Proof
The existence of $\xi$ is immediate from the hypothesis $[\oldUpsilon_1]\preceq[\oldUpsilon_2]$. To prove uniqueness, we must show that if
two embeddings $\xi,\xi':\oldUpsilon_1\to\oldUpsilon_2$, are isotopic in $\otheroldLambda$ then they are isotopic in $\oldUpsilon_2$. This follows from Corollary \ref{i guess} if we let $\oldUpsilon_1$ and $\oldUpsilon_2$ play the respective roles of $\oldUpsilon$ and $\otheroldLambda_0$ in that lemma.
\EndProof
%\oldXi\frakX

\Number\label{my remark}
Suppose that $\otheroldLambda_1$ and $\otheroldLambda_2$ are negative, closed,
orientable $2$-orbifolds. We will denote by
$\calm=\calm(\otheroldLambda_1,\otheroldLambda_2)$ the set of all triples
$(\oldUpsilon_1,\newf ,\oldUpsilon_2)$, where $\oldUpsilon_i$ is an
element of $\Theta_-(\otheroldLambda_i)$ for $i=1,2$, and $\newf
:\oldUpsilon_1\to\oldUpsilon_2$ is a homeomorphism. We define an
equivalence relation $
\equiv_{\otheroldLambda_1,\otheroldLambda_2}$ on the set $\calm$ (to be denoted simply by
$\equiv$ when the $\otheroldLambda_i$ are understood) 
 by stipulating that $(\oldUpsilon_1,\newf ,\oldUpsilon_2)\equiv (\oldUpsilon_1',\newf ',\oldUpsilon_2')$ if and only if there are homeomorphisms $\eta_i:\oldUpsilon_i\to\oldUpsilon_i'$ for $i=1,2$ such that (i) $\eta_i$, regarded as an embedding of $\oldUpsilon_i$ in $\otheroldLambda_i$, is isotopic to the inclusion, and (ii) the homeomorphisms $\newf '$ and $\eta_2\circ \newf \circ\eta_1^{-1}$ from $\oldUpsilon_1'$ to $\oldUpsilon_2'$ are isotopic. The set of equivalence classes under the relation $\equiv$ on $\calm$ will be denoted by  $\barcalm=\barcalm(\otheroldLambda_1,\otheroldLambda_2)$. The equivalence class of an element 
$(\oldUpsilon_1,\newf ,\oldUpsilon_2)$ of $\calm$ will be denoted by
$[\oldUpsilon_1,\newf ,\oldUpsilon_2]\in\barcalm$, or simply by $[\newf ]$ when there
is no ambiguity. 
%f

If $(\oldUpsilon_1,\newf,\oldUpsilon_2)\in\calm$ is given, then since $\oldUpsilon_i\in\Theta_-(\otheroldLambda)$, we have $Y_i:=[\oldUpsilon_i]\in\barcaly_-(\otheroldLambda_i)$ for $i=1,2$. It follows from the definition of the equivalence relation $\equiv_{\otheroldLambda_1,\otheroldLambda_2}$ that $Y_1$ and $Y_2$ depend only on the equivalence class $\emm:=[\oldUpsilon_1,\newf,\oldUpsilon_2]$ of $(\oldUpsilon_1,\newf,\oldUpsilon_2)$. We will denote $Y_1$ and $Y_2$ by $\dom \emm$ and $\range \emm$ respectively. 

If $(\oldUpsilon_1,\newf,\oldUpsilon_2)\in\calm$ is given, then
$(\oldUpsilon_2,\newf^{-1},\oldUpsilon_1)$ is also an element of $\calm$. It is
clear that $[\oldUpsilon_2,\newf^{-1},\oldUpsilon_1]\in\barcalm$ depends only on the
equivalence class $\emm=[\oldUpsilon_1,\newf,\oldUpsilon_2]\in\barcalm$. We will
denote the element $[\oldUpsilon_2,\newf^{-1},\oldUpsilon_2]\in\barcalm$ by
$\emm^{-1}$. Note that for any $\emm\in\barcalm$ we have $\dom \emm^{-1}=\range
\emm$, $\range \emm^{-1}=\dom \emm$, and $(\emm^{-1})^{-1}=\emm$.

If $\otheroldLambda$ is a If closed, orientable $2$-orbifold, and
$\oldUpsilon$ is an element of $\Theta_-(\otheroldLambda)$,
then
$(\oldUpsilon,{\id}_\oldUpsilon,\oldUpsilon)$ is also an element of $\calm(\otheroldLambda,\otheroldLambda)$, and it is
clear that $[\oldUpsilon,{\id}_\oldUpsilon,\oldUpsilon]\in\barcalm$ depends only on the
isotopy class $Y:=[\oldUpsilon]\in\barcaly_-(\otheroldLambda)$. We will
denote the element $[\oldUpsilon,\id_\oldUpsilon,\oldUpsilon]\in\barcalm(\otheroldLambda,\otheroldLambda)$ by
$1_Y$. Note that for any $Y\in\barcaly_-(\otheroldLambda)$ we have $\dom 1_Y=\range1_Y=Y$.

\EndNumber

\Number\label{zoltan}
Suppose that $\otheroldLambda_1$, $\otheroldLambda_2$ and $\otheroldLambda_3$ are negative, closed, orientable $2$-orbifolds, and that $\emm_1\in\barcalm(\otheroldLambda_1, \otheroldLambda_2)$ and $\emm_2\in\barcalm(\otheroldLambda_2$, $\otheroldLambda_3)$ are elements such that $\range\emm_1\preceq\dom\emm_2$. Let us choose elements $(\oldUpsilon_1,\newf_1,\oldUpsilon_2)\in\calm(\otheroldLambda_1,\otheroldLambda_2)$ and $(\oldUpsilon_2',\newf_2,\oldUpsilon_3)\in\calm(\otheroldLambda_2,\otheroldLambda_3)$ such that
$[\oldUpsilon_1,\newf_1,\oldUpsilon_2]=\emm_1$ and $[\oldUpsilon_2',\newf_2,\oldUpsilon_3]=\emm_2$. The hypothesis $\range\emm_1\preceq\dom\emm_2$ means that $[\oldUpsilon_2]\preceq[\oldUpsilon_2']$. Hence
 by Proposition \ref{unique embedding}
there is an embedding $\xi:\oldUpsilon_2\to\oldUpsilon_2'$, unique up to isotopy in $\oldUpsilon_2'$, such that $\xi$, regarded as an embedding of $\oldUpsilon_2$ in $\otheroldLambda_2$, is isotopic to the inclusion  map $\oldUpsilon_2\to\otheroldLambda_2$. 
Then $\newf_2\circ\xi\circ \newf_1$ is a well-defined embedding of $\oldUpsilon_1$ in $\oldUpsilon_3$, and the uniqueness property of $\xi$ implies that the element $[\oldUpsilon_1,\newf_2\circ\xi\circ \newf_1, 
\newf_2(\xi(\oldUpsilon_2))]$
%\newf_2\circ\xi\circ \newf_1(\oldUpsilon_1)]$ \redcomment{Could I write 
%$\newf_2\circ\xi(\oldUpsilon_2)$ in place of $\newf_2\circ\xi\circ %\newf_1(\oldUpsilon_1)$? It seems a little better}
of $\barcalm(\otheroldLambda_1,\otheroldLambda_3)$ depends only on $\emm_1$ and $\emm_2$. It will be denoted by $\emm_2\circ\emm_1$. Note that $\dom(\emm_2\circ\emm_1)=\dom\emm_1$ and that $\range(\emm_2\circ\emm_1)\preceq\range\emm_2$.

In the special case where $\dom\emm_2=\range\emm_1$, we may take $\oldUpsilon_2=\oldUpsilon_1$ in the above definition, and we may take $\xi$ to be the identity map of $\oldUpsilon_1$. It follows that in this case we have $(\emm_2\circ\emm_1)^{-1}=\emm_1^{-1}\circ\emm_2^{-1}$.

If the elements $\emm_1,\emm_2$ of $\barcalm(\otheroldLambda_1,\otheroldLambda_2)$ and $\barcalm(\otheroldLambda_2,\otheroldLambda_3)$ are given in the form 
$\emm_1=[\oldUpsilon_1,\newf_1,\oldUpsilon_2]$ and
$\emm_2=[\oldUpsilon_2',\newf_2,\oldUpsilon_3]$, where the elements 
$(\oldUpsilon_1,\newf_1,\oldUpsilon_2)\in\calm(\otheroldLambda_1,\otheroldLambda_2)$ and
$(\oldUpsilon_2',\newf_2,\oldUpsilon_3)\in\calm(\otheroldLambda_2,\otheroldLambda_3)$ satisfy
 $\oldUpsilon_2\subset\oldUpsilon_2'$, then in the construction given in
 \ref{zoltan} we may take $\xi:\oldUpsilon_2\to\oldUpsilon_2'$ to be the
 inclusion map. 
Hence in this case we have
 $\emm_2\circ\emm_1=
[\oldUpsilon_1,\newf_2\circ \newf_1,\newf_2(\oldUpsilon_2)]$.
\EndNumber

%%[\oldUpsilon_1,\newf_2\circ \newf_1,\newf_2\circ \newf_1(\oldUpsilon_1)]$. \redcomment{Or

\Number\label{more zoltan}
Suppose that $\otheroldLambda_1$, $\otheroldLambda_2$, $\otheroldLambda_3$ and $\otheroldLambda_4$ are
negative, closed, orientable $2$-orbifolds, that $\emm_i$ is an
element of $\barcalm(\otheroldLambda_i, \otheroldLambda_{i+1})$ for $i=1,2,3$, and
that $\range\emm_i\preceq\dom\emm_{i+1}$ for $i=1,2$. Then
$
\emm_2\circ\emm_1$ and $\emm_3\circ\emm_2$ are defined, and we have
$\range(\emm_2\circ\emm_1)\preceq\range\emm_2\preceq\dom\emm_3$ and
$\range\emm_1\preceq\dom\emm_2=\dom(\emm_3\circ\emm_2)$, so that
$\emm_3\circ(\emm_2\circ\emm_1)$ and $(\emm_3\circ\emm_2)\circ\emm_1$
are defined.
 We claim that
\Equation\label{easy associativity}
\emm_3\circ(\emm_2\circ\emm_1)=(\emm_3\circ\emm_2)\circ\emm_1.
\EndEquation

To prove this, for $i=1,2,3$ we choose elements
$(\oldUpsilon_i,\newf_1,\oldUpsilon_2)\in\calm(\otheroldLambda_1,\otheroldLambda_2)$,
$(\oldUpsilon_2',\newf_2,\oldUpsilon_3)\in\calm(\otheroldLambda_2,\otheroldLambda_3)$, and 
$(\oldUpsilon_3',\newf_3,\oldUpsilon_4)\in\calm(\otheroldLambda_3,\otheroldLambda_4)$ such that
$[\oldUpsilon_1,\newf_1,\oldUpsilon_2]=\emm_1$, $[\oldUpsilon_2',\newf_2,\oldUpsilon_3=\emm_2$, and 
$[\oldUpsilon_3',\newf_3,\oldUpsilon_4]=\emm_3$.
The hypothesis
$\range\emm_i\preceq\dom\emm_{i+1}$ for $i=1,2$ means that $[\oldUpsilon_i]\preceq[\oldUpsilon_i']$ for $i=2,3$. Hence
 by Proposition \ref{unique embedding}, for $i=2,3$
there exists an embedding $\xi_i:\oldUpsilon_i\to\oldUpsilon_{i}'$, unique up to isotopy in $\oldUpsilon_i'$, such that $\xi_i$, regarded as an embedding of $\oldUpsilon_i$ in $\otheroldLambda_i$, is isotopic to the inclusion  map $\oldUpsilon_i\to\otheroldLambda_i$. By the definition given in \ref{zoltan}, we have 
$$\emm_3\circ(\emm_2\circ\emm_1)=
[\newf_3\circ\xi_3\circ(\newf_2\circ\xi_2\circ \newf_1)]
=[(\newf_3\circ\xi_3\circ \newf_2)\circ\xi_2\circ \newf_1]
=
(\emm_3\circ\emm_2)\circ\emm_1,$$
which proves \ref{easy associativity}.
\EndNumber

\Number\label{restriction}
Let $\otheroldLambda_1$ and $\otheroldLambda_2$ be negative, closed, orientable $2$-orbifolds, and
suppose that we are given an element $\emm$ of
$\barcalm(\otheroldLambda_1,\otheroldLambda_2)$, and an element $Y$ of
$\barcaly_-(\otheroldLambda_1)$ with $Y\preceq\dom\emm$. Since $\range 1_Y=Y$,
it follows from \ref{zoltan} that $\emm\circ1_Y$ is a well-defined
element of $\barcalm(\otheroldLambda)$. It will be denoted by
$\emm|Y$. According to \ref{zoltan} we have
$\dom(\emm|Y)=\dom1_Y=Y$, and
$\range(\emm|Y)\preceq\range\emm$. The element $\range(\emm|Y)$ of $\Theta_-(\otheroldLambda_1)$
will be denoted by $\emm(Y)$.

Note that $\emm\,|\dom\emm=\emm$.

If we write $\emm=[\oldUpsilon_1,\newf,\oldUpsilon_2]$ for some
$(\oldUpsilon_1,\newf,\oldUpsilon_2)\in\calm(\otheroldLambda_1,\otheroldLambda_2)$, and if
$Y\in\barcaly_-(\otheroldLambda_1)$ is given in the form $[\oldUpsilon]$ for some $\oldUpsilon\in\Theta_-(\otheroldLambda_1)$
%compact suborbifold $\oldUpsilon$ of $\otheroldLambda_1$ 
such that
$\oldUpsilon\subset\oldUpsilon_1$, then it follows from \ref{zoltan} (with $1_Y$ and $\emm$ playing the roles of $\emm_1$ and
$\emm_2$) that
$\emm|Y=[\oldUpsilon,\newf|\oldUpsilon,\newf(\oldUpsilon)]$. Hence
$\emm(Y)=[\newf(\oldUpsilon)]$.

If we fix negative, orientable,
compact orbifolds $\otheroldLambda_1$ and $\otheroldLambda_2$, and an element $\emm$ of
$\barcalm(\otheroldLambda_1,\otheroldLambda_2)$, the assignment $Y\mapsto \emm(Y)$ is an order-preserving mapping from the
set $\{Y:Y\preceq\dom\emm\}$ to the set
$\{Z:Z\preceq\range\emm\}$. This is because, if $Y\preceq
Y'\preceq\dom\emm$, we can write $Y=[\oldUpsilon]$ and $Y'=[\oldUpsilon']$
for some $\oldUpsilon,\oldUpsilon'\in\Theta_-(\otheroldLambda)$ with
$\oldUpsilon\subset\oldUpsilon'\subset\oldUpsilon_1$. We then have
$\newf(\oldUpsilon)\subset \newf(\oldUpsilon')$ and hence
$\emm(Y)=[\newf(\oldUpsilon)]\preceq[ \newf(\oldUpsilon')]=\emm(Y')$. Note also
that this order-preserving map is bijective, as it has the two-sided
inverse $Z\mapsto \emm^{-1}(Z)$. 
\EndNumber

\Number\label{special composition}
Let $\otheroldLambda_1$, $\otheroldLambda_2$ and $\otheroldLambda_3$  be
negative, closed, orientable $2$-orbifolds, let $\emm_1\in
\barcalm(\otheroldLambda_1, \otheroldLambda_{2})$ and $\emm_2
\in\barcalm(\otheroldLambda_2, \otheroldLambda_{3})$ be given. Suppose that $\range\emm_1\preceq\dom\emm_2$, so that $\emm_2\circ\emm_1$ is defined. Then
for 
any element $Y$  of
$\barcaly_-(\otheroldLambda_1)$, we have 
\Equation\label{that's why}
\emm_2\circ\emm_1(Y)=\emm_2(\emm_1(Y)).
\EndEquation
 (Indeed, if we write $Y=[\oldUpsilon]$ for some 
$\oldUpsilon\in\Theta_-(\otheroldLambda_1)$,
 and if for $i=1,2$ we write
 $\emm_i=[\oldUpsilon_i,\newf_i,\oldUpsilon_{i+1}]$
for some $(\oldUpsilon_i,\newf_i,\oldUpsilon_{i+1})\in
\calm(\otheroldLambda_i, \otheroldLambda_{i+1})$, then
$\emm_2\circ\emm_1(Y)=[\newf_2\circ\newf_1(\oldUpsilon)]=[\newf_2(\newf_1(\oldUpsilon))]=\emm_2([\newf_1(\oldUpsilon)]=
\emm_2(\emm_1(Y))$.)
\EndNumber

\Number\label{image and restriction}
We will sometimes encounter a more general situation than the one considered in \ref{special composition}. Let $\otheroldLambda_1$, $\otheroldLambda_2$ and $\otheroldLambda_3$  be
negative, closed, orientable $2$-orbifolds, let $\emm_1\in
\barcalm(\otheroldLambda_1, \otheroldLambda_{2})$ and $\emm_2
\in\barcalm(\otheroldLambda_2, \otheroldLambda_{3})$ be given,
%elements such
%that $\range\emm_1\preceq\dom\emm_{2}$, 
and let 
$Y$ be
an element of
$\barcaly_-(\otheroldLambda_1)$ such that $Y\preceq\dom\emm_1$ and $\emm_1(Y_1)\preceq\dom\emm_2$. 
%According to
%\ref{easy associativity}, we have
%$(\emm_2\circ\emm_1)\circ1_Y=\emm_2\circ(\emm_1\circ1_Y)$, i.e.
%\Equation\label{without range}
%$$(\emm_2\circ\emm_1)|Y=\emm_2\circ(\emm_1|Y).$$
%\EndEquation
%In particular, the two sides of (\ref{first equation}) have the same
%range, i.e.
%\Equation\label{within range}
We claim that
\Equation\label{number one}
\emm_2\circ(\emm_1|Y)=(\emm_2|\emm_1(Y))\circ(\emm_1|Y)
\EndEquation
and
\Equation\label{number three}
\range(\emm_2\circ(\emm_1|Y))=\emm_2(\emm_1(Y)).
\EndEquation
To prove (\ref{number one}) and (\ref{number three}), note that since $Y\preceq\dom\emm_1$, we may write $Y=\oldUpsilon$ and $\emm_1= [\oldUpsilon_1,\newf_1,\oldUpsilon_2]$, where $\oldUpsilon$ is an element of $\Theta_-(\otheroldLambda_1$), and $(\oldUpsilon_1,\newf_1,\oldUpsilon_2)$ is an element of $\calm(\otheroldLambda_1,\otheroldLambda_2)$ such that $\oldUpsilon_1\supset\oldUpsilon$. Since $\oldUpsilon\subset\oldUpsilon_1$, we have $\emm_1|Y=[\oldUpsilon,\newf_1|\oldUpsilon,\newf_1(\oldUpsilon)]$, and hence $\emm_1(Y)=[\newf_1(\oldUpsilon)]$. 
%\redmissingref{Wasn't this proved in the preceding subsection?} 
Since
$\emm_1(Y)\preceq\dom(\emm_2)$, we may write
$\emm_2= [\oldUpsilon_2',\newf_2,\oldUpsilon_3]$, where $(\oldUpsilon_2',\newf_2,\oldUpsilon_3)$ is an element of $\calm(\otheroldLambda_2,\otheroldLambda_3)$ such that $\oldUpsilon_2'\supset \newf_1(\oldUpsilon_1)$. Since $
\newf_1(\oldUpsilon)\subset\newf_1(\oldUpsilon_1)\subset
\oldUpsilon_2' $, we have 
\Equation\label{another number}
\emm_2|\emm_1(Y)=[\newf_1(\oldUpsilon),\newf_2|\newf_1(\oldUpsilon)
%\emm_1(\oldUpsilon_1)
,\newf_2(\newf_1(\oldUpsilon))].
\EndEquation
Hence
\Equation
\label{yeah, this one too}
\begin{aligned}
(\emm_2|\emm_1(Y))\circ(\emm_1|Y)&=[\oldUpsilon,(\newf_2|\newf_1(\oldUpsilon ))\circ (\newf_1|\oldUpsilon ),\newf_2(\newf_1(\oldUpsilon))]\\
&=[\oldUpsilon,\newf_2
%\newf_1(\oldUpsilon ))
\circ (\newf_1|\oldUpsilon ),\newf_2(\newf_1(\oldUpsilon))]=
\emm_2\circ(\emm_1|Y),
\end{aligned}
\EndEquation
which proves (\ref{number one}). To prove (\ref{number three}), note that by definition we have $\emm_2(\emm_1(Y))=\range(\emm_2|\emm_1(Y))$, so that (\ref{another number}) gives $\emm_2(\emm_1(Y))
=[\newf_2(\newf_1(\oldUpsilon))]$. But the final equality of (\ref{yeah, this one too}) shows that $[\newf_2(\newf_1(\oldUpsilon_1))]=
\range(\emm_2\circ(\emm_1|Y))$, and hence
$\emm_2(\emm_1(Y))=\range(\emm_2\circ(\emm_1|Y))$, as
required. 
%\redcomment{I think this all makes sense given the observations
  %in the preceding subsection. Maybe I should emphasize that I'm using
  %them freely. Or is that silly because I just said all those things?}

We also record the special case of (\ref{number one}) in which $Y=\dom\emm_1$: if $\emm_1\in
\barcalm(\otheroldLambda_1, \otheroldLambda_{2})$ and $\emm_2
\in\barcalm(\otheroldLambda_2, \otheroldLambda_{3})$ are elements such that $\range\emm_1\preceq\dom\emm_2$, then 
\Equation\label{lunch}
\emm_2\circ\emm_1=(\emm_2|\range\emm_1)\circ\emm_1.
\EndEquation
\EndNumber

\Number
Now suppose that $\otheroldLambda_1$, $\otheroldLambda_2$ and $\otheroldLambda_3$ are negative, closed, orientable $2$-orbifolds, and that arbitrary elements $\emm_1\in\barcalm(\otheroldLambda_1, \otheroldLambda_2)$ and $\emm_2\in\barcalm(\otheroldLambda_2$, $\otheroldLambda_3)$ are given. 
By \ref{my remark}, $\emm_1^{-1}$ is a well-defined element of $\barcalm(\otheroldLambda_2,\otheroldLambda_1)$, and $\dom\emm_1^{-1}=\range\emm_1$. Since $\range\emm_1\wedge\dom\emm_2\preceq\range\emm_1$, \ref{restriction} gives a well-defined element $\emm_1^{-1}|(\range\emm_1\wedge\dom\emm_2)$ of $\barcalm(\otheroldLambda_2,\otheroldLambda_1)$, with
$\dom(\emm_1^{-1}|(\range\emm_1\wedge\dom\emm_2))
=
\range\emm_1\wedge\dom\emm_2$. Another application of \ref{my remark} gives an element 
$(\emm_1^{-1}|(\range\emm_1\wedge\dom\emm_2))^{-1}$ of $\barcalm(\otheroldLambda_1,\otheroldLambda_2)$, with
$\range((\emm_1^{-1}|(\range\emm_1\wedge\dom\emm_2))^{-1})
=\dom (\emm_1^{-1}|(\range\emm_1\wedge\dom\emm_2))
=\range\emm_1\wedge\dom\emm_2\preceq\dom\emm_2$. Now \ref{zoltan} gives an element $\emm_2\circ(\emm_1^{-1}|(\range\emm_1\wedge\dom\emm_2))^{-1}$ of $\barcalm(\otheroldLambda_1,\otheroldLambda_3)$. We define
\Equation\label{display it}
\emm_2\diamond\emm_1=\emm_2\circ(\emm_1^{-1}|(\range\emm_1\wedge\dom\emm_2))^{-1}
\EndEquation
for any $\emm_1\in\barcalm(\otheroldLambda_1,\otheroldLambda_2)$ and
$\emm_1\in\barcalm(\otheroldLambda_2,\otheroldLambda_3)$.

According to \ref{zoltan} we have
$$
%\begin{aligned}
\dom(\emm_2\diamond\emm_1)
%&
=
\dom(\emm_2\circ(\emm_1^{-1}|(\range\emm_1\wedge\dom\emm_2))^{-1})
%\preceq
=\dom((\emm_1^{-1}|(\range\emm_1\wedge\dom\emm_2))^{-1}),
$$
%\\&
which by \ref{my remark} means that $\dom(\emm_2\diamond\emm_1)=\range (\emm_1^{-1}|(\range\emm_1\wedge\dom\emm_2))$. By \ref{restriction} this may be written in the form
\Equation\label{yes oi do}
\dom(\emm_2\diamond\emm_1)= \emm_1^{-1}(\range\emm_1\wedge\dom\emm_2)
%\preceq\range (\emm_1^{-1})
.
%\end{aligned}$$
\EndEquation
%
%so that
%\Equation\label{what you thought}
%\dom(\emm_2\diamond\emm_1)\preceq\dom\emm_1.
%\EndEquation
\EndNumber

\Number\label{return of zoltan}
Suppose that $\otheroldLambda_1$, $\otheroldLambda_2$ and $\otheroldLambda_3$ are negative, closed, orientable $2$-orbifolds, and that $\emm_1\in\barcalm(\otheroldLambda_1, \otheroldLambda_2)$ and $\emm_2\in\barcalm(\otheroldLambda_2$, $\otheroldLambda_3)$ are elements such that $\range\emm_1\preceq\dom\emm_2$; thus by \ref{zoltan}, 
$\newf_2\circ \newf_1$ is defined. In this case, we have $\newf_2\diamond \newf_1=\newf_2\circ \newf_1$. (Indeed, 
since $\range\emm_1\preceq\dom\emm_2$, we have
$\range\emm_1\wedge\dom\emm_2=\range\emm_1$; thus 
$\emm_1^{-1}|(\range\emm_1\wedge\dom\emm_2))^{-1}=
\emm_1^{-1}|(\range\emm_1)=\emm_1^{-1}$, and hence
$\emm_2\diamond\emm_1=\emm_2\circ(\emm_1^{-1}|(\range\emm_1\wedge\dom\emm_2))^{-1}=\emm_2\circ(\emm_1^{-1})^{-1}=\emm_2\circ\emm_1
$.)
\EndNumber

 \Lemma\label{before associativity}Let $\otheroldLambda_1$, $\otheroldLambda_2$ and $\otheroldLambda_3$ be negative, closed, orientable $2$-orbifolds, and let $\emm_1\in\barcalm(\otheroldLambda_1, \otheroldLambda_2)$ and $\emm_2\in\barcalm(\otheroldLambda_2$, $\otheroldLambda_3)$ be given. Let $Y$ be an element of $\barcaly_-(\otheroldLambda_1)$. Then the following conditions are equivalent:
 \begin{enumerate}
 \item$Y\preceq\dom(\emm_2\diamond\emm_1)$.
 \item $Y\preceq\dom\emm_1$, and $\emm_1(Y)\preceq\dom\emm_2$ (where $\emm(Y)$ is defined as in \ref{restriction}).
 \end{enumerate}
 Furthermore, if (1) (or (2)) holds, then $(\emm_2\diamond\emm_1)|Y=\emm_2\circ(\emm_1|Y)$ (where  $\emm_2\circ(\emm_1|Y)$ is defined because $\range(\emm_1|Y)=\emm_1(Y)\preceq\dom\emm_1$ by (2)).
 \EndLemma

 \Proof
Choose an element $(\oldUpsilon_1,\newf_1,\oldUpsilon_2)$ of $\barcalm(\otheroldLambda_1,\otheroldLambda_2)$ such that $[\oldUpsilon_1,\newf_1,\oldUpsilon_2]=\emm_1$. Then $[\oldUpsilon_2]=\range\emm_1$. According to \ref{my remark}, we have
$[\oldUpsilon_2,\newf_1^{-1},\oldUpsilon_1]=\emm_1^{-1}$. 
Since $\range\emm_1\wedge\dom\emm_2\preceq\range\emm_1$, there is a suborbifold $\oldGamma$ of $\oldUpsilon_2$ such that $[\oldGamma]=\range\emm_1\wedge\dom\emm_2$. According to \ref{restriction} we have 
\Equation\label{is this it?}
\emm_1^{-1}|(\range\emm_1\wedge\dom\emm_2)=[\oldGamma,\newf_1^{-1}|\oldGamma,\newf_1^{-1}(\oldGamma)].
\EndEquation
Equating the ranges of the two sides of (\ref{is this it?}) gives
$\emm_1^{-1}(\range\emm_1\wedge\dom\emm_2)=
[\newf_1^{-1}(\oldGamma)]$. With (\ref{yes oi do}) this gives
\Equation\label{how about you?}
 \dom(\emm_2\diamond\emm_1)=[\newf_1^{-1}(\oldGamma)].
\EndEquation

Now suppose that (1) holds. Then by (\ref{how about you?}) we have $Y\preceq[\newf_1^{-1}(\oldGamma)]$. Hence there is a suborbifold $\oldUpsilon$ of $\newf_1^{-1}(\oldGamma)$ such that $\oldUpsilon\in\Theta_-(\otheroldLambda)$ and $[\oldUpsilon]=Y$. Since $\newf_1^{-1}(\oldGamma)\subset\oldUpsilon_1$, we have $Y\preceq[\oldUpsilon_1]=\dom\emm_1$. Furthermore, according to \ref{restriction}, we have $\emm_1|Y=[\oldUpsilon,\newf_1|\oldUpsilon,\newf_1(\oldUpsilon)]$, so that $\emm_1(Y)=\range(\emm_1|Y)=[\newf_1(\oldUpsilon)]$. Since
$\newf_1(\oldUpsilon)\subset\oldGamma
%\subset\oldUpsilon_2
$, we have $\emm_1(Y)\preceq[\oldGamma]=
\range\emm_1\wedge\dom\emm_2\preceq
\dom\emm_2$, which gives (2).

Continuing to assume that (1) holds, let us prove that $(\emm_2\diamond\emm_1)|Y=\emm_2\circ(\emm_1|Y)$. By (\ref{display it}) we have
$\emm_2\diamond\emm_1=\emm_2\circ(\emm_1^{-1}|(\range\emm_1\wedge\dom\emm_2))^{-1}$. But since $\emm_1^{-1}|(\range\emm_1\wedge\dom\emm_2)=[\oldGamma,\newf_1^{-1}|\oldGamma,\newf_1^{-1}(\oldGamma)]$ by (\ref{is this it?}), it follows from the definitions in \ref{my remark} that
$(\emm_1^{-1}|(\range\emm_1\wedge\dom\emm_2))^{-1}=
[\newf_1^{-1}(\oldGamma),(\newf_1^{-1}|\oldGamma)^{-1},\oldGamma]=
[\newf_1^{-1}(\oldGamma),\newf_1|\newf_1^{-1}(\oldGamma),\oldGamma]$.
On the other hand, since 
%$\dom\emm_2=[\oldUpsilon_2]$,
$[\oldGamma]=
\range\emm_1\wedge\dom\emm_2\preceq
\dom\emm_2$, we may write $\dom\emm_2$ in the form $[\oldUpsilon_2']$ for
some $\oldUpsilon_2'\in\Theta_-(\otheroldLambda_2)$ such that $\oldUpsilon_2'\supset\oldGamma$. We may then write
$\emm_2$ as $[\oldUpsilon'_2,\newf_2,\oldUpsilon_3]$, where $\newf_2$ is a
homeomorphism from $\oldUpsilon_2$ to some $\oldUpsilon_3\in\Theta_-(\otheroldLambda_3)$. We now have
$$
%\begin{aligned}
\emm_2\diamond\emm_1
%&
=\emm_2\circ(\emm_1^{-1}|(\range\emm_1\wedge\dom\emm_2))^{-1}
%\\
%&
=[\oldUpsilon'_2,\newf_2,\oldUpsilon_3]\circ[\newf_1^{-1}(\oldGamma),\newf_1|\newf_1^{-1}(\oldGamma),\oldGamma].
%\end{aligned}
$$

But since $
%\newf_1(\newf_1^{-1}(\oldGamma)=
\oldGamma\subset\oldUpsilon_2'$, it follows from
\ref{zoltan} that 
$[\oldUpsilon'_2,\newf_2,\oldUpsilon_3]\circ[\newf_1^{-1}(\oldGamma),\newf_1|\newf_1^{-1}(\oldGamma),\oldGamma]=[
\newf_1^{-1}(\oldGamma),\newf_2\circ (\newf_1|\newf_1^{-1}(\oldGamma)),\newf_2(\oldGamma)]$.
Hence $\emm_2\diamond\emm_1=[
\newf_1^{-1}(\oldGamma),\newf_2\circ (\newf_1|\newf_1^{-1}(\oldGamma)),\newf_2(\oldGamma)]$, and
since $\oldUpsilon\subset \newf_1^{-1}(\oldGamma)$, \ref{restriction}
gives 
$$(\emm_2\diamond\emm_1)|Y=[
\oldUpsilon,\newf_2\circ (\newf_1|\oldUpsilon),\newf_2( \newf_1(\oldUpsilon)]=\emm_2\circ(\emm_1|Y),$$
as required.

It remains to prove that (2) implies (1). If (2) holds then
$Y\preceq\dom\emm_1=[\oldUpsilon_1]$, and hence $Y=[\oldUpsilon]$ for
some $\oldUpsilon\in\Theta_-(\otheroldLambda_1)$ with $\oldUpsilon\subset\oldUpsilon_1$. By \ref{restriction},
since $\oldUpsilon\subset\oldUpsilon_1$, we
have $\emm_1|Y=[\oldUpsilon,\newf_1|\oldUpsilon,\newf_1(\oldUpsilon)]$, and hence
$\emm_1(Y)=\range(\emm_1|Y)=[\newf_1(\oldUpsilon)]$. According to (2) we have
$\emm_1(Y)\preceq\dom\emm_2$, and hence
$[\newf_1(\oldUpsilon)]\preceq\dom\emm_2$. But since $\newf_1(\oldUpsilon)\subset
\newf_1(\oldUpsilon_1)$, we have $[\newf_1(\oldUpsilon)]\preceq
[\newf_1(\oldUpsilon_1)]=\range\emm_1$. Since $\range\emm_1\wedge\dom\emm_2$
is by definition the infimum of $\range\emm_1$ and $\dom\emm_2$, it
follows that $[\newf_1(\oldUpsilon)]\preceq
\range\emm_1\wedge\dom\emm_2$. 
%Hence there is an taut, negative,
%compact suborbifold $\oldGamma'$ of $\otheroldLambda_2$ such that
%$[\oldGamma']=\range\emm_1\wedge\dom\emm_2$ and $\oldGamma'\supset
%\newf_1(\oldUpsilon)$. Thus $Y=[\oldUpsilon]\preceq [\newf_1^{-1}(\oldGamma')]$. 
%But since $\oldGamma^{-1}\subset it follows from \ref{restrict} that $\emm_1^{-1}|[\Ga%%mma']=
This means that
$\newf_1(\oldUpsilon)$ is isotopic in $\otheroldLambda_2$ to a suborbifold $\oldDelta$ of
$\oldGamma$. We have $\oldDelta\subset\oldGamma\subset\oldUpsilon_2$
%\redcomment{That seems wrong. I haven't said that
  %$\oldGamma\subset\oldUpsilon_2$. In fact in the step before this
  %one, I chose a representative $\oldUpsilon_2'$ of $\dom\emm_2$ such
  %that $\oldGamma\subset\oldUpsilon_2'$. I may need to do the same
  %thing here} 
and
$\newf_1(\oldUpsilon)\subset \newf_1(\oldUpsilon_1)=\oldUpsilon_2$. Since 
$\newf_1(\oldUpsilon)$ and $\oldDelta$ are isotopic in $\otheroldLambda_2$, and are both contained in the taut,
negative compact suborbifold  $\oldUpsilon_2$ of $\otheroldLambda_2$, it follows
from Corollary \ref{i guess} that they are isotopic in $\oldUpsilon_2$. Since $\newf_1$
is a homeomorphism from $\oldUpsilon_1$ to $\oldUpsilon_2$, it now follows
that $\oldUpsilon$ and $\newf_1^{-1}(\oldDelta)$ are isotopic in $\oldUpsilon_1$,
and therefore in $\otheroldLambda_1$. Hence
$Y=[\oldUpsilon]=[\newf_1^{-1}(\oldDelta)]\preceq[\newf_1^{-1}(\oldGamma)]$. But by
\ref{restriction} we have 
$\emm_1^{-1}|[\oldGamma]=[\oldGamma,\newf_1^{-1}|\oldGamma, \newf_1^{-1}(\oldGamma)]$, so
that $\range(\emm_1^{-1}|[\oldGamma])=[ \newf_1^{-1}(\oldGamma)]$. Thus
$Y\le\range(\emm_1^{-1}|[\oldGamma])=\range(\emm_1^{-1}(
\range\emm_1\wedge\dom\emm_2))$. According to (\ref
{yes oi do}), this means 
$Y\le\dom(\emm_2\diamond\emm_1)$, which is Condition (1).
 \EndProof

\Proposition\label{associativity}
Suppose that $\otheroldLambda_1$, $\otheroldLambda_2$  and $\otheroldLambda_3$ are
negative, closed, orientable $2$-orbifolds, and that
$\emm_i\in\barcalm(\otheroldLambda_i, \otheroldLambda_{i+1})$ is given for
$i=1,2$. Then
$(\emm_2\diamond\emm_1)^{-1}=\emm_1^{-1}\diamond\emm_1^{-1}$. Furthermore,
if  $\otheroldLambda_4$ is another negative, closed, orientable $2$-orbifold, and
$\emm_3\in\barcalm(\otheroldLambda_3, \otheroldLambda_{4})$ is given, then
$(\emm_3\diamond\emm_2)\diamond\emm_1=\emm_3\diamond(\emm_2\diamond\emm_1)$.
\EndProposition

\Proof

We prove the second assertion first.
Set $Y=\dom(\emm_3\diamond(\emm_2\diamond\emm_1))$ and $Y'=\dom((\emm_3\diamond\emm_2)\diamond\emm_1)$.
According to Lemma \ref{before associativity} we have $Y\preceq\dom(\emm_2\diamond\emm_1)$ and
$(\emm_2\diamond\emm_1)(Y)\preceq\dom\emm_3$; and in addition
$\emm_3\circ((\emm_2\diamond\emm_1)|Y)=(\emm_3\diamond(\emm_2\diamond\emm_1))|Y$,
which in view of \ref{restriction} means that 
$\emm_3\circ((\emm_2\diamond\emm_1)|Y)=\emm_3\diamond(\emm_2\diamond\emm_1)$.
Since $Y\preceq\dom(\emm_2\diamond\emm_1)$, another application of Lemma \ref{before associativity} gives $Y\preceq\dom\emm_1$ and $\emm_1(Y)\preceq\dom\emm_2$, and in addition that $(\emm_2\diamond\emm_1)|Y=\emm_2\circ(\emm_1|Y)$. Hence $$(\emm_2\diamond\emm_1)(Y)=\range(\emm_2\diamond\emm_1|Y)=\range(\emm_2\circ(\emm_1|Y))=\emm_2(\emm_1(Y)),$$
where the last equality follows from (\ref{number three}).
Thus if we set $Z=\emm_1(Y)$, we have $Z\preceq\dom\emm_2$ and
$\emm_2(Z)=(\emm_2\diamond\emm_1)(Y) \preceq\dom\emm_3$. The first
assertion of Lemma \ref{before associativity} now gives
$Z\preceq\dom(\emm_3\diamond\emm_2)$. Furthermore, since
$Y\preceq\dom\emm_1$ and
$\emm_1(Y)=Z\preceq\dom(\emm_3\diamond\emm_2)$, the first assertion of
Lemma \ref{before associativity} gives
$Y\preceq\dom((\emm_3\diamond\emm_2)\diamond\emm_1)$, i.e. $Y\preceq
Y'$. The second assertion of Lemma \ref{before associativity} gives
$(\emm_3\diamond\emm_2)|Z=\emm_3\circ(\emm_2|Z)$ and
$((\emm_3\diamond\emm_2)\diamond\emm_1)|Y=(\emm_3\diamond\emm_2)\circ(\emm_1|Y)$. Hence
$$
\begin{aligned}
((\emm_3\diamond\emm_2)\diamond\emm_1)|Y&= (\emm_3\diamond\emm_2)\circ(\emm_1|Y)=((\emm_3\diamond\emm_2)|Z)\circ(\emm_1|Y) \quad\text{(by (\ref{number one}))}
\\&=(\emm_3\circ(\emm_2|Z)) \circ(\emm_1|Y)\\
&=\emm_3\circ((\emm_2|Z) \circ(\emm_1|Y)) \quad\text{(by \ref{easy associativity})}\\
&=\emm_3\circ(\emm_2 \circ(\emm_1|Y)) 
\quad\text{(by (\ref{number one}))}\\
& =\emm_3\circ((\emm_2 \diamond\emm_1)|Y)=\emm_3\diamond(\emm_2 \diamond\emm_1).
\end{aligned}
$$

We now claim that $Y=Y'$, so that $((\emm_3\diamond\emm_2)\diamond\emm_1)|Y=
(\emm_3\diamond\emm_2)\diamond\emm_1$. With the equality
$((\emm_3\diamond\emm_2)\diamond\emm_1)|Y=\emm_3\diamond(\emm_2 \diamond\emm_1)$ established above, this will complete the proof of the second assertion of the lemma. As we have shown that $Y\preceq Y'$, in view of Proposition \ref{new partial order} it is enough to show that $Y'\preceq Y$. To this end, note that by Lemma \ref{before associativity} we have $Y'\preceq\dom\emm_1$ and $\emm_1(Y')\preceq\dom(\emm_3\diamond\emm_2)$. Set $Z'=\emm_1(Y')$, so that $Z'\preceq\dom(\emm_3\diamond\emm_2)$; then Lemma \ref{before associativity} gives $Z'\preceq\dom(\emm_2)$, and $\emm_2(Z')\preceq\dom(\emm_3)$. Since $Y'\preceq\dom\emm_1$ and  $\emm_1(Y')=Z'\preceq\dom(\emm_2)$, another application of Lemma \ref{before associativity} gives $Y'\preceq\dom(\emm_2\diamond\emm_1)$, and $(\emm_2\diamond\emm_1)|Y'=\emm_2\circ(\emm_1|Y')$.
Hence
\Equation\label{why not}
\begin{aligned}
\range((\emm_2\diamond\emm_1)|Y')&=\range(\emm_2\circ(\emm_1|Y'))=\emm_2(\emm_1(Y')) \quad\text{(by (\ref{number three}))}
\\&=\emm_2(Z')\preceq\dom\emm_3.
\end{aligned}
\EndEquation
By definition (see \ref{restriction}) we have 
$((\emm_2\diamond\emm_1))(Y')=\range((\emm_2\diamond\emm_1)|Y')$, which
with (\ref{why not}) gives $((\emm_2\diamond\emm_1))(Y')\preceq\dom\emm_3$.

A final application of Lemma \ref{before associativity} now shows that
$Y'\preceq\dom(\emm_3\diamond(\emm_2\diamond\emm_1))=Y$, as
required. Thus the second assertion is proved.

To prove the first assertion, we begin with the definition $\emm_2\diamond\emm_1=\emm_2\circ(\emm_1^{-1}|(\range\emm_1\wedge\dom\emm_2))^{-1}$. Since $\range (\emm_1^{-1}|(\range\emm_1\wedge\dom\emm_2))^{-1}=\dom
(\emm_1^{-1}|(\range\emm_1\wedge\dom\emm_2))=
\range\emm_1\wedge\dom\emm_2$, 
we may apply (\ref{lunch}), with $(\emm_1^{-1}|(\range\emm_1\wedge\dom\emm_2))^{-1}$ playing the role of $\emm_1$, to rewrite the definition as
\Equation\label{fort lauderdoodle}
\emm_2\diamond\emm_1=(\emm_2|
(\range\emm_1\wedge\dom\emm_2))
\circ(\emm_1^{-1}|(\range\emm_1\wedge\dom\emm_2))^{-1}.
\EndEquation
Applying the final observation of \ref{zoltan}, we then obtain
$$(\emm_2\diamond\emm_1)^{-1}=
(\emm_1^{-1}|(\range\emm_1\wedge\dom\emm_2))
\circ
(\emm_2|
(\range\emm_1\wedge\dom\emm_2))^{-1}.$$
Since $\range\emm_1=\dom(\emm_1^{-1})$ and
$\dom\emm_2=\range(\emm_2^{-1})$, this may be rewritten as
\Equation\label{fort noodle by the poodle}
(\emm_2\diamond\emm_1)^{-1}=
(\emm_1^{-1}|(\range(\emm_2^{-1})\wedge\dom(\emm_1^{-1})))
\circ
(\emm_2|(\range(\emm_2^{-1})\wedge\dom(\emm_1^{-1})))
^{-1}.
\EndEquation
On the other hand, substituting $\emm_2^{-1}$ and $\emm_1^{-1}$,
respectively, for $\emm_1$ and $\emm_2$ in (\ref{fort lauderdoodle})
shows that the right hand side of (\ref{fort noodle by the poodle}) is equal to $\emm_1^{-1}\diamond\emm_2^{-1}$. This proves the first assertion of the proposition.
\EndProof

\Number\label{more associativity}
In view of the second assertion of Proposition \ref{associativity}, if $\otheroldLambda$ is a
negative, closed, orientable $2$-orbifold, then for any elements
$\emm_1,\cdots,\emm_k$ of $\barcalm(\otheroldLambda,\otheroldLambda)$, where $k\ge0$, we may
unambiguously define an element $\emm_k\diamond\cdots\diamond\emm_1$ of
$\barcalm(\otheroldLambda,\otheroldLambda)$ (which is interpreted as the identity if $k=0$); and for
$0\le  j\le k$ we have
$(\emm_k\diamond\cdots\diamond\emm_{j+1})\diamond(\emm_j\diamond\cdots\diamond\emm_1)=\emm_k\diamond\cdots\diamond\emm_1$. 
The first assertion of Proposition \ref{associativity} then implies by
induction that 
\Equation\label{snicker snoodle}
(\emm_k\diamond\cdots\diamond\emm_1)^{-1}=\emm_1^{-1}\diamond\cdots\diamond\emm_k^{-1}.
\EndEquation

If
$\emm\in\barcalm(\otheroldLambda,\otheroldLambda)$ is given, then for every 
$\emm\in\barcalm(\otheroldLambda,\otheroldLambda)$ and every $k\ge0$, we set
$\emm^{\diamond k}=\emm\diamond\cdots\diamond\emm$, where $\emm$
appears $k$ times in the latter expression. We also define
$\emm^{\diamond(-m)}=(\emm^{\diamond m})^{-1}$, which by (\ref{snicker snoodle}) is the same as $(\emm^{-1})^{\diamond m}$.
\EndNumber

We conclude this section with a result will be needed in Section \ref{clash section}:

\Proposition\label{steinbeck}
Suppose that $\otheroldLambda$ is a
negative, closed, orientable $2$-orbifold, that
$\frakJ$ is a taut suborbifold of $\otheroldLambda$, and that a finite group $G$ acts on $\otheroldLambda$ by homeomorphisms in such a way that $g\cdot\frakJ$ is isotopic to $\frakJ$ for every $g\in G$. Then there is a suborbifold of $\otheroldLambda$ which is isotopic to $\frakJ$ and is invariant under the action of $G$.
\EndProposition

The basic strategy for proving  Proposition \ref{steinbeck} is to transition
to the smooth category, and then isotop $\frakJ$ to a suborbifold
whose boundary curves are geodesic in a $G$-invariant
metric. Technical issues arise, first because some components of
$\partial\frakJ$ may be isotopic to curves that doubly cover geodesics
homeomorphic to $[[0,1]]$, and second because different components of
$\partial\frakJ$ may be isotopic to one another. The latter issue
becomes particularly complicated if $\frakJ$ has annular components.

\Proof[Proof of Proposition \ref{steinbeck}]
Since $\otheroldLambda$ is negative, the orbifold $\otheroldLambda/G$ is also
negative. Let us fix a smooth structure on $\otheroldLambda/G$ that is
compatible with the PL structure that it inherits from
$\otheroldLambda$. We pull this smooth structure back to
$\otheroldLambda$ to obtain a smooth structure which is $G$-invariant
and compatible with the PL structure of $\otheroldLambda$; we denote
by $\oldLambda$ the orbifold $\otheroldLambda$ endowed with
this PL structure. It is then natural to denote the orbifold
$\otheroldLambda/G$, endowed with its aforementioned smooth structure,
by $\oldLambda/G$.

Let $\frakJsmooth$ be a smooth suborbifold of $\oldLambda$ which is piecewise smoothly isotopic to $\frakJ$.
Set $B=|\partial\frakJsmooth|$. Since $\partial\frakJsmooth$ is a smooth, locally
separating, closed $1$-suborbifold of the orientable $2$-orbifold
$\oldLambda$, the components of $B$ are smooth
simple closed curves in the smooth $2$-manifold
$|\oldLambda|-\fraks_\otheroldLambda$. 
The tautness of $\frakJ$
implies that these curves are all $\pi_1$-injective in the orbifold
$\oldLambda$, and therefore in the manifold $|\oldLambda|-\fraks_\oldLambda$. Let us define an equivalence relation on the finite
set $\calc(B)$ by declaring two components  of $B$ to be equivalent if
and only if they are homotopic in $|\otheroldLambda|-\fraks_\otheroldLambda$. According to \cite[Lemma 2.4]{epstein}
(translated into the smooth category), two components of $B$ are
equivalent if and only if they cobound a weight-$0$ smooth annulus in $|\oldLambda|-\fraks_\otheroldLambda$. We let $E$ denote the set of all equivalence classes in $\calc(B)$. For any  $c\in E$, let
%component $\frakC$ of $\frakB$, let 
$n_c$ denote the cardinality of the equivalence class $c$.
%the equivalence class of $\frakC$. 
%Let us choose a closed $1$-submanifold of $\frakB_0$ of $\frakB$ such that $\calc(\frakB_0)$ contains exactly one representative of each equivalence class in $\calc(\frakB)$. 
For each  $c\in E$, we denote by $W_c$ the union of the curves in the equivalence class $c$, and we fix a
 weight-$0$ smooth annulus
 $R_c\subset|\oldLambda|-\fraks_\oldLambda$ which is a
 common tubular neighborhood of all the curves in
 $c$; we take the family $(R_c)_{c\in E}$ to be pairwise
 disjoint. Let us also fix an element $C_c$ of the equivalence class
 $c$ for each $c\in E$, and set $Y=\bigcup_{c\in E}C_c$; thus $Y$ is a
 union of components of $B$, and is in particular a
 $\pi_1$-injective closed $1$-manifold in $\oldLambda$.

% as $c$ ranges over the equivalence classes in $\frakR_c$.
% $\frakR_\frakC$ denotes the component of $\frakR$ containing $\frakC$, element \oldLambda$ such that, if $\frakC$ is any component of $\frakB_0$, every component of $\frakB$ belonging to the equivalence class of $\frakC$  is contained in the component of $\frakR$ containing $\frakC$, and is a core curve of that component.
%We will denote by $\frakR_\frakC$ the component of $\frakR$ containing $\frakC$, for every $\frakC\in\frakB_0$.

By construction, the (pairwise disjoint) simple closed curves in the family  $(C_c)_{c\in E}$ are pairwise non-homotopic in $|\oldLambda|-\fraks_\oldLambda$, and hence  the  $\obd(C_c)$ are  pairwise non-isotopic in $\oldLambda$.

%By the equivalence of the PL and smooth categories for $2$-manifolds, there is an isotopy of $|\otheroldLambda|-\fraks_\otheroldLambda|$, constant outside a compact subset of   $|\otheroldLambda|-\fraks_\otheroldLambda|$, which carries $\frakB$ and each $\frakR_c|$ onto smooth submanifolds of  $|\oldLambda|-\fraks_\otheroldLambda|$. Since this smoothing isotopy is constant outside a compact subset of   $|\otheroldLambda|-\fraks_\otheroldLambda|$, it extends to an orbifold isotopy of $\otheroldLambda$. We denote the (smooth) images of $\frakB$, $C_c$ and $\frakR_c$ under this isotopy by $\frakBsmooth$, $\Ccsmooth$ and $\frakRcsmooth$. Then the smooth counterparts of the assertions made above in the PL category about $\frakB$, $C_c$ and $\frakR_c$ hold for  $\frakBsmooth$, $\Ccsmooth$ and $\frakRcsmooth$ respectively. 

The negativity of $\otheroldLambda/G$ implies that $\oldLambda/G$
% denote the and $\other, and hence
admits a metric under which each component is a complete hyperbolic  orbifold without cusps. It follows
that $\oldLambda$ admits a $G$-invariant metric in which each
component is a complete hyperbolic orbifold without cusps. %Fix an element $\frakJ$ of
%$\Theta(\otheroldLambda)$ such that $[\frakJ]=J$. Since $\otheroldLambda$ is
%closed, and since the definition of $\Theta(\otheroldLambda)$ implies that

Since $\oldLambda$ is negative, it is very good.
Let $p:\toldLambda\to\oldLambda$ be a finite-sheeted smooth regular
covering such that $\toldLambda$ is a smooth $2$-manifold. Let $N$
denote the covering group of this regular covering. We equip
$\toldLambda$ with the componentwise hyperbolic metric inherited
from $\oldLambda$, so that $N$ acts on $\toldLambda$ by
isometries. We denote by $\tG$ the group of all self-diffeomorphisms
of $\toldLambda$ that are lifts of elements of $G$. Then $\tG$ is
finite, and there is a homomorphism $q:\tG\to G$ such that for each
$\tg\in\tG$ we have $p\circ\tg=q(g)\circ p$. We may choose the covering
$\toldLambda$ to be characteristic in the sense that every
self-diffeomorphism of $\oldLambda$ admits a lift to
$\toldLambda$. This implies that $q$ is surjective.

Set $\tY=p^{-1}(Y)\subset\toldLambda$. Since $Y$ is a
 $\pi_1$-injective closed $1$-manifold in $\oldLambda$, the
$N$-invariant closed $1$-manifold $\tY$ is $\pi_1$-injective in
 $\toldLambda$. Hence for each component $\tC$ of $\tY$, the
 inclusion map $i_{\tC}: \tC\to\toldLambda$ is homotopic to a
 geodesic immersion $\gamma_{\tC}:\tC\to\toldLambda$, which is
 unique (up to reparametrization preserving the orientation of $\tC$). Since
 $i_{\tC}$ is injective for each component $\tC$ of $\tY$, and
 $i_{\tC}(\tC)\cap i_{\tC'}(\tC')=\emptyset$ for any distinct
 components $\tC$ and $\tC'$ of $\tY$, 
the immersion $i_{\tC}$ is injective for each component $\tC$ of $\tY$, and
the simple closed geodesics $\gamma_{\tC}(\tC)$ and
$\gamma_{\tC'}(\tC')$ are either disjoint or equal for any distinct
 components $\tC$ and $\tC'$ of $\tY$ (see, for example, \cite{otherFHS}). 
Hence
$\tYgeo:=\bigcup_{\tC\in\calc(\tY)}\gamma_{\tC}(\tC)$ is a
$1$-submanifold of $\toldLambda$ whose components are simple
closed geodesics. The uniqueness (up to reparametrization) of the geodesic immersion 
$\gamma_{\tC}$ homotopic to
$i_{\tC}$ implies that for each component $\tC$ of $\tY$, and each
$g\in N$, we have $\gamma_{g\cdot\tC}=g\cdot\gamma_{\tC}$ (up to reparametrization). Hence
$\tYgeo$ is $N$-invariant. If we set
$\Ygeo=p(\tYgeo)$, it now follows that
$\obd(\Ygeo)$ is a geodesic $1$-suborbifold
of $\oldLambda$; each component of  $\obd(\Ygeo)$ is
(orbifold-)diffeomorphic either to $\SSS^1$ or to $[[0,1]]$, and $|p|^{-1}(\Ygeo)=\tYgeo$.

We claim:
\Claim\label{tohu vavohu}
 The $1$-orbifold $\obd(\Ygeo)$ is $G$-invariant. 
\EndClaim

To prove \ref{od tohu vavohu}, it suffices to show that $\Ygeo$ is
$G$-invariant. Let $g\in G$ be given, and let $\Cgeo$ 
be a component of $\Ygeo$. Since $q:\tG\to G$ is surjective, we
may fix an element $\tg\in\tG$ with $q(\tg)=g$. Choose a component
$\tCgeo$ of $p^{-1}(\Cgeo)$. Then $\tCgeo$ is in particular a
component of $\tYgeo$, and hence $\tCgeo=\gamma_{\tC}(\tC)$
for some component $\tC$ of $\tY$. Thus $C:=|p|(\tC)$ is a component of
$Y$, and hence of $B=|\partial\frakJsmooth|$. By hypothesis,
$g\cdot\frakJ$ is PL orbifold-isotopic to $\frakJ$ in
$\otheroldLambda$. Hence $g\cdot\obd( C)$ is smoothly isotopic in
$\oldLambda$ to
some component  of $\obd(B)$. Since each component of $B$ is, by the
definition of $Y$, homotopic in $|\oldLambda|-\fraks_\oldLambda$
to some component 
%$C'$ 
of $Y$, it follows that $g\cdot\obd(C)$ is isotopic in $\oldLambda$ to
some component of $\obd(Y)$, and hence   that the component $\tg\cdot\tC$ of $p^{-1}(g\cdot\obd( C))$ is isotopic
in $\toldLambda$
to some component $\tC'$ of $p^{-1}(\obd(Y))=\tY$. The simple closed curve
$\tC'$ is in turn homotopic in $\toldLambda$ to the simple closed
geodesic $\gamma_{\tC'}(\tC')$, which is a component of $\tYgeo$. But
since $\tCgeo$ is the unique simple closed geodesic
homotopic to $\tC$,  the unique simple closed geodesic
homotopic to $\tg\cdot\tC$
 is $\tg\cdot\tCgeo$; it follows that
$\tg\cdot\tCgeo=\gamma_{\tC'}(\tC')\subset\tYgeo$. Hence
$g\cdot\Cgeo=|p|(\tg\cdot\tCgeo)\subset |p|(\tYgeo)=\Ygeo$. This
establishes \ref{tohu vavohu}.

Let us write $\Ygeo=\Ygeo^{(1)}\discup\Ygeo^{(2)}$, where each component of $\obd(\Ygeo^{(1)})$ is
diffeomorphic to $\SSS^1$ and each component of $\obd(\Ygeo^{(2)})$ is
diffeomorphic to $[[0,1]]$. From \ref{tohu vavohu} we immediately
deduce:

\Claim\label{den sum}
Each of the $1$-orbifolds $\obd(\Ygeo^{(1)})$ and $\obd(\Ygeo^{(2)})$ is
$G$-invariant. 
\EndClaim

Let $\oldPi$ denote a tubular neighborhood
of $\obd(\Ygeo^{(2)})$ in $\oldLambda$ which is disjoint from
$\obd(\Ygeo^{(1)})$
. By \ref{den sum} we may take $\oldPi$ to be $G$-invariant.
% (so that $\oldPi$ is well-defined up to isotopy in
% $\otheroldLambda-\obd(\Ygeo^{(1)})$). 
Set $Y^*=\Ygeo^{(1)}\discup\partial|\oldPi|$, so
that $Y^*$ is a closed $1$-submanifold of
$|\oldLambda|-\fraks_\oldLambda$. Then \ref{den sum} and the
$G$-invariance of $\oldPi$ imply:

\Claim\label{od tohu vavohu}
 The $1$-manifold $Y^*$ is $G$-invariant. 
\EndClaim

We claim: 

\Claim\label{bored}
There is a closed $1$-submanifold of
$|\oldLambda|-\fraks_\oldLambda$   which is
isotopic to $Y$ in $|\oldLambda|-\fraks_\oldLambda$ and is
disjoint from $Y^*$.
\EndClaim

To prove \ref{bored}, let $Y_0$ be a closed $1$-submanifold of
$|\oldLambda|-\fraks_\oldLambda$   which is
isotopic to $Y$ in $|\oldLambda|-\fraks_\oldLambda$ and transverse to
$Y^*$, and which is chosen, among all $1$-submanifolds of
$|\oldLambda|-\fraks_\oldLambda$ that are
isotopic to $Y$ in $|\oldLambda|-\fraks_\oldLambda$ and transverse to
$Y^*$, in such a way as to minimize the cardinality of its intersection
with $Y^*$. Note that since $Y_0$ is  isotopic to $Y$ in
$|\oldLambda|-\fraks_\oldLambda$, the $1$-manifold
$\tY_0:=|p|^{-1}(Y_0)$ is isotopic to $\tY=|p|^{-1}(Y)$ in
$\toldLambda$. We must show that $Y_0\cap Y^*=\emptyset$.

Assume to the contrary that $Y_0\cap Y^*\ne\emptyset$, so that
$\tY_0$ has non-empty intersection with $\tY^*:=|p|^{-1}(Y^*)$. Then
there are components $\tC_0$ and $\tC_1$ of $\tY_0$ and $\tY^*$ respectively
such that $\tC_0\cap \tC_1\ne\emptyset$. It follows from the
definitions of $Y^*$ and
$\tY^*$ that $\tC_1$ is homotopic in $\toldLambda$ to a curve disjoint
from $\tY$, and hence (since $\tY$ and $\tY_0$ are isotopic) to a curve
disjoint from $\tY_0$, and in particular from $\tC_0$. It then follows
from Proposition \ref{beats epstein} that $\tC_1$ and $\tC_0$ have a
degenerate crossing in the $2$-manifold $\toldLambda$. (In
\ref{nondeg cross}, degenerate crossings and disks of degeneracy were
defined only in the context of the PL category; however, the
definition, and the proof of Proposition \ref{beats epstein}, go through in the smooth category without change, except
that a disk of degeneracy is only a piecewise smooth disk.) In particular $\tY_0$ and $\tY^*$ have a degenerate
crossing. Let  $D$ be a disk of degeneracy for $\tY_0$ and
$\tY^*$ which is minimal with respect to inclusion among all  disks
of degeneracy for $\tY_0$ and $\tY^*$. Then $\partial
D=D\cap(\tY_0\cup \tY^*)=\ta_0\cup \ta_1$, where 
%\subset|\otheroldLambda|-\fraks_\otheroldLambda$ such that $\partial D$ has
%the form $\ta_1\cup \ta_2$, where
 $\ta_0$ is an arc contained in a component $\tC_0'$ of $\tY_0$, and $\ta_1$
 is an arc contained in a component $\tC_1'$ of $\tY^*$. It follows
 that $\ta_0\cap\tY^*=\partial \ta_0$ and that
 $\ta_1\cap\tY_0=\partial \ta_1$. 
%YC

For $i=0,1$, set $C_i'=|p|(\tC_i')$. Then $p_i:=|p||\tC_i':\tC_i'\to C_i'$ is
a covering map for $i=0,1$. Since $\tC_0'\cap\tC_1'\ne\emptyset$, we
have $C_0'\cap C_1'\ne\emptyset$. 
Since $\ta_0\cap\tY^*=\partial \ta_0$, we have
$|p|(\inter \ta_0)\cap Y^*=|p|((\inter \ta_0)\cap\tY^*)=\emptyset$.
If $p_0|\inter \ta_0:\inter \ta_0\to C_0'$ is
surjective, then $C_0'\cap
C_1'\subset |p|(\inter \ta_0)\cap Y^*=\emptyset$,
a contradiction. Hence $p_0|\inter \ta_0:\inter \ta_0\to C_0'$ is not
surjective. If $p_0|\ta_0:\ta_0\to C_0'$ is surjective, it now follows that
$p_0(\partial \ta_0)$ consists of a single point $P$ of $C_0'$;
since $|p|(\inter \ta_0)\cap Y^*=\emptyset$ and $C_0'\cap
C_1'\ne\emptyset$, we have $P\in C_1'$. By orientability,  $P$ is then a
non-transverse point of intersection of $C_0'$ and $C_1'$, a
contradiction. Hence $p_0|\ta_0:\ta_0\to C_0'$ is not surjective. Since
$p_0:\tC_0'\to C_0'$ is a covering map, it follows that $p_0|\ta_0$ is
one-to-one. The same argument, with the roles of $\tY^*$ and $\tY_0$
interchanged, shows that $p_1|\ta_1$ is one-to-one. 
Since $\ta_0\cap |p|^{-1}(Y^*)=\ta_0\cap\tY^*=\partial \ta_0$, it follows that
$|p|\big|(\partial D)=|p|\big|(\ta_0\cup \ta_1)$ is one-to-one. 

Thus $L:=|p|(\partial D)$ is a piecewise smooth simple closed curve in $|\oldLambda|-\fraks_\oldLambda$,
and $|p|$ maps $\partial D$ homeomorphically onto $L$. Since
$\partial D$ represent the trivial element of
$H_1(\toldLambda;\FF_2 )$, the curve $L$ represent the trivial
element of $H_1(|\oldLambda|;\FF_2 )$, and therefore separates
$|\oldLambda|$. Since $D\cap(\tY_0\cup \tY^*) =\partial D$, the
map $|p|$ must take $\inter D$ into a single component $R$ of
$|\oldLambda|-L$. Then $\overline{R}$ is a topological $2$-submanifold of
$|\oldLambda|$, and $|p|\big| D:D\to\overline{R}$ is a
boundary-preserving map. But
$r:=|p|:|\toldLambda|\to|\oldLambda|$ is a branched covering
map with branch locus $\fraks_\oldLambda$, and $Y_0$ and $Y^*$
are disjoint from $\fraks_\oldLambda$. Hence the
boundary-preserving map $r:D\to\overline{R}$ is also a branched
covering. Since $r|\partial D=|p|\big|\partial D$ is one-to-one, the
degree of the branched covering $r$ must be $1$, i.e. $r$ is a
homeomorphism. This implies that $R$ is a disk and that
$R\cap\fraks_\oldLambda=\emptyset$. Since $|p|$ maps $\partial D$
homeomorphically onto $\partial R$, we may write
$\partial R=a_0\cup a_1$, where $a_i:=|p|(\ta_i)$ is an arc for
$i=0,1$; we have $a_0\cap a_1=\partial a_0=\partial a_1$, and $a_0$
and $a_1$ are respectively contained in $Y_0$ and $Y^*$. Thus $R$ is a
disk of degeneracy for $Y_0$ and $Y^*$ in the smooth $2$-manifold
$|\oldLambda|-\fraks_\oldLambda$. This implies
(cf. \ref{nondeg cross}) that  $Y_0$ is isotopic in
$|\oldLambda|-\fraks_\oldLambda$, via a smooth isotopy constant
outside a small neighborhood of $R$, to a $1$-manifold $Y_0'$ meeting
$Y^*$ transversally, with $ \card(Y_0'\cap Y^*)<\card(Y_0\cap
Y^*)$. This contradicts the minimality property of $Y_0$, and thus
(\ref{bored}) is proved.

According to \ref{bored}, we may fix a closed $1$-submanifold $Y_0$
of
$|\oldLambda|-\fraks_\oldLambda$,
isotopic to $Y$ in $|\oldLambda|-\fraks_\oldLambda$ (and
therefore $\pi_1$-injective), such that
$Y_0\cap Y^*=\emptyset$. In particular $Y_0$ is disjoint from $\partial|\oldPi|$, so
that each component of $Y_0$ is either contained in $\inter|\oldPi|$
or disjoint from $|\oldPi|$. Since each component of $\oldPi$
is annular, it follows from \ref{cobound} that each component of $Y_0$ contained in $\inter|\oldPi|$ must cobound  a weight-$0$ annulus with
a boundary component of $|\oldPi|$. Hence after possibly modifying
$Y_0$ by a further isotopy in
$|\oldLambda|-\fraks_\oldLambda$, we may assume that
\Equation\label{energized}
Y_0\cap|\oldPi|=\emptyset.
\EndEquation

Since $Y_0$ is $\pi_1$-injective in
$\oldLambda$, the $1$-submanifold $\tY_0$ of $\toldLambda$
is $\pi_1$-injective. Since $\Ygeo\subset\Ygeo^{(1)}\cup|\oldPi|\subset
Y^*\cup|\oldPi|$, and since $Y_0$ is by definition disjoint from $Y^*$,
and is disjoint from $|\oldPi|$ by \ref{energized}, we have
$Y_0\cap\Ygeo=\emptyset$.
Setting  $\tY_0:=|p|^{-1}(Y_0)$, we deduce that
\Equation\label{reenergized}
\tY_0\cap\tYgeo=\emptyset.
\EndEquation

We claim:

\Claim\label{blague}
If two components of $\tY_0$ are homotopic simple closed
  curves in $\toldLambda$, then they are in the same $N$-orbit.
\EndClaim

To prove \ref{blague}, suppose that $\tC_0$ and $\tC_0'$  are homotopic components of $\tY_0$. By
\cite[Lemma 2.4]{epstein}, $\tC_0$
and $\tC_0'$ cobound an annulus $\tA \subset \toldLambda$. We will show by
induction on the number $n$ of components of $\tY_0$ contained in
$\inter \tA $ that $\tC_0$ and $\tC_0'$ lie in the same $N$-orbit. If $n=0$
then $\tA $ is the closure of a component of $\toldLambda-\tY_0=|p|^{-1}(|\oldLambda|-Y_0)$;
since $p$ is an orbifold covering, and since the hypothesis that
$\otheroldLambda$ is negative implies that no component of
$\oldLambda$ is toric, it follows that $p|\obd(\tA)$ is an orbifold
covering of some suborbifold $\frakA$ of $\oldLambda$. Since
$\tA$ is an annulus, $\frakA$ is annular. If $|\frakA|$ is a
weight-$0$ annulus then its two boundary components are components of
$Y_0$ that are homotopic to each other in $|\oldLambda|-\fraks_\oldLambda$; this implies that $Y$ has two
distinct components that are homotopic to each other, a contradiction
to the definition of $Y$. Hence $|\frakA|$ is a weight-$0$ disk. In
particular, $\partial\frakA$ is connected, which implies that $\tC_0$
and $\tC_0'$ have the same image under $|p|$ and therefore lie in the
same $N$-orbit.

Now suppose that $n>0$, and choose a component $\tC_0''$ of $\tY_0$ contained in
$\inter \tA $.  Since $\tY_0$ is $\pi_1$-injective in
$\toldLambda$, the component $\tC_0''$ of $\tY_0$ cobounds
sub-annuli of $\tA$ with $\tC_0$ and
$\tC_0'$. As the interior of each of these annuli contains fewer than
$n$ components of $\tY_0$, the induction hypothesis implies that the
$N$-orbit of $\tC_0''$ contains both $\tC_0$ and $\tC_0'$. This completes
the proof of \ref{blague}.
%A

Let us denote by $\scrQ$ the set of all
annuli in the $2$-manifold $\toldLambda$ that have one boundary
component in $\tY_0$
and one in $\tYgeo$. The importance of this definition lies in the following fact:

\Claim\label{cafe du monde}
Every component of $Y_0$ is the image under $|p|$ of a boundary component of an element of $\scrQ$.
\EndClaim

To prove \ref{cafe du monde},  suppose that $C_0$ is a component of
$Y_0$. Since $Y_0$ and $Y$ are isotopic, $C_0$ is isotopic to a
component $C$ of $Y$. Choose a component $\tC$ of
$|p|^{-1}(C)$. Then $\gamma_{\tC}(\tC)$ is a simple 
closed geodesic which is homotopic to $\tC$ and therefore to a
component $\tC_0$ of $|p|^{-1}(C_0)$. Since $\gamma_{\tC}(\tC)\subset\tYgeo$
and $\tC_0\subset\tY_0$ are disjoint by (\ref{reenergized}), it now
follows from \cite[Lemma 2.4]{epstein}  that $\gamma_{\tC}(\tC)$ and $\tC_0$
cobound an annulus $\tQ\subset\toldLambda$, which is by
definition an element of $\scrQ$. Since $C_0=|p|(\tC_0)$, this proves \ref{cafe du monde}.

Now we claim:

\Claim\label{i hardly know you}
For any two distinct annuli $\tQ_0,\tQ_1\in\scrQ$, either $\tQ_0\cap \tQ_1=\emptyset$, or
$\tQ_0\cap \tQ_1$ is a component of $\tYgeo$.
\EndClaim

To prove \ref{i hardly know you}, assume that $\tQ_0\cap
\tQ_1\ne\emptyset$. By the definition of $\scrQ$, we may write
 $\partial \tQ_i=\tC_i\cup\tCgeo_i$ for $i=1,2$, where the
 $\tC_i$  are components of $\tY_0$, and the $\tCgeo_i$ are components
 of $\tYgeo$. Since $\tY_0\cap\tYgeo=\emptyset$ but $\tQ_0\cap
\tQ_1\ne\emptyset$, one of the $\tQ_i$ must contain a boundary component
of the other. But $\partial \tQ_i$ is $\pi_1$-injective in
$\toldLambda$  for $i=1,2$, since both $\tY_0$ and $\tYgeo$ are
$\pi_1$-injective. It therefore follows that the core curves of $\tQ_0$
and $\tQ_1$ are isotopic to each other. By hyperbolicity, a given
homotopy class of closed curves in $\toldLambda$ can contain only
one geodesic. Hence $\tCgeo_0=\tCgeo_1$. This implies that either (i)
$\tQ_0\subset \tQ_1$, (ii) $\tQ_1\subset \tQ_0$, (iii) $\tQ_0\cap \tQ_1=\tCgeo_0$,
or (iv) $\tQ_0\cup \tQ_1=\toldLambda$, and $\toldLambda$ is a
torus. Alternative (iv) contradicts the hypothesis that
$\otheroldLambda$ is negative. Now suppose that (i) holds. Since $\tQ_0$
and $\tQ_1$ have isotopic
core curves, $\tC_0$ and $\tC_1$ are homotopic; in view of
\ref{blague}, $\tC_0$ and $\tC_1$ are in the same $N$-orbit. If we fix
$g\in N$ such that $g\cdot\tC_0=\tC_1$, then $g$ leaves the homotopy
class of $\tC_0$ invariant; since $\tCgeo_0$ is the unique closed
geodesic in this homotopy class, $g$ also leaves $\tCgeo_0$
invariant. Hence $\partial(g\cdot
\tQ_0)=(g\cdot\tC_0)\cup(g\cdot\tCgeo_0)=\tC_1\cup\tCgeo_0=\tC_1\cup\tCgeo_1=\partial
\tQ_1$. Again using that $\toldLambda$ is not a torus, we deduce
that the annuli $g\cdot \tQ_0$ and $\tQ_1$, which have the same
boundary, must coincide. Alternative (i) therefore implies that
$\tQ_0\subset g\cdot \tQ_0$, and since $\tQ_0\ne \tQ_1$ the inclusion is
proper. This is a contradiction since $g$ has finite order. If (ii)
holds we obtain the same contradiction; hence (iii) must hold, and
\ref{i hardly know you} is proved.

Set $\tcalq=\bigcup_{\tQ\in\scrQ}\tQ$. It follows from \ref{i hardly know
  you} that we may write $\tcalq$ as a disjoint union
$\tcalq^{(1)}\discup\tcalq^{(2)}$, where each component of $\tcalq^{(1)}$
is an element of $\scrQ$, and each component of $\tcalq^{(2)}$ has the
form $\tQ\cup \tQ'$ for some $\tQ,\tQ'\in\scrQ$ such that $\tQ\cap \tQ'$ is a
component of $\tYgeo$. Note also that since $\tY_0$
and  $\tYgeo$ are $N$-invariant, the collection of annuli $\scrQ$ is
$N$-invariant, and the sets $\tcalq^{(1)}$ and $\tcalq^{(2)}$ are
therefore $N$-invariant. Hence if we set $\calq=|p|(\tcalq)$, we may
write $\calq=\calq^{(1)}\discup\calq^{(2)}$, where
$\calq^{(i)}=|p|(\tcalq^{(i)})$ for $i=1,2$,  The components of
$\calq^{(i)}$ (for $i=1.2$) are precisely the sets of the form
$|p|(\tQ)$, where
$\tQ$ is a component
of $\tcalq^{(i)}$. For each component $\tQ$ of $\tcalq^{(i)}$, the
component $Q=|p|(\tQ)$ of $\calq^{(i)}$ is canonically identified with the
quotient of $\tQ$ by its stabilizer in $N$; in particular
$p|\obd(\tQ):\obd(\tQ)\to\obd(Q)$ is an orbifold covering, and hence
$\obd(Q)$ is annular. Furthermore, for each
component $Q$ of $\calq^{(i)}$, the components of $\tcalq^{(i)}$ that
are mapped onto $Q$ by $|p|$ constitute a single $N$-orbit of
components of $\tcalq^{(i)}$.

Consider an arbitrary component $Q$ of $\calq^{(1)}$, and choose a
component $\tQ$ of $\tcalq^{(1)}$ with $|p|(\tQ)=Q$. According to the
definitions we have $\tQ\in\scrQ$, which means that $\tQ$ is an
annulus with one boundary component in $\tY_0=|p|^{-1}(Y_0)$
and one in $\tYgeo=|p|^{-1}(\Ygeo
)$. Since $Y_0\cap\tYgeo=\emptyset$
by  (\ref{reenergized}), the annular orbifold $\obd(Q)$ has two boundary
components, one in $\obd(Y_0)$ and one in $\obd(\Ygeo)$; the latter
component must be diffeomorphic to $\SSS^1$ since $\otheroldLambda$ is
orientable, and is therefore a component of $\obd(\Ygeo^{(1)})$. Hence $Q$ is
a weight-$0$ annulus. We denote its boundary components by $C_0(Q)$
and $\Cgeo(Q)$, where $C_0(Q)$
and
$\Cgeo(Q)$ are respectively components of $Y_0$ and $\Ygeo^{(1)}$.

Now consider an arbitrary component $Q$ of $\calq^{(2)}$, and choose a
component $\tQ$ of $\tcalq^{(2)}$ with $|p|(\tQ)=Q$.  According to the
definitions we may write  $\tQ=\tQ\cup \tQ'$, where $\tQ,\tQ'\in\scrQ$ and $\tQ\cap \tQ'$ is a
component $\tCgeo$ of $\tYgeo$. The definition of $\scrQ$  implies that
each $\tQ_i$ is an
annulus with one boundary component in $\tY_0$
and one in $\tYgeo$. Hence $\tQ$ is an annulus whose boundary
components $\tC$ and $\tC'$ are contained in $\tY_0$. Since 
$\tC$ and $\tC'$ are homotopic in $\toldLambda$, it follows from
\ref{blague} that they are in the same $N$-orbit. Hence the orientable annular
orbifold $\obd(Q)$ has connected boundary, so that $Q$ is a weight-$2$
disk, and the points of $Q\cap\fraks_\oldLambda$ are of order
$2$. In particular $\partial Q$ is a single component of $Y_0$. Now let us fix an
element $g\in N$ such that $g\cdot \tC=\tC'$. 
%Since $\otheroldLambda$ is orientable, 
%$\partial\obd( Q)$ is
%diffeomorphic to $\SSS^1$, and is therefore a component of $\oldUpsilon^{(1)}$.
In particular we have
$(g\cdot\tQ)\cap\tQ\ne\emptyset$, and since $\tQ$ is a component of
the $N$-invariant set $\tcalq^{(2)}$, it must itself be $N$-invariant;
thus $g$ interchanges $\tC$ and $\tC'$. Since the collection $\scrQ$
is also $N$-invariant, it now follows that $g$ interchanges $\tQ$
and $\tQ'$ and therefore leaves $\tCgeo$ invariant. But since
$\otheroldLambda$ is orientable, $g$ preserves the orientation of
$\tQ$, and must therefore reverse the orientation of $\tCgeo$. This
implies that $p(\obd(\tCgeo))\subset\inter\obd(Q)$ is diffeomorphic to
$[[0,1]]$, so that $p(\obd(\tCgeo))$ is a component of
$\obd(\Ygeo^{(2)})$. In this case we set $C_0(Q)=\partial Q$ and
$\Cgeo(Q)=|p|(\tCgeo)$, so that  $C_0(Q)$
and
$\Cgeo(Q)$ are respectively components of $Y_0$ and
$\Ygeo^{(2)}$. The fact that the components of $\tcalq^{(2)}$ that
are mapped onto $Q$ by $|p|$ constitute a single $N$-orbit of
components of $\tcalq^{(2)}$ implies that $\Cgeo(Q)$ is well-defined,
i.e. does not depend on the choice of $\tQ$.

If $Q$ is a component of $\calq^{(2)}$, then since
$\Cgeo(Q)\subset\inter Q$ is a component of $\Ygeo^{(2)}$, the
definitions imply that $\Cgeo(Q)$ is contained in a component of
$|\oldPi|$, which we denote by $P_Q$. Since $\partial Q=C_0(Q)\subset
Y_0$ is disjoint from $P_Q\subset\oldPi$ by (\ref{energized}), we have
$P_Q\subset\inter Q$. We denote by $T_Q$ the annulus $Q-\inter P_Q$. Since $\obd(\Cgeo)$ is diffeomorphic to $[[0,1]]$, we
have $\wt\Cgeo=2=\wt Q$, and hence $\wt T_Q=0$. Set $C^*(Q)=\partial
P_Q\subset\partial|\oldPi|\subset Y^*$. The boundary
components of $T_Q$ are $C_0(Q)$ and $C^*(Q)$.

If $Q$ is a component of $\calq^{(1)}$, we have seen that $Q$ is a
weight-$0$ annulus with boundary components $C_0(Q)\subset Y_0$ and
$\Cgeo(Q)\subset\Ygeo^{(1)}\subset Y^*$;
in this case we set $T_Q=Q$ and $C^*(Q)=\Cgeo(Q)$.

Thus we have shown:

\Claim\label{youth in asia}
For every component $Q$ of $\calq$, the set $T_Q$ is a weight-$0$
annulus in $|\oldLambda|$, whose boundary components are
$C_0(Q)\subset Y_0$ and $C^*(Q)\subset Y^*$. Furthermore, we have
$T_Q\subset Q$ for every component $Q$ of $\calq$, so that the family
$(T_Q)_{Q\in\calc(\calq)}$ is pairwise disjoint.
\EndClaim

Now consider an arbitrary element $\tQ$ of $\scrQ$. The definition of
$\scrQ$ implies that one component $\tC_0$
of $\partial\tQ$ is contained
in $\tY_0$, and that its other component $\tCgeo$ is contained
in $\tYgeo$. The definitions of $\tcalq^{(1)}$ and $\tcalq^{(2)}$ imply
that $\tQ$ is contained in either
$\tcalq^{(1)}$ or $\tcalq^{(2)}$. If $\tQ\subset\tcalq^{(1)}$ then
$\tQ$ is a component of $\tcalq^{(1)}$, so that $Q:=|p|(\tQ)$ is a
component of $\calq$. The components $|p|(\tC_0)$ and $|p|(\tC^*)$ of $Q$
are respectively components of $Y_0$
and $\Ygeo^{(1)}$, so 
%that
%where $C_0(Q)$
%and
%$\Cgeo(Q)$ are respectively components of $Y_0$ and $\Ygeo^{(1)}$
 that $C_0(Q)=|p|(\tC_0)$ and $\Cgeo(Q)=|p|(\tCgeo)$. 
If $\tQ\subset\tcalq^{(2)}$ then
the component $\tQ$  of $\tcalq^{(2)}$ containing $\tQ$ has the form
$\tQ\cup\tQ'$, where $\tQ'\in\scrQ$, and $\tQ\cap\tQ'$ is a component
of  $\tYgeo$, which must be
$\tCgeo$ since $\tC_0$ is disjoint from $\tYgeo$ by
(\ref{reenergized}). It then follows from the discussion above that 
$Q:=|p|(\tQ\cup\tQ')$ 
is a
component of $\calq$,  that $|p|(\tC_0)=\partial Q=C_0(Q)$, and that
$|p|(\tCgeo)=\Cgeo(Q)$. Thus:

\Claim\label{pasketti}
If $\tQ$ is
%\subset\toldLambda$, which is by
%definition 
an element of $\scrQ$, then $|p|$ maps the component of $\partial\tQ$ contained
in $\tY_0$ onto $C_0(Q)$, and maps the component of $\partial\tQ$ contained
in $\tYgeo$ onto $\Cgeo(Q)$, for some component $Q$ of $\calq$.
\EndClaim

We claim:

\Claim\label{petit monde}
Every component of $Y_0$ has the form $C_0(Q)$ for some component $Q$
of $\calq$; every component of $\Ygeo$ has the form $\Cgeo(Q)$ for some component $Q$
of $\calq$; and every component of $Y^*$ has the form $C^*(Q)$ for some component $Q$
of $\calq$.
\EndClaim

To prove \ref{petit monde}, first suppose that $C_0$ is a component of
$Y_0$. According to \ref{cafe du monde}, there exist an element $\tQ$ of $\scrQ$ and a component $\tC_0$ of $\partial\tQ$ such that $|p|(\tC_0)=C_0$. It then follows from \ref{pasketti}
that $C_0=|p|(\tC_0)=C_0(Q)$ for some component $Q$ of $\calq$.

Next suppose that $\Cgeo$ is a component of
$\Ygeo$. Choose a component $\tCgeo$ of
$|p|^{-1}(\Cgeo)$. We may write $\tCgeo=\gamma_{\tC}(\tC)$ for some
component $\tC$ of $\tY$. 
%Since $Y_0$ and $Y$ are isotopic, $\tC$ is isotopic to a
%component $\tC_0$ of $Y_0$. 
Then $\tCgeo$ 
%is a simple 
%closed geodesic which
 is homotopic to $\tC$ and therefore to a
component $\tC_0$ of $|p|^{-1}(Y_0)$. Since $\tCgeo\subset\tYgeo$
and $\tC_0\subset\tY_0$ are disjoint by (\ref{reenergized}), it now
follows from \cite[Lemma 2.4]{epstein}  that $\tCgeo$ and $\tC_0$
cobound an annulus $\tQ\subset\toldLambda$, which is by
definition an element of $\scrQ$. It then follows from \ref{pasketti}
that $\Cgeo=|p|(\tCgeo)=\Cgeo(Q)$ for some component $Q$ of $\calq$.

To prove the final assertion of \ref{petit monde}, suppose that $C^*$ is a component of $Y^*$. Then $C^*$ is a
component of either $\Ygeo^{(1)}$ or $\partial|\oldPi|$. In the
case where $C^*$ is a
component of  $\Ygeo^{(1)}$, then in particular $C^*$ is a
component of  $\Ygeo$, and by the part of \ref{petit monde}
already proved, it has the form $\Cgeo(Q)$ for some component $Q$ of
$\calq$. If  $Q$ were a component of $\calq^{(2)}$ 
we would have  $\wt Q=2$, a contradiction since
$C^*\subset\Ygeo^{(1)}$. Hence $Q$ is a component of
$\calq^{(1)}$, so that $C^*(Q)=\Cgeo(Q)=C^*$. 

Now consider the
case where $C^*$ is a
component of  $\partial|\oldPi|$. Then $C^*=\partial P$ for some
component $P$ of $|\oldPi|$, and $\obd(P)$ is a tubular neighborhood
of $\obd(\Cgeo)$ for some component $\Cgeo$ of $\Ygeo^{(2)}$.
By the part of \ref{petit monde}
already proved, we may write $\Cgeo=\Cgeo(Q)$ for some component $Q$ of
$\calq$. If  $Q$ were a component of $\calq^{(1)}$,
we would have  $\wt Q=0$, a contradiction since
$C^*\subset\Ygeo^{(2)}$. Hence $Q$ is a component of
$\calq^{(2)}$, and the component $P$ of $|\oldPi|$ containing
$\Cgeo=\Cgeo(Q)$ is by definition equal to $P_Q$. Thus we have
$C^*(Q)=\partial P_Q=\partial P=C^*$. This completes the proof of
\ref{petit monde}.

According to \ref{youth in asia}, $(T_Q)_{Q\in\calc(\calq)}$ is a
pairwise disjoint family of 
%For every component $Q$ of $\calq$, the set $T_Q$ is a 
weight-$0$
annuli, and $C_0(Q)$ and $C^*(Q)\subset Y^*$ are the
boundary components of $T_Q$ for each $Q\in\calq$. It follows from
\ref{youth in asia} and \ref{petit monde} that
$\bigcup_{Q\in\calc(\calq)}C_0(Q)= Y_0$ and that
$\bigcup_{Q\in\calc(\calq)}C^*(Q)= Y^*$. This implies that $Y_0$ and
$Y^*$ are isotopic in $|\oldLambda|-\fraks_\oldLambda$. Since $Y_0$ is by definition isotopic
to $Y$ in $|\oldLambda|-\fraks_\oldLambda$, we deduce:

\Claim\label{very old-fashioned}
The $1$-manifolds $Y$ and
  $Y^*$ are isotopic in $|\oldLambda|-\fraks_\oldLambda$. 
\EndClaim

According to \ref{very old-fashioned}, we may fix a(n
orbifold) self-diffeomorphism $h$ of $\oldLambda$, smoothly isotopic to
the identity, such that $|h|(Y)=Y^*$. Then we have $Y=Y_1\cup Y_2$,
where $Y_1=|h|^{-1}(\Ygeo^{(1)})$ and
$Y_2=|h|^{-1}(\partial|\oldPi|)$. Hence we may write $E=E_1\discup E_2$
  in such a way that for each $i\in\{1,2\}$ and each $c\in E_i$ we have
$C_c\subset Y_i$. Then for each $c\in E_1$, the  homeomorphism $|h|$
maps $C_c$ 
 onto a component of $\Ygeo^{(1)}$, which we denote by
$\Ccgeo$. For each $c\in E_2$, it 
maps $C_c$ onto $|\partial\oldPi_c|$ for some component $\oldPi_c$ of $\oldPi$; the
unique component of $\Ygeo^{(2)}$, contained in $|\oldPi_c|$ will be denoted by $\Ccgeo$.  Thus
$\Ygeo$ may be written as a disjoint union $\bigcup_{c\in
  E}\Ccgeo$. 
%The definition of $C_c$ for any $c\in E$ is independent
%of the choice of $\oldPi$ within its isotopy class in $\otheroldLambda-\oldUpsilon^{(1)}$.
%\frakB
%.  
For each $\epsilon>0$ and each $c\in E$, let $\frakN_{\epsilon,c}$
denote the closed $\epsilon$-neighborhood of $\obd(\Ccgeo)$. Let us fix an
$\epsilon$ which is small enough to guarantee that the family
$(\frakN_{\epsilon,c})_{c\in E}$ is pairwise disjoint, and that
$\frakN_{\epsilon,c}$ is a geometric tubular neighborhood of $\obd(\Ccgeo)$
for each $c\in E$;  the latter condition means that the nearest point
map  $\frakN_{\epsilon,c}\to\obd(\Ccgeo)$ is well-defined and is an
$I$-fibration. This $I$-fibration is trivial if $c\in E_1$ and is
nontrivial if $c\in E_2$. For each $c\in E$, let us set $k_c=\lfloor
n_c/2\rfloor$ if $c\in E_1$, and $k_c=n_c$ if $c\in E_2$.
For each $c\in E$ and each positive real number $s$, let $K_{c,s}$ denote the set of all
points in $\frakN_{\epsilon,c}$ whose minimum distance from $\obd(\Ccgeo)$ is equal to $s$. Let us define a subset $X_c$ of
$[0,\epsilon)$ by setting $X_c=\{\epsilon/(k_c+1),
2\epsilon/(k_c+1),\ldots, k_c\epsilon/(k_c+1)\}$ if $n_c$ is even or
$c\in E_2$, and $X_c=\{0,\epsilon/(k_c+1),
2\epsilon/(k_c+1),\ldots, k_c\epsilon/(k_c+1)\}$ if $n_c$ is odd and
$c\in E_1$. 
For each $c\in E$, set $W^!_c=\bigcup_{s\in
  X_c}K_{c,s}$. Note that if $c\in E$ and $s\in(0,\epsilon)$, then
$K_{c,s}$ intersects each fiber of the $I$-fibration
$\frakN_{\epsilon,c}\to\obd(\Ccgeo)$ in exactly two points; whereas for any
$c\in E$,
the intersection of
$K_{c,0}$ with each fiber of 
$\frakN_{\epsilon,c}$ consists of exactly one point.
Hence, given any $c\in E$ and any $s\in[0.\epsilon)$, we have
$\compnum(K_{c,s})=2$ if the $I$-fibration of $\frakN_{\epsilon,c}$ is
trivial and $s>0$, while $\compnum(K_{c,s})=1$ if the $I$-fibration of $\frakN_{\epsilon,c}$ is
nontrivial or $s=0$. It follows that for any $c\in E$ we have
$\compnum( W^!_c)=2\card X_c$ if
$c\in E_1$
%I$-fibration of $\frakN_{\epsilon,c}$ is
%trivial
 and $0\notin X_c$, while $\compnum( W^!_c)=2(\card X_c)-1$ if
$c\in E_1$
% the $I$-fibration of $\frakN_{\epsilon,c}$ is
%trivial
 and $0\in X_c$, and $\compnum( W^!_c)=\card X_c$ if
$c\in E_2$ 
% the $I$-fibration of $\frakN_{\epsilon,c}$ is
%trivial and
(in which case $0\notin X_c$). In view of the definitions of $k_c$ and
$X_c$, it follows that 
\Equation\label{bug it}
\compnum(W^!_c)=n_c \text{ for every }c\in E.
\EndEquation

Now fix a real number $\delta$ such that $0<\delta<\epsilon/(k_c+1)$
for every $c\in E_2$. Set $R_c^!=\frakN_{\epsilon,c}$ for every $c\in
E_1$, and $R_c^!=\overline{\frakN_{\epsilon,c}-\frakN_{\delta,c}}$ for
  every $c\in E_2$. Then for each $c\in E$, the submanifold $R_c^!$ of
  $|\oldLambda|-\fraks_\oldLambda$ is a
 common tubular neighborhood of all the components of $W^!_c$. 

For each $c\in E_2$, a core curve of the annulus $R_c^!$ is the
boundary of the underlying manifold of a tubular neighborhood of
$\Ccgeo$;  the core curve $h(C_c)$ of $h(R_c)$
also bounds the underlying manifold of a tubular 
neighborhood of $\Ccgeo$, namely $\oldPi_c$. Likewise, for each $c\in
E_2$, both $h(R_c)$ and $R_c^!$ are annuli having $\Ccgeo$ as a core curve.
Furthermore, if we set $R^+_c=R_c$ for each $c\in E_1$, and $R^+_c=R_c\cup|\oldPi_c|$ for each $c\in E_2$, by taking $\epsilon$ sufficiently small we can guarantee
that $R_c^!$ lies in an arbitrarily small neighborhood of $h(R^+_c)$ for
each $c\in E$. Since $(h(R^+_c))_{c\in E}$ is a disjoint family, it
follows that there is a smooth isotopy of the orbifold $\oldLambda$
which
% is constant on $\bigcup_{c\in E_1}R_c^!$ and 
carries $h(R_c)$ onto $R_c^!$ for every
 $c\in E$. Using this isotopy we may modify $h$ to obtain a
 self-diffeomorphism $h'$ of $\oldLambda$, smoothly isotopic to the
 identity, such that $h'(R_c)=R_c^!$ for every $c\in E$.

 For every $c\in E$, the
components of $B$ belonging to the equivalence class $c$ are core
curves of the weight-$0$ annulus $R_c$, while the
components of $W^!_c$ are core curves of $R_c^!$; furthermore, we have 
$\card c=n_c=\compnum( W^!_c)$ by (\ref{bug it}). We may therefore
choose $h'$ so that $|h'|(\bigcup_{C\in c}C)=W^!_c$ for every $c\in
E$. If we set $B^!=\bigcup_{c\in E}W^!_c$, it follows that
$|h'|(B)=B^!$. Hence:

\Claim\label{occasional cortex} 
The  $1$-manifolds
$B$ and $B^!$ are (smoothly) orbifold-isotopic in $\oldLambda$.
\EndClaim

Now we claim:
\Claim\label{allons enfants}
 The $1$-manifold $B^!$ is $G$-invariant.
\EndClaim

To prove \ref{allons enfants}, consider an arbitrary component $C^!$
of $B^!$, and an arbitrary element $g$ of $G$. Then $C^!$ is a
component of $W^!_c$ for some $c\in E$. Let $j\in\{1,2\}$ denote the
index such that $c\in E_j$. Then $C_c\subset Y_j$, and hence
$\Ccgeo\subset\Ygeo^{(j)}$. It then follows from \ref{den sum}
that $g\cdot\Ccgeo$ is a component of $\Ygeo^{(j)}$. We may therefore
write $g\cdot\Ccgeo=\Ccprimegeo$ for some $c'\in E_j$. If $j=1$, we have
$h(\obd(C_c))=\obd(\Ccgeo)$; if $j=2$ then
$h(\obd(C_c))=\partial\oldPi_c$. Likewise, if $j=1$ we have
$h(\obd(C_{c'}))=\obd(\Ccprimegeo)=g\cdot\obd(\Ccgeo)$, and if $j=2$ we have
$h(\obd(C_{c'}))=\partial\oldPi_{c'}=g\cdot\partial\oldPi_c$, where the last
equality follows from the fact that $\oldPi_c$ and $\oldPi_{c'}$ are
the boundaries of the components of the $G$-invariant
orbifold $\oldPi$ that contain   $\obd(\Ccgeo)$ and
$\obd(\Ccprimegeo)$ respectively. Hence in any case we have
$h(\obd(C_{c'}))=g\cdot h(\obd(C_c))$. But since $h$ is orbifold-isotopic to the
identity of $\oldLambda$, the curve $|h|(C_{c})$ is homotopic to
$C_c$ in
$|\oldLambda|-\fraks_\oldLambda$, and $|h|(C_{c'})$ is
homotopic to $C_{c'}$ in
$|\oldLambda|-\fraks_\oldLambda$. Hence $C_{c'}$ is
homotopic to $g\cdot C_c$ in $|\oldLambda|-\fraks_\oldLambda$.

Now $n_{c'}$ is the number of components of $B$ that are
homotopic in
$|\oldLambda|-\fraks_\oldLambda$ to $C_{c'}$, or
equivalently to $g\cdot C_c$; hence $n_{c'}$ is also the number of
components of $g^{-1}\cdot B$ that are
homotopic in
$|\oldLambda|-\fraks_\oldLambda$ to $C_c$. Since by
hypothesis $g^{-1}\cdot\frakJ$ is orbifold-isotopic to $\frakJ$
in
$\otheroldLambda$, the $1$-manifold  $g^{-1}\cdot B=\partial|g^{-1}\cdot\frakJsmooth|$ is manifold-isotopic in
$|\oldLambda|-\fraks_\oldLambda$ to
$B=\partial|\frakJsmooth|$. It follows that $n_{c'}$ is the number of components of $B$ that are
homotopic in
$|\oldLambda|-\fraks_\oldLambda$ to $C_c$,
i.e. $n_{c'}=n_c$. Since \ref{den sum} implies that $c$ and $c'$ are either both in $E_1$ or both
in $E_2$, it follows from the definitions that $k_{c'}=k_c$ and that
$X_{c'}=X_c$.

Since $G$ acts by isometries and $g\cdot\obd(\Ccgeo)=\obd(\Ccprimegeo)$, we have
$g\cdot K_{c,s}=K_{c',s}$ for each real number $s\ge0$. Hence 
$g\cdot
C^!\subset g\cdot W^!_c=g\cdot(\bigcup_{s\in
  X_c} K_{c,s})=\bigcup_{s\in X_c}K_{c',s}$. Since $X_{c'}=X_c$, we obtain
$g\cdot
C^!\subset \bigcup_{s\in X_{c'}}K_{c',s}=W^!_{c'}\subset B^!$, and
\ref{allons enfants} is proved.

In view of \ref{allons enfants}, $\obd(B^!)/G$ is a well-defined
smooth $1$-suborbifold of $\oldLambda/G$.  Let $\frakV$ be a
PL $1$-suborbifold of $\otheroldLambda/G$ that is piecewise smoothly
isotopic to $\obd(B^!)/G$. Then the preimage of $\frakV$ under the
quotient map $\otheroldLambda\to\otheroldLambda/G$ is a $G$-invariant
PL $1$-suborbifold $\frakB^!$  of $\otheroldLambda$. It is
piecewise-smoothly isotopic in $\otheroldLambda$ to $\obd(B^!)$, which
in turn piecewise-smoothly isotopic to $\obd(B)$ by \ref{occasional
  cortex}. Hence $\frakB^!$ and $\obd(B)=\partial\frakJ$ are PL
isotopic in $\otheroldLambda$. This implies that $\frakJ$ is PL
isotopic in $\otheroldLambda$ to a $2$-suborbifold $\frakJ^!$ with $\partial\frakJ^!=\frakB^!$.  To prove the proposition it therefore remains only to prove that $\frakJ^!$ is $G$-invariant.

%We now turn to the proof that $\frakJ^!$ is $G$-invariant, which will complete the %proof of the proposition. 
It suffices to show that if $g\in G$ is given,  and  $\frakZ$ is any component of $\otheroldLambda$, then $g\cdot (\frakJ^!\cap\frakZ)\subset \frakJ^!$. We  set $\frakZ'=g\cdot\frakZ$, $\oldUpsilon=\frakJ^!\cap\frakZ$, and
$\oldUpsilon'=\frakJ^!\cap\frakZ'$. 
Since $\frakB^!$ is $G$-invariant, we have $\partial( g\cdot\oldUpsilon)=
\partial (g\cdot(\frakJ^!\cap\frakZ))=(g\cdot\partial\frakJ^!)\cap\frakZ'=
(g\cdot \frakB^!)\cap\frakZ'
=\frakB^!\cap\frakZ'$.
%$g\cdot(\frakB^!)=\frakB^!$. 
But since $\frakB^!\cap\frakZ'=\partial\oldUpsilon'$
%(\frakJ^!\cap\frakZ')$ 
is a two-sided closed suborbifold of the connected orbifold $\frakZ'$, 
%$\otheroldLambda$, \redcomment{Oops! It was not assumed connected} 
the only suborbifolds of $\frakZ'$ having boundary $\frakB^!\cap\frakZ'$ are 
$\oldUpsilon'$ and $\overline{\frakZ'-\oldUpsilon'}$.
%$\frakJ^!\cap\frakZ'$ and $\overline{\frakZ'\setminus \frakJ^!}$. 
If $g\cdot\oldUpsilon=\oldUpsilon'$, then
%$g\cdot(\frakJ^!\cap\frakZ)=
%\frakJ^!\cap\frakZ'$, then 
$g\cdot \oldUpsilon\subset \frakJ^!$,
%(\frakJ^!\cap\frakZ')\subset \frakJ^!$, 
as required. We will assume that $g\cdot \oldUpsilon=\overline{\frakZ'-\oldUpsilon'}$, %$g\cdot(\frakJ^!\cap\frakZ)=
%\overline{\frakZ'\setminus \frakJ^!}$
%\overline{\frakB^!-\frakJ}$ 
and obtain a contradiction.

Since $\frakJ^!$ is isotopic in $\otheroldLambda$ to $\frakJ$, the suborbifold $\oldUpsilon=\frakJ^!\cap\frakZ$ is isotopic in $\frakZ$ to $\frakU:=\frakJ\cap\frakZ$. Likewise, $\oldUpsilon'=\frakJ^!\cap\frakZ'$ is isotopic in $\frakZ'$ to $\frakU':=\frakJ\cap\frakZ'$. 
Since $\oldUpsilon$ is isotopic to $\frakU$, the suborbifold $g\cdot\oldUpsilon$ is isotopic in $\frakZ'=g\cdot\frakZ$ to $g\cdot\frakU$.
On the other hand, the hypothesis that $g\cdot\frakJ$ is isotopic in $\otheroldLambda$ to $\frakJ$ implies that $g\cdot\frakU$ is isotopic in $\frakZ'$ to $\frakU'$.
It now follows that $g\cdot\oldUpsilon$ is isotopic to $\oldUpsilon'$. The assumption $g\cdot \oldUpsilon=\overline{\frakZ'-\oldUpsilon'}$ therefore implies that $\oldUpsilon'$ is isotopic to $\overline{\frakZ'-\oldUpsilon'}$.

Since $\frakJ$ is taut, every component of $\frakJ^!$, and in particular every component of $\oldUpsilon'$, has non-positive Euler characteristic. First consider the case in which every component of $\oldUpsilon'$ has Euler characteristic $0$. Since $\oldUpsilon'$ is isotopic to $\overline{\frakZ'-\oldUpsilon'}$, every component of 
$\overline{\frakZ'-\oldUpsilon'}$ also
has Euler characteristic $0$. But then $\chi(\frakZ')=\chi(\oldUpsilon')+\chi(
\overline{\frakZ'-\oldUpsilon'})=0$. This contradicts the hypothesis that $\otheroldLambda$ is negative.

Now consider the case in which some component $\oldGamma$ of $\oldUpsilon'$ has negative Euler characteristic. Since $\oldUpsilon'$ is isotopic to $\overline{\frakZ'-\oldUpsilon'}$, the component $\oldGamma$ is isotopic to some suborbifold $\oldGamma_1$ of $\frakZ'$ such that $\oldGamma_1\cap\oldGamma=\emptyset$.
It then follows from Corollary \ref{nafta} that $[\oldGamma_1]\wedge[\oldGamma]=[\emptyset]$. On the other hand, since $\oldGamma_1$ and $\oldGamma$ are isotopic, we have $[\oldGamma_1]=[\oldGamma]$. 
Hence 
$[\oldGamma_1]\wedge[\oldGamma]=[\oldGamma]\wedge[\oldGamma]=[\oldGamma]$. 
It
now follows that $[\oldGamma]=[\emptyset]$, which is impossible since $\oldGamma$ is a connected (and hence non-empty) negative suborbifold of $\otheroldLambda$. 
\EndProof
%\otheroldLambda\oldUspilon_0Y_0^+'\frakJ_0\oldDelta
%fffffffffffff[C'j_1\frakD\oldDelta\frakZ\oldGamma\frakG\Lambda\frakB_1\frakR_c\frakZ\oldDelta\oldUpsilon\otheroldLambda
%h_1(\frakJ)\frakJ'\frakB'\frakV\frakB\frakC\frakD\oldDelta\fra ,h ^0
%\frakR $p Y \frakZ \frakQ X\frakC \oldPi \frakP \frakQ $N$ G_1 q
%\compnum Y_0 Y_1 %\frakB_1 \obd(B^!) \frakB \frakJ_1 D L a_ B E Q
%\frakA Q' g( p\oldUpsilon_i\oldUpsilon_2\oldUpsilon_i\oldUpsilon_j
%\frakB \frakC |\oldUpsilon|\tYgeo\oldUpsilon^{(1)}
%\oldUpsilon^{(2)}\oldUpsilon^{(i)}\oldUpsilon^{(j)} diffeo homeo
%\oldLambda \toldLambda \frakJ \frakC\oldPi $T R_Q X W \frakW B^! B_

%\gamma
%\label\beta\gamma

\section{Higher characteristic $2$-orbifolds}\label{higher section}

\Number\label{what's iota?}
If $\oldPsi$ is an orientable $3$-orbifold which is componentwise
strongly 
\simple\ and componentwise boundary-irreducible, we will denote by $\oldSigma^-=\oldSigma^-(\oldPsi)$
the union of all components of
$\oldSigma(\oldPsi)$ that have strictly negative (orbifold) Euler
characteristic, and by $\oldPhi^-=\oldPhi^- (\oldPsi)$ the union of
all components of $\oldPhi(\oldPsi)$ that have strictly negative
 Euler characteristic. Since every \bindinglike\ connected
\Ssuborbifold\ of $\oldPsi$ has Euler characteristic $0$ by Lemma \ref{when a tore a fold}, every
component of $\oldSigma^-$ is a \pagelike\ \Ssuborbifold\ of
$\oldPsi$. Note also  that by \ref{tuesa day}, if $\oldLambda$ is any \Ssuborbifold\ of $\oldPsi$ then each component of $\oldLambda\cap\partial\oldPsi$ has either the same Euler characteristic as $\oldLambda$, or twice the Euler characteristic of $\oldLambda$. Hence we have
$\oldSigma^-\cap\partial\oldPsi=\oldPhi^-$.

In view of the definition of a \pagelike\ \Ssuborbifold,
$\oldSigma^-$
may be equipped with an $I$-fibration $q:\oldSigma^-\to \frakB$ over a (possibly disconnected) $2$-orbifold in
such a way that $\partialh\oldSigma^-=\oldPhi^-$. By
\ref{fibered stuff}, $q|\oldPhi^-:\oldPhi^-\to \frakB$ is a two-sheeted covering map. Its non-trivial deck transformation is 
%restriction of the 
%we may define an involution
an involution of $\oldPhi^-$ which will be denoted by $\iota_\oldPsi$. While this definition of $\iota_\oldPsi$  depends on the $I$-fibration $q$, it follows from Corollary \ref{flubadub o'connor} that the strong equivalence class of $\iota_\oldPsi$ in $\oldPhi^-$ is independent of the choice of the $I$-fibration $q$. Since $\oldPhi^-$ is itself well defined up to isotopy in $\partial\oldPsi$, the involution $\iota_\oldPsi$  is well defined up to strong equivalence in $\partial\oldPsi$.

The definition of $\iota_\oldPsi$ may be paraphrased by saying that it is the involution of $\oldPhi^-$ which interchanges the endpoints of each generic fiber of $\oldSigma^-$, i.e. each fiber which is orbifold-homeomorphic to $[0,1]$. Hence the intersection of $\partial\oldPsi$ with each saturated suborbifold of $\oldSigma^-$ is $\iota_\oldPsi$-invariant; in particular $\frakA(\oldPsi)\cap\oldSigma^-\cap\partial\oldPsi$, which is the intersection of $\partial\oldPsi$ with the vertical boundary of $\oldSigma^-$, is $\iota_\oldPsi$-invariant. Conversely, if a suborbifold $\frakC$ of $\oldPhi^-$ is $\iota_\oldPsi$-invariant, then the union of all fibers of $\oldSigma^-$ that have their endpoints in $\frakC$ is a saturated suborbifold of $\oldSigma^-$ whose intersection with $\oldPsi^-$ is $\frakC$.

Note also that if $\frakU$ is a component of $\oldSigma^-$ whose intersection with $\partial\oldPsi$ has two components (so that these two components are the components of $\partialh\frakU$ under the fibration $q$), then the components of 
$\frakU\cap\partial\oldPsi$
are interchanged by  $\iota_\oldPsi$.
\EndNumber

\Number\label{in duck tape}
Now let $\Mh$ be
a closed, hyperbolic, orientable $3$-orbifold, and set $\oldOmega=(\Mh)\pl$. Let
 $\oldTheta$ be an incompressible, closed $2$-suborbifold
of $\oldOmega$.
Set
$\oldPsi=\oldOmega\cut\oldTheta$ and ${\toldTheta}=\partial\oldPsi$. Since $\oldTheta$ is incompressible, 
$\oldPsi$ is componentwise strongly 
\simple\ and componentwise
boundary-irreducible by \ref{oops lemma}. Hence 
by
\ref{oldSigma def} and \ref{tuesa day}, $\oldSigma(\oldPsi)$, 
$\oldPhi(\oldPsi)$,
$\kish(\oldPsi)$, $\book(\oldPsi)$ and
$\frakA(\oldPsi)$ are well defined. According to \ref{boundary is negative}, $\toldTheta$ is negative. These facts---the 
componentwise strong 
\simple ity and componentwise
boundary-irreducibility of
$\oldPsi$, the  well definedness of $\oldSigma(\oldPsi)$, 
$\oldPhi(\oldPsi)$,
$\kish(\oldPsi)$, $\book(\oldPsi)$ and
$\frakA(\oldPsi)$, and the negativity of $\toldTheta$---will be used, often without being mentioned explicitly, not only in this subsection but in its applications throughout the remainder of the monograph.

Set 
$\tau=\tau_\oldTheta$. Thus $\tau$ is an involution of ${\toldTheta}$.
% which
%restricts to an involution of ${\toldTheta}$. 
Set
$\oldPhi^-=\oldPhi^-(\oldPsi)$. According to
\ref{tuesa day}, $\oldPhi^-$
is a taut suborbifold  of ${\toldTheta}$; it is negative by the
definition of $\oldPhi^-$, and hence belongs to
$\Theta_-({\toldTheta})$ (see \ref{chidef}).  Set $\iota=\iota_\oldPsi$, 
which
according to \ref{what's iota?} is an involution of $\oldPhi^-$, well defined up to strong equivalence. It follows from the definition of strong equivalence that %According to an observation made in Definition \ref{strong equivalence},  
$\iota$ is in particular
%of 
%$\calm({\toldTheta},{\toldTheta})$
%$(\oldPhi^-,\iota,\oldPhi^-)$ is in particular 
well-defined up to isotopy in $\toldTheta$, so that $[\iota]=[\oldPhi^-,\iota,\oldPhi^-]$ is a well-defined element of $\barcalm({\toldTheta},{\toldTheta})$ (see \ref{my remark}). We also have an element $[\tau]=[{\toldTheta},\tau,{\toldTheta}]$ of $\barcalm({\toldTheta},{\toldTheta})$. 
%Thus
%and $({\toldTheta},\tau,{\toldTheta})$ are elements of 
%(see \ref{my remark}), so that and are elements of
%\redcomment{I want to say that $[\iota]$ is well defined because the strong equivalence class determines the isotopy class. State this in some decent way.} 
Since $[\tau]$ and $[\iota]$ are involutions of
their respective domains, we have $[\tau]^{-1}=[\tau]$ and
$[\iota]^{-1}=[\iota]$.

For every integer $n\ge1$, we will set
$\noodge_n=\noodge_n^{\oldOmega,\oldTheta}=([\iota]\diamond[\tau])^{\diamond
  (n-1)}\diamond[\iota]\in\barcalm({\toldTheta},{\toldTheta})$. In view of \ref{more associativity}, we may
write $\noodge_n=[\iota]\diamond([\tau]\diamond[\iota])^{\diamond (n-1)}$.
By (\ref{snicker snoodle}), and
the equalities $[\tau]^{-1}=[\tau]$ and
$[\iota]^{-1}=[\iota]$, we have
$$\noodge_n^{-1}=
([\iota]^{-1}\diamond[\tau]^{-1})^{\diamond(n-1)}\diamond[\iota]^{-1}
=
([\iota]\diamond[\tau])^{\diamond(n-1)}\diamond[\iota]=j_n.$$
%([\iota]\diamond[\tau])^{\diamond k}\diamond
%\iota\diamond([\tau]\diamond[\iota])^{\diamond k}
%=^{-1}
%Note that if $n$ is odd and $k=(n-1)/2$, then by (\ref{snicker snoodle}), and
%the equalities $[\tau]^{-1}=[\tau]$ and
%$[\iota]^{-1}=[\iota]$, we have
%$$\noodge_n=
%([\iota]\diamond[\tau])^{\diamond k}\diamond
%\iota\diamond([\tau]\diamond[\iota])^{\diamond k}
%=
%([\iota]\diamond[\tau])^{\diamond k}\diamond
%\iota\diamond([\iota]\diamond[\tau])^{\diamond(- k)}, 
%$$
%which with (\ref{snicker snoodle}) gives
%$\noodge_n^{-1}=
%([\iota]\diamond[\tau])^{\diamond k}\diamond
%\iota^{-1}\diamond([\iota]\diamond[\tau])^{\diamond(-k)}=\noodge_n$.
%
%If $n=2k$ is even,  we have
%$$
%\begin{aligned}
%\noodge_n&=
%(([\iota]\diamond[\tau])^{\diamond (k-1))}\diamond[\iota])\diamond
%\tau\diamond
%([\iota]\diamond([\tau]\diamond[\iota]) ^{\diamond (k-1)})\\
%&=
%(([\iota]\diamond[\tau])^{\diamond (k-1))})\diamond[\iota])\diamond
%\tau\diamond
%(([\iota]\diamond[\tau])^{\diamond (k-1))})\diamond[\iota])^{-1},
%\end{aligned}
%$$
%so that
%$\noodge_n^{-1}=
%(([\iota]\diamond[\tau])^{\diamond (k-1))})\diamond[\iota])\diamond
%\tau^{-1}\diamond
%(([\iota]\diamond[\tau])^{\diamond (k-1))})\diamond[\iota])^{-1}=\noodge_n$.
This shows:
\Equation\label{yes virginia}
\noodge_n^{-1}=\noodge_n\qquad\text{for every }n\ge1.
\EndEquation
%\redcomment{I doubt that the division into odd and even cases is really needed. Rethink this.}
According to \ref{my remark} we have
$\dom\noodge_n=\range\noodge_n^{-1}=\range\noodge_n$. We shall write 
$V_n=V_n^{\oldOmega,\oldTheta}=\dom\noodge_n^{\oldOmega,\oldTheta}=\range\noodge_n^{\oldOmega,\oldTheta}\in\barcaly_-(\toldTheta)$.

Note that the definition of the $\noodge_n$, with \ref{more associativity}, implies that for every $n\ge1$ we have $\noodge_{n+1}=([\iota]\diamond[\tau])\diamond\noodge_n$. It follows from Lemma \ref{before associativity} that  
$\dom(([\iota]\diamond[\tau])\diamond\noodge_n)\preceq\dom\noodge_n$, i.e. $V_{n+1}\preceq V_n$. Hence
\Equation\label{i thought so too}
V_{n'}\preceq V_n\qquad\text {whenever }n'\ge n\ge1.
\EndEquation
The definition also implies that 
\Equation\label{and so did you}
\noodge_1=[\iota], \text{ and hence }V_1=[\oldPhi^-].
\EndEquation
\EndNumber
%taut

\Lemma\label{j-lemma}
Let $\Mh$ be a closed, hyperbolic, orientable $3$-orbifold, and let $\oldTheta$ be an incompressible, closed $2$-suborbifold
of $\oldOmega:=(\Mh)\pl$.
Set
$\oldPsi=\oldOmega\cut\oldTheta$, ${\toldTheta}=\partial\oldPsi$, and
$\tau=\tau_\oldTheta$; and for every $n\ge1$ set
$\noodge_n=\noodge_n^{\oldOmega,\oldTheta}$ and
$V_n=V_n^{\oldOmega,\oldTheta}$. Let $Z\in\barcaly_-(\oldTheta)$ be
given, and let $m$ and $k$ be integers with $k>m>1$. Then the following conditions are equivalent:
\begin{enumerate}[(a)]
\item $Z\preceq V_k$;
\item $Z\preceq V_m$ and $[\tau](\noodge_m(Z))\preceq V_{k-m}$.
\end{enumerate}
Furthermore, if (a) (or (b)) holds, then 
$$\noodge_k|Z=\noodge_{k-m}\circ [\tau]
\circ(\noodge_m|Z),$$
where the threefold composition on the right hand side is unambiguously defined in view of 
(\ref{easy associativity}).
\EndLemma

\Proof
Recall that by definition (see \ref{in
  duck tape}) we have
$\noodge_k=([\iota]\diamond[\tau])^{\diamond
  (k-1)}\diamond[\iota]$, where $\iota=\iota_\oldPsi$.  In view of
\ref{more associativity}, we may write 
\Equation\label{chico}
\noodge_k=([\iota]\diamond[\tau])^{\diamond
  (k-m-1)}\diamond[\iota]\diamond[\tau]\diamond ([\iota]\diamond[\tau])^{\diamond
 ( m-1)}\diamond[\iota]=\noodge_{k-m}\diamond[\tau]\diamond\noodge_{m}.
\EndEquation

Now apply Lemma \ref{before associativity}, taking $Y=Z$, $\emm_1=\noodge_m$, and $\emm_2=\noodge_{k-m}\diamond[\tau]$. Since $\noodge_k= (\noodge_{k-m}\diamond[\tau])\diamond\noodge_{m}$ by (\ref{chico}), and since $V_k=\dom\noodge_k$ by \ref{in duck tape}, 
this application of Lemma \ref{before associativity} shows:
\Claim\label{boychik}
 Condition (a) of the present Lemma 
holds if and only if we have
%with these choices of $Y$, $\emm_1$ and $\emm_2$. Hence we have 
$Z\preceq\dom\noodge_m=V_m$ and 
$\noodge_m(Z)\preceq\dom (\noodge_{k-m}\diamond[\tau])$. Furthermore, if
Condition (a) does hold, then $
%\noodge_k=
\noodge_{k}|Z=(\noodge_{k-m}\diamond[\tau])\circ(\noodge_m|Z)$.
\EndClaim

Now, under the assumption that $Z\preceq V_m$, let us make a second application of Lemma \ref{before associativity}, this time taking $Y=\noodge_m(Z)$, and letting $[\tau]$ and $\noodge_{k-m}$ play the respective roles of 
$\emm_1$ and $\emm_2$. In this setting,
Condition (2) of Lemma \ref{before associativity}  asserts that
$\noodge_m(Z) \preceq\dom[\tau]$ and that
$[\tau](\noodge_m(Z))\preceq\dom \noodge_{k-m}=V_{k-m}$. But the relation
$\noodge_m(Z) \preceq\dom[\tau]$ holds automatically because $\dom[\tau]={\toldTheta}$. 
Hence this second application of Lemma \ref{before associativity} shows:
\Claim\label{bapchuk}
Under the assumption that $Z\preceq V_m$, we have $\noodge_m(Z)\preceq\dom (\noodge_{k-m}\diamond[\tau])$ if and only if $[\tau](\noodge_m(Z))\preceq V_{k-m}$; and furthermore, 
if these equivalent conditions hold, then 
$
(\noodge_{k-m}\diamond[\tau])|(\noodge_m(Z) )=\noodge_{k-m}\circ( [\tau]|
(\noodge_m(Z))$. 
\EndClaim

Taken together, the first assertion of \ref{boychik} and the first assertion of \ref{bapchuk} immediately imply that Conditions (a) and (b) of the present lemma are equivalent. The second assertion of \ref{boychik} and the second assertion of \ref{bapchuk} imply that if (a) (or (b)) holds then
$
\noodge_k|Z=\noodge_{k-m}\circ( [\tau]|
(\noodge_m(Z) )\circ(\noodge_m|Z)$, which in view of (\ref{number one})  may be rewritten as
$
\noodge_k|Z=\noodge_{k-m}\circ [\tau]
\circ(\noodge_m|Z)$. 
%\redcomment{I'm implicitly using (\ref{easy
    %associativity}); cf. the comment after the statement of the
  %lemma.} 
This establishes the final assertion of the lemma.
\EndProof
%\noodge_k

\Number\label{oldXi}
Suppose that
$\Mh$ is
a closed, hyperbolic, orientable $3$-orbifold, that
 $\oldTheta$ is an incompressible, closed $2$-suborbifold
of $\oldOmega:=(\Mh)\pl$, and that $\oldTheta'$ is a union of components of $\oldTheta$. 
Set
$\oldPsi=\oldOmega\cut\oldTheta$ and ${\toldTheta}=\partial\oldPsi$, and set
$\oldPsi'=\oldOmega\cut{\oldTheta'}$ and ${\toldTheta}'=\partial\oldPsi'$, so that
${\toldTheta}'$ is canonically identified with a union of components of
${\toldTheta}$. For each integer $n\ge1$ we will denote by 
 $\oldfrakX_n=\oldfrakX_n(\oldOmega,\oldTheta,\oldTheta')$
the set of all components (see
\ref{chidef}) $U$ of
$V_n:=V_n^{\oldOmega,\oldTheta}$ with the property that
$\noodge_n(U)\preceq[{\toldTheta}']$. We will set
$X_n=X_n(\oldOmega,\oldTheta,\oldTheta')=\bigvee_{U\in\oldfrakX_n}U$ (see \ref{distribute}).
For each $U\in\Xi_n$, we have $U\preceq
V_n$; by the definition of a supremum it follows that $X_n\preceq V_n$.
By (\ref{i thought so too}) and (\ref{and so did you}) it follows that $X_n\preceq V_1=[\oldPhi^-(\oldPsi)]$.

If we write $V_n=[\frakV_n]$ for some $\frakV_n\in\Theta_-({\toldTheta})$, the elements of $\oldfrakX_n$ may be
written as $[\oldXi_{n,1}],\ldots,[\oldXi_{n,p}]$ for certain components
$\oldXi_{n,1},\ldots,\oldXi_{n,p}$ of $\frakV_n$; in particular the
$\oldXi_{n,i}$ are pairwise disjoint, which implies by \ref{join and meet} that
$X_n=[\oldXi_{n,1}\cup\cdots\cup\oldXi_{n,p}]$. 
It follows that each component of $X_n$ is an
  element of $\oldfrakX_n$, and in particular a component of $V_n$.

We set $X=X(\oldOmega,\oldTheta,\oldTheta')=\bigvee_{n\ge1}X_n(\oldOmega,\oldTheta,\oldTheta')\in\barcaly_-({\toldTheta})$. We set $X'=X'(\oldOmega,\oldTheta,\oldTheta')=X(\oldOmega,\oldTheta,\oldTheta')\wedge[\toldTheta']$ and  $X''=X''(\oldOmega,\oldTheta,\oldTheta')=X(\oldOmega,\oldTheta,\oldTheta')\wedge[\toldTheta-\toldTheta']$. Note that
$X'$ and
$X''$ may be regarded as elements of
$\barcaly_-(\toldTheta')$ and $\barcaly_-(\toldTheta-\toldTheta')$
respectively. Since $X_n\preceq[\oldPhi^-(\oldPsi)]$ for each $n\ge1$, the definition of a supremum shows that $X\preceq [\oldPhi^-(\oldPsi)]$.
\EndNumber
%\frakV\oldPhi

\Lemma\label{here goethe}
Suppose that $\Mh$ is a closed, hyperbolic, orientable $3$-orbifold, $\oldTheta$ is an incompressible, closed $2$-suborbifold
of $\oldOmega:=(\Mh)\pl$, and that
$\oldTheta'$ is a union of components of $\oldTheta$. If $U$ and $W$ are respectively components of $X_m(\oldOmega,\oldTheta,\oldTheta')$ and $X_n(\oldOmega,\oldTheta,\oldTheta')$ for some (not
necessarily distinct) positive integers $m$ and $n$, then 
either (a) $U\preceq W$, (b) $W\preceq U$, or (c) there exist elements $\frakU,\oldUpsilon\in\Theta_-(\toldTheta)$ such that $[\frakU]=U$, $[\oldUpsilon]=W$, and $\frakU\cap\oldUpsilon=\emptyset$. Furthermore, each component of
$X(\oldOmega,\oldTheta,\oldTheta')$ is a component of
$X_n(\oldOmega,\oldTheta,\oldTheta')$ for some positive integer $n$.
\EndLemma

\Proof
To prove the first assertion, note that by \ref{oldXi}, $U$ and $W$
are respectively elements of $\oldfrakX_m$ and $\oldfrakX_n$ and are therefore
components of $V_m$ and $V_n$ respectively. Consider the case in which $m\le n$. In this case  we have $V_n\preceq
V_m$ by \ref{i thought so too}. In particular $W\preceq V_m$. Hence there are elements $\oldUpsilon$ and $\frakV_m$ of $\Theta_-(\oldTheta)$ such that
 $[\oldUpsilon]=W$, $[\frakV_m]=V_m$, and $\oldUpsilon\subset\frakV_m$. Since $W$ is connected, $\oldUpsilon$ is contained in a component $\oldUpsilon'$ of $\frakV_m$. Then $W':=[\oldUpsilon']$ is a component of $V_m$, and $W\preceq W'$.  On the other hand, since $U$ is a component of $V_m$, there is a component (see \ref{components}) $\frakU$ of $\frakV_m$ such that $[\frakU]=U$. If $\frakU=\oldUpsilon'$ then $W\preceq U$, and if $\frakU\ne\oldUpsilon'$ then $\frakU\cap\oldUpsilon\subset \frakU\cap\oldUpsilon'=\emptyset$. This proves that in the case 
% so that $U\wedge W=[\emptyset]$  by Assertion (2) of Corollary \ref{negative lattice}. Thus in the case
 $m\le n$, one of the alternatives (a) or (c) of the first assertion of the lemma holds.
% we have $U\preceq W$ or $U\wedge W=\emptyset$.
%and $W$ are 
%is connected, it is
%$\frakU$ and $\oldUpsilon$ are connected and is therefore 
%contained in a component $\oldUpsilon'$ of denotes the component of $\frakV_m$ containing 
%either $U\preceq W$, or there are
%disjoint elements $\frakU$ and $\oldUpsilon$ of $\Theta_-(\oldTheta)$ such
%that $[\frakU]=U$ and $[\oldUpsilon]=W$; in the latter case we have $U\wedge
%W=[\emptyset]$. 
The same argument shows that if $n\le m$ then one of the alternatives (b) or (c) holds, and so the first assertion is proved.
%either
%$W\preceq U$ or $U\wedge W=[\emptyset]$. This proves the first assertion.

To prove the second
assertion, let ${\mathscr U}_n$ denote the set of components (see \ref{components}) of $X_n$ for each $n\ge1$, and set 
$\mathscr U=\bigcup_{n=1}^\infty \mathscr U_n$. 
According to \ref{components}, we have $X_n=\bigvee_{U\in {\mathscr U}_n}U$ for each $n\ge1$, and the definition of $X$ gives $X=\bigvee_{n\ge1}X_n$. Hence $X=\bigvee_{U\in{\mathscr U}}U$.

We claim: 

\Claim\label{it stops}
There are elements $U_1,\cdots,U_N$ of $\mathscr U$, for some integer $N\ge0$, such that $X=U_1\vee\cdots\vee U_N$.
\EndClaim

To prove \ref{it stops}, note that by construction the set $\mathscr U$ is either finite or countably infinite. If $\mathscr U$ is finite, the assertion follows from the equality 
$X=\bigvee_{U\in{\mathscr U}_n}U$. If
$\mathscr U$ is infinite we may write $\mathscr U=\{U_1,U_2,\ldots\}$; setting $Z_k=U_1\vee\cdots\vee U_k$ for each $k\ge1$, we have $Z_1\preceq Z_2\preceq\cdots$, and it follows from Proposition \ref 
{ascending}
that there exists an index $N\ge1$ such that $Z_k=Z_N$ for every $k\ge N$. Hence $X=U_1\vee U_2\vee\cdots=Z_N=U_1\vee\cdots\vee U_N$, and the proof of \ref{it stops} is complete.

Fix elements $U_1,\cdots,U_N$ of $\mathscr U$  having the property stated in \ref{it stops}. After possibly replacing $\{U_1,\cdots,U_N\}$ by a subcollection, we may assume that for any distinct indices $i,j\in\{1,\ldots,N\}$ we have $U_i\not\preceq U_j$. Since $U_i\in\mathscr U$ for each $i\in\{1,\ldots,N\}$, it follows from the first assertion of the present lemma, which has already been proved, that
for any distinct indices $i,j\in\{1,\ldots,N\}$ there exist elements $\frakU$ and $\oldUpsilon$ of $\Theta_-(\toldTheta)$ such that $[\frakU]=U_i$, $[\oldUpsilon]=U_j$, and $\frakU\cap 
\oldUpsilon =\emptyset$. It then follows from Proposition \ref{all at once} that there are pairwise disjoint elements $\frakU_1,\ldots,\frakU_N$ of $\Theta_-(\toldTheta)$ such that $[\frakU_i]=U_i$ for $i=1,\ldots,N$. According to an observation made in \ref{join and meet}, it now follows that $U_1,\ldots,U_N$ are the components of $U_1\vee\cdots\vee U_N=X$. Since $U_1,\ldots,U_N\in\mathscr U$ by \ref{it stops}, this shows that each component of $X$ is an element of $\mathscr U$, which by the definition of $\mathscr U$ means that each component of $X$ is a component of $X_n$ for some $n$.
\EndProof
%VW\frakV\oldUpsilon

The importance of $X'(\oldOmega,\oldTheta,\oldTheta')$ arises from the following result:

\Proposition\label{herman had burped}
Let $\Mh$ be a closed, hyperbolic, orientable $3$-orbifold, and let $\oldTheta$ be an incompressible, closed $2$-suborbifold
of $\oldOmega:=(\Mh)\pl$. 
Set
$\oldPsi=\oldOmega\cut\oldTheta$ and ${\toldTheta}=\partial\oldPsi$.
%, and
%$\tau=\tau_\oldTheta$.
Let $\oldTheta'$ be a union of components of $\oldTheta$, and set
$\oldPsi'=\oldOmega\cut{\oldTheta'}$ and ${\toldTheta}'=\partial\oldPsi'$.
Then the relation
$[\oldPhi^-(\oldPsi')]\preceq
X'(\oldOmega,\oldTheta,\oldTheta')$ holds in $\barcaly_-(\toldTheta')$
(where $[\oldPhi^-(\oldPsi')]$ and $X'(\oldOmega,\oldTheta,\oldTheta')$,  are well-defined element of $\barcaly_-(\toldTheta')$ by \ref{in duck tape} and 
\ref{oldXi} respectively).
\EndProposition

\Remark
We believe that the conclusion of Proposition \ref{herman had burped} may be replaced by the stronger conclusion $[\oldPhi^-(\oldPsi')]=
X'(\oldOmega,\oldTheta,\oldTheta')$, although we have not written down a proof of this. The conclusion stated in Proposition \ref{herman had burped}, and proved below, will be sufficient for the applications in this monograph.
\EndRemark

\Proof[Proof of  Proposition \ref{herman had burped}]
% It may be that instead of showing 
%$[\oldPhi^-(\oldPsi')]=
%X'(\oldOmega,\oldTheta,\oldTheta')$, I will have to settle for 
Set $\oldTheta'':=\oldTheta-\oldTheta'$. Then $\oldTheta''$
is
canonically identified with a closed, incompressible suborbifold of
$\inter\oldPsi'$. Furthermore, $\oldPsi$ is canonically identified with $(\oldPsi')\cut{\oldTheta''}$. We also 
identify
${\toldTheta}'$  with a union of components of
${\toldTheta}$, and
set $\toldTheta'':=\toldTheta-\toldTheta'$; under our identifications we have 
$\toldTheta''=\rho_{\oldTheta''}^{-1}(\oldTheta'')$.
We set
$\rho=\rho_{\oldTheta''}:\oldPsi\to\oldPsi'$ and $\tau=\tau_{\oldTheta''}:\toldTheta''\to\toldTheta''$.

According to Proposition \ref{mystery}, it is enough to prove that if $\frakK$ is any component of $\oldPhi^-(\oldPsi')$ then $[\frakK] \preceq
X'(\oldOmega,\oldTheta,\oldTheta')$. Let $\oldLambda$ denote the component of $\oldSigma(\oldPsi')$ containing $\frakK$. Since $\chi(\frakK)<0$, we have $\chi(\oldLambda)<0$ by \ref{tuesa day}, and it then follows from  Lemma \ref{when a tore a fold} that $\oldLambda$ is a \pagelike\ \Ssuborbifold\ of $\oldPsi'$. We first apply Assertion (1) of
Proposition \ref{newer no disks lemma}, letting $\oldPsi'$, $\oldTheta''\subset\inter\oldPsi'$,
and $\Fr_{\oldPsi'}\oldLambda$ play the respective roles of $\oldPsi$,
$\oldPi$, and $\frakB$. This shows that, after modifying
$\oldSigma(\oldPsi')$ within its isotopy class in $\oldPsi'$, we may assume that
$\Fr_{\oldPsi'}\oldLambda$ has reduced intersection (see Definition \ref{reduced intersection}) with $\oldTheta''\subset\inter\oldPsi'$. Next, we apply Assertion (3) of Proposition \ref{newer no disks lemma} to conclude that each component of the
  suborbifold $\oldLambda\cap\oldTheta''$ of $\oldLambda$ is
parallel in the pair $(\oldLambda,\Fr\oldLambda)$
%parallel \redmissingref{rel something?} 
to $\frakK$.  Hence if $\oldGamma$ is any component of
$\oldLambda\cap\oldTheta''$, there is a unique connected suborbifold
$\frakD_\oldGamma$ of $\oldLambda$ such that $\Fr_\oldLambda
\frakD_\oldGamma=\oldGamma$ and $\frakD_\oldGamma\supset\frakK$. Furthermore, $\frakD_\oldGamma$ admits a trivial
$I$-fibration under which its the components of its horizontal
boundary are $\oldGamma$ and $\frakK$. The collection
$\{\frakD_\oldGamma:
\oldGamma\in\calc(\oldLambda\cap\oldTheta'')\}$ consists of
suborbifolds of $\oldLambda$ whose frontiers are connected and pairwise disjoint, and
they all contain $\frakK\subset\partial\oldLambda$; hence this collection is linearly ordered
by inclusion. We may therefore index the components of
$\oldLambda\cap\oldTheta''$ as $\oldGamma_1,\ldots,\oldGamma_m$, for
some $m\ge0$, in such 
a way that, if we set $\frakD_k=\frakD_{\oldGamma_k}$ for $1\le k\le m$, we have
$\frakD_1\supset\cdots\supset \frakD_k$. We set $\frakD_0=\oldLambda$ and
$\frakD_{m+1}=\emptyset$, and we set
$\frakU_k=\overline{\frakD_k-\frakD_{k+1}}$ for $k=0,\ldots,m$.
%, and we set
  %$\frakU_0=\overline{\oldLambda-\frakD_1}$ and
  %$\frakU_n=\frakD_n$. \redmissingref{I guess the alternative would be to set
Let us also set $\oldGamma_{m+1}=\frakK$. Then for $k=1,\ldots,m$,
there is a trivial $I$-fibration of 
$\frakU_k$ for which the
components of $\partialh\frakU_k$ are $\oldGamma_k$ and
$\oldGamma_{k+1}$. We have $\oldGamma_k\subset\oldTheta''$ for $1\le
k\le m$, and $\oldGamma_{m+1}\subset\partial\oldPsi'$.

The pair
$(\frakU_0,\oldGamma_1)$ is homeomorphic to $(\oldLambda,\frakK)$;
hence $\frakU_0$ admits a $I$-fibration such
that $\oldGamma_1$ is a component of $\partialh\frakU_0$. This
fibration is trivial if and only if the connected \pagelike\ \Ssuborbifold\
$\oldLambda$ of $\oldPsi'$ is untwisted.
%\redcomment{the original thing was
 % untwisted. Use this to replace refs. to triviality of 
  %throughout the proof.} 
Thus we
have $\partialh\frakU_0=\oldGamma_1$ if the connected \pagelike\ \Ssuborbifold\ $\oldLambda$ of $\oldPsi'$ is
twisted;
%the fibration of $\frakU_0$
%is non-trivial; 
and if it is untwisted, $\partialh\frakU_0$
consists of $\oldGamma_1$ and a second component $\oldGamma_0$. In the latter case, $\oldGamma_0$ is the component of $\oldLambda\cap\partial\oldPsi'$ distinct from $\frakK$, so that in particular $\oldGamma_0\subset\toldTheta'\subset\toldTheta$. In all cases we have $\frakU_{k-1}\cap\frakU_k=\oldGamma_k$ for $k=1,\ldots,m$. 

For $k=0,\ldots,m$ we have
$\frakU_k-\partialh\frakU_k\subset\oldPsi'-\oldTheta''$. Hence there
is an embedding
$\eta_k:\frakU_k\to\oldPsi=(\oldPsi')\cut{\oldTheta''}$ such that
$\rho\circ\eta_k$ is the inclusion map $\frakU_k\to\oldPsi'$. 

Set
$\toldGamma_k=\eta_{k}(\oldGamma_k))$ for $0\le k\le m$ if $\oldLambda$ is untwisted, and for  $1\le k\le m$ if $\oldLambda$ is twisted.
Set $\toldGamma'_k:=\eta_{k-1}(\oldGamma_k))$ for $1\le k\le m+1$.
For $k=1,\ldots,m$, the maps
$\rho\circ(\eta_{k-1}|\oldGamma_k)$ and
$\rho\circ(\eta_{k}|\oldGamma_k)$ are both equal to the inclusion map
$\oldGamma_k\to\oldPsi'$. Since
$\frakU_{k-1}\cap\frakU_k=\oldGamma_k$, the orbifolds 
$\toldGamma_k$ and
$\toldGamma'_k$
%$\eta_{k-1}$ and
%$\eta_k$ 
are contained in distinct sides of the component of
$\toldTheta''$ containing $\oldGamma_k$. 
%\redcomment{Small problem here:
  %If this is the def. of $\toldGamma'_k$, it has been given only for
  %$1\le k\le m$. But I need $\toldGamma'_{m+1}$ later. Tweak
  %required.} 
Hence
$\eta_{k-1}|\oldGamma_k=\tau\circ(\eta_{k}|\oldGamma_k)$, so that
%. It follows
%that if we set  for $i=1,\ldots,m$, we have
\Equation\label{why not another}
\toldGamma'_k=\tau(\toldGamma_k) \text{ for } k=1,\ldots,m.
\EndEquation

Note also that since $\oldGamma_{m+1}=\frakK$, we may use our
identification of $\toldTheta'$ with a union of components of
$\toldTheta$ to write $\toldGamma_{m+1}'=\frakK$. Likewise, in the case
where $\oldLambda$ is untwisted, $\toldGamma_0$ is identified with
$\oldGamma_0$, so that in particular
$[\toldGamma_0]\preceq[\toldTheta']$. 
%\redcomment{The following sentence doesn't seem to make
  %much sense. Should it be $\toldGamma_0 $  and $\toldGamma_{m+1}$
  %instead of $\oldGamma_0 $  and $\oldGamma_{m+1}$, and suborbifolds
  %instead of components? Is it quoted anywhere?} In particular, $\oldGamma_0 $  and $\oldGamma_{m+1}=\frakK$ may be regarded as components of $\toldTheta' $. \redmissingref{ I was confused about this passage. Check again that it makes sense. An earlier comment said to recheck the redcomment below.}

For $k=0,\ldots,m$, set $\tfrakU_k=\eta_k(\frakU_k)$. Then $\tfrakU_k$
inherits an $I$-fibration from $\frakU_k$ via the homeomorphism
$\eta_k:\frakU_k\to\tfrakU_k$. For $k=1,\ldots,m$, the $I$-fibration
of $\tfrakU_k$ is trivial, and the components of $\partialh\tfrakU_k$
are $\toldGamma_k$ and $\toldGamma'_{k+1}$. 
We
have $\partialh\tfrakU_0=\toldGamma_1'$ if the connected \pagelike\
\Ssuborbifold\ $\oldLambda$ of $\oldPsi'$ is twisted; and if it is untwisted, $\partialh\tfrakU_0$ has two components, $\toldGamma'_1$ and $\toldGamma_0$. We claim:
\Claim\label{sycamore slew}
For $k=0,\ldots,m$, the orbifold $\tfrakU_k$ is an \Ssuborbifold\ of $\oldPsi$. (In particular, each component of $\Fr_\oldPsi\tfrakU_k$ is essential in $\oldPsi$.)
\EndClaim
To prove \ref{sycamore slew}, first note that we have
$\partialh\tfrakU_k=\toldGamma_k\cup\toldGamma'_{k+1}=\tfrakU_k\cap\partial\oldPsi$
if $0<k\le m$, or if $k=0$ and $\oldLambda$ is untwisted. Furthermore, if $\oldLambda$ is twisted, we have
$\partialh\tfrakU_0=\toldGamma'_{1}=\tfrakU_0\cap\partial\oldPsi$. Hence
% in the
%case where
 the pair ($\tfrakU_k,\tfrakU_k\cap\partial\oldPsi')$ is
always an \spair. According to Definition \ref{S-pair def}, in order to show that $\tfrakU_k$ is an \Ssuborbifold, it remains only to show that $\Fr_\oldPsi\tfrakU_k$ is essential in $\oldPsi$. 

Let $\tfrakV$ be any component of $\Fr_\oldPsi\tfrakU_k$. Then $\rho$ maps $\tfrakV$ homeomorphically onto a suborbifold $\frakV$ of some component $\frakV^+$ of $\Fr_{\oldPsi'}\oldLambda$.
Since the discussion above shows that each component of $\partialh\frakU_k$ is either equal to $\frakK$ or
parallel  to $\frakK$ in the pair $(\oldLambda,\Fr\oldLambda)$,
%parallel \redmissingref{rel something?} 
each boundary component of $|\frakV|$ either cobounds a weight-$0$ annulus with, or is equal to, some common boundary component of %there exist a component \supset\frakV$ of $\Fr_{\oldPsi'}\oldLambda$, and a simple closed curve component $C$ of 
$|\frakV^+|$ and $|\frakK|$. 
It follows that the annular orbifold $\frakV$ is $\pi_1$-injective in the annular orbifold
$\frakV^+$.
But the component $\frakV^+$ of $\frakA$ is essential in $\oldPsi'$ by \ref{tuesa day}; in particular, $\frakV^+$ is $\pi_1$-injective in $\oldPsi'$. It follows that
$\frakV$ is $\pi_1$-injective in $\oldPsi'$, and hence that $\tfrakV$ is $\pi_1$-injective in $\oldPsi$.  %follows that 
 %$ such that each component of $|\partial\frakV|$ is the frontier in $|\frakV^+|$ of a  weight-$0$ annulus in $|\frakV^+|$ whose intersection with $\partial\frakV^+$ is $C$.
%In particular $\frakV$ is $\pi_1$-injective in $\frakV^+$. Since the component $\frakV^+$ of $\frakA(\oldPsi')$ is essential, and hence $\pi_1$-injective, in $\oldPsi'$, it follows that $\frakV$ is $\pi_1$-injective in $\oldPsi'$, and hence that $\tfrakV$ is $\pi_1$-injective in $\oldPsi$. \redcomment{I may need to rewrite this to provide notation for the second part of the proof of essentiality, but I'm no longer sure of this. Make sure I've dealt with the possibility that $|\frakV|$ is a weight-$2$ disk.}

Now suppose that $\tfrakV$ is parallel in $\oldPsi$ to a $2$-suborbifold of $\partial\oldPsi$. Then $\frakV$ is parallel in $\oldPsi'$ to a $2$-suborbifold $\frakZ$ of $\toldTheta'\cup\oldTheta''\subset\oldPsi'$. Since $\frakZ$ is homeomorphic to $\frakV$, it is annular, and in particular connected; hence either $\frakZ\subset\toldTheta'$ or $\frakZ\subset\oldTheta''$.
If $\frakZ\subset\oldTheta''$, then conditions (i)---(iii) of Definition \ref{reduced intersection} hold with $\oldPsi'$, $\oldTheta''$, $\Fr_{\oldPsi'}\oldLambda$, and $\frakZ$ playing the respective roles of
$\oldPsi$, $\frakP$, $\frakB$, and $\frakC$ in that definition. This contradicts the fact that $\Fr_{\oldPsi'}\oldLambda$ has reduced intersection (see Definition \ref{reduced intersection}) with $\oldTheta''\subset\inter\oldPsi'$. Now suppose that
$\frakZ\subset\toldTheta'$. Then we have $\partial\frakV\subset\frakV^+\cap\toldTheta'=\partial\frakV^+$, and hence $\frakV=\frakV^+$. (In the notation above, we have $m=0$.) Hence $\frakV^+$ is parallel in $\oldPsi'$ to a $2$-suborbifold of $\partial\oldPsi'$. This contradicts the essentiality of 
$\frakV^+$ in $\oldPsi'$. Thus we have shown that
$\tfrakV$ is not parallel in $\oldPsi$ to a $2$-suborbifold of $\partial\oldPsi$. This completes the proof that $\tfrakV$ is essential in $\oldPsi$, and \ref{sycamore slew} is thereby proved.

It follows from \ref{sycamore slew} and
Proposition \ref{slide it} that $\tfrakU_k$ is isotopic in $\oldPsi$ to
a suborbifold $\tfrakU_k^0$ of $\oldPhi(\oldPsi)$. Note also that since the
components of $\partialh\tfrakU_k$ are parallel to $\frakK$ in
$\oldLambda$, they are negative; since $\Fr_\oldPsi\tfrakU_k$ is
$\pi_1$-injective by \ref{sycamore slew}, the component of $\oldPhi(\oldPsi)$ containing
$\tfrakU_k^0\cap\partial\oldPsi$ must be negative,
i.e. $\tfrakU_k^0\subset\oldSigma^-(\oldPsi)$. 
The discussion in \ref{what's iota?} shows that $\oldSigma^-(\oldPsi)$ is a \pagelike\ \Ssuborbifold pf $\oldPsi)$, and as in \ref{what's iota?} we fix
an $I$-fibration $q:\oldSigma^-\to \frakB$ over a  $2$-orbifold in
such a way that $\partialh\oldSigma^-=\oldPhi^-$.

%Now by  Proposition \ref{new characteristic} and Definition \ref{oldSigma def}, each component of $\oldSigma(\oldPsi)$ is an \Ssuborbifold\ of $\oldPsi$, and by Lemma \ref{when a tore a fold}, no component of $\oldSigma(\oldPsi)$ admits an $\SSS^1$-fibration. 
%\redmissingref{a
  %cross-ref., and I really think there should be a universal
  %convention about this. After all, the very definition of
  %$\iota_\oldPsi$ involves an $I$-fibration of $\oldSigma^-$ which is
  %canonical up to some sort of isotopy. Also, $\oldSigma^-(\oldPsi)\cap\partial\oldPsi$ is usually called $\oldPhi^-(\oldPsi)$.} 
%Hence there is an $I$-fibration of
%$\oldSigma^-(\oldPsi)$ for which
%$\partialh\oldSigma^-(\oldPsi)=\oldSigma^-(\oldPsi)\cap\partial\oldPsi$; by \ref{tuesa day}, the latter set is equal to $\oldPhi^-(\oldPsi)$. 
Since
the (annular) components of $\Fr_\oldPsi\tfrakU_k^0$ are essential in $\oldPsi$ by \ref{sycamore slew}, they are in particular $\pi_1$-injective in $\oldSigma^-(\oldPsi)$, and none of them is parallel in the pair $(\oldSigma^-(\oldPsi),\oldPhi^-(\oldPsi))$  to a suborbifold of $\oldPhi^-(\oldPsi)$; hence Proposition \ref{when vertical}, applied with $\oldSigma^-(\oldPsi)$ and $\Fr_\oldPsi\tfrakU_k^0$ playing the roles of $\oldLambda$ and $\frakC$ respectively, allows us to choose
%containing $\tfrakU_k^0$ is a
%\pagelike\ \Sssuborbifold by \redmissingref{a cross-ref. that may
  %disappear when the comment below is dealt with}; 
$\tfrakU_k^0$ within its isotopy class in $\oldSigma^-(\oldPsi)$ so
that $\Fr\tfrakU_k^0$, and hence $\tfrakU_k^0$, is saturated in the $I$-fibration of
$\oldSigma^-(\oldPsi)$. It now follows from the discussion in \ref{what's iota?} that
$\tfrakU_k^0\cap\partial\oldPsi$, which has either one or two
components, is invariant under $\iota:=\iota_\oldPsi$, and that if it
has two components then they are interchanged by $\iota$. This proves:
\Claim\label{iota be able to do this}
We have   $[\toldGamma_{k}]\preceq [\oldPhi^-(\oldPsi)]$ for
$k=1,\ldots,m$ and
$[\toldGamma'_{k}]\preceq [\oldPhi^-(\oldPsi)]$ for $k=1,\ldots,m+1$. If
the connected \pagelike\ \Ssuborbifold\ $\oldLambda$ of $\oldPsi'$ is untwisted,
we also have $[\toldGamma_{0}]\preceq
[\oldPhi^-(\oldPsi)]$. Furthermore, we have $[\iota]([\toldGamma_{k+1}'])=[\toldGamma_{k}]$ for
$k=1,\ldots,m$, and for $k=0$ if $\oldLambda$ is untwisted. In the
case where $\oldLambda$ is twisted, we have
$[\iota]([\toldGamma_{1}'])=[\toldGamma_{1}']$.
\EndClaim

Now consider any index $k$ with $1\le k\le
m$.
%, that $[\toldGamma_k']\preceq V_k$, and that
%$\noodge_k([\toldGamma_k'])\preceq[\toldGamma_0]$. 
By \ref{iota be able
  to do this} we have
$[\toldGamma'_{k+1}]\preceq [\oldPhi^-(\oldPsi)]=\dom[\iota]=\dom([\tau]\diamond[\iota])$
and $[\iota]([\toldGamma_{k+1}'])=[\toldGamma_{k}]$. By
(\ref{why not another}) we have
$\toldGamma'_{k}=\tau(\toldGamma_k)$ and hence:
$([\tau\circ\iota])([\toldGamma_{k+1}'])=[\toldGamma_{k}']$. Using \ref{return of zoltan} to rewrite $[\tau\circ\iota]$ as $[\tau]\diamond[\iota]$, we deduce:

\Claim\label{a poor thing but mine own}
For $k=1,\ldots,m$ we have
$[\toldGamma'_{k+1}]\preceq\dom([\tau]\diamond[\iota])$
and 
$([\tau]\diamond[\iota])([\toldGamma_{k+1}'])=[\toldGamma_{k}']$.
\EndClaim

Next, for each $n\ge1$, set
$\noodge_n=\noodge_n^{\oldOmega,\oldTheta}$ and $V_n=
V_n^{\oldOmega,\oldTheta}=\dom\noodge_n$. We claim:
\Claim\label{now here's the thing}
If $\oldLambda$ is untwisted, then for $k=1,\ldots,m+1$ we have
$[\toldGamma_k']\preceq V_k$ and
$\noodge_k([\toldGamma_k'])\preceq[\toldGamma_0]$. 
If $\oldLambda$ is twisted, then for $k=1,\ldots,m+1$ we have
$[\toldGamma_k']\preceq V_{2k-1}$ and
$\noodge_{2k-1}([\toldGamma_k'])\preceq[\toldGamma_k']$. 
\EndClaim

To prove the first assertion of \ref{now here's the thing}, we argue
by induction on $k$. For the proof when $k=1$, observe that since
$\oldLambda$ is untwisted, \ref{iota be able to do this} and (\ref{and
  so did you}) give
$[\toldGamma_{1}']\preceq
[\oldPhi^-(\oldPsi)]=V_1$ and 
$\noodge_1
([\toldGamma_{1}'])=[\iota]([\toldGamma_{1}'])=[\toldGamma_{0}]$,
which implies the assertion in this case. Now suppose that $1\le k\le
m$, that $[\toldGamma_k']\preceq V_k$, and that
$\noodge_k([\toldGamma_k'])\preceq[\toldGamma_0]$. By \ref{a poor thing but mine own} we have
$[\toldGamma'_{k+1}]\preceq\dom([\tau]\diamond[\iota])$
and 
$([\tau]\diamond[\iota])([\toldGamma_{k+1}'])=[\toldGamma_{k}']$. Since we
also have $[\toldGamma_k']\preceq V_k=\dom\noodge_k$, 
Condition (2) of Lemma \ref{before associativity} holds with
$Y=[\toldGamma'_{k+1}]$, $\emm_1=[\tau]\diamond[\iota]$, and
$\emm_2=\noodge_k$. Furthermore, the definitions of the $\noodge_n$ (see \ref{in duck tape}), with
\ref{more associativity}, imply that
%$\noodge_k\diamond[\tau\diamond\iota]=\noodge_{k+1}$, which by \ref{return of zoltan} may be rewritten as 
$\noodge_k\diamond([\tau]\diamond[\iota])=\noodge_{k+1}$. Lemma \ref{before associativity} therefore implies
that
%$[\toldGamma'_{k+1}]\preceq\dom(\noodge_k\diamond([\tau]\diamond[\iota]))$,
%and that 
%$(\noodge_k\diamond([\tau]\diamond[\iota]))|[\toldGamma'_{k+1}]=
%\noodge_k\circ([\tau\circ\iota]|[\toldGamma'_{k+1}])$.
%But
% Thus we have
$[\toldGamma'_{k+1}]\preceq\dom\noodge_{k+1}=V_{k+1}$
and that
$\noodge_{k+1}|[\toldGamma'_{k+1}]=
\noodge_k\circ(([\tau]\diamond[\iota])\big|[\toldGamma'_{k+1}])$. 
Hence
$\noodge_{k+1}([\toldGamma'_{k+1}])
=(\noodge_k\circ([\tau]\diamond[\iota]))([\toldGamma'_{k+1}]))$.
But by \ref{special composition}, 
$(\noodge_k\circ([\tau]\diamond[\iota]))([\toldGamma'_{k+1}]))=\noodge_k(([\tau]\diamond[\iota])([\toldGamma'_{k+1}]))$.
We therefore have
$$\noodge_{k+1}([\toldGamma'_{k+1}])
=\noodge_k(([\tau]\diamond[\iota])([\toldGamma'_{k+1}]))
=\noodge_k([\toldGamma'_{k}])\preceq[\toldGamma_0].$$
This completes the induction. 

To prove the second assertion of \ref{now here's the thing}, we again argue
by induction on $k$. For the proof when $k=1$, observe that since
$\oldLambda$ is twisted, \ref{iota be able to do this} and (\ref{and
  so did you}) give
$[\toldGamma_{1}']\preceq
[\oldPhi^-(\oldPsi)]=V_1$ and 
$\noodge_1
([\toldGamma_{1}'])=[\iota]([\toldGamma_{1}'])=[\toldGamma_{1}']$,
which implies the assertion in this case. Now suppose that $1\le k\le
m$, that $[\toldGamma_k']\preceq V_{2k-1}$, and that
$\noodge_{2k-1}([\toldGamma_k'])\preceq[\toldGamma_k']$. 
By \ref{a poor thing but mine own} we have
$[\toldGamma'_{k+1}]\preceq\dom([\tau]\diamond[\iota])$
and 
$([\tau]\diamond[\iota])([\toldGamma_{k+1}'])=[\toldGamma_{k}']$.
Since we
also have $[\toldGamma_k']\preceq V_{2k-1}=\dom\noodge_{2k-1}$, 
Condition (2) of Lemma \ref{before associativity} holds with
$Y=[\toldGamma'_{k+1}]$, $\emm_1=[\tau]\diamond[\iota]$, and
$\emm_2=\noodge_{2k-1}$. Furthermore, the definitions of the $\noodge_n$, with
\ref{more associativity}, imply that
$\noodge_{2k-1}\diamond[\tau\diamond\iota]=\noodge_{2k}$.
%, which by \ref{return of zoltan} may be rewritten as
%$\noodge_{2k-1}\diamond([\tau]\diamond[\iota])=\noodge_{2k}$.
Lemma \ref{before associativity} therefore implies
that
%$[\toldGamma'_{k+1}]\preceq\dom(\noodge_k\diamond([\tau]\diamond[\iota]))$,
%and that 
%$(\noodge_k\diamond([\tau]\diamond[\iota]))|[\toldGamma'_{k+1}]=
%\noodge_k\circ(([\tau]\diamond[\iota])|[\toldGamma'_{k+1}])$.
%But
% Thus we have
$[\toldGamma'_{k+1}]\preceq\dom\noodge_{2k}$
and 
$\noodge_{2k}|[\toldGamma'_{k+1}]=
\noodge_{2k-1}\circ(([\tau]\diamond[\iota])|[\toldGamma'_{k+1}])$.
Hence
$\noodge_{2k}([\toldGamma'_{k+1}])=
\noodge_{2k-1}(([\tau]\diamond[\iota])([\toldGamma'_{k+1}]))
=\noodge_{2k-1}([\toldGamma'_{k}])\preceq [\toldGamma'_{k}]$. 

Since $[\iota]^{-1}=[\iota]$ and $[\tau]^{-1}=[\tau]$ (see \ref{in duck tape}),
it follows from Proposition \ref{associativity} that
 $[\iota]\diamond[\tau]=([\tau]\diamond[\iota])^{-1}$. In view of
\ref{my remark}, and the relations 
$[\toldGamma'_{k+1}]\preceq\dom([\tau]\diamond[\iota])$
and 
$([\tau]\diamond[\iota])([\toldGamma_{k+1}'])=[\toldGamma_{k}']$, it follows
 that
$[\toldGamma'_{k}]\preceq\dom([\iota]\diamond[\tau])$
and that
$([\iota]\diamond[\tau])([\toldGamma_{k}'])=[\toldGamma_{k+1}']$.  Condition (2) of Lemma \ref{before associativity} therefore holds with $Y=\toldGamma_{k+1}'$, $\emm_1=\noodge_{2k}$ and $\emm_2=[\iota]\diamond[\tau]$. Furthermore, the definitions of the $\noodge_n$, with
\ref{more associativity}, also imply that
$([\iota]\diamond[\tau])\diamond\noodge_{2k}=\noodge_{2k+1}$. Lemma \ref{before associativity} therefore implies
that
$[\toldGamma'_{k+1}]\preceq\dom\noodge_{2k+1}$
and that
$\noodge_{2k+1}|[\toldGamma'_{k+1}]=([\iota]\diamond[\tau])\circ
(\noodge_{2k}|[\toldGamma'_{k+1}])$.
In view of \ref{special composition} it follows that
$\noodge_{2k+1}([\toldGamma'_{k+1}])=
([\iota]\diamond[\tau])(
\noodge_{2k}([\toldGamma'_{k+1}]))$. Since 
$\noodge_{2k}([\toldGamma'_{k+1}])\preceq [\toldGamma'_{k}]\preceq\dom([\iota]\diamond[\tau])$, the order-preserving property stated in \ref{restriction} gives $([\iota]\diamond[\tau])(
\noodge_{2k}([\toldGamma'_{k+1}]))\preceq
([\iota]\diamond[\tau])(
[\toldGamma'_{k}])=[\toldGamma'_{k+1}]$. 
It follows that
$\noodge_{2k+1}([\toldGamma'_{k+1}])\preceq 
[\toldGamma'_{k+1}])$,
and again the induction is complete. Thus \ref{now here's the thing} is proved.

%$Y\preceq\dom(\emm_2\diamond\emm_1)$.
% \item $Y\preceq\dom\emm_1$, and $\emm_1(Y)\preceq\dom\emm_2$ (where $\emm(Y)$ is 

%$[\toldGamma'_{k}]\preceq [\oldPhi^-(\oldPsi)]$ for $k=1,\ldots,m$. If
%the \pagelike\ \Ssuborbifold\ $\oldLambda$ of $\oldPsi'$ is untwisted,
%we also have $[\toldGamma_{0}]\preceq
%[\oldPhi^-(\oldPsi)]$. Furthermore, we have $[\iota]([\toldGamma_{k+1}'])=[\toldGamma_{k}]$

 %\redproofsummary{See how much of the notation I got right
  %there. It seems better to use an equality rather than a $\preceq$ in
  %the second assertion of \ref{now here's the thing}. The proof of each of the two assertions of \ref{now here's the thing} is an
  %induction using (\ref{why not another}) and \ref{iota be able to do this}.
%The conclusion should then follow easily. 

To prove that 
%$\oldPhi^-(\oldPsi')$ then
 $\frakK \preceq
X'(\oldOmega,\oldTheta,\oldTheta')$, which has been seen to imply the conclusion of the proposition, first consider the case in
which $\oldLambda$ is untwisted. In this case, applying the first
assertion of \ref{now here's the thing} with $k=m+1$, and recalling that
$\toldGamma_{m+1}'=\frakK$, we find that $[\frakK]\preceq
V_{m+1}$ and that
$\noodge_{m+1}([\frakK])\preceq[\toldGamma_0]\preceq[\toldTheta']$. 
%\redcomment{I'm
  %pretty uncertain about the notation in that sentence. The issue is
  %related to the comment above that begins ``The following 
%sentence doesn't seem to make
  %much sense.''}
%\redcomment{The
 % fact that $[\toldGamma_0]\preceq[\toldTheta']$ should have been
  %mentioned earlier. Also the fact that
  %$[\frakK]\preceq[\toldTheta']$. But these depend on either the same implicit identifications mentioned in an earlier comment, or similar ones.} 
If $U$ denotes the component of
  $V_{m+1}$ such that $[\frakK]\preceq U$ (see \ref{components}), it follows that $\noodge_{m+1}(U)\preceq[{\toldTheta}']$. In view of the definitions
  given in \ref{oldXi}, this means that
  $U\in\oldfrakX_{m+1}(\oldOmega,\oldTheta,\oldTheta')$, so that $U\preceq
  X_{m+1}(\oldOmega,\oldTheta,\oldTheta')$ and hence $[\frakK]\preceq
  X_{m+1}(\oldOmega,\oldTheta,\oldTheta')\preceq
  X(\oldOmega,\oldTheta,\oldTheta')$. Since $[\frakK]\preceq\toldTheta'$, it follows that 
$[\frakK]\preceq
  X'(\oldOmega,\oldTheta,\oldTheta')$, as required.  

Now consider the case in
which $\oldLambda$ is twisted. In this case, applying the second
assertion of \ref{now here's the thing} with $k=m+1$,  and again using that
$\toldGamma_{m+1}'=\frakK$,  we find that $[\frakK]\preceq
V_{2m+1}$ and that
$\noodge_{2m+1}([\frakK])\preceq[\frakK]\preceq\toldTheta'$. If in
this case $U$ denotes the component of
  $V_{2m+1}$ such that $\frakK\preceq U$, it follows that $\noodge_{2m+1}(U)\preceq[{\toldTheta}']$. This means that
  $U\in\oldfrakX_{2m+1}(\oldOmega,\oldTheta,\oldTheta')$, so that $U\preceq
  X_{2m+1}(\oldOmega,\oldTheta,\oldTheta')$ and hence $\frakK\preceq
  X_{2m+1}(\oldOmega,\oldTheta,\oldTheta')\preceq
  X(\oldOmega,\oldTheta,\oldTheta')$. Since $\frakK\preceq\toldTheta'$, it follows that 
$\frakK\preceq
  X'(\oldOmega,\oldTheta,\oldTheta')$, as required.  
\EndProof
%,n \le n+1_ndefinitions\oldLambda\oldTheta
%\Theta\oldSigma(a)\frakE\frakB\frakC\frakD\oldPi\frakP\frakU\tfrakU\\tfrakR \frakA_\oldPi\frakP\frakVC\frakZ
%[\iota\circ\tau]{newer

\Corollary\label{burpollary}
Let $\Mh$ be a closed, hyperbolic, orientable $3$-orbifold, and let $\oldTheta$ be an incompressible, closed $2$-suborbifold
of $\oldOmega:=(\Mh)\pl$. 
Set
$\oldPsi=\oldOmega\cut\oldTheta$ and ${\toldTheta}=\partial\oldPsi$, and
$\tau=\tau_\oldTheta$.
Let $\oldTheta'$ be a union of components of $\oldTheta$, and set
$\oldPsi'=\oldOmega\cut{\oldTheta'}$ and ${\toldTheta}'=\partial\oldPsi'$, so that
${\toldTheta}'$ is canonically identified with a union of components of
${\toldTheta}$. Then 
$[\oldPhi^-(\oldPsi')]\preceq
[\oldPhi^-(\oldPsi)]$.
\EndCorollary

\Proof
According to \ref{oldXi}, we have
$X(\oldOmega,\oldTheta,\oldTheta')\preceq [\oldPhi^-(\oldPsi)]$. Since Proposition \ref{herman had burped} gives
$[\oldPhi^-(\oldPsi')]\preceq
X'(\oldOmega,\oldTheta,\oldTheta')\preceq
X(\oldOmega,\oldTheta,\oldTheta')$, the corollary
follows. 
\EndProof

\Lemma\label{before veggie}
Let $\Mh$ be a closed, hyperbolic, orientable $3$-orbifold, and let $\oldTheta$ be an incompressible, closed $2$-suborbifold
of $\oldOmega:=(\Mh)\pl$. Set
$\oldPsi=\oldOmega\cut\oldTheta$ and ${\toldTheta}=\partial\oldPsi$, 
$\tau=\tau_\oldTheta$, and $\noodge_n=\noodge_n^{\oldOmega,\oldTheta}$ for every $n\ge1$.
Let $\oldTheta'$ be a union of components of $\oldTheta$. Set
$\oldPsi'=\oldOmega\cut{\oldTheta'}$ and ${\toldTheta}'=\partial\oldPsi'$, so that
${\toldTheta}'$ is canonically identified with a union of components of ${\toldTheta}$. Let $m>0$ be an integer, and let $W$ be a component of $V_m^{\oldOmega,\oldTheta}$ such that $\noodge_k(W)\preceq[\oldTheta-\oldTheta']$ for $k=1,\ldots,m$. Then 
$$
W\wedge X(\oldOmega,\oldTheta,\oldTheta')=\noodge_m(\noodge_m(W)\wedge
[\tau](X(\oldOmega,\oldTheta,\oldTheta')).
$$
\EndLemma

Note that the right hand side of the displayed formula in Lemma
\ref{before veggie} makes sense, because
$V_m=\dom\noodge_m=\range\noodge_m$ by \ref{in duck tape}, and
therefore $\noodge_m(W)\wedge
[\tau](X(\oldOmega,\oldTheta,\oldTheta'))\preceq \noodge_m(W)\preceq
V_m$.

\Proof[Proof of Lemma \ref{before veggie}]
Set $X=X(\oldOmega,\oldTheta,\oldTheta')$, and set $V_n=V_n^{\oldOmega,\oldTheta}$ and
$X_n=X_n(\oldOmega,\oldTheta,\oldTheta')$ for every $n\ge1$. 
 We claim:
\Claim\label{main point}
For each $k$
with $1\le k\le m$ we have $W\wedge X_k=[\emptyset]$, and
for each $k>m$ we have
$W\wedge X_k=\noodge_m(\noodge_m(W)\wedge
[\tau](X_{k-m}))$.
\EndClaim

To prove \ref{main point}, consider any index $k\ge1$.
Suppose that $Z$ is a component of
$W\wedge X_k$. In particular we have $Z\preceq
X_k$, and since $Z$ is connected, an
observation made in \ref{components} implies that $Z\preceq
U$ for some component $U$ of $X_k$. According to an observation made in  \ref{oldXi}, each component of $X_k$ is an element of $\oldfrakX_k(\oldOmega,\oldTheta,\oldTheta')$. Hence we have
$U\in\oldfrakX_k(\oldOmega,\oldTheta,\oldTheta')$; that is, $U$ is a component of $V_k$ and 
$\noodge_k(U)\preceq[{\toldTheta}']$. 
Since 
$\noodge_k(Z)\preceq\noodge_k(U)$ (by the order-preserving property
observed in \ref{restriction}), we have $\noodge_k(Z)\preceq[{\toldTheta}']$.

Consider the case in which $ k\le m$.
%Since 
%$\noodge_k(Z)\preceq\noodge_k(U)$ (by the order-preserving property
%observed in \ref{restriction}), \redmissingref{awkward phrase} we have
%$\noodge_k(Z)\preceq[{\toldTheta}']$. 
%On the other hand, s
Since $Z\preceq W$,
we have $\noodge_k(Z)\preceq\noodge_k(W)$, which with the hypothesis
gives $\noodge_k(Z)\preceq[{\toldTheta}-{\toldTheta}']$, which
contradicts $\noodge_k(Z)\preceq[{\toldTheta}']$. This
shows that $W\wedge X_k=[\emptyset]$ when $1\le k\le m$, which is the first assertion of \ref{main point}.

Now consider the case in which $k>m$. 
%Let $Z$ be any component of $W\wedge
%X_k$. In particular we have $Z\preceq
%X_k$, and hence $Z\preceq
%U$ for some component $U$ of $X_k$. It again follows from \ref{oldXi}
%that
%$U\in\oldfrakX_n(\oldOmega,\oldTheta,\oldTheta')$; that is, $U$ is a component of $V_k$ and 
%$\noodge_k(U)\preceq[{\toldTheta}']$.
%that $U$ is an element of $\oldfrakX_n component of $V_k=\dom\noodge_k$ \redmissingref{I'm adding a
  %comment to \ref{oldXi} saying that this needs to be made explicit
  %there}  and that
%$\noodge_k(U)\preceq[{\toldTheta}']$. 
%In particular we have 
Since $Z\preceq U\preceq V_k$, it follows from
% and
%$\noodge_k(Z)\preceq[{\toldTheta}']$. According to
 Lemma \ref{j-lemma}
%, the relation $Z\preceq V_k$ implies 
that $Z\preceq V_m$, that $[\tau](\noodge_m(Z))\preceq V_{k-m}$, and that
$\noodge_k|Z=\noodge_{k-m}\circ [\tau]
\circ(\noodge_m|Z)$. With \ref{special composition} this implies that 
$\noodge_k(Z)=\noodge_{k-m} ([\tau](\noodge_m|Z))$.
Hence the relation 
$\noodge_k(Z)\preceq[{\toldTheta}']$ may be rewritten as 
%\emm_2\diamond
$\noodge_{k-m}([\tau](\noodge_m(Z))\preceq[{\toldTheta}']$. If we set
$T=[\tau](\noodge_m(Z))\preceq 
%[\tau](\noodge_m(V_k))\preceq
 V_{k-m}$, we therefore have
$\noodge_{k-m}(T)\preceq[{\toldTheta}']$. If $R$ denotes the component
of $V_{k-m}$ such that $T\preceq R$ (see \ref{components}), then since ${\toldTheta}'$ is a
union of components of ${\toldTheta}$ we must have
$\noodge_{k-m}(R)\preceq[{\toldTheta}']$. By definition this means that
$R$ is an element of $\oldfrakX_{k-m}$, so that $R\preceq X_{k-m}$; hence $T\preceq X_{k-m}$,
i.e. $[\tau](\noodge_m(Z))\preceq X_{k-m}$. The order-preserving
property pointed out in \ref{restriction} then gives 
%Since $[\tau]^{-1}=[\tau]$
%and $\noodge_m^{-1}=\noodge_m$ by \redmissingref{give cross-refs. for
  %these, and I think I need another fact which may or may not have
  %been written down}, it follows that 
$[\tau]([\tau](\noodge_m(Z
% X_{k-m}
)))
\preceq  [\tau]( X_{k-m})$, i.e. $\noodge_m(Z)
%=[\tau]([\tau]( X_{k-m}))
\preceq  [\tau]( X_{k-m})$.

%\redmissingref{I'd said $Z\preceq \noodge_m( [\tau](
%X_{k-m}))$, which seems irrelevant.} 
%\redcomment{Think about what principle I am using here. Has
  %it ever been stated? Is it the same principle mentioned in another
  %comment? Maybe not, because this one involves inverses.}

%Since this holds for every component $Z$ of $W\wedge X_k$, it follows from \redcomment{a lemma that I need to write} that $W\wedge X_k\preceq \noodge_m( [\tau]( X_{k-m}))$.
% Hence $\noodge_m(Z)\preceq  [\tau]( X_{k-m})$. 
Since we also have $\noodge_m(Z)\preceq\noodge_m( W\wedge X_k)\preceq  \noodge_m( W)$,
%$\noodge_m( W\wedge X_k)\preceq  
%$\noodge_m( W)$, 
it follows that
$\noodge_m( Z)\preceq \noodge_m( W)\wedge  [\tau]( X_{k-m})$.
The order-preserving
property pointed out in \ref{restriction} then gives  
%We therefore have
$ Z=\noodge_m( \noodge_m(Z))
\preceq \noodge_m( \noodge_m( W)\wedge  [\tau](
X_{k-m}))$. 
Since this holds for every component $Z$ of $W\wedge X_k$, it now follows from Proposition \ref{mystery}  that 
\Equation\label{sagetimes}
W\wedge X_k\preceq
\noodge_m( \noodge_m( W)\wedge  [\tau]( X_{k-m})).
\EndEquation
%$\noodge_m( W\wedge X_k)\preceq \noodge_m( W)  [\tau]( X_{k-m})$.
%$\noodge_m( W\wedge X_k)\preceq \noodge_m( W)  [\tau]( X_{k-m})$.

In view of (\ref{sagetimes}) and Proposition \ref{new partial order}, to complete the proof of the second assertion of \ref{main point}, it remains only to show that $\noodge_m(
  \noodge_m( W)\wedge [\tau]( X_{k-m}))\preceq W\wedge X_k$. To this end, suppose that $H$ is a component of
$\noodge_m( \noodge_m( W)\wedge [\tau]( X_{k-m}))$. Then
$H\preceq \noodge_m( \noodge_m( W)\wedge [\tau]( X_{k-m}))\preceq\range\noodge_m=V_m$, and
the order-preserving
property pointed out in \ref{restriction} gives 
$\noodge_m(H)\preceq\noodge_m(\noodge_m(
\noodge_m( W)\wedge [\tau](
X_{k-m})
%\preceq[\tau](X_{k-m})
))
=
\noodge_m( W)\wedge [\tau](
X_{k-m})\preceq[\tau](X_{k-m})$.
%\redcomment{I had written ``$\noodge_m(H)$ is a component of
%  $\noodge_m( W)\wedge [\tau](
  %X_{k-m})\preceq[\tau](X_{k-m})$,'' but I think it's irrelevant. I
  %think I should be quoting \ref{before associativity} to justify the
  %statements just made.}
  %\redmissingref{I think there may more or less be cross-refs for these
    %things.} 
Hence $[\tau](\noodge_m(H))\preceq X_{k-m}$. If $Q$
  denotes the component of $X_{k-m}$ such that
  $[\tau](\noodge_m(H))\preceq Q$ (see \ref{components}), then it follows from \ref{oldXi} that
 $Q\in\oldfrakX_{k-m}(\oldOmega,\oldTheta,\oldTheta')$, so that $Q$ is a component of $V_{k-m}$ and 
  $\noodge_{k-m}(Q)\preceq{\toldTheta}'$. Since $H\preceq V_m$ and
  $[\tau](\noodge_m(H))\preceq Q\preceq V_{k-m}$, it follows from Lemma
  \ref{j-lemma} that $H\preceq V_k$,
% implies that $H\preceq V_m$, that
  %$[\tau](\noodge_m(H))\preceq V_{k-m}$, 
and that
  $\noodge_k|H=\noodge_{k-m}\circ [\tau] \circ(\noodge_m|H)$.  Hence, by \ref{special composition}, we have
  $\noodge_k(H)=\noodge_{k-m}(
[\tau](\noodge_m(H)))\preceq\noodge_{k-m}(Q)\preceq{\toldTheta}'$. 
  Since ${\toldTheta}'$ is a union of components of ${\toldTheta}$,
  the component $P$ of $V_m$ such that $H\preceq P$ must satisfy
  $\noodge_k(P)\preceq{\toldTheta}'$. The definition of $X_k$ now
  implies that
$P\in\oldfrakX_{k}(\oldOmega,\oldTheta,\oldTheta')$, so that
 $P\preceq X_k$ and hence $H\preceq X_k$.  But we have
  $H\preceq\noodge_m( \noodge_m( W)\wedge [\tau]( X_{k-m}))
  \preceq\noodge_m( \noodge_m( W))=W$, where we have again used the
  order-preserving property. 
%\redcomment{This needs
    %justification, as do a lot of other similar implications.}
  Since $H\preceq X_k$ and $H\preceq W$, we have
  $H\preceq W\wedge X_k$. As this holds for every component $H$ of
  $\noodge_m( \noodge_m( W)\wedge [\tau]( X_{k-m}))$, it follows from
Proposition \ref{mystery} that 
$$ \noodge_m(
  \noodge_m( W)\wedge [\tau]( X_{k-m}))\preceq W\wedge X_k,
$$
and the proof of  \ref{main point} is complete.

Now, using \ref{main point}, the definition of
$X$, and \ref{distribute}, we find:
\Equation\label{nasty}
\begin{aligned}
W\wedge X&=W\wedge\bigvee_{k=1}^\infty
X_k=\bigvee_{k=1}^\infty (W\wedge
X_k)\\
&=\bigg(\bigvee_{k=1}^m\emptyset\bigg)\vee\bigg(  \bigvee_{k=m+1}^\infty 
\noodge_m(\noodge_m(W)\wedge[\tau](X_{k-m}))\bigg)\\
&=\bigvee_{k=1}^\infty 
\noodge_m(\noodge_m(W)\wedge[\tau](X_{k})).
\end{aligned}
\EndEquation
Since $\dom\noodge_m=\range\noodge_m=V_m$, it follows from the
discussion in \ref{restriction} that $Y\mapsto\noodge_m(Y)$ is an
automorphism of the partially ordered set $\{Y:Y\preceq
V_m\}\subset\barcaly_-({\toldTheta})$. Likewise, 
$Y\mapsto[\tau](Y)$ is an
automorphism of the partially ordered set $\barcaly_-({\toldTheta})$. Using
the observation that
an automorphism of a partially ordered set commutes with the formation
of supremum and infima, and another application of \ref{distribute}, we obtain
$$
\begin{aligned}
\bigvee_{k=1}^\infty 
\noodge_m(\noodge_m(W)\wedge[\tau](X_{k}))
&=
\noodge_m\bigg(\bigvee_{k=1}^\infty
(\noodge_m(W)\wedge[\tau](X_{k}))\bigg)\\
&=
\noodge_m\bigg(
\noodge_m(W)\wedge \bigvee_{k=1}^\infty[\tau](X_{k}))\bigg)\\
&=
\noodge_m\bigg(
(\noodge_m(W)\wedge [\tau](\bigvee_{k=1}^\infty X_{k}))\bigg)\\
&=\noodge_m(
(\noodge_m(W)\wedge [\tau](X)),
\end{aligned}
$$
which with (\ref{nasty}) gives the conclusion of the lemma.
\EndProof
%YRQPH

\Lemma\label{veggie prop}
Let $\Mh$ be a closed, hyperbolic, orientable $3$-orbifold, and let $\oldTheta$ be an incompressible, closed $2$-suborbifold
of $\oldOmega:=(\Mh)\pl$. Set
$\oldPsi=\oldOmega\cut\oldTheta$ and ${\toldTheta}=\partial\oldPsi$, and
$\tau=\tau_\oldTheta$.
Let $\oldTheta'$ be a union of components of $\oldTheta$. Set
$\oldPsi'=\oldOmega\cut{\oldTheta'}$ and ${\toldTheta}'=\partial\oldPsi'$, so that
${\toldTheta}'$ is canonically identified with a union of components of ${\toldTheta}$.
Suppose that
$\chibar(X''(\oldOmega,\oldTheta,\oldTheta')\wedge
[\tau](X(\oldOmega,\oldTheta,\oldTheta')))<
%\max (
\chibar(X''(\oldOmega,\oldTheta,\oldTheta'))$.
%,\chibar(
%\tau(\oldXi'')))$. 
Then $\chibar(X'(\oldOmega,\oldTheta,\oldTheta'))<\chibar(\oldPhi^-(\oldPsi)\cap{\toldTheta}')$.
\EndLemma

\Proof
 Set $\oldPhi^-=\oldPhi^-(\oldPsi)$,
$X=X(\oldOmega,\oldTheta,\oldTheta')$,  $X'=X'(\oldOmega,\oldTheta,\oldTheta')$  and  $X''=X''(\oldOmega,\oldTheta,\oldTheta')$. For each integer $n\ge1$ set $X_n=X_n (\oldOmega,\oldTheta,\oldTheta')$, $V_n=V_n (\oldOmega,\oldTheta,\oldTheta')$, and $\noodge_n=\noodge_n^{\oldOmega,\oldTheta,\oldTheta'}$. Let $Z^{(1)},\ldots,Z^{(p)}$ denote the components of $X''$ (where $p\ge0$). Applying Lemma \ref{babi also}, with $\otheroldLambda={\toldTheta}$, $Y=[\tau](X)$, $T=X''$ and $T_i=Z^{(i)}$, we  find that $\chibar(X''\wedge[\tau](X))= \sum_{i=1}^p\chibar(Z^{(i)}\wedge[\tau](X))$.  But since the $Z^{(i)}$ are the components of $X''$, we also have 
$\chibar(X'')=\sum_{i=1}^p\chibar(Z^{(i)})$. The hypothesis
$\chibar(X''\wedge
[\tau](X''))<
\chibar(X'')$ can therefore be rewritten as 
$\sum_{i=1}^p\chibar(Z^{(i)}\wedge[\tau](X))<\sum_{i=1}^p\chibar(Z^{(i)})$. It follows that for some index $i_0$ we have
$\chibar(Z^{(i_0)}\wedge[\tau](X))<\chibar(Z^{(i_0)})$. Set $Z=Z^{(i_0)}$, so that 
\Equation\label{Y Z}
\chibar(Z\wedge[\tau](X))<\chibar(Z).
\EndEquation

Since $Z$ is a component of $X''$, it is in particular a component of $X$. Hence by Lemma \ref{here goethe}, $Z$ is a component of $X_{m_1}$ for some ${m_1}\ge1$. By \ref{oldXi} this implies that $Z\in\oldfrakX_{m_1}(\oldOmega,\oldTheta,\oldTheta')$, so that $Z$ is a component of $V_{m_1}$ and $\noodge_{m_1}(Z)\preceq[{\toldTheta}']$. Since
(\ref{i thought so too}) gives
$V_{m_1}\preceq V_n$ for every $n\le {m_1}$, we have $Z \preceq V_n=\dom \noodge_n$ for every $n\le {m_1}$. Let $m\le m_1$ denote the smallest strictly positive index for which $\noodge_{m}(Z)\preceq[{\toldTheta}']$. We claim:

\Claim\label{still crazy}
$Z$ is a component of $X_m$.
\EndClaim

To prove \ref{still crazy}, note that since $m\le m_1$ we have $Z\preceq V_m$. Let $Z'$ denote the component of $V_m$ such that $Z\preceq Z'$ (see \ref{components}). We have $\noodge_m(Z)\preceq\noodge_m(Z')$ (by the order-preserving property
observed in \ref{restriction}). We also 
have $\noodge_m(Z)\preceq[{\toldTheta}']$. Since $Z'$ is connected, and
${\toldTheta}'$ is a union of components of ${\toldTheta}$, it follows
that $\noodge_m(Z')\preceq[{\toldTheta}']$; hence $Z'\in\oldfrakX_m (\oldOmega,\oldTheta,\oldTheta')$, so that $Z'\preceq
X_m\preceq X$. But $Z$ has been seen to be a component of $X$, and $Z'$ is connected and satisfies $Z\preceq Z'$. This implies that $Z'=Z$, and \ref{still crazy} is proved.

It was pointed out in \ref{oldXi} that every component of $X_m$ is a
component of $V_m$. Thus \ref{still crazy} implies that $Z$ is a
component of $V_m$. Hence 
$W:=\noodge_m(Z)$ is a component of $j_m(V_m)$. But since $\range\noodge_m=\dom\noodge_m=V_m$ (see
\ref{in duck tape}), we have $j_m(V_m)=V_m$; thus  $W$ is a
component of $V_m$, and hence $W\preceq V_k$ for $k=1,\ldots,m$. We claim:

\Claim\label{some day}
For $k=1,\ldots,m$ we have $\noodge_k(W)\preceq[\oldTheta-\oldTheta']$.
\EndClaim

To prove \ref{some day}, first consider the case $k=m$. Since
$W=\noodge_m(Z)$, and since $\noodge_m^{-1}=\noodge_m$ by \ref{in duck
  tape}, and $\noodge_m^{-1}(\noodge_m(Z))=Z$ by \ref{special composition}, we have $\noodge_m(W)=Z$. But since $Z$ is a component of $X''$, the definition of $X''=X''(\oldOmega,\oldTheta,\oldTheta')$ (see \ref{oldXi}) gives $Z\preceq[\oldTheta-\oldTheta']$, which proves \ref{some day} in this case.

Next consider the case in which $1\le k<m$. In this case we apply Lemma \ref{j-lemma}, taking $Z$ as above and letting $m-k$ and $m$ play the respective roles of $m$ and $k$ in Lemma \ref{j-lemma}. Since $Z\preceq V_m$, Lemma \ref{j-lemma} implies that $Z\preceq V_{m-k}$, that  $[\tau](\noodge_{m-k}(Z))\preceq V_k$, and that
$\noodge_m|Z=\noodge_{k}\circ [\tau]
\circ(\noodge_{m-k}|Z)$. Thus if we set $P=[\tau](\noodge_{m-k}(Z))$, we
have $\noodge_{k}(P)=W$. 
Since $\noodge_k^{-1}=\noodge_k$ by \ref{in duck
  tape}, and $\noodge_k^{-1}(\noodge_k(P))=P$ by \ref{special composition}, we have $\noodge_k(W)=P$. On the other hand,
since we took $m$ to be the smallest strictly positive index for which
$\noodge_{m}(Z)\preceq[{\toldTheta}']$, we have
$\noodge_{m-k}(Z)\not\preceq[{\toldTheta}']$; as $Z$ is connected, this
implies that $\noodge_{m-k}(Z)\preceq[{\toldTheta-\toldTheta}']$. Since $\tau$ leaves ${\toldTheta-\toldTheta}'=\rho_\oldTheta^{-1}({\oldTheta-\oldTheta'})$ invariant, it follows that $P=[\tau](\noodge_{m-k}(Z))\preceq[{\toldTheta-\toldTheta}']$. This shows that $\noodge_{k}(W)\preceq[\toldTheta-{\toldTheta}']$, and completes the proof of \ref{some day}.

It follows from \ref{some day} that the hypothesis of Lemma \ref{before veggie} is satisfied with the choice of $W$ made above. Hence Lemma \ref{before veggie} guarantees that $W\wedge X=\noodge_m(\noodge_m(W)\wedge
[\tau](X))$, i.e.
\Equation\label{main component}
W\wedge X=\noodge_m(Z\wedge
[\tau](X)).
\EndEquation
Now, using (\ref{Y Z}) and (\ref{main component}), we find that
%\Equation\label{barclay's}
$\chibar(W\wedge X)=\chibar(\noodge_m(Z\wedge [\tau](X))=
\chibar(Z\wedge[\tau](X))<\chibar(Z)=\chibar(\noodge_m(W))$,
i.e.
\Equation\label{vasistas}
\chibar(W\wedge X)<\chibar(W).
\EndEquation

For each $n\ge1$ we will set $V_n'=V_n\wedge[{\toldTheta}']$. Since
${\toldTheta}'$ is a union of components of ${\toldTheta}$, the
components of $V_n'$ are simply the components of $V_n$ that are
$\preceq[{\toldTheta}']$. For any indices $n,p$ with $n\ge p\ge 1$, we
have $V_n\preceq V_p$ by (\ref{i thought so too}); by \ref{join and
  meet} and lattice theory, it follows that  $V_n'\preceq V_p'$. We will also set $\oldPhi'=\oldPhi^-\cap{\toldTheta}'$. By (\ref{and so did you}) we have $V_1=[\oldPhi^-]$, and hence $V_1'=[\oldPhi']$.
%_i_j

We have observed that $W$ is a component of $V_m$. On the other hand,  by our choice of $m$ we have $W=\noodge_m(Z)\preceq[{\toldTheta}']$. Hence $W$ is a component of $V_m'$.

Let us now fix an element $\oldUpsilon$ of $\Theta_-({\toldTheta})$ such that $[\oldUpsilon]=W$. Since $W\preceq V_m'\preceq V_1'=[\oldPhi']$, 
%$\oldUpsilon\subset\inter(\oldPhi^-.\cap{\toldTheta}')$, so that . and 
%by \ref{in duck tape} we have $V_m'\preceq V_1=[\oldPhi^-\cap{\toldTheta}']$, 
we may choose $\oldUpsilon$ within its isotopy class so that $\oldUpsilon\subset\inter\oldPhi'$. 
Let $\frakU$ denote the union of all negative components of the orbifold $\oldPhi'-\inter\oldUpsilon$, so that $\frakU\in\Theta_-({\toldTheta})$. Set $U=[\frakU]$. We claim:

\Claim\label{trouble in mind}
For every component $Q$ of $X'$, at least one of the following alternatives holds: (i) $Q\preceq W$, or (ii) $Q\preceq U$.
\EndClaim

To prove \ref{trouble in mind}, first note that any component $Q$ of
$X'$ is in particular a component of $X$; hence by Lemma \ref{here
  goethe}, $Q$ is a component of $X_q$ for some
$q\ge1$, which by  \ref{oldXi} implies that it is
a component of $V_q$. Since $Q$ is a component of $X'$ we also have $Q\preceq[{\toldTheta}']$, and hence $Q$ is a component of $V_q'$. By (\ref{i thought so too}) we have $V'_q\preceq V'_m\preceq V'_1=[\oldPhi']$ if $q\ge m$, and $V'_m\preceq V'_q\preceq [\oldPhi']$ if $q\le m$. 
Hence there are elements $\oldGamma_m$ and $\oldGamma_q$ of $\Theta_-({\toldTheta})$ such that $[\oldGamma_m]=V'_m$, $[\oldGamma_q]=V'_q$,  and either $\oldGamma_q\subset\oldGamma_m\subset\inter\oldPhi'$ or $\oldGamma_m\subset\oldGamma_q\subset\inter\oldPhi'$. Since $W$ is a
  component of $V'_m$, there is a component $\oldGamma^0_m$ of
  $\oldGamma_m$ such that $W=[\oldGamma^0_m]$. Since
  $[\oldGamma^0_m]=W=[\oldUpsilon]$, the suborbifolds $\oldGamma^0_m$ and
  $\oldUpsilon$ of $\oldPhi'$ are isotopic in ${\toldTheta}$. It therefore
  follows from Corollary \ref{i guess} that they are isotopic in
  $\oldPhi'$. Hence after possibly modifying $\oldGamma^0_m$ and
  $\oldGamma_q$ by a (single) isotopy in $\oldPhi'$ we may assume that
  $\oldGamma^0_m=\oldUpsilon$, i.e. that $\oldUpsilon$ is a component of
  $\oldGamma_m$.

Now since $Q$ is a component of $V_q'$, there is a component $\frakQ$
of $\oldGamma_q$ such that $Q=[\frakQ]$. Since $\frakQ$ and $\oldUpsilon$
are respectively components of $\oldGamma_q$ and $\oldGamma_m$, and
since we have either
$\oldGamma_q\subset\oldGamma_m\subset\inter\oldPhi^-$ or
$\oldGamma_m\subset\oldGamma_q\subset\inter\oldPhi^-$, we must have
either (a) $\frakQ\subset\oldUpsilon$, or (b)
$\frakQ\cap\oldUpsilon=\emptyset$, or (c) $\oldUpsilon\subset\frakQ$. 

If Alternative (c) holds then $W\preceq Q\preceq X'\preceq X$, and
hence $W\wedge X=W$. This is impossible, since $\chibar(W\wedge
X)<\chibar(W)$ by (\ref{vasistas}). Alternative (a) immediately
implies Alternative (i) of the conclusion of \ref{trouble in mind}. If
Alternative (b) holds, then since $\frakQ$ is negative and taut,
the component of $\oldPhi'-\inter\oldUpsilon$ containing $\frakQ$ must be negative, and is therefore by definition a component of $\frakU$; this implies Alternative (ii) of the conclusion, and the proof of \ref{trouble in mind} is thus complete.
%n

To complete the proof of the lemma, let $Q_1,\ldots,Q_s$ denote the
components of $X'$. (Here $s\ge0$.) We may take the $Q_i$ to be
indexed in such a way that $Q_i\preceq W$ when $1\le i\le t$, and
$Q_i\not\preceq W$ when $t< i\le s$, where $t$ is an integer with
$0\le t\le s$. According to \ref{trouble in mind}, we have $Q_i\preceq
U$ when $t<i\le s$. Set $G=\bigvee_{i\le t}Q_i$ and
$G^*=\bigvee_{i>t}Q_i$. Since we have $Q_i\preceq W$ when $i\le t$ and
$Q_i\preceq U$ when $i>t$, it follows from Proposition \ref{mystery}
that $G\preceq W$ and $G^*\preceq U$. Since $G$ is a join of components
of $X'$ we have $G\preceq X'\preceq X$ and hence $G\preceq X\wedge W$; by Proposition \ref{new partial order} and (\ref{vasistas}) it follows that $\chibar(G)\le\chibar(X\wedge W)<\chibar(W)$. Since $G^*\preceq U$, Proposition \ref{new partial order} implies that $\chibar(G^*)\le\chibar(U)$. But since $Q_1,\ldots,Q_s$ are the components of $X'$, the discussion in \ref{join and meet} shows that an element $\oldXi$ of $\Theta_-({\toldTheta})$ with $[\oldXi]=X'$ can be written as a disjoint union $\oldDelta\discup\oldDelta^*$ where $[\oldDelta]=G$ and $[\oldDelta^*]=G^*$. Hence
\Equation\label{last i think}
\chibar(X')=\chibar(G)+\chibar(G^*)<\chibar(W)+\chibar(U).
\EndEquation
Since $W=[\oldUpsilon]$, and since $U=[\frakU]$ where $\frakU$ is the union of all negative components of the taut orbifold $\oldPhi'-\inter\oldUpsilon$,  we have $\chibar(W)+\chibar(U)=\chibar(\oldUpsilon
)+\chibar(\frakU)=\chibar(\oldPhi')$.
Thus (\ref{last i think}) becomes $\chibar(X')<\chibar(\oldPhi')=\chibar(\oldPhi^-\cap{\toldTheta}')$, which is the conclusion of the lemma.
\EndProof
%\oldXi\frakX\frakZ\calq\eta\oldPhi

%\mu

\chapter{Hyperbolic $3$-orbifolds with (possibly) reducible underlying manifolds  }\label{general chapter}

The goal of this chapter is to prove the results that were referred to
in the Introduction as Theorems B, C, and D (Theorems \ref{manifold
  homology bound}, \ref{orbimain} and \ref{other orbimain}). The
reader who has examined the sketch of the proof of Theorem B in
the Introduction will recognize from the section titles that Sections \ref{dandy
  section}, \ref{semi-dandy section}, \ref{dandy existence section},
and \ref{clash section} are devoted to concepts and results that
were mentioned in that sketch. Section \ref{sphere section} gives some
background about spheres in $3$-manifolds, related for example to the
arguments of \cite{milnorprime}, that is applied in later sections to
spheres which are the underlying surfaces of suborbifolds of a given $3$-orbifold. Section
\ref{structure section} establishes the crucial fact that the graph
$\calf$ described in the sketch of the proof of Theorem B is a forest,
along with some other important properties of the forest mentioned in
the sketch. These ingredients are assembled in Section \ref{manifold
  hom section} to prove Theorem B, and in Section \ref{main theorem
  section} we combine Theorem B with Proposition A (Proposition
\ref{my little sony}) to prove Theorems C and D.

\section{Spheres in $3$-manifolds}\label{sphere section}

\Definition\label{complete def}
A {\it system of spheres} in an orientable $3$-manifold $M$ is a (possibly empty) closed $2$-manifold $\cals\subset \inter M$
such that every component of $\cals$ is a $2$-sphere.
A system of spheres $\cals$ will be termed {\it complete} if
every component of $M\cut\cals$ is $+$-irreducible (see
Definition \ref{P-stuff}).
\EndDefinition

\Lemma\label{that's true tooz}
If $X$ is a $+$-irreducible, compact, orientable $3$-manifold whose boundary components are all $2$-spheres, and $\cals$ is a system of spheres in $X$, then $X\cut\cals$ is $+$-irreducible.
\EndLemma

\Proof
Since $X^+$ is irreducible, every $2$-sphere in $\inter X^+$ separates $X^+$. But $X$ is homeomorphic to a submanifold of $X^+$, and hence every $2$-sphere in $\inter X$ separates $X$. In particular each component of $\cals$ separates $X$; it follows that if $Y$ is a component of $X\cut\cals$, then $Y$ is homeomorphic to the closure of a component of $X-\cals$, and in particular to a submanifold $Y'$ of $X^+$. It suffices to prove that $Y'$ is $+$-irreducible. Each component $S$ of $\partial Y'$ is a $2$-sphere and therefore, in view of the irreducibility of $X^+$, bounds a $3$-ball $B_S\subset X^+$. For each component $S$ of $\partial Y'$ we have either $B_s\cap Y'=S$ or $B_S\supset Y'$. If $B_S\cap Y'=S$ for every component $S$ of $\partial Y'$, then by the connectedness of $X^+$ we have $X^+=Y'\cup\bigcup_{S\in\calc(\partial Y')}B_S$; thus $(Y')^+$ is homeomorphic to $X^+$ and is therefore irreducible. If for some component $S_0$ of $\partial Y'$ we have $B_{S_0}\supset Y'$, then $Y'$ is a submanifold of a $3$-ball. Since the components of $Y'$ are $2$-spheres, $Y'$ is a $3$-sphere-with-holes (see \ref{P-stuff}), so that $(Y')^+$ is homeomorphic to $\SSS^3$ and is therefore irreducible.
\EndProof
%XN

\Number\label{explain why}
If a system
of spheres $\cals_1$ in a closed, orientable $3$-manifold contains a complete system of spheres $\cals$,
then $\cals_1$ is itself complete. Indeed, every component $Y$
of $M\cut{\cals_1}$ is obtained by cutting a component $X$ of $M\cut\cals$
along a system of spheres;  since $X$ is $+$-irreducible, it follows from Lemma \ref{that's true tooz} that $Y$ is $+$-irreducible. 
\EndNumber
%X_1

%\redcommentB{Strangely, \ref{explain why} is not quoted, although it provides context for the next def. Is it used implicitly?

%I have deleted Prop ``min consequence,'' because it was not quoted. It can be found in newester7.tex. It said that if $\cals$ is a  minimal complete system of spheres in a connected, closed, orientable $3$-manifold $M$, and some component of $M\cut\cals$ is a $3$-sphere-with-holes, then $M\cut\cals$ is connected.

%Likewise, I've removed Proposition ``ambient sum'' because it doesn't seem to be quoted anywhere. It can be found in never.tex.
 %}

\Definition\label{min def}
A complete system of spheres $\cals$ in a closed, orientable $3$-manifold $M$ is said to be {\it minimal} if there is no complete system of spheres in $M$ which is the union of a proper subset of the components of $\cals$.
\EndDefinition

\Proposition\label{kneser}
Every closed, orientable $3$-manifold contains a complete system of
spheres.
\EndProposition

\Proof
According to
 Kneser's prime decomposition theorem
\cite[Theorem 3.15]{hempel},
%}
any closed, orientable $3$-manifold $M$ may be written as a connected
sum $M_1\#\cdots\#M_n$ of prime manifolds. We may take the $M_i$ to be
indexed in such a way that, for some $k$ with $0\le k\le n$, the
manifold $M_i$ is homeomorphic to $\SSS^2\times \SSS^1$ whenever $1\le i\le
k$ and is irreducible whenever $k<i\le n$. We may write $M=X_1\cup\cdots\cup
X_n$, where $X_{i}\cap X_{j}$ is empty whenever $i,j\in\{1,\ldots,n\}$ differ by at least $2$, and is a sphere $S_i$ whenever $1\le i<n$ and $j=i+1$; the
$X_i$ have pairwise disjoint interiors; and $ X_i^+$ is homeomorphic
to $M_i$. For $1\le i\le k$, let $t_i\subset \inter X_i$ be a sphere
that is mapped to $\SSS^2\times\{{\rm pt}\}$ under a homeomorphism from $\widehat
X_i$ to $\SSS^2\times \SSS^1$. Then $\bigcup_{1\le i<n}S_i\cup \bigcup_{1\le
  i\le k}T_i$ is a complete system of spheres in $M$.
%; if $n=1$ and $k=1$, then
%$M$ may be homeomorphically identified with $\SSS^2\times\\SSS^1$, and
%$\SSS^2\times{\rm pt}$ is a complete system of spheres. If $n>1$, there exist spheres
%$S_1,\ldots,S_{n1}$ such that $\SSS^1$ bounds a submanifold homeomorphic
%to the complement of the interior of a ball in 
\EndProof

\Proposition\label{noble}
If $\cals$ is a minimal complete system of spheres in a closed, orientable $3$-manifold $M$, then no component of $\cals$ bounds a ball in $M$.
\EndProposition

\Proof
Suppose that some component $S_0$ of $\cals$ bounds a ball $B\subset M$. Among all $3$-balls in $M$ that are bounded by components of $\cals$, we may take $B$ to be minimal with respect to inclusion; then $\cals\cap\inter B=\emptyset$. Let 
$\tB$ denote the component of $M\cut\cals$ which is mapped homeomorphically onto $B$ by $\rho_\cals$. Then $\tS_0:=\partial\tB$ is a component of $\rho_\cals^{-1}(S_0)$, and the other component of $\rho_\cals^{-1}(S_0)$, which we shall denote by $\tS_1$, is a boundary component of some component $X_1$ of $M\cut\cals$. Now if we set $\cals'=\cals-S$, then $M\cut{\cals'}$ has a component which is obtained from the disjoint union $X_1\cup B_0$ by gluing $\tS_1$ to $\tS_0$; and every other component of $M\cut{\cals'}$ is homeomorphic to a component of $M\cut\cals$. Hence every component of $(M\cut{\cals'})^+$ is homeomorphic to a component of $(M\cut\cals)^+$, and is therefore irreducible by the completeness of  $\cals$. This shows  that $\cals'$ is complete, a contradiction to the minimality of $\cals$.
\EndProof

\Definition\label{involvement def}
Let $\cals$ be a minimal complete system of spheres in a closed, orientable $3$-manifold $M$. Set $N=M\cut\cals$ and $\rho=\rho_\cals$. Let $S$ be a component of $\cals$, and let $T\subset M$ be a $2$-sphere with $T\cap\cals=\emptyset$. We will say that $T$ {\it involves} $S$ if there is a submanifold $Z$ of $N$ such that
\begin{itemize}
\item $Z$ is a $3$-sphere-with-holes,
\item $\Fr_NZ=\rho^{-1}(T)$, and
\item $\partial Z$ contains exactly one component of
  $\rho^{-1}(S)$. Note that the truth this condition depends on the
  system $\cals$ as well as the spheres $T$ and $S$; when it is
  necessary to specify the system in question, we will say that $T$
  involves the component $S$ of the system $\cals$.
\end{itemize}
\EndDefinition

\Remark In the special case when $M$ is a connected sum of copies of
$\SSS^2\times \SSS^1$, a minimal complete system of spheres $\cals$ in $M$ defines a
basis of $H_2(M;\FF_2 )$, consisting of the elements $[S]$ for
$S\in\cals$. In this case, if $T\subset M$ is a $2$-sphere with
$T\cap\cals=\emptyset$, it is not hard to show that $T$ involves a sphere
$S\in\calc(\cals)$ if and only if $[S]$ has coefficient
$1$ in the expansion of $[T]$ in the basis defined by
$\cals$. Although this fact will not be used, it is the source of the
term ``involves'' and provides motivation for some of the results in
this section. For example, Proposition \ref{just do it} below is the
analogue of a standard replacement principle in linear algebra.
\EndRemark

\Proposition\label{just do it}
Let $\cals$ be a complete system of spheres in a closed, orientable $3$-manifold $M$. 
%Set $N=M\cut\cals$ and $\rho=\rho_\cals$. 
Let $S$ be a component of $\cals$, and let $T\subset M$ be a $2$-sphere with
$T\cap\cals=\emptyset$. Suppose that $T$ involves the component $S$ of
$\cals$. Then
$(\cals-S)\cup T$ is  a complete system of spheres in
$M$. \abstractcomment{\tiny I had the sentence ``Furthermore, if $T$ is connected [which it now is anyway], then $\cals\glurk$ is
minimal.'' I'm almost certain I don't need this
  sentence. If it's re-instated, make it clear that $\cals\glurk$ means $(\cals-S)\cup T$.}
\EndProposition

\Proof
Set $N=M\cut\cals$ and $\rho=\rho_\cals:N\to M$.
Let $T^*$ denote the $2$-sphere $\rho^{-1}(T)\subset\inter N$. 
%If we set $\cals\smurk=\cals\cup T$, 
If we set $\cals\smurk=\cals\cup T$, then
$N\cut {T^*}$ is canonically homeomorphic to $M\cut{\cals\smurk}$. 

We may also write $\cals\smurk$ as the disjoint union of $S$ with $\cals'':=(\cals-S)\cup T$.  Hence if we set 
$X=M\cut{\cals\glurk}$,
% and $\sigma=\rho_{\cals\glurk}:X\to M$, 
and let $S^*$ denote the $2$-sphere 
$\rho_{\cals\glurk}^{-1}(S)\subset\inter X$, 
%If we set $\cals\smurk=\cals\cup T$, 
%If we set $\cals\smurk=\cals\cup T$, 
then
$X\cut {S^*}$ is canonically homeomorphic to $M\cut{\cals\smurk}$. Thus
$X\cut {S^*}$ and $N\cut {T^*}$ are canonically homeomorphic. If we set $\alpha=\rho_{T^*}:N\cut{T^*}\to N$ and $\beta=\rho_{S^*}:X\cut {S^*}\to X$, the canonical homeomorphism from $X\cut {S^*}$ to $N\cut {T^*}$ carries $\beta^{-1}(S^*)$ onto $\alpha^{-1}(\rho^{-1}(S))$. Hence 
$X$ is homeomorphic to the manifold obtained from $N\cut {T^*}$ by gluing together the two components of 
$\alpha^{-1}(\rho^{-1}(S))$, which are boundary components of $N\cut {T^*}$.

%But we can also write $\cals\smurk$ as the disjoint union of $S$ with $\cals'':=(\cals-S)\cup T$.  
%\cals\cup T}$. If we set $\cals\glurk=(\cals-S)\cup T$, then 
%\cals\smurk}$
%It therefore follows from \redmissingref{a cross-ref that I thought was mentioned elsewhere, but I can't %find it} that every component of  %$M\cut{\cals\smurk}$ is $+$-irreducible. 

Since $T$ involves $S$, there is
a submanifold $Z$ of $N$, homeomorphic to a $3$-sphere-with-holes, such that
$\Fr_NZ=T^*$, and
$\partial Z$ contains exactly one component of $\rho^{-1}(S)$, say $\tS_1$. 
Set $W=\overline{N-Z}$. Then the remaining component of $\rho^{-1}(S)$, say $\tS_2$, 
is contained in $\partial W$. Since $\Fr_NZ=T^*$, the disjoint union $Z\discup W$ is canonically homeomorphic to $N\cut{T^*}$. 
%this shows that $N\cut_{T^*}$ has a component $Z_1$ which is a $3$-sphere-with-holes, that one component of 
%$\alpha{-1}(\rho^{-1}(S))$ is contained in $Z_1$, and that the other component of 
%$\alpha{-1}(\rho^{-1}(S))$ is contained in $W_1:=$N\cut_{T^*}-Z_1$.  
The canonical homeomorphism between them maps the components of $\alpha^{-1}(\rho^{-1}(S))$ onto $\tS_1$ and $\tS_2$. 
Hence 
$X$ is homeomorphic to the manifold obtained from $Z\discup W$ by gluing $\tS_1$ to $\tS_2$.

The proposition asserts that $\cals''$ is  a complete system of spheres in
$M$, which is tantamount to saying that every component of $X=M\cut{\cals\glurk}$ is $+$-irreducible. By hypothesis, $\cals$ is  a complete system of spheres in
$M$, i.e. every component of $N=M\cut\cals$ is $+$-irreducible. It therefore follows from Lemma \ref{that's true tooz} that every component of  $N\cut{T^*}$, and hence every component of $ W$, is $+$-irreducible. Now let $W_0$ denote the component of $W$ containing $\tS_2$. Then each component of $X$ is  either homeomorphic to a component of $W$ distinct from $W_0$, and therefore $+$-irreducible, or
% manifold obtained from $Z\discup W$ by gluing $\tS_1$ to $\tS_2$.
homeomorphic to the manifold $V$ obtained from $Z\discup W_0$ by gluing $\tS_1$ to $\tS_2$. Since $Z$ is a $3$-sphere-with-holes, $\plusV$ is homeomorphic to $\plusw$; since $W$ is $+$-irreducible, it follows that $V$ is $+$-irreducible. 
This shows that each component of $X$ is $+$-irreducible, as required. 
\EndProof

%' ---> \glurk\sigma
%'' ---> \smurk
%glurk ---> ''
\Proposition\label{shmazer}
Let $\cals$ be a minimal complete system of spheres  in a closed, orientable $3$-manifold $M$. 
Let $Y\subset M$ be a $3$-sphere with holes, and suppose that
$\cals_0:=\cals\cap Y$ is a union of components of $\partial Y$. Set
$\calt=\partial Y-\cals_0$. Then for every component $S$ of $\cals_0$,
there is a component $T$ of $\calt$ which involves the component $S$
of $\cals$.
\EndProposition

\Proof
The statement is vacuously true if $\cals_0=\emptyset$. Now assume that $\cals_0\ne\emptyset$ and that $S$ is a component of $\cals_0$.
Let us 
index the components of $\cals_0$ as $S_1,\ldots,S_m$, where $m\ge1$
and $S_1=S$. Let us 
index the components of $\calt$ as $T_1,\ldots,T_n$, where a priori we
have $n\ge0$. Then $S_1,\ldots,S_m,T_1,\ldots,T_n$ are the components
of $\partial Y$.

Set $N=M\cut \cals$.
Since $\cals_0=\cals\cap Y$ is a union of components of $\partial Y$, there is a natural
identification of $\inter Y$ with a subset of $N$, and if
$\tY$ denotes the closure of $\inter Y$ in $N$, then
$\rho:=\rho_\cals$ maps $\tY$ homeomorphically onto $Y$. Hence we may 
index the components of $\partial\tY$ as
$\tS_1,\ldots,\tS_m,\tT_1,\ldots,\tT_n$, where $\rho$ maps $\tS_i$
homeomorphically onto $S_i$ for $1\le i\le m$, and maps $\tT_j$
homeomorphically onto $t_j$ for $1\le j\le n$.

Set $N=M\cut \cals$.
Since $\cals_0=\cals\cap Y$ is a union of components of $\partial Y$, there is a natural
identification of $\inter Y$ with a subset of $N$, and if
$\tY$ denotes the closure of $\inter Y$ in $N$, then
$\rho:=\rho_\cals$ maps $\tY$ homeomorphically onto $Y$. Hence we may 
index the components of $\partial\tY$ as
$\tS_1,\ldots,\tS_m,\tT_1,\ldots,\tT_n$, where $\rho$ maps $\tS_i$
homeomorphically onto $S_i$ for $1\le i\le m$, and maps $\tT_j$
homeomorphically onto $t_j$ for $1\le j\le n$.

Since $\rho|\tY$ is one-to-one, we have $\rho^{-1}(S_1)\cap\tY=\tS_1$. Hence, if 
$\tS_1'$ denotes the side of
$S$ distinct from $\tS_1$, we have
\Equation\label{that's all folks}
\tS_1'\cap\tY=\emptyset.
\EndEquation

I claim:
\Claim\label{ain't nothing}
 We have $n\ge 1$, i.e. $\calt\ne\emptyset$.
\EndClaim 
To prove \ref{ain't nothing}, assume to the contrary that
$\calt=\emptyset$.  Then $\tY$ is a component of $N$, having $\tS_1$
as one boundary component. It follows from (\ref{that's all folks}) that
%Since $\rho|\tY$ is one-to-one, $\rho$ maps
%the components of $\partial\tY$ distinct from $\tS_1$ onto components
%of $\cals$ distinct from $S=S_1$. In particular, 
$\tS_1'$  lies in a component $W$ of $N$ distinct from
$\tY$.
 Now set $\cals'=\cals-S$ and $N'=M\cut{\cals'}$. Then  $N'$ has
a component $W'$ homeomorphic to the manifold obtained from the disjoint
union $\tY\discup W$ by gluing $\tS_1$ to $\tS_1'$. Since $\tY$ is
homeomorphic to $Y$, the hypothesis implies that $\tY$ is a $3$-sphere
with holes; hence $(W')^+$ is homeomorphic to $W^+$. But
$W^+$ is irreducible since $\cals$ is complete, and it now follows that $W'$ is $+$-irreducible. Every
component of $N'$ distinct from $W'$ is homeomorphic to a component
of $N$ and is therefore $+$-irreducible. Hence $\cals'$ is
complete. This contradicts the minimality of $\cals$, and so
\ref{ain't nothing} is proved.

Now let $N_0$ denote the component of $N$ containing $\tY$. Since
$\cals$ is complete, $N_0$ is $+$-irreducible. Hence for each index $j$ with $1\le j\le n$, the sphere
$\tT_j$ separates $N_0$. It follows that for $j=1,\ldots,n$ there exists a unique component $X_j$ of
$\overline{N_0-\tY}$ such that $X_j\cap\tY=\tT_j$. Since $\tY$ is
homeomorphic to $Y$, and is therefore a $3$-sphere with holes, the
manifold $N_0^+$ is homeomorphic to the connected sum
$X_1^+\#\cdots\#X_n^+$. But $N_0^+$ is irreducible since $N_0$ is
$+$-irreducible, and so at most one of the summands $X_j^+$ can fail to be a
$3$-sphere. This proves:
\Claim\label{a poor thing}
There is at most one index $j\in\{1,\ldots,n\}$ such that $X_j$ is not
a $3$-sphere with holes.
\EndClaim

It follows from \ref{ain't nothing} and \ref{a poor thing} that there
is an index $j_0\in\{1,\ldots,n\}$ such that $X_j$ is a $3$-sphere
with holes for every $j\ne j_0$. After reindexing we may assume that
$j_0=1$, i.e.
\Claim\label{but my own}
 $X_j$ is a $3$-sphere
with holes for every $j$ with $1<j\le n$.
\EndClaim

(Of course \ref{but my own}  is a vacuous statement if $n=1$.)

Set $Z_1=\overline{N_0-\inter X_1}=\tY\cup\bigcup_{1<j\le n}X_j$. 
We have observed that $\tY$ is
homeomorphic to $Y$, and is therefore a $3$-sphere with holes; by
\ref{but my own}, $X_j$ is also a $3$-sphere with holes for $1<j\le n$. Our choice of $X_j$ gives
$X_j\cap\tY=T_j$ for each $j>1$. For
any distinct indices $j$ and $j'$, the components $X_j$ and $X_{j'}$ of 
$\overline{N_0-\tY}$ are distinct since their intersections $T_j$ and $T_{j'}$ with $\tY$ are distinct, and hence
$X_j\cap X_{j'}=\emptyset$. It follows that:
\Claim\label{what i say}
$Z_1$ is
a $3$-sphere with holes.
\EndClaim

We now distinguish two cases. First consider the case in which
$\tS_1'\not\subset\partial Z_1$. In this case, $\tS_1$ is the unique
component of $\rho^{-1}(S)$ contained in $\partial Z_1$, and it follows from
Definition \ref{involvement def}, taking $Z$ to be the
$3$-sphere-with-holes $Z_1$ and noting that $\rho^{-1}(T_1)=\tT_1=\Fr_NZ_1$, that $T_1$
involves $S$. 

Now consider the case in which
$\tS_1'\subset\partial Z_1$.  The definition of $Z_1$ implies that every component of $\partial Z_1$ is a component of 
$\partial\tY$ or of $\partial X_j$ for some $j$ with $1<j\le n$. By (\ref{that's all folks}), $\tS_1'$ cannot be a component of $\partial\tY$. Hence
 there is an index $j_1>1$ such
that $\tS_1'\subset\partial X_{j_1}$. By \ref{but my own}, $X_{j_1}$
is a $3$-sphere-with-holes. Since $\tS_1\subset\tY$, and since $X_{j_1}\cap\tY=\tT_{j_1}$ is disjoint from $\tS_1$, we have
$\tS_1\not\subset X_{j_1}$; hence $\tS_1'$ is the unique
component of $\rho^{-1}(S)$ contained in $\partial X_{j_1}$. It therefore follows from
Definition \ref{involvement def}, taking $Z$ to be the
$3$-sphere-with-holes $X_{j_1}$
and noting that $\rho^{-1}(T_{j_1})=\tT_{j_1}=\Fr_N X_{j_1}$, that $T_{j_1}$
% that $t_{j_1}=\Fr_NX_{j_1}$
involves $S$. 
\EndProof

\abstractcomment{\tiny I'm planning to replace the following result.

\Proposition\label{holester}
Let 
$\cals$ be a complete system of spheres in a closed, orientable
hyperbolic $3$-manifold $M$. Let $S^!$ be a sphere contained in
$M-\cals$, and let $X$ denote the component of $M-\cals$ containing
$S^!$. Let
t$Y$ be a component of $(X-S^!) \,\hat\empty$, and let $\tS$ be a component of
$\partial\hat X\subset (X-S^!) \,\hat\empty$
%\subset \partial(X-S^!)^+$
which is {\it not} contained in
$\partial Y$. Suppose that either (a)
$\grock_\cals(\tS)$ is a
component of $\partial Y$, 
(b) $Y$ is {\it not} a $3$-sphere with holes. Then 
%Definition \ref{P-stuff}). 
%such that
%$S^{!}\cap\cals=\emptyset$, and such that,
$\cals^!:=(\cals\cup S^!)-\rho(\tS)$ is a complete system of spheres
in $M$.
\EndProposition

\Proof
The proof summary said ``A straightforward cut and paste argument which I think
  will be the main ingredient in the proof of the next prop.''
\EndProof
%+\plus

}

%\redcomment{The following result will be a corollary of an new proposition if it survives. The
%proof will be replaced.

\Corollary\label{taser}
Let $\cals$ be a complete system of
spheres in a closed, orientable $3$-manifold $M$. Let $D\subset M$ be
a disk with $D\cap\cals=\partial D$, let $S$ denote the component of
$\cals$ containing $\partial D$, let $D_1$ and $D_2$ denote the
closures of the components of $S-\partial D$, and for $i=1,2$ set
$S_i=D\cup D_i$ and $\cals_i=(\cals-S)\cup S_i$. Then either $\cals_1$
or $\cals_2$ is a complete system of spheres in $M$.
\EndCorollary

\Proof
Let $Y_0$ be a regular neighborhood of $S\cup D$, disjoint from $\cals-S$. Then $Y_0$ is a
$3$-sphere with three holes, and its
boundary components are isotopic in $M':=M-(\cals-S)$ to $S$, $S_1$ and
$S_2$. Hence $Y_0$ is isotopic in $M'$ to a $3$-sphere with three
holes $Y$ having $S$ as a boundary component; the other two boundary
components of $Y$, which are isotopic in $M':=M-(\cals-S)$ to $S_1$ and
$S_2$, will be denoted respectively by  $S_1'$ and $S_2'$. 
%We have
%$\tS\cap Y=S$. 
Thus $\cals$ and $Y$ satisfy the hypothesis of
Proposition \ref{shmazer}, 
%\redcomment{Not quite. It's stated with a minimality hypothesis. Either add that hyp. here (and check it in apps.), or get rid of it in Proposition \ref{shmazer}.} 
and in the notation of that proposition we
have $\cals_0=S$ and $\calt=S_1'\cup S_2'$. Hence Proposition
\ref{shmazer} implies that $S_m'$ involves the component $S$ of
$\cals$ for some
$m\in\{1,2\}$. We may therefore apply Proposition \ref{just do it},
with $T=S_m'$, to deduce that $(\cals-S)\cup S_m'$ is a complete
system of spheres in $M$. Since $(\cals-S)\cup S_m'$ is isotopic to
$(\cals-S)\cup S_m=\cals_m$, the conclusion follows.
\EndProof

\abstractcomment{\tiny I'm in the process of reworking Section 11 using these results and the cheat sheet.}

%%minimal proofreadingnote comment missingref comment

\section{Admissible systems of spheres}\label{dandy section}

While the contents of  Sections \ref{dandy
    section}---\ref{structure section} were sketched in the
  Introduction and at the beginning of this chapter, it seems worth repeating that the goal of these sections to establish machinery that
  allows one to bound the number of spheres in a suitably chosen complete system of
spheres in a manifold of the form $|(\Mh)\pl|$, where $\Mh$ is a hyperbolic $3$-orbifold. In Sections
\ref{manifold hom section} and \ref{main theorem section}, this machinery will be combined with the results of Section \ref{irr-M
  section}  to prove the main results of the monograph.

\Definition
Let $\Mh$ be a closed, orientable hyperbolic $3$-orbifold, and set $\oldOmega=(\Mh)\pl$. An
{\it $\oldOmega$-admissible system of spheres} in $|\oldOmega|$ is a compact $2$-manifold $\cals\subset
|\oldOmega|$ such that every component of $\cals$ is a sphere,
$\cals$ is transverse to $\fraks_\oldOmega$, and
$\obd(\cals)$ is incompressible in $\oldOmega$. Note that an admissible
system of spheres may be empty, and it may be connected (i.e. consist
of a single sphere). We shall say ``admissible'' in place of
$\oldOmega$-admissible'' when it is clear which orbifold is involved.
\EndDefinition

\Number\label{it's ok}
Let $\Mh$ is a closed, orientable, hyperbolic $3$-orbifold. Set $\oldOmega=(\Mh)\pl$, and let $\cals$ be an
$\oldOmega$-admissible system of spheres in $M:=|\oldOmega|$. Then $\oldOmega$ is strongly \simple\ and $\obd(\cals)$ is incompressible in $\oldOmega$. Hence, if $X$ is any union of
connected components of $M-\cals$, then according to  \ref{t-defs} the invariants $t_\oldOmega(X)$, $s_\oldOmega(X)$, $s'_\oldOmega(X)$, $q_\oldOmega(X)$ and $y_\oldOmega(X)$ are well defined.
\EndNumber

\begin{remarksnotationdefinitions}\label{doublemint}
Let $\Mh$ be a closed, orientable hyperbolic $3$-orbifold, and set $\oldOmega=(\Mh)\pl$.
Let $\cals$ be an
$\oldOmega$-admissible system of spheres in $M:=|\oldOmega|$. Set $\oldTheta=\obd(\cals)$,
$\oldPsi=\oldOmega\cut\oldTheta$, and $\tcals=\partial
M\cut\cals=\partial |\oldPsi|$. The $\oldOmega$-admissibility of
$\cals$ means that $\oldTheta$ is incompressible, so that all the observations and conventions laid out in \ref{in duck tape} apply. In particular,
$\oldPsi$ is componentwise strongly 
\simple\ and componentwise
boundary-irreducible; $\obd(\tcals)$ is negative; $\oldSigma(\oldPsi)$, 
$\oldPhi(\oldPsi)$,
$\kish(\oldPsi)$, $\book(\oldPsi)$ and
$\frakA(\oldPsi)$ are well defined; and $\noodge_n^{\oldOmega,\oldTheta}$ and $V_n^{\oldOmega,\oldTheta}$ are defined for every $n\ge1$. Furthermore, if $\cals'$ is a subsystem of $\cals$, and if we set $\oldTheta'=\obd(\cals')$, then $X(\oldOmega,\oldTheta,\oldTheta')$, $X'(\oldOmega,\oldTheta,\oldTheta')$, $X''(\oldOmega,\oldTheta,\oldTheta')$, and 
$\xi_n(\oldOmega,\oldTheta,\oldTheta')$
$X_n(\oldOmega,\oldTheta,\oldTheta')$ for each $n\ge1$, are defined by \ref{oldXi}. All this structure will be invoked in this section and the succeeding ones when we are dealing with an admissible system $\cals$. In particular, the componentwise strong \simple ity of $\oldPsi:=\oldOmega\cut\oldTheta$, the consequence that $\oldSigma(\oldPsi)$, 
$\oldPhi(\oldPsi)$,
$\kish(\oldPsi)$, $\book(\oldPsi)$ and
$\frakA(\oldPsi)$ are well defined,  and the negativity of $\obd(\tcals)$ will often be used without necessarily being explicitly mentioned.

When an admissible system $\cals$ is given, 
a component $\tS$ of  $\tcals$ will be
termed  {\it \full} (relative to $\cals$) if we have
$\obd(\tS)\subset\book(\oldPsi)$ (cf. \ref{doublemint}), or equivalently if $\obd(\tS)\cap\kish(\oldPsi)=\emptyset$.
A component of $\cals$
%and let $\tS$ and $\tS'$ denote
%the sides (see \ref{nbhd stuff}) of $S$. The sphere $S$ 
will be termed
{\it \doublefull} (relative to $\cals$) if both its sides (see \ref{nbhd stuff}) are \full.
\end{remarksnotationdefinitions}

\Remark\label{special ouenelitte}
Let $\Mh$ be a closed, orientable hyperbolic $3$-orbifold, and set $\oldOmega=(\Mh)\pl$.
Let $\cals$ be an
$\oldOmega$-admissible system of spheres in $M:=|\oldOmega|$. Set $\oldPsi=\oldOmega\cut{\obd(\cals)}$, and $\tcals=\partial
M\cut\cals=\partial |\oldPsi|$. It follows immediately from Corollary \ref{abu gnu} that a component $\tS$ of $\tcals$ is \full\ if and only if every component of $\overline{\obd(\tS)\setminus\oldPhi(\oldPsi)}$ is an annular orbifold.
\EndRemark

\Proposition\label{hugh manatee prop}
Let $\Mh$ be a closed, orientable
hyperbolic $3$-orbifold, and set $\oldOmega=(\Mh)\pl$. Suppose that 
$\oldOmega$ contains no  embedded negative
turnovers. Set $M=|\oldOmega|$.  Let $\cals$ be an
$\oldOmega$-admissible system of spheres in $M$. Suppose that $X$ is a union of  components of $M-\cals$. Then
%Set $\tcals=\partial M\cut\cals$. 
 the number of non-\full\ 
components of $\partial\hatX$ is bounded above by
$s'_\oldOmega(X)$ (see \ref{t-defs}).  
% is bounded below by the number of non-\full\ 
%components of $\partial M\cut\cals$.
In particular,
%Set $\tcals=\partial M\cut\cals$. 
 the number of non-\full\ 
components of
$\partial M\cut\cals$
  is bounded above by
$s'_\oldOmega(M-\cals)$.
\EndProposition

\tinymissingref{\tiny  I have removed Prop. ``oh THAT hugh manatee!'' from this
  version. It can be found in eclipse4.tex. 
}

\Proof
 Set
 $\oldPsi=\omega(\hatX)$. 
Since 
$\oldOmega$ contains no  embedded negative
turnovers, no component of $\partial\oldPsi$ is a negative turnover.

Set $\oldGamma=\partial\oldPsi\cap\kish(\oldPsi)$ (cf. \ref{doublemint}). It follows from Corollary \ref{abu gnu}  that $\oldGamma$ is the union of all non-annular components of
$\overline{\partial\oldPsi-\oldPhi(\oldPsi)}$; since $\partial\oldPsi$ has no toric components by \ref{boundary is negative}, and all components of 
$\overline{\partial\oldPsi-\oldPhi(\oldPsi)}$  have non-positive Euler characteristic by \ref{tuesa day},
$\oldGamma$ is the union of all components of
$\overline{\partial\oldPsi-\oldPhi(\oldPsi)}$ that have strictly negative
Euler characteristic.
 In particular, $\oldGamma$ is negative, and it is not a negative turnover since no component of $\partial\oldPsi$ is a negative turnover.
Proposition \ref{at least a sixth} therefore gives
\Equation\label{de shmoozer}
\compnum(|\oldGamma|)\le6
\chibar(\oldGamma ).
\EndEquation

If it happens that  $\lambda_\oldPsi=2$, then according to
Proposition \ref{fortunately},
for every
integer $d>1$, the number of points of order $d$ in 
$\fraks_{\overline{\partial\oldPsi-\oldPhi(\oldPsi)}}$ is even. But every component of 
$\overline{\partial\oldPsi-\oldPhi(\oldPsi)}$ not contained in $\oldGamma$ is annular; since such a component is automatically orientable, its underlying surface is therefore either an annulus with no
singular points, or a disk with exactly two singular points, each of
which has order $2$. Hence for every
integer $d>1$, the number of points of order $d$ in $\fraks_\oldGamma$ has the
same parity as the number of points of order $d$ in 
$\fraks_{\overline{\partial\oldPsi-\oldPhi(\oldPsi)}}$, and is
therefore even. 
%\redcomment{I think I get it.} 
In
particular the number of points of order $2$ in $\oldGamma$ is even. Since
no component of $\partial\oldPsi$ is a negative turnover, no component of $\oldGamma$
is a negative turnover.  Thus in the case $\lambda_\oldPsi=2$,  it follows from 
Proposition \ref{at least a sixth} that
\Equation\label{do i need all this?}
\compnum(|\oldGamma|)\le2\lfloor
3\chibar(\oldGamma )\rfloor.
\EndEquation

On the other hand, since every component of 
$\overline{\partial\oldPsi-\oldPhi(\oldPsi)}$ not contained in $\oldGamma$ has
Euler characteristic $0$, we have 
$\chibar(\oldGamma)=\chibar(\overline{\partial\oldPsi-\oldPhi(\oldPsi)})$. Since
by \ref{tuesa day} we have
$2\chibar(\kish\oldPsi))=
2\chibar(\overline{(\oldPsi-\oldSigma(\oldPsi)})=
\chibar(\overline{\partial\oldPsi-\oldPhi(\oldPsi)})$, it follows that
\Equation\label{Will Haggle}
\chibar(\oldGamma)=2\chibar(\kish(\oldPsi)).
\EndEquation
Now  (\ref{de shmoozer})  and (\ref{Will Haggle}) give
$
\compnum(|\oldGamma|)\le12
\chibar(\kish(\oldPsi) )$, which in view of the definition of $\sigma(\oldPsi)$ given in \ref{t-defs}, and the integrality of $\compnum(\oldGamma)$, may be written as 
$\compnum(|\oldGamma|)\le\lfloor\sigma(\oldPsi)\rfloor$.
If $\lambda_\oldPsi=2$, then 
 (\ref{do i need all this?}) and (\ref{Will Haggle}) give
$
\compnum(|\oldGamma|)\le2\lfloor
6\chibar(\kish(\oldPsi) )\rfloor=2\lfloor
\sigma(\oldPsi)/2\rfloor$.
Recalling from \ref{t-defs} that we have $\sigma'(\oldPsi)
%=\sigma(\oldPsi)$ if $\lambda_\oldPsi=1$, 
 %and
%$\sigma'(\oldPsi)
=\lambda_\oldPsi\lfloor{\sigma(\oldPsi)}/\lambda_\oldPsi\rfloor$, we deduce that 
$\compnum(|\oldGamma|)\le
\sigma'(\oldPsi)$ in all cases.
 But according to the definitions in
\ref{t-defs}, we have $
s'_\oldOmega(X)=
\sigma'(\oldPsi)$. Hence
\Equation\label{maan toor}
\compnum(|\oldGamma|)\le s'_\oldOmega(X).
\EndEquation

Finally, if $\tS$ is a non-\full\ component of $\partial \oldPsi$,
then by definition we have
$\partial\oldPsi\cap\kish(\oldPsi)\ne\emptyset$. In view of the definition of
$\oldGamma$, this means that 
%\redcomment{I have changed the def. of $\oldGamma$. Rewrite the
%  following stuff accordingly.} By \redmissingref{cross-refs}, every component of 
%$\tS\cap\kish(\oldPsi)$ is a component of
%$\overline{(\partial\oldPsi)-\oldPhi(\oldPsi)}$, and has strictly negative
%Euler characteristic; thus such a component is a component of
%$\oldGamma$. Thus
 every non-\full\
component of $\partial\oldPsi$ contains at least one component of $|\oldGamma|$. Hence the number of
 non-\full\
components of $\partial\oldPsi$ is at most $\compnum(|\oldGamma|)$. With (\ref{maan
  toor}), this gives the first conclusion of the proposition. The second conclusion is simply the special case $X=M-\cals$.
\EndProof

\Corollary\label{hugh manatee cor}
Let $\Mh$ be a closed, orientable
hyperbolic $3$-orbifold, and set $\oldOmega=(\Mh)\pl$. Suppose that 
$\oldOmega$ contains no  embedded negative
turnovers. Set $M=|\oldOmega|$.  Let $\cals$ be an
$\oldOmega$-admissible system of spheres in $M$.
Then the number of non-\doublefull\ 
components of $\cals$ is bounded above by
$s'_\oldOmega(M-\cals)$.
% is bounded below by the number of non-\full\ 
%components of $\partial M\cut\cals$.
\EndCorollary

\Proof
Since by definition every non-\doublefull\ component of $\cals$ has a non-\full\ side, the number of non-\doublefull\ 
components of $\cals$ is bounded above by the number of non-\full\ 
components of $\tcals:=\partial M\cut\cals$. The assertion therefore follows from the second assertion of Proposition \ref{hugh manatee prop}.
\EndProof

\Corollary\label{lose lose} Let $\Mh$ be a closed, orientable,
hyperbolic $3$-orbifold. Set $\oldOmega=(\Mh)\pl$ and $M=|\oldOmega|$. Assume that  $\oldOmega$
contains no  embedded negative turnovers. Let $\cals$ be an
$\oldOmega$-admissible system of spheres in $M$, and let $X$ be a 
connected component of $M-\cals$. If $\partial \hatX$ has at least one
component which is not \full, then $s_\oldOmega(X)\ge\lambda_\hatX$ (see \ref{it's ok}). 
\EndCorollary

\Proof[Proof of Corollary \ref{lose lose}]
It follows from Proposition \ref{hugh manatee prop} that  
$s'_\oldOmega(X)$ is bounded above by
 the number of non-\full\ 
components of $\partial\hatX$; as the latter number is strictly positive by hypothesis, we have
$s'_\oldOmega(X)>0$. But it was observed in \ref{t-defs} that 
$\sigma'(\oldPsi)$ is divisible by $\lambda_\oldPsi$ for any componentwise strongly \simple, componentwise boundary-irreducible $3$-orbifold $\oldPsi$; thus
$s'_\oldOmega(X)=\sigma'(\hatX)$ is divisible by $\lambda_\hatX$, and since $s'_\oldOmega(X)>0$ it follows that $s'_\oldOmega(X)\ge \lambda_\hatX$.
\EndProof

\Proposition\label{loster lemster} Let $\Mh$ be a closed, orientable,
hyperbolic $3$-orbifold, and set $\oldOmega=(\Mh)\pl$.
Let
$\cals$ be a non-empty $\oldOmega$-admissible system of spheres in $M:=|\oldOmega|$, and let $X$ be a
connected component of $M-\cals$. Suppose that $\hatX$ is $+$-irreducible, and that
every component of $\partial \hatX$ is \full. Then $\pi_1(X)$ is cyclic. 
\EndProposition

\Proof Set $N=\hatX$, $\tcals=\partial N$, and $\oldPsi=\obd(N)$.  The
hypothesis that every component of $\tcals$ is \full\ means, by
definition (see \ref{doublemint}), that 
$\tcals\subset|\book(\oldPsi)|$. Since $M$ is connected and $\cals\ne\emptyset$, we have $\tcals\ne\emptyset$ and hence $|\book(\oldPsi)|\ne\emptyset$.

Thus if $|\book(\oldPsi)|$ is a proper submanifold of the connected manifold $N$, the frontier $F$ of $|\book(\oldPsi)|$ in $N$ is non-empty.
Now according to \ref{tuesa day}, $F$ is a union of components of
$|\frakA(\oldPsi)|$, and hence
the components of $\obd(F)$ are essential annular orbifolds
in $\oldPsi$. It follows that the components of $F$ are
%if $A$ is a component of $\Fr_{N}
%|\book(\oldPsi)|$, then $A$ is a
%properly embedded annuli and disks in $N$
%if $A$ is a component of $\Fr_{N}
%|\book(\oldPsi)|$, then $A$ is a
therefore properly embedded annuli and disks in $N$; hence the boundary curves of $F$ are contained in $\Fr_{\tcals}(\tcals\cap|\book(\oldPsi)|)$. But the latter set is empty since
$\tcals\subset|\book(\oldPsi)|$. This shows that $|\book(\oldPsi)|=N$.
In view of the definitions given in \ref{oldSigma def} and \ref{tuesa day}, it follows that:
\Equation\label{bert and harry}
\oldSigma(\oldPsi)=\oldSigma_1(\oldPsi)=\overline{\oldPsi-\frakH(\oldPsi)}.
\EndEquation

Recall from \ref{oldSigma def} that $\frakH(\oldPsi) $ is a strong regular neighborhood of the two-sided annular suborbifold
$\frakQ:=\frakQ(\oldPsi)$ of $\oldPsi$. We claim: 
%(defined up to isotopy) given by Proposition \ref{new characteristic}, and we will denote by 
%Since $|\book(\oldPsi)|=N$, 
%
%each component $\frakC$ of $\obd(\overline{N-|\oldSigma(\oldPsi)|})=
%\overline{\book(\oldPsi)-\oldSigma(\oldPsi)}$
%is a component of $\frakH(\oldPsi)$ such that %$\frakC\cap\oldSigma(\oldPsi)=\Fr_\oldPsi\frakC$.
%In view of Proposition \ref{new characteristic} and Definition \ref{oldSigma def},  it follows that 
%$\obd(\overline{N-|\oldSigma(\oldPsi)|})$
%is a regular neighborhood in $\oldPsi$ of a properly embedded annular suborbifold of %$\oldPsi$; hence $\overline{
%|\book(\oldPsi)|-|\oldSigma(\oldPsi)|}=
%\overline{N-|\oldSigma(\oldPsi)|}
%$ is a regular neighborhood in $N$ of a properly embedded submanifold $\calb$ of $N$ whose components are annuli and disks. We claim:
\Claim\label{yes it do} Every component of $|\frakQ|$ separates $N$.
\EndClaim

 To prove this, let $B$ be any component of $|\frakQ|$. Since $\obd(B)$ is annular and orientable, $B$ is either an annulus or a disk; we first consider the case in which $B$ is an annulus. Since $B$ is properly embedded in $N$, and every boundary component of $N$ is a sphere, the boundary components $C_1$ and $C_2$ of $B$ bound disjoint disks $D_1$ and $D_2$ in $\tcals$. (The disks $D_1$ and $D_2$ are uniquely determined  in the subcase where $C_1$ and $C_2$ lie in the same component of $\tcals$; otherwise there are two choices for each of them.) The sphere $T:=D_1\cup B\cup D_2$ is contained in $N$ and hence in $N^+$. Since $N$ is $+$-irreducible by hypothesis, $N^+$ is irreducible; thus $T$ bounds a ball in $N^+$, and hence $B$ separates $N$.

Similarly, in the case where $B$ is a disk, $\partial B$ lies in a sphere component of $\tcals$ and hence bounds a disk $D\subset\tcals$. The sphere $T:= B\cup D\subset N\subset N^+$ bounds a ball in the irreducible manifold $N^+$, and hence $B$ separates $N$. This completes the proof of \ref{yes it do}.

It follows from \ref{yes it do} that $\pi_1(N)$ is generated by the
images under inclusion of the fundamental groups of the components of
$N-|\frakQ|$. Since $|\frakH(\oldPsi)|$ is a regular
neighborhood  of $|\frakQ|$, and since  $\overline{N-|\frakH(\oldPsi)|}=|\oldSigma(\oldPsi)|$ by (\ref{bert and harry}), each component of $N-|\frakQ|$ deform-retracts
to a component of $|\oldSigma(\oldPsi)|$. Hence, if for each component $\frakU$ 
of $\oldSigma(\oldPsi)$, we denote by $H_\frakU$ the image of the inclusion
homomorphism $\pi_1(|\frakU|)\to\pi_1(N)$, then:

\Claim\label{we know it now}
The group $\pi_1(N)$ is generated by the subgroups $H_\frakU$, where $\frakU$ ranges over the  components of $\oldSigma(\oldPsi)$.
\EndClaim

We claim:

\Claim\label{when is a group not a group}
For every component $\frakU$ of $\oldSigma(\oldPsi)$, the group $H_\frakU$ is cyclic.
\EndClaim

To prove \ref{when is a group not a group}, consider an arbitrary
 component $\frakU$ of $\oldSigma(\oldPsi)$. By (\ref{bert and harry}),
$\frakU$ is a
 component of $\oldSigma_1(\oldPsi)$, and is therefore an
%; hence, according to
 %\redmissingref{cross-reference. Maybe I should say it's an
   \Ssuborbifold. If $\frakU$ is a \bindinglike\ \Ssuborbifold, then 
 by Lemma \ref{when a tore a fold} it is a \torifold; hence by
   Proposition \ref{three-way equivalence} $|\frakU|$ is either a solid torus or a
   ball, and hence $\pi_1$ is cyclic. In particular $H_\frakU$ is cyclic.
%it's a solid
  % torus},  \redproofreadingnote{Actually that missingref (this is a
 %  proofreadingnote) confuses me. I want to say that aa component of
   %$\oldSigma(\oldPsi)$ is an \Ssuborbifold, that if it's Is there something else going on?
   %At the moment I don't even remember what $\oldSigma_1$ is} $|\frakU|$ is either a solid torus or a \pagelike\ \Ssuborbifold\ of
% $\frakU$ may be equipped with an $I$-fibration over a $2$-orbifold in such a way that it is standardly embedded in
% $\oldPsi$. If $|\frakU|$ is a solid torus then $\pi_1(|\frakU|)$ is cyclic, and hence

 Now consider the case in which $\frakU$ is a \pagelike\ \Ssuborbifold\ of $\oldPsi$.
%equipped with an $I$-fibration over a $2$-orbifold and is standardly embedded in $\oldPsi$. 
Fix a component $V$ of $|\frakU|\cap\tcals$. Then by Proposition \ref{doesn't exist}, the image of the inclusion homomorphism $\pi_1(V)\to\pi_1(|\frakU|)$ has index at most $2$ in $\pi_1(|\frakU|)$. Hence if $L$ denotes the image of the inclusion homomorphism $\pi_1(V)\to\pi_1(N)$, we have $|H_\frakU:L|\le2$. On the other hand, since the component of $\tcals$ containing $F$ is a $2$-sphere, the group $L$ is trivial. Thus $H_\frakU$ has order at most $2$, and is therefore cyclic in this case as well. This completes the proof of \ref{when is a group not a group}.

Next, we claim:

\Claim\label{a kind of binding}
If $\frakU$ is any component of $\oldSigma(\oldPsi)$, then $|\frakU|$ is contained in a connected submanifold $W$ of $N$ whose
boundary components are all spheres, and such that the inclusion homomorphism $\pi_1(|\frakU|)\to\pi_1( W)$ is surjective.
\EndClaim

To prove \ref{a kind of binding},  first note that if $\partial |\frakU|\subset\tcals$, then since $\partial|\frakU|$ is closed, every component of $\partial|\frakU|$ is a sphere, and \ref{a kind of binding} follows if we set $W=|\frakU|$. We may therefore assume that $\partial |\frakU|\not\subset\tcals$, i.e. that $\Fr_N|\frakU|\ne\emptyset$. According to \ref{tuesa day}, every component of $\obd(\Fr_N|\frakU|)\subset\frakA(\oldPsi)$ is an annular orbifold, and so every component of $\Fr_N|\frakU|$ has non-empty boundary. Hence there is a component $\tS$ of $\tcals$ such that $\calv:=|\frakU|\cap\tS$ is a non-empty, proper subsurface of $\tS$.

For every component $C$ of $\partial \calv$, let $D_C$ denote the disk
which is contained in the sphere $\tS$, has boundary $C$, and contains
the component of $\calv$ having $C$ as a boundary curve. 
The set of disks of the form $D_C$, where $C$ ranges over the
components of $\partial \calv$, is partially ordered by inclusion and
hence has a minimal element $D_0$. We have $D_0=D_{C_0}$ for some
component $C_0$ of $\partial \calv$. Let $V$ denote the component of
$\calv$ containing $C_0$. Since $D_0=D_{C_0}$, we have $V\subset
D_0$. Let $C_1,\ldots,C_m$ denote the components of $\partial V$
distinct from $C_0$, where $m\ge0$. For $0<i\le m$, let
$D_i\subset\inter D_0$ denote the disk in $\tS$ such that $D_i\cap
V=\partial D_i=C_i$. The minimality  of $D_0$ implies that
$|\frakU|\cap D_i=\calv\cap D_i=C_i$ for $i=1,\ldots,m$. Hence
$Z:=|\frakU|\cup(D_1\cup\cdots\cup D_m)$ is obtained from $|\frakU|$
by attaching disks to $C_1,\ldots,C_m$ along their boundaries. In
particular, $Z$ is connected, and the inclusion homomorphism
$\pi_1(|\frakU|)\to\pi_1(Z)$ is surjective. On the other hand,
$\frakU$ is an \Ssuborbifold\ of $\oldPsi$ by the definition of
$\oldSigma(\oldPsi)$, and it therefore follows from Proposition \ref{doesn't exist}
the inclusion homomorphism $H_1(V;\QQ)\to H_1(|\frakU| ;\QQ)$ is surjective. This implies that the inclusion homomorphism $H_1(V;\QQ)\to H_1(Z;\QQ)$ is surjective; since $V\subset D_0\subset Z$, and $D_0$ is a disk, we deduce that
%\redcomment{Is this a real comment? `Fix the passage from here. Fragment: ``and that the component of $\oldPhi$ in question here is a disk'''} 
that 
$H_1(Z;\QQ)=0$. Hence if $W$ denotes a regular neighborhood of $Z$ in $N$, then $W$ is connected, the inclusion homomorphism $\pi_1(|\frakU|)\to\pi_1(W)$ is surjective, and $H_1(W;\QQ)=0$. It is a standard consequence of Poincar\'e-Lefschetz duality that if a compact orientable $3$-manifold $W$ satisfies $H_1(W;\QQ)=0$, then every boundary component of $W$ is a sphere. This establishes \ref{a kind of binding}.

To establish the conclusion of the proposition, it suffices to prove
that $\pi_1(N)$ is cyclic. For this purpose we will distinguish two cases. In the case where $H_\frakU$ is trivial for every component $\frakU$ of $\oldSigma(\oldPsi)$, it follows from \ref{we know it now} that $\pi_1(N)$ is trivial. 

Now consider the case in which there is a component $\frakU$ of
$\oldSigma(\oldPsi)$ such that $H_\frakU$ is non-trivial. In this case, according
to \ref{a kind of binding},
$|\frakU|$ is contained in a connected submanifold $W$ of $N$ such that (i) the
boundary components of $W$ are all spheres, and (ii) the inclusion homomorphism $\pi_1(|\frakU|)\to\pi_1( W)$ is surjective. After possibly modifying $W$ by a non-ambient isotopy, we may assume that $W\subset\inter N$. It then follows from (i) that
the inclusion homomorphism $\pi_1(W)\to\pi_1(N)$ is injective.
According to (ii) and the definition of $H_\frakU$, it now follows that $\pi_1(W)$ is
isomorphic to $H_\frakU$. In particular $\pi_1(W)$ is non-trivial. On the
other hand, (i) implies that
%$W$ is non-ambiently isotopic to a submanifold of $\inter N$ whose frontier components are spheres, 
the manifold $W^+$ is a connected summand of $N^+$. Since $N^+=(\hatX)^+$ is irreducible by hypothesis, any connected summand of $N^+$ which is non-simply connected must be homeomorphic to $N^+$. Hence $N^+$ is homeomorphic to $W^+$, and so $\pi_1(N)$ is isomorphic to $H_\frakU$. In view of \ref{when is a group not a group}, it follows that $\pi_1(N)$ is cyclic.
\EndProof
%U\calb_1

\Proposition\label{new get loster}
Let $\Mh$ be a closed, orientable, hyperbolic $3$-orbifold. Set $\oldOmega=(\Mh)\pl$, and assume that $\Mh$
 contains no  embedded negative turnovers. Suppose that $\vol(\Mh)\le3.44$. 
Let
$\cals$ be a non-empty $\oldOmega$-admissible system of spheres in $M:=|\oldOmega|$, and let $X$ be a connected component of $M-\cals$. Suppose that $\hatX$ is
$+$-irreducible. Then at least one of the
following alternatives holds (where $s'_{\oldOmega}(X)$, $y_{\oldOmega}(X)$, and $q_{\oldOmega}(X)$ are well defined by \ref{it's ok}):
\Alternatives
\item $0\le h(X)\le1$;
\item $2\le h(X)\le3$, and 
$s'_{\oldOmega}(X)\ge\lambda_\hatX $; 
\item $4\le h(X)\le7$, while $s'_{\oldOmega}(X)\ge \lambda_\hatX $ and $y_{\oldOmega}(X)\ge 6$; 
\item $4\le h(X)\le5$, while $s'_{\oldOmega}(X)\ge \lambda_\hatX $  and 
$q_{\oldOmega}(X)>4$; or
\item  $6\le h(X)\le7$, while $s'_{\oldOmega}(X)\ge \lambda_\hatX $  and $q_{\oldOmega}(X)>10$.
\EndAlternatives
\EndProposition

\abstractcomment{\tiny I hope I've corrected all the places where $Y$
  appeared instead of $\hatY$. }

\Proof 
Set $\oldPsi=\obd(\hatX)$, and 
$N=|\oldPsi |=\hatX$, so that $h(X)=h(N)$; and $s'_\oldOmega(X)$, $y_\oldOmega(X)$ and
$q_\oldOmega(X)$ are by definition equal to $\sigma(\oldPsi)$, $\delta(\oldPsi)$ and
$\theta(\oldPsi)$ respectively (see \ref{t-defs}). 

First consider the case in which
$s'_{\oldOmega}(X)<\lambda_\hatX $. Since  $\oldOmega$
contains no  embedded negative  turnovers, it follows from Corollary \ref{lose lose} that every component of
$\partial \hatX$ is \full. Hence by Proposition \ref{loster lemster}, $\pi_1(X)$ is
cyclic, and in particular $h(X)\le1$. Thus Alternative (i) of the conclusion holds in this case.

There remains the case in which 
$s'_{\oldOmega}(X)\ge\lambda_\hatX $. In this case, if $h(X)\le3$, either Alternative (i)
or Alternative (ii) of the conclusion holds. For the rest of the
proof, we will assume that $s'_{\oldOmega}(X)\ge\lambda_\hatX $ and that $h(X)\ge4$.

Note that since the components of $\cals$ are spheres, every
component of $\partial N$ is a sphere. Since $\Mh$ is hyperbolic and orientable
and $\cals$ is $\oldOmega$-admissible, it follows from
Lemma \ref{oops lemma}  that $\oldPsi$ is  componentwise strongly \simple\ and
boundary-irreducible. In particular $\oldPsi$ is very good, so that $\smock_0(\oldPsi)$ is defined. By
the additivity of $\smock_0$ over components (see \ref{t-defs}),
Lemma \ref{frobisher},
Corollary \ref{smockollary}, and the hypothesis of the present proposition, we have 
$\smock_0(\oldPsi)\le\smock_0(\oldOmega\cut{\obd(\cals)})\le 
\smock_0(\oldOmega)=\vol\Mh\le3.44$. Thus $\oldPsi$
satisfies the hypotheses of Proposition \ref{new get lost}. 

In view of
Conclusion (1)  of Proposition \ref{new get lost}, we have $
h(X)\le7$. In view of Conclusion (2)  of Proposition \ref{new get
  lost}, we have either $y_{\oldOmega}(X)\ge6$ or
$q_{\oldOmega}(X)>4$. Hence in the subcase $4\le h(X)\le5$, either Alternative
(iii) or Alternative (iv) of the conclusion of the present proposition
holds.

Finally, consider the subcase $6\le h(X)\le7$. In view of Conclusion (3)  of Proposition \ref{new get
  lost}, we have either $y_{\oldOmega}(X)\ge6$ or
$q_{\oldOmega}(X)>10$ in this subcase. Hence either Alternative
(iii) or Alternative (v) of the conclusion of the present proposition holds.
\EndProof

%s_ 2\oldOmega
\abstractcomment{\tiny Old comment: ``The point of 
  the statement of Theorem \ref{darts theorem}. In the latter theorem
  I talk directly about cutting an orbifold along a suborbifold.''
  Completely out of date, as far as I can see.}

\Corollary\label{new get loster corollary}
Let $\Mh$ be a closed, orientable, hyperbolic $3$-orbifold. Set $\oldOmega=(\Mh)\pl$. Suppose that $\oldOmega$
 contains no  embedded negative turnovers, and that $\volorb(\oldOmega)\le3.44$. 
Let
$\cals$ be a non-empty $\oldOmega$-admissible system of spheres in $M:=|\oldOmega|$, and let $X$ be a connected component of $M-\cals$. Suppose that $\hatX$ is
$+$-irreducible. Then $h(X)\le t_\oldOmega(X)+2/\lambda_\oldOmega$ (see \ref{it's ok}).
\EndCorollary

\Proof
We set $\lambda=\lambda_\oldOmega$. According to \ref{lambda thing}, we have $\lambda_\hatX\ge\lambda$. The invariants 
$t_\oldOmega(X)$, $s_\oldOmega(X)$, $s'_\oldOmega(X)$, $q_\oldOmega(X)$ and $y_\oldOmega(X)$ are well defined by \ref{it's ok}. 

The hypotheses of the corollary are the same as those Proposition \ref{new get loster}, and hence one of the alternatives of the latter proposition must hold.

First consider the case in which Alternative (ii) holds. Thus $h(X)\le3$ and $s'_\oldOmega(X)\ge\lambda_\hatX\ge\lambda$. According to an observation made in \ref{t-defs}  we have $s'_\oldOmega(X)\le s_\oldOmega(X)$, so that $s_\oldOmega(X)\ge\lambda$. Hence 
$s_\oldOmega(X)+2/\lambda\ge\lambda+2/\lambda$. But since $\lambda\in\{1,2\}$ we have $\lambda+2/\lambda=3$. Thus 
$s_\oldOmega(X)+2/\lambda\ge3\ge h(X)$. Since
Corollary \ref{lollapalooka} gives $s_\oldOmega(X)\le t_\oldOmega(X)$,  the conclusion follows in this case.

In the remaining cases I will show that $h(X)\le t_\oldOmega(X)+1$, which implies the conclusion since $1\le2/\lambda$.

If Alternative (i) of Proposition \ref{new get loster} holds, we have
$h(X)\le1\le  t_\oldOmega(X)+1$. If Alternative (iii) holds we have 
$h(X)\le1+y_\oldOmega(X)$, which with Corollary \ref{lollapalooka}
implies that $h(X)\le1+t_\oldOmega(X)$. If Alternative (iv) or (v) holds we have 
$h(X)\le1+q_\oldOmega(X)$, which with Corollary \ref{lollaheplooka} 
implies that $h(X)\le1+t_\oldOmega(X)$.
\EndProof
%_\hatX\vol
%oldPsi Y

\section{Dandy systems of spheres}\label{semi-dandy section}

\begin{definitionnotation} \label{size def}
Let $\oldTheta$ be a $2$-orbifold such that every component of $|\oldTheta|$ is a
$2$-sphere. 
If $C\subset |\oldTheta|$ is a simple closed curve such that
$C\cap\fraks_\oldTheta=\emptyset$, we will define the {\it size} of
$C$, denoted $\size C$, to be
$\min(\wt_\oldTheta(\scrd_1),\wt_\oldTheta(\scrd_2) )$,
where $\scrd_1$ and
  $\scrd_2$ denote the components of the complement of $C$ relative to
  the component of $|\oldTheta|$ containing $C$. 
\end{definitionnotation}
%S

\begin{definition}\label{splinter def}
Let $\oldTheta$ be a $2$-orbifold such that $S:=|\oldTheta|$ is a
$2$-sphere, and let $F\subset |\oldTheta|$  be a compact, connected surface with
$\partial F\cap\fraks_\oldTheta=\emptyset$. We will say that $F$
{\it splinters} $\oldTheta$, or respectively {\it semi-splinters} $\oldTheta$, if for each component $\Delta$ of $\overline{S-F}$ we have
$\wt_\oldTheta\Delta<\wt_\oldTheta (S) /2$, or respectively
$\wt_\oldTheta\Delta\le\wt_\oldTheta (S) /2$. (Note that these conditions hold vacuously if $F=|\oldTheta|$.)
%\redmissingref{That seems pointless
%  since all I need to do is add an $\obd$, but fixing it will require a
  %lot of subtitutions. Actually leaving it this way may have some advantages, although I can't put them into words yet. In any event, there are probably several places
%where the term ``semi-splinter'' ought to be used.}
\end{definition}

\Lemma\label{negative splinter}
Let $\oldTheta$ be a negative $2$-orbifold such that $S:=|\oldTheta|$ is a
$2$-sphere, and let $F\subset |\oldTheta|$  be a compact, connected surface with
$\partial F\cap\fraks_\oldTheta=\emptyset$. Suppose that $\obd(F)$
is taut (see \ref{praxis}) in $\oldTheta$ and that $F$
 splinters $\oldTheta$. Then $\chi(\obd(F))<0$.
\EndLemma

\Proof
Since $\obd(F)$ is taut and $\oldTheta$ is negative, we have $\chi(\obd(F))\le0$. Suppose that $\chi(\obd(F))=0$, so that $\obd(F)$ is annular. If $F$ is a weight-$0$ annulus then $|\oldTheta|-\inter F$ has two components $D_1$ and $D_2$, each of weight strictly less than $(\wt S)/2$; but $\wt S=\wt D_1 + \wt D_2$, a contradiction. There remains the possibility that $F$ is a weight-$2$ disk; in this case, $D=|\oldTheta|-\inter F$ is connected, and $\wt S-2=\wt D<(\wt S)/2$, so that $\wt S<4$ and $\wt D\le1$. The latter inequality is impossible since $\obd(F)$ is taut and $\partial F=\partial D$.
\EndProof
%C

\begin{definition}\label{what's a belt?}
Let $\Mh$ be a closed, orientable hyperbolic $3$-orbifold, and set $\Mh=(\oldOmega)\pl$.
Let $\cals$ be an
$\oldOmega$-admissible system of spheres in $M:=|\oldOmega|$.  Set
$\oldPsi=\oldOmega\cut{\obd(\cals)}$,
$N=|\oldPsi|= M\cut\cals$, and
$\tcals=\partial N$. Let  $\tS$ be a component of $\tcals$.
According to Lemma \ref{oops lemma}, $\oldPsi$ is componentwise strongly \simple\ and componentwise boundary-irreducible, so that $\oldSigma(\oldPsi)$ and $\oldPhi(\oldPsi$ are defined. By a {\it
belt} for $\tS$ (relative to $\cals$) we shall mean a component $G$ of $|\oldPhi(\oldPsi)|$, contained in $\tS$, such that 
\begin{itemize}
\item $\obd(G)$ is an annular orbifold, and
\item $\wt_\oldPsi D\le\wt_\oldPsi(\tS)/2$ for every component $D$ of $\tS-\inter
  G$.
\end{itemize}

%\Definition
%Let $\oldOmega$ be a closed, orientable hyperbolic $3$-orbifold, and
%let $\cals$ be an
%$\oldOmega$-admissible system of spheres in $M:=|\oldOmega|$.  Set
%$\oldPsi=\oldOmega\cut{\obd(\cals)}$,
%$N=|\oldPsi|= M\cut\cals$, and
%$\tcals=\partial N$. 
A {\it pseudo-belt} (relative to $\cals)$ for a component
$\tS$ of  $\tcals$ is defined to be 
a component $E$ of
  $|\oldPhi(\oldPsi)|$, contained in $\tS$, such that
\begin{itemize}
\item  $E$
  splinters $\obd(\tS)$, and
\item
there is a component $\frakR$ of $\oldSigma(\oldPsi)$ such
  that $E$ is the full intersection of $|\frakR|$ with $\tcals$.
\end{itemize}
%\end{definition} 

Note that if $G$ or $E$ is, respectively, a belt or pseudo-belt, then $\obd(G)$ or $\obd(E)$ is, respectively, a component of $\oldPhi(\oldPsi)$, and is therefore taut in $\partial\oldPsi$ and $\pi_1$-injective in $\oldPsi$ by \ref{tuesa day}.

\end{definition}

\Lemma\label{how do i recognize a belt?}
Let $\Mh$ be a closed, orientable hyperbolic $3$-orbifold, and set $\oldOmega=(\Mh)\pl$.
Let $\cals$ be an
$\oldOmega$-admissible system of spheres in $M:=|\oldOmega|$.  Set
$\oldPsi=\oldOmega\cut{\obd(\cals)}$ and
$N=|\oldPsi|= M\cut\cals$. Let $\tS$ be a component of
$\tcals:=\partial N$, and let $G$ be a component of $|\oldPhi(\oldPsi)|$, contained in $\tS$, such that $\obd(G)$ is an annular orbifold. Then the following conditions are equivalent:
\begin{enumerate}
\item $G$ is a
belt for $\tcals$;
\item $\wt_\oldPsi D=\wt_\oldPsi(\tS)/2$ for every component $D$ of $\tS-\inter
  G$;
\item $\size C=\wt_\oldPsi(\tS)/2$ for some component $C$ of $\partial
  G$;
\item $\size C=\wt_\oldPsi(\tS)/2$ for every component $C$ of $\partial
  G$.
\end{enumerate}

In particular, if there is a belt for $\tS$ then $\wt_\oldPsi\tS$ is even.
\EndLemma

\Proof
It is trivial that (2) implies (1). To prove the converse, suppose that $G$ is a belt. Since $\obd(G)$ is annular (and orientable), $G$ is either an annulus of weight $0$ or a disk of weight $2$. If $G$ is an annulus of weight $0$, then since $\tS$ is a sphere, $\tS-\inter G$ has exactly two components, say $D_1$ and $D_2$. For $i=1,2$ we have $\wt D_i\le\wt(\tS)/2$ by the definition of a belt, but since $\wt G=0$ we have $\wt D_1+\wt D_2=\wt(\tS)$. Hence $\wt D_1=\wt D_2=\wt(\tS)/2$, which proves (2) in this case. If $G$ is  a disk of weight $2$, then $D:=\tS-\inter G$ is a disk, and $\wt D=\wt(\tS)-2$. But by the definition of a belt, we have $\wt D\le\wt(\tS)/2$; it follows that $\wt(\tS)\le4$ and that $\wt D\le2$.
%If $\wt(\tS)<4$, then $\wt(D)<2$, so that $\obd(D)$ is a discal %orbifold. But this is impossible, because 
On the other hand,
$\obd(\partial D)=\obd(\partial G)\subset\partial\oldPhi(\oldPsi)$ is
$\pi_1$-injective in $\oldPhi$, and hence in $\tS$, by \ref{tuesa day}; this implies that the orientable $2$-orbifold $D$ is not  discal,  and therefore that $\wt D\ge2$. Hence we must have $\wt D=2$ and therefore $\wt(\tS)=4$, so that $\wt D=\wt(\tS)/2$, and (2) holds in this case as well. Thus (1) and (2) are equivalent.

The implication (4) $\Rightarrow$ (3) holds because an annular orbifold has a non-empty boundary. To prove the converse in the case where $G$ is a weight-$0$ annulus, we need only note that the components of $\partial G$ are parallel in $\tS\setminus\fraks_\oldPsi$, and hence have the same size. In the case where $G$ is a weight-$2$ disk, $\partial G$ has only one component, so that the implication (3) $\Rightarrow$ (4) is immediate. Thus (3) and (4) are equivalent.

It remains to prove that (2) is equivalent to (4). For this purpose, consider any component $C$ of $\partial G$, and let $D$ denote the component of $\tS-\inter G$ bounded by $C$. We shall prove:
\Claim\label{sharrop}
We have  $\size C=\wt(\tS)/2$ if and only if $\wt D=\wt(\tS)/2$
\EndClaim
The equivalence of (2) and (4) will follow from \ref{sharrop}, since every component of $\partial G$ bounds a component of $\tS-\inter G$, and each component of $\tS-\inter G$ is bounded by a component of $\partial G$. 

To prove \ref{sharrop}, set $D'=\tS-\inter G$. The definition of size gives $\size C=\min(\wt D,\wt D')=\min(\wt D,\wt\tS-\wt D)$. But for any integers $m$ and $n$ with $0\le m\le n$, we have $\min(m,n-m)=n/2$ if and only if $m=n/2$.
\EndProof

\abstractcomment{\tiny I had this:
``Note that in the setting of Definition xxx we have $\fraks_\oldTheta=\tS\cap\fraks_\oldPsi$. Hence
the condition that $G$ be a balanced annulus may be rewritten as
  $\card (\scrd_1\cap\fraks_\oldPsi)=\card (\scrd_2\cap\fraks_\oldPsi)$.'' I think this will come out in the wash.}

\Lemma\label{it's-a unique}
Let $\Mh$ be a closed, orientable hyperbolic $3$-orbifold, and set $\oldOmega=(\Mh)\pl$. Let $\cals$ be an
$\oldOmega$-admissible system of spheres in $M:=|\oldOmega|$. 
%If $S$ 
Set $\oldPsi=\oldOmega\cut{\obd(\cals)}$, 
$N=M\cut\cals
%=|\oldPsi|
$,
and $\tcals=\partial N$. If $\tS$ is a component of $\tcals$, 
there
can be at most one belt for $\tS$, and $\tS$ can contain at most one
component of $|\oldPhi(\oldPsi)|$ which splinters $\tS$. Furthermore, if $\tS$  has a belt, then no component of $|\oldPhi(\oldPsi)|$  splinters $\tS$. 

In particular,
there
can be at most 
%one belt for $\tS$, and there can be at most one
pseudo-belt for $\tS$. Furthermore, $\tS$ cannot have both a belt and
a pseudo-belt.
\EndLemma

\Proof
The final sentence of the conclusion follows from the preceding ones, because the definition of pseudo-belt implies that any pseudo-belt splinters $\obd(\tS)$. It therefore suffices to prove the  assertions before the final sentence.

If some component $F$ of  $|\oldPhi(\oldPsi)|$ is a belt, then by definition  $\chi(\obd(F))=0$; since $\obd(F)$ is taut by \ref{tuesa day}, it follows from Lemma \ref{negative splinter} that $F$ cannot splinter $\obd(\tS)$. It remains 
%, the  definition of splintering gives
%$\wt( D)<\wt(\tS) /2$ for every component $D$ of $\tS-\inter
%F$, while the assumption that $F$ is a belt, together with Lemma \ref{how do i recognize a belt?}, implies that $\wt( D)=\wt(\tS) /2$ for every component $D$ of $\tS-\inter
%F$. It follows that $\tS-\inter
%F$ has no components, i.e. that $F=\tS$. But the definition of a belt
%implies that $\obd(F)$ is an annular orbifold, so that $F$ is a disk
%or an annulus, not a $2$-sphere. This contradiction shows that no
 %belt can splinter $S$. Hence in order to prove the lemma, it
%suffices
 to show that  $|\oldPhi(\oldPsi)|$ can not have two distinct components 
contained in $\tS$, each of which either is a belt or splinters $\tS$.
Let us assume that   $F_0$ and $F_1$ are two such
  distinct components.

In particular $F_0$ and $F_1$ are disjoint, compact, connected subsurfaces of
  the sphere $\tS$. Hence there is a component $V$ of
% Furthermore, we have $F_0,F_1\subset \tS$. Since $F_0$ and $F_1$ are
% disjoint annuli on the $2$-sphere $\tS$, the components of 
$\tS-\inter(F_0\cup F_1)$ which is an annulus sharing one
boundary component with $F_0$ and one with  $ F_1$. 
%Set
%$C_i=F_1\cap V$ 
For $i=0,1$, let $D_i$ denote the  component of $\tS-\inter V$
containing $F_i$.

The disk $V\cup D_0$ is a component of $\tS-\inter F_1$. If $F_1$ 
splinters $\obd(\tS)$, we have
$\wt(V\cup D_0)<\wt(\tS) /2$. If $F_1$ is a belt, then by Lemma
\ref{how do i recognize a belt?} we have
$\wt(V\cup D_0)=\wt (\tS) /2$. Since these observations remain
true if the indices $0$ and $1$ are interchanged, we deduce:
\Claim\label{ought to buy now}
For $i=0,1$ we have $\wt(V\cup D_i)\le\wt (\tS) /2$, with equality if
and only if $F_{1-i}$ is a belt.
\EndClaim

%If either $F_0$ or $F_1$ is a pseudo-belt, 
It follows from \ref{ought
  to buy now} that
\Equation\label{scuzz}
\wt(V\cup D_0)+\wt(V\cup D_1)\le\wt (\tS),
\EndEquation
with equality only if $F_0$ and $F_1$ are both belts.
But
since $D_0$, $V$, and $D_1$ have mutually disjoint interiors, and their
union is $\tS$, we have $\wt(\tS)=\wt(D_0)+\wt(V)+\wt(D_1)\le
\wt(D_0)+2\wt(V)+\wt(D_1)=\wt(V\cup D_0)+\wt(V\cup D_1)$. Hence
\Equation\label{fuzz}
\wt(\tS)\le
\wt(V\cup D_0)+\wt(V\cup D_1),
\EndEquation
with equality only if $\wt V=0$. Taken together, (\ref{scuzz}) and
(\ref{fuzz}) give a contradiction unless both are equalities,
i.e. unless $\wt V=0$ and both the $F_i$ are belts. But in this case, $V$ is an annulus component of
$\partial|\oldPsi|-\inter|\oldPhi(\oldPsi)|$ having weight $0$, so that $\obd(V)$ is an annular orbifold; and each component of
$\partial V$ is contained in a component of $|\oldPhi(\oldPsi)|$ which
is an annular orbifold. 
This contradicts the uniqueness assertion of Proposition \ref{lady edith}.
\EndProof
%oldTheta N

\begin{remarksnotationdefinitions}\label{stop beeping}
Let $\Mh$ be a closed, orientable hyperbolic $3$-orbifold, and set $\oldOmega=(\Mh)\pl$. Let $\cals$ be an
$\oldOmega$-admissible system of spheres in $M:=|\oldOmega|$. Set
$\oldPsi=\oldOmega
\cut 
{\obd(\cals)}
$, and set 
$N=M\cut\cals
=|\oldPsi|
$, $\tcals=\partial N$, and $\tau=\tau_\cals$.

We will say that a component 
$\tS$  of $\tcals$ is {\it belted relative to $\cals$},  or {\it pseudo-belted
 relative to $\cals$}, if there
exists a belt, or, respectively, a pseudo-belt for $\tS$. (We will say
simply that $\tS$ is belted or pseudo-belted when it is clear which
admissible system is involved.) According to Lemma \ref{it's-a
  unique}, a component of $\tcals$ cannot be both belted and pseudo-belted. If $\tS$ is belted, its belt, which is
unique by Lemma \ref{it's-a
  unique} (once $\oldPhi(\oldPsi)$ has been fixed within its isotopy class), will be
denoted $G_\tS^\cals$, or simply by $G_\tS$ when it is clear from
the context which system of spheres is involved. If $\tS$ is pseudo-belted, its pseudo-belt, which is
unique by Lemma \ref{it's-a
  unique} (once $\oldPhi(\oldPsi)$ has been fixed within its isotopy class), will be
denoted $E_\tS^\cals$, or simply by $E_\tS$ when it is clear which system of spheres is involved. If  $\oldPhi(\oldPsi)$ has not been fixed within its isotopy class, then $\obd(G_\tS)$ and $\obd(E_\tS)$ are well defined up to isotopy in $\tS$, in the belted and pseudo-belted cases respectively.

%\frakU

If $\tS$ is a pseudo-belted component of $\tcals$, 
then since $E_\tS$ is by definition a component of
$|\oldPhi(\oldPsi)|$, there is a unique component $\frakR$ of $\oldSigma(\oldPsi)$ such
that $E_\tS\subset|\frakR|$. We will set $\frakR_\tS=\frakR$. The definition of a
pseudo-belt implies that
$|\frakR_\tS|\cap\tcals=E_\tS$. 
The definition also implies that
$E_\tS$ splinters $\obd(\tS)$; since $\obd(\tS)$ is negative by \ref{doublemint}, it then follows from Lemma \ref{negative splinter} that $\chi(\obd(F))<0$, i.e. $E_\tS$ is a component of
$\oldPhi^-(\oldPsi)$. 
%But the definition of a pseudo-belt also implies that the connected orbifold $\obd(E_\tS)$ is the full intersection of $\frakR_\tS$ with $\partial\oldPsi$ (so that by \ref{S-pair def} the connected \pagelike\ \Ssuborbifold\ $\oldLambda$ is \twisted). 
By \ref{what's iota?}, 
%$\Sigma^-(\oldPsi)$, 
$\frakR_\tS$ is a component of $\oldSigma^-(\oldPsi)$. In particular  
%the
%connected \Ssuborbifold\ 
$\frakR_\tS$ must be saturated in a fibration $q$ of $\oldSigma^-(\oldPsi)$ of the kind described in \ref{what's iota?}. It therefore follows from the discussion in \ref{what's iota?} that
%compatible with $\oldPhi^-(\oldPsi)$ \redproofreadingnote{Which is compatible with which? And should this fibration, which is defined up to isotopy, have a name? And should the term ``isotopy'' be defined for fibrations? I skirted the issue in stating Prop. \ref{unique fibration}}, so that by \ref{what's iota?}, 
% Since the component $\frakR_\tS$ of $\oldSigma^-(\oldPsi)$ \redmissingref{ justify the next clause as needed; it's tied up with $\iota$ issues}, which implies that saturated
$\obd(E_\tS)$ is invariant under $\iota_\oldPsi$. (The fibration $q$ is non-trivial since $\obd(E_\tS)=\partialh\frakR_\tS$ is connected.) The restriction of
$\iota_\oldPsi$ to $\obd(E_\tS)$ is an involution of $\obd(E_\tS)$ which will be
denoted by  $\epsilon_\tS^\cals$, or simply by $\epsilon_\tS$ when it
is clear which admissible system is involved. According to  \ref{what's iota?},  the involution $\iota_\oldPsi$ of $\oldPhi^-(\oldPsi)$ is 
well defined up
  to strong equivalence in $\partial\oldPsi=\obd(\tcals)$ 
once $\oldPhi(\oldPsi)$ has been fixed within its isotopy class. Hence 
the involution $\epsilon_\tS$ of
$\obd(E_\tS)$ is 
well defined up
  to strong equivalence in $\obd(\tS)$, even if $\oldPhi(\oldPsi)$ has not been fixed within its isotopy class.

If a component 
$\tS$  of $\tcals$ is belted, then it follows from Lemma \ref{how do i recognize a belt?} that  $\wt_\oldOmega\tS$ is
even. Furthermore if the belted component $\tS$ has weight $2d$, and if $G_\tS$
is an annulus, then it follows from Lemma \ref{how do i recognize a belt?} that a core curve for $G_\tS$ has size $d$.

If $\tS$ is a belted component 
of $\tcals$, then since $G_\tS$ is by definition a component of $|\oldPhi(\oldPsi)|$, it is contained in a unique component of $|\oldSigma(\oldPsi)|$, which
% containing
%$G_\tS$ 
will be denoted $L_\tS$.  Since the definition of a belt implies that
$\obd(G_\tS)$ is an annular $2$-orbifold, it follows from Proposition
\ref{what does it say?} that $\obd(L_\tS)$ is a \torifold, and that
every component of $\obd(L_\tS\cap\partial|\oldPsi|)$ is an annular
$2$-orbifold. By \ref{tuesa day}, $\obd(L_\tS\cap\partial|\oldPsi|)$ is $\pi_1$-injective in $\oldPsi$. 
%Proposition \ref{what does it say?},
  %and it could be put back. Added 4/7/18: I'm confused about this, because I think the general fact is mentioned in \ref{tuesa day}. Or am I misunderstanding what the issue is? This seems like a very long comment on a complete triviality. ??

\end{remarksnotationdefinitions}
%\torifold

\Lemma\label{need that too}
Let $\Mh$ be a closed, orientable hyperbolic $3$-orbifold, and set $\oldOmega=(\Mh)\pl$. Let $\cals$ be an
$\oldOmega$-admissible system of spheres in $M:=|\oldOmega|$. Set
$N=M\cut\cals$, , and let $\tS$ be a belted component of 
$\tcals:=\partial N$. Then for every component $X$ of $\tS-\inter G_\tS$ we have $\chi(\obd(X))<0$.
\EndLemma

\Proof
Set $\oldTheta=\obd(\tS)$ and $\oldPsi=\oldOmega\cut\oldTheta$, so that $N=|\oldPsi|$.
Since $\tS$ is a sphere and $G_\tS$ is connected, the component  $X$
of $\tS-\inter G_\tS$ is a disk. Since the belt $G_\tS$ is by
definition a component of $|\oldPhi(\oldPsi)|$, the boundary of
$\obd(G_\tS)$ is $\pi_1$-injective in $\oldTheta$ by \ref{tuesa day}, and hence $\chi(\obd(X))\le0$. Let us assume that $\chi(\obd(X))=0$, i.e. that $\obd(X)$ is an annular orbifold. Let $V\subset X$ denote the component of
$\partial|\oldPsi|-\inter|\oldPhi(\oldPsi)|$ containing $\partial X$. By \ref{tuesa day},  $\partial \obd(V)$ is $\pi_1$-injective in $\oldPsi$ and hence in $\obd(X)$; since $\obd(X)$ is  annular, it then follows from \ref{cobound} that $\obd(V)$ is annular. Hence by Proposition \ref{lady edith},  $V$ is a weight-$0$ annulus. Let $Y\subset X$ denote the component of $|\oldPhi(\oldPsi)|$ whose boundary contains the simple closed curve $(\partial V)-(\partial D)$. Since $\obd(X)$ is annular, and $\Fr_{\obd(X)}\obd(Y)$ is $\pi_1$-injective, it follows from \ref{cobound} 
% fAgain applying \redrealmissingref{same two cross-refs}, we deduce 
that $\obd(Y)$ is annular. But now $G_\tS$ and $Y$ are the components of $|\oldPhi(\oldPsi)|$ sharing boundary components with $V$, and since $\obd(G_\tS)$ and $\obd(Y)$ are both annular, we have a contradiction to the uniqueness assertion of Proposition \ref{lady edith}. 
\EndProof

\begin{definitionremark}\label{strict def}
Let $\Mh$ be a closed, orientable hyperbolic $3$-orbifold, and set $\oldOmega=(\Mh)\pl$. Let $\cals$ be an
$\oldOmega$-admissible system of spheres in $M:=|\oldOmega|$. Set
$\oldPsi=\oldOmega\cut{\obd(\cals)}$ and $\tS=\partial M\cut \cals=\partial|\oldPhi|$.
A component 
$\tS$  of $\tcals$ will be said to be {\it \central} (relative to the system $\cals$) if no component of $\tS\cap|\kish(\oldPsi)|$ semi-splinters $\obd(\tS)$.

A component of $\cals$ will be termed {\it \doublecentral} if its sides
  are both \central\ in $\cals$.
% $\tS$ is \full\ 
%and contains no bad component (see \ref{stop beeping})  of $|\oldPhi(\oldPsi)|$.

Note that if $\tS$ is a \full\ component of $\tS$ then in particular
$\tS$ is \central, since by definition we have 
$\tS\cap\kish(\oldPsi)=\emptyset$ in this case. 
\end{definitionremark}

\Lemma\label{corpulent porpoises}
Let $\Mh$ be a closed, orientable hyperbolic $3$-orbifold, and set $\oldOmega=(\Mh)\pl$. Let $\cals$ be an
$\oldOmega$-admissible system of spheres in $M:=|\oldOmega|$. Set
$\oldPsi=\oldOmega\cut{\obd(\cals_)}$ and $\tS=\partial M\cut \cals=\partial|\oldPhi|$.
Let 
$\tS$ be a component of $\tcals$ which is both \central\ and
belted. Let $V$ be a component of $\tS\setminus\inter|\oldPhi(\oldPsi)|$ such
that $\partial V\cap\partial G_\tS\ne\emptyset$. Then
$\obd(V)$ is an annular orbifold.
\EndLemma

\Proof
Set $G=G_\tS$. Since $\obd(G)$ is an annular orbifold by the
definition of a belt (and is orientable), $G$ is either an annulus of weight $0$ or a disk
of weight $2$. 

In the case where $G$ is a disk of weight $2$, then $\wt(\tS-\inter G)
=\wt(\tS)-2$; but it follows from Lemma
\ref{how do i recognize a belt?} that $\wt(\tS-\inter
G)=\wt(\tS)/2$. Hence $\wt(\tS)=4$ and $\wt(\tS-\inter
G)=2$, so that $\obd(\tS-\inter
G)$ is an annular orbifold. 
Now since $V\cap G\ne\emptyset$, we have $\partial G\subset V$. On the other hand, since $V$ be a component of $\tS\setminus\inter|\oldPhi(\oldPsi)|$, there is a component $\kappa$ of $\tS\cap\kish(\oldPsi)$ with $V\subset|\kappa|$. We now have $\partial G\subset V\subset|\kappa|$, so that any component of $\overline{\tS-\kappa}$ must be contained in either $G$ or
$\tS-\inter
G$. Hence every component of $\overline{\tS-\kappa}$ has weight at most $2=\wt(\tS)/2$. By definition this means that $\kappa$ semi-splinters $\tS$, a contradiction to the hypothesis that $\tS$ is \central. Thus this case cannot occur.

Now consider the case in which $G$ is a weight-$0$ annulus. Since $G$ and $V$  are components of $|\oldPhi(\oldPsi)|$
and $\tS\setminus|\oldPhi(\oldPsi)|$ respectively, their intersection is a union of
common boundary components, which by hypothesis is non-empty. Since $\tS$ is a
sphere, $G\cap V$ is a single simple closed curve $C_0$. Let $C_1$
denote the other boundary component of $G$. For $i=0,1$, let $D_i$
denote the component of $\tS-\inter G$ whose boundary is $C_i$. The
definition of a belt gives $\wt D_i\le\wt(\tS)/2$ for $i=0,1$.

According to \ref{tuesa day} we have $\chi(\obd(V))\le0$. We
need to show that $\chi(\obd(V))=0$. If
we assume that $\chi(\obd(V))<0$, then
the component
$\obd(V)$ of $\obd(\tS)\setminus\inter\oldPhi(\oldPsi)|$ 
is not contained in $\frakH(\oldPsi)$; hence 
for some
component $\kappa$ of $\tS\cap|\kish(\oldPsi)|$ we have
$\obd(V)\subset\kappa$, i.e.
$V\subset|\kappa|$.  Since $\tS$  is a \central\ component 
of $\tcals$, the component $\kappa$ of $\tS\cap|\kish(\oldPsi)|$ cannot semi-splinter $\obd(\tS)$;  there therefore exists a component
$D$ of
$\tS-\inter|\kappa|$ such that $\wt_\oldPsi D>\wt_\oldPsi(\tS)/2$. If $E$
denotes the component of $\tS-\inter V$ containing $D$, we have
$\wt_\oldPsi E\ge\wt_\oldPsi D>\wt_\oldPsi(\tS)/2$.

Note that 
$D_1\cup A$ is one component of $\tS-\inter V$, and that
every other component of $\tS-\inter V$ is contained in $D_0$. If
$E=D_1\cup A$, we have $\wt E=\wt D_1+\wt A=\wt D_1\le\wt(\tS)/2$,
a contradiction. If $E\subset D_0$, we have $\wt E\le\wt
D_0\le\wt(\tS)/2$, which again gives a contradiction.
\EndProof

\Definition\label{bad def}
Let $\Mh$ be a  closed, orientable hyperbolic $3$-orbifold, and set $\oldOmega=(\Mh)\pl$. Let $\cals$ be an
$\oldOmega$-admissible system of spheres in $|\oldOmega|$. Set $\oldPsi=\oldOmega\cut{\obd(\cals)}$.
Let  $F$ be a component of $|\oldPhi(\oldPsi)|$, 
and let $\tS$ denote the component of $\tcals$ containing $F$. We will
say that $F$ is {\it \bad} (relative to $\cals$) if (i) $F$ splinters $\obd(\tS)$, and (ii) the component of $|\oldSigma(\oldPsi)|$ containing $F$
meets $\tau(\tS)$.
\EndDefinition

\Definition \label{semidandy def}Let $\Mh$ be a closed, orientable hyperbolic $3$-orbifold, and set $\oldOmega=(\Mh)\pl$. A {\it$\oldOmega$-\dandy} system of spheres in $|\oldOmega|$ is an
$\oldOmega$-admissible system of spheres $\cals$ in $|\oldOmega|$ such that the following conditions hold when we set 
$\oldTheta=\obd(\cals)$, $\oldPsi=\oldOmega\cut\oldTheta$, 
$N=M\cut\cals
=|\oldPsi|
$, $\tcals=\partial N$, and ${\toldTheta}=\obd(\tcals)$:
\begin{itemize}
\item for any weight-$0$
annulus $A$, properly embedded in $N$, the two components of
  $\partial A$ have the same size (see \ref{size def}) in ${\toldTheta}$; 
\item 
for every component $F$
% of $\tcals$ such that $\tS$ contains no
%bad component 
of $|\oldPhi(\oldPsi)|$ which splinters $\obd(\tS)$, where $\tS$ denotes  the component of $\tcals$ containing
$F$, either
%. Let
%$E$ be a component of $|\oldPhi(\oldPsi)|$ contained in $\tS$. Suppose that
$F$ is a pseudo-belt for $\tS$, or $F$ is \bad; and
\item for
every properly embedded disk $D$ in $N$, transverse to $\fraks_\oldPsi$, such that $\obd(D)$ is
annular and $\pi_1$-injective in $\oldPsi$,
%The stronger
%condition seems to be needed for the proof of  Proposition\ref{dandy
  %subsystem}. It should not affect other refs. to the def., but I
%should check this. Of course I need a slight rewording of the proof of Proposition \ref{semidandies exist}, and
%I'll have to say a word when \dandy ness is applied later in the paper
there is a disk $E\subset\tcals$ such that $\obd(E)$ is an annular orbifold and $\partial E=\partial D$.
\end{itemize}

We will say ``\dandy'' in place of ``$\oldOmega$-\dandy'' when
 there is no ambiguity about which orbifold $\oldOmega$ is involved.
\EndDefinition

\Remark
In navigating the logic of this monograph, it may be useful to note that
the second condition in the definition given above is used in a
crucial way in the proof of Lemma \ref{and
  the foos just keep on fooing}, which is in turn quoted in the
proof of Propositions \ref{dandy subsystem} and  \ref{three or better}. \EndRemark

\abstractcomment{\tiny I had written ``The second condition in this definition originally said that 
every \strict\ component $\tS$ of $\tcals$ is either  belted or
pseudo-belted. Perhaps there should be a lemma to the effect that
\dandy\ systems have that property, but it should be stronger: it
should assert not only that every \strict\ component $\tS$ of
$\tcals$ is either  belted or pseudo-belted, but that components that
``look \strict\ around the middle'' are either  belted or
pseudo-belted.'' This seems to have been dealt with in the statements
of Lemmas \ref{and the foos just keep on fooing} and the now-superseded ``also because
strict,'' and also the comment after Lemma \ref{and the foos just keep on fooing}.}

\tinymissingref{\tiny See pseudo-belt.tex for the old def. of ``\dandy'' and a
  lot of stuff about how to prove that \dandy\ systems in the old
  sense exist. Also for all the stuff about \eclipsing\ spheres.}

%oldTheta

The rest of this section is devoted to exploring properties of \dandy\ systems. The question of the existence of (useful) \dandy\ systems will be addressed in the next section. Our first result on \dandy\ systems, Lemma \ref{and the foos just keep on fooing}, will depend on the following two preliminary results.

\Lemma\label{decomp}Let $\oldTheta$ be a $2$-orbifold such that $S:=|\oldTheta|$ is a $2$-sphere,  and let $\frakC$ be a compact (but possibly disconnected) $2$-suborbifold of $S$. Then either 
\begin{itemize}
\item there is a component  of $|\frakC|$ which splinters $\oldTheta$,or
\item
there is a component of $|\overline{\oldTheta-\frakC}|$ which semi-splinters $\oldTheta$.
\end{itemize}
\EndLemma

\Proof
Set $X=\fraks_\oldTheta$, a finite subset of $S$.
Since $\oldTheta$ must be orientable, $F:=|\frakC|$ is a
% compact $2$-suborbifold $\frakC$ must have the form $\obd(F)$, where $T$ is a 
compact $2$-submanifold of $S$ with $\partial T\cap X=\emptyset$.

If $F=S$ then $S=|\frakC|$ splinters $\oldTheta$, and if $F=\emptyset$ then $S=|\oldTheta-\frakC|$ splinters (and hence semi-splinters) $\oldTheta$. We may therefore assume that $\partial F\ne\emptyset$.

Let $\cald$ denote the set of all disks $D\subset S$ such that (i) $\partial D\subset\partial F$ and (ii) $\card(D\cap X)\le\card(X)/2$. Note that if $D\subset S$ is any disk bounded by a component of $\partial F$, then either $D\in\cald$ or $S-\inter D\in\cald$. Hence $\cald\ne\emptyset$. Since $\cald$ is obviously finite, there is an element $D_0$ of $\cald$ which is maximal with respect to inclusion. Since $D_0\in\cald$ we have $\partial D_0\subset\partial F$, and hence there is a component $V$ of either $F$ or $\overline{S-F}$ such that $V\cap D_0=\partial D_0$. 
We claim:
\Claim\label{moscoot}
For every component $L$ of $S-V$, we have $\card(X\cap L)\le\card(X)/2$.
\EndClaim
To prove \ref{moscoot}, first note that the required inequality holds when $L$ is the component $\inter D_0$ of $S-V$, because $D_0\in\cald$. Now suppose that $L$ is a component of $S-V$ distinct from $\inter D_0$. 
Since $S$ is a $2$-sphere and the subsurface $V$ is connected, the component $L$ of $S-V$ is the interior of a disk in $S$, and hence $D_1:=S-L$ is a closed disk. Since $L\ne\inter D_0$, the disk $D_1$ properly contains $D_0$. By the maximality of $D_0$ we have $D_1\notin\cald$; since $\partial D_1=\partial \overline{L}\subset\partial F$, it follows that $\card(D_1\cap X)>\card(X)/2$. Hence $\card(X\cap L)=\card(X)-\card(X\cap D_1)<\card(X)/2$. This proves \ref{moscoot}.

By construction, $V$ is either a component of $F$ or a component of $\overline{S-F}$. If  $V$ is a component of  $\overline{S-F}$, then $\obd(V)$ is  a component of $\overline{\oldTheta-\frakC}$, and \ref{moscoot} says 
by definition that this component semi-splinters $\oldTheta$; thus the second alternative of the conclusion of the lemma holds.
%follows from \ref{moscoot}.
% upon setting $T=V$. 
Now suppose that $V$ is a component of $F$. 
According to \ref{moscoot}, either we have
(a) $\card(X\cap L)<\card(X)/2$
for every component $L$ of $S-V$, or (b)
$\card(X\cap L_0)=\card(X)/2$
for some component $L_0$ of $S-V$. If (a) holds, then 
by definition $\obd(V)$ is  a component of $\overline{\frakC}$ which splinters $\oldTheta$, so that the first alternative of the conclusion of the lemma holds.
If (b) holds,
%the first alternative of the conclusion of the lemma holds with $T=V$, or there is a component $L_0$ of $S-V$ such that $\card(X\cap L_0)=\card(X)/2$. In the latter case, 
then since $S$ is a $2$-sphere and the subsurface $V$ is connected, the component $L_0$ of $S-V$ is the interior of a disk $E$ in $S$. Set $c=\partial E$. Since $V$ is a component of $F$, and the disk $E$ is a component of $\overline{S-V}$ bounded by $c\subset\partial V$, there is a component $T$ of $\overline{S-F}$ such that $c\subset T\subset E$. 
Now $\frakE
:=\obd(T)$ is a component of $\overline{\oldTheta-\frakC}$; we claim that $\frakE$ semi-splinters $\oldTheta$, which will imply the second alternative of the conclusion and thus complete the proof. Thus it suffices 
%In order to establish the second alternative of the conclusion in this case, it now suffices 
to prove that for every component $K$ of $S-T$ we have $\card(X\cap K)\le\card(X)/2$. For the component $S-E$ of $S-T$ we have $\card(X\cap(S-E))=\card X-\card(X\cap L_0)=\card(X)/2$. If $K$ is a component of $S-T$ distinct from $S-E$, we have $K\subset E$, and hence $\card(X\cap K)\le\card(X\cap L_0)=\card(X)/2$.
%X
\EndProof
%\frakE(a) (b)\frakC\oldGamma\frakE V

\Lemma\label{can't have it both ways}
Let $\oldTheta$ be a $2$-orbifold such that $|\oldTheta|$ is a
$2$-sphere, and let $F$ and $Q$   be compact, connected subsurfaces of $|\oldTheta|$ whose boundaries are disjoint from
$\fraks_\oldTheta$. Suppose that $F$ splinters $\oldTheta$ and that $Q$ semi-splinters $\oldTheta$. Then $F\cap Q\ne\emptyset$.
\EndLemma

\Proof
Assume that $F\cap Q=\emptyset$. Then $F$ and $Q$ are disjoint,
compact, connected subsurfaces of the sphere $S:=|\oldTheta|$, and hence there are
disjoint disks $D,E\subset S$ such that $F\subset D$, $\partial
D\subset\partial F$, $Q\subset E$, and $\partial E\subset\partial
F$. Let us set $n=\wt S=\card\fraks_\oldTheta$. Since $S- D$ is a
component of $S- F$, and since $F$ splinters $\oldTheta$, we have
$\wt(S-D)<n/2$; hence $\wt D>n/2$. Since $S- E$ is a component of $S-
Q$, and since $Q$ semi-splinters $\oldTheta$, we have $\wt(S-E)\le
n/2$; thus $\wt E\ge n/2$. Hence we have $\wt D+\wt E>n=\wt S$, which is impossible since $D$ and $E$ are disjoint.
\EndProof

\Lemma\label{and the foos just keep on fooing}
Let $\Mh$ be a closed, orientable hyperbolic $3$-orbifold, and set $\oldOmega=(\Mh)\pl$. Set
$M=|\oldOmega|$, and let $\cals$ be an $\oldOmega$-\dandy\ system of spheres
in $M$. Set $N=M\cut\cals$ and
$\tcals=\partial N$. Then for every \central\ component $\tS$ of
$\tcals$,  either $\tS$ is
belted or pseudo-belted, or 
$\tS$ contains an \bad\ component (see Definition \ref{bad def}) of $|\oldPhi(\oldPsi)|$. \EndLemma

\Proof
Set 
$\oldPsi=\oldOmega\cut{\obd(\cals)}$,
so that $N=|\oldPsi|=M\cut\cals$, and
set $\oldPhi=\oldPhi(\oldPsi)$. 
%Set
%$\rho=\rho_\cals$ and $S=\rho(\tS)$. 

Let $\tS$ be any \central\ component of $\tcals$. Set ${\toldTheta}=\obd(\tS)$.
We apply Lemma \ref{decomp}, letting 
$\toldTheta$ play the role of $\oldTheta$ (so that $\tS$ plays the role of $S$), and letting $\oldPhi\cap\obd(\tS)$ play the role of $\frakC$. This shows that either
\begin{enumerate}[(i)]
\item there is a component $\frakE\subset \obd(\tS)$ of $\oldPhi$ such that $|\frakE|$ splinters $\obd(\tS)$, or
\item 
there is a component $\frakE$ of $\obd(\tS)\setminus\inter\oldPhi$ such that $|\frakE|$ semi-splinters $\obd(\tS)$.
\end{enumerate}

If (i) holds, 
it follows from 
the second condition in the definition of a \dandy\ system (see \ref{semidandy def}) that $|\frakE|$ is either a pseudo-belt for $\tS$ or an \bad\ component of $|\oldPhi|$. 
Hence the conclusion of the lemma holds in this subcase.

The rest of the proof will be devoted to the case in which (ii) holds. We fix a component $\frakE$ of $\obd(\tS)\setminus\inter\oldPhi$ such that $V:=|\frakE|$ semi-splinters $\obd(\tS)$.
We claim:
\Claim\label{it's annular} The orbifold $\frakE$ is annular.
\EndClaim

To prove \ref{it's annular}, 
assume that $\frakE$ is not annular. Since $\frakE$ is a component of $\partial\oldPsi-\inter\oldPhi(\oldPsi)$, it  follows from Corollary \ref{abu gnu} and the non-annularity of $\frakE$ that $\frakE$ is a
 component 
of 
$\obd(\tS)\cap\kish(\oldPsi)$. Thus $V$ is a component of $\tS\cap|\kish(\oldPsi)|$ which semi-splinters $\obd(\tS)$. But in view of Definition \ref{strict def}, this contradicts the hypothesis that $\tS$ is \central. Thus \ref{it's annular} is proved.

It follows from \ref{it's annular} and  Proposition \ref{lady edith} that
$V$ is an weight-$0$ annulus, and that one component $C_0$ of
$\partial V$ is contained in a component $G$ of $|\oldPhi(\oldPsi)|$ such that $\obd(G)$ is an (orientable) annular orbifold. In particular $G$ is a disk or an annulus. Let $C_1$ denote the component of $\partial V$ distinct from $C_0$, and for $i=0,1$, let $L_i$  denote the component of $\tS- V$ bounded by $C_i$.

We claim that $G$ is a belt for $\tS$; this will imply that $\tS$ is belted, and complete the proof. According to the definition, showing that $G$ is a belt is tantamount to showing that
$\wt_\oldPsi D\le\wt_\oldPsi(\tS)/2$ for every component $D$ of $\tS-\inter
  G$. Note that $D_1:=\overline{L_1}\cup V$ is a component  of $\tS-\inter
  G$, and that if $\tS-\inter
  G$ has another component (i.e. if $G$ is an annulus) then this second component is contained in $L_0$. We have $\wt D_1=\wt L_1+\wt V=\wt L_1$. But since $L_1$ is a component of   
$\tS- V$, and $V$ semi-splinters $\obd(\tS)$, we have
$\wt L_1 \le \wt_\oldPsi(\tS)/2$. This gives the required inequality $\wt D_1\le \wt_\oldPsi(\tS)/2$ for the component $D_1$ of  $\tS-\inter
  G$.  If $D$ is a component  of $\overline{\tS-
  G}$ contained in $L_0$, we have
$\wt D\le\wt L_0\le\wt_\oldPsi(\tS)/2$, where the last inequality holds because $L_0$ is a component of   
$\tS- V$ and $V$ semi-splinters $\obd(\tS)$.
Thus $G$ is indeed a belt for $\tS$.
\EndProof
%T splinter|\oldGamma|\frakC\frakE\oldGamma$F$V

\Lemma\label{lemmanade}Let $\Mh$ be a closed, orientable hyperbolic $3$-orbifold, and set $\oldOmega=(\Mh)\pl$. Set
$M=|\oldOmega|$, and let $\cals$ be an $\oldOmega$-\dandy\ system of spheres
in $M$. Set $\oldPsi=\oldOmega\cut{\obd(\cals)}$, $N=|\oldPsi|=M\cut\cals$, $\toldTheta =\partial\oldPsi$, and
$\tcals=|\toldTheta|=\partial N$. 
%Set },{\toldTheta})$.
 %$\oldPhi^-=\oldPhi^- (\oldPsi)$. 
Let  $n\ge1$ be an integer, 
%and (in the notation of \ref{in duck tape}) 
and set $\noodge_m=\noodge_m^{\oldOmega,\toldTheta}\in\barcalm({\toldTheta},{\toldTheta})$
and 
$V_m=V_m^{\oldOmega,\toldTheta}\in\barcaly_-(\toldTheta)$ for $m=1,\ldots,n$
%$\iota=\iota_\oldPsi$ 
(see \ref{in duck tape}). Let $Z$ be a connected element of $\barcaly_-(\toldTheta)$ such that $Z\preceq
V_n$ (so that for $1\le m\le n$ we have $Z\preceq V_m$ by \ref{i
  thought so too}).
Then there exist an element $\frakZ$ of $\Theta_-(\toldTheta)$ with $[\frakZ]=Z$,
and embeddings
%$m=1,\ldots,n$, 
%let 
$\otherdge_1,\ldots,\otherdge_n:\frakZ\to\toldTheta$   such that
% element of $\calm(\toldTheta,\toldTheta)$ such that 
$[\frakZ,\otherdge_m,\otherdge_m(\frakZ)]=\noodge_m|Z$ for
$1\le m\le n$. Furthermore, if we 
denote by $\tS_0$ the component of $\tcals$ containing $|\frakZ|$, and
by $\tS_m$ 
the component of $\tcals$ containing
$|\otherdge_m(\frakZ)|$ for $m=1,\ldots,n$, and if we assume  that
 $|\frakZ|$ semi-splinters $\obd(\tS_0)$ (see \ref{splinter
  def}) and that $\wt\tS_m=\wt\tS_0$ for $m=1,\ldots,n$, then the
following conclusions hold:
%such that $[\frakZ]+Z$ and \preceq
%V_n$
%. \redcomment{If I keep this general form of the
%  statement I need to define $Z=[\frakZ]$} and $\noodge_m|Z\in\calm(\toldTheta,\toldTheta)$).
%(and hence $\noodge_m|[\frakZ]$ is an element of $\calm(\toldTheta,\toldTheta)$). 
%Let $\tS_0$ denote the component of $\tcals$ containing $|\frakZ|$. For $m=1,\ldots,n$, 
%let $\otherdge_m:\frakZ\to\toldTheta$ be an embedding such that
% element of $\calm(\toldTheta,\toldTheta)$ such that 
%$[\frakZ,\otherdge_m,\otherdge_m(\frakZ)]=\noodge_m|[\frakZ]$; and
%let $\tS_m$ denote
%the component of $\tcals$ containing
%$|\otherdge_m(\frakZ)|$. 
 %Suppose that $|\frakZ|$ semi-splinters $\obd(\tS_0)$ (see \ref{splinter
  %def}), and that $\wt\tS_m=\wt\tS_0$ for $m=1,\ldots,n$. Then:
\begin{enumerate}
%Let $j:\frakZ\to\toldTheta$ be an element of $\calm(\toldTheta,\toldTheta)$ such that $[j]=\noodge_n|[\frakZ]$.   
\item If $C$ is any component of $\partial\frakZ$, and if $D$ and $D'$ denote the components of $\tS_0-\inter|\frakZ|$ and $\tS_n-\inter \otherdge_n(|\frakZ|)$ bounded by $C$ and $\otherdge_n(C)$ respectively, we have $\wt D=\wt D'$. 
\item The subsurface $|\otherdge_n(\frakZ)|$ semi-splinters $\obd(\tS_n)$; and if $|\frakZ|$ splinters $\obd(\tS_0)$, then $|\otherdge_n(\frakZ)|$ splinters $\obd(\tS_n)$.
\end{enumerate}
\EndLemma
%\tS_0

\Proof
Let us fix an element $\frakZ\in\Theta_-(\toldTheta)$ such that
$[\frakZ]=Z$. Now given any $m$ with $1\le m\le n$, fix an element
$(\frakZ_m,\otherdge^0_m,\frakZ_m')$ of $\calm(\toldTheta,
\toldTheta)$ such that
$[\frakZ_m,\otherdge^0_m,\frakZ_m']=j_m|Z$. Since
$[\frakZ_m]=\dom(j_m|Z)=Z=[\frakZ]$, there is a homeomorphism
$h_m:\toldTheta\to\toldTheta$, isotopic to the identity, such that
$h_m(\frakZ)=\frakZ_m$. Then $\otherdge_m:=\otherdge^0_m\circ h_m:\frakZ\to\toldTheta$
is an
embedding, $(\frakZ,\otherdge_m,\otherdge_m(\frakZ))\in \calm(\toldTheta,
\toldTheta)$, and $[\frakZ,\otherdge_m,\otherdge_m(\frakZ)]=[\frakZ_m,\otherdge_m^0,\otherdge_m^0(\frakZ)]
=[\frakZ_m,\otherdge^0_m,\frakZ_m']
=\noodge_m|Z$.
%$m=1,\ldots,n$, 
%let 
%$\otherdge_1,\ldots,\otherdge_n:\frakZ\to\toldTheta$   such that
% element of $\calm(\toldTheta,\toldTheta)$ such that 
%$[\frakZ,\otherdge_m,\otherdge_m(\frakZ)]=\noodge_m|Z$ for
%$1\le m\le n$.

We must now show that if $\Mh$, $\cals$, $n$ and $Z$ satisfy the
general hypotheses of the lemma, if we are given $\frakZ\in \Theta_-(\toldTheta)$ with $[\frakZ]=Z$,
and embeddings
%$m=1,\ldots,n$, 
%let 
$\otherdge_1,\ldots,\otherdge_n:\frakZ\to\toldTheta$   such that
% element of $\calm(\toldTheta,\toldTheta)$ such that 
$[\frakZ,\otherdge_m,\otherdge_m(\frakZ)]=\noodge_m|Z$ for
$1\le m\le n$,  if we define $\tS_0$, and $\tS_m$ for $1\le m\le n$ as in the
statement of the lemma, and  if in addition we assume
that
 $|\frakZ|$ semi-splinters $\obd(\tS_0)$ and that $\wt\tS_m=\wt\tS_0$
 for $m=1,\ldots,n$, then Assertions (1) and (2) hold.

Let us first show that for a given
$\Mh$, $\cals$, $n$, $Z$,  
$\frakZ$
and
%$m=1,\ldots,n$, 
%let 
$\otherdge_1,\ldots,\otherdge_n$
satisfying all these assumptions, Assertion (1) implies
Assertion (2). 
Suppose that (1) holds, and consider an arbitrary  component $D$ of $\tS_n-\inter \otherdge_n(|\frakZ|)$. If we set $C=\partial D$, $C'=\otherdge_n^{-1}(C)$, and denote by $D'$ the component of $\tS_0-\inter|\frakZ|$ bounded by $C'$, then  Assertion (1) gives $\wt D=\wt D'$. By hypothesis we have $\wt S_n=\wt S_0$, and since $|\frakZ|$ semi-splinters $\obd(\tS_0)$ we have $\wt D'\le(\wt\tS_0)/2$. Hence
$\wt D\le(\wt\tS_n)/2$. Since $D$ is an arbitrary  component of
$\tS_n-\inter \otherdge_n(|\frakZ|)$, this means that
 $|\otherdge_n(\frakZ)|$ semi-splinters $\obd(\tS_n)$. Likewise, 
if $|\frakZ|$ splinters $\obd(\tS_0)$ we have $\wt D'<(\wt\tS_0)/2$, so that
$\wt D<(\wt\tS_n)/2$. This shows that
 $|\otherdge_n(\frakZ)|$ splinters $\obd(\tS_n)$, and Assertion (2) is
 established.
%C' cobble C D' D dobble

In order to prove Assertion (1), we first observe:
\Claim\label{duly noted}
The truth of Assertion (1) is independent of the
choice of $\frakZ$ and of $\otherdge_m$, 
subject to the conditions $[\frakZ]=Z$ and
$[\otherdge_m]=[\frakZ,\otherdge_m,\otherdge_m(\frakZ)]=\noodge_m|Z$.
\EndClaim

Assertion (1) will be proved by induction on $n$. For the case $n=1$ note that, according to
(\ref{and so did you}),  we have $V_1=[\oldPhi^-]$ and $\noodge_1=[\iota]$, where $\oldPhi^-=\oldPhi^-(\oldPsi)$ and $\iota=\iota_\oldPsi$.
In view of \ref{duly noted}, it follows that for $n=1$ we may take $\frakZ$ to be a suborbifold of $\oldPhi^-$, and $\otherdge_1$ to be $\iota|\frakZ$.

Set $F=|\frakZ|\subset\tS_0$.
%, and let $\tS$ denote the component of $\tcals$
%containing $F$ . 
Set $D_0=D$, and let $D_1,\ldots,D_k$ denote the
remaining components of $\tS_0-\inter F$, where $k\ge0$. Set
$C_i=\partial D_i$ for $i=1,\ldots,k$. Then $C_0,\ldots,C_k$ are the
boundary components of $F$, and hence the boundary components of $F':=|\iota|(F)$
may be indexed as $C_0',\ldots,C_k'$ where $C_i'=|\iota|(C_i)$
%$C_i'=h(C_i\times\{1\})$ 
for $i=1,\ldots,k$. According to \ref{what's iota?}, each component of
$\frakA(\oldPsi)\cap\oldSigma^-(\oldPsi)$ is invariant under $\iota$; hence
for each $i$ with $1\le i\le k$, the curves $\obd(C_i)$ and $\obd(C_i')$ are contained in the same component $\frakB_i$ of $\frakA(\oldPsi)$. Here $\frakB_i$ is an
essential annular suborbifold of $\oldPsi$. If $|\frakB_i|$ is a weight-$0$ annulus, it
follows from the first condition in the definition of a \dandy\ system that $C_i$ and $C_i'$ have the same size; if $|\frakB_i|$ is a disk we have $C_i=C_i'$. Thus in any case we have
$\size(C_i)=\size(C_i')=s_i$, say, for $i=0,\ldots,k$. Let $D_i'$
denote the component of $\tS_1-\inter F'$ bounded by $C_i'$, for
$i=0,\ldots,k$, so that $D_0'=D'$.

Set $p=\wt\tS_0$, so that by hypothesis we have $\wt\tS_1=p$. 
Since $F$ semi-splinters $\obd(\tS_0)$, and $p=\wt\tS_0$, we have $\wt D_i\le(\wt\tS_0)/2=p/2$ for
$i=0,\ldots,k$. In view of the definition of size, this means that
$\wt D_i=s_i$. We are required to prove that $\wt D_0=\wt D_0'$,
i.e. that $\wt D_0'=s_0$. Since $\size C_0'=s_0$, and $p=\wt\tS_1$, we must have either
$\wt D_0'=s_0$, or $s_0<p/2$ and $\wt D_0'=p-s_0$. We shall assume the
latter alternative and obtain a contradiction, thus completing the proof of Assertion (1) in the case $n=1$.

%k

Since
$D_0,\ldots,D_k$ are the components of $\tS_0-\inter F$, we have $p=\wt\tS_0=\wt F+\sum_{i=0}^k\wt D_i$.
Likewise, since
$D_0',\ldots,D_k'$ are the components of $\tS_1-\inter F'$, we have $p=\wt\tS_1=\wt F'+\sum_{i=0}^k\wt D_i'$.
% Since
%$\wt D_0'=p-s_0/2> p/2$, and since
 %$D_i'\cap D_0'=\emptyset$ for $i\ge1$, 
In particular we have $\wt D_i'\le p-\wt
 D_0'=s_0< p/2$ for $1\le i\le k$. Hence $\wt D_i'=\size C_i'=s_i$ for
 $1\le i\le k$. Note also that since $\obd(F)$ and $\obd(F')$ are
 homeomorphic, we have $\wt F=\wt F'$. It follows that 
%$$p-s_0
%=\wt D_0'
%=p-(\wt F'+\sum_{i=1}^k\wt D_i')=p-(\wt F'+\sum_{i=1}^ks_i)
%=p-(\wt F+\sum_{i=1}^k\wt D_i)=\wt D_0=s_0.$$
$$s_0=
p-\wt D_0'
=\wt F'+\sum_{i=1}^k\wt D_i'=\wt F'+\sum_{i=1}^ks_i
=\wt F+\sum_{i=1}^k\wt D_i=p-\wt D_0=p-s_0.$$
This is a contradiction, since $s_0<p/2$, and the case $n=1$ of (1) is proved.
%\wt D_0'\ge p/2$ and $\wt D_0<p/2$. Thus the second assertion is proved.

Now suppose that Assertion (1) is true for a given $n$, and suppose that
$\frakZ$ and $\otherdge_{n+1}$ satisfy the hypotheses with $n+1$ in
place of $n$.  
Let $C$ be any component of $\partial\frakZ$, and let $D$ and $D'$ denote the components of $\tS_0-\inter|\frakZ|$ and $\tS_{n+1}-\inter \otherdge_{n+1}(|\frakZ|)$ bounded by $C$ and $\otherdge_{n+1}(C)$ respectively. We must show that $\wt D=\wt D'$.

%Fix

Since $\wt\tS_m=\wt\tS_0$ for $m=1,\ldots,n+1$, we have in particular  $\wt\tS_m=\wt\tS_0$ for $m=1,\ldots,n$. 

According to the definitions (see \ref{in duck tape}), with \ref{more associativity},  we
have $\noodge_{n+1}=[\iota]\diamond[\tau]\diamond\noodge_n$. 
Since
$V_n=\dom\noodge_n$ and $V_{n+1}=\dom\noodge_{n+1}$, the hypothesis
$Z\preceq V_{n+1}$ may be rewritten as $ Z\preceq\dom
([\iota]\diamond[\tau]\diamond\noodge_n)$. Lemma \ref{before
  associativity} then gives that $ Z\preceq\dom\noodge_n=V_n$, 
that
$\noodge_n( Z)\preceq\dom([\iota]\diamond[\tau])$, and that
$\noodge_{n+1}| Z=([\iota]\diamond[\tau])\circ(\noodge_n|Z)$. A
second application of
Lemma \ref{before
  associativity} then gives that
$\noodge_n( Z)\preceq\dom
%([\iota]\diamond
[\tau]$, that
$[\tau](\noodge_n( Z))\preceq\dom[\iota]=[\oldPhi^-]$, and that
$([\iota]\diamond[\tau])|\noodge_n( Z)=[\iota]\circ([\tau]|\noodge_n( Z)$.
%[\frakZ]
Hence there is an element $\otherdge_n^*$ of
 $\calm(\toldTheta,\toldTheta)$ with $[\otherdge^*_n]=\noodge_n$, such that
 $\frakZ_n:=\otherdge^*_n(\frakZ)$ satisfies
 $\tau(\frakZ_n)\subset\oldPhi^-$;
%There is awkwardness in the statement, which begins with the fact that if I started with $Z\preceq V_n$ and defined $(\frakZ, \otherdge_m,\frakZ_m)$ to be a representative of the equivalence class $\noodge_m|Z$, it would not be a priori true that all these reos (for different values of $m$) would have the same domain.}
 and we then have
 $[\iota\circ\tau\circ(\otherdge^*_n|\frakZ)]=\noodge_{n+1}|\frakZ$. 
%\otherdge
In view of \ref{duly noted} it follows that, for the purpose of this
induction step. we may
 assume that $\otherdge_n$ and
 $\otherdge_{n+1}$ have been chosen within their equivalence
 classes in $\calm(\toldTheta,\toldTheta)$ so that $\otherdge_{n+1}=\iota\circ\tau\circ(\otherdge_n|\frakZ)$. 
Set $\frakZ_{n+1}=\otherdge_{n+1}(\frakZ)$, so that $\iota(\tau(\frakZ_n))=\frakZ_{n+1}$. 

Set $C''=\otherdge_n(C)$. Since $\otherdge_{n+1}=\iota\circ\tau\circ(\otherdge_n|\frakZ)$, we have 
$\otherdge_{n+1}(C)=\iota(\tau(C''))=\otherdge_1(\tau(C''))$.

Let $D''$ denote the component of 
$\tS_n-\inter |\frakZ_n|$  bounded
by $C''=\otherdge_n(C)$.
The induction
hypothesis gives $\wt D=\wt D''$. Since we have shown that (1) implies
(2), the
 induction hypothesis also implies that $|\frakZ_n|$ semi-splinters $\obd(\tS_n)$.

Since $|\frakZ_n|$ semi-splinters $\obd(\tS_n)$, the surface $|\tau(\frakZ_n)|$ semi-splinters $\tau(\obd(\tS_n))$. Furthermore, the components $\tau(\tS_n)$ and $\tS_{n+1}$ containing $|\tau(\frakZ_n)$ and $|\iota(\tau(\frakZ_n)|=|\otherdge_1(\tau(\frakZ_n)|$ have the same weight; $\tau(C'')$ is a component of $\partial\tau(\frakZ_n)$; and $\tau(D'')$ and $D'$ are the components of $\tS_n-\inter|\tau(\frakZ_n)|$ and $\tS_{n+1}-\inter \otherdge_{n+1}(|\frakZ|)$ bounded by $\tau(C'')$ and $\otherdge_{n+1}(C)=\otherdge_1(\tau(C'')$ respectively. Thus the hypotheses of Assertion (1) continue to hold if we let $1$, $\tau(\frakZ_n)$, $\tau(\tS_n)$, $\tS_{n+1}$, $\tau(C'')$, $\tau(D'')$ and $D'$ play the respective roles of $n$, $\frakZ$, $\tS_0$, $\tS_n$, $C$, $D$ and $D'$. Since the case $n=1$ of (1) has already been proved, it follows that $\wt\tau(D'')=\wt(D')$. Since we have seen that $\wt D=\wt D''=\wt\tau(D'')$, we may now conclude that $\wt D=\wt D'$, and
the
induction is complete. Thus Assertion (1) is proved.
\EndProof
%$h*_0

\Lemma\label{i'm glad you mentioned foos}
Let $\Mh$ be a closed, orientable hyperbolic $3$-orbifold, and set $\oldOmega=(\Mh)\pl$. Set
$M=|\oldOmega|$, and let $\cals$ be an $\oldOmega$-\dandy\ system of spheres
in $M$. Set $N=M\cut\cals$ and
$\tcals=\partial (M\cut\cals)$. If a component $\tS$ of
$\tcals$ contains an \bad\ component (see Definition \ref{bad def}) $F$ of $|\oldPhi(\oldPsi)|$, and if $U$ denotes the component of $|\oldSigma(\oldPsi)|$ containing $F$, then $U\cap\tcals=F\cup F'$ for some \bad\ component $F'$ of $|\oldPhi(\oldPsi)|$ contained in $\tau_\cals(\tS)$.
\EndLemma

\Proof
The definition of an \bad\ component implies that $F$ splinters $\obd(\tS)$, and that the component $U$ of $|\oldSigma(\oldPsi)|$ containing $F$
meets $\tS':=\tau_\cals(\tS)$. Since $F$ splinters $\obd(\tS)$, and since $\obd(\tS)$ is negative by \ref{doublemint}, we have
$\chi(\obd(F))<0$ by Lemma \ref{negative
  splinter}. Hence by
Lemma \ref{when a tore a fold}, $\obd(U)$ is a \pagelike\
\Ssuborbifold\ of $\oldPsi$, so that $U\cap\tcals$ has at most two
components. Since $U\cap\tcals\ne\emptyset$, there must be exactly one
component $F'$ of $U\cap\tS'$ distinct from $F$, and we must have
$F'\subset\tS'$. It follows from \ref{what's iota?} that
$\obd(F')$ is isotopic to $\iota_\oldPsi(\obd(F))$, which with (\ref{and so did
  you}) implies that in
$\barcaly_-(\obd(\tcals)$ we have $[\obd(F')]=\noodge_1([\obd(F)])$, where 
$\noodge_1=\noodge_1^{\oldOmega,\toldTheta}$.  We now apply Lemma
\ref{lemmanade}, taking $n=1$ and 
%$\frakZ=\obd(F)$
$Z=[\obd(F)]$, so that $\tS_0=\tS$, $\tS_1=\tS'$. Since $\tS'=\tau_\cals(\tS)$, we have $\wt\tS'=\wt\tS$, as required for the application of Lemma \ref{lemmanade}.
Since $[\obd(F)]=Z$ and 
$\noodge_1([\obd(F)])=[\obd(F')]$, 
we may take the element  $\frakZ$ of $\Theta_-(\obd(\tcals)$ given by Lemma
\ref{lemmanade} to be $\obd(F)$, and we may choose the
embedding $\otherdge_1:\frakZ\to\obd(\tcals)$ given by Lemma
\ref{lemmanade} in such a way that $\otherdge_1(\frakZ)=\obd(F')$.
Since $F$ splinters $\obd(\tS)$, it
%the 
%$[, and $\noodge_1(\frakZ)=\obd(F')$. It
%therefore
 follows from  Assertion (2) of Lemma \ref{lemmanade} 
that $F'$ splinters $\obd(\tS')$.  Since the component $U$ of $|\oldSigma(\oldPsi)|$ containing $F'$
meets $\tS=\tau(\tS')$, the component $F'$ of $|\oldPhi(\oldPsi)|$ is
by  definition \bad, and the conclusion is proved.
\EndProof

%DEsubcaseL_i=0

\Definition\label{subsystem}
Let $\oldOmega$ be a closed,
orientable $3$-orbifold, and let $\cals$ be an $\oldOmega$-admissible system of
spheres in $M:=|\oldOmega|$. Note that if $\cals'$ is any union of
components of $\cals$, then $\cals'$ is itself
an $\oldOmega$-admissible system of spheres in
$M$. We will call such a system $\cals'$ a {\it subsystem} of $\cals$.
\EndDefinition

\Lemma\label{before dandy subsystem}
Let $\Mh$ be a closed,
orientable hyperbolic $3$-orbifold, and set $\oldOmega=(\Mh)\pl$. Let $\cals$ be an $\oldOmega$-admissible system of
spheres in $M:=|\oldOmega|$. Let $\cals_1$ be a subsystem of $\cals$. Set $\oldPsi=\oldOmega\cut{\obd(\cals)}$ and $\oldPsi_1=\oldOmega\cut{\obd(\cals_1)}$, set $\tcals=\partial\oldPsi$ and $\tcals_1=\partial\oldPsi_1$, and identify $\tcals_1$ in the natural way with a union of components of $\tcals$.  Then:
\begin{enumerate}
 
\item if $Z$ is an \bad\ component 
of $|\oldPhi^-(\oldPsi)|$, and if the component $\tS$   of $\tcals$ containing $Z$ is contained in $\tcals_1$, 
% \bad\ , and if $Z$ is either a pseudo-belt for $\tS$ relative  to  $\cals$  or is \bad\ relative  to  $\cals$,
then $\obd(Z)$  is isotopic in $\obd(\tS)$ to a suborbifold %$\frakR$ 
of some component $\oldUpsilon$ of $\oldPhi^-(\oldPsi_1)$. Furthermore, if $\oldLambda$ denotes the component of $\oldSigma^-(\oldPsi_1)$ containing $\oldUpsilon$ (cf. \ref{what's iota?}), then  
%the inclusion homomorphism $\pi_1(\frakR)\to\pi_1(\oldLambda)$ is non-surjective in the case where $Z$ is a pseudo-belt relative to $\cals$; and 
$|\oldLambda|$ has non-empty intersection with $\tau_{\cals_1}(\tS)$.
% in the case where $Z$ is \bad\ relative to $\cals$.
\item if $\tS$ is a component of $\tcals_1$ which is pseudo-belted relative to $\tcals$, then there is a connected, $\iota_{\oldPsi_1}$-invariant
%of $|\oldPhi^-(\oldPsi)|$, if $\tS$ denotes the component of $\tcals$ containing $Z$, and if $Z$ is either a pseudo-belt for $\tS$ relative  to  $\cals$  or is \bad\ relative  to  $\cals$,
%then $\obd(Z)$  is isotopic in $\tS$ to a 
suborbifold
 $\frakV$ of 
%some component $\frakR$ of
 $\oldPhi^-(\oldPsi_1)\cap\obd(\tS)$  such that the involution $\iota_{\oldPsi_1}|\frakV$ of $\frakV$ is strongly equivalent (see \ref{strong equivalence}) in $\obd(\tS)$ to the involution $\epsilon_\tS^{\cals}$ of $\obd(E_\tS^\cals)$.
% Furthermore, if $\oldLambda$ denotes the component of $\oldSigma^-(\oldPsi_1)$ containing $\oldUpsilon$ (cf. \ref{what's iota?}), then  
%the inclusion homomorphism $\pi_1(\frakR)\to\pi_1(\oldLambda)$ is non-surjective in the case where $Z$ is a pseudo-belt relative to $\cals$; and 
%$|\oldLambda|$ has non-empty intersection with $\tS':=\tau_{\cals_1}(\tS)$ in the case where $Z$ is \bad\ relative to $\cals$.\frakR\oldUpsilon
 \item 
If $\cals$ is \dandy, and if $F$ is a component 
of $|\oldPhi(\oldPsi_1)|$ which splinters $\obd(\tS)$, where $\tS$ denotes  the component of $\tcals_1$ containing
$F$, then either 
(a) $F$ is \bad\ relative to $\cals_1$, and
$|\oldPhi(\oldPsi)|$ has an \bad\ component $Z\subset\tS$ relative to $\cals$
such that $\obd(F)$ and $\obd(Z)$ are isotopic in $\tS$, 
or 
(b) $\tS$
is  pseudo-belted relative to both $\cals$ and $\cals_1$; furthermore, 
$\obd(F)$ is isotopic in $\obd(\tS)$ both to $\obd(E_\tS^\cals)$ and to $\obd(E_\tS^{\cals_1})$, and $\epsilon _\tS^\cals$ is strongly equivalent to $\epsilon _\tS^{\cals_1}$ in $\obd(\tS)$.
 %and $\tS$ has a pseudo-belt $Z$ relative to $\tcals$ such that %$\obd(F)$ and $\obd(Z)$ are isotopic, 
\end{enumerate}
\EndLemma

Note that, under the hypotheses of Proposition \ref{before dandy subsystem}, if a  component $\tS$  of $\tcals_1$ is pseudo-belted relative to both
  $\cals_1$ and $\cals$, then  
$\epsilon_\tS^\cals$ and $\epsilon_\tS^{\cals_1}$  are defined up to strong equivalence in $\tS$ by \ref{stop beeping}; and for any
component $\tS$  of $\tcals_1$, the set of all components of $\oldPhi(\oldPsi_1)$ contained in $\obd(\tS)$ is well-defined up to isotopy by \ref{oldSigma def}. This shows that the three conclusions in the statement of Proposition \ref{dandy subsystem} are well formulated.

\Proof[Proof of Lemma \ref{before dandy subsystem}]
Set $\calt=\rho_{\cals_1}^{-1}(\cals-\cals_1)\subset \inter N_1$. The components of $\calt$ are $2$-spheres, and $\rho_{\cals_1}:N\to M$ maps $\calt$ homeomorphically onto $\cals-\cals_1$.

To prove (1) and (2), 
suppose that
$Z$ is a component 
of $|\oldPhi^-(\oldPsi)|$ contained in a component $\tS$  of $\tcals_1$, and that $Z$ is either a pseudo-belt for $\tS$ relative  to  $\cals$  or is \bad\ relative  to  $\cals$.

Set $\frakZ=\obd(Z)$, and let $\frakE$ denote the component of $\oldSigma^-(\oldPsi)$ containing $\frakZ$ (see \ref{what's iota?}). Set $\tS'=\tau_{\cals}(\tS)=\tau_{\cals_1}(\tS)$. 
Note that according to the definition of a pseudo-belt and the definition of \bad\ component, we have:

\Claim\label{NEW MAD}
If $Z$ is a pseudo-belt relative to $\cals$, then $\frakE\cap\obd(\tcals)=\frakZ$.
If $Z$ is \bad\ relative to $\cals$, then $\frakE\cap\obd(\tcals)$ has two components, one of which is $
\frakZ\subset\obd(\tS)$, while the other is contained in $\obd(\tS')$. 
\EndClaim

%\Claim\label{MAD}
%Either (a)
%$\frakE\cap\obd(\tcals)=\frakZ\subset\obd(\tS)$, or (b) %$\frakE\cap\obd(\tcals)$ has two components, of which one is $
%\frakZ\subset\obd(\tS)$ and the other is contained in $\tS'$.
%\EndClaim

In either case we have $\frakE\cap\obd(\tcals)\subset\obd(\tS)\cup\obd(\tS')\subset\obd(\tcals_1)$. Hence $\rho_{\obd(\calt)
}:\oldPsi\to\oldPsi_1$ maps $\frakE$ homeomorphically onto a
suborbifold $\frakE_1$ of $\oldPsi_1$, and
%. Note also that since 
%$\frakE\cap\obd(\tcals)\subset\obd(\tcals_1)$, we have
$\frakE_1\cap\obd(\calt)=\emptyset$.

Recall that $\tcals_1=\partial|\oldPsi_1|$ has been canonically identified with a union of components of $\tcals=\partial|\oldPsi|$. From this point of view, the map
$\rho_{\obd(\calt)}:\oldPsi\to\oldPsi_1$ restricts to the identity map on $\obd(\tcals_1)$. In particular, $\obd(\tS),\obd(\tS')\subset\obd(\tcals_1)$ may be regarded as boundary components of $\oldPsi_1$, so that $\frakZ$ is a suborbifold of the boundary of $\oldPsi_1$. Now \ref{NEW MAD} implies:

\Claim\label{MAD ONE}
If $Z$ is a pseudo-belt relative to $\cals$, we have
$\frakE_1\cap\obd(\tcals_1)=\frakZ\subset\obd(\tS)$. If $Z$ is \bad\ relative to $\cals$, then $\frakE_1\cap\obd(\tcals_1)$ has two components, of which one is $
\frakZ\subset\obd(\tS)$ and the other is contained in $\obd(\tS')$.
\EndClaim

Since $\frakE$ is a component of $\oldSigma(\oldPsi)$, the pair $(\frakE,\frakE\cap\obd(\tcals))$ is an \spair. Since $\chi(\frakE))<0$, it follows from Lemma \ref{when a tore a fold} that the \spair\ $(\frakE,\frakE\cap\obd(\tcals))$ is \pagelike.
The homeomorphism $\rho_{\obd(\calt)}| \frakE:\frakE\to\frakE_1$ may be regarded as a homeomorphism between the pairs
$(\frakE,\frakE\cap\obd(\tcals))$ and
$(\frakE_1,\frakE_1\cap\obd(\tcals_1))$. Hence
$(\frakE_1,\frakE_1\cap\obd(\tcals_1))$ is a \pagelike\ \spair.

%the relation $\frakZ\subset\frakE\cap\obd(\tcals)\subset\obd(\tS)\cup\obd(\tS')$ in $\oldPsi$ may be written as a relation $\frakZ\subset\frakE_1\cap\obd(\tcals)\subset\obd(\tS)\cup\obd(\tS'))$ in $\oldPsi_1$. 
%Furthermore, if we endow the component $\frakE$ of $\oldSigma^-(\oldPsi)$ with the $I$-fibration described in \ref{what's iota?}, then $\frakE_1$ inherits an $I$-fibration from $\frakE$ via the homeomorphism $\rho_{\obd(\calt)}| \frakE:\frakE\to\frakE_1$. Since $\partialh\frakE=\frakE\cap\obd(\tcals)\subset\obd(\tS)\cup\obd(\tS')\subset\obd(\tcals_1)$, we have $\partialh\frakE_1=\frakE_1\cap\obd(\tcals_1)\subset\obd(\tS)\cup\obd(\tS')$. This says that the pair $(\frakE_1,\frakE_1\cap\obd(\tcals_1))$ is a \pagelike\ \spair. In particular we have $\partialv\frakE_1=\Fr_{\oldPsi_1}\frakE_1$.

We claim;

\Claim\label{ok, ok}
Every component of  $\Fr_{\oldPsi_1}\frakE_1$ is an essential annular suborbifold of $\oldPsi_1$.
\EndClaim

To prove this, first note that $\Fr_{\oldPsi_1}\frakE_1$ is the image of $\Fr_\oldPsi\frakE$ under the homeomorphism $\rho_{\obd(\calt)}| \frakE:\frakE\to\frakE_1$. Since $\frakE$ is a component of $\oldSigma(\oldPsi)$, its frontier components are components of $\frakA(\oldPsi)$, and are therefore essential annular suborbifolds of $\oldPsi$ by \ref{tuesa day}. In particular, $\Fr_{\oldPsi}\frakE$ is $\pi_1$-injective in $\oldPsi$. But the map $\rho_{\obd(\calt)}:\oldPsi\to\oldPsi_1$ is $\pi_1$-injective by the admissibility of $\cals$. Hence
$\Fr_{\oldPsi_1}\frakE_1$ is $\pi_1$-injective in $\oldPsi_1$.

To complete the proof 
%that every component of $\partialv\frakE_1=\Fr_{\oldPsi_1}\frakE_1$ is an essential annular suborbifold of $\oldPsi_1$, 
of \ref{ok, ok}, it remains to show that no component of $\Fr_{\oldPsi_1}\frakE_1$ is parallel in $\oldPsi_1$ to a suborbifold of $\tcals_1$. If a component $\oldDelta$ of $\Fr\frakE_1$ is parallel in $\oldPsi_1$ to a suborbifold $\oldDelta'$ of $\tcals_1$, there is a suborbifold $\frakK$ of $\oldPsi$ admitting a trivial $I$-fibration such that $\partialv\frakK\subset\tcals_1$, and $\oldDelta$ and $\oldDelta'$ are the components of $\partialh\frakK$. Since $\oldDelta$ is annular, $\frakK$ is a \torifold. We have $\Fr_{\oldPsi_1}\frakK=\oldDelta$. We have observed that $\calt\subset\inter N_1$ (i.e. $\obd(\calt)\subset\inter\oldPsi_1$), and we have $\obd(\calt)\cap \Fr_{\oldPsi_1}\frakK
=\obd(\calt)\cap\oldDelta
\subset\obd(\calt)\cap\frakE_1=\emptyset$. Thus if $\obd(\calt)\cap\frakK\ne\emptyset$, then $\obd(T)\subset\inter\frakK$ for some component $T$ of $\calt$. Since the admissibility of $\cals$ implies that $\obd(\calt)$ is $\pi_1$-injective in $\oldPsi_1$, it then follows that $\obd(T)$ is $\pi_1$-injective in $\frakK$. This is impossible, since $\frakK$ is a \torifold, and $\obd(T)$ is negative by \ref{in duck tape}. Hence $\obd(\calt)\cap\frakK=\emptyset$. This implies that $\rho_{\obd(\calt)}$ maps some suborbifold $\tfrakK$ of $\oldPsi$ homeomorphically onto $\frakK$. The homeomorphism $\rho_{\obd(\calt)}|\tfrakK:\tfrakK\to\frakK$ maps some suborbifolds $\toldDelta$ and $\toldDelta'$ onto $\oldDelta$ and $\oldDelta'$ respectively; furthermore, pulling back the $I$-fibration of $\frakK$ via 
$\rho_{\obd(\calt)}|\tfrakK:\tfrakK\to\frakK$, we obtain an $I$-fibration of $\tfrakK$ such that $\partialv\tfrakK\subset\tcals$, and the components of $\partialh\tfrakK$ are $\toldDelta$ and $\toldDelta'$. This means that $\toldDelta$ and $\toldDelta'$ are parallel in $\oldPsi$. But $\toldDelta$ is a component of $\Fr_\oldPsi\frakE$, and $\toldDelta'\subset\tcals$. This contradicts the essentiality of the components of $\Fr_\oldPsi\frakE$ in $\oldPsi$. Thus the proof of \ref{ok, ok} is complete.

We have shown that the pair $(\frakE_1,\frakE_1\cap\obd(\tcals_1))$ is a \pagelike\ \spair, and according to \ref{ok, ok}, every component of 
$\Fr_{\oldPsi_1}\frakE_1$ is an essential annular suborbifold of $\oldPsi_1$, By definition this means that $\frakE_1$ is a \pagelike\ \Ssuborbifold\ of $\oldPsi_1$. Hence by Proposition \ref{slide it}, $\frakE_1$ is isotopic in $\oldPsi_1$ to a suborbifold of $\oldSigma(\oldPsi_1)$. Thus there is a homeomorphism $h:\oldPsi\to\oldPsi$, isotopic to the identity, such that  %and therefore to a suborbifold 
$\frakE_1':=h(\frakE_1)\subset\oldLambda-\Fr_\oldPsi\oldLambda$ for some component  $\oldLambda$ of $\oldSigma(\oldPsi_1)$.
%, with $\frakE_1'\cap \Fr_\oldPsi\oldLambda=\emptyset$. 
%In particular, $\frakZ\subset\obd(\tS)$, which by \ref{MAD ONE} is a component of $\frakE_1\cap\obd(\tcals_1)$, is isotopic to a component  $\frakR\subset\obd(\tS)$ of $\frakE_1'\cap\obd(\tcals_1)\subset\oldLambda\cap\obd(\tcals_1)$. 
The component $\oldUpsilon$ of $\oldLambda\cap\obd(\tcals_1)$
containing $\frakV:=h(\frakZ)$ is in particular a component of $\oldSigma(\oldPsi_1)\cap\obd(\tcals_1)=\oldPhi(\oldPsi_1)$. %and therefore to a suborbifold of $\frakV$ for some component $\frakV$ of %$\oldPhi(\oldPsi_1)$. 
Since $\frakZ$ is taut by \ref{tuesa day} and is negative by the definition of $\oldPhi^-(\oldPsi)$, the orbifold $\oldUpsilon$ is also negative, and is therefore a component of $\oldPhi^-(\oldPsi_1)$. 
The discussion in \ref{what's iota?} then shows that $\oldLambda$ is a component of $\oldSigma^-(\oldPsi_1)$.

In the case where $Z$ is \bad\ relative to $\cals$, it follows from \ref{MAD ONE} that $|\frakE_1|\cap\tS'\ne\emptyset$. Since $\frakE_1$ and $\frakE_1'$ are isotopic, we have $|\frakE_1'|\cap\tS'\ne\emptyset$, and since $\oldLambda\supset\frakE_1'$, we have $|\oldLambda|\cap\tS'\ne\emptyset$. This proves Assertion (1).

Now consider the case where $Z$ is a pseudo-belt relative to $\cals$. 
Since $\oldLambda$ is a component of $\oldSigma^-(\oldPsi_1)$, the discussion in \ref{what's iota?} gives an $I$-fibration $q:\oldLambda\to\frakB$, where $\frakB$ is some $2$-orbifold, such that 
$\partialh\oldLambda=\oldLambda\cap\tcals_1$. The discussion in \ref{what's iota?} also shows that we may choose $\iota_{\oldPsi_1}$ within its strong equivalence class in $\tcals_1$ so that $\iota_{\oldPsi_1}|\partialh\oldLambda$ is the non-trivial deck transformation of the two-sheeted cover $q|\partialh\oldLambda:\partialh\oldLambda\to\frakB$; in particular $\partialh\oldLambda$ is $\iota_{\oldPsi_1}$-invariant.
Since  \ref{ok, ok} implies that every component of  
$\Fr_{\oldLambda}\frakE_1'=
\Fr_{\oldPsi_1}\frakE_1'$ is an essential annular suborbifold of $\oldPsi_1$, no component of
$\Fr_{\oldPsi_1}\frakE_1'$ is 
parallel in the 
 pair $(\oldLambda,\oldLambda\cap\tcals_1)=(\oldLambda, \partialh\oldLambda)$ to a suborbifold of $\partialh\oldLambda=\oldLambda\cap\tcals_1$. It therefore follows from 
Proposition \ref{when vertical} that $\Fr_{\oldPsi_1}\frakE_1'$ is
isotopic, via an isotopy of the pair
$(\oldLambda,\partialh\oldLambda)$, to a saturated annular suborbifold of $\oldLambda$.
Hence we may suppose the homeomorphism $h$ to  have been chosen in such a way that  $\frakE_1'$ is a saturated  suborbifold of $\oldLambda$. Hence there is a $2$-suborbifold $\frakB_0$ of $\frakB$ such that $\frakE_1'=q^{-1}(\frakB)$, and $q_0:=q|\frakE_1':\frakE_1'\to\frakB_0$ is an $I$-fibration. In terms of this $I$-fibration we have $\partialh\frakE_1'=\frakE_1'\cap\partialh\oldLambda=\frakE_1'\cap\tcals_1$.
%E_1 cov

Since $Z$ is a pseudo-belt relative to $\cals$, it follows from \ref{MAD ONE} that $\frakE_1\cap\obd(\tcals_1)=\frakZ$. Hence $\partialh\frakE_1'=\frakE_1'\cap\tcals_1=\frakV$.
Thus the pre-image of $\frakB_0$ under  the two-sheeted cover $q|\partialh\oldLambda:\partialh\oldLambda\to\frakB$ is the connected suborbifold $\frakV$ of $\oldLambda\cap\tcals_1$. 
%This implies that 
 %the two-sheeted cover $q|\oldLambda\cap\tcals_1:\oldLambda\cap\tcals_1\to\frakB$admits no section, and hence that 
%$\oldLambda\cap\tcals_1$ is connected, i.e. $\oldLambda\cap\tcals_1=\oldUpsilon$.
This shows that  $\frakV$ is $\iota_{\oldPsi_1}$-invariant, and that
$\iota_{\oldPsi_1}|\frakV$ is the non-trivial deck transformation of
the two-sheeted cover $q|\partialh\oldLambda$. 

Since $q_0:\frakE_1'\to\frakB_0$ is an $I$-fibration, and $h\circ\rho_{\obd(\calt)}$ maps $\frakE$ homeomorphically onto $\frakE_1'$, the map $r:=q_0\circ h\circ\rho_{\obd(\calt)}:\frakE\to\frakB_0$ is an $I$-fibration. In terms of this $I$-fibration, recalling that $\rho_{\obd(\calt)}$ restricts to the identity map on $\obd(\tcals_1)$ we find that $\partialh \frakE=h^{-1}(\partialh\frakE_1')=h^{-1}(\frakV)=\frakZ$. Thus the $I$-fibration $r$ is compatible (in the sense of \ref{S-pair def}) with $\frakZ$. The discussion in \ref{what's iota?} then shows that, up to strong equivalence, the non-trivial deck transformation of the two-sheeted covering $r|\frakZ:\frakZ\to\frakB_0$ is equal to $\iota_\oldPsi|\frakZ$, which may be denoted $\epsilon^{\cals}_\tS$ since $Z=|\frakZ|$ is  a pseudo-belt for $\tS$ relative  to  $\cals$. But using the definition of $r$ and the fact that  $\iota_{\oldPsi_1}|\frakV$ is the non-trivial deck transformation of the two-sheeted cover $q|\partialh\oldLambda$,
 and again using
that $\rho_{\obd(\calt)}$ restricts to the identity map on $\obd(\tcals_1)$, we find that the non-trivial deck transformation of  $r|\frakZ$ is $(h|\frakZ)^{-1}\circ (\iota_{\oldPsi_1}|\frakV) \circ (h|\frakZ)$. Since $h$ is isotopic to the identity, this shows that $\epsilon^{\tcals}_\tS$ is strongly equivalent to 
$\iota_{\oldPsi_1}|\frakV$, and the proof of Assertion (2) is complete.

To prove Assertion (3), 
we assume that $\cals$ is \dandy, and we consider a component $F$
% of $\tcals$ such that $\tS$ contains no
%bad component 
of $|\oldPhi(\oldPsi_1)|$ which splinters $\obd(\tS)$, where $\tS$ denotes the component $\tS$ of $\tcals_1$ containing
$F$. Since $F$ splinters $\obd(\tcals)$, Lemma \ref{negative splinter} asserts that $\chi(\obd(F))<0$, so that $\obd(F)$ is a component of $\oldPhi^-(\oldPsi_1)$. According to Corollary \ref{burpollary}, we have $[\oldPhi^-(\oldPsi_1)]\preceq
[\oldPhi^-(\oldPsi)]$. Hence $\oldPhi^-(\oldPsi)$ may be chosen within its isotopy class in $\obd(\tS)$ so that $\oldPhi^-(\oldPsi_1)\subset
\oldPhi^-(\oldPsi)$. Let $\frakZ$ denote the component of $\oldPhi^-(\oldPsi)$ containing $\obd(F)$. Since $F$ splinters $\obd(\tS)$, so does $Z:=|\frakZ|$. Since $\cals$ is \dandy\ and therefore satisfies the second condition of Definition \ref{semidandy def}, either
$Z$ is a pseudo-belt for $\tS$ relative to $\cals$, or $Z$ is \bad\ relative to $\cals$. 

If $Z$ is \bad\ relative to $\cals$, we fix a component $\oldUpsilon$ of $\oldPhi^-(\oldPsi_1)$ having the properties stated in Assertion (1), and we fix a suborbifold $\frakV$ of $\oldUpsilon$ which is isotopic to $\frakZ=\obd(Z)$   in $\obd(\tS)$. We may take $\frakV$ to be contained in $\inter\oldUpsilon$. If $Z$ is a pseudo-belt for $\tS$ relative to $\cals$, we fix 
a connected, $\iota_{\oldPsi_1}$-invariant
suborbifold
 $\frakV$ of 
 $\oldPhi^-(\oldPsi_1)\cap\obd(\tS)$  
having the properties stated in Assertion (2). In this case we let $\oldUpsilon\subset\obd(\tS)$ denote the component of $\oldPhi^-(\oldPsi)$ containing $\frakV$; we may again take $\frakV$ to be contained in $\inter\oldUpsilon$. Since $\iota_{\oldPsi_1}|\frakV$  is strongly equivalent  in $\obd(\tS)$ to $\epsilon_\tS^{\tcals}$, in particular $\frakV$ is isotopic in $\obd(\tS)$ to  $\obd(E_\tS^\cals)=\obd(Z)=\frakZ$ (see \ref{strong equivalence}). Thus in either case, $\oldUpsilon$ is a component of $\oldPhi^-(\oldPsi)$ whose interior contains $\frakV$, and $\frakV$ is isotopic in $\obd(\tS)$ to $\frakZ$.

%It now follows from Assertion (1) that $\frakZ:=\obd(Z)$  is isotopic in $\tS$ to a suborbifold $\frakV$ of some component $\oldUpsilon$ of $\oldPhi^-(\oldPsi_1)$. Furthermore, if $\oldLambda$ denotes the component of $\oldSigma^-(\oldPsi_1)$ containing $\oldUpsilon$ (cf. \ref{what's iota?}), then  either (A) $Z$ is a pseudo-belt relative to $\cals$, and
%the inclusion homomorphism $\pi_1(\frakV)\to\pi_1(\oldLambda)$ is non-surjective; or (B) $Z$ is \bad\ relative to $\cals$, and 
%$|\oldLambda|$ has non-empty intersection with $\tS':=\tau_{\cals_1}(\tS)$.

Since $\obd(F)\subset\frakZ$, and $\frakZ$ is isotopic to  $\frakV\subset\oldUpsilon$, we have the relation $[\obd(F)]\preceq[\frakZ]=[\frakV]\preceq[\oldUpsilon]$ in $\barcaly_-(\obd(\tcals_1))$. Since $[\obd(F)]\preceq[\oldUpsilon]$, we have $[\obd(F)]\wedge[\oldUpsilon]=[\obd(F)]$. But $\obd(F)$ and $\oldUpsilon$ are both components of $\oldPhi^-(\oldPsi_1)$, and are therefore either equal or disjoint. If $\obd(F)\cap\frakZ=\emptyset$, it follows from Corollary \ref{nafta} that $[\obd(F)]\wedge[\oldUpsilon]=\emptyset$; this now implies $[\obd(F)]=[\emptyset]$, which is impossible since $F$ is connected and hence non-empty. We must therefore have $\oldUpsilon=\obd(F)$. We now obtain $[\obd(F)]\preceq[\frakV]\preceq[\obd(F)]$, and since $\preceq$ is a partial order by Proposition \ref{new partial order}, it follows that $[\frakV]=[\obd(F)]$.
% Equivalently, we have $[\frakV]=[\obd(F)]$, which 
This means that the
inclusion map $i:\frakV\to\oldUpsilon=\obd(F)$ is isotopic in
$\obd(\tS)$ to a homeomorphism $j:\frakV\to\obd(F)$. Since $i$ and
 $j$ are isotopic in $\obd(\tS)$, it follows from Corollary \ref{i
 guess} that they are (non-ambiently) isotopic in $\oldUpsilon$; that
 is, the inclusion map from $\frakV$ to $\obd(F)$ is isotopic in
 $\obd(F)$ to a homeomorphism of $\frakV$ onto $\obd(F)$. In
 particular, the inclusion $\frakV\to\obd(F)$ is a homotopy
 equivalence.
Since $\frakV\subset\inter\oldUpsilon$, it then follows that $\oldUpsilon=\obd(F)$ is a strong regular neighborhood of $\frakV$ in $\obd(\tS)$.

Consider the case in which $Z$ is \bad\ relative to $\cals$.  Let
$\oldLambda$ denote the component of $\oldSigma^-(\oldPsi_1)$
containing $\oldUpsilon$. 
%We will equip $\oldLambda$ with the
%$I$-fibration described in \ref{what's iota?}. We then have
%$\partialh\oldLambda=\oldLambda\cap\obd(\tcals_1)$. 
%In particular, $\oldUpsilon\subset\tS$ is a component of
%$\partialh\oldLambda$. But 
According to Assertion (1), $|\oldLambda|$ 
%also
 has non-empty intersection with $\tS':=\tau_{\tcals_1}(\tS)=\tau_{\tcals}(\tS)$.  
%Thus $\oldLambda\cap \obd(\tS)$  and $\oldLambda\cap \obd(\tS')$ are both non-empty. %intersection with both $\obd(\tS)$ and $\obd(\tS')$.
%, and since $\frakE_1$ and $\frakE_1'$ are isotopic, $\frakE_1'$ also meets both $\obd(\tS)$ and $\obd(\tS')$. 
%Hence $\partialh\oldLambda=\oldLambda\cap\obd(\tcals_1)$ meets both
%$\obd(\tS)$ and $\obd(\tS')$, and is therefore not connected. This
%shows that the $I$-fibration of $\oldLambda$ is trivial, and that one
%component of $\partialh\oldLambda$ is contained in $\obd(\tS)$ while
%the other is contained in $\obd(\tS)$. Since $\oldUpsilon$ is a
%component of $\partialh\oldLambda$, it must be the unique component
%contained in $\obd(\tS)$. 
Since $F=|\oldUpsilon|$ splinters $\obd(\tS)$, 
%and  $\oldLambda$ is a component of $\oldSigma^-(\oldPsi_1)$, 
it 
%now 
follows from the definition that $F$ is \bad\ relative to $\cals_1$. But we have seen that $[\obd(F)]=[\frakV]=[\frakZ]=[\obd(Z)]$, i.e.  $\obd(F)$ and $\obd(Z)$ are isotopic in $\tS$. This gives Alternative (a) of the conclusion of Assertion (3) in the case where  $Z$ is \bad\ relative to $\cals$.

Now consider the case in which $Z$ is a pseudo-belt for $\tS$ relative to $\cals$, so that $\tS$ is pseudo-belted relative to $\cals$ and $\frakZ$ is isotopic to $\obd(E^{\cals}_{\tS})$. Since we have seen that  $[\obd(F)]=[\frakZ]=[\obd(Z)]$, it follows that $\obd(F)$ is isotopic to $\obd(E^{\cals}_{\tS})$. 

According to Assertion (2), $\frakV$ is invariant under the involution
$\iota_{\oldPsi_1}$ of $\oldPhi^-(\oldPsi_1)$. Since
$\oldUpsilon\subset\tS$ is the component of $\oldPhi^-(\oldPsi_1)$
containing $\frakV$, the suborbifold $\oldUpsilon$ is also
$\iota_{\oldPsi_1}$-invariant. Now if $\oldLambda$ denotes the
component of $\oldSigma^-(\oldPsi_1)$ containing $\oldUpsilon$, it
follows from \ref{what's iota?} that $\oldLambda\cap\tcals_1$ has
either one or two components, and that if it has two then they are
interchanged by $\iota_{\oldPsi_1}$. Since the component $\oldUpsilon$
of $\oldLambda\cap\tcals_1$  is $\iota_{\oldPsi_1}$-invariant,
$\oldLambda\cap\tcals_1$ must be the connected suborbifold
$\oldUpsilon$. Furthermore, we have seen that
$[\oldUpsilon]=[\obd(F)]=[\frakZ]$, i.e. $\oldUpsilon$ is isotopic to
$\frakZ$ in $\tS$; since $Z=|\frakZ|$ is a pseudo-belt for $\tS$
relative to $\cals$, it follows from Definition \ref{what's a belt?}
that $F=|\oldUpsilon|$ splinters $\obd(\tS)$. Again applying
Definition \ref{what's a belt?}, we now deduce that $F$ is a
pseudo-belt for $\tS$ relative to $\cals_1$ (so that $\tS$ is in
particular pseudo-belted relative to $\cals_1$). The definition of
$\epsilon^{\cals_1}_\tS$ (see \ref{stop beeping}) now gives that
$\iota_{\oldPsi_1}|\oldUpsilon$ is strongly equivalent to
$\epsilon^{\cals_1}_\tS$. But since $\frakV$ is
$\iota_{\oldPsi_1}$-invariant, and since $\oldUpsilon$ has been seen
to be a strong regular neighborhood of $\frakV$, the involution
$\iota_{\oldPsi_1}|\oldUpsilon$ is also strongly equivalent to
$\iota_{\oldPsi_1}|\frakV$, which is in turn strongly equivalent to
$\epsilon^{\cals}_\tS$ by Assertion (2). Hence $\epsilon^{\cals}_\tS$
and $\epsilon^{\cals_1}_\tS$ are strongly equivalent. Furthermore, the
strong equivalence of $\iota_{\oldPsi_1}|\oldUpsilon$,
$\epsilon^{\cals}_\tS$
and $\epsilon^{\cals_1}_\tS$ implies by definition that $\oldUpsilon=\obd(F)$, $\obd(E^{\cals}_\tS)$
and $E^{\cals_1}_\tS$ are all isotopic in $\obd(\tS)$.
This gives Alternative (b) of the conclusion of Assertion (3) in the case where  $Z$ is a pseudo-belt for $\tS$ relative to $\cals$.
%Now suppose that (A) holds. Then the inclusion homomorphism $\pi_1(\frakV)\to\pi_1(\oldLambda)$ is non-surjective.
%Since we have shown that the inclusion $\frakV\to\obd(F)$ is a homotopy equivalence, it now follows that 
%the inclusion homomorphism from $\pi_1(\obd(F))=\pi_1(\oldUpsilon)$ to $\pi_1(\oldLambda)$ is not surjective. Since $\obd(F)=\oldUpsilon$ is a component of $\partialh(\oldLambda)$, it now follows that the $I$-fibration of $\oldLambda$ is non-trivial, and hence $\oldLambda\cap\tcals_1=\partialh(\oldLambda)=\obd(F)$. Since $\oldLambda$ is a component of $\oldSigma^-(\oldPsi_1)$, and $F$ splinters $\obd(\tS)$, it now follows from the definition that $F$ is a pseudo-belt for $\tS$. 
\EndProof
%\frakZ'\frakV'\frakC\frakB(2) F' (A) (B)

\Proposition\label{dandy subsystem}
Let $\Mh$ be a closed,
orientable hyperbolic $3$-orbifold, and set $\oldOmega=(\Mh)\pl$. Let $\cals$ be an $\oldOmega$-dandy system of
spheres in $M:=|\oldOmega|$. Let $\cals_1$ be a subsystem of $\cals$. Then $\cals_1$ is also
$\oldOmega$-\dandy. Furthermore, if we set $\oldPsi=\oldOmega\cut{\obd(\cals)}$, $\oldPsi_1=\oldOmega\cut{\obd(\cals_1)}$, $\tcals=\partial|\oldPsi|$ and $\tcals_1=\partial|\oldPsi_1|$, and identify $\tcals_1$ in the natural way with a union of components of $\tcals$, then the following implications hold: 
\begin{itemize}
%%\redcomment{The assumption that
    %$\tS$ is a \central\ component of $\tcals$ is new. Fix
    %cross-refs. going forward. I have revised the statement, but not
    %the proof (yet), of Lemma \ref{no wonder 2}. I also need to fix the app. of Lemma \ref{no wonder 2} to Theorem \ref{manifold homology bound}.}
\item If a  component $\tS$  of $\tcals_1$ is belted relative to $\cals_1$, and
  is \central\ relative to  $\cals$, then $\tS$ is belted relative to $\cals$, and  $\obd(G_\tS^\cals)$ and $\obd(G_\tS^{\cals_1})$ are isotopic in $\obd(\tS)$. 
\item If a component $\tS$ of $\tcals_1$ contains a component $F_1$ of $|\oldPhi(\oldPsi_1)|$ which is \bad\ relative to $\cals_1$, then 
$\tS$ contains a component $F_1$ of $|\oldPhi(\oldPsi)|$ which is \bad\ relative to $\cals$; furthermore, $\obd(F)$ and $\obd(F_1)$ are isotopic in $\obd(\tS)$. 
\item If a  component $\tS$  of $\tcals_1$ is pseudo-belted relative to
  $\cals_1$ then $\tS$ is pseudo-belted relative to $\cals$. Furthermore, 
$\epsilon_\tS^\cals$ and $\epsilon_\tS^{\cals_1}$ (see
  \ref{stop beeping}) are strongly equivalent (see \ref{strong equivalence}), so that in particular
  $\obd(E_\tS^\cals)$ and $\obd(E_\tS^{\cals_1})$ are  isotopic.  
\end{itemize}
\EndProposition

Note that, under the hypotheses of Proposition \ref{dandy subsystem},
if  a  component $\tS$  of $\tcals_1$ is belted relative to both
$\cals_1$ and $\cals$, then  $\obd(G_\tS^\cals)$ and
$\obd(G_\tS^{\cals_1})$ are defined up to isotopy in $\obd(\tS)$
according to \ref{stop beeping}. Likewise, as we pointed out above in
the context of the statement of Lemma \ref{before dandy subsystem}, if a  component $\tS$  of $\tcals_1$ is pseudo-belted relative to both
  $\cals_1$ and $\cals$, then  
$\epsilon_\tS^\cals$ and $\epsilon_\tS^{\cals_1}$  are defined up to strong equivalence in $\tS$ by \ref{stop beeping}; and for any
component $\tS$  of $\tcals_1$, the set of all components of
$\oldPhi(\oldPsi_1)$ contained in $\obd(\tS)$ is well-defined up to
isotopy by \ref{oldSigma def}. This shows that the three bulleted
assertions 
in the statement of Proposition \ref{dandy subsystem} are well formulated.

\Proof[Proof of Proposition \ref{dandy subsystem}]
Set $M=\oldOmega$, $N=|\oldPsi|=M\cut\cals$,
$N_1=|\oldPsi_1|=M\cut{\cals_1}$, $\tcals=\partial N$, and
$\tcals_1=\partial N$.  Set $\calt=\rho_{\cals_1}^{-1}(\cals-\cals_1)\subset \inter N_1$. The components of $\calt$ are $2$-spheres, and $\rho_{\cals_1}:N\to M$ maps $\calt$ homeomorphically onto $\cals-\cals_1$.

In order to verify
that $\cals_1$ satisfies the first condition in the definition of a
\dandy\ system (see \ref{semidandy def}), suppose that $A$ is a weight-$0$ annulus, properly embedded in $N_1$. We must show that the boundary curves $C$ and $C'$ of $A$ have the same size. First consider the case in which $\obd(A)$ is not $\pi_1$-injective in
%Since $A\cap\fraks_\oldPsi=\emptyset$, the boundary curves
% $C$ and $C'$ are both homotopically trivial in
$\oldPsi_1$. 
%\redmissingref{Make sure it's clear what ``homotopically
   %trivial'' means.} 
Since the
\dandy\ system $\cals$ is by definition admissible, the $2$-orbifold $\obd
(\tcals_1)$ is $\pi_1$-injective in $\oldPsi_1$, and hence $C$ and $C'$ bound
discal suborbifolds of $\obd(\tcals_1)$. 
Since $\obd(\tcals)$ is orientable, this means that $C$ and $C'$ bound
  disks $E,E'\subset\cals_1$ such that
  $\wt E,\wt E'\le1$. Hence $\size C,\size C'\le1$. If $\size C\ne\size C'$, we may assume by symmetry that $\size C=1$ and that $\size C=0$, so that $C$ and $C'$ bound disks $E_1,E_1'\subset\tcals_1$ with $\wt E_1=1$ and $\wt E_1'=0$. But then the sphere $X:=E_1\cup A\cup E_1'$ has weight $1$; this means that $\obd(X)$ is a bad $2$-orbifold. The map $\rho_{\obd(\cals_1)}|\obd(X)$ is an immersion of this bad $2$-orbifold in $\oldOmega$, a contradiction to the hyperbolicity of $\Mh$. 

There remains the case in which 
 $\obd(A)$ is $\pi_1$-injective in
$\oldPsi_1$. In this case we apply Proposition \ref{newer no disks lemma}, letting $\oldPsi_1$ play the role of $\oldPsi$ in that lemma, and taking $\oldPi=\obd(\calt
%\rho_{\cals_1}^{-1}(\cals-\cals_1)
)\subset \inter \oldPsi_1$. The incompressibility of $\oldPi$ follows from the hypothesis that $\cals$ is \dandy, and in particular admissible. According to Assertion (1) of Proposition \ref{newer
 no disks lemma}, applied with $\frakB=\obd(A)$, we may suppose $A$ to be chosen within its isotopy class in $N_1-\fraks_{\oldPsi_1}$ in such a way that $\obd(A)$
%$\obd(A)$ is isotopic rel $\partial \frakB$  to a suborbifold which
has reduced intersection with $\frakP$. 
It then follows from Assertion (2) of Proposition \ref{newer no disks lemma} 
that
no component of
$\obd(A\cap\calt)$ bounds a discal suborbifold of $\obd(A)$. Since $A$ has weight $0$, this implies that no component of
$A\cap\calt$ bounds a disk in $A$. Hence
%  that $C$ and $C'$ bound a properly embedded weight-$0$ annulus $A_0\subset \oldPsi_1$ such that
%no component of $A_0\cap\calt$ bounds a disk in $A_0$. (In the notation of Proposition \ref{newer no disks lemma} we have $A_0=|\frakB_0|$.)
%This implies that
 there exist an integer $n\ge1$ and a homeomorphism
 $\eta:\SSS^1\times[0,n]\to A$ such that
 $\eta(\SSS^1\times(\{1,\ldots,n-1\}))=A\cap\calt$. We may suppose $\eta$
 to be chosen in such a way that $\eta(\SSS^1\times\{0\})=C$ and
 $\eta(\SSS^1\times\{n\})=C'$. Set $C_i=\eta(\SSS^1\times\{i\})$ for
 $i=0,\ldots,n$ and $B_i=\eta(\SSS^1\times[i-1,i])$ for
 $i=1,\ldots,n$. Then $B_i\subset N_1$ is an annulus whose boundary
 components $C_{i-1}$ and $C_i$ lie in $\tcals_1\cup\calt$, and whose
 interior is disjoint from $\tcals_1\cup\calt$. Hence $B_i$ is the
 image under $\rho_\calt:N\to N_1$ of a properly embedded annulus
 $\tB_i$ in $N$. Since $B_i\subset A$, we have $\wt B_i=0$, and hence
 $\wt \tB_i=0$. Let us label the components of $\partial \tB_i$ as
 $\tC_i$ and $\tC_i'$, where $\rho_\calt(\tC_i)=C_i$ and
 $\rho_\calt(\tC_i')=C_{i-1}$. Then $\size(\tC_i)=\size(C_i)$ and
 $\size(\tC_{i}')=\size(C_{i-1})$. But since  $\tB_i$ is a weight-$0$ annulus, properly embedded in $N$, and since $\cals$ is by hypothesis a \dandy\ system, we have $\size(\tC_i')=\size(\tC_i)$. Hence $\size(C_{i-1})=\size(C_i)$. Since this holds for $i=1,\ldots,n$, it follows that $\size(C_{0})=\size(C_n)$, i.e. that $\size(C)=\size(C')$. This completes the verification that $\cals_1$ satisfies the first condition in the definition of a
\dandy\ system.

The second condition in the definition of a \dandy\ system asserts that if $F$ is a component
of $|\oldPhi(\oldPsi_1)|$ which splinters $\obd(\tS)$, where $\tS$ denotes  the component of $\tcals$ containing
$F$, then
%. Let
%$E$ be a component of $|\oldPhi(\oldPsi)|$ contained in $\tS$. Suppose that
either $F$ is a pseudo-belt for $\tS$ relative to $\cals_1$, or $F$ is \bad\ relative to $\cals_1$. Since $\cals$ is \dandy, this is included in
Assertion (3) of Lemma \ref{before dandy subsystem}.

We now turn to the verification that $\cals_1$ satisfies the third condition in the def of a \dandy\ system. We will use the following preliminary observation:

\Claim\label{prelims}
Let $\Delta\subset|\oldPsi_1|$ be a disk (not necessarily properly embedded), transverse to $\fraks_{\oldPsi_1}$, and suppose that the orientable $2$-orbifold $\obd(\Delta)$ is annular, or equivalently that $\wt\Delta=2$ and that the points of $\Delta\cap\fraks_{\oldPsi_1}=\fraks_{\obd(\Delta)}$ are of order $2$. Then the inclusion homomorphism $\pi_1(\obd(\Delta))\to\pi_1(\oldPsi_1)$ is injective if and only if the conjugacy class in $\pi_1(\oldPsi_1)$ defined by (an orientation of) $\partial\Delta$ consists of elements of infinite order.
\EndClaim

To prove \ref{prelims}, note that $\pi_1(\obd(\Delta))$ is an infinite dihedral group, and that the subgroup which is the image of the inclusion homomorphism $\pi_1(\partial(\obd(\Delta)))\to\pi_1(\obd(\Delta))$ is the index-$2$ infinite cyclic subgroup of the infinite dihedral group. A homomorphism from the infinite dihedral group to another group is injective if and only if its restriction to the index-$2$ infinite cyclic subgroup is injective. This proves \ref{prelims}.

Now suppose that $D$ is a properly embedded disk  in $N$,
transverse to $\fraks_\oldPsi$, such that (i) $\obd(D)$ is annular and
(ii) $\obd(D)$ is $\pi_1$-injective in $\oldPsi$. We must show that
$\gamma:=\partial D$ is the boundary of a disk $E\subset\tcals$ such
that $\obd(E)$ is an annular orbifold. After modifying $D$ by a small
general-position isotopy, we may assume that (iii) $D$ is transverse
to $\calt$, and $D\cap\calt\cap\fraks_{\oldPsi_1}=\emptyset$. Among
all disks satisfying (i)---(iii), we may assume that $D$ is chosen so
as to minimize the number of components of $D\cap\calt$. We claim:

\Claim\label{albel}
For any component  $\beta$ of $D\cap\calt$, the following conditions are equivalent:
\begin{enumerate}[(a)]
\item the conjugacy class in $\pi_1(\oldPsi_1)$ defined by (an orientation of) $\beta$ consists of elements of finite order;
\item the size of the simple closed curve $\beta\subset\calt_0$, where $\calt_0$ denotes the component of $\calt$ containing $\beta$, is at most $1$;
\item the subdisk of $D$ bounded by $\beta$ has weight at most $1$.
\end{enumerate}
\EndClaim

To prove \ref{albel}, we will first show that (a) and (b) are equivalent. If (b) holds then $\beta$ bounds a disk $\Delta\subset\calt_0$ with $\wt\Delta\le1$. It follows that $\obd(\Delta)\subset\oldPsi_1$ is discal, and hence (a) holds. Now suppose that (b) does not hold. Then $\beta$ divides $\calt_0$ into two disks of weight at least $2$, and hence defines an element of $\pi_1(\obd(\calt_0))$ consisting of elements of infinite order. Since the inclusion homomorphism $\pi_1(\obd(\calt))\to\pi_1(\oldPsi_1)$ is injective by admissibility, it now follows that (a) does not hold.

Now let us show that  (a) and (c) are equivalent. Let $\Delta_1$ denote the subdisk of $D$ bounded by $\beta$. If (c) holds, i.e. if $\wt\Delta_1\le1$, then $\obd(\Delta_1)\subset\oldPsi_1$ is discal, and hence (a) holds. Now suppose that (c) does not hold, i.e. that $\wt(\Delta_1)>1$. Since $\wt(D)=2$, the curve $\beta$ cobounds a weight-$0$ annulus with $\gamma$. Hence $\beta$ and $\gamma$ (suitably oriented) define the same conjugacy class in $\pi_1(\obd(\oldPsi_1)$. But since the inclusion homomorphism $\pi_1(\obd(D))\to\pi_1(\oldPsi_1)$ is injective, it follows from \ref{prelims} that the conjugacy class in $\pi_1(\oldPsi_1)$ defined by (an orientation of) $\gamma$ consists of elements of infinite order. Hence (a) does not hold. This completes the proof of \ref{albel}.

Next, we claim:

\Claim\label{dumbbell}
If  $\beta$ is a  component of $D\cap\calt$ such that the size of the simple closed curve $\beta\subset\calt_0$, where $\calt_0$ denotes the component of $\calt$ containing $\beta$, is at most $2$, then the size of $\beta$ in $\calt_0$ is equal to the weight of the subdisk of $D$ bounded by $\beta$.
\EndClaim

To prove \ref{dumbbell}, let $s$ denote the size of $\beta\subset\calt_0$, and set $w=\wt\Delta_1$, where $\Delta_1$ denotes the subdisk of $D$ bounded by $\beta$. We have $w\le\wt D=2$, and by the hypothesis of \ref{dumbbell} we have $s\le2$. According to \ref{albel} we have $s\le1$ if and only if $w\le1$. Hence we need only rule out the possibility that either $s=1$ and $w=0$, or that $s=0$ and $w=1$. Suppose that one of these situations occurs. Let $\Delta\subset\calt_0$ be a disk such that $\wt D=s$, and let $h$ be an immersion of $\SSS^2$ in $|\oldPsi|$ which maps the upper hemisphere to $\Delta$ and the lower hemisphere to $\Delta_1$. Then there exist an orbifold $\calr$ such that $|\calr|=\SSS^2$, and an orbifold immersion $\frakh:\calr\to\oldPsi_1$ such that $|\frakh|=h$. We have $\wt\calr=w+s=1$, and since $|\calr|$ is a $2$-sphere, $\calr$ is a bad orbifold. But $\oldPsi_1$ is strongly \simple\ (by \ref{doublemint}) and therefore very good, and a bad orbifold cannot be immersed in a good one. Thus \ref{dumbbell} is proved.

We now claim:

\Claim\label{dumber bell}
If $D\cap\calt\ne\emptyset$, there is a component  $\beta$ of
$D\cap\calt$ such that the size of the simple closed curve $\beta$ in $\calt_0$, where $\calt_0$ denotes the component of $\calt$ containing $\beta$, is at most $2$.
\EndClaim

To prove \ref{dumber bell}, note that any component of $D\cap\calt$ bounds a subdisk of $D$. If $D\cap\calt\ne\emptyset$, then among all subdisks of $D$ bounded by components of $D\cap\calt$, we may choose one, say, $\Delta_0$, which is minimal with respect to inclusion. Set $\beta=\partial\Delta_0$. Let $s$ denote the size of $\beta\subset\calt_0$, and set $w=\wt\Delta_0$. We have $w\le\wt D=2$. If $w\le1$, then by \ref{albel} we have $s\le1$. Now suppose that $w=2$. The minimality of $\Delta_0$ implies that $\Delta_0\cap\calt=\beta$. Hence there is a properly embedded disk $\tDelta_0\subset|\oldPsi$ such that $\rho_{\obd(\calt)}:\oldPsi\to\oldPsi_0$ maps $\obd(\tDelta_0)$ homeomorphically onto $\obd(\Delta_0)$. We have $\wt\tDelta_0=\wt\Delta_0=2$. Now since $w>1$, it  follows from \ref{albel} that the conjugacy class in $\oldPsi_1$ defined by $\beta$ consists of elements of infinite order; hence by \ref{prelims}, $\obd(\Delta_0)$ is $\pi_1$-injective in $\oldPsi_1$. It then follows that
$\obd(\tDelta_0)$ is $\pi_1$-injective in $\oldPsi$. Since
 $\cals$ is \dandy\ and therefore satisfies the third condition of
 Definition \ref{semidandy def}, $\partial\tDelta_0$ must bound a
 weight-$2$ disk $\tE_0\subset\tcals$. Since
 $\partial\tE_0=\partial\tDelta_0\subset\rho_{\calt}^{-1}(\calt_0)$,
 the disk $\tE_0$ must be contained in a component of
 $\rho_{\calt}^{-1}(\calt_0)$. The image of $\tE_0$ under $\rho_\calt$
 is therefore a disk of weight $2$ contained in $\calt_0$ and bounded
 by $\beta$. This shows that $\size\beta\le2$, and \ref{dumber bell} is proved.

Next, we claim:

\Claim\label{tinkerbell}
We have $D\cap\calt=\emptyset$.
\EndClaim

To prove \ref{tinkerbell}, assume that $D\cap\calt\ne\emptyset$. Then
by \ref{dumber bell}, there is a component  $\beta$ of $D\cap\calt$
such that the size $s$ of the simple closed curve
$\beta\subset\calt_0$, where $\calt_0$ denotes the component of
$\calt$ containing $\beta$, is at most $2$. The definition of size
also guarantees that $s\le\wt(\calt_0)/2$, and that $\beta$ bounds a
disk $\Delta\subset\calt_0$ with $\wt\Delta=s$. Among all subdisks of
$\Delta$ (including $\Delta$ itself) which are bounded by components
of $D\cap\calt_0$, choose one, say $\Delta_0$, which is minimal with
respect to inclusion. Set $\beta_0=\partial\Delta_0$. We have
$\wt\Delta_0\le\wt\Delta=s\le\wt(\calt_0)/2$, and hence
$\wt\Delta_0=\size\beta_0$. But $\wt\Delta_0\le\wt\Delta=s\le2$, and
hence $\size\beta_0\le2$. It therefore follows from \ref{dumbbell}
that, if $\Delta_1$ denotes the subdisk of $D$ bounded by $\beta_0$,
then $\wt\Delta_1=\size\beta_0=\wt\Delta_0$. Now the minimality of
$\Delta_0$ implies that $\Delta_0\cap D=\beta_0$. Hence $D_1:=(D-\Delta_1)\cup\Delta_0$ is a properly embedded disk in $|\oldPsi_1|$. Since $\wt\Delta_1=\wt\Delta_0$, we have $\wt D_1=\wt D=2$. Modifying $D_1$ by a small isotopy we obtain a weight-$2$ properly embedded disk $D_2$ in $|\oldPsi_1|$, transverse to $\calt$, such that $\partial D_2=\partial D_1=\gamma$, $D_2\cap\calt\cap\fraks_{\oldPsi_1}=\emptyset$ and $D_2\cap\calt=(D_1\cap\calt)-\Delta_0$. Since $\obd(D)$ is $\pi_1$-injective and $\partial D_2=\gamma=\partial D$, it follows from \ref{prelims} that $\obd(D_2)$ is $\pi_1$-injective in $\oldPsi_1$. Thus Conditions (i)---(iii) above hold with $D_2$ in place of $D$. But $D_2\cap\calt$ is a union of components of $D\cap\calt$, not including the component $\beta_0$. This shows that $D_2\cap\calt$ has fewer components than $D\cap\calt$, and our choice of $D$ is contradicted. Thus \ref{tinkerbell} is proved.

We are now in a position to show that $\gamma$ is the boundary of a  disk $E\subset\tcals$ such that $\obd(E)$ is an annular orbifold. 

It follows from \ref{tinkerbell} that there exists a  properly embedded disk $\tD$ in $N=|\oldPsi|$ such that $\rho_\calt$ maps $\tD$ homeomorphically onto $D$. Then $\obd(\tD)$ is annular since $\obd(D)$ is annular, and is $\pi_1$-injective in $\oldPsi$ since $\obd(D)$ is $\pi_1$-injective in $\oldPsi_1$. Since $\cals$ is dandy and therefore satisfies the third condition of Definition \ref{semidandy def}, there is a disk $\tE\subset\cals$ such that $\partial\tE=\partial\tD$, and $\obd(\tE)$ is annular. If we now set $E=\rho_\calt(\tE)$, it follows that $E$ has the required properties. This
completes the proof that $\cals_1$ satisfies  the third condition of Definition \ref{semidandy def}. This completes the proof that $\cals_1$ is \dandy.

We now turn to the proofs of three bulleted assertions in the statement. We prove the last two of these first. Suppose that 
$\tS$ is a
component  of $\tcals_1$, and that  $F\subset\tS$ is a component of $|\oldPhi(\oldPsi_1)|$ which is either a pseudo-belt of $\tS$ relative to
  $\cals_1$, or a \bad\ component relative to $\cals_1$. In either case it follows from the definitions that $F$ splinters $\obd(\tS)$. Hence one of the alternatives (a) and (b) of  Assertion (3) of Lemma \ref{before dandy subsystem} must hold. 

If $F\subset\tS$ were an \bad\ component of $|\oldPhi(\oldPsi_1)|$
relative to $\cals_1$, but Alternative (b) of Assertion (3) of Lemma
\ref{before dandy subsystem}  held, then both $\obd(F)$ and
$\obd(E^{\tcals_1}_{\tS})$ would be components of
$|\oldPhi(\oldPsi_1)|$ contained in $\tS$. These components of
$|\oldPhi(\oldPsi_1)|$ would have to be distinct, because the
component of $\oldSigma(\oldPsi_1)$ containing $\obd(F)$ would meet $\tau_{\cals_1}(\obd(\tS))$, while the component of $\oldSigma(\oldPsi_1)$ containing $\obd(E^{\tcals_1}_{\tS})$ would not. Hence 
$|\oldPhi(\oldPsi_1)|$ would have two distinct components which splinter
$\obd(\tS)$, a contradiction to Lemma \ref{it's-a unique}. This shows that if 
$F$ is a \bad\ component of $|\oldPhi(\oldPsi_1)|$ relative to $\cals_1$, then Alternative (a) of Assertion (3) of Lemma \ref{before dandy subsystem} must hold. The same argument shows that if $F\subset\tS$ is a 
pseudo-belt of $\tS$ relative to
  $\cals_1$
then Alternative (b) of Assertion (3) of Lemma \ref{before dandy subsystem}  must hold.

If $F$ is an \bad\ component of $|\oldPhi(\oldPsi_1)|$ relative to $\cals_1$, then since Alternative (a) holds, we have in particular that $\tS$ contains a component $F_1:=Z$ of $|\oldPhi(\oldPsi)|$ which is \bad\ relative to $\cals$, and, $\obd(F)$ and $\obd(F_1)$ are isotopic in $\obd(\tS)$. This establishes the second bulleted assertion of the present proposition. If $F$ is a 
pseudo-belt of $\tS$ relative to
  $\cals_1$,
then since Alternative (b) holds, we have in particular that $\tS$ is pseudo-belted relative to $\cals$, and that 
$\epsilon_\tS^\cals$ and $\epsilon_\tS^{\cals_1}$  are strongly equivalent in $\obd(\tS)$. This establishes the third bulleted assertion of the present proposition.

To prove the first bulleted assertion, suppose that $\tS$ is a component  of $\tcals_1$ which is belted relative to $\cals_1$ and
  is \central\ relative to  $\cals$. Assume that $\tS$ is not belted relative to $\cals$. 
Then according to
Lemma \ref{and the foos just keep on fooing}, either $\tS$ is
 pseudo-belted relative to $\cals$, or 
$\tS$ contains an \bad\ component of $|\oldPhi(\oldPsi)|$ relative to $\cals$. 
Let
$Z$ denote an \bad\ component of $|\oldPhi(\oldPsi)|$ relative to $\cals$ if one exists,
and if $\tS$ is
 pseudo-belted relative to $\cals$ set $Z=E_\tS^\cals$.
%of $|\oldPhi^-(\oldPsi)|$ which either is a pseudo-belt for $\tS$
%relative  to  $\cals$  or is \bad\ relative  to  $\cals$. 
In both
cases it follows from the definitions that $Z$ splinters $\obd(\tS)$.
It  follows from Assertion (1) of Lemma \ref{before dandy
  subsystem} (in the case where $\tS$ contains an \bad\ component of
$|\oldPhi(\oldPsi)|$ relative to $\cals$) or from 
Assertion (2) of Lemma \ref{before dandy
  subsystem} and the definition of strong equivalence (in the case where $\tS$ is pseudo-belted relative to
$\cals$) that $\obd(Z)$ is isotopic to a suborbifold of some component $\oldUpsilon\subset\obd(\tS)$ of $\oldPhi^-(\oldPsi_1)$.
%\redcomment{This doesn't work, of course. But I think it may
%be easy to get what I want if I directly quote Corollary
%\ref{burpollary} to get $[\oldPhi^-(\oldPsi_1)]\preceq
%[\oldPhi^-(\oldPsi)]$; since $Z$ splinters $\obd(\tS)$ (the component of
%$\obd(\tcals)$ containing it, and I hope I have the conventions for
%splintering right), the component of } that, after suitably choosing $\oldPhi(\oldPsi)$ within its isotopy class, we may assume that 
%$\frakZ=\obd(Z)$  is  a suborbifold of some component $\oldUpsilon$
%of $\oldPhi^-(\oldPsi_1)$. 
Since $Z$ splinters $\obd(\tS)$,
%and $Z\subset|\oldUpsilon|$, 
the surface $|\oldUpsilon|$ also splinters $\obd(\tS)$. But since $\tS$ is belted relative to $\cals_1$, it follows from the first assertion of Lemma \ref{it's-a unique} that no component of $|\oldPhi(\oldPsi_1)|$ can splinter $\tS$. This contradiction shows that $\tS$ is belted relative to $\tS$. It remains to show that $\obd(G_\tS^\cals)$ and $\obd(G_\tS^{\cals_1})$  are isotopic in $\obd(\tS)$.

To prove this, note that since 
Corollary \ref{burpollary} gives $[\oldPhi^-(\oldPsi_1)]\preceq
[\oldPhi^-(\oldPsi)]$, we may assume $\oldPhi^-(\oldPsi_1)$ and $\oldPhi^-(\oldPsi)$ to be chosen within their isotopy classes so that $\oldPhi^-(\oldPsi_1)\subset\inter
\oldPhi^-(\oldPsi)$. Since $G_\tS^\cals$ and $G_\tS^{\cals_1}$  are respectively components of 
$|\oldPhi^-(\oldPsi)|$ and
$|\oldPhi^-(\oldPsi_1)|$, we must have either 
$G_\tS^{\cals_1}\subset
\inter G_\tS^\cals$ or
$G_\tS^{\cals_1}\cap
G_\tS^\cals=\emptyset$. In either case we have
\Equation\label{mumblety peg}
(\partial G_\tS^{\cals_1})\cap
(\partial G_\tS^\cals)=\emptyset.
\EndEquation

Let us choose components $C$ and $C_1$ of
$\partial G_\tS^\cals$
and $\partial G_\tS^{\cals_1}$ respectively. By (\ref{mumblety peg}) we have $C\cap C_1
=\emptyset$. Since $\tS$ is a sphere, $C$ and $C_1$ respectively bound disks $D$ and $D_1$ with $D_1\subset\inter D$. But by Lemma \ref{how do i recognize a belt?}, we have $\size C=\size C_1=\wt(\tS)/2$, which implies that $\wt D=\wt D_1=\wt(\tS)/2$. Hence $\wt(\overline{D-D_1})=\wt(D)-\wt(D_1)=0$, so that $\overline{D-D_1}$ is a weight-$0$ annulus whose boundary curves are $C$ and $C_1$. This implies:
\Claim\label{i didn't see nothin'}
The $1$-suborbifolds $\obd(C)$ and $\obd(C_1)$ are isotopic in $\tS$.
\EndClaim

Now since the orbifolds $\obd(G_\tS^\cals)$
and $\obd(G_\tS^{\cals_1})$  are annular (and orientable), each of the manifolds 
$G_\tS^\cals$
and $G_\tS^{\cals_1}$ 
is either a weight-$0$ annulus or a weight-$2$ disk. If $G_\tS^\cals$
and $G_\tS^{\cals_1}$ are both weight-$0$ annuli, it follows from \ref{i didn't see nothin'} that
 $\obd(G_\tS^\cals)$
and $\obd(G_\tS^{\cals_1})$  are isotopic in $\tS$, as required. If $G_\tS^\cals$
and $G_\tS^{\cals_1}$ are both weight-$2$ disks, it follows from \ref{i didn't see nothin'} that
 $\obd(G_\tS^\cals)$ is isotopic in $\tS$ either to
 $\obd(G_\tS^{\cals_1})$ or to $\obd(\tS-\inter G_\tS^{\cals_1})$. But since  $\obd(G_\tS^\cals)$
and $\obd(G_\tS^{\cals_1})$ are both annular, the latter alternative would imply that $\obd(\tS)$ is toric, a contradiction since $\obd(\tS)$ is negative according to \ref{doublemint}. Hence in this case $\obd(G_\tS^\cals)$
and $\obd(G_\tS^{\cals_1})$  are  again isotopic in $\tS$.

There remains the possibility that one of the manifolds $G_\tS^\cals$
and $G_\tS^{\cals_1}$ is a weight-$0$ annulus while the other is a
weight-$2$ disk. Thus we may write
$\{\obd(G_\tS^\cals),\obd(G_\tS^{\cals_1})\}=\{\oldGamma,\oldGamma'\}$,
where $\oldGamma$ and $\oldGamma'$ are annular,
$|\oldGamma|$ is a weight-$0$ annulus, and $|\oldGamma'|$ is a
weight-$2$ disk. Set $C_0=C$ if $\oldGamma=\obd(G_\tS^\cals)$, and set
$C_0=C_1$ if $\oldGamma=\obd(G_\tS^{\cals_1})$; thus $C_0$ is a
boundary component of $\oldGamma$. It follows from \ref{i didn't see
  nothin'} that $C_0$ bounds a suborbifold $\frakZ_0$ of $\obd(\tS)$ homeomorphic to
$\oldGamma'$. Since $\oldGamma'$ has connected boundary, either (i) $\frakZ_0$ is a component of
$\overline{\obd(\tS)-\oldGamma}$, or (ii) 
$\frakZ_0\supset\oldGamma$. We define a suborbifold $\frakZ$ of
$\obd(\tS)$ by setting $\frakZ=\frakZ_0$ if (i) holds, and
$\frakZ=\overline{\frakZ_0-\oldGamma}$ if (ii) holds. Thus in any event $\frakZ$ is a component of
$\tS-\inter\oldGamma$. Furthermore, $\frakZ$ is homeomorphic to
$\frakZ_0$; this is trivial if (i) holds, and if (ii) holds it follows
from the fact that $\oldGamma$ is a weight-$0$ annulus. In particular,
$\frakZ$ is annular.

Now set $\oldPsi_0=\oldPsi$ if $\oldGamma=\obd(G_\tS^\cals)$, and set
$\oldPsi_0=\oldPsi_1$ if $\oldGamma=\obd(G_\tS^{\cals_1})$. In either
case, $\oldPsi_0$ is a componentwise strongly \simple, componentwise
boundary-irreducible, orientable $3$-orbifold (see \ref{doublemint}, and 
%(equal to either $\oldPsi$ or $\oldPsi_1$) and an annular component 
$\oldGamma$ 
is an annular
  component 
of $\oldPhi(\oldPsi_0)$. The component of $|\partial\oldPsi_0|$
containing $|\oldGamma|$ is the $2$-sphere $\tS$, and the component
$\frakZ$ of $(\partial\oldPsi_0)-\oldGamma$ is annular. This
contradicts Corollary \ref{mangenfeffer}. 
% (equal, respectively, to $\obd(G_\tS^\cals)$ or $\obd(G_\tS^{\cals_1})$), such that the component $Z:=\tS$ of $|\partial\oldPsi_0|$ containing $\oldGamma$ is a $2$-sphere, $|\oldGamma|$ is a weight-$0$ annulus, and $|\partial\oldGamma|$  has a component $C_0$  (equal, respectively, to $C$ or to $C_1$) such that $\obd(C_0)$ is isotopic in $\obd(Z)$ to the boundary of a suborbifold of an annular suborbifold of $\obd(Z)$ (namely   $\obd(G_\tS^{\cals_1})$ or $\obd(G_\tS^{\cals})$ respectively) whose underlying surface is a weight-$2$ disk. It then follows that some component $E$ of $Z-\inter|\oldGamma|$ is a weight-$2$ disk such that $\obd(E)$ is annular. If $E$ is a component of $Z\setminus\inter\oldPhi(\oldPsi_0)$ then we have a contradiction to Corollary \ref{mangenfeffer}. If $E$ is not a component of $Z\setminus\inter\oldPhi(\oldPsi_0)$, then it follows from \ref{cobound} that $E$ contains a weight-$0$ annulus which shares one boundary component with $|\oldGamma|$ and the other with $|\oldGamma'|$, where $|\oldGamma'|$ is another component of $\oldPhi(\oldPsi_0)$. This contradicts the uniqueness assertion of Proposition \ref{lady edith}. Thus the proof that $\obd(G_\tS^\cals)$ and $\obd(G_\tS^{\cals_1})$  are isotopic in $\obd(\tS)$ is complete.
\EndProof
%\cals'\tcals'N'\oldPsi'\frakR
%F{newer\frakK\frakZ\oldGamma\frakE\frakB\frakC\frakD\oldDelta\oldGamma\oldUpsilon\frakW\oldDelta\frakD\frakK\oldLambda\frakL\oldPi\frakP\frakR\frakZ
%(a) (i) \frakR \frakV 

\Proposition\label{kavanaugh}
Let $\Mh$ be a closed,
orientable hyperbolic $3$-orbifold, and set $\oldOmega=(\Mh)\pl$. Let $\cals$ be an $\oldOmega$-dandy system of
spheres in $M:=|\oldOmega|$. Set $\oldPsi=\oldOmega\cut{\obd(\cals)}$ and
$\tcals=|\partial\oldPsi|=\partial M\cut\cals$. Suppose that $\tS$ is a pseudo-belted component of $\tcals$. If $C\subset E_\tS$ is a simple closed curve such that
$C\cap\fraks_{\obd(\tS)}=\emptyset$, then $\size|\epsilon_\tS|(C)=\size C$.
\EndProposition

\Proof
We may assume that $C\subset\inter E_\tS$. Set $\epsilon=\epsilon_\tS$.
Set $n=\wt\tS$, and set $s=\size C$. Then $s\le n/2$, and $C$ bounds a disk $\Delta\subset\tS$ of weight $s$. Since $E_\tS$ is connected, each component of $\tS-\inter E_\tS$ is a disk. Let $D_1,\ldots,D_m$ denote the components of $\tS-\inter E_\tS$ contained in $\Delta$. Then $G:=\overline{\Delta-(D_1\cup\cdots\cup D_m)}$ is a connected (planar) surface. Set $C_i=\partial D_i$ for $i=1,\ldots,m$. Then $C,C_1,\ldots,C_m$ are the boundary components of $G$, and hence the boundary components of $G':=|\epsilon|(G)$ are $C':=|\epsilon|(C)$, and $C_i':=|\epsilon|(C_i)$ for $i=1,\ldots,m$. For $i=1,\ldots,m$, the curve $C_i$ is a boundary component of $E_\tS$, and hence $C_i'$ is also a boundary component of $E_\tS$. Let $D_i'$ denote the component of $\tS-\inter E_\tS$ bounded by $C_i'$.

By the definition of a pseudo-belt, $E_\tS$ semi-splinters $\tS$, and hence for $i=1,\ldots,m$, each of the components $D_i$ and $D_i'$ of $\tS-\inter E_\tS$ has weight at most $n/2$. It follows that $\size C_i=\wt D_i$ and that $\size C_i'=\wt D_i'$. Now recall from \ref{stop beeping} that $\frakR_\tS$ is a component of $\oldSigma^-(\oldPsi)$ and
  that $\obd(E_\tS)=\frakR\cap\obd(\tcals)$. Thus for $i=1,\ldots,m$
  we have
  $\obd(C_i)\subset\partial\obd(E_\tS)\subset\partial\frakA(\oldPsi)$. Let
  $\frakB_i$ denote the component of $\frakA(\oldPsi)$ containing
  $\obd(C_i)$. By \ref{what's iota?}, we may assign an $I$-fibration
  to $\oldSigma^-(\oldPsi)$ in such a way that
  $\oldPhi^-(\oldPsi)=\partialh\oldSigma^-(\oldPsi)$. Then
  $\partialv\oldSigma^-(\oldPsi)$ is a union of components of $\frakA(\oldPsi)$. In particular the annular orbifold $\frakB_i$ is saturated in the fibration, and is therefore $\iota_\oldPsi$-invariant by \ref{what's iota?}. Hence $\obd(C_i')= \epsilon(\obd(C_i))=\iota_\oldPsi(\obd(C_i))$ is a component of $\partial\frakB_i$. 

The orientable annular orbifold $\frakB_i$ has at most two boundary components. Hence if $C_i'\ne C_i$, then $\partial \frakB_i$ has exactly the two boundary components
$\obd(C_i)$ and $\obd(C_i')$. In this case $|\frakB_i|$ is a weight-$0$ annulus with boundary curves $C_i$ and $C_i'$. The first condition in the definition of a \dandy\ system (\ref{semidandy def}) then implies that $\size(C_i)=\size(C_i')$. The latter equality is of course trivial in the case where $C_i=C_i'$, and therefore holds in all cases. But we have seen that $\size C_i=\wt D_i$ and that $\size C_i'=\wt D_i'$. Hence $\wt D_i=\wt D_i'$ for $i=1,\ldots,m$. 

On the other hand, $\epsilon$ maps $\obd(G)$ homeomorphically onto $\obd(G')$, and therefore $\wt G=\wt G'$. If $\Delta'$ denotes the disk $G'\cup D_1'\cup\cdots\cup D_m'$, then $\wt\Delta'=\wt G'+\sum_{i=1}^m\wt D_i'=\wt G+\sum_{i=1}^m\wt D_i=\wt\Delta=s$. Since $C':=\epsilon(C)=\partial\Delta'$ and $s\le n/2$, it follows that $\size C'=s$, which is the conclusion of the proposition.
\EndProof

%proofread

\section{Existence of complete dandy systems}\label{dandy existence section}

As we indicated in the Introduction and at the beginning of this
chapter, the object of this section is to prove that for every closed, orientable hyperbolic $3$-orbifold  $\Mh$, there exists an $(\Mh)\pl$-\dandy\ system of spheres $\cals\subset
|(\Mh)\pl|$ which is a complete system of spheres (see Definition \ref{complete
  def}) in $|(\Mh)\pl|$. This is included in  Proposition \ref{semidandies exist} below, which is stated in a slightly stronger form with an eye to possible future applications. In order to state Proposition \ref{semidandies exist} we need a few more definitions.

\begin{notationdefinitions}\label{spell it out}
Let $\Mh$ be a closed, orientable hyperbolic $3$-orbifold, and set $\oldOmega=(\Mh)\pl$. We will
denote by $\Gamma_\oldOmega$  the set of all complete systems
of spheres in $M:=|\oldOmega|$ which are transverse to
$\fraks_\oldOmega$. Let $\cals\in\Gamma_\oldOmega$ be given, and set
$\oldPsi=\oldOmega\cut{\obd(\cals)}$,
$N=|\oldPsi|= M\cut\cals$, and
$\tcals=\partial N$.
%We will define $r(\cals)=\wt_\oldOmega\cals$ and $q(\cals)=\card \calc(\cals)$. \redcomment{%This hardly seems worth introducing new letters for now that I have $\wt$ and $\comp%num$ (which I forgot to use here).}
Let us denote by $\nu(\cals)\in\NN$ the number of points of $\cals\cap\fraks_\oldOmega $ that have order strictly greater than $2$. We shall set $\mu(\cals)=(
\wt_\oldOmega\cals,\compnum(\cals),\nu(\cals))\in\NN^3$, and regard the set $\NN^3$ as a well-ordered set under the lexicographical ordering.

An element $\cals\in\Gamma_\oldOmega$ will be termed {\it$\mu$-minimal} if
the inequality
$\mu(\cals)\le\mu(\cals')$ holds in  $\NN^3$ for every
$\cals'\in\Gamma_\oldOmega$. 
\end{notationdefinitions}

%i

\Lemma\label{it is what it is}
If $\Mh$ is any closed, orientable hyperbolic $3$-orbifold, and if we set $\oldOmega=(\Mh)\pl$, then $\Gamma_\oldOmega$ has a $\mu$-minimal element. Furthermore, any $\mu$-minimal element of $\Gamma_\oldOmega$ is a minimal complete system of spheres in $|\oldOmega|$, in the sense of Definition \ref{min def}.
\EndLemma

\Proof
It follows from
Proposition \ref{kneser} that $|\oldOmega|$ contains a complete system of
spheres. A small
isotopy then gives a complete system of spheres which is transverse to
$\fraks_\oldOmega$. This shows that $\Gamma_\oldOmega\ne\emptyset$. Since 
$\NN^3$ is well ordered, the first assertion of the lemma now follows. To prove the second assertion, note that if
$\cals$ is a $\mu$-minimal element of $\Gamma_\oldOmega$, and
 $\cals'$ is a complete system of spheres in $M:=|\oldOmega|$ which is the union of a proper subset of the components of $\cals$, then $\cals'\in\Gamma_\oldOmega$; and since 
$\wt_\oldOmega\cals'\le \wt_\oldOmega\cals$ and
$\compnum(\cals'))<\compnum(\cals))$, we have
$\mu(\cals')<\mu(\cals)$, a contradiction to $\mu$-minimality.
\EndProof

\Proposition\label{semidandies exist} Let $\Mh$ be a closed, orientable hyperbolic $3$-orbifold, and set $\oldOmega=(\Mh)\pl$. Then any $\mu$-minimal element of $\Gamma_\oldOmega$ is an $\oldOmega$-dandy system of spheres.
In particular, there exists an $\oldOmega$-\dandy\ system of spheres $\cals\subset
|\oldOmega|$ which is a complete system of spheres in $|\oldOmega|$.
%such that every component of $M\cut\cals$ is $+$-irreducible 
\EndProposition 

The material in Subsections \ref{why admissible}---\ref{fluorescent foo-foo} consists of preliminaries for the proof of Proposition
\ref{semidandies exist}.

\Lemma\label{why admissible}
Let $\Mh$ be any closed, orientable hyperbolic $3$-orbifold, and set $\oldOmega=(\Mh)\pl$. Let $\cals$ be any $\mu$-minimal element of  $\Gamma_\oldOmega$. Then 
$\cals$ is an $\oldOmega$-admissible system of spheres,
i.e. $\obd(\cals)$ is incompressible in
$\oldOmega$. 
\EndLemma

\Proof
We must show that $\obd(\cals)$ is $\pi_1$-injective in $\oldOmega$, and that $\cals$ has no
component $S$ such that
$\obd(S)$ is a spherical orbifold. If  $\cals$ does have a
component $S$ such that
$\obd(S)$ is a spherical orbifold, then since $\Mh$ is hyperbolic 
and therefore irreducible, and $\obd(S)$ is two-sided, $\obd(S)$ bounds a discal $3$-suborbifold
 $\oldLambda\subset\oldOmega$. In
particular, $|\oldLambda|$ is a ball. 
But according to Lemma \ref{it is what it is}, $\cals$ is a minimal complete system of spheres in $M:=|\oldOmega|$, and hence by Proposition \ref{noble}, no component of $\cals$  can bound a ball in $M$. This contradiction shows that no component of $\obd(S)$ is spherical.

Now suppose that $\obd(\cals)$ is not $\pi_1$-injective in $\oldOmega$. Then by Proposition \ref{kinda dumb} (more specifically the implication (c)$\Rightarrow$(a)), there is an orientable  discal $2$-suborbifold $\frakD$ of $\oldOmega$
%, in general position with respect to $\fraks_\oldOmega$, with $\frakD\cap\obd(\cals)=\partial\frakD$, 
such that $\partial\frakD$ does not bound a discal suborbifold of $\obd(\cals)$. (The hypothesis in Proposition \ref{kinda dumb} that the ambient manifold is very good follows from the hyperbolicity of $\Mh$.) 
%general position  implies that $|\frakD|\cap\fraks_\oldPsi$ is finite, so that $\frakD$ is orientable, and discality then implies that $\wt|\frakD|\le1$. 
Thus there is a disk $D\subset M$ with
$D\cap\cals=\partial D$, such that $D$ meets $\fraks_\oldOmega$
transversally in at most one point, but
such that any disk contained in $\cals$ and bounded by $\partial D$
meets $\fraks_\oldOmega$
in at least two points. Let $S$ denote the component of
$\cals$ containing $\partial D$, let $D_1$ and $D_2$ denote the
closures of the components of $S-\partial D$ (so that $\wt D_i\ge2$ for $i=1,2$), and for $i=1,2$ set
$S_i=D\cup D_i$ and $\cals_i=(\cals-S)\cup S_i$. 
Then according to
Proposition \ref{taser}, either $\cals_1$
or $\cals_2$ is a complete system of spheres in $M$. After relabeling
the $D_i$ if necessary, we may assume that $\cals_1$ is a complete
system of spheres, i.e. that $\cals_1\in\Gamma_\oldOmega$. We have 
$$\wt\nolimits_\oldOmega(\cals)-\wt\nolimits_\oldOmega(\cals_1)=
\wt\nolimits_\oldOmega(D_2)-\wt\nolimits_\oldOmega(D)>0$$
since $\wt_\oldOmega D\le1$ and
$\wt_\oldOmega D_2\ge2$.
%\wt
 Hence
$\wt_\oldOmega\cals_1<\wt_\oldOmega\cals$, so that 
$\mu(\cals_1)<\mu(\cals)$. This again contradicts the $\mu$-minimality of
$\cals$, and so  the lemma is proved.
\EndProof

\Number\label{yup they are}Let $\Mh$ be a closed, orientable
hyperbolic $3$-orbifold, and set $\oldOmega=(\Mh)\pl$. If  $\cals$ is a $\mu$-minimal element of
$\Gamma_\oldOmega$, and if we set $\oldPsi=\oldOmega\cut{
  \obd(\cals)}$, then it follows from Lemmas \ref{why admissible} and \ref{oops lemma} that
the components of
 $\oldPsi$ are strongly \simple\ and boundary-irreducible (see  \ref{what, no soap?}); in particular, $ \oldSigma(\oldPsi)$, $ \oldPhi(\oldPsi)$, and 
 $\frakA(\oldPsi)$ are defined in view of Definition \ref{oldSigma def} and \ref{tuesa day}. 
These objects will be referred to freely in the lemmas in this section.
According to the discussion in \ref{tuesa day}, every component of $\frakA(\oldPsi)$ is an orientable annular orbifold, and hence every component of $|\frakA(\oldPsi)|$
is a disk or annulus. 
\EndNumber
%oldTheta

\Lemma\label{six of one} Let $\Mh$ be a closed, orientable
hyperbolic $3$-orbifold, and set $\oldOmega=(\Mh)\pl$. Let $\cals$ be a $\mu$-minimal element of
$\Gamma_\oldOmega$.  Set $\oldPsi=\oldOmega\cut {\obd(\cals)}$,
$N=M\cut\cals =|\oldPsi| $, and $\tcals=\partial N$.  Let $A$ be a weight-$0$
annulus properly embedded in $N$. Then 
the two components of $\partial A$ have the same size (see \ref{size
  def}) in $\obd(\tcals)$.
\EndLemma

\Proof
Let $\tC$ and $\tC'$ denote the boundary components of $A$, and set $s=\size \tC$ and $s'=\size \tC'$. By symmetry it is enough to prove that $s\le s'$. Let $\tS$ and $\tS'$ denote the components of
 $\tcals$ containing $\tC$ and $\tC'$ respectively. (It is possible
 that $\tS'=\tS$ or that $\tS'=\tau_\cals(\tS)$.)
% this does not affect
 %the argument.)
%\redcomment{It will, slightly. I think the issue is that $\tC$ could lie
  %inside the smaller disk bounded by $\tC'$ in $\tcals$, in which case
  %embeddedness of $\tD$ would seem to fail. But then the inequality
  %$s\le s'$, which is what I'm actually proving here, seems to be
  %trivial. I should double-check this. If I'm right, a very small
  %reorganization is needed: I should begin by saying that I'm proving
  %$s\le s'$, and that this is enough by symmetry. Then I'll have a
  %hard case (the one that's written down now) and an easy case.}
%\redcomment{Maybe so, but it still confuses me.}) 
Then $\tC'$ bounds a disk
 $\Delta\subset \tS'$ such that $\wt_\oldPsi \Delta=s'$. 

Consider first the case in which $\tC\subset\Delta$ (so that in particular $\tS=\tS'$). In this case, $\tC$ bounds a subdisk of $\Delta$, which must have weight at most $\wt\Delta=s'$; it then follows from the definition of size that $s=\size\tC\le s'$, as required. 

For the rest of the proof, we assume that $\tC\not\subset\Delta$. Since $\tC$ and $\tC'$ are the boundary components of the properly embedded annulus $A\subset N$, they are disjoint. Hence we must have $\tC\cap\Delta=\emptyset$. It follows that $A\cup\Delta$ is a disk contained in $N$ whose
 boundary is $\tC$.
Since $\wt A=0$, we have $\wt(A\cup\Delta)=\wt\Delta=s'$.
%hypothesis we have
 %$A\cap\fraks_\oldPsi=\emptyset$. Hence the disk $A\cup \Delta$
%, whose
 %boundary is $\tC$, 
%meets $\fraks_\oldPsi$ in $s'$ points. 
Let
 $\tD$ be a disk, properly embedded in $N$ and transverse to $\fraks_\oldPsi$, obtained from $A\cup \Delta$ by a
 small non-ambient isotopy, such that $\partial \tD=\tC$ and
 $\wt_\oldPsi\tD=s'$.  If we set $\rho=\rho_\cals$,
 $S=\rho(\tS)$, $D=\rho(\tD)$ and $C=\rho(\tC)$, then 
we have
$D\cap\cals=\partial D =C\subset S$. Let $D_1$ and $D_2$ denote the
closures of the components of $S-C$, and for $i=1,2$ set
$S_i=D\cup D_i$ and $\cals_i=(\cals-S)\cup S_i$. Then according to Corollary \ref{taser}, either $\cals_1$
or $\cals_2$ is a complete system of spheres in $M$. After possibly
re-indexing, we may assume that $\cals_1$ is a complete system. Since
$\cals_1$ is clearly transverse to $\fraks_\oldOmega$, we have
$\cals_1\in\Gamma_\oldOmega$.

The $\mu$-minimality of $\cals$ implies that $\mu(\cals_1)\ge
\mu(\cals)$. Hence, in particular, $\wt_\oldOmega\cals_1\ge
\wt_\oldOmega\cals$. Since $\cals_1=(\cals-S)\cup S_1$, it
follows that $\wt_\oldOmega S_1\ge
\wt S$. Since $S=D_1\cup D_2$ and $S_1=D\cup D_1$,
this implies that 
\Equation\label{toothloose}
\wt\nolimits_\oldOmega D_2\le
\wt\nolimits_\oldOmega D=\wt\nolimits_\oldPsi\tD=s'.
\EndEquation

On the other hand, we have
$s=\size C=\min(\wt_\oldOmega(D_1),\wt_\oldOmega(D_2))$, so that
in particular $\wt_\oldOmega D_2\ge s$. With
(\ref{toothloose}) this implies that $s\le s'$, as required. 
%The same argument,
%with the roles of $\tC$ and $\tC'$ reversed, shows that $s'\le s$, and the
%proof is complete.
\EndProof
%oldTheta\oldDeltaE\tS

\abstractcomment{\tiny I had written: ``Question about the proof
  above: Does it matter whether $\rho(\tC)$ and $\rho(\tC)$ lie in the
  same component of $\cals$?'' Apparently it does not.}

\Lemma\label{disks are ok}
Let $\Mh$ be a closed, orientable hyperbolic $3$-orbifold,  and set $\oldOmega=(\Mh)\pl$. Let
$\cals$ be a $\mu$-minimal element of  $\Gamma_\oldOmega$.  Set
$\oldPsi=\oldOmega\cut{
\obd(\cals)}$, $N=|\oldPhi|$, and $\tcals=\partial N$.  Then
for
every properly embedded disk $D$ in $N$, transverse to $\fraks_\oldPsi$, such that $\obd(D)$ is
annular and $\pi_1$-injective in $\oldPsi$,
%The stronger
%condition seems to be needed for the proof of  Proposition\ref{dandy
  %subsystem}. It should not affect other refs. to the def., but I
%should check this. Of course I need a slight rewording of the proof of Proposition \ref{semidandies exist}, and
%I'll have to say a word when \dandy ness is applied later in the paper
there is a disk $E\subset\tcals$ such that $\obd(E)$ is an annular orbifold and $\partial E=\partial D$.
\EndLemma

\Proof
Set $M=|\oldOmega|$, $N=|\oldPsi|$, and $\rho=\rho_\cals:N\to M$. Let
$\tS$ denote the component of $\tcals$ containing $\partial D$, and
set $S=\rho(\tS)$.
Since $D$ is a disk and $\obd(D)$ is annular (and orientable), we have $\wt_\oldPsi D=2$, and  both points of $D\cap\fraks_\oldPsi$ are of order $2$. Hence $\wt_\oldOmega\rho( D)=2$, and both points of $\rho(D)\cap\fraks_\oldOmega$ are of order $2$. Let $\tE_1$ and $\tE_2$ denote the
closures of the components of $\tS-\partial  D$, and for $i=1,2$ set $E_i=\rho(\tE_i)$,
$S_i=\rho(D)\cup E_i$ and $\cals_i=(\cals-S)\cup S_i$. Applying
Corollary \ref{taser}, with $\rho(D)$ playing the role of $D$ in the corollary, we deduce that either $\cals_1$
or $\cals_2$ is a complete system of spheres in $M$. After possibly
re-indexing we may assume that $\cals_1$ is a complete system of
spheres. Since 
$\cals_1$ is clearly transverse to 
$\fraks_\oldOmega$, we have  $\cals_1\in\Gamma_\oldOmega$. The
 $\mu$-minimality of $\cals$ therefore implies that
$\mu(\cals_1)\ge \mu(\cals)$. In particular we have
$\wt_\oldOmega\cals_1\ge 
\wt_\oldOmega\cals$, or equivalently  
$\wt_\oldPsi D=\wt_\oldOmega \rho(D)\ge 
\wt_\oldOmega E_2$.
Hence $\wt_\oldOmega E_2\le2$. If $\wt_\oldOmega E_2\le1$ then $E_2$ is discal; but since $\obd(D)$ is $\pi_1$-injective in $\oldPsi$ by hypothesis, $\partial\obd(D)=\partial\obd(E_2)$ cannot bound a discal suborbifold of $\oldPsi$. 
%But by \redmissingref{I don't think it's true that $\wt_\oldOmega E_2$ cannot be less than $2$. I think there should be another alternative in the conclusion, that $E$ is discal, and the apps. have to be adjusted. I had written ``Or do I need something like $\pi_1$-injectivity of $\obd(D)$ in $\oldOmega$?'' which I think is an equivalent modification of the statement},
%$\wt_\oldOmega E_2$ cannot be less than $2$. 
We must therefore have $\wt_\oldOmega E_2=2$. It follows that
$\wt_\oldOmega\cals_1= 
\wt_\oldOmega\cals$; and since we clearly have $\compnum(\cals_1)=\compnum(\cals)$, the inequality
$\mu(\cals_1)\ge \mu(\cals)$ implies that 
$\nu(\cals_1)\ge \nu(\cals_1)$. This means that $\fraks_{\rho(D)}$ has
at least as many points of order greater than $2$ as $\fraks_{E_2}$
has. Since both points of 
$\fraks_{\rho(D)}$
 are of order  $2$, both points of $\fraks_{E_2}$ are also of order
 $2$. Thus $\obd(E_2)$ is annular, so that $\obd(\tE_2)$ is annular, and the conclusion of the lemma is true with $E=\tE_2$. 
\EndProof
%oldThetaD_}

\abstractcomment{\tiny I have removed subsection ``bookkeeping'' and Lemma
``yes it do''. They can be found in not-doubly2.tex.  

 I am removing Lemmas ``no she don't'' and ``twerpentine.'' They can be found in the file doubleclash1.tex. I am also removing Lemma ``snooing the bedoodly out of it'' because I no longer believe it.}

\Lemma\label{fluorescent foo-foo}
Let $\Mh$ be a closed, orientable hyperbolic $3$-orbifold, and set $\oldOmega=(\Mh)\pl$. Set
$M=|\oldOmega|$, and let $\cals\subset M$ be a $\mu$-minimal element of
$\Gamma_\oldOmega$. Set $\oldPsi=\oldOmega\cut{
\obd(\cals)}$,  $N=|\oldPsi|=M\cut\cals$, and $\tcals=\partial N$.
%, and $\tau=\tau_\cals$.
Let $F$ be a component
% of $\tcals$ such that $\tS$ contains no
%bad component 
of $|\oldPhi(\oldPsi)|$ which splinters $\obd(\tS)$, where $\tS$ denotes the component $\tS$ of $\tcals$ containing
$F$. Then either
%. Let
%$F$ be a component of $|\oldPhi(\oldPsi)|$ contained in $\tS$. Suppose that
$F$ is a pseudo-belt for $\tS$, or $F$ is \bad.
\EndLemma

\Proof
Set  $\tau=\tau_\cals$,  $\rho=\rho_\cals$, and $S=\rho(\tS)$.
Let $\frakU$ denote the component of $\oldSigma(\oldPsi)$ such that $F\subset|\frakU|$.

Since $F$ splinters $\obd(\tS)$, and $\obd(\tS)$ is negative by \ref{doublemint}, it follows from Lemma \ref{negative splinter} that $\chi(\obd(F)<0$. Hence by Lemma \ref{when a tore a fold}, the connected \Ssuborbifold\ $\frakU$ of $\oldPsi$ must be \pagelike. It follows (cf. \ref{S-pair def} that $|\frakU|\cap\tcals$ has at
most two components. 

If $|\frakU|\cap\tcals$ is connected, so that $|\frakU|\cap\tcals=F$, then by definition $F$
is a pseudo-belt for $\tS$; thus the
conclusion of the lemma holds if  $|\frakU|\cap\tcals$ is connected. For
the rest of the argument we shall assume that $|\frakU|\cap\tcals$ has
exactly two components, (i.e. $\frakU$ is untwisted, cf. \ref{S-pair def}) and we shall complete the proof by showing
that, under this assumption, $F$ is ungainly. Let $F'$ denote the component of $|\frakU|\cap\tcals$  distinct from $F$, and let $\tS'$ denote the component of $\tcals$ containing $F'$.

Since $(\frakU,\obd(F\discup F')$ is an untwisted \pagelike\ \spair\ (with $\frakU$ connected), the
components of $\Fr_N|\frakU|$ are weight-$0$ annuli 
$A_1,\ldots,A_n$, where a priori we have $n\ge0$; and  each $A_i$ has one boundary curve, $C_i$, in $F$, and one, $C_i'$, in $F'$. 
We have $\partial F=C_1\cup\cdots\cup C_n$ and $\partial
F'=C_1'\cup\cdots\cup C_n'$. For $i=1,\ldots,n$, let $\Delta_i$ denote
the component of $\tS-\inter F$ bounded by $C_i$, and let $\Delta_i'$ denote
the component of $\tS'-\inter F'$ bounded by $C_i'$.

For $i=1,\ldots,n$, set $k_i=\wt_\oldPsi\Delta_i$ and
$k_i'=\wt_\oldPsi\Delta_i'$. The hypothesis that $F$ splinters
$\obd(\tS)$ means that $k_{i}<\wt_\oldPsi (\tS)/2$ for $i=1,\ldots,n$, i.e.
\Equation\label{forget the alamo}
k_{i}<\wt\nolimits_\oldOmega (S)/2 \qquad\text{for }
i=1,\ldots,n.
\EndEquation 

 Since $\wt A_i=0$,  the sphere $\Delta_i\cup A_i\cup\Delta_i'$ meets 
$\fraks_\oldPsi$ transversally in $k_i+k_i'$ points.

For each $i$, since $A_i$ is a weight-$0$ annulus, Lemma \ref{six of one} implies
that the boundary curves $C_i$ and $C_i'$ of $A_i$ have the same size. 
%Let us denote their common size by $s_i$. 

Since the components of $\tS-C_i$ are $\inter\Delta_i$ and $\tS-\Delta_i$, the definition of size gives
$$
\size C_i=\min(\wt\nolimits_\oldPsi(\Delta_i),\wt\nolimits_\oldPsi(\tS-\Delta_i))=\min(k_i, \wt\nolimits_\oldOmega(S)-k_i),
$$
which with (\ref{forget the alamo}) implies
\Equation\label{peel me a grape}
\size C_i=k_i.
\EndEquation
The definition of size also gives
\Equation\label{two ships}
\size C_i'=\min(\wt\nolimits_\oldPsi(\Delta_i'),\wt\nolimits_\oldPsi(\tS'-\Delta_i'))=\min(k_i', \wt\nolimits_\oldOmega(S')-k_i').
\EndEquation
Since $\size C_i=\size C_i'$, it follows from (\ref{peel me a grape}) and (\ref{two ships}) that
\Equation\label{consort}
k_i=\min(k_i', \wt\nolimits_\oldOmega(S')-k_i').
\EndEquation

Now note  that Since $(\frakU,\obd(F\discup F')$ is an untwisted \pagelike\ \spair\ (with $\frakU$ connected), the
orbifolds $\obd(F)$ and $\obd(F')$ are homeomorphic, and hence
 $\wt_\oldPsi F=\wt_\oldPsi F'$; let $\ell$ denote the common value of  $\wt_\oldPsi F$ and $\wt_\oldPsi F'$. 
Since $\Delta_1',\ldots,\Delta_n'$ are the components of $\tS'-\inter F'$, we have
\Equation\label{i yam that i yam}
\wt\nolimits_\oldOmega S'=\wt\nolimits_\oldPsi \tS'=\ell+k_1'+\cdots k_n'.
\EndEquation
 Similarly,
\Equation\label{you yam that you yam}
\wt\nolimits_\oldOmega S=\ell+k_1+\cdots k_n.
\EndEquation

The proof will be divided into three cases. It will turn out that the
first two cases cannot occur, and in the third case it will be proved
that $F$ is ungainly, as required.

{\bf Case I.} The sphere $\tS'$ is distinct from both $\tS$ and $\tau(\tS)$.

Let $\tY$ be a small regular neighborhood of $|\frakU|\cup\tS \cup\tS'$ in $N$.
The boundary components of $\tY$ are $\tS$, $\tS'$, and spheres $\tT_1,\ldots,\tT_n$, where $\tT_i$ is a nearby parallel copy of $\Delta_i\cup A_i\cup\Delta_i'$. We may choose $\tY$ in such a way that 
\Equation\label{how to pick 'em}
\wt\nolimits_\oldPsi\tT_i=\wt\nolimits_\oldPsi(\Delta_i\cup A_i\cup\Delta_i')=k_i+k_i'
\EndEquation
for $i=1,\ldots, n$.

%T

Since $(\frakU,\obd(F\discup F')$ is an untwisted \pagelike\ \spair, we have a homeomorphism from $|\frakU|$ to $F\times[0,1]$ that maps $F$ and $F'$ to $F\times\{0\}$ and $F'\times\{1\}$ respectively. This extends to a homeomorphism from $|\frakU|\cup\tS \cup\tS'$ to $(F\times[0,1])\cup(\tS\times\{0,1\})\subset \tS\times[0,1]$. Hence 
$\tY$ is a $3$-sphere-with holes. 

Set $\rho=\rho_\cals$.
Since $\tS'$ is distinct from $\tau(\tS)$ as well as from $\tS$, the spheres $S:=\rho(\tS)$ and $S':=\rho(\tS')$ are distinct components of $\cals$, and $\rho$ maps $\tY$ homeomorphically onto a $3$-sphere-with holes $Y\subset M$ whose boundary components are $S$, $S'$, and $T_1,\ldots,T_n$, where $T_i=\rho(\tT_i)$. The hypotheses of 
Proposition \ref{shmazer} now hold if we make this choice of $Y$ and
take $\cals$ as above. In the notation of Proposition \ref{shmazer},
we have $\cals_0=S\cup S'$ and $\calt=T_1\cup\cdots\cup T_n$. It now
follows from Proposition \ref{shmazer} that there is an index
$i\in\{1,\ldots,n\}$ such that $T_i$ involves the component $S'$ of $\cals$. Hence $n\ge1$, and after possibly re-indexing we may assume that $T_n$ involves $S'$.
Hence by Proposition \ref{just do it}, $(\cals-S')\cup T_n$ is a complete system of spheres in $M$. As
$(\cals-S')\cup T_n$ is clearly transverse to $\fraks_\oldOmega$, it is an element of $\Gamma_\oldOmega$. 
The $\mu$-minimality of $\cals$ therefore implies that $\mu((\cals-S')\cup T_n )\ge\mu(\cals)$. In particular we have  
$\wt_\oldOmega((\cals-S')\cup T_n )\ge
\wt_\oldOmega\cals$. This means that
\Equation\label{getting there}
\wt\nolimits_\oldOmega T_n \ge
\wt\nolimits_\oldOmega S'.
\EndEquation

To interpret the two sides of (\ref{getting there}), first note that $\wt_\oldOmega T_n=\wt_\oldPsi\tT_n=k_n+k_n'$ by (\ref{how to pick 'em}), and
$\wt_\oldOmega S'=\ell+k_1'+\cdots k_n'$ by (\ref{i yam that i yam}).
Thus (\ref{getting there}) gives
$k_n+k_n'\ge \ell+k_1'+\cdots k_n'$, or
\Equation\label{roller coaster}
k_n\ge \ell+k_1'+\cdots k_{n-1}'.
\EndEquation

Now (\ref{consort}) implies that $k_i'\ge k_i$ for each $i$, which with 
(\ref{roller coaster}) gives
\Equation\label{oh boy}
k_n\ge \ell+k_1+\cdots k_{n-1}.
\EndEquation
But by (\ref{you yam that you yam}), we have
$\wt_\oldOmega S=\ell+k_1+\cdots k_n$.
Hence (\ref{oh boy}) may be rewritten as $k_n\ge\wt_\oldOmega( S)-k_n$, i.e. $k_n\ge\wt_\oldOmega (S)/2$. This contradicts (\ref{forget the
  alamo}), and we have shown that Case I cannot occur.
%But according  to (\ref{not a prime}), for each $i\in\{1,\ldots,n-1\}$ we have either 
%$s_i=k_i$ or $s_i= \wt_\oldOmega (S)-k_i\le k_i$. Consider the subcase in which 
%$s_i=k_i$ for every $i\in\{1,\ldots,n-1\}$. In this subcase, it follows from (\ref{oh boy}) that 
%$$k_n\ge \ell+k_1+\cdots k_{n-1}=\wt\nolimits_\oldPsi(\tS)-k_n %=\wt\nolimits_\oldOmega(S)-k_n,$$
%so that $k_n\ge\wt_\oldOmega (S)/2$; this contradicts (\ref{forget the
  %alamo}).
% implies that $F$ does

%There remains the subcase in which
%$s_{i_0}= \wt_\oldOmega(S)-k_{i_0}\le k_{i_0}$ for some $i_0\in\{1,\ldots,n-1\}$. In this case, we have 
%$k_{i_0}\ge\wt_\oldOmega(S)/2$, and again we have a contradiction to (\ref{forget the
  %alamo}). Thus we have shown that Case I cannot occur.

{\bf Case II.} We have $\tS'=\tS$.

In this case, $F$ and $F'$ are disjoint subsurfaces of $\tS$, and hence each of them has a non-empty boundary; thus $n>0$. Since $F$ and $F'$ are connected and $\tS$ is a $2$-sphere, there are disjoint disks $D,D'\subset\tS$, containing $F$ and $F'$ respectively, such that $\partial D$ and $\partial D'$ are components of $\partial F$ and $\partial F'$ respectively. After possibly re-indexing we may assume that $\partial D=C_n$. We have $\partial D'=C'_m$ for some $m$ with $1\le m\le n$.

Note that 
$$D=F\cup\bigcup_{1\le i<n}\Delta_{i} \qquad \text{and}\qquad D'=F'\cup\bigcup_{i\ne m}\Delta_{i}',$$
so that
$$\wt\nolimits_\oldPsi D=\ell+\sum_{1\le i<n}k_{i} \qquad \text{and}\qquad\wt\nolimits_\oldPsi D'=\ell+\sum_{i\ne m}
k_{i}'. $$
We also have $\Delta'_m=\tS-\inter D'\supset D$,  so that $k'_m\ge \wt_\oldPsi D$; similarly, we have $k_n\ge \wt_\oldPsi D'$. Hence
\Equation\label{you're going crazy}
k'_m\ge\ell+\sum_{1\le i<n}k_{i}\qquad \text{and}\qquad
k_n\ge\ell+\sum_{i\ne m}k_{i}'
.
\EndEquation

The second inequality of (\ref{you're going crazy}) implies that for every $i\ne m$ we have $k'_{i}\le k_n$. Since $k_n<\wt_\oldOmega(S)/2$
%=\wt_\oldOmega(\tS)/2$ 
by (\ref{forget the alamo}), and since we are in the case where $\tS=\tS'$ and hence $S=S'$, it follows that $k_i'<\wt_\oldOmega(S')/2$ for every $i\ne m$. In view of (\ref{consort}) it follows that
\Equation\label{buy a quack}
k_i'=k_i\qquad\text{for every }i\ne m.
\EndEquation

Consider the subcase in which $m< n$, which implies in particular that
$n>1$. In this subcase, the sum on the right hand side of the second
inequality of (\ref{you're going crazy}) includes the term $k_n'$,
which by (\ref{buy a quack}) is equal to $k_n$ since $n\ne m$. Hence
the inequality implies that $\ell=0$ and that $k_i'=0$ whenever $i\ne
m,n$. Since 
%$0$, we have 
%$F\cap\fraks_\oldPsi=\emptyset$, i.e. 
$\wt F=\ell=0$, and since we have
assumed that $\chi(\obd(F)<0$, we must have
$n\ge3$. Thus there is an index $j\ne m, n$ in $\{1,\ldots,n\}$, and
we have $k_j'=0$. But this is impossible, since the component $A_j$  of $|\frakA(\oldPsi)|$ is an essential  annular orbifold in $\oldPsi$ by \ref{tuesa
  day}).

Now consider the subcase in which $m=n$. Then for $1\le i<n$ we have $k_i'=k_i$ by (\ref{buy a quack}). Thus the second inequality of (\ref{you're going crazy}) becomes $k_n\ge\ell+\sum_{1\le i<n}k_{i}$. But by \ref{you yam that you yam} we have
$\wt_\oldOmega S
%=\wt_\oldPsi \tS
=\ell+k_1+\cdots k_n$. Hence $k_n\ge \wt_\oldOmega(S)/2$, a contradiction to (\ref{forget the alamo}). Thus Case II cannot occur.

{\bf Case III.} We have $\tS'=\tau(\tS)$. 

In this case, according to the definition of $\tS'$, we have
$F'\subset \tau(\tS)$, so that $|\frakU|$ meets $\tau(\tS)$.  Since $F$ splinters $\obd(\tS)$ by hypothesis, it now follows from the definition (see
\ref{stop beeping}) that $F$ is \bad, as asserted by the lemma.
\EndProof
%oldTheta %\oldDelta\tau\Delta\tDelta $t

\abstractcomment{\tiny For the summary of the old version of the proof above, see doubleclash1.tex.}

\Proof[Proof of Proposition \ref{semidandies exist}] Let $\Mh$ be a
closed, orientable hyperbolic $3$-orbifold, and set $\oldOmega=(\Mh)\pl$.
Suppose that $\cals$ is a $\mu$-minimal element 
of $\Gamma_\oldOmega$.
According to Lemma \ref{why
  admissible}, $\cals$ is an $\oldOmega$-admissible system of spheres. Lemma 
\ref{six of one} asserts that $\cals$ satisfies the first bulleted
condition in the definition (\ref{semidandy def}) of a \dandy\ system,
 Lemma \ref{fluorescent foo-foo}
%\ref{and the foos just keep on fooing} 
asserts that $\cals$ satisfies the second bulleted condition, and Lemma \ref{disks are ok} asserts that $\cals$ satisfies the third bulleted condition. This proves the first assertion of the proposition. To prove the second assertion, we need only note that $\Gamma_\oldOmega$ has a $\mu$-minimal element  $\cals$ according to Lemma \ref{it is what it is}, that $\cals$ is a
complete system of spheres in $|\oldOmega|$
 according to the definition of $\Gamma_\oldOmega$, and that $\cals$ is \dandy\ by the first assertion of the present proposition.
\EndProof

\section{Spheres of incompatibility and clash components}\label{clash section}

\begin{definitionremark}\label{squerkimer}
Let $\Mh$ be a closed,
orientable, hyperbolic $3$-orbifold, and set $\oldOmega=(\Mh)\pl$. Let $\cals$ be an
$\oldOmega$-admissible system of spheres in $M:=|\oldOmega|$. 
Set
$\rho=\rho_\cals$ and $\tau=\tau_\cals$,
 %and $\oldPsi=\oldOmega\cut{ \obd(\cals)}$.
Let $S$ be a   component of $\cals$, and set
$X''=X''(\oldOmega,\obd(\cals),\obd(\cals-S))$ (in the notation of
\ref{oldXi}).  
We will call $S$ an {\it $\oldOmega$-\clashsphere} for $\cals$  (or simply
a \clashsphere\ for $\cals$ when it is clear which orbifold is involved) if we have
$\chibar(X''\wedge
[\tau](X''))<
\chibar(X'')$ (in the notation of Chapter \ref{higher chapter}).

Note that according to the definitions in \ref{oldXi} we have $X''=X\wedge[\tcals_0]$, where $X''=X''(\oldOmega,\obd(\cals),\obd(\cals-S))$ and $\tcals_0=\rho^{-1}(S)$. Since $\tau(\tcals_0)=\tcals_0$, we have $X''\wedge\tau(X'')=(X\wedge[\tcals_0])\wedge (\tau(X)\wedge[\tcals_0])=X \wedge[\tcals_0]\wedge \tau(X)=X''\wedge\tau(X)$. Hence the inequality in the above definition may be rewritten in the equivalent form $\chibar(X''\wedge
[\tau](X))<
\chibar(X'')$.
\end{definitionremark}

\Number\label{dots trrivial}
Note that if $\Mh$ is a closed,
orientable, hyperbolic $3$-orbifold, if we set $\oldOmega=(\Mh)\pl$, and if $\cals\subset |\oldOmega|$ is an
$\oldOmega$-admissible system of spheres which is connected,
i.e. consists of a single sphere $S=\cals$, then $S$ cannot be a
\clashsphere\ for $\cals$; for in this
case we have
$X'':=X''(\oldOmega,\obd(\cals),\obd(\cals-S))\preceq X'(\oldOmega,\obd(\cals),\obd(\cals-S))=\emptyset$, so that $\chibar(X''\wedge
[\tau](X''))$ and
$\chibar(X'')$ are both equal to $0$.
\EndNumber

\abstractcomment{\tiny It said ``This def. of \clashsphere\ may work, but it makes me
  uncomfortable that it is not made clear how it is related to the
  def. of an \eclipsing\ sphere. The reason that an \eclipsing\ sphere
  may fail to be a \clashsphere\ is that it may not be
  \doublefull. The reason that a \clashsphere\ 
  may fail to be  an \eclipsing\ sphere is that it may not have
  strictly greater weight than the other spheres involved in the
  'clash.''' This will be obsolete because I plan to remove the
  def. of \eclipsing\ sphere, possibly replacing it with something else.}

\Proposition\label{fine and dandies exist} Let $\Mh$ be a closed,
orientable hyperbolic $3$-orbifold, and set $\oldOmega=(\Mh)\pl$. Set $M=|\oldOmega|$.
Let $\cals$  be an $\oldOmega$-admissible system of spheres in $M$, and let $S$
be a \doublefull\ component of $\cals$ which is an $\oldOmega$-\clashsphere\ for $\cals$. Set
$\cals':=\cals-S$. Then
$$s'_{{\oldOmega}}(
M-\cals'
%X^*
)\ge\lambda_\oldOmega+
%\sum_{Z\in\calc(X^*-S)}
s'_{\oldOmega}(M-\cals
%X^*-S
)$$ 
(see \ref{it's ok}).
\EndProposition

\Proof
%Set $\lambda=\lambda_\oldOmega$. According to \ref{lambda thing}, we have %$\lambda_\oldPsi=\lambda_{\oldPsi'}=\lambda$. 
Set $\oldTheta=\obd(\cals)$. Since $\cals$ is in particular an $\oldOmega$-admissible system of spheres, $\oldTheta$ is an incompressible suborbifold of $\oldOmega$. Set $\oldTheta'=\obd(\cals')$, so that $\oldTheta'$ is a union of components of $\oldTheta$. Set
$\oldPsi=\oldOmega\cut{\oldTheta}$ and ${\toldTheta}=\partial\oldPsi$,
$\oldPsi'=\oldOmega\cut{\oldTheta'}$, and ${\toldTheta}'=\partial\oldPsi'$, so that
${\toldTheta}'$ is canonically identified with a union of components of ${\toldTheta}$. Set ${\toldTheta}''={\toldTheta}-{\toldTheta}'=\rho_\oldTheta^{-1}(\obd(S))$.
Set $\tau=\tau_\oldTheta$.  It follows from Corollary \ref{burpollary}  that
 $[\oldPhi^-(\oldPsi')]\preceq[\oldPhi^-(\oldPsi)]$. Hence
we may suppose $\oldPhi^-(\oldPsi')$ to be chosen within its orbifold-isotopy class in such a way that $\oldPhi^-(\oldPsi')\subset\inter(\oldPhi^-(\oldPsi) \cap{\toldTheta'})
$. 

Set $\lambda=\lambda_\oldOmega$. According to \ref{lambda thing}, we have $\lambda_\oldPsi=\lambda_{\oldPsi'}=\lambda$. 

Set $X=X(\oldOmega,\oldTheta,\oldTheta')$, $X'=X'(\oldOmega,\oldTheta,\oldTheta')$ and $X''=X''(\oldOmega,\oldTheta,\oldTheta')$. 
By hypothesis $S$ is an $\oldOmega$-\clashsphere\ for $\cals$, which by
\ref{squerkimer} means that
$\chibar(X''\wedge
[\tau](X))<
\chibar(X'')$.
%[\tau](X''))<
%\chibar(X'')$. 
It therefore follows from Lemma \ref{veggie prop} 
that $\chibar(X')<\chibar(\oldPhi^-(\oldPsi)\cap{\toldTheta}')$.
Since Proposition \ref{herman had burped} asserts that
$[\oldPhi^-(\oldPsi')]\preceq
X'$, we have
$\chibar(\oldPhi^-(\oldPsi'))\preceq\chibar(
X')$ by Lemma \ref{old partial order}, and
 it now follows that 
\Equation\label{lovely display}
\chibar(\oldPhi^-(\oldPsi'))<\chibar(\oldPhi^-(\oldPsi) \cap{\toldTheta}').
\EndEquation
Set $\oldGamma=(\oldPhi^-(\oldPsi) \cap\toldTheta')-\inter\oldPhi^-(\oldPsi')$. Then (\ref{lovely display}) gives $\chibar(\oldGamma)=\chibar(\oldPhi^-(\oldPsi) \cap{\toldTheta}')-\chibar(\oldPhi^-(\oldPsi'))>0$, i.e. 
\Equation\label{still less than nothing}
\chi(\oldGamma)<0.
\EndEquation
It  follows from Proposition \ref{fortunately} (and from the componentwise strong \simple ity and
  componentwise boundary-irreducibility of $\oldPsi$, see \ref{doublemint}) that the number of points of order $2$ in $\fraks_{{\toldTheta}-\inter\oldPhi(\oldPsi)}$, which we shall denote by $a$, is divisible by $\lambda_\oldPsi=\lambda$. The same proposition, %\redmissingref{or ``the same lemma'' if things are changed as I'm suggesting}, 
applied with $\oldPsi'$ in place of $\oldPsi$, shows that the number of points of order $2$ in 
$\fraks_{{\toldTheta}'-\inter\oldPhi(\oldPsi')}$, which we shall denote by $a'$, is divisible by $\lambda_{\oldPsi'}=\lambda$. 
On the other hand, since $S$ is \doublefull\ by hypothesis, it follows from Remark \ref{special ouenelitte} that  every component of $|\overline{\toldTheta''\setminus\oldPhi(\oldPsi)}|$ is either a weight-$0$ annulus, or a disk containing exactly two points of $\fraks_\oldPsi$, each of order $2$. Hence the number of points of order $2$ in $\fraks_{\overline{{\toldTheta}''\setminus\oldPhi(\oldPsi)}}$, which we shall denote by $b$, is even, and in particular is divisible by $\lambda\in\{1,2\}$. But the definitions of $\oldGamma$ and $\oldTheta''$ imply that
\Equation\label{is this what you wanted?}
({\toldTheta}'-\inter\oldPhi(\oldPsi'))\cup \overline{{\toldTheta}''\setminus\oldPhi(\oldPsi)}=({\toldTheta}-\inter\oldPhi(\oldPsi))\cup\oldGamma.
\EndEquation
Each side of (\ref{is this what you wanted?}) is exhibited as a union
of two orbifolds which meet in a (possibly empty) union of boundary
components, and these boundary components contain no singular
points. Hence the numbers of singular points of order $2$ in  the
orbifolds appearing on the left
and right hand sides of (\ref{is this what you wanted?})
are respectively equal to $a'+b$ and to $a+c$, where $c$ denotes the number of points of order $2$ in $\fraks_{\oldGamma}$. Thus $a'+b=a+c$, which implies that $\lambda|c$. We have shown:
\Claim\label{piled higher}
The number of points of order $2$ in $\fraks_{\oldGamma}$ is divisible by $\lambda$.
\EndClaim

Let $\oldGamma^-$ denote the union of all negative components of $\oldGamma$. It follows from (\ref{still less than nothing}) that 
\Equation\label{not really nothing}
\oldGamma^-\ne\emptyset.
\EndEquation
Each component of $\oldGamma-\oldGamma^-$ has Euler characteristic $0$, and is therefore either an annulus containing no singular points, or a disk containing exactly two singular points, each of order $2$. In view of \ref{piled higher}, this implies:
\Claim\label{and deeper}
The number of points of order $2$ in $\fraks_{\oldGamma^-}$ is divisible by $\lambda$.
\EndClaim
%Since the \dandy\ system $\cals$ is in particular admissible,
%$\oldTheta$ is incompressible. 
Since by hypothesis $\oldOmega$
contains no  embedded negative turnovers, no component of
$\oldGamma^-$ is a negative turnover. The definition of $\oldGamma^-$ implies that
$\oldGamma^-$ is negative. These observations show that $\oldGamma^-$  the hypotheses of the first assertion of Proposition \ref{at least a sixth} hold with $\oldGamma^-$ in place of $\oldGamma$; by \ref{and deeper}, if $\lambda=2$, the hypotheses of the second assertion also hold with $\oldGamma^-$ in place of $\oldGamma$. Hence 
$\compnum(|\oldGamma^-|)\le6\chibar(\oldGamma^- )$, and if $\lambda=2$ we have
 $\compnum(|\oldGamma^-|)\le2\lfloor
3\chibar(\oldGamma^- )\rfloor$.
But according to (\ref{not really nothing}) we have $\compnum(|\oldGamma^-|)>0$. Hence $\chibar(\oldGamma^-)$ is bounded below by $1/6$, and by $1/3$ if $\lambda=2$; that is,
\Equation\label{THAT's what I wanted!}
\chibar(\oldGamma^-)\ge\lambda/6.
\EndEquation

We have observed that each side of (\ref{is this what you wanted?}) is exhibited as a union of two orbifolds which meet at most in a (possibly empty) union of boundary components. Hence
\Equation\label{crooks}
\chibar({\toldTheta}'-\inter\oldPhi(\oldPsi'))+\chibar(\overline{{\toldTheta}''\setminus\oldPhi(\oldPsi)})=\chibar({\toldTheta}-\inter\oldPhi(\oldPsi))+\chibar(\oldGamma).
\EndEquation
Since $S$ is \doublefull\ by hypothesis, it follows from \ref{special ouenelitte} that every component of $\overline{{\toldTheta}''\setminus\oldPhi(\oldPsi)}$ has Euler characteristic $0$. We have observed that every component of $\oldGamma-\oldGamma^-$ has Euler characteristic $0$, and hence $\chibar(\oldGamma)=\chibar(\oldGamma^-)$. Thus (\ref{crooks}) becomes
$\chibar({\toldTheta}'-\inter\oldPhi(\oldPsi'))=\chibar({\toldTheta}-\inter\oldPhi(\oldPsi))+\chibar(\oldGamma^-)$, which with (\ref{THAT's what I wanted!}) gives
\Equation\label{more crooks}
\chibar({\toldTheta}'-\inter\oldPhi(\oldPsi'))\ge\chibar({\toldTheta}-\inter\oldPhi(\oldPsi))+\lambda/6.
\EndEquation
Since
by \ref{tuesa day} we have
$2\chibar(\kish\oldPsi))=
2\chibar(\overline{(\oldPsi-\oldSigma(\oldPsi)})=
\chibar(\overline{{\toldTheta}-\oldPhi(\oldPsi)})$, and similarly $2\chibar(\kish\oldPsi'))=
%2\chibar(\overline{(\oldPsi-\oldSigma(\oldPsi)})=
\chibar(\overline{{\toldTheta}'-\oldPhi(\oldPsi')})$, it follows from (\ref{more crooks}) that
\Equation\label{crooks mit crooks}
\chibar(\kish(\oldPsi'))\ge\chibar(\kish(\oldPsi))+\lambda/12.
\EndEquation

Finally, according to the definitions in \ref{t-defs}, we have 
$s'_\oldOmega(M-\cals)=\sigma'(\oldPsi)=\lambda
%\bigg
\lfloor12\chibar(\kish(\oldPsi))/\lambda\rfloor$, and similarly
$s_\oldOmega(M-\cals')=\lambda
%\bigg
\lfloor12\chibar(\kish(\oldPsi'))/\lambda\rfloor$.
But it follows from (\ref{crooks mit crooks}) that
$12\chibar\kish(\oldPsi')/\lambda\ge (12\chibar\kish(\oldPsi')/\lambda)+1$,
which implies that 
$\lfloor12\chibar\kish(\oldPsi')/\lambda\rfloor\ge \lfloor12\chibar\kish(\oldPsi')/\lambda\rfloor+1$.
Hence
$s'_{{\oldOmega}}(
M-\cals'
)\ge
s'_{\oldOmega}(M-\cals
)+\lambda$, as required. 
\EndProof

\Corollary\label{toss 'em out}
Let 
$\Mh$ be a closed, orientable hyperbolic $3$-orbifold, and set $\oldOmega=(\Mh)\pl$ and 
$M=|\oldOmega|$. Let
$\cals$ be an $\oldOmega$-admissible
system of spheres in $M$. Then there is a subsystem $\cals^*$
of $\cals$ such that (i)
$s'_\oldOmega(M-\cals^*)\ge s'_\oldOmega(M-\cals)+\lambda_\oldOmega\compnum(\cals-\cals^*)$ (see \ref{it's ok}), and (ii) %$\cals^*$
%has
 no \doublefull\ component of $\cals^*$ is an $\oldOmega$-\clashsphere\ for $\cals^*$.
%\ which is a.
\EndCorollary

\Proof
Set $\lambda=\lambda_\oldOmega$.
Assume that $\cals$ has no subsystem $\cals^*$ satisfying (i) and
(ii). We shall construct an infinite sequence $\cals_0,\cals_1,\ldots$
of subsystems of $\cals$ such that for each $n\ge0$ we have
\Equation\label{compnum equation}
\compnum(\cals-\cals_n)=n \qquad\text{and}\qquad
s'_\oldOmega(M-\cals_n)\ge s'_\oldOmega(M-\cals)+\lambda n.
\EndEquation
The first condition of (\ref{compnum equation}) is obviously impossible when $n>\compnum(\cals)$. Hence the proof will be complete when the sequence $\cals_0,\cals_1,\ldots$ has been constructed.

The $\cals_n$ will be constructed recursively. We set $\cals_0=\cals$;
it is then obvious that (\ref{compnum equation}) holds for
$n=0$. Suppose that for a given $n$ the system $\cals_n$ has been
constructed and satisfies  (\ref{compnum equation}). It follows from
(\ref{compnum equation}) that $s'_\oldOmega(M-\cals_n)\ge
s'_\oldOmega(M-\cals)+\lambda\compnum(\cals-\cals_n)$. Thus Condition (i) of the
statement of the corollary holds with $\cals_n$ in place of $\cals^*$. 
Since we have assumed that $\cals$ has no subsystem satisfying Conditions (i) and
(ii), Condition (ii) cannot hold with $\cals_n$ in place of $\cals^*$;
that is, $\cals_n$ must have a \doublefull\ component $S_n$ which is a \clashsphere\ for $\cals_n$. Now set
$\cals_{n+1}=\cals_n-S_n$. Then $\cals_{n+1}$ is a subsystem of
$\cals_n$ and hence of $\cals$. We have
$\compnum(\cals-\cals_{n+1})=\compnum(\cals-\cals_n)+1=n+1$. On the other
hand, we may apply Proposition 
\ref
{fine and dandies exist}, with $\cals_n$ and $S_n$ playing the roles
of $\cals$ and $S$, to deduce that 
$s'_{{\oldOmega}}(
M-\cals_{n+1}
%X^*
)\ge\lambda+
%\sum_{Z\in\calc(X^*-S)}
s'_{\oldOmega}(M-\cals_n
%X^*-S_n
)$.
With (\ref{compnum equation}), this implies that 
$s'_\oldOmega(M-\cals_{n+1})\ge s'_\oldOmega(M-\cals)+\lambda n+\lambda$. Thus
(\ref{compnum equation}) holds with $n+1$ in place of $n$, and the
recursive construction of the $\cals_n$ is complete. As observed above, this completes the proof of the corollary.
\EndProof
%s_

\Lemma\label{tell me why}
Let $\Mh$ be a closed, orientable
hyperbolic $3$-orbifold, and set $\oldOmega=(\Mh)\pl$.
Let $\cals$ be an
$\oldOmega$-admissible system of spheres
in $M:=|\oldOmega|$, and let $S$ be a \doublefull\ component of $\cals$. If 
$0<\chibar(X''(\oldOmega,
\obd(\cals),
\obd(\cals-S)))<2\chibar(\obd(S))$,
then $S$
is a \clashsphere\ for $\cals$.
\EndLemma

\Proof
Set $\oldPsi=\oldOmega\cut{\obd(\cals)}$,
$\tcals=\partial|\oldPsi|$,
$\rho=\rho_{\obd(\cals)}:\oldPsi\to\oldOmega$,
$\tau=\tau_{\obd(\cals)}:\partial\oldPsi\to\partial\oldPsi$, and $\oldPhi^-=\oldPhi^-(\oldPsi)$. Set $\tcals_0=|\rho|^{-1}(S)$. Set $V_n=V_n^{\oldOmega,\obd(\cals)}$ and $\noodge_n=\noodge_n^{\oldOmega,\obd(\cals)}$ for every $n\ge1$. Set $X=X(\oldOmega,
\obd(\cals),
\obd(\cals-S))$, $X_n=X_n(\oldOmega,
\obd(\cals),
\obd(\cals-S))$ for each $n\ge1$, and $X''=X''(\oldOmega,
\obd(\cals),
\obd(\cals-S))$. 
Set $U=X''\wedge[\tau](X'')\in\barcaly_-(\tcals)$.

Let us assume that $S$ is not a \clashsphere\ for $\cals$. By definition this means that we do not have
$\chibar(U)<\chibar(X'')$. Since
$U\preceq X''$, it follows from Proposition \ref{new partial order} that we do have 
$\chibar(U)\le\chibar(X'')$, and hence 
\Equation\label{the basics}
\chibar(U)=\chibar(X'').
\EndEquation

Fix an element $\oldXi$ of $\Theta_-(\obd(\tcals))$ such that
$[\oldXi]=X$. Since $X\preceq V_1=[\oldPhi^-]$ by
\ref{oldXi}, \ref{i thought so too} and (\ref{and so did you}),  we may take $\oldXi$ to be contained in
$\inter\oldPhi^-$. Set $\oldXi''=\oldXi\cap\obd(\tcals_0)$. 
Since $\obd(\tcals_0)$ is a union of components of $\obd(\tcals) $, the suborbifolds $\oldXi$ and $\obd(\tcals_0)$ of $\obd(\tcals)$ are in standard position, and $\oldXi''$ is negative. Hence by (a rather degenerate case of) Assertion (2) of Proposition \ref{negative lattice}, we have $[\oldXi'']=X''$.

According to 
(\ref
{the basics}), we have
$\chibar(X''\wedge[\tau](X'') )=\chibar(X'')$.
Thus the hypotheses of Lemma \ref{inchworm lemma} hold with $X''$,  $[\tau](X'')$, $\frakX''$ and $\tau(\frakX'')$ playing the respective roles of $Z_1$, $Z_2$, $\frakZ_1$ and $\frakZ_2$. Hence if $\oldUpsilon$ denotes the union of $\frakX''$ with all components of $\obd(\tcals_0)-\inter\oldXi''$ that are annular orbifolds, 
then Lemma \ref{inchworm lemma} implies that $[\oldUpsilon]=[\tau(\oldUpsilon)]$, i.e.:

\Claim\label{glocca morra}
The suborbifolds $\oldUpsilon$ and $\tau(\oldUpsilon)$ of $\obd(\tcals_0)$ are isotopic. 
\EndClaim

Since $\obd(\tcals_0)$ is invariant under the involution $\tau$, it follows that if we set $\oldUpsilon^*=\obd(\tcals_0)-\inter\oldUpsilon$, then:
\Claim\label{yes they ah}
The suborbifolds $\oldUpsilon^*$ and $\tau(\oldUpsilon^*)$ of $\obd(\tcals_0)$ are isotopic. 
\EndClaim

Note that the definition of $\oldUpsilon$  directly implies:

\Claim\label{directly} No component of $\oldUpsilon^*=\obd(\tcals_0)-\inter\oldUpsilon$ is annular.
\EndClaim

Next we claim:

\Claim\label{ankle}
Each of the suborbifolds $\oldUpsilon$ and $\oldUpsilon^*$ of $\obd(\tcals_0)$ is non-empty, negative and taut.
\EndClaim

To prove \ref{ankle}, first note that by the definition of $\oldUpsilon$, each component of $\oldUpsilon$ is a union of some non-empty set of components of $\oldXi''$ and some (possibly empty) set of  annular orbifolds, each of which has non-empty intersection with $\oldXi''$; and that any two of the orbifolds whose union defines $\oldUpsilon$ meet precisely in a union of common boundary components. But $\oldXi''$ is negative and taut since 
$\oldXi\in\Theta_-(\obd(\tcals))$, and it follows that
 $\oldUpsilon$ is negative and taut. Since $\oldUpsilon$ is taut and $\partial\oldUpsilon^*=\partial\oldUpsilon$, the suborbifold $\oldUpsilon^*$ is also taut. Furthermore, since $\oldUpsilon^*$ has no annular component by \ref{directly}, and has no toric component since $\partial\oldPsi$ is negative by \ref{doublemint}, the suborbifold $\oldUpsilon^*$ is negative.

To complete the proof of \ref{ankle}, it remains to show that $\oldUpsilon$ and $\oldUpsilon^*$ are non-empty. The definition of $\oldUpsilon$ implies that $\chi(\oldUpsilon)=\chi(\oldXi'')$. With the hypothesis of the lemma, this gives $0<\chibar(\oldUpsilon)<2\chibar(\obd(S))=\chibar(\obd(\tcals_0))$. Thus $\chibar(\oldUpsilon)>0$, while $\chibar(\oldUpsilon^*)=\chibar(\tcals_0)-\chibar(\oldUpsilon)>0$. This shows that $\oldUpsilon$ and $\oldUpsilon^*$ are non-empty, and \ref{ankle} is proved.

Now let $\frakV$ denote the union of all negative components of $(\obd(\tcals_0)\cap\oldPhi^-)-\inter\oldXi''$. We claim:

\Claim\label{my little frotzy}
We have $\frakV\subset\oldUpsilon^*$. Furthermore, every component $\frakB$ of $\frakD:=\overline{\oldUpsilon^*-\frakV}$ is an annular orbifold, each of whose boundary components is a component of either $\partial\frakV$ or of $\partial\oldUpsilon=\partial\oldUpsilon^*$; and at least one component of $\partial\frakB$ is a component of $\partial\frakV$.
\EndClaim

To prove \ref{my little frotzy}, first note that if $\frakK$ is any component of $\frakV$, then by definition $\frakK$ is negative, and $\frakK\subset (\obd(\tcals_0)\cap\oldPhi^-)-\inter\oldXi''
\subset \obd(\tcals_0)-\inter\oldXi''
$. Furthermore, $\frakK$ is
 taut, since each component of $\partial\frakK$ is a boundary component of one of the taut suborbifolds  
$\oldPhi^-$ and $\oldXi''$.
%\subset \obd(\tcals_0)-\inter\oldXi''$. 
Since $\frakK$ is negative and taut, the component of $
\obd(\tcals_0)-\inter\oldXi''$ containing $\frakK$ cannot be annular; hence $\frakK$, in addition to being disjoint from $\inter\oldXi''$, is disjoint from all annular components of $\obd(\tcals_0)-\inter\oldXi''$. In view of the definition of $\oldUpsilon$, it follows that $\frakK$ is disjoint from $\inter\oldUpsilon$. 
%\redcomment{This seems wrong, because $\oldUpsilon$ was defined in terms of $\obd(\tcals_0)-\inter\oldXi''$, not $(\obd(\tcals_0)\cap\oldPhi^-)-\inter\oldXi''$. Should $\frakV$ have been defined in terms of $\obd(\tcals_0)-\inter\oldXi''$ also?} 
This shows that $\frakV\subset\oldUpsilon^*$.

To prove the remaining assertions of \ref{my little frotzy}, first note that according to the definitions we have 
$\partial\oldUpsilon^*=\partial\oldUpsilon\subset\partial\oldXi''$ and $\partial\frakV\subset\partial\oldPhi^-\cup\partial\oldXi''$. Since $\oldXi\subset\inter\oldPhi^-$, the intersection of the $1$-manifolds $\partial\frakV$ and $\partial\oldUpsilon^*$ is a union of common components.
Hence 
$\frakD=\overline{\oldUpsilon^*-\frakV}$ is a suborbifold of $\obd(\tcals_0)$, each of whose boundary components is contained in either $\partial\frakV$ or $\partial\oldUpsilon^*=\partial\oldUpsilon$.

This argument also shows that
 each component of  $\partial\frakD$ is a component either of $\partial\oldPhi^-$ or of $\partial\oldXi''\subset\inter\oldPhi^-$. 
%Since $\oldPhi^-$ and $\frakX''$ are taut, it follows that $\frakD$ is taut in $\obd(\tcals)$. 
%Since each component of  $\partial\frakD$ is a component  of $\partial\oldPhi^-$ or of $\partial\oldXi''$, 
We may therefore write
$\frakD=\frakD_0\cup\frakD_1$, where $\frakD_0$ is a suborbifold of $\oldPhi^-$ meeting $\partial\oldPhi^-$ in a union of components of $\partial\oldPhi^-$, and $\frakD_1$ is a union of components of $\tcals_0\setminus\inter\oldPhi^-$ which also meets  $\partial\oldPhi^-$ in a union of components of $\partial\oldPhi^-$. Since $S$ is \doublefull, each component of $\tcals_0\setminus\inter\oldPhi^-$ is annular (see \ref {special ouenelitte}). 
%Since $\frakD$ is taut, it follows via \ref{cobound} that 
In particular, each component of $\frakD_1$ is annular. Now suppose that $\frakB_0$ is a component of $\frakD_0$. We have $\frakB_0\subset\frakD_0\subset \frakD\cap\oldPhi^-=\overline{\oldUpsilon^*-\frakV}\cap\oldPhi^-$. Since $\oldUpsilon^*$ is the complement relative to $\obd(\tcals_0)$ of $\inter\oldUpsilon
%\s, so that $\frakB_0$ is disjoint from %$\inter\oldUpsilon
\supset\tcals_0\cap\inter\oldXi''$, it follows that $\frakB_0\subset
(\obd(\tcals_0)\cap\oldPhi^-)-\inter\oldXi''$. Since $\frakB_0$ is disjoint from $\inter\frakV$, the definition of $\frakV$ implies that the component $\oldPi$ of
$(\obd(\tcals_0)\cap\oldPhi^-)-\inter\oldXi''$ containing $\frakB_0$
%The component $\oldPi$ of %$\obd(\tcals_0)-\inter\oldXi''$ containing %$\frakB_0$ is disjoint from %$\inter\oldUpsilon$, and hence 
cannot be negative. But $\oldPi$ is taut, since each component of $\partial\oldPi$ is a boundary component of one of the taut suborbifolds $\oldPhi^-$ and $\oldXi''$ of $\obd(\tcals)$. Hence $\chi(\oldPi)\le0$, and since $\oldPi$ cannot be negative we have
% in view of the definition of %$\inter\oldUpsilon$. 
%\redcomment{Is this any better? If so, fix up %the argument from this point.}
%follows at all!! I think the information I need %is that $\frakB_0$ is disjoint from %$\inter\frakV$.} 
%By tautness \redproofreadingnote{possible %cross-ref to \ref{ankle}}, we must have
 $\chi(\oldPi)=0$.  Since $\partial\oldPsi$ is negative by \ref{doublemint}, $\oldPi$ cannot be toric, and  must therefore be annular. It then follows from \ref{cobound} that $\frakB_0$ is annular. Thus each component of $\frakD_0$, as well as each component of $\frakD_1$, is annular; this implies that each component of $\frakD$ has Euler characteristic $0$. 
But the  negativity of  $\partial\oldPsi$ also implies that no component of $\frakD$ can be toric, and hence
each component of $\frakD$ is annular. %\redproofreadingnote{The part about
%each of its boundary components being a component of %either $\partial\frakV$ or of $\partial\oldUpsilon$,  %should have been made clearer than it was.}

To complete the proof of \ref{my little frotzy}, it remains only to prove that  $\frakD=\overline{\oldUpsilon^*-\frakV}$ has no component $\frakB$ such that $\partial\frakB\subset\partial\oldUpsilon$. If $\frakB$ is such a component, it must be a component of $\oldUpsilon^*$; but by what we have already proved, the component $\frakD$ of $\frakB$ must be annular. This contradicts \ref{directly}. Thus \ref{my little frotzy} is proved.

Assertions \ref{very latest}---\ref{shuck not} below are preparation for the proof of
%The next key point in the proof of the present lemma is 
\ref{blech}, which is equivalent to the assertion that the orbifold $\frakV\subset\oldPhi^-$ is isotopic in $\obd(\tcals)$ to its image under the involution $\iota_\oldPsi$ of $\oldPhi^-$ (see \ref{what's iota?}). 
The strategy of the proof of the lemma will then be to 
%prove ), and then try to give a hint of how the rest of the proof will work: I will prove \ref{blech} (give some version of the statement here), then 
combine \ref{blech} with \ref{yes they ah}, \ref{my little frotzy}, and \ref{ankle} 
to get a contradiction to the hyperbolicity of $\Mh$.

%\redcomment{The next couple of pages are devote to the proof of \ref{blech}. Revise them as needed to make them fit with what comes before.}

According to the definitions we have $X''\preceq[\obd(\tcals_0)]$ and hence $U\preceq[\obd(\tcals_0)]$. Recalling (see \ref{in duck tape}) that $\dom[\iota]=[\oldPhi^-]$, we claim:
\Claim\label{very latest}
We have $U\preceq [\oldPhi^-]$, and $[\iota](U)\preceq X$.
\EndClaim
To prove \ref{very latest}, first note that $U\preceq X''$ by the definition of $U$, while $X''\preceq X\preceq[\oldPhi^-]$ by an observation made in \ref{oldXi}. It follows that $U\preceq[\oldPhi^-]$. 

To prove that $[\iota](U)\preceq X$, according to Proposition
\ref{mystery}, it suffices to show that for every component $G$ of $U$
we have $[\iota](G)\preceq X$. If $G$ is a component of $U$, then since
$G\preceq U\preceq [\tau](X'')\preceq [\tau](X)$, and since $[\tau]^{-1}=[\tau]$ by \ref{in duck tape}, we deduce via \ref{special composition} that $[\tau](G)\preceq X$. If $Y$ denotes the component of $X$ such that
$[\tau](G)\preceq Y$ (see \ref{components}), then by 
Lemma \ref{here goethe}, $Y$ is a component of $X_n$ for some
$n\ge1$. In particular we have $[\tau](G)\preceq X_n\preceq
V_n$ (where the relation $X_n\preceq V_n$ holds by \ref{oldXi}).

Set $G'=[\iota](G)$. Since $[\iota]^{-1}=[\iota]$ by \ref{in duck tape}, we have $G'\preceq 
\dom[\iota]=
[\oldPhi]$, and in view of \ref{special composition} we have $[\iota](G')=G$. Since we also have $G\preceq[\obd(\tcals)]=\dom[\tau]$, Lemma \ref{before associativity} then implies that $G'\preceq\dom([\tau]\diamond[\iota])$ and that $([\tau]\diamond[\iota])|G'=[\tau]\circ([\iota]|G')$. In particular, using \ref{special composition}, we have 
$([\tau]\diamond[\iota])(G')=[\tau]([\iota](G'))=[\tau](G)\preceq V_n=\dom\noodge_n$. 
A second application of Lemma \ref{before associativity} then gives that $G'\preceq\dom
(\noodge_n\diamond
([\tau]\diamond[\iota]))
$, and that $(\noodge_n\diamond
([\tau]\diamond[\iota]))|G'=
\noodge_n\circ
([\tau]\diamond[\iota])|G')$.
But it follows from the Definition of $\noodge_n$ and $\noodge_{n+1}$
(see \ref{in duck tape}) and from \ref{more associativity} that 
%$\noodge_n=[\iota]\diamond([\tau]\diamond[\iota])^{\diamond (n-1)}
$\noodge_n\diamond
([\tau]\diamond[\iota]))=[\iota]\diamond([\tau]\diamond[\iota])^{\diamond (n-1)}\diamond
[\tau]\diamond[\iota]=\noodge_{n+1}$. 
Hence $G'\preceq\dom\noodge_{n+1}=V_{n+1}$, and $\noodge_{n+1}(G')=\noodge_n(([\tau]\diamond[\iota])(G'))=
\noodge_n(([\tau]\circ[\iota])(G'))=
\noodge_n([\tau](G))$, by \ref{special composition} and the equality $[\iota](G')=G$. Since $[\tau](G)\preceq X_n$, and since $X_n\in\oldfrakX_n(\oldOmega,
\obd(\cals),
\obd(\cals-S))$ by \ref{oldXi}, we have $\noodge_n([\tau](G))
\preceq \noodge_n(X_n)
\preceq[\obd(\tcals-\tcals_0)]$, i.e.
$\noodge_{n+1}(G')\preceq[\obd(\tcals-\tcals_0)]$. If $P$ denotes the component of $V_{n_1}$ such that $G'\preceq P$ (see \ref{components}), then since
$\tcals-\tcals_0$ is a union of components of $\tcals$, we have 
$\noodge_{n+1}(P)\preceq[\obd(\tcals-\tcals_0)]$. By definition this means that $P\in
\oldfrakX_{n+1}(\oldOmega,
\obd(\cals),
\obd(\cals-S))$, so that $P\preceq X_{n+1}$. In particular we have $G'\preceq X_{n+1}\preceq X$, i.e. $[\iota](G)\preceq X$, and the proof of \ref{very latest} is complete.

Let us
fix an element $\frakU$ of $\Theta_-(\obd(\tcals))$ such that $[\frakU]=U$. Since $U\preceq X''$, and since $\chibar(U)=\chibar(X'')$ by \ref{the basics},  it follows from Lemma \ref{before inchworm} (applied with $X''$ and $U$ playing the respective roles of $Z$ and $Z'$) that we may choose $\frakU$ in such a way that (1) $\frakU\subset\oldXi''$, (2) $\partial\oldXi''\subset\partial\frakU$, and (3) each component of $\overline{\oldXi''-\frakU}$ is an annular orbifold contained in $\inter\oldXi''$.
This implies:
\Claim\label{beppo} Every component of $\oldPhi^--\inter\frakU$ is either a component of $\oldPhi^--\inter\oldXi''$ or an annular orbifold.
\EndClaim
On the other hand, it follows from \ref{very latest} that the
restriction to $\frakU$ of the homeomorphism $\iota:\oldPhi^-\to\oldPhi^-$ is isotopic in $\obd(\tcals)$ to an embedding $f$
having domain $\frakU$ and having image contained in $\oldXi \subset\inter\oldPhi^-$. Since
$\iota|\frakU$ and $f$ are isotopic in $\obd(\tcals)$ and both their
images are contained in $\inter\oldPhi^-$, it follows from Corollary \ref{i guess} that they are isotopic in $\oldPhi^-$. Hence 
there is a homeomorphism $\iota':\oldPhi^-\to\oldPhi^-$, isotopic
%\redmissingref{I had the phrase ``}
in $\oldPhi^-$
%\redmissingrefcontinued{'' here but I don't
%  think it belongs} 
to $\iota$, such that
$\iota'(\frakU)\subset\oldXi$. (It is not obvious whether $\iota'$ can be taken to be an involution, and the issue will not affect the argument.) We have
% to there is an element $\iota'$ of  $\calm(\obd(\tcals))$ such that
 $[\iota']=[\iota] \in\calm(\obd(\tcals))$. Set
 $\iota''=(\iota')^{-1}:\oldPhi^-\to\oldPhi^-$. Then
 $\iota''\in\calm(\obd(\tcals))$ and
 $[\iota'']=[\iota']^{-1}=[\iota]^{-1}=[\iota]$ (where the last equality holds by \ref{in duck tape}).

The definition of $\frakV$ implies that $\frakV$ is negative, and that each of its boundary components is a boundary component of one of the taut orbifolds   $\oldXi''$ and $\oldPhi^-$. Hence $\frakV$ is taut, and so $\frakV\in\Theta_-(\obd(\tcals))$. Set $Q=[\frakV]\in\barcaly_-(\obd(\tcals))$. Then $Q\preceq[\oldPhi^-]=\dom[\iota]$. We claim:
\Equation\label{not an iota}
[\iota](Q)\preceq[\obd(\tcals_0)].
\EndEquation
To prove (\ref{not an iota}), assume to the contrary that
$[\iota](Q)\not\preceq[\obd(\tcals_0)]$. Then according to Proposition \ref{mystery}, for some component $K$
of $Q$ we have
$[\iota](K)\not\preceq[\obd(\tcals_0)]$, and therefore
$[\iota](K)\preceq[\obd(\tcals-\tcals_0)]$. 
By (\ref{and so did you}) this means $[\noodge_1](K)\preceq[\obd(\tcals-\tcals_0)]$. If $K^+$ denotes the component of $V_1=\oldPhi^-$ (cf. \ref{and so did you}) such that $K\preceq K^+$ (see \ref{components}), then by 
\ref{oldXi},  $K^+$ is an element of $\oldfrakX_{1}(\oldOmega,
\obd(\cals),
\obd(\cals-S))$, and  therefore $K^+\preceq X_1$. Hence $K\preceq X_1\preceq X$. Since
$K\preceq Q=[\frakV]\preceq[\obd(\tcals_0)]$, it follows that $K\preceq X''$,
and therefore that $K\wedge X''=K$. But the definition of $\frakV$
implies that $\frakV$ is isotopic to an element of
$\Theta_-(\obd(\tcals))$ which is disjoint from $\oldXi''$,
and we may
therefore write $K=[\frakK
]$ for some $\frakK\in\Theta_-(\obd(\tcals))$
with $\frakK\cap\oldXi''=\emptyset$. It then follows from the first assertion of Corollary \ref{nafta} that $K\wedge X''=[\emptyset]$; hence $K=[\emptyset]$, which is impossible. Thus (\ref{not an iota}) is proved.

Note that since $\oldXi''=\oldXi\cap\obd(\tcals_0)$, we have
\Equation\label{go shuck yourself} \frakV\subset(\obd(\tcals_0)\cap\oldPhi^-)-\inter\oldXi''\subset \oldPhi^--\inter\oldXi.
\EndEquation
Next, note that
since $\iota'(\frakU)\subset\oldXi$ and $\iota'(\oldPhi^-)=\oldPhi^-$, we have $\iota''(\oldPhi^--\inter\oldXi)\subset\oldPhi^--\inter\frakU$. With (\ref{go shuck yourself}), this gives
\Equation\label{shuck not}
\iota''(\frakV)\subset\oldPhi^- -\inter\frakU.
\EndEquation

If $\oldGamma$ is any  component of $\frakV$,
it follows from (\ref{shuck not}) and \ref{beppo} that either 
$\iota''(\oldGamma)\subset\oldPhi^- -\inter\oldXi''$, or
$\iota''(\oldGamma)$ is contained in an annular suborbifold of
$\oldPhi^-$. But since $\iota''$ is a self-homeomorphism of
$\oldPhi^-$, and since $\oldGamma$ is an element of 
$\Theta_-(\obd(\tcals))$, the suborbifold $\iota''(\oldGamma)$ is
negative and $\pi_1$-injective in $\oldPhi^-$, and therefore cannot be
contained in an annular suborbifold of $\oldPhi^-$. Hence
$\iota''(\oldGamma)\subset\oldPhi^- -\inter\oldXi''$. On the other
hand, since $[\iota'']=[\iota]$, it follows from (\ref{not an iota})
that
$[\iota''](Q)\preceq[\obd(\tcals_0)]$, so that $\iota''(\oldGamma)\subset\obd(\tcals_0)$.
 Combining this with the inclusion $\iota''(\oldGamma)\subset\oldPhi^- -\inter\oldXi''$, we deduce that $\iota''(\oldGamma)\subset
(\obd(\tcals_0)\cap\oldPhi^-) -\inter\oldXi''$, which in view of the negativity and tautness of $\iota''(\oldGamma)$ implies that $\iota''(\oldGamma)\subset\frakV$. As this holds for every component $\oldGamma$ of $\frakV$, we have $\iota''(\frakV)\subset\frakV$. Hence $[\iota](Q)\preceq Q$. By the order-preserving property observed in \ref{restriction}, and \ref{special composition},
it follows that $Q\preceq[\iota]^{-1}(Q)=[\iota](Q)$. Since $\preceq$ is a partial order by Proposition \ref{new partial order}, we now deduce:
\Equation\label{blech}
[\iota](Q)=Q.
\EndEquation

We now claim:
\Claim\label{marigold}
There is a compact $3$-suborbifold $\frakZ$ of $\oldPsi$ such that $\frakZ\cap\obd(\tcals)=
\oldUpsilon^*$, and every component of $\Fr_\oldPsi\frakZ$ has Euler characteristic $0$. 
\EndClaim

%\redproofsummary{Revise the material below to give a proof of this. The argument as it's given now gives the corresponding result in terms of $\frakV^+$ (I think it's called) instead of $\oldUpsilon^*$. The adjustment should involve using \ref{my little frotzy} to say that $\oldUpsilon^*$ is isotopic to something having similar properties to what is called $\frakV^+$, which means that it's enough to prove \ref{marigold} with this other thing in place of $\oldUpsilon^*$. The argument below should do that. The role of what is called $\frakD_!$ below will be played by $\frakD=\overline{\oldUpsilon^*-\frakV}$, which I  should remind the reader has been named canonically in the statement of \ref{my little frotzy}.}

To prove \ref{marigold}, first note that by (\ref{blech}) and Proposition \ref{steinbeck}, there is an element $\frakV_! $ of $\Theta_-(\obd(\tcals))$ such that $[\frakV_! ]=Q$ and $\iota(\frakV_! )=\frakV_! $. Now according to \ref{what's iota?},
$\oldSigma^-:=\oldSigma^-(\oldPsi)$
may be equipped with an $I$-fibration $q|\oldSigma^-:\oldSigma^-\to \frakB$  over a $2$-orbifold in
such a way that $\partialh\oldSigma^-=\oldSigma^-\cap\obd(\tcals)=\oldPhi^-$, and $\iota$ is the non-trivial deck transformation of the
two-sheeted covering map
$q|\oldPhi^-:\oldPhi^-\to \frakB$. 
Since $\iota(\frakV_! )=\frakV_! $, there is a $2$-suborbifold 
$\frakC$ of $\frakB$ such that $(q|\oldPhi^-)^{-1}(\frakC)=\frakV_! $. If we set $\frakZ^0=q^{-1}(\frakC)$, it follows that $q|\frakZ^0:\frakZ^0\to\frakC$ is an $I$-fibration, and that with respect to this $I$-fibration we have $\partialh\frakZ^0=\frakV_! =\frakZ^0\cap\obd(\tcals)$. It follows that $\Fr_\oldPsi\frakZ^0=\partialv\frakZ^0$, so that the components of $\Fr_\oldPsi\frakZ^0$ are annular orbifolds.
%\reddelete{In particular, $\frakZ^0$ is a \pagelike\ \Ssuborbifold\ of $\oldPsi$. 
%
%(1) For a \pagelike\ \Ssuborbifold, don't I need to argue that the frontier components are essential? But I guess after the revision this will be irrelevant. (2)}

Since $[\frakV_!]=Q$, there is an (orbifold) isotopy $(h_t)_{0\le t\le1}$ of $\tcals$ such that $h_1(\frakV)=\frakV_!$. 
It follows from \ref{my little frotzy} that $h_1(\oldUpsilon^*)$ is a suborbifold $\oldUpsilon^*_!$  containing $\frakV_!$; and that every component  of $\frakD_!:=\overline{\oldUpsilon_!^*-\frakV_!}$ is an annular orbifold, each of whose boundary components is a component of either $\partial\frakV_!$ or of $
%\partial\oldUpsilon_!=
\partial\oldUpsilon_!^*$, and  at least one of whose boundary components  is a component of $\partial\frakV_!
$. 
%\redcommentB{Do I need all that? The last bit (``at least one''), and the corresponding part of \ref{my little frotzy}, may not matter.}

In particular, $\frakD_!\cap\frakZ^0=\frakD_!\cap\frakV_!$ is a union of components of $\partial\frakD_!$. 
%Let $\frakD_!$ denote the union of all components of $\obd(\tcals_0)-\inter\frakV_!$ that are annular orbifolds, and set ${\frakV_+}_!=\frakV_! \cup\frakD_!$. Comparing this definition with the definitions of $\frakD$ and $\frakV_+$ (and using the fact that $\frakV_!$ and $\frakV$ are isotopic since $[\frakV_!]=Q=[\frakV]$), we find that ${\frakV_+}_!$ is isotopic to $\frakV_+$. 
%It follows from the definition of $\frakD_!$ that
%$\frakD_!\subset\obd(\tcals_0)$ and that
%$\frakD_!\cap\frakV_! =\partial\frakD_!$. Since $\frakV_!
%=\frakZ^0\cap\obd(\tcals)$, it then follows that
%$\frakD_!\cap\frakZ^0=\partial\frakD_!$. 
Now let $\frakZ^1$ denote the
union of $\frakZ^0$ with a strong (orbifold) regular neighborhood of
$\frakD_!$ in $\overline{\oldPsi-\frakZ^0}$. Then $\Fr_\oldPsi\frakZ^1$
is homeomorphic to an orbifold obtained from the disjoint union
$(\Fr_\oldPsi\frakZ^0)\discup\frakD_!$ by gluing certain boundary
components of $\frakD_!$ to certain boundary components of
$\Fr_\oldPsi\frakZ^0$; each  boundary component of
$\Fr_\oldPsi\frakZ^0$ is glued to at most one boundary component of
$\frakD_!$, and vice versa. Since the components of
$\frakD_!$
% are annular orbifolds, 
and 
%since  $\partialh\frakZ^0=\frakV_! =\frakZ\cap\obd(\tcals)$,
%$\partialh\oldSigma^-=\oldPhi^-$,
%$\frakZ^0$ is standardly embedded in $\oldPsi$ with respect to its $I$-fibration \redcomment{Rewrite that somehow in the language of \Ssuborbifold s}, 
%the components of
  $\Fr_\oldPsi\frakZ^0$ are  annular
orbifolds, each component of $\Fr_\oldPsi\frakZ^1$ has Euler
characteristic $0$. But the definition of $\frakZ^1$ implies that
$\frakZ^1\cap\obd(\tcals)$ is the union of $\frakV_!$ with a strong regular neighborhood of $\frakD_!$ in $\obd(\tcals_0)-\inter\frakV_!$. Hence $\frakZ^1\cap\obd(\tcals)$ is isotopic in $\tcals_0$ to $\frakV_!\cup\frakD_!=\oldUpsilon^*_!$, and  therefore to $\oldUpsilon^*$. 
%At the moment it looks to me as if 
%$\frakZ^1\cap\obd(\tcals)$ is equal to $\frakV_!\cup\frakD_!=\oldUpsilon^*_!$
%rather than being a regular neighborhood of it. No, it's more complicated than that. I think 
%$|\frakZ^1|\cap\tcals$ is obtained from 
%$|\oldUpsilon^*_!|$ by adding weight-$0$ annuli to some, not necessarily all, boundary components.}
It follows that
$\frakZ^1$ is isotopic to a suborbifold $\frakZ$ such that
$\frakZ\cap\obd(\tcals)=\oldUpsilon^*$. This proves \ref{marigold}.

We are now ready to establish the final contradiction that will
complete the proof of the present lemma. 

Let us fix a component $\tS$ of $\tcals_0$. If $\tS\subset|\oldUpsilon|$, then  it follows from \ref{glocca morra} that $|\tau|(\tS)\subset|\oldUpsilon|$, so that $\tcals_0=\tS\cup|\tau|(\tS)\subset|\oldUpsilon|$; this is impossible since $\tcals_0-\inter|\oldUpsilon|=|\oldUpsilon ^*|\ne\emptyset$ by \ref{ankle}. Hence $\tS\not\subset\oldUpsilon$. The same argument, with the roles of $\oldUpsilon$ and $\oldUpsilon^*$ reversed, and \ref{yes they ah} used in place of \ref{glocca morra},  shows that $\tS\not\subset|\oldUpsilon^*|$, i.e. $\tS\cap|\oldUpsilon|\ne\emptyset$. Since $\tS$ is connected, it follows that $\Fr_\tS(\tS\cap|\oldUpsilon|)\ne\emptyset$; we fix a component $C$ of $\Fr_\tS(\tS\cap|\oldUpsilon|)$. Then $C$ is a common boundary component of some component $W$ of $\tS\cap|\oldUpsilon|$ and of some component $W^*$ of $\tS\cap|\oldUpsilon^*|$.

Set $\tS'=|\tau|(\tS)$.
Let $\frakZ$ be an orbifold having the properties stated in \ref{marigold}. According to \ref{yes they ah}, 
$\oldUpsilon^*=\frakZ\cap\obd(\tcals)
$ is isotopic in $\tcals_0$ to 
$\tau(\oldUpsilon^*)=\tau(\frakZ\cap\obd(\tcals))$. In particular, $\oldUpsilon^*\cap\obd(\tS')=\frakZ\cap\obd(\tS')$ is isotopic in $\obd(\tS')$ to  $\tau(\oldUpsilon^*)\cap\obd(\tS')=\tau(\frakZ\cap\obd(\tS))$. Hence $\frakZ$ is isotopic in $\oldPsi$, rel $\obd(\tS)$,  to
%Since $\tau$ interchanges the two components of $\obd(\tcals_0)$, it follows  that $\frakZ$ is isotopic to 
a suborbifold $\frakZ_1$ such that 
$\frakZ_1\cap\obd(\tS')=\tau(\frakZ\cap\obd(\tS))$. We now have
$\frakZ_1\cap\obd(\tS')= 
\tau(\frakZ_1\cap\obd(\tS))$.
%, so  
%that
$\frakZ_1\cap\obd(\tcals_0)$ is $\tau$-invariant.
% Since  $\oldUpsilon$ is invariant under $\tau$, the orbifold 
%$\oldUpsilon':=\obd (\tcals_0)-\inter\oldUpsilon$
%is also invariant under $\tau$; since \ref{marigold} gives that
%$\frakZ\cap\obd(\tcals)=\oldUpsilon'$, 
It follows that 
$\rho(\frakZ_1)$ is a suborbifold of $\oldOmega$, that
$\frakE:=\partial(\rho(\frakZ_1))=\rho(\Fr_\oldPsi\frakZ_1)$, and that
up to homeomorphism the $2$-orbifold $\frakE$ is obtained from $\Fr_\oldPsi\frakZ_1$ by gluing the boundary components of $\Fr_\oldPsi\frakZ$ in pairs. But it follows from  \ref{marigold} that each component of $\Fr_\oldPsi\frakZ$ has Euler characteristic $0$, and hence each component of $\frakE$ is toric.

Since $\frakZ_1\cap\tS=\frakZ\cap\tS=\oldUpsilon^*\cap\tS$, we have $\Fr_\tS(\oldUpsilon^*\cap\tS)\subset\Fr_\oldPsi\frakZ_1$. Since the definition of $C$ implies that $C\subset\Fr_\tS(\oldUpsilon^*\cap\tS)$, it follows that $C\subset\Fr_\oldPsi\frakZ_1$ and hence that $\rho(C)\subset\rho(\Fr_\oldPsi\frakZ_1)=\frakE$.
Let $\frakE_0$ denote the component of $\frakE$ containing $\rho(C)$. Then $\frakE_0$ is toric, and is two-sided since $\frakE:=\partial(\rho(\frakZ_1))$.

%Added 6/17/18: Maybe it is true if, as suggested above, I use what was called $\frakZ^1$ to  play the role of 
%$\frakZ$ 
%} $C\subset\frakE$. 

By Proposition \ref{preoccupani}, $\frakE_0$ is the boundary of a
suborbifold $\frakJ$ of $\oldOmega$ such that the inclusion homomorphism $\pi_1(\frakJ)\to\pi_1(\oldOmega)$ has a virtually cyclic image. In particular, $|\frakE_0|$ separates the manifold $|\oldOmega|$, and the components of its complement are $\inter|\frakJ|$ and $|\oldOmega|-|\frakJ|$. Note also that the $2$-manifold $|\Fr_\oldPsi\frakZ_1|$ is properly embedded in $|\oldPsi|$, and hence that the closed $2$-manifolds $S$ and $|\frakE_0|$ meet transversally (in the PL sense) in $|\oldOmega|$. Since $\rho$ maps $\tS$ homeomorphically onto $S$, we have connected subsurfaces $\rho(W)$  and $\rho(W^*)$ of $S$ which have disjoint interiors and share the boundary component $\rho(C)$. It follows that $\rho(W)$  and $\rho(W^*)$ lie in distinct components of 
$|\oldOmega|-|\frakE_0|$. In particular,
% separates the manifold 
%\redcomment{Explain why} 
either $\rho(W)$ or $\rho(W^*)$ is contained in $\frakJ$.

Hence at least one of the homomorphisms $(\rho|W)_\sharp:\pi_1(W)\to\pi_1(\oldOmega)$ 
and $(\rho|W^*)_\sharp:\pi_1(W^*)\to\pi_1(\oldOmega)$ has a virtually cyclic image.
But it follows from \ref{ankle} that $W$ and $W^*$ are  $\pi_1$-injective in $\oldPsi$, and since $\cals$ is admissible it follows that the homomorphisms $(\rho|W)_\sharp$ and $(\rho|W^*)_\sharp$ are injective. Hence either $\pi_1(W)$ or $\pi_1(W^*)$ is virtually cyclic. This is impossible because $\chi(W)$ and $\chi(W^*)$ are negative by \ref{ankle}.
\EndProof
%\frakWU$P\frakW_1\frAKW^*_1\oldUpsilon_1\oldUpsilon^*_1\tau\Πρ

\Lemma\label{dun da-dun dun}
Let $\Mh$ be a closed, orientable
hyperbolic $3$-orbifold, and set $\oldOmega=(\Mh)\pl$.
Let $\cals$ be an
$\oldOmega$-admissible system of spheres
in $M:=|\oldOmega|$, and let $S$ be a \doublefull\ component of $\cals$. If 
$X''(\oldOmega,\obd(\cals),\obd(\cals-S))=[\emptyset]$, then $\cals=S$.
\EndLemma

\Proof
Set $\oldPsi=\oldOmega\cut{\obd(\cals)}$,  $N=|\oldPsi|=M\cut\cals$,
$\tcals=\partial N$, and $\rho=\rho_{\obd(\cals)}:\oldPsi\to \oldOmega$, so that $|\rho|=\rho_{\cals}:N\to M$. Set
$\oldSigma=\oldSigma(\oldPsi)$, $\oldSigma^-=\oldSigma^-(\oldPsi)$, 
 $\oldPhi=\oldPhi(\oldPsi)$, $\oldPhi^-=\oldPhi^-(\oldPsi)$, and $\iota=\iota_\oldPsi$. Set $\tcals_0=|\rho|^{-1}(S)\subset\tcals$. Let $\frakZ$ denote the union of all components $\frakU$ of 
$\oldSigma^-$ such that $\frakU\cap\obd(\tcals_0)\ne\emptyset$. 
Since by \ref{what's iota?} we have $\oldSigma^-\cap\tcals=\oldPhi^-$,
we may alternatively describe
 $\frakZ$ as the union of all components $\frakU$ of 
$\oldSigma^-$ such that such that $\frakU\cap\oldPhi^-\cap\obd(\tcals_0)\ne\emptyset$. We claim:
\Equation\label{karsh}
|\frakZ|\cap\tcals\subset\tcals_0.
\EndEquation

To prove (\ref{karsh}), let $\frakU$ be any component of $\frakZ$. We
must prove that $|\frakU|\cap\tcals\subset\tcals_0$. 
%By definition we
%have $\frakU\cap\oldPhi^-\cap\obd(\tcals_0)\ne\emptyset$. 
According to the (alternative) description of $\frakZ$, the orbifold
 $\frakU$ is a component of $\oldSigma$ which contains a component
$\oldPi$ of $\oldPhi^-$ such that $|\oldPi|\subset\tcals_0$.
% and
%$\chi(\oldPi)<0$. By Lemma \ref{when a tore a fold}
According to the discussion in \ref{what's iota?}, $\frakU$ admits an $I$-fibration such that $\frakU\cap\obd(\tcals)=\partialh\frakU$. Thus $\oldPi$ is a component of $\partialh\frakU$.
%it follows that $\frakU$
%is a \pagelike\ \Ssuborbifold\ of $\oldPsi$. 
%If the connected \pagelike\
%\Ssuborbifold\ $\frakU$ is twisted, then 
If the $I$-fibration of $\frakU$ is non-trivial, then
$|\frakU|\cap\tcals=|\oldPi|$, so that in particular
$|\frakU|\cap\tcals\subset\tcals_0$, as required. 

There remains the case in which the fibration of $\frakU$ is trivial. Then 
%$\frakU$ admits a trivial
%$I$-fibration for which $\oldPi$ is a component of
%$\partialh\frakU$. It follows that
%$\oldPi=\partialh\frakU$, then properly contained in
%$\partialh\frakU$. Then by \redmissingref{cross-references}, 
$\partialh\frakU$ has two components, $\oldPi$ and
$\oldPi':=\iota(\oldPi)$ (see \ref{what's iota?}). It suffices to
prove that $|\oldPi'|\subset\tcals_0$. Assume to the contrary that
$|\oldPi'|\subset\tcals-\tcals_0$.
%=|\rho|^{-1}(\cals-S)$. 
Set $P=[\oldPi]\in\barcaly_-(\obd(\tcals))$, and
$P'=[\oldPi']=[\iota](P)$. Then $P$ and $P'$ are components of
$[\oldPhi^-]$, which by (\ref{and so did you})  is equal to $V_1:=V_1^{\oldOmega,\obd(\tcals)}$, and (\ref{and so did you}) also gives  
$[\iota]=\noodge_1:=\noodge_1^{\oldOmega,\obd(\tcals)}$. Hence
$P'=\noodge_1(P)$. Since $|\oldPi|\subset\tcals_0
% =|\rho|^{-1}(S)
$, and
since we have assumed that $|\oldPi'|\subset\tcals-\tcals_0$,
%$|\oldPi'|\subset|\rho|^{-1}(\cals-S)$,
we
have $P\preceq[\obd(\tcals_0)]$ and
$\noodge_1(P)=P'\preceq[\obd(\tcals-\tcals_0)]$. According
to  \ref{oldXi}, the relation $\noodge_1(P)\preceq[\obd(\tcals-\tcals_0)]$  means that 
$P\in\oldfrakX(\oldOmega,\obd(\cals),\obd(\cals-S))$, and therefore that
$P\preceq X(\oldOmega,\obd(\cals),\obd(\cals-S))$. Hence $P\preceq
X(\oldOmega,\obd(\cals),\obd(\cals-S))
\wedge
[\obd(\tcals_0)]
=X''(\oldOmega,\obd(\cals),\obd(\cals-S))$. Since by hypothesis we have
$X''(\oldOmega,\obd(\cals),\obd(\cals-S))=[\emptyset]$, this is a contradiction, and thus (\ref{karsh}) is proved.

Next we claim:
\Claim\label{after karsh}
If $\oldGamma$ is any component of $\partial(\overline{\oldPsi-\frakZ})$ such that $\oldGamma\cap\frakZ\ne\emptyset$, then (i) $\oldGamma$ is toric, (ii) $|\oldGamma|\cap\tcals \ne\emptyset$, and (iii) $|\oldGamma|\cap\tcals\subset\tcals_0$. 
\EndClaim

To prove \ref{after karsh}, write
$\oldGamma=\oldGamma_1\cup\oldGamma_2$, where
$\oldGamma_1=\oldGamma\cap\obd(\tcals)$ and
$\oldGamma_2=\oldGamma\cap\frakZ\subset\Fr_\oldPsi\frakZ$. 
Then
$\oldGamma_2$ is a union of components of $\frakA(\oldPsi)$, so
that by \ref{tuesa day}, every component of $\oldGamma_2$ is
an annular orbifold. Since $\oldGamma$ is by definition a closed $2$-orbifold, $|\oldGamma_2|$ is a proper subset of $|\oldGamma|$, and hence $|\oldGamma_1|\ne\emptyset$. This is Assertion (ii) of \ref{after karsh}.

We have
$\oldGamma_1\cap\oldGamma_2=\partial\oldGamma_1=\partial\oldGamma_2$. According
to the hypothesis of \ref{after karsh}, we have
$\oldGamma_2\ne\emptyset$; since $\oldGamma$ is connected, it follows
that every component of $\oldGamma_1$ has non-empty boundary. Hence
every component of $|\oldGamma_1|$ has non-empty intersection with
$|\frakZ|\cap\tcals$. But by (\ref{karsh}) we have
$|\frakZ|\cap\tcals\subset\tcals_0$, and hence every component of
$|\oldGamma_1|$ is contained in $\tcals_0$. This means that
$|\oldGamma|\cap\tcals=
|\oldGamma_1|\subset\tcals_0$, which is  Assertion (iii) of
\ref{after karsh}.

To prove Assertion (i) of \ref{after karsh}, first note that
we have already observed that 
%$\oldGamma_2$ is a union of components of $\frakA(\oldPsi)$, so
%that by \ref{tuesa day}, 
every component of $\oldGamma_2$ is
an annular orbifold. It therefore suffices to prove that every
component of $\oldGamma_1$ has Euler characteristic $0$ . If $\oldPi$
is a component of $\oldGamma_1$, we may write $\oldPi
=\oldPi_1\cup\oldPi_2$, where $\oldPi_1=\oldPi\cap(\overline{\oldOmega-\oldSigma})$ and $\oldPi_2=\oldPi\cap(\oldSigma-\frakZ)$. Since $\frakZ$ is a union of components of $\oldSigma$, we have $\partial\oldPi\subset\partial\oldPi_1$ and $\oldPi_1\cap\oldPi_2=\partial\oldPi_1-\partial\oldPi=\partial\oldPi_2$. 
Since $|\oldGamma_1|\subset\tcals_0$ by Assertion (iii), each
component of $|\oldPi_1|$ is a component of
$\overline{\tcals_0\setminus|\oldPhi|}$. According to Remark
\ref{special ouenelitte}, the hypothesis that $S$ is \doublefull\
implies that every component of
$\overline{\obd(\tcals_0)\setminus\oldPhi}$ is an annular orbifold,
and hence every component of
$\oldPi_1$ is  annular. Now suppose that $\frakK$ is a component of
$\oldPi_2$. Then $\frakK$ is in particular a component of $\oldPhi$,
so that $\chi(\frakK)\le0$ by \ref{tuesa day}.  If
$\chi(\frakK)<0$, then since
$|\frakK|\subset|\oldGamma_1|\subset\tcals_0$, we have
$\frakK\subset\obd(\tcals_0)\cap\oldPhi^-$, which by definition
implies $\frakK\subset\frakZ$. This is impossible, since
$\frakK\subset\oldPi_2\subset\oldSigma-\frakZ$. Hence
$\chi(\frakK)=0$. Since the components of $\oldPi_1$ and $\oldPi_2$
have Euler characteristic $0$, and meet in a union of common boundary
components, we have $\chi(\oldPi)=0$, and the proof of \ref{after karsh} is complete.

Now we claim:
\Claim\label{on accouna occupani}
If $\frakJ$ is any component of $\overline{\oldPsi-\frakZ}$ such that $\frakJ\cap\obd(\tcals_0)\ne\emptyset$, then $\frakJ$ is a \torifold\ and $\frakJ\cap\obd(\tcals)\subset\obd(\tcals_0)$.
\EndClaim

To prove \ref{on accouna occupani}, first choose a component $\tS$ of $\tcals_0$ such that $\frakJ\cap\obd(\tS)\ne\emptyset$. Note that $\tS$ is \full\ since $S$ is
\doublefull. It therefore follows from 
Remark \ref{special ouenelitte}
that every component of $\overline{\obd(\tS)\setminus\oldPhi}$ is an annular orbifold; since $\chi(\obd(\tS))<0$ by \ref{boundary is negative} (and the componentwise strong simplicity and componentwise boundary-irreducibility of $\oldPsi$, see \ref{doublemint}), it follows that
$\obd(\tS)$
contains at least one component of $\oldPhi^-$, and hence
$\frakZ\cap\obd(\tS)\ne\emptyset$. 
%Since $|\frakJ|\cap\tS\ne\emptyset$, 
It follows that the subsurface $|\frakJ|\cap\tS$ of $\tS$ is proper as well as non-empty, and must therefore have non-empty frontier
 since $\tS$ is connected.
This implies that
the component $\frakJ$ of $\overline{\oldPsi-\frakZ}$ has non-empty
intersection with ${\frakZ}$. Fix a component $\oldGamma$ of
$\partial\frakJ$ such that $\oldGamma\cap\frakZ\ne\emptyset$. By
\ref{after karsh}, $\oldGamma$ is toric and
$\oldGamma\cap\obd(\tcals)\subset\obd(\tcals_0)$. 
On the other hand, since $\oldGamma$ is a component of $\partial(\overline{\oldPsi-\frakZ})$ whose intersection with $\frakZ$ is non-empty, $\oldGamma$ must contain a component of $\Fr_\oldPsi\frakZ$, which is in particular a component of $\frakA(\oldPsi)$, and is therefore annular and $\pi_1$-injective by \ref{tuesa day}; this implies that the inclusion homomorphism $\pi_1(\oldGamma)\to\pi_1(\oldPsi)$ has infinite image. We also know that $\oldGamma$ is two-sided since it is a component of $\partial(\overline{\oldPsi-\frakZ})$.
Hence
by 
Corollary \ref{preoccucorollary}
%Proposition \ref{preoccupani} 
(and the componentwise strong simplicity of $\oldPsi$),  $\oldGamma$ bounds a
\torifold\ $\frakJ'$ in $\oldPsi$. According to \ref{after karsh}, $|\oldGamma|$ has non-empty intersection with $\tcals=\partial|\oldPsi|$. Hence we cannot have $\frakJ'\cap\frakJ=\oldGamma$; we must therefore have
%We have $\partial|\frakJ'|=|\oldGamma|\subset|\partial|\frakJ|$, and since %$|\frakJ|\cap\partial|\oldPsi|\supset|\oldGamma|\cap\tcals$, it follows from \ref{after karsh} that %$|\frakJ|\cap\partial|\oldPsi|\ne\emptyset$. Hence 
$|\frakJ|\subset|\frakJ'|$, and any component of $\partial|\frakJ|$ distinct from $|\oldGamma|=\partial|\frakJ'|$ must be contained in $\inter|\frakJ'|\subset\inter|\oldPsi|$. But it follows from \ref{after karsh} that every component of $\partial|\frakJ|$ has non-empty intersection with $\tcals=\partial|\oldPsi|$. Hence we must have $\partial|\frakJ|=|\oldGamma|$ and therefore 
$\frakJ'=\frakJ$. Thus $\frakJ$ is a \torifold. Since
$\partial\frakJ=\oldGamma$, we have $\frakJ\cap\obd(\tcals)=\oldGamma\cap\obd(\tcals)\subset\obd(\tcals_0)$, and \ref{on accouna occupani} is proved.

Now let $\frakZ^*$ denote the union of all components of $\overline{\oldPsi-\frakZ}$ which meet $\obd(\tcals_0)$. 
According to \ref{on accouna occupani}, we have
$\frakZ^*\cap\obd(\tcals)\subset\obd(\tcals_0)$, and according to (\ref{karsh}) we have
$\frakZ\cap\obd(\tcals)\subset\obd(\tcals_0)$.
 We have $\Fr_\oldPsi\frakZ^*\subset\Fr_\oldPsi\overline{\oldPsi-\frakZ}=\Fr_\oldPsi\frakZ$. On the other hand, if $\frakB$ is any component of $\Fr\frakZ$, then $\frakB$ is in particular a
%Every component of $\Fr_\oldPsi\frakZ$ or $\Fr_\oldPsi\frakZ^*$ is a 
component of $\frakA(\oldPsi)$, and is therefore annular. We have
%If $\frakB$ is any component of $\Fr\frakZ^*$, and if $\frakJ$ denotes the component of $\overline{\oldPsi-\frakZ}$ containing $\frakB$, we have 
 $\emptyset\ne\partial\frakB\subset
\frakZ\cap\obd(\tcals)\subset\obd(\tcals_0)$; thus the component of $\overline{\oldPsi-\frakZ}$ containing $\frakB$ is contained in $\frakZ^*$. This shows that
$\Fr_\oldPsi\frakZ^*=\Fr_\oldPsi\frakZ$. Since $\frakZ$ and $\frakZ^*$ have disjoint interiors, it follows that 
$\Fr_\oldPsi(\frakZ^*\cup\frakZ)=\emptyset$. Thus $\oldPsi':=
\frakZ^*\cup\frakZ$ is a union of components of $\oldPsi$. We have 
$\oldPsi'\cap\obd(\tcals)=
(\frakZ^*\cap\obd(\tcals))\cup(\frakZ\cap\obd(\tcals))\subset\obd(\tcals_0)$. But the definition of  $\frakZ^*$ gives $\obd(\tcals_0)=(\obd(\tcals_0)\cap\frakZ)\cup(\obd(\tcals_0)\cap\overline{\oldPsi-\frakZ})\subset\frakZ\cup\frakZ^*=\oldPsi'$. Hence
$\oldPsi'\cap\obd(\tcals)=\obd(\tcals_0)$.

Since  $\oldPsi'$ is a union of components of $\oldPsi$, and since
$\tcals_0=|\oldPsi'|\cap\tcals$ is invariant under
$|\tau|=\tau_\cals$, the manifold  $|\rho(\oldPsi')|\subset|\oldOmega|$ is closed. Since $\oldOmega$ is connected it follows that $\rho_\cals(\oldPsi')=\oldOmega$ and hence that $\oldPsi'=\oldPsi$. In particular, $\obd(\tcals)=\oldPsi'\cap\obd(\tcals)=\obd(\tcals_0)$, and therefore $\cals=S$.
\EndProof
%\tcals_0\tau\oldGamma\frakG\frakB\frakC\frakD\frakE\oldPhi^-\oldSigma(i\oldPi
%$P<mmmmmm \rho cals_0

\begin{definition}\label{self-clash def}
Let $\Mh$ be a closed, orientable hyperbolic $3$-orbifold, and set $\oldOmega=(\Mh)\pl$. Let $\cals$ be an
$\oldOmega$-admissible system of spheres in $M:=|\oldOmega|$. Set
$\tau=\tau_{\obd(\cals)}$, $\oldPsi=\oldOmega\cut{\obd(\cals_)}$ and $\tS=\partial M\cut \cals=\partial|\oldPhi|$.
A component 
$S$  of $\cals$ will be called a {\it \centralclashsphere} if it satisfies at least one of the following conditions:
\begin{enumerate}[(i)]
\item one of the sides of $S$ contains an \bad\ component (see \ref{stop beeping}) of $|\oldPhi(\oldPsi)|$;
\item both sides of $S$ are pseudo-belted;
\item the sides $\tS$ and $\tS'$ of $S$ are both belted, and $\obd(G_\tS)$ and $\tau (\obd(G_{\tS'}))$ are not (orbifold)-isotopic in $\obd(\tS)$; or
\item $S$ has a pseudo-belted side $\tS$ and a belted side $\tS'$,  and either (a) $\tau (\obd(G_{\tS'}))$
is not isotopic in $\obd(\tS)$ to a suborbifold of 
$\obd(E_\tS)$, or (b) 
$\tau (\obd(G_{\tS'}))$
is isotopic in $\obd(\tS)$ to a suborbifold $\frakC$ of $\obd(E_\tS)$, and
$\frakC$ is not isotopic in $\obd(E_\tS)$ to
$\epsilon_\tS(\frakC)$.
\end{enumerate}
\end{definition}

\Proposition\label{no wonder 1}
Let $\Mh$ be a closed, orientable
hyperbolic $3$-orbifold, and set $\oldOmega=(\Mh)\pl$.
Let $\cals$ be an
$\oldOmega$-\dandy\ system of spheres
in $M:=|\oldOmega|$ having at least two components.
%$\cals_1$ be a (not necessarily proper) subsystem of $\cals$, and set
%$\oldPsi_1=\oldOmega\cut{\obd(\cals_1)}$.  
Then every \doublefull\ 
%Let $S$ be a    component
%of $\cals$ which is a
 \centralclashsphere\ of $\cals$
%. Then $S$
is a \clashsphere\ for $\cals$.
\EndProposition

\Proof
Set
$\oldPsi=\oldOmega\cut{\obd(\cals)}$, $N=|\oldPsi|$, $\tcals=\partial
N$,
$\rho=\rho_\cals:N\to M$,
$\oldSigma=\oldSigma(\oldPsi)$, $\oldPhi=\oldPhi(\oldPsi)$ and $\oldPhi^-=\oldPhi^- (\oldPsi)$ (see \ref{what's iota?}).
Set $\tau=\tau_{\obd(\cals)}:\obd(\tcals)\to \obd(\tcals)$, so that $|\tau|=\tau_\cals:\tcals\to\tcals$.

Suppose that $S$ is a \doublefull\ \centralclashsphere\ of $\cals$. We
are required to prove that $S$ is a \clashsphere\ for $\cals$. Set
$n=\wt S$, so that each side of $S$ has weight $n$. Set
$\tcals_0=\rho^{-1}(S)$.

Set
$X=X(\oldOmega,\obd(\cals),\obd(\cals-\tcals_0))$ and
 $X''=X''(\oldOmega,\obd(\cals),\obd(\cals-\cals_0))$
(in the notation of Section \ref{vegematic section}). By definition we have 
$X''=X\wedge[\obd(\tcals_0)]$.  For each $n\ge1$, set $X_n=X_n(\oldOmega,\obd(\cals),\obd(\cals-S))$.

Since $\cals$ has at least two components, and $S$ is \doublefull, it follows from Lemma
\ref{dun da-dun dun} that $X''\ne[\emptyset]$. Since by definition we
have $X''\in\barcaly_-((\obd(\tcals_0))$, it follows that
$\chi(X'')<0$, i.e.
\Equation\label{not nothing}
\chibar(X'')>0.
\EndEquation

Since $S$ is a \centralclashsphere, one of the alternatives (i)---(iv) of Definition \ref{self-clash def} holds. We will divide the proof into two cases: the case in which  one of the alternatives (i) or (ii) holds, and the case in which  one of the alternatives (iii) or (iv) holds.

{\bf Case I: One of the alternatives (i) or (ii) of Definition \ref{self-clash def} holds.}

In this case we claim:
\Claim\label{gejumpt}
Each side $\tS$ of $S$ contains a component $F$ of $|\oldPhi|$ such that (a) $F$ splinters $\obd(\tS)$, and (b) the component $P$ of $|\oldSigma|$ containing $F$ satisfies $P\cap\tcals\subset\tcals_0$. \EndClaim

To prove \ref{gejumpt}, first notice that if Alternative (ii) of Definition \ref{self-clash def} holds, then each side $\tS$ of $S$ is pseudo-belted. If we set $F=E_\tS$, then by definition $F$ is a component of
  $|\oldPhi|$, contained in $\tS$ and
  splintering $\obd(\tS)$, and the component $P=\frakR_\tS$ of $|\oldSigma|$ containing $F$ satisfies
$P\cap\tcals=F\subset\tS\subset\tcals_0$ (cf. \ref{stop beeping}). Now suppose that Alternative (i) holds, and let $F_1$ be an \bad\ component of $|\oldPhi|$ contained in a component $\tS_1$ of $\tcals_0$. Let $P$ denote the component of $|\oldSigma(\oldPsi)|$ containing $F_1$. 
Then according to Lemma \ref{i'm glad you mentioned foos}, applied with $\tS_1$ and $F_1$ playing the roles of
$\tS$ and $F$,  we have $P\cap\tcals=F_1\cup F_2$ for some \bad\
component $F_2$ of $|\oldPhi|$ contained in
$\tS_2:=|\tau|(\tS_1)$.
The definition of a \bad\ component then
implies that $F_i$ splinters $\obd(\tS_i)$ for $i=1,2$. Each side
$\tS$ of $S$ has the form $\tS_i$
for some $i\in\{1,2\}$, and the assertions of \ref{gejumpt} then follow upon setting $F=F_i$. Thus \ref{gejumpt} is proved.

To prove the conclusion of the lemma in Case I, we will show that $\chibar(X'')<\chibar(\obd(\tcals_0))=2
\chibar(\obd(S))$. Since we have $\chibar(X'')>0$ by (\ref{not
  nothing}), and $S$ is \doublefull, it will then follow from 
Lemma \ref{tell me why} that $S$ is a \clashsphere, as required.

Since $X''=X\wedge[\obd(\tcals_0)]
%According to the definition (see \ref{oldXi} we have  %$X''=X''(\oldOmega,\obd(\cals),\obd(\cals-\tcals_0))
\preceq
[\obd(\tcals_0)]$, Proposition \ref{new partial order} gives $\chibar(X'')
\le \chibar(\obd(\tcals_0))
=2\chibar(\obd(S))$. For the rest of the argument in Case I we will assume 
$\chibar(X'')=\chibar(\obd(\tcals_0))$, and the goal will be to obtain
a contradiction. 
%Let
%$\oldXi''\in\Theta_-(\obd(\tS))$ denote a representative of the
%equivalence class
%$X''\wedge \tS\in\barcaly_-(\obd(\tcals_0))$.
Let
$\oldXi''\in\Theta_-(\obd(\tcals_0))$ denote a representative of the
equivalence class
$X''\in\barcaly_-(\obd(\tcals_0))$.

By Lemma \ref{old partial 
  order}, the assumption $\chibar(X'')=\chibar(\obd(\tcals_0))$
implies:
\Claim\label{there's the rub}
Every component of $\obd(\tcals_0-\inter\oldXi'')$ is annular.
\EndClaim

%setting
%$n=\wt S$, and 
% for each
%side $\tS$ of $S$, \redcomment{That's pretty confusing, since $X''$ is
%  defined (I think) to be $X\cap(\tS\cup\tau(\tS))$ in the first
%  place. The notation can stand some improvement} 
We claim:
\Claim\label{abagli} 
For each side $\tS$ of $S$, there is a component $\frakZ$ of
$\oldXi''$, contained in $\obd(\tS)$, such that $|\frakZ|$ semi-splinters $\obd(\tS)$.
\EndClaim

To prove \ref{abagli}, we 
%set $n=\wt S=\wt\tS$, and we
 apply  
Lemma \ref{decomp}, letting $\obd(\tS)$ and $\tS$ play the respective roles of $\oldTheta$ and $S$, and taking $\frakE=\obd(\tS)\setminus\oldXi''$.
%$X=\tS\cap\fraks_\oldPsi$ \redproofreadingnote{$X$ is a bad letter %here. Try to change the notation in Lemma \ref{decomp} and %$F=|\tS-\inter\oldXi''(\tS)|$. 
According to that
lemma, either the assertion of \ref{abagli} is true 
%(i.e. there is a component $\frakZ$ of $\oldXi''(\tS)$ such that every
%component of $\tS-\inter|\frakZ|$ has weight at most $n/2$), 
or 
%Let $X$ be a finite subset of a $2$-sphere $S$,  and let $F$ be a
%%compact (but possibly disconnected) $2$-dimensional submanifold of
%%%$S$ whose boundary is disjoint from $X$. Then either 
%\begin{itemize}
 there is a component $T$ of $\tS\setminus\inter|\oldXi''| $ which splinters $\obd(\tS)$. 
If the latter alternative holds, it follows from  Lemma \ref{negative splinter} (and \ref{doublemint}) that $\chi(\obd(T))<0$. (The tautness of the component $\obd(T)$ of $\obd(\tcals)-\inter\oldXi''$, which is required for the application of Lemma \ref{negative splinter}, follows from the tautness of $\oldXi''$, which is in turn included in the fact that $\oldXi''$ is an element of $\Theta_-(\obd(\tcals))$.)  But $\obd(T)$ is annular by \ref{there's the rub}. This contradiction proves \ref{abagli}.
%Suppose
 %that there is a there is a component $T$ of
%% $|\tS-\inter\oldXi''(\tS)| $ having the latter property. It now
 %follows from \ref{there's the rub}
%Since  every component of $\obd(\tcals_0-\inter X'')$ is annular,
%But the condition $\chibar(X'')=\chibar(\obd(S))$ implies that every
%component of %$\tS-\inter\oldXi''(\tS) $ is annular, so
% that
 %$T$ is either a weight-$0$ annulus or a
%weight-$2$ disk. If $T$ is a weight-$0$ annulus then $\tS-\inter T$ %has two components, the sum of whose weights is $n$; thus they %cannot each have weight strictly less than $n/2$. If $T$ is a
%weight-$2$ disk, then $\tS-\inter T$ is a single disk $D$ of weight
%$n-2$. If $n-2<n/2$ then $n\le3$ and hence $\wt D\le1$, which is
%impossible since $\partial D$ is a boundary component of
%$\oldXi''(\tS)\in\Theta_-(\obd(\tcals))$, which is  a taut
%suborbifold of $\obd(\tS)$ by the
%definition of $\Theta_-(\obd(\tcals))$. Thus \ref{abagli} is proved.

Now choose a side $\tS_0$ of $S$. By \ref{abagli}, we may fix a component $\frakZ$ of
$\oldXi''$, contained in $\obd(\tS_0)$, such that $|\frakZ|$ semi-splinters
$\obd(\tS_0)$. Since $[\frakZ]\preceq[\oldXi'']=X''\preceq X$, and
since $[\frakZ]$ is connected, it follows from \ref{components} that
$[\frakZ]\preceq Y$ for some component $ Y$ of
$X=X(\oldOmega,\obd(\cals),\obd(\cals-S))$. By  Lemma \ref{here
  goethe}, $ Y$ is a component of
$X_n=X_n(\oldOmega,\obd(\cals),\obd(\cals-S))$ for some integer
$n\ge1$. 
%By \ref{oldXi} it follows that $Y$ is a component of $V_n$. 
%In particular we have $[\frakZ]\preceq X_n$. 
Let  $n_0$ denote the smallest positive integer such that $Y\preceq
X_{n_0}$. The definition of
$X_{n_0}(\oldOmega,\obd(\cals),\obd(\cals-S))$ (see \ref{oldXi})
implies that $Y\preceq V_{n_0}:=V_{n_0}^{\oldOmega,\cals}$, which by
(\ref{i thought so too}) implies that for every $m$ with $0<m\le n_0$
we have $Y\preceq V_{m}:= V_{m}^{\oldOmega,\cals}=\dom \noodge_{m}$,
where $\noodge_{m}=\noodge_{m}^{\oldOmega,\cals}$. We claim:
\Claim\label{goose egg}
For any index $m$ with $0<m<n_0$, we have $\noodge_m(Y)\preceq[\obd(\tcals_0)]$.
\EndClaim
To prove \ref{goose egg}, let $Y^+$ denote
the component of $V_m$ such that $Y\preceq Y^+$ (see
\ref{components}). Note that since $Y^+$ is connected, we
have either $\noodge_{m}(Y^+)\preceq[\obd(\tcals_0)]$ or
$\noodge_{m}(Y^+)\preceq[\obd(\tcals-\tcals_0)]$. If
$\noodge_{m}(Y^+)\preceq[\obd(\tcals-\tcals_0)]$, then by \ref{oldXi} we have $Y^+\in\oldfrakX_m
(\oldOmega,\cals,\cals-S)$ and hence $Y\preceq Y^+\preceq X_{m}$, a
contradiction to the minimality of $n_0$. Hence
$\noodge_m(Y)\preceq\noodge_m( Y^+)\preceq[\obd(\tcals_0)]$, and \ref{goose egg} is proved. 
%\tcals_0

%if $n>1$. 
Let us set $Q=[\tau](\noodge_{n_0-1}([\frakZ]))\in\barcaly_-(\obd(\tcals))$ if $n_0>1$, and $Q=[\frakZ]$ if $n_0=1$. Setting $\iota=\iota_\oldPsi$, we claim:
\Claim\label{easy does it}
We have $Q\preceq[\obd(\tcals_0)]$, and if $\frakU$ is an element of $\Theta_-(\obd(\tcals))$ with $[\frakU]=Q$, then $|\frakU|$ semi-splinters the component of $\obd(\tcals_0)$ containing it. 
Furthermore, we have $Q\preceq[\oldPhi^-]$, and $[\iota](Q)\preceq[\obd(\tcals-\tcals_0)]$. 
\EndClaim

When $n_0=1$, since $Q=[\frakZ]$, the first sentence of \ref{easy does it} is immediate from our choice of $\frakZ$. To prove the second sentence when $n_0=1$, note that according to
\ref{oldXi} and (\ref{and so did you}) we have $Q=[\frakZ]\le Y\preceq X_{n_0}=X_1\preceq V_1=[\oldPhi^-]$. The order-preserving property of \ref{restriction} then gives
$[\iota](Q)=\noodge_1([\frakZ])\preceq\noodge_1(X_1)$. If $R$ denotes the component of $X_1$ such that 
$[\iota](Q)\preceq\noodge_1(R)$ (see \ref{components})
then \ref{oldXi} gives $R\in
\oldfrakX_1(\oldOmega,\cals,\cals-S)$ and hence
%we have the definition of %$X_1=X_1(\oldOmega,\cals,\cals-S)\in guarantees that 
$\noodge_1(R)\preceq
\obd(\tcals-\tcals_0)$.
We therefore have $[\iota](Q)\preceq \obd(\tcals-\tcals_0)$, and \ref{easy does it} is proved in this case.

To prove \ref{easy does it} when $n_0>1$, first notice that  the case
$m=n_0-1$ of \ref{goose egg} gives 
$\noodge_{n_0-1}([\frakZ])\preceq \noodge_{n_0-1}(Y)
\preceq[\obd(\tcals_0)]$, and hence $Q\preceq[\obd(\tcals_0)]$.
% since $\frakZ$ is connected, we have either $\noodge_{n_0-1}([\frakZ])\preceq[\obd(\tcals_0)]$ or %$\noodge_{n_0-1}([\frakZ])\preceq[\obd(\tcals-\tcals_0)]$. If %$\noodge_{n_0-1}([\frakZ])\preceq[\obd(\tcals-\tcals_0)]$, 
%the definition of $X_{n_0-1}=X_{n_0-1}(\oldOmega,\cals,\cals-S)$ implies that $[\frakZ]\preceq %X_{n_0-1}$, a contradiction to the minimality of $n$. Hence 
 %$\noodge_{n_0-1}([\frakZ]) \preceq[\obd(\tcals_0)]$, and therefore $Q=[\tau](\noodge_{n_0-1}([\frakZ])) %\preceq[\obd(\tcals_0)]$.
Next,
for $m=1,\ldots,n$, 
let $\otherdge_m:\frakZ\to\obd(\tcals)$ be an embedding such that
% element of $\calm(\obd(\tcals),\obd(\tcals))$ such that 
$[\frakZ,\otherdge_m,\otherdge_m(\frakZ)]=\noodge_m|[\frakZ]$. After modifying the element $\frakU$  of $\Theta_-(\obd(\tcals))$ which appears in the statement of \ref{easy does it}  within its isotopy class, we may assume that $\frakU=\tau(\otherdge_{n_0-1}(\frakZ))$.
%for
%some $\otherdge_{n_0-1}\in\calm(\tcals,\tcals)$ with
%$[\otherdge_{n_0-1}]=\noodge_{n_0-1}$. 
In order to show that $\frakU$ semi-splinters the component of
$\obd(\tcals_0)$ containing it, we apply Lemma \ref{lemmanade},
letting $n_0-1$ play the role of $n$ in that lemma, and letting $Q$
play the role of $Z$. Since, in the notation introduced above, we have
$[\frakZ]=Q$ and $[\otherdge_m]=\noodge_m$,
we may take the element  $\frakZ$ of $\Theta_-(\obd(\tcals)$ and the
embedding $\otherdge_m:\frakZ\to\obd(\tcals)$ given by Lemma
\ref{lemmanade} to be the ones introduced above. The sphere denoted
$\tS_0$ in Lemma \ref{lemmanade} coincides with the sphere denoted
$\tS_0$ above. As in the statement of Lemma \ref{lemmanade}, we denote
by $\tS_m$ 
the component of $\tcals$ containing
$|\otherdge_m(\frakZ)|$ for $m=1,\ldots,n_0-1$. The subsurface $|\frakZ|$ semi-splinters  $\obd(\tS_0)$ by our choice of $\frakZ$. It follows from \ref{goose egg} that $\tS_1,\ldots,\tS_{n_0-1}$ are all components of $\tcals_0$, i.e. are sides of $S$, and therefore have the same weight as $\tS_0$, as required for the application of Lemma \ref{lemmanade}. 
%=\tS$, $\tS_1=\tS'$, and 
%$\noodge_1(\frakZ)=\obd(F')$ 
%\redcomment{Where does this go? ``Since $Q\preceq[\obd(\tcals_0)]$, the orbifold $\frakU$ is contained in a component $\tS$ of $\tcals_0$.''}
According to Assertion (2) of
Lemma \ref{lemmanade}, $|\otherdge_{n_0-1}(\frakZ)|$ semi-splinters 
%the
%component 
$\obd(\tS_{n_0-1})$.
% of $\tcals$ containing
%it. 
%\redcomment{I guess the hypothesis in Lemma \ref{lemmanade} that
%  $\wt\tS_m=\wt\tS_0$ for $m=1,\ldots,n$ (I think it's $n-1$ here, not
  %$n$) does hold in this context, because the sphere are all the
  %same. I need to point that out, though. It requires generalizing the
%argument used to prove $Q \preceq[\obd(\tcals_0)]$.} 
Since $\tau$ is an involution of $\obd(\tcals)$, it follows that
$|\frakU|=|\tau(\otherdge_{n_0-1}(\frakZ))|$ semi-splinters $\tau(\obd(\tS_{n_0-1}))$.
% the component of $\tcals$ containing it.

Now note that according to the definitions (see \ref{in duck tape}),
and \ref{more associativity}, we
have $\noodge_{n_0}=[\iota]\diamond[\tau]\diamond\noodge_{n_0-1}$. Since
%$V_{n_0-1}=\dom\noodge_{n_0-1}$ and 
$V_{n_0}=\dom\noodge_{n_0}$, the relation
 $[\frakZ]\preceq V_{n_0}$ may be rewritten as $[\frakZ]\preceq\dom
([\iota]\diamond[\tau]\diamond\noodge_{n_0-1})$. Lemma \ref{before
  associativity} then gives 
%that $[\frakZ]\preceq \dom\noodge_{n_0-1}$,
that 
$\noodge_{n_0-1}([\frakZ])\preceq\dom([\iota]\diamond[\tau])$, and that
$\noodge_{n_0}|[\frakZ]=([\iota]\diamond[\tau])\circ(\noodge_{n_0-1}|\frakZ)$. A
second application of
Lemma \ref{before
  associativity} then gives that
%$\noodge_{n_0-1}([\frakZ])\preceq\dom
%([\iota]\diamond
%[\tau]$, that
$[\tau](\noodge_{n_0-1}([\frakZ]))\preceq\dom[\iota]=[\oldPhi^-]$, and that
$([\iota]\diamond[\tau])|\noodge_{n_0-1}([\frakZ])=[\iota]\circ([\tau]|\noodge_{n_0-1}([\frakZ])$.
In view of the definition of $Q$, it follows that
$Q\preceq[\oldPhi^-]$ and that
$[\iota](Q)=\noodge_{n_0}([\frakZ])$. But we have $[\frakZ]\preceq Y\preceq
X_{n_0}$. If $Z^+$ denotes the component of $X_{n_0}$ such that $[\frakZ]\preceq Z^+$ (see \ref{components}), then by \ref{oldXi} we have $Z^+\in\oldfrakX_{n_0}(\oldOmega,\cals,\cals-S)$, so that
%  and the definition of %$X_{n_0}=X_{n_0}(\oldOmega,\cals,\cals-S)$ %guarantees that 
$\noodge_{n_0}(X_{n_0})\preceq
\obd(\tcals-\tcals_0)$. Hence $[\iota](Q)\preceq \obd(\tcals-\tcals_0)$, and the proof of \ref{easy does it} is complete.

Now fix an element $\frakU$ of $\Theta_-(\tcals)$ with $[\frakU]=Q$,
and let $\tS$ denote the component of $\tcals$ containing
$|\frakU|$. It follows from \ref{easy does it} that
$\tS\subset\tcals_0$. Hence, according to \ref{gejumpt}, $\tS$
contains a component $F$ of $|\oldPhi|$ such that (a) $F$ splinters
$\obd(\tS)$, and (b) the component $P$ of $|\oldSigma|$ containing $F$
satisfies $P\cap\tcals\subset\tcals_0$. It follows from Lemma
\ref{negative splinter} that $\chi(\obd(F))<0$, so that $\obd(F)$ is a
component of $\oldPhi^-$. Hence by \ref{what's iota?}, $\obd(P)$ is a
component of $\oldSigma^-$. In particular $\obd(P)$ is saturated in a
fibration $q$ of the kind described in \ref{what's iota?}, so that
according to the discussion in \ref{what's iota?}, $
\obd(P\cap\tcals)$ is $\iota$-invariant. It follows that
$|\iota|(F)\subset P\cap\tcals
%$, and hence
%\redmissingref{basically the def. of $\iota$ and the $I$-fibration of $\obd(U)$. Another $\iota$ issue} it follows that
%$|\iota|(F)
\subset\tcals_0$. But \ref{easy does it} implies that $|\frakU|$ is contained in a component $F'$ of $|\oldPhi|$ such that $|\iota|(F')\subset\tcals-\tcals_0$, and that $F'$ semi-splinters $\obd(\tS)$. The components $F$ and $F'$ of $|\oldPhi|$ are distinct since
$|\iota|(F)\subset\tcals_0$ and $|\iota|(F')\subset\tcals-\tcals_0$. Hence $\frakU\cap F=\emptyset$. Since $F$ splinters $\obd(\tS)$ and $Q$ semi-splinters $\obd(\tS)$, this contradicts Lemma \ref{can't have it both ways}, and the proof of the present lemma is complete in Case I.
%\iotaV_nY\oldUpsilon\frakC

{\bf Case II. One of the alternatives (iii) or (iv) of Definition \ref{self-clash def} holds.}

%\redcomment{I need to fix this up to cover both alternatives, and i think the difference between them will be fairly small. Old version: ``[Suppose that Condition (iii) of Definition \ref{self-clash def} holds]: thus the sides $\tS$ and $\tS'$ of $S$ are both belted, and $\obd(G_\tS)$ and $\obd (\tau (G_{\tS'}))$  are not orbifold-isotopic in $\obd(\tS)$. Fix boundary components $C$ and $C'$ of $G_\tS$ and $\tau(G_{\tS'})$ respectively. We claim:''}

We shall refer to the subcases in which Alternative (iii), Sub-alternative (a) of Alternative (iv), and Sub-alternative (b) of Alternative (iv) of Definition \ref{self-clash def} hold as Subcases (iii), (iv-a), and (iv-b) respectively.

If we are in
 Subcase (iv-b), we may assume $\oldSigma$ to have been chosen within
 its (orbifold)-isotopy class so that $\tau(\obd(G_{\tS'}))\subset
\inter  \obd(E_\tS)$, and so that $\tau (\obd(G_{\tS'}))$
is not isotopic in $\obd(E_\tS)$ to 
$\epsilon_\tS (\tau(\obd  (G_{\tS'})))$. We claim:

\Claim\label{margaret dumont}
In Subcase (iv-b), the suborbifolds $\tau (\obd(G_{\tS'}))$
and $\epsilon_\tS (\tau(\obd  (G_{\tS'})))$ of $\obd(E_\tS)$ are not isotopic in $\obd(\tS)$.
\EndClaim

To prove \ref{margaret dumont}, assume that $\tau (\obd(G_{\tS'}))$ and $\epsilon_\tS (\tau(\obd  (G_{\tS'})))$  are  isotopic in $\obd(\tS)$. Then according to Proposition \ref{before i guess}, either (A) 
 $\tau (\obd(G_{\tS'}))$ and $\epsilon_\tS (\tau(\obd  (G_{\tS'})))$
%$\epsilon_\tS (\obd  (G_{\tS'}))$ of $\obd(E_\tS)$ 
are  isotopic in $\obd(E_\tS)$, or (B) each of the manifolds 
$|\tau| (G_{\tS'})$ and 
$|\epsilon_\tS| (|\tau|(G_{\tS'}))$
is a weight-$0$ annulus, and has a core curve which cobounds a
weight-$0$ annulus with some boundary component of $E_\tS$. We have
observed that (A) cannot occur in Case (iv-b). Now suppose that (B)
holds. Fix a component $L$ of $\partial E_\tS$ which cobounds a
weight-$0$ annulus in $E_\tS$ with a core curve of $|\tau|
(G_{\tS'})$, and a component $L'$ of $\partial E_\tS$ which cobounds a
weight-$0$ annulus in $E_\tS$ with a core curve of $|\epsilon_\tS|
(|\tau|(G_{\tS'}))$. Since  $\tau (\obd(G_{\tS'}))$ and $\epsilon_\tS (\tau(\obd
(G_{\tS'})))$ are isotopic in $\obd(\tS)$, the $1$-suborbifolds
$\obd(L)$ and $\obd(L')$ are isotopic in $\obd(\tS)$; hence $L$ and
$L'$ are isotopic in the $2$-manifold $\tS-\fraks_{\obd(\tS)}$. (It is worth
mentioning that $\obd(L)$ and $\obd(L)$ are technically the same objects as $L$ and $L'$. We use the notation $L$ and $L'$ when we think of them as $1$-submanifolds of $\tS-\fraks_{\obd(\tS)}$, and we use $\obd(L)$ and $\obd(L')$ when we think of them as suborbifolds of $\obd(\tS)$ which are $1$-manifolds.) If $L=L'$ then  $\tau (\obd(G_{\tS'}))$ and $\epsilon_\tS (\tau(\obd  (G_{\tS'})))$
%$\epsilon_\tS (\obd  (G_{\tS'}))$ of $\obd(E_\tS)$ 
are  isotopic in $\obd(E_\tS)$, and we again have a contradiction. Now suppose that $L\ne L'$. Since $L$ and $L'$ are disjoint and are isotopic in $\tS-\fraks_{\obd(\tS)}$, it follows from  \cite[Lemma 2.4]{epstein} that they cobound an annulus $A\subset \tS-\fraks_{\obd(\tS)}$. 

On the other hand, since $\tS$ is a $2$-sphere, there are disk components $D$ and $D'$ of $\tS-\inter E_\tS$ with $\partial D=L$ and $\partial D'=L'$. The unique annulus contained in $\tS$ and bounded by $L\cup L'$ is $\tS-\inter(D\cup D')$, and hence $\tS-\inter(D\cup D')=A$. In particular we have $\wt(\tS-\inter(D\cup D'))=0$. If $L''$ is a component of $\partial E_\tS$ distinct from $L$ and $L'$, and if $D''$ denotes the disk component of $\tS-\inter E_\tS$ bounded by $L''$, then $D''\subset \tS-\inter(D\cup D')$ and hence $\wt D''=0$; this implies that $\obd(L'')$ is not $\pi_1$-injective in $\obd(\tS)$, a contradiction to the tautness of $\obd(E_\tS)$ (see \ref{what's a belt?}). Hence we must have $\partial E_\tS=L\cup L'$ and therefore $E_\tS=
\tS-\inter(D\cup D')$. This shows that $E_\tS$ is a weight-$0$
annulus,  so that $\chi(\obd(E_\tS)=0$. But since the pseudo-belt
$E_\tS$ splinters $\obd(\tS)$ by definition, we have
$\chi(\obd(E_\tS)<0$ by Lemma \ref{negative splinter}. This
contradiction completes the proof of \ref{margaret dumont}.
%C

In Case II we may label the sides of $S$ as $\tS$ and $\tS'$ in such a way that $\tS'$ is belted. We fix a boundary component $C'$ of $G_{\tS'}$. 
Let us define a 
closed, possibly disconnected $1$-manifold $C\subset\tS$,
disjoint from $\fraks_\oldOmega$ and $\pi_1$-injective in $\obd(\tS)$, as follows.  If we are in Subcase
(iii), so that $\tS$ is belted, we take $C$ to be a boundary component
of $G_\tS$. In  Subcase (iv-a) we set $C=\partial
E_\tS$. If we are in
 Subcase (iv-b), as we have chosen $\oldSigma$ to have been chosen
 within its (orbifold)-isotopy class so that
 $\tau(\obd(G_{\tS'}))\subset \obd(E_\tS)$, we  may define $C$ to be a boundary component of $|\epsilon_\tS|(|\tau|(G_{\tS'}))$.

In Subcases (iii) and (iv-b), the $1$-manifold $C$ is by definition
connected, so that $\size C$ (as well as $\size C'$) is defined.
We claim:

\Claim\label{aznavour} In Subcases (iii)
and (iv-b), we have $\size C= \size C'=n/2$.
\EndClaim 

To prove \ref{aznavour}, first note that Lemma \ref{how do i recognize a belt?} implies
that $\size C'=n/2$ because $C'$ is a boundary component of the belt
$G_{\tS'}$. If we are in Subcase (iii), $C$ is a boundary component of
$G_\tS$, and hence Lemma \ref{how do i recognize a belt?} implies that
$\size C=n/2$.  We now turn to Subcase (iv-b). In this subcase, $C$ is a boundary
component of $|\epsilon_\tS|(|\tau|(G_{\tS'})$, so that
$C=|\epsilon_\tS|(|\tau|(C_0))$ for some component $C_0$ of $\partial
G_{\tS'}$. We have $\size C_0=n/2$ by Lemma \ref{how do i recognize a
  belt?}. By Proposition \ref{kavanaugh},
we have
%$C=|\tau|(C_0)$, so that 
$\size C=\size|\epsilon_\tS|(|\tau|(C_0))
=\size|\tau|(C_0)
=\size C_0=n/2$. This completes the proof of \ref{aznavour}.

In all three subcases, we claim:

\Claim\label{can't split 'em}
The $1$-manifolds $C$ and $|\tau|(C')$ meet essentially (see \ref{meet me in st. louis}) in $\obd(\tS)$.
\EndClaim

To prove \ref{can't split 'em}, first consider Subcase (iv-a).
Suppose that \ref{can't split 'em} is false, i.e. that $C=\partial
E_\tS$ and $|\tau|(C')$ do not meet essentially. Then after
modifying $\oldSigma$ by an isotopy which is constant on all
components of $\oldSigma$ except $\frakR_\tS$, we may assume that
$\partial E_\tS$ and $|\tau|(C')$ are disjoint. 
Since $|\tau|(C')$ is a boundary component of  $|\tau|(G_{\tS'})$, and since $\tau( \obd(G_{\tS'}))$ is annular (and orientable), it follows from \ref{cobound} (and the $\pi_1$-injectivity  of $\partial E_\tS$, see \ref{tuesa day}) that each component of $\partial E_\tS$ either is disjoint from $|\tau|(G_{\tS'})$ or cobounds a weight-$0$ annulus with $\partial(|\tau|(G_{\tS'})$. Hence,
after an additional (orbifold) isotopy, we may assume that 
$\obd(\partial E_\tS)$ and $\tau(\obd(G_{\tS'}))$ are disjoint. But according to Sub-alternative (a) of Alternative (iv), $\tau(\obd( G_{\tS'}))$ is not isotopic in $\tS$ to a suborbifold of $\obd(E_\tS)$. Hence we must have $|\tau|(G_{\tS'})\subset  \tS-E_\tS$.
%, and in particular $\tau(C')\subset\tS- E_\tS$. 
Now
%In this subcase, $\obd(\tau(G_{\tS'}))$ is not isotopic in $\obd(\tS)$ to a %suborbifold of 
%$\obd(E_\tS)$. Hence if B
by definition the pseudo-belt $E_\tS$ splinters $\obd(\tS)$; hence if $D$ denotes the component of $\tS-\inter E_\tS$ containing $|\tau|(G_{\tS})$, we have $\wt D<n/2$. Since $|\tau|(C')\subset|\tau|(G_{\tS})\subset D$, and since $D$ is a disk, it follows that $\size |\tau|(C')<n/2$. But $\size|\tau|(C')=\size C'$, and since $C'$ is a boundary component of the belt $G_\tS$, it follows from Lemma \ref{how do i recognize a belt?} that $\size C'=n/2$. This contradiction establishes  \ref{can't split 'em} in this subcase.

We now turn to the proof of 
 \ref{can't split 'em} in Subcases (iii) and (iv-b). In these subcases, $C$, as well as $C'$, is connected, and we have 
$\size C= \size C'=n/2$ by \ref{aznavour}.

Assume that $C$ and $|\tau|(C')$ do not meet essentially. Then $C$ is isotopic in $\tS\setminus\fraks_\oldOmega$ to a simple closed curve $C_1$ such that $C_1\cap |\tau|(C')=\emptyset$. Hence $C_1$ and $|\tau|(C')$ respectively bound disks $D,D'\subset\tS$ such that $D'\subset\inter D$. Since 
$\size C_1=\size C$ and $\size C'=\size|\tau|(C')$ are both equal to $n/2=\wt(\tS')/2$, the weights of  $ D$ and $ D'$ are both equal to $n/2$. It follows that the annulus $\overline{D-D'}$ has weight $0$, and hence that  $\obd(C)$ and $\tau(\obd(C'))$ are orbifold-isotopic in $\obd(\tS)$.

Now define a subsurface $T$ of $\tS$ by setting $T=G_\tS$ in Subcase
(iii), and setting $T=|\epsilon_\tS|(|\tau|(G_{\tS'}))$ in Subcase
(iv-b). Thus $\obd(T)$ is homeomorphic to either $\obd(G_\tS)$ or
$\obd(G_{\tS'})$, and in either case the definition of a belt
implies that $\obd(T)$ is an (orientable) annular orbifold. Likewise, $\tau (\obd (G_{\tS'}))$ is annular (and orientable). Hence
each of the surfaces $T$ and $|\tau|(G_{\tS'})$ is either a weight-$0$ annulus or a weight-$2$ disk. According to our construction, $C$ and $|\tau|(C')$ are boundary components of $T$ and $|\tau|(G_{\tS'})$ respectively. 

In Subcase (iii), $\obd(G_\tS)$ and $\tau( \obd(G_{\tS'}))$ are not
orbifold-isotopic in $\obd(\tS)$. In Subcase (iv-b), according to \ref{margaret dumont},
$\tau(\obd(G_{\tS'}))$ and $\epsilon_\tS (\tau(\obd (G_{\tS'})))$ are
not orbifold-isotopic in $\obd(\tS)$. 
Thus in either subcase, the definition of $T$ implies that $\obd(T)$ and $\tau(\obd(G_{\tS'}))$ are not orbifold-isotopic in $\obd(\tS)$.

If $T$ and $|\tau|(G_{\tS'})$ are both weight-$0$ annuli, then since $\obd(C)$
and $\tau(\obd(C'))$ are orbifold-isotopic in $\obd(\tS)$, the suborbifolds $\obd(T)$ and $\tau(\obd(G_{\tS'}))$ are also isotopic, a contradiction.

If $T$ and $G_{\tS'}$ are both weight-$2$ disks, then in particular
their boundaries are connected, and hence
$C=\partial T$ and $C'=\partial G_{\tS'}$. Since $\obd(C)=\partial\obd( T)$ and $\tau(\obd(C')) =\partial\tau(\obd(G_{\tS'}))$ are isotopic in
$\obd(\tS)$, the suborbifold $\tau( \obd (G_{\tS'}))$ of
$\obd(\tS)$ is isotopic either to $\obd(T)$ or to
$\obd(\tS-\inter T)$. If $\tau (\obd (G_{\tS'}))$ is isotopic
to $\obd(T)$, we have a contradiction. If
$\tau (\obd (G_{\tS'}))$ is isotopic to $\obd(\tS-\inter T)$, then both 
$\obd(T)$ and $\obd(\tS-\inter T)$ are annular, and hence
$\chi(\obd(S))=
\chi(\obd(\tS))=\chi(\obd(T)+\chi(\obd(\tS-\inter T))=0$; since the admissibility of $\cals$ implies that $\obd(S)$ is incompressible in the $3$-orbifold $\oldOmega$, we have a contradiction to the hyperbolicity of $\Mh$.

Next consider the possibility that $T$ is a weight-$0$ annulus and
$G_{\tS'}$ is a weight-$2$ disk. This situation cannot arise in
Subcase (iv-b), since in that subcase
$T=|\epsilon_\tS|(|\tau|(G_{\tS'}))$ is homeomorphic to
$G_{\tS'}$. Hence we must be in Subcase (iii), and we then have
$T=G_\tS$. Since the component $\obd(C)$ of $\partial \obd( G_\tS)$ is isotopic in
$\obd(\tS)$ to $\tau(\obd(C'))=\partial\tau(\obd(G_{\tS'}))$, there is a disk
 $\Delta \subset\tS$ with $\partial \Delta =C$ such that $\obd(\Delta )$ is isotopic to
$\tau(\obd(G_{\tS'}))$. Since $\obd(G_{\tS'})$ is annular, the
two points of $\Delta \cap\fraks_\oldPsi$ are both of order $2$. Since $\partial
\Delta =C$, we have either $\Delta \cap G_\tS=C$ or $G_\tS\subset \Delta $.
Define a set
$\Delta _1\subset\tS$ by setting $\Delta _1=\Delta $ if $\Delta \cap G_\tS=C$, and
$\Delta _1=\overline{\Delta -G_\tS}$ if $G_\tS\subset \Delta $. Then  $\Delta _1$ is a disk containing exactly two
points of $\fraks_\oldPsi$, both of order $2$. (This is obvious if $\Delta _1=\Delta $, and if $\Delta _1=\overline{\Delta -G_\tS}$ it follows from the fact that
$T=G_\tS$ is a weight-$0$ annulus.)
In particular,
$\obd(\Delta _1)$ is annular.  Furthermore,
$C_1:=\partial \Delta _1=\Delta _1\cap G_\tS$ is a component of $\partial
G_\tS$. (If  $\Delta \cap G_\tS=C$, we have $C_1=C$.) In either subcase we have $\Delta _1\cap G_\tS=C_1$. 
Hence the  disk $\Delta _1$ is a component of $\tS-\inter G_\tS$. 
Since $G_\tS$ is a component of $|\oldPhi(\oldPsi)|$ and $\obd(\Delta _1)$ is annular, this contradicts Corollary \ref{mangenfeffer}, with $\obd(G_\tS)$ playing the role of $\oldGamma$.

There remains the possibility that  $G_{\tS'}$ is a weight-$0$ annulus
and $T$ is a weight-$2$ disk. Again this situation cannot arise in
Subcase (iv-b), since in that subcase
$T=|\epsilon_\tS|(|\tau|(G_{\tS'}))$ is homeomorphic to
$G_{\tS'}$. Hence again we must be in Subcase (iii), and we have
$T=G_\tS$. 
Since the $1$-manifolds $\obd(C)$ and $\tau(\obd(C'))$ are isotopic in
$\obd(\tS)$, the $1$-manifolds $\obd(C')$ and $\tau(\obd(C))$ are isotopic in
$\obd(\tS')$. Hence we may apply the same argument that was used in
the situation where $T$ is a weight-$0$ annulus
and $G_{\tS'}$ is a weight-$2$ disk, but with the roles of $\tS$ and
$\tS'$ reversed, and with the roles of $C$ and
$C'$ reversed, to obtain a contradiction.
%If this occurs, then again we cannot be
%in Subcase
%we obtain the same contradiction as in the situation where
  %$T$ is a weight-$0$ annulus and $G_{\tS'}$ is a weight-$2$
  %disk, by reversing the roles of
  %$T$ and $G_{\tS'}$.
The proof of \ref{can't split 'em} is thus complete.

%If $Y$ is any component of $X''$, then by Lemma \ref{here goethe},    $Y$ is a component of $X_n$ for some $n\ge1$. By \ref{oldXi}, (\ref{i thought so too}), and (\ref{and so did you}), 
%we have $X_n\preceq V_n\preceq V_1=\oldPhi^-$; hence $Y\preceq\oldPhi^-$ for every component $Y$ of $X''$, which by Proposition\ref{mystery} implies that $X''\preceq[\oldPhi^-]$.

Let us now  write $X''=[\oldXi'']$ for some
$\oldXi''\in\Theta_-(\obd(\tcals_0))$. Since, by an observation made in \ref{oldXi}, we have $X''\preceq X\preceq[\oldPhi^-]$, we may take $\oldXi''$ to be contained in $\inter\oldPhi^-$.
We claim:
\Claim\label{zeroth claim}The simple closed curve $C'$ is disjoint from
$|\oldXi''|$, and the closed $1$-manifold $C$ is isotopic in
$\tcals_0-\fraks_{\obd(\tcals_0)}$ to a $1$-manifold disjoint from
$|\oldXi''|$. 
\EndClaim

To prove \ref{zeroth claim}, first note that since $\obd(G_{\tS'})$ is an annular orbifold and a component of $\oldPhi$, the surface $G_{\tS'}$
is disjoint from $|\oldPhi^-|$. Since $C'\subset\partial G_{\tS'}$ and
$\oldXi''\subset\oldPhi^-$, in particular $C'$ is disjoint from
$|\oldXi''|$, and the first assertion is proved. In Subcase (iii),
since $C\subset\partial G_\tS$, the argument used to prove the first
assertion shows that $C$ is disjoint from
$|\oldXi''|$, which is stronger than the second assertion. 
% follows from the same argument as the first. 
In Subcase (iv-a) we have $C=\partial E_\tS$; since $E_\tS$ is a component of $|\oldPhi^-|$ by \ref{stop beeping}, we have
$C\subset|\partial\oldPhi^-|$. Since
%According to \redmissingref{cross-ref---Is this part of the stuff that goes in \ref{oldXi}?} we have 
$\oldXi''\subset\inter\oldPhi^-$, in particular $C$ is  disjoint from $|\oldXi''|$, and the second assertion is true in this subcase as well.

To prove the second assertion of \ref{zeroth claim} in Subcase (iv-b),
recall that in this subcase, $C$ is a boundary component of
$B:=|\epsilon_\tS|(|\tau|(G_{\tS'}))\subset \inter E_\tS$. 
%Let $B$ denote the component of $|\epsilon_\tS|(|\tau|(G_{\tS'}))$ containing $C$. 
It suffices to show that $\oldXi''$ is (orbifold)-isotopic in $\obd(\tcals_0)$ to a suborbifold of $\obd(\tcals_0-\inter 
B)$. In the notation of Section \ref{vegematic section}, this means showing that $X''\preceq
%%be a regular neighborhood of $C$ in $\inter E_\tS$ which is disjoint
%from $\fraks_\oldPhi$. \redmissingref{I do know $C$ is disjoint from
%  that, right? I wish I knew what I was talking about.} 
%We must show
%that $X''=[\oldXi'']\preceq
[\tcals_0-\inter B]$ in $\barcaly(\obd(\tcals_0))$. We have
$X''\in\barcaly_-(\obd(\tcals_0))$. Hence, according to
Corollary \ref{mystery corollary}, it suffices to show that for every
component $Y$ of $X''$ we have $Y\preceq[\obd(\tcals_0-\inter
B)]$. 
%\tcals_0-\inter B 

Let  $Y$ be any component  of $X''$. Then $Y$ is in particular a component of $X$, and according to an observation in \ref{oldXi} we have $X\preceq\oldPhi^-$; thus
%According to Lemma \ref{here goethe},
%$Y$ is a component of $X_n$ for some $n\ge1$. Since $X_n\preceq V_n$ by \ref{oldXi}, we have 
%$Y\preceq V_n=\dom\noodge_n$. By (\ref{i thought so too}) and (\ref{and so did you}), we have
%$V_{n}\preceq V_1=[\oldPhi^-$], so that 
$Y\preceq{\oldPhi^-}$.
%By \redmissingref{same cross-ref. as above} we have %$X''\preceq[\oldPhi^-]$. 
Hence (see \ref{components})
%for any component $Y$ of $X''$ 
there is a component $W$ of $[\oldPhi^-]$ such that $Y\preceq W$. If
$W$ is distinct from the component $[\obd(E_\tS)]$ of $[\oldPhi^-]$,
then
$W\subset\obd(\tcals_0-E_\tS)\subset\obd(\tcals_0-\inter B)$, and hence
 $Y\preceq W\preceq
%[\obd(\tcals_0-E_\tS)]\preceq
[\obd(\tcals_0-\inter B)]$, as required. We may therefore assume that $W=[\obd(E_\tS)]$, so that $Y\preceq[\obd(E_\tS)]$.

Since $Y$ is  a component of $X$, it follows from Lemma \ref{here goethe} that
$Y$ is a component of $X_n$ for some $n\ge1$. It then
%By \ref{oldXi}, we have  Since $X_n\preceq V_n$, 
%$Y\preceq V_n=\dom\noodge_n$. It also 
follows from \ref{oldXi} that 
%the
%According to Lemma \ref{here goethe},
%\redmissingref{I think this is something I want to %mention
%  in \ref{oldXi}}, 
%$Y$ is a component of $X_n$ for some $n\ge1$. Since %$X_n\preceq V_n$ by \ref{oldXi}, we have 
%The
%definition of %$X_n=X_n(\oldOmega,\obd(\cals),\obd(\cals-S))$
%implies \redmissingref{I think} that 
%$Y\preceq V_n=\dom\noodge_n$.
%It follows from \ref{oldXi} that the 
%component $Y$ of $X_n$ 
 $Y\in\xi_n$, so that $Y$ is a component of $V_n=\dom\noodge_n$, and $\noodge_n(Y)\preceq[\obd(\tcals-\tcals_0)]$. On the other hand, since $j_1=]\iota_\oldPsi]$  by (\ref{and so did you}),  and 
since
$E_\tS$ is contained in the domain of $\iota_\oldPsi$, and is $\iota_\oldPsi$-invariant, by \ref{stop beeping},
we may invoke
the order-preserving property of \ref{restriction} to deduce that $\noodge_1(Y)=[\iota_\oldPsi](Y)
\preceq
%\noodge_1([\obd(E_\tS)])=
[\iota_\oldPsi(\obd(E_\tS))]=[\obd(E_\tS)]
%$, and it follows that $\noodge_1(Y)\preceq\bd(E_\tS)
% \redmissingref{other
%cross-refs. presumably referring to the twisted $I$-fibration} we have $\noodge_1(Y)\preceq\noodge_1([\obd(E_\tS)])=[\iota(\obd(E_\tS))]=[\epsilon_\tS( \obd (E_\tS))]=[\obd(E_\tS)]
\preceq[\tcals_0]$. Since $\noodge_n(Y)\preceq[\obd(\tcals-\tcals_0)]$ and 
$\noodge_1(Y)\preceq[\tcals_0]$,
 we must have $n>1$. 

Using the definition of $\noodge_n$ (see \ref{in duck tape}) and \ref{more associativity}, we may then write $\noodge_n=([\iota]\diamond[\tau])^{\diamond
(n-1)}\diamond[\iota]=(\noodge_{n-1}\diamond[\tau])\diamond[\iota]$.
Hence we have $Y\preceq\dom(
(\noodge_{n-1}\diamond[\tau])\diamond[\iota])$, and
%\redcomment{Use this: the relation
%$[\frakZ]\preceq V_{n}$ may be rewritten as $[\frakZ]\preceq\dom
%([\iota]\diamond[\tau]\diamond\noodge_{n-1})$. 
Lemma \ref{before
  associativity} then gives that
%$Y\preceq\dom([\iota]=\oldPhi^-$ \redmissingref{I think that last equality is OK}, and that 
$[\iota](Y)\preceq\dom (\noodge_{n-1}\diamond[\tau])$. 
A 
second application of
Lemma \ref{before
  associativity} then gives that  $[\tau]([\iota](Y))\preceq\dom\noodge_{n-1}=V_{n-1}$. Since $n>1$, it follows from (\ref{i thought so too}) and (\ref{and so did you}) that $V_{n-1}\preceq V_1=[\oldPhi^-]$. Hence $[\tau]([\iota](Y))\preceq[\oldPhi^-]$. 

Since $Y\preceq[\obd(E_\tS)]$, we may write $Y=[\oldUpsilon]$ for some taut negative suborbifold $\oldUpsilon$ of $\obd(E_\tS)$. Since $\epsilon_\tS$ is by definition the restriction of $\iota_\oldPsi$ to $\obd(E_\tS)$, the relation
$[\tau]([\iota](Y)])\preceq[\oldPhi^-]$ implies that
%, and we have $[\iota(\oldUpsilon)]\preceq\dom (\noodge_{n-1}\diamond[\tau])$. 
%This means that
 $\tau(\epsilon_\tS(\oldUpsilon))$
is (orbifold-)isotopic in $\obd(\tS')$ to a suborbifold $\frakJ$ of $\oldPhi^-$. The component of $\oldPhi^-$ containing $\frakJ$ is by definition a negative component of $\oldPhi$, and is therefore disjoint from the annular component $\obd(G_{\tS'})$ of $\oldPhi$. Thus we have
 $\frakJ\subset\obd(\tS'-\inter G_{\tS'})$, and $\epsilon_\tS(\oldUpsilon)$ is isotopic in $\obd(\tS)$ to the suborbifold $\tau(\frakJ)$ of $\obd(\tS)-\inter\tau(\obd(G_{\tS'}))$. 

Now since the component $\obd(G_{\tS'})$ of $\oldPhi$ is taut  (see
\ref{tuesa day}), $\obd(\tS)-\inter\tau(\obd(G_{\tS'}))$ is also taut. In particular, each component of $\obd(\tS)-\inter\tau(\obd(G_{\tS'}))$ has non-positive Euler characteristic. According to Corollary \ref{mangenfeffer}, applied with the annular orbifold $\tau(\obd (G_{\tS'}))$ playing the role of $\oldGamma$, no component of
$\obd(\tS)-\inter\tau(\obd(G_{\tS'})$ can be annular; furthermore, no component can be toric, since $\obd(\tcals)$ is negative (see \ref{doublemint}). Hence $\obd(\tS)-\inter\tau(\obd(G_{\tS'}))$ is negative as well as being taut, and is therefore an element of $\Theta_-(\obd(\tcals))$. The fact that $\epsilon_\tS(\oldUpsilon)$ is isotopic in $\obd(\tS)$ to a suborbifold of $\obd(\tS)-\inter\tau(\obd(G_{\tS'}))$ can now be expressed in $\barcaly_-(\obd(\tcals))$ by the relation $[
\epsilon_\tS(\oldUpsilon)]\preceq[\obd(\tS)-\inter\tau(\obd(G_{\tS'}))]$. 
But since $\epsilon_\tS(\oldUpsilon)\subset\epsilon_\tS(\obd(E_\tS))=\obd(E_\tS)$, we also have $[\epsilon_\tS(\oldUpsilon)]\preceq[\obd(E_\tS)]$. Hence
$[
\epsilon_\tS(\oldUpsilon)]\preceq[\obd(\tS)-\inter\tau(\obd(G_{\tS'}))]\wedge[\obd(E_\tS)]$.

Since we are in Case (iv-b), we have $\tau(\obd(G_{\tS'}))\subset \inter\obd(E_\tS)$. In particular, $\partial \tau(\obd(G_{\tS'}))=\partial(\obd(\tS)-\inter\tau(\obd(G_{\tS'}))$ is disjoint from $\partial\obd(E_\tS)$. It therefore follows from Corollary \ref{nafta} that
$[\obd(\tS-\inter\tau(G_{\tS'})]\wedge[\obd(E_\tS)]=[\oldGamma_-]$, where $\oldGamma_-$ denotes the union of all negative components of $(\obd(\tS)-\inter\tau(\obd(G_{\tS'})))\cap\obd(E_\tS)=
\obd(E_\tS)-\inter\tau(\obd(G_{\tS'}))$. Thus we have
$[(
\epsilon_\tS(\oldUpsilon)]\preceq[\oldGamma_-]$; in particular, $
\epsilon_\tS(\oldUpsilon)$ is isotopic in $\obd(\tS)$ to a suborbifold of 
$\obd(E_\tS)-\inter\tau(\obd(G_{\tS'}))$. Since
$\epsilon_\tS(\oldUpsilon)$ is negative and is contained in 
$\obd(E_\tS)$, it now follows from Corollary \ref{i guess} that
$\epsilon_\tS(\oldUpsilon)$ is isotopic in $\obd(E_\tS)$ to a suborbifold of 
$\obd(E_\tS)-\inter\tau(\obd(G_{\tS'}))$. Since $\epsilon_\tS$ is an involution of $\obd(E_\tS)$,  it follows that $\oldUpsilon$ is isotopic in $\obd(E_\tS)$---and in particular in $\obd(\tS)$---to a suborbifold of $\obd(E_\tS)$ disjoint from
 $\epsilon_\tS(\tau (\obd( G_{\tS'}))=\obd(B)$. 
This shows that $Y\preceq[\obd(\tS-\inter B)]\preceq[\obd(\tcals_0-\inter B)]$, as required. Thus the proof of \ref{zeroth claim} is complete.

Now we claim:
\Claim\label{first and second claims}
If  $\chibar(X'')
=\chibar(\obd(\tcals_0))$, then each component of each of the $1$-manifolds 
%curves
$C$ and $C'$ is isotopic in
$\tcals_0-\fraks_{\obd(\tcals_0)}$ to a component of $|\partial\oldXi''|$. 
\EndClaim

%\redcomment{I've added the $\chibar$ hypothesis! It seems crucial, but
  %the point is that it holds in the app. below.}

To prove \ref{first and second claims}, note that by the second assertion of \ref{zeroth claim} $C$ is isotopic in $\tcals_0-\fraks_{\tcals_0}$ to a $1$-manifold $C_0$ which is disjoint from $|\oldXi''|$. If $c$ is any component of $C_0$, we then have 
$c\subset|\inter\frakB|$ for some component $\frakB$ of $\obd({\tcals_0})-\inter \oldXi''$. 
%\redcomment{Not exactly right,
  %since $C$ need not be connected. It should be OK if we replace $C$
  %in this statement by a component of $C$, but then I would need to go
  %from a statement about components to a statement about the whole
  %thing. The alternative would be to replace the statement of the
  %claim by a componentwise one, but then I would need to fiddle with
  %the app. below.} 
Since $\chibar(X'')
=\chibar(\obd(\tcals_0))$, it follows from  Lemma \ref{old partial
  order} that the orientable $2$-orbifold $\frakB$ is annular. It then follows from
%particular we have $c\subset\inter\frakB$, and by \redmissingref{cross-ref}, $c$ is $\pi_1$-injective in $\oldOmega$ and hence in $\frakB$. It therefore follows from 
\ref{cobound} that $c$ cobounds a weight-$0$ annulus in $\frakB$ with some component of $\partial|\frakB|\subset\partial\oldXi''$. This
% Hence $c$ is isotopic in $\tcals_0-\fraks_{\tcals_0}$ to a submanifold of some regular neighborhood of $\partial|\frakB|=\partial|\oldXi''|$ in $|\oldXi''|\setminus\fraks_{\tcals_0}$, which 
implies the assertion for $C$.
%component of $\partial|\frakB|\subset\partial|\oldXi''|$. 
The proof of the assertion for $C'$ is the same if one uses the first assertion of \ref{zeroth claim} in place of the second. Thus 
\ref{first and second claims} is proved.
%C

We now proceed to the proof that $S$ is a \clashsphere\ for $\cals$ in Case II. According to Lemma \ref{old partial order} we have $\chibar(X'')\le
\chibar(\obd(\tcals_0))=
2\chibar(\obd(S))$ If it happens that $\chibar(X'')
<2\chibar(\obd(S))$, then since we have $\chibar(X'')>0$ by (\ref{not
  nothing}),  it follows from 
Lemma \ref{tell me why} (and the hypothesis that $S$ is \doublefull) that $S$ is a \clashsphere. We may therefore restrict our attention to the subcase where 
$\chibar(X'')
=2\chibar(\obd(S))=\chibar(\obd(\tcals_0))$. In view of the definition of a \clashsphere, we are required to show that $\chibar(X''\wedge
[\tau](X''))< 
\chibar(X'')$.
%, or equivalently that 
%$\chibar(X''\wedge
%[\tau](X''))<
%\chibar(\obd(\tcals_0))$. 
%If we set $Y=\tS\wedge[\tau](X)=[\tau](\tS')$, then $X''\wedge
%[\tau](X)=X''\wedge Y$, so that it suffices to prove that $\chibar(X''\wedge
%Y)<
%\chibar(\obd(S))$. 
%Since $X''\wedge
%[\tau](X'')\preceq [\tau](X'')$, we have
%$\chibar(X''\wedge
%[\tau](X''))\le\chibar([\tau](X''))$ by Proposition \ref{new partial order}; thus if  %$\chibar([\tau](X''))<\chibar(\obd(\tS\cup\tS'))$, the conclusion is true. We may %therefore assume that $\chibar([\tau](X''))=\chibar(\obd(\tS\cup\tS'))$. 

Since
 $\chibar(X'')
=\chibar(\obd(\tcals_0))$, it follows from
\ref{first and second claims} that each component of $C$ or
$\tau(C')$ is orbifold-isotopic in $\obd(\tS)$ to a component of $\partial\oldXi''$ or
$\partial\tau(\oldXi'')$ respectively. With
 \ref{can't split 'em}, this implies that
$\partial\oldXi''$ and $\tau(\partial\oldXi'')$ meet essentially in
$\obd(\tcals_0)$. Applying  Proposition \ref{tpp}, with
$\oldXi''$ and $\tau(\oldXi'')$ playing the roles of $\oldUpsilon_1$
and $\oldUpsilon_2$, we deduce that $\chibar(X''\wedge
[\tau](X''))<\min(\chibar(X''),\chibar([\tau](X'')))\le\chibar(X'')
%=\chibar(\obd(\tcals_0))
$. This
completes the proof of the lemma in Case II. 
\EndProof
%{tell\frakC _n %m $pP k_{n\tau\frakE\frakB\oldGamma\frakQ\frak$UYWG
%curve BVWXYZR\oldGamma(A) {dun R P C $Z^+ Y \frakC \noodge j X $Z Q D
%\Delta $Y $L \toldTheta

%\redmissingref{Make it clear that the \doublefull\ condition is needed for the application of Lemma \ref{tell me why}, which is quoted in the proof.}

\Lemma\label{no wonder 2}
Let $\Mh$ be a closed, orientable
hyperbolic $3$-orbifold, and set $\oldOmega=(\Mh)\pl$. Let $\cals$ be an
$\oldOmega$-\dandy\ system of spheres
in $M:=|\oldOmega|$. Let
$\cals_1$ be a subsystem of $\cals$.
Let $S$ be a \centralclashsphere\ of $\cals_1$ which is a
\doublecentral\ component of $\cals$. Then $S$
is a \centralclashsphere\ of $\cals$.
 \EndLemma

\Proof
Set $\tcals=\partial M\cut\cals$ and  $\tcals_1=\partial M\cut{\cals_1}$. Then $\tcals_1$ is canonically identified with a union of components of $\tcals$. Set $\tau=\tau_\cals:\tcals\to\tcals$. Then $\tau$ leaves $\tcals_1$ invariant, and under the canonical identification we have $\tau|\tcals_1=\tau_{\cals_1}$. 
If $S$ is a \centralclashsphere\ of $\cals_1$, one of the conditions (i)---(iv) of Definition \ref{self-clash def} holds with $\cals_1$ in place of $\cals$.
If (i) holds, let $\tS$ be a side of $S$ which contains an \bad\  component of $|\oldPhi(\oldOmega\cut{\obd(\tcals_1)})|$ (relative to $\cals_1$). According to Proposition \ref{dandy subsystem}, 
$\tS$  contains an \bad\  component of $|\oldPhi(\oldOmega\cut{\obd(\tcals)})|$ (relative to $\cals$). Hence $S$
is a \centralclashsphere\ of $\cals$, which is the required
conclusion. Now suppose that (ii) holds, so that the two sides $\tS$
and $\tS'$ of $S$ are pseudo-belted relative to $\cals_1$.  By Proposition
\ref{dandy subsystem}, $\tS$ and $\tS'$  are pseudo-belted relative to $\cals$,
and again the conclusion follows. 

Next suppose that (iii) holds, so
that the two sides $\tS$ and $\tS'$ of $S$ are belted relative to $\cals_1$,
and $\obd(G_\tS^{\cals_1})$ and $\tau(\obd(G_{\tS'}^{\cals_1}))$ are not isotopic. According to the hypothesis of the present
lemma, $\tS$ and $\tS'$ are \central\ relative to the system $\cals$. It therefore follows
from Proposition \ref{dandy subsystem} that $\tS$ and $\tS'$  are
belted relative to $\cals$, and that 
$\obd(G_\tS^\cals)$ and $\obd(G_{\tS'}^\cals)$ are respectively isotopic to $\obd(G_\tS^{\cals_1})$  and $\obd(G_{\tS'}^{\cals_1})$.
% are also isotopic.
%$\oldPhi(\oldPsi)$ and $\oldPhi(\oldPsi_1)$ 
%may be chosen within their isotopy classes so that $G_\tS^\cals=G_\tS^{\cals_1}$ and $G_{\tS'}^\cals=G_{\tS'}^{\cals_1}$.  
Hence $\obd(G_\tS^{\cals})$ and $\tau(\obd(G_{\tS'}^{\cals}))$ are not isotopic, and again the conclusion follows. 

Now suppose that (iv) holds. Thus the sides of $S$ may be labeled
as $\tS$ and $\tS'$, where $\tS$ is pseudo-belted side in $\cals_1$
and $\tS'$ is belted relative to $\cals_1$, and one of the sub-alternatives (a)
or (b) holds. 
Again applying Proposition  \ref{dandy subsystem} (and
again using that $\tS'$ is \central\ relative to $\cals$), we deduce that $\tS'$ is belted relative to $\cals$, 
and that   $\obd(G_{\tS'}^\cals)$ and $\obd(G_{\tS'}^{\cals_1})$ are isotopic; while
 $\tS$  is pseudo-belted relative to $\cals$,  %$\oldPhi^-(\oldPsi_1)$, and their $I$-fibrations, may be chosen within their isotopy classes so that
%,
and $\epsilon_\tS^\cals$ and $\epsilon_\tS^{\cals_1}$ are strongly
equivalent.
%(so that in particular
%$\obd(E_\tS^\cals)$ and $\obd(E_\tS^{\cals_1})$ are isotopic).
Hence we may assume 
$\obd(G_{\tS'}^{\cals_1})$ to have been chosen within its isotopy
class so that
$G_{\tS'}^\cals=G_{\tS'}^{\cals_1}$,  and we may assume
 $\epsilon_\tS^{\cals_1}$  to have been chosen within its strong
 equivalence 
class so that $E_\tS^\cals=E_\tS^{\cals_1}$ and
$\epsilon_\tS^\cals=\epsilon_\tS^{\cals_1}$.
If (a) holds with $\cals_1$ in place of $\cals$, i.e. if $\tau(\obd(G_{\tS'}^{\cals_1}))$ is not isotopic in $\obd(\tS)$ to a suborbifold of 
$\obd(E_\tS^{\cals_1})$, it follows that $\tau(\obd(G_{\tS'}^{\cals}))$ is not isotopic in $\obd(\tS)$ to a suborbifold of 
$\obd(E_\tS^{\cals})$; this is Sub-alternative (a) with respect to the system $\cals$, and it again follows that 
$S$
is a \centralclashsphere\ of $\cals$. Finally, if (b) holds with $\cals_1$ in place of $\cals$, i.e. if $\tau(\obd(G_{\tS'}^{\cals_1}))$ is orbifold-isotopic in $\obd(\tS)$ to a suborbifold $\frakC$ of $\obd(E_\tS^{\cals_1})$, and 
$\frakC$ is not orbifold-isotopic in $\obd(E_\tS^{\cals_1})$ to 
$\epsilon_\tS^{\cals_1}(\frakC)$, it follows that $\tau(\obd(G_{\tS'}^{\cals}))$ is orbifold-isotopic in $\obd(\tS)$ to a suborbifold $\frakC$ of $\obd(E_\tS^{\cals})$, and that
$\frakC$ is not orbifold-isotopic in $\obd(E_\tS^{\cals})$ to 
$\epsilon_\tS^{\cals}(\frakC)$. This is Sub-alternative (b) with respect to the system $\cals$, and again the conclusion follows.
\EndProof
%\tau\tcals

\section{The structure of dandy systems without clash components}\label{structure section}

\abstractcomment{\tiny
I'm removing Lemma
``in the beginning,'' which can be found in calL.tex, and
replacing it by Lemma \ref{new beginning} below. I'm not sure I need
that one either!}

\tinymissingref{\tiny I have removed Lemmas ``because strict'' and``also because strict.'' They can be found in newest.tex.}

\Definition
Let $\oldOmega$ be a closed,
orientable $3$-orbifold. An admissible system of
spheres  $\cals$  in $M:=|\oldOmega|$ will be said to be {\it balanced}
if all its components have the same weight. 
\EndDefinition

%\begin{remarknotation}\label{e-sub-m}
%Let $\oldOmega$ be a closed,
%orientable $3$-orbifold, and let $\cals$ be an admissible system of
%spheres in $M:=|\oldOmega|$. Set $\tcals=\partial M\cut\cals$. For each $m\in\NN$, we %will denote by
%$\sea_m(\cals)\subset\cals$ the union of all components $S$ of $\cals$
%having weight $m$. We
%will set 
%d
%$E_m(\cals)=\bigcup_{0\le d\le m}\sea_d(\cals)$. \redcomment{I'm guessing
  %I won't need $E_m(\cals)$.}
%\Lemma
%Let $\oldOmega$ be a closed,
%orientable $3$-orbifold, and let $\cals$ be an admissible system of
%spheres in $M:=|\oldOmega|$. 
%Note that $\sea_m(\cals)$ is an balanced subsystem of $\cals$ for every $m\in\NN$.
%We will also set
%$E_m(\cals)=\bigcup_{0\le d\le m}\sea_d(\cals)$,
%$\tsea_d(\cals)=\rho_\cals^{-1}(\sea_d(\cals))\subset\tcals$, \redmissingref{Not
  %sure whether this is used.} 
%and $\tE_m(\cals)=\rho_\cals^{-1}(E_m(\cals))=
%\bigcup_{0\le d\le m}\tsea_d(\cals)$
 %for every $m\in\NN$. 
%\end{remarknotation}

\begin{remarksnotationdefinitions}\label{old trier}
 Let $\Mh$ be a closed,
orientable hyperbolic $3$-orbifold, and set $\oldOmega=(\Mh)\pl$ and $M=|\oldOmega|$. 
Let $\cals$  be an admissible
% \dandy\
%comment{seems like enough for
  %now??? how natural is $\calL$ if the system isn't even
  %semi-\dandy?} 
system of spheres in $M$.
% having no \clashspheres. 
Set 
$\oldPsi=\oldOmega\cut
{\obd(\cals)}$ and $N=|\oldPsi|=M\cut\cals$.
%$N=M\cut\cals=|\oldPsi|$,
%$\tcals=\partial N$, $\tau=\tau_\cals$, and $\rho=\rho_\cals$.  
We will define $\calL\oldomecals\subset N$ to be the union of all sets of
the form  $L_\tS$, where $\tS$ ranges over all belted components of
$\tcals$ (see \ref{stop beeping}). Since each set $L_\tS$ is by
definition a component of
$|\oldSigma(\oldPsi)|$, any two distinct sets in the union defining  $\calL\oldomecals$
are disjoint, and the components of $\calL\oldomecals$ are
precisely the sets of the form $L_\tS$, where $\tS\in\calc(\tcals)$ is
belted.
Since $\Fr_N\calL\omecals$ is a union of components of
$|\frakA(\oldPsi)|$, it follows from the discussion in \ref{tuesa day}
that every component of  $\obd(\Fr_N\calL\oldomecals)$ is an essential annular $2$-suborbifold of $\oldPsi$. Furthermore, according to
 \ref{stop
    beeping},  if  $L$ is any component of
$\calL\oldomecals$, then $\obd(L)$ is a \torifold, and
 $\obd(L\cap\tcals)$ is a non-empty disjoint union of annular
 $2$-orbifolds, contained in
 $\obd(\partial L)$ and $\pi_1$-injective in $\obd(L)$. 
Thus for each component $Q$ of $L\cap\tcals$,  the inclusion homomorphism
 $\pi_1(\obd(Q))\to\pi_1(\obd(L))$ is an injection with infinite domain; furthermore, $\pi_1(\obd(L))$ is virtually cyclic since $\obd(L)$ is a \torifold. Hence the 
 index $n_Q$ of the image of this inclusion homomorphism
 $\pi_1(\obd(Q))\to\pi_1(\obd(L))$ is finite. We shall set
 $w_L=\min_{Q\in\calc(L\cap\tcals)}n_Q$.
\end{remarksnotationdefinitions}

\Lemma\label{new beginning}
 Let $\Mh$ be a closed,
orientable hyperbolic $3$-orbifold, and set $\oldOmega=(\Mh)\pl$. Set $M=|\oldOmega|$ and
$\tcals=\partial M\cut\cals$.
Let $\cals$  be a  \dandy, balanced system of spheres in $M$.
%having no \clashspheres.
%, and let
%$m\in\NN$ be given. 
If $\tS$
is a component of $\tcals$ such that $\calL\oldomecals\cap\tS\ne\emptyset$, then
$\tS$ is belted and $\calL\oldomecals\cap\tS=G_\tS$.
\EndLemma

\Proof
Set $\oldPsi=\oldOmega\cut{ \obd(\cals)}$ and $N=|\oldPsi|=M\cut\cals$, so that $\tcals=\partial N$.

Let $L$ be any component of $\calL:=\calL\oldomecals$ such that $L\cap\tS\ne\emptyset$, and let $G$ be any component of $L\cap\tS$. It suffices to show that $\tS$ is belted and that $G=G_\tS$. 

According to the definition of $\calL$, we have $L=L_{\tT}$ for some belted component $\tT$ of
$\tcals$. 
It  was observed in \ref{old trier} that
$\obd(L)=\obd(L_{\tT})$ is a \torifold, and that every component
of $\obd(L\cap\tcals)=\obd(L_{\tT}\cap\partial|\oldPsi|)$ is an
annular $2$-orbifold.  In particular $\obd(G)$ is an annular $2$-orbifold.

I claim:

\Claim\label{same size}All the boundary components of the $2$-manifold $L\cap\tcals$ have the same size (see
\ref{size def}) in $\tcals$.
\EndClaim

To prove \ref{same size}, note that since
$\obd(L)=\obd(L_{\tT})$ is a \torifold, the orientable  orbifold $\obd(\partial L)$ is
toric, and hence either (i) $\partial L$ is a weight-$0$
torus, or (ii) $\partial L$ is a weight-$4$ sphere.

First consider the case in which (i) holds. 
%Since $\obd(\partial L)$ is a toric orbifold disjoint from $\fraks_\oldPsi$, the $2$-manifold $\partial L$ is a torus. Likewise, t
The components of
$L\cap\tcals$ are disjoint from $\fraks_\oldPsi$ and are underlying
surfaces of orientable annular orbifolds; these components are therefore annuli
in the torus $\partial L$. 
%Here is one from Lemma ``new beginning'': ``But by (missingref:cross-reference to characteristic suborbifold section), $\obd(L\cap\tcals)$, as a union of components of $\oldPhi(\oldPsi)$, is $\pi_1$-injective in $\oldPsi$, and is therefore also $\pi_1$-injective in $\partial L\subset\oldPsi$.''
%
But since the components of $\frakA(\oldPsi)$ are essential in
$\oldPsi$ by \ref{tuesa day},
% by \redmissingref{cross-reference to characteristic suborbifold
% section},
 $\obd(L\cap\tcals)$, as a union of components of $\oldPhi(\oldPsi)$, is $\pi_1$-injective in $\oldPsi$, and is therefore also $\pi_1$-injective in $\partial L\subset\oldPsi$. Hence the components of
$L\cap\tcals$  are homotopically non-trivial annuli in $\partial
L$. It follows that  if we set $m=\compnum(L\cap\tcals)$, we may index the components of $L\cap\tcals$ as %$B_0,\lcdots,B_{n_1}$, where $0,\ldots,n-1$ are interpreted as elements of $\ZZ/m\ZZ $,
$(G_i)_{i\in\ZZ/m\ZZ }$, label the components of $G_i$ as $C_i$ and $C_i'$, and index the components of $\Fr L=\partial L-\inter(L\cap\tcals)$ as
$(A_i)_{i\in\ZZ/m\ZZ }$, in such a way that $C_{i}'$ and $C_{i+1}$ are the components of $\partial A_i$ for each $i\in\ZZ/m\ZZ $. 

For each $i\in\ZZ/m\ZZ $, the annulus $A_i$ is properly embedded in $N$, and
has weight $0$ since $\partial L\supset A_i$ has weight $0$.
According to the first condition in the definition of a \dandy\ system, it then follows that the two components of
  $\partial A_i$ have the same size; that is, we have $\size C_i'=\size C_{i+1}$ for each $i\in\ZZ/m\ZZ $. But since the annulus $G_i\subset\partial L$ has weight $0$, the simple closed curves $C_i$ and $C_i'$ are isotopic in $\tcals-\fraks_\oldPsi$, are therefore have the same size. Hence $\size C_i=\size C_{i+1}$ for each $i\in\ZZ/m\ZZ $. It follows that all the $C_i$ have the same size, and \ref{same size} is proved in this case.

Now suppose that (ii) holds.
Each component of
$\obd(L\cap\tcals)$ is an orientable annular orbifold, and hence each component
of $L\cap\tcals$ is either a weight-$0$ annulus  or a weight-$2$
disk. Since $\wt L=4$, exactly two
components of $L\cap\tcals$, say $G^+$ and $G^-$, are disks, and
$\tcals\cap\fraks_\oldPsi\subset G^+\cup G^-$. 
%\redcomment{That sentence
  %is wrong. We don't know that $\card(\tcals\cap\fraks_\oldPsi)=4$. We
  %know that $\card((\partial L)\cap\fraks_\oldPsi)=4$. The rest of the
  %paragraph needs to be redone. I think the conclusion ``we may index the components of %$L\cap\tcals$ as $(G_i)_{1\le i\le{m-2}}$, label the components of $G_i$ as $C_i^-$ and $C_i^+$, set $C_0^+=\partial G^-$ and $C_{m-1}^-=\partial G^+$, and index the components of $\Fr L=(\partial L)-\inter(L\cap\tcals)$ as
%$(A_i)_{1\le i\le{m-1}}$, in such a way that 
%$C_{i-1}^+$ and $C_{i}^-$ are the components of $\partial A_i$, for
%$1\le i\le m-1$,'' or something like it, may be true unless $\partial
%L$ consists of one disk in $\partial\oldPsi$ and one frontier disk, in
%which case \ref{same size} is trivial.} 
The remaining components of $L\cap\tcals$ are annuli. By \ref{tuesa day}, $\obd(L\cap\tcals)$, as a union of components of $\oldPhi(\oldPsi)$, is $\pi_1$-injective in $\oldPsi$. Hence no annulus component of $L\cap\tcals$ can have a boundary component that bounds a disk in $\partial L\setminus\fraks_\oldPsi$. This means that every 
annulus component of $L\cap\tcals$ separates $G^+$ from $G^-$ in
$\tS$.  
Hence if $m=\compnum(L\cap\tcals)$, then $m\ge2$, and we may index the annulus components of $L\cap\tcals$ as $(G_i)_{1\le i\le{m-2}}$, label the components of $\partial G_i$ as $C_i^-$ and $C_i^+$, set $C_0^+=\partial G^-$ and $C_{m-1}^-=\partial G^+$, and index the components of $\Fr L=\partial L-\inter(L\cap\tcals)$ as
$(A_i)_{1\le i\le{m-1}}$, in such a way that 
$C_{i-1}^+$ and $C_{i}^-$ are the components of $\partial A_i$, for $1\le i\le m-1$.

For each $i$ with $1\le i\le m-1$, the annulus $A_i$ is properly
embedded in $N$, and has weight $0$ since $A_i$ is disjoint from
$\inter G^+$ and $\inter G^-$. According to the first condition in the definition of a \dandy\ system, it then follows that the two components of
  $\partial A_i$ have the same size; that is, 
$\size C_{i-1}^+=\size C_{i}^-$ for $1\le i\le m-1$.
But since $G_i\subset\partial L$ is disjoint from $\fraks_\oldPsi$, the simple closed curves $C_i^-$ and $C_i^+$ are isotopic in $\tcals-\fraks_\oldPsi$, so that $\size C_i^-=\size C_{i}^+$ for 
$1\le i\le{m-2}$. It follows that all the $C_i^-$ and $C_i^+$ have the same size, and the proof of \ref{same size} is complete.

By hypothesis $\cals$ is balanced, which by definition means that its
components all have the same weight in $\oldOmega$, or equivalently
that all components of $\tcals$ have the same weight in $\oldPsi$;
we shall denote this common weight by $n$. Since $\tT$ is belted and
$L=L_\tT$, the surface $G_\tT$ is a component of $L\cap\tcals$. Since
$G_\tT$ is a belt for $\tT$, Lemma \ref{how do i recognize
  a belt?}   implies that each boundary component of
$G_\tT$ has size $(\wt\tT)/2=n/2$. It therefore follows from
\ref{same size} that every boundary component of $L\cap\tcals$ has
size $n/2$. In particular each component of $\partial G$ has size
$n/2$. But since $\tS$ is a component of $\tcals$ we have
$\wt_\oldPsi\tS=n$, so that each boundary component of $L\cap\tcals$
has size $(\wt_\oldPsi\tS )/2$. Furthermore, the component $L$ of $\calL\oldomecals$ is by definition a component of $|\oldSigma(\oldPsi)$, and hence $G$ is a component of $|\oldPhi(\oldPsi)$. As we have observed that $\obd(G)$ is an
annular $2$-orbifold, it now follows from Lemma \ref{how do i recognize
  a belt?} that $G$ is a belt for $\tS$. In view of the uniqueness of belts (see Lemma \ref{it's-a unique}) and the definition of $G_\tS$, this shows that $\tS$ is belted and that $G=G_\tS$, as required.
\EndProof
%4

\begin{remarksnotationdefinitions}\label{new trier} 
%%\redcomment{Adapt the following to get a def. of $\calf_d\oldomecals$, which is the for%est defined by belts of components of $\tsea_d(\cals)$, and by the corresponding things %for non-belted, \full\ components of $\tsea_d(\cals)$.}

\tinymissingref{\tiny I'm starting a new version here with a revised def. of
  $\calf_m$. For the old def. see eclipse3.tex. }

Let $\Mh$ be a closed hyperbolic $3$-orbifold $\oldOmega$, and set $\oldOmega=(\Mh)\pl$. Given an admissible %$\oldOmega$-\dandy,
%balanced
 system of spheres $\cals\subset M:=|\oldOmega|$,
% without
%\centralclashspheres \redmissingref{do these conditions matter for the def.? and does it matter whether they matter? One point is that if they are removed here, the sentence ``It will follow from Proposition \ref{new enchanted forest} that
%  there are no loops in $\calf$'' will have to be revised so as to include these conditions}, 
we will describe a finite graph
$\calf=\calf\oldomecals$ canonically associated with $\oldOmega$ and $\cals$. For
the purpose of the description, set $\tcals=\partial M\cut\cals$ and $\calL=\calL\oldomecals$.

Let $\tcalsu$ denote the union of all unbelted components of
$\tcals$. The vertex set $\calv=\calv\oldomecals$ is defined to be a
bijective copy of the disjoint union $\calc(\calL)\discup\calc(\tcalsu)$; we
fix a bijection $X\to v_X$ from $\calc(\calL)\discup\calc(\tcalsu)$ to
$\calv$. Thus every element of $\calv$ can be written uniquely in
exactly one of the forms $v_L$, where $L$ is a component of $\calL$,
or $v_\tS$, where $\tS$ is an unbelted component of $\tcals$. Vertices
of the former type will be termed {\it material}, while vertices
of the latter type will be termed {\it ideal}. 

The edge set $\cale=\cale\oldomecals$ is defined to be a
bijective copy of $\calc(\cals)$; we
fix a bijection $S\to e_S$ from $\calc(\cals)$ to
$\cale$. Thus every element of $\cale$ can be written uniquely in the
form $e_S$, where $S$ is a component of $\cals$.

If $S$ is a component of $\cals$, with each side $\tS$ of $S$ we shall
associate an element $v$ of $\calv$, by defining $v$ to be the
material vertex $v_{L_\tS}$ if $\tS$ is belted, and defining it to be
the ideal vertex $v_\tS$ if $\tS$ is unbelted. The vertices incident
to the edge $e_S$ will be defined to be the two elements of $\calv$
associated with the sides of $S$. Thus an arbitrary edge $e=e_S$ of
$\calf$ is incident to exactly two
distinct vertices of $\calf$ unless it happens that the elements of $\calv$
associated with the two sides of $S$ coincide, in which case $e$ is a
loop. (It will follow from Proposition \ref{new enchanted forest} that
in the case where $\cals$ is $\oldOmega$-\dandy\ and
balanced, and has no
\centralclashspheres,  there are no loops in $\calf$.)

It follows from the definitions given above that a material vertex
$v=v_L$ of $\calf$ is incident to at least one edge, because the
component $L$ of $\calL$ can by definition be written in the form
$L=L_\tS$ for some component $\tS$ of $\tcals$ (see \ref{old trier}),
and $e_{\rho(\tS)}$ is then by definition incident to $v$. It also
follows from the definitions given above that an ideal vertex $v_\tS$
is incident to exactly one edge of $\calf$, namely the edge
$e_{\rho(\tS)}$. Thus $\calf$ has no isolated vertices, and every
ideal vertex of $\calf$ has valence $1$.

We define a {\it special material vertex} of $\calf$ to be a material
  vertex of the form $v_L$, where $L$ is a component
of $\calL$ such that $w_L>1$ (where $w_L$ is defined as in \ref{old
  trier}).
% or $L\cap\fraks_\oldPsi\ne\emptyset$. \redcomment{That last 
  %condition is not the right one. The model example of what I have in
  %mind is a solid torus with a singular core curve. Actually with the
  %new def. of $w_L$, the additional condition may not be needed---but
%  I'm a little uncertain about the right def.}  
We define a {\it special ideal vertex} of $\calf$ to be an ideal
  vertex of the form $v_\tS$, where $\tS$ is a component of $\tcals$ which is
pseudo-belted relative to $\cals$. (According to Lemma \ref{it's-a unique}, a
pseudo-belted component 
$\tS$  of $\tcals$ is not belted, and $v_\tS$ is therefore indeed
a well-defined ideal vertex.) By a {\it special vertex} of $\calf$ we
shall mean simply a vertex which is either a special material vertex
or a special ideal vertex.
\end{remarksnotationdefinitions}
%isolated

\Proposition\label{new enchanted forest} Let $\Mh$ be a closed,
orientable hyperbolic $3$-orbifold, and set $\oldOmega=(\Mh)\pl$. Set $M:=|\oldOmega|$. Let $\cals$ be a \dandy, balanced system
of spheres in $M$ having no \centralclashspheres. Then the graph
$\calf\oldomecals$ is a (possibly empty) forest, that is, each of its
components is a tree. Furthermore, each component of $\calf\oldomecals$
contains at most one special vertex. 
\EndProposition

\Proof
Set $\oldPsi=\oldOmega\cut{\obd(\cals)}$, $N=|\oldPsi|=M\cut\cals$,
$\tcals=\partial N$, 
$\calL=\calL\oldomecals$,
$\rho=\rho_{\obd(\cals)}$, and $\tau=\tau_{\obd(\cals)}$. Then we have $|\rho|=\rho_\cals$ and $|\tau|=\tau_\cals$.

For $m=0$ and for $m=2$, we let $\cals_m$ denote the union of all components of $\cals$ that have exactly $m$ belted sides. We denote by $\cals_2^{!}$ the union of all components of $\cals$ that have
a pseudo-belted side and a belted side, and we denote by $\cals_1$ the union of all components of $\cals$ that have one belted side and one side that is neither belted nor pseudo-belted. According to Lemma \ref{it's-a unique}, no
component of $\tcals$ can be both belted and pseudo-belted; it follows that $\cals$ is the disjoint union of $\cals_0$, $\cals_1$, $\cals_2$ and $\cals_2^!$, and that each component of $\cals_2^{!}$ has a
unique belted side and a unique pseudo-belted side.

If $S$ is any   component of $\cals_2^{!}$, the hypothesis implies that $S$ is not a \centralclashsphere\ of $\cals$; in particular, it does not satisfy Alternative (iv) of Definition \ref{self-clash def}. Hence if $\tS$ and $\tS'$ denote
 respectively the pseudo-belted side and the belted side of $S$, then
% $G_{\tS'}$ 
%may be chosen within its (orbifold-)isotopy class so that
 $\obd(|\tau|(G_{\tS'}))
=\tau(\obd(G_{\tS'}))
$
%\subset 
 is orbifold-isotopic in $\obd(\tS)$ to a suborbifold of $
\obd(E_\tS)$, which will here be denoted by $\frakC_\tS$; and
$\frakC_\tS$
 is orbifold-isotopic in $\obd(E_\tS)$ to 
$\epsilon_\tS(\frakC_\tS)$. 
It then follows from Prop. \ref{steinbeck}
that, after modifying $\frakC_\tS$ within its isotopy class, we may assume that $\frakC_\tS$ is invariant under  the involution $\epsilon _\tS$, which by the definition of $\epsilon_\tS$ (see \ref{stop beeping})means that it is invariant under $\iota_\oldPsi$. By an observation made in \ref{what's iota?},
%In view of the definition of $\iota_\oldPsi$ (see \ref{what's iota?}), 
we may therefore
fix an $I$-fibration of $\frakR_\tS$, which is compatible with $\frakR_\tS$ and therefore satisfies $\partialh \frakR_\tS=\frakR_\tS\cap\partial
\oldPsi=\frakR_\tS\cap\obd(\tcals)$, and a suborbifold $\frakD_\tS$ of
  $\frakR_\tS$ which is saturated and therefore inherits an $I$-fibration
  over a $2$-orbifold, 
% the structure of an
 %$I$-bundle 
  %over a $2$-orbifold, 
such that 
%\redproofreadingnote{The fact that there is an $I$-fibration of $\frakU$ satisfying the first (three-way) equality says that $(\frakU,\frakR_\tS\cap\partial\oldPsi)$ is an \spair, which is true because $\frakU$ is an \Ssuborbifold. The equality $\frakC_\tS 
%=\partialh \frakD_\tS$ stated below is the consequence of the $\epsilon$-invariance. It should be easier to write this passage correctly now that all the issues about $\epsilon$ are cleared up} and  
$\frakC_\tS 
=\partialh \frakD_\tS$.
Thus 
 $|\frakC_\tS|\subset
E_\tS\subset\tS$. Now $\obd(G_{\tS'})$ is an annular orbifold, since $G_{\tS'}$ is a belt for
  $\tS'$ (see \ref{what's a belt?}). Hence
%, the
  %orbifold. 
$\frakC_\tS$ is an annular $2$-orbifold. 
%with; and $\frakC_\tS$
% is disjoint from $\fraks_\oldPsi$
  %(since $G_{\tS'}$ is), and
 %is (orbifold)-isotopic to $\obd(\tau (G_{\tS'}))$ in
%  $\obd(\tS)$.
%\setminus\fraks_\oldPsi$; 
We have $|\frakD_\tS|\cap\tcals=|\frakD_\tS|\cap(|\frakR_\tS|\cap\tcals)=|\frakD_\tS\cap\partialh\frakR_\tS|=|\partialh\frakD_\tS|$, since $\frakD_\tS$ inherits its fibration from $\frakR_\tS$; thus $|\frakD_\tS|\cap\tcals=|\frakC_\tS|$.

Since $\frakD_\tS$ is equipped with an $I$-fibration and its horizontal boundary is an annular
orbifold,  $\frakD_\tS$ is a
\torifold. It follows from \ref{fibered stuff} that the inclusion homomorphism $\pi_1(\frakC_\tS)\to\pi_1(\frakD_\tS)$ has index
$2$ in $\pi_1(\frakD_\tS)$.

To summarize: 
\Claim\label{replaces readin' john o'hara}
For every  component $S$  of $\cals_2^{!}$, if $\tS$ and $\tS'$ denote
 respectively the pseudo-belted side and the belted side of $S$, the annular $2$-orbifolds $\frakC_\tS$ and $\obd(|\tau| (G_{\tS'}))
=\tau (\obd (G_{\tS'}))
$ are
isotopic in
  $\obd(\tS)$. We have $|\frakC_\tS|\subset
  %contained in
 %$ 
E_\tS\subset\tS$ and
$|\frakD_\tS|\cap\tcals=|\frakC_\tS|$. Furthermore, $\frakD_\tS$ is a
\torifold, and the inclusion homomorphism $\pi_1(\frakC_\tS)\to\pi_1(\frakD_\tS)$ has index
$2$ in $\pi_1(\frakD_\tS)$.
\EndClaim

If $S$ is any   component of $\cals_2$, the hypothesis implies that $S$ is not a \centralclashsphere\ of $\cals$; in particular, it does not satisfy Alternative (iii) of Definition \ref{self-clash def}. Hence:

\Claim\label{replaces kimbecile}
If $S$ is any   component of $\cals_2$,  and if $\tS$ and $\tS'$
denote the sides of $S$, then, $\obd(G_{\tS})$ and
$\obd(|\tau|(G_{\tS'}))
=\tau (\obd (G_{\tS'}))
$ are isotopic in $\obd(\tS)$. 
\EndClaim 

The manifold $\calL$ is a union of components of $|\oldSigma(\oldPsi)|$, where $\oldSigma(\oldPsi)$ is well-defined up to (orbifold) isotopy.
% rel $\fraks_\oldOmega$. 
For each belted component $\tS$ of $\tcals$, Lemma \ref{new beginning} implies that $\calL\cap\tS=G_\tS$. It follows from \ref{replaces readin' john o'hara} and \ref{replaces kimbecile} that, after possibly modifying  $\oldSigma(\oldPsi)$ within its (orbifold) isotopy class, we may assume that  the following condition holds:
\Claim\label{fixed it up}
For every  component $S$  of $\cals_2^{!}$, if $\tS$ and $\tS'$ denote
respectively the pseudo-belted side and the belted
 side of $S$, we have 
$|\frakC_\tS|=|\tau| (G_{\tS'})$;
and
for every  component $S$  of $\cals_2$, if $\tS$ and $\tS'$ denote the
 sides of $S$, we have $G_{\tS}=|\tau|(G_{\tS'})$.
\EndClaim

If $S$ and $T$ are distinct components of $\cals_2^{!}$, and $\tS$ and
$\tT$ denote their pseudo-belted sides, then
 according to the definitions
in \ref{stop beeping}, each of the sets $|\frakR_\tS|$ and
$|\frakR_\tT|$ is a component of $|\oldSigma(\oldPsi)|$, and we have
$|\frakR_\tS|\cap\cals=E_\tS\subset\tS$, while
$|\frakR_\tT|\cap\cals=E_\tT\subset\tT$. Since $S$ and $T$ are distinct
components of $\cals$, we have $\tS\ne\tT$, so that 
$|\frakR_\tS|$ and
$|\frakR_\tT|$ are distinct components of $|\oldSigma(\oldPsi)|$. This proves:

\Claim\label{U ain't Uself}
If $S$ and $T$ are distinct components of $\cals_2^{!}$, and $\tS$ and
$\tT$ denote their pseudo-belted sides, then
$|\frakR_\tS|\cap|\frakR_\tT|=\emptyset$.
\EndClaim

Now let $\calr$ denote the union of all sets of
the form $|\frakR_\tS|$, where
$\tS$ ranges over the pseudo-belted sides of all components of
$\cals_2^{!}$; according to \ref{U ain't Uself}, this is a disjoint
union. Let 
$\cald\subset\calr$ denote the union of all sets of
the form $|\frakD_\tS|$, where
$\tS$ ranges over the pseudo-belted sides of all components of $\cals_2^{!}$.
According to \ref{replaces readin' john o'hara}, each component of $\obd(\cald)$ is a 
\torifold. We claim:

\Claim\label{U ain't L}
We have $\calr\cap\calL=\emptyset$. In particular, $\cald\cap\calL=\emptyset$.
\EndClaim

To prove \ref{U ain't L}, it suffices to show that if 
 $\tS$ is a belted component of $\cals$, and $\tT$ is the pseudo-belted side
 of a component of $\cals_2^{!}$,  then $L_\tS\cap|\frakR_\tT|=\emptyset$.  According to the definitions
in \ref{stop beeping}, each of the sets $L_\tS$ and
$|\frakR_\tT|$ is a component of $|\oldSigma(\oldPsi)|$. Furthermore, we have
$|\frakR_\tT|\cap\tcals=E_\tT\subset\tT$, so that in particular $|\frakR_\tT|\cap\tT\ne\emptyset$. On the other hand, Lemma \ref{new
  beginning} implies that $L_\tS\subset\calL$ can meet only belted
components of $\tcals$; since $\tT$ is pseudo-belted, Lemma
\ref{it's-a unique} implies that it cannot be belted. Hence
$L_\tS\cap\tT=\emptyset$. This shows that 
$L_\tS$ and
$|\frakR_\tT|$ are distinct components of $|\oldSigma(\oldPsi)|$, and the assertion
follows. Thus \ref{U ain't L} is proved.

Now set $\calp=\calL\cup\cald$; by \ref{U ain't L}, this is a
disjoint union. Hence each component of $\calp$ is a component of either $\calL$ or $\cald$. According to \ref{old trier}, $\obd(L)$ is a \torifold\ for each component $L$ of $\calL$; and, also by \ref{old trier}, we have $L=L_\tS$ for some component $\tS$ of $\tcals$, so that in particular $L\cap\tcals\supset G_\tS$. According to the definition of $\cald$, each component of $\cald$ has the form $|\frakD_\tS|$ for some component $\tS$ of $\tcals$, and by \ref
{replaces readin' john o'hara}, $\frakD_\tS$ is a \torifold\ and 
\Equation\label{just fits}
|\frakD_\tS|\cap\tcals= |\frakC_\tS|.
\EndEquation
In particular:

\Claim\label{they're torifolds}
For each component $P$ of $\calp$, the orbifold $\obd( P)$ is a \torifold, and $P\cap\cals\ne\emptyset$.
\EndClaim

Let us set $\calq=\calp\cap\tcals$. We claim:
\Claim\label{what is calg?}
For any component $\tS$ of $\tcals$, we have
\begin{itemize}
\item $\calp\cap\tS=G_\tS$ if $\tS$ is belted,
\item $\calp\cap\tS=|\frakC_\tS|$ if $\tS$ is the pseudo-belted side of a
  component of $\cals_2^{!}$, and
\item $\calp\cap\tS=\emptyset$ otherwise.
\end{itemize}
In particular, for each component $\tS$ of $\tcals$, the set $\calp\cap\tS=\calq\cap\tS$ is either the
empty set or the underlying surface of an annular $2$-orbifold, $\pi_1$-injective in $\oldPsi$. Hence the components
of $\obd(\calq)$ are all annular $2$-orbifolds, $\pi_1$-injective in $\oldPsi$.
\EndClaim

To prove \ref{what is calg?}, first consider a belted component $\tS$
of $\tcals$. According to Lemma \ref{new beginning} we have
$\calL\cap\tS=G_\tS$. On the other hand, any component of $\calr$ has
the form $|\frakR_\tT|$ for some pseudo-belted component of $\tcals$; according to the definitions
in \ref{stop beeping}, we have
$|\frakR_\tS|\cap\cals=E_\tT\subset\tT$. Since $\tS$ is belted while $\tT$
is pseudo-belted, Lemma \ref{it's-a unique} implies that $\tS$ and
$\tT$ are distinct components of $\tcals$. Hence
$|\frakR_\tT|\cap\tS=\emptyset$. This shows that
$\calr\cap\tS=\emptyset$. It follows that
$\calp\cap\tS=(\calL\cap\tS)\cup(\calr\cap\tS)=G_\tS\cup\emptyset=G_\tS$,
which proves the first bulleted assertion of \ref{what is calg?}.

Next consider a component $\tS$ of $\tcals$ which is the pseudo-belted side of a
  component of $\cals_2^{!}$.
Since $\tS$ is pseudo-belted, it follows from Lemma \ref{it's-a
  unique} that it is not belted. Hence by Lemma \ref{new beginning} we have
$\calL\cap\tS=\emptyset$. On the other hand, any component of $\cald$ has
the form $|\frakD_\tT|$ where $\tT$ is the pseudo-belted side of some component of $\cals_2^{!}$; according to \ref{replaces readin' john o'hara}, we have
$|\frakD_\tT|\cap\cals=|\frakC_\tT|\subset\tT$. We therefore have
$|\frakD_\tT|\cap\tS=|\frakC_\tS|$ if $\tT=\tS$, and $|\frakD_\tT|\cap\tS=\emptyset$ if
$\tT\ne\tS$. This shows that
$\cald\cap\tS=|\frakC_\tS|$. It follows that
$\calp\cap\tS=(\calL\cap\tS)\cup(\cald\cap\tS)=\emptyset\cup |\frakC_\tS|=
|\frakC_\tS|$,
which proves the second bulleted assertion of \ref{what is calg?}. 

Now consider a component $\tS$ of $\tcals$ which is neither a
belted component, nor the pseudo-belted side of a
  component of $\cals_2^{!}$.
Since $\tS$ is not belted, Lemma \ref{new beginning} implies that
$\calL\cap\tS=\emptyset$. On the other hand, any component of $\calr$ has
the form $|\frakR_\tT|$ where $\tT$ is the pseudo-belted side of some component
of $\cals_2^{!}$; according to the definitions
in \ref{stop beeping}, we have
$|\frakR_\tT|\cap\cals=E_\tT\subset\tT$. Our assumption about $\tS$ implies
that $\tS$ and $\tT$ are distinct components of $\tcals$, so that
$|\frakR_\tT|\cap\tS=\emptyset$. We therefore have
$\calr\cap\tS=\emptyset$. It follows that
$\calp\cap\tS=(\calL\cap\tS)\cup(\calr\cap\tS)=\emptyset\cup\emptyset=\emptyset$.
This completes the proof of the third bulleted assertion of \ref{what is calg?}. 

The remaining assertions of \ref{what is calg?} follow from the  three bulleted assertions in view of the annularity of the $\obd(G_\tS)$ and $\frakC_\tS$ and their $\pi_1$-injectivity. (The $\pi_1$-injectivity of the $\obd(G_\tS)$  was pointed out in \ref{what's a belt?}. Each suborbifold of the form $\frakC_\tS$ is by construction isotopic to $\tau(\obd(G_{\tau(S)}))$, which is $\pi_1$-injective by virtue of the $\pi_1$-injectivity of $\obd(G_{\tau(S)})$ and the admissibility of $\cals$.)

Next, we claim:
\Claim\label{how does calg move?}
If $Q$ is any component of $\calq$, we have $Q= \calp\cap\tS$ for some component $\tS$ of $|\rho|^{-1}(\cals_1\cup\cals_2\cup\cals_2^!)\subset\tcals$. Furthermore, we have $\calp\cap|\tau|(\tS)=|\tau|(Q)$ if $\tS\subset|\rho|^{-1}(\cals_2\cup\cals_2^!)$, and $\calp\cap|\tau|(\tS)=\emptyset$ if $\tS\subset|\rho|^{-1}(\cals_1)$. In particular, $|\tau|(Q)$
either is a component of $\calq$ or is disjoint from $\calq$.
\EndClaim

To prove the first assertion, consider any component $Q$ of $\calq=\calp\cap\tcals$. Since $Q$ is connected it is contained in a component $\tS$ of $\tcals$ and is therefore a component of $\calp\cap\tS$. But it follows from \ref{what is calg?} that $\calp\cap\tS$ has at most one component, and hence $Q= \calp\cap\tS$. 
Since $\calp\cap\tS\ne\emptyset$, \ref{what is calg?} also implies
that either $\tS$ is the pseudo-belted side of a component of $\cals_2^!$,
or that $\tS$ is belted;  it then follows from the
definitions of $\cals_1$, $\cals_2$ and $\cals_2^!$ that $\tS$ is a
side of a component of $\cals_1\cup\cals_2\cup\cals_2^!$.
This proves the first assertion of \ref{how does calg move?}.

To prove the second assertion, set $\tS'=|\tau|(\tS)$ and $S=|\rho|(\tS)$. First consider the case in which
$S$ is a component of $\cals_2$, so that $\tS$ and $\tS'$ are belted.
It follows from \ref{what is calg?} and the first assertion that
$Q=\calp\cap\tS=G_{\tS}$, which with
\ref{fixed it up} gives
$|\tau| (Q)= |\tau| (G_{\tS})=
G_{\tS'}$. It then follows from \ref{what is calg?} that
$\calp\cap\tS'=G_{\tS'}=|\tau|(Q)$, as required for the second assertion. %$\calp\cap|\tau|(\tS)=|\tau|(Q)$.

In the case where $S$ is a component of $\cals_1$, the sphere $S$ has a belted side of and a side which is neither belted nor pseudo-belted. Since $\calp\cap\tS\ne\emptyset$, it follows from \ref{what is calg?} that $\tS$ is the belted side of $S$.
It also follows from \ref{what is calg?} that $\tS'=|\tau|(\tS)$ is disjoint from
$\calp$, which implies
the second assertion of \ref{how does calg
  move?} in this case.

There remains the case in which $S$
  is a component of $\cals_2^{!}$. Consider the subcase in which $\tS$ is the pseudo-belted side of $S$, so that  $\tS'$ is the belted side. By \ref{what is calg?} and the first assertion of \ref{how does calg move?} we have
$Q=\calp\cap\tS=|\frakC_\tS|$; and
%where of $\cals_2^{!}$.
%It follows from \ref{what is calg?} that $\calp\cap\tS=|\frakC_\tS|=Q$, which gives the %first assertion in this case. 
% Let   $\tS'=|\tau|(\tS)$ denote the belted side of $S$. 
according to \ref{fixed it up}, we have
$|\tau| (Q)= |\tau| (|\frakC_\tS|)=G_{\tS'}
$. It follows from \ref{what is calg?} that
$\calp\cap\tS'=G_{\tS'}=|\tau|(Q)$, so that the second assertion holds in this subcase. 

Now consider the subcase in which $\tS$ is the belted side of $S$, so that  $\tS'$ is the  pseudo-belted side. By \ref{what is calg?} and the first assertion of \ref{how does calg move?} we have
$Q=\calp\cap\tS=G_\tS$; and
%where of $\cals_2^{!}$.
%It follows from \ref{what is calg?} that $\calp\cap\tS=|\frakC_\tS|=Q$, which gives the %first assertion in this case. 
% Let   $\tS'=|\tau|(\tS)$ denote the belted side of $S$. 
according to \ref{fixed it up}, applied with the roles of $\tS$ and $\tS'$ reversed, we have
$|\tau| (Q)= |\tau| (G_\tS)=|\frakC_{\tS'}|
$. It follows from \ref{what is calg?} that
$\calp\cap\tS'=|\frakC_{\tS'}|=|\tau|(Q)$, so that the second assertion
holds in this subcase as well.

To prove the final assertion of \ref{how does calg move?}, assume that $|\tau|(Q)\cap\calq\ne\emptyset$. Then according to the first two assertions of \ref{how does calg move?}, there is a component $\tS$ of $|\rho|^{-1}(\cals_2\cup\cals_2^!)\subset\tcals$ such that
 $Q= \calp\cap\tS$ and $\calp\cap|\tau|(\tS)=|\tau|(Q)$. In particular the connected set $|\tau|(Q)$ is contained in $\calp\cap\tcals=\calq$, so that $|\tau|(Q)\subset Q'$ for some component $Q'$ of $\calq$. Again applying the first assertion of \ref{how does calg move?}, this time with $Q'$ playing the role of $Q$ in that assertion, we deduce that  $Q'= \calp\cap\tT$ for some component $\tT$ of $\tcals$. The non-empty set $|\tau|(Q)$ is contained in $\tT$ since $|\tau|(Q)\subset Q'$, and is also contained in $|\tau|(\tS)$; hence the components $\tT$ and $|\tau|(\tS)$ of $\tcals$ must coincide. We now have $|\tau|(Q)\subset Q'=\calp\cap\tT=\calp\cap|\tau|(\tS)=|\tau|(Q)$, so that $|\tau|(Q)=Q'$; in particular, $|\tau|(Q)$ is a component of $\calq$. This completes the proof of \ref{how does calg
  move?}.

%The first two assertions of \ref{how does calg move?} obviously imply the third.
%To prove the first two assertions of \ref{how does calg move?}, first note that for each %component
%$Q$ of $\calq$, it follows from \ref{what is calg?} that $Q$ has
%either the form $G_\tS$ for some belted component  $\tS$ of
%$\tcals$, 
%or the form
%$|\frakC_\tS|$ where $\tS$ is the pseudo-belted side of a
  %component of $\cals_2^{!}$. Consider first the case in which  $Q=G_\tS$
  %for some belted component  $\tS$ of $\tS$. Set $S=|\rho|(\tS)$ and
  %$\tS'=|\tau|(\tS)$, so that $\tS$ and $\tS'$ are the sides of $S$. It follows from %\ref{what is calg?} that $\calp\cap\tS=G_\tS=Q$, which gives the first assertion in %this case. 

We will also need the following fact:

\Claim\label{baby one more fact}
If $\tS$ is a component of $\tcals$ such that $\tS\cap\calq\ne\emptyset$, then for every component $X$ of $\tS\setminus\inter\calq$ we have $\chi(\obd(X))<0$.
\EndClaim

To prove \ref{baby one more fact}, note that by \ref{what is calg?}, either (a) $\tS$ is belted and $\calp\cap\tS=G_\tS$, or (b) $\tS$ is the pseudo-belted side of a
  component $S$ of $\cals_2^{!}$ and $\calp\cap\tS=|\frakC_\tS|$. If (a) holds, then $X$ is a component of $\tS-\inter G_\tS$, and it follows from Lemma \ref{need that too} that $\chi(\obd(X))<0$. If (b) holds, 
then $X$ is a component of $\tS-\inter |\frakC_\tS|$, and 
it follows from \ref{fixed it up} that $|\tau|(|\frakC_\tS|)=G_{\tS'}$, where $\tS'=|\tau|(\tS)$. Hence $X':=|\tau|(X)$ 
is a component 
of $\tS'-\inter G_{\tS'}$. Applying
Lemma \ref{need that too}, with the belted component $\tS'$ playing the role of $\tS$, we deduce that $\chi(\obd(X'))<0$. Since $\tau|\obd(\tS')$ is an orbifold homeomorphism from $\obd(\tS')$ to $\obd(\tS)$, we have
$\chi(\obd(X))=\chi(\obd(X'))<0$, and the proof of \ref{baby one more fact} is complete.

Let us set  $\calw=|\rho|(\calp)\subset M$.
Recall that two distinct points $x$ and
$x'$ of $N$ have the same image
under the surjective map $|\rho|:N\to M$ if and only if they both
lie in $\tcals$ and are interchanged by the involution $|\tau|$ of
$\tcals$. Hence if $\barp$ denotes the quotient space of $\calp$ obtained by
identifying  $x$ with $x'$ for each pair of points
$x,x'\in\calq=\calp\cap\tcals$ that are interchanged by $|\tau|$,
%\redcomment{The bar notation is no good. I've used bar throughout the
  %paper to mean closure. Find another notation}
and if $q:\calp\to\barp$ denotes the quotient map,
then the surjection $|\rho|\big|\calp:\calp\to\calw$ factors as $\barrho
\circ q$,
where $\barrho:\barp\to\calw$ is  a homeomorphism.

It follows from \ref{how does calg move?} that we have $\calq\subset |\rho|^{-1}(\cals_1\cup\cals_2\cup\cals_2^!)$. Hence we may write $\calq$
as the disjoint union of $\calq_1:=\calq\cap |\rho|^{-1}(\cals_1)$ and
$\calq_2:=\calq\cap |\rho|^{-1}(\cals_2\cup\cals_2^!)$. 
It follows from \ref{how does calg move?} that we have
\Equation\label{oh all right}
|\tau|(\calq_2)=\calq_2\quad\text{ and }\quad
|\tau|(\calq_1)\cap\calq=\emptyset.
\EndEquation
 Since  the involution $|\tau|$
leaves no component of $\tcals$
invariant, for each component $Q$ of $\calq_2$ the set $|\tau|(Q)$, which by (\ref{oh all right}) is a
component of $\calq_2$, is distinct from $Q$. With
(\ref{oh all right}) this implies:
\Claim\label{carp} The manifold  $\barp$ is obtained
from $\calp$ by gluing together, via the restriction of $|\tau|$, each
pair of components of $\calq_2$ that are interchanged by $|\tau|$. 
\EndClaim
It follows from \ref{carp} and (\ref{oh all right}) that
$\barp$ is a compact $3$-manifold, and that
$\barq_1:=q(\calq_1)$ is a compact submanifold of $\partial\calw$,
while 
$\barq_2:=q(\calq_2)$
%$\calb_2$ 
is a properly embedded $2$-manifold in $\calw$. The
surjection $q:\calp\to\barp$ maps $\calq_1$ homeomorphically onto
$\barq_1$, maps each component of $\calq_2$ homeomorphically onto a
component of $\barq_2$, and maps $\calp-\calq\supset\inter\calp$ homeomorphically onto $\barp-\barq$, where $\barq=\barq_1\cup\barq_2$. We have $\barq_1\cap\barq_2=\emptyset$. Since
$\barrho:\barp\to\calw$ is a  homeomorphism, we deduce:

\Claim\label{still sick}
The set $\calw$ is a compact $3$-manifold, and
$\calb_1:=|\rho|(\calq_1)$ is a compact submanifold of $\partial\calw$
while $\calb_2:=|\rho|(\calq_2)$
is a properly embedded $2$-manifold in $\calw$. We have $\calb_1\cap\calb_2=\emptyset$. The
surjection $|\rho||\calp:\calp\to\calw$ maps $\calq_1$ homeomorphically
to $\calb_1$, maps each component of $\calq_2$ homeomorphically to
a component of $\calb_2$, and maps $\calp-\calq\supset\inter\calp$
homeomorphically onto $\calw-\calb$, where
$\calb=|\rho|(\calq)=\calb_1\cup\calb_2\subset\calw$. Furthermore, the pre-image under $|\rho|$ of each component of $\calb_2$ is the union of two components of $\calb$.
\EndClaim

Since $\rho:\oldPsi\to\oldOmega$ is an orbifold immersion, its restriction to any suborbifold of $\oldPsi$, such as $\obd(\calp)$ or any component of $\obd(\calq)$, is an immersion. Now observe that, if $g:\oldXi_1\to\oldXi_2$ is an immersion of orbifolds such that $|g|:|\oldXi_1|\to|\oldXi_2|$ is a homeomorphism of spaces, then $g$ is an orbifold homeomorphism. Hence \ref{still sick} implies:

\Claim\label{sicker still}
The orbifold $\obd(\calw)$ is  compact, and
$\obd(\calb_1)=\rho(\obd(\calq_1))$ is a compact submanifold of $\partial(\obd(\calw)$,
while $\obd(\calb_2)=\rho((\obd(\calq_2))$
is a properly embedded $2$-orbifold in $\calw$. We have $\obd(\calb_1)\cap\obd(\calb_2)=\emptyset$. The
surjective orbifold immersion $\rho|\obd(\calp):\obd(\calp)\to\obd(\calw)$ maps $\obd(\calq_1)$ homeomorphically
to $\obd(\calb_1)$, maps each component of $\obd(\calq_2)$ homeomorphically to
a component of $\obd(\calb_2)$, and maps $\obd(\calp-\calq\supset\inter\calp)$
homeomorphically onto $\obd(\calw-\calb)$. (Here
$\obd(\calb)=\obd(\calb_1)\cup\obd(\calb_2)$.) Furthermore, the pre-image under $\rho$ of each component of $\obd(\calb_2)$ is the union of two components of $\obd(\calq_2)$.
\EndClaim

Next we claim:

\Claim\label{zusammen}
Each component of $\obd(\calb)$ is an annular $2$-orbifold, and is $\pi_1$-injective in $\obd(\calw)$.
\EndClaim

To prove \ref{zusammen}, let $B$ be any component of $\calb$. By \ref{sicker still}, there is a component $Q$ of $\calq$ such that $\rho$ maps $\obd(Q)$  homeomorphically onto $\obd(B)$. Since  $\obd(Q)$ is annular by \ref{what is calg?}, it follows that  $\obd(B)$ are  annular. Now \ref{what is calg?} also asserts that $\obd(Q)$ is $\pi_1$-injective in $\oldPsi$. But since $\tcals$ is \dandy, and in particular admissible, $\rho:\oldPsi\to\oldOmega$ is $\pi_1$-injective. Hence $\rho|\obd(Q)$, regarded as a map from $\obd(Q)$ to $\oldOmega$, is $\pi_1$-injective. Since $\rho$ maps $\obd(Q)$  homeomorphically onto $\obd(B)$, it follows that $\obd(B)$ is $\pi_1$-injective in $\oldOmega$. This completes the proof of \ref{zusammen}.

We now claim:

\Claim\label{my mcblister}
Every component of $\partial W$ contains the image under $|\rho|$ of some component  of $\Fr_N\calp$.
\EndClaim

To prove \ref{my mcblister}, let $T$ be a component of $\partial W$. It follows from \ref{still sick} that $\calb\cap T$ is a union of components of $\partial\calb_1$ and of $\calb_2$, and that $(\partial\calb_1)\cap\calb_2=\emptyset$. For each component $B$ of $\calb_2$, the orbifold $\obd(B)$ is annular by \ref{zusammen}. Hence each component of $\calb_2$ is an annulus or a disk. Since $T$ is a closed surface, it cannot be a disjoint union of simple closed curves, annuli and disks; hence $T\not\subset\calb$. Choose a component $Z$ of $T\setminus\calb$. Then $Z$ is in particular a component of $\calw-\calb$; hence by \ref{still sick}, $Z$ is the homeomorphic image under $|\rho|$ of some component $T$ of $\calp-\calq=\Fr_N\calp$. In particular we have $|\rho|(Y)\subset T$. Thus \ref{my mcblister} is proved.

We also claim:
\Claim\label{anemone}For every component $T$ of
  $\partial\calw$, we have $\chi(\obd(T))=0$.
\EndClaim

To prove \ref{anemone}, first recall that by \ref{they're torifolds},
%from \ref{old trier}
%\redmissingref{be sure that all this stuff is stated there} that 
%for each
%belted component $\tS$ of $\tcals$, the orbifold $\obd(L_\tS)$ is a
%\torifold, and $\obd(G_\tS)$ is an annular $2$-orbifold contained in %$\obd(\partial
%L_\tS)$  and $\pi_1$-injective in $L_\tS$. As we have observed that each component of $\cald$ is a 
%\torifold, it now follows that
 each component of $\calp$ is a
\torifold; and that by \ref{what is calg?}, each component of $\obd(\calq)$---and in particular
each component of $\obd(\calq_2)$---is an annular $2$-orbifold,
%Euler characteristic $0$,
and is $\pi_1$-injective in $\obd(\calp)$. Hence each component of $\obd(\overline{\partial\calp-\calq_2})$
%  an annular $2$-orbifold, and therefore 
has
  Euler characteristic $0$. But it follows from \ref{sicker still} that $\obd(\calw)$ is homeomorphic to an orbifold obtained from $\obd(\calp)$
  by gluing together the components of $\obd(\calq_2)\subset\partial\obd(\calp)$
  in pairs; hence $\obd(\partial\barp)$ is homeomorphic to an orbifold obtained from
  $\obd(\overline{\partial\calp-\calq_2})$ by gluing the     boundary components of 
  $\obd(\overline{\partial\calp-\calq_2})$ in pairs. It follows that each
    component of $\obd(\partial\calw)$ has Euler characteristic $0$, and
\ref{anemone} is proved.

Next we claim:
\Claim\label{P-claim}
If $Y$ is any component of $\Fr_N\calp$, then $Y\cap\tcals\ne\emptyset$. 
\EndClaim

To prove \ref{P-claim}, let $P$ denote the component of $\calp$ containing $Y$. 
%We have   $P\cap\tcals\ne\emptyset$. 
According to \ref{they're torifolds}, $\obd(P)$ is a \torifold; by Proposition \ref{three-way equivalence}, this implies that $P$ is either a ball or a solid torus, and in particular $\partial P$ is connected. On the other hand, it follows from \ref{they're torifolds} that $P\cap\tcals\ne\emptyset$. Since by \ref{what is calg?} every component of $P\cap\tcals$ has non-empty boundary, we have 
$(\partial P)\cap\tcals\ne\emptyset$. In view of the connectedness of $\partial P$, it follows that every component of
$(\partial P)\setminus \tcals$ has non-empty frontier relative to $\partial P$. Since $Y$ is the closure of a component of 
$(\partial P)\setminus \tcals$, it follows that $Y$ has non-empty boundary and therefore meets $\tcals$.
Thus \ref{P-claim} is proved.

We also claim: 
\Claim\label{just put it together}
Let $W$ be any
component of
$ \calw$. Then for every
component $T$ of
$\partial W$,
there exist a negative, compact, connected $2$-orbifold  $\frakZ$  with connected boundary, and an (orbifold) immersion $f:\frakZ\to\obd(\overline{M-W})$ such that $|f|^{-1}(\partial W)=|f|^{-1}(T)=\partial\frakZ$, and such that $f$, regarded as an (orbifold) immersion of $\frakZ$ in $\oldOmega$, is $\pi_1$-injective.
\EndClaim

To prove \ref{just put it together}, note that by \ref{my mcblister}, the component $T$ of $\partial W$ contains $|\rho|(Y)$ for some component $Y$ of $\Fr_N\calp$. %sincee $\calw$ is by definition equal to $|\rho|(\calp)$, there is a component $P$ of %$\calp$ such that $P\subset W$. 
According to \ref{P-claim} we have $Y\cap\tcals\ne\emptyset$. Fix a component $c$ of $Y\cap\tcals=\partial Y$. Then $c$ is a component of $\partial Q$ for some
% denote so that $P$ contains at least one 
component $Q$ of $\calq=\calp\cap\tcals$. According to \ref{how
 does calg move?},
we have $Q=\calp\cap\tS$ for some component $\tS$ of $\tcals$, and we
have either (1) $\calp\cap|\tau|(\tS)=|\tau|(Q)$, or (2)
$\calp\cap|\tau|(\tS)=\emptyset$. 
In either case, it follows that if we set
$S=|\rho|(\tS)=|\rho|(|\tau|(\tS))$, then $|\rho|$ maps $Q$ homeomorphically
onto a component $B$  of $\calw\cap\cals$, contained in the component
$S=|\rho|(\tS)$ of $\cals$; and that $B=|\rho|(\calp)\cap S=\calw\cap
S$. 
%In particular, $|\rho|$ maps $c$ homeomorphically onto a component
%$c$ of $\partial B$. Since $c\subset Y$, we have $c\subset %T$. Since
%$c\subset\partial W\cap\cals$, we have
%$c\cap\fraks_\oldOmega=\emptyset$. \redcomment{The  last two sentences, and maybe more, may be unnecessary.

Since $\tS$ is a sphere and $Q$ is connected, the  simple closed curve $c$ bounds a disk component of $\tS-\inter Q$, which may be written as $|\frakZ|$ for some connected suborbifold $\frakZ$ of $\obd(\tS)$. According to \ref{baby one more fact}, we have $\chi(\frakZ)<0$. Since $Q=\calp\cap\tS$, we have
$|\inter\frakZ|\cap\calp\subset|\inter\frakZ|\cap Q=\emptyset$.

If Alternative
(1) above holds, i.e. if $\calp\cap|\tau|(\tS)=|\tau|(Q)$, 
then $|\tau(\inter\frakZ)|\subset|\tau|(\tS)-|\tau|(Q)=|\tau|(\tS)\setminus\calp$, so that $|\tau(\inter\frakZ)|\cap\calp=\emptyset$.
If Alternative (2) holds, i.e. if
$\calp\cap|\tau|(\tS)=\emptyset$, then in particular $|\tau(\inter\frakZ)|\cap\calp=\emptyset$. Thus in any event we have $|\tau(\inter\frakZ)|\cap\calp=\emptyset$. Since we have also seen that
$|\inter\frakZ|\cap\calp=\emptyset$, it follows that $|\rho(\inter\frakZ)|\cap\calw=\emptyset$. 
%\redmissingref{Get from this---and make sure the whole thing is clear---to asserting that}
%$|\rho|$
%maps
%the disk $\frakZ$  homeomorphically onto a properly embedded disk $X\subset\overline{M-\calw}$. 
But 
%$|\rho(\partial\frakZ)|\subset |\rho|(\partial Q)=\partial B\subset\partial \calw$ (where the last equality and inclusion follow from \ref{still sick}), and 
$|\rho|(\partial\frakZ)=|\rho|(c)\subset T\subset\partial\calw $. Hence if we set $f=|\rho\big|\frakZ|$, we have $|f|^{-1}(\partial\calw)=|f|^{-1}(T)=\partial\frakZ$. 
%It also follows that $X$  is a component of $S-\inter Z$, We have $\chi(\obd(X))=\chi(\obd(\frakZ))<0$.

By \ref{what is calg?}, $\obd(Q)$ is a $\pi_1$-injective annular
$2$-orbifold in $\tS$; hence $\obd(B)$ is annular and $\pi_1$-injective in $\obd(S)$. Since $\obd(B)$ is in particular not discal, the $\pi_1$-injectivity of $\obd(B)$ in $\obd(S)$ implies that $\obd(\partial B)$ is $\pi_1$-injective in $\obd(S) $. Hence  $\obd(\overline{ S- B})$ is $\pi$-injective in $\obd(S) $, and in particular its component $\frakZ$ is  $\pi$-injective in $\obd(S)$. But $S$ is a component of $\cals$, which is \dandy\ and in particular admissible, so that $\rho|\obd(\tS)$ is  $\pi$-injective in $\oldOmega$. Hence $f=\rho|\frakZ:\frakZ\to\oldOmega$ is $\pi_1$-injective, and
% This proves the $\pi_1$-injectivity of $\obd(X)$ in $\oldOmega$, and  
the proof of \ref{just put it together} is complete.
%Z \tX still

From \ref{anemone}, \ref{just put it together} and 
Lemma \ref{occupani}, we immediately deduce:
\Claim\label{purty solid}
If  $W$ is a component of $\calw$, and $T$ is any component of
$\partial W$, then $W$  is contained in a
 submanifold $J$
of $M$ such that $\obd(J)$ is a \torifold, and
$\partial J=T$.
\EndClaim

Now we claim:

\Claim\label{don't die ET}
For every component $W$ of $\calw$, the orbifold $\obd(W)$ is a \torifold.
\EndClaim

To prove \ref{don't die ET}, first note that
if $W$ is a component of $\calw$ whose boundary is connected, then
it is immediate from \ref{purty solid} that $\obd(W)$ is a torifold. If $W$ is a component of $\calw$ whose boundary is not connected, fix distinct boundary components $t_1$, $t_2$ of $W$. According to \ref{purty solid}, for $i=1,2$ there is a
%If  $W$ is component of $\calw$, and $T$ is any component of $\partial W$, then $W$  is contained in a 
submanifold $J_i\subset M$ with $W\subset J_i$ and $\partial J_i=T_i$,
such that $\obd(J_i)$ is a \torifold. In particular we have $\partial J_1=T_1\subset\inter J_2$ and $\partial J_2=T_2\subset\inter J_1$. Hence  $\partial(J_1\cup J_2)=\emptyset$.
Since $M$ is connected, it
follows that $M=\inter J_1\cup\inter J_2$. Since $\obd(J_i)$ is a \torifold\ for each $i$, each of the $J_i$ is either a solid torus or a ball. The closed surface $t_1$ is contained in $\inter J_2$, and since $J_1$ is a ball or solid torus, $t_1$ separates $J_2$; similarly, $t_2$ separates $J_1$. Hence $J_1\cap J_2$ is a connected $3$-manifold whose boundary components are $t_1$ and $t_2$. For $i=1,2$, since $\obd(J_i)$ is a \torifold, the inclusion homomorphism $\pi_1(\obd(\partial J_i)\to\pi_1(\obd( J_i))$ is surjective; in particular, the inclusion homomorphism $\pi_1(\obd(J_1\cap J_2))\to\pi_1(\obd( J_i))$ is surjective for $i=1,2$. By van Kampen's theorem for orbifolds, the inclusion homomorphism $\pi_1(\obd(J_1\cap J_2))\to\pi_1(\obd( J_1\cup J_2))=\pi_1(\oldOmega)$ is surjective, and in particular so is the inclusion homomorphism $\pi_1(\obd(J_1))\to\pi_1(\oldOmega)$. But $\pi_1(\obd(J_1))$ is virtually cyclic since $\obd(J_1)$ is a \torifold, and  so $\pi_1(\oldOmega)$ is virtually cyclic.
This contradicts the hyperbolicity of $\Mh$, and so
\ref{don't die ET} is proved in all cases.

Now we claim:

\Claim\label{comic cockroach}
Each component of $\calb_2$ is a  separating surface in the component of $\calw$ containing it.
\EndClaim

To prove \ref{comic cockroach}, let $B$ be any component of $\calb_2$,
and let $W$ denote the component of $\calw$ containing $B$. According
to \ref{still sick}, $B$ is properly embedded in $W$. According to
\ref{don't die ET}, $\obd(W)$ is a \torifold. According to
\ref{zusammen}, 
$\obd(B)$ is an annular $2$-orbifold and is $\pi_1$-injective in
$\obd(W)$. Now since $\obd(W)$ is a \torifold, it follows from
Proposition \ref{three-way equivalence} that either (a) $W$ is a ball or (b) $W$ is a
solid torus and $\pi_1(\obd(W))$ is isomorphic to
$\ZZ\times(\ZZ/m\ZZ)$ for some $m\ge1$. Every properly embedded
surface in a $3$-ball is separating; hence if (a) holds then the
conclusion of \ref{comic cockroach} holds. Now suppose that (b)
holds. Since  $\obd(B)$ is annular and orientable, $B$ is
either a weight-$0$ annulus, or a weight-$2$
disk whose singular points are of order $2$. In the latter subcase,
$\pi_1(\obd(B))$ is an infinite dihedral group, and is therefore
non-abelian. Since $\pi_1(\obd(W))$ is abelian in Case (b), we have a contradiction to the
$\pi_1$-injectivity of $\obd(B)$  in $\obd(W)$. Hence if (b) holds, $\obd(B)$ is 
a weight-$0$ annulus. It follows that $\pi_1(\obd(B))$ is infinite cyclic, and in particular torsion-free, so that the image of the inclusion homomorphism $\pi_1(\obd(B))\to\pi_1(\obd(W))$ has trivial intersection with the torsion subgroup of $\pi_1(\obd(W))\cong\ZZ\times(\ZZ/m\ZZ)$, which is the kernel of the natural homomorphism $\pi_1(\obd(W))\to\pi_1(W)$. It follows that $B$ is $\pi_1$-injective in $W$. Since every $\pi_1$-injective, properly embedded annulus in a solid torus $W$ is a separating surface in $W$, the conclusion of  \ref{comic cockroach} holds in this case as well, and the proof of \ref{comic cockroach} is complete.

We also claim:

\Claim\label{one more}
We have $\calb\subset \cals$. Furthermore, $S\cap\calb$ is connected for every component $S$ of $\cals_1\cup\cals_2\cup\cals_2^!$, while $S\cap\calb=\emptyset$ for every component $S$ of $\cals_0$. Hence each component of $\calb$ is contained in a component of $\cals_1\cup\cals_2\cup\cals_2^!$, and each component of $\cals_1\cup\cals_2\cup\cals_2^!$ contains a unique component of $\calb$. 
\EndClaim

To prove \ref{one more}, first note that the definition of $\calq$ gives $\calq\subset\tcals$, so that $\calb=|\rho|(\calq)\subset|\rho|(\tcals)=\cals$, and the first assertion is proved. Next note that if $S$ is a component of $\cals_1\cup\cals_2\cup\cals_2^!$, then by definition $S$ has at least one belted side; let us fix a belted side $\tS$ of $S$. By \ref{what is calg?}, $Q:=\calq\cap\tS$ is connected. The other side of $S$ is $\tS':=|\tau|(\tS)$, and by \ref{how does calg move?} we have either $\calq\cap\tS'=|\tau|(Q)$ or $\calq\cap\tS'=\emptyset$. In either case, 
$S\cap\calb=|\rho|((\calq\cap\tS)\cup(\calq\cap\tS'))=|\rho|(Q)$, so that $S\cap\calb$ is connected. Finally, if $S$ is a component of $\cals_0$, and $\tS$ and $\tS'$ denote the sides of $S$, then \ref{what is calg?} implies that $\calq\cap\tS=\calq\cap\tS'=\emptyset$, and hence $S\cap\calb=\emptyset$. This proves the second assertion of \ref{one more}; since the final assertion follows from the first two, the proof of \ref{one more} is complete.

Let us now construct a $3$-manifold $\calw_1$ from the disjoint union of $W$ with $\calb_1\times[0,1]$ by gluing $\calb_1\subset\partial\calw$ to $\calb_1\times\{0\}\subset\calb_1\times[0,1]$ by the homeomorphism $x\to(x,0)$. Thus $\calw$ is canonically  identified with a submanifold of $\calw_1$.
For each component $B$ of $\calb_1$, the component $B\times[0,1]$ of $\calb_1\times[0,1]$ is likewise canonically identified with a submanifold of $\calw_1$, which we will denote by $C_B$.
According to \ref{still sick},
$\calb_1:=|\rho|(\calq_1)$ is a compact submanifold of $\partial\calw$,
while $\calb_2:=|\rho|(\calq_2)$
is a properly embedded $2$-manifold in $\calw$, and $\calb_1\cap\calb_2=\emptyset$. Hence $\calb=\calb_1\cup\calb_2\subset\calw\subset\calw_1$ is properly embedded in $\calw_1$. 
Since, by \ref{comic cockroach}, each component of $\calb_2$ is a  separating surface in the component of $\calw$ containing it, we deduce:
\Claim\label{cockroach redux}
Each component of $\calb$ is a  separating surface in the component of $\calw_1$ containing it.
\EndClaim

We now turn to the proof that $\calf=\calf\oldomecals$ is a forest. For each component $S$ of $\cals_0$, the sides $\tS$ and $\tS'$ are by definition unbelted, so that $v_\tS$ and $v_{\tS'}$ are distinct ideal vertices incident to $e_S$. Since an ideal vertex has valence $1$, the subgraph $I_S$ consisting of the edge $e_S$ and the vertices $v_{\tS}$ and $v_{\tS'}$ is a component  of $\calf$, which is clearly a tree. Hence to prove that $\calf$ is a forest, it suffices to prove that $\calf^*:=\calf-\bigcup_{S\in\calc(\cals_0)}I_S$ is a forest.

For this purpose, consider the dual graph (see \ref{dual graph}) $\calg$ of the two-sided surface $\calb$ in the orientable
$3$-manifold $\calw_1$. 
%By definition, the vertices of $\calg$ are in
%bijective correspondence with the components of $\calw_1-\calb$; the
%edges of $\calg$ are in bijective correspondence with the components
%of $\calb$; and a vertex $v$ of $\calg$ is incident to an edge $e$ if
%and only if the component of $\calb$ corresponding to $e$ is contained
%in the closure of the component of $\calw_1-\calb$ corresponding to
%$v$. 
It follows from \ref{cockroach redux} that $\calg$ is a
forest. We will show that $\calf^*$ is isomorphic to $\calg$, which
will prove that $\calf^*$ is a forest. We must define a bijection
$K\mapsto v^K$ from $\calc(\calw_1-\calb)$ to the vertex set
$\calv^*\subset\calv\oldomecals$ of $\calf^*$, and a bijection
$B\mapsto e^B$ from $\calc(\calb)$ to the edge set
$\cale^*\subset\cale\oldomecals$ of $\calf^*$, such that for each
$K\in\calc(\calw_1-\calb)$ and each $B\in\calc(\calb)$, the vertex
$v^K$ is incident to the edge $e^B$ if and only if
$B\subset\overline{K}$.

To construct the bijection $B\mapsto e^B$, note that according to \ref{one more}, each component of $\calb$ is contained in a component of $\cals_1\cup\cals_2\cup\cals_2^!$, and each component of $\cals_1\cup\cals_2\cup\cals_2^!$ contains a unique component of $\calb$. Furthermore, the definition of $\calf^*$ implies that $\cale^*=\{e_S:S\in\calc (\cals_1\cup\cals_2\cup\cals_2^!\}$. Hence we may define a bijection $B\mapsto e^B$ from $\calc(\calb)$ to $\cale^*$   by setting $e^B=e_S$, where $S$ is the component of $\cals$ containing $B$. 

The construction of the bijection $K\mapsto v^K$ is somewhat more involved. First note that, according to the construction of $\calw_1$, each $K\in\calc(\calw_1-\calb)$ has exactly one of two forms: either 
$K=C_B-B$ for a unique component $B$ of $\calb_1$, or 
$K$ is a component of $\calw-\calb_2$. 
For a component $K$ of  the latter type, it follows from \ref{still sick} that we have $K=|\rho|(P\setminus\calq)\supset|\rho|(\inter P)$ for a unique component $P$ of $\calp=\calL\discup\cald$. 
Thus we may write $\calc(\calw_1-\calb)$ as a disjoint union $\calk_1\discup\calk_2\discup\calk_3$, where $\calk_1$, $\calk_2$ and $\calk_3$ consist, respectively, of the components of $\calc(\calw_1-\calb)$ such that
(1) $K=C_B-B$ for a (unique) component $B$ of $\calb_1$, (2)
$K=|\rho|(L\setminus\calq) \supset|\rho|(\inter L)$ for a (unique) component $L$ of $\calL$, or
(3) $K=|\rho|(R\setminus\calq) \supset|\rho|(\inter R)$ for a (unique) component $R$ of
$\cald$. 

On the other hand, the definition of $\calf^*$ gives that $\calv^*$ consists of (i) all material vertices of $\calf$, and (ii) all ideal vertices  of the form $v_\tS$ where $\tS$ is an unbelted side of a component $S$ of $\cals$ that is not contained in $\cals_0$. A component $S$ of $\cals$ of the kind mentioned in (ii) cannot be contained in $\cals_2$, as a component of $\cals_2$ has two belted sides; hence it is contained in $\cals_1$ or in $\cals_2^!$. We may therefore write $\calv^*$ as a disjoint union $\calv_1^*\discup\calv_2^*\discup\calv_3^*$, where $\calv_1^*$, $\calv_2^*$ and $\calv_3^*$ consist, respectively, of the vertices $v$ having the form
(1) $v=v_\tS$ for a (unique) unbelted component $\tS$ of $|\rho|^{-1}(\cals_1)$, (2) $v=v_L$ for a (unique) component $L$ of $\calL$, or (3) $v=v_\tS$ for a (unique) unbelted component $\tS$ of $|\rho|^{-1}(\cals_2^!)$.

%Let us define a bijection from $\calk_1$ to $\calv_1$ by $C_B-B\mapsto
If $K\in\calk_1$ is given, and if we write $K=C_B-B$ with $B\in\calc(\calb_1)$, then in view of the definitions of $\calq_1$ and $\calb_1$, we have
$B\subset\calb_1=\rho(\calq_1)\subset\cals_1$.
%\ref{still sick} we have $B=|\rho|(Q)$ for a unique component $Q$ of 
%$\calq_1=\calq\cap |\rho|^{-1}(\cals_1)$; in particular we have
%$B\subset\cals_1$. 
We set $v^K=v_\tS$, where $\tS$ is the unique
unbelted side of the component of $\cals_1$ containing $B$. By
definition we have $v^K\in\calv_1^*$. 

On the other hand, if
$v\in\calv_1^*$ is given, and if we write $v=v_\tS$, where $\tS$ is
the unbelted side of a component $S$ of $\cals_1$, then the definition
of $\cals_1$ says that $S$ has a unique belted side $\tS'$ and has no
pseudo-belted side. It follows from \ref{what is calg?} that
$Q:=G_{\tS'}$ is the unique component of $\calq$ 
that meets
$|\rho|^{-1}(S)$; and since $S\subset\cals_1$, we have $Q\subset\calq_1$ by definition. It now follows that
%contained in $\tS$. Furthermore, $Q$ is a component of $\calq_1$ since $S\subset\cals_1$. If we set
$B:=|\rho|(Q)$, 
which by \ref{still sick} is a component of $\calb_1$, 
is the unique component of $\calb$ contained in $S$, and hence 
%it now follows that
$K:=C_B-B$ is the unique element of $\calk_1$ with $v^K=v$. This shows that the assignment $K\mapsto v^K$ is a bijection from $\calk_1$ to $\calv^*_1$.

If $K\in\calk_2$ is given, and if we write 
$K=|\rho|(L\setminus\calq)$  with $L\in\calc(\calL)$, then $L$ is uniquely determined by $K$ according to \ref{still sick}, and $v_L$ is by
definition an element of $\calv_2^*$. We set $v^K=v_L$. 

On the other hand, if an
arbitrary element $v$ of $\calv_2^*$ is given, and if we write
$v=v_L$, where $L$ is a component of $\calL$, then
$K:=|\rho|(L\setminus\calq)$ is by definition the unique element of
$\calk_2$ with $v^K=v$. Thus the assignment $K\mapsto v^K$ is a
bijection from $\calk_2$ to $\calv^*_2$.

If $K\in\calk_3$ is given, and if we write 
$K=|\rho|(R\setminus\calq)$ with $R\in\calc(\cald)$, then by definition
we have $R=|\frakD_\tS|$, where $\tS$ is the pseudo-belted side of a component of $\cals_2^{!}$. By \ref{still sick}, $R$ is uniquely determined by $K$, and $\tS$ since 
$|\frakD_\tS|\cap\tcals=|\frakC_\tS|$ by 
\ref{replaces readin' john o'hara}, $\tS$ is also uniquely determined by $K$.
We set $v^K=v_\tS$; by definition we have $v^K\in\calv_3^*$. 

On the other hand, if an arbitrary element $v$ of $\calv_3^*$ is given, and if we write $v=v_\tS$, where $\tS$ is the pseudo-belted side of a component of $\cals_2^!$, then $K:=|\rho|(|\frakD_\tS|\setminus\calq)=|\rho|(|\frakD_\tS|-|\frakC_\tS|)$ is by definition the unique element of $\calk_3$ with $v^K=v$. Thus the assignment $K\mapsto v^K$ is a bijection from $\calk_3$ to $\calv^*_3$.

Since we have expressed 
$\calc(\calw_1-\calb)$ and $\calv^*$  as disjoint unions $\calk_1\discup\calk_2\discup\calk_3$ and $\calv_1^*\discup\calv_2^*\discup\calv_3^*$, the observations about bijectivity made above show that the assignment $K\mapsto v^K$ is a bijection from $\calc(\calw_1-\calb)$ to $\calv^*$.

To complete the proof that $\calf$ is a forest, it remains to show:
\Claim\label{what's left}
Given $K\in\calc(\calw_1-\calb)$ and $B\in\calc(\calb)$, the vertex $v^K$ is incident to the edge $e^B$ if and only if $B\subset\overline{K}$. 
\EndClaim

To prove \ref{what's left}, let $S$ denote the component of $\cals$
containing $B$, so that $e^B=e_S$. First consider the case in which
$K\in\calk_1$, so that $K=C_{B'}-B'$ for some component $B'$ of
$\calb_1$. In this case we have $v^K=v_{\tS_1}$, where $\tS_1$ is the
unique unbelted side of the component $S_1$ of $\cals_1$ containing
$B'$. The definition of $\calf$ implies that $v_{\tS_1}$ is incident
to $e_S$ if and only if $\tS_1$ is a side of $S$, i.e. if and only if
$S_1=S$. 
On the other hand, the definition of $\calw_1$ implies that
we have $B\subset\overline K=C_{B'}$ if and only if $B=B'$.  The final
assertion of \ref{one more} implies that $B=B'$ if and only if
$S=S_1$, and thus \ref{what's left} is proved in this case.

Next consider the case in which $K\in\calk_2$, so that
$K=|\rho|(L\setminus\calq)$ for some component $L$ of $\calL$.  In this
case, we have $v^K=v_L$. The definition of $\calf$ implies that
$v_{L}$ is incident to $e_S$ if and only if $L=L_\tS$ for some belted
side $\tS$ of $S$.   Thus to prove \ref{what's left} in this case, we must show that $B\subset\overline K=|\rho|(L)$ if and only if
$L=L_\tS$ for some belted side $\tS$ of $S$. If $L=L_\tS$ for some
belted side $\tS$ of $S$, we have $\tS\cap\calq=G_\tS$ by \ref{what is
  calg?}, and it follows from \ref{how does calg move?} that
$|\tau|(\tS)\cap\calq$ is either $|\tau|(G_\tS)$ or the empty set. Hence
$|\rho|(G_\tS)$ is the unique component of $\calb=|\rho|(\calq)$ contained
in $S$, and therefore $B=|\rho|(G_\tS)\subset
|\rho|(L_\tS)=|\rho|(L)$. Conversely, if $B\subset|\rho|(L)$, then in
particular $S\cap|\rho|(L)\ne\emptyset$, so that $\tS\cap L\ne\emptyset$
for some side $\tS$ of $S$. According to Lemma \ref{new beginning}, this
implies that $\tS$ is belted and that $\calL\cap\tS=G_\tS$; in
particular, $L_\tS$ is the only component of $\calL$ having non-empty
intersection with $\tS$, and hence $L=L_\tS$. Thus \ref{what's left}
is proved in this case.

The final case of \ref{what's left} is the one in which $K\in\calk_3$, so that $K=|\rho|(R\setminus\calq)$ for some component $R$ of $\cald$. 
We then have $R=|\frakD_{\tS_1}|$, where ${\tS_1}$ is the pseudo-belted side
of a component $S_1$ of $\cals_2^{!}$, and $v^K=v_{\tS_1}$. 
The definition of $\calf$ implies that $v_{{\tS_1}}$ is incident to $e_S$ if and only if ${\tS_1}$ is a side of $S$, i.e. if and only if $S_1=S$. Thus we must show that $B\subset \overline K=|\rho|(R)$ if and only if $S_1=S$. First suppose that $S_1=S$. We have ${\tS_1}\cap\calq=|\frakC_{\tS_1}|$ by \ref{what is calg?}, and it follows from \ref{how does calg move?} that $|\tau|({\tS_1})\cap\calq$ is either $|\tau|(|\frakC_{\tS_1}|)$ or the empty set. Hence $|\rho|(|\frakC_{\tS_1}|)$ is the unique component of $\calb=|\rho|(\calq)$ contained in $S_1=S$. We therefore have $B=|\rho|(|\frakC_{\tS_1}|)\subset |\rho|(|\frakD_{\tS_1}|)=|\rho|(R)$. Conversely, if $B\subset|\rho|(R)$, then in particular $S\cap|\rho|(R)\ne\emptyset$, so that ${\tS}\cap R\ne\emptyset$ for some side ${\tS}$ of $S$. But by \ref{just fits}, we have
$|\frakD_{\tS_1}|\cap\tcals= |\frakC_{\tS_1}|$, so that in particular $\tS_1$ is the only component of
$\tcals$ that meets $R=|\frakD_{\tS_1}|$. Hence $\tS_1=\tS$, and
therefore $S_1=S$. Thus \ref{what's left} is proved in all cases; hence $\calf$ is a forest, and the first assertion of the proposition
is established.
%\tS

To prove the second assertion of the proposition, suppose that some component $\calt$ of $\calf$ has at least two special vertices. By the definition of $\calf^*$, every component of $\calf$ either (a) is a component of $\calf^*$, or (b) consists of the edge $e_S$ for some $S\in\calc(\calc_0)$ and the two ideal vertices incident to $e_S$. In each case, we shall obtain a contradiction.

If the component $\calt$ satisfies (b), and the sides of $S$ are denoted by $\tS$ and $\tS'$,
then both $v_{\tS}$ and $v_{\tS'}$ are special ideal vertices. By
definition this means that $\tS$ and $\tS'$ are pseudo-belted sides of
$S$. Thus $S$ satisfies Alternative (ii) of Definition \ref{self-clash
  def} and is therefore a \centralclashsphere\ of $\cals$. This
contradicts the hypothesis of the proposition.

The rest of the proof will be devoted to obtaining a contradiction in
Case (a), in which $\calt$ is a component of $\calf^*$.

We have shown that the maps $K\mapsto v^K$ and $B\mapsto e^B$ define
an isomorphism of the dual graph of $\calb$ in $\calw_1$ to
$\calf^*$. Hence there is a component $W_1$ of $\calw_1$ such that
these maps give an isomorphism from the dual graph of $\calb\cap W_1$ in
$W_1$ to $\calt$. Fix distinct components $K_1$ and $K_2$ of
$W_1\setminus\calb$ such that $v_i:=v^{K_i}$ is a special vertex for
$i=1,2$. We claim:

\Claim\label{one or two} For each $i\in\{1,2\}$, the suborbifold
$\obd(\overline{K_i})$  of $\oldOmega$ is a \torifold. Furthermore, for every component $B$ of
$\Fr_{W_1}\overline{K_i}$, the inclusion homomorphism
$\pi_1(\obd(B))\to\pi_1(
\obd(\overline{K_i}))$ is non-surjective.
\EndClaim

To prove \ref{one or two}, first note that $K_i\in
\calc(\calw_1-\calb)=\calk_1\discup\calk_2\discup\calk_3$. If
$K_i\in\calk_1$, the definition of the assignment $K\mapsto v^K$
implies that $v_i=v^{K_i}\in\calv_1^*$, which by definition means that
$v_i=v_\tS$ for some unbelted component $\tS$ of
$|\rho|^{-1}(\cals_1)$. In particular $v_i$ is ideal, and by the
definition of a special ideal vertex, $\tS$ must be pseudo-belted. But
according to the definition of $\cals_1$ (and Lemma \ref{it's-a unique}), a component of $\cals_1$ has
no pseudo-belted side. Hence we cannot have $K_i\in\calk_1$.

If $K_i\in\calk_2$, then
$K_i=|\rho|(L\setminus\calq)$  for some $L\in\calc(\calL)$. If
$K_i\in\calk_3$, then $K=|\rho|(R\setminus\calq)$ for some component $R$
of $\cald$.  Thus in the case $K_i\in\calk_2\cup\calk_3$, we have $K_i=|\rho|(P\setminus\calq)$ for some component $P$ of $\calp$. Let us set
 $\calq_P=P\cap\tcals$, so that $\calq_P$ is a union of components of
 $\calq$ and $K_i=|\rho|(P-\calq_P)$. It follows from \ref{still sick}
 that $|\rho|$ maps $P-\calq_P$ homeomorphically onto $K_i$, and maps
 each component of $\calq_P$ homeomorphically onto a component of
 $\calb$. Since by \ref{cockroach redux} each component
 of $\calb$ is separating in the component of
 $\calw_1$ containing it, 
%\redcomment{Double check that I have
  % that right} 
$|\rho|$ maps $P$ homeomorphically
 onto $\overline{K_i}$. Since $\rho$ is an orbifold
 immersion of $\oldPsi$ in $\oldOmega$, the map $\rho|\obd(P)$ is an
 orbifold homeomorphism from $\obd(P)$ to
 $\obd(\overline{K_i})$. Since $\obd(P)$ is a
 \torifold\ by \ref{they're torifolds},
 $\obd(\overline{K_i})$ is a \torifold, which is the
 first assertion of \ref{one or two} in this case. In addition to knowing that 
$\rho|\obd(P):
\obd(P)\to
 \obd(\overline{K_i})$ is an
 orbifold homeomorphism, we know by \ref{sicker still} that 
each component of
$\Fr_{W_1}\overline{K_i}\subset\obd(\calb)$
is the homeomorphic image under the restriction of $\rho$ of a component of $\calq_P$. Hence in order to prove
 the second assertion of \ref{one or two}, that for every component $B$ of $\Fr_{W_1}\overline{K_i}$ the inclusion homomorphism $\pi_1(\obd(B))\to\pi_1(\obd(\overline{K_i}))$ is non-surjective, it suffices to show that 
for every component $Q$ of $\calq_P$, the inclusion homomorphism $\pi_1(\obd(Q))\to\pi_1(\obd(P))$ is non-surjective.

In the subcase where $K_i\in\calk_2$, so that $L:=P\in\calc(\calL)$, the definition of the assignment $K\mapsto v^K$ implies that  $v^{K_i}=v_L$.
Since the material vertex $v_L=v^{K_i}$ is a special material vertex, by definition we have $w_L>1$. The definition of $w_L$ given in \ref{old trier} now implies that
for every component $Q$ of $\calq_L=L\cap Q$,
the image of the inclusion homomorphism 
%from Thus
 $\pi_1(\obd(Q))\to\pi_1(\obd(L))$ has index $w_L>1$ in
 $\pi_1(\obd(L))$; in particular, this inclusion homomorphism is
 not surjective.
 Now consider the subcase in which $K_i\in\calk_3$, so that $P\in\calc(\cald)$, and hence 
$P=|\frakD_\tS|$, where
$\tS$ is the pseudo-belted side of some component of $\cals_2^{!}$. According to \ref{replaces readin' john o'hara}, we have $
\calq_R=
|\frakD_\tS|\cap\tcals=|\frakC_\tS|$, a connected set.
Again by \ref{replaces readin' john o'hara}, the image of the inclusion homomorphism 
%from Thus
 $\pi_1(\frakC_\tS)\to\pi_1(\frakD_\tS)$ has index $2$ in
 $\pi_1(\frakD_\tS)$; in particular, this inclusion homomorphism is
 not surjective. This completes the proof of \ref{one or
   two}. 
%\redcomment{Check that this makes sense after it's fixed up.}
% is disjoint from $\fraks_\oldPsi$
  %(since $G_{\tS'}$ is), and

For $i=1,2$, it follows from \ref{cockroach redux} that every
component of $\Fr_{W_1}\overline{K_i}$ separates $W_1$. Hence, if
$\overline{K_1}\cap \overline{K_2}=\emptyset$, there is a component
$A$ of $W_1-(K_1\cup K_2)$ such that $A\cap \overline{K_i}$ is a
single component of $\Fr_{W_1}\overline{K_i}$ for $i=1,2$. In this
subcase we set $K_2'=K_2\cup A$, and denote by $B_1$ the component
$\overline{K_1}\cap A$ of $\Fr_{W_1}\overline{K_1}$.  If
$\overline{K_1}\cap \overline{K_2}\ne\emptyset$, then again because
each component of each $\Fr_{W_1}\overline{K_i}$ separates $W_1$, the
intersection of $\overline{K_1}$ and $\overline{K_2}$ is a single
component of $\Fr_{W_1}\overline{K_1}$; in this case we set $K_2'=K_2$
and $B_1=\overline{K_1}\cap \overline{K_2}$. Note that in each case,
$B_1$ is a component of $\Fr \overline{K_1}$, and $\overline{K_1}\cap
\overline{K_2'}=B_1$. We claim:

\Claim\label{ha ha}
The inclusion homomorphisms $\pi_1(\obd(B_1))\to \pi_1(\obd(\overline{K_1}))$ and $\pi_1(\obd(B_1))\to \pi_1(\obd(\overline{K_2'}))$ are non-surjective.
\EndClaim

To prove \ref{ha ha}, first note that the non-surjectivity of $\pi_1(\obd(B_1))\to \pi_1(\obd(\overline{K_1}))$ follows immediately from \ref{one or two}. If $\overline{K_1}\cap \overline{K_2}\ne\emptyset$, so that $K_2'=K_2$, the non-surjectivity of $\pi_1(\obd(B_1))\to \pi_1(\obd(\overline{K_2'}))$ also follows from \ref{one or two}. Now suppose that
$\overline{K_1}\cap \overline{K_2}=\emptyset$, so that $K_2'=K_2\cup
A$, and $\overline{K_2}\cap A$ is a single component $B_2$ of $\Fr
\overline{K_2}$. 
%Since $\obd(B_2)$ is $\pi_1$-injective in $\obd(W_1)$ by
It follows from \ref{zusammen} that
%{cross-ref, but maybe the statement should be that
  $\obd(\Fr\overline {K})$ is $\pi_1$-injective in $\obd(\overline{K})$ for every
  component $K$ of $\calw_1\setminus\calq$; 
%in any event I have added lots of $\obd$'s in this passage that I think were missing, and should check them---same thing in the next paragraph}, 
we may therefore use van Kampen's
theorem for orbifolds to identify
$\pi_1(\obd(B_2))$, $\pi_1(\obd(\overline{K_2}))$ and
$\pi_1(\obd((A))$ with subgroups of $\pi_1(\obd(\overline{K_2'}))$
(using a common base point in $B_2$) in such a way that
$\pi_1(\obd(\overline{K_2'}))$ is a free product with amalgamation
$\pi_1(\obd(\overline{K_2}))\star_{\pi_1(\obd(B_2))}\pi_1(\obd((A))$. In
particular, under this identification we have
$\pi_1(\obd(\overline{K_2}))\cap\pi_1(\obd((A))={\pi_1(\obd(B_2))}$. But
since $\pi_1(\obd(B_2))\to \pi_1(\obd(\overline{K_2}))$ is
non-surjective by \ref{one or two}, $\pi_1(\obd(B_2))$ is identified
with a proper subgroup of $\pi_1(\obd(\overline{K_2}))$, and hence
$\pi_1(\obd(A))$ is identified with a proper subgroup of
$\pi_1(\obd(\overline{K_2'}))$. This means that the inclusion
homomorphism $\pi_1(\obd(A))\to\pi_1(\obd(\overline{K_2'}))$ is
non-surjective; in particular, since $B_1\subset A$,  the inclusion homomorphism $\pi_1(\obd(B_1))\to\pi_1(\obd(\overline{K_2'}))$ is non-surjective, and the proof of \ref{ha ha} is complete.

Now note that since $\obd(B_1)$ is $\pi_1$-injective in $\obd(W_1)$ by \ref{zusammen}, we may use van Kampen's theorem for orbifolds to identify $\pi_1(\obd(B_1))$, $\pi_1(\obd(\overline{K_1}))$ and $\pi_1(\obd((\overline{K_2'}))$ with subgroups of $\pi_1(\obd(\overline{K_1}\cup \overline{K_2'}))$ (using a common base point in $B_1$) in such a way that $\pi_1
(\obd(\overline{K_1}\cup \overline{K_2'}))$ is a free product with amalgamation  $\pi_1(\obd(\overline{K_1}))\star_{\pi_1(\obd(B_1))}\pi_1(\obd((\overline{K_2'}))$. Now $\pi_1(\obd(B_1))\to \pi_1(\obd(\overline{K_1}))$ and $\pi_1(\obd(B_1))\to \pi_1(\obd(\overline{K_2'}))$ are non-surjective by \ref{one or two}. Furthermore, since $\obd(\overline{K_1})$ and $\pi_1(\obd(\overline{K_2}))$ are \torifold s by \ref{one or two}, $\pi_1(\obd(\overline{K_1}))$ and $\pi_1(\obd(\overline{K_2}))$ are infinite; since it follows from \ref{zusammen} that $\obd(\overline{K_2})$ is $\pi_1$-injective in $\obd(\overline{K_2'})$, the group $\pi_1(\obd(\overline{K_2'}))$ is also infinite. Thus $\pi_1
(\obd(\overline{K_1}\cup \overline{K_2'}))$ is a free product with amalgamation  of two infinite proper subgroups. But \ref{zusammen} also implies that $\obd(\overline{K_1}\cup \overline{K_2'})$ is $\pi_1$-injective in $\obd(W_1)$; and since the construction of $\obd(\calw_1)$ immediately implies that it is (orbifold-)homeomorphic to $\obd(\calw)$, it follows from \ref{don't die ET} that $\obd(W_1)$ is a \torifold\ and therefore has a virtually cyclic fundamental group. Hence $\pi_1
(\obd(\overline{K_1}\cup \overline{K_2'}))$ is virtually cyclic. But a virtually cyclic group cannot be expressed as a free product with amalgamation  of two infinite proper subgroups. This contradiction completes the proof of the second assertion of the proposition.
\EndProof
%\oldomecals\calb\cals_\calv\calp\calL\barl\calqQ*PPPPAAAK_1K_2K_iH_}' Z
%\cals_1\calg\eta\cale\calv (a)(i)\torifold KKKCCC\calk\etaTTTTTTTTT\caltYYYCCCXXX\oldGamma\rho\tau balanced(a)(i)(1)(A)\oldPi\frakP\oldPi$r \frakZ \tc R \frakrZ H E$R\frakB\frakC proofread\calb\calc\cald\frakD\calu\calr\frakB \oldXi

\Proposition\label{three or better}
Let $\Mh$ be a closed, orientable hyperbolic $3$-orbifold, and set $\oldOmega=(\Mh)\pl$ and
$M=|\oldOmega|$. Let $\cals$
be a \dandy, balanced system of spheres in $M$ having no \centralclashspheres. Let $v$ be a non-special
vertex of $\calf\oldomecals$ having valence at most $2$ in $\calf\oldomecals$. Then each edge incident to $v$ has the form $e_S$, where $S$ is a component of $\cals$ which is not \doublecentral. 
\EndProposition

\Proof
Set $\oldPsi=\oldOmega\cut{\obd(\cals)}$ and $\tcals=\partial(M\cut\cals)=\partial|\oldPsi|$.
Set $\calf=\calf\oldomecals$.

First consider the case in which $v$ is an ideal vertex. By definition
we have $v=v_\tS$ for some unbelted component  $\tS$ of $\tcals$. It
was pointed out in \ref{new trier} that the unique edge  incident to
$v$ is 
$e_S$, where $S=\rho(\tS)$. We are required to show that $e_S$ is not
\doublecentral. If we assume that  $e_S$ is 
\doublecentral, then in particular $\tS$ is \central. Since 
$\cals$ is a \dandy\ system of spheres in $M$, and since $\tS$ is not belted, it then follows from Lemma \ref{and
  the foos just keep on fooing} that either $\tS$ is
pseudo-belted or $\tS$ contains an \bad\ component of
$|\oldPhi(\oldPsi)|$. If $\tS$ contains an \bad\ component  of $|\oldPhi(\oldPsi)|$,
then $S$ satisfies Alternative (i) of Definition \ref{self-clash def}
and is therefore a  \centralclashsphere, a contradiction to the
hypothesis that $\cals$  has no
\centralclashspheres.  If  $\tS$ is
pseudo-belted, then by definition (see \ref{new trier}),
$v=v_\tS$ is a special ideal vertex, and we have a contradiction to
the hypothesis that $v$ is non-special.

Now consider the case in which $v$ is a material vertex. By definition
(see \ref{new trier})
we have $v=v_L$ for some component  $L$ of $\calL\oldomecals$. 
According to the definition of a special vertex given in \ref{new
  trier}, the hypothesis that $v$ is non-special means that $w_L=1$.
%and that $L\cap\fraks_\oldPsi=\emptyset$. \redmissingref{The condition
  %$L\cap\fraks_\oldPsi=\emptyset$ is the wrong one, and may not be needed
  %after $w_L$ is redefined properly.}

Let ${\tS_1},\ldots,{\tS_m}$ denote the components of $\tcals$ that
meet $L$, where $m\ge0$.
It follows from
Lemma \ref{new beginning} that each $\tS_i$ is belted, and that
$L\cap{\tS_i}=G_{\tS_i}$ for $1\le i\le m$. Hence
$G_{\tS_1},\ldots,G_{\tS_m}$ are the components of
$L\cap\tcals$. Since $L_{\tS_i}$ is by definition the component
of $|\oldSigma(\oldPsi)|$ containing $G_{\tS_i}$ (see \ref{stop beeping}),
and $L$ is a component of
$|\oldSigma(\oldPsi)|$, (see \ref{old trier}) we have $L=L_{\tS_i}$ for $i=1,\ldots,m$.

The definition of the incidence relation on $\calf$ implies
that the edges incident to $v$ are precisely those of the form
$e_{\rho(\tS)}$, where $\tS$ is a component of $\tcals$ such that
$L_\tS=L$. 
%\redcomment{This is wrong! Edges are associated with
  %components of $\cals$, not components of $\tcals$. Everything in
  %this passage needs fixing. I guess I should be using
  %$e_{S_1},\ldots,e_{S_m}$, where $S_i=\rho(\tS_i)$. These edges are
  %distinct because $\calf$ is a forest by the preceding prop. (and I
  %should make sure the hyps. hold).} 
Hence if we set $S_i=\rho(\tS_i)$ for $i=1,\ldots,m$, the edges $e_{S_1},\ldots,e_{S_m}$ are incident
to $v$. Conversely, if an edge $e_S$ is incident to $v$, then some side $\tS$ of $S$
we have
% of $S=\rho(\tS)$ for some 
$L_\tS=L$ and hence $\tS\cap L\ne\emptyset$; thus for
some $i$ we have $\tS=\tS_i$ and therefore $S=S_i$. This shows that $e_{S_1},\ldots,e_{S_m}$ are the only
edges incident
to $v$. These edges are pairwise distinct, since $\calf$ has no loops by Proposition \ref{new enchanted forest}. Hence $m$ is the valence of $v$.
The hypothesis then implies that $m\le2$. We have $m\ge1$
since, as observed in \ref{new trier}, $\calf\oldomecals$ has no isolated
vertices.

We claim:
\Claim\label{trrrivial}
The orbifold $\obd(L)$ is an \untwisted\ \pagelike\ \Ssuborbifold\ (see \ref{S-pair def}) of
%can be equipped with a trivial $I$-fibration over a $2$-orbifold in such a way that it is standardly embedded in
 $\oldPsi$.
\EndClaim

To prove \ref{trrrivial}, first note that by \ref{old trier}, 
%cross-ref using the annularity of the components of $\obd(L\cap\tcals)$---this is almost Proposition \ref{what does it say?}, but I need to know that $L\cap\tcals\ne\emptyset$}, 
$\obd(L)$ is a \torifold, and that by \ref{what's a belt?}, the annular orbifold
$\obd(G_{\tS_i})$ is $\pi_1$-injective for $1\le i\le m$. We have
seen that $m$ is either $1$ or $2$ and that $w_L=1$. The definition of
$w_L$ (see \ref{old trier}) gives $1=w_L=\min(n_1,n_2)$ if $m=2$ and
$1=n_1$ if $m=1$, where $n_i$ is the index of the image of the
inclusion homomorphism
$\pi_1(\obd(G_{\tS_i}))\to\pi_1(\obd(L))$. By symmetry we may
assume that $n_1=1$, i.e. that
$\pi_1(\obd(G_{\tS_1}))\to\pi_1(\obd(L))$ is an isomorphism. We now apply
Proposition \ref{one for the books}, letting
$\obd(L)$ and $\obd(G_{\tS_1})$
play the respective roles of $\oldUpsilon$ and $\frakZ$;
the hypotheses in Proposition \ref{one for the books} that $\oldUpsilon$ is very good and has no spherical boundary components hold here because $\obd(L)$ is a \torifold. Proposition \ref{one for the books} now asserts that the orbifold pair
$(\obd(L),\obd(G_{\tS_1}))$ is homeomorphic to
$(\oldGamma \times[0,1]),\oldGamma\times\{0\})$ for some
$2$-orbifold $\oldGamma$, which is homeomorphic to $\obd(G_{\tS_1})$
and therefore annular. 

If $m=1$ then $L\cap\partial\oldPsi=G_{\tS_1}$. Since $\obd(L)$ is a \torifold,
and $G_{\tS_1}$ is annular and $\pi_1$-injective, $\Fr_NL=\overline{\partial L
-G_{\tS_1}}$ is also annular. Thus the hypotheses of Lemma \ref{why the role}
hold with $A=\Fr_NL$, $\oldLambda=\obd(L)$, and $\oldUpsilon=G_{\tS_1}$. But
since $\pi_1(\obd(G_{\tS_1}))\to\pi_1(\obd(L))$ is an isomorphism,
we have a contradiction to the final assertion of Lemma \ref{why the role}.
Hence $m=2$.

Since
%We may identify
$(\obd(L),\obd(G_{\tS_1}))$ is homeomorphic to
%with
$(\oldGamma\times[0,1],\oldGamma\times\{0\})$,
% in such a way that
the orbifold $\frakK:=\obd( (\partial L)-\inter G_{\tS_1})$ is homeomorphic to $\oldGamma$ and is therefore annular.
% $G_{\tS_2}\subset|\oldGamma|\times\{1\}$. 
Since $\obd(G_{\tS_2})$ is an annular
suborbifold of $\inter\frakK$,
% the annular orbifold $\oldGamma\times\{1\}$, 
and is
$\pi_1$-injective by \ref{what's a belt?}, it follows from \ref{cobound} that either (a)
$\frakK$ is an orbifold regular neighborhood of
$\obd(G_{\tS_2})$ in $\obd(\partial L)$, or (b) $|\frakK|$ is a weight-$2$ disk and both points
of $\fraks_\frakK$ have order $2$, while $G_{\tS_2}$ is a weight-$0$ annulus.
% and
%$\obd(G_{\tS_2})$ has no singular points. 
Alternative (a) implies the
conclusion of \ref{trrrivial}. 

Now suppose that (b) holds. Since $\obd(G_{\tS_2})$ is
$\pi_1$-injective in $\frakK$, it follows from \ref{cobound} that the complement of the open annulus
$\inter G_{\tS_2}$ in the disk $|\frakK|$ must have two components,
an annulus $B$ of weight $0$ and a disk $D$ of weight
$2$. The disk $D$ is a component of
$|\frakA(\oldPsi)|$, and 
is therefore a properly embedded disk in $N$, transverse to $\fraks_\oldPsi$, with $\obd(D)$ 
annular and $\pi_1$-injective in $\oldPsi$. The third condition in the definition of a
  \dandy\ system (see \ref{semidandy def}) then implies that there is a
  disk $E\subset\tcals$ such that $\obd(E)$ is an annular orbifold
  and $\partial E=\partial D$.  Now the hypotheses of Lemma \ref{why the role}
hold with $A=D$, $\oldLambda=\obd(L)$, and $\oldUpsilon=\obd(E)$. The first
assertion of Lemma \ref{why the role} then implies that $D$ is the
full frontier of $L$ in $N$. But this is false, since the annulus $B$
is disjoint from $D$ and is contained in $\Fr_N L$. Thus Alternative
(b) cannot occur, and the proof of \ref{trrrivial} is
complete. 
%\redcomment{This is another argument that I believe is correct,
  %although I'm not 100\% comfortable with it. It may be a matter of
  %understanding Lemma \ref{why the role} better.}

%Now set $S_i=\rho(\tS_i)$ for $i=1,2$. \oldGamma
To establish the conclusion of
the proposition we must show that
neither of the components $S_1,S_2$
of $\cals$ is \doublecentral; by symmetry it suffices to prove that $S_1$
is not \doublecentral. Let us assume that $S_1$ is \doublecentral. Then in
particular $\tS_1$ is \central. 

Let us fix a component $C$ of $\partial G_{\tS_1}$, and let $V$ denote
the component of 
$\tS\setminus\inter|\oldPhi(\oldPsi)|$ 
%$\tS_1-\inter  G_{\tS_1}$ 
containing $C$. According
to  Lemma \ref{corpulent porpoises}, $\obd(V)$ is an annular orbifold. Thus
the hypotheses of  Proposition \ref{lady edith} hold with this choice
of $V$, and the component $F$ of $|\oldPhi(\oldPsi)|$ given by
Proposition \ref{lady edith} must be equal to $G_{\tS_1}$. Since $L=L_{\tS_1}$ is the component of
$|\oldSigma(\oldPsi)|$ containing $G_{\tS_1}$, the last sentence of
Proposition \ref{lady edith} asserts that $L$
%But $V$ shares the boundary component $C$ with $V$, which by \ref{trrrivial} 
is not a \pagelike\ \Ssuborbifold\ of
% can be equipped with a trivial $I$-fibration in such a way that it is standardly embedded in 
$\oldPsi$. But this contradicts  \ref{trrrivial}.
\EndProof
%{new beginning}

\Proposition\label{alt twice times lemma} 
 Let $\Mh$ be a closed,
orientable hyperbolic $3$-orbifold, and set $\oldOmega=(\Mh)\pl$ and $M=|\oldOmega|$.
Let $\cals$  be a non-empty, $\oldOmega$-\dandy, balanced system of spheres in $M$ having no \centralclashspheres. Set $\tcals=\partial M\cut\cals$, and let $n$ denote the number of  non-\central\ components of $\tcals$.
Then 
we have $\compnum(\cals)\le2n-1$. 
\EndProposition

\Proof 
Set $\tcals=\partial M\cut\cals$, and set
$\calf=\calf\oldomecals$. Set $\calv=\calv\oldomecals$ and
$\cale=\cale\oldomecals$. Given $v\in\calv$ and $e\in\cale$, we shall
write $e>v$ if $e$ is incident to $v$. Let us denote by $\cald$ the
set of all ordered pairs  $(v,e)\in\calv\times\cale$ such that $e>v$.

If $(v,e)\in\cald$ is given, the definition of $\calf$ implies that there is a unique component $S$ of $\cals$ such that $e=e_S$. Furthermore, if $v$ is a material vertex, 
there is a unique
component $L$  of $\calL\oldomecals$ such that $v=v_L$; we have $L=L_\tS$ for
some belted side $\tS$ of $S$. Since $\calf$ has no loops by Proposition \ref{new enchanted forest}, the side $\tS$ is uniquely determined by $(v,e)$. If $v$ is an ideal vertex, 
then the definition of $\calf$ implies that $v=v_{\tS}$ for a unique unbelted side  $\tS$ of $S$. In each case, the unique component $\tS$ of $\tcals$ so defined will be denoted by $\tS_{(v,e)}$. 

Let us define an integer-valued function $\phi$ on $\cald$ by setting $\phi(v,e)=1$ if $\tS_{(v,e)}$ is a \central\ component  of $\tcals$, and $\phi(v,e)=3$ otherwise.

If $v$ is a non-special vertex of $\calf$, then by Proposition \ref{three or better}, either $v$ has valence at least $3$, or $\phi(v,e)=3$ for every edge incident to $v$. Since $v$ is not isolated (see \ref{new trier}), it follows that
\Equation\label{three or better claim}
\sum_{e>v}\phi(v,e)\ge3 \text{ for every non-special vertex }v\text{ of }\calf.
\EndEquation
On the other hand, since $\calf$ has no isolated vertices (see \ref{new trier}), and since the values of $\phi$ are strictly positive integers by definition, we have 
\Equation\label{three or worser claim}
\sum_{e>v}\phi(v,e)\ge1 \text{ for every vertex }v\text{ of }\calf.
\EndEquation

%Since each edge of $\calf$ is incident to exactly two vertices, we have
%\Equation\label{label it}
%2\sum_{e\in\cale}\phi(e)=\sum _{v\in\vns} \sum_{e>v}\phi(e)+\sum _{v\in\vs} \sum_{e>v}\phi(e),
%\EndEquation
If $\vns$ and $\vs$ denote, respectively, the set of all non-special
vertices and of all special vertices of $\calf$, then (\ref{three or
  better claim}) and (\ref{three or worser claim}) assert that
$\sum_{e>v}\phi(v,e)$ is bounded below by $3$ when $v\in\vns$ and by
$1$ when $v\in\vs$; 
 hence
%where $\vns$ and $\vs$ denote, respectively, the set of all non-special vertices and of all special vertices of $\calf$. According to \ref{three or better claim}, we have $\sum_{e>v}\phi(e)\ge3$ for every $v\in\vns$. Since $\calf$ has no isolated vertices by Corollary \ref{three or better cor}, and since the values of $\phi$ are strictly positive integers by definition, we have $\sum_{e>v}\phi(e)\ge1$ for every $v\in\vs$. Thus (\ref{label it}) gives
\Equation\label{that too}
\sum_{(v,e)\in\cald}\phi(v,e)=\sum_{v\in\calv}\sum_{e>v}\phi(v,e)\ge3\card\vns+\card\vs.
\EndEquation

Let $\caldc$ denote the set of all elements $(v,e)$ of $\cald$ such that $\tS_{(v,e)}$ is \central, and set $\caldnc=\cald-\caldc$. The definition of $\phi$ implies that
$$\sum_{(v,e)\in\cald}\phi(v,e)=\card\caldc+3\card\caldnc,$$
which with (\ref{that too}) gives
\Equation\label{i'm weary}
\card\caldc+3\card\caldnc
\ge3\card\vns+\card\vs.
\EndEquation
We shall obtain a lower bound for the right hand side of (\ref{i'm weary}) in terms of edges of $\calf$. For this purpose, we consider an arbitrary component $T$ of $\calf$, and let $\calv^T\subset\calv$ and $\cale^T\subset\cale$ denote, respectively, the set of vertices and of edges of $T$. Since $T$ is a tree by Proposition \ref{new enchanted forest}, we have $\card\cale^T=\card\calv^T-1$.  On the other hand, 
%if $\delta^T$ denotes the number of special vertices of $T$, 
Proposition \ref{new enchanted forest} also asserts that 
%$\delta^T\le1$. Hence
$T$ contains at most one special vertex, i.e. $\card(\vs\cap\calv^T)\le1$. Hence
\Equation\label{if only}
\card\cale^T\le\card\calv^T-\card(\vs\cap\calv^T)=\card(\vns\cap\calv^T),
\EndEquation
and equality can hold in (\ref{if only}) only if $T$ contains a special vertex. Summing (\ref{if only}) over the components of $\calf$, we obtain
\Equation\label{summed if only}
\card\cale\le\card\vns,
\EndEquation
and equality can hold in (\ref{summed if only}) only if each component of $\calf$ contains a special vertex. The hypothesis $\cals\ne\emptyset$ implies that $\calf$ has at least one component; hence if equality holds in (\ref{summed if only}) then $\card\vs\ge1$. Thus the inequality $3\card\vns\ge3\card\cale^T$, which follows from (\ref{summed if only}), is strict unless $\card\vs>0$. Hence in any event we have
$3\card\vns+\card\vs>3\card\cale$, which with (\ref{i'm weary}) gives
\Equation\label{sick o' tryin'}
\card\caldc+3\card\caldnc>3\card\cale.
\EndEquation

Next note that we have $\card\caldc+\card\caldnc=\card\cald$, so that $\card\caldc+3\card\caldnc=\card\cald+2\card\caldnc$; and that since each edge of $\calf$ is incident to exactly two vertices, we have $\card\cale=(\card\cald)/2$. Hence (\ref{sick o' tryin'}) may be rewritten in the form $\card\cald+2\card\caldnc>3 (\card\cald)/2$, i.e. $\card\cald<4\card\caldnc$. Again using that
$\card\cale=\card\cald/2$, we deduce
$\card\cale<2\card\caldnc$, and therefore
\Equation\label{requires interpretation}
\card\cale\le2\card\caldnc-1.
\EndEquation

To interpret the left hand side of (\ref{requires interpretation}), note that by the definition of $\calf$, the assignment $S\to e_S$ is a bijection from $\calc(\cals)$ to $\cale$. Hence
\Equation\label{what's on the left}
\card\cale=\compnum(\cals).
\EndEquation
In order to estimate the right hand side of (\ref{requires interpretation}), we will show:
\Claim\label{another bloomin' injection}
The assignment $(v,e)\mapsto\tS_{(v,e)} $ is an injection from $\cald$ to $\calc(\tcals)$.
\EndClaim

To prove \ref{another bloomin' injection}, suppose that $(v_1,e_1)$ and $(v_2,e_2)$ are elements of $\cald$ such that $\tS_{(v_1,e_1)}=\tS_{(v_2,e_2)}=\tS$, say. If $i\in\{1,2\}$ is given, 
there is a unique component $S_i$ of $\cals$ such that $e_i=e_{S_i}$. Furthermore, according to the definition of $\tS_{(v_i,e_i)}$, one of the following alternatives holds for each $i\in\{1,2\}$:
\Alternatives
\item $v_i$ is a material vertex, 
so that $v_i=v_{L_i}$ for some
component $L_i$  of $\calL\oldomecals$; moreover, $\tS=\tS_{(v_i,e_i)}$ is a belted side of $S_i$, and $L_i=L_\tS$; or 
\item  $\tS=\tS_{(v_i,e_i)}$ is an unbelted side of $S_i$, and $v_i$ is the ideal vertex  $v_{\tS}$.
\EndAlternatives
In any event, for $i=1,2$, we have $S_i=\rho(\tS)$. Hence $S_1=S_2$, so that
$e_1=e_{S_1}=e_{S_2}=e_2$. It remains to show that $v_1=v_2$. If $\tS$
is belted, then (ii) obviously cannot hold for $i=1$ or for
$i=2$. Hence (i) must hold for $i=1,2$, so that $L_i=L_\tS$. Hence
$L_1=L_2$, so that $v_1=v_{L_1}=v_{L_2}=v_2$. If $\tS$ is unbelted, then (i) obviously
cannot hold for $i=1$ or for $i=2$. Hence (ii) must hold for $i=1,2$,
so that $v_i=v_\tS$. Again it follows that $v_1=v_2$, and the proof of
\ref{another bloomin' injection} is  complete. 

In view of the definition of $\caldnc$, the injection given by \ref{another bloomin' injection} restricts to an injection from $\caldnc$ to the set of  non-\central\ components of $\tcals$. Hence if, as in the statement of the proposition, we denote by $n$ denote the number of  non-\central\ components of $\tcals$, we have
\Equation\label{caramba}
\card\caldnc\le n.
\EndEquation
The conclusion of the proposition follows immediately from (\ref{requires interpretation}), (\ref{what's on the left}), and (\ref{caramba}).
\EndProof

\Corollary\label{brand new corollary} 
 Let $\Mh$ be a closed,
orientable hyperbolic $3$-orbifold, and set $\oldOmega=(\Mh)\pl$ and $M=|\oldOmega|$.
Let $\cals$  be a non-empty, $\oldOmega$-\dandy, balanced system of spheres in $M$ having no \centralclashspheres. 
Then 
we have 
$\compnum(\cals)\le2s'_\oldOmega(M-\cals)-1$ (see \ref{it's ok}).
\EndCorollary

\Proof
Set $\tcals=\partial M\cut\cals$, and let $n$ denote the number of components of $\tcals$ that are not
\central\ relative to $\cals$.  
According to Proposition \ref{alt twice times lemma} we have $\compnum(\cals)\le2n-1$. Since every \full\ component of $\tcals$ is
\central\ (see \ref{strict def}), $n$ is bounded above by the number
of components of $\tcals$ that are not 
\full\ relative to $\tcals$, which is in turn bounded above by $s'_\oldOmega(M-\cals)$ according to Proposition \ref{hugh manatee prop}. This proves the corollary.
\EndProof

\Corollary\label{just one sphere}
 Let $\Mh$ be a closed,
orientable hyperbolic $3$-orbifold, and set $\oldOmega=(\Mh)\pl$ and $M=|\oldOmega|$.
% homeomorphic to $\SSS^2\times \SSS^1$ or to a (possibly trivial) lens space. 
Let $\cals$  be an $\oldOmega$-\dandy\ system of spheres in $M$ which is
connected, i.e. consists of a single sphere $S=\cals$, and suppose
that $S$ is not a \centralclashsphere\ of $
\cals$. Then $S$ has at least one side which is not \central.
%no component of which is an
%$\oldOmega$-\clashsphere . 
%at most half
%at most half the\ components of $\cals$ are \doublefull.
%bounded above by the number of components of $\tcals$ that are not \full.
% strictly pless than $\compnum(\cals)/2$.
%55 are
%Furthermore, 
\EndCorollary

\Proof
Since $\cals$ has exactly one component, it is non-empty and
balanced, and thus satisfies the hypotheses of Proposition
\ref{alt twice times lemma}. Hence if
$n$ denotes the number of components of 
$\tcals:=\partial M\cut\cals$ 
 that are not
\central\ in $\cals$, we have
$\compnum(\cals)\le2n-1$. Since $\compnum(\cals)=1$, it follows that
$n\ge1$, which gives the conclusion.
\EndProof

\section{Orbifold volume and manifold homology}\label{manifold hom section}

\tinymissingref{\tiny 
Lemma ``merium olso'' is being removed. It can be found in
new-march.tex. So can the old version of the proof of Theorem
\ref{manifold homology bound}.
}

The main result of this section, Theorems \ref{manifold
  homology bound}, was stated in the Introduction as Theorem B.

\Lemma\label{4 goes to 6}
Let $\Mh$ be a hyperbolic
$3$-orbifold, and set $\oldOmega=(\Mh)\pl$,
% such that $\fraks_\oldOmega$ is a link, \redcomment{That doesn't look as if it belongs there, but I need to think about it!!!! I notice the expression for $\sigma_\oldOmega$ in the second sentence of the proof is wrong (but would be correct without the floor). Study the statement and proof here. Should $s_\oldOmega$ perhaps be $s_\oldOmega'$?? Whatever is going on, this is a serious issue, and I don't even know whether the current statement of the lemma is correct} 
Let $\cals$ be
an $\oldOmega$-admissible system of spheres in $M:=|\oldOmega|$. Suppose
that $\cals$ has a component $S$ such that (i) $\wt_\oldOmega(S)>4$, and (ii)
$S$ has at least one side which is not \central\ relative to $\cals$. Then
$s_\oldOmega(M-\cals)\ge3\lambda_{\oldOmega}$ (see \ref{it's ok}).
\EndLemma
%\toldTheta

\Proof
Set $\tcals=
\partial M\cut\cals$, $\oldPsi=\oldOmega\cut{\obd(\cals)}$ and ${\toldTheta}=\partial\oldPsi$, so that
$\tcals=|{\toldTheta}|$. Set $\lambda=\lambda_\oldOmega$, so that by \ref{lambda thing} we have $\lambda_\oldPsi=\lambda$. According to \ref{t-defs}, 
we have $s_\oldOmega(M-\cals)=\sigma(\oldPsi)=
12
\chibar(\kish(\oldPsi))
$. Hence in order to show that $s_\oldOmega(M-\cals)\ge3\lambda$, it
suffices to show that $\chibar(\kish(\oldPsi))\ge\lambda/4$. Since by
\ref{tuesa day} we have
$\chibar(\kish(\oldPsi))=\chibar(\kish(\oldPsi)\cap{\toldTheta})/2$, it is enough
to show that $\chibar(\kish(\oldPsi)\cap{\toldTheta})\ge\lambda/2$. 
%Since the
%hypothesis implies that $\cals\ne\emptyset$, 
It follows from
Corollary \ref{after abu gnu}
that every component of
$\kish(\oldPsi)\cap{\toldTheta}$ has negative Euler characteristic. Hence it is in fact enough to prove:

\Claim\label{shoulda woulda}
For some union of components $\frakK$ of
$\kish(\oldPsi)\cap{\toldTheta}$ we have $\chibar(\frakK)\ge\lambda/2$.
\EndClaim

{\bf Case I: Every component of $|\kish(\oldPsi)|\cap\tS$ has  weight divisible by $\lambda$.} In this case, we use the hypothesis that $\cals$ has a component $S$ such that (i) $\wt_\oldOmega(S)>4$, and (ii)
$S$ has at least one side, say $\tS$, which is not \central\ in
$\tcals$; by definition (see \ref{strict
  def}) this means
there is a 
component $\kappa$ of $\tS\cap|\kish(\oldPsi)|$ which semi-splinters $\tS$, i.e. such that for every component
$\Delta$ of
$\tS-\inter\kappa$, we have $\wt \Delta\le\wt (\tS)/2$. Since
we are in Case I, $\wt\kappa$ is divisible by $\lambda$. We will complete the proof in
this case by showing that $\chibar(\obd(\kappa))\ge\lambda/2$, thus verifying
\ref{shoulda woulda} with $\frakK=\obd(\kappa)$.

Consider first the subcase in which $\kappa$ is a disk. Then
$\tS-\inter\kappa$ is connected, and therefore
$\wt (\tS-\inter\kappa)\le\wt (\tS)/2$. Since
$\wt _\oldPsi(\tS)=\wt_\oldOmega(S)>4$, it follows that
$\wt (\kappa)>2$. Since $\wt (\kappa)$ is divisible by $\lambda\in\{1,2\}$, it follows that $\wt (\kappa)\ge2+\lambda$. If $p_1,\ldots,p_m$ are the orders of
the distinct points of $\fraks_{\obd(\kappa)}$, we have $m=\wt\kappa\ge2+\lambda$ and
$p_i\ge2$ for $i=1,\ldots,m$. Hence
$\chi(\obd(\kappa))=1-\sum_{i=1}^m(1-1/p_i)\le1-m/2\le-\lambda/2$, so that
$\chibar(\obd(\kappa))\ge\lambda/2$. 

Next consider the subcase in which $\kappa$ is an annulus. By
 Corollary \ref{after abu gnu} 
we have $\chi((\obd(\kappa))<0$. Thus
$\wt (\kappa)>0$. Since $\wt (\kappa)$ is divisible by $\lambda$, it follows that $\wt \kappa\ge\lambda$. If $p_1,\ldots,p_m$ are the orders of
the distinct points of $\fraks_{\obd(\kappa)}$, we have $m=\wt\kappa\ge\lambda$ and
$p_i\ge2$ for $i=1,\ldots,m$. Hence
$\chi(\obd(\kappa))=-\sum_{i=1}^m(1-1/p_i)\le-m/2\le-\lambda/2$, so that
$\chibar(\obd(\kappa))\ge\lambda/2$.

If $\kappa$ is neither a disk nor an annulus, it is a planar surface
with at least three boundary components. In this subcase we have
$\chi(\obd(\kappa))\le \chi(\kappa)\le-1$, so that
$\chibar(\obd(\kappa))\ge1\ge\lambda/2$. This completes the proof in Case I.

{\bf Case II: Some component of $|\kish(\oldPsi)|\cap\tS$ has  weight not divisible by $\lambda$.} We must have $\lambda=2$ in this case. Thus $\lambda_\oldPsi=2$, so that 
%For the argument in this case, note that since
%$\fraks_\oldOmega$ is a link,  
$|\kish(\oldPsi)|\cap\fraks_\oldPsi$ is a $1$-manifold properly embedded in $|\oldPsi|$. We have $|\kish(\oldPsi)|\cap\tcals \cap\fraks_\oldPsi=\partial 
(|\kish(\oldPsi)|\cap\fraks_\oldPsi)
$, and hence $ \wt_\oldPsi(|\kish(\oldPsi)|\cap\tcals)=\card
(|\kish(\oldPsi)|\cap\tcals \cap\fraks_\oldPsi)$
is even. Hence if $q$ denotes the number of components of $
  |\kish(\oldPsi)|\cap\tcals$ that have odd weight, then $q$ is even. Since we are in 
  Case II, and since $\lambda=2$, the set 
$
  |\kish(\oldPsi)|\cap\tS$ has an odd-weight component, which is in
  particular an odd-weight component of
$
  |\kish(\oldPsi)|\cap\tcals$; it follows that
we have $q>0$, and hence $q\ge2$. We may therefore fix
distinct  components $\kappa_1$ and $\kappa_2$ of $
  |\kish(\oldPsi)|\cap\tcals$ such that
  $m_j:=\wt \kappa_j$ is odd for $i=1,2$. It now suffices
  to show that $\chibar(\obd(\kappa_j)\ge1/2$ for $i=1,2$, as this
  will verify \ref{shoulda woulda} with
  $\frakK=\obd(\kappa_1\cup\kappa_2)$ (and $\lambda=2$). 

Let
  $j\in\{1,2\}$ be given. If $\kappa_j$ is a disk, then since $\chi(\obd(\kappa_j))<0$ by Corollary \ref{after abu gnu}, we have $m_j>1$; since $m_j$ is odd, it follows that
 $m_j\ge3$. If $p_{1,j},\ldots,p_{m_j,j}$ are the orders of
the distinct points of $\fraks_{\obd(\kappa_j)}$, we have $p_{i,j}\ge2$
for $i=1,\ldots,m_j$, so that
$\chi(\obd(\kappa_j))=1-\sum_{i=1}^{m_j}(1-1/p_{i,j})\le1- m_j/2$. Since $m_j\ge3$ it follows
that 
$\chi(\obd(\kappa_j))\le-1/2$, i.e. $\chibar(\obd(\kappa_j))\ge1/2$. 

If $\kappa_j$ is not a disk, we note that since $m_j$ is odd we have
$m_j\ge1$; if $p_{1,j},\ldots,p_{m_j,j}$ are the orders of
the distinct points of $\fraks_{\obd(\kappa_j)}$, 
%we have $p_{i,j}\ge2$
%for $i=1,\ldots,m_j$. 
then since $\kappa_j$ is not a disk we have
$\chi(\obd(\kappa_j))\le-\sum_{i=1}^{m_j}(1-1/p_{i,j})$. Since $m_j\ge1$, and
since $p_{i,j}\ge2$ for $i=1,\ldots,m_j$, it follows
that
$\chi(\obd(\kappa_j))\le-1/2$,
i.e. $\chibar(\obd(\kappa_j))\ge1/2$. 
\EndProof
%kappa)\wt_\oldPsi\tS\kappa\tcals

%\redcomment{Point out
%Fix $\voct$ stuff from this point. And try to remember
 %why the restriction $\vol\oldOmega<1.83$ is so damn
 %important for the following proof.
%, and
%. OK, I'm starting to get that. One small point is that I should
 %point out where it's important that $1.83<3.44$.
%}

\Theorem\label{manifold homology bound} Let $\Mh$ be a closed,
orientable hyperbolic $3$-orbifold, and set $\oldOmega=(\Mh)\pl$. Suppose that $\vol\Mh<0.915\lambda_\oldOmega$. 
%  $\voct/2=1.83\ldots$ is the best I can do. Fix the proof and the
  %application, where $1.72$ will be replaced by $\voct/4=0.9159\ldots$.} 
Suppose
 that  $\oldOmega$ contains no  embedded negative turnovers. Then
$$h(|\Mh|)
\le\max\bigg(3,\bigg\lfloor\frac {\vol(\Mh)}{0.305}\bigg\rfloor+(3-\lambda_\oldOmega)
\max\bigg( (3\lambda_\oldOmega+1)
\bigg\lfloor
\frac{\vol(\Mh)}{0.305\lambda_\oldOmega}\bigg\rfloor
,\bigg\lfloor\frac{\vol(\Mh)}{0.305\lambda_\oldOmega}\bigg\rfloor+2\bigg)\bigg).
$$
%(In particular, if $\lambda_\oldOmega=2$, i.e. if $\fraks_{\Mh}$ is a link and $\vol\Mh<1.83$, we have $\dim
%H_1(|\Mh|;\FF_2 )=h(|\Mh|)\le40$; and if $\lambda_\oldOmega=1$, i.e. if $\fraks_{\Mh}$ is not a link and $\vol\Mh<0.915$, we have $\dim
%H_1(|\Mh|;\FF_2 ) =h(|\Mh|)\le18$.)
 \EndTheorem

The results of Theorem \ref{manifold homology bound} are expressed in the following table (where $\Mh$ is a 
closed, orientable, hyperbolic 3-orbifold such that  $(\Mh)\pl$
contains no embedded negative turnovers):

\begin{center}
 \begin{tabular}{||c | c | c||} 
 \hline
If $\lambda_\oldOmega=$ & and $\vol(\Mh)<$ & then $h(|\Mh|)\le$ \\ [0.5ex] 
 \hline\hline
1 & 0.305  & 4  \\ 
 \hline
 1 & 0.61 & 9  \\
 \hline
 1 & 0.915 & 18  \\
 \hline
 2 & 0.305 & 3  \\
 \hline
 2 & 0.61 & 8  \\ [1ex] 
 \hline
 2 & 0.915 & 16  \\ [1ex] 
 \hline
 2 & 1.22 & 24  \\ [1ex] 
 \hline
 2 & 1.525 & 32  \\ [1ex] 
 \hline
 2 & 1.83 & 40 \\ [1ex] 
 \hline

\end{tabular}
\end{center}

\Proof[Proof of Theorem \ref{manifold homology bound}]
Set $M=|\oldOmega|$. Set $V=\vol(\Mh)$, $\lambda=\lambda_\oldOmega$, and $\alpha=\lambda\lfloor V /(0.305\lambda)\rfloor$. 
Note that since  $V< 0.915\lambda$ by hypothesis, we have 
%$\alpha\le V /0.305<3\lambda$. Since $\alpha$ is an integer, and is even if $\lambda=2$, we have
\Equation\label{spillane}
\alpha\le2\lambda.
\EndEquation

According to Proposition \ref{semidandies exist} \tinymissingref{\tiny note
  corrected cross-ref here}, there exists an $\oldOmega$-\dandy\ system of spheres $\cals\subset
M$ which is a complete system of spheres  in $M$. 

If $\cals=\emptyset$ then the completeness of $\cals$ implies that $M$
is irreducible. It then follows from Proposition \ref{lost corollary}
that if $\lambda=2$ then $h(M)\le5$, and that if
  $ V <0.915$ (and in particular if $\lambda=1$), then $h(M)\le3$. This immediately implies the conclusion of the theorem.

For the rest of
  the proof we will assume that $\cals\ne\emptyset$. According to the definition of  a \dandy\ system of spheres,
$\cals$ is in particular admissible; hence by  \ref{doublemint}, 
 $\oldPsi:=\oldOmega\cut\cals$ is componentwise strongly \simple\ and componentwise boundary-irreducible. According to  \ref{lambda thing}, we have
$\lambda_\oldPsi=\lambda$.

We claim:
\Claim\label{less nor 6}
For every subsystem $\calt$ of $\cals$, we have 
$$s'_\oldOmega(M-\calt)\le\alpha.$$
\EndClaim

To prove \ref{less nor 6}, 
note that by the definitions of 
$s'_\oldOmega(M-\calt)$ and $\sigma'(
\oldOmega\cut{\obd(\calt)}
%\obd(M\cut\calt)
)$ (see \ref{t-defs}), we have 
\Equation\label{justino}
s'_\oldOmega(M-\calt)=\sigma'(
\oldOmega\cut{\obd(\calt)}
)=\lambda\lfloor\sigma(
\oldOmega\cut{\obd(\calt)}
)/\lambda\rfloor.
\EndEquation
By Corollary \ref{bloody hell}
we have
$\smock(\oldOmega)\ge\delta(\oldOmega)$. On the other hand, the definition of the invariant $\delta$ (see \ref{t-defs}) gives
$\delta(\oldOmega)=\sup_\oldPi\sigma(\oldOmega\cut\oldPi)$, where  $\oldPi$ ranges
over all (possibly empty) incompressible, closed $2$-suborbifolds of $\inter\oldOmega$. In particular, $\zeta(\oldOmega)\ge\sigma(\oldOmega\cut{\obd(\calt)})
%=\sigma(\oldPsi)
$,
 %\redcomment{Fix this. I have now defined $\oldPsi$ to be $\oldOmega\cut\cals$. I should probably go back to defining it to be $\oldOmega\cut\calt$.} 
%It follows that $\sigma(
%\oldOmega\cut{\obd(\calt)})\le\smock(\oldOmega)$, 
which with (\ref{justino}) gives
% and by (\ref{kitsch}) we have
%$\delta(\oldPsi)\ge\sigma(\oldPsi)$. Hence 
$s'_\oldOmega(M-\calt)\le
\lambda\lfloor\smock(\oldOmega)/\lambda\rfloor$. In view of the
definition of $\smock(\oldPsi)$ (see \ref{t-defs}), this means that
$s'_\oldOmega(M-\calt)\le
\lambda\lfloor\smock_0(\oldOmega)/(0.305\lambda)\rfloor$. 
 Since Corollary \ref{smockollary} asserts
that $\smock_0(\oldOmega)=\vol\Mh=V$, the inequality
$s'_\oldOmega(M-\calt)\le
\lambda\lfloor V /(0.305\lambda)\rfloor=\alpha$
now follows. Thus \ref{less nor 6} is proved.

According to  Corollary \ref{toss 'em out}, there is a subsystem
$\cals^*$ of $\cals$ such that no \doublefull\ component of $\cals^*$ is an $\oldOmega$-\clashsphere\ for $\cals^*$, and
$s'_\oldOmega(M-\cals^*)\ge
s'_\oldOmega(M-\cals)+\lambda\compnum(\cals-\cals^*)$. In particular we have
%sphere in the system $\cals^*$, and (ii)
\Equation
\label{barnery your googles}
\compnum(\cals-\cals^*)\le
s'_\oldOmega(M-\cals^*)/
%2
\lambda.
\EndEquation
%\redcomment{Am I throwing away too much info there? At this point (9/1/18) I think not.}
According to Proposition
  \ref{dandy subsystem} (applied with $\cals^*$ playing the role of $\cals_1$ in that proposition), $\cals^*$ is \dandy.

Now note that by Corollary \ref{hugh manatee cor},  the number of non-\doublefull\ 
components of $\cals^*$ is bounded above by 
$s'_\oldOmega(M-\cals^*)$. Hence \ref{less nor 6}, applied to $\calt=\cals^*$, implies:

\Claim\label{great-grandfather 1}
The number of components of $\cals^*$ that are not \doublefull\ (relative to $\cals^*$) is bounded above by $\alpha$.
\EndClaim

Next, we claim:

\Claim\label{great-grandfather 2}
Either $\compnum(\cals^*)\le 1$, or  $\cals^*$ has no \doublefull\ component of weight strictly greater than $4$.
\EndClaim

To prove \ref{great-grandfather 2}, assume that  $\compnum(\cals^*)\ge2$ and that some \doublefull\ component  $S$ of $\cals^*$ has weight strictly greater than $4$. %If $S$ is a maximum-weight component of $\cals^*_1$ then $\wt S>4$. 
Since we chose $\cals^*$ in such a way that no \doublefull\ component of $\cals^*$ is a \clashsphere\  for $\cals^*$, the component $S$ 
%is a \doublefull\ component 
of $\cals^*$ is not a \clashsphere\ for $\cals^*$.
Since $\cals^*$ has 
%no
%\clashspheres\ and has
 at least two components, and since $S$ is \doublefull, it follows from Proposition \ref{no wonder 1}, applied with the \dandy\ system $\cals^*$ playing the role of $\cals$, that
 $S$ is not a \centralclashsphere\ of $\cals^*$. 

Now regard $S$ as a subsystem $\cals^*_1$ of $\cals^*$. Since  $S$ is not a \centralclashsphere\ of $\cals^*$, but is a \doublefull\ and hence \doublecentral\ component of $\cals^*$,
we may apply Lemma
\ref{no wonder 2} , with $\cals^*$ and $\cals_1^*$ playing the
respective roles of $\cals$ and $\cals_1$, to deduce that $S$ is not a
\centralclashsphere\  of $\cals_1^*$, i.e. that $\cals_1^*$ has no
\centralclashsphere. By Proposition \ref{dandy subsystem} (applied with $\cals^*$ and $\cals_1^*$ playing the respective roles of $\cals$ and $\cals_1$ in that proposition), $\cals_1^*$ is \dandy. 
%, and by \ref{dots trrivial} it has no \clashspheres. 
It now follows from Corollary \ref{just one sphere}, applied with $\cals^*_1$ playing the role of
$\cals$,
that $S$ has at least one side which is not \central\ relative to the system
$\cals^*_1$.
Since in addition we have
$\wt_\oldOmega(S)>4$, the hypotheses of 
Lemma \ref{4 goes to 6} hold with $\cals^*_1$ playing the role of
$\cals$. It therefore follows from Lemma \ref{4 goes to 6} that
$s_\oldOmega(M-\cals^*_1)\ge3\lambda$. But by
\ref{less nor 6} (applied with $\calt=\cals^*_1$) and (\ref{spillane}) we have 
$s_\oldOmega(M-\cals^*_1)\le\alpha\le2\lambda$. 
This contradiction establishes \ref{great-grandfather 2}.

Now let $\cals^{!}$ denote the union of all components of $\cals^*$
that are \doublefull\ relative to $\cals^*$. By definition, $\cals^!$ is a subsystem of $\cals^*$. 
Hence Proposition \ref{dandy subsystem} (applied with $\cals^*$ and $\cals^!$ playing the respective roles of $\cals$ and $\cals_1$ in that proposition) gives:
\Claim\label{easy stuff}
$\cals^!$ is \dandy.
\EndClaim

It follows from \ref{great-grandfather 1} and the definition of $\cals^!$ that
\Equation\label{dream}
\compnum(\cals^*-\cals^!)\le\alpha.
\EndEquation

According to \ref{great-grandfather 2}, either $\compnum(\cals^*)\le 1$, or
every component of $\cals^!$ has weight
at most $4$. Now since the $\oldOmega$-\dandy\ system $\cals^*$ is in particular admissible
by definition, and since by hypothesis  $\oldOmega$ contains no
 embedded negative turnovers, no component of  $\obd(\cals^*)$ is a negative turnover. Since $\obd(\cals^*)$ is negative by \ref{doublemint}, this means that no component of $\cals^*$ has weight $3$. If a 
component $T$ of $\cals^*$ had weight strictly less than $3$, then $\obd(T)$ would be spherical, which is impossible since $\obd(\cals)$ is incompressible (by admissibility) in the irreducible $3$-orbifold $\oldOmega$. Hence:
\Claim\label{when I say exactly I mean exactly}
Either $\compnum(\cals^*)\le 1$, or every component of $\cals^!$ has weight exactly
$4$. 
\EndClaim

If  every component of $\cals^!$ has weight exactly
$4$, then in particular, $\cals^!$ is balanced. If $\compnum(\cals^*)\le 1$, then in particular we have $\compnum(\cals^!)\le 1$, which trivially implies that $\cals^!$ is balanced. Thus we have shown: 

\Claim\label{uncle}
The system $\cals^!$ is balanced.
\EndClaim

We claim:
\Claim\label{cathedral} Either $\compnum(\cals^*)\le1$, or $\cals^!$ has no \centralclashsphere.
\EndClaim

To prove \ref{cathedral}, assume that $\compnum(\cals^*)\ge2$ and that  $\cals^!$ has  a \centralclashsphere\ $S$. 
The definition of $\cals^!$ implies that $S$ is a \doublefull\ component of $\cals^*$, so that it is in particular a \doublecentral\ component of $\cals^*$. 
By applying Lemma \ref{no wonder 2}, with $\cals^*$ and
$\cals^!$ playing the respective roles of $\cals$ and $\cals_1$,
we deduce that $S$ is a \centralclashsphere\ of
$\cals^*$. Since $\compnum(\cals^*)\ge2$, and since $S$ is \doublefull\ in $\cals^*$, we may therefore apply Proposition \ref{no wonder 1}, with the \dandy\ system $\cals^*$ playing the role of $\cals$, to deduce that $S$ is a \clashsphere\ for $\cals^*$. This is a contradiction, since $\cals^*$ has no \doublefull\ \clashspheres. Thus \ref{cathedral} is proved.

If we assume that $\compnum(\cals^*)\ge2$, then according to \ref{easy stuff}, \ref{uncle}, and \ref{cathedral}, 
%To summarize, 
$\cals^!$  is an $\oldOmega$-\dandy, balanced system of spheres in $M$, which 
% by \ref{when I say exactly I mean exactly}, 
has no \centralclashspheres.
Hence if 
in addition we assume $\cals^!\ne\emptyset$, 
we can apply Corollary \ref{brand new corollary}, with $\cals^!$ playing the role of $\cals$, to deduce that
$\compnum(\cals^!)\le2s'_\oldOmega(M-\cals^!)-1$. 
This proves:

\Claim\label{listen to the mocking bird}
Either $\compnum(\cals^*)\le1$, or $\cals^!=\emptyset$, or 
$\compnum(\cals^!)\le2s'_\oldOmega(M-\cals^!)-1$.
\EndClaim

In all cases, we now claim:
\Equation\label{all right, have it your way. You heard a Dixie melody}
\compnum(\cals)\le\max\bigg(
2s'_\oldOmega(M-\cals^!)+\frac1{\lambda}s'_\oldOmega(M-\cals^*)+\alpha-1
%2s_\oldOmega(M-\cals^!)+\frac12s_\oldOmega(M-\cals^*)+3
,\frac1\lambda s'_\oldOmega(M-\cals^*)+\max(\alpha,1)
\bigg).
\EndEquation

To prove (\ref{all right, have it your way. You heard a Dixie melody}), first note that
\Equation\label{dixie shmixie}
\compnum(\cals)= \compnum(\cals^!) +
\compnum(\cals^*-\cals^!)+
\compnum(\cals-\cals^*).
\EndEquation
In the case where the third alternative of \ref{listen to the mocking bird} holds, 
it follows from (\ref{dixie shmixie}),  (\ref{barnery your googles}), and (\ref{dream})  that 
$$\compnum(\cals)\le2s'_\oldOmega(M-\cals^!)+\frac1{\lambda}s'_\oldOmega(M-\cals^*)+\alpha-1,$$
which implies (\ref{all right, have it your way. You heard a Dixie melody}).
In the case where the second alternative of \ref{listen to the mocking bird} holds, 
it follows from (\ref{dixie shmixie}),  (\ref{barnery your googles}), and (\ref{dream})  that 
$$\compnum(\cals)\le\frac1{\lambda}s'_\oldOmega(M-\cals^*)+\alpha,$$
which implies (\ref{all right, have it your way. You heard a Dixie melody}).
In the case where the first alternative of \ref{listen to the mocking bird} holds, 
we have
$\compnum(\cals)= 
\compnum(\cals^*)+
\compnum(\cals-\cals^*)\le \compnum(\cals-\cals^*)+1$, which with (\ref{barnery your googles}) gives
$$\compnum(\cals)\le
\frac1\lambda s'_\oldOmega(M-\cals^*)+1.$$
This again implies (\ref{all right, have it your way. You heard a Dixie melody}), which is thereby proved in all cases.

From (\ref{all right, have it your way. You heard a Dixie melody}) and (\ref{less nor 6}), it follows that
\Equation\label{cranbagels}
\begin{aligned}\compnum(\cals)&\le\max\bigg(\bigg(3+\frac1\lambda\bigg)\alpha-1,\frac\alpha\lambda+\max(\alpha,1)\bigg)\\
&=\max\bigg(\bigg(3+\frac1\lambda\bigg)\alpha-1,\frac\alpha\lambda+1\bigg).
\end{aligned}
\EndEquation

If $X$ is any component of $M-\cals$, then $\hatX$ is $+$-irreducible
by the definition of completeness (\ref{complete def}). Since in
addition we have that $ V 
%=\voct/2
<1.83<3.44$ and that
$\cals\ne\emptyset$, it follows from
Corollary \ref{new get loster corollary} that
$h(X)\le t_\oldOmega(X)+2/\lambda$. Summing over
the components of $M-\cals$, we obtain
\Equation\label{shillelagh law}
\begin{aligned}h(M-\cals)&=\sum_{X\in \calc(M-\cals)}h(X)\le\sum_{X\in
\calc(M-\cals)}\bigg(t_\oldOmega(X)+\frac2\lambda\bigg)\\&=(3-\lambda)\compnum(M-\cals)+t_\oldOmega(M-\cals),\end{aligned}
\EndEquation
where in the last step we have rewritten $2/\lambda$ as $3-\lambda$, since $\lambda\in\{1,2\}$.
%so that
%\Equation\label{shillelagh law}
%h(M-\cals)-\frac2\lambda\compnum(M-\cals)\le t_\oldOmega(M-\cals).
%\EndEquation

%\redcomment{The coefficient $2/\lambda$ will make that formula trickier to apply in the next step. Decide how to deal with this.

On the other hand, Proposition \ref{dual homology} asserts that $h(M)\le h(M-\cals)+\beta$, where $\beta$ denotes the first betti number of the dual graph of $\cals$ in $M$. The dual graph has $\compnum(M-\cals)$ vertices and $\compnum(\cals)$ edges, and hence $\beta=1-\compnum(M-\cals)+\compnum(\cals)$. We therefore have \Equation\label{new mabel mabel} 
%\begin{aligned}
h(M)
%\le1+\compnum(\cals )+h((M-\cals )\to M)
%-\compnum(M-\cals )\\
\le1+\compnum(\cals )+h(M-\cals )
-\compnum(M-\cals ).
%\end{aligned}
\EndEquation
Combining (\ref{shillelagh law}) with (\ref{new mabel mabel}), we find that
\Equation\label{new getting there}
h(M)\le1+\compnum(\cals)+t_\oldOmega(M-\cals)
+(2-\lambda)\compnum(M-\cals).
%, which with
%(\ref{cultural appropriation day}) gives
%\Equation\label{getting there}
%h(M)\le\frac52t_\oldOmega(M-\cals^!)+t_\oldOmega(M-\cals)+4.
\EndEquation

Since $M$ is connected, the dual graph of $\cals$ in $M$ is also connected, and therefore has Euler characteristic at most $1$. Since $\compnum(M-\cals)$ and $\compnum(\cals)$ are respectively the number of vertices and the number of edges of the dual graph, it follows that
$\compnum(M-\cals)\le\compnum(\cals)+1$; hence (\ref{new getting there}) implies
\Equation\label{shmoller shmoaster}
h(M)
\le
t_\oldOmega(M-\cals)
+(3-\lambda)(1+\compnum(\cals)).
\EndEquation

Now Lemma \ref{frobisher}  asserts that
$\smock(\oldOmega)\ge\smock(\oldOmega\cut{\obd(\cals)})$, 
%\redcomment{Should this be $\oldPsi^*$??}
\tinymissingref{\tiny I am eliminating
  references to $\smock'$, and all things prime, here. If it turns out
  I need them after all, they are in new-march.tex.} which in view of the definition
of $t_\oldOmega(M-\cals)$ (see \ref{t-defs} and \ref{it's ok}) says that
$t_\oldOmega(M-\cals)\le \smock(\oldOmega)$. 
By definition (again see \ref{t-defs}) we have $\smock(\oldOmega)= 
%6
\smock_0(\oldOmega)/
0.305
%\voct
$, 
%\redcomment{That's wrong. There is no floor in the def. of $\smock$. Fix from here.}
%\redcomment{
%Fix the proof from this point. Some of the occurrences of $t_\oldOmega$ should be replaced by $t'_\oldOmega$. The expression for $\zeta$ in terms of floors is not consistent with what is now in \ref{t-defs}. This may be partially corrected by the use of $t'$, but it's not that simple because some of the inequalities really do seem to involve $t_\oldOmega$. It's not clear to me at this point that I will end up with a stronger statement of the theorem than the one that's there now, except that $\alpha$ will enter in and probably give better homology bounds for more restrictive volume bounds}
and
by Corollary \ref{smockollary} we have $
\smock_0(\oldOmega)=
V$.
Thus
\Equation\label{no star}
t_\oldOmega(M-\cals)
\le\smock(\oldOmega)
= 
%6}{\voct}
\smock_0(\oldOmega)/
{0.305}
%\rfloor
=
%2\lfloor
%6}{\voct}
 V /
{0.305}
%\rfloor.
\EndEquation
%The same argument, with $\cals$ replaced by $\cals^!$ and $\cals^*$ respectively, shows that
%\Equation\label{yes dagger}
%t_\oldOmega(M-\cals^!)
%\le
%2\bigg\lfloor\frac{
%6}{\voct}
 %V }
%{0.61}
%\bigg\rfloor
%\EndEquation
%and
%\Equation\label{yes star}
%t_\oldOmega(M-\cals^*)
%\le
%2\bigg\lfloor\frac{
%6}{\voct}
 %V }
%{0.61}
%\bigg\rfloor.
%\EndEquation

From (\ref{shmoller shmoaster}), (\ref{no star}), (\ref{cranbagels}), 
% (\ref{no you didn't}), (\ref{yes dagger}) and (\ref{yes star}), 
it follows that
%\Equation\label{the jinkles}
$$
\begin{aligned}
h(M)&\le\frac V{0.305}+(3-\lambda)(1+\compnum(\cals))\\
&\le\frac V{0.305}+(3-\lambda)
\max\bigg(\bigg (3+\frac1\lambda\bigg)\alpha,\frac\alpha\lambda+2\bigg).,
\end{aligned}
$$
According to the definition of $\alpha$, it follows that
$$h(|\Mh|)
\le\max\bigg(3,\frac {\vol(\Mh)}{0.305}+(3-\lambda_\oldOmega)
\max\bigg( (3\lambda_\oldOmega+1)
\bigg\lfloor
\frac{\vol(\Mh)}{0.305\lambda_\oldOmega}\bigg\rfloor
,\bigg\lfloor\frac{\vol(\Mh)}{0.305\lambda_\oldOmega}\bigg\rfloor+2\bigg)\bigg),
$$
which  implies the conclusion of the theorem since $h(\Mh)$ is an integer.
\EndProof

Note that in the proof of Theorem \ref{manifold homology bound}, the
hypothesis $\vol\oldOmega<0.915\lambda$ is needed to establish Claim \ref{less
  nor 6}, which is in turn needed, via the intermediate claims 
\ref{great-grandfather 2} and \ref{when I say exactly I mean exactly},
to prove Claim \ref{uncle}. Since \ref{uncle} is crucial for the
application of Corollary \ref{brand new corollary}, it is not clear
how one can get a result similar to Theorem \ref{manifold homology
  bound}---say in the case $\lambda_\oldOmega=2$---without the hypothesis $\vol\oldOmega<1.83$. Apart from this
essential step, the 
argument would go through under the weaker hypothesis
$\vol\oldOmega<3.44$, giving an only slightly weaker bound on $\dim H_1(|\oldOmega|,\FF_2 )$.

\section{Hyperbolic volume and orbifold homology: the general case}\label{main theorem section}

In this section we will prove the following two theorems, which were stated in the introduction as Theorems C and D:

\Theorem\label{orbimain}Let $\Mh$ be a 
closed, orientable, hyperbolic 3-orbifold such that $V:=\vol(\Mh)$ is strictly less than
$
%\voct/4=
0.915$. Suppose that $\fraks_{\Mh}$ is a link, and that $\pi_1(\Mh)$
contains no hyperbolic triangle group. Then 
$$\begin{aligned}
\dim H_1(\Mh;\FF_2 )
\le1 &+\max\bigg(3,\bigg\lfloor\frac {V}{0.1525}\bigg\rfloor+
\max\bigg( 7
\bigg\lfloor
\frac{V}{0.305}\bigg\rfloor
,\bigg\lfloor\frac{V}{0.305}\bigg\rfloor+2\bigg)\bigg)\\&+
\max\bigg(3,\bigg\lfloor\frac {V}{0.305}\bigg\rfloor+
\max\bigg( 7
\bigg\lfloor
\frac{V}{0.61}\bigg\rfloor
,\bigg\lfloor\frac{V}{0.61}\bigg\rfloor+2\bigg)\bigg).
\end{aligned}
$$
In particular, $\dim H_1(\Mh;\FF_2 )\le29$. 
\EndTheorem

\Theorem\label {other orbimain}Let $\Mh$ be a 
closed, orientable, hyperbolic 3-orbifold such that $V:=\vol(\Mh)$ is strictly less than
$0.22875$. Suppose that  $\pi_1(\Mh)$
contains no hyperbolic triangle group. Then
$$\dim H_1(\Mh;\FF_2 )\le7+4
  \max\bigg( 9
\bigg\lfloor
\frac{V}{0.07625}\bigg\rfloor
,3\bigg\lfloor\frac{V}{0.07625}\bigg\rfloor+4\bigg) .
  $$
%.305 4V
\EndTheorem

The content of Theorems \ref{orbimain} and \ref{other orbimain} are summarized in the following table (where $\Mh$ is a 
closed, orientable, hyperbolic 3-orbifold such that  $\pi_1(\Mh)$
contains no hyperbolic triangle group).

\begin{center}
 \begin{tabular}{||c | c | c||} 
 \hline
If $\lambda_{(\Mh)\pl}=$ & and $\vol(\Mh)<$ & then $\dim H_1(\Mh;\FF_2 )\le$ \\ [0.5ex] 
 \hline\hline
1 & 0.07625  & 23  \\ 
 \hline
 1 & 0.1525 & 43  \\
 \hline
 1 & 0.22875 & 79  \\
 \hline
 2 & 0.305 & 7  \\
 \hline
 2 & 0.4575 & 13  \\ [1ex] 
 \hline
 2 & 0.61 & 14  \\ [1ex] 
 \hline
 2 & 0.7625 & 28  \\ [1ex] 
 \hline
 2 & 0.915 & 29  \\ [1ex] 
 \hline

\end{tabular}
\end{center}

\Proof[Proof of Theorem \ref{orbimain}]
Set $\oldOmega=(\Mh)\pl$. Then $\fraks_\oldOmega$ is a link in $|\oldOmega|$. 
%Since by \ref{tsar}, $\Psi'_\oldOmega$ is contained in $\fraks_\oldOmega$ and is the underlying set of a trivalent graph, we have $\Psi'_\oldOmega=\emptyset$. 
It therefore follows from the second assertion of Proposition
\ref{boogie-woogie bugle boy} 
%\redcomment{It's a lemma now, and this will have to be rewritten. The sentence before this one should be useful in doing the proof of the new ``Proposition A''} 
that
$\Mh$ is covered with degree at most $2$ by an
orbifold $\tMh$ such that 
$$\dim H_1(\Mh;\FF_2 )\le1+ h(|{\tMh}|)+h(|{\Mh}|).$$
Since $V<0.915$, we have
$\vol(\tMh)=2V<1.83$.  Since $\Mh$ has a link as its singular set,
so does $\tMh$; that is, $\lambda_{\toldOmega}=\lambda_\oldOmega=2$. If $\oldOmega$ or $\toldOmega$ contained
an  embedded negative turnover, it would be $\pi_1$-injective by Corollary \ref{injective hamentash}, a contradiction to the hypothesis that  $\pi_1(\Mh)$
contains no hyperbolic triangle group. Hence neither $\Mh$ nor $\tMh$ contains
any  embedded negative turnover. As $\vol\Mh<0.915<1.83$ and $\vol\tMh<1.83$, we may apply
Theorem \ref{manifold homology bound} to $\tMh$ and to 
$\Mh$. This gives 
$$h
(|\tMh|)
%$$h(M)
%&\le\frac V{0.305}+(3-\lambda)(1+\compnum(\cals))\\
%&
\le\max\bigg(3,\bigg\lfloor\frac {2V}{0.305}\bigg\rfloor+
\max\bigg( 7
\bigg\lfloor
\frac{2V}{0.61}\bigg\rfloor
,\bigg\lfloor\frac{2V}{0.61}\bigg\rfloor+2\bigg)\bigg)
$$
and
$$h
(|\Mh|)
%$$h(M)
%&\le\frac V{0.305}+(3-\lambda)(1+\compnum(\cals))\\
%&
\le\max\bigg(3,\bigg\lfloor\frac {V}{0.305}\bigg\rfloor+
\max\bigg( 7
\bigg\lfloor
\frac{V}{0.61}\bigg\rfloor
,\bigg\lfloor\frac{V}{0.61}\bigg\rfloor+2\bigg)\bigg).
$$
Hence
$$\begin{aligned}
\dim H_1(\Mh;\FF_2 )&\le 1+ h(|{\tMh}|)+h(|{\Mh}|)\\
&\le1+\max\bigg(3,\bigg\lfloor\frac {2V}{0.305}\bigg\rfloor+
\max\bigg( 7
\bigg\lfloor
\frac{2V}{0.61}\bigg\rfloor
,\bigg\lfloor\frac{2V}{0.61}\bigg\rfloor+2\bigg)\bigg)\\&\qquad+
\max\bigg(3,\bigg\lfloor\frac {V}{0.305}\bigg\rfloor+
\max\bigg( 7
\bigg\lfloor
\frac{V}{0.61}\bigg\rfloor
,\bigg\lfloor\frac{V}{0.61}\bigg\rfloor+2\bigg)\bigg).
%1+\frac {3V}{0.305}+
%\max\bigg( 7
%\bigg\lfloor
%\frac{V}{0.305}\bigg\rfloor
%,\bigg\lfloor\frac{V}{0.305}\bigg\rfloor+2\bigg)
%+
%\le\frac {2V}{0.305}+
%\max\bigg( 7
%\bigg\lfloor
%\frac{V}{0.305}\bigg\rfloor
%,\bigg\lfloor\frac{V}{0.305}\bigg\rfloor+2\bigg).
\end{aligned}
$$

%&\le 1+
%\max\bigg(5,4+7\bigg\lfloor\frac{
%6}{\voct}
%\vol(\Mh)}
%{0.305}
%\bigg\rfloor\bigg)+
%\max\bigg(5,4+7\bigg\lfloor\frac{
%12}{\voct}
%\vol(\Mh)}
%{0.305}
%\bigg\rfloor\bigg)
%.
%\end{aligned}
%$$

%1+7\bigg(\bigg\lfloor\frac{12}{\voct}\vol(\Mh)\bigg\rfloor+\bigg\lfloor\frac{6}{\voct}\vol%(\Mh)\bigg\rfloor\bigg).

%7+2\bigg\lfloor\frac{20}{3}V\bigg\rfloor+2\bigg\lfloor\frac{12}{\voct}V\bigg\rfloor
%+\max\bigg(1,2\bigg\lfloor\frac{12}{\voct}V\bigg\rfloor+2\bigg\lfloor\frac{6}{\voct}V\bigg\rfloor\bigg).
%
%6+2\bigg\lfloor\frac{20}{3}V\bigg\rfloor+
%4\bigg\lfloor\frac{12}{\voct}V\bigg\rfloor
%+2\bigg\lfloor\frac{6}{\voct}V\bigg\rfloor.
The right hand side is strictly bounded above by $29$  since $V
<0.915$. 
%\redcomment{The logic is correct, but if some of the formulae
  %involving $\voct$ are changed slightly, it will need to be
  %reconsidered---possibly by strengthening the hypothesized volume
  %bound very slightly.}
%\oldOmega
\EndProof

\Proof[Proof of Theorem \ref{other orbimain}]
Set $\oldOmega=(\Mh)\pl$. It follows from Proposition \ref{my little sony} that either $\oldOmega$ has a covering $\toldOmega$ of degree at most $2$ such that
$\dim H_1(\oldOmega;\FF_2 )\le1+ h(|{\toldOmega}|)+h(|{\oldOmega}|)$, or there is a $((\ZZ/2\ZZ)\times(\ZZ/2\ZZ))$ cover $\toldOmega$ of $\oldOmega$ such that
%orbifold $\tMh$ such that
$\dim H_1(\oldOmega;\FF_2 )\le 3+h(|\oldOmega|)+4h(|\toldOmega_2|)$.
Thus in either case, $\Mh$ has a regular covering $\tMh$, with covering group isomorphic to $(\ZZ/2\ZZ )^d$ for some $d\le2$, such that 
\Equation\label{statistical thai}
\dim H_1(\Mh;\FF_2 )\le 3+h(|\Mh|)+4h(|\tMh|).
\EndEquation

Since  $\pi_1(\Mh)$
contains no hyperbolic triangle group, it follows from Corollary \ref{injective hamentash}
that neither $\oldOmega$ nor $\toldOmega$ contains
any embedded negative turnover. Since $\vol\Mh<0.22875$, and since the covering $\toldOmega$ of $\oldOmega$ have degree at most $4$, the volume $V$ of $\Mh$ and the volume $4V$ of $\tMh$ are both strictly bounded above by $4\times0.22875=0.915$.
We may therefore apply
Theorem \ref{manifold homology bound} to $\tMh$ and to 
$\Mh$. For $\tMh$ this gives 
\Equation\label{house of mouse}
\begin{aligned}
h
(|\tMh|)
%$$h(M)
%&\le\frac V{0.305}+(3-\lambda_{\toldOmega})(1+\compnum(\cals))\\
%&
&\le\max\bigg(3,\bigg\lfloor\frac {4V}{0.305}\bigg\rfloor+
(3-\lambda_{\toldOmega})\max\bigg( (3\lambda_{\toldOmega}+1)
\bigg\lfloor
\frac{4V}{0.305}\bigg\rfloor
,\bigg\lfloor\frac{4V}{0.305}\bigg\rfloor+2\bigg)\bigg).
%\\
%&=\max\bigg(3,
%\max\bigg( 8
%\bigg\lfloor
%\frac{4V}{0.305}\bigg\rfloor
%,2\bigg\lfloor\frac{4V}{0.305}\bigg\rfloor+2\bigg)\bigg)\lambda
\end{aligned}
\EndEquation
For $\lambda_{\toldOmega}=1$, (\ref{house of mouse}) becomes
  $$
  \begin{aligned}
    h(|\tMh|)&\le
  \max\bigg(3,\bigg\lfloor\frac {4V}{0.305}\bigg\rfloor+
2\max\bigg( 4
\bigg\lfloor
\frac{4V}{0.305}\bigg\rfloor
,\bigg\lfloor\frac{4V}{0.305}\bigg\rfloor+2\bigg)\bigg) \\
&=
\max\bigg( 9
\bigg\lfloor
\frac{4V}{0.305}\bigg\rfloor
,3\bigg\lfloor\frac{4V}{0.305}\bigg\rfloor+4\bigg),
\end{aligned}
$$
while for $\lambda_{\toldOmega
}=2$, (\ref{house of mouse}) becomes
$$
\begin{aligned}
    h(|\tMh|)&\le
  \max\bigg(3,\bigg\lfloor\frac {4V}{0.305}\bigg\rfloor+
\max\bigg( 7
\bigg\lfloor
\frac{4V}{0.305}\bigg\rfloor
,\bigg\lfloor\frac{4V}{0.305}\bigg\rfloor+2\bigg)\bigg)\\
&=\max\bigg(3,
\max\bigg( 8
\bigg\lfloor
\frac{4V}{0.305}\bigg\rfloor
,2\bigg\lfloor\frac{4V}{0.305}\bigg\rfloor+2\bigg)\bigg).
\end{aligned}
$$
Hence in all cases we have
\Equation\label{playing with lyre}
h(|\tMh|)\le
\max\bigg( 9
\bigg\lfloor
\frac{4V}{0.305}\bigg\rfloor
,3\bigg\lfloor\frac{4V}{0.305}\bigg\rfloor+4\bigg),
\EndEquation
Applying
Theorem \ref{manifold homology bound} to $\Mh$, and observing that $\lfloor V/0.305\rfloor\le\lfloor0.22875/0.305\rfloor=0$, we obtain
$$
\begin{aligned}
h
(|\Mh|)&
%$$h(M)
%&\le\frac V{0.305}+(3-\lambda_\oldOmega)(1+\compnum(\cals))\\
%&
\le\max\bigg(3,\bigg\lfloor\frac {V}{0.305}\bigg\rfloor+
(3-\lambda_\oldOmega)\max\bigg( (3\lambda_\oldOmega+1)
\bigg\lfloor
\frac{V}{0.305}\bigg\rfloor
,\bigg\lfloor\frac{V}{0.305}\bigg\rfloor+2\bigg)\bigg)\\
&=\max(3,2(3-\lambda_\oldOmega))
.
%&=\redcomment{what, and give the reason}\lambda
\end{aligned}
$$
Since $\lambda_\oldOmega\in\{1,2\}$, it follows that
\Equation\label{game of clones}
h(|\Mh|)\le4.
\EndEquation
%\redcomment{which is $3$ if $\lambda_\Mh=2$ and $4$ if $\lambda_\Mh=1$.}
% 2V 61
%where the last equality holds
Now by (\ref{statistical thai}), (\ref{playing with lyre}) and (\ref{game of clones}), we have
%it follows \redcomment{in the case $\lambda_{\tMh}=1$ that 
$$
\begin{aligned}
  \dim H_1(\tMh;\FF_2 )&\le
3+h(|\Mh|)+4h(|\tMh|)\\
&\le  7+4
  \max\bigg( 9
\bigg\lfloor
\frac{4V}{0.305}\bigg\rfloor
,3\bigg\lfloor\frac{4V}{0.305}\bigg\rfloor+4\bigg) .
\end{aligned}
$$
\EndProof

\bibliographystyle{plain}

\end{document}

%% file: newbookmacros.tex
\usepackage[dvips,text={6.5truein,9truein},left=1truein,top=1truein]{geometry}
\usepackage{amssymb,amsmath,amscd,enumerate,
}
\usepackage[pdftex]{graphicx}
\usepackage{url}

\usepackage[colorlinks,linkcolor=blue,citecolor=blue,pdfstartview=FitH]{hyperref}
\usepackage{color}
\usepackage{mathrsfs}

\newcommand\nonessentialproofreadingnote[1]{\empty}
\newcommand\tinymissingref[1]{\empty}
\newcommand\abstractcomment[1]{\empty}

\usepackage{amsthm}
\theoremstyle{definition}
\newtheorem{para}{}[section]

\newtheorem{remark}[para]{Remark}
\newtheorem{remarks}[para]{Remarks}
\newtheorem{notation}[para]{Notation}

\newtheorem{convention}[para]{Convention}
\newtheorem{definition}[para]{Definition}
\newtheorem{definitions}[para]{Definitions}
\newtheorem{definitionnotation}[para]{Definition and Notation}
\newtheorem{remarksnotation}[para]{Remarks and Notation}
\newtheorem{remarksnotationdefinitions}[para]{Remarks, Notation, and
  Definitions}
\newtheorem{remarksdefinitions}[para]{Remarks and Definitions}
\newtheorem{remarknotation}[para]{Remark and Notation}
\newtheorem{remarksdefinition}[para]{Remarks and Definition}
\newtheorem{notationremark}[para]{Notation and Remark}
\newtheorem{notationremarks}[para]{Notation and Remarks}
\newtheorem{definitionsnotation}[para]{Definitions and Notation}
\newtheorem{notationdefinitions}[para]{Notation and Definitions}
\newtheorem{reviewdefinition}[para]{Review and Definition}
\newtheorem{definitionremark}[para]{Definition and Remark}
\newtheorem{definitionremarks}[para]{Definition and Remarks}
\newtheorem{definitionsremarks}[para]{Definitions and Remarks}
\newcommand\Alternatives{\begin{enumerate}[(i)]}
\newcommand\EndAlternatives{\end{enumerate}}
\newcommand\Conditions{\begin{enumerate}[(1)]}
\newcommand\EndConditions{\end{enumerate}}

\theoremstyle{plain}
\newtheorem{theorem}[para]{Theorem}
\newtheorem{lemma}[para]{Lemma}
\newtheorem{proposition}[para]{Proposition}
\newtheorem*{propA}{Proposition A}

\newtheorem*{TheoremB}{Theorem B}
\newtheorem*{TheoremC}{Theorem C}
\newtheorem*{TheoremD}{Theorem D}

\newtheorem*{propE}{Proposition E}
\newtheorem*{propF}{Proposition F}
\newtheorem*{propG}{Proposition G}
\newtheorem{corollary}[para]{Corollary}
\newtheorem{conjecture}[para]{Conjecture}

\newtheorem{claim}[equation]{}
\numberwithin{equation}{para}
\numberwithin{figure}{section}
\numberwithin{section}{chapter}

\newcommand\Number{\begin{para}}
\newcommand\EndNumber{\end{para}}
\newcommand\Definition{\begin{definition}}
\newcommand\EndDefinition{\end{definition}}
\newcommand\Definitions{\begin{definitions}}
\newcommand\DefinitionsNotation{\begin{definitionsnotation}}
\newcommand\EndDefinitionsRemarks{\end{definitionsremarks}}
\newcommand\DefinitionsRemarks{\begin{definitionsremarks}}
\newcommand\RemarksNotation{\begin{remarksnotation}}
\newcommand\EndRemarksNotation{\end{remarksnotation}}
\newcommand\RemarkNotation{\begin{remarknotation}}
\newcommand\EndRemarkNotation{\end{remarknotation}}
\newcommand\NotationRemark{\begin{notationremark}}
\newcommand\NotationRemarks{\begin{notationremarks}}
\newcommand\EndNotationRemark{\end{notationremark}}
\newcommand\EndNotationRemarks{\end{notationremarks}}
\newcommand\EndDefinitionsNotation{\end{definitionsnotation}}
\newcommand\DefinitionNotation{\begin{definitionnotation}}
\newcommand\EndDefinitionNotation{\end{definitionnotation}}
\newcommand\ReviewDefinition{\begin{reviewdefinition}}
\newcommand\EndReviewDefinition{\end{reviewdefinition}}
\newcommand\EndDefinitions{\end{definitions}}
\newcommand\Theorem{\begin{theorem}}
\newcommand\EndTheorem{\end{theorem}}
\newcommand\Conjecture{\begin{conjecture}}
\newcommand\EndConjecture{\end{conjecture}}
\newcommand\Remark{\begin{remark}}
\newcommand\EndRemark{\end{remark}}
\newcommand\Remarks{\begin{remarks}}
\newcommand\EndRemarks{\end{remarks}}
\newcommand\Convention{\begin{convention}}
\newcommand\EndConvention{\end{convention}}
\newcommand\Notation{\begin{notation}}
\newcommand\EndNotation{\end{notation}}
\newcommand\Lemma{\begin{lemma}}
\newcommand\EndLemma{\end{lemma}}
\newcommand\Proposition{\begin{proposition}}
\newcommand\EndProposition{\end{proposition}}
\newcommand\Corollary{\begin{corollary}}
\newcommand\EndCorollary{\end{corollary}}
\newcommand\Claim{\begin{claim}}
\newcommand\EndClaim{\end{claim}}
\newcommand\Proof{\begin{proof}}
\newcommand\EndProof{\end{proof}}
\newcommand\Equation{\begin{equation}}
\newcommand\EndEquation{\end{equation}}
\newcommand\NoProof{{\hfill$\square$}}
\newcommand\Bullets{\begin{itemize}}
\newcommand\EndBullets{\end{itemize}}

\newcommand\pl{_{\rm PL}}
\newcommand\PL{^{\rm PL}}
\newcommand\smooth{_{\rm smooth}}
\newcommand\oldOmega{{\mathfrak M}} 
\newcommand\Mh{{\oldOmega_{\rm h}}} 
\newcommand\Nh{{N_{\rm h}}} 
\newcommand\ph{{p_{\rm h}}} 
\newcommand\tMh{{\toldOmega_{\rm h}}} 
\newcommand\Ends{{\rm Ends}}
\newcommand\bends{{\rm ends}}
\newcommand\oldomecals{{(\oldOmega,\cals)}}
\newcommand\oldGamma{{\mathfrak G}}
\newcommand\oldDelta{{\mathfrak D}}
\newcommand\oldTheta{{\mathfrak T}}
\newcommand\oldLambda{{\mathfrak L}}   
\newcommand\otheroldLambda{{\mathfrak O}}

\newcommand\oldXi{\mathfrak X} 
\newcommand\frakJsmooth{{\frakJ^{\rm smooth}}}

\newcommand\oldPi{\mathfrak P} 
\newcommand\oldSigma{\mathfrak S} 
\newcommand\oldUpsilon{\mathfrak W} 
\newcommand\oldPhi{\mathfrak F} 
\newcommand\oldPsi{\mathfrak N}
\newcommand\frakA{\mathfrak A} 
\newcommand\frakO{\mathfrak O} 
\newcommand\frakB{\mathfrak B}
\newcommand\frakD{\mathfrak D}
\newcommand\frakC{\mathfrak C}
\newcommand\frakE{\mathfrak E}
\newcommand\frakI{\mathfrak I}
\newcommand\frakJ{\mathfrak J}
\newcommand\frakK{\mathfrak K}

\newcommand\frakN{\mathfrak N}
\newcommand\frakZ{\mathfrak Z}
\newcommand\frakU{\mathfrak U}
\newcommand\tfrakU{{\widetilde{\mathfrak U}}}

\newcommand\tfrakK{{\widetilde{\mathfrak K}}}
\newcommand\tfrakV{{\widetilde{\mathfrak V}}}

\newcommand\tL{{\widetilde L}}

\newcommand\tfrakR{{\widetilde{\mathfrak R}}}

\newcommand\tfrakC{{\widetilde{\mathfrak C}}}
\newcommand\tfrakE{{\widetilde{\mathfrak E}}}

\newcommand\tfrakZ{{\widetilde{\mathfrak Z}}}
\newcommand\frakV{\mathfrak V}
\newcommand\frakP{\mathfrak P}
\newcommand\frakR{\mathfrak R}
\newcommand\tfrakB{\widetilde{\mathfrak B}}

\newcommand\toldOmega{\widetilde{\oldOmega}} 

\newcommand\toldGamma{\widetilde{\mathfrak G}}
\newcommand\tcalv{\widetilde{\calv}}
\newcommand\tv{\widetilde{v}}
\newcommand\tsigma{\widetilde{\sigma}}
\newcommand\ts{\widetilde{s}}
\newcommand\tg{\widetilde{g}}
\newcommand\toldDelta{\widetilde{\mathfrak D}}
\newcommand\toldTheta{\widetilde{\mathfrak T}}
\newcommand\toldLambda{\widetilde{\mathfrak L}}   
\newcommand\toldXi{\widetilde{\mathfrak X}}
\newcommand\toldPi{\widetilde{\mathfrak P}}

\newcommand\toldUpsilon{\widetilde{\mathfrak W}}

\newcommand\toldPsi{\widetilde{\oldPsi}}
\newcommand\Ccgeo{C_c^{\rm geo}}
\newcommand\Cgeo{C^{\rm geo}}
\newcommand\Ygeo{Y_{\rm geo}}
\newcommand\Ccprimegeo{C_{c'}^{\rm geo}}
\newcommand\tCgeo{{\tC}^{\rm geo}}
\newcommand\tYgeo{{\tY}_{\rm geo}}

\newcommand\compnum{{\frak c}}
\newcommand\newgamma{\gamma}
\newcommand\newalpha{\alpha}
\renewcommand\epsilon{\varepsilon}

\newcommand\Fann{{F}_{\rm ann}}
\newcommand\Fannin{{F}_{\rm ann}^{\rm in}}
\newcommand\Fannout{{F}_{\rm ann}^{\rm out}}
\newcommand\Dandy{Dandy}
\newcommand\spair{S-pair}
\newcommand\Ssuborbifold{S-suborbifold}
\newcommand\Asuborbifold{A-suborbifold}
\newcommand\pagelike{page-like}
\newcommand\bindinglike{binding-like}
\newcommand\untwisted{untwisted}
\newcommand\twisted{twisted}

\newcommand\eclipsing{eclipsing}

\newcommand\simple{atoral}
\newcommand\dandy{dandy}
\newcommand\bad{ungainly}

\newcommand\torifold{solid toric orbifold}

\newcommand\book{\mathop{\rm book}}

\newcommand\kish{\mathop{\rm kish}}

\newcommand\discup{\mathbin{\rotatebox[origin=c]{90}{$\vDash$}}}
\newcommand\FF{{\bf F}}
\newcommand\RR{{\bf R}}
\newcommand\SSS{{\bf S}}
\newcommand\TTT{{\bf T}}
\newcommand\PP{{\bf P}}
\newcommand\DD{{\bf D}}

\newcommand\OO{{\bf O}}
\newcommand\fraks{{\Sigma}}
\newcommand\frakg{{\mathfrak g}}

\newcommand\hateta{\widehat\eta}

\newcommand\hatC{{\widehat C}}
\newcommand\hatX{{\widehat X}}
\newcommand\hatY{{\widehat Y}}
\newcommand\plusX{{ X^+}}
\newcommand\plusV{{ V^+}}
\newcommand\plusN{{ N^+}}
\newcommand\plusw{{ W^+}}

\newcommand\chibar{\overline\chi}

\newcommand\inter{\mathop{\rm int}}
\newcommand\size{\mathop{\rm size}}

\newcommand\thh{\widetilde h}
\newcommand\tH{\widetilde H}
\newcommand\tM{\widetilde M}
\newcommand\tC{\widetilde C}
\newcommand\tcalq{\widetilde{\calq}}

\newcommand\tX{\widetilde X}

\newcommand\tZ{\widetilde Z}
\newcommand\tY{\widetilde Y}

\newcommand\tT{{\widetilde T}}
\newcommand\tU{{\widetilde U}}
\newcommand\tV{{\widetilde V}}
\newcommand\tW{{\widetilde W}}
\newcommand\tR{{\widetilde R}}
\newcommand\tQ{\widetilde Q}
\newcommand\tq{\widetilde q}

\newcommand\tDelta{\widetilde \Delta}

\newcommand\tphi{{\widetilde \phi}}

\newcommand\tcals{\widetilde \cals}

\newcommand\barp{\Pi}
\newcommand\overp{\overline{p}}

\newcommand\barq{\Lambda}

\newcommand\barrho{\overline{\rho}}

\newcommand\tS{{\widetilde S}}
\newcommand\ta{{\widetilde a}}

\newcommand\tA{{\widetilde A}}

\newcommand\txi{{\widetilde\xi}}
\newcommand\tB{{\widetilde B}}
\newcommand\tD{{\widetilde D}}
\newcommand\tE{{\widetilde E}}

\newcommand\tG{{\widetilde G}}
\newcommand\tGee{{\widetilde G}}

\newcommand\tj{{\widetilde j}}

\newcommand\tx{{\widetilde x}}
\newcommand\tu{{\widetilde u}}

\newcommand\calc{{\mathcal C}}

\newcommand\cala{{\mathscr A}}
\newcommand\scrt{{\mathscr T}}

\newcommand\calp{{\mathcal P}}

\newcommand\ch{\calc_{\rm h}}
\newcommand\cald{{\mathcal D}}
\newcommand\frakH{{\mathfrak H}}
\newcommand\frakf{{\mathfrak f}}
\newcommand\frakh{{\mathfrak h}}
\newcommand\frakX{{\mathfrak X}}

\newcommand\oldfrakX{{\Xi}}

\newcommand\calf{{\mathcal F}}

\newcommand\calL{{\mathscr L}}
\newcommand\calK{{\mathcal K}}
\newcommand\calm{{\Phi}}
\newcommand\newf{{\mathfrak f}}
\newcommand\emm{{f}}

\newcommand\otherdge{{\mathfrak j}}
\newcommand\noodge{j}
\newcommand\frakQ{{\mathfrak Q}}

\newcommand\dom{\mathop{\rm dom}}
\newcommand\range{\mathop{\rm range}}
\newcommand\barcalm{{\overline{\calm}}}

\newcommand\calq{{\mathcal Q}}
\newcommand\calt{{\mathcal T}}
\newcommand\calb{{\mathcal B}}

\newcommand\partialh{\partial_{\rm h}}
\newcommand\partialv{\partial_{\rm v}}
\newcommand\calg{{\mathcal G}}
\newcommand\cale{{\mathcal E}}

\newcommand\calx{{\mathcal X}}
\newcommand\caly{{\mathcal Y}}
\newcommand\barcaly{{\overline{\Theta}}}
\newcommand\calw{{\mathcal W}}
\newcommand\calv{{\mathscr V}}

\newcommand\vs{{\mathcal V}_{\rm s}}
\newcommand\vns{{\mathcal V}_{\rm ns}}

\newcommand\caldc{\cald_{\rm c}}
\newcommand\caldnc{\cald_{\rm nc}}
\newcommand\scrQ{\mathscr Q}

\newcommand\tcalsu{{\tcals}_{\rm u}}

\newcommand\grock{\tau}

\newcommand\id{{\rm id}}
\newcommand\smock{\zeta}
\newcommand\smurk{'}
\newcommand\glurk{''}

\newcommand\Isom{{\rm Isom}}
\newcommand\area{\mathop{{\rm area}}}

\newcommand\cut{^\dagger_}

\newcommand\ZZ{{\bf Z}}

\newcommand\NN{{\bf N}}
\newcommand\CC{{\bf C}}
\newcommand\QQ{{\bf Q }}
\newcommand\HH{{\bf H}}

\newcommand\calk{{\mathcal K}}
\newcommand\cals{{\mathscr S}}

\newcommand\calr{{\mathscr R}}
\newcommand\omecals{{(\Omega,\cals)}}

\newcommand\scrd{{\mathscr D}}

\newcommand\wt{\mathop{{\rm wt}}}

\newcommand\obd{\mathop{\omega}}

\newcommand\silv{\mathop{\rm silv}}

\newcommand\vol{\mathop{\rm vol}}
\newcommand\volG{\mathop{\rm vol_G}}
\newcommand\volorb{\mathop{\rm vol_{ASTA}}}
\newcommand\voct{V_{\rm oct}}
\newcommand\vtet{V_{\rm tet}}

\newcommand\Fix{\mathop{{\rm Fix}}}
\newcommand\card{\mathop{{\rm card}}}
\newcommand\rank{\mathop{{\rm rank}}}
\newcommand\image{\mathop{{\rm Im}}}

\newcommand\piggletwo{{\rm PGL}_2}

\newcommand\Fr{\mathop{\rm Fr}}

\newcommand\full{fully booked}
\newcommand\central{centrally booked}

\newcommand\clashspheres{spheres of incompatibility}
\newcommand\clashsphere{sphere of incompatibility}
\newcommand\centralclashsphere{clash component}
\newcommand\centralclashspheres{clash components}

\newcommand\strict{strictly booked}

\newcommand\doublefull{doubly fully booked}
\newcommand\doublecentral{doubly centrally booked}

\DeclareFontFamily{U}{rcjhbltx}{}
\DeclareFontShape{U}{rcjhbltx}{m}{n}{<->rcjhbltx}{}
\DeclareSymbolFont{hebrewletters}{U}{rcjhbltx}{m}{n}

\let\aleph\relax\let\beth\relax
\let\gimel\relax\let\daleth\relax

\DeclareMathSymbol{\aleph}{\mathord}{hebrewletters}{39}
\DeclareMathSymbol{\beth}{\mathord}{hebrewletters}{98}
\DeclareMathSymbol{\gimel}{\mathord}{hebrewletters}{103}
\DeclareMathSymbol{\daleth}{\mathord}{hebrewletters}{100}

\DeclareMathSymbol{\lamed}{\mathord}{hebrewletters}{108}
\DeclareMathSymbol{\mem}{\mathord}{hebrewletters}{109}
\DeclareMathSymbol{\ayin}{\mathord}{hebrewletters}{96}
\DeclareMathSymbol{\tsadi}{\mathord}{hebrewletters}{118}
\DeclareMathSymbol{\qof}{\mathord}{hebrewletters}{114}
\DeclareMathSymbol{\shin}{\mathord}{hebrewletters}{152}

%% file: orbifolds1.bbl
\begin{thebibliography}{10}

\bibitem{twocusps}
Ian Agol.
\newblock The minimal volume orientable hyperbolic 2-cusped 3-manifolds.
\newblock {\em Proc. Amer. Math. Soc.}, 138(10):3723--3732, 2010.

\bibitem{rankfour}
Ian Agol, Marc Culler, and Peter~B. Shalen.
\newblock Dehn surgery, homology and hyperbolic volume.
\newblock {\em Algebr. Geom. Topol.}, 6:2297--2312, 2006.

\bibitem{last}
Ian Agol, Marc Culler, and Peter~B. Shalen.
\newblock Singular surfaces, mod 2 homology, and hyperbolic volume. {I}.
\newblock {\em Trans. Amer. Math. Soc.}, 362(7):3463--3498, 2010.

\bibitem{ast}
Ian Agol, Peter~A. Storm, and William~P. Thurston.
\newblock Lower bounds on volumes of hyperbolic {H}aken 3-manifolds.
\newblock {\em J. Amer. Math. Soc.}, 20(4):1053--1077 (electronic), 2007.
\newblock With an appendix by Nathan Dunfield.

\bibitem{accs}
James~W. Anderson, Richard~D. Canary, Marc Culler, and Peter~B. Shalen.
\newblock Free {K}leinian groups and volumes of hyperbolic {$3$}-manifolds.
\newblock {\em J. Differential Geom.}, 43(4):738--782, 1996.

\bibitem{atkinson}
Christopher~K. Atkinson.
\newblock Two-sided combinatorial volume bounds for non-obtuse hyperbolic
  polyhedra.
\newblock {\em Geom. Dedicata}, 153:177--211, 2011.

\bibitem{bp}
Riccardo Benedetti and Carlo Petronio.
\newblock {\em Lectures on hyperbolic geometry}.
\newblock Universitext. Springer-Verlag, Berlin, 1992.

\bibitem{bbmbp}
Laurent Bessi{\`e}res, G{\'e}rard Besson, Sylvain Maillot, Michel Boileau, and
  Joan Porti.
\newblock {\em Geometrisation of 3-manifolds}, volume~13 of {\em EMS Tracts in
  Mathematics}.
\newblock European Mathematical Society (EMS), Z\"urich, 2010.

\bibitem{blp}
Michel Boileau, Bernhard Leeb, and Joan Porti.
\newblock Geometrization of 3-dimensional orbifolds.
\newblock {\em Ann. of Math. (2)}, 162(1):195--290, 2005.

\bibitem{bmp}
Michel Boileau, Sylvain Maillot, and Joan Porti.
\newblock {\em Three-dimensional orbifolds and their geometric structures},
  volume~15 of {\em Panoramas et Synth\`eses [Panoramas and Syntheses]}.
\newblock Soci\'et\'e Math\'ematique de France, Paris, 2003.

\bibitem{other-bs}
F.~Bonahon and L.~Siebenmann.
\newblock The classification of {S}eifert fibred {$3$}-orbifolds.
\newblock In {\em Low-dimensional topology ({C}helwood {G}ate, 1982)},
  volume~95 of {\em London Math. Soc. Lecture Note Ser.}, pages 19--85.
  Cambridge Univ. Press, Cambridge, 1985.

\bibitem{bonahon-siebenmann}
F.~Bonahon and L.~C. Siebenmann.
\newblock The characteristic toric splitting of irreducible compact
  {$3$}-orbifolds.
\newblock {\em Math. Ann.}, 278(1-4):441--479, 1987.

\bibitem{borel}
A.~Borel.
\newblock Commensurability classes and volumes of hyperbolic {$3$}-manifolds.
\newblock {\em Ann. Scuola Norm. Sup. Pisa Cl. Sci. (4)}, 8(1):1--33, 1981.

\bibitem{bcsz}
Steve Boyer, Marc Culler, Peter~B. Shalen, and Xingru Zhang.
\newblock Characteristic subsurfaces and {D}ehn filling.
\newblock {\em Trans. Amer. Math. Soc.}, 357(6):2389--2444, 2005.

\bibitem{bredon}
Glen~E. Bredon.
\newblock {\em Introduction to compact transformation groups}.
\newblock Academic Press, New York-London, 1972.
\newblock Pure and Applied Mathematics, Vol. 46.

\bibitem{canmac}
Richard~D. Canary and Darryl McCullough.
\newblock Homotopy equivalences of 3-manifolds and deformation theory of
  {K}leinian groups.
\newblock {\em Mem. Amer. Math. Soc.}, 172(812):xii+218, 2004.

\bibitem{cao-m}
Chun Cao and G.~Robert Meyerhoff.
\newblock The orientable cusped hyperbolic {$3$}-manifolds of minimum volume.
\newblock {\em Invent. Math.}, 146(3):451--478, 2001.

\bibitem{Cao-Zhu}
Huai-Dong Cao and Xi-Ping Zhu.
\newblock A complete proof of the {P}oincar\'e and geometrization
  conjectures---application of the {H}amilton-{P}erelman theory of the {R}icci
  flow.
\newblock {\em Asian J. Math.}, 10(2):165--492, 2006.

\bibitem{chk}
Daryl Cooper, Craig~D. Hodgson, and Steven~P. Kerckhoff.
\newblock {\em Three-dimensional orbifolds and cone-manifolds}, volume~5 of
  {\em MSJ Memoirs}.
\newblock Mathematical Society of Japan, Tokyo, 2000.
\newblock With a postface by Sadayoshi Kojima.

\bibitem{hodad}
Marc Culler, Jason Deblois, and Peter~B. Shalen.
\newblock Incompressible surfaces, hyperbolic volume, {H}eegaard genus and
  homology.
\newblock {\em Comm. Anal. Geom.}, 17(2):155--184, 2009.

\bibitem{hyperhaken}
Marc Culler and Peter~B. Shalen.
\newblock Volumes of hyperbolic {H}aken manifolds. {I}.
\newblock {\em Invent. Math.}, 118(2):285--329, 1994.

\bibitem{inject}
Marc Culler and Peter~B. Shalen.
\newblock Betti numbers and injectivity radii.
\newblock {\em Proc. Amer. Math. Soc.}, 137(11):3919--3922, 2009.

\bibitem{lastplusone}
Marc Culler and Peter~B. Shalen.
\newblock Singular surfaces, mod 2 homology, and hyperbolic volume, {II}.
\newblock {\em Topology Appl.}, 158(1):118--131, 2011.

\bibitem{fourfree}
Marc Culler and Peter~B. Shalen.
\newblock 4-free groups and hyperbolic geometry.
\newblock {\em J. Topol.}, 5(1):81--136, 2012.

\bibitem{after-hodad}
Jason DeBlois and Peter~B. Shalen.
\newblock Volume and topology of bounded and closed hyperbolic 3-manifolds.
\newblock {\em Comm. Anal. Geom.}, 17(5):797--849, 2009.

\bibitem{epstein}
D.~B.~A. Epstein.
\newblock Curves on {$2$}-manifolds and isotopies.
\newblock {\em Acta Math.}, 115:83--107, 1966.

\bibitem{otherFHS}
Michael Freedman, Joel Hass, and Peter Scott.
\newblock Closed geodesics on surfaces.
\newblock {\em Bull. London Math. Soc.}, 14(5):385--391, 1982.

\bibitem{FHS}
Michael Freedman, Joel Hass, and Peter Scott.
\newblock Least area incompressible surfaces in {$3$}-manifolds.
\newblock {\em Invent. Math.}, 71(3):609--642, 1983.

\bibitem{gromov}
Michael Gromov.
\newblock Volume and bounded cohomology.
\newblock {\em Inst. Hautes \'Etudes Sci. Publ. Math.}, 56:5--99 (1983), 1982.

\bibitem{hatcherbook}
Allen Hatcher.
\newblock {\em Algebraic topology}.
\newblock Cambridge University Press, Cambridge, 2002.

\bibitem{hempel}
John Hempel.
\newblock {\em 3-manifolds}.
\newblock AMS Chelsea Publishing, Providence, RI, 2004.
\newblock Reprint of the 1976 original.

\bibitem{illman}
S\"oren Illman.
\newblock Smooth equivariant triangulations of {$G$}-manifolds for {$G$} a
  finite group.
\newblock {\em Math. Ann.}, 233(3):199--220, 1978.

\bibitem{jaco-rubinstein}
William Jaco and J.~Hyam Rubinstein.
\newblock P{L} equivariant surgery and invariant decompositions of
  {$3$}-manifolds.
\newblock {\em Adv. in Math.}, 73(2):149--191, 1989.

\bibitem{js}
William~H. Jaco and Peter~B. Shalen.
\newblock Seifert fibered spaces in {$3$}-manifolds.
\newblock {\em Mem. Amer. Math. Soc.}, 21(220):viii+192, 1979.

\bibitem{Jo}
Klaus Johannson.
\newblock {\em Homotopy equivalences of {$3$}-manifolds with boundaries},
  volume 761 of {\em Lecture Notes in Mathematics}.
\newblock Springer, Berlin, 1979.

\bibitem{jones-reid}
Kerry~N. Jones and Alan~W. Reid.
\newblock Minimal index torsion-free subgroups of {K}leinian groups.
\newblock {\em Math. Ann.}, 310(2):235--250, 1998.

\bibitem{kapovich}
Michael Kapovich.
\newblock {\em Hyperbolic manifolds and discrete groups}.
\newblock Modern Birkh\"auser Classics. Birkh\"auser Boston, Inc., Boston, MA,
  2009.
\newblock Reprint of the 2001 edition.

\bibitem{kawasaki}
Tetsuro Kawasaki.
\newblock The signature theorem for {$V$}-manifolds.
\newblock {\em Topology}, 17(1):75--83, 1978.

\bibitem{kleiner-lott}
Bruce Kleiner and John Lott.
\newblock Notes on {P}erelman's papers.
\newblock {\em Geom. Topol.}, 12(5):2587--2855, 2008.

\bibitem{lange}
Christian Lange.
\newblock Equivariant smoothing of piecewise linear manifolds.
\newblock 1. arXiv:1507.02395.

\bibitem{meeks-scott}
William~H. Meeks, III and Peter Scott.
\newblock Finite group actions on {$3$}-manifolds.
\newblock {\em Invent. Math.}, 86(2):287--346, 1986.

\bibitem{milnorprime}
J.~Milnor.
\newblock A unique decomposition theorem for {$3$}-manifolds.
\newblock {\em Amer. J. Math.}, 84:1--7, 1962.

\bibitem{miyamoto}
Yosuke Miyamoto.
\newblock Volumes of hyperbolic manifolds with geodesic boundary.
\newblock {\em Topology}, 33(4):613--629, 1994.

\bibitem{morgan-poincare}
John Morgan and Gang Tian.
\newblock {\em Ricci flow and the {P}oincar\'e conjecture}, volume~3 of {\em
  Clay Mathematics Monographs}.
\newblock American Mathematical Society, Providence, RI; Clay Mathematics
  Institute, Cambridge, MA, 2007.

\bibitem{Morgan-Tian}
John Morgan and Gang Tian.
\newblock {\em The geometrization conjecture}, volume~5 of {\em Clay
  Mathematics Monographs}.
\newblock American Mathematical Society, Providence, RI; Clay Mathematics
  Institute, Cambridge, MA, 2014.

\bibitem{smithbook}
John~W. Morgan and Hyman Bass, editors.
\newblock {\em The {S}mith conjecture}, volume 112 of {\em Pure and Applied
  Mathematics}.
\newblock Academic Press, Inc., Orlando, FL, 1984.
\newblock Papers presented at the symposium held at Columbia University, New
  York, 1979.

\bibitem{mukherjee}
Amiya Mukherjee.
\newblock {\em Differential topology}.
\newblock Hindustan Book Agency, New Delhi; Birkh\"{a}user/Springer, Cham,
  second edition, 2015.

\bibitem{rafalski-turnover}
Shawn Rafalski.
\newblock Immersed turnovers in hyperbolic 3-orbifolds.
\newblock {\em Groups Geom. Dyn.}, 4(2):333--376, 2010.

\bibitem{schoen-yau}
R.~Schoen and Shing~Tung Yau.
\newblock Existence of incompressible minimal surfaces and the topology of
  three-dimensional manifolds with nonnegative scalar curvature.
\newblock {\em Ann. of Math. (2)}, 110(1):127--142, 1979.

\bibitem{prelim}
Peter~B. Shalen.
\newblock Volume and homology for hyperbolic 3-orbifolds, {I}.
\newblock arXiv:1709.07413.

\bibitem{arithmetic}
Peter~B. Shalen.
\newblock Homology and volume of hyperbolic {$3$}-orbifolds, and enumeration of
  arithmetic groups.
\newblock {\em Ann. Fac. Sci. Toulouse Math. (6)}, 24(5):1147--1156, 2015.

\bibitem{smale}
Stephen Smale.
\newblock A {V}ietoris mapping theorem for homotopy.
\newblock {\em Proc. Amer. Math. Soc.}, 8:604--610, 1957.

\bibitem{soma}
Teruhiko Soma.
\newblock The {G}romov invariant of links.
\newblock {\em Invent. Math.}, 64(3):445--454, 1981.

\bibitem{thurstonnotes}
William~P. Thurston.
\newblock The geometry and topology of $3$-manifolds.
\newblock Lecture notes, 1978.

\bibitem{Waldhausen}
Friedhelm Waldhausen.
\newblock On irreducible {$3$}-manifolds which are sufficiently large.
\newblock {\em Ann. of Math. (2)}, 87:56--88, 1968.

\end{thebibliography}
